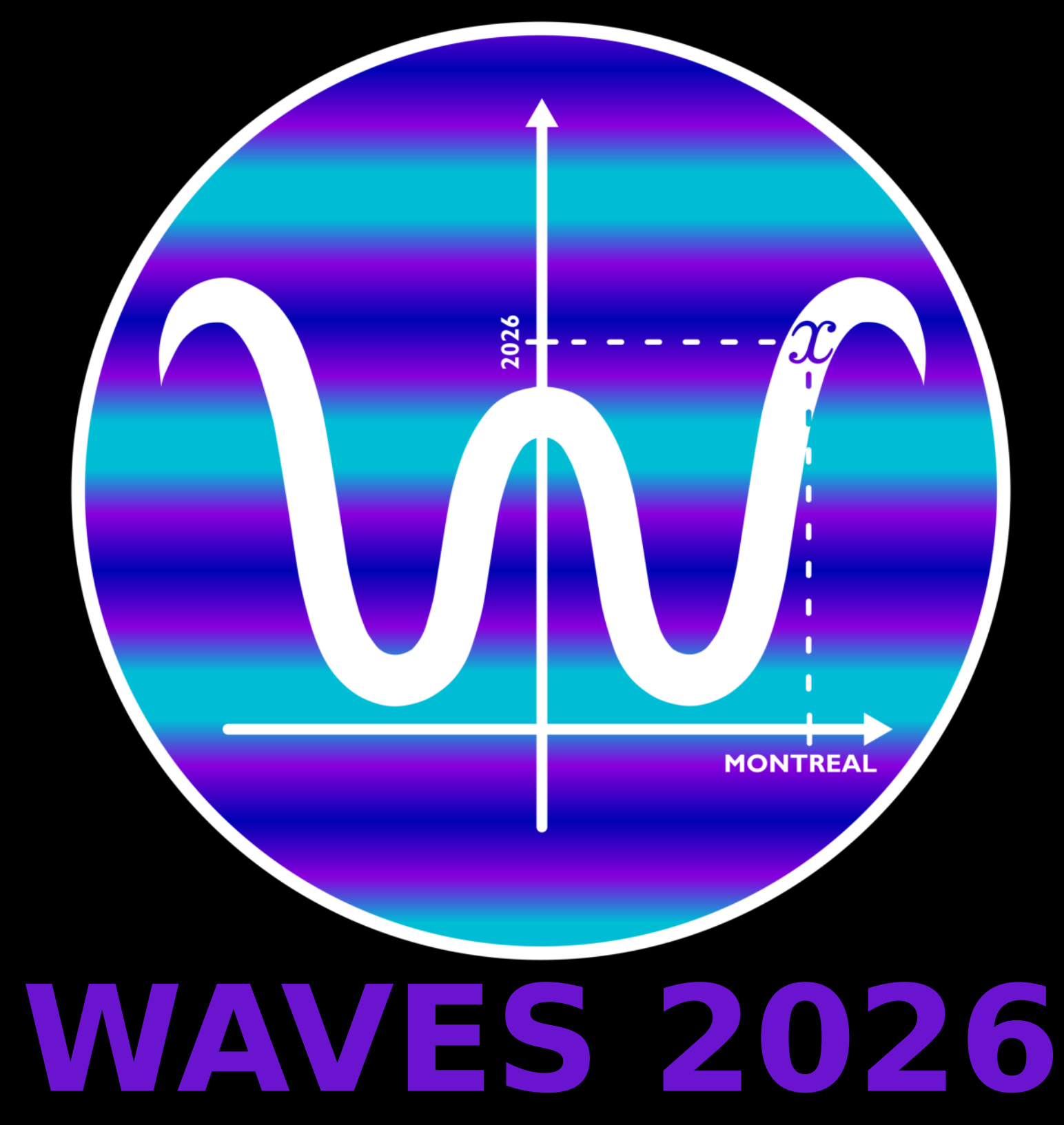

# WAVES 2026

The 17th International Conference on

*Mathematical and Numerical Aspects of Wave Propagation*

# BOOK OF ABSTRACTS

**John Molson School of Business**

Concordia University

Montreal, Canada

**June 22 – 26, 2026**

## Organizing committee

**Lise-Marie Imbert-Gerard, Chrysoula Tsogka, Daniel Appelö, Thomas Hagstrom, Yann-Meing Law, Serge Prudhomme, Behrooz Yousefzadeh**

## Scientific committee

**Daniel Appelo, USA, Tilo Arens, Germany, Lehel Banjai, UK, Alex Barnett, USA, Hélène Barucq, France, Eliane Bécache, France, Anne-Sophie Bonnet-Ben Dhia, France, Marc Bonnet, France, Fioralba Cakoni, USA, Maxence Cassier, France, Stephanie Chaillat, France, Simon Chandler-Wilde, UK, Lucas Chesnel, France, Xavier Claeys, France, Bérengère Delourme, France, Bruno Després, France, Florian Faucher, France, Sonia Fliss, France, Damien Fournier, Germany, Adrianna Gillman, USA, Laure Giovangigli, France, Dan Givoli, Israel, Laurent Gizon, Germany, Roland Griesmaier, Germany, Marcus Grote, Switzerland, Bojan Guzina, USA, Houssem Haddar, France, Thomas Hagstrom, USA,eMartin Halla, Germany, Christophe Hazard, France, Dave Hewett, UK, Ralf Hiptmair, Switzerland, Marlis Hochbruck, Germany, Thorsten Hohage, Germany, Lise-Marie Imbert-Gerard, USA, Sébastien Imperiale, France, Maryna Kachanovska, France, Manfred Kaltenbacher, Austria, Olivier Lafitte, France, Yann-Meing Law, Canada, Christoph Lehrenfeld, Germany, Frédérique Le Louër, France, Bruno Lombard, France, Wangtao Lu, China, Paul Martin, USA, Axel Modave, France, Andrea Moiola, Italy, Peter Monk, USA, Lothar Nannen, Austria, Tram Nguyen, Germany, Vanja Nikolic, The Netherlands, Nilima Nigam, Canada, Lauri Oksanen, Finland, Konstantin Pankrashkin, Germany, Ilaria Perugia, Austria, Ha Pham, France, Jerónimo Rodríguez García, Spain, Olof Runborg, Sweden, Tonatiuh Sánchez-Vizuet, USA, Claire Scheid, France, Euan Spence, UK, Antoine Tonnoir, France, Sebastien Tordeux, France, Chrysoula Tsogka, USA, Jorn Zimmerling, Sweden**

# Wavebeam propagation in heterogeneous media

**Guillaume Bal**[1,*], **Anjali Nair**[1]
[1]CCAM, University of Chicago, Chicago IL, USA
*Email: guillaumebal@uchicago.edu

## Abstract

Wavebeam models concern the propagation of time-harmonic waves in heterogeneous media along a privileged direction of propagation, whereby neglecting back-scattering along that axis while still retaining transverse dispersion and refraction effects. Such effects include the formation of speckle, a phenomenon of coherent interferences leading to rapid and random fluctuations in the wavebeam intensity.

We consider asymptotic models of such speckle formation, and in particular present recent results on the resolution of the so-called Gaussian conjecture, stating that wavefields are approximately complex-circular Gaussian, for paraxial and Itô-Schrödinger models in the weak-coupling regime. We also justify the use of phase screen method, a form of splitting algorithm, to simulate such wavebeam propagation.

## 1 Introduction

This abstract presents recent results in the theoretical and computational analyses of laser light propagation through turbulent atmospheres. Such analyses find applications in open space communication; see right panel of Fig. 1.

A general feature of laser light propagation in heterogeneous media is the apparition of speckle, a granular structure where light intensity oscillates between large and small values in a random pattern; see Fig. 1 and Fig. 2. Speckle is a well known phenomenon of coherent interference [18]. Heuristically, it may be interpreted as a random convergence and divergence of caustics; see Fig. 4. As the strength of the random turbulent atmosphere increases, one expects speckle patterns to display smaller and smaller scales; see Fig. 2. Yet quantitative descriptions of the phenomenon remain difficult to derive. Until recently, no mathematical model fully explained how a deterministic incident beam develops such a speckle pattern as the distance of propagation through the turbulent medium increases.

This problem was referred to as the Gaussian conjecture for wavebeam propagation [16]; see also [12,13]. This conjecture, implying that a complex-valued scalar wavebeam model would be asymptotically given by a complex circular Gaussian law, was recently solved in the weak-coupling regime and in the Itô-Schrödinger model of wave propagation in [5], following earlier results in [16] showing that fourth-order statistical moments of the field were compatible with the Gaussian conjecture. Note that other regimes of propagation displaying non-Gaussian statistics are also possible [10,13]. These results were extended in the weak-coupling regime for a more general paraxial model in [7]. Generalizations of these techniques were also applied to the propagation of partially coherent sources in [6]. In this setting, incident beams are themselves modeled statistically and the latter paper shows that the speckle field is in fact no longer necessarily Gaussian. Limiting distributions are then characterized by their statistical moments.

Beyond confirming that the wave-field becomes complex Gaussian after long-distance propagation, the papers [5–7] identify a regime of propagation distances where the correlation function, which fully characterizes the wavefield asymptotically, solves an anomalous diffusion equation. In that regime, the scintillation index, a gauge of the randomness of the wavefield intensity defined more precisely as the intensity variance rescaled by the square of the average intensity, is asymptotically unity. This implies that at any point, the wavefield intensity satisfies an exponential distribution, which is consistent with speckle observations in, e.g., Fig.2.

This result is not good news for communication purposes. Indeed, strong randomness in intensity measurements could potentially significantly degrade communication capabilities. At the same time, the aforementioned diffusion model also predicts that speckle dots become smaller as turbulence strength or propagation distance increase. This observation, which is consistent with the simulations in Fig.2, shows that local intensity averages might still significantly mitigate the effect of speckle fluctuations while preserving coherent features of the propagating signal. Recent results on wavefield classifications in [24] quantify, and partially confirm, this observation.

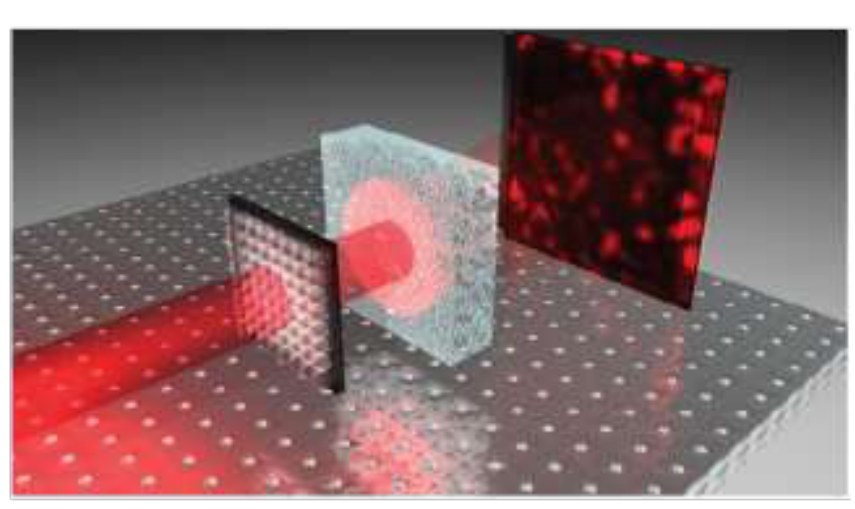

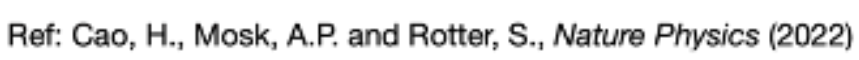

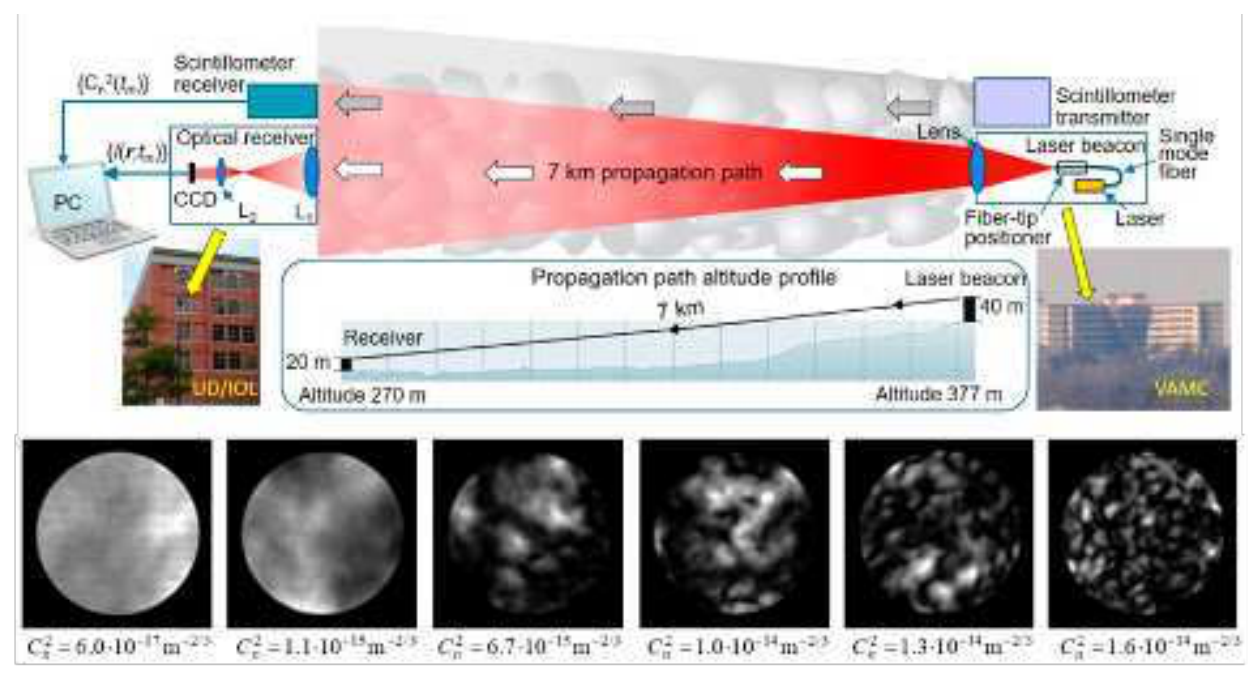

Figure 1: Examples of speckle generation for laser light propagating through heterogeneous media. Left: Light propagation through slab of random material. Right: application in open space communication where light speckle appears as turbulence strength increases (from left to right).

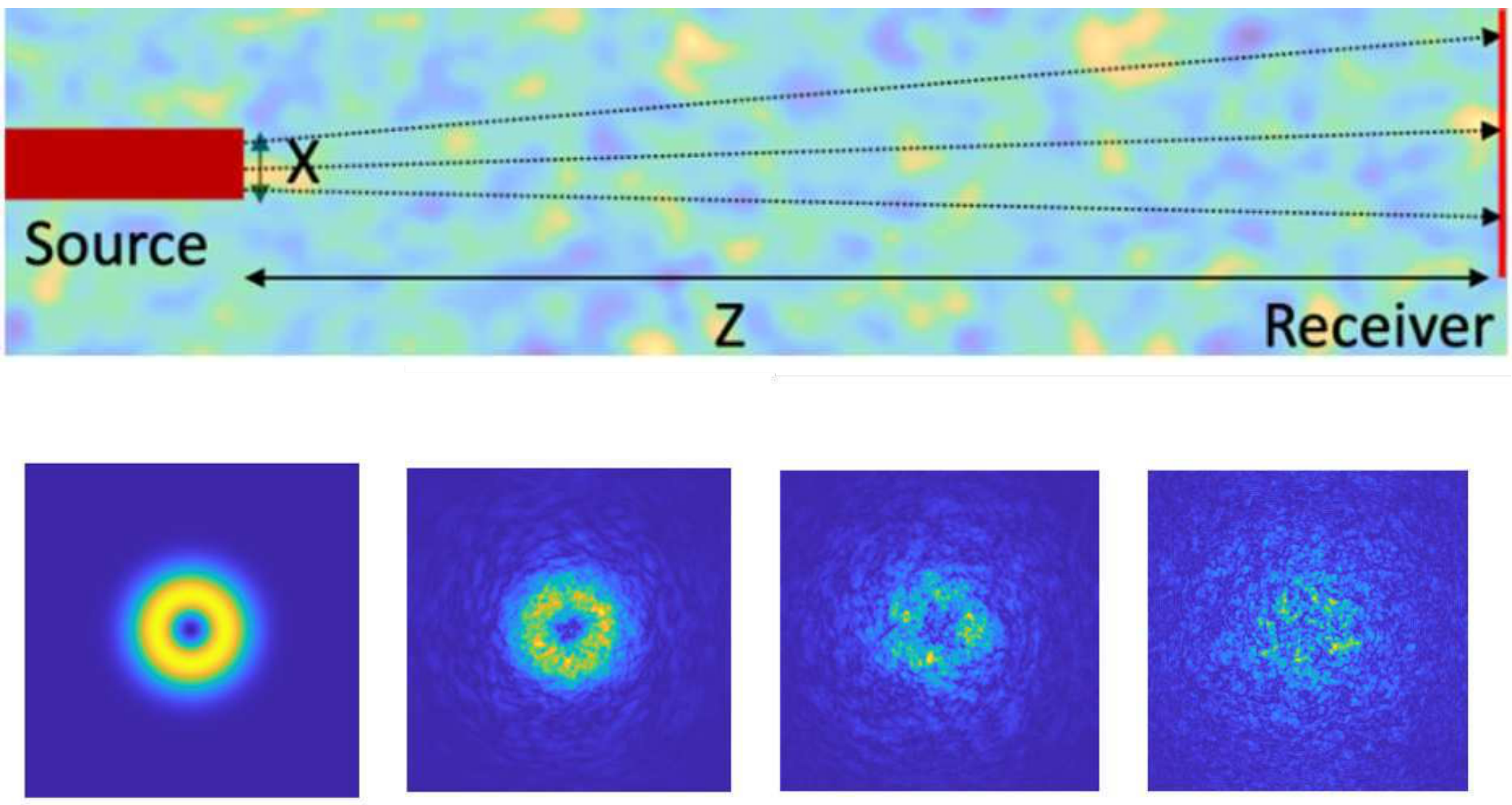

Figure 2: Numerical simulation by phase screen method (a splitting method) of laser light propagating through media with heterogeneous refractive indices. Speckle formation increases as turbulence strength increases (bottom from left to right: snapshots at a given axial distance $Z$).

Due to their practical relevance, wavebeam propagation in turbulent media is routinely solved numerically. Such simulations are challenging because the typical wavelength, on the order of micrometers, and the medium correlation length, on the order of millimeters, are significantly smaller than the typical overall distance of propagation, of order of kilometers in the open space communication application. This implies that the spatial discretization of (time-harmonic) wavebeams must be much larger than the aforementioned small parameters of the model.

One of the few numerical methods of choice is the phase screen method. It is a splitting method, alternatively integrating a dispersive term, which is easily simulated in Fourier variables, and a local refractive term modeling interactions with the underlying medium, which is easily simulated in physical variables. Surprisingly, this method seems to provide accurate results even when the splitting step is chosen much larger than the wavelength $\lambda$ of the propagating beam and the correlation length $\ell$ of the medium. In other regimes of wave propagation, such expectations would simply be incorrect [8]. In wavebeam propagation, however, the paraxial model and the Itô-Schrödinger model, its stochastic limit, have close-by solutions statistically. This heuristically justifies why the phase screen method stands a chance at providing quantitative estimates. In the preprint [3], we prove convergence results as the splitting step $\Delta z$ converges to 0 even when it is significantly larger than $\ell$, which is the case of practical interest. More precisely, we show that a centered (Strang) splitting is limited to first-order accuracy when convergence is measured in a root mean square sense. However, it is indeed second-order for statistical moments of the wavefield. Full spatial discretizations of these equations are also considered and numerical simulations [3] confirm the theoretical predictions.

## 2 Paraxial and Itô-Schrödinger models

We start with a Helmholtz scalar model for time-harmonic wave fields

$$\Big(\partial_z^2 + \Delta_x + k_0^2 n^2(z,x)\Big)p(z,x) = 0$$

posed for $z > 0$ coordinate on the main axis of propagation and for transverse variables $x \in \mathbb{R}^d$, and augmented with $p(0,x) = u_0(x)$ at $z = 0$. Here, $k_0 = \omega_0/c_0$ with $\omega_0$ frequency of the laser beam and $c_0$ light speed, while $n^2(z,x) = \frac{c^2(z,x)}{c_0^2} = 1 + \nu(z,x)$ is the index of refraction modeling interactions with a turbulent atmosphere. We will assume $\nu$ to be Gaussian in what follows. Consider the high frequency long propagation scaling:

$$z \to \frac{\eta z}{\varepsilon\theta},\ x \to x,\ k_0 \to \frac{k_0}{\theta},\ \nu \to \frac{\varepsilon^{\frac{1}{2}}\theta^{\frac{3}{2}}}{\eta^{\frac{3}{2}}}\nu,$$

with $0 < \theta, \varepsilon \ll 1$ and $\eta = O(1)$ in the *kinetic* regime and $\eta = o(1)$ in the *diffusive* regime. The medium correlation length $\ell$ (mm in practice) is scaled to be of order $O(1)$. Thus, $\theta = \lambda/\ell$ is of order $10^{-3}$. The parameter $\varepsilon$ plays two roles: it models the weak coupling regime, stating that random fluctuations are of order $O(\varepsilon)$ but integrated over distances at least $O(\varepsilon^{-1})$; it also models the width of the incident beam with $\varepsilon = \ell/w_0$ so that the incident beam is of the form $u_0(\varepsilon x)$.

We now introduce the standard paraxial envelope ansatz $u^\theta(z,x) = p\big(\frac{\eta z}{\varepsilon\theta}, x\big)e^{-\frac{i\eta k_0 z}{\varepsilon\theta^2}}$. After formally ignoring the backscattering term $(\frac{\varepsilon\theta}{\eta})^2\partial_z^2 u^\theta$, we arrive at the paraxial approximation

$$2ik_0\partial_z u^\theta + \frac{\eta}{\varepsilon}\Delta_x u^\theta + \frac{k_0^2}{(\eta\varepsilon\theta)^{\frac{1}{2}}}\nu\Big(\frac{\eta z}{\varepsilon\theta}, x\Big)u^\theta = 0,$$

with incident beam profile $u^\theta(0,x) = u_0^\theta(x)$. Justifying the paraxial equation from the Helmholtz equation is difficult; see [1, 15]. However, it is valid only when $\theta = \lambda/\ell \ll 1$. Otherwise, large-angle scattering leading for instance to radiative transfer models is unavoidable [2]. And in the regime $\theta \ll 1$, a central limit approximation shows that

$$\frac{1}{\tau^{\frac{1}{2}}}\nu(\frac{z}{\tau}, x)dz \approx B(dz, x)$$

with $B$ a Gaussian field. Therefore, as soon as a paraxial wavebeam model is valid, then heuristically so is its stochastic Stratonovich-Schrödinger limit

$$du^\varepsilon = \frac{i\eta}{2k_0\varepsilon}\Delta_x u^\varepsilon dz + \frac{ik_0}{2\eta}u^\varepsilon \circ dB$$

where the singular product $u^\varepsilon \circ dB$ needs to be understood in a Stratonovich sense.

Define the statistics

$$C(z,x) = \mathbb{E}\nu(z,x)\nu(0,0),\ R(x) = \int_{\mathbb{R}} C(z,x)dz$$

so that $\mathbb{E}B(z,x)B(z',x') = \min(z,z')R(x-x')$. After standard correction, we thus obtain the Itô-Schrödinger (SPDE) model

$$du^\varepsilon = \frac{i\eta}{2k_0\varepsilon}\Delta_x u^\varepsilon \mathrm{d}z - \frac{k_0^2 R(0)}{8\eta^2}u^\varepsilon \mathrm{d}z + \frac{ik_0}{2\eta}u^\varepsilon dB,$$

with incident conditions $u^\varepsilon(0,x) = u_0(\varepsilon x)$ and where the singular product $u^\varepsilon dB$ is now understood in the Itô sense. The above SPDE is analyzed mathematically in [10] and its derivation in, e.g., [11].

The convergence results for $u^\varepsilon$ obtained in [5] were generalized for $u^\theta$ in [7] essentially by showing that statistical moments of $u^\theta$ were appropriately close to the statistical moments of $u^\varepsilon$. We now focus on the latter. The Itô formula for Hilbert-valued processes [16] shows that the statistical moments for $p, q \geq 0$

$$\mu_{p,q}^\varepsilon(z,X,Y) = \mathbb{E}\prod_{j=1}^{p} u^\varepsilon(z,x_j)\prod_{j=1}^{q} u^{\varepsilon *}(z,y_j)$$

satisfy equations $\partial_z \mu^\varepsilon_{p,q} = \mathcal{L}^\varepsilon_{p,q}\mu^\varepsilon_{p,q}$ with

$$\mathcal{L}^\varepsilon_{p,q} = \frac{i\eta}{2k_0\varepsilon}(\sum_{j=1}^{p}\Delta_{x_j} - \sum_{j=1}^{q}\Delta_{y_j}) + \frac{k_0^2}{4\eta^2}\mathcal{U}_{p,q}(X,Y),$$

for $\mathcal{U}$ a potential involving sums of terms of the form $R(a-b)$ for $a,b$ variables in $(X,Y)$.

## 3 Convergence Results

We now describe the convergence process to a limiting Gaussian field. The convergence result is multiscale in nature; see Fig. 3. Let $(z,r)$ be *macroscopic* variables and $x$ (transverse) *microscopic* variables and introduce the field:

$$x \mapsto \phi^\varepsilon(z,r,x) := u^\varepsilon\Big(z, \frac{r}{\varepsilon} + \eta x\Big).$$

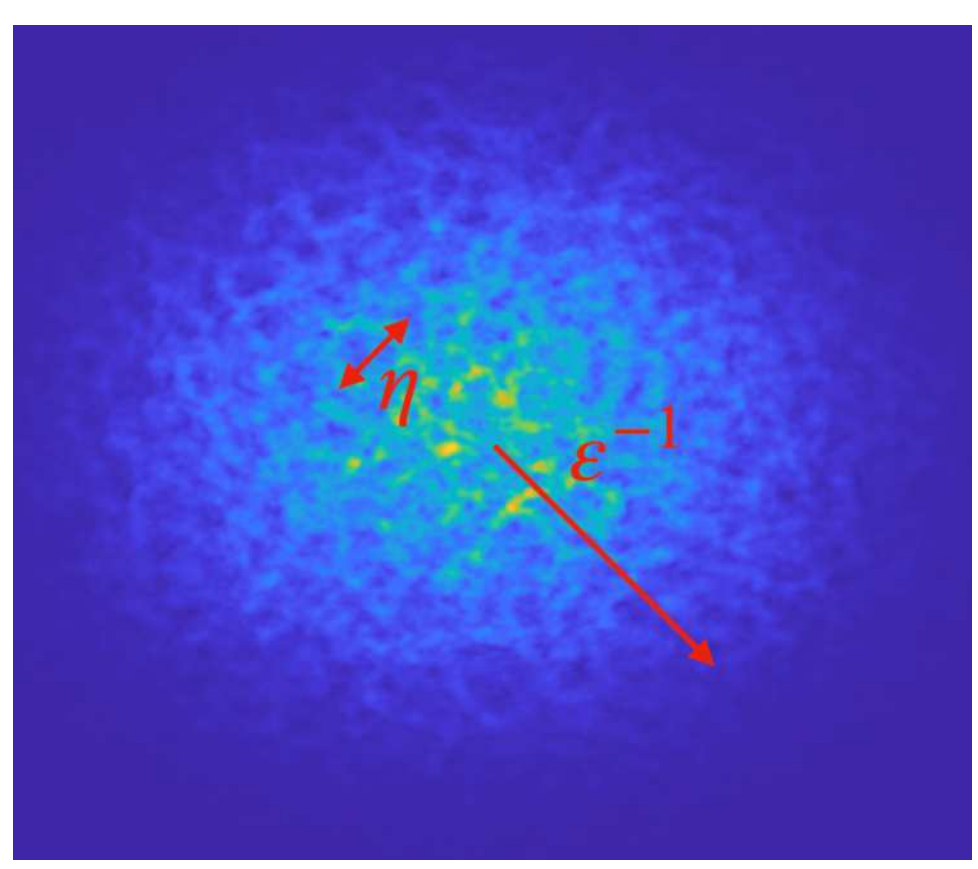

Figure 3: Multiscale convergence result with an incident beam scaling at the scale $\varepsilon^{-1} \gg 1$ and speckle grains of size $\eta \ll 1$.

Let $\eta = 1$ in the *kinetic* regime and $\eta \approx (\ln\ln\varepsilon^{-1})^{-1}$ in the *diffusive* regime. We need this separation of scales for technical reasons. Then we have the following result:

**Theorem** [Gaussian conjecture diffusive regime] At *fixed* $z > 0$ and *fixed* $r \in \mathbb{R}^d$, $x \mapsto \phi^\varepsilon$ converges in law in the space of continuous functions as $\varepsilon \to 0$ to a mean-zero complex circular Gaussian field. The limiting Gaussian process $\phi(x)$ is characterized by $\mathbb{E}\{\phi(x)\phi(y)\} = 0$ and by $\mathcal{C}(x,y) = \mathbb{E}\{\phi(x)\phi^*(y)\}$ given by

$$\begin{aligned}\mathcal{C}(x,y) \;\; = \;\; & e^{\frac{k_0^2}{32}z(y-x)^t\Gamma(y-x)}e^{-i\frac{3k_0}{2z}(y-x)\cdot r} \\ & G(z^3,r) * \Big[e^{i\frac{3k_0}{2z}(y-x)\cdot r}|u_0|^2(r)\Big]\end{aligned}$$

where $G$ solves the diffusion equation

$$(\partial_t + \frac{1}{24}\nabla_r\cdot\Gamma\nabla_r)G(t,r) = 0,$$

with $G(0,r) = \delta_0(r)$, and where the negative definite $\Gamma$ is given by

$$\Gamma = \nabla^2 R(0).$$

A similar result holds for the paraxial solution $u^\theta$ [7] and holds for both $u^\varepsilon$ and $u^\theta$ with appropriate modifications in the kinetic regime where $\eta = 1$. The latter is characterized by a coherent part $\mathbb{E}u^\varepsilon(z,x)$ decaying exponentially with $z$ with non mean-zero limit as $\varepsilon \to 0$.

The tensor $\Gamma$ is the only term that depends on the random medium. In the scalar case, it is proportional to $-\int |k|^2\hat{R}(k)dk < 0$ with $\hat{R}(k)$ the Fourier transform of transversal correlation function $R(x)$.

The above diffusion is *anomalous* in the sense that $z^3 \approx |x|^2$. This scaling may be explained as follows: in the absence of randomness, ballistic transport is modeled by $z^2 \approx |x|^2$ while transverse scattering is responsible for an additional standard factor $z$ in the mean distance $|x|^2$.

For incident plane waves [5], the correlation function $e^{\frac{k_0^2}{8}z(y-x)^t\Gamma(y-x)}$ describes speckle size of order $\mathfrak{s}^2 \approx \frac{4}{k_0^2|\Gamma|z}$ and hence inversely proportional to $\sigma\sqrt{z}$ where $\sigma$ is random fluctuations strength and $z$ is distance of propagation.

In this regime, the scintillation index $S = \frac{\mathbb{E}I^2-(\mathbb{E}I)^2}{(\mathbb{E}I)^2}$, with $I(x) = |u|^2(x)$, is approximately 1. In fact, in the limit $\varepsilon \to 0$, $|u|^2$ satisfies an exponential law. However, any average of the intensity over a domain that is large compared to $\mathfrak{s}\eta$ results in a statistically stable quantity.

**Elements of derivation.** The convergence result is based on an analysis of finite dimensional distributions of $\phi^\varepsilon$. This is carried out by analyzing (rescaled versions of) the moments $\mu^\varepsilon_{p,q}$. The equations they satisfy cannot be solved explicitly. However, it turns out that their solution operator $U^\varepsilon_{p,q}$ satisfies $U^\varepsilon_{p,q} = N^\varepsilon_{p,q} + E^\varepsilon_{p,q}$, where $N^\varepsilon_{p,q}$ is an operator that propagates moments as if the underlying process was Gaussian. A combinatorial analysis of the remaining terms shows that $\|E^\varepsilon_{p,q}\| \to 0$ as $\varepsilon \to 0$. Here, the norm $\|\cdot\|$ is that of bounded measures $\mathcal{M}_B(\mathbb{R}^{d(p+q)})$ after passing to the Fourier domain. This implies a uniform convergence in the physical variables and explains how we can obtain a convergence in law in the space of continuous functions.

This convergence of finite dimensional distributions needs to be augmented with a tightness

(compactness) argument. Using a similar machinery to the finite dimensional convergence result, we obtain the following (sufficient) stochastic continuity result

$$\mathbb{E}|\phi^\varepsilon(s,r,x+h)-\phi^\varepsilon(s,r,x)|^{2n} \le C(z,\alpha_0)|h|^{2\alpha_0 n},$$

uniformly in $s \in [0,z]$ and $h \in B(0,1) \subset \mathbb{R}^d$ for some $\alpha_0 > 0$.

These results then suffice to show that $\phi^\varepsilon$ converges to a Gaussian process, which is mean-zero in the diffusive regime. It remains to characterize the limit by computing its correlation function $\mathcal{C}(x,y)$. This is done explicitly by a direct analysis of the moment $\mu_{1,1}$.

## 4 Numerical simulations

Splitting techniques called phase screen methods are routinely used in the numerical simulations of paraxial wave propagation in random media [21–23]. They allow one to solve the paraxial model, with $\varepsilon = \eta = 1$ to simplify but still $\theta \ll 1$, and $\kappa_{1,2} > 0$,

$$\partial_z u^\theta = i\kappa_1 \Delta_x u^\theta + i\kappa_2 \frac{1}{\sqrt{\theta}} \nu\Big(\frac{z}{\theta}, x\Big) u^\theta,$$

with $u^\theta(0,x) = u_0(x)$. For $\gamma \in [0,1)$ with $\gamma = \frac{1}{2}$ for centered (Strang) splitting, define

$$\tau_\gamma(z) := \Delta z \sum_{n\ge 0} \delta(z-(n+\gamma)\Delta z).$$

The splitting scheme may then be defined as the solution to

$$\partial_z u^{\theta\Delta} = i\tau_\gamma(z)\kappa_1 \Delta_x u^{\theta\Delta} + i\kappa_2 \nu^\theta(z,x) u^{\theta\Delta}$$

with $u^{\theta\Delta}(0,x) = u_0(x)$. This is a convenient way to write the splitting scheme, with an equation on $(z_n, z_n+\Delta z)$ with $z_n = (n+\gamma)\Delta z$ solved with $\kappa_1 \equiv 0$ while the equation on $(z_n-0, z_n+0)$ is solved by applying $\mathcal{F}^{-1} e^{-i\kappa_1|\xi|^2\Delta z}\mathcal{F}$, where $\mathcal{F} = \mathcal{F}_{x\to\xi}$ is Fourier transform.

In practice, we need to choose $\Delta z \gg \theta$ since $\Delta z \ll \theta$ for $Z/\theta \approx 10^6$ is prohibitively expensive. Splitting schemes when $\Delta z$ is not the smallest scale have no reason to converge to the right solution a priori [8, 9]. However, since the paraxial model is well approximated by an Itô-Schrödinger limit for $u$, it may be expected that a splitting scheme for the SPDE limit indeed converges [19, 20]. Such convergence results for both the SPDE and the paraxial models have been obtained recently in [3]. The main strong convergence result [3, Theorem 2.1] is that for $\mathsf{v} \in \{u, u^\theta\}$,

$$|||\mathsf{v} - \mathsf{v}^\Delta||| \le C\Delta z$$

where $|||u||| := \sup_{0\le s\le Z}(\mathbb{E}\|u(s,\cdot)\|^2_{L^2(\mathbb{R}^d)})^{\frac{1}{2}}$. We thus obtain first-order convergence even for the centered splitting algorithm $\gamma = \frac{1}{2}$ in the root-mean-square sense. However, statistical moments enjoy better stability results [3, Theorem 2.2] as

$$\|\mu_{p,q}[\mathsf{v}] - \mu_{p,q}[\mathsf{v}^\Delta]\|_{L^\infty(\mathbb{R}^{(p+q)d})} \le C(\Delta z)^2$$

for $\mathsf{v} \in \{u, u^\theta\}$ when $\theta \le \Delta z$. See also [3] for convergence results of fully discretized equations using a spectral method for the transverse variables. Numerical simulations in [3] fully confirm these orders of convergence.

The convergence results to a Gaussian field shown in the preceding section involve a correlation function $\mathcal{C}(x,y)$ at a fixed value of $z$ with rapid Gaussian decay in $(x-y)$. One may more generally ask about correlations in $z$, as well as correlations across frequencies $\omega_0$ [14, 17]. This problem is analyzed in [4], where it is shown that fields at close-by values of $z$ are indeed jointly Gaussian with a correlation function that now decays much more slowly as $z^{-\frac{d}{2}}$. This is consistent with the numerical simulations presented in Fig.4 in transverse dimension $d=1$, where we see much longer correlation in the axial direction than in the transverse direction(s).

**Acknowledgment.** The first author acknowledges funding from ONR Grant N00014-26-1-2017 and NSF Grant DMS-2306411.

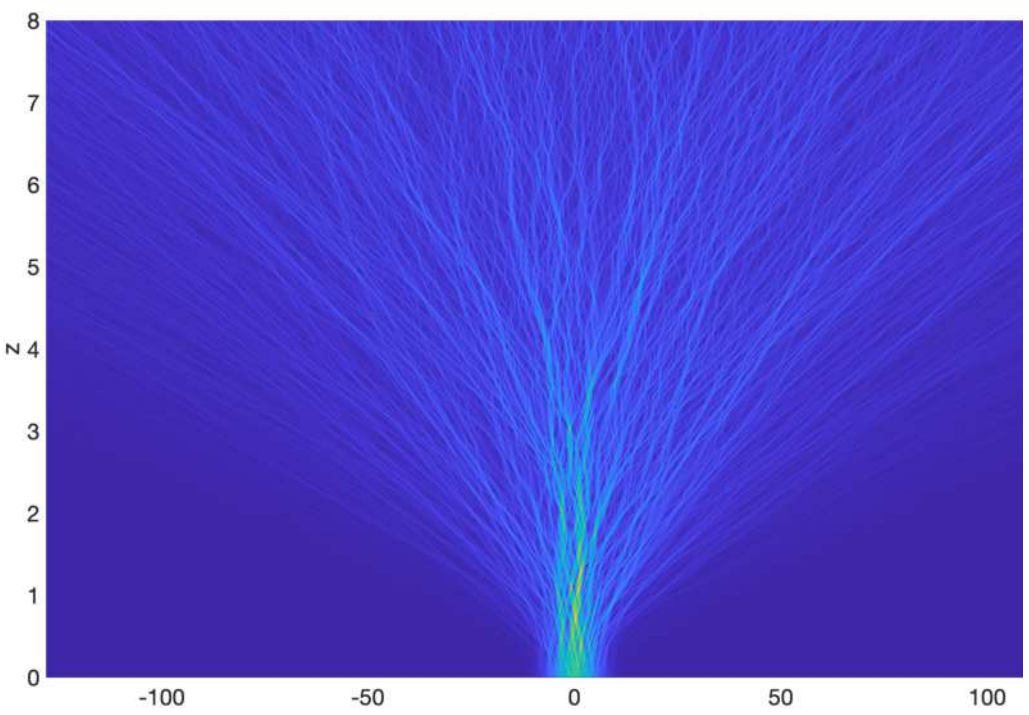

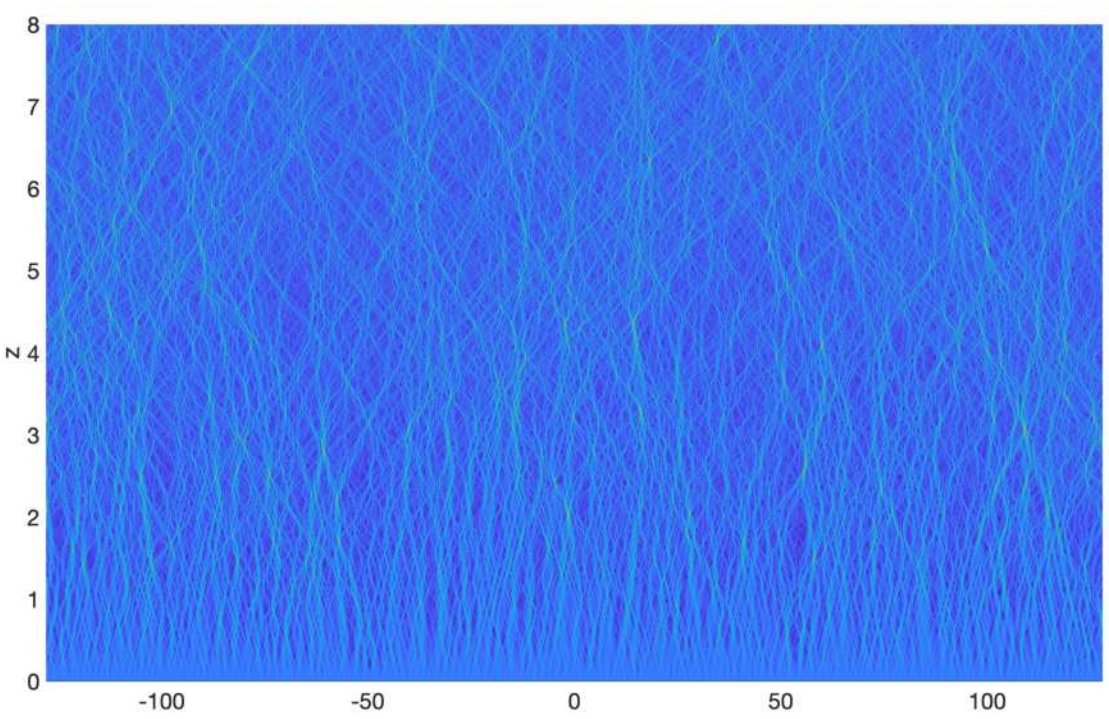

Figure 4: Long-distance propagation with $d = 1$ showing random convergence and divergence of caustics (branch flows) along axial direction. Left: localized incident beam; Right: plane-wave incident beam.

# Modal analysis of open acoustic waveguides

**Christophe Hazard**[1,*]
[1]POEMS, CNRS, Inria, ENSTA, Institut Polytechnique de Paris, 91120 Palaiseau, France
*Email: christophe.hazard@ensta.fr

## Abstract

This presentation overviews the results of a long-standing collaboration with my colleagues and friends Anne-Sophie Bonnet-Ben Dhia and Sonia Fliss, also including Sarah Al Humaikani, Lahcène Chorfi, Ghania Dakhia, Benjamin Goursaud and Antoine Tonnoir. It concerns the mathematical analysis of perturbations or junctions of acoustic open waveguides, more specifically two issues related to time-harmonic propagation in such devices: (i) the statement of an appropriate *radiation condition* and the well-posedness of the associated radiation problem; (ii) the possible existence of *trapped* modes. The results we have obtained over the last two decades are all based on a *modal approach.* Roughly speaking, this consists in using separation of spatial coordinates which allows us to represent the acoustic field as a superposition of evanescent and propagating modes.

## 1 Introduction

A *uniform waveguide* is a device made of a propagation medium whose physical features are invariant in one space direction, called here *longitudinal*, so that waves can propagate in this direction while remaining confined to a bounded region in the orthogonal *transverse* direction (for clarity, we restrict ourselves to the 2D case). Such a waveguide is called *closed* if its transverse cross-section is bounded; the confinement is then simply due to this boundedness. Here, we are rather interested in *open* waveguides, i.e., when the cross-section can be considered as unbounded, as for instance optical fibers, immersed pipes, or other guiding devices embedded into a propagating matrix. In such cases, the possible transverse confinement results from particular arrangements of the constituent materials.

We consider here the case of acoustic waves, modeled by the usual *scalar* wave equation. We are concerned with time-harmonic propagation in perturbations or junctions of such uniform waveguides, as illustrated in fig. 1. We assume that each configuration of interest is described by a bounded positive wavenumber function $k : \mathbb{R}^2 \to \mathbb{R}$, so that the acoustic pressure field $u$ satisfies the Helmholtz equation

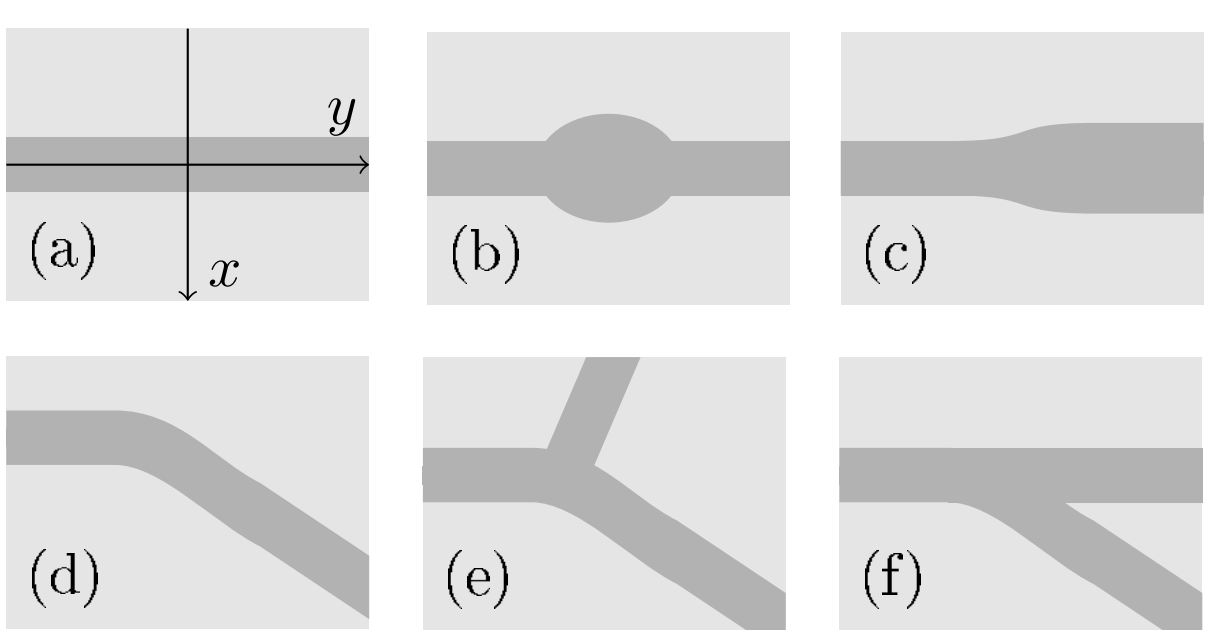

Figure 1: Typical configurations of interest. The two shades of gray represent different materials associated with different constant values of $k$ in (1). (a): uniform waveguide. (b): local perturbation of (a). (c): straight junction of two semi-infinite waveguides. (d): bend junction. (e) and (f): multiple junctions.

$$-\Delta u - k^2\, u = f \quad \text{in } \mathbb{R}^2, \tag{1}$$

where $f$ represents a compactly supported excitation. The questions we address are twofold:
(i) If $f \neq 0$, can we state an *outgoing radiation condition* expressing that $u$ radiates towards infinity? If so, can we prove the existence and uniqueness of an outgoing solution $u$ to (1)?
(ii) If $f = 0$, are there solutions $u \in L^2(\mathbb{R}^2)\setminus\{0\}$ to (1)? Such $u$'s are often called *trapped modes*, since they do not radiate towards infinity (their existence implies non-uniqueness in (i)).

We have explored a *modal* perspective on these issues. We review here the basic ideas of our approach and the main results we have obtained. In short, as regards (i), we can answer positively only in cases (a), (b) and (c) of fig. 1. Concerning (ii), we claim that there are no trapped modes in cases (a) to (e).

## 2 Modal representations

We consider here an infinite uniform waveguide filling $\mathbb{R}^2$ endowed with coordinates $(x, y)$ where $x$ refers to the *transverse* direction and $y$, to the *longitudinal* one. It is characterized by a function $k$ which only depends on the transverse co-ordinate, i.e., $k = k(x)$. We assume that there exist positive constants $k_{\inf}$, $k_{\sup}$, $k_\infty$ and $R$ such that for all $x \in \mathbb{R}$,

$$k_{\inf} \le k(x) \le k_{\sup} \quad \text{and} \tag{2}$$
$$|x| > R \implies k(x) = k_\infty. \tag{3}$$

In order to construct a *modal* representation of solutions $u$ to (1), our parti pris is to *diagonalize* the $x$-dependent part of the Helmholtz operator, i.e., the unbounded *selfadjoint* operator $A : L^2(\mathbb{R}) \to L^2(\mathbb{R})$ with domain $H^2(\mathbb{R})$ defined by

$$A := -\partial_x^2 - k^2. \tag{4}$$

This leads us to a separation of variables process (hereafter noted as SVP) first presented below for a homogeneous medium, then for a general open waveguide.

### 2.1 SVP for a homogeneous medium

We assume here that $k(x) = k_\infty$ for all $x \in \mathbb{R}$. Consider the usual Fourier transform $\mathcal{F}_\infty$ defined by

$$\forall \xi \in \mathbb{R}, \quad \mathcal{F}_\infty \varphi(\xi) := \int_{\mathbb{R}} \varphi(x)\, \overline{\Phi_\xi^\infty(x)}\, \mathrm{d}x, \tag{5}$$
$$\text{where} \quad \Phi_\xi^\infty(x) := \mathrm{e}^{\mathrm{i}\xi x}/\sqrt{2\pi}.$$

It is well-known that $\mathcal{F}_\infty$ is unitary from $L^2(\mathbb{R}_x)$ to $L^2(\mathbb{R}_\xi)$. Its inverse expresses as follows:

$$\forall x \in \mathbb{R}, \quad \mathcal{F}_\infty^{-1} \widehat{\varphi}(x) = \int_{\mathbb{R}} \widehat{\varphi}(\xi)\, \Phi_\xi^\infty(x)\, \mathrm{d}\xi. \tag{6}$$

Moreover, $\mathcal{F}_\infty$ *diagonalizes* $A_\infty := -\partial_x^2 - k_\infty^2$ in the sense that

$$\mathcal{F}_\infty A_\infty = \lambda_\xi\, \mathcal{F}_\infty \text{ where } \lambda_\xi := \xi^2 - k_\infty^2, \tag{7}$$

which means that $\mathcal{F}_\infty$ converts the action of $A_\infty$ into the multiplication by $\lambda_\xi$. This property is related to the fact that for all $\xi \in \mathbb{R}$, $\Phi_\xi^\infty$ satisfies the eigenvalue equation

$$(-\partial_x^2 - k_\infty^2)\, \Phi_\xi^\infty = \lambda_\xi\, \Phi_\xi^\infty \quad \text{in } \mathbb{R}. \tag{8}$$

As $\Phi_\xi^\infty \notin L^2(\mathbb{R})$, it is called a *generalized eigenfunction* of $A_\infty$ associated with $\lambda_\xi$. When $\xi$ varies in $\mathbb{R}$, the spectral parameter $\lambda_\xi$ covers the spectrum $[-k_\infty^2, +\infty)$ of $A_\infty$, which is purely *essential* spectrum. The family $\{\Phi_\xi^\infty \mid \xi \in \mathbb{R}\}$ can be seen as a generalized spectral basis of $A_\infty$. From (5) and (6), we see that $\mathcal{F}_\infty$ appears as a *decomposition* operator on this family, which furnishes the spectral components of a given function $\varphi$, whereas $\mathcal{F}_\infty^{-1}$ acts as a *recomposition* operator, which allows us to reconstruct a function $\varphi$ from its spectral components.

Suppose that $u$ satisfies (1) and consider the partial Fourier transform of $u$ and $f$ in the $x$-direction, i.e., $\widehat{u}_\xi(y) := \mathcal{F}_\infty[u(\cdot\,, y)](\xi)$ and $\widehat{f}_\xi(y) := \mathcal{F}_\infty[f(\cdot\,, y)](\xi)$. As $\Delta = \partial_x^2 + \partial_y^2$, we deduce from (7) that for almost every $\xi \in \mathbb{R}$,

$$-\partial_y^2 \widehat{u}_\xi + \lambda_\xi\, \widehat{u}_\xi = \widehat{f}_\xi \quad \text{in } \mathbb{R}. \tag{9}$$

When $f = 0$, the general solution to (9) is a linear combination of $\mathrm{e}^{-\sqrt{\lambda_\xi}\, y}$ and $\mathrm{e}^{+\sqrt{\lambda_\xi}\, y}$. Using the convention

$$\sqrt{\lambda_\xi} := -\mathrm{i}\, \sqrt{-\lambda_\xi} \quad \text{if } \lambda_\xi < 0, \tag{10}$$

we say that $\mathrm{e}^{-\sqrt{\lambda_\xi}\, y}$ is *outgoing* at $+\infty$ since it is either evanescent (i.e., tends to 0 as $y \to +\infty$) if $\lambda_\xi > 0$ (i.e., $|\xi| > k_\infty$), or propagating towards[1] $+\infty$ if $\lambda_\xi < 0$ (i.e., $|\xi| < k_\infty$). Similarly, $\mathrm{e}^{+\sqrt{\lambda_\xi}\, y}$ is *outgoing* at $-\infty$.

When $f \neq 0$, a particular solution to (9) is given by $\gamma_\xi * \widehat{f}_\xi$, i.e.,

$$\widehat{u}_\xi(y) = \int_{\mathbb{R}} \gamma_\xi(y - y')\, \widehat{f}_\xi(y')\, \mathrm{d}y', \tag{11}$$

where $\gamma_\xi(y) := \mathrm{e}^{-\sqrt{\lambda_\xi}\, |y|}/(2\sqrt{\lambda_\xi})$ satisfies (9) with right-hand side $\delta$ (Dirac measure at $y = 0$). Formula (11) actually provides the only solution to (9) which is *outgoing on both sides* $\pm\infty$ in the sense that

$$\widehat{u}_\xi(y) = \widehat{\alpha}_\xi^\pm\, \mathrm{e}^{\mp\sqrt{\lambda_\xi}\, y} \quad \text{if } \pm y > \rho, \tag{12}$$
$$\text{where } \widehat{\alpha}_\xi^\pm := \int_{-\rho}^{\rho} \frac{\mathrm{e}^{\pm\sqrt{\lambda_\xi}\, y'}}{2\sqrt{\lambda_\xi}}\, \widehat{f}_\xi(y')\, \mathrm{d}y'$$

and $\rho > 0$ is such that $\operatorname{supp}(f) \subset \{|y| < \rho\}$. Using the inverse Fourier transform then yields

$$u = \int_{\mathbb{R}} \widehat{\alpha}_\xi^\pm\, w_\xi^{\infty,\pm}\, \mathrm{d}\xi \quad \text{in } \{\pm y > \rho\}, \tag{13}$$
$$\text{where } w_\xi^{\infty,\pm}(x, y) := \Phi_\xi^\infty(x)\, \mathrm{e}^{\mp\sqrt{\lambda_\xi}\, y}.$$

This is a *modal representation* of $u$: in each half-plane $\{\pm y > \rho\}$, $u$ appears as a continuous superposition of *$y$-outgoing modes* $w_\xi^{\infty,\pm}$ (which satisfy (1) in $\mathbb{R}^2$ for $f = 0$), with density $\widehat{\alpha}_\xi^\pm$.

[1] This interpretation follows from an implicit time-dependence in $\mathrm{e}^{-\mathrm{i}\omega t}$, where $\omega$ is the angular frequency.

### 2.2 Case of a general waveguide

In order to obtain a similar representation under the general assumptions (2)-(3), the price to pay is to introduce a generalization $\mathcal{F}$ of the Fourier transform $\mathcal{F}_\infty$, which diagonalizes $A$ (see (4)) in the sense that (7) becomes

$$\mathcal{F}\, A = \lambda_\xi\, \mathcal{F}, \tag{14}$$

*with the same definition of* $\lambda_\xi := \xi^2 - k_\infty^2$. If we are able to do so, the Helmholtz equation (1) in the waveguide will be transformed into the same equation (9), so that the same SVP will apply, despite the inhomogeneities of the medium.

The idea is to search for $\mathcal{F}$ in a form analogous to (5) by constructing a family of generalized eigenfunctions of $A$ which will be denoted

$$\{\Phi_\xi \mid \xi \in \mathbb{R} \cup \mathbb{G}\}. \tag{15}$$

The $\Phi_\xi$'s for $\xi \in \mathbb{R}$ are associated with the *essential* spectrum $[-k_\infty^2, +\infty)$ of $A$, which is the same as $A_\infty$ since both operators differ from a relatively compact perturbation. But unlike $A_\infty$, $A$ may have a *discrete* spectrum. The latter is precisely related to the possible existence of *guided* waves, which justifies our use of the word *waveguide* and the notation $\mathbb{G}$, which is made precise below. The $\Phi_\xi$'s for $\xi \in \mathbb{G}$ are eigenfunctions associated with the eigenvalues of $A$.

The construction of the $\Phi_\xi$'s for $\xi \in \mathbb{R}$ is based on a perturbation argument related to scattering theory. Considering each $\Phi_\xi^\infty$ as an *incident wave*, we search for a *scattered wave* $\Phi_\xi^{\rm scat}$ such that $\Phi_\xi := \Phi_\xi^\infty + \Phi_\xi^{\rm scat}$ satisfies $A\, \Phi_\xi = \lambda_\xi\, \Phi_\xi$. Using (4) and (8), this means that

$$-\partial_x^2 \Phi_\xi^{\rm scat} - (\lambda_\xi + k^2)\, \Phi_\xi^{\rm scat} = (k^2 - k_\infty^2)\, \Phi_\xi^\infty$$

in $\mathbb{R}$ (recall that $k^2 - k_\infty^2$ is compactly supported, see (3)). Its solution is uniquely defined if we further assume that $\Phi_\xi^{\rm scat}$ is outgoing (i.e., $\Phi_\xi^{\rm scat}(x) = C_\xi^\pm\, {\rm e}^{{\rm i}|\xi||x|}$ as $x \to \pm\infty$).

The remaining $\Phi_\xi$'s are associated with the discrete spectrum $\Lambda_{\rm disc}$ of $A$, which is a (possibly empty) finite set. More precisely,

$$\begin{array}{ll} \Lambda_{\rm disc} \subset [-k_{\sup}^2, -k_\infty^2) & \text{if } k_{\sup} > k_\infty, \\ \Lambda_{\rm disc} = \varnothing & \text{otherwise.} \end{array}$$

The min-max principle provides sufficient conditions for the existence of eigenvalues. For every $\lambda \in \Lambda_{\rm disc}$, one can choose $\xi \in {\rm i}\mathbb{R}$ such that $\lambda = \lambda_\xi$. We denote $\mathbb{G}$ the finite subset of ${\rm i}\mathbb{R}$ thus obtained, and $\{\Phi_\xi \mid \xi \in \mathbb{G}\}$ a family of normalized eigenfunctions.

For conciseness, we use the somewhat abusive notation

$$\int_{\mathbb{R}\cup\mathbb{G}} F(\xi)\, {\rm d}\xi := \int_{\mathbb{R}} F(\xi)\, {\rm d}\xi + \sum_{\xi\in\mathbb{G}} F(\xi)$$

for any function $F : \mathbb{R} \cup \mathbb{G} \to \mathbb{C}$, and denote by $L^2(\mathbb{R} \cup \mathbb{G})$ the associated space of square integrable functions. The following theorem can be extended to numerous scattering problems [13].

**Theorem 1** (i) *The generalized Fourier transform defined by*

$$\forall \xi \in \mathbb{R} \cup \mathbb{G},\ \mathcal{F}\varphi(\xi) := \int_{\mathbb{R}} \varphi(x)\, \overline{\Phi_\xi(x)}\, {\rm d}x$$

*is unitary from* $L^2(\mathbb{R})$ *to* $L^2(\mathbb{R} \cup \mathbb{G})$.
(ii) *It diagonalizes* $A$ *in the sense of* (14).
(iii) *Its inverse expresses as follows:*

$$\forall x \in \mathbb{R},\ \mathcal{F}^{-1}\widehat{\varphi}(x) = \int_{\mathbb{R}\cup\mathbb{G}} \widehat{\varphi}(\xi)\, \Phi_\xi(x)\, {\rm d}\xi$$

Considering a solution $u$ to (1) and setting $\widehat{u}_\xi(y) := \mathcal{F}[u(\cdot\,, y)](\xi)$ and $\widehat{f}_\xi(y) := \mathcal{F}[f(\cdot\,, y)](\xi)$, we can follow the same SVP as in §2.1. Formula (11) provides again the only $y$-*outgoing* solution to (9). Using the above expressions of $\mathcal{F}$ and $\mathcal{F}^{-1}$, as well as Fubini's theorem, we infer that the only $y$-*outgoing* solution to (1) is given by the integral representation

$$u(M) = \mathcal{G}f(M) := \int_{\mathbb{R}^2} G(M, M') f(M') {\rm d}M' \tag{16}$$

for all $M \in \mathbb{R}^2$, where $G$ is the *outgoing Green's function* of the waveguide defined for all $M = (x, y)$ and $M' = (x', y')$ by

$$G(M, M') := \int_{\mathbb{R}\cup\mathbb{G}} \frac{{\rm e}^{-\sqrt{\lambda_\xi}\,|y-y'|}}{2\sqrt{\lambda_\xi}}\, \Phi_\xi(x)\, \overline{\Phi_\xi(x')}\, {\rm d}\xi.$$

Moreover, on both sides of the support of $f$, (12) remains valid so that (13) becomes

$$u = \int_{\mathbb{R}\cup\mathbb{G}} \widehat{\alpha}_\xi^\pm\, w_\xi^\pm\, {\rm d}\xi \quad \text{in } \{\pm y > \rho\}, \tag{17}$$
$$\text{where } w_\xi^\pm(x, y) := \Phi_\xi(x)\, {\rm e}^{\mp\sqrt{\lambda_\xi}\, y}.$$

This modal representation now includes a discrete sum of *guided* modes $w_\xi^\pm$ (for $\xi \in \mathbb{G}$) which do not decay as $|y| \to \infty$.

The rigorous justification of the above formal approach raises a technical obstacle related to the generalized Fourier transform $\mathcal{F}$. Indeed, for a given $y \in \mathbb{R}$, $u(\cdot, y) \notin L^2(\mathbb{R})$ in general, while Theorem 1 defines $\mathcal{F}$ in $L^2(\mathbb{R})$. So we need to extend $\mathcal{F}$ to a larger space. A natural way to do so is to mimic the duality argument used to extend the usual Fourier transform to the Schwartz space of distributions, in order to build a new distribution space adapted to the operator $A$ (see (4)). We **conjecture** that this can be achieved, which is done in [11] in the particular case where $k$ is two-valued. The general case remains an open question.

## 3 Perturbation of a uniform waveguide

Consider the configuration (b) of fig. 1, i.e., a local perturbation of a uniform waveguide (for which we keep the notations of §2). It is represented by a bounded function $k_* = k_*(x, y)$ such that $k_* - k$ is compactly supported, say in $\{|y| < \rho\}$. The Helmholtz equation

$$-\Delta u - k_*^2\, u = f \quad \text{in } \mathbb{R}^2 \tag{18}$$

can be rewritten equivalently

$$-\Delta u - k^2\, u = (k_*^2 - k^2)\, u + f \quad \text{in } \mathbb{R}^2, \tag{19}$$

where the right-hand side is still compactly supported. Hence the modal approach of §2 applies, which suggest to use (17) or its equivalent form (12) as a *radiation condition*: we will say that a solution to (18) is *y-outgoing* if there exists two functions $\mathbb{R} \cup \mathbb{G} \ni \xi \mapsto \widehat{\alpha}_\xi^\pm$ such that

$$\widehat{u}_\xi(y) = \widehat{\alpha}_\xi^\pm\, \mathrm{e}^{-\sqrt{\lambda_\xi}\,|y|} \quad \text{for } \pm y > \rho. \tag{20}$$

**Theorem 2** ( **[2]**) *Under the above-mentioned conjecture, for $f \in L^2(\mathbb{R}^2)$ compactly supported, (18) has a unique solution $u \in H^1_{\mathrm{loc}}(\mathbb{R}^2)$ which is y-outgoing in the sense of (20).*

*Proof.* Choose a *bounded* domain $\Omega \subset \mathbb{R}^2$ that contains the supports of $k_* - k$ and $f$, and denote by $\mathcal{K} : L^2(\Omega) \to L^2(\Omega)$ the operator defined by $\mathcal{K}u = \mathcal{G}\big((k_*^2 - k^2)u\big)|_\Omega$, where $\mathcal{G}$ is given in (16). It is then clear that (19)-(20) is equivalent to the Lippmann–Schwinger equation

$$(I - \mathcal{K})\, u = (\mathcal{G}f)|_\Omega \quad \text{in } L^2(\Omega). \tag{21}$$

It can be shown that $\mathcal{K}$ is compact, so that the existence of $u$ follows from its uniqueness.

To prove uniqueness, the first step is based on an *energy* argument. Assume $f = 0$ in (18). Multiplying this equation by $\overline{u}$, integrating over the strip $\mathbb{R} \times (-\rho, \rho)$, using *formally* Green's formula and taking the imaginary part yields

$$E^+(u) + E^-(u) = 0 \quad \text{where} \tag{22}$$
$$E^\pm(u) := \operatorname{Im} \int_{\mathbb{R}} \frac{\partial u}{\partial |y|}(x, \pm\rho)\, \overline{u(x, \pm\rho)}\, \mathrm{d}x$$

represents the outgoing energy flux across the transverse section $\{y = \pm\rho\}$. As $\mathcal{F}$ is unitary, we can transform this integral by a Plancherel-like identity. From (20), we obtain

$$E^\pm(u) = \int_{\mathbb{R}\cup\mathbb{G}} \operatorname{Im}\big(-\sqrt{\lambda_\xi}\big)\, \Big|\widehat{\alpha}_\xi^\pm\, \mathrm{e}^{-\sqrt{\lambda_\xi}\,\rho}\Big|^2\, \mathrm{d}\rho.$$

As $\operatorname{Im}\big(-\sqrt{\lambda_\xi}\big)$ vanishes if $\lambda_\xi > 0$ and is positive otherwise (see (10)), we deduce that $E^\pm(u) \geq 0$. So (22) implies that $E^\pm(u) = 0$, thus $\widehat{\alpha}_\xi^\pm = 0$ for $|\xi| < k_\infty$ and $\xi \in \mathbb{G}$. In other words, *the spectral components of $u$ associated with all y-propagating modes vanish on both sides.* By (20), we have proved that $\widehat{u}_\xi(y) = 0$ for $\pm y > \rho$ and $\xi \in (-k_\infty, k_\infty) \cup \mathbb{G}$.

To prove that this results holds true for all $\xi \in \mathbb{R}$ with $|\xi| > k_\infty$, we use an ingenious *analyticity* argument inspired from [17]. We start from the integral representation (16), which writes equivalently in the Fourier variable as

$$\widehat{u}_\xi(y) = \int_{\mathbb{R}^2} \frac{\mathrm{e}^{-\sqrt{\lambda_\xi}\,|y-y'|}}{2\sqrt{\lambda_\xi}}\, \overline{\Phi_\xi(x')}\, f_*(x', y')\, \mathrm{d}x'\mathrm{d}y',$$

where $f_* := (k_*^2 - k^2)\, u$. On the one hand, the function $\mathbb{R} \ni \xi \mapsto \sqrt{\lambda_\xi}$ has an analytic extension in $\mathbb{C} \setminus \mathbb{B}$ where $\mathbb{B}$ is the union of two branch cuts with endpoints $\pm k_\infty$. On the other hand, the theory of scattering resonances [9] tells us that for any given $x' \in \mathbb{R}$, the function $\mathbb{R} \ni \xi \mapsto \overline{\Phi_\xi(x')}$ has a *meromorphic* extension in $\mathbb{C} \setminus \mathbb{B}$. As $f_*$ is compactly supported, we infer that $\widehat{u}_\xi(y)$ has a *meromorphic* extension in $\mathbb{C} \setminus \mathbb{B}$ for all $y \in \mathbb{R}$. As we already know that $\widehat{u}_\xi(y) = 0$ for $\xi \in (-k_\infty, k_\infty)$, we deduce from the principle of isolated zeros that $\widehat{u}_\xi(y) = 0$ for $\xi \in \mathbb{R} \cup \mathbb{G}$ and $\pm y > \rho$. Applying $\mathcal{F}^{-1}$ yields $u(x, y) = 0$ for $x \in \mathbb{R}$ and $\pm y > \rho$, so $u = 0$ in $\mathbb{R}^2$ by the unique continuation principle, which completes the uniqueness proof. □

No more conjecture is needed if we address the question of the possible existence of trapped modes. Indeed, we have:

**Theorem 3** ( [12]) *For $f = 0$, the only solution $u \in L^2(\mathbb{R})$ to* (18) *is $u = 0$.*

*Proof.* As we assume that $u \in L^2(\mathbb{R})$, we know that for a.e. $y \in \mathbb{R}$, $u(\cdot, y) \in L^2(\mathbb{R})$, so its generalized Fourier transform $\widehat{u}_\xi(y)$ is now well-defined. As $\mathcal{F}$ is unitary, we have

$$\|u\|^2_{L^2(\mathbb{R}^2)} = \int_{\mathbb{R}} \|\widehat{u}_\bullet(y)\|^2_{L^2(\mathbb{R}\cup\mathbb{G})}\,\mathrm{d}y.$$

Hence Fubini's theorem tells us that for a.e. $\xi \in \mathbb{R}\cup\mathbb{G}$, $\widehat{u}_\xi(\cdot) \in L^2(\mathbb{R})$. Applying the SVP decribed in §2, we infer from this Plancherel-type argument that if $|y| > \rho$, the modal components of $u$ associated with *propagating* modes must vanish (i.e., $\widehat{u}_\xi(y) = 0$ for $\xi \in (-k_\infty, k_\infty) \cup \mathbb{G}$), which is exactly what the *energy* argument gave us in the proof of Theorem 2. The *analyticity* argument can then be repeated, which completes the proof. □

## 4 Junctions of waveguides

Consider now the case (c) of fig. 1, i.e., a straight junction of two uniform waveguides characterized by two functions $k^\pm = k^\pm(x)$. On each side of the junction, the associated generalized Fourier transform $\mathcal{F}^\pm : L^2(\mathbb{R}) \to L^2(\mathbb{R} \cup \mathbb{G}^\pm)$ provides a radiation condition as in (20). In [5], we proved the well-posedness of the corresponding radiation problem again in the case where both functions $k^\pm$ are two-valued (because of the conjecture of §2). The general ideas of the proof are the same as for Theorem 2, but their implementation is far more intricate, because of the simultaneous use of $\mathcal{F}^-$ and $\mathcal{F}^+$. The generalization of our approach to the other configurations of junctions in fig. 1 remains open.

On the other hand, the use of modal representation in half-planes proves to be a powerful tool to address the question of trapped modes in such junctions. We focus below on the case (d) of fig. 1 for which the absence of trapped modes follows from the following theorem (see [3]) and the unique continuation principle. Recently, we extended our approach to case (e) of fig. 1 (see [1]). But this approach no longer applies in case (f), since we can no longer use modal representations as soon as the angle between two waveguides in smaller than $\pi/2$.

**Theorem 4** *Let $\Omega_\theta := \{(x,y) \in \mathbb{R}^2 \,|\, y > -|x| \tan\theta\}$ for a given $\theta \in (0, \pi/2)$. If $u \in L^2(\Omega_\theta)$ satisfies $\Delta u + k_\infty^2\, u = 0$ in $\Omega_\theta$, then $u = 0$.*

*Proof.* In the limit case $\theta = 0$, the SVP of §2 yields the following modal representation

$$u = \int_{|\xi|>k_\infty} \widehat{\alpha}_\xi\, w_\xi^{\infty,+}\,\mathrm{d}\xi \quad \text{in } \Omega_0, \tag{23}$$

where the $w_\xi^{\infty,+}$'s are defined as in (13) and $\widehat{\alpha}_\xi := \mathcal{F}_\infty[u(\cdot, 0)](\xi)$ must vanish for $|\xi| < k_\infty$ thanks to the same Plancherel-type argument as in the proof of Theorem 3. In other words, (23) only involves evanescent modes.

As $\theta > 0$, we can prove in addition that $\widehat{\alpha}_\xi = 0$ also if $|\xi| > k_\infty$ using again an *analyticity* argument. This will complete the proof, since (23) will imply that $u = 0$ in $\Omega_0$, then also in $\Omega_\theta$ by the unique continuation principle. The idea is to rewrite the definition of $\widehat{\alpha}_\xi$ as

$$\widehat{\alpha}_\xi := \int_{\mathbb{R}} u(x,0)\,\overline{\Phi_\xi^\infty(x)}\,\mathrm{d}x = \int_{-\infty}^{0} \ldots + \int_{0}^{+\infty} \ldots$$

and express $u(x, 0)$ in both integrals using modal representations similar to (23) in half-planes rotated of angle $\pm\theta$, which yields

$$u(x,0) = \int_{|\eta|>k_\infty} \widehat{\beta}_\eta^\pm\, w_\eta^{\infty,+}(x\cos\theta, \pm x \sin\theta)\,\mathrm{d}\eta$$

if $\pm x > 0$, for some functions $\eta \mapsto \widehat{\beta}_\eta^\pm$. A simple computation then shows that $\widehat{\alpha}_\xi$ extends to a analytic function in a complex vicinity of the real $\xi$-axis. As $\widehat{\alpha}_\xi = 0$ for $\xi \in (-k_\infty, k_\infty)$, the conclusion follows. □

## 5 Related works

Some of the results presented above hold true for higher dimensions, namely the absence of trapped modes in cases (a), (b) and (d) of fig. 1 (see [3, 12]). But mainly because of the conjecture mentioned at the end of §2, we are unable to extend our results for the radiation problem.

Other radiation conditions for perturbed uniform waveguides can be found in the literature [6, 7, 14, 19]. They are based on a wave splitting which amounts to distinguishing the discrete and continuous parts in our modal representation (17). The discrete part is a superposition of outgoing guided modes which can be dealt with separately, whereas the remaining *free wave* associated with the continuous part is assumed to satisfy a usual Sommerfeld-type radiation condition. We conjecture that these radiation conditions are equivalent to ours. Unfortunately, the uniqueness proofs proposed in [6, 19] are false. An alternative very powerful technique for tackling the same issues as those considered here is the commutator method developed in mathematical physics in the early 80's (see, e.g., [8, 18, 20]). Let us mention finally the recent work [10] that considers a 2D straight junction as in §4 and exploits the approach of [16], which seems to provide very general results for waveguide networks using intricate tools of microlocal analysis.

Modal representations such as (17) play a key role in the so-called half-space matching method [4] developed in POEMS lab for solving scattering problems in unbounded media. From a numerical point of view, the idea is to couple half-space representations in overlapping half-spaces with a Finite Element discretization in a bounded region. The most recent development (in the framework of Sarah Al Humaikani's thesis) concerns precisely the junction of waveguides (see also [15]).

# The DFDM: a Promising Approach for Seismic Wavefield Computations in Global Seismology

Barbara Romanowicz[1,*], Yuan Tian[1], Muaaz Al Gawan[2], Dorian Soergel[1], Chao Lyu[3], Jack Deslippe[2]

[1]Department of Earth and Planetary Science, University of California at Berkeley, USA
[2]Lawrence Berkeley National Laboratory, NERSC, Berkeley (CA), USA
[3]University of Chinese Academy of Sciences, China

*Email: barbarar@berkeley.edu

**Abstract**

Imaging of the earth's interior using seismic waves presents specific challenges, due to the very uneven distribution of sources (natural earthquakes) and receivers. In the past forty years, significant progress has been made in the deployment of high dynamic range, low noise digital broadband and very broadband sensors around the world, and the development of publicly accessible on-line archives. Owing to this, and to access to increasingly powerful super-computers, attention has progressively shifted from classical "travel time" tomography to the use of full three-component waveforms. This requires the computation of the predicted wavefield in a complex 3D earth, which, in the last 15 years has benefited from the availability of codes based on the Spectral Element Method (SEM). While the SEM has led to significant gains in resolution of the earth's deep structure, it also has some limitations. We present the status of our on-going efforts to implement and optimize a promising new approach, the Distributional Finite Difference Method, for the computation of the seismic wavefield in complex 3D media, in particular in the context of global seismic imaging.

## 1 Background and Motivation

Global imaging of Earth's mantle elastic and anelastic structure is particularly challenging because of the non-uniform distribution of sources (earthquakes) and seismic stations around the globe. Traditionally, global seismic tomography has relied on the measurement of secondary observables: 1) travel time measurements of a small number of seismic "body wave phases", corresponding to signals in the records that are well-separated from others, often exhibit a high signal-to-noise ratio, and at large, "teleseismic", distances, illuminate the earth's lower mantle and core e.g. [1]; 2) group and phase velocity dispersion measurements of surface waves, that allow the illumination of the upper mantle in vast areas without earthquake sources or stations, in particular in the oceans. These secondary observables are typically related to the earth model through the sensitivity kernels based on high frequency approximations (e.g. Ray theory for body waves, "path average approximation" for surface waves). This computationally very efficient approach led to the development of several generations of 3-dimensional (3D) earth models e.g. [2,3] that informed our present knowledge of prominent large scale features of whole mantle structure.

Full waveform inversion (FWI), was introduced early on in global seismology [4]. It allows one to exploit more of the information contained in seismic records, but necessitates the computation of the complete forward seismic wavefield in earth models in which structure varies not only with depth, but also laterally (i.e. 3D models). Early work relied on first order normal mode first order perturbation theory, and more specifically on an asymptotic approximation which only considers along-branch mode coupling by lateral heterogeneity, and is only appropriate for surface waves and for the determination of long wavelength smooth structure. Later, across-branch mode coupling was added, albeit also asymptotically ( [5] , figure 1). This allowed a more accurate treatment of body wave sensitivity and led to the development of the first global whole mantle shear velocity models based entirely on waveform inversion, within the limitations of normal mode first order perturbation theory [6, e.g.].

More accurate predictions of the seismic wavefield in a realistic, complex, 3D earth model require the numerical integration of the wave equation, which becomes computationally challenging at increasing frequencies, especially for targets in the deep mantle, because of the long paths involved with sources and stations confined near Earth's surface.

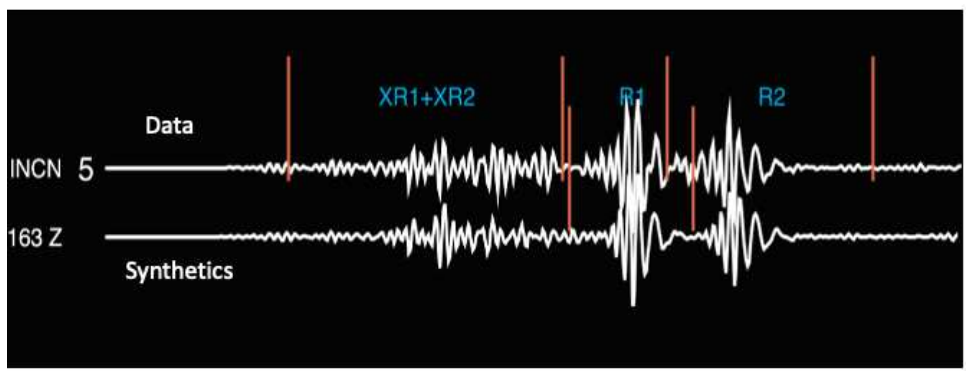

Figure 1: Comparison of observed (top) and predicted (bottom) vertical component records, showing overtone (XR1+XR2) and fundamental mode minor (R1) and major (R2) arc surface waves for a $\sim$ M 7 earthquake observed at an epicentral distance distance of $\sim 130^o$. The waveforms have been low pass filtered with a cut-off at 60 s period. Red bars mark limits of windows that are used for weighting different parts of the record. Here the synthetics have been computed using normal mode perturbation theory including cross-branch coupling [5] in model $SAW24B16$ ( [6]. At these periods, normal mode perturbation theory provides a conveniently efficient approach for predicting the seismic wavefield.

Over the last 15 years, the spectral element method (SEM) has become the method of choice for forward modeling of the wavefield at teleseismic distances, because of its relatively manageable computational cost (due to a diagonal mass matrix), and its accuracy, in particular for the modeling of surface waves, essential in global seismic tomography. SEM was first introduced to global seismology almost 25 years ago [7, 8]. Since then, it has been developed into a popular tool used widely in the earthquake seismology community [9] and has led to higher resolution mantle images at global and continental scales. In particular, the long-debated existence of deeply rooted mantle "plumes" as the cause of mid-plate volcanism has now been confirmed. These plumes extend quasi-vertically from the core-mantle boundary across the lower mantle, making their way to the earth's surface in the vicinity of "hotspot" volcanoes, such as Hawaii in the Pacific ( [10], Figure 2) or Cape Verde in the Atlantic.

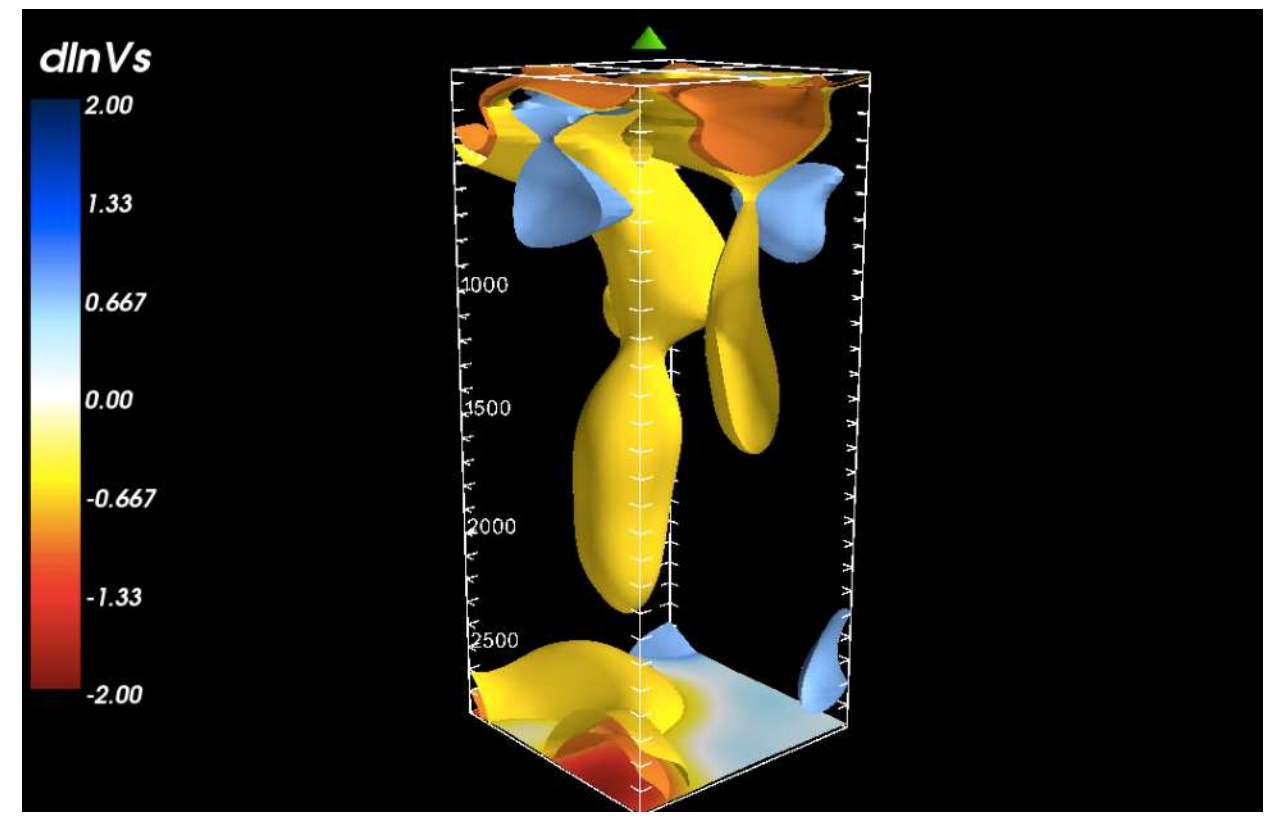

Figure 2: 3D view of the Hawaian plume in model $SEMUCB_WM1$ [10].Relative variations of shear wave velocity (dlnVs, in percent) with respect to global radially symmetric average from 50 km depth to the core-mantle boundary at 2900 km depth. The location of the active Hawaian hotspot volcano is shown by the green cone. Only contours corresponding to shear velocity variations larger than 0.5% (in absolute value) are shown. Low shear velocities indicate hotter than average material.

As increasingly higher spatial resolution is sought to investigate the nature of remote, small scale objects of geodynamic interest, the complexity of meshing, accuracy of the wavefield, and soaring computational time at high frequencies represent increasingly significant challenges, justifying continued efforts at designing the next generation of numerical wave propagation codes, beyond the SEM. There is a clear need for further development of flexible methodologies that take advantage of the variable resolution required in different locations of an earth 3D model, while keeping the wavefield computation accurate, including at fluid-solid boundaries, and computational costs competitive with the most efficient codes presently available, such as the SPECFEM suite [9].

The recently proposed "Distributional Finite Difference Method" (DFDM) [12–14] presents some advantages over both the SEM and classical finite-differences. Here we present this method and our efforts to implement it in a scalable HPC environment for applications in global seis-mology and beyond.

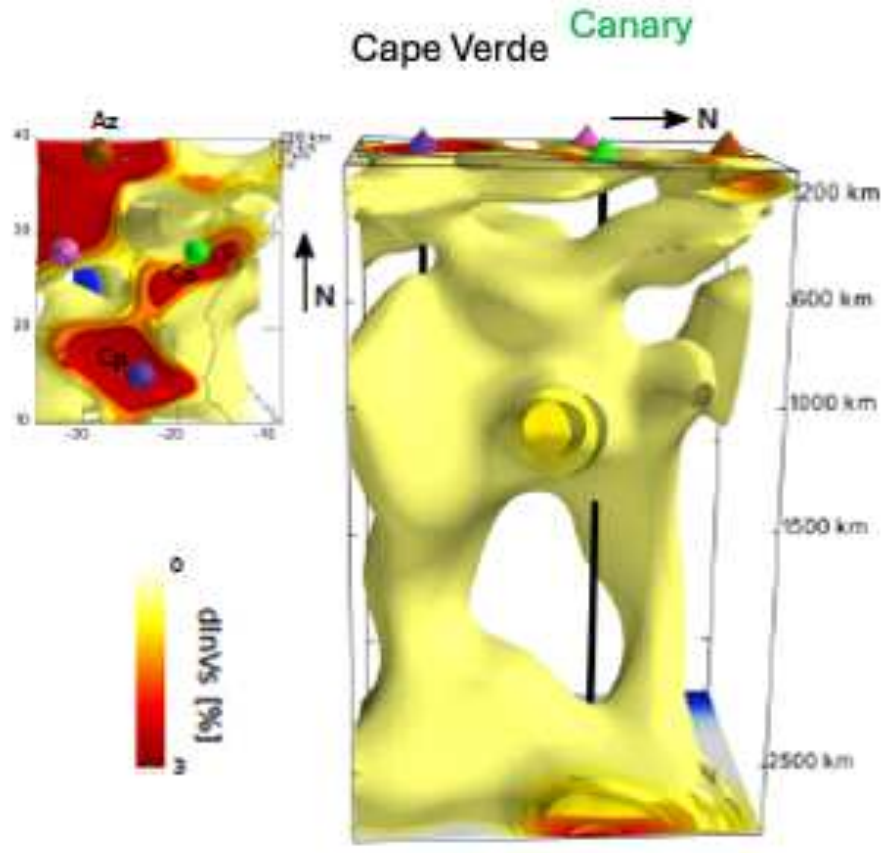

Figure 3: 3D view of the Cape Verde and Canary Plumes in model $SEMATL_23$ [11]. Relative variations of shear wave velocity (dlnVs, in percent) with respect to global radially symmetric average, from 100 km depth to the core-mantle boundary at 2900 km depth. The location of the Cape Verde, Canary and Azores volcanoes are shown by color cones, also in the left panel, showing a map view from the top at 130 km depth. Only contours corresponding to negative shear velocity variations of absolute value larger than 1.0% are represented.

## 2 The Distributed FD Method

The structure of the DFDM algorithm is very similar to that of finite-differences (FD) [15], with added versatility, in that a unique scheme allows to account for domain decomposition, heterogeneity, curvilinear mesh, anisotropy, non-conformal interfaces, or discontinuous grid and fluid-solid interfaces. The DFDM algorithm accurately and naturally accounts for the free surface and a solid-fluid interface with topography, such as the earth's surface or the core-mantle boundary (CMB). The discrete wavefield is discontinuous between neighboring elements, as in the Discontinuous Galerkin Method, which provides flexibility in specifying boundary conditions. Also, mass and stiffness matrices are band diagonal and accurate. In DFD, the accuracy is independent of the size of the elements, it is adjusted by varying the number of grid points inside individual elements. Thus, in regions where the model is smooth, the elements can be very large. An important difference with FD is that, in order to evaluate partial derivatives, the usual discrete operator (here in 1D):

$$Df(x+\frac{h}{2}) = \lim_{h\to 0} \frac{f(x+h)-f(x)}{h} \quad (1)$$

is replaced by the discrete equivalent of the distributional derivative:

$$Df(x) = -\int_{-\infty}^{+\infty} f(y)\delta'(y-x)dx \quad (2)$$

which leads to the following expression for the d-th derivative:

$$\begin{aligned} f^{(d)}(x) &= \int_{-\infty}^{+\infty} f^{(d)}(y)\delta(y-x)dy \quad (3)\\ &= (-1)^d \int_{-\infty}^{+\infty} f(y)\delta^{(d)}(y-x)dy \end{aligned}$$

Otherwise, the algorithmic structure is similar to a classical FD, with 1D operators constructed from basis functions with compct support, leading to computations that involve operations on band-diagonal matrices. The basis functions are staggered according to the so-called Lebedev grid, which has 2nd order accuracy in space, and which removes non-physical parasitic modes that appear in centered FD operations [13, 16]. Field variables at a given time are expanded in specific orthonormal bases with local support. The 1D orthonormal basis functions are built from the staggered sets of B-spline basis functions of polynomial order p, although other choices can be made [12]. This is in contrast to SEM, which uses Legendre polynomials with global support. The use of orthogonal basis functions maintains the self-adjoint antisymmetric nature of the wave equation after discretization. The extension from 1D to 2D and 3D is performed via tensor products, using the same implementation of the tensor product as in the SEM.

## 3 Preliminary benchmarks and current status

Preliminary benchmarks of the DFDM algorithm itself (on a single node, before parallelization and optimization) have been performed extensively in 2D global elastic models [14] as well as 3D elastic regional models [17] and have shown that this approach exhibits comparable accuracy to the SEM in the solid domain and superior accuracy in the fluid domain, reduced memory usage, and comparable computational efficiency. Benchmarks show that the superior accuracy of the band-diagonal mass and stiffness matrices in DFDM enables fewer points per wavelength, approaching the spectral limit, and compensates for the increased computational burden due to four Lebedev staggered grids. Figure 4 compares the global mesh required by SEM for the accurate computation of the seismic wavefield in a global 2D elastic reference earth model (in which structure varies only with depth), and the corresponding mesh in DFD.

DFDM offers multiple opportunities for parallelization and optimization, such as leveraging the diagonal banded nature of most operator matrices to reduce the overall computational requirements, independent processing of elements, and efficiently handles the physical quantities exchange between connected elements with simple to solve boundary conditions. Similarly, the compute-bound nature of the overall workflow makes it an ideal candidate for GPU architectures that can tremendously improve the performance, thus making way for even larger simulations.

DFD therefore represents an attractive alternative to SEM-based approaches, once parallelized and optimized. It may be particularly advantageous, possibly in combination with hybrid modeling approaches such as “Box Tomography”, for imaging small targets located in the deep mantle, such as, for example, Ultra Low Velocity Zones (ULVZs), patches of extreme elastic properties detected at the earth’s core-mantle boundary, in particular at the roots of major mantle “plumes”.

Beyond global seismology, once fully developed and optimized, DFDM may be of interest for modeling acoustic and elastic structures at various scales, such as in geophysical prospecting, medical imaging or to predict ground motions generated by large earthquakes.

Recently, we have been working on implementing a C++, parallelized version of elastic DFDM in 2D and 3D Our preliminary benchmarks in a 2D elastic global model, on 4 cores on a Macpro laptop, against the optimized SPECFEM2D [19] SEM code are encouraging (figure 5. Without much effort at optimization, the DFDM run obtains a $\sim 25\%$ wall time reduction compared to the $SEM$ in a standard setup $SEM-1$, at comparable accuracy.

We will present the scientific motivation that drives our work, the potential advantages of the

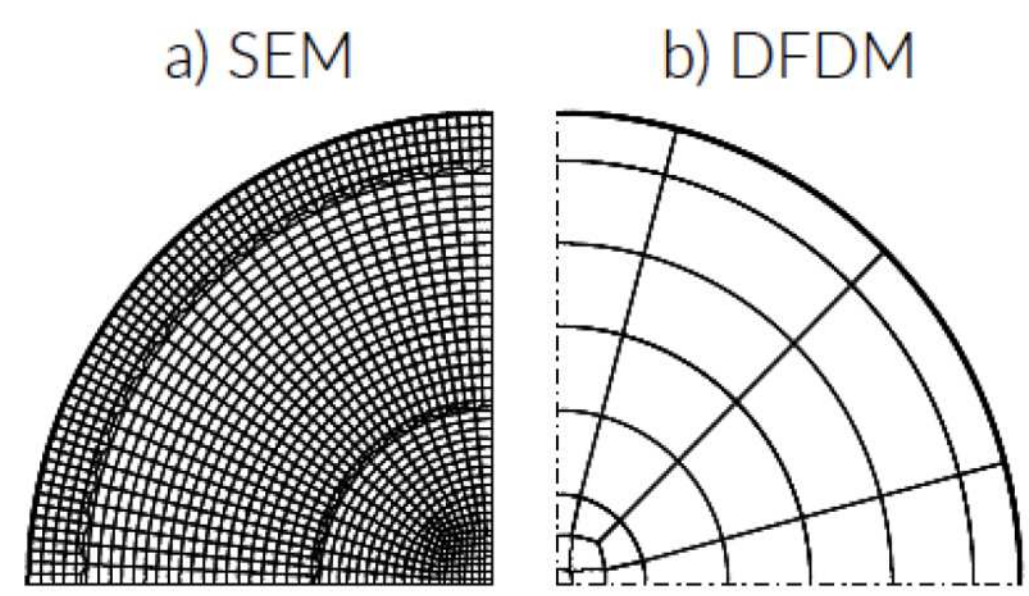

Figure 4: Numerical meshes employed for the SEM (left) and DFDM (right) simulations for a global 2D elastic benchmark in earth model AK135 [18]. In both cases, the mesh matches the solid-solid discontinuities in the model (crust/mantle boundary, 670 km discontinuity) as well as the solid/fluid discontinuities (core-mantle boundary and inner core boundary). In the SEM simulation, the mesh has two doubling layers (under the 670 and in the middle of the outer core), and the number of elements is adjusted according to the desired accuracy. For the DFDM, the number of elements is independent of the desired accuracy. Because DFDM allows for non-conformal interfaces, the number of grid points is adjusted independently inside individual elements. Note that the SEM grid has actually been simplified for clarity of the figure. After [14]

.

DFDM methodology, and provide an update on our efforts to 1) optimize the DFD computation on an HPC platform at NERSC (National Energy Research Scientific Computing Center) in the framework of the NESAP (NERSC Science Acceleration Program) and 2)produce a DFDM-based wave propagation solver ready for application to comparisons with real data. We discuss the next steps and illustrate potential applications of interest, in particular for the imaging of target objects in the deep earth.

Table 1: Simulation parameters and error metrics with respect to SEM-Ref. Ne= number of elements; p = polynomial degree; ppw = min. points per wavelength; $\Delta t$ = timestep; NSTEP = number of time steps in computation; $\sigma = \frac{std(u-u_{ref}}{std(u_{ref}}$; $= \frac{max|u-u_{ref}|}{max|u_{ref}|}$; Pearson correlation coef. $\rho = \frac{cov(u,u_{ref})}{(std(u)std(u_{ref}))}$

| | SEM-Ref | SEM-1 | DFDM | DFDM-hires |
|---|---|---|---|---|
| $N_e$ | 322560 | 81408 | 93 | 93 |
| $p$ | 4 | 4 | 4 | 6 |
| ppw | - | - | 3 | 5 |
| $\Delta t$ (s) | 0.0700 | 0.0700 | 0.0700 | 0.0700 |
| NSTEP | 72000 | 72000 | 72000 | 72000 |
| Time (s) | 3055.0 | 690.0 | 531.0 | 871.9 |
| $\sigma$ | - | 0.2338 | 0.2335 | 0.2309 |
| $\varepsilon$ | - | 0.1773 | 0.2178 | 0.3015 |
| $\rho$ | - | 0.9723 | 0.9725 | 0.9730 |
| Amp diff (%) | - | 4.89 | 4.54 | 1.53 |

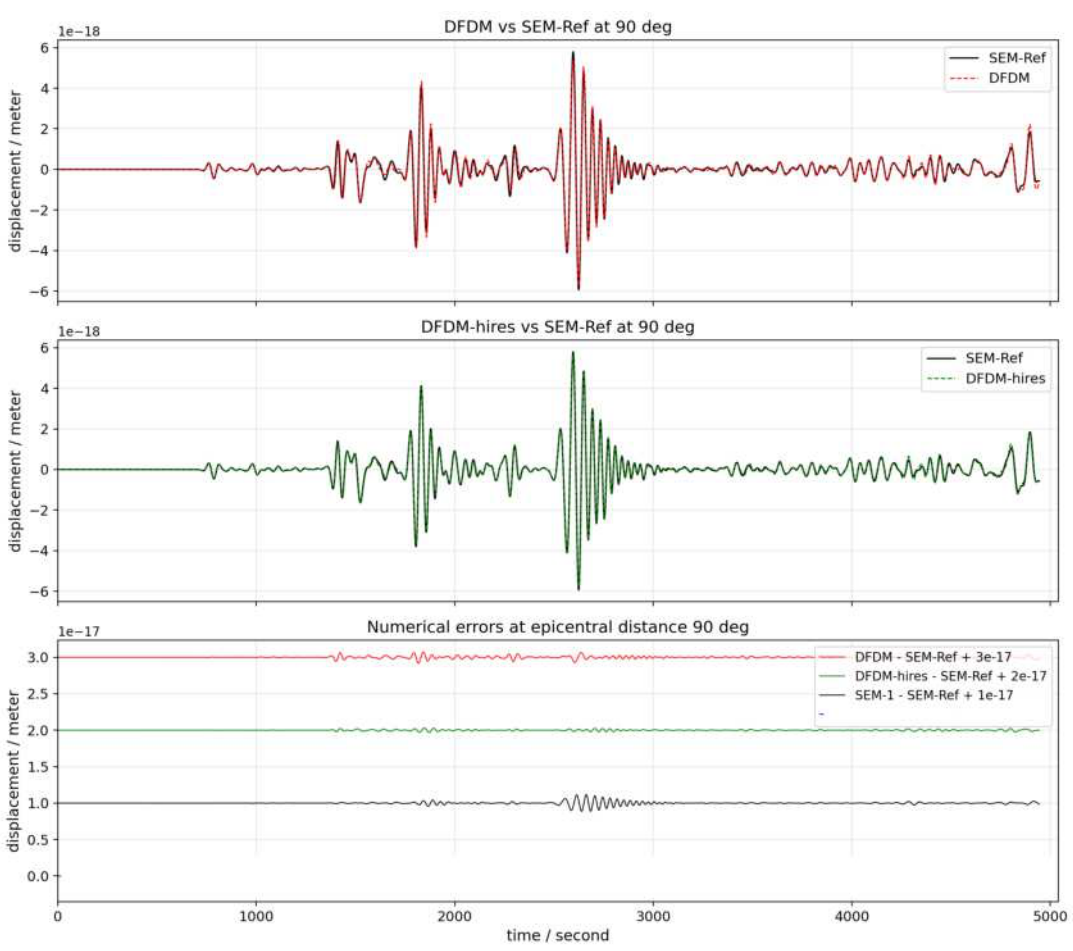

Figure 5: Comparison of vertical component displacement seismograms recorded at a distance of $90^o$ from a moment-tensor source at a depth of 50 km (Ricker wavelet with peak frequency of 0.0125 Hz), computed in a 2D elastic global earth model modified from model $AK135$ [18]. The 2 top panels show one-to-one comparisons of a reference SEM high accuracy computation (SEM-ref) with, respectively, two computations with DFDM (top: DFDM; second row: DFDM-hires). The bottom panel shows the differences between the SEM-ref and 3 other computations: DFDM, DFDM-hires and SEM-1. SEM-ref (broken lines in top panels) is the reference, most accurate SEM computation (see Table 1), while SEM-1 is the standard set-up used for global computations. The parameters of each computation are given in Table 1. While acceptable in practice for comparison with real data, SEM-1 shows the largest error. DFD exhibits slightly larger "errors" than SEM-1 in the body waves, but smaller "errors" in the fundamental mode surface waves. The "higher resolution" DFD-hires exhibits the smallest difference with SEM-ref.

# The extended sampling method for inverse scattering revisited

**Isaac Harris**[1,*]
[1]Department of Mathematics, Purdue University, West Lafayette, USA
*Email: harri814@purdue.edu

## Abstract

We revisit the extended sampling method (ESM) for inverse scattering and provide an alternative formulation with a more rigorous analytical foundation based on the factorization method. Using range characterizations for translated far–field operators associated with a fixed radius sampling disk, we derive practical imaging functionals that can recover the unknown scatterer. In contrast to the original ESM analysis, which is closely tied to linear sampling arguments and a sound–soft (Dirichlet) sampling operator, the proposed approach yields two related indicators built from either sound–soft (Dirichlet) or sound–hard (Neumann) sampling operator.

## 1 Introduction

The inverse scattering problem of determining an unknown scatterer from time–harmonic acoustic and electromagnetic waves has received a lot of attention. This is due to its wide applications in non–destructive testing. Here, we will revisit one of the most recently developed *qualitative reconstruction methods* applied to shape reconstruction. In some cases, it is optimal to use qualitative methods for two reasons: first, iterative methods may require a priori information that may not be readily available, such as the number of regions to be recovered, and second, qualitative methods are computationally cheap to implement. Qualitative methods seek to recover the shape of an object from little a priori information by relating the support of the object to the range of a known operator. One of the trade–offs is that many qualitative methods require multi–static data collected all around the region of interest.

The method we will study is the extended sampling method (ESM) which was first introduced in [4]. This method allows one to recover the location of an unknown scatterer from the measured far–field pattern from one (or multiple) incident directions. For a radiation solution to the Helmholtz equation denoted $u$, the far–field pattern is defined via the asymptotic behavior

$$u(x) = \frac{\mathrm{e}^{\mathrm{i}\pi/4}}{\sqrt{8\pi k}} \cdot \frac{\mathrm{e}^{\mathrm{i}k|x|}}{\sqrt{|x|}} u^{\infty}(\hat{x}) + \mathcal{O}\left(\frac{1}{|x|^{3/2}}\right)$$

as $|x| \to \infty$ in $\mathbb{R}^2$ with $\hat{x} = x/|x|$ where $u^{\infty}$ is the far–field pattern. Since the ESM can recover the location of the scatterer with only one incident direction, this alleviates a main drawback of qualitative reconstruction methods. Since the introduction of the ESM, this imaging technique has been applied to other scattering problems, see for e.g., [2, 5]. The ESM is implemented by (numerically) asking when the measured far–field pattern is in the range of the far–field operator for a sound–soft sampling disk $B_z$ with fixed radius $R > 0$ centered at $z$. By sweeping through sampling points $z \in \mathbb{R}^2$ this allows one to recover the location of the scatterer. However, the standard analysis uses the ideas from the linear sampling method whereas we present an ESM based on the more analytically rigorous factorization method. We also consider using both sound–soft (Dirichlet) and sound–hard (Neumann) sampling disks $B_z$ in our new implementation.

## 2 Extended Sampling Revisited

In this section, we will revisit the ESM for the case of an acoustic scatterer in $\mathbb{R}^2$ for simplicity. We begin by providing a result that we can use to derive a new ESM. The following result is essentially the linchpin for the linear sampling and factorization methods which have been more widely studied [1, 3]. Now, assume that there are two open bounded regions $D_1$ and $D_2$ subsets of $\mathbb{R}^2$. Let the radiating scattered fields $u_{D_1}$ and $u_{D_2}$ satisfy

$$\Delta u_{D_j} + k^2 u_{D_j} = 0 \text{ in } \mathbb{R}^2 \setminus \overline{D}_j$$

and

$$\lim_{r\to\infty} \sqrt{r}(\partial_r u_{D_j} - \mathrm{i}k u_{D_j}) = 0$$

for $j = 1, 2$. We have that, for a fixed wavenumber $k > 0$, if the far–field patterns are non–trivial

$$u^{\infty}_{D_1}(\hat{x}) = u^{\infty}_{D_2}(\hat{x}) \text{ for all } \hat{x} \in \mathbb{S} \implies D_1 \cap D_2 \neq \emptyset$$

where $\mathbb{S}$ denotes the unit circle. Note that this result was proven in [5] for electromagnetic scatterers. The proof is similar for our case associated acoustic scattering.

Now, assume that we have the measured far–field pattern $u^{\infty}_{D}(\cdot\,; d)$ which is produced by a time–harmonic incident plane wave, denoted by

$$u^{\text{inc}}(x) = \mathrm{e}^{\mathrm{i}kx\cdot d}$$

for some fixed $d \in \mathbb{S}$. To fix a model problem for the purpose of exposition, we let the unknown scatterer $D$ be an anisotropic region such that the associated scattered field $u^{\text{scat}}_{D}$ satisfies

$$\nabla \cdot A \nabla u^{\text{scat}}_{D} + k^2 n u^{\text{scat}}_{D} = \nabla \cdot (I - A)\nabla u^{\text{inc}} + k^2(1-n)u^{\text{inc}} \quad \text{in} \quad \mathbb{R}^2$$

along with the radiation condition

$$\lim_{r\to\infty} \sqrt{r}(\partial_r u^{\text{scat}}_{D} - \mathrm{i}k u^{\text{scat}}_{D}) = 0.$$

The coefficient matrix $A \in C^1(D, \mathbb{R}^{2\times 2})$ is symmetric uniformly positive definite in $D$ and the scalar function $n \in C(D)$, denotes the material parameters of the $D$ where we assume that $I - A$ and $1 - n$ have support equal to $D$.

From the direct scattering problem for $u^{\text{scat}}_{D}$, we have that it is a radiating solution to the Helmholtz equation $\mathbb{R}^2 \setminus \overline{D}$. Now, the idea presented in [4] is to use the far–field equation

$$(F^{\text{Dir}}_{B_z} g)(\hat{x}) = u^{\infty}_{D}(\hat{x}, d), \quad \hat{x} \in \mathbb{S}.$$

Here, we define far–field operator

$$(F^{\text{Dir}}_{B_z} g)(\hat{x}) = \int_{\mathbb{S}} \mathrm{e}^{-\mathrm{i}kz\cdot(\hat{x}-\hat{y})} U^{\infty}_{\text{Dir}}(\hat{x}, \hat{y}) g(\hat{y}) \mathrm{d}s(\hat{y})$$

where $U^{\infty}_{\text{Dir}}$ is the far–field pattern for the sound–soft (Dirichlet) obstacle $B$ i.e., the disk centered at the origin with radius $R > 0$. By using separation of variables we have that

$$U^{\infty}_{\text{Dir}}(\hat{x}, \hat{y}) = -\frac{4}{\mathrm{i}}\Big[J_0(kR)/H^{(1)}_0(kR) + 2\sum_{p=1}^{\infty} J_p(kR)/H^{(1)}_p(kR)\cos\big(p(\theta_{\hat{x}} - \theta_{\hat{y}})\big)\Big].$$

By solving the far–field equation for a sampling disk $B_z$ with radius $R > 0.5\,\text{diam}(D)$, we have that the minimal value of the solution norm $\|g_z\|_{L^2}$ happen when $D \subset B_z$. This allows one to approximate the location of the scatterer.

Here, we propose an alternative ESM where the analysis is based on the factorization method. Indeed, from the analysis in Chapter 1 of [3] we have that, if $-k^2$ is not a Dirichlet eigenvalue of the Laplacian in $B_z$ then

$$\text{Range}\left(\left(F^{\text{Dir}*}_{B_z} F^{\text{Dir}}_{B_z}\right)^{1/4}\right) = \text{Range}\left(G^{\text{Dir}}_{B_z}\right).$$

The operator $G^{\text{Dir}}_{B_z}$ is defined by

$$G^{\text{Dir}}_{B_z} f = u^{\infty}_{B_z}(\hat{x})$$

where

$$\Delta u_{B_z} + k^2 u_{B_z} = 0 \text{ in } \mathbb{R}^2 \setminus \overline{B}_z$$

with

$$u_{B_z}\big|_{\partial B_z} = f \quad \text{and} \quad \lim_{r\to\infty} \sqrt{r}(\partial_r u_{B_z} - \mathrm{i}k u_{B_z}) = 0.$$

From this, we see that

$$u^{\infty}_{D} \in \text{Range}\left(\left(F^{\text{Dir}*}_{B_z} F^{\text{Dir}}_{B_z}\right)^{1/4}\right) \implies \overline{B}_z \cap \overline{D} \neq \emptyset$$

since this would imply that the $u^{\infty}_{D} = u^{\infty}_{B_z}$. This can be tested using Picard's Criteria since $F^{\text{Dir}}_{B_z}$ is injective, compact and normal. Therefore, if we define

$$W_{\text{Dir}}(z; d) = \left[\sum_{j=1}^{\infty} \frac{\left|\left(u^{\infty}_{D}(\cdot\,, d), \varphi^{z}_{j}\right)_{L^2(\mathbb{S})}\right|^2}{|\lambda^{z}_{j}|}\right]^{-1}$$

where $\varphi^{z}_{j}$ and $\lambda^{z}_{j}$ are the eigenvalues and vectors of $F^{\text{Dir}}_{B_z}$. With this, Picard's Criteria gives that

$$W_{\text{Dir}}(z; d) > 0 \implies \overline{B}_z \cap \overline{D} = \emptyset.$$

Thus, one can recover the scatterer (or at least its location) by plotting the imaging functional $W_{\text{Dir}}(z; d)$. Notice, that if you have data from a finite number of incident directions $d_l$ for $l = 1, \cdots, N_d$ then the imaging functional can be given as

$$z \longmapsto \sum_{l=1}^{N_d} W_{\text{Dir}}(z; d_l).$$

The support of the imaging functional should coincide with the unknown scatterer $D$.

From this analysis, we notice that sound–soft (Dirichlet) obstacle assumption for the far–field operator is not necessary. Therefore, we also propose using a sound–hard (Neumann) obstacle assumption for the far–field operator. To this end, just as before we can define

$$(F_{B_z}^{\text{Neu}} g)(\hat{x}) = \int_{\mathbb{S}} \mathrm{e}^{-\mathrm{i}kz\cdot(\hat{x}-\hat{y})} U_{\text{Neu}}^{\infty}(\hat{x},\hat{y}) g(\hat{y}) \mathrm{d}s(\hat{y})$$

where $U_{\text{Neu}}^{\infty}$ is the far–field pattern for the sound hard (Neumann) obstacle $B$ i.e., the disk centered at the origin with radius $R > 0$. By using separation of variables we have that

$$U_{\text{Neu}}^{\infty}(\hat{x},\hat{y}) = -\frac{4}{\mathrm{i}}\Big[J_0'(kR)/H_0^{(1)'}(kR) + 2\sum_{p=1}^{\infty} J_p'(kR)/H_p^{(1)'}(kR)\cos\big(p(\theta_{\hat{x}}-\theta_{\hat{y}})\big)\Big].$$

Now, we can define another imaging functional using the factorization method analysis in Chapter 1 of [3]. Just as in the previous case, we define

$$W_{\text{Neu}}(z;d) = \left[\sum_{j=1}^{\infty} \frac{\left|\big(u_D^{\infty}(\cdot\,,d),\psi_j^z\big)_{L^2(\mathbb{S})}\right|^2}{|\mu_j^z|}\right]^{-1}$$

where $\psi_j^z$ and $\mu_j^z$ are the eigenvalues and vectors of $F_{B_z}^{\text{Neu}}$ provided $-k^2$ is not a Neumann eigenvalue of the Laplacian in $B_z$. Again, by appealing to Picard's Criteria we have that

$$W_{\text{Neu}}(z;d) > 0 \implies \overline{B}_z \cap \overline{D} = \emptyset$$

and for a finite number of incident directions $d_l$ for $l = 1,\cdots,N_d$ this imaging functional can be given as

$$z \longmapsto \sum_{l=1}^{N_d} W_{\text{Neu}}(z;d_l).$$

Here as before, the support of the imaging functional should coincide in some sense with the unknown scatterer $D$.

From this, we now have two new imaging functionals for recovering an unknown scatterer. Here, we have presented this for the case when the scatterer is an anisotropic medium. Notice that the analysis here would work for other scattering regions. This shows that this new variant of the ESM can be applied to many types of problems. One question for future work is whether this can be replicated for other scattering problems, e.g., electromagnetic, elastic, and flexural wave scattering.

## 3 Numerical Examples

In this section, we provide some numerical examples to illustrate the effectiveness of the new ESM imaging functionals derived in the previous section. Here, we will use the synthetic far–field data from an anisotropic region denoted $D$. In all of the examples given we let

$$A = \begin{bmatrix} 4 & 1 \\ 1 & 4 \end{bmatrix} \quad \text{and} \quad n = 0.5$$

be the material parameters for the scatterer with $k = 4$ being fixed. The boundary of the scatterer is given by

$$\partial D = 0.5\big(1 + 0.25\sin(3\theta)\big)\big(\cos\theta, \sin\theta\big) - (1,1)$$

for $\theta \in [0, 2\pi)$. This region is a pear–shaped scatterer that we wish to recover.

In order to evaluate the imaging functionals, we need to compute the synthetic far–field data associated with the pear–shaped scatterers. By appealing to the Lippmann–Schwinger integral equation for the scattered field the Born approximation is given by

$$u_D^{\infty}(\hat{x},d) \approx k^2 \int_D \hat{x}\cdot(I-A)d\,\mathrm{e}^{-\mathrm{i}k\omega\cdot(\hat{x}-d)}\,\mathrm{d}\omega + k^2\int_D (n-1)\mathrm{e}^{-\mathrm{i}k\omega\cdot(\hat{x}-d)}\,\mathrm{d}\omega.$$

Next, we need an approximation of the far–field operators $F_{B_z}^{\text{Dir}}$ and $F_{B_z}^{\text{Neu}}$ where we fix $R = 0.75$. This is done by truncating the series used to define the kernel of the integral operators and using a standard Riemann sum approximation for the integrals. Since in many applications one only has access to noisy data, in the following reconstructions we add 2% or 5% random noise to the synthetic far–field data.

We now consider the case when one only has the far–field pattern $u_D^{\infty}(\hat{x},d)$ produced from a single incident direction $d$. For this set of examples, we take the incident direction given by

$$d = \big(\cos\theta, \sin\theta\big) \text{ for } \theta = \frac{3\pi}{2}.$$

To stabilize the reconstruction in the presence of noisy data, the imaging functionals are computed via a spectral cut–off such that we stop the sum when eigenvalues are less than $10^{-3}$.

From Figures 1 and 2, we see that even with only one incident direction that the imaging functionals can recover the scatterer. Now we give the reconstruction using four incident directions. With this, we use the computed far–field data with incident directions given by

$$d_l = \big(\cos\theta_l, \sin\theta_l\big) \text{ for } \theta_l = 0, \frac{\pi}{2}, \pi, \text{ and } \frac{3\pi}{2}.$$

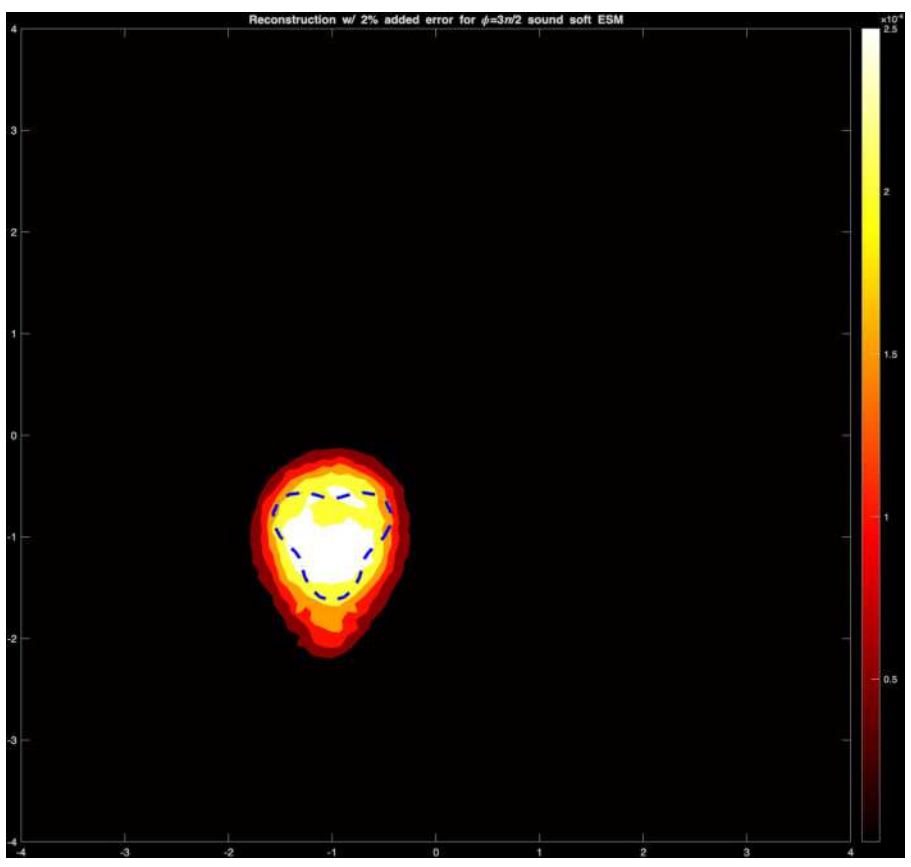

Figure 1: Reconstruction of the scatterer using the sound–soft $F^{\mathrm{Dir}}_{B_z}$ far–field operator with four incident directions with 2% added noise. The dotted lines are the boundary of the scatterer.

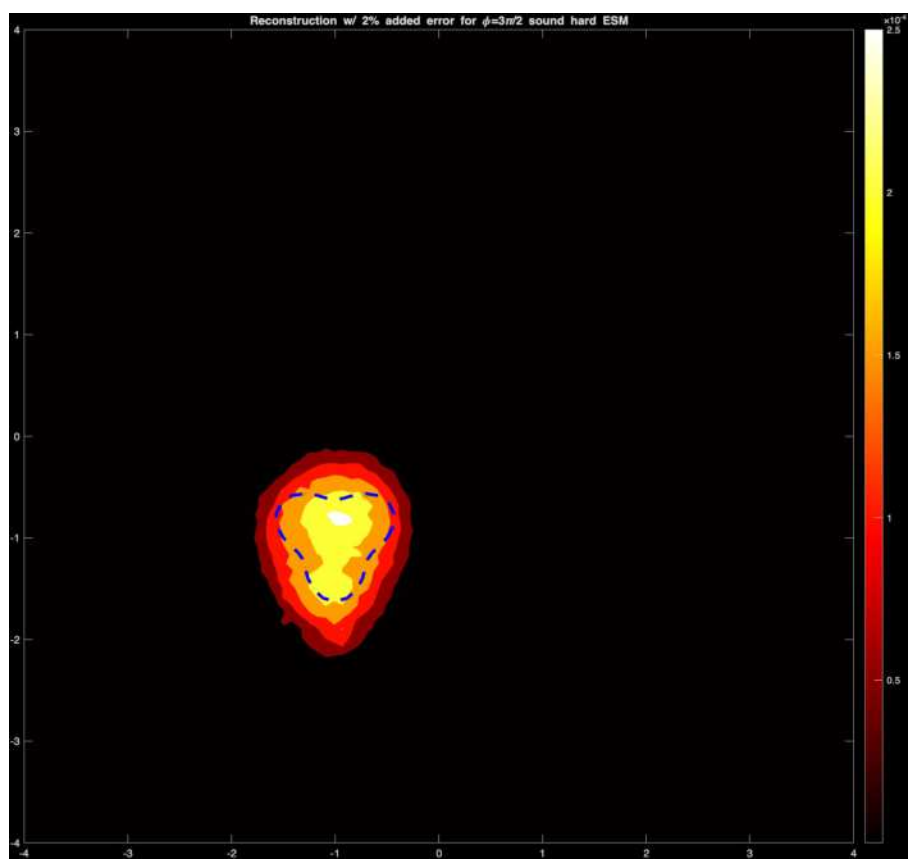

Figure 2: Reconstruction of the scatterer using the sound–hard $F^{\mathrm{Neu}}_{B_z}$ far–field operator with only one incident direction with 2% added noise. The dotted lines are the boundary of the scatterer.

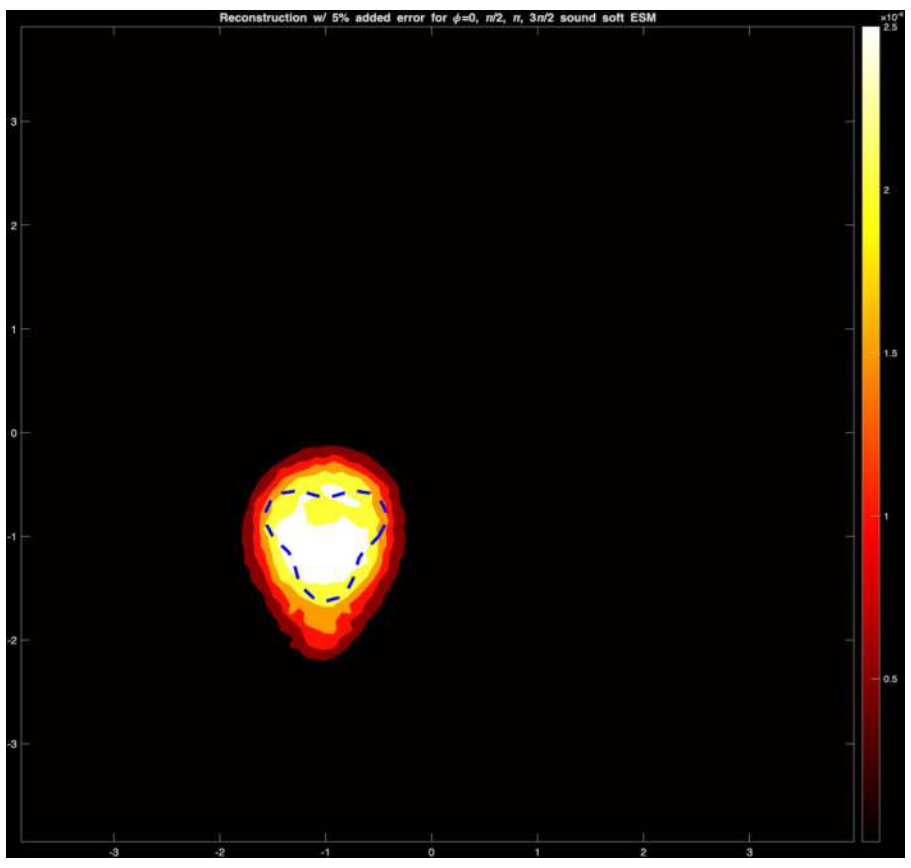

Figure 3: Reconstruction of the scatterer using the sound–soft $F^{\mathrm{Dir}}_{B_z}$ far–field operator with only one incident direction with 5% added noise. The dotted lines are the boundary of the scatterer.

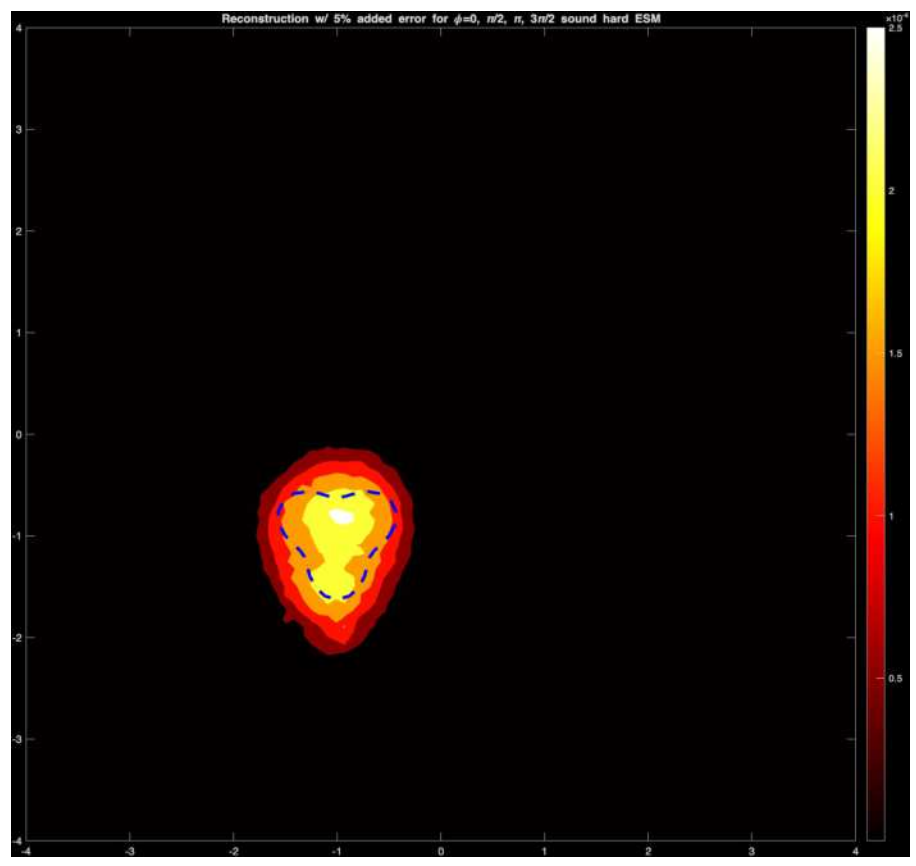

Figure 4: Reconstruction of the scatterer using the sound–hard $F^{\mathrm{Neu}}_{B_z}$ far–field operator with four incident directions with 5% added noise. The dotted lines are the boundary of the scatterer.

Again, in Figures 3 and 4 we see that the scatterer is recovered in this case as well. Note that in our numerical experiments, this method works best for $R \approx 0.5\,\mathrm{diam}(D)$.

# Trefftz methods with evanescent plane waves

Andrea Moiola[1,*], Nicola Galante[2], Emile Parolin[2]
[1]Department of Mathematics "F. Casorati", University of Pavia, Pavia, Italy
[2]Sorbonne Université, Université Paris Cité, CNRS, INRIA, Laboratoire Jacques-Louis Lions, LJLL, EPC ALPINES, Paris, F-75005, France
*Email: andrea.moiola@unipv.it

**Abstract**

Classical Trefftz methods approximate Helmholtz solutions using propagative plane waves and are subject to strong numerical instabilities. Evanescent plane wave bases can substantially mitigate this phenomenon. We propose a simple recipe to select such basis functions. We show that the numerical results obtained by the Ultraweak Variational Formulation (UWVF) greatly improve thanks to this choice. More details and examples will soon be available in [6].

## 1 Trefftz methods and plane waves

Exactly 100 years ago, at the International Congress of Applied Mechanics that took place in Zürich in 1926, Erich Trefftz proposed the approximate solution of a Laplace boundary value problem (BVP) by piecewise-harmonic functions [12]. Decades later, his work gave the name to a wide class of numerical methods applicable to any linear PDE: **Trefftz methods** approximate BVP solutions with discrete spaces made of piecewise solutions of the underlying PDE. Trefftz methods often take the form of discontinuous Galerkin schemes (DG).

Trefftz methods have been extensively used to approximate the **Helmholtz equation** [8]

$$\Delta u + \kappa^2 u = 0 \tag{1}$$

on a subdomain of $\mathbb{R}^n$, where $\kappa > 0$ is the wavenumber. Classical finite element and DG methods approximate solutions of (1) by piecewise polynomials. Since there is no polynomial solution of (1) other than the constant zero, Trefftz methods for this PDE resort to non-polynomial bases. Most Trefftz schemes (e.g. [3, 7, 9]) for (1) use **propagative plane wave** (PPW) basis functions:

$$\mathbf{x} \mapsto e^{\imath\kappa\mathbf{d}\cdot\mathbf{x}} \quad \text{with} \quad \mathbf{d} \in \mathbb{S}^{n-1}, \tag{2}$$

where $\mathbb{S}^{n-1} := \{\mathbf{v} \in \mathbb{R}^n : |\mathbf{v}| = 1\}$. Discrete spaces spanned by PPWs with different propagation directions $\mathbf{d}$ approximate general Helmholtz solutions with fast convergence rates [10]. However, it is widely acknowledged that PPW-based Trefftz implementations suffer from strong numerical instabilities.

## 2 A negative PPW approximation result

The instabilities of PPW-based Trefftz schemes have been described and addressed in several ways, mostly in terms of near-dependence of the PPWs, ill-conditioning of the linear systems, and linear algebra techniques (e.g. [9, §4.2], [2], [4], [8, §4.3]). We follow a different approach, inspired by [1], and report the result of [11, Lemma 4.2] (2D) and [5, Lemma 4.4] (3D).

**Theorem 1** *Let $B$ be the unit ball in $\mathbb{R}^n$, $n \in \{2, 3\}$, and $\kappa > 0$. Take any*

- *target relative accuracy $0 < \delta < 1$,*
- *large threshold value $M > 1$.*

*Then, there exists $u \in C^\infty(\mathbb{R}^n)$, solution of (1) on $\mathbb{R}^n$, normalised as $\|u\|_{H^1(B)} = 1$, such that, if $u$ is approximable by a linear combination of a finite number of PPWs with accuracy $\delta$, i.e.*

$$\left\| u - \sum_{p=1}^{P} \mu_p e^{\imath\kappa\mathbf{d}_p\cdot\mathbf{x}} \right\|_{H^1(B)} \le \delta \qquad \begin{array}{c}\text{for some } P \in \mathbb{N},\\ \boldsymbol{\mu} \in \mathbb{C}^P,\\ \{\mathbf{d}_p\}_{p=1}^P \subset \mathbb{S}^{n-1},\end{array}$$

*then the Euclidean norm of the coefficient vector is larger than $M$: $\|\boldsymbol{\mu}\|_{\mathbb{C}^P} \ge M$.*

The function $u$ is simply a high-index circular or spherical wave. The $H^1(B)$ norm in the statement could be replaced by any finite-order Sobolev norm.

Theorem 1 shows that, for some target Helmholtz solutions, **accurate PPW approximations necessarily require huge coefficients**. If the function we seek to approximate is the $u$ provided by the theorem for large $M$ (e.g.

$M$ larger than the reciprocal of machine precision), then numerical cancellation is inevitable and there is no way to stably represent its PPW approximation in computer arithmetic.

This is the fundamental issue underlying all PPW-Trefftz instabilities. Linear algebra techniques (such as preconditioning or matrix decompositions) can not help finding an approximation if an approximation that is representable in computer arithmetic does not exist. Recall also that [1] shows that, roughly speaking, stable numerical approximations are possible if and only if small-coefficient approximation exists.

Theorem 1 states that accurate approximations and small-coefficient representations are mutually exclusive for *some* Helmholtz solutions (the "evanescent modes"). Other Helmholtz solutions can be stably approximated by PPWs, as witnessed by the numerous successful uses of Trefftz schemes [8].

## 3 Evanescent plane waves

To deal with the PPW approximation instability, following [5,11] we propose the use of evanescent plane waves (EPW):

$$\mathbf{x} \mapsto e^{\imath\kappa\mathbf{d}\cdot\mathbf{x}}, \quad \mathbf{d} \in \mathbb{C}^n, \quad \mathbf{d}\cdot\mathbf{d} = \sum_{j=1}^{n} d_j^2 = 1. \quad (3)$$

EPWs allow complex propagation vectors $\mathbf{d}$ and include PPWs as special case ($\mathbf{d} \in \mathbb{R}^n \subset \mathbb{C}^n$). The condition $\mathbf{d}\cdot\mathbf{d} = 1$ is equivalent to

$$|\Re(\mathbf{d})|^2 - |\Im(\mathbf{d})|^2 = 1, \quad \Re(\mathbf{d})\cdot\Im(\mathbf{d}) = 0. \quad (4)$$

An EPW is determined by three parameters:

$$\text{propagation direction} \quad \mathbf{p} := \frac{\Re(\mathbf{d})}{|\Re(\mathbf{d})|} \in \mathbb{S}^{n-1},$$

$$\text{evanescence direction} \quad \mathbf{e} := \frac{\Im(\mathbf{d})}{|\Im(\mathbf{d})|} \in \mathbb{S}^{n-1},$$

$$\text{evanescence strength} \quad \eta := |\Im(\mathbf{d})| \in [0, +\infty).$$

Note that $\mathbf{e}$ has to be orthogonal to $\mathbf{p}$ by (4), and that if $\eta = 0$ then the EPW is a PPW and $\mathbf{e}$ is undefined. We also introduce a second measure of the evanescence strength:

$$\zeta := |\Re(\mathbf{d})| = \sqrt{1+\eta^2} \in [1, +\infty).$$

Then the evanescent plane wave is

$$\mathrm{EW}_{[\mathbf{p},\mathbf{e},\eta]}(\mathbf{x}) := e^{\imath\kappa\mathbf{d}\cdot\mathbf{x}} = e^{\imath\kappa\zeta\mathbf{p}\cdot\mathbf{x}} e^{-\kappa\eta\mathbf{e}\cdot\mathbf{x}}. \quad (5)$$

This EPW oscillates and propagates in direction $\mathbf{p}$ with apparent wavenumber $\kappa\zeta$, and decays with exponential rate $\kappa\eta$ in direction $\mathbf{e}$; see Figure 1.

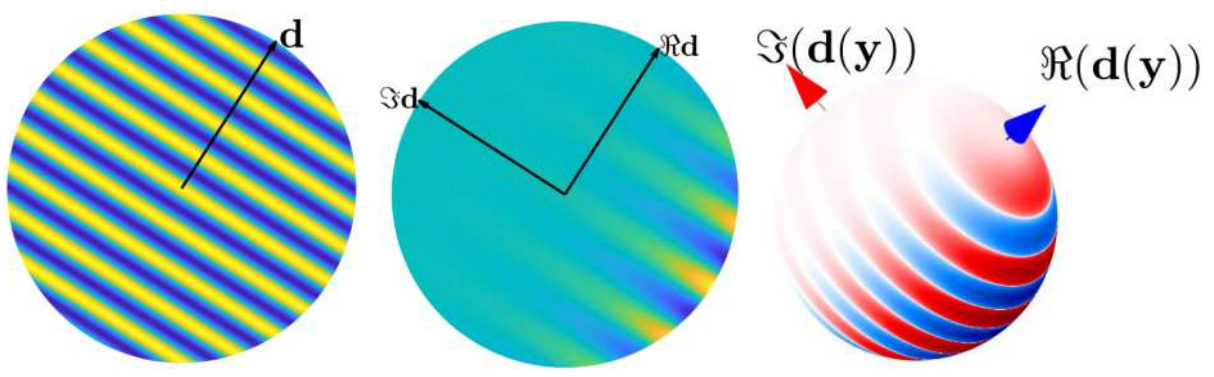

Figure 1: A PPW and an EPW on the disc in $\mathbb{R}^2$, an EPW on the ball in $\mathbb{R}^3$.

Why do we expect linear combinations of the EPWs (3) to give accurate small-coefficient approximations? The reason is in [11, Thm. 6.7] (2D) and [5, Thm. 3.9] (3D): any Helmholtz solution $u \in H^1(B)$, with $B$ the unit disc/ball, is a continuous EPW superposition. The coefficient of this superposition, measured in a weighted $L^2$ norm, is bounded by the $H^1(B)$ norm of $u$. No such a representation can hold for PPWs. This "Herglotz transform" is a step towards a convergence and stability theory for EPW-based schemes; further work is needed to extend it to more general domains and discrete EPW spaces.

## 4 EPW selection

Implementing a PPW-based Trefftz method, one typically selects a $P$-dimensional discrete space by choosing PPW directions $\{\mathbf{d}_p\}_{p=1}^{P}$ that are roughly equispaced on the unit sphere $\mathbb{S}^{n-1}$.

When implementing an EPW-based method, how should one choose the EPW parameters $\mathbf{p}, \mathbf{e}, \eta$? The propagation directions $\mathbf{p}$ can be taken uniformly on $\mathbb{S}^{n-1}$ and the evanescence directions $\mathbf{e}$ uniformly in $\{\mathbf{e} \in \mathbb{S}^{n-1} : \mathbf{e}\cdot\mathbf{p} = 0\} \sim \mathbb{S}^{n-2}$. The choice of the evanescence strengths $\eta$ of the basis functions is less obvious but is key to ensure the quality of the method.

We describe in detail a simple concrete strategy to choose a discrete EPW space in 2D ($n = 2$). We refer to [6] for the 3D case.

In 2D, the propagation direction $\mathbf{p}$ is parametrised by the polar angle $\theta \in [0, 2\pi)$ and the decay direction $\mathbf{e}$ by $\varphi \in \{\pm 1\}$, as $\mathbf{e}$ can only be one of the two unit vectors orthogonal to $\mathbf{p}$. So (5) reads

$$\mathrm{EW}_{[\mathbf{p},\mathbf{e},\eta]}(\mathbf{x}) = e^{\imath\kappa\zeta\binom{\cos\theta}{\sin\theta}\cdot\mathbf{x}} e^{-\kappa\eta\varphi\binom{-\sin\theta}{\cos\theta}\cdot\mathbf{x}}.$$

### 4.1 EPW sampling recipe in 2D

Assume we are given a "DOF budget" $P \in \mathbb{N}$ and we want to construct a $P$-dimensional EPW Trefftz space on a polygonal domain $K \subset \mathbb{R}^2$ (which in §5 will represent a mesh element). We define a "Fourier truncation level"

$$L := \frac{P}{4},$$

borrowing the name from [11, §7.3], [5, p. 462]. We sample $P$ points $\{\mathbf{y}_p = (\theta_p, \varphi_p, \xi_p)\}_{p=1}^P$ from the uniform probability measure on

$$Y := [0, 2\pi) \times \{\pm 1\} \times [0, 1].$$

This probability is the Lebesgue measure scaled by $\frac{1}{2\pi}$ in the first component $\theta_p$, equal probability $\frac{1}{2}$ for $\varphi_p = 1$ and $\varphi_p = -1$, and the Lebesgue measure in the last component $\xi_p$. We define the sampled evanescence parameters

$$\zeta_p := \max\left\{1, \frac{2L}{\kappa\,\mathrm{diam}(K)}\xi_p\right\} \text{ and } \eta_p^2 := \zeta_p^2 - 1.$$

The discrete Trefftz space is

$$\mathbb{T}_P(K) := \mathrm{span}\left\{\frac{\mathrm{EW}_{\mathbf{y}_p}}{\|\mathrm{EW}_{\mathbf{y}_p}\|_{L^\infty(K)}}\right\}_{p=1}^P, \text{ with}$$
$$\mathrm{EW}_{\mathbf{y}_p}(\mathbf{x}) := e^{\imath\kappa\zeta_p \binom{\cos\theta_p}{\sin\theta_p}\cdot\mathbf{x}} e^{-\kappa\eta_p\varphi_p\binom{-\sin\theta_p}{\cos\theta_p}\cdot\mathbf{x}}. \quad (6)$$

The norm $\|\mathrm{EW}_{\mathbf{y}_p}\|_{L^\infty(K)}$ can be computed by evaluating the EPW at each vertex of the cell $K$.

### 4.2 Comments on the sampling recipe

The choices made in §4.1 might seem arbitrary: they are motivated in [6] from an asymptotic approximation of the density appearing in the Herglotz integral representation [11] (see also [5, eq. (5.8)]) mentioned in §3.

The sampling of the set $Y$ can be done in several different ways: random, quasi-random (e.g. with Sobol sequences, as in §6), deterministic (e.g. equispaced in $\theta_p$ and/or $\xi_p$), with or without a tensor-product structure.

If $2L \le \kappa\,\mathrm{diam}(K)$, then all $\zeta_p = 1$ and all basis functions are PPWs. In this case the DOF budget $P$ is low, and all the approximation effort is devoted to the approximation of the propagative components of the solution. For larger values of $P$, an increasing fraction $\frac{2L}{\kappa\,\mathrm{diam}(K)}$ of the $P$ DOFs will be EPWs.

In 3D, we propose a similar recipe with $L = \sqrt{\frac{P}{8}}$ and $\zeta_p = \max\{1, \frac{2L}{\kappa\,\mathrm{diam}(K)}\sqrt{\xi_p}\}$. In this case the uniform sampling of the orthogonal directions $\mathbf{p}$ and $\mathbf{e}$ is more delicate and we refer to [6] for the details.

## 5 Ultraweak variational formulation

The EPW discrete space (6) can be used within any of the Trefftz variational formulations in [8, §2] (LS, TDG, UWVF, VTCR, WBM,...). Here we focus on the UWVF [3,9], which can be recast as a special instance of the TDG [6–8].

Consider the Helmholtz impedance BVP

$$\begin{aligned} \Delta u + \kappa^2 u = 0 \quad &\text{on } \Omega \subset \mathbb{R}^n, \\ \partial_{\mathbf{n}} u - \imath\kappa\sigma u = g \quad &\text{on } \partial\Omega, \end{aligned} \qquad (7)$$

with $0 < \sigma \in L^\infty(\partial\Omega)$ and $g \in L^2(\partial\Omega)$. Let $\Omega$ be Lipschitz, bounded, partitioned in a mesh $\mathcal{T} = \{K\}$, and let $\mathcal{F} := (\bigcup_{K\in\mathcal{T}} \partial K) \setminus \partial\Omega$ be the internal mesh skeleton. Define the Trefftz space

$$T(\mathcal{T}) := \prod_{K\in\mathcal{T}} T(K),$$
$$T(K) := \big\{v \in H^1(K) : \Delta v + \kappa^2 v = 0, \ \partial_{\mathbf{n}_K} v \in L^2(\partial K)\big\}.$$

Extend $\sigma$ to a positive $L^\infty(\partial\Omega \cup \mathcal{F})$ and define two Robin traces:

$$\gamma_\pm^K : T(K) \to L^2(\partial K), \ \gamma_\pm^K v := \pm\partial_{\mathbf{n}_K} v - \imath\kappa\sigma v.$$

For a finite dimensional subspace $V_h$ of $T(\mathcal{T})$, the classical UWVF [3,9] of BVP (7) consists in finding $u_h \in V_h$ such that

$$\begin{aligned} &\sum_{K\in\mathcal{T}} \int_{\partial K} \sigma^{-1}\, \gamma_-^K u_h\, \overline{\gamma_-^K v_h} \\ &- \sum_{K_1\in\mathcal{T}} \sum_{K_2\in\mathcal{T}} \int_{\partial K_1\cap\partial K_2} \sigma^{-1}\, \gamma_-^{K_1} u_h\, \overline{\gamma_+^{K_2} v_h} \\ &= \sum_{K\in\mathcal{T}} \int_{\partial K\cap\partial\Omega} \sigma^{-1}\, g\, \overline{\gamma_+^K v_h} \qquad \forall v_h \in V_h. \end{aligned} \qquad (8)$$

See [6] for the extension to piecewise-constant parameters and general mixed (Dirichlet, Neumann, impedance) boundary conditions. The UWVF (8) is consistent with the BVP (7), well-posed, coercive and quasi-optimal in skeleton norms [6–8].

The linear system of (8) can be written as $(\mathbf{D} - \mathbf{C})\mathbf{u} = \mathbf{b}$, where $\mathbf{D}$ and $\mathbf{C}$ are the matrices corresponding to the two integrals at the left-hand side. The matrix $\mathbf{D}$ is Hermitian, positive-definite, and block-diagonal with each block corresponding to a mesh element, so the system is equivalent to $(\mathbf{I}-\mathbf{D}^{-1}\mathbf{C})\mathbf{u}=\mathbf{D}^{-1}\mathbf{b}$ as in [3, eq. (2.32)]. The inversion of $\mathbf{D}$ is a local (elementwise) operation.

If the mesh facets are flat and EPW or PPW bases are used, then the entries of the matrices $\mathbf{D},\mathbf{C}$ can be computed analytically.

**Oversampled regularised UWVF** We use the UWVF (8) with the sampled EPW discrete space from (6):

$$V_h=\left\{v\in L^2(\Omega):\ v|_K\in\mathbb{T}_P(K)\ \ \forall K\in\mathcal{T}\right\}.$$

The sampling of §4.1 can be done independently for each element $K$ or, if the element sizes are comparable, once for all elements.

However, following [1], to ensure a stable solution even with basis functions that are close to linearly dependent, we need to use **oversampling** (10% oversampling is used in §6) and **regularisation**. The regularisation relies on the truncated SVD of the $\mathbf{D}$ matrix, which can be done in parallel (elementwise) as in [2].

We proceed as follows. For each $K\in\mathcal{T}$, let $\mathbb{T}_{N_K^{\text{trial}}}(K)$ and $\mathbb{T}_{N_K^{\text{test}}}(K)$ be two EPW spaces in the form (6), with $N_K^{\text{test}}\geq N_K^{\text{trial}}$, and define $N^\bullet:=\sum_{K\in\mathcal{T}}N_K^\bullet$ for $\bullet\in\{\text{trial,test}\}$. Denote $\mathbf{D}_K\in\mathbb{C}^{N_K^{\text{test}}\times N_K^{\text{trial}}}$ the tall matrix corresponding to the diagonal block associated to the element $K$ of the Galerkin matrix $\mathbf{D}\in\mathbb{C}^{N^{\text{test}}\times N^{\text{trial}}}$ of the first integral in (8). Let $\mathbf{D}_K=\mathbf{U}_K\mathbf{\Sigma}_K\mathbf{V}_K^*$ be its SVD and $\sigma_1\geq\sigma_2\geq\cdots\geq\sigma_{N_K^{\text{trial}}}\geq 0$ be the singular values. For a regularisation parameter $\epsilon\in(0,1)$ ($\epsilon=10^{-14}$ in §6), larger than the machine precision, let $\mathbf{\Sigma}_{K,\epsilon}^\dagger$ be the truncated pseudo-inverse of $\mathbf{\Sigma}$: the diagonal matrix with $q$-th diagonal entry $\sigma_q^{-1}$ if $\sigma_q\geq\epsilon\sigma_1$ and 0 otherwise. We then define $\mathbf{D}_\epsilon^\dagger\in\mathbb{C}^{N^{\text{trial}}\times N^{\text{test}}}$ the matrix with diagonal blocks $\mathbf{V}_K\mathbf{\Sigma}_{K,\epsilon}^\dagger\mathbf{U}_K^*$. Denoting by $\mathbf{C}\in\mathbb{C}^{N^{\text{test}}\times N^{\text{trial}}}$ the Galerkin matrix of the second integral in (8), the linear system $(\mathbf{I}-\mathbf{D}^{-1}\mathbf{C})\mathbf{u}=\mathbf{D}^{-1}\mathbf{b}$ is approximated by the $N^{\text{trial}}\times N^{\text{trial}}$ system

$$(\mathbf{I}-\mathbf{D}_\epsilon^\dagger\mathbf{C})\mathbf{u}_\epsilon=\mathbf{D}_\epsilon^\dagger\mathbf{b}. \tag{9}$$

The regularized linear system (9) is sparse and can be solved by any direct or iterative method.

## 6 Numerical experiments

We finally illustrate the differences in the approximation capabilities of PPWs (taking $L=0$ in the sampling recipe in §4.1) and EPWs.

**Point source** The domain is the square $\Omega:=(0,1)\times(-\frac{1}{2},\frac{1}{2})$, discretized into an unstructured quasi-uniform mesh with 41 triangles. We consider the numerical approximation of the fundamental solution $u(\mathbf{x})=\frac{\imath}{4}H_0^{(1)}(\kappa|\mathbf{x}-\mathbf{s}|)$ using the UWVF (8) with a suitable datum $g$. The wavenumber is $\kappa=128$ and the singularity $\mathbf{s}=(-\frac{\pi}{5\kappa},0)$ is located a tenth of a wavelength away from the left boundary.

Figure 2 reports the real part of the solution $u$, the convergence of the $\kappa$-weighted relative $H^1(\Omega)$-error with respect to the number of degrees of freedom, and the pointwise error in the domain (for the largest approximation space, $P=815$) when using either PPWs or EPWs. Note the different scales in the error colorbars.

We observe a first regime of fast convergence, obtained for both types of waves since the two approximation spaces are actually of similar nature (propagative) when the number of degrees of freedom is moderate. Then the PPW instability makes convergence stall and prevents from obtaining higher accuracy, whereas EPWs incur a much smaller error. The ratio between the EPW and the PPW errors reaches values below $10^{-6}$. The convergence is not monotone since the approximation spaces are not nested.

**Influence of the wavenumber** In the same configuration, we now study the influence of the wavenumber $\kappa$ when the size of the approximation space is proportional to $\kappa$ itself. More precisely, the number of PWs per cell is fixed to $P=4\kappa$. As before, the singularity is placed at $\mathbf{s}=(-\frac{\pi}{5\kappa},0)$, lying at a distance from the boundary proportional to the wavelength.

Figure 3 reports the convergence of the $\kappa$-weighted relative $H^1(\Omega)$-error with respect to the number of degrees of freedom and the wavenumber. The error achieved using EPWs improves with the wavenumber as the size of the approximation space increases linearly with $\kappa$.

**Scattering problem** We consider the scattering of an incident plane wave $u_i(\mathbf{x}):=e^{\imath\kappa\mathbf{d}\cdot\mathbf{x}}$ for $\mathbf{d}=(\frac{1}{2},-\frac{\sqrt{3}}{2})$ and wavenumber $\kappa=16$ by a sound-soft (i.e. Dirichlet) non-convex scatterer of diameter 2. The UWVF formulation in the total field features a non-homogeneous Robin boundary condition on a truncation boundary of diameter 4 surrounding the obstacle. The mesh contains 64 triangles. Since the solution to this problem is not known analytically, we use as reference solution the EPW approximation with twice as many waves as the largest approximation space. This test case is more challenging because of the presence of the cavity and the corner singularities at each vertex of the polygonal scatterer.

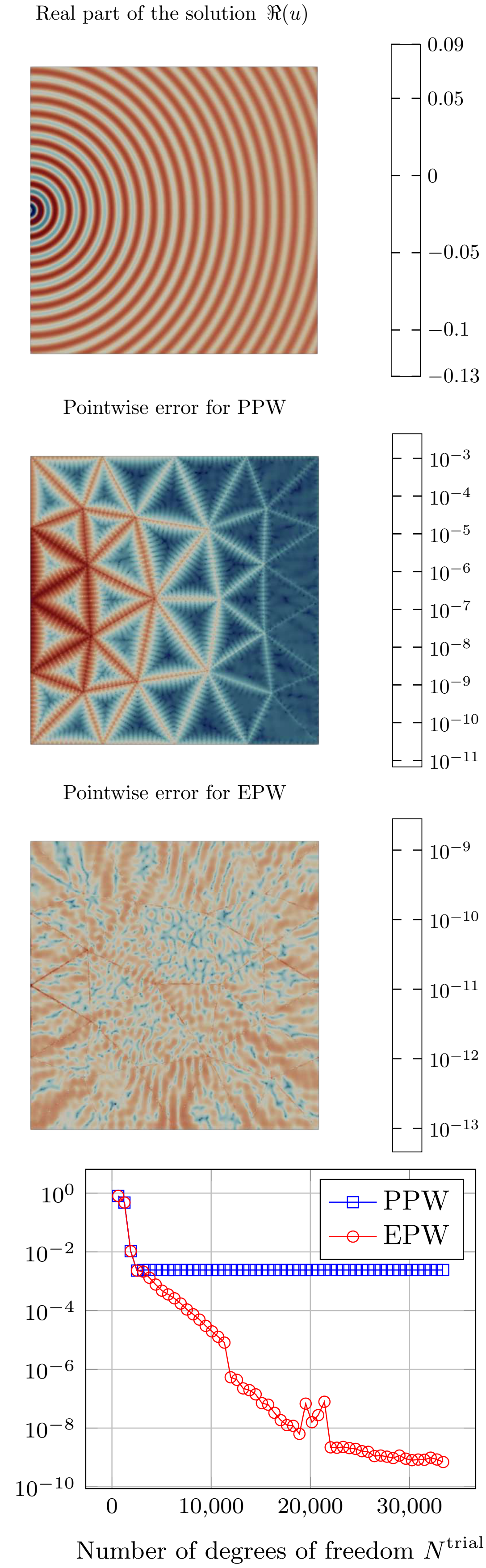

Figure 2: Fundamental solution: $\Re(u)$ (top), pointwise error (middle) for PPWs and EPWs, convergence history (bottom).

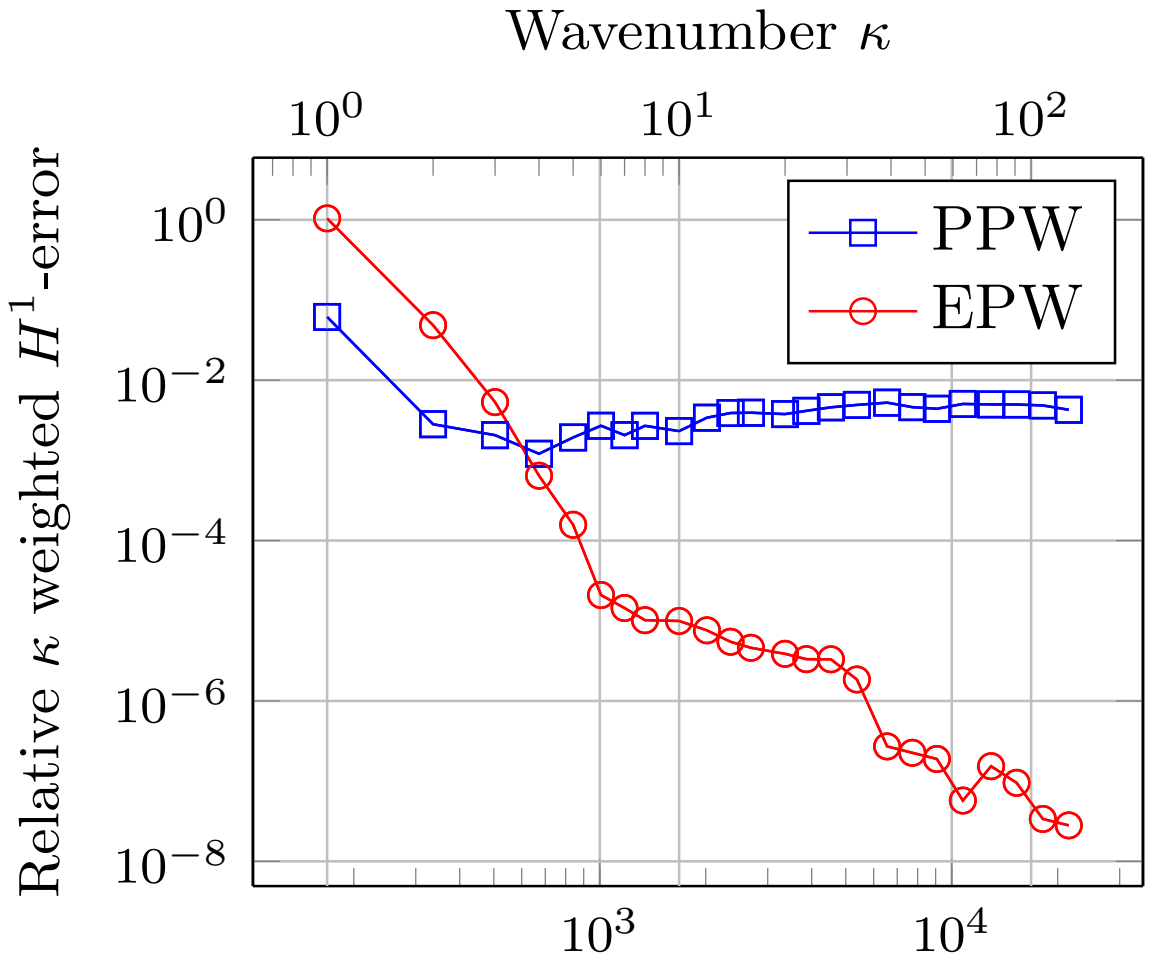

Figure 3: Fundamental solution: convergence result for increasing $\kappa$ using an approximation space of size increasing linearly in $\kappa$.

Figure 4 reports the real part of the solution, the convergence history of the $\kappa$-weighted relative $H^1(\Omega)$-error with respect to the number of degrees of freedom and the pointwise error in the domain (for the largest approximation space, $P = 2080$) when using either PPWs or EPWs.

Again we observe that the convergence of the UWVF using PPWs stagnates whereas EPWs yield much better accuracy.

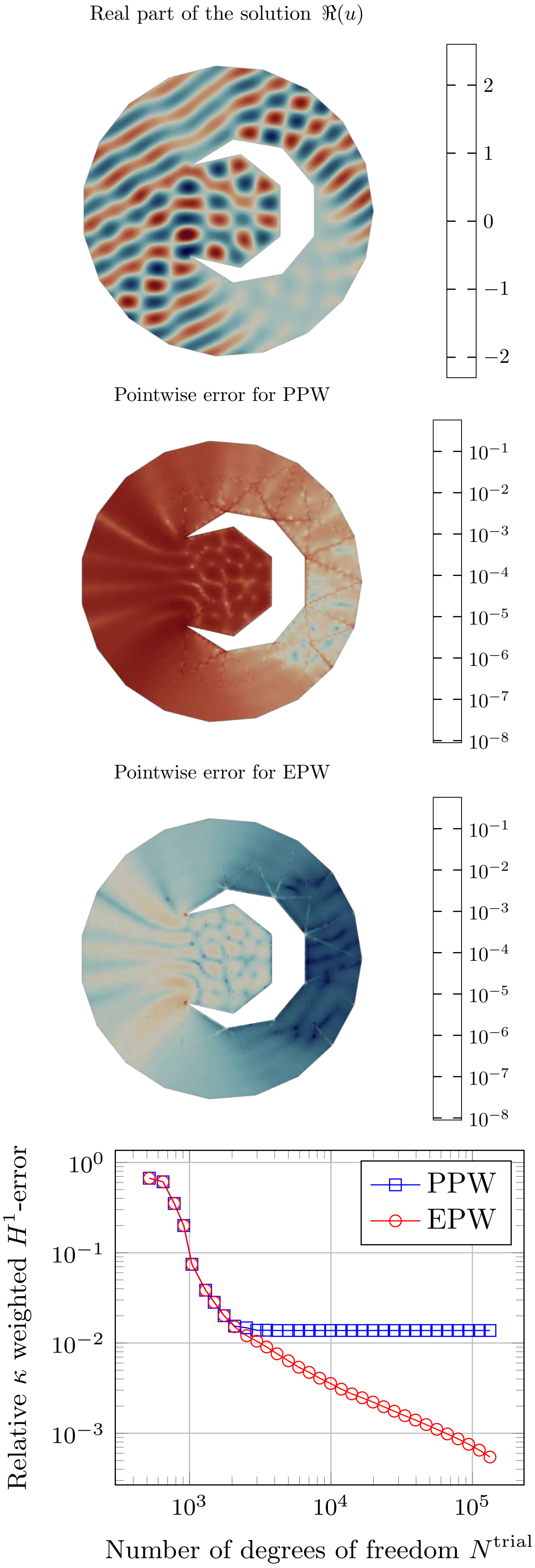

Figure 4: Scattering problem: $\Re(u)$ (top), pointwise error (middle) for PPWs and EPWs, convergence history (bottom).

# Moment limiter framework for the discontinuous Galerkin method

**Lilia Krivodonova**[1,*]
[1]Department of Applied Mathematics, University of Waterloo, Waterloo, Canada
*Email: lgk@uwaterloo.ca

## Abstract

Solutions of hyperbolic conservation laws develop discontinuities in finite time even with smooth initial data. High-order numerical methods have difficulty approximating such solutions because they develop oscillations near discontinuities of exact solutions, making computations unstable or unreliable. To control oscillations and to guarantee convergence of a numerical scheme, we apply limiters. The main challenge in designing limiters is to preserve high-order accuracy while stabilizing numerical solutions.

In this lecture, we present a framework for constructing moment limiters for the discontinuous Galerkin method. The limiter acts directly on solution coefficients by reconstructing slopes along directions in which the moments decouple. This leads to compact algorithms that are applicable across mesh types and computational settings, and that integrate naturally with adaptive and high-performance implementations.

## 1 Introduction

The development of effective limiters is central to the practical use of high-order methods for problems with discontinuities. While the discontinuous Galerkin method provides a natural framework for high-order approximation, its behavior near nonsmooth features remains a limiting factor in many applications.

A large body of work has been devoted to limiter design. However, the tension between stability and accuracy persists: methods that are robust tend to degrade the underlying approximation, while methods that preserve high-order accuracy often fail in the presence of strong gradients. This motivates approaches that act directly on the structure of the numerical solution.

## 2 Model

We consider systems of hyperbolic conservation laws written in the form

$$u_t + \nabla \cdot F(u) = 0,$$

where $u$ denotes the vector of conserved quantities and $F(u)$ is the flux.

A primary example is given by the compressible Euler equations, written in conservative form as

$$\frac{\partial}{\partial t}\begin{pmatrix}\rho \\ \rho \mathbf{v} \\ E\end{pmatrix} + \nabla \cdot \begin{pmatrix}\rho \mathbf{v} \\ \rho \mathbf{v} \otimes \mathbf{v} + pI \\ (E+p)\mathbf{v}\end{pmatrix} = 0,$$

where $\rho$ is the density, $\mathbf{v}$ is the velocity, $E$ is the total energy, and $p$ is the pressure.

The pressure is determined through an equation of state, typically taken in the form

$$p = (\gamma - 1)\left(E - \tfrac{1}{2}\rho|\mathbf{v}|^2\right),$$

with $\gamma$ denoting the ratio of specific heats.

## 3 Framework for moment limiters

We construct moment limiters that operate directly on the coefficients of the DG approximation. The main idea is to reconstruct slopes along a set of directions in which the moments decouple. In this representation, the multidimensional limiting problem is reduced to a sequence of one-dimensional ones.

This perspective leads to a simple and coherent construction. The limiter is local, uses compact stencils, and does not rely on large-scale reconstruction. At the same time, it remains consistent with the high-order structure of the approximation and avoids unnecessary modification of smooth solutions.

A key aspect of the framework is its suitability for parallel implementation. Because the reconstruction is performed on compact stencils and along decoupled directions, the algorithm admits straightforward parallelization. In particular, all examples are computed on graphics processing units, where the locality of the method allows efficient execution. A significant portion of the reconstruction is shifted to a preprocessing step on the mesh, reducing runtime cost and leading to a fast limiter in practice.

An important feature of the framework is its generality. The same construction applies to triangular, quadrilateral, and tetrahedral meshes without change in the underlying principles.

## 4 Properties and numerical results

The limiter was first developed for nonuniform triangular grids. For second order approximations, the limiter satisfies a local maximum principle, providing a clear stability mechanism and ensuring boundedness of the numerical solution in terms of neighboring values.

For higher order approximations, we adopt a hierarchical strategy in which higher order moments are limited first. Lower order components are modified only when necessary, reducing overlimiting and preserving accuracy.

**Quadrilateral meshes.** We next consider quadrilateral meshes. The limiter applies directly in this setting using the same directional reconstruction framework. However, it is more involved because the Jacobian of mapping a physical element to a computational element is nonlinear. The results show that the limiter preserves high-order accuracy while controlling oscillations in a manner consistent with other mesh types (Fig. 1)

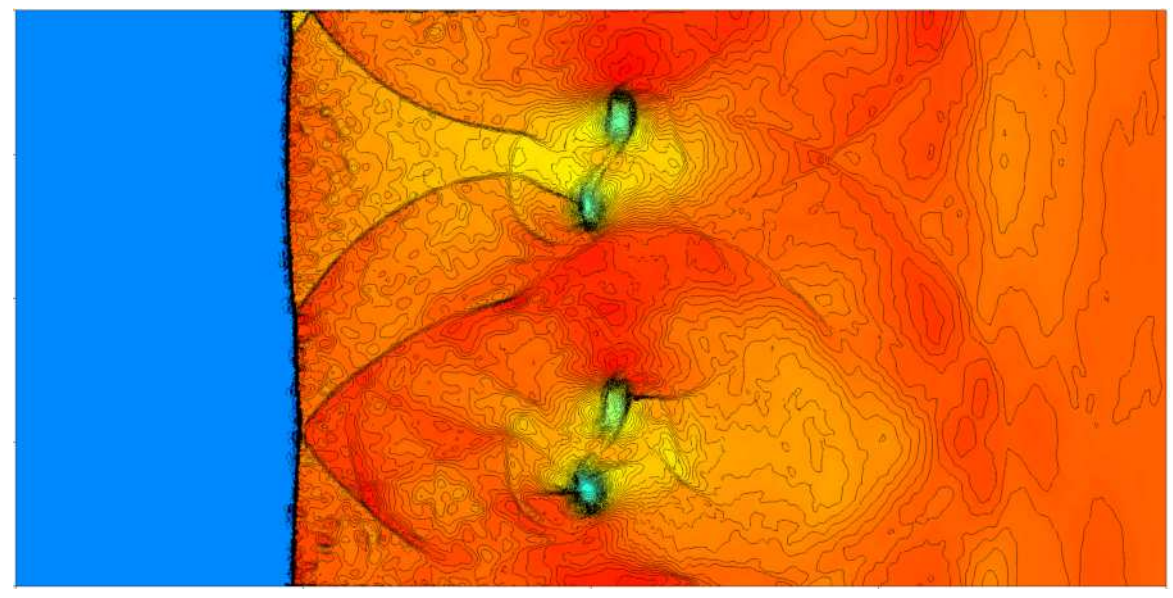

Figure 1: Shock vortex interaction example on an unstructured mesh quads with $p = 3$.

**Tetrahedral meshes.** We construct high-order moment limiters on tetrahedral meshes. The method is applied to representative three-dimensional problems, including bubble-shock interaction (Fig. 2), demonstrating stability in the presence of strong gradients and the ability to resolve fine-scale structures.

**Curved elements.** We also consider problems posed on curved elements. Such configurations arise naturally in applications involving complex geometries, where geometric representation plays a central role in achieving high-order accuracy. The limiter extends directly to curved elements through the same moment-based formulation. We show the performance on an example of a shock-cylinder interaction, with third order approximation in Fig. 3.

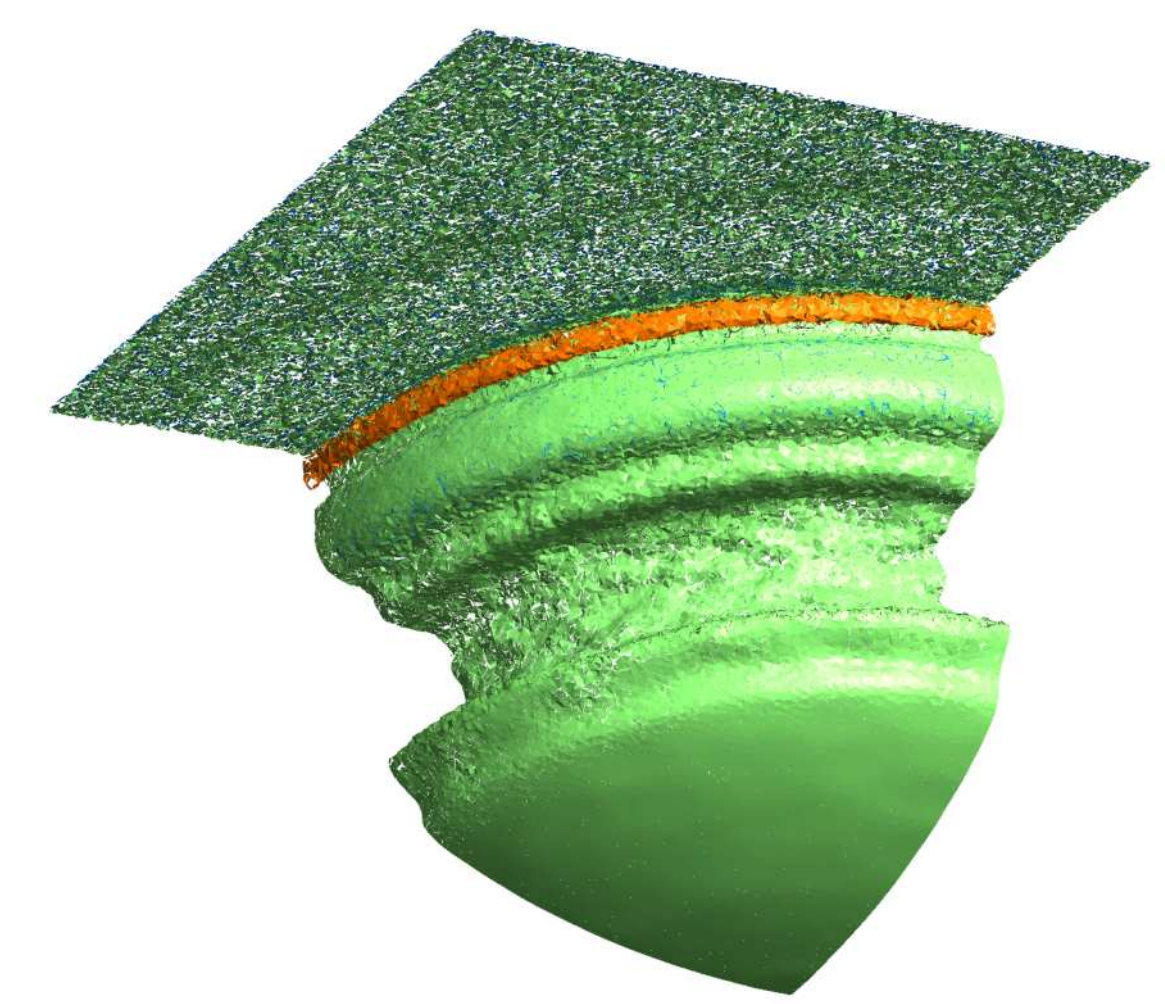

Figure 2: Bubble shock interaction case using $p = 2$.

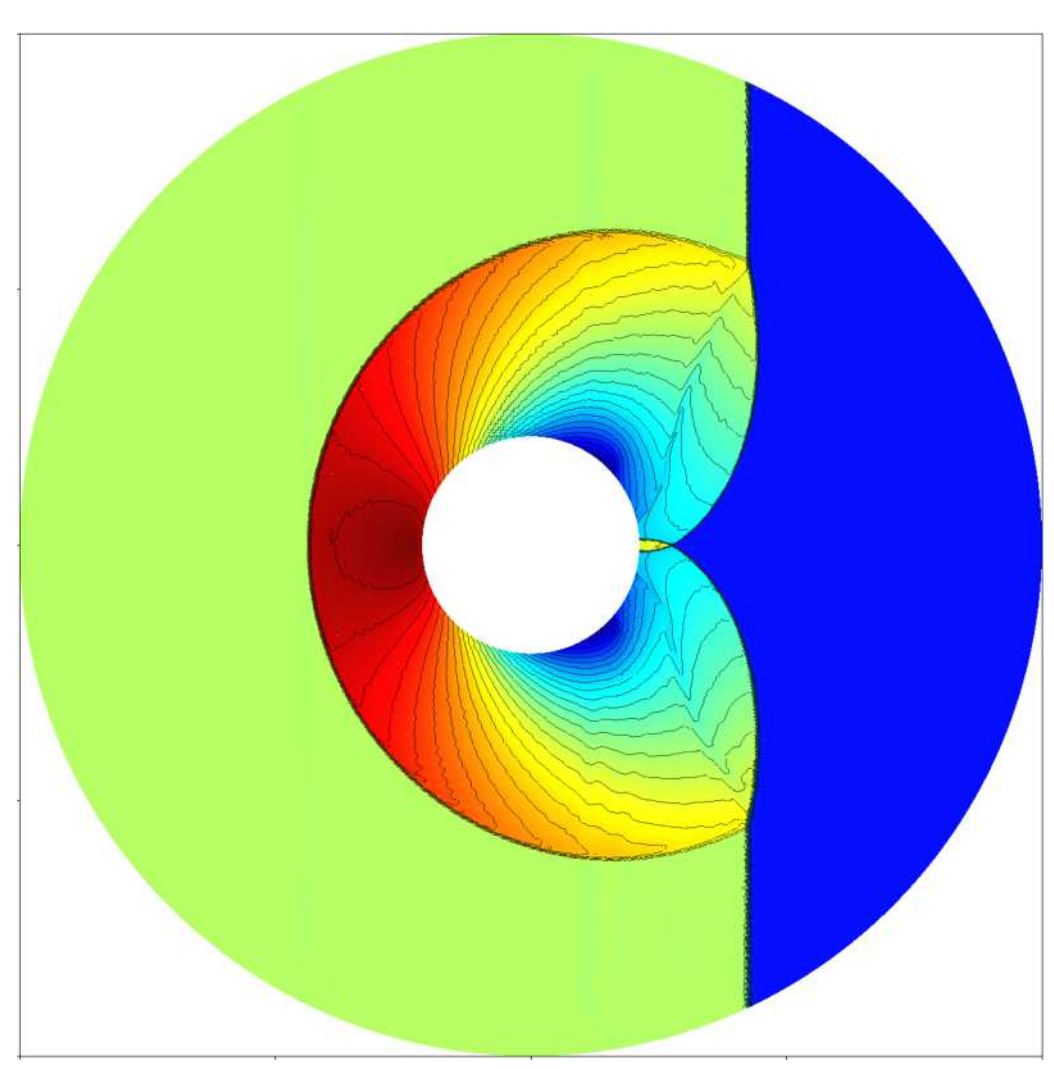

Figure 3: Shock cylinder interaction.

**Adaptive computations.** In many applications, the use of adaptive mesh refinement is essential for resolving multiscale features of the solution efficiently. Regions containing shocks, contacts, or fine-scale structures require increased resolution, while smooth regions can be represented on coarser meshes. This leads naturally to nonconforming meshes.

Within the present framework, these situations are handled without modification of the underlying limiter. The reconstruction is based on local connectivity and does not depend on mesh uniformity, allowing the same construction to be applied across refinement levels.

By combining adaptive refinement with a moment limiter, the method provides both resolution and stability. Fine-scale structures are resolved where needed, while oscillations are controlled uniformly across the mesh, including across interfaces between different refinement levels.

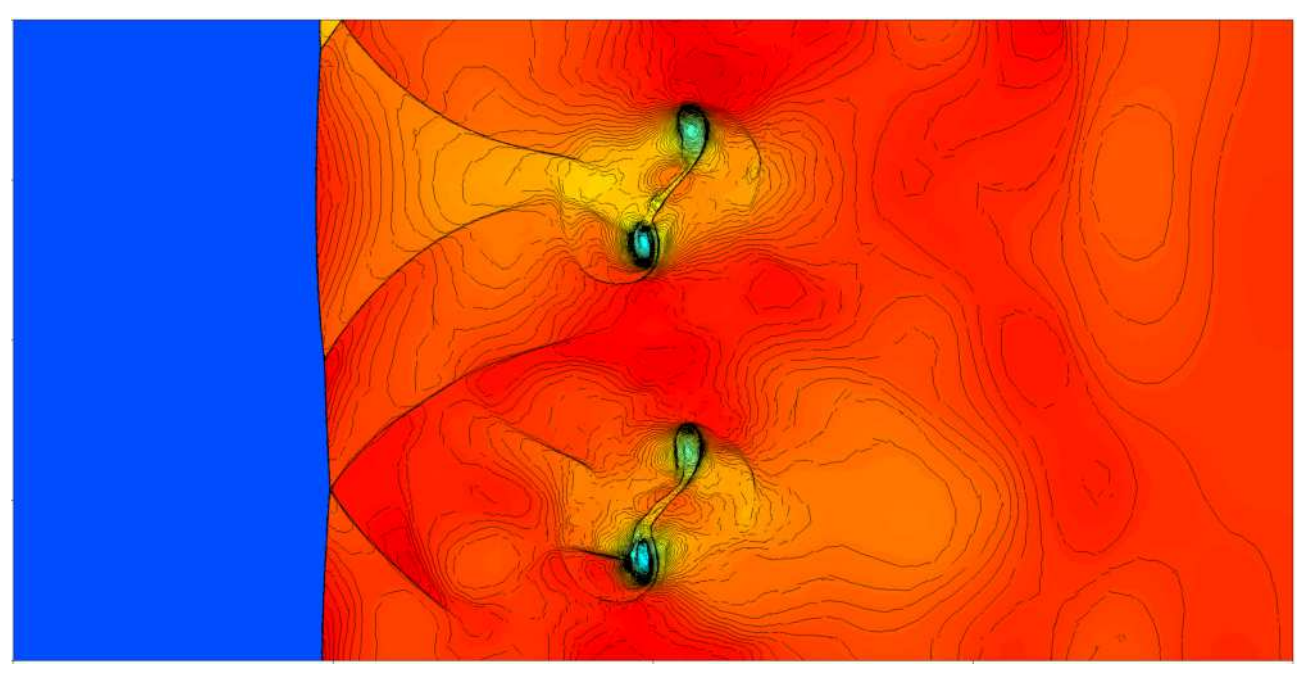

(a) Density isolines.

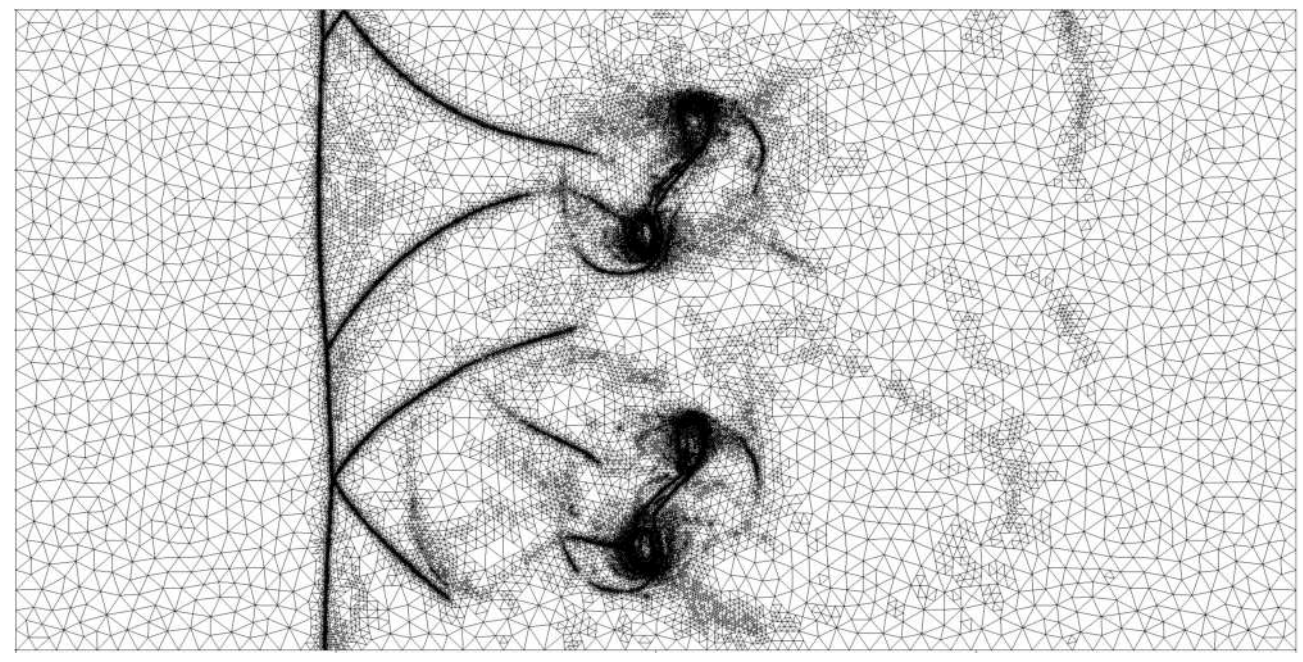

(b) Final adaptively refined mesh.

Figure 4: Shock vortex test case, obtained using a second order approximation.

# Physics-driven generative methods for inverse scattering

**Ethan Hanold**[1], **Nicholas Boffi**[2], **Leonardo Zepeda-Núñez**[3], , **Qin Li**[4]
[1]University of Wisconsin-Madison
[2]Carnegie Mellon University
[3]Google Research
[4]University of Wisconsin-Madison

**Abstract**

Inverse scattering is ubiquitous across science and engineering, yet classical methods based on PDE-constrained optimization often suffer from severe non-convexity and high computational costs. To address these challenges, we propose to leverage stochastic interpolants, a generative framework, with a physics-based initialization for the reference distribution. By anchoring the reference distribution within the generative mapping, our method reduces the transport distance to the target distribution. This yields more accurate reconstructions and requires substantially fewer integration steps than standard generative models.

## 1 Introduction

Wave-based inverse problems aim to reconstruct the properties (typically the refractive index) of an unknown medium by probing it with impinging waves and measuring the scattered impulse response at the boundary. This task arises in biomedical imaging [18], synthetic aperture radar [6], nondestructive testing [15], and geophysics [17]. The inverse problem is inherently ill-posed [9] and extremely sensitive to noise, often leading to spurious reconstructions [22].

Traditionally, hand-crafted asymptotic approximations [7] provided qualitative reconstructions [4, 8]. However, modern high-performance computing shifted the paradigm toward PDE-constrained optimization [16, 21, 23] yielding high-quality results [24]. Despite their success these methods are plagued by local minima, they require complex algorithmic pipelines [5, 13, 20], and they incur prohibitive computational costs at high resolutions.

Fortunately, over the past decade, machine learning has emerged as exciting tool for inverse problem, seeking to approaximate the inverse map *directly*. Most of these approaches are *deterministic* and they seek to leverage knowledge of the underlying physic into to network architecture to highly improve the quality of reconstruction at a fraction of the computational cost [10, 12, 27, 30]. On the other hand, recent *generative* methods adopt a probabilistic perspective [19] to reconstruct the media. Although successful in other domains, generative models for inverse scattering are hindered by spectral bias [25, 26], making the recovery of fine-grained details difficult. In this work, we incorporate physical and analytical knowledge, specifically the filtered back-projection formula, into a flow-based model to design an initial distribution and network for accurate and efficient reconstructions.

## 2 Model

We consider the two-dimensional inverse scattering problem with constant-density acoustic physics modeled by the Helmholtz equation. We seek to recover a compactly supported perturbation of the refractive index, $m$, against a homogeneous background.

**Forward problem:** The medium is illuminated by monochromatic plane waves $e^{i\omega s\cdot x}$ for incidence directions $s \in \mathrm{S}^1$. The total field decomposes as $u^{\mathrm{tot}} = e^{i\omega s\cdot x} + u^{\mathrm{sc}}$, where the scattered wave satisfies

$$\Delta u^{\mathrm{sc}}(x) + \omega^2(1+m(x))u^{\mathrm{sc}}(x) = -\omega^2 m(x) e^{i\omega s\cdot x},$$

subject to the Sommerfeld radiation condition. Measurements given by $y_\omega(r,s) = u^{\mathrm{sc}}(Rr;s)$ are collected on a circle enclosing $m$ with receivers $r \in S^1$, defining the forward map $\mathcal{F}_\omega(m) = y_\omega$.

**Inverse problem:** We aim to reconstruct $m$ from wideband data $\{y_\omega\}$, requiring the inversion of $\mathcal{F}_\omega$. Direct inversion is infeasible; thus, classical optimization minimizes the data misfit $\min_m \sum_\omega \|\mathcal{F}_\omega(m) - y_\omega\|^2$. Although highly successful [16], this approach is computationally

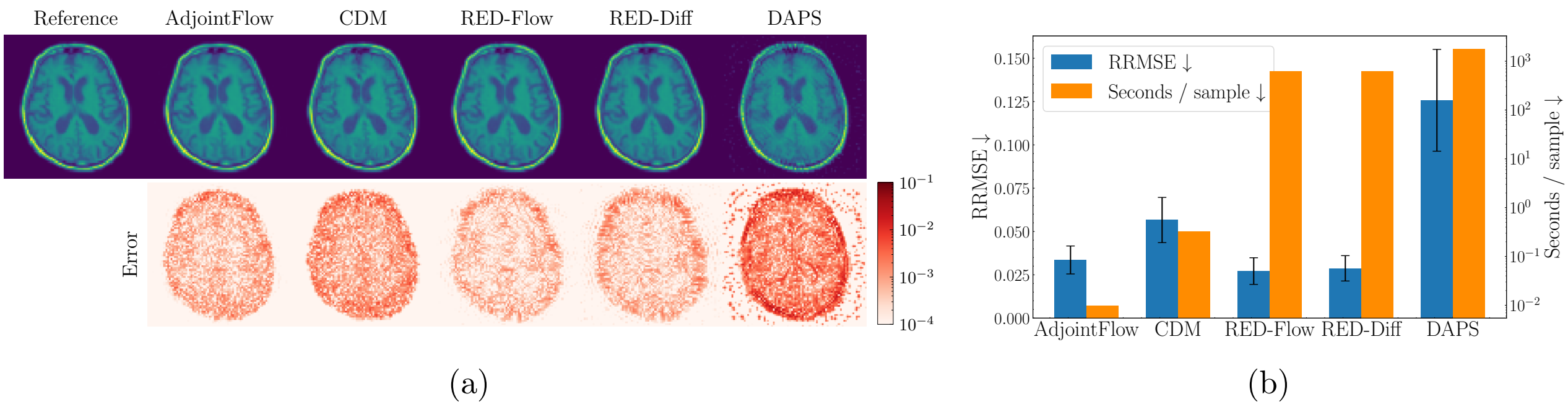

(a) (b)

Figure 1: (a) Qualitative display of AdjointFlow performance alongside conditional diffusion model (CDM) [29] as well as with RED-Flow and other PnP-DP methods [31] on the brain scan dataset. (b) Efficiency–accuracy tradeoff between methods (Metrics calculated using the same NVIDIA RTX A6000 GPU).

demanding, prone to local minima, and cost-prohibitive for posterior sampling.

## 3 Methods and Results

Our goal is to learn a generative model that *efficiently* samples the posterior $p(m \mid y)$. To this end, we leverage the stochastic interpolant (SI) framework [1] for conditional generative modeling. For each training pair $(y, m)$, we construct the interpolant $m_t = \alpha_t m_0 + \beta_t m$, where $m_0 = m_0(y)$ is a measurement-dependent base point, and $\alpha_t, \beta_t, \gamma_t$ satisfy boundary conditions $\alpha_0 = \beta_1 = 1$ and $\alpha_1 = \beta_0 = 0$. This couples the noisy measurement $y$ to the target $m$ [2].

To sample from the posterior, we learn a velocity field $\hat{b}_t(m, y)$ that transports samples along the probability flow $\dot{m}_t = b_t(m_t; y)$ initialized at $m_{t=0} = m_0(y)$, where $b_t(m; y) = \mathbb{E}[\dot{m}_t \mid m_t = m; y]$. Sampling thus reduces to numerically integrating this ODE from $t = 0$ to $t = 1$.

**Initial condition:** Generative models typically initialize $m_0$ with an uninformative Gaussian noise, thus, requiring many integration steps [29]. Instead, we use the back-projection formula to construct a physics-based initial condition $m_0 = F^* y$, where $F$ is the linearization of $\mathcal{F}$ given by the Born approximation. Anchoring the generative mapping with this physics-based initialization significantly reduces the transport distance to the target distribution, boosting efficiency and reconstruction quality.

**Results:** We consider the challenging brain scans dataset [11], which comprises MRI data. We solve the Helmholtz equation discretized with 6th-order finite differences, using a sparse multifrontal solver for the resulting linear systems. We treat the brain reconstructions from [11] as the target perturbation $m$. We simulate measurements $y$ for 5,000 distinct perturbations. For each $y$, we evaluate the initial condition $m_0(y) = F^* y$, computing the adjoint [3]. Overall, we obtain 5,000 triplets $(m, y = \mathcal{F}(m), m_0(y) = F^* y)$, which we use to train $\hat{b}_t(m, y)$ following [1].

To prevent the inverse crime during inference, we use a different discretization to compute the initial condition $m_0(y) = F^* y$. The state at $t = 1$ yields the reconstructed parameter. We evaluate the model on 1,000 hold-out samples against several state-of-the-art baselines. Figure 1(a) demonstrates that the proposed method accurately recovers the perturbation. Table 1 shows that our method outperforms competing methodologies across multiple metrics, and Figure 1(b) highlights its improved computational efficiency, which can produce samples in tens of milliseconds. This approach naturally combines with inference-time scaling techniques, such as methods based on variational techniques [14] (AdjointRED), which further improve accuracy at the expense of additional inference cost.

Table 1: Comparison of methods using RRMSE (%), (mean energy log ratio) MELR ($\times 10^{-2}$), and Sinkhorn divergence (SD) metrics.

| Method | RRMSE ↓ | MELR ↓ | SD↓ |
|---|---|---|---|
| **PnP Diffusion Prior** | | | |
| DAPS [28] | 12.572 | 18.700 | 0.278 |
| RED-Diff [14] | 2.873 | 7.370 | 0.107 |
| **Conditional Diffusion** | | | |
| EquiNet CDM [29] | 5.661 | **3.494** | 0.157 |
| **Ours** | | | |
| AdjointFlow | 3.356 | 4.776 | 0.115 |
| AdjointRED | **2.719** | 5.922 | **0.098** |

# Recovery of acoustic properties in the Westervelt equation using high-frequency large-amplitude waves

**Sebastian Acosta**[1,*], **Benjamin Palacios**[2]
[1]Texas Children's Hospital and Baylor College of Medicine, Houston, USA
[2]Department of Mathematics, Pontificia Universidad Catolica de Chile, Santiago, Chile
*Email: sacosta@bcm.edu

## Abstract

The Westervelt equation models the propagation of nonlinear ultrasound waves. For medical imaging purposes, it is needed to estimate the parameters in the equation that model wave speed, diffusivity, and nonlinearity. We explore the well-posedness of the forward and inverse problems for quadratic and general models of nonlinearity, by constructing time-periodic solutions driven at a sufficiently high frequency and large amplitude.

## 1 Introduction

Nonlinear ultrasound has become a cornerstone of modern medical technology, offering enhanced visualization of internal structures, blood perfusion in both organs and tumors. Beyond diagnostic imaging, it enables therapeutic interventions such as the ablation of pathological tissues via high-intensity focused ultrasound and the targeted delivery of genetic or pharmacological agents using ultrasound responsive microagents. The inherent portability of these systems is particularly advantageous in surgical environments, allowing for the real-time, intraoperative monitoring of both the clinical procedure and patient hemodynamics.

Recent mathematical analysis of nonlinear ultrasound has focused on identifying the source of nonlinearity through various inverse problem frameworks. Efforts have established local uniqueness based on linearized models, utilizing infinitesimally small special solutions such as Gaussian beams and wave packets for non-diffusive models, and complex geometrical optics solutions for time-periodic diffusive models. For a comprehensive review of these developments in nonlinear acoustics, see [1].

## 2 Quadratic Nonlinearity

We consider the propagation of nonlinear ultrasound waves governed by the following Westervelt equation posed for a bounded connected domain $\Omega \subset \mathbb{R}^d$, for dimension $d = 2, 3$, with smooth boundary $\partial\Omega$ and time interval $(0, T)$,

$$\left(u - \alpha u^2\right)_{tt} - \nabla \cdot (\gamma \nabla u + \beta \nabla \partial_t u) = 0, \qquad (1)$$

with boundary condition on $(0, T) \times \partial\Omega$,

$$\gamma u_\nu + \beta u_{t\nu} + \lambda u_t + \eta u = g\, e^{-i\omega t} \qquad (2)$$

and periodicity conditions

$$u|_{t=0} = u|_{t=T} \quad \text{and} \quad u_t|_{t=0} = u_t|_{t=T} \qquad (3)$$

where $T = 2\pi/\omega$ is the period, and $\omega$ is the angular frequency of oscillation for the boundary source with amplitude $g : \partial\Omega \to \mathbb{C}$. The variable $u$ denotes the fluctuations in acoustic pressure, where $\alpha(x)$ represents the nonlinearity coefficient, $\beta(x)$ the diffusivity, and $\gamma(x)$ corresponds to the squared wave speed. The expression $(u - \alpha u^2)$ within (1) originates from the constitutive relation between pressure and density variations in fluid or biological media. As this term is derived from a second-order Taylor approximation, it remains accurate only for moderate pressure amplitudes. Consequently, the model is constrained by the physical requirement that $2\alpha u < 1$. Exceeding this threshold causes the governing equation to degenerate and lose its well-posedness. To navigate this constraint, it is common to establish $L^\infty$ bounds by assuming that either the nonlinearity $\alpha$ or the boundary excitation $g$ are sufficiently small.

To build a solution for (1)-(3), the strategy followed in [2] considers the following ansatz,

$$u^K(t, x) = \sum_{k=0}^{K} u_k(x) e^{-ik\omega t} \qquad (4)$$

where the terms $u_k$ are defined recursively and the limit as $K \to \infty$ of the sequence $\{u_K\}$ is shown to converge to a function $u$ that satisfies (1)-(3) in the appropriate normed space provided that the following assumption

$$\|\alpha\|_{L^\infty(\Omega)}\|g\|_{H^{1/2}(\partial\Omega)} < C\,\omega^{1/2} \tag{5}$$

is satisfied.

Boundary excitations and corresponding boundary measurements are encoded in the so-called Robin-to-Dirichlet map given by

$$\Lambda_{\alpha,\beta,\gamma} : g \mapsto u|_{(0,T)\times\partial\Omega} \tag{6}$$

In [2] we showed that knowledge of the map $\Lambda_{\alpha,\beta,\gamma}$ is sufficient to determine the coefficient of nonlinearity $\alpha$, diffusivity $\beta$ and squared wave speed $\gamma$ uniquely (under appropriate regularity and positivity conditions). Note that the assumption (5) does not require $g$ to be infinitesimally small and that it may be allowed to be large as the frequency $\omega$ increases. However, there is still a limitation to how large the amplitude of the waves can be for this quadratic nonlinear model.

## 3 Generalized Nonlinearity

Our current efforts are aimed at generalizing the model of nonlinearity to allow waves of much larger amplitude. In practice, linear ultrasound utilizes oscillatory amplitudes in the pressure field around 10-100 kPa, while nonlinear ultrasound applications require amplitudes around 10-100 MPa. Hence, there is an unresolved need for a mathematical model that can handle solutions with very large amplitude.

One way to generalize the modeling of nonlinearity is to start from the basic laws governing fluid dynamics. Conservation of mass dictates

$$\rho_t + \nabla\cdot(\rho\boldsymbol{v}) = 0, \tag{7}$$

and balance of momentum

$$(\rho\boldsymbol{v})_t + \nabla\cdot(\rho\boldsymbol{v}\otimes\boldsymbol{v}) = -\nabla p + \nabla\cdot\boldsymbol{\tau} + \mathbf{f}, \tag{8}$$

where $\boldsymbol{v}$ denotes fluid velocity, $\rho$ is the mass density, and the stress tensor $\boldsymbol{\tau}$ satisfies $\nabla\cdot\boldsymbol{\tau} = (2\mu+\lambda)\nabla(\nabla\cdot\boldsymbol{v})$ for irrotational flows. Here the shear viscosity $\mu$ and bulk viscosity $\lambda$ are assumed to have negligible variations. When ignoring the convective terms, these equations reduce to

$$\rho_t + \nabla\cdot(\rho\boldsymbol{v}) = 0, \tag{9}$$
$$(\rho\boldsymbol{v})_t = -\nabla p - \beta\nabla\rho_t + \mathbf{f}, \tag{10}$$

where $\beta = (2\mu+\lambda)/\rho_{\text{ref}}$. Now combining the divergence of (10) with the time derivative of (9), we obtain

$$\rho_{tt} - \nabla\cdot(\beta\nabla\rho_t) - \Delta p = f \tag{11}$$

where $f = -\nabla\cdot\mathbf{f}$ augmented by the boundary condition

$$\partial_\nu p + \beta\partial_\nu\rho_t + \lambda\rho_t = g. \tag{12}$$

The density $\rho$ and pressure $p$ are related by an equation of state in differential form,

$$\nabla p = G(\rho,x)\nabla\rho. \tag{13}$$

The function $G : \mathbb{R}\times\overline{\Omega}\to\mathbb{R}$ satisfies

$$0 < G_{\min} \le G(\rho,x) \quad \text{for } (\rho,x)\in\mathbb{R}\times\overline{\Omega}, \tag{14}$$
$$|\partial_\rho G(\rho,x)| \le L \quad \text{for } (\rho,x)\in\mathbb{R}\times\overline{\Omega}, \tag{15}$$

for some constants $G_{\min}$ and $L \ge 0$. In the physical context, $G(\rho,x)$ translates into the instantaneous wave speed (squared) at a density value of $\rho$. Hence, for instance if $\partial_\rho G > 0$, $G(\rho,x)$ increases as $\rho$ increases, that is, as the medium is compressed, which is the primary characteristic of nonlinear acoustics.

Assuming apriori bounds of the form

$$\|\varrho\|_{L^\infty((0,T)\times\Omega)} \le C\,\omega^{1/2} \tag{16}$$
$$\|\varrho_t\|_{L^\infty((0,T)\times\Omega)} \le C\,\omega^{2} \tag{17}$$

where $\varrho = \rho - \overline{\rho}$ and $\overline{\rho}$ is the time-average density field, one can arrive to this stability estimate

$$\|\nabla\overline{\rho}\|^2_{L^2(\Omega)} \le CS$$
$$\omega^2\|\nabla\varrho\|^2_{L^2((0,T)\times\Omega)} \le CS$$
$$\|\nabla\varrho_t\|^2_{L^2((0,T)\times\Omega)} + \|\varrho_t\|^2_{L^2((0,T)\times\partial\Omega)} \le CS$$

where $S = \|f\|^2_{L^2((0,T)\times\Omega)} + \|g\|^2_{L^2((0,T)\times\partial\Omega)}$. These estimates reveal the role played by the frequency $\omega$ and the amplitude of the waves for forward and inverse problems. These advances will be expounded upon at the conference talk.

# Fast high-order solvers for the Lippmann-Schwinger equation in piecewise-smooth heterogeneous media

Thomas G. Anderson[1,*], Juan Burbano-Gallegos[2], Luiz M. Faria[3], Carlos Pérez-Arancibia[2]

[1]CMOR, Rice University, Houston, USA
[2]Applied Mathematics, University of Twente, Netherlands
[3]POEMS, INRIA, ENSTA Paris, Institut Polytechnique de Paris, France

*Email: thomas.anderson@rice.edu

**Abstract**

Inhomogeneous scattering problems have long been formulated as Lippmann-Schwinger (LS) or similar volume integral equations, but high-order accuracy typically requires that the volume integral operator (VIO) numerical schemes depend on unstructured domain discretizations, i.e. using (possibly curved) triangular or tetrahedral elements. These discretizations, unlike commonly-used schemes, remain high-order accurate even when, as is often the case in applications, the contrast jumps discontinuously across complex boundaries. This work proposes a class of preconditioners for high-order unstructured grid discretizations using preconditioners on structured grids to reap the positive aspects of each.

## 1 Introduction

Given a solution $u^{inc}$ to the free-space Helmholtz equation with frequency $k > 0$, we are interested in the scattering problem due to an inhomogeneity $n(x) > 0$,

$$\Delta u^t + k^2 n(x) u^t = 0 \quad \text{in} \quad \mathbb{R}^d \setminus \Gamma,$$

with the resulting scattered field $u^{sc} = u^t - u^{inc}$ satisfying the Sommerfeld radiation condition. The LS formulation of this is: seek $u^t : \mathbb{R}^d \to \mathbb{C}$ satisfying for all $x \in \Omega$,

$$u^t(x) + \mathcal{V}_{k,\Omega}[(k^2 - k^2 n)u^t](x) = u^{inc}(x), \quad (1)$$

where $\Omega$ denotes the region over which $1 - n$ does not vanish and where $\mathcal{V}_k$ denotes the VIO

$$\mathcal{V}_{k,\Omega}[f](x) = \int_\Omega G_k(x,y) f(y) \, \mathrm{d}y. \quad (2)$$

Recently proposed [1] discretization methods for eq. (2) yield high-order accurate Nyström integral equation methods using triangulations that conform to possible discontinuities in $n(x)$.

The performance of iterative solvers using these discretizations quickly degrades for high-frequency problems—typically resulting in, at least, linear growth in iteration counts as the wavenumber increases (cf. Tables 1-2). Meanwhile, successful preconditioners have been recently proposed [6] whose performance, relying on concepts such as wave sweeping or polarized traces, is impressive—showing moderate growth in, or even bounded, iteration counts as the frequency grows. However, all such methods are based on structured grid discretizations that are always low-order accurate in the presence of discontinuities in $n(x)$. Indeed, while a variety of fast solvers have been provided [2–5] in recent decades, all revert to low-order accuracy in the presence of discontinuities in $n(x)$.

## 2 Hybrid Structured-Unstructured Grid Methods

The proposed method introduces a structured grid on the rectangular (we restrict discussion to 2D, $d = 2$) domain $\Omega_R \supset \Omega$, discretized with mesh size $h_C$ so that $kh_C = \mathcal{O}(1)$. We assume also an unstructured triangular mesh of size $h_T$, $kh_T = \mathcal{O}(1)$. We solve the LS equation

$$u^t(x) + \mathcal{V}_{k,\Omega_R}[(k^2 - \kappa_R^2 n)u^t](x) = u^{inc}(x), \quad (3)$$

on this domain, where here

$$\kappa_R(x) = \begin{cases} \kappa(x) := k(n(x))^{1/2}, & x \in \Omega, \\ k, & x \in \Omega_R \setminus \Omega. \end{cases}$$

The hybrid mesh contains the $N_H$ quadrature nodes that arise from the union of quadrature nodes of the curved triangulation of $\Omega$ and the regular Cartesian grid vertices on $\Omega_R \setminus \Omega$. The Cartesian grid on $\Omega_R$ has $N_C < N_H$ quadrature/evaluation nodes (typically, $N_C \ll N_H$ for high order discretizations of $\mathcal{V}_{k,\Omega_R}$).

The preconditioning scheme proceeds by exchanging information at iterate $k$ between the grids. We call $A \in \mathbb{R}^{N_H \times N_C}$ the matrix that maps, using bilinear interpolation, from data on structured mesh vertices to the hybrid mesh quadrature nodes in $\Omega$ (and is the identity on the $N_{C,\text{out}}$ points in the Cartesian grid outside $\Omega$). Conversely, we let $C \in \mathbb{R}^{N_C \times N_H}$ denote the matrix that maps from $N_q$ quadrature nodes in $\Omega$, using Galerkin projection, to function values on triangle vertices and then to interior Cartesian grid points using barycentric interpolation ($C$ acts as an identity on hybrid mesh points outside $\Omega$). For efficiency and theoretical simplicity we use mass lumping in constructing $C$.

We denote by $B \in \mathbb{C}^{N_C \times N_C}$ some structured-grid approximation to the inverse to the LS solution operator. It follows from elementary rank theorems that a first attempt at a preconditioner,

$$L = ABC,$$

is not invertible; empirically, it is not effective as a preconditioner. We instead propose an alternative methodology that can lead to a provably-invertible preconditioner and is in practice effective.

**Theorem 1** *Let $\mu \in \mathbb{C} \setminus \{0\}$ and assume*

$$\|I_n - CA\| \left\|\mu B^{-1} - I\right\| < 1.$$

*Then*

$$P = \mu I + A(B - \mu I)C$$

*is invertible.*

Sufficient conditions are provided that ensure $\left\|\mu B^{-1} - I\right\| < 1$ in terms of the contrast $n(x)$; likewise, sufficient conditions are provided that ensure $\|I - CA\| < 1$ in terms of the relative mesh sizes $h_T$ and $h_C$.

## 3 Numerical Results

We solve a problem with smooth refractive index for which we can compare accuracy between two high-order accurate (4th-order for the Cartesian-grid preconditioner, 5th-order for unstructured) solvers, and then demonstrate the iteration counts for a discontinuous media problem. In the smooth media case we keep $kh_C = 1/8$ while for the discontinuous media case we let $kh_C = 1/4$; in all cases $h_T = 2/3h_C$ and $\Omega$ is the unit ball $B_1(0)$. The GMRES iterative solver is chosen and the requested tolerance is set to $10^{-6}$.

| $k$ | # Struct. Its | # Hybrid Its | Accuracy |
|---|---|---|---|
| 2 | 3 | 3 | $1.7e-10$ |
| 4 | 3 | 3 | $7.3e-10$ |
| 8 | 4 | 4 | $5.2e-11$ |
| 16 | 4 | 4 | $2.9e-9$ |
| 32 | 5 | 5 | $4.5e-10$ |
| 64 | 4 | 5 | $3.6e-9$ |

Table 1: Smooth media solver statistics; $n(x) = 1 - 0.7e^{-40(x^2+y^2)}$.

| $k$ | # Struct. Its | # Hybrid Its |
|---|---|---|
| 2 | 6 | 8 |
| 4 | 6 | 8 |
| 8 | 7 | 9 |
| 16 | 7 | 11 |
| 32 | 8 | 13 |
| 64 | 12 | 19 |

Table 2: Discontinuous media solver statistics; $n(x) = 2$ for $x^2 + y^2 < 1$.

# Time dependent acoustic scattering from penetrable objects with nonlinearities

**Paschalis Angelidis**[1,*], **Alexander Rieder**[1]
[1]Institute of Analysis and Scientific Computing, TU Wien, Vienna, Austria
*Email: paschalis.angelidis@tuwien.ac.at

**Abstract**

In this work we consider the propagation of acoustic waves with nonlinearities in full space. A Time Domain Boundary Integral Formulation (TDBIE) is used for the exterior unbounded domain and is coupled with the semi-linear wave equation in the interior domain. For the spatial discretization we use the Finite Element Method for the variational formualtion of the interior and the Boundary Element Method (BEM) for the TDBIE. For the temporal discretization we propose a fully implicit scheme, namely we use an Implicit Euler scheme in the interior and Convolution Quadrature (CQ) for the convolution integrals arising from the TDBIE. We prove stability of the fully discrete scheme by energy methods and use this estimate along with the convergence properties of the nonlinear term and the properties of the boundary potentials to prove error bounds of the fully discrete scheme.

## 1 Introduction

For Lipschitz domains $U \subset \mathbb{R}^d$, we denote $H^1_\Delta(U)$ such that $H^1_\Delta(U) := \{u \in H^1(U) : \Delta u \in L^2(U)\}$. We consider the following scattering problem for $u(t) \in H^1_\Delta(\mathbb{R}^d \setminus \Gamma)$:

$$\begin{cases} \ddot{u} - \Delta u = f(u) & \text{in } \Omega \times (0,T] \\ \ddot{u} - \Delta u = 0 & \text{in } \Omega^{\text{ext}} \times (0,T] \\ [\![\gamma u]\!] = -\gamma^{\text{ext}} u^{\text{inc}} & \text{on } \Gamma \times (0,T] \\ [\![\partial_n u]\!] = -\partial_n^{\text{ext}} u^{\text{inc}} & \text{on } \Gamma \times (0,T] \\ u(0) = u_0,\ \dot{u}(0) = v_0 & \text{in } \Omega \\ u(0) = \dot{u}(0) = 0 & \text{in } \Omega^{\text{ext}} \end{cases} \tag{1}$$

where $\Omega \subset \mathbb{R}^d$ is an open, bounded, Lipschitz domain, $\Gamma = \partial\Omega$, $\Omega^{\text{ext}} := \mathbb{R}^d \setminus \overline{\Omega}$ and $u^{\text{inc}}$ is an incident wave that satisfies the homogeneous wave equation in the exterior. The nonlinear function $f(u) = -\frac{dG}{du}(u) = -G'(u)$ is considered to be the negative antiderivative of a function $G(u)$ with respect to the unknown $u$ where we assume the following for $d = 3$:

(i) $G \in C^2(\mathbb{R})$ and $G(\eta) \geq 0,\ \forall \eta \in \mathbb{R}$.

(ii) $|G'(\eta)| \leq C(1 + |\eta|^q)$ where $2 \leq q < 5$.

(iii) $|G''(\eta)| \leq C(1 + |\eta|^p)$, where $p < 4$.

One example that satisfies these conditions is $G(\eta) = 1 - \cos\eta$, $\eta \in \mathbb{R}$, i.e. the sine-Gordon equation.

## 2 Coupled Domain-Boundary Integral Formulation

To derive a weak formulation of the interior equation, one has to test with a test function $v \in H^1(\Omega)$ and integrate by parts in the spatial variable as follows:

$$(\ddot{u}, v)_\Omega + (\nabla u, \nabla v)_\Omega - \langle \partial_n^{\text{int}} u, \gamma^{\text{int}} u \rangle_\Gamma + (G'(u), v)_\Omega = 0 \tag{2}$$

for all $v \in H^1(\Omega)$. The solution now in the exterior can be written in terms of layer potentials on the boundary using the following representation formula:

$$u = -\mathcal{S}(\partial_t)\partial_n^{\text{ext}} u + \mathcal{D}(\partial_t)\gamma^{\text{ext}} u \tag{3}$$

where $\mathcal{S}(s) : H^{-1/2}(\Gamma) \to H^1_\Delta(\Omega^{\text{ext}})$ and $\mathcal{D}(s) : H^{1/2}(\Gamma) \to H^1_\Delta(\Omega^{\text{ext}})$ are the single and double layer potential operators respectively. Denoting the boundary densities as $\psi = \gamma^{\text{ext}}\dot{u}$ and $\phi = -\partial_n^{\text{ext}} u$ and using the properties of the layer potentials along with the transmission conditions on the boundary we arrive to the following coupled system: find $u \in H^1(\Omega)$, $\phi \in H^{-1/2}(\Gamma)$ and $\psi \in H^{1/2}(\Gamma)$ such that

$$\begin{aligned} (\ddot{u}, v)_\Omega + (\nabla u, \nabla v)_\Omega + \langle \phi, \gamma^{\text{int}} v \rangle_\Gamma \\ + (G'(u), v)_\Omega = \langle \partial_n^{\text{ext}} u^{\text{inc}}, \gamma^{\text{int}} v \rangle_\Gamma \\ \begin{pmatrix} -\gamma^{\text{int}}\dot{u} \\ 0 \end{pmatrix} + B(\partial_t) \begin{pmatrix} \phi \\ \psi \end{pmatrix} = \begin{pmatrix} -\gamma^{\text{ext}}\dot{u}^{\text{inc}} \\ 0 \end{pmatrix} \end{aligned} \tag{4}$$

$v \in H^1(\Omega)$, and $B(s) : H^{-1/2}(\Gamma) \times H^{1/2}(\Gamma) \to H^{1/2}(\Gamma) \times H^{-1/2}(\Gamma)$ is the acoustic Calderon operator that collects Dirichlet and Neumann traces of the potentials.

## 3 Numerical scheme and results

We introduce a finite element and boundary element spaces such as $V_h \subset H^1(\Omega)$, $X_h \subset H^{1/2}(\Gamma)$ and $Y_h \subset H^{-1/2}(\Gamma)$ for the spatial approximation. In time, for a fixed $\Delta t > 0$ we use the abbreviation $t_n = n\Delta t$ for $n = 0, \dots, T/\Delta t$. It is also useful to introduce the approximation of the first order time derivative of an arbitrary scalar function $\varphi$ as follows:

$$d_t\varphi^n = \frac{\varphi^n - \varphi^{n-1}}{\Delta t} \tag{5}$$

and second order time derivative:

$$d_t^2\varphi^n = \frac{\varphi^n - 2\varphi^{n-1} + \varphi^{n-2}}{\Delta t^2} \tag{6}$$

The time approximation of the nonlinear term follows similarly:

$$d_u G(\varphi^n) = \begin{cases} \frac{G(\varphi^n) - G(\varphi^{n-1})}{\varphi^n - \varphi^{n-1}}, & \varphi^n \neq \varphi^{n-1} \\ -f(\varphi^n), & \varphi^n = \varphi^{n-1} \end{cases}$$

We also look to approximate convolution integrals of the form:

$$K(\partial_t)\varphi := \mathscr{L}^{-1}\left\{K(\cdot)\mathscr{L}\left\{\varphi\right\}\right\}$$

where $K$ is the Laplace transform of a kernel $k$, that satisfies $|K(s)| \lesssim |s|^\mu$, $\mu \in \mathbb{R}$. This can be done by replacing $K(\partial_t)\varphi(t)$ with a discrete convolution such that:

$$K(\partial_n^{\Delta t})\varphi(t_n) = \sum_{j=0}^{n} \omega_{n-j}(K)\varphi(t_j)$$

and follow the methods described in [3]. Putting everything together we have the following discretization of (1): find $u_h^n \in V_h$, $\phi_h^n \in Y_h$ and $\psi_h^n \in X_h$ such that:

$$\begin{aligned} &(d_t^2 u_h^n, v)_\Omega + (\nabla u_h^n, \nabla v)_\Omega + \langle \phi_h^n, \gamma^{\text{int}} v\rangle_\Gamma \\ &+(d_u G(u_h^n), v)_\Omega = \langle \partial_n^{\text{ext}} u^{\text{inc}}(t_n), \gamma^{\text{int}} v\rangle_\Gamma \\ -&\langle \gamma^{\text{int}}(d_t u_h^n), \xi\rangle_\Gamma + \left\langle B(\partial_t^{\Delta t}) \begin{pmatrix} \phi_h \\ \psi_h \end{pmatrix}(t_n), \begin{pmatrix} \xi \\ \eta \end{pmatrix} \right\rangle_\Gamma \\ &= -\langle \gamma^{\text{ext}} \dot{u}^{\text{inc}}(t_n), \xi\rangle_\Gamma \end{aligned} \tag{7}$$

$v \in V_h$, $\xi \in Y_h$ and $\eta \in X_h$. The following theorems state the stability of the fully discrete scheme with repsect to the energy functional:

$$E_h^n = \frac{1}{2}\|d_t u_h^n\|^2_{L^2(\Omega)} + \frac{1}{2}\|u_h^n\|^2_{H^1(\Omega)} + \int_\Omega G(u_h^n)\, dx$$

and also the convergence of the numerical methods towards the continuous solutions $(u, \phi, \psi)$.

**Theorem 1** *Under sufficient regularity of $u^{\text{inc}}$, the system (7) admits the following stability estimate for $n = 1, \dots, N$:*

$$E_h^n \leq C(T)\left(E_h^0 + \max_n(R_h^n)\right) \tag{8}$$

*where $C(T)$ is a constant that scales exponentially with the end time $T$.*

To state the theorem for the error bounds we need to define the following: $e_u = R_h u - u_h$, $e_\phi = \Pi_h\phi - \phi_h$ and $e_\psi = Q_h\psi - \psi_h$, where $R_h$ is the elliptic projection and $\Pi_h$, $Q_h$ are the $L^2$ projections into the boundary element spaces. Moreover we need to define the product norm:

$$|\!|\!|(\phi, \psi)|\!|\!|^2_\Gamma = \|\phi\|^2_{H^{-1/2}(\Gamma)} + \|\psi\|^2_{H^{1/2}(\Gamma)} \tag{9}$$

**Theorem 2** *Under the assumption that the exact solution $u$ and the boundary densities $\phi, \psi$ have sufficient smoothness, the temporal discretization of the boundary integral operators is $A$-stable of first order and the growth condition of $G'$ is $q \leq 3$, then*

$$\|d_t e_u^n\|_{L^2(\Omega)} + \|e_u^n\|_{H^1(\Omega)} \leq C(\Delta t + h) \tag{10}$$

*and*

$$\Delta t \sum_{n=0}^{\infty} \left|\!\left|\!\left| \begin{pmatrix} (\partial_t^{\Delta t})^{-1} e_\phi \\ (\partial_t^{\Delta t})^{-1} e_\psi \end{pmatrix}(t_n) \right|\!\right|\!\right|_\Gamma \leq C(\Delta t + h) \tag{11}$$

# Numerical Study of Eigenvector Deflation to Accelerate the WaveHoltz Method

**Daniel Appelö**[1,*], **William D. Henshaw**[2], **Zhichao Peng**[3]
[1]Department of Mathematics, Virginia Tech., Blacksburg, USA
[2]Department of Mathematical Sciences, Rensselaer Polytechnic Institute, Troy, NY 12180, USA
[3]Department of Mathematics, The Hong Kong University of Science and Technology, Hong Kong
*Email: appelo@vt.edu

**Abstract**

## 1 Introduction

In this work we study how deflation techniques can be used to accelerate the WaveHoltz method with energy-conserving boundary conditions. For such boundary conditions the WaveHoltz method, by itself, converges in approximately $O(\omega^d)$ iterations, where $d$ is the dimension of the problem.

We are interested in finding numerical approximations to the solution $u = u(\mathbf{x}) \in \mathbb{R}$ of a Helmholtz boundary value problem (BVP)

$$\begin{aligned} \nabla \cdot (c(\mathbf{x})^2 \nabla) u + \omega^2 u = f(x), \quad & \mathbf{x} \in \Omega, \\ \mathcal{B} u(\mathbf{x}) = 0, \quad & \mathbf{x} \in \partial\Omega. \end{aligned} \tag{1}$$

Here the boundary condition operator $\mathcal{B}$ is taken here to be a Dirichlet condition.

The WaveHoltz method is an iterative solver to find solutions to (1) that is built upon a time-domain solver for the associated wave equation initial boundary value problem (IBVP)

$$\begin{aligned} & w_{tt} = \nabla \cdot (c(\mathbf{x})^2 \nabla) w - f(\mathbf{x}) \cos(\omega t), \\ & w(\mathbf{x}, t) = 0, \quad \mathbf{x} \in \partial\Omega, \\ & w(\mathbf{x}, 0) = v^{(k)}(\mathbf{x}), \\ & w_t(\mathbf{x}, 0) = 0. \end{aligned} \tag{2}$$

## 2 Deflated WaveHoltz

The original WaveHoltz iteration with deflation additions can be found in Figure 1. The algorithm describes how the WaveHoltz algorithm can be accelerated using some subset $\phi_m \in \mathcal{D}$ of the eigenfunctions of the Laplacian, where $\mathcal{D}$ is called the *deflation set*. We will denote the dimension of this set by $N_{\mathcal{D}}$. There are two key changes to the original algorithm. The first is on line 5 where the forcing $f$ is adjusted by removing the components of $f$ along the eigenfunctions in the deflation set. Here $(f, \phi_m)_\Omega$ denotes the inner product of $f$ and $\phi_m$. The second key change is on line 13 where the components of the solution in the deflation set are added to the WaveHoltz solution. Line 10 contains an additional change where the WaveHoltz iterate $v^{(k+1)}$ is deflated after each iteration. This line is useful to add when approximate eigenvectors are being used. It is however not added in the computational results presented here.

## 3 Numerical Example

We consider the wave equation with constant speed of sound equal to one and with homogenous Dirichlet boundary conditions discretized on a single curvilinear grid. The curvilinear grid is described by a mapping $x(r, s) = r + 0.1 \sin(2s)$, $y(r, s) = s - 0.3 \sin(2r)$, $[r, s] \in [-1, 1]^2$. The Laplacian is discretized using a conservative formulation and the metric coefficients and derivatives are approximated using summation by parts finite difference operators. Boundary conditions are enforced by projection. We use fourth order accurate SBP operators and determine the number of grid points according to the simplified rule of thumb

$$\mathrm{PPW}_p = \pi \left[ \frac{N_\Lambda}{\epsilon} \right]^{1/p},$$

with $\epsilon = 10^{-4}$. We use WaveHoltz with $N_p = 3$ periods and solve the time dependent wave equation using an implicit method with 10 timesteps per period. The forcing is a point source at $x_0 = 0.1, y_0 = 0.2$.

Choosing $\omega = 40$ results in a grid with $756 \times 756$ unknowns. We use 100 eigenpairs for the deflation, computed to a tolerance of $10^{-14}$. These 100 eigenpairs pairs are found using 593 wave-solves. The location of the deflated eigenvalues on the $\beta$-curve (see [2]) are displayed in Figure 1 where the ACR can also be read off to be $\approx 0.85$.

Figure 1 also display the residuals for this problem when using the standard CG accelerated WaveHoltz with with and without the forc-

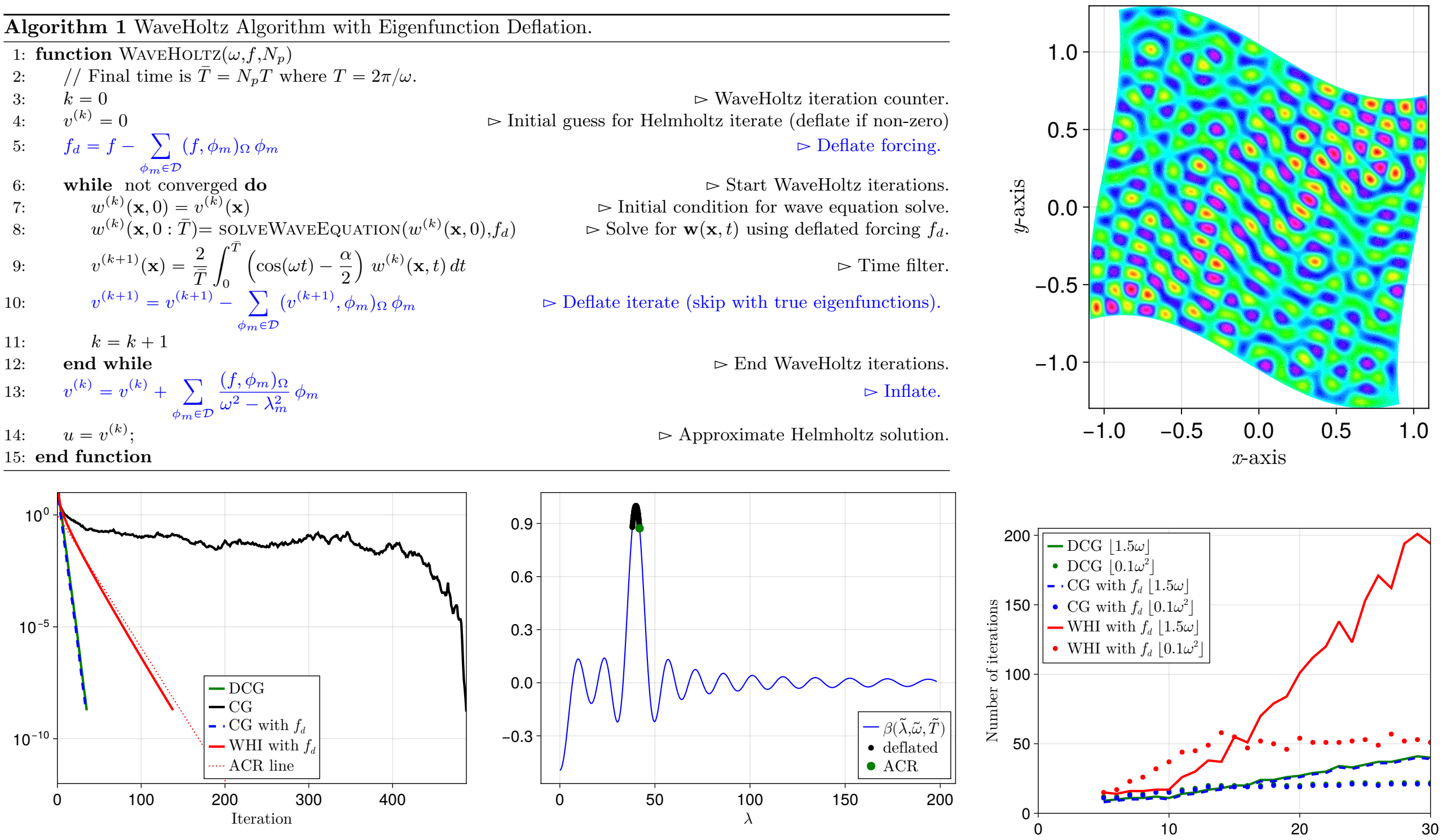
**Algorithm 1** WaveHoltz Algorithm with Eigenfunction Deflation.
1: **function** WaveHoltz($\omega$,$f$,$N_p$)
2: // Final time is $\bar{T} = N_p T$ where $T = 2\pi/\omega$.
3: $k = 0$ ▷ WaveHoltz iteration counter.
4: $v^{(k)} = 0$ ▷ Initial guess for Helmholtz iterate (deflate if non-zero)
5: $f_d = f - \sum_{\phi_m \in \mathcal{D}} (f, \phi_m)_\Omega \, \phi_m$ ▷ Deflate forcing.
6: **while** not converged **do** ▷ Start WaveHoltz iterations.
7: $w^{(k)}(\mathbf{x}, 0) = v^{(k)}(\mathbf{x})$ ▷ Initial condition for wave equation solve.
8: $w^{(k)}(\mathbf{x}, 0 : \bar{T})$= SolveWaveEquation($w^{(k)}(\mathbf{x}, 0)$,$f_d$) ▷ Solve for $\mathbf{w}(\mathbf{x}, t)$ using deflated forcing $f_d$.
9: $v^{(k+1)}(\mathbf{x}) = \frac{2}{\bar{T}} \int_0^{\bar{T}} \left(\cos(\omega t) - \frac{\alpha}{2}\right) w^{(k)}(\mathbf{x}, t) \, dt$ ▷ Time filter.
10: $v^{(k+1)} = v^{(k+1)} - \sum_{\phi_m \in \mathcal{D}} (v^{(k+1)}, \phi_m)_\Omega \, \phi_m$ ▷ Deflate iterate (skip with true eigenfunctions).
11: $k = k + 1$
12: **end while** ▷ End WaveHoltz iterations.
13: $v^{(k)} = v^{(k)} + \sum_{\phi_m \in \mathcal{D}} \frac{(f, \phi_m)_\Omega}{\omega^2 - \lambda_m^2} \phi_m$ ▷ Inflate.
14: $u = v^{(k)}$; ▷ Approximate Helmholtz solution.
15: **end function**
$y$-axis
$x$-axis
Iteration
DCG
CG
CG with $f_d$
WHI with $f_d$
ACR line
$\lambda$
$\beta(\tilde{\lambda}, \tilde{\omega}, \bar{T})$
deflated
ACR
Number of iterations
DCG $\lfloor 1.5\omega \rfloor$
DCG $\lfloor 0.1\omega^2 \rfloor$
CG with $f_d$ $\lfloor 1.5\omega \rfloor$
CG with $f_d$ $\lfloor 0.1\omega^2 \rfloor$
WHI with $f_d$ $\lfloor 1.5\omega \rfloor$
WHI with $f_d$ $\lfloor 0.1\omega^2 \rfloor$

Figure 1: WaveHoltz Algorithm with Eigenfunction Deflation.

ing being deflated (denoted "CG" and "CG with $f_d$"), the plain WaveHoltz iteration applied to the deflated forcing (denoted "WHI with $f_d$"), and the deflated conjugate gradient method by Saad (denoted "DCG") applied to the problem with the non-deflated forcing. As can be seen CG with $f_d$ and DCG have the lowest iteration count, followed by WHI with $f_d$ (whoose rate of convergence is well predicted by the ACR). The standard CG method is significantly slower.

Ignoring the cost of the deflation and inflation we can estimate the cost to solve for $N_f$ right hand sides. If, for the lower frequency we take the number of wave solves to reach convergence to be $\sim 20$ for the deflated (D)CG method and $\sim 500$ for the standard CG method we can estimate the cost (in number of wave solves) to be "Cost deflated WaveHoltz" $\sim 20N_f + 593$, "Cost standard WaveHoltz" $\sim 500N_f$. From this we can conclude that (in terms of number of wave solves) the deflated method can be faster already with as few as two right hand sides.

Finally, in the lower right graph in Figure 1 we display the number of iterations needed to reach a residual reduction my 10 orders of magnitude for different and different choices for the size of the deflation space. We chose this dimension to grow either linearly with as $\frac{3\omega}{2}$ , or growing quadratically with as $\frac{\omega^2}{10}$ . As can be seen these result in a linear growth in the number of iterations or constant number of iteratons when the dimension of the deflation space is linear and quadratic, respectively.

# A fast direct solver based on FEM-BEM coupling for time-harmonic electromagnetic scattering problems

**Emanuele Arcese**[1,*], **Francis Collino**[2]
[1]CEA CESTA, Le Barp, France
[2]Electromagnetics consultant, Paris, France
*Email: emanuele.arcese@cea.fr

### Abstract

We develop a fast direct solver based on the coupling of finite element and boundary element methods (FEM-BEM) for transmission problems of time-harmonic electromagnetic scattering involving complex 3-D objects, which can be made up of anisotropic and/or highly inhomogeneous components. This solver combines the Schur complement technique with hierarchical matrix ($\mathcal{H}$-matrix)-based compression algorithms.

## 1 Introduction

The FEM-BEM coupling under consideration is the formulation '(C)' (as presented in (1)) from Levillain's doctoral research [1]. It can be viewed as a domain decomposition technique that couples in a symmetric fashion the Rumsey's reaction principle-based boundary integral formulation and two weak formulations of the second-order Maxwell's equations in both the electric and the magnetic fields. Although this formulation introduces twice as many FEM-related unknowns as classical coupling approaches (as, e.g., the Costabel symmetric formulation [2]), it becomes particularly advantageous when dealing with an inhomogeneous material confined to localized regions within a geometrically complex object. This is primarily due to the formulation's high accuracy in far-field analysis and its ease of implementation into existing PMCHWT-type BEM solvers. A notable characteristic of the linear system resulting from this combined discretization is its dense-sparse structure. For its numerical solution, compression and preconditioning techniques have been extensively studied; however, significant challenges remain when addressing real-life problems. In this presentation, we will discuss recent advancements in designing a fast hierarchical direct solver for the aforementioned coupled problem in a high performance computing context.

## 2 Problem setting

We consider the transmission problem derived from time-harmonic Maxwell's equations with variable material parameters. The scattering object occupies a bounded domain $\Omega_2 \subset \mathbb{R}^3$ having a non-smooth boundary $\Gamma$, which is equipped with a unit normal vector $\mathbf{n}$ directed from $\Omega_2$ into $\Omega_1 := \mathbb{R}^3 \backslash \overline{\Omega}_2$. For simplicity, the object is assumed to consist of a material with parameters $\varepsilon = \varepsilon(\mathbf{x})$ and $\mu = \mu(\mathbf{x})$, $\mathbf{x} \in \Omega_2$, and to be illuminated by an incident plane wave whose fields are $\mathbf{E}^I$ and $\mathbf{H}^I$. For the derivation of Levillain variational formulation of this scattering problem, the general idea is to simply use the transmission conditions on $\Gamma$ : $\gamma_\pi^+\mathbf{E} = \gamma_\pi^-\mathbf{E}$ and $\gamma_\pi^+\mathbf{H} = \gamma_\pi^-\mathbf{H}$; where $\mathbf{E} = \mathbf{E}(\mathbf{x})$ and $\mathbf{H} = \mathbf{H}(\mathbf{x})$ are the total electric and (normalized) magnetic fields respectively and $\gamma_\pi^\pm\mathbf{U}$ stands for the tangential components trace of a vector $\mathbf{U}$, $\mathbf{n} \times (\mathbf{U} \times \mathbf{n})$, onto $\Gamma$ from outside and inside $\Omega_2$, respectively. Particularly, the formulation reads: seek $(\mathbf{E}_2, \mathbf{H}_2, \mathbf{J}, \mathbf{M}) \in \boldsymbol{V}$, provided that the wave number $\kappa > 0$, such that

$$\begin{aligned} &\mathrm{i}\kappa \int_\Gamma \mathcal{R}_\kappa \mathbf{J} \cdot \mathbf{J}' d\gamma + \int_\Gamma \mathcal{Q}_\kappa \mathbf{M} \cdot \mathbf{J}' d\gamma \\ &+ \int_\Gamma \mathcal{Q}_\kappa \mathbf{J} \cdot \mathbf{M}' d\gamma - \mathrm{i}\kappa \int_\Gamma \mathcal{R}_\kappa \mathbf{M} \cdot \mathbf{M}' d\gamma \\ &- \frac{1}{2\mathrm{i}\kappa} \left( a_E \left( \mathbf{E}_2, \mathbf{E}_2' \right) - a_H \left( \mathbf{H}_2, \mathbf{H}_2' \right) \right) \\ &= \int_\Gamma \gamma_\pi^+ \mathbf{E}^I \cdot \mathbf{J}' - \gamma_\pi^+ \mathbf{H}^I \cdot \mathbf{M}' d\gamma, \end{aligned} \tag{1}$$

for all $\left(\mathbf{E}_2', \mathbf{H}_2', \mathbf{J}', \mathbf{M}'\right) \in \boldsymbol{V}$. Here, $\mathbf{J} = \mathbf{n} \times \mathbf{H}_2$ and $\mathbf{M} = \mathbf{E}_2 \times \mathbf{n}$ represent the surface currents with notation $\mathbf{U}_2 = \mathbf{U}|_{\Omega_2}$, the function space is $\boldsymbol{V} = [\boldsymbol{H}(\mathbf{curl}, \Omega_2)]^2 \times \left[\boldsymbol{H}^{-1/2}(\mathsf{div}_\Gamma, \Gamma)\right]^2$, the boundary integral operators $\mathcal{R}_\kappa$ and $\mathcal{Q}_\kappa$ are obtained by applying $\gamma_\pi$ to Maxwell single and double layer potentials respectively, $a_E$ and $a_H$ are the sesqui-linear forms defined as follows:

$$a_E(\mathbf{E}_2, \mathbf{E}_2') := \kappa^2 \left(\varepsilon \mathbf{E}_2, \mathbf{E}_2'\right)_0 - \left(\mu^{-1} \mathbf{curl}\, \mathbf{E}_2, \mathbf{curl}\, \mathbf{E}_2'\right)_0,$$
$$a_H(\mathbf{H}_2, \mathbf{H}_2') := \kappa^2 \left(\mu \mathbf{H}_2, \mathbf{H}_2'\right)_0 - \left(\varepsilon^{-1} \mathbf{curl}\, \mathbf{H}_2, \mathbf{curl}\, \mathbf{H}_2'\right)_0,$$

with $(\cdot,\cdot)_0$ being the $\boldsymbol{L}^2(\Omega_2)$-inner product.

## 3 Methods and numerical results

We discretize system (1) by means of a standard lowest-order Galerkin method combining Nédélec's **curl**-conforming finite elements on tetrahedral meshes of the interior domain with Raviart-Thomas' div-conforming boundary elements on triangular meshes of the coupling surface. The 4-D singular integrals involving the boundary integral operators are evaluated numerically via the integration scheme described in [3]. The resulting linear system features a matrix with a dense-sparse structure. For its numerical solution, we propose a fully direct method consisting mainly in computing two Schur complements, arising from the interior problems, and assembling them on the fly into the $\mathcal{H}$-matrix associated to integral operators, which is then factorized in a $LDL^t$ form. The proposed method is implemented in the CEA in-house BEM code, which relies on a task-based hierarchical direct solver using both shared- and distributed-memory parallelisms.

We compare the designed approach with two existing coupling formulations: a 3-D solver (noted as DDMRIK) [4] and a 2-D axis-symmetric solver (assumed here to provide a reference solution); on a radar cross section (RCS) test case from the 2023 EM ISAE workshop. This deals with a cylinder with hemispherical ends in which a small toroidal region is filled by an inhomogenous material ($\kappa = 104.8\ \mathrm{m}^{-1}$). Both existing formulations address in a classic manner the problem by employing boundary elements for the radiation condition and finite elements for the whole interior problem. Instead, the proposed approach enables the use of finite elements solely for the inhomogeneous domain while treating the homogeneous problems with boundary elements (see fig. 1). Figure 2 and table 1 demonstrate the ability of the present solver to obtain accurate results while controlling the problem size.

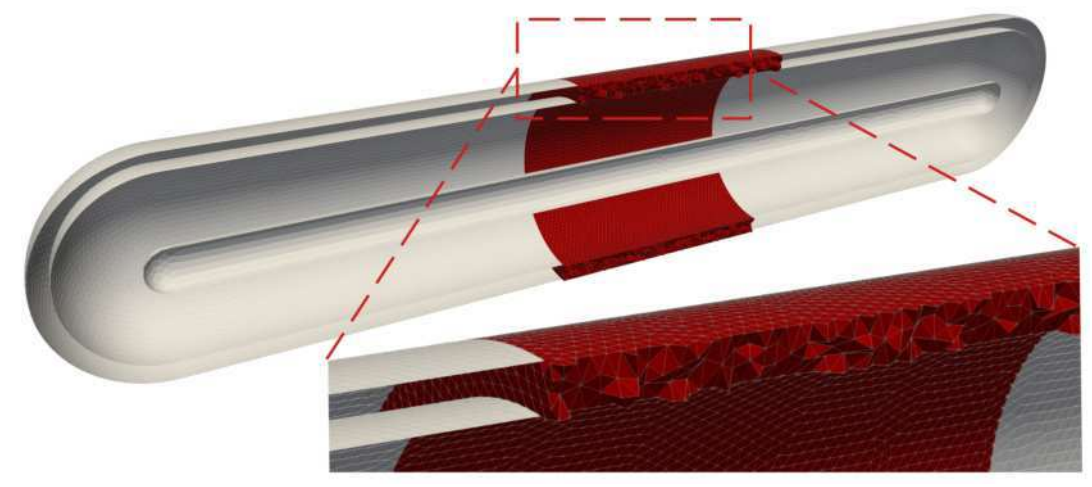

Figure 1: A half-object mesh handled by the present solver. The red region is occupied by the inhomogeneous material.

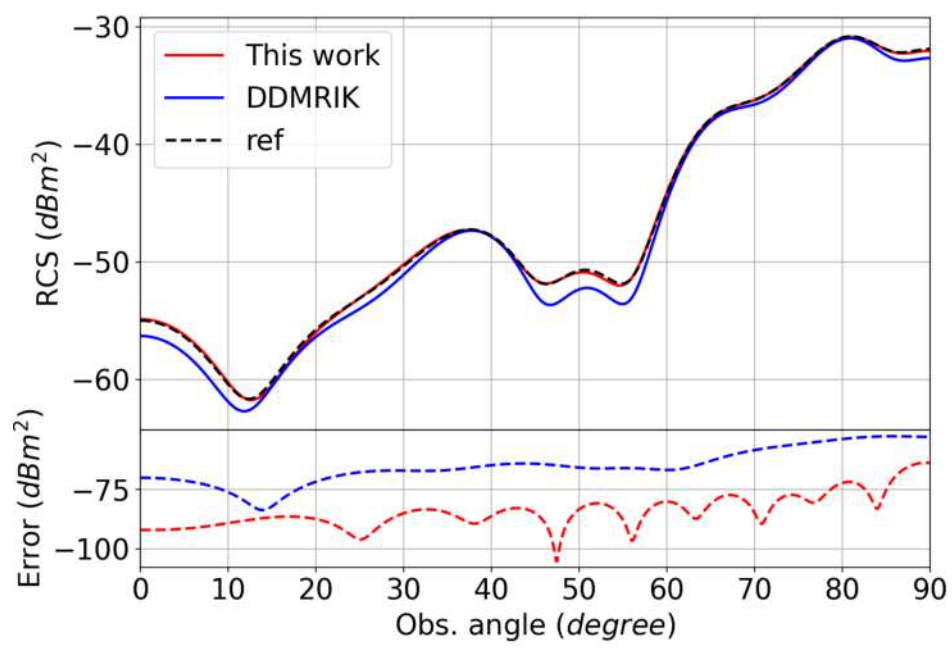

Figure 2: RCS and its error (with respect to the reference) versus monostatic angle.

| | This work ($\epsilon_{\text{threshold}} = 10^{-6}$) | DDMRIK |
|---|---|---|
| $h_{mean}$ (m) | $\pi/(24\kappa)$ | $\pi/(20\kappa)$ |
| # dofs | 549 571 | 2 318 992 |

Table 1: Simulation data.

# Optimization of EM-Chiral Scatterers

**Tilo Arens**[1,*], **Marvin Knöller**[2], **Raphael Schurr**[1]
[1]Department of Mathematics, Karlsruhe Institute of Technology (KIT), Germany
[2]Department of Mathematics and Statistics, University of Helsinki, Finland
*Email: tilo.arens@kit.edu

**Abstract**

Electromagnetic chirality is a quantity that describes a scatterer's different interaction with electromagnetic waves of opposite circular polarization. We optimize the shape of scatterers in order to maximize this quantity. We consider a specific class of scatterers, defined by a spine curve and a radius function along that curve. The shape optimization is performed by a BFGS scheme based on the domain derivative of the far field operator associated with the scatterer. We prove that this derivative exists and numerically calculate it and the corresponding derivative of the measure of EM chirality. In each iteration step, a large number of Maxwell problems have to be solved. We present a number of examples of optimized scatterers with non-intuive shapes that attain up to 90% of the theoretical maximum value the measure of EM chirality.

## 1 Introduction

Chirality is a phenomenon encountered in many applications. Geometrically, an object is called *chiral,* if it cannot be superimposed on its mirror image. In an electromagnetic scattering problem, the vectorial nature of the wave field implies that under a reflection the electromagnetic field changes its helicity, i.e. its direction of circular polarization. *Electromagnetic chirality (em-chirality)* is the phenomenon that the interaction of an object with fields of one helicity cannot be reproduced by the interaction of its mirror image with fields of the opposite helicity.

It turns out that there is a natural way to quantify the discrepancy between both configurations via the scattering or far field operator associated with the object, leading to so-called *measures of em-chirality.* Such measures are bounded by the Hilbert-Schmidt norm of the far field operator. Objects with an em-chirality measure very close to its maximum value are nearly invisible to fields of one helicity, making them desirable in many applications. We will demonstrate how the shape of a certain class of obstacles may be optimized in order to maximize their measure of em-chriality.

## 2 Long tubular objects

Previous research in Physics focussed on metallic scatterers of a size much smaller than the wavelength. In this case, helical shapes have turned out to produce (close to) maximal em-chiral objects. Thus, it seems natural to also consider helicies of a diameter similar to the wavelength as candidates for objects of maximal em-chirality. Although helices of very high em-chirality measure have been found, these objects turn out to be exceedingly hard to fabricate due to the very thin wire that forms the helix [1, 3, 4]. Hence, we consider tubular objects of larger and longitudinally varying thickness [2, 5]. Precisely, given a $C^1$-smooth curve $\boldsymbol{z} : [0,1] \to \mathbb{R}^3$, we denote by $(\boldsymbol{t}, \boldsymbol{n}, \boldsymbol{b})$ a continuous accompanying orthogonal frame. Also given a function $r : [0,1] \to \mathbb{R}_{>0}$, we consider curved tubes paramterized by

$$\begin{aligned} \boldsymbol{x}_{\text{body}}(\tau,\varphi) &= \boldsymbol{z}(\tau) + r(\tau)\,\boldsymbol{\zeta}(\tau,\varphi) \\ \text{where } \boldsymbol{\zeta}(\tau,\varphi) &= \cos(\varphi)\,\boldsymbol{n}(\tau) + \sin(\varphi)\,\boldsymbol{b}(\tau)\,, \end{aligned}$$

for $\tau \in (0,1)$ and $\varphi \in (-\pi,\pi]$. At both ends, the tube is closed by a cap that renders the entire surface $C^1$-smooth. A careful analysis yields explicit expressions for the first order effect on this parametrization caused by perturbing $\boldsymbol{z}$ and $r$. This forms the basis of computable domain derivatives discussed in the next section [2].

## 3 The domain derivative of the far field

We here limit the presentation to perfectly conducting objects but remark that homogeneous penetrable objects may be handled analogously. The scattered field $(\boldsymbol{E}^s, \boldsymbol{H}^s)$ produced by an incident field $(\boldsymbol{E}^i, \boldsymbol{H}^i)$ hitting a perfectly conducting obstacle $D$ with outward drawn unit normal $\boldsymbol{\nu}$ can be computed as the weak solution to the exterior boundary value problem

$$\left.\begin{aligned} \mathbf{curl}\,\boldsymbol{E}^s - \mathrm{i}\omega\mu_0\,\boldsymbol{H}^s &= 0 \\ \mathbf{curl}\,\boldsymbol{H}^s + \mathrm{i}\omega\varepsilon_0\,\boldsymbol{E}^s &= 0 \end{aligned}\right\}\quad \text{in } \mathbb{R}^3\setminus\overline{D}\,,$$

$$\boldsymbol{E}^s\times\boldsymbol{\nu} = -\boldsymbol{E}^i\times\boldsymbol{\nu}\qquad \text{on } \partial D$$

together with the Silver-Müller radiation condition. By the asymptotics of the scattered field, there is an associated far field pattern $\boldsymbol{E}^\infty \in L^2_t(\mathbb{S}^2)$. Given some compactly supported $C^1$-smooth vector field $\boldsymbol{h}$ one obtains a perturbed domain $D_{\boldsymbol{h}}$ and an associated far field $\widehat{\boldsymbol{E}}^\infty_{\boldsymbol{h}}$. The domain derivative of the far field is a function $(\boldsymbol{E}')^\infty \in L^2_t(\mathbb{S}^2)$, such that

$$\left\|\widehat{\boldsymbol{E}}^\infty_{\boldsymbol{h}} - \boldsymbol{E}^\infty - (\boldsymbol{E}')^\infty\right\|_{L^2_t(\mathbb{S}^2)} \le C\,\|\boldsymbol{h}\|^2_{1,\infty}\,\|\boldsymbol{E}^i\|_{H(\mathbf{curl},\Omega)}\,.$$

It is the far field of the solution of the above boundary value problem, but with the perfect conductor boundary condition replaced by

$$\boldsymbol{E}'\times\boldsymbol{\nu} = \mathrm{i}\,k\,h_{\boldsymbol{\nu}}\,\boldsymbol{\nu}\times(\boldsymbol{H}\times\boldsymbol{\nu}) - \mathbf{Curl}\,(h_{\boldsymbol{\nu}}E_{\boldsymbol{\nu}})$$

Here, the index $\boldsymbol{\nu}$ indicates normal components of the corresponding vector fields.

We prove that the domain derivative of the far field for restricted to long tubular domains is obtained by replacing $\boldsymbol{h}$ in this boundary condition with the corresponding first order pertubations obtained in the previous section.

## 4 Shape optimization

We use the chirality measure defined in [1] which is based on spectral properties of the far field operator $\mathcal{F}$ associated with the scatterer $D$. One of the central results in [5] is the Fréchet differentiability of $\mathcal{F}$ in the Hilbert-Schmidt norm with respect to perturbations of $D$:

$$\|\mathcal{F}_{D_{\boldsymbol{h}}} - \mathcal{F}_D - \mathcal{F}'_D\boldsymbol{h}\|_{\mathrm{HS}} \le C\,\|\boldsymbol{h}\|^2_{1,\infty}.$$

Based on this property, shape derivatives of the em-chirality functional may be computed and form the basis of a corresponding shape optimization scheme. Suitable regularization terms have had to be formulated to stabilize the optimization. Using the BFGS algorithm, a number of concrete examples were carried out, one of which is shown in Figure 1. Other examples will be shown in the presentation.

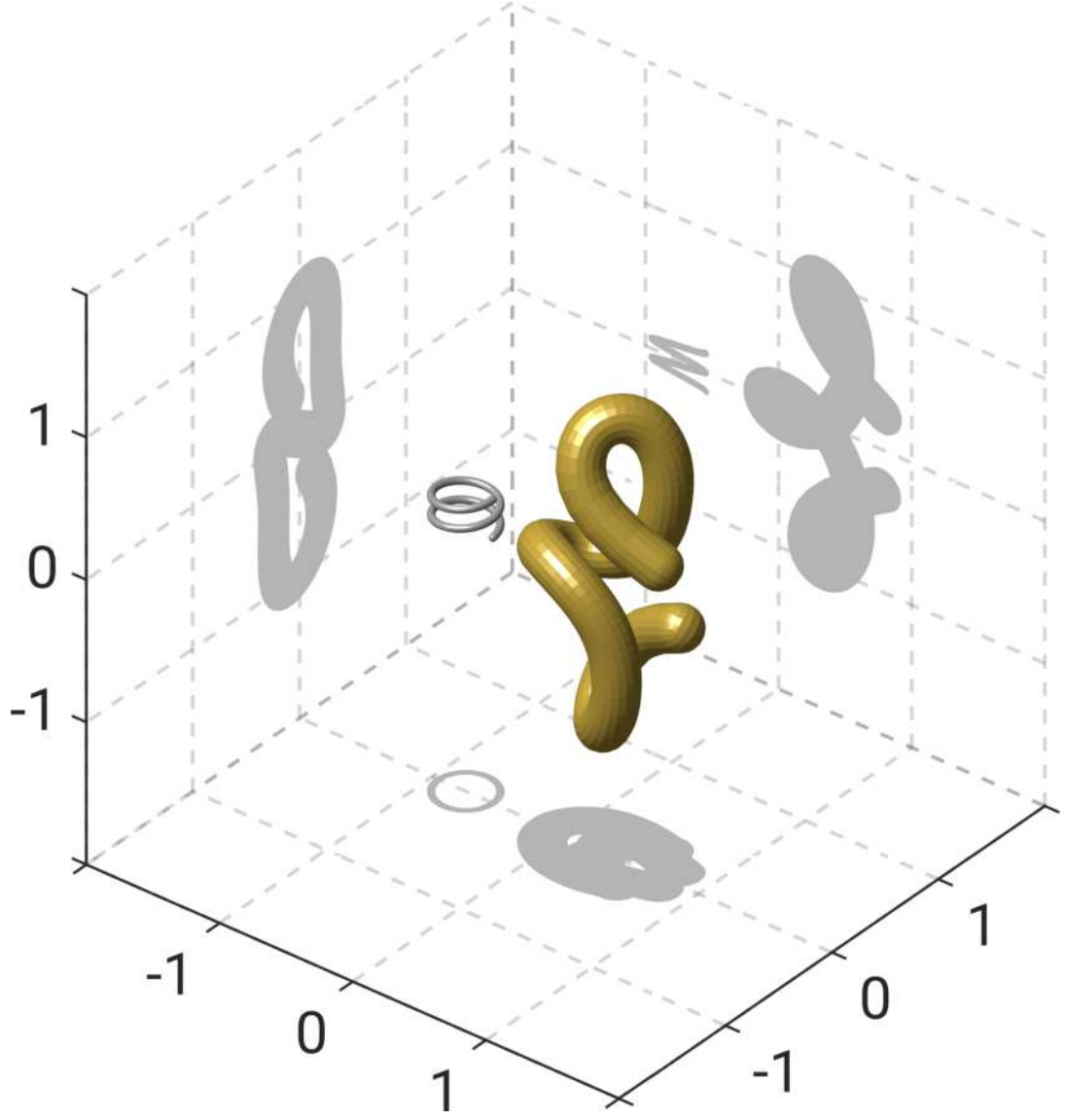

Figure 1: In scale comparison of a highly em-chiral scatterer produced by our proposed algorithm and a silver helix with similar properties.

In each iteration step, it is required to numerically solve a large number of Maxwell problems. This is done using the boundary element method. We spend particular effort on efficiency when solving the boundary integral equations. A discussion of some of the implementation steps will be included in the presentation.

# Fast assembly of H-matrix and FMM-based boundary operators for time-harmonic Maxwell equations in a frequency-sweep context

**Clément Aumonier**[1,*], **David Levadoux** [1]

[1]DTIS, ONERA, Université de Toulouse, Toulouse 31000, France

*Email: clement.aumonier@onera.fr

**Abstract**

We present a method to construct discrete boundary integral operators for time-harmonic Maxwell equations in a frequency-sweep context, ensuring compatibility with both FMM and $\mathcal{H}$-matrix compression. First, the frequency-independent operations to build the elemental matrix are factorized. We then use the EIM [1] combined with a non-intrusive procedure [2] to subsample the required number of frequencies while ensuring error control. The result is a significant acceleration of BEM operators constructions for multiple frequency points simultaneously. Its performance is further validated through an integration into an RBM framework.

## 1 Introduction

Time-harmonic electromagnetic simulations commonly use BEM with RBM for frequency sweeping [2]. However, repeatedly assembling the fully populated linear system represents a high CPU cost, so we aim to address this bottleneck. Similar objectives have been pursued in [4–6] with a different techniques from ours. Consider a wave with wavenumber $\kappa$ propagating in free space and scattered by an object $\Gamma$. The scattered field is obtained by solving an integral equation, which formulation with electric field and a perfect electrical conductor material for $\Gamma$ is given by $T(\kappa)J = -E^{inc}_{|\tan}$ with :

$$\begin{aligned}\forall x \in \Gamma,\, T(\kappa)J(x) &= \int_\Gamma G(x,y)J(y)\, d\Gamma_y \\ &+ \frac{1}{\kappa^2}\nabla_x \int_\Gamma G(x,y)\nabla_y \cdot J(y)\, d\Gamma_y\end{aligned} \tag{1}$$

and $G(x,y) = e^{i\kappa||x-y||}/4\pi||x-y||$ the Green kernel. A suitable mesh $\mathcal{T}_h$ is built on $\Gamma$, and is associated with a set of basis functions $\{\phi_i\}_{i=1}^{\mathcal{N}}$. Using a Galerkin method, the fully populated linear system $\mathbf{A}(\kappa)\mathbf{x}(\kappa) = \mathbf{b}(\kappa)$ is built by taking the entries of the matrix $\mathbf{A}(\kappa)$ as :

$$\mathbf{A}_{ij}(\kappa) = \langle T(\kappa)\phi_i, \phi_j\rangle \tag{2}$$

This discrete representation of the BEM operator is compressed using the FMM, allowing for fast convolution, or with $\mathcal{H}$-matrices for both fast convolution and direct inversion. These compression methods involve an expensive assembly step.

## 2 Multi-build frequency

We seek to reach the BEM operator $\forall\kappa \in \mathcal{D}$, where $\mathcal{D} \subset \mathbb{R}_+$ is an arbitrary interval. The multi-build under frequency consists in a first step to factorize the different operations required for computing the Galerkin matrix entries. For a discrete surrogate $\Xi = \{\kappa_1, ..., \kappa_{N_f}\}$ of $\mathcal{D}$, the data are represented using a third dimension (the frequency), and the internal procedure building a Galerkin matrix take the form : $\mathbf{A}_{i,j,f} \leftarrow \mathbf{A}_{i,j,f} + \epsilon^{K,L}_{m,n} EM^{K,L}_{m,n,f} \quad 1 \le f \le N_f$. Here $EM$ is the elemental matrix associated with the two elements $K$ and $L$, $m$ and $n$ refer to the local numeration of global unknowns $i$ and $j$. A significant CPU saving is realized by using a minimum time every finite element connection, frequency independent operations, and mainly, the calculus of Green kernel can be vectored. This strategy improves computational efficiency by avoiding the independent assembly of BEM operators at each frequency. This first consideration can be applied to FMM with :

- radiation and reconstruction operator

$$\mathcal{T}^{rad}_{\mathcal{P}a} J(s) = \int_{\Gamma\cap\mathcal{P}a} e^{i\kappa s\cdot(x - C_{\mathcal{P}a})} J(x)dx$$

$$\mathcal{T}^{rec}_{\mathcal{P}a} \mathcal{G}(y) = \int_{\mathbb{S}^2} e^{i\kappa s\cdot(y - C_{\mathcal{P}a})} \mathcal{G}(s)ds$$

- near-field interactions $\mathbf{A}^{close}$.

It also can be applied to $\mathcal{H}$-matrices with :

- rows $(b^l)^T$ and columns $(a^l)$ selected for low rank approximation $M \approx \sum_{l=1}^{k} a^l (b^l)^T$ of admissible blocks $M \in P^+_{I\times J}$
- non-admissible full blocks $M \in P^-_{I\times J}$.

In the second step, the Empirical Interpolation Method [1] is used to obtain a tensorized representation of the Green kernel (or of a variant that gathers all the $\kappa$-dependence in Eq. 1). This tensorization separates the variables between space, $r = ||x - y|| \in \mathcal{R}$, and frequency, $\kappa \in \mathcal{D}$, by selecting iteratively the so-called magic points. Once reflected on the whole BEM operator, it leads to an interpolation said weakly intrusive. Therefore, the EIM is coupled with the non-intrusive procedure [2] to obtain an interpolation in the form of Eq. 3, for all $\kappa \in \mathcal{D}$.

$$\mathbf{A}(\kappa) \approx \sum_{p}^{Q} \theta_p(\kappa)\mathbf{A}(\hat{\kappa}_p) \quad (3)$$

When integrated into the core FMM procedures, this approach yields fewer interpolation points due to the reduced domain spaces $\mathcal{R}$, as detailed below. We denote by $R_L$ the radius of the the leaf level of the FMM octree. So the space intervals are set as:

- rad/rec operator, $\mathcal{R} = [0; 2\sqrt{3}R_L]$,
- close part, $\mathcal{R} = [0; 6\sqrt{3}R_L]$.

In a same way, reduced space intervals can be considered for the assembly of non-admissible and admissible blocks of an $\mathcal{H}$-matrix.

## 3 Numerical results

We compare the two steps of the proposed method and analyze the speed-up obtained with respect to assembling all desired frequencies individually. The figure 1 shows the speed-up versus the number of frequencies, when assembling the close part of FMM on the unit sphere of 127008 unknowns. This experiment shows the savings made in the elemental matrix $EM$ builder. The three curves compared are the first step of the method without kernel computation vectorization (MBF I), with kernel computation vectorization (MBF II), and adding the EIM interpolation (MBF-EI). Each curve follows a trend that confirms their respective asymptotic complexities. Of course, the number of computed

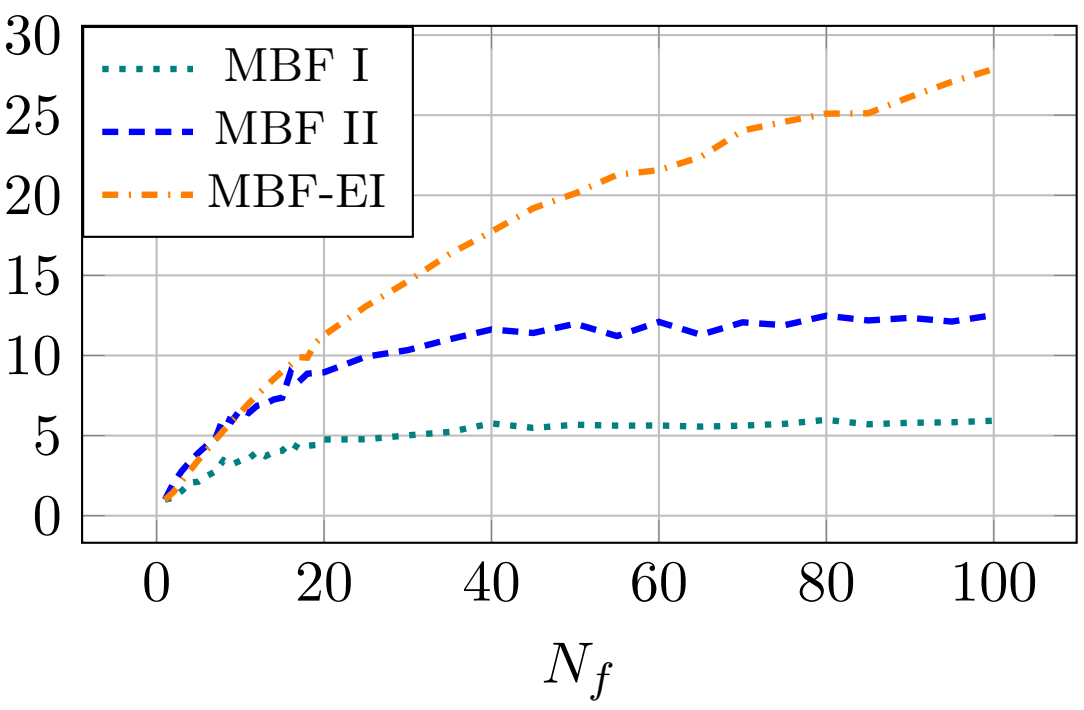

Figure 1: Speed-up of the FMM near-field on the unit sphere.

frequencies varies for each test case, and in practice, a speed-up of 20 is never achieved.

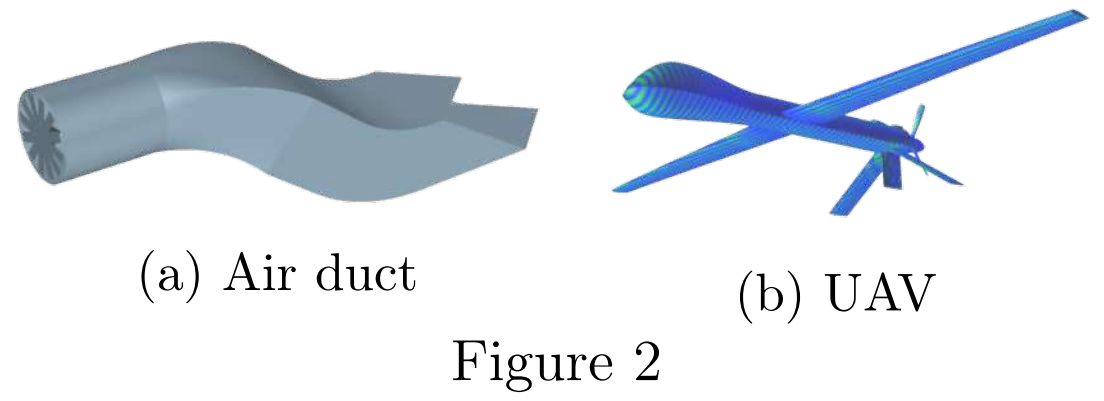

(a) Air duct (b) UAV

Figure 2

However more results will be discussed in the presentation where speed-ups of about 4-5 are reached. We will test it on an air duct (see fig. 2a) representative of industrial stealthiness problems issued from [3], and the relevancy of this approach in a RBM context for a drone target (fig. 2b).

# Convergence of the WaveHoltz iteration on $\mathbb{R}^d$

Elliot Backman[1,*], Olof Runborg[1]
[1]Department of Mathematics, Royal Institute of Technology, Stockholm, Sweden
*Email: elliotba@kth.se

**Abstract**

The WaveHoltz method is an iterative time-domain method for numerically solving the Helmholtz equation. We define and study the WaveHoltz iteration on the entirety of $\mathbb{R}^d$. In the constant coefficients case we show that the difference between the $n$th WaveHoltz iterate and the outgoing solution to the Helmholtz problem again satisfies a Helmholtz equation with a modified forcing. Careful analysis of this equation produces a frequency-explicit estimate of the error in the $n$th WaveHoltz iterate against the real part of the outgoing solution. In particular, this error decays as $n^{-1/2}$.

## 1 Introduction

Finding numerical solutions to the Helmholtz equation is a challenging problem, especially in the high-frequency case. Discretising the Helmholtz equation generally produces large indefinite linear systems which are costly to solve with standard Krylov methods. One approach to this problem is to consider time-domain methods [1, 2], which have the benefits of being implementable in a parallel and memory-lean manner with relative ease. The WaveHoltz method is one such iterative time-domain method. Previous work on the WaveHoltz method has shown that the method converges for interior problems with Dirichlet or Neumann boundary conditions [2], and numerical evidence suggests that the method also converges in more general cases for which proofs have yet to be found. In this contribution we present a convergence result for the open problem set in the entire space $\mathbb{R}^d$

$$\Delta u + \omega^2 u = f, \qquad x \in \mathbb{R}^d, \tag{1}$$

subject to the radiation condition

$$\lim_{|x|\to\infty} |x|^{\frac{d-1}{2}}\big(\partial_r - i\omega\big)u(x) = 0.$$

## 2 The WaveHoltz Iteration

The WaveHoltz method revolves around the operator $\Pi$, defined for compactly supported smooth functions $\varphi$ as

$$\Pi\varphi = \frac{2}{T}\int_0^T \Big(\cos(\omega t) - \frac{1}{4}\Big) w(x,t)\,dt,$$

where $T = 2\pi/\omega$ and the function $w(x,t)$ is the solution of the associated initial-value problem

$$\begin{aligned} &\partial_t^2 w = \Delta w - f(x)\cos(\omega t), && x \in \mathbb{R}^d, t > 0, \\ &w(x,0) = \varphi(x), && x \in \mathbb{R}^d, \\ &\partial_t w(x,0) = 0, && x \in \mathbb{R}^d. \end{aligned} \tag{2}$$

We then define the WaveHoltz method by the fixed-point iteration

$$u^{n+1} = \Pi u^n, \ n \geq 0, \quad u^0 \equiv 0,$$

motivated by the fact that this operator at least formally has the solution to (1) as a fixed point. Note that the WaveHoltz iterates are real by construction, but that the outgoing solution to (1) in general is complex-valued. To analyse the problem set in $\mathbb{R}^d$ we need to extend $\Pi$ to the weighted $L^p$ space in which the solution to (1) lies, [3]. Such spaces, $L^p_s(\mathbb{R}^d)$, are defined for $s \in \mathbb{R}$ as the space of all functions in $L^p(\mathbb{R}^d)$ for which the norm $\|f\|_{L^p_s(\mathbb{R}^d)} = \|\langle x\rangle^s f(x)\|_{L^p(\mathbb{R}^d)}$ is finite, where $\langle x\rangle = (1+|x|^2)^{\frac{1}{2}}$. The extension of $\Pi$ into $L^p_s(\mathbb{R}^d)$ makes use of a density argument and the finite speed of propagation of wave solutions. It holds for any choice of $s$ and positive integer $p$, provided that $f$ is in $L^2_t(\mathbb{R}^d)$ for some $t > 1/2$.

It can then be shown that the differences $e^n := u^n - u$ between the iterates and the outgoing solution $u$ to (1) are themselves outgoing solutions to the problems

$$\Delta e^n + \omega^2 e^n = -\mathcal{S}^n f, \tag{3}$$

where $\mathcal{S}$ is defined analogously to the operator $\Pi$, but with $f = 0$ in (2). Moreover, this operator can be represented by use of the Fourier transform $\mathcal{F}$. We have

$$(\mathcal{F}\mathcal{S}^n\varphi)(\xi) = \beta^n(|\xi|/\omega)(\mathcal{F}\varphi)(\xi), \tag{4}$$

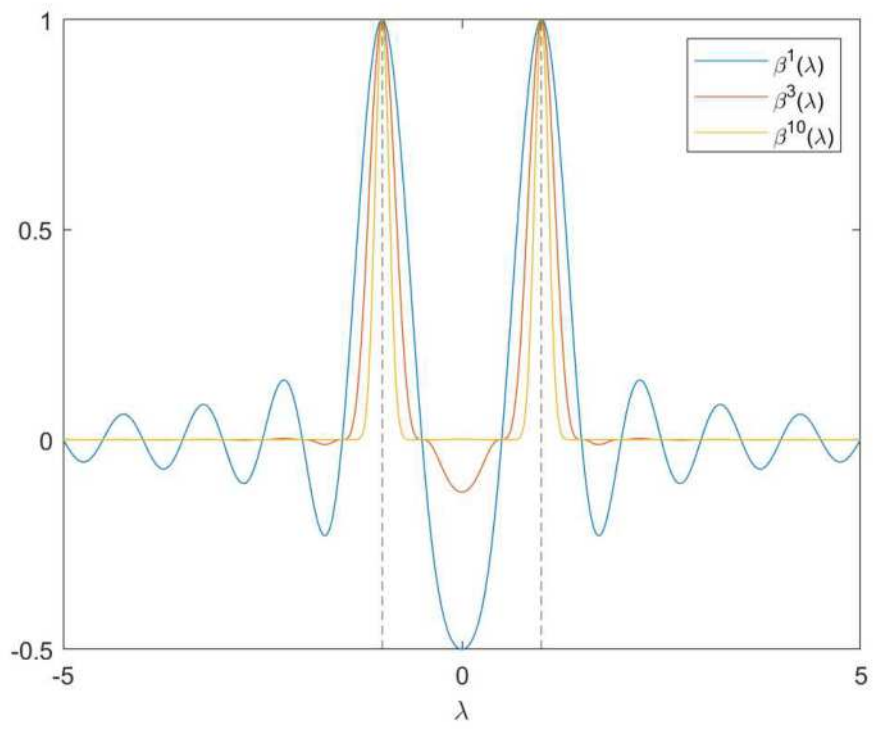

Figure 1: Some powers of the function $\beta(\lambda)$, the lines $\lambda = 1$ and $\lambda = -1$ being shown in black.

where $\beta$ is given by

$$\beta(\lambda) = \frac{1}{\pi}\int_0^{2\pi}\left(\cos(t) - \frac{1}{4}\right)\cos(\lambda t)dt.$$

The proof of our convergence result relies on the properties of the function $\beta$, perhaps most important of which for us being that $\beta^n(\lambda)$ approaches zero for every $\lambda$ away from $\pm 1$ as $n$ grows larger. Our approach is now to estimate in a suitable norm the solution to (3) in terms of another norm of the forcing $-\mathcal{S}^n f$, and to use the pointwise decay of the forcing to demonstrate that the norm of the forcing converges to zero with the iteration number $n$. It turns out that the standard estimates of the solution in terms of the forcing found for instance in [3] are not quite sharp enough for our analysis, which motivates a closer investigation.

## 3 A Convergence Result

We analyse the damped Helmholtz equation

$$\Delta u + (\omega^2 + i\omega\alpha)u = G \tag{5}$$

where $\alpha > 0$. The limiting-absorption principle guarantees that the solution to the damped problem approaches the outgoing solution to (1) as $\alpha \to 0^+$; see [3]. Inspired by (4), we assume that the forcing $G$ satisfies

$$(\mathcal{F}G)(\xi) = h(|\xi|/\omega)(\mathcal{F}g)(\xi),$$

for some functions $h \in L^1(\mathbb{R}^d) \cap L^\infty(\mathbb{R}^d)$, and $g \in L^2_s(\mathbb{R}^d)$ with some $s > 3/2$. Moreover we require the technical assumption that

$$\int_0^\delta \left|\frac{h(1+x) - h(1-x)}{2x}\right| dx < \infty \tag{6}$$

for some $\delta \in (0,1)$. Under these assumptions we show that the solution to (5) satisfies

$$\|\mathfrak{Re}\{u\}\|_{L^2_{-s}(\mathbb{R}^d)} \leq C\omega^{2s-2}|||h|||\,\|f\|_{L^2_s(\mathbb{R}^d)}. \tag{7}$$

for a particular norm $|||\cdot|||$ involving a term of the type (6), and some constant $C$ that is independent of $\alpha$. This result is achieved by proving the estimate for $\omega = 1$ and then rescaling it to find (7) for $\omega > 1$. Moreover, as the constants in our estimate are independent of $\alpha$ we can use the limiting-absorption principle to extend (7) to outgoing solutions to (1) as well ($\alpha = 0$).

The estimate (7) is finally applied to (3) with $h = \beta$. We verify that (6) holds and find that $|||\beta^n||| \sim n^{-1/2}$. This results in

**Theorem 1** *Suppose $s > \frac{3}{2}$ and $f \in L^2_s(\mathbb{R}^d)$. Then for $\omega \geq 1$ the real part of the WaveHoltz iterates $u^n$ converge in $L^2_{-s}(\mathbb{R}^d)$-norm to the real part of the outgoing solution $u$ to (1), with the error $e^n$ at the nth iteration satisfying*

$$\|\mathfrak{Re}\{e^n\}\|_{L^2_{-s}(\mathbb{R}^d)} \leq C\omega^{2s-2}n^{-\frac{1}{2}}\|f\|_{L^2_s(\mathbb{R}^d)},$$

*where $C$ is independent of $f$, $u$, $\omega$ and $n$.*

In other words, under the assumptions of the theorem we find that the WaveHoltz iterates converge to the real part of the outgoing solution to (1). To also find the imaginary part of the solution, another type of iteration is needed.

# Revisiting Artificial Boundary Conditions Through Modal Expansion

**Hélène Barucq**[1,*], **Florian Faucher**[1], **Olivier Lafitte**[2], **Ha Pham**[1]
[1]Team Makutu, Inria Bordeaux, University of Pau and Pays de l'Adour, CNRS UMR 5142, France.
[2]University Sorbonne Paris Nord, LAGA, CNRS UMR 7539, France.
*Email: helene.barucq@inria.fr

**Abstract**

We revisit the question regarding the contribution of the transition and glancing region in approximating the exterior Dirichlet-to-Neumann operator (DtN) of the Helmholtz operator for circular curves. This symmetry provides, via Fourier series transform, analytical expressions for the DtN operator, its nonlocal approximations via Debye and Airy-type expansions, as well as input and output traces. Our region-specific modal error analysis compares these nonlocal approximations with low-order transparent boundary conditions containing no tangential derivatives commonly employed in numerical computations with finite elements.

## 1 Introduction

In computing outgoing-at-infinity solutions to the Helmholtz equation with finite element methods, a transparent boundary condition (TBC) imposed on an artificial boundary is commonly employed to truncate the infinite domain. These conditions are constructed by approximating the exterior DtN, and in most cases, only the hyperbolic region (HR) in the phase space on local coordinate charts is considered, ignoring the transition (TR) and evanescent region (ER), cf. [2].

In this work, we revisit the question regarding the contribution of these regions for circular symmetry. This setting provides, via Fourier series transform, analytical expressions in terms of Bessel functions of the input traces, the DtN, the exact outgoing solutions, as well as those obtained with TBCs. We also employ approximations coming from the Debye and Airy-type expansion, cf. [1]. Our numerical investigation first studies the contribution of each region in input Dirichlet traces. Secondly, the modal analysis of error in each region is investigated for the hard-scattering of plane waves by circular obstacles. We compare among the first-order Debye and Airy-type TBCs, as well as the Sommerfeld and the first-order RBC.

The construction of TBCs treating the glancing set and evanescent region were carried out in [3, 4]. These conditions contain tangential derivatives [3] and fractional ones [4]. Our work focuses on those with no tangential derivatives.

## 2 DtN on a circle and its approximations

Consider the Helmholtz operator,

$$-\Delta - \mathsf{k}^2, \qquad \text{with constant } \mathsf{k} > 0,$$

on a disc centered at the origin $\mathbb{B}_\mathfrak{r}$ of radius $\mathfrak{r}$, with boundary $\Gamma$. The exterior Dtn on $\Gamma^+$ maps a Dirichlet trace to the Neumann trace of an outgoing solution: for $\mathfrak{t} \in H^{1/2}(\Gamma)$, $\Lambda^+(\mathfrak{t}) = \partial_\nu u_+$, where $u_+$ is outgoing and uniquely solves:

$$(-\Delta - \mathsf{k}^2)u^+ = 0 \text{ on } \mathbb{R}^2 \setminus \overline{\mathbb{B}_\mathfrak{r}}, \text{ with } u^+|_\Gamma = \mathfrak{t}.$$

The Fourier series on $\Gamma$ is denoted as,

$$\mathfrak{t}(\theta) = \sum_{p\in\mathbb{Z}} [\mathfrak{t}]_p \, e^{\mathrm{i}p\theta}, \text{ where } [\mathfrak{t}]_p = \frac{1}{2\pi}\int_0^{2\pi} \mathfrak{t}(\tilde\theta) e^{-\mathrm{i}p\tilde\theta} d\tilde\theta.$$

We work with special functions and their derivative, cf. [1, Ch. 9, 10]: the Hankel $\mathrm{H}_p^{(1)}$, the Bessel $\mathrm{J}_p$ functions, and $\mathcal{A}(t) := \frac{\mathrm{Ai}'(t) - \mathrm{i}\,\mathrm{Bi}'(t)}{\mathrm{Ai}(t) - \mathrm{i}\,\mathrm{Bi}(t)}$ with Airy functions, $\mathrm{Ai}(\cdot), \mathrm{Bi}(\cdot)$.

Introduce variables $\tau$ and $\zeta$ with,

$$\tau := p/(\mathsf{k}\mathfrak{r}), \qquad \zeta(|\tau|^{-1}) \text{ see [1, line 10.20]}. \quad (1)$$

The phase space $(\mathsf{k}\mathfrak{r}, p)$ is separated into purely imaginary and real square roots of $\zeta$ by the turning points of the Bessel equation, we have

$$(\mathrm{HR}) : \{\mathsf{k}\mathfrak{r} > |p|\} = \{|\tau| < 1\} = \{\zeta < 0\},$$
$$(\mathrm{ER}) : \{\mathsf{k}\mathfrak{r} < |p|\} = \{|\tau| > 1\} = \{\zeta > 0\},$$

separated by the glancing set

$$\{p = \mathsf{k}\mathfrak{r}\} = \{|\tau| = 1\} = \{\zeta = 0\}.$$

The Fourier transform conjugates the DtN into a diagonal multiplication operator, with

$$[\Lambda\,\mathfrak{t}]_p = \Lambda_p\,[\mathfrak{t}]_p, \quad \text{where } \Lambda_p := \mathsf{k}\frac{\mathrm{H}_p^{(1)\prime}(\mathsf{k}\mathfrak{r})}{\mathrm{H}_p^{(1)}(\mathsf{k}\mathfrak{r})}. \quad (2)$$

Commonly employed in literature are local low-order approximations of $\Lambda$ such as, cf. [2],

$$\mathcal{I}^{\text{som}} = \mathsf{ik}, \qquad \mathcal{I}^{\text{c}} = \mathsf{ik} - 1/(2\mathfrak{r}). \tag{3}$$

Away from the glancing set, the leading part of the Debye expansion, cf. [1, 10.19], is given by,

$$\mathcal{I}_p^{\text{Deb}} = \mathsf{k} \begin{cases} \mathrm{i}\sqrt{1-\tau^2}, & |\tau| < 1, \\ (-1)\sqrt{\tau^2 - 1}, & |\tau| > 1. \end{cases} \tag{4}$$

We thus retrieve the principal symbol of $\Lambda$ (away from glancing), in the class of pseudo-differential operators on $\Gamma$, cf. [3, Eq. (2.19)].

On the other hand, employing $\zeta$ which is smooth in $\tau^{-1}$ and well-defined across the glancing set, the leading term in the uniform expansion [1, 10.20] inherits this property,

$$\mathcal{I}_p^{\text{UO}} = -\mathsf{k}\, \frac{2\,\tau}{p^{1/3}} \left(\frac{1-\tau^{-2}}{4\zeta}\right)^{1/2} \mathcal{A}(p^{2/3}\zeta). \tag{5}$$

## 3 Numerical experiments

We consider the scattering of a planewave $\mathfrak{u}_{\text{pw}} = e^{\mathsf{i}\mathsf{k}\mathbf{x}\cdot(\cos\theta_{\text{inc}},\sin\theta_{\text{inc}})}$ by a hard circular obstacle $\mathbb{B}_{\text{obs}}$. The truncated problem for the scattered wave writes as:

$$(-\Delta - \mathsf{k}^2)u^{\mathcal{I}} = 0 \text{ on } \mathbb{B}_{\mathfrak{r}} \setminus \overline{\mathbb{B}_{\text{obs}}},$$
$$u^{\mathcal{I}} = -\mathfrak{u}_{\text{pw}} \text{ on } \partial\mathbb{B}_{\text{obs}},\ (\partial_n - \mathcal{I})u^{\mathcal{I}}_{\text{sc}} = 0 \text{ on } \Gamma.$$

The solutions are written as Fourier expansions. The $p$-th coefficients on $\Gamma$ of the reference solution $u^+$ and the TBC $u^{\mathcal{I}}$ are

$$[u^+]_p = -\mathrm{i}^p\, \frac{\mathrm{J}_p(\mathsf{k}\,\delta_{\text{obs}})}{\mathrm{H}_p^{(1)}(\mathsf{k}\,\delta_{\text{obs}})}\, \mathrm{H}_p^{(1)}(\mathsf{kr})\; e^{-\mathrm{i}p\,\theta_{\text{inc}}},$$
$$[u^{\mathcal{I}}]_p = -\mathrm{i}^p \left(a_p\, \mathrm{H}_p^{(1)}(\mathsf{kr}) \;+\; b_p\, \mathrm{H}_p^{(2)}(\mathsf{kr})\right) e^{-\mathrm{i}p\theta_{\text{inc}}},$$

where coefficient $a_p$, $b_p$ are constrained by the boundary conditions on $\partial\mathbb{B}_{\text{obs}}$ and $\Gamma$.

For the experiment, the obstacle is centered at the origin with size $\delta_{\text{obs}} = 0.25$, the wavenumber is $\mathsf{k} = 30$ and the trace is taken at $\mathsf{r} = 1$, see Fig. 1. The magnitude of the coefficients $[u^+]_p$ is shown in Fig. 2 for $\mathsf{kr} = 30$ and $\mathsf{kr} = 50$.

We then compare the relative difference $\mathfrak{e}_p$ between the analytic coefficients and the approximate impedance ones such that,

$$\mathfrak{e}_p \;=\; |\mathcal{I}^{\text{ref}} - \mathcal{I}^{\bullet}|/|\mathcal{I}^{\text{ref}}|\;, \tag{6}$$

where $\mathcal{I}^{\text{ref}}$ is the condition for the DtN and $\mathcal{I}^{\bullet}$ the approximate condition. In Fig. 3, we see that the condition $\mathcal{I}_p^{UO}$ is the most accurate near glancing, while $\mathcal{I}_p^{Deb}$ becomes the most accurate in the evanescent regions at higher $p$. Near $p = 0$, it is $\mathcal{I}_p^{c}$ that is the most accurate. From these results, we will construct TBC to be applied in a discretization scheme.

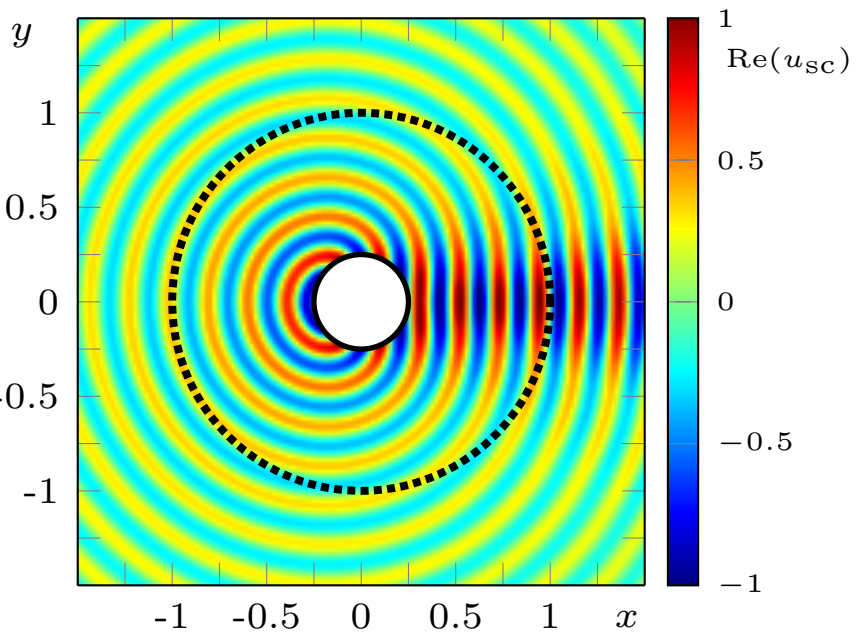

Figure 1: Scattered field of a planewave by a hard obstacle of radius $\delta_{\text{obs}} = 0.25$ with $\mathsf{k} = 30$, $\mathsf{r} = 1$ (dotted line).

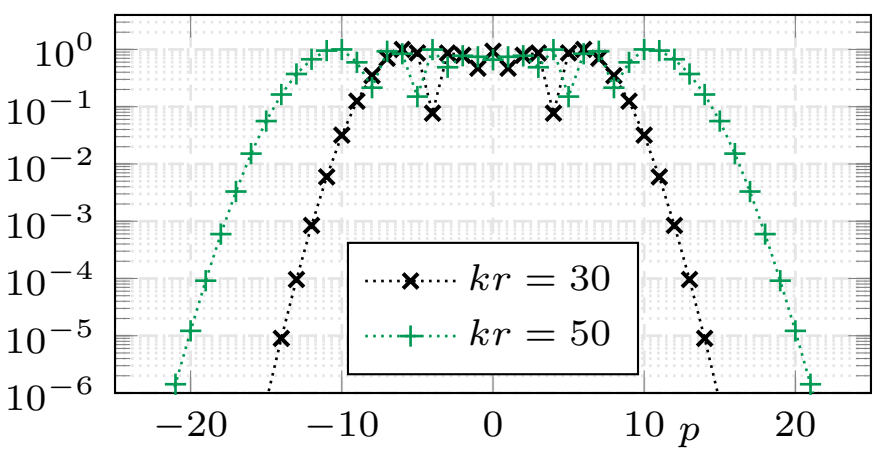

Figure 2: Normalized amplitude of the analytic modal coefficients $[u^+]_p$ for $\mathsf{kr} = 30$ and $\mathsf{kr} = 50$.

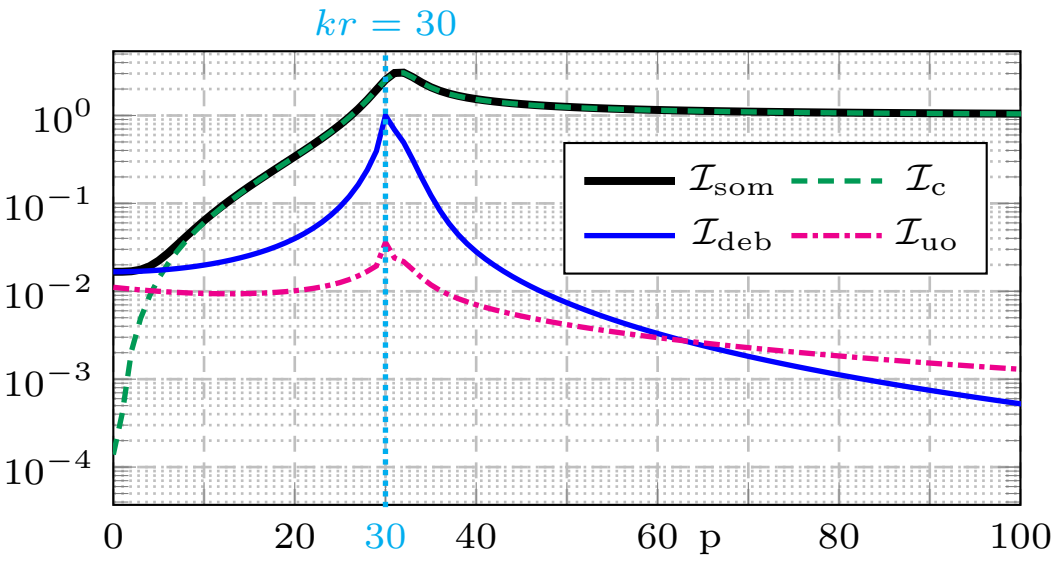

Figure 3: Relative difference $\mathfrak{e}_p$ between the analytic and approximated coefficients, $\mathsf{kr} = 30$.

# Effective models for Helmholtz's equation in a corner with a thin layer

**Cédric Baudet**[1,*], **Sonia Fliss**[1], **Patrick Joly**[1]
[1]POEMS, CNRS, Inria, ENSTA, Institut Polytechnique de Paris, 91120 Palaiseau, France
*Email: cedric.baudet@ensta.fr

**Abstract**

We consider the scattering by a (possibly locally perturbed) Dirichlet angular sector covered on one side by a thin layer. We provide effective models of order 3 and 5 w.r.t. the thickness of the layer when it tends to 0. These models replace the layer by effective boundary conditions and the corner by a transparent boundary condition taking into account the singular behavior of the solution.

## 1 The original problem

Let $\Omega := \{(r\cos\theta, r\sin\theta) \mid r \in \mathbb{R}_+^*, \theta \in (0,\Theta)\}$ with $\Theta \in (0, 2\pi)$ and $L := \mathbb{R}_+^* \times (-1, 0]$. The domain $\Omega_1 \subset \mathbb{R}^2$ is an open set that coincides with $\Omega \cup L$ outside of the disc $B(0, R_c)$ for some $R_c > 0$. Then let $\mu, \rho \in L^\infty(\Omega_1)$ be two functions bounded from below by positive constants, that are constant in $\Omega \backslash B(0, R_c)$ and $L \backslash B(0, R_c)$. See Figure 1. Now let $\varepsilon > 0$. The physical domain is $\Omega_\varepsilon := \{(x, y) \in \mathbb{R}^2 \mid (\frac{x}{\varepsilon}, \frac{y}{\varepsilon}) \in \Omega_1\}$ and the physical coefficients are given by $\mu_\varepsilon : (x, y) \in \Omega_\varepsilon \mapsto \mu(\frac{x}{\varepsilon}, \frac{y}{\varepsilon})$ and $\rho_\varepsilon : (x, y) \in \Omega_\varepsilon \mapsto \rho(\frac{x}{\varepsilon}, \frac{y}{\varepsilon})$. We consider $u_\varepsilon \in H_0^1(\Omega_\varepsilon)$ the unique solution of

$$\operatorname{div}(\mu_\varepsilon \nabla u_\varepsilon) + \omega^2 \rho_\varepsilon u_\varepsilon = f_s \quad \text{in } \Omega_\varepsilon, \qquad (1)$$

where $\Im(\omega) \neq 0$ and $f_s \in H^{-1}(\Omega_\varepsilon)$.

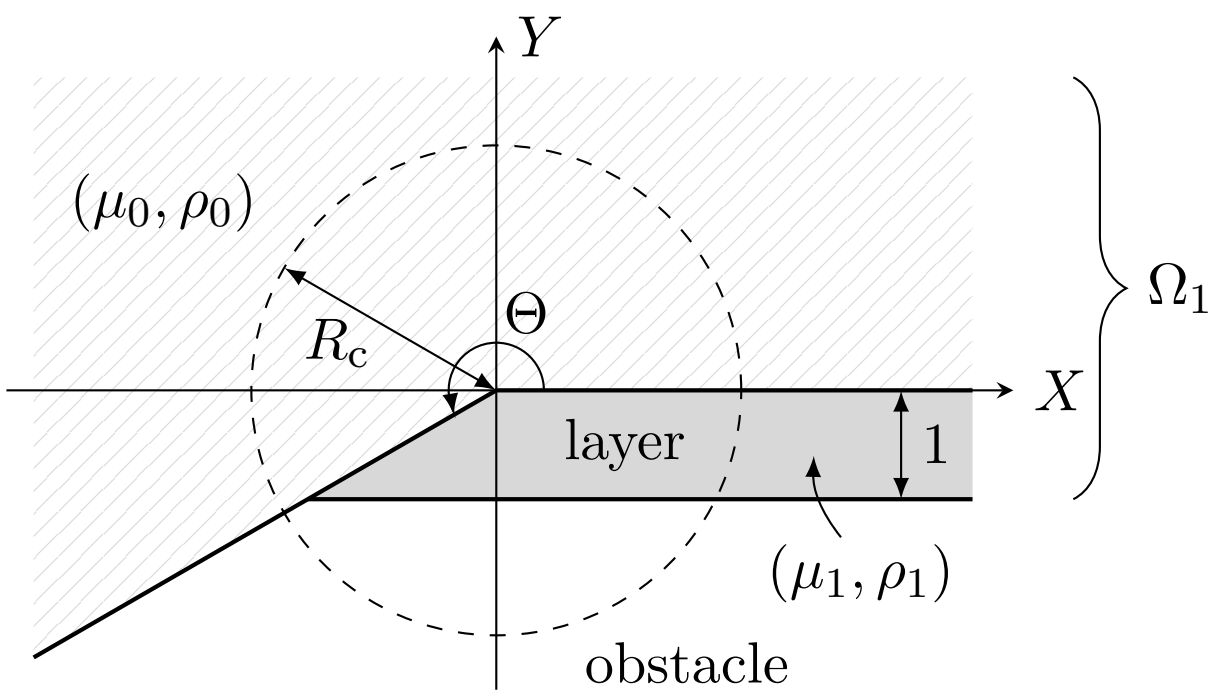

Figure 1: An example for $\Omega_1$

## 2 The effective models

We build two effective models that allow computing approximate solutions of (1) of order $O(\varepsilon^3)$, resp. $O(\varepsilon^5)$. The models are posed in a domain $\Omega^\circ := \Omega \setminus \overline{B}(0, R^\circ)$ for some fixed $R^\circ > 0$. The sets $\Gamma^\circ$, $\Sigma^\circ$ and $\Lambda^\circ$ and the point $P^\circ$ denote parts of $\partial\Omega^\circ$ as in the following figure.

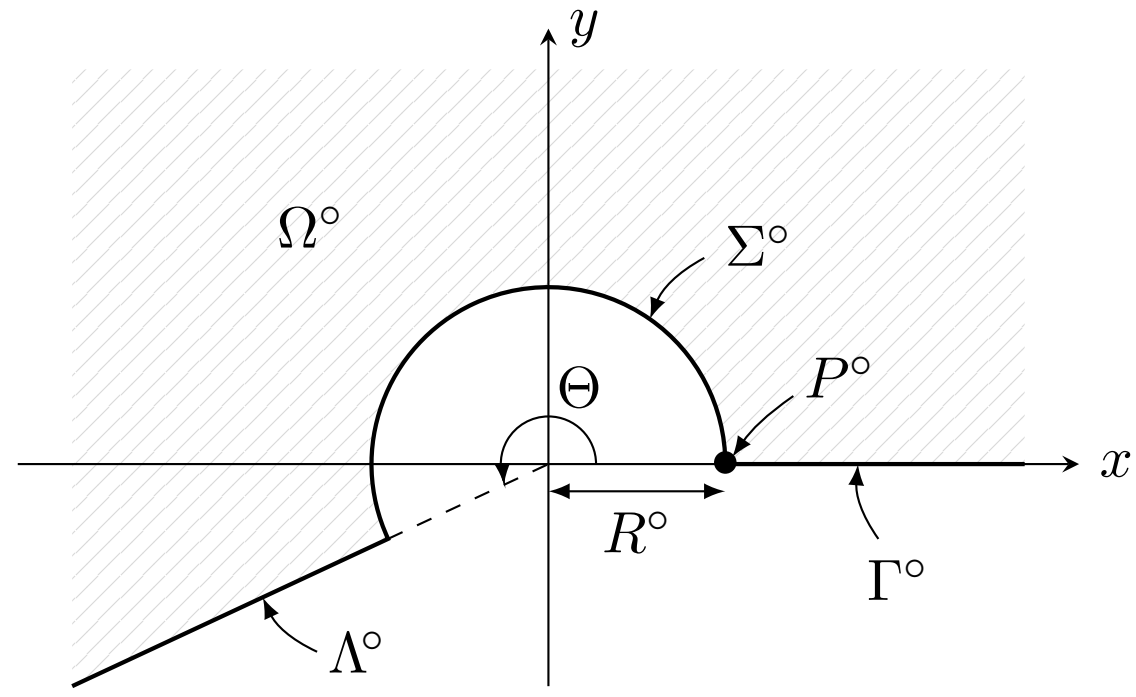

The approximate solution of order $m \in \{3, 5\}$ satisfies $\mu_0 \Delta u_\varepsilon^{(m)} + \omega^2 \rho_0 u_\varepsilon^{(m)} = f_s$ in $\Omega^\circ$ with a Dirichlet condition $u_\varepsilon^{(m)} = 0$ on $\Lambda^\circ$. Moreover the layer is replaced by the following effective boundary conditions on $\Gamma^\circ$:

$$\begin{aligned} \mu_0 \partial_n u_\varepsilon^{(3)} &= -\tfrac{\mu_1}{\varepsilon} u_\varepsilon^{(3)}, \\ \mu_0 \partial_n u_\varepsilon^{(5)} &= -\tfrac{\mu_1}{\varepsilon} u_\varepsilon^{(5)} + \tfrac{\varepsilon}{3}(\mu_1 \partial_x^2 + \omega^2 \rho_1) u_\varepsilon^{(5)}. \end{aligned} \qquad (2)$$

These conditions are classical in presence of an infinite layer. The difficulty here lies in the vicinity of the corner, where the exact solution $u_\varepsilon$ has a highly intricate behavior [1]. We propose a Dirichlet-to-Neumann condition (DtN) on the artificial boundary $\Sigma^\circ$:

$$\partial_n u_\varepsilon^{(m)} = \mathbf{S}_\varepsilon^{(m)}(u_{\varepsilon|\Sigma^\circ}^{(m)}) \quad \text{on } \Sigma^\circ \qquad (3)$$

where $\mathbf{S}_\varepsilon^{(m)} \in \mathcal{L}(L^2(\Sigma^\circ))$ is an operator of finite rank $N$, which is an additional parameter of the method that is not included in the notation for simplicity. The construction of $\mathbf{S}_\varepsilon^{(m)}$ is explained in Section 3. Condition (3) ensures that $u_\varepsilon^{(m)}$ has a behavior around the corner that is similar to $u_\varepsilon$.

Moreover, when $m = 5$, (2) involves the derivative $\partial_x^2$ on $\Gamma^\circ$. This requires an additional boundary condition at the point $P^\circ$, given by $\partial_x u_\varepsilon^{(5)}(P^\circ) = \mathbf{P}_\varepsilon(u^{(5)}_{\varepsilon|\Sigma^\circ})$ with $\mathbf{P}_\varepsilon \in \mathcal{L}(L^2(\Sigma^\circ), \mathbb{C})$. The construction of $\mathbf{P}_\varepsilon$ is similar to $\mathbf{S}_\varepsilon^{(m)}$, so in the next part we focus on $\mathbf{S}_\varepsilon^{(m)}$.

## 3 The DtN operator

The construction of $\mathbf{S}_\varepsilon^{(m)}$ is based on a previous work [1] that established an asymptotic expansion of $u_\varepsilon$ of the form

$$u_\varepsilon = \sum_{p\in\mathbb{N}+\frac{\pi}{\Theta}\mathbb{N},\, p\leqslant m} \sum_{\ell=0}^{\ell_p} \varepsilon^p \ln^\ell \varepsilon \cdot u_{p,\ell} + o(\varepsilon^m) \quad (4)$$

far from the corner, using the method of matched asymptotic expansion. The fields $u_{p,\ell}$ are built inductively w.r.t. $p$ and have an intricate singular behavior near the corner. We show that this behavior can be represented as a modal decomposition of the form

$$u_{p,\ell} = \sum_{j=0}^{\lfloor p \rfloor} \sum_{d\in\frac{\pi}{\Theta}\mathbb{Z}^*,\, d\geqslant j-p} \sigma_d(u_{p-j,\ell})\, \Phi_{d,j}. \quad (5)$$

Here $\Phi_{d,j}$ is a computable function that behaves like $r^{d-j}$ at the corner and depends only on $(\Theta, \mu_0, \mu_1, \rho_0, \rho_1, \omega)$. In particular, $\Phi_{d,0}$ is proportional to $J_d(k_0 r)\sin(d\theta)$ where $J_d$ is the Bessel function of the first kind and $k_0 := \omega\sqrt{\rho_0/\mu_0}$ is the wave number. Hence $\Phi_{d,0}$ is a classical mode of Helmholtz's equation in a corner with Dirichlet conditions. Moreover, in the asymptotic analysis [1], using the matching conditions and the description of the solution near the corner, one obtains the following inductive relation for $d < 0$

$$\sigma_d(u_{p,\ell}) = \sum_{0<d'\leqslant p+d} \left(c_{d,d',p',\ell'}\right) \underset{(p',\ell')}{*} \left(\sigma_{d'}(u_{p',\ell'})\right) \quad (6)$$

where $*$ denotes a convolution product. The coefficient $c_{d,d',p',\ell'}$ depends on algebraic matching coefficients and coefficients coming from corner profiles. Combining (4) and (5), we show

$$u_\varepsilon = \sum_{d\in\frac{\pi}{\Theta}\mathbb{Z}^*,\, d\geqslant -m} \sigma_d(u_\varepsilon)\, \Phi_d^\varepsilon + o(\varepsilon^m) \quad (7)$$

with
$$\begin{cases} \sigma_d(u_\varepsilon) := \displaystyle\sum_{p\in\mathbb{N}+\frac{\pi}{\Theta}\mathbb{N},\, p\leqslant m} \sum_{\ell=0}^{\ell_p} \varepsilon^p \ln^\ell \varepsilon \cdot \sigma_d(u_{p,\ell}) \\ \Phi_d^\varepsilon := \displaystyle\sum_{j=0}^{m} \varepsilon^j \Phi_{d,j}. \end{cases} \quad (8)$$

Then, using (6) and writing $d$ as $n\frac{\pi}{\Theta}$, we deduce

$$u_\varepsilon = \sum_{n=1}^{\infty} \sigma_{n\frac{\pi}{\Theta}}(u_\varepsilon)\, \Psi_n^\varepsilon + o(\varepsilon^m) \quad (9)$$

where $\Psi_n^\varepsilon := \Phi^\varepsilon_{n\frac{\pi}{\Theta}} + \sum_{-m\leqslant d<0} c^\varepsilon_{n,d} \Phi_d^\varepsilon$ with:

$$c^\varepsilon_{n,d} := \sum_{\substack{p\in\mathbb{N}+\frac{\pi}{\Theta}\mathbb{N} \\ 1\leqslant p\leqslant m}} \sum_{\ell=0}^{\ell_p} \varepsilon^p \ln^\ell \varepsilon \cdot c_{d,n\frac{\pi}{\Theta},p,\ell} \quad (10)$$

To build $\mathbf{S}_\varepsilon^{(m)}$, we exploit the series (9) truncated at a rank $N$, which is a parameter of the model. $\mathbf{S}_\varepsilon^{(m)}$ is designed so that

$$\mathbf{S}_\varepsilon^{(m)}\left(\sum_{n=1}^{N} \lambda_n \Psi_n^\varepsilon\right) = \sum_{n=1}^{N} \lambda_n\, \partial_n \Psi_n^\varepsilon \quad (11)$$

for any $(\lambda_n) \in \mathbb{C}^N$. Denoting $u_\lambda := \sum \lambda_n \Psi_n^\varepsilon$, one has $\left[\langle u_\lambda, \Psi_n^\varepsilon\rangle_{L^2(\Sigma^\circ)}\right]_{n\in[\![1,N]\!]} = G_\varepsilon \left[\lambda_n\right]_{n\in[\![1,N]\!]}$ where $G_\varepsilon := \left[\langle \Psi_n^\varepsilon, \Psi_{n'}^\varepsilon\rangle_{L^2(\Sigma^\circ)}\right]_{n,n'\in[\![1,N]\!]}$ is the Gram matrix of the modes $\Psi_n^\varepsilon$. Note that (8) and (10) imply $\Psi_n^\varepsilon = \Phi_d^\varepsilon + O(\varepsilon) = \Phi_{d,0} + O(\varepsilon)$ with $d := n\frac{\pi}{\Theta}$. Therefore, when $\varepsilon$ is small, the modes $\Psi_n^\varepsilon$ are almost orthogonal in $L^2(\Sigma^\circ)$ and $G_\varepsilon$ can be proven to be invertible. Hence the actual definition of $\mathbf{S}_\varepsilon^{(m)}$ is: $\forall u \in L^2(\Sigma^\circ)$,

$$\mathbf{S}_\varepsilon^{(m)}(u) := \Big[\langle u, \Psi_n^\varepsilon\rangle_{L^2(\Sigma^\circ)}\Big]_n^{\mathsf{T}} G_\varepsilon^{-1} \Big[\partial_n \Psi^\varepsilon_{n|\Sigma^\circ}\Big]_n.$$

## 4 Final results

Using again $\Psi_n^\varepsilon = \Phi_{d,0} + O(\varepsilon)$, one can show that $\mathbf{S}_\varepsilon^{(m)}$ is a perturbation of a positive DtN operator, which allows to prove the well-posedness of the effective models for sufficiently small $\varepsilon$. Let us assume that $R_\mathrm{s} := \mathrm{dist}(0, \mathrm{supp}\, f_\mathrm{s})$ is positive. We prove the error estimate

$$\|u_\varepsilon^{(m)} - u_\varepsilon\|_{H^1(\Omega^\circ)} \lesssim \left(\frac{\varepsilon}{R^\circ}\right)^m N^2 + \left(\frac{R^\circ}{R_\mathrm{s}}\right)^{N\frac{\pi}{\Theta}} \quad (12)$$

where $\lesssim$ denotes an inequality valid up to a constant independant of $\varepsilon$, $R^\circ$ and $N$. The second term of (12) is due to the truncation of the modal decomposition (9). Numerical validations will be presented during the talk.

# Continuous and discrete stability of dispersive perturbation of hyperbolic initial boundary value problems on the half-line

**Geoffrey Beck**[1,*], **Nicolas Crouseilles**[1], **Ludovic Martaud**[1]
[1]Univ Rennes, IRMAR UMR 6625, Centre Inria de l'Université de Rennes (MINGuS), France
*Email: geoffrey.a.beck@inria.fr

**Abstract**

This work concerns the theoretical study and numerical approximation of a class of dispersive perturbation of hyperbolic initial boundary value problems which describes the wave propagation in dispersive media in the context of thin or micro-structured waveguides. The number of boundary conditions is not necessarily the same for the hyperbolic and dispersive systems. Given a numerical flux adapted to the hyperbolic case without boundary, we explain how to modify it to take into account the dispersive perturbation and the boundary.

## 1 Dispersive perturbation of hyperbolic initial boundary probems

In the framework of gravity water-wave dynamics, the Boussinesq-type systems are obtained from the free-surface Euler equations through an asymptotic expansion in a shallow-water regime. When the shallowness parameter $\kappa$ vanishes, one recovers the so-called Saint-Venant equation which is a hyperbolic system. A more accurate description, incorporating dispersive phenomena, is achieved by retaining terms of order $\kappa^2$ in the asymptotic expansion. In the linear regime, this more accurate model takes the following form

$$\big(I_p - (\kappa\pi\partial_x)^2\big)\partial_t u + \partial_x A u = 0, \qquad (1)$$

with unknown $u(x,t) \in \mathbb{R}^p$, $t \geq 0, x \geq 0$, $I_p$ the identity matrix, $\pi$ is a projector and $A$ is a $\mathbb{R}$-diagonalizable matrix satisfying the *Kawashima condition* $\ker(A) \cap \ker(\pi) = \{0\}$.

One would like to study (1) on the *half-line* $x \in \mathbb{R}^+$. In the case of the Saint-Venant equations, it is well established that relevant boundary conditions consist in prescribing either the discharge or the free-surface elevation. It is therefore natural to expect a similar structure for dispersive extensions of the Saint-Venant system. However, dispersive corrections are not unique. As illustrated by the so-called abcd Boussinesq systems, several distinct perturbations can be constructed. Certain abcd Boussinesq systems require supplementary boundary conditions, whose physical interpretation may be unclear. This is why, one considers quite general type of boundary conditions:

$$\mathsf{N}(u, \kappa\pi\partial_x u)^{\mathsf{T}} = g. \qquad (2)$$

At the discrete level, some difficulties arise from the boundary conditions that do not provide a complete knowledge of the solution at the boundary, whereas any standard finite volume numerical scheme requires this knowledge. A hidden set of ordinary differential equations on the trace of $u$ enables one to provide the missing boundary conditions and help use to establish the fully continuous and discrete stability.

## 2 How many boundary conditions ?

The best way to understand which boundary conditions (2) are compatible with the linear problem (1) is to solve explicitly the system when the only data are the boundary conditions. More precisely, applying the Laplace transform $\mathcal{L}$ : $u \mapsto \hat{u}(s)$ to (1,2) leads to the following ODE

$$\begin{cases} s\hat{U}_\kappa + \mathsf{A}_\kappa(s)\partial_x \hat{U}_\kappa = 0, \\ \mathsf{N}\hat{U}_\kappa(x=0) = \hat{g}, \end{cases} \qquad (3)$$

where

$$\hat{U}_\kappa := \begin{pmatrix} \hat{u} \\ \kappa\pi\partial_x\hat{u} \end{pmatrix}, \quad \mathsf{A}_\kappa(s) = \begin{pmatrix} A & -s\kappa\pi \\ -s\kappa\pi & 0 \end{pmatrix},$$

for all $s$ such that $\Re(s) > 0$.

One first shows that $\mathsf{A}_\kappa(s)$ is $\mathbb{C}$-diagonalizable thanks to the Kawashima condition. One denotes by $p^+(s)$ the number of eigenvalues $\mu_i(s)$ such that $\mathfrak{Re}\big(s\overline{\mu_j(s)}\big) > 0$ and by $E_e^+(s)$ the set of associated eigenvectors. This are assotiated to exponentially decreasing (with respect to $x$) solutions to (3). Secondly, when $A$ admits a block-diagonal Friedrichs symmetrizer along $\pi$, one can show that

$$p^+(s) = \text{Rank } (\pi) + N_+(A),$$

where $N_+(A)$ is the the number of positive eigenvalue of $A$. We say that (1,2) satisfy the *Dispersive Uniform Kreiss Lopatinskii (DUKL) condition if for all s such that* $\Re(s) > 0$, $p^+ = \mathit{Rank}\ (N)$ *and* $E_e^+(s) \not\subset \ker \mathsf{N}$. In that case, there exists a unique solution $\hat{U}_\kappa(s) \in L^2_x(\mathbb{R}^+)$ to (3) for all s such that $\Re(s) > 0$. Thus $p^+$ is the number of mandatory boundary conditions.

## 3 Hidden trace equations

The dispersive system (1) can be reformulated under a *non-local conservative form*:

$$\partial_t u + \partial_x\big(\pi A(\mathcal{R}u) + \pi^* Au + \kappa\pi\dot{u}(0,t)e^{-\frac{x}{\kappa}}\big) = 0,$$

where $\pi^* = Id - \pi$ and the operator $\mathcal{R}$ is the inverse operator of $1-(\kappa\partial_x)^2$ with homogeneous Neumann boundary conditions. If one moreover assumes some technical assumptions on $\mathsf{N}$ and $A$ more restrictive than the (DUKL) condition, then there exists a set of *hidden ODEs*

$$\kappa\dot{w}(t) = Kw(t) + L(\mathcal{R}\pi u)(0,t) + Q\mathsf{g}(t),$$

satisfied by the boundary conditions

$$w(t) = (\kappa\dot{u}, u, \mathcal{R}\pi^* u)^\mathsf{T}(0,t),$$

where $K, L, Q$ are explicit matrices and $\mathsf{g}$ is related to the source $g$. This set of hidden ODEs is the keystone to find an energy estimate and show well-posedness.

## 4 Numerical methods based on hidden ODEs

Given a consistent flux $\mathfrak{F}[u_L, u_R] = f(A, u_L, u_R)$ adapted to the hyperbolic case ($\kappa = 0$) without boundary, one defines the flux

$$\mathcal{F}_\kappa[u_L, u_R] = \tilde{\mathfrak{F}}_\kappa[u_L, u_R] - \frac{\nu}{2}\left((\mathcal{R}^{\Delta x}u)_R - (\mathcal{R}^{\Delta x}u)_L\right).$$

where $\tilde{\mathfrak{F}}_\kappa[u_L, u_R] := f(\pi A\mathcal{R}^{\Delta x} + \pi^* A, u_L, u_R)$ and $\mathcal{R}^{\Delta x}$ stands for one spatial discretization to the elliptic problem that define the operator $\mathcal{R}$. The equation under a non-local conservative form and the hidden ODEs are discretized through the following scheme

$$\begin{cases} \frac{u_i^{n+1}-u_i^n}{\Delta t} + \frac{\mathcal{F}_\kappa[u_i^n, u_{i+1}^n] - \mathcal{F}_\kappa[u_{i-1}^n, u_i^n]}{\Delta x} = \pi w_2^n \gamma_i \\ w^{n+1} = w^n + \frac{\Delta t}{\kappa}[Kw^n + L(\mathcal{R}^{\Delta x}\pi u)_0^n + Q\mathsf{g}^n], \end{cases}$$

where $u_i^n$ denotes the finite volume numerical unknown. It is worth mentioning that $\gamma_i$ is a numerical approximation of $e^{-x/\kappa}$ obtained with the same discretization chosen for $\mathcal{R}^{\Delta x}$.

We have shown that for some choice of $\mathfrak{F}$ and $\mathcal{R}^{\Delta x}$ there is $\nu > 0$ such that the considered scheme is stable.

## 5 Numerical methods based on augmented hyperbolic problem after relaxation

The dispersion in the flux $\mathcal{F}_\kappa$ induces the computation of non-local terms related to $\mathcal{R}$. We can replace the dispersive term $\kappa^2\partial_x\partial_t\pi u$ by $p = \dot{J}$ where $J$ solves $\alpha^2\ddot{J}+J = \kappa^2\partial_x\pi u$ and $\alpha \ll 1$ is the relaxation parameter. Then the augmented variable $(u, J, p)^\mathsf{T}$ satisfy a characteristic hyperbolic problem plus an additional stiff term in $\alpha^{-2}$. Formally, when $\alpha$ vanishes, one recovers the dispersive system (1). This is formally, well-prepared initial data and compatibility conditions are mandatory to avoid some pathological cases. This method can be applied on the half-line only for certain boundary conditions when $\pi^* A\pi^* = 0$. Discrete boundary procedure computations with Riemann invariants and an artificial viscosity are proposed to correct any three-points finite volume schemes (see [1]). We use this procedure to solve numerically the relaxed system on $(u, J, p)^\mathsf{T}$. We compare simulations performed using the dispersive and relaxed schemes with an exact solution in the case where (1) is the classical Boussinesq-Abbott system. When $\kappa$ is small, the dispersive scheme gives a better approximation of the solution near the boundary.

**Acknowledgement**

The first author is supported by the BOURGEONS project, grant ANR-23-CE40-0014-01 of the French National Research Agency (ANR).

# A Flux-Based Domain Decomposition Method for Time-Harmonic Wave Problems

**Nawfel Benatia**[1,*], **Christophe Geuzaine**[1]
[1]Department of Electrical Engineering and Computer Science, University of Liège, Belgium
*Email: nawfel.benatia@uliege.be

**Abstract**

We propose a non-overlapping domain decomposition method for the time-harmonic Helmholtz equation, incorporating Robin-type transmission conditions directly into the numerical fluxes. We study the well-posedness of the local subproblems and establish that the associated scattering operator is strictly contractive, thereby guaranteeing convergence of the global iteration scheme. Numerical experiments with a high-order $H^1/H(\text{div})$ finite element discretization support the theoretical analysis.

## 1 Introduction

Domain Decomposition Methods (DDMs) offer an effective framework for solving large-scale partial differential equations by partitioning the domain into smaller subdomains, naturally enabling parallel computation. A central difficulty lies in enforcing suitable transmission conditions to ensure interface compatibility. Classical Schwarz methods with Dirichlet or Neumann conditions fail to converge for wave propagation problems even with overlap, and break down entirely in non-overlapping settings. The introduction of Robin transmission conditions by Lions [1] and Després [2] restored convergence for elliptic and Helmholtz problems, forming the basis of optimized Schwarz methods [3] and marking a major theoretical and computational advance.

In this work, we combine the flexibility of Robin transmission conditions with the advantages of hybridizable discontinuous Galerkin discretizations. We formulate the DDM through numerical fluxes rather than explicit boundary conditions, thereby constructing an interface operator that is strictly contractive in a suitable energy norm. This is equivalent to the Characteristic Hybridizable Discontinuous Galerkin method (CHDG) [4], albeit at the level of a subdomain instead of each mesh element.

## 2 Model

We consider the time-harmonic scalar wave problem on a Lipschitz polygonal domain $\Omega \subset \mathbb{R}^2$ formulated as a mixed first-order system : Find $u : \Omega \to \mathbb{C}$ and $\mathbf{q} : \Omega \to \mathbb{C}^2$ such that

$$\begin{cases} \boldsymbol{\nabla}\cdot\mathbf{q} + \kappa^2 u = 0, & \text{in } \Omega, \\ -\mathbf{q} + \boldsymbol{\nabla} u = \mathbf{0}, & \text{in } \Omega, \\ \mathbf{n}\cdot\mathbf{q} - \imath\kappa u = g, & \text{on } \partial\Omega, \end{cases}$$

where $\kappa > 0$ is the wavenumber, $\mathbf{n}$ the outward unit normal, and $g$ the boundary data.

Let $\mathscr{T}_h$ be a conforming, shape-regular triangulation of $\Omega$, such that

$$\overline{\Omega} \simeq \overline{\Omega}_h = \bigcup_{T\in\mathscr{T}_h} \overline{T}, \qquad h = \max_{T\in\mathscr{T}_h} \operatorname{diam}(T).$$

The discrete domain $\Omega_h$ is partitioned into non-overlapping subdomains $\{\Omega_i\}_{i=1}^{N_{\text{sub}}}$ with disjoint interiors. Interfaces between adjacent subdomains are defined by

$$\Sigma_{ij} \coloneqq \partial\Omega_i \cap \partial\Omega_j, \quad i < j,$$

with outward normals satisfying $\mathbf{n}_i = -\mathbf{n}_j$ on each $\Sigma_{ij}$. The skeleton of the decomposition is

$$\Sigma \coloneqq \bigcup_{i=1}^{N_{\text{sub}}} \partial\Omega_i, \quad \partial\Omega_i = \Gamma_i \cup \Big( \bigcup_{j\in\mathfrak{N}_i} \Sigma_{ij} \Big),$$

where $\Gamma_i \coloneqq \overline{\Omega}_i \cap \partial\Omega$, and $\mathfrak{N}_i \coloneqq \{j : \Sigma_{ij} \neq \emptyset\}$ denotes the set of neighbors of $\Omega_i$.

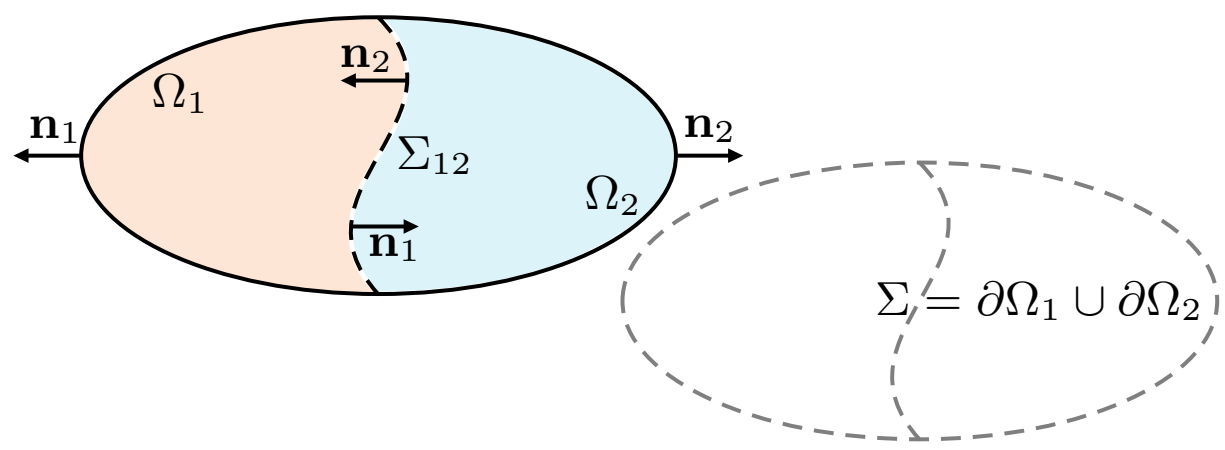

Figure 1: Decomposition into two subdomains $\Omega_1$ and $\Omega_2$ with interface $\Sigma_{12}$ and skeleton $\Sigma$.

## 3 Methods

On each subdomain $\Omega_i$, we consider the conforming finite element space $V_h(\Omega_i) \subset H^1(\Omega_i)$

of continuous piecewise $\mathbb{P}_p$ polynomials and the Raviart–Thomas space $\mathbf{W}_h(\Omega_i) \subset H(\text{div}, \Omega_i)$. We denote by $V_h(\partial\Omega_i) := \text{tr}(V_h(\Omega_i))$ the trace space of piecewise polynomial interface unknowns, which are *a priori* discontinuous across interfaces.

For each subdomain $\Omega_i$, we solve a local weak mixed problem in which physical boundary conditions are replaced by numerical fluxes. On an interface $\Sigma_{ij}$, characteristic fluxes are defined by

$$\widehat{u}_i = -\frac{\lambda_i + \lambda_j}{2\imath\kappa}, \qquad \mathbf{n}_i \cdot \widehat{\mathbf{q}}_i = \frac{\lambda_j - \lambda_i}{2},$$

where the transmission traces are given by

$$\lambda_i = -\mathbf{n}_i \cdot \mathbf{q}_i - \imath\kappa u_i.$$

These fluxes implicitly enforce Robin transmission conditions. These are the same numerical fluxes as in CHDG, which implicitly embed Robin transmission conditions at interfaces. To decouple the local solves, we introduce a global multi-trace variable $\mu_h \in \mathbb{V}_h := \prod_i V_h(\partial\Omega_i)$, with local restrictions $\mu_i := \mu_h|_{\partial\Omega_i}$. The fluxes become

$$\widehat{u}_i(\mu_i) = -\frac{\lambda_i + \mu_i}{2\imath\kappa}, \qquad \mathbf{n}_i \cdot \widehat{\mathbf{q}}_i(\mu_i) = \frac{\mu_i - \lambda_i}{2},$$

where $\mu_i$ represents the incoming trace (equal to $\lambda_j$ on $\Sigma_{ij}$ and to $g$ on $\Gamma_i$). Each subdomain thus depends only on its local interface data $\mu_i$. Global coupling is enforced through the swap operator $\Pi : \mathbb{V}_h \to \mathbb{V}_h$,

$$(\Pi(\xi_h))|_{\partial\Omega_i} = \begin{cases} \xi_j, & \text{on } \Sigma_{ij}, \\ 0, & \text{on } \Gamma_i, \end{cases}$$

and the scattering operator $\mathsf{S} : \mathbb{V}_h \longrightarrow \mathbb{V}_h$,

$$(\mathsf{S}(\mu_h))|_{\partial\Omega_i} := -\mathbf{n}_i \cdot \mathbf{q}_i(\mu_i) - \mathrm{i}\kappa\, u_i(\mu_i).$$

Let $b$ the global right-hand side

$$b|_{\partial\Omega_i} = \begin{cases} 0 & \text{on } \Sigma_{ij}, \\ g & \text{on } \Gamma_i, \end{cases}, \qquad \text{and} \qquad b_h := \mathcal{P}_h b \in \mathbb{V}_h,$$

where $\mathcal{P}_h$ is the $L^2(\Sigma)$-projection onto $\mathbb{V}_h$. The global skeleton problem reads:

$$\begin{cases} \text{Find } \mu_h \in \mathbb{V}_h \text{ such that} \\ \quad (\mathsf{I} - \Pi \circ \mathsf{S})\mu_h = b_h. \end{cases}$$

The composed operator $\Pi \circ \mathsf{S}$ is strictly contractive: $\|\Pi \circ \mathsf{S}\,(\mu_h)\|_\Sigma < \|\mu_h\|_\Sigma$, for all $\mu_h \neq 0$, ensuring geometric convergence with rate $\rho = \|\Pi \circ \mathsf{S}\| < 1$. After convergence, local solutions $(u_i(\mu_i), \mathbf{q}_i(\mu_i))$ are obtained independently and assembled to form $(u_h, \mathbf{q}_h)$.

## 4 Results

The flux-based DDM is assessed on a unit disk partitioned into 16 subdomains using a 45° plane wave as the analytical solution at $\kappa = 16\pi$. The method is formulated in $H^1/H(\text{div})$ spaces with polynomial degrees $p = 4/5$ on a mesh of size $h = 0.025$. Although $\Pi \circ \mathsf{S}$ is contractive and guarantees convergence of the fixed-point iteration $\mu_h^{(\ell+1)} = (\Pi \circ \mathsf{S})(\mu_h^{(\ell)}) + b_h$, the global system $(\mathsf{I} - \Pi \circ \mathsf{S})\mu_h = b_h$ is solved in practice using GMRES, requiring 120 iterations. Numerical eigenvalue analysis confirms the theoretical prediction: the skeleton operator is strictly contractive with spectral radius $\rho = 0.9999993860$, and all 4,408 eigenvalues satisfy $|\lambda_i| < 1$.

The full contribution will provide detailed analysis and the application to larger scale problems.

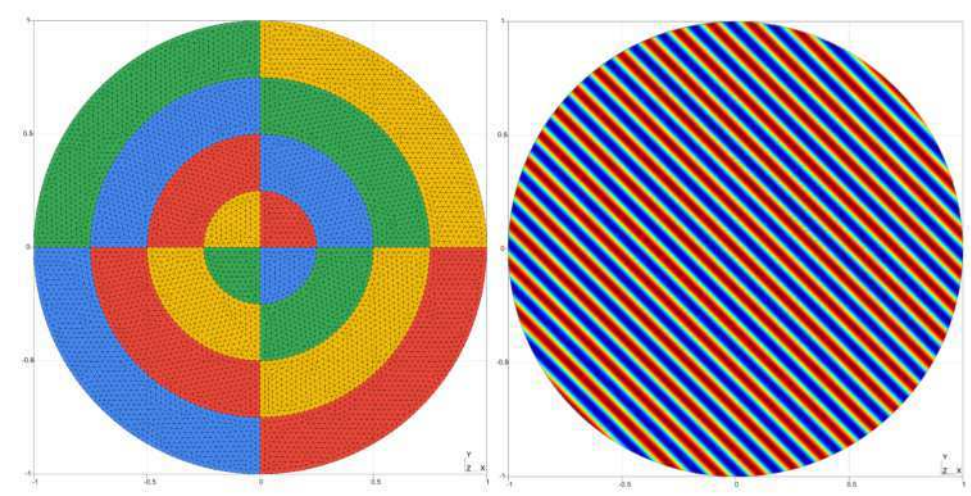

Figure 2: Mesh decomposition of the unit disk for $h = h_1$ (left) and the real part of the solution field, $\Re[u(x, y)]$, at $\kappa = 16\pi$ for $h = 0.01$ (right).

# A Robust Model Order Reduction Method to Accelerate Full Waveform Inversion via Fréchet Derivatives

Julien Besset[1,*], Hélène Barucq[1], Victor Martins Gomes[1]
[1]Makutu, Inria, TotalEnergies, University of Pau et des Pays de l'Adour, France
*Email: julien.besset@totalenergies.com

## Abstract

Subsurface imaging allows the characterization of underground structures using geophysical methods. It plays a crucial role in energy production (O&G and geothermal), civil engineering, and CO2 storage. Full Waveform Inversion (FWI) is a widespread imaging technique that relies on repeated solutions of the wave equation, resulting in a high computational cost. Seeking to accelerate the acoustic FWI process, we employ a Model Order Reduction (MOR) approach to reduce the dimensionality of the systems to be solved while preserving their essential dynamics. This method, called MOR-T$_L$, incorporates Fréchet derivatives of the propagation operator with respect to the model parameter (velocity) to allow the reduced bases to remain valid despite successive parameter updates during the inversion. Our approach is successfully validated for 2D acoustic FWI problems.

## 1 Model

We seek for $u : \Omega \times (0,T) \longrightarrow \mathbb{R}$ satisfying

$$\begin{cases} \dfrac{1}{c^2}\dfrac{\partial^2 u}{\partial t^2}(\mathbf{x},t) - \Delta u(\mathbf{x},t) = f(\mathbf{x},t) & \text{in } \Omega \times (0,T) \\ u(\mathbf{x},t) = 0 & \text{on } \Gamma_1 \times (0,T) \\ \dfrac{\partial u}{\partial \mathbf{n}}(\mathbf{x},t) + \dfrac{1}{c}\dfrac{\partial u}{\partial t}(\mathbf{x},t) = 0 & \text{on } \Gamma_2 \times (0,T) \\ u(\mathbf{x},0) = 0,\ \dfrac{\partial u}{\partial t}(\mathbf{x},0) = 0 & \text{in } \Omega \end{cases} \tag{1}$$

We consider $\Omega \subset \mathbb{R}^2$ an open bounded domain with boundary $\partial\Omega = \Gamma_1 \cup \Gamma_2$, where $\Gamma_1$ is the top and $\Gamma_2$ the remainder of the boundary. We assume $\Omega$ is convex (generally a cuboid) with Lipschitz boundary. The parameter $c(\mathbf{x}) > 0$, $\forall \mathbf{x} \in \Omega$, is the acoustic wave speed in the domain, $u$ is the pressure wave field, $f$ is the source term which is independent of $c$ and taken here as a Ricker wavelet. We impose a free surface condition on $\Gamma_1$, and an Absorbing Boundary Condition (ABC) on $\Gamma_2$. Here, $\mathbf{n}$ is the outward-pointing unit normal vector on $\partial\Omega$.

## 2 Methodology

Let $\Theta$ be the space of parameters with $\theta = \dfrac{1}{c^2(\mathbf{x})} \in \Theta$. $\theta_0$ is an initial parameter and $\delta\theta \in \Theta$ a small perturbation. Suppose that $u \in \mathcal{C}^L(\Theta)$, $L \in \mathbb{N}^*$. We introduce the Taylor polynomial of order $L$ centered at $\theta_0$:

$$Tu(\theta_0) = \sum_{\ell=0}^{L} \frac{(\alpha\delta\theta)^\ell}{\ell!} D^\ell u(\theta_0)$$

with $\alpha \in \mathbb{R}^*$.
We seek a basis $\{\phi_i^*\}_{i=1}^N$ satisfying:

$$\underset{\{\phi_i\}_{i=1}^N}{\arg\min}\ \frac{1}{2}\|u(\theta_0 + \alpha\delta\theta) - \mathcal{P}_\phi Tu(\theta_0)\|_{L^2}^2, \tag{2}$$

where $\mathcal{P}_\phi$ denotes the projection onto the reduced space spanned by $\{\phi_i\}_{i=1}^N$. It follows that the reduced basis $\{\phi_i^*\}_{i=1}^N$ incorporates variations of $\theta_0$ in the direction $\alpha\delta\theta$.

Let $v^\ell(\theta_0)$ denote the $\ell$-th order Fréchet derivative of $u$ at $\theta_0$ in the direction $\delta\theta$, denoted by $D^\ell u(\theta_0)$, for $\ell = 0, \ldots, L$. $v^\ell$ satisfies (1) with different right-hand sides that depend on $\ell$ [1].

We use the Spectral Element Method (SEM) [2] to obtain a high-fidelity representation of $v^\ell$. For $\theta = \theta_0$,

$$v^\ell(\theta_0) \coloneqq v^\ell_{\mathbf{SEM}}(\theta_0) = \sum_{i=1}^{N_x} v_i^\ell(t)\psi_i(\mathbf{x}), \tag{3}$$

where $N_x$ denotes the degrees of freedom for the finite element discretization and $\{\psi_i\}_{i=1}^{N_x}$ are the Lagrange basis functions defined at the Gauss-Lobatto-Legendre points.

Following [1], we replace the $\alpha$-dependent minimization problem (2) with the search for (up to a controlled approximation error in $\alpha$ and $L$):

$$\underset{\{\phi_i\}_{i=1}^N}{\arg\min}\ \frac{1}{2}\|v^\ell_{\mathbf{SEM}}(\theta_0) - \mathcal{P}_\phi v^\ell_{\mathbf{SEM}}(\theta_0)\|_{L^2}^2, \tag{4}$$

## 3 Implementation

We aim to construct a reduced basis $\{\phi_i\}_{i=1}^N$ that approximates the Fréchet derivatives $v^\ell(\theta_0)$ while minimizing the projection error (4). In standard data-driven MOR [3], the reduced basis is obtained via a truncated SVD of the snapshot matrix, which collects high-fidelity SEM solutions $v^\ell_{\mathbf{SEM}}(\theta_0)$ over time. For high-dimensional problems, this approach is costly and memory-prohibitive, as it lacks a systematic strategy for selecting snapshots. We mitigate these issues by employing a truncated QR decomposition based on a Gram-Schmidt (GS) orthonormalization process. The reduced basis is constructed incrementally from the SEM snapshots of each Fréchet derivative. This avoids the need to store the full snapshot matrix, providing a crucial advantage for large-scale applications. Each new snapshot is orthogonalized against the current basis vectors and retained only if its residual norm exceeds a prescribed truncation threshold $\epsilon$ [1].

The final reduced basis $\{\phi_i\}_{i=1}^N$ is obtained by concatenating and reordering the partial bases of each $\ell$-th order derivative, followed by a second GS orthonormalization. The QR-based approach offers superior memory efficiency and scalability compared to SVD-based methods. Its effectiveness depends on the truncation parameter $\epsilon$, which balances computational accuracy and cost.

Once the basis has been constructed, the solution for a perturbed model $\theta_0 + \alpha\delta\theta$ can be represented as:

$$u(\theta_0+\alpha\delta\theta) \coloneqq u_{\mathbf{ROM}}(\theta_0+\alpha\delta\theta) = \sum_{i=1}^{N} a_i(t)\phi_i(\mathbf{x}) \tag{5}$$

with $N \ll N_x$ and where $\{a_i\}_{i=1}^N$ are computed as the solution of the Reduced Order Model (ROM) equation by plugging (5) into (1).

### 3.1 Numerical Results

Our ROM-FWI algorithm, which integrates MOR-$\mathrm{T}_L$ directly within the inversion loop, successfully accelerates the classical FWI workflow. At each iteration, the computed gradient defines a perturbation direction $\delta\theta$, used to generate the reduced basis via the solutions of the corresponding Fréchet derivatives. The resulting ROM-system enables the line search, that is, the evaluation of velocity parameter updates along the perturbation direction $\delta\theta$ (i.e., $\theta_0 + \alpha\delta\theta$), to be performed entirely in the reduced space. Compared to classical FWI, our approach significantly reduces computational cost while maintaining similar parameter uncertainty, as shown in Table 1.

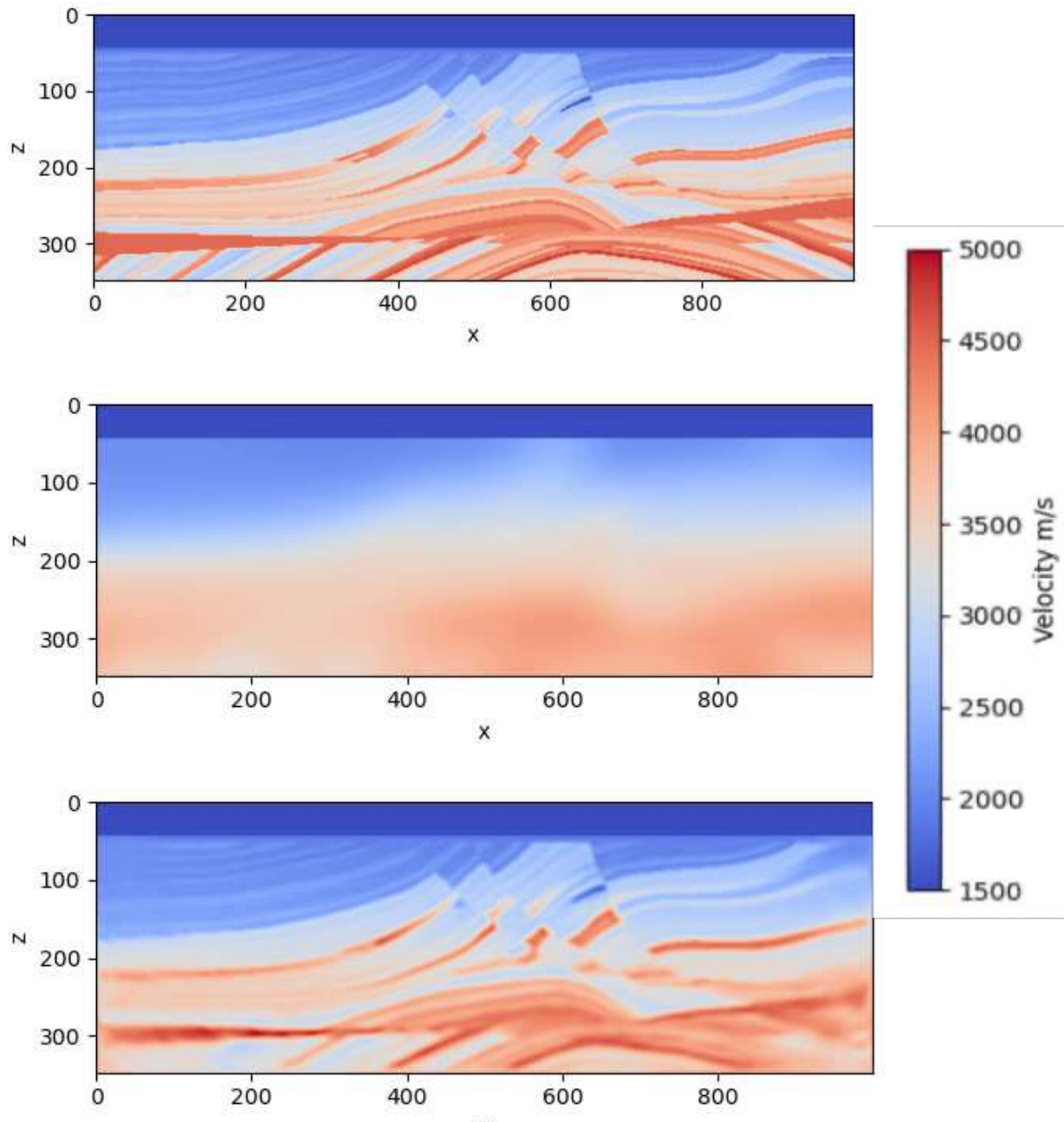

Figure 1: Target (top), Starting (center) and Final (bottom) model.

| Method | LS CCT | Total CCT | $\frac{\lVert c_f - c_{\text{tar}} \rVert_2}{\lVert c_{\text{tar}} \rVert_2}$ |
|---|---|---|---|
| ROM-FWI | 760.95h | 1650.63h | 5.72% |
| FWI | 1599.98h | 2356.33h | 5.75% |

Table 1: Line Search (LS) and Total Cumulated Computational Time (CCT) for classical FWI and ROM-FWI. Rightmost column shows the relative L2-norm error between calculated and target models. Both FWI ran with 21 clients and 48 MPIs on a model containing 1,394,697 degrees of freedom.

# Analysis and Simulation of 2D Maxwell's equations in Cold Plasma

Manaswinee Bezbaruah[1,*], Patrick Ciarlet[1], Maryna Kachanovska[1]
[1]POEMS (ENSTA Paris), Palaiseau, France
*Email: manaswinee.bezbaruah@ensta.fr

## Abstract

This work is concerned with hybrid resonances in cold magnetised plasma, a phenomenon that plays a key role in plasma heating. Such problems are associated with the singular solution of Maxwell's equations with frequency-dependent, sign-changing dielectric tensors. In simplified settings, we arrive at a PDE with smooth sign-changing coefficients vanishing on a curve inside the domain. This PDE falls outside of the scope of standard elliptic theory and requires non-classical variational formulations. In particular, such variational formulations have to be consistent with the limiting absorption principle and must allow for singular solutions.

Our recent work has been focused on the numerical implementation of two formulations of such boundary value problems. We concentrate on the 2D case, and re-interpret the corresponding problem as a transmission problems, in the spirit of the recent PhD thesis of É. Peillon [4]. We present two numerical methods and discuss the related computational implementation.

## 1 Introduction

Plasma is recognised as the fourth state of matter and is one of the more prevalent form of matter in the universe. It consists of charged particles and ions, with concentrations and densities that can change both spatially and temporally. In practice, plasma is mainly used for industrial purposes, for example in electric energy production via fusion nuclear reactors. This, in turn, motivates large facets of academic and industrial research. One of the major challenges of achieving stable nuclear fusion is plasma heating, and it is resolved by sending electromagnetic waves at specific frequencies and directions depending on the characteristics of the plasma. We will investigate *hybrid resonance* in cold plasma - which leads to a *sign-changing degenerate* PDE.

## 2 The Model

Consider $\Omega = (-a, a) \times [0, L] \subset \mathbb{R}^2$ with left and right boundaries $\Gamma_n, \Gamma_p$ respectively, and top and bottom boundary $\Gamma_+, \Gamma_-$ respectively.

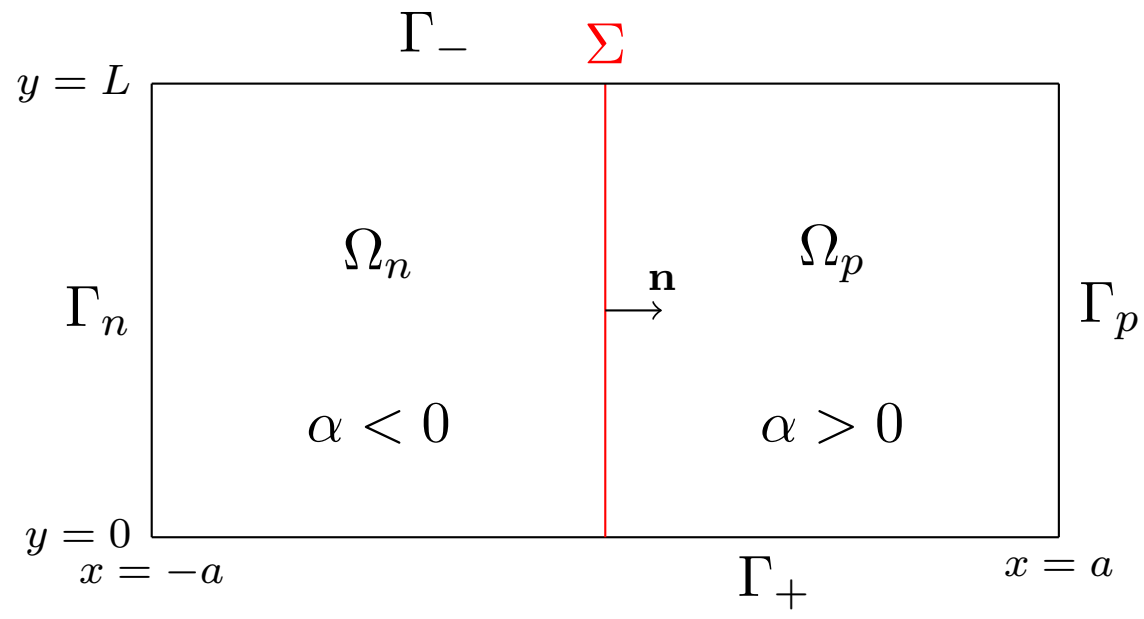

Let the function $\alpha \in C^2(\overline{\Omega})$, be such that $\alpha > 0$ in $\Omega_p$ and $\alpha < 0$ in $\Omega_n$. Moreover, $\alpha = x \cdot r(y) + \mathcal{O}(x^2)$ in the vicinity of $\Sigma = \{(x, y) \in \Omega : x = 0\}$.

We are interested in finding $u^+$ that is the $L^2$-limit of $\lim_{\nu \to 0^+} u^\nu \in H^1(\Omega)$ s.t.,

$$\operatorname{div}((\alpha + i\nu)\nabla u^\nu) = f \quad \text{in } \Omega, \quad u^\nu\big|_{\partial\Omega} = 0.$$

The direct simulation of the above for small $\nu > 0$ is computationally inefficient as a very high mesh refinement is required to ensure convergence with respect to $\nu$. However, with a sufficiently regular $f$, $u^+$ can be represented as

$$u^+ = g^+(y)\left(\log|x| + i\pi \mathbf{1}_{\{x<0\}}\right) + u^+_{\text{reg}}.$$

where $g^+ \in H^1(\Sigma)$ and $u^+_{\text{reg}} \in H^1(\Omega)$ (i.e. *globally* $H^1$) [1, 3]. The challenge in designing a numerical method to compute $u^+$ is the fact that $L^2$-limit $u^- := \lim_{\nu \to 0^-} u^\nu \neq u^+$ and thus the method must be able to distinguish between the two solutions. One such variational formulation was proposed in [2], and further analysed and improved upon in [1], and simplified in [4].

In this work, we focus on the analysis and implementation of the formulation in [4] as well as propose a new alternate numerical method based on *non-local operators*.

## 3 Methods and Results

### 3.1 Variational Formulation from [4]

The key idea here is to use the splitting of $u^+$ into a singular and a regular part. Let us introduce the auxiliary notation.

For $j \in \{p, n\}$, $H^1_{1/2}(\Omega_j) := \overline{C^\infty(\overline{\Omega_j})}^{\|\cdot\|_{|\alpha|^{1/2}}}$, with the norm:

$$\|v\|^2_{|\alpha|^{1/2}} = \int_{\Omega_j} |v|^2 + \int_{\Omega_j} |\alpha||\nabla v|^2.$$

Define also,

$$H^1_{1/2,0}(\Omega) := \{u \in H^1_{1/2}(\Omega) \text{ s.t. } u\big|_{\partial\Omega} = 0\}.$$

Let $\mathrm{S}(x) := \log|x| + i\pi\mathbf{1}_{x<0}$ and for $g \in L^2(\Sigma)$, $s_g := S(x)g(y)$. We look for $u^+$ in the form

$$u^+ = u^+_{\text{reg}} + u^+_{\text{sing}},$$

where the pair $(u^+_{\text{reg}}, u^+_{\text{sing}})$ is such that

$$u^+_{\text{reg}} \in H^1_{1/2}(\Omega_{p,n}),$$

and

$$u^+_{\text{sing}}(x, y) = g^+(y)\,\mathrm{S}(x).$$

To set up the variational formulation we introduce the bilinear formulations $b_{\text{reg}} : H^1_{1/2}(\Omega) \times H^1_{1/2,0}(\Omega) \to \mathbb{C}$ and $b_{\text{sing}} : H^1_0(\Sigma) \times H^1_{1/2,0}(\Omega) \to \mathbb{C}$ as follows,

$$b_{\text{reg}}(u, v) := \int\limits_{\Omega} (\alpha\nabla u) \cdot \overline{\nabla v}\,\mathrm{d}x,$$

$$b_{\text{sing}}(g, v) := \int\limits_{\Omega} (\alpha\partial_y s_g)\overline{\partial_y v}\,\mathrm{d}x - \int\limits_{\Omega} \partial_x(\alpha\partial_x s_g)\overline{v}\,\mathrm{d}x.$$

Furthermore, we introduce for a piecewise regular $v \in H^1(\Omega\backslash\Sigma)$, $[v]_\Sigma := v(0+, y) - v(0-, y)$, $y \in \Sigma$. For $v \in H^1_{1/2}(\Omega)$, the jump is understood variationally as in [1].

We look for the solution $(u^+_{\text{reg}}, g^+) \in H^1_{1/2}(\Omega) \times H^1_0(\Sigma)$ such that for $f \in L^2(\Omega)$, for all $v \in H^1_{1/2,0}(\Omega)$ and $k \in H^1_0(\Sigma)$,

$$\begin{aligned} b(\mathbf{u}, \mathbf{v}) := b_{\text{reg}}(u^+_{\text{reg}}, v) + b_{\text{sing}}(g^+, v) \\ + \langle [u^+_{\text{reg}}]_\Sigma, k\rangle_{H^{-1}(\Sigma), H^1_0(\Sigma)} = \int\limits_{\Omega} f v \mathrm{d}x, \\ u^+_{\text{reg}} + g^+(y)\left(\log|x| + i\pi\mathbf{1}_{\{x<0\}}\right)\big|_{\partial\Omega} = 0. \end{aligned}$$

This formulation was shown to be injective in [4].

### 3.2 Non-local Formulation

An alternative approach to finding $u^+$ is to compute the flux $g^+$ first and then recover $u^+ = u^+_{\text{reg}} + u^+_{\text{sing}}$ by solving a well-posed (still sign-changing) problem for $u^+_{\text{reg}} \in H^1_{1/2}(\Omega)$.

$$\begin{aligned} \operatorname{div}(\alpha\nabla u^+_{\text{reg}}) &= -\operatorname{div}(\alpha\nabla(g^+\log|x| + i\pi\mathbf{1}_{\{x<0\}})) + f, \\ u^+_{\text{reg}}\big|_{\Gamma_p\cup\Gamma_n} &= -g^+(y)\left(\log|x| + i\pi\mathbf{1}_{\{x<0\}}\right), \\ u^+_{\text{reg}}\big|_{\Gamma_+\cup\Gamma_-} &= 0. \end{aligned}$$

The equation for $g^+$ below is inspired by the Green's formula from [1–3] and involves two operators.
First, consider $\mathcal{S} \in \mathcal{L}(H^1_0(\Sigma), H^{-1}(\Sigma))$ that maps $k \mapsto [u_k]_\Sigma$, where $u_k \in H^1_{1/2}(\Omega)$ solves,

$$\begin{aligned} \operatorname{div}(\alpha\nabla u_k) &= -\operatorname{div}(\alpha\nabla(k\log|x|)), \\ u_k\big|_{\Gamma_p\cup\Gamma_n} &= -k\log|x|\big|_{\Gamma_p\cup\Gamma_n}, \\ u_k\big|_{\Gamma_+\cup\Gamma_-} &= 0. \end{aligned}$$

Then, consider $\mathcal{F} \in \mathcal{L}(L^2(\Omega), H^{1/2}(\Sigma))$ that maps $f \mapsto [U_f]_\Sigma$ where $U_f \in H^1_{1/2}(\Omega)$ solves:

$$\begin{aligned} \operatorname{div}(\alpha\nabla U_f) &= f, \\ U_f\big|_{\partial\Omega} &= 0. \end{aligned}$$

It can be shown [3] that the singular part $g^+ \in H^1_0(\Sigma)$ satisfies the following equation:

$$(-i\pi + \mathcal{S}^*)g^+ = \mathcal{F}f.$$

The accompanying talk will expand on the computational methodologies and results.

# A note on the discrete Unmapped Tent Pitching for the heterogeneous wave equation

Marcella Bonazzoli[1,*], Gabriele Ciaramella[2], Ilario Mazzieri[2]
[1]Inria, Unité de Mathématiques Appliqueés, ENSTA, Institut Polytechnique de Paris, 91120 Palaiseau, France
[2]MOX, Dipartimento di Matematica, Politecnico di Milano, 20133, Milano, Italy
*Email: marcella.bonazzoli@inria.fr

**Abstract**

The Unmapped Tent Pitching (UTP) algorithm is a space–time domain decomposition method for the parallel solution of wave-type problems. We have recently extended UTP to heterogeneous settings and compared, at the continuous level, the computational cost of different space–time decompositions. In this note, we report discrete-level observations showing that the optimal decomposition strategy may differ from the one predicted by the continuous analysis.

## 1 Introduction

The UTP algorithm was designed for the parallel solution of the homogeneous wave equation in [2]. It is inspired by the Mapped Tent Pitching (MTP) algorithm [3]. While in MTP each physical tent is mapped onto a space-time Cartesian box in nD (rectangle in 1D) where local problems are solved, UTP avoids the mapping by computing the solution directly on a physical space-time box containing the tent. In [1], we extended UTP to the heterogeneous case:

$$\partial_{tt}u(x,t) = c(x)^2\partial_{xx}u(x,t),\ (x,t)\in\Omega\times(0,T),$$

coupled with homogeneous Dirichlet conditions on $\{0,L\}\times(0,T)$ and non-null initial conditions. Here, $\Omega = (0,L)$, $T>0$, $c(x) = c_1 > 0$ for $x\in(0,L/2]$ and $c(x) = c_2 > 0$ for $x \in (L/2,L)$. For simplicity, we take $c_1 > c_2$ and $T = L/(2c_1)$.

## 2 UTP algorithm and discrete setting

UTP is based on a *red-black* decomposition of $\Omega\times(0,T)$ into space-time subdomains. On $\Omega$ we set a grid $\{x_{\mathrm{R},j}\}_j$ having $m_1$ and $m_2$ uniformly distributed points in $(0,L/2)$ and $(L/2,L)$, resp. Thus, the *red* subintervals $I_{\mathrm{R},j} := (x_{\mathrm{R},j}, x_{\mathrm{R},j+1})$ have length $L_1 = L/(2m_1)$ and $L_2 = L/(2m_2)$ in $(0,L/2)$ and $(L/2,L)$, resp. The *black* subintervals $I_{\mathrm{B},j} := (x_{\mathrm{B},j}, x_{\mathrm{B},j+1})$, $j = 1,\dots,m_1 + m_2 - 1$, are defined with $x_{\mathrm{B},j} := x_{\mathrm{R},j} + L_1/2$, $j = 1,\dots,m_1$, $x_{\mathrm{B},j} = x_{\mathrm{R},j} + L_2/2$, $j = m_1+1,\dots,m_1+m_2$. While [1] also considers strategies with different time heights $H_1$ and $H_2$ in the two regions, here we assume $H_1 = H_2 = H$, with $H = L_1/(2c_1)$ determined by the smaller characteristic slope $|\pm 1/c_1|$. More precisely, if $k$ is the iteration index, we define the time subintervals $\tau_1 := (0,H)$ and $\tau_k := (t_{k,a}, t_{k,b}) := \big((k-2)H, \min(kH,T)\big)$ for $k\geq 2$. Space-time rectangles $\mathcal{T}_j^k$ are pitched alternately on the red and black subintervals: $\mathcal{T}_j^k := I_k\times\tau_k$, with $I_k := I_{\mathrm{R},j}$ for $k$ odd, $I_k := I_{\mathrm{B},j}$ for $k$ even. See, e.g., Fig. 1, where $c_1 = 2c_2$. Thus, given an approximation $u^{k-1}$ in $\Omega\times(0,T)$ satisfying initial and boundary conditions, the next $u^k$ is computed by solving the subproblems

$$\begin{aligned}&\partial_{tt}u_j^k(x,t) = c(x)^2\partial_{xx}u_j^k(x,t),\ (x,t)\in\mathcal{T}_j^k,\\ &u_j^k = u^{k-1},\ \text{on}\ \partial I_k\times\tau_k,\\ &u_j^k = u^{k-1},\ \partial_t u_j^k = \partial_t u^{k-1},\ \text{on}\ I_k\times\{t_{k,a}\},\end{aligned}$$

and then setting $u^k = u_j^k$ in $\mathcal{T}_j^k$ for the solved $j$ subdomains and $u^k = u_j^{k-1}$ elsewhere. The algorithm terminates when the exact solution is computed in the entire $\Omega\times(0,T)$; see the hatched regions in Fig. 1.

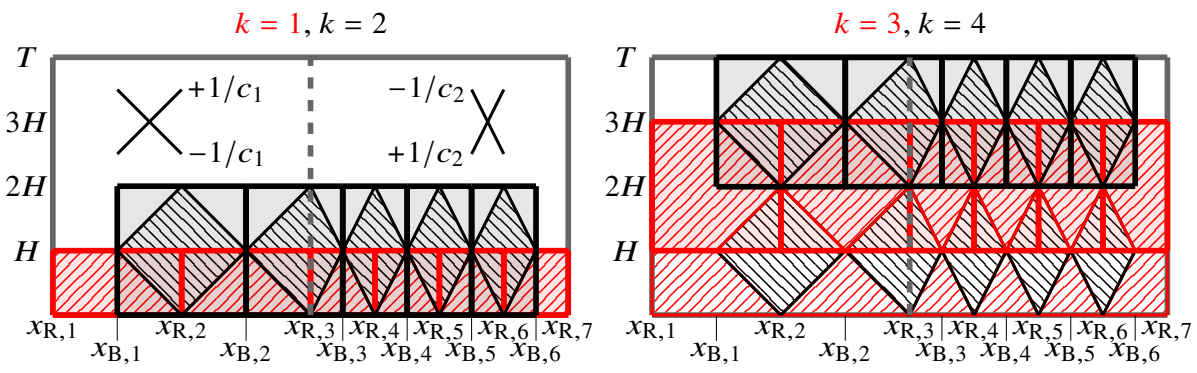

Figure 1: First 4 iterations for $4 = m_2 > m_1 = 2$ with $m_2/m_1 = L_1/L_2 = c_1/c_2$. Red and black boxes denote space–time subdomains; hatched regions indicate exact solution areas.

We discretize $\overline{\Omega}\times[0,T]$ by a uniform time grid of size $\Delta t$, and a piecewise-uniform space grid of size $h_1$ in $[0,L/2]$ and $h_2$ in $[L/2,L]$. The discrete space–time subdomains ${}^h\mathcal{T}_j^k$ are defined on this grid; see Fig. 2. The wave equa-

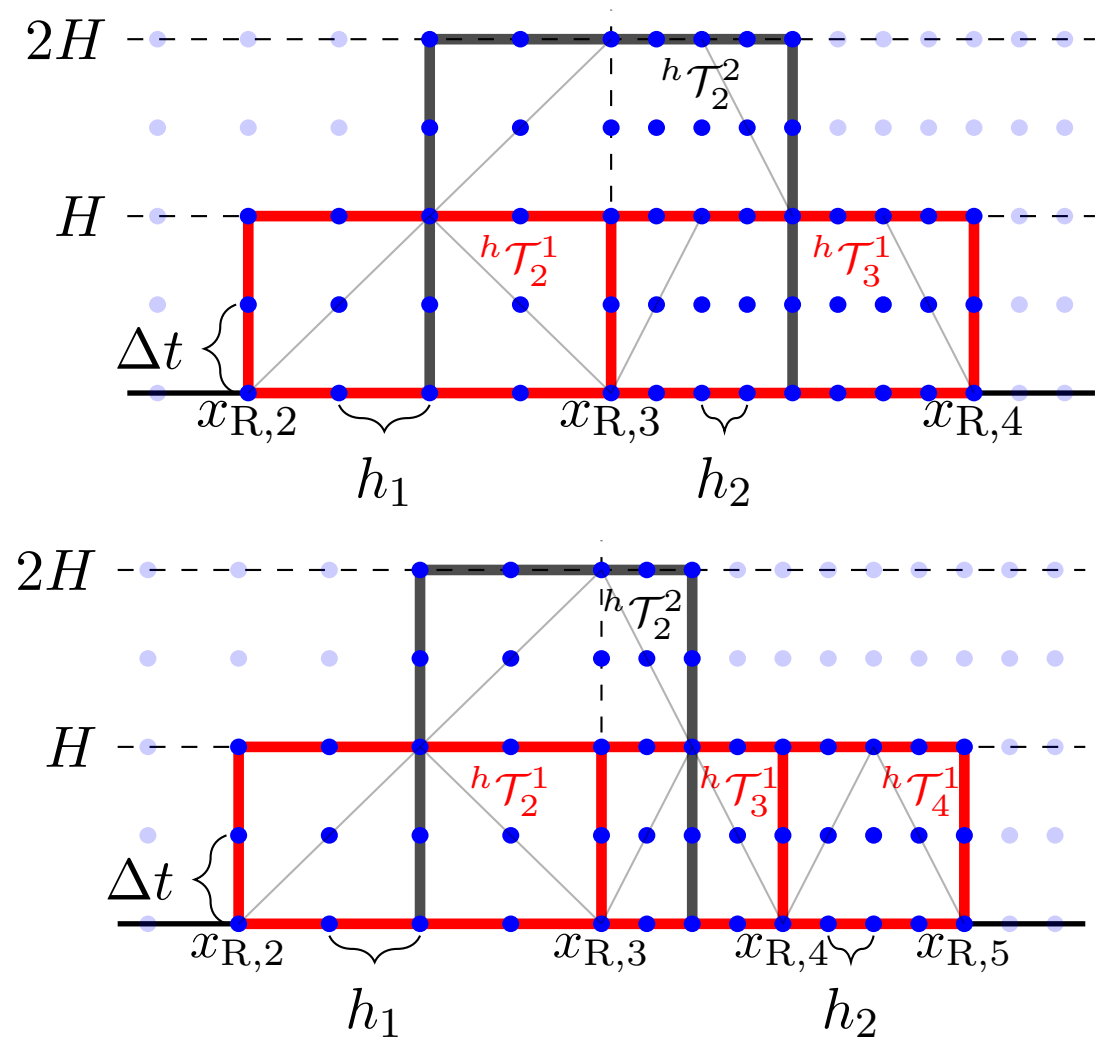

Figure 2: Discrete subdomains for $L_1 = L_2$ (top) and $m_2/m_1 = L_1/L_2 = c_1/c_2$ (bottom).

tion is discretized by the leapfrog scheme

$$u_{j,k}^{n+1,\ell}=2(1-\alpha)u_{j,k}^{n,\ell}-u_{j,k}^{n-1,\ell}+\alpha(u_{j,k}^{n,\ell-1}+u_{j,k}^{n,\ell-1}),$$

where $n$ and $\ell$ are time and space grid indices, $\alpha := \frac{c^2\Delta t^2}{h^2}$ is the CFL number, and $j$ and $k$ are the space-time subdomain indices.

## 3 Computational cost analysis

In [1], we considered $m_1+m_2$ parallel processes, one per subdomain solve, and modeled the cost of solve $(j,k)$ by the subdomain measure $|\mathcal{T}_j^k|$. In that continuous framework, the decomposition choice $H_1 = H_2$ and $m_1 = m_2$ (equivalently, $L_1 = L_2$) is optimal in terms of computational time and load balancing. At the discrete level, however, the situation changes due to the mismatch between the continuous characteristics (used in [1]) and the numerical ones. The latter stem from the explicit scheme, depend on the CFL condition, and are obtained by tracking stencil-based information propagation.

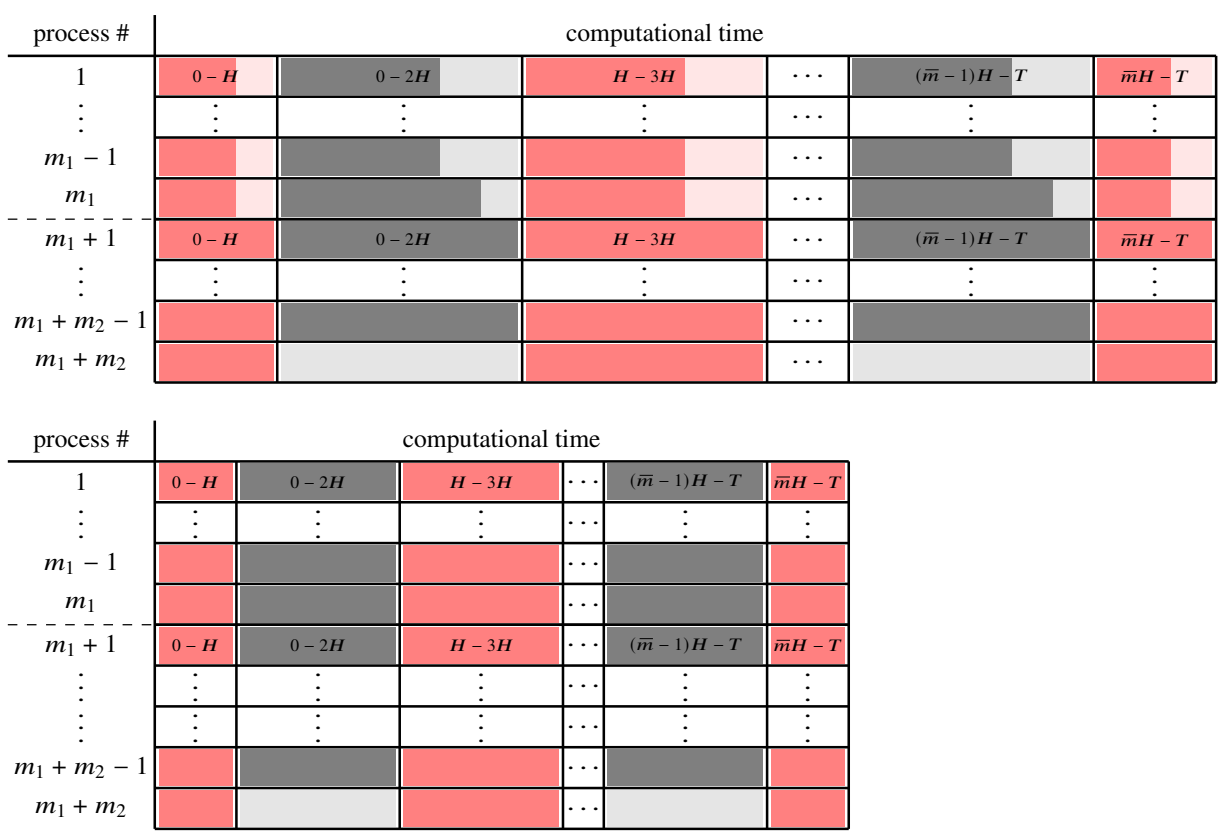

Figure 3: Pipelines for $L_1 = L_2$ (top) and $m_2/m_1 = L_1/L_2 = c_1/c_2$ (bottom).

In our discrete setting, the cost of solve $(j,k)$ is modeled by $|{}^h\mathcal{T}_j^k|$: the number of space–time grid points in ${}^h\mathcal{T}_j^k$, namely the number of variables computed by the solve in the discrete subdomain $(j,k)$, see Fig. 2. In [1], we showed that the height of a subdomain must not exceed the height of the solution front obtained after the first solve on that subdomain. Otherwise, additional iterations of the corresponding subproblem are required, slowing down the overall solution process. This remains true in the discrete setting, and we assume $H_1 = H_2$ (same number of time steps within the discrete subdomains). Moreover, to avoid dissipation and dispersion in the numerical solution, we assume $\alpha = 1$, guaranteeing that numerical and continuous characteristics coincide. This implies that $h_1 \neq h_2$, since $c_1 > c_2$. Now, only two cases are possible: $L_1 = L_2$ and $L_1/L_2 = c_1/c_2$. In the first case $|\mathcal{T}_j^k| = |\mathcal{T}_i^k|$ while $|{}^h\mathcal{T}_j^k| < |{}^h\mathcal{T}_i^k|$ for $j < m_1$ and $i > m_1$. In contrast, in the second case $|\mathcal{T}_j^k| > |\mathcal{T}_i^k|$ while $|{}^h\mathcal{T}_j^k| = |{}^h\mathcal{T}_i^k|$ for $j < m_1$ and $i > m_1$. See also Fig. 2. This means that, while in the continuous setting the subdomain costs are balanced for $L_1 = L_2$, the situation is opposite in the discrete setting, where the balance occurs for $L_1/L_2 = c_1/c_2$. This fact leads to the pipelines in Fig. 3, where the case $L_1/L_2 = c_1/c_2$ is the optimal one (at the expense of using more processes).

# A well-posed effective strain-gradient model for elastic waves in arbitrary periodic media

**Nicolas Auffray**[1], **Marc Bonnet**[2],[*], **Rémi Cornaggia**[2], **Giuseppe Rosi**[3]
[1]Sorbonne Université, CNRS, UMR 7190, Institut Jean Le Rond ∂'Alembert, F-75005 Paris, France
[2]POEMS (CNRS, Inria, ENSTA Paris), France
[3]CNRS, MSME, Univ Paris Est Creteil, F-94010 Creteil, France
[*]Email: mbonnet@ensta.fr

**Abstract** Following up on our Waves 2024 communication, we present an updated version of our development of a well-posed effective strain-gradient model for elastic waves in arbitrary periodic media and its use for both (frequency-domain) dispersion analysis and (transient) propagation simulation. Our presentation will include discussion of special cases (centrosymmetry, homogeneous elasticity or mass density) and numerical results on several 2D periodicity cell configurations with comparisons to results from Floquet-Bloch analyses and actual microstructures.

**1 Strain-gradient model.** We address waves propagating in an unbounded $d$-dimensional elastic medium whose elasticity tensor $\boldsymbol{C}$ and mass density $\rho$ are both spatially periodic and otherwise arbitrary ($L^\infty$, positive-definite). The displacement $\boldsymbol{u}$ satisfies the wave equation

$$\boldsymbol{\nabla}\cdot\{\boldsymbol{C}(\boldsymbol{X}):(\boldsymbol{e}[\boldsymbol{u}](\boldsymbol{X},t))\}=\rho(\boldsymbol{X})\boldsymbol{u}''(\boldsymbol{X},t) \tag{1}$$

where $\boldsymbol{e}[\boldsymbol{u}] := \frac{1}{2}(\boldsymbol{\nabla}\boldsymbol{u}+\boldsymbol{\nabla}^{\text{T}}\boldsymbol{u})$ is the linearizsed strain tensor. We aim at providing an *effective* model that reproduces key features of elastic wave propagation (anisotropy, dispersion ...) in the long-wavelength regime $\lambda \gg \ell$, for a periodicity cell $\ell\mathcal{Y}$ (linear size $\ell$, unit cell $\mathcal{Y}$).

Phenomenological *strain-gradient* (SG) models [4] provide effective wave equations such as

$$\begin{aligned}C^0_{ijk\ell}U_{k,j\ell}+M_{ijk\ell m}U_{k,\ell mj}-A_{ijk\ell mn}U_{k,j\ell mn}\\ =\rho^0U''_i+K_{ijk}U''_{j,k}-J_{ijk\ell}U''_{k,j\ell},\end{aligned} \tag{2}$$

(in component form and using Einstein's summation convention), where $\boldsymbol{U}$ is a macroscopic field approximating the mean motion, $\boldsymbol{C}^0$ and $\rho^0$ are the classical leading-order homogenized elasticity tensors and mass density [2], $\boldsymbol{M}$ and $\boldsymbol{K}$ are coupling tensors, and $\boldsymbol{A}$ and $\boldsymbol{J}$ are second-order elasticity and inertia tensors. Thanks to its higher-order terms, the model (2) has a wider range of application than the leading-order homogenized model based only on $(\boldsymbol{C}^0,\rho^0)$, in particular accounting for anisotropic behavior that using solely the tensor $\boldsymbol{C}^0$ may fail to capture [4].

**2 Asymptotic homogenization.** We define and compute the effective tensors in (2) using two-scale asymptotic homogenization in terms of the ratio $\varepsilon=\ell/\lambda \ll 1$. Following classical steps [1, 2], we define "slow" and "fast" normalized coordinates $\boldsymbol{x}=\boldsymbol{X}/\lambda$ and $\boldsymbol{y}=\boldsymbol{X}/\ell=\boldsymbol{x}/\varepsilon$, introduce a formal expansion

$$\boldsymbol{u}(\boldsymbol{X},t)=\textstyle\sum_{n=0}^{3}\varepsilon^n\boldsymbol{u}_n(\boldsymbol{x},\boldsymbol{y},t)+o(\varepsilon^3) \tag{3}$$

of $\boldsymbol{u}(\boldsymbol{X},t)$ in the wave equation (1), apply the chain rule $\boldsymbol{\nabla}_{\boldsymbol{X}}=\lambda^{-1}(\boldsymbol{\nabla}_{\boldsymbol{x}}+\varepsilon^{-1}\boldsymbol{\nabla}_{\boldsymbol{y}})$ to the differential operators, and finally obtain for $n\le 3$ separated variable solutions

$$\begin{aligned}\boldsymbol{u}_n(\boldsymbol{x},\boldsymbol{y},t)&=\boldsymbol{U}_n(\boldsymbol{x},t)\\&+\textstyle\sum_{j=1}^{n}\boldsymbol{\nabla}^j_{\boldsymbol{x}}\boldsymbol{U}_{n-j}(\boldsymbol{x},t)\bullet^j\boldsymbol{\chi}^j(\boldsymbol{y})\\&+\rho^0\omega^2\textstyle\sum_{j=0}^{n-2}\boldsymbol{\nabla}^j_{\boldsymbol{x}}\boldsymbol{U}_{n-2-j}(\boldsymbol{x},t)\bullet^j\boldsymbol{\zeta}^{j+2}(\boldsymbol{y})\end{aligned}$$

in terms of macroscopic fields $\boldsymbol{U}_n$, *cell functions* $\boldsymbol{\chi}^j$ that satisfy recursive elastostatic problems

$$-\boldsymbol{\nabla}_{\boldsymbol{y}}\cdot(\boldsymbol{C}:\boldsymbol{e}_{\boldsymbol{y}}[\boldsymbol{\chi}^n])=\boldsymbol{g}_n+\boldsymbol{\nabla}_y\cdot\boldsymbol{h}_n \tag{4}$$

posed on $\mathcal{Y}$ with $(\boldsymbol{g}_n,\boldsymbol{h}_n)$ depending on $\boldsymbol{C},\rho$ and previous cell functions $(\boldsymbol{\chi}^j)_{j<n}$, and similarly-defined *inertial cell functions* $\boldsymbol{\zeta}^j$ for which $\boldsymbol{g}_n,\boldsymbol{h}_n=\boldsymbol{0}$ if $\rho=\rho^0$. This formal procedure, supplemented by a "Boussinesq trick" [1] that features a tunable weight $\vartheta$, provides an one-parameter family of effective PDEs. The latter, of the form (2) and with tensors $(\boldsymbol{M},\boldsymbol{K})=O(\varepsilon)$ and $(\boldsymbol{A}_\vartheta,\boldsymbol{J}_\vartheta)=O(\varepsilon^2)$, are verified within a $o(\varepsilon^2)$ residual, for any $\vartheta$, by the second-order approximation $\boldsymbol{U}_\varepsilon:=\boldsymbol{U}_0+\varepsilon\boldsymbol{U}_1+\varepsilon^2\boldsymbol{U}_2$ of the macroscopic field. Notable special cases include that of $\mathcal{Y}$ centrosymmetric, for which $\boldsymbol{M}=\boldsymbol{0}$, $\boldsymbol{K}=\boldsymbol{0}$.

**3 Reciprocity and its exploitation.** The derivation and systematic use of *reciprocity identities* between cell problems, by means of combinations of their weak forms for suitably chosen test functions, allow to express tensors $(\boldsymbol{M},\boldsymbol{K})$ and $(\boldsymbol{A}_\vartheta,\boldsymbol{J}_\vartheta)$ explicitly in terms of cell functions $(\boldsymbol{\chi}^1,\boldsymbol{\chi}^2)$ (i.e. *not* $\boldsymbol{\chi}^3$) and $\boldsymbol{\zeta}^2$ (i.e. *not* $\boldsymbol{\zeta}^3$) and to thereby show that they have the symmetry properties expected of a SG model [4].

Setting the resulting $\vartheta$-dependent PDE (2) in the form $\mathcal{Q}_\vartheta\boldsymbol{U}+\mathcal{R}_\vartheta\boldsymbol{U}''=0$, the bilinear forms associated with the differential operators $\mathcal{Q}_\vartheta$ (of order 4) and $\mathcal{R}_\vartheta$ (of order 2) are respectively $H^2(\mathbb{R}^d)$-coercive and $H^1(\mathbb{R}^d)$-elliptic for any $\vartheta > \vartheta_0$, where $\vartheta_0$ does not depend on $\varepsilon$; in particular *positive* tensors $(\boldsymbol{A}_\vartheta, \boldsymbol{J}_\vartheta)$ are obtained.

Moreover, the threshold $\vartheta_0$ is easily computable by solving two small eigenvalue problems defined in terms of the "raw" tensor coefficients, of dimension $d^2(d+1)/2$ (resp. $d^2$), that yield the positivity threshold of $\boldsymbol{A}_\vartheta$ (resp. $\boldsymbol{J}_\vartheta$).

**4 Generalized Christoffel equation.** On trying a time-harmonic plane wave of the form $\boldsymbol{U}(\boldsymbol{x},t) = \boldsymbol{V} \exp\big(\mathrm{i}(\kappa \boldsymbol{n}\cdot\boldsymbol{x} - \omega t)\big)$, with given $\kappa$ and (unit) $\boldsymbol{n}$, in the effective SG PDE, propagation may only occur at angular frequencies $\omega = \omega(\boldsymbol{n},\kappa)$ such that $\mu := \rho^0\omega^2/\kappa^2$ and $\boldsymbol{V}$ solve

$$\begin{aligned}\big[\boldsymbol{Q}^0 + \mathrm{i}\delta\boldsymbol{Q}^1 + \delta^2\boldsymbol{Q}^2_\vartheta\big]\cdot\boldsymbol{V}\\ -\,\mu\big[\boldsymbol{I} + \mathrm{i}\delta\boldsymbol{R}^1 + \delta^2\boldsymbol{R}^2_\vartheta\big]\cdot\boldsymbol{V} = \boldsymbol{0},\end{aligned} \quad (5)$$

(with $\delta := \varepsilon\kappa\lambda = \ell\kappa$ and the matrices in (5) depending on $\boldsymbol{n}$), with polarizations $\boldsymbol{V}$ that are corresponding eigenvectors. The above-outlined selection rule for $\vartheta$ ensures that both matrices in the eigenvalue problem (5) are SPD. With $\kappa = 0$, (5) reduces to the standard Christoffel equation for the leading-order homogenized medium.

**Proposition 1** *Let $\mu_0 > 0$ be a standard Christoffel eigenvalue (solving (5) with $\kappa = 0$), and $\boldsymbol{V}_0$ a corresponding eigenvector. We have $\mu(\delta) = \mu_0 + O(\delta^2)$ unless $\mathcal{Y}$ is non-centrosymmetric and $\mu_0$ has multiplicity 2, in which case $\mu(\delta) = \mu^0 + O(\delta)$. Polarizations verify $\boldsymbol{V}(\delta) = \boldsymbol{V}_0 + O(\delta^2)$ for centrosymmetric $\mathcal{Y}$, $\boldsymbol{V}(\delta) = \boldsymbol{V}_0 + O(\delta)$ otherwise.*

**5 Solvability of initial-value problems.** We then consider initial-value problems

$$\begin{aligned}\mathcal{Q}_\vartheta\boldsymbol{U} + \mathcal{R}_\vartheta\boldsymbol{U}'' = 0 \quad \text{in } \mathbb{R}^d\times[0,T]\\ \boldsymbol{U}(0) = \boldsymbol{U}_0, \quad \boldsymbol{U}'(0) = \boldsymbol{V}_0\end{aligned} \quad (6)$$

for the effective strain-gradient wave equation. The previously-established symmetry, sign and coercivity properties of $\mathcal{Q}_\vartheta, \mathcal{R}_\vartheta$ allow to establish the well-posedness of problem (6) in finite-energy spaces, by applying the Hille-Yosida theorem [3] to a first-order reformulation of (6):

**Theorem 1** *Let $\boldsymbol{V}_0 \in \boldsymbol{H}^2(\mathbb{R}^d)$ and $\boldsymbol{U}_0 \in \boldsymbol{H}^2_{\mathcal{Q}}(\mathbb{R}^d)$, with $\boldsymbol{H}^2_{\mathcal{Q}}(\mathbb{R}^d) := \big\{\, \boldsymbol{W} \in \boldsymbol{H}^2(\mathbb{R}^d),\, \mathcal{Q}_\vartheta\boldsymbol{W} \in \boldsymbol{H}^{-1}(\mathbb{R}^d) \,\big\}$. Then, problem (6) has a unique solution* $\boldsymbol{U} \in C^1([0,T];\boldsymbol{H}^2(\mathbb{R}^d)) \cap C^0([0,T];\boldsymbol{H}^2_{\mathcal{Q}}(\mathbb{R}^d))$.

A similar well-posedness result is obtained for a forced-response version of problem (6).

**6 Numerical illustration.** Dispersion using the generalized Christoffel equation (compared with Floquet-Bloch analysis and leading-order homogenization) and a transient forced response (compared with leading-order homogenization and actual microstructure) based on the effective SG model are both illustrated in Fig. 1 on a $D_6$ (centrosymmetric) periodicity cell.

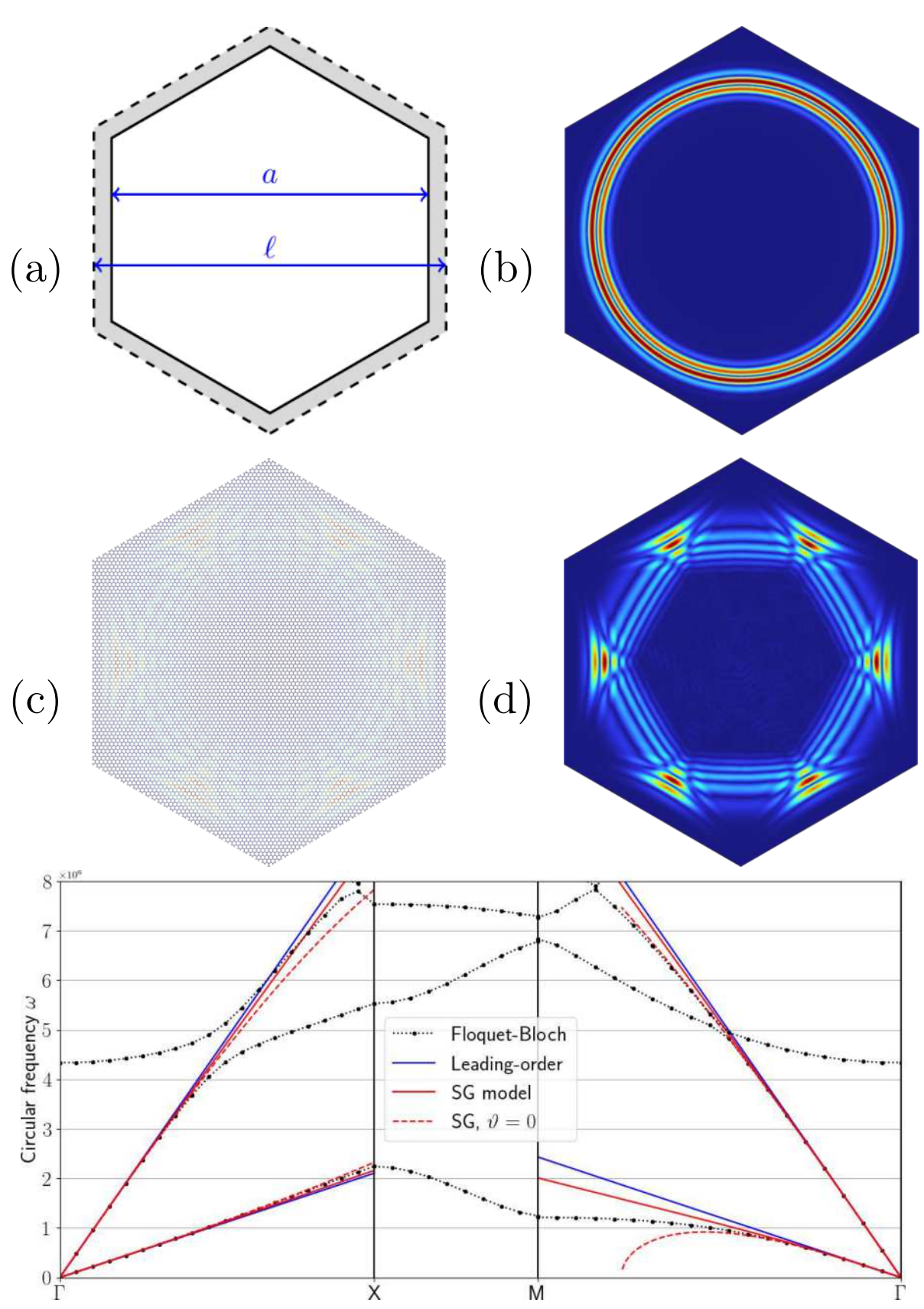

Figure 1: *Waves in microstructure generated by* $D_6$ *periodicity cell (a). Shear wavefront (displacement magnitude) at $t = 0.7\mu s$ (central frequency 80kHz) in (b) the leading-order homogenized Cauchy medium, (c) the microstructure and (d) the well-posed SG medium. Dispersion diagram and dispersion branches (bottom)*

## References.

# Green tensor for strain gradient elastodynamics

Marc Bonnet[1,*]
[1]POEMS (CNRS, Inria, ENSTA Paris), France
*Email: mbonnet@ensta.fr

**Abstract** The strain gradient elastic (SGE) model is an enriched form of linear elasticity featuring internal length scales to account in a continuum model for microstructure effects, and leads to a 4th order linear PDE. Green tensors for the general SGE PDE have been investigated over the last decade or so for elastostatics. We study here the 3D elastodynamic SGE Green tensor in the Laplace and frequency domains, obtaining general properties for anisotropic SGE media as well as closed-form expressions and radiation conditions for isotropic and entrosymmetric media.

**1 Anisotropic case.** We consider a 3D unbounded, homogeneous medium endowed with strain gradient elastic (SGE) properties [4, 5]. The SGE model is an enriched form of linear elasticity featuring internal length scales, to account for microstructure effects associated e.g. with architectured materials via a phenomenological continuum model [3]. Moreover, second-order two-scale homogenization of periodic heterogeneous elastic media produces effective 4th-order PDEs of SGE type, see our companion communication. Green tensors for the general SGE PDE have been investigated only over the last decade or so and for elastostatics [1, 2]. In this work, we investigate the SGE Green tensor for 3D elastodynamics.

Displacement fields $\boldsymbol{u}$ in linear SGE media are governed by PDEs of the form

$$-C^0_{ijk\ell}u_{k,j\ell} - M_{ijk\ell m}u_{k,\ell mj} + A_{ijk\ell mn}u_{k,j\ell mn} + \rho^0 u''_i + K_{ijk}u''_{j,k} - J_{ijk\ell}u''_{k,j\ell} = f_i \quad (1)$$

(in component form), where $f_i$ is a body force density, $\boldsymbol{C}$ is the classical elasticity tensor, $\rho$ is the mass density, $\mathbf{M}$ and $\mathbf{K}$ are coupling tensors, and $\mathbf{A}$ and $\mathbf{J}$ are second-order elasticity and inertia tensors. In addition to the usual positive-definiteness of $\boldsymbol{C},\rho$, the additional tensors must be such that $\boldsymbol{A}-\boldsymbol{M}^{\text{T}}\!:\!\boldsymbol{C}^{-1}\!:\!\boldsymbol{M} \in S^2\big(S^2(\mathbb{R})\otimes\mathbb{R}\big)$ and $\boldsymbol{J}-\frac{1}{\rho}\boldsymbol{K}^{\text{T}}\!\cdot\!\boldsymbol{K} \in S^2(\mathbb{R}\otimes\mathbb{R})$ are positive definite.

The sought Green tensor is $\boldsymbol{G} = \boldsymbol{G}^p \otimes \boldsymbol{e}_p$ where each $\boldsymbol{G}^p$ $(p = 1, 2, 3)$ solves (1) with $f_i = \delta(t)\delta(\boldsymbol{x})\delta_{ip}$. Applying the spatial Fourier transform and the Fourier-Laplace transform in time to (1) yields the representation

$$\boldsymbol{G}(\boldsymbol{x},\omega) = \frac{1}{8\pi^3}\int_H\int_{\mathbb{R}} \big[\xi^2\mathbf{Q} - \omega^2\mathbf{R}\big]^{-1} e^{\mathrm{i}\xi\hat{\boldsymbol{\alpha}}\cdot\boldsymbol{x}}\xi^2 \,\mathrm{d}\xi\,\mathrm{d}\hat{\boldsymbol{\alpha}}, \quad (2)$$

with $\boldsymbol{\xi} = \xi\hat{\boldsymbol{\alpha}}$ $(|\hat{\boldsymbol{\alpha}}| = 1)$, $H := \big\{\hat{\boldsymbol{\alpha}} : \hat{\boldsymbol{\alpha}}\boldsymbol{x} \geq 0\big\}$, $\omega = \omega + \mathrm{i}\gamma$ (so that $s = -\mathrm{i}\omega$ is the classical Laplace variable) and the Fourier symbol matrices

$$\mathbf{Q}(\xi,\hat{\boldsymbol{\alpha}}) := \mathbf{C}(\hat{\boldsymbol{\alpha}}) + \xi\mathbf{M}(\hat{\boldsymbol{\alpha}}) + \xi^2\mathbf{A}(\hat{\boldsymbol{\alpha}}),$$
$$\mathbf{R}(\xi,\hat{\boldsymbol{\alpha}}) = \rho\mathbf{I} + \xi\mathbf{K}(\hat{\boldsymbol{\alpha}}) + \xi^2\mathbf{J}(\hat{\boldsymbol{\alpha}}).$$

($\mathbf{C},\mathbf{A},\mathbf{J}$ real symmetric, $\mathbf{M},\mathbf{K}$ imaginary skew-symmetric). The matrix $\big[\xi^2\mathbf{Q}-\omega^2\mathbf{R}\big]^{-1}$ has six pairs $\pm\xi_i$ of poles with (say) $\text{Im}(\xi_i) > 0$, which for any non-real $\omega$ are all non-real by the constitutive assumptions. Due to its $O(\xi^{-2})$ decay, $\big[\xi^2\mathbf{Q}-\omega^2\mathbf{R}\big]^{-1}\xi^2 \in \boldsymbol{L}^1(\mathbb{R}^3)$, in marked contrast with the classical-elastodynamics case. More-explicit (closed-form in some cases) expressions of (2) are found by applying the residue theorem with poles $\xi_i$. Representation (2) or its residue-theorem allow to establish some general properties of $\boldsymbol{G}$, which largely echo previously known ones for the anisotropic elastostatic case [2], notably its *non-singular* nature:

**Proposition 1** *Assume $\boldsymbol{C},\boldsymbol{M},\boldsymbol{A}$ general SGE elastic tensors, $\boldsymbol{K},\boldsymbol{J}$ general SGE inertial tensors (all sign-suitable). Let $Im(\omega) > 0$. Then:*

*(a) $\boldsymbol{G}(\boldsymbol{\cdot},\omega)$ is continuous in $\mathbb{R}^3$, and in particular non-singular at $\boldsymbol{x} = \boldsymbol{0}$;*

*(b) $\boldsymbol{\nabla G}$ is finite (not uniquely defined) at $\boldsymbol{x} = \boldsymbol{0}$;*

*(c) $\boldsymbol{G}(\boldsymbol{x},\omega)$ decays exponentially as $|\boldsymbol{x}| \to \infty$.*

The time-harmonic Green tensor $\boldsymbol{G}_{\omega_0}(\boldsymbol{\cdot})$ (with prescribed angular velocity $\omega_0$) is obtained using the limiting-amplitude principle as

$$\boldsymbol{G}_{\omega_0}(\boldsymbol{\cdot}) = \boldsymbol{G}(\boldsymbol{\cdot},\omega_0).$$

The residue-theorem form of (2) with $\omega = \omega_0$ allows to obtain the following general properties:

**Proposition 2** *Let $\boldsymbol{C},\boldsymbol{M},\boldsymbol{A}$ and $\boldsymbol{K},\boldsymbol{J}$ as in Prop. 1. The time-harmonic Green tensor $\boldsymbol{G}_{\omega_0}(\boldsymbol{\cdot})$ has properties (a), (b) of Prop. 1. In addition: $\boldsymbol{G}(\boldsymbol{x},\omega_0) = O(|\boldsymbol{x}|^{-1})$ as $|\boldsymbol{x}| \to \infty$.*

**2 Isotropic and centrosymmetric case.** In this case, $\boldsymbol{A}$ and $\boldsymbol{J}$ are determined by up to 5 and 3 material parameters, respectively (two of each only remaining active in the governing PDE) while $\boldsymbol{M} = \mathbf{0}, \boldsymbol{K} = \mathbf{0}$. The spatial-Fourier matrix symbols have explicit expressions

$$\mathbf{Q}(\xi, \hat{\boldsymbol{\alpha}}) = \mu(1+\ell_{\mathrm{S}}^2\xi^2)(\mathbf{I} - \hat{\boldsymbol{\alpha}}\otimes\hat{\boldsymbol{\alpha}}) + (\lambda+2\mu)(1+\ell_{\mathrm{P}}^2\xi^2)\hat{\boldsymbol{\alpha}}\otimes\hat{\boldsymbol{\alpha}}$$

$$\mathbf{R}(\xi, \hat{\boldsymbol{\alpha}}) = \rho(1+\lambda_{\mathrm{S}}^2\xi^2)(\mathbf{I} - \hat{\boldsymbol{\alpha}}\otimes\hat{\boldsymbol{\alpha}}) + \rho(1+\lambda_{\mathrm{P}}^2\xi^2)\hat{\boldsymbol{\alpha}}\otimes\hat{\boldsymbol{\alpha}};$$

(where the internal lengths $\ell_\nu$, resp. $\lambda_\nu$ are material parameters for $\boldsymbol{A}$, resp. $\boldsymbol{J}$), and so does $\big[\xi^2\mathbf{Q} - \omega^2\mathbf{R}\big]^{-1}$. The latter has poles $\pm\xi_{\nu a}, \pm\xi_{\nu b}$ ($\nu = \mathrm{P}, \mathrm{S}$) given by

$$\xi_{\nu a}^2 = \frac{1}{2\ell_\nu^2}\big[\lambda_\nu^2\kappa_\nu^2 - 1 + \sqrt{\Delta_\nu}\big], \quad \xi_{\nu b}^2 = \frac{-\kappa_\nu^2}{(\xi_{\nu a}\ell_\nu)^2}$$

with $\Delta_\nu = (1-\lambda_\nu^2\kappa_\nu^2)^2 + 4\ell_\nu^2\kappa_\nu^2$ and $\kappa_{\mathrm{S}}^2 := \omega^2\rho/\mu$, $\kappa_{\mathrm{P}}^2 := \omega^2\rho/(\lambda+2\mu)$. This allows evaluating the spatial-Fourier representation (2) in closed form, to obtain

$$\boldsymbol{G}(\boldsymbol{x}) = A(\boldsymbol{x})\mathbf{I} + \frac{1}{\kappa_{\mathrm{S}}^2}\boldsymbol{\nabla\nabla}\big[B_{\mathrm{S}}(\boldsymbol{x}) - B_{\mathrm{P}}(\boldsymbol{x})\big]$$

with

$$A(\boldsymbol{x}) = \frac{1}{4\pi\mu x\sqrt{\Delta_{\mathrm{S}}}}\big[e^{\mathrm{i}x\xi_{\mathrm{S}a}} - e^{\mathrm{i}x\xi_{\mathrm{S}b}}\big],$$
$$B_\nu(\boldsymbol{x}) = \frac{\ell_\nu^2}{4\pi\mu x\sqrt{\Delta_\nu}}\big(\xi_{\nu a}^2 e^{\mathrm{i}x\xi_{\nu b}} - \xi_{\nu b}^2 e^{\mathrm{i}x\xi_{\nu a}}\big) \tag{3}$$

In the limit of vanishing strain-gradient tensors, at any $\boldsymbol{x} \neq \mathbf{0}$, we have

$$\big|\boldsymbol{G}(\boldsymbol{x}, \omega) - \boldsymbol{G}^0(\boldsymbol{x}, \omega)\big| \to 0 \quad \text{as } \sqrt{\ell_\nu^2 + \lambda_\nu^2} \to 0$$

($\boldsymbol{G}^0$: Green tensor of classical elastodynamics).

**Time-harmonic case** The limiting values of the residue theorem-relevant poles are

$$\xi_{\nu a} = (\sqrt{2}\ell_\nu)^{-1}\sqrt{\lambda_\nu^2\kappa_\nu^2 - 1 + \sqrt{\Delta_\nu}} \quad \in \mathbb{R},$$
$$\xi_{\nu b} = (\sqrt{2}\ell_\nu)^{-1}\mathrm{i}\sqrt{1 - \lambda_\nu^2\kappa_\nu^2 + \sqrt{\Delta_\nu}} \quad \in \mathrm{i}\mathbb{R}.$$

with $\kappa_\nu, \Delta_\nu$ evaluated for $\omega = \omega_0$. The far-field leading behavior of $\boldsymbol{G}$ then stems from

$$A(\boldsymbol{x}) = \frac{e^{\mathrm{i}x\xi_{\nu a}}}{4\pi\mu x\sqrt{\Delta_{\mathrm{S}}}} + o\Big(\frac{1}{x}\Big),$$
$$\boldsymbol{\nabla\nabla}B_\nu(\boldsymbol{x}) = -\frac{k_\nu^2\, e^{\mathrm{i}x\xi_{\nu a}}}{4\pi\mu x\sqrt{\Delta_\nu}}\hat{\boldsymbol{x}}\otimes\hat{\boldsymbol{x}} + o\Big(\frac{1}{x}\Big) \tag{4}$$

(with $x = |\boldsymbol{x}|$, $\hat{\boldsymbol{x}} = \boldsymbol{x}/x$). While equations (3), (4) are analogous in form to corresponding expressions for $\boldsymbol{G}^0$, the far-field behavior (4) notably differs from that of $\boldsymbol{G}^0$ (due to factors $\Delta_\nu^{-1/2}$ and differing wavenumbers except if $\ell_\nu = \lambda_\nu$). By contrast, for elastostatics, the far-field leading term of the SGE Green tensor was found in [2] to coincide with the classical anisotropic Green tensor.

The far-field approximations (4) allow to show that $\boldsymbol{G}$ satisfies the following generalization of Kupradze's elastodynamic radiation conditions:

**Proposition 3** *Each* $\boldsymbol{G}^p$ *verifies*

$$\begin{aligned}
&\big[\big(\boldsymbol{\Sigma}^p(\boldsymbol{x}) - \mathit{div}\,\boldsymbol{T}^p(\boldsymbol{x})\big)\cdot\hat{\boldsymbol{x}} \\
&\quad - (\lambda+2\mu)\mathrm{i}\xi_{\mathrm{P}a}(1+\ell_{\mathrm{P}}^2\xi_{\mathrm{P}a}^2)\boldsymbol{G}^p(\boldsymbol{x})\big]\cdot\hat{\boldsymbol{x}} = o(\tfrac{1}{x}),\\
&\big[\boldsymbol{T}^p(\boldsymbol{x})\cdot\hat{\boldsymbol{x}} - (\lambda+2\mu)\mathrm{i}\xi_{\mathrm{P}a}\boldsymbol{E}^p\big]\cdot\hat{\boldsymbol{x}} = o(\tfrac{1}{x}),\\
&\big[\big(\boldsymbol{\Sigma}^p(\boldsymbol{x}) - \mathit{div}\,\boldsymbol{T}^p(\boldsymbol{x})\big)\cdot\hat{\boldsymbol{x}} \\
&\quad - \mu\mathrm{i}\xi_{\mathrm{S}a}(1+\ell_S^2\xi_{\mathrm{S}a}^2)\boldsymbol{G}^p(\boldsymbol{x})\big]\cdot\boldsymbol{P}(\hat{\boldsymbol{x}}) = o(\tfrac{1}{x}),\\
&\big[\boldsymbol{T}^p(\boldsymbol{x})\cdot\hat{\boldsymbol{x}} - 2\mu\mathrm{i}\xi_{\mathrm{S}a}\boldsymbol{E}^p\big]\cdot\boldsymbol{P}(\hat{\boldsymbol{x}}) = o(\tfrac{1}{x}).
\end{aligned}$$

*in the far-field limit, with* $\boldsymbol{P}(\hat{\boldsymbol{x}}) := \boldsymbol{I} - \hat{\boldsymbol{x}}\otimes\hat{\boldsymbol{x}}$ *and* $\boldsymbol{E}^p$, $\boldsymbol{\Sigma}^p$, $\boldsymbol{T}^p := \boldsymbol{A}\therefore\boldsymbol{\nabla}\boldsymbol{E}^p$*: linearized strain, Cauchy stress, hyper-stress for* $\boldsymbol{G}^p$*.*

**3 Additional topics.** The oral presentation will include additional topics, in particular (a) the isotropic non-centrosymmetric case (involving nonzero coupling tensors $\boldsymbol{M}, \boldsymbol{K}$, each carrying one specific material parameter), (b) integral form of radiation conditions, following [6] and (c) mapping properties of volume potentials. In addition, regarding the transient case, we expect to show that $\boldsymbol{G}(\cdot, t)$ is the sum of a propagative term (with finite wave speed) and an evanescent term.

**References.**

# A Model Reduction Approach for Wave Phenomena using dG-Trefftz Methods

**Tobias Born**[1,*], **Karsten Urban**[1]
[1]Institute for Numerical Mathematics, Ulm University, Germany
*Email: tobias.born@uni-ulm.de

**Abstract**

Transport and wave-type problems are known to be hard problems for model reduction. Standard linear approaches cannot lead to an efficient reduction. Discontinuous Galerkin (dG) Trefftz methods offer variational formulations of parameterized problems whose spaces are inseparably interwoven with the parameter. By means of a traveling wave model problem, we investigate whether dG-Trefftz could be used for nonlinear model reduction in order to try to outperform standard linear approaches.

## 1 Model reduction of a traveling wave

Following [2], we consider a traveling wave problem: find $u$, such that for $(t,x)\in(0,1)^2$

$$\mu^{-2}\frac{\partial^2}{\partial t^2}u(t,x;\mu) - \frac{\partial^2}{\partial x^2}u(t,x;\mu) = 0 \qquad (1)$$

as well as initial and boundary conditions

$$u(0,x;\mu) = 0, \qquad \frac{\partial}{\partial t}u(0,x;\mu) = 0,$$
$$u(t,0;\mu) = \tanh(5t)^3, \quad u(t,1;\mu) = 0,$$

are fulfilled weakly. The parameter $\mu\in[0.3,2] =: \mathcal{P}$ denotes the wave speed. This describes a wave traveling from $x=0$ towards $x=1$ within a time period of 1. For $\mu<1$, the wave is too slow to reach $x=1$ in the given time period, for $\mu>1$, the wave is reflected at $x=1$, see Figure 1. Although seemingly simple, this model problem contains all relevant challenges for model reduction, namely a traveling wave phenomenon and completely different scenarios for different parameter ranges.

This model problem is an example of a parameterized PDE (PPDE). If one wishes to solve a PPDE for different values of the parameter often, very fast or on restricted computing devices, model reduction is a must. Linear reduction approaches, such as the Reduced Basis Method, [6], establish a low-dimensional, parameter-independent basis and approximate the solution of the PPDE for a given $\mu\in\mathcal{P}$ as superposition of these basis elements. The best approximation of linear reduction approaches is bounded by the Kolmogorov $n$-width $d_n(\mathcal{P})$, [6]. While $d_n(\mathcal{P})$ decreases exponentially for elliptic and parabolic problems, for transport and wave-type problems a slow polynomial decay must be expected, [3].

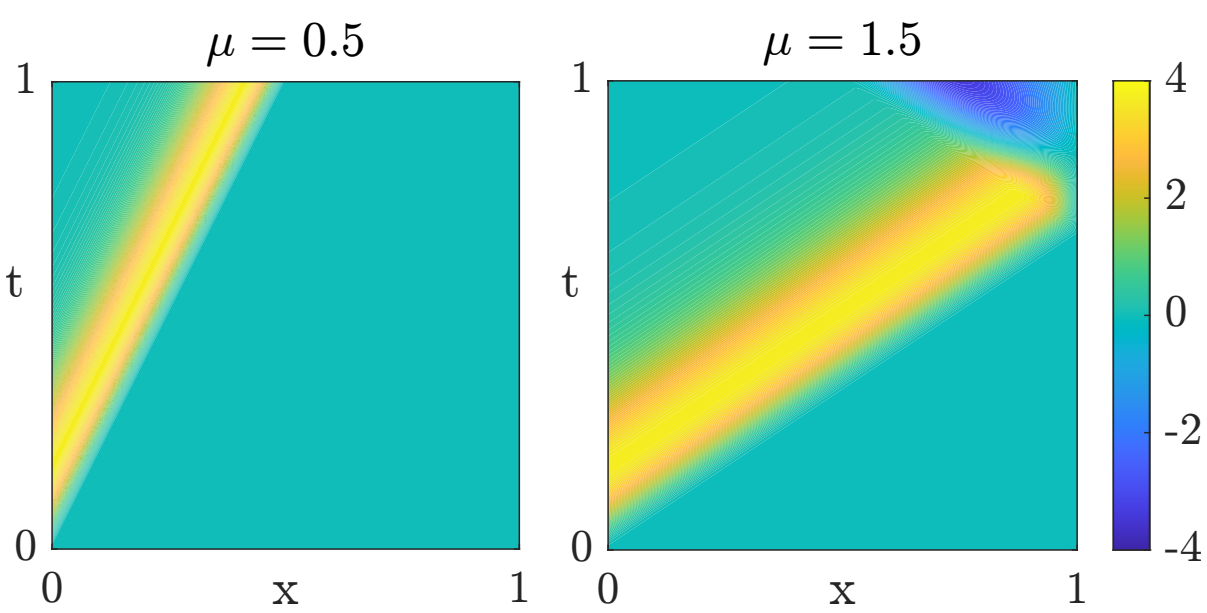

Figure 1: Plot of $\frac{\partial}{\partial x}u$ for different wave speeds.

Thus, for a reduction approach to be able to overcome the shortcomings of $d_n(\mathcal{P})$, it must not yield superpositions of fixed, parameter-independent basis elements. Instead, it has to involve nonlinearities. The question is what these must be to actually achieve faster convergence than $d_n(\mathcal{P})$. Our approach is to impose nonlinearity through a low-dimensional basis that depends on the parameter using a dG-Trefftz method.

## 2 A reduction approach using dG-Trefftz

The Trefftz space $\mathbb{T}(\mu)$ incorporates all functions that on each single mesh element are polynomial of a predefined degree and solve the equivalent first-order formulation of (1) exactly. Between mesh elements, Trefftz functions can be discontinuous, thus $\mathbb{T}(\mu)$ consists of many local bases inherently depending on the wave speed $\mu$.

In [5], a space-time dG formulation of the model problem is constructed, namely finding $(v,\sigma)\in\mathbb{T}(\mu)$ such that for all $(\omega,\tau)\in\mathbb{T}(\mu)$

$$\begin{aligned}&\mathcal{A}\big(\big[v(t,x;\mu),\, \sigma(t,x;\mu)\big],\, \big[\omega,\tau\big];\, \mu\big)\\ &\quad = \ell\big(\,[\omega,\tau]\,;\, \mu\big).\end{aligned} \qquad (2)$$

The parameter dependency of $\mathbb{T}(\mu)$ offers opportunities for nonlinear reduction. Take the basis expansion of the solution to (2) for $\mu \in \mathcal{P}$,

$$(v,\sigma)(\mu) = \sum\nolimits_k r_k(\mu)\Big(\varphi_k^{(v)}(\mu),\varphi_k^{(\sigma)}(\mu)\Big) \in \mathbb{T}(\mu).$$

The idea is to transform it from $\mathbb{T}(\mu)$ to $\mathbb{T}(\tilde{\mu})$ by keeping the coefficients but replacing the basis,

$$(v,\sigma)_{\mu\to\tilde{\mu}} = \sum\nolimits_k r_k(\mu)\Big(\varphi_k^{(v)}(\tilde{\mu}),\varphi_k^{(\sigma)}(\tilde{\mu})\Big) \in \mathbb{T}(\tilde{\mu}).$$

The resulting function of course differs from the solution to (2) for $\tilde{\mu} \in \mathcal{P}$, but we could hope that it is a good approximation of the latter.

Using the strong greedy method, $n$ samples $\mu^{(i)} \in \mathcal{P}$ are determined. For $\mu \in \mathcal{P}$, the reduced solution is defined as superposition of the transformations into $\mathbb{T}(\mu)$ of the solutions to (2) regarding the sample values $\mu^{(i)} \in \mathcal{P}$,

$$(v_n,\sigma_n)(\mu) = \sum\nolimits_{i=1}^{n} \eta_i(\mu)\,(v,\sigma)_{\mu^{(i)}\to\mu} \in \mathbb{T}(\mu).$$

The coefficients $\eta_i(\mu)$ are determined by moving (2) into the reduced space $U_n(\mu) \subset \mathbb{T}(\mu)$,

$$U_n(\mu) = \operatorname{span}\Big\{(v,\sigma)_{\mu^{(1)}\to\mu},\ \ldots,\ (v,\sigma)_{\mu^{(n)}\to\mu}\Big\}.$$

A key point is that for each $\mu \in \mathcal{P}$, the reduced solution is found in a different reduced space. This makes the reduction approach nonlinear and the Kolmogorov $n$-width no longer decisive.

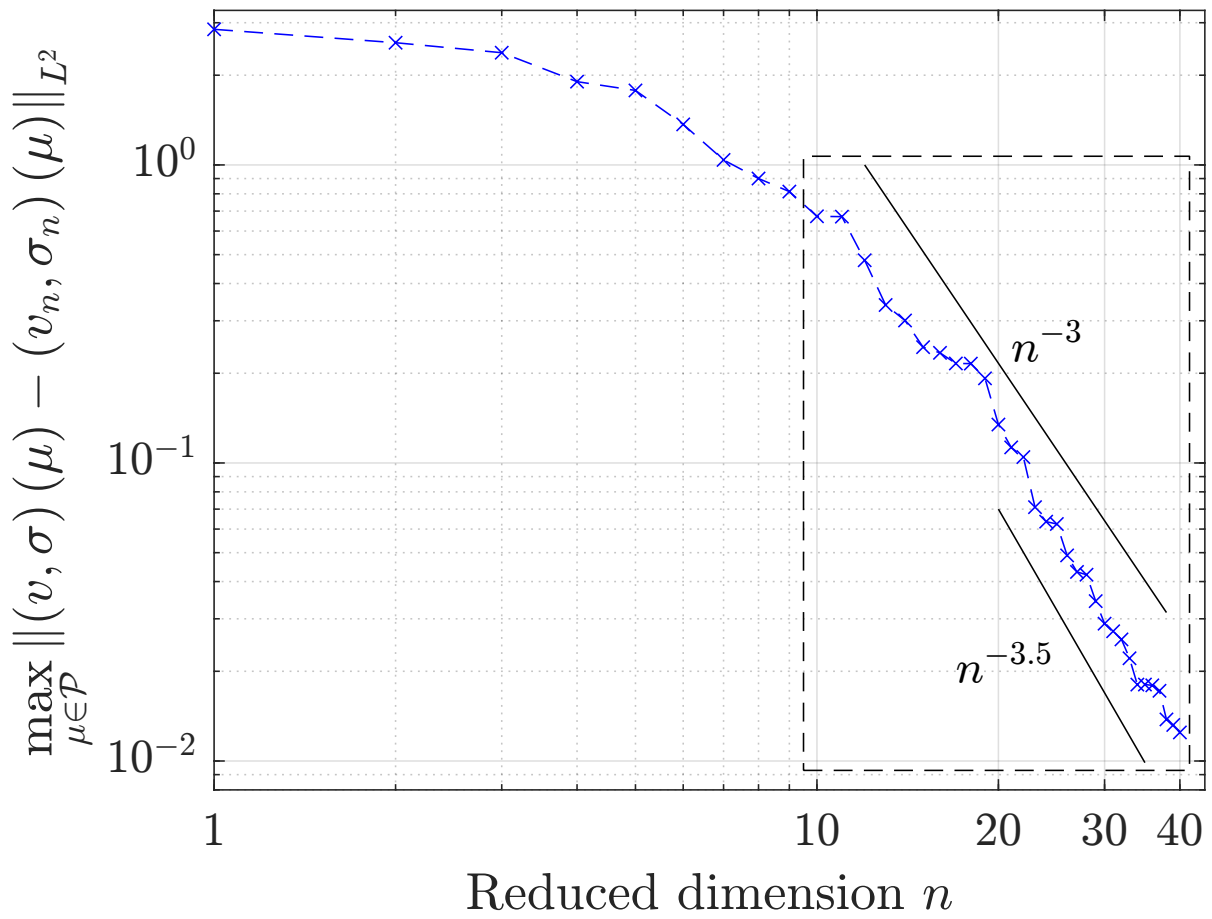

Figure 2: Strong greedy error of the nonlinear dG-Trefftz reduction approach.

Figure 2 shows the largest reduction error across the entire parameter set $\mathcal{P} = [0.3, 2]$ of the nonlinear dG-Trefftz reduction approach for different reduced space dimensions $n$. It seems to form a plateau and then converge at a polynomial rate between $n^{-3}$ and $n^{-3.5}$. In [2], a reduction approach for the traveling wave model problem was introduced that is linear and thus bounded by $d_n(\mathcal{P})$. Numerically, it exhibited a convergence rate of $n^{-3.5}$, which is in line with our nonlinear approach.

## 3 Conclusion and perspectives

The numerical results indicate that the nonlinear dG-Trefftz reduction approach, just as it is, is not able to outperform linear approaches and the shortcomings of the Kolmogorov $n$-width. Nevertheless, it still seems very promising to use dG-Trefftz methods within the model reduction chain, since their parameter-dependent spaces inherently create nonlinearity in the reduction process. Therefore, we are further investigating in what ways dG-Trefftz could be used for nonlinear model reduction.

The question inevitably arises, how nonlinearity in the reduction process could be learned, see [1,4]. Furthermore, it might be easier to reduce not the full state of a hyperbolic PPDE, but functional evaluations of the state. These could be averages, evaluations on boundaries, or similar.

# Diffraction by a rough thin layer on an arbitrary shaped object: the periodic and random cases

**Pierre Boulogne**[1,2,*], **Sonia Fliss**[2], **Laure Giovangigli**[2], **Justine Labat**[1]
[1]CEA-CESTA, Le Barp, France
[2]POEMS, CNRS, Inria, ENSTA, Institut Polytechnique de Paris, 91120 Palaiseau, France
*Email: pierre.boulogne@ensta.fr

**Abstract**

We are interested in the time-harmonic scattering by a bounded regular object coated with a thin rough layer. The layer's thickness is small compared to the object's characteristic size and the incident wavelength, which makes the computational cost of a direct simulation prohibitive. This scale separation can be exploited to derive effective models that avoid meshing the thin layer, reducing computational costs.
We study two cases: periodic and random stationary roughness.

## 1 Introduction and model

We consider a $C^\infty$ bounded connected domain of $\mathbb{R}^2$, denoted $\mathcal{D}_0$, illuminated by an incident plane wave $u_\mathrm{i}$. The equation governing the problem is

$$-\nabla\cdot(\mu_\delta\nabla u_\delta) - k^2\rho_\delta u_\delta = 0 \quad \text{in } \mathcal{D} := \mathbb{R}^2\setminus\overline{\mathcal{D}_0}$$

with a Dirichlet boundary condition

$$u_\delta = 0 \quad \text{on } \partial\mathcal{D}_0$$

and $u_\delta - u_\mathrm{i}$ satisfying the Sommerfeld radiation condition.
To describe the layer let us introduce the local coordinates tied to $\partial\mathcal{D}_0$, $(s,\nu)$, $s$ being the curvilinear abscissa along $\partial\mathcal{D}_0$ and $\nu$ the distance along the normal. The layer around the object is modelled *via* two reference functions $\mu$ and $\rho$, defined on the half-plane $\mathcal{B} := \mathbb{R}\times\mathbb{R}_+$. They are both equal to 1 for heights $Y > H > 0$. We consider two different settings: in (case 1) $X \mapsto (\mu,\rho)(X,\cdot)$ is $2\pi$-periodic; in (case 2) $(X,\omega) \mapsto (\mu_\omega,\rho_\omega)(X,\cdot)$ is stationary ergodic, where $\omega$ represents a realisation. We omit the $\omega$ notation in the following for readability. Let $\delta > 0$ be a small parameter. We finally have $(\mu_\delta,\rho_\delta)(s,\nu) := (\mu,\rho)\left(\frac{s}{\delta},\frac{\nu}{\delta}\right)$. For all $\alpha \geq 0$ let us define $\Gamma_\alpha := \{x\in\mathcal{D}, \mathrm{dist}(x,\mathcal{D}_0) = \alpha\}$, $\mathcal{D}^\infty_\alpha := \{x\in\mathcal{D}, \mathrm{dist}(x,\mathcal{D}_0) > \alpha\}$ and $\mathcal{D}^\alpha := \{x\in\mathcal{D}, \mathrm{dist}(x,\mathcal{D}_0) < \alpha\}$. The $(s,\nu)$ variables may not be defined far from $\Gamma_0$. In order to perform the asymptotic analysis, let us introduce $\eta > 0$ such that $(s,\nu)$ exists in $\mathcal{D}^\eta$.

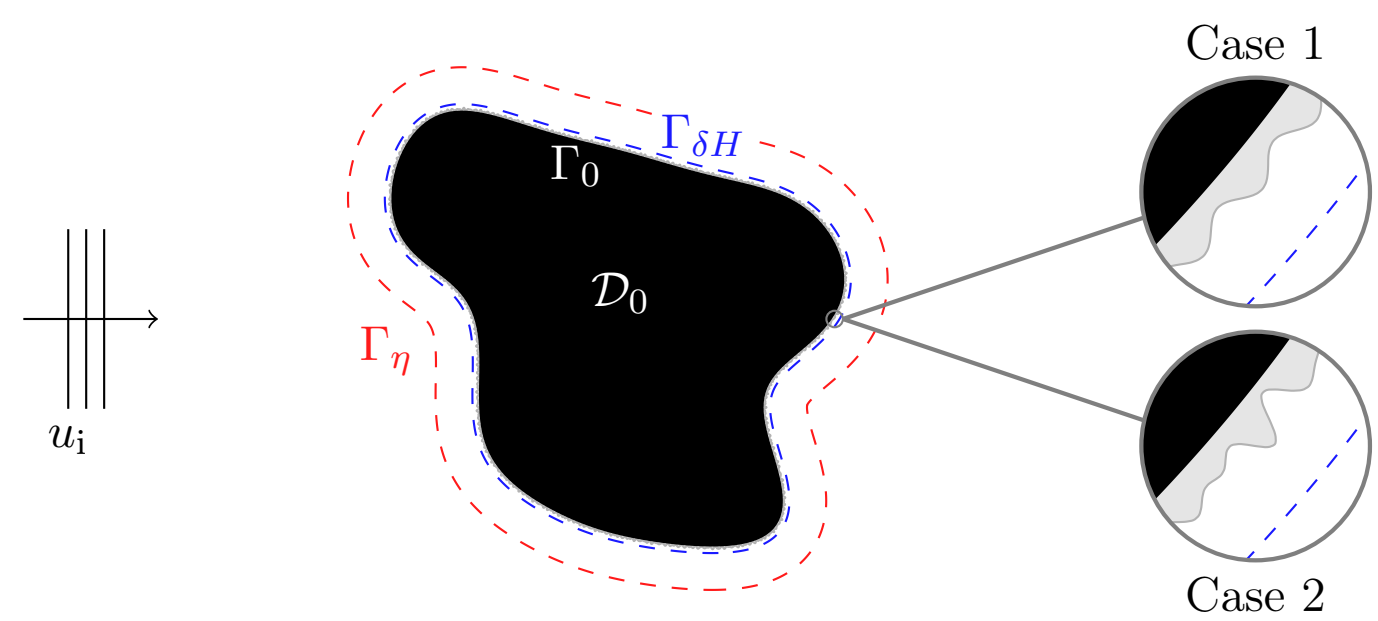

Figure 1: Representation of the model.

## 2 Formal asymptotic expansion

We assume that $u_\delta(s,\nu)$ can be written when $(s,\nu)\in\mathcal{D}^\eta$ as:

$$\begin{cases} \displaystyle\sum_{n=0}^{\infty} \delta^n\, u_n^{NF}\left(s;\frac{s}{\delta},\frac{\nu}{\delta}\right) & \text{if } \nu < \delta H, \\ \displaystyle\sum_{n=0}^{\infty} \delta^n\left[u_n^{NF}\left(s;\frac{s}{\delta},\frac{\nu}{\delta}\right) + u_n^{FF}(s,\nu)\right] & \text{if } \eta > \nu > \delta H. \end{cases} \tag{1}$$

The so-called near-field terms $u_n^{NF}$ depend on the curvilinear abscissa and microscopic variables $(X,Y) := \left(\frac{s}{\delta},\frac{\nu}{\delta}\right)$ to capture the quick variations of the solution near the layer. Depending on the case, they must satisfy for all $s$:

$$\text{(case 1) } X \mapsto u_n^{NF}(s;X,\cdot) \ 2\pi\text{-periodic}, \tag{2}$$
$$\text{(case 2) } (X,\omega) \mapsto u_{n,\omega}^{NF}(s;X,\cdot) \ \text{stationary}. \tag{3}$$

Furthermore, since $u_\delta$ does not oscillate far from the layer, the near-field terms should satisfy: $u_n^{NF}(\cdot;\cdot,Y) \underset{Y\to\infty}{\to} 0$. The so-called far-field terms $u_n^{FF}$ depend only on the macroscopic variables and capture the behaviour of the solution far from the object. By substituting the ansatz of the solution (1) in the problem satisfied by $u_\delta$ we get that the far-field terms satisfy:

$$-\Delta u_n^{FF} - k^2 u_n^{FF} = 0 \quad \text{in } \mathcal{D}_{\delta H}^{\infty}$$

and $u_n^{FF} - \delta_{n,0} u_{\mathrm{i}}$ is outgoing. The boundary condition on $\Gamma_{\delta H}$ that is missing is provided by the near-field terms. These latter satisfy Laplace-type problems, parametrized by the macroscopic variable $s$, of the form

$$\begin{cases} -\nabla_{X,Y}\cdot\left(\mu\nabla_{X,Y}u_n^{NF}\right) = F_{n-1}^{NF}(s) \text{ in } \mathcal{B}, \\ u_n^{NF} = 0 \quad \text{on } \Gamma := \mathbb{R}\times\{0\}, \\ [\![u_n^{NF}]\!]_{Y=H} = -u_n^{FF}|_{\Gamma_{\delta H}}(s), \\ [\![\partial_Y u_n^{NF}]\!]_{Y=H} = -\partial_\nu u_{n-1}^{FF}|_{\Gamma_{\delta H}}(s), \\ u_n^{NF} \text{ satisfies (2) or (3)}, \end{cases} \tag{4}$$

where $F_{n-1}^{NF}$ depends on the previous near-field terms.

## 3 Periodic case

In the first case, we show that (4) admits one and only one solution satisfying $\int_0^{2\pi}\int_0^{\infty}\frac{|v|^2}{1+Y^2} + |\nabla v|^2 \mathrm{d}Y\mathrm{d}X < \infty$ because $F_{n-1}^{NF}$ decreases fast enough. Furthermore, the solution converges exponentially fast to a constant when $Y\to\infty$. By imposing this constant to be 0 we obtain a compatibility condition which relates the traces of the far-field terms:

$$u_n^{FF}|_{\Gamma_{\delta H}} = \sum_{j=0}^{n-1}\sum_{i=0}^{n-j-1} c_{n-j-1,i}(s)\partial_\nu\partial_s^i u_j^{FF}|_{\Gamma_{\delta H}}, \tag{5}$$

where $c_{i,j}$ are functions depending both on the layer properties and the geometry of the object. The exponential decay of the near-field terms ensures that (4) is well-posed for every $n$.

**Theorem:** For every compact $K$ at a strictly positive distance from $\mathcal{D}_0$, for every $n\in\mathbb{N}$,

$$\left\| u_\delta - \sum_{j=0}^{n}\delta^j u_j^{FF}\right\|_{H^1(K)} = O\left(\delta^{n+1}\right)$$

## 4 Stationary case

In the second case, similarly, (4) admits one and only one solution satisfying $\mathbb{E}\left[\int_0^{+\infty}\frac{|v|^2}{1+Y^2} + |\nabla v|^2\right] < +\infty$ for $n=0$ and 1. The solutions still converge to a constant as $Y$ tends to infinity. Imposing that this constant equals 0 we obtain boundary conditions for the far-field terms.

For $n=0$, the solution of (4) is shown to be trivial, giving the condition $u_0^{FF}|_{\Gamma_{\delta H}} = 0$.

For $n=1$, $F_0^{NF} = 0$. Let us introduce the corrector $\widetilde{W}_1$ as the unique solution of

$$\begin{cases} -\nabla\cdot\left(\mu\nabla\widetilde{W}_1\right) = 0 \quad \text{in } \mathcal{B} \\ \widetilde{W}_1 = 0 \quad \text{on } \Gamma \\ [\![\widetilde{W}_1]\!]_{Y=H} = 0, \quad [\![\partial_Y\widetilde{W}_1]\!]_{Y=H} = -1 \\ (X,\omega)\mapsto\widetilde{W}_1(X,\cdot) \text{ is stationary} \end{cases}$$

such that $\mathbb{E}\left[\int_0^{+\infty}\frac{|v|^2}{1+Y^2} + |\nabla v|^2\right] < +\infty$. Let us denote $c_1 := \lim_{Y\to\infty}\widetilde{W}_1$ that we show to be independant of $\omega$. We have $u_1^{NF} \underset{Y\to\infty}{\to} 0$ if and only if $u_1^{FF} = \partial_\nu u_0^{FF} c_1$ on $\Gamma_{\delta H}$ and then $u_1^{NF} = \left(\widetilde{W}_1 - c_1\chi_{Y>H}\right)\partial_\nu u_0^{FF}|_{\Gamma_{\delta H}}$. The following error estimate under some quantitative assumptions on $\mu_\omega$ is proved in a close setting [1]:

$$\sqrt{\mathbb{E}\left[\left\|u_0^{FF} + \delta u_1^{FF} - u_\delta\right\|_{H^1(K)}^2\right]} = O\left(\delta^{\frac{3}{2}}\left|\log(\delta)\right|^{\frac{1}{2}}\right). \tag{6}$$

Since $u_1^{NF}\to 0$ slowly, $F_1^{NF}$ does not decay fast enough so that $u_2^{NF}$ is not well defined. In consequence, we cannot derive a boundary condition for $u_2^{FF}$ and we cannot improve the effective model.

## 5 Numerical simulation

We perform numerical simulations of effective solutions. We observe the expected convergence rates (5) for the periodic case and slightly better than (6) for the stationary one.

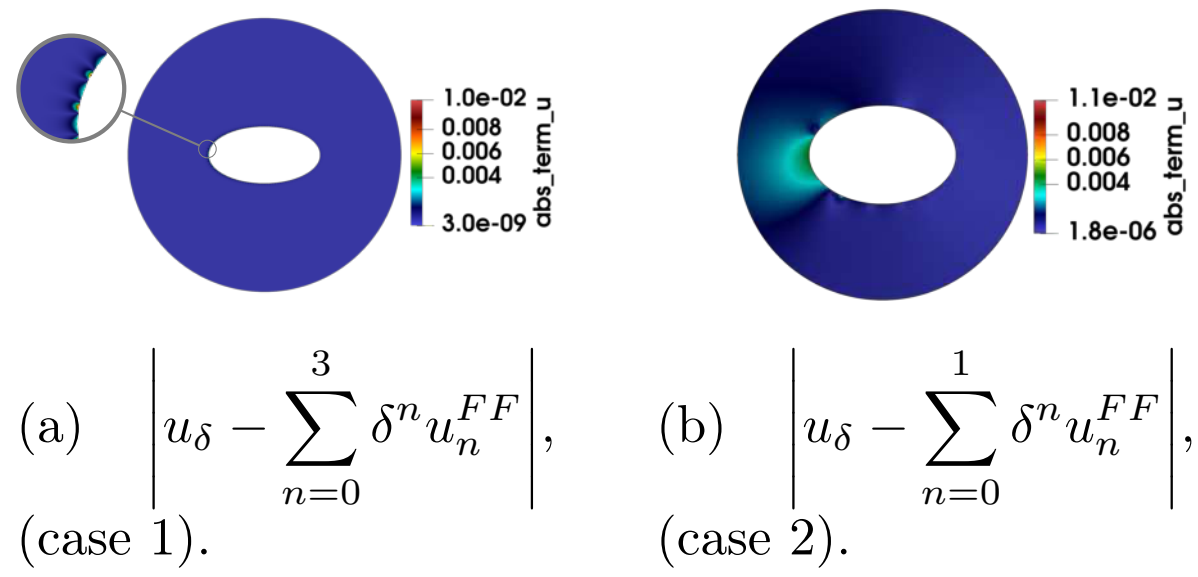

(a) $\left|u_\delta - \sum_{n=0}^{3}\delta^n u_n^{FF}\right|$, (case 1). (b) $\left|u_\delta - \sum_{n=0}^{1}\delta^n u_n^{FF}\right|$, (case 2).

# Localized Time-Periodic Patterns in a 1D Lattice

**Jason J. Bramburger**[1,*], **Eric Berglund**[2], **Björn Sandstede**[3]
[1]Department of Mathematics and Statistics, Concordia University, Canada
[2]One Health Trust, Washington, DC, USA
[3]Division of Applied Mathematics, Brown University, USA
*Email: jason.bramburger@concordia.ca

## Abstract

Motivated by numerical continuation studies of coupled mechanical oscillators, this talk presents recent results on branches of localized time-periodic solutions in one-dimensional chains of coupled oscillators. We consider Ginzburg–Landau equations with nonlinearities of Lambda–Omega type and establish the existence of spatially localized synchrony patterns in the regime of weak coupling and weak amplitude dependence of the intrinsic oscillator periods. Depending on the form and strength of the coupling, these localized states organize either into a discrete stack of isolas or along a single connected snaking branch, revealing distinct bifurcation structures underlying the emergence of localized, breather-type wave patterns.

## 1 Introduction

Localized oscillatory patterns of synchrony arise in many natural and engineered bistable systems, where individual units admit two coexisting stable states and coupling mediates transitions between them. In a typical setting, identical bistable elements—each capable of resting in one of (at least) two attractors or undergoing state-dependent oscillatory excursions—are arranged on a network or spatial domain. When uncoupled, each unit evolves independently toward one of its stable states; however, coupling can induce coordinated switching and sustained oscillations that remain spatially confined, forming coherent localized patches of synchrony rather than global phase locking.

Such dynamics appear in neural tissue and intracellular signaling, and are closely related to breather solutions in lattice and continuum models from mathematical physics, where nonlinearity and bistability generate time-periodic, spatially localized oscillations. From a wave propagation perspective, these states can be viewed as pinned, time-periodic wave packets arising from the interplay of nonlinearity, coupling, and lattice effects.

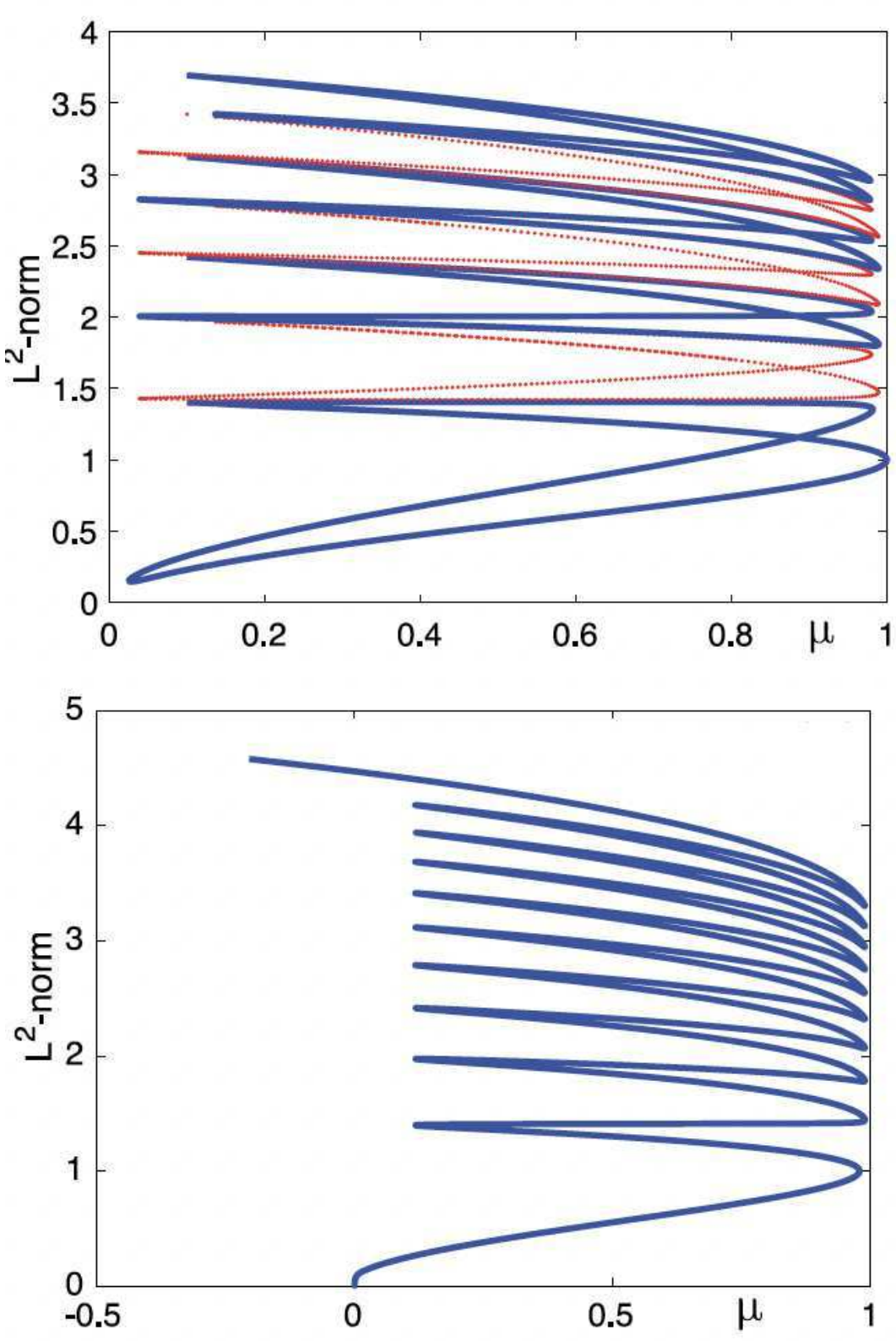

Figure 1: Numerical continuation of isola (top) and snaking (bottom) branches of (1) with conservative ($c = \mathrm{i}$) and dissipative ($c = 1$) coupling, respectively, for $N = 10$ nodes with coupling strength $\varepsilon = 0.01$. Images originally appear in [2].

## 2 Model

Complicated branches of localized synchrony patterns were identified in [1] using a one-mode Fourier truncation within a harmonic balance approximation of weakly coupled mechanical oscillators with linear springs and nonlinear dampers. The resulting reduced model takes the form

$$\dot{Z}_n = f(|Z_n|, \mu)Z_n + \varepsilon c(Z_{n+1} + Z_{n-1} - 2Z_n), \quad (1)$$

where $Z_n \in \mathbb{C}$ for $n = 1, \ldots, N$, the coupling strength satisfies $0 < \varepsilon \ll 1$, and $f$ is a complex-valued nonlinearity depending on a parameter $\mu$. In the mechanical setting of [1], the coupling constant is $c = \mathrm{i}$, corresponding to conservative coupling. When $c = 1$, the coupling is dissipative and (1) can be viewed as a generalized discrete Ginzburg–Landau equation.

A prototypical choice of nonlinearity is $f(r, \mu) = -\mu + 2r^2 - r^4$, the normal form for a subcritical Hopf bifurcation. Numerical continuation (Figure 1) shows that this system supports localized synchrony patterns whose bifurcation structure depends strongly on the coupling: for $c = \mathrm{i}$, solutions organize into discrete stacks of closed curves (isolas), whereas for $c = 1$, they lie on connected snaking branches. The aim of this work is to explain this distinct organization of localized states in parameter space.

## 3 Methods and Results

Our analysis begins in the weak-coupling regime $0 < \varepsilon \ll 1$, where the system can be viewed as a perturbation of decoupled bistable oscillators undergoing a subcritical Hopf bifurcation. In this limit, spatially uniform equilibria and periodic orbits are readily characterized, and localized synchrony patterns correspond to spatially inhomogeneous steady states of the lattice dynamical system obtained in a suitable rotating frame. We exploit the small parameter $\varepsilon$ to derive reduced equations governing the amplitude profile and to identify the leading-order balance between onsite bistability and discrete diffusion that enables spatial localization.

To establish existence, we combine spatial dynamics techniques with perturbative arguments. Viewing the lattice index $n$ as a discrete spatial variable, localized states are constructed as homoclinic orbits to small-amplitude backgrounds in an associated reversible map for conservative coupling and in a gradient-like setting for dissipative coupling. In the conservative case $c = \mathrm{i}$, the absence of dissipation leads to a constrained phase–amplitude structure and isolates families of solutions that persist as closed curves in parameter space. In contrast, for dissipative coupling $c = 1$, variational structure and pinning mechanisms give rise to snaking behavior, with solution branches oscillating back and forth across a narrow parameter interval.

These analytical results clarify the bifurcation organization observed numerically. In particular, we show how the interplay between subcriticality, weak coupling, and the nature of the coupling constant selects either disconnected isolas or a single connected snaking branch. The theory captures the scaling of solution amplitude and spatial width with $\varepsilon$, provides conditions for persistence of localization, and explains the transition between conservative and dissipative regimes within a unified framework.

## 4 Discussion and Outlook

The results presented here provide a unified explanation for the emergence and organization of localized synchrony patterns in weakly coupled bistable oscillator chains. By linking the bifurcation structure to the underlying coupling mechanism, we clarify how conservative interactions favor disconnected families of localized states, while dissipative coupling promotes snaking and parameter pinning. This distinction highlights the fundamental role of the phase-amplitude structure and the energy balance in shaping the global solution landscape.

Several directions remain open. An important next step is to move beyond the weak-coupling regime and determine the robustness of these bifurcation structures for moderate coupling strength and in higher-dimensional lattices. It would also be natural to investigate the stability and dynamical selection mechanisms for the localized states identified here, as well as their persistence under heterogeneity or external forcing. Such extensions would further connect the present theory with experimentally observable localized oscillations and breather-type phenomena in mechanical and physical lattice systems.

# Fast Singular-Kernel Convolution on General Non-Smooth Domains via Truncated Fourier Filtering

**Jinghao Cao**[1,*], **Oscar Bruno**[1]

[1]Computing and Mathematical Sciences, California Institute of Technology, Pasadena, USA

*Email: jinghao.cao@caltech.edu

## Abstract

This paper [1] develops a fast, high-order numerical method for evaluating singular convolution integrals

$$\int_\Omega K(x-y)\,\phi(y)\,\mathrm{d}y, \qquad x\in\Omega\subset\mathbb{R}^2,$$

on *general non-smooth domains* (including corners), focusing on the logarithmic kernel $K(z) = \log|z|$. The method extends the recently introduced Truncated Fourier Filtering (TFF) framework—originally designed for smooth kernels—to singular kernels, while preserving a grid-based, FFT-compatible structure and achieving *super-algebraic convergence.*

## 1 Introduction

The rapid and accurate evaluation of convolutions with singular kernels plays a central role in numerous scientific and engineering applications.

Existing approaches, while highly effective in specialized settings, often rely on geometry-adapted discretizations, local analytic expansions, recursive refinements, or assumptions of globally smooth densities, which limit robustness on complex or non-smooth domains. As a result, achieving high-order accuracy typically entails significant geometric or analytic complexity.

The proposed method instead employs truncated Fourier expansions of the domain characteristic function together with Fourier representations of products of characteristic and singular factors, yielding a simple FFT-compatible scheme that attains high-order accuracy on arbitrary geometries without boundary-fitted meshes or local singularity treatments.

**Theorem 1 (Truncated Fourier Filtering)** *Let $f$ be a smooth $2\pi$-periodic function and let $g$ be a $2\pi$-periodic function that may be discontinuous or weakly singular, with Fourier coefficients $\widehat{g}_k$. For $F\in\mathbb{N}$, let*

$$g_F(x) = \sum_{|k|\le F} \widehat{g}_k e^{ikx}$$

*denote the truncated Fourier series of $g$, and let $Q_N$ denote the trapezoidal rule with $N$ equispaced nodes on $[0,2\pi]$, with $N\ge bF$, for $b>1$.*

*Then, for every $p\in\mathbb{N}$, there exists a constant $C_p$ (independent of $F$ and $N$) such that*

$$\left|\int_0^{2\pi} f(x)g(x)\,\mathrm{d}x - Q_N[f\,g_F]\right| \le C_p F^{-p}.$$

The result shows that the rapid decay of Fourier coefficients of the smooth factor compensates for Gibbs oscillations introduced by truncation of the non-smooth factor.

## 2 Window decomposition

The method presented in this paper utilizes a smooth radial window $W_1$ to split the convolution integral into

$$\int_\Omega \log|x-y|\phi(y)\,\mathrm{d}y = I_1(x) + I_2(x),$$

where $I_1$ captures the localized singularity near $y=x$, and $I_2$ has a globally smooth integrand. The two terms are treated by different but complementary TFF strategies.

### 2.1 Localized singular integral $I_1$

Using polar coordinates centered at $x$,

$$I_1(x) = \int_{-\pi}^{\pi}\int_0^{d(\theta)} r\log|r|\,\varphi(r,\theta)W_1(r)\,\mathrm{d}r\,\mathrm{d}\theta,$$

where $d(\theta)$ is the distance to $\partial\Omega$ along a ray.

Three geometric regimes arise: (i) the window lies fully inside $\Omega$, (ii) it intersects $\partial\Omega$ smoothly, and (iii) it intersects a corner. In cases (i)–(ii), careful symmetrization and truncation of the Fourier series of $|r|\log|r|$ times a characteristic function lead to a smooth angular integrand, and standard trapezoidal rules in $(r, \theta)$ converge superalgebraically. In case (iii), angular smoothness breaks down; the remedy is to apply TFF also in the angular variable by replacing angular characteristic functions with their truncated Fourier series. This restores superalgebraic convergence even in the presence of corners.

### 2.2 Acceleration near the boundary

When $x$ approaches $\partial\Omega$, the polar distance function $d(\theta)$ develops sharp peaks. Although convergence remains superalgebraic, efficiency is improved by parameterizing the angular integral along the boundary and applying a graded parametrization as seen in [2]. This clusters quadrature nodes near the closest boundary point and significantly reduces the required resolution. Higher order convergence is depicted in Figure 1. We can prove that the derivatives of the graded parametrization are uniformly bounded in terms of the distance parameter $\varepsilon$ near the singularity.

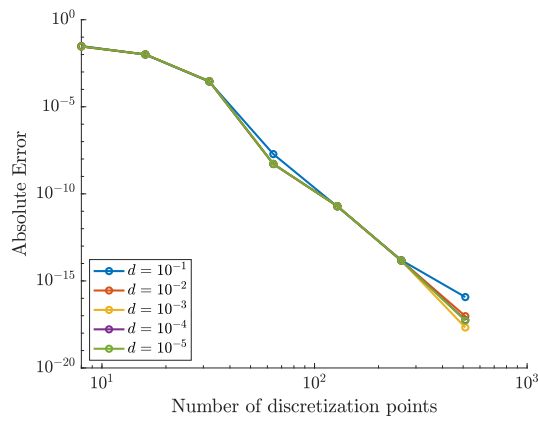

(a) Near a smooth boundary segment

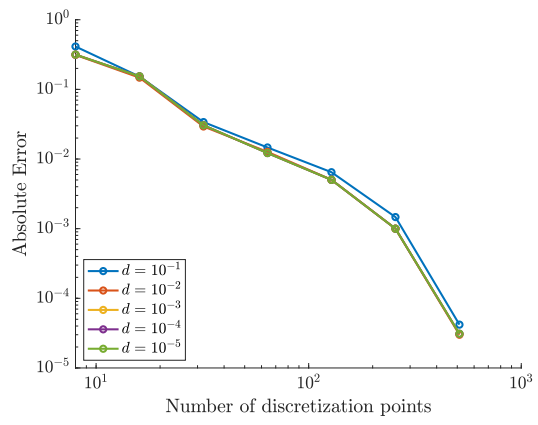

(b) Near a non-smooth boundary segment

Figure 1: Convergence of the accelerated TFF method for points at varying proximity to boundary features.

### 2.3 Global smooth integral $I_2$

After windowing, $I_2$ has a smooth integrand. Embedding $\Omega$ in a rectangle $R$ and approximating $\chi_\Omega$ by a truncated Fourier series $\chi_\Omega^F$, the integral reduces to a discrete convolution

$$I_2(x_{mn}) \approx \sum_{k,\ell} f^F(y_{k\ell})\, g(x_{mn} - y_{k\ell}),$$

which is evaluated efficiently using FFTs with zero padding. Existing TFF theory guarantees superalgebraic convergence in $F$.

## 3 Representative numerical results

Figure 2 shows typical convergence behavior for the three geometric regimes, confirming superalgebraic decay of the error. Figure 3 illustrates computed solutions on non-smooth domains, with errors around $10^{-10}$ and runtimes of a few seconds on a single CPU core using MATLAB for thousands of evaluation points.

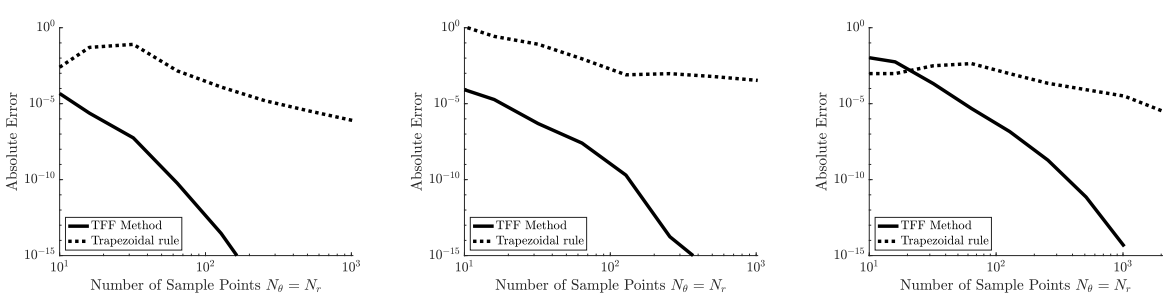

(a) Interior point (b) Smooth boundary (c) Corner case

Figure 2: Convergence of the TFF method for the three window–boundary regimes.

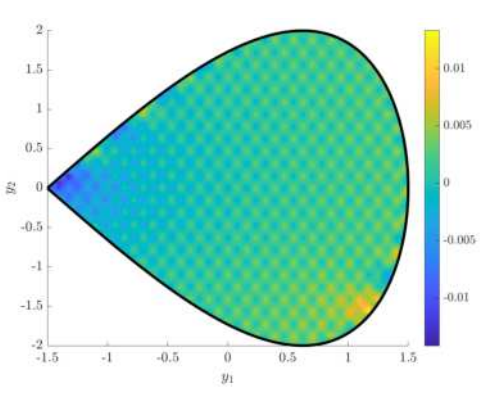

(a) Drop-shaped domain

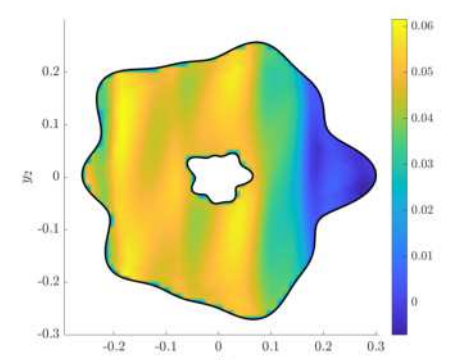

(b) Double-ring domain

Figure 3: Computed convolution values on representative non-smooth domains for $\phi(y) = \exp(\mathrm{i}40y_1 - \mathrm{i}20y_2)$.

## 4 Summary

We have developed a robust FFT-based method that combines smooth windowing, truncated Fourier series of singularity, and designed quadrature rules in radial and angular variables. The resulting scheme achieves superalgebraic convergence uniformly across the domain, including for geometries with corners and for evaluation points arbitrarily close to the boundary.

Beyond its immediate 2D implications, this work establishes the analytical and numerical foundation for an ongoing project aimed at the solution of three-dimensional scattering problems. In particular, the present framework provides the structural ingredients required to treat hypersingular and weakly singular kernels in 3D, while preserving high-order accuracy and computational efficiency.

# Quantitative ultrasound imaging in layered media

**Josselin Garnier**[1], **Laure Giovangigli**[2], **Pierre Millien**[3], **Sofia Suarez Casabiell**[1], [*]
[1]CMAP, CNRS, Ecole polytechnique, Institut Polytechnique de Paris, 91120 Palaiseau, France
[2]POEMS, CNRS, Inria, ENSTA, Institut Polytechnique de Paris, 91120 Palaiseau, France
[3]Institut Langevin, ESPCI Paris, PSL University, CNRS, 1 rue Jussieu, 75005 Paris, France
[*]Email: sofia.suarez-casabiell@polytechnique.edu

**Abstract**

In the context of medical imaging, we consider an acoustic two-layer medium with a planar interface. Each layer is made of small heterogeneities. A transducer array transmits pulses that are backscattered by the medium and recorded by the array. The objective is to estimate from the recorded signals the effective speeds of sound of the layers and the depth of the interface, in the presence or absence of a strong point-like reflector in the bottom layer. We propose a quantitative method based on the characterization of the wave fluctuations in the stochastic homogenization regime.

## 1 Introduction

Ultrasound images are formed by emitting pulses from a linear array and applying delay-and-sum beamforming assuming a sound speed. If this assumed speed is wrong, delays misalign, causing destructive interference and a blurred, distorted image, so estimating the effective speed of sound is essential.

## 2 Model

We model a tissue-mimicking medium by an incompressible homogeneous one in which there are numerous contrasted heterogeneities. Let $D = D_1 \cup D_2 \subset \mathbb{R}^2$ be an open bounded domain with at least Lipschitz boundary regularity, separated by a planar interface $\{y = d^\star\}$, $d^\star > 0$. In each layer $D_i$ there is a random collection of small heterogeneities $(S_j^{(i)})_j$. We suppose that the outer medium, the backgrounds and each heterogeneity are homogeneous so that the square index of refraction $n$ is given by (see Figure 1):

$$n = n_m \, \mathbb{1}_{\mathbb{R}^2 \setminus \overline{D}} + \sum_{i=1}^{2} \left( n_i \, \mathbb{1}_{D_i} + \sum_{j \in [1, \mathbb{N}^{(i)}]} (n_j^{(i)} - n_i) \, \mathbb{1}_{S_j^{(i)}} \right),$$

with $n_m, n_i, n_j^{(i)} \in [n_{\min}, n_{\max}]$.

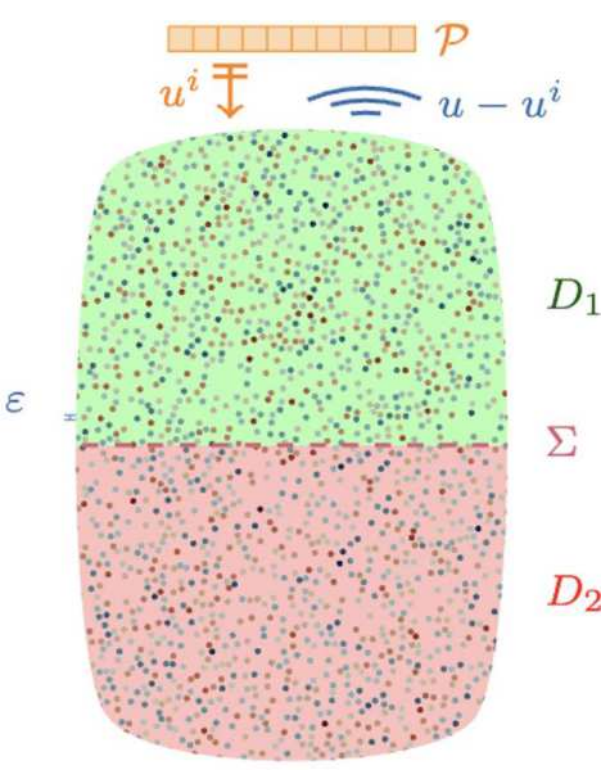

Figure 1: Model setup

The probe $\mathcal{P} = [-L, L] \times \{0\}$, $L > 0$, is discretised into $N_p$ elements at positions $\{\mathbf{x_e}\}_{e=1}^{N_p}$ (emitters) or $\{\mathbf{x_r}\}_{r=1}^{N_p}$ (receivers). We work in the frequency domain at angular frequency $\omega$ within the probe bandwidth $\mathcal{BW} := [\omega_c - \mathcal{B}, \omega_c + \mathcal{B}]$ where $\omega_c$ is the central frequency and $\mathcal{B} < \omega_c$. The incident plane wave with amplitude $U_i \in \mathbb{C}$ and direction $\theta \in \mathbb{S}^1$, is $u^i(\mathbf{x}) = U_i e^{ik_m \theta \cdot \mathbf{x}}$, $k_m = \omega \sqrt{n_m}$. The scattering problem by $D$ is modeled by the Helmholtz equation in $\mathbb{R}^2$, with the Sommerfeld radiation condition.

## 3 Homogenized problem and effective parameters

With stationarity and ergodicity assumptions on the random distribution of scatterers, it is possible to apply stochastic homogenization theory when the diameter of the scatterers becomes much smaller than the wavelength, [1]:
The unique solution $u \in H^1(B_R)$ of the scattering problem converges a.s. weakly in $H^1(B_R)$

to the unique solution $u^\star$ of a deterministic homogenized problem where each layer is modeled by a piecewise-constant effective index $n_i^\star$, so the microstructure can be replaced by effective layer parameters. We denote by $G^{\mathbf{p}}$ the Green's function of the two-layer medium with homogeneous parameters $\mathbf{p} = (c_1, c_2, d)$. $\mathbf{p}^\star = (c_1^\star, c_2^\star, d^\star)$ stands for the true homogenized parameters, to be estimated from the backscattered signals.

## 4 Imaging Functional

As in [2] for $\omega \in \mathcal{BW}$, define the reflection matrix $M(\mathbf{x_e}, \mathbf{x_r}, \omega)$, whose elements are the scattered signals recorded by $\mathbf{x_r}$ when $\mathbf{x_e}$ is used as a source. We can reconstruct the image by the following functional.

**Definition 1** (Confocal imaging functional)**.**

$$I^{\mathbf{p}}(\mathbf{z}) := \int_{\mathcal{P}\times\mathcal{P}\times\mathcal{BW}} \Big( \overline{M(\mathbf{x_e}, \mathbf{x_r}, \omega)} \, G^{\mathbf{p}}(\mathbf{z}, \mathbf{x_e}) \, G^{\mathbf{p}}(\mathbf{z}, \mathbf{x_r}) \Big) \, d\mathbf{x_e} \, d\mathbf{x_r} \, d\omega$$

*Using the 1st-order corrector in [3],*

$$I^{\mathbf{p}}(\mathbf{z}) = \int_D \big(n(\mathbf{y}) - n^\star(\mathbf{y})\big) \, F^{\mathbf{p}}(\mathbf{z}, \mathbf{y}) \, d\mathbf{y},$$

*where the PSF (point spread function) is*

$$F^{\mathbf{p}}(\mathbf{z}, \mathbf{y}) = \int_{\mathcal{BW}} \omega^2 \left( \int_{\mathcal{P}} G^{\mathbf{p}}(\mathbf{z}, \mathbf{x}) \overline{G^{\mathbf{p}^\star}(\mathbf{y}, \mathbf{x})} d\mathbf{x} \right)^2 d\omega.$$

## 5 Focal spot for $\mathbf{p} = (c_1, c_2, d)$

In the paraxial regime it is possible to find the center $\varphi_{\mathbf{p}}(\mathbf{y})$ and sizes of the focal spot on the image when $\mathbf{p} = (c_1, c_2, d)$ and there is a strong point-like reflector at $\mathbf{y} = (\eta^{\frac{1}{2}} y^\perp, y^{\shortparallel})$.

**Theorem 1** (Center of the focal spot)**.** *For* $\mathbf{y}$ *in the paraxial approximation,*

$$\varphi_{\mathbf{p}}(\mathbf{y}) = \left( \eta^{\frac{1}{2}} \frac{dc_1 + c_2(\varphi_{\mathbf{p}}(y)^{\shortparallel} - d)}{d^\star c_1^\star + c_2^\star (y^{\shortparallel} - d^\star)} y^\perp, \; d + \frac{c_2}{c_2^\star}\big(y^{\shortparallel} - d^\star\big) + \; c_2 \left( \frac{d^\star}{c_1^\star} - \frac{d}{c_1} \right) \right).$$

*Of course,* $\varphi_{\mathbf{p}^\star}(\mathbf{y}) = \mathbf{y}$.

**Theorem 2** (Broadband paraxial PSF widths)**.** *In the paraxial regime,*

-Axial width*:* $\Delta^{\shortparallel} \sim \sqrt{\frac{3}{2} \frac{c_2}{\mathcal{B}}}$.

-Lateral width*:*

$$\Delta^\perp \sim 2\sqrt{\frac{3}{2} \frac{c_2\big(\varphi_{\mathbf{p}}(\mathbf{y})^{\shortparallel} - d\big) + dc_1}{\omega_c L}}.$$

*(see Figure 2).*

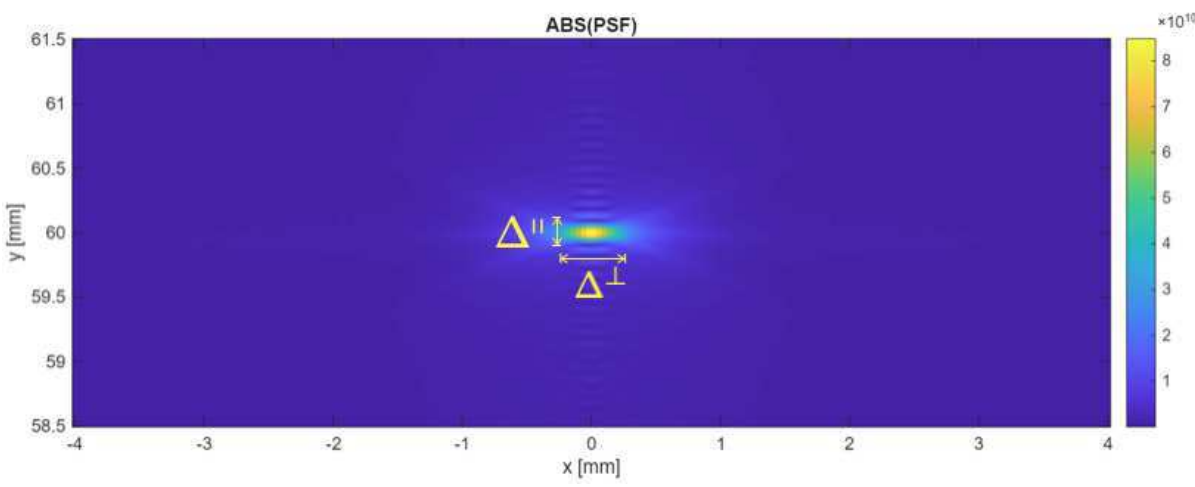

Figure 2: Shape of the PSF when $\mathbf{p} = \mathbf{p}^\star = (1450, 1580, 0.03)$ and $\mathbf{y} = (0, 0.06)$.

## 6 Estimating $\mathbf{p}^\star$ in the speckle

When there is no strong reflector, by exploiting the stationarity and ergodicity of the microstructure, we extend the focusing criterium used for a single reflector (maximizing $F^{\mathbf{p}}$ with respect to $\mathbf{p}$) to the speckle-only regime by averaging locally the imaging functional, as in [1]. Then, we define $N_r$ virtual guide stars in the speckle, that are $N_r$ local reference points $y_i$ (not true isolated reflectors) obtained through local spatial averaging of the imaging functional. They play the role of surrogate focal targets, and the estimator becomes

$$\widehat{\mathbf{p}} = \arg\max_{\mathbf{p}\in\mathcal{T}} \sum_{i=1}^{N_r} F^{\mathbf{p}}\big(\varphi_{\mathbf{p}}(\mathbf{y}_i), \mathbf{y}_i\big),$$

where $\mathcal{T}$ encodes the physical constraints on $\mathbf{p} = (c_1, c_2, d)$.

**Perspectives** A natural perspective is to study how the method behaves when the interface is not planar, since this may lead to biased estimates and poorer focusing.

# Analysis of the geometric spectral properties of electromagnetic waveguides

**Philippe Briet**[1], **Maxence Cassier**[2,*], **Thomas Ourmières-Bonafos**[3], **Michele Zaccaron** [2,4]

[1]Université de Toulon, Aix Marseille Univ, CNRS, CPT, Toulon, France
[2]Aix Marseille Univ, CNRS, Centrale Med, Institut Fresnel, Marseille, France
[3]Aix-Marseille Univ, CNRS, I2M, Marseille, France
[4]Inria, UMA, ENSTA Paris, Institut Polytechnique de Paris, 91120 Palaiseau, France

*Email: maxence.cassier@fresnel.fr

## Abstract

We consider a homogeneous isotropic electromagnetic waveguide with simply connected cross-section, embedded in a perfect conductor. We study how bending and twisting affect the real spectrum of the associated self-adjoint Maxwell operator. We first give geometric conditions ensuring preservation of the essential spectrum, which is symmetric, has a gap around the origin, and coincides with the spectrum of a reference waveguide with constant twist and zero curvature. We then provide in the case of purely bended waveguides sufficient conditions for the appearance of discrete eigenvalues in the gap (associated with trapped modes) and analyze them both theoretically and numerically.

## 1 Introduction

Constructing lossless waveguides that support trapped waves is a longstanding problem in physics, as waves typically propagate to infinity along the guiding direction. Mathematically, this amounts to identifying geometries for which the associated self-adjoint operator admits non-zero eigenvalues whose eigenfunctions (bound states or trapped modes) are spatially localized within the waveguide. A central question is how geometric perturbations of the waveguide can induce such spectral properties. In this context, we restrict our attention to the existence of a discrete spectrum consisting of isolated eigenvalues of finite multiplicity.

For the Dirichlet Laplacian in quantum waveguides, bending creates bound states [2], while some twisting prevents them. In contrast, the vectorial electromagnetic case remains largely open. As noted in [2], "the existence and nature of bound states depend more delicately on shape and curvature than in the scalar case". In [1], we partially address this issue.

## 2 Mathematical model

The *cross-section* $\omega$ of the considered waveguides is a bounded, simply connected Lipschitz domain in $\mathbb{R}^2$, with Cartesian coordinates $\mathbf{y} = (y_2, y_3) \in \omega$. In particular, a *straight waveguide* is defined as $\Omega_0 := \mathbb{R} \times \omega$.

The base curve $\Gamma \subset \mathbb{R}^3$ of the waveguide is assumed to admit a $\mathscr{C}^2$-injective arc-length parametrization $\gamma : \mathbb{R} \to \mathbb{R}^3$ whose curvature $\kappa : s \in \mathbb{R} \mapsto |\gamma''(s)|$ is bounded. One needs now to define an appropriate moving frame along $\Gamma$ as e.g. the Frenet frame. However, the Frenet frame is not defined when the curvature $\kappa$ vanishes, which makes it unsuitable for straight or locally straight waveguides. To overcome these limitations, we equip $\Gamma$ with a $\mathscr{C}^1$-smooth relatively parallel adapted frame which generalizes the Frenet frame. Such frame is uniquely defined by setting, for $s \in \mathbb{R}$, $\mathbf{e}_1(s) := \gamma'(s)$ and by choosing two vectors $\mathbf{e}_2(s)$, $\mathbf{e}_3(s)$ which are relatively parallel, i.e. satisfying for all $s \in \mathbb{R}$:

$$\begin{cases} \dfrac{\mathrm{d}\mathbf{e}_1}{\mathrm{d}s}(s) &= k_1(s)\mathbf{e}_2(s) + k_2(s)\mathbf{e}_3(s), \\ \dfrac{\mathrm{d}\mathbf{e}_2}{\mathrm{d}s}(s) &= -k_1(s)\mathbf{e}_1(s), \\ \dfrac{\mathrm{d}\mathbf{e}_3}{\mathrm{d}s}(s) &= -k_2(s)\mathbf{e}_1(s). \end{cases}$$

Here, $k_1, k_2$ are $\mathscr{C}^0$-functions (defined via $\gamma$) satisfying $k_1(s)^2 + k_2(s)^2 = \kappa(s)^2$ on $\mathbb{R}$ and the initial conditions $(\gamma'(0), \mathbf{e}_2(0), \mathbf{e}_3(0))$ form a positively oriented orthonormal basis of $\mathbb{R}^3$. (See the appendix section A of [1] for more details.)

Let $\theta : \mathbb{R} \to \mathbb{R}$ be a $\mathscr{C}^1$-smooth map and introduce the vectors $\mathbf{e}_1^\theta(s) := \mathbf{e}_1(s)$ and

$$\begin{pmatrix} \mathbf{e}_2^\theta(s) \\ \mathbf{e}_3^\theta(s) \end{pmatrix} := \begin{pmatrix} \cos(\theta(s))\,\mathbf{e}_2(s) - \sin\big(\theta(s)\big)\,\mathbf{e}_3(s) \\ \sin\big(\theta(s)\big)\,\mathbf{e}_2(s) + \cos\big(\theta(s)\big)\,\mathbf{e}_3(s) \end{pmatrix}.$$

Hence, $(\mathbf{e}_1^\theta, \mathbf{e}_2^\theta, \mathbf{e}_3^\theta)$ is nothing but a rotation in the transverse plane of angle $\theta$ of the relatively parallel adapted frame $(\mathbf{e}_1, \mathbf{e}_2, \mathbf{e}_3)$. Thus, its derivative $\theta'$ encodes the twist of the waveguide.

The perturbed waveguide $\Omega$ is defined as the image of the straight waveguide $\Omega_0 = \mathbb{R} \times \omega$ via

$$\Phi(s, y_2, y_3) = \gamma(s) + y_2 \mathbf{e}_2^\theta(s) + y_3 \mathbf{e}_3^\theta(s),$$

where $\Phi$ is a $\mathscr{C}^1$-diffeomorphism. This holds if $\Phi$ is globally injective and $b\|\kappa\|_{L^\infty} < 1$, with $b = \sup_{\mathbf{y}\in\omega} |\mathbf{y}|$ (which ensures that the waveguide does not locally self-intersect). For $\theta = \beta s$ with $\beta \in \mathbb{R}$ and $\kappa = 0$, $\Omega_\beta := \Phi(\Omega_0)$ is the reference waveguide with constant twist $\theta' = \beta$.

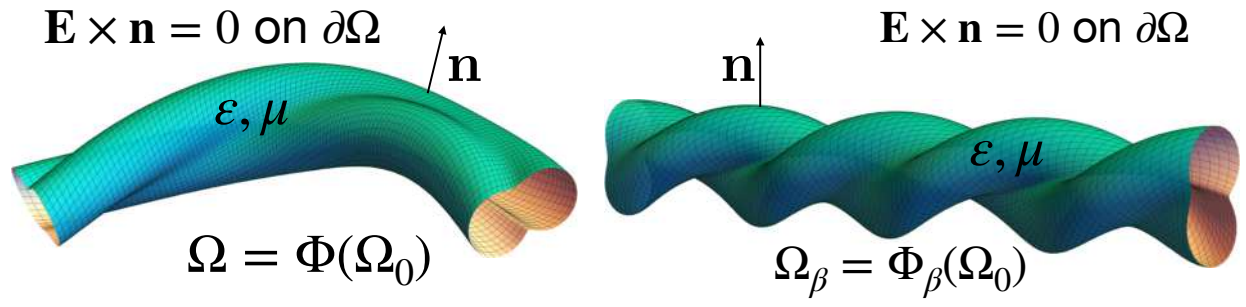

Figure 1: Left: general perturbed waveguide $\Omega$; Right: reference waveguide $\Omega_\beta$.

The waveguide $\Omega$ is enclosed in a perfect conductor and filled with an isotropic, homogeneous dielectric of permittivity $\varepsilon > 0$ and permeability $\mu > 0$. In $\Omega$, the evolution of the electric field $\mathbf{E}$ and magnetic field $\mathbf{H}$ is governed by the Maxwell's equations rewritten here as a conservative first order evolution equation:

$$\frac{\mathrm{d}\mathbf{U}}{\mathrm{d}t} + \mathrm{i}\mathcal{A}\mathbf{U} = \mathbf{F}, \ \forall t > 0 \ \text{ with } \ \mathbf{U}(0) = \mathbf{U}_0, \quad (1)$$

where the operator $\mathcal{A}$ is unbounded and self-adjoint on the Hilbert space $\mathbf{L}^2_{\varepsilon,\mu}(\Omega) = \mathbf{L}^2(\Omega)^2$ equipped with the following inner product:

$$(\mathbf{U}, \mathbf{U}')_{\mathbf{L}^2_{\varepsilon,\mu}(\Omega)} = \int_\Omega \left(\varepsilon \mathbf{E}_1 \cdot \overline{\mathbf{E}}_2 + \mu \mathbf{H}_1 \cdot \overline{\mathbf{H}}_2\right) \mathrm{d}\mathbf{x},$$
$$\forall \mathbf{U} = (\mathbf{E}, \mathbf{H})^\top, \quad \mathbf{U}' = (\mathbf{E}', \mathbf{H}')^\top \in \mathbf{L}^2_{\varepsilon,\mu}(\Omega).$$

In (1), $\mathcal{A} : D(\mathcal{A}) \subset \mathbf{L}^2_{\varepsilon,\mu}(\Omega) \to \mathbf{L}^2_{\varepsilon,\mu}(\Omega)$ is given by

$$\mathcal{A}\mathbf{U} = \mathcal{A}(\mathbf{E}, \mathbf{H}) := \mathrm{i}\big(\varepsilon^{-1}\, \mathbf{curl}\ \mathbf{H}, -\mu^{-1}\, \mathbf{curl}\ \mathbf{E}\big)$$

with domain $D(\mathcal{A}) := H_0(\mathbf{curl}\,, \Omega) \times H(\mathbf{curl}\,, \Omega)$ (we refer to [1] for the definition of these functional spaces), $\mathbf{U}_0 = (\mathbf{E}_0, \mathbf{H}_0)$ are the initial fields and $\mathbf{F} = (-\mathbf{J}/\varepsilon, 0)^\top$ is the source term defined via the current density $\mathbf{J}$. When $\Omega = \Omega_\beta$, the operator $\mathcal{A}$ is denoted by $\mathcal{A}_\beta$.

## 3 Main results and method

Our first result states that, under assumptions (2) on the curvature and twist at infinity, the essential spectrum $\sigma_{\mathrm{ess}}(\mathcal{A})$ of $\mathcal{A}$ in $\Omega$ coincides with the spectrum (and also the essential spectrum) of $\mathcal{A}_\beta$ in $\Omega_\beta$ (which is deducible from existing results). The proof consists, after applying two appropriate unitary transformations to straighten the waveguide and flatten the metric, in establishing a resolvent idendity between the perturbed and reference operators, and then to apply a meromorphic Fredholm theorem.

**Theorem 1** *If there exists $\beta \in \mathbb{R}$ such that*

$$\lim_{|s|\to\infty} \kappa(s) = 0 \ \textit{and} \ \lim_{|s|\to\infty} (\theta'(s) - \beta) = 0, \qquad (2)$$

*then there exists $a_\beta > 0$ such that*

$$\sigma_{\mathrm{ess}}(\mathcal{A}) = \sigma(\mathcal{A}_\beta) = (-\infty, -a_\beta] \cup \{0\} \cup [a_\beta, +\infty),$$

*(with $\sigma(\mathcal{A}_\beta) = \sigma_{\mathrm{ess}}(\mathcal{A}_\beta)$). Moreover, if $\beta = 0$, then $a_0 = \sqrt{\lambda_2^N(\omega)}\, c$, with $c = (\varepsilon\mu)^{-1/2}$ and $\lambda_2^N(\omega)$ the first non-zero Neumann eigenvalue in $\omega$.*

Let $\psi$ be a normalized mode associated to the first non-zero Neumann eigenvalue $\lambda_2^N(\omega)$. We define $\mathbf{X}$ and $\mathbf{Y}^\theta \in \mathbb{R}^2$ (for $\theta \in \mathbb{R}$) by

$$\mathbf{X} = \int_{\partial\omega} \mathbf{n}(\mathbf{y})|\psi(\mathbf{y})|^2 \mathrm{d}\sigma(\mathbf{y}) \text{ and } \mathbf{Y}^\theta = \mathfrak{R}_\theta \int_{\mathbb{R}} \begin{pmatrix} k_1(s) \\ k_2(s) \end{pmatrix} \mathrm{d}s$$

where $\mathbf{n}$ is the unit outward normal vector to $\partial\omega$ and $\mathfrak{R}_\theta$ the counterclockwise rotation of angle $\theta$.

**Theorem 2** *Assume that $\theta' = 0$, $\kappa \in L^1(\mathbb{R})$ and $\lim_{|s|\to\infty} \kappa(s) = 0$. If*

$$\mathbf{X} \cdot \mathbf{Y}^\theta > 2\lambda_2^N(\omega)\, \frac{b^2\|\kappa\|_{L^\infty(\mathbb{R})}}{1 - b\|\kappa\|_{L^\infty(\mathbb{R})}} \|\kappa\|_{L^1(\mathbb{R})} \qquad (3)$$

*then $\sigma_{\mathrm{dis}}(\mathcal{A})$ is not empty (and symmetric).*

Theorem 2 gives sufficient conditions for the existence of discrete eigenvalues in the gap of $\mathcal{A}$ around 0. Its proof uses min-max theory on the quadratic form of $\mathcal{A}^2$ restricted to divergence-free functions in $\ker(\mathcal{A})^\perp$. The geometric conditions $\mathbf{X} \neq 0$ and $\mathbf{Y}^\theta \neq 0$, related to the cross-section and base curve, are required for the inequality (3). Since $\theta$ is constant in $\mathbf{Y}^\theta$, this angle has to be chosen in the frame $(\mathbf{e}_1^\theta, \mathbf{e}_2^\theta, \mathbf{e}_3^\theta)$ to maximize the left-hand side of (3), reflecting the vectorial nature of Maxwell's equations and the selection of an optimal bending orientation. In [1], we analyze the geometric conditions $\mathbf{X} \neq 0$ and $\mathbf{Y}^\theta \neq 0$ theoritically and numerically to show that many waveguides satisfying (3) (and thus $\sigma_{\mathrm{dis}}(\mathcal{A}) \neq \emptyset$) can be constructed.

**Funding**

This work received support from the french government under the France 2030 investment plan, as part of the Initiative d'Excellence d'Aix-Marseille Université-A*MIDEX AMX-21-RID-012.

# Convergence analysis for Numerical Steepest Descent contour deformation methods

**Thomas Caussade**[1,*], **David P. Hewett**[1]

[1]Department of Mathematics, University College London, United Kingdom

*Email: thomas.caussade.23@ucl.ac.uk

**Abstract**

The evaluation of highly oscillatory integrals is an important topic in many areas of computational wave propagation. The Numerical Steepest Descent (NSD) method is a powerful approach to computing such integrals, which combines complex contour deformation with quadrature rules. However, for fixed frequencies NSD may lose accuracy when stationary points coalesce with each other, or with endpoints of the integration contour. It was recently demonstrated in [1] that this can be dealt with by using standard quadrature inside neighbourhoods of stationary points, in which the number of oscillations is bounded and small, combined with NSD away from stationary points. To date, a rigorous analysis of such an approach has not been performed. In this talk we present a framework for such an analysis, and show for a number of model functions how this delivers provably high-order accuracy, even in presence of branch point singularities and coalescing stationary points.

## 1 Introduction

Consider the numerical evaluation of integrals of the form

$$\int_\Gamma f(z)e^{i\omega g(z)}\,\mathrm{d}z, \tag{1}$$

where $\Gamma$ is a complex contour in $\mathbb{C}$, $f$ is the amplitude, $g$ is the phase, and $\omega > 0$ is a frequency parameter. Such integrals arise in numerous applications, particularly in the study of high-frequency wave propagation and quantum mechanics.

These are typically difficult to evaluate when $\omega$ is large, because standard quadrature rules usually need a large number of quadrature points to achieve high accuracy. However, there are many approaches to effectively evaluate integrals of the form (1), which exploit the oscillatory behaviour of the integral to produce numerical schemes where accuracy typically improves as frequency increases.

In this talk, we focus on one such method, known as the Numerical Steepest Descent (NSD) method, which involves a combination of contour deformation techniques and quadrature rules.

## 2 Numerical Steepest Descent method

In its simplest form [3], in the NSD method one deforms $\Gamma$ onto a union of contours on which $\Re g(z)$ is (piecewise) constant so that the integrand is no longer oscillatory. These contours coincide with the steepest descent (SD) curves of $-\Im g$, and they connect endpoints of $\Gamma$, valleys at infinity (which correspond to angular sectors where the integrand decays rapidly as $|z| \to \infty$), and stationary points of $g$, which are points $\xi \in \mathbb{C}$ where $g'(\xi)$ vanishes. Then, one applies a suitable numerical quadrature to evaluate the contribution of each SD contour to the value of the integral.

However, implementation of the NSD usually requires specific knowledge of the phase because evaluating integrals along SD contours that emanate from or near stationary points is fraught with numerical difficulties, especially when stationary points are close to another stationary point or an endpoint.

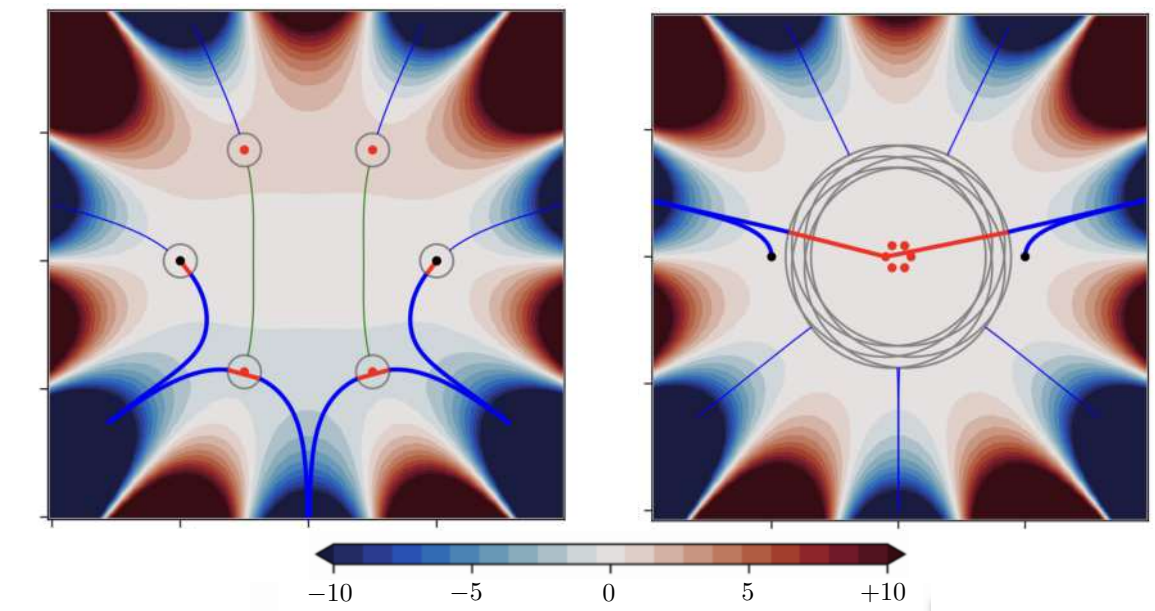

Figure 1: Quasi-SD deformation of $[-1, 1]$ by following SD contours (in blue) and straight lines (in red) inside the non-oscillatory balls centred at stationary points (red dots), for $g(z) = z^7/7 - \varepsilon^6 z$ with $\varepsilon = 1$ (left) and $\varepsilon = 0.1$ (right). The colorbar shows $-\Im g(z)$.

Algorithmically, this can be handled by in-troducing non-oscillatory regions, as in the PathFinder algorithm described in [1,2]. The algorithm exploits the fact that the integrand is locally slowly varying at stationary points, and introduces bounded "non-oscillatory balls" containing each stationary point, within which the integrand provably undergoes at most fixed (and small) number of oscillations. Then, the deformed contour, which we call the "quasi-SD deformation", is the union of straight-line contours inside non-oscillatory balls, where standard quadrature rules are efficient, and SD contours outside, where we proceed with NSD techniques. We call such a method Regularised NSD (RNSD) method.

## 3 Quadrature error analysis

Existing error analysis for NSD tends to focus on proving asymptotic convergence as $\omega \to \infty$, and it has already been demonstrated that if $\omega$ is sufficiently large coalescence phenomena dissapears, and extreme high accuracy can be attained [3].

For fixed $\omega$, we expect numerical convergence as the number of quadrature points increases, albeit with a rate that depends on $\omega$ and other parameters of the phase. However, there is little rigorous theory for numerical convergence of NSD in the literature, because of the difficulties with stationary point coalescence mentioned above.

In this talk, we show that RNSD methods are provably high-order methods with computational cost essentially independent of frequency. We derive rigorous quadrature error bounds for fixed $\omega$ and increasing number of quadrature points, showing that convergence can be achieved uniformly in $\omega$, even when coalescence and branch-point singularities are involved.

Our analysis relies on the estimation of the size of neighbourhoods of analyticity around each contour of the quasi-SD deformation, and the total quadrature error is the combination of error contributions from each traversed contour.

On finite contours, we apply Gauss-Legendre quadrature rules, and the error analysis follows from classical estimates with Bernstein ellipses [4], giving exponential convergence as the number of quadrature points increases. In the presence of near-singularities, e.g. branch points, we resort to a composite quadrature on a graded mesh instead, which also provides provably high-order accuracy.

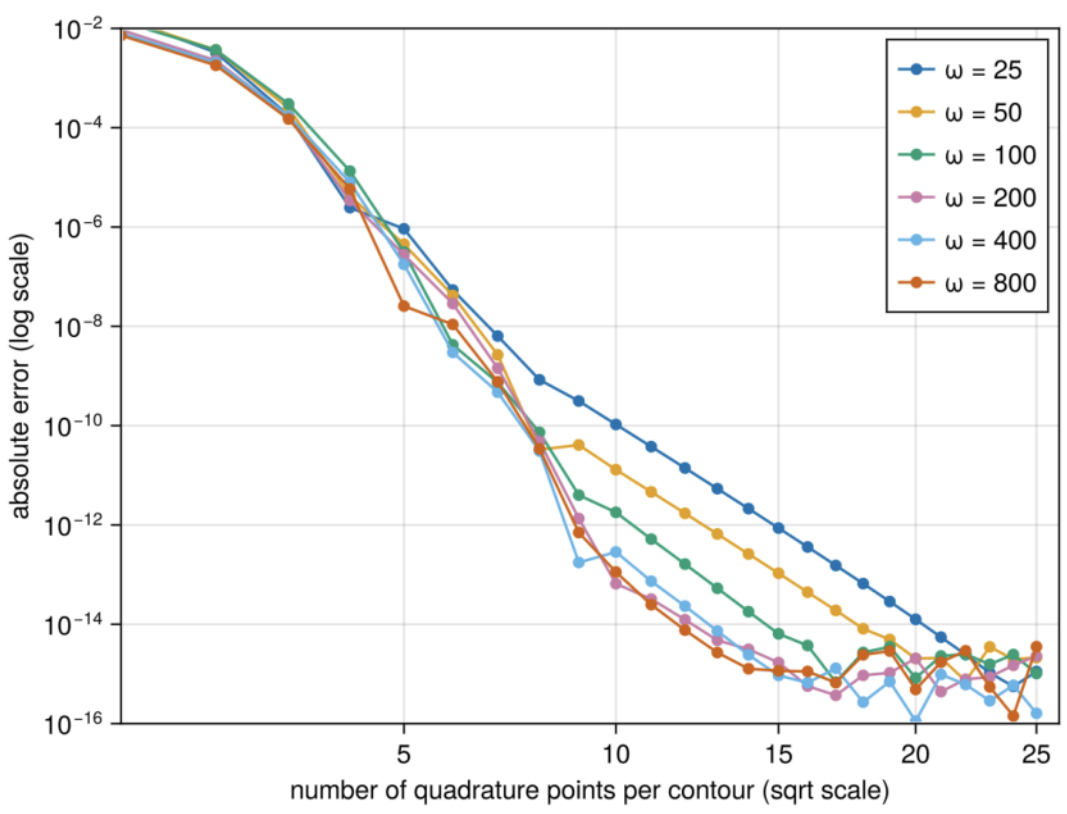

Figure 2: Absolute self-convergence error (log scale) as the number of quadrature points (sqrt scale) increases for different fixed values of frequency, with $g(z) = z^7/7 - z$.

On infinite SD contours, we apply Gauss-Laguerre quadrature rules. The error analysis relies on more recent estimates [5], which show that error decays with root-exponentially as the number of quadrature points increase. The error estimate depends on the size of a parabolic region around the SD contour and the growth of the integrand at infinity inside this region.

We present RNSD quadrature error estimates for a number of example phase functions that are explicit in $\omega$ and parameters of the phase. These include a monomial phase, a cubic phase exhibiting coalescence, and a model singular phase with branch points.

# Efficient 3D elastodynamic multiple scattering using the Method of Reflections: convergence dependence on boundary integral formulation

Stéphanie Chaillat[1,*], Marion Darbas[2]
[1]Laboratoire POEMS (CNRS-INRIA-ENSTA), Institut Polytechnique de Paris, Palaiseau, France
[2]LAGA CNRS UMR 7539, Université Sorbonne Paris Nord, Villetaneuse, France
*Email: stephanie.chaillat@ensta.fr

**Abstract**

Boundary Element Methods have reached maturity for solving problems in homogeneous domains with both efficiency and accuracy. The next challenge lies in developing extensive domain decomposition strategies tailored for parallelization and iterative solutions of very large-scale problems. As an initial focus, we consider multiple scattering problems, where natural boundary domain decomposition methods can be derived using the Method of Reflections.

In previous works, we developed new variants of this approach for 2D Helmholtz and elastodynamic problems. In particular, numerical simulations show that the convergence of the standard Method of Reflections is affected in the presence of closed-space obstacles. We successfully introduced domain decomposition techniques such as obstacle overlap to improve performance. In this study, we investigate the influence of the integral formulation on each scatterer. Our goal is to develop efficient preconditioners for 3D multiple scattering problems.

## 1 Fast BEMs for waves

Wave propagation can be simulated with methods like finite elements, finite differences, spectral elements, and Boundary Element Methods (BEMs). In this approach, the exterior problem is reformulated as a boundary integral equation using the Green's function, determining its efficiency and limitations. It is particularly efficient in homogeneous, unbounded domains, where the free-space Green's function satisfies radiation conditions. A major drawback is that BEM leads to dense, non-symmetric linear systems, resulting in a high computational cost for large scale problems. These constraints limit its application to realistic geometries or heterogeneous media.

BEMs have become well-established for solving problems in homogeneous domains offering both efficiency and accuracy thanks to the development of fast BEMs (based on $\mathcal{H}$-matrices or Fast Multipole Method). These approaches have demonstrated their potential for modeling not only acoustic wave propagation but also seismic waves in the ground [1, 4, 6]. A major challenge is to design scalable domain decomposition strategies compatible with parallellization and iterative solutions of very large-scale problems.

## 2 The Method of Reflections

We focus on 3D elastic multiple scattering problems. Boundary domain decomposition methods can be derived easily using the Method of Reflections. We consider a collection of $I$ disjoint bounded open sets $\{O_i\}_{1\leq i\leq I}$ of $\mathbb{R}^3$ with smooth boundaries $\{\Gamma_i\}_{1\leq i\leq I}$ and note $\Omega = \mathbb{R}^3 \setminus \bigcup_{i=1}^I \bar{O}_i$ the exterior domain (see Fig. 1 for the notations). In $\Omega$, the propagation of time-harmonic waves in an isotropic and homogeneous elastic medium is governed by the Navier equation

$$\begin{aligned} &\mu\Delta\mathbf{u} + (\lambda+\mu)\nabla(\nabla\cdot\mathbf{u}) + \rho\omega^2\mathbf{u}, && \text{in } \Omega,\\ &\mathbf{u}|_{\Gamma_i} = \mathbf{g}_i = -\mathbf{u}_{\text{inc}}|_{\Gamma_i}, && i=1,\dots,I,\\ &\text{Radiation conditions,} \end{aligned} \tag{1}$$

where $\mathbf{u}$ is the displacement, $\lambda$, $\mu$ and $\rho$ the elastic parameters and $\omega$ the angular frequency. It can be formulated as a system of boundary integral equations and solved using BEMs. For instance, if we consider a formulation based solely on the single-layer potential,

$$\begin{aligned} \mathbf{u}(\mathbf{x}) = \mathbf{S}(\boldsymbol{t})(\mathbf{x}) &= \sum_{i=1}^{I} \mathbf{S}_i(\boldsymbol{t}_i)(\mathbf{x})\\ &= \sum_{i=1}^{I} \int_{\Gamma_i} \mathbf{\Phi}(\mathbf{x},\mathbf{y})\,\boldsymbol{t}_i(\mathbf{y})\,ds(\mathbf{y}), \quad \mathbf{x}\notin\partial\Omega_i, \end{aligned} \tag{2}$$

where $\mathbf{\Phi}(\mathbf{x},\mathbf{y})$ is the fundamental solution of the Navier equation, we obtain the system:

$$\begin{pmatrix} S_{11} & S_{12} & \cdots & S_{1I}\\ \vdots & \vdots & \ddots & \vdots\\ S_{I1} & S_{I2} & \cdots & S_{II} \end{pmatrix} \begin{pmatrix} q_1\\ \vdots\\ q_I \end{pmatrix} = \begin{pmatrix} g_1\\ \vdots\\ g_I \end{pmatrix}. \tag{3}$$

Each block $S_{ij}$ represents the discretized interaction between obstacles $i$ and $j$. The vector $q_i$ approximates the unknown traction $\mathbf{t}_i$ on obstacle $i$, while $g_i$ corresponds to the discretized boundary data associated with obstacle $i$. Formulating the linear system is intuitive, but its resolution is computationally intensive due to the non-symmetric and dense block matrix.

The Method of Reflections [5,7] is a domain decomposition approach for multiple scattering problems. It decomposes the scattered field generated by the multi-body configuration into contributions from individual obstacles, solving each scattering problem iteratively. For each obstacle, the incident field is defined as the superposition of the original wave and the fields scattered by the other obstacles at previous iterations. Since the approach is independent of the underlying numerical method, and due to the homogeneous nature of the exterior problem, it is well suited to BEMs.

## 3 Improvement of the convergence

Practical implementation of the Method of Reflections is straightforward. At the discrete level, the sequential variant can be directly obtained by applying the block Gauss-Seidel method to the full system, while the parallel version follows naturally from the block Jacobi method. The latter not only simplifies implementation but also provides an inherent way to parallelize BEMs for multiple scattering problems. These formulations are not equivalent, and their performance varies significantly. What is truly compelling is analyzing how these methods behave as a function of obstacle spacing, geometry, and frequency, as these factors critically influence their computational efficiency. In [2], we study the Method of Reflections for 2D Helmholtz problems in detail. Our results show that both the sequential and parallel forms of the Method of Reflections, when based on boundary integral equations, provide an efficient solver. However, we observe that small distances between scatterers can significantly slow down the method, and in some cases prevent the parallel version from converging.

To address this issue, we introduce in [2] domain decomposition techniques, such as overlap and coarse spaces, within the Method of Reflections for solving the block EFIE system (3), which substantially improve performance.

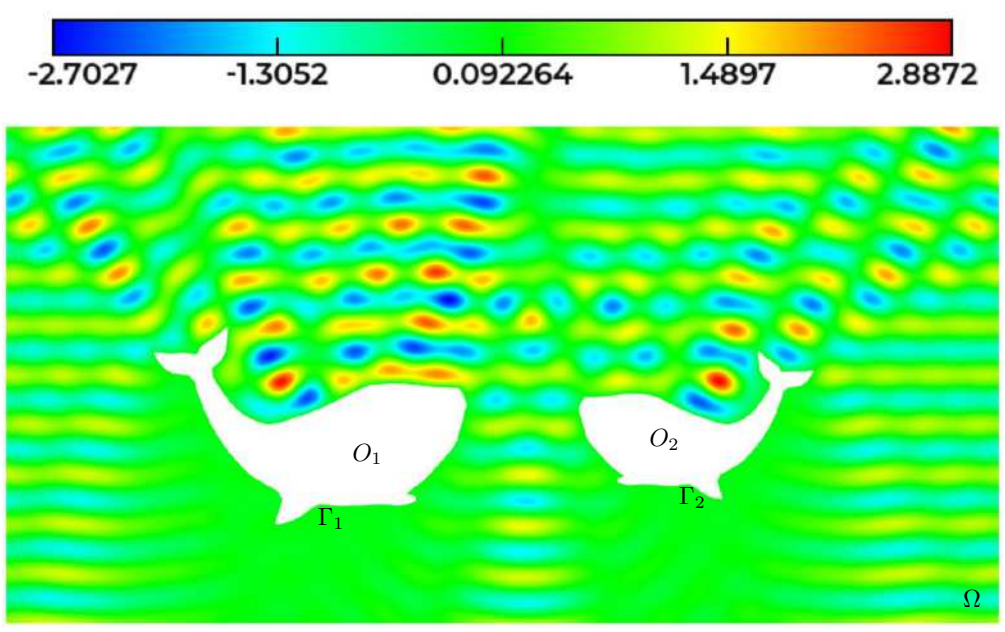

Figure 1: Illustration of a multiple scattering problem (incident plane wave, 2D Helmholtz).

For instance, nearby obstacles are grouped into overlapping subdomains and treated together within the Method of Reflections. The main strength of the new approach lies in its effectiveness as a preconditioner for Krylov methods. In this contribution, we extend the approach to 3D elastodynamics (with FM-BEM) and consider the question of the most robust boundary integral formulations. We will compare approaches based on single- or double-layer potentials and assess the efficiency of well-conditioned combined formulations using analytic preconditioners [3].

## A neural operator framework for solving inverse scattering problems.

**Victor Chenu**[1], **Houssem Haddar**[1], **Hadrien Montanelli**[1]
[1]Inria, ENSTA, UMA, Institut Polytechnique de Paris, Palaiseau, France

**Abstract**

We present a neural operator framework for solving inverse scattering problems [3]. A neural operator produces a preliminary indicator function for the scatterer, which, after appropriate rescaling, is used as a regularization parameter within the Linear Sampling Method to validate the initial reconstruction. The neural operator is implemented as a DeepONet with a fixed radial-basis-function trunk, while the noise level required for rescaling is estimated using a dedicated neural network. Two-dimensional numerical experiments demonstrate the effectiveness of our approach.

### 1 Introduction

Inverse acoustic scattering problems consist in reconstructing the geometry of an unknown obstacle from measurements of scattered waves. At fixed frequency, sampling methods such as the Linear Sampling Method (LSM) [1] provide robust qualitative reconstructions, but they require the solution of many ill-posed linear systems together with a suitable choice of regularization parameters. In particular, Morozov's discrepancy principle, while reliable, is computationally expensive and requires prior knowledge of the noise level.

In this work, we propose a hybrid scientific machine learning strategy combining neural operators with the LSM. Rather than replacing classical inversion techniques, our approach enhances them by injecting learned information into the regularization procedure. A neural operator produces a preliminary indicator function encoding the approximate location and size of the obstacle. This indicator is then rescaled using a predicted noise level and incorporated into the LSM as a spatially varying regularization function.

The neural operator is implemented as a radial-basis-function (RBF) DeepONet, based on the DeepONet architecture [2]. A second neural network estimates the noise level directly from the singular values of the far-field matrix. The resulting method preserves the robustness of the LSM while significantly reducing computational cost.

### 2 Model

We consider time-harmonic acoustic wave propagation in $\mathbb{R}^2$ in the presence of a sound-soft obstacle $D \subset \mathbb{R}^2$. For a given incident plane wave

$$u^{\mathrm{i}}(x, d) = e^{ikx\cdot d}, \qquad d \in \mathbb{S}^1,$$

the scattered field $u^{\mathrm{s}}$ satisfies the exterior Helmholtz problem

$$\Delta u^{\mathrm{s}} + k^2 u^{\mathrm{s}} = 0 \quad \text{in } \mathbb{R}^2 \setminus \overline{D}, \tag{1}$$

with Dirichlet boundary condition

$$u^{\mathrm{s}} = -u^{\mathrm{i}} \quad \text{on } \partial D, \tag{2}$$

and the Sommerfeld radiation condition at infinity.

The scattered field admits the far-field expansion

$$u^{\mathrm{s}}(x) = \frac{e^{ik|x|}}{\sqrt{|x|}} \left( u^{\infty}(\hat{x}, d) + \mathcal{O}(|x|^{-1}) \right), \quad |x| \to \infty,$$

where $u^{\infty}(\hat{x}, d)$ denotes the far-field pattern.

Given $n$ incident directions and $m$ observation directions, the measurements are collected into a far-field matrix

$$F_{ij} = u^{\infty}(\hat{x}_i, d_j).$$

The inverse problem consists in reconstructing the geometry of the obstacle $\partial D$ from $F$. This problem is nonlinear and severely ill-posed.

The Linear Sampling Method transforms the inverse problem into a family of linear systems

$$F g_z = \phi_z,$$

indexed by sampling points $z$ in a probing domain. The norm of the solution $g_z$ defines an indicator function that characterizes the obstacle. In practice, the system must be regularized, typically via Tikhonov regularization with a parameter usually chosen using Morozov's discrepancy principle.

## 3 Methods and Results

Our approach replaces Morozov's discrepancy principle with a learned regularization strategy. The workflow consists of three components.

First, a neural operator maps the noisy far-field matrix $F_\delta$ to a normalized indicator function $I_\theta(z)$. We implement this operator as a DeepONet, where the trunk net is fixed with Gaussian radial basis functions.

Second, the noise level $\delta = \|F - F_\delta\|_2$ is estimated using a separate network trained on the singular values of the far-field matrix $F_\delta$. This choice is motivated by the super-exponential decay of singular values of the far-field operator and the appearance of a noise-dependent plateau.

Third, the predicted indicator and noise level are combined into a spatially varying regularization function

$$\alpha_\theta(z) = \delta_\theta I_\theta(z),$$

which is used within the LSM. This produces a refined indicator that validates and improves the neural prediction.

Two-dimensional numerical experiments (see 1 )demonstrate that the neural operator accurately recovers the location and approximate geometry of single obstacles, and provides reasonable initial guesses for multiple obstacles outside the training regime. When incorporated into the LSM, the reconstruction quality is comparable to the standard Morozov-based approach while significantly reducing computational .

These results highlight the complementary roles of neural operators and classical inversion methods: the neural network provides a fast, data-driven indicator, while the LSM restores robustness and theoretical consistency.

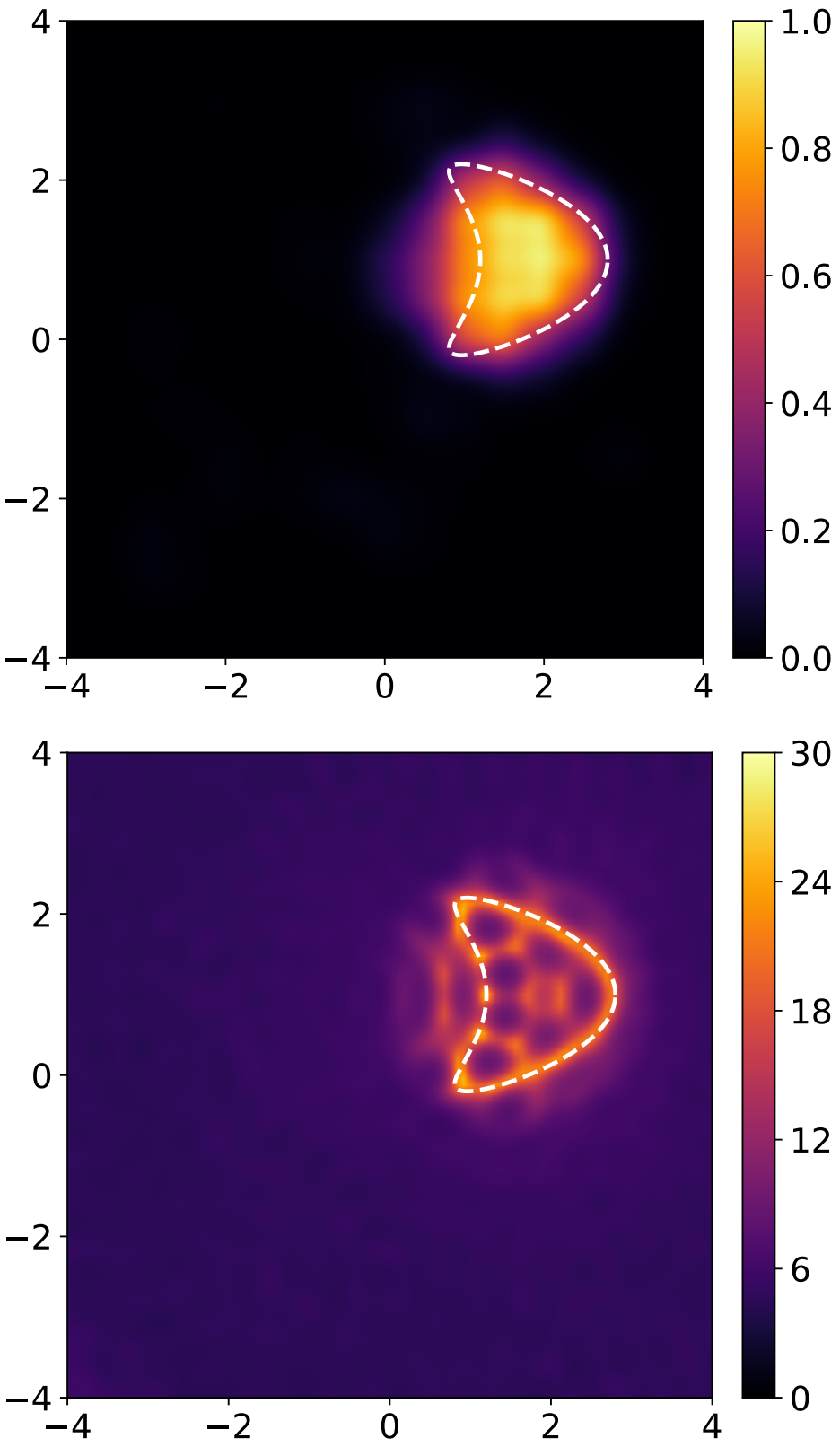

Figure 1: Reconstruction of a kite using neural network (top) and using the LSM (bottom).

# Reconstruction of the Shear Modulus in Nearly Incompressible Media via Adaptive Spectral Inversion and the Virtual Fields Method

**Nagham Chibli**[1,2,*], **Martin Genet**[1], **Sébastien Imperiale**[2]
[1]LMS, Ecole Polytechnique, CNRS, Institut Polytechnique de Paris
[2]Inria, Inria-Saclay - CMAP, Ecole Polytechnique, CNRS, Institut Polytechnique de Paris
*Email: nagham.chibli@inria.fr

**Abstract**

In time-harmonic elastography, the shear modulus is reconstructed from full-field displacement data by solving an inverse problem governed by the elastodynamic equation. When the parameter distribution is partially known, for instance as a heterogeneous region with a known location, the Virtual Fields Method (VFM) provides an efficient reconstruction framework. This method relies on the principle of virtual work, that is, the variational formulation of the direct problem assuming the displacement field is known. In nearly incompressible media, the formulation involves both the shear modulus and an unknown pressure term, which cannot be directly measured and may strongly affect stability. To address this, suitable virtual displacement fields, i.e., test functions, are chosen to generate independent equations for the shear modulus while canceling the pressure contribution. An optimal class of virtual fields is introduced to minimize the relative error on reconstructed parameters. Numerical experiments show improved robustness against noise. In practice, the functional basis representing the shear modulus is often unknown. We propose using the Adaptive Spectral Inversion (ASI) method, which iteratively builds a low-dimensional basis from the eigenfunctions of a suitable operator. The idea is to combine the ASI method and the VFM approach by integrating the adaptive basis construction of ASI with the pressure elimination and robustness of VFM. Numerical results confirm that this strategy effectively reconstructs the shear modulus in a stable manner.

## 1 Introduction

This work focuses on the reconstruction of a heterogeneous shear modulus $\mu$ in the incompressible linear elastodynamics equation in the frequency domain

$$\nabla \cdot \big(2\mu\, \boldsymbol{\varepsilon}(\mathbf{u})\big) + \nabla p + \rho\,\omega^2 \mathbf{u} = \mathbf{0}, \quad \text{in } \Omega, \quad (1)$$

where $\Omega \subset \mathbb{R}^d$ $(d = 1, 2, 3)$ with boundary $\partial\Omega = \Gamma_D \cup \Gamma_N$, $\Gamma_D \cap \Gamma_N = \varnothing$. The full displacement field $\mathbf{u}$ is assumed known, the density $\rho$ and the angular frequency $\omega$ are given, and $p$ denotes the hydrostatic pressure, which acts as a Lagrange multiplier and remains unknown. The boundary conditions are

$$\begin{cases} \mathbf{u} = \mathbf{u}_D & \text{on } \Gamma_D, \\ \big(2\mu\boldsymbol{\varepsilon}(\mathbf{u}) - p\mathbf{I}\big)\mathbf{n} = \mathbf{t} & \text{on } \Gamma_N. \end{cases}$$

One of the difficulties compared to the compressible case is the presence of the unknown pressure field $p$, which cannot be measured. Moreover, the traction $\mathbf{t}$ is assumed to be unknown.

## 2 Methods and Results

At first, we assume that the shear modulus distribution can be expressed in a given basis. In other words, we seek coefficients $\alpha_j$ such that

$$\mu = \sum_{j=1}^{n} \alpha_j\, \phi_j, \quad (2)$$

where $\phi_1, \dots, \phi_n$ are known basis functions. The Virtual Fields Method (VFM) is based on the principle of virtual work [1], which corresponds to the weak form of the equilibrium equation, that reads, for all $\mathbf{v} \in \mathbf{V}_0$,

$$\sum_{j=1}^{n} \alpha_j\, a(\phi_j; \mathbf{u}, \mathbf{v}) + b(p, \mathbf{v}) = \ell(\mathbf{v}). \quad (3)$$

Here $\mathbf{V}_0 = \boldsymbol{H}^1_{0,\Gamma_D}(\Omega)$ and

$$a(\phi_j; \mathbf{u}, \mathbf{v}) = \int_\Omega 2\phi_j\, \boldsymbol{\varepsilon}(\mathbf{u}) : \boldsymbol{\varepsilon}(\mathbf{v})\, d\Omega,$$

$$b(p, \mathbf{v}) = \int_\Omega p\, \nabla \cdot \mathbf{v}\, d\Omega,$$

$$\ell(\mathbf{v}) = \int_{\Gamma_N} \mathbf{t} \cdot \mathbf{v}\, d\Gamma + \int_\Omega \rho\omega^2 \mathbf{u} \cdot \mathbf{v}\, d\Omega.$$

We select exactly $n$ virtual fields (test functions), chosen to eliminate the dependence on $p$ and on $\mathbf{t}$ ($\mathbf{T} \in L^\infty(\Gamma_N)^{d\times d}$ is introduced to remove the traction contribution). Specifically, each virtual field $\mathbf{v}_i$, $1 \le i \le n$, belongs to

$$K_0 = \left\{ \mathbf{v} \in V_0 \;\middle|\; \begin{array}{ll} b(p,\mathbf{v}) = 0, & \forall p \in L^2(\Omega), \\ \mathbf{T}\mathbf{v} = \mathbf{0}, & \text{on } \Gamma_N \end{array} \right\}.$$

With this choice, by Eq. (3), we obtain

$$\boldsymbol{A}\,\boldsymbol{\alpha} = \boldsymbol{\lambda},$$

where $\boldsymbol{\alpha} = (\alpha_1, \dots, \alpha_n)^T$ and, for $1 \le i, j \le n$,

$$A_{ij} = a(\phi_j; \mathbf{u}, \mathbf{v}_i), \quad \lambda_i = \ell(\mathbf{v}_i).$$

We propose to choose the virtual fields so that the matrix $\boldsymbol{A}$ is diagonal, namely,

$$\mathbf{v}_i \in H_i := \left\{ \mathbf{v} \in K_0 \;\middle|\; a(\phi_j, \mathbf{u}, \mathbf{v}) = 0, \ \forall j \neq i \right\}.$$

Moreover, each $\mathbf{v}_i$ is chosen to maximize the corresponding diagonal term, more precisely,

$$\mathbf{v}_i = \underset{\mathbf{v}\in H_i,\ \|\mathbf{v}\|_V \le 1}{\operatorname{argmax}} \ a(\phi_i, \mathbf{u}, \mathbf{v}).$$

With this choice, we establish a stability estimate with respect to measurement noise. Let $\mathbf{u}^\delta$ denote the noisy measurement and $\mathbf{u}^0$ the noiseless one, satisfying (1) with the true parameter $\mu^0$ of the form (2). If $\mu^\delta$ denotes the parameter reconstructed using the proposed virtual fields, we can show that, for some $p \in [1, +\infty]$, if there exists a constant $C_p$ such that

$$\sum_{i=1}^{n} \frac{\|\phi_i\|_{L^p(\Omega)}}{|a(\phi_i, \mathbf{u}^\delta, \mathbf{v}_i^\delta)|} \le C_p,$$

then the following stability estimate holds

$$\|\mu^\delta - \mu^0\|_{L^p(\Omega)} \le C(C_p, \rho, \omega, \mu^0)\, \|\mathbf{u}^\delta - \mathbf{u}^0\|_{H^1(\Omega)}.$$

Numerical experiments show that these virtual fields significantly improve the reconstruction of the shear modulus, even in the presence of noise and partially known boundary conditions. Figure 1 illustrates an example of such a reconstruction, compared with other types of virtual fields proposed in the literature [2].

However, in practice, the search space is rarely known *a priori*. The stability estimate is valid only if $\mu^0$ is expressed in the chosen basis and, in particular, shows that as $n$ increases, the stability constant also grows, indicating that one cannot arbitrarily increase the dimension of the search space for a distributed parameter $\mu$.

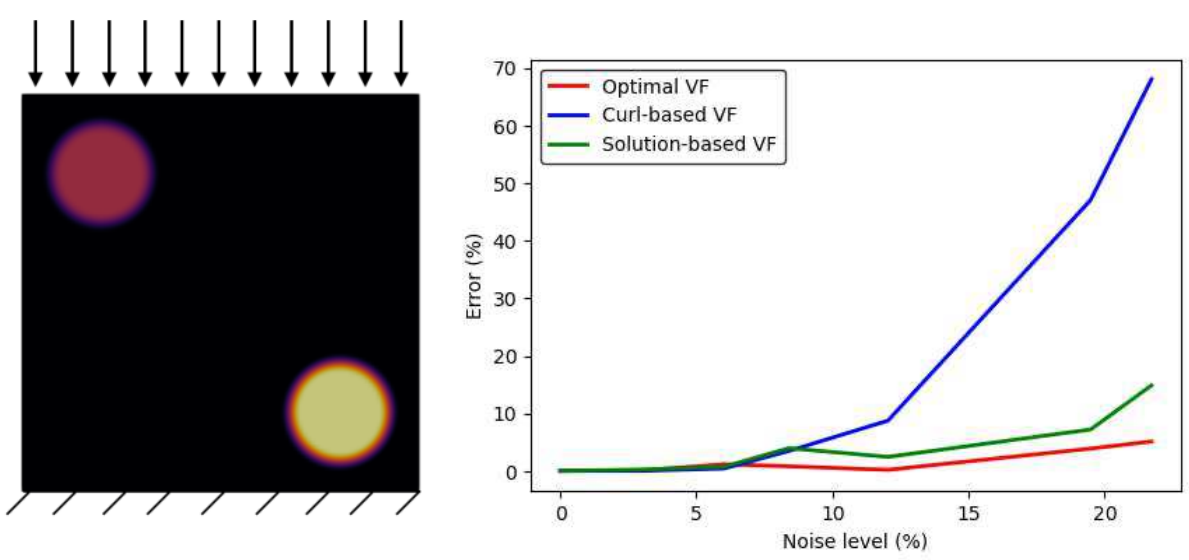

Figure 1: Shear modulus $\mu$ and comparison of the reconstruction with different virtual fields.

The *Adaptive Spectral Inversion* (ASI) method addresses this limitation by iteratively constructing a low-dimensional basis from the eigenfunctions of a suitable elliptic operator $L_\varepsilon[w]$ that depends on the medium [3]:

$$L_\varepsilon[w]v = -\nabla \cdot \big(q_\varepsilon[w]\nabla v\big),$$
$$q_\varepsilon[w](x) = \frac{1}{\sqrt{|\nabla w(x)|^2 + \varepsilon^2}}.$$

We propose a novel *ASI–VFM approach*, which combines the adaptive basis construction of ASI with the VFM. In other words, at each iteration $m$, given the previous estimate $\mu^{m-1}$, we construct the elliptic operator $L_\varepsilon[\mu^{m-1}]$ and search for $\mu^m$ in the basis formed by the eigenfunctions of $L_\varepsilon[\mu^{m-1}]$ using the VFM method. This coupling effectively mitigates pressure-related noise amplification. Numerical results confirm the efficiency of the proposed ASI–VFM method in nearly incompressible media.

# The return to the equilibrium of a floating cylinder in the Boussinesq regime.

**Geoffrey Beck**[1] **Ewan Contentin**[1,*], **Ludovic Martaud**[1]

[1]Univ Rennes, IRMAR UMR 6625, centre INRIA de l'université de Rennes (MINGuS

*Email: ewan.contentin@univ-rennes.fr

**Abstract**

In this work we study the interactions between waves at the surface of a sea and a partially immersed moving floating cylinder in the Boussinesq regime. The modelisation of the full waves-structure system can be written as a system of PDEs combined to a system of ODEs. Here we have an interest in the case of the return to the equilibrium in the linear regime, that is to say the solid is dropped in a fluid at rest without any speed. In this case we show that the solid returns to an equilibrium position in infinite time, but slower than $\mathcal{O}(t^{-1/2})$. Also we intend to give two numerical schemes for this system.

## 1 Modelisation

We study the case of a freely floating cylinder of radius $R > 0$ and of mass $m$ in a fluid. The cylinder is supposed to be allowed to move only vertically. We consider that there is no swirl in the fluid or at its surface, and that there is an axisymmetry with respect to the central axis of the cylinder.

With the previous considerations, the 2D situation can be written as a 1D problem with radial terms and we reduce the study to a slice of the problem, that is to say on $[0, +\infty)$. The behavior of the waves is governed by the dimensionless Boussinesq-Abbott system:

$$\begin{cases} \partial_t \zeta + \mathrm{d_r}\, q = 0, \\ D_\kappa \partial_t q + \partial_r \zeta = 0, \\ \mathrm{d_r} = \partial_r + \frac{d-1}{r}, \quad D_\kappa := (1 - \kappa^2 \partial_r \mathrm{d_r}) \end{cases} \tag{1}$$

where $d \in \{1, 2\}$ is the dimension of the problem. In fact we can treat both the 1D and the 2D situations simultaneously (see [3] for the 1D problem, and [1] for the 2D problem). The shallowness is denoted by $\kappa$. Mathematically $\kappa$ describes the dispersion. The variables are $\zeta$, for the elevation of the waves and $q$ for the discharge.

It was shown in [1, 3] that (1) can only be studied in the domain $r > R$ where the wave's elevation is free, with boundary condition $\underline{q} = -\frac{R}{2}\dot{\delta}$, where $\delta$ is the distance of the center of mass to its equilibrium and $\underline{q} = \lim_{\eta \to 0^+} q(R + \eta)$. Also $\delta$ is described by the Newton's equation.

## 2 Some theoritical results

Here we study the case of the return to the equilibrium in the linear case. That is to say at the initial time $t = 0$ the solid is dropped without any speed in a fluid at rest. All the following results are showed in [3] for $d = 1$ and in [1] for $d = 2$.

By inversing the dispersive operator $D_\kappa$ we obtain the *hidden equation*:

$$\kappa^2 \ddot{\underline{\zeta}} + \underline{\zeta} = \underline{\mathcal{R}\zeta} + B(r/\kappa)(\underline{\mathcal{R}\zeta} - \delta), \quad r \geq R. \tag{2}$$

where $B$ is given by decreasing exponential (resp. Bessel functions) if $d = 1$ (resp. if $d = 2$). The operator $\mathcal{R}$ is the dispersive operator such that for any $u \in L^2((R, +\infty))$ we define $\mathcal{R}u \in H^2((R, +\infty))$ by solving the elliptic problem:

$$(1 - \kappa^2 \mathrm{d_r} \partial_r)\mathcal{R}u = u, \quad \underline{\partial_r \mathcal{R}u} = 0. \tag{3}$$

The Newton's equation contains the pressure force exerted by the fluid under the solid. We can obtain an expression of this pressure term depending on $\zeta$. If we replace it in the Newton's equation and combine the result with the hidden equation it leads to:

$$M_{R,\kappa}\ddot{\delta} + \delta = \underline{\mathcal{R}\zeta} \tag{4}$$

where $M_{R,\kappa} = m + m_{add}$ where $m_{add} > 0$ is the *added mass*, representing the inertia of the solid.

We showed [1, 3] that the model is equivalent to the *augmented formula*, satisfied in the domain $r > R$. It is formulated with the *augmented variables* $u := (\zeta, q)$ for the fluid's variables and the boundary variable $Z := (\delta, \underline{\zeta})$:

$$\begin{cases} \partial_t u + \partial_r \mathcal{F}_\kappa[u] + \frac{1}{r}A[u] = \kappa(\underline{\mathcal{R}\zeta} - \delta)\gamma, \; r > R, \\ M\ddot{Z} + TZ = \mathfrak{F}(\underline{\mathcal{R}\zeta}), \\ \underline{u} = (\underline{\zeta}, -\frac{R}{2}\dot{\delta}), \quad (u, Z)_{|t=0} = (u_0, Z_0). \end{cases} \tag{5}$$

The terms $\mathcal{F}_\kappa$, $A$, $\gamma$, $\mathfrak{F}$, $M$ and $T$ are explicited in [1] and [3]. $\gamma$ is a dispersive boundary layer exponentially decreasing. Since $\underline{\mathcal{R}\zeta}$ requires to know $\zeta$ on $(R, +\infty)$, the PDEs system and the ODEs system are coupled.

All the equations (2), (4) and (5) have an expression in the nonlinear regime. In particular in (5) $\mathcal{F}_\kappa$ is still a flux in the nonlinear-regime. With the augmented formula we show that the linear system is globally wellposed and that the nonlinear system is locally wellposed.

In the return to the equilibrium the interesting quantity is $\delta$ and one wants to avoid computing the other variables. We showed in [1] that $\delta$ satisfies a Cummin's type equation both in 1D and 2D. It provides an expression of the Laplace transform of $\delta$ depending only on $\delta_{|t=0}$. Sharp estimates on this Laplace transform allows to show that $\delta$ lies in $H^2(\mathbb{R}^+)$ and we showed that $\delta$ and $\dot{\delta}$ decay to 0 slower than $t^{-1/2-\rho}$ for any $\rho \in (0, 1/2)$, *i.e* for large $t$, there exists $c > 0$ such that $|\delta(t), \dot{\delta}(t)| > ct^{-1/2-\rho}$.

## 3 Numerical simulations

The difficult part of the numerical simulations is the boundary $r = R$: the initial problem (1) gives only one boundary condition on $q$, given by $\delta$, coupled with the Newton's equation, but the classical finite volume method requires to know all the quantities at the boundary. We overpass this problem by using the hidden equation (2) in $r = R$. The theorical analysis provides several methods to make numerical simulations of the problem. A preliminary step consists in discretizing (3) as an elliptic problem.

### 3.1 The flux method.

The first method is based on (5) and is an adaptation of [4] in the 2D case. The simulations are done through a finite volume method, where the left boundary cell is filled with the order 2 ODE on $Z$, which gives $\underline{u}$. This scheme is explicit in time. In 1D, such a scheme needs adding an artificial viscosity to make the scheme stable. It is quite usual to do it, and in [4] the simulations were conclusive. Therefore in 2D we added viscosity with the same philosophy. The interesting point of this method is that it is adaptable to the non-linear regime because the nonlinearities just modify the flux. For this method, a computation of simultaneously $\zeta$ and $q$ is necessary on the contrary to the second one.

This method does not fully use the hidden equation (2). Indeed we only use it at $r = R$ whereas it is satisfied on the entier domain $r \geq R$. It leads to the second method.

### 3.2 The ODE method.

This one uses the ODE formulation, which is the combination of (2) and (4). Indeed thanks to the regularizing properties of $\mathcal{R}$, (2) is an ODE on each space's point. Thanks to this formulation, the only importance of the space discretization is contained in the discretization of (3). Therefore the only discussion is on the time discretization. The coupling terms ($\underline{\mathcal{R}\zeta}$ and $\delta$ in the hidden equation) are explicited. However the nonlinearities involve a coupling between $\zeta$ and $q$. It might be difficult to adapt the scheme to the nonlinear regime.

This scheme opens several perspectives: we almost get rid of the space discretization, that impacts the scheme only through $\mathcal{R}\zeta$ and the CFL can be independent on the space step. Also we can get rid of the first cell in $\zeta$ and avoid to compute it. Furthermore higher order schemes are reachable with this method.

**Aknowledgments**

The authors are supported by the BOURGEONS project, grant ANR-23-CE40-0014-01 of the French National Research Agency (ANR) and the Simons Foundation Award ID: 651475.

# A deep learning approach for optical cloaking device design

Elsie Cortes[1,*], Camille Carvalho[2], Symeon Papadimitropoulos[1], Chrysoula Tsogka[1]
[1]University of California Merced, Applied Mathematics Department, Merced CA 95434 USA
[2]INSA Lyon, ICJ UMR5208, CNRS, Ecole Centrale de Lyon, Universite Claude Bernard Lyon 1, Université Jean Monnet, 69621 Villeurbanne, France
*Email: ecortes7@ucmerced.edu

**Abstract**

We propose a deep learning framework based on an encoder–decoder architecture for the design of optical cloaking devices. The cloaking strategy relies on concentric layers surrounding the object to be cloaked, whose geometry and material parameters must be determined to minimize scattering. To reduce the number of unknowns, we restrict the cloak design to either object-fitted layers or circular layers, so that the unknown parameters are limited to the layer radii and material properties. We consider wave propagation in 2D multilayered media modeled by the Helmholtz equation. Datasets are generated using a boundary element method that, together with a close evaluation treatment, enables accurate field computation on closely spaced interfaces without refined discretization, significantly reducing computational costs and enabling efficient dataset generation. We demonstrate successful cloak design and show that object-fitted cloaking layers generally outperform circular geometries within the same NN architecture. The proposed framework is flexible and can be extended to more complex cloaking scenarios.

## 1 Problem Formulation

We study optical cloaking by analyzing wave propagation through $N+1$ concentric material regions with smooth, nested boundaries $\Gamma_j$ of radii $r_j$ and permittivities $\varepsilon_j$, for $j = 0, \dots, N$ at angular frequency $\omega$ (see Figure 1). We assume non-magnetic media, i.e. $\mu_j = 1$ in all regions. In each region, the wavenumber $k_j = \omega\sqrt{\varepsilon_j}$. The cloaking design problem consists in determining the characteristic parameters of all the concentric layers encapsulating an object characterized by radius $r_{\text{ob}}$ and permittivity $\varepsilon_{\text{ob}}$ so as to minimize the scattered field. In each layer, the total field $u_j$ satisfies the Helmholtz equation $\nabla^2 u_j + k_j^2 u_j = 0$. The incident field is a plane wave $u_0^{\text{in}} = e^{ik_0 \mathbf{a}\cdot\mathbf{x}}$, where $\mathbf{a} = (\cos\alpha, \sin\alpha)$ is the incidence direction. Across each interface $\Gamma_j$, the field satisfies transmission conditions s.t. $[u_j]_{\Gamma_j} = 0$ and $[\varepsilon_j^{-1}\partial_n u_j]_{\Gamma_j} = 0$ and Sommerfeld radiation condition is imposed in the unbounded region.

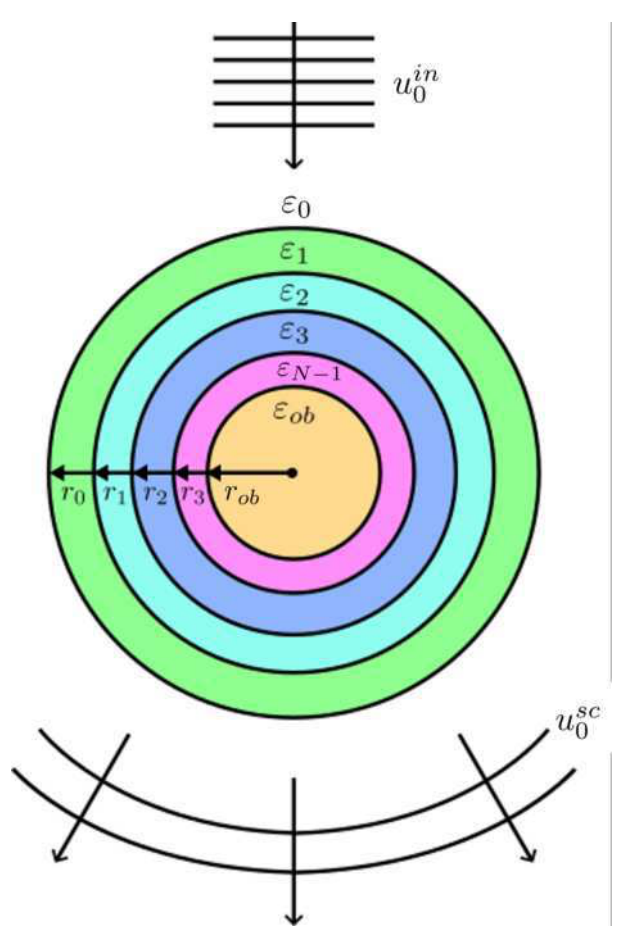

Figure 1: Schematic and notations of the problem.

## 2 Deep Learning Approach with BEM

Neural network (NN) approaches have been explored for cloaking design, and encoder–decoder architectures have shown promise in determining parameters for cloaking [1]. In contrast to traditional inverse scattering methods which require expensive, repeated PDE solves, the NN approach shifts the computational cost to the initial training of the model for fast inference. We propose the following general encoder–decoder framework (Figure 2). The decoder solves the forward problem, mapping the vectors $\Delta\mathbf{r} = (r_j - r_{j+1})_{j=0}^{N-2}$ and $\mathbf{e} = (\varepsilon_j)_{j=1}^{N-1}$ containing the layer thickness and permittivities of the cloak, respectively, to the field measured at a set of randomly placed points outside the cloaking device. The encoder solves the inverse scattering problem by mapping the field to the corresponding cloaking parameters. The decoder is trained first and the encoder is subsequently trained so that the combined encoder–decoder acts as the identity.

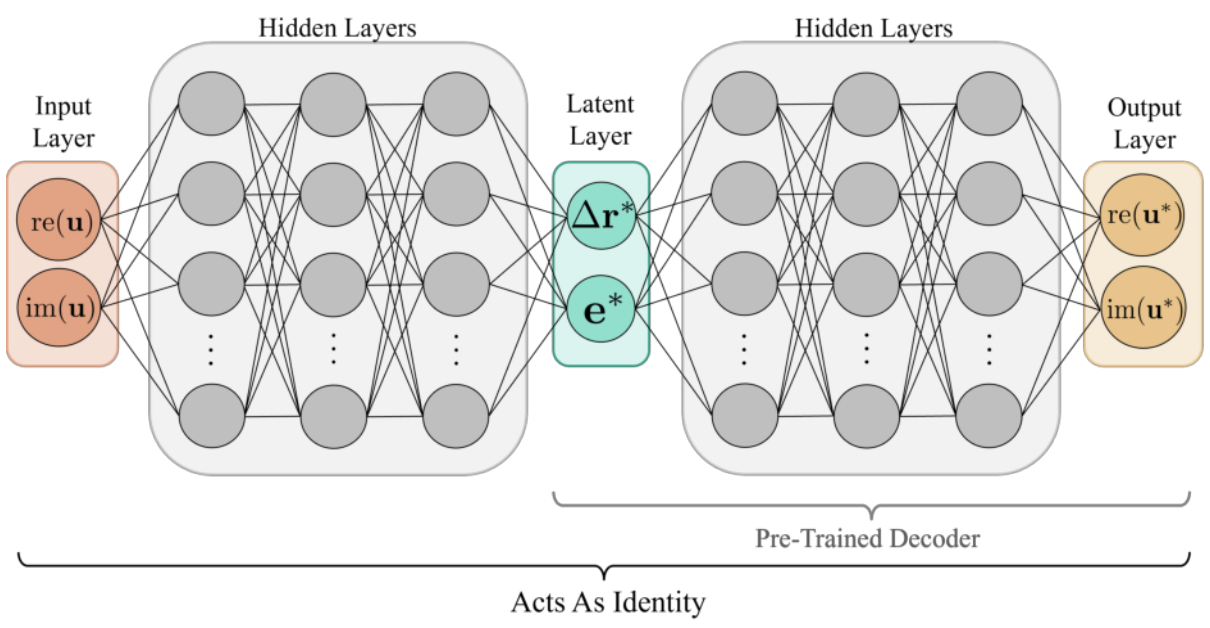

Figure 2: Combined encoder-decoder architecture.

Training data generation for the NN is computationally demanding. We adopt a Boundary Element Method (BEM), which is computationally efficient because the unknowns are defined only on the interfaces, allowing field evaluation at arbitrary points in the domain [2]. However, for thin layers, this approach can exhibit close-evaluation errors arising from nearly singular interactions [3]. We mitigate these errors using a modified boundary integral representation method (BE-BRIEF) [4].

## 3 Assessing Design Results

We consider an object with a star-shaped geometry. A natural question in cloaking design is how the surrounding layers should be shaped. Circular layers are often preferred because their simpler design is commonly assumed to be optimal for arbitrary objects. However, one may ask whether an object-fitted cloak, whose layers conform to the star-shaped object boundary, can achieve superior scattering reduction. In the following we fix $N = 5$, the object parameters $r_{\text{ob}} = 1$, $\varepsilon_{\text{ob}} = 1$, and the permittivity of the background medium $\varepsilon_0 = \sqrt{2}$. The parameters of the circular/star-shaped cloaking layers, $(r_{j-1}, \varepsilon_j)_{j=1}^{N-1}$, are determined by the NN.

To systematically quantify and compare the cloaking performance of various cloaking designs, we compute the two-dimensional radar cross section (RCS) that measures the object's visibility to the incident wave: larger RCS values correspond to stronger scattering and higher detectability. The RCS is evaluated on a circle of radius $R = m\lambda_0 = m\frac{2\pi}{k_0}$, where $\lambda_0$ is the background wavelength and $m$ denotes the observation distance (we use $m = 5$). The 2D RCS is defined as $\sigma(\theta) = 2\pi R\frac{|u_0(R,\theta)-u^{\text{in}}(R,\theta)|^2}{|u^{\text{in}}(R,\theta)|^2}$.

The NN cloaking design outputs for both

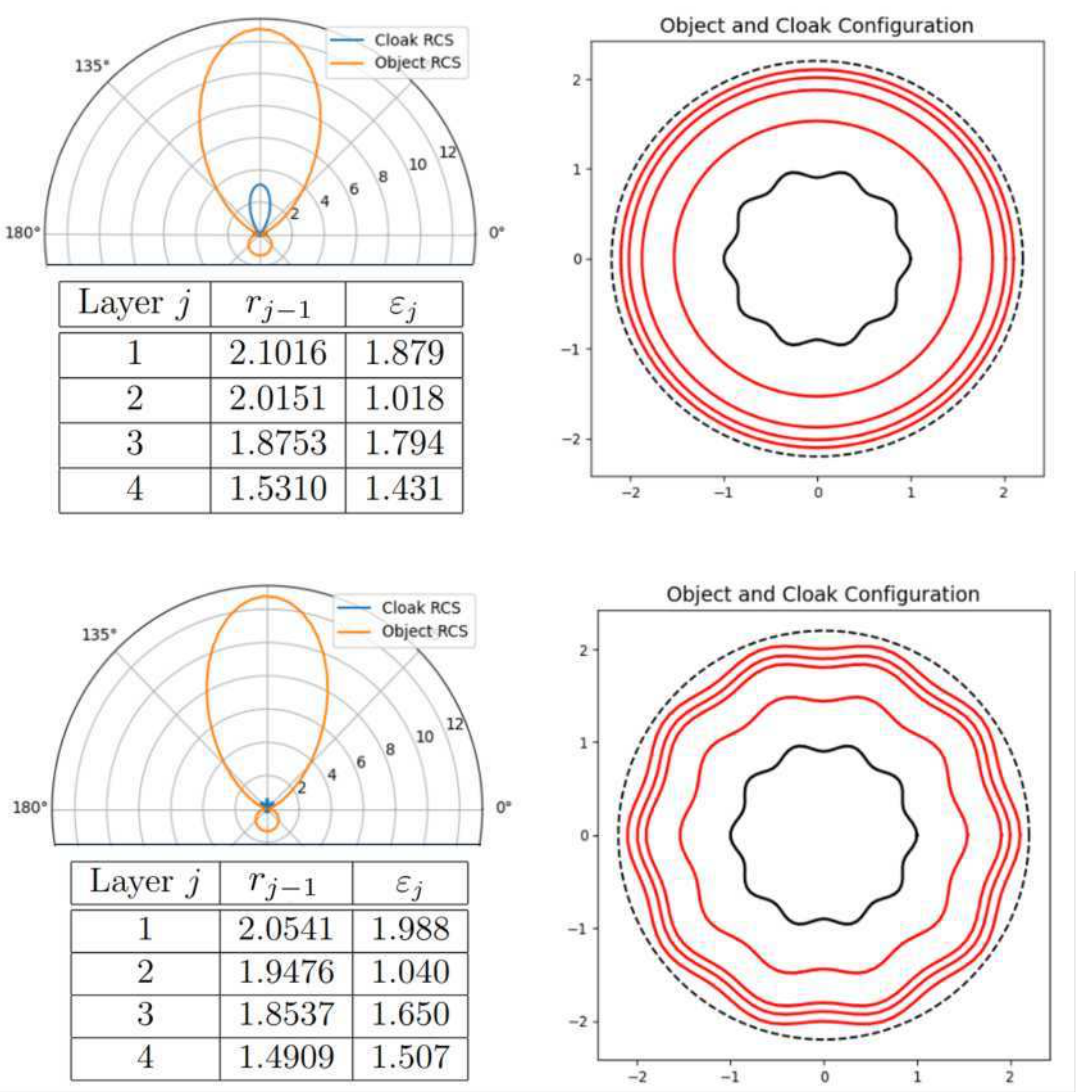

| Layer $j$ | $r_{j-1}$ | $\varepsilon_j$ |
|---|---|---|
| 1 | 2.1016 | 1.879 |
| 2 | 2.0151 | 1.018 |
| 3 | 1.8753 | 1.794 |
| 4 | 1.5310 | 1.431 |

| Layer $j$ | $r_{j-1}$ | $\varepsilon_j$ |
|---|---|---|
| 1 | 2.0541 | 1.988 |
| 2 | 1.9476 | 1.040 |
| 3 | 1.8537 | 1.650 |
| 4 | 1.4909 | 1.507 |

Figure 3: RCS plots and cloaking parameters (left) for given cloaking device configurations (right) with object interface in black and cloak interfaces in red.

types of cloak interface geometries are shown in Figure 3. The predicted parameters are similar in both cases. However, the star-shaped cloak achieves improved scattering reduction, suggesting object-fitted layer geometries may yield more effective cloaking performance.

The proposed NN architecture has been used to effectively design cloaks for circular, star-shaped, and kite-shaped objects, including cases in which orientation and geometric concavity influence scattering behavior. Based on these results, we believe that the proposed framework can serve as a useful tool for designing and assessing cloaking configurations across a wide range of object geometries and physical settings.

# Preconditioning for the DG plane wave method.

**Joseph Coyle**[1,*], **Nilima Nigam**[2]
[1]Department of Mathematics, Monmouth University,West Long Branch, USA
[2]Department of Mathematics, Simon Fraser University, Burnaby, Canada
*Email: jcoyle@monmouth.edu

**Abstract**

We present recent results on the preconditioning of linear systems arising from plane wave DG methods. Inspired by previous work on preconditioning on a single element, we consider block-diagonal preconditioning.

## 1 Introduction

Plane wave Trefftz methods (PWB) can offer significant advantages over approaches based on polynomials in terms of accuracy [5], at the expense of the poor conditioning of the system matrices. One fix is the use of evanescent waves [6] or the use of preconditioners, [1]. In previous work [3], we investigated the behaviour of the PWB for the Helmholtz equation focusing on the mass and system matrices on a *single* element, extending the work of [7]. Based on insights gained from circular elements and cyclic polygons, we designed several circulant preconditioners and examined their performance. We now report on the behaviour of block-diagonal preconditioners based on these strategies.

In a Trefftz method, we seek $u_p \in V_p(\mathcal{T})$ where we set $\mathcal{T}$ to be a mesh of elements $K$,

$$V_p := \left\{ v \in L^2(\mathcal{T}) \, : v(x)|_K = \textstyle\sum_{j=1}^{p} \alpha_j e^{i\omega d_j \cdot (x - x_c)} \right\},$$

and $d_j$ are the plane wave directions in a given element. The standard DG formulation leads to a linear system, denoted by $\mathtt{SyS}\,\mathtt{u} = \mathtt{f}$.

## 2 Element-wise preconditioning

On a single element $K$, the local interaction leads to three types of matrices: mass, cross, and stiffness matrices whose entries are given by

$$[\mathtt{M}]^K_{m,\ell} := \int_{\partial K} \varphi_m \overline{\varphi_\ell} ds,$$

$$[\mathtt{S}]^K_{m,\ell} := \int_{\partial K} \nabla \varphi_m \cdot \mathbf{n}\, \overline{\varphi_\ell} ds,$$

and

$$[\mathtt{D}]^K_{m,\ell} := \int_{\partial K} \nabla \varphi_m \cdot \mathbf{n}\, \overline{\nabla \varphi_\ell} ds,$$

respectively. The conditioning of the system matrix $\mathtt{SyS}$ on $K$ deteriorates sharply and unavoidably with the number of plane wave $p$ and the size of the element $h$.

**Theorem 1** *Let $K$ be a disk of radius $h$ centered at the origin, and consider the cross matrix resulting from $p$ (odd) evenly spaced plane waves. $\mathtt{S}^{disk}$ is circulant, and the following hold:*

1. *The generating function of the cross matrix is*

$$F^{disk}_{\mathtt{S}} := 2\pi\kappa h^2 \sin(x/2)\, J_1(2\kappa h \sin(x/2)).$$

2. *There is a simple eigenvalue,*

$$\lambda_0(\mathtt{S}^{disk}) = \sum_{j=0}^{p-1} F^{disk}_{\mathtt{S}}(2j\pi/p) \approx p\kappa h^2 \int_0^{2\pi} \sin(x/2)\, J_1(2\kappa h \sin(x/2))\, dx.$$

3. *The eigenvalues $\lambda_L(\mathtt{S}^{disk}), L = 1, ..., (p-1)/2$ of $\mathtt{S}^{disk}$ have multiplicity 2 and are given by*

$$\lambda_L(\mathtt{S}^{disk}) = \sum_{j=0}^{p-1} F^{disk}_{\mathtt{S}}(2j\pi/p)\, e^{-i2jL\pi/p} \approx -2ph \int_0^{2\kappa h} \left(\sqrt{4\kappa^2 h^2 - z^2}\right) J''_{2L}(z)\, dz.$$

From this we see the condition number is approximated by

$$|\lambda_{max}/\lambda_{min}| \approx C\, \left((\Gamma(p-1)/p\kappa h^2)(\kappa h)^{1-p}\right).$$

Similar results were obtained for the mass and stiffness matrices as well.

On a general element, these matrices are neither Toeplitz nor circulant. We observe, however, that as the number of plane waves $p$ increases, these matrices *tend* towards Toeplitz.

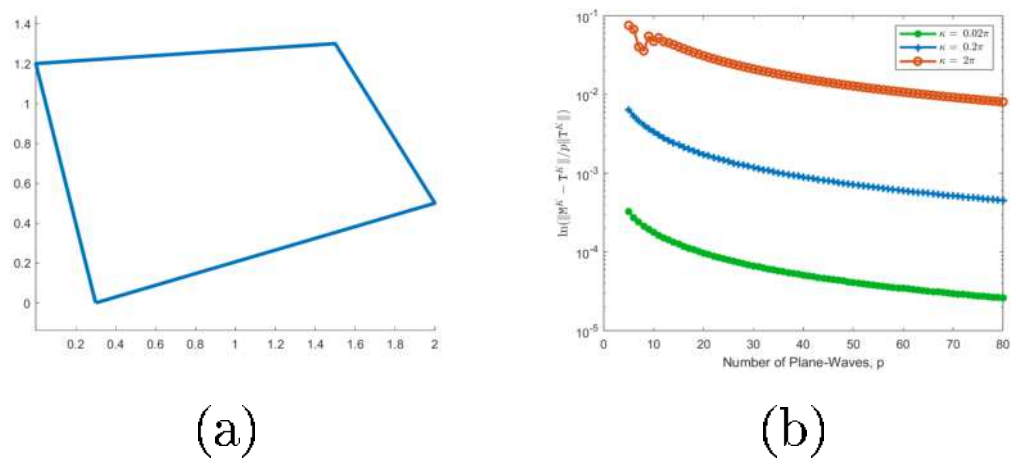

Figure 1: Plot (b) of $\ln\left(\|\mathtt{M}^K - \mathtt{T}^K\|/p\|\mathtt{T}^K\|\right)$ as a function of the number of plane-waves, $p = 0.02\pi, 0.2\pi, 2\pi$. K is shown in (a).

For instance, we can construct a Toeplitz matrix $T^K$ using the average of the diagonals of $\mathtt{M}^K$. A measure of how close to Toeplitz the mass matrix $\mathtt{M}^K$ is depicted in Figure 1. This suggests the use of *circulant* preconditioners.

We highlight three preconditioners from the several considered in [3]. The first is based on the system matrix on a disk of the same size $h$ as $K$,

$$\mathtt{P}^K_{disk} := \mathtt{U}(\Lambda(\mathtt{SyS}^{disk}))^{-1/2}\mathtt{U}^*$$

and the second is a circulant matrix generated from the average value of the diagonals of the mass matrix on the element $K$:

$$\mathtt{P}^K_{best} := \mathtt{U}(\kappa^2\Lambda\left(\mathtt{C}^K_{best}\right)^{-1/2})\mathtt{U}^*,$$

where $\mathtt{C}^K_{best} := \mathrm{circ}[c_0, c_1, ..., c_{p-1}]$ and

$$c_{i-j} = (1/n) \sum_{p-q=i-j(mod\, n)} [\mathtt{M}^K]_{pq}.$$

Finally, we consider a singular preconditioner [4]

$$\mathtt{P}^K_{reg} := (1/\kappa)\mathtt{U}\,\mathtt{diag}(1/\mu_1, ..., 1/\mu_i, 0, 0, ...0)^{1/2}\mathtt{U}^*,$$

where the $\mu_i$ are the eigenvalues of the mass matrix on the disk $\mathtt{M}^{disk}$ so that $\mu_i \geq \delta$ in magnitude for some $\delta > 0$. We note $\mathtt{U}$ is a unitary Fourier matrix, and since both $\mathtt{SyS}^{disk}$ and $\mathtt{C}^K_{best}$ are circulant, their eigenvalues are computed efficiently in closed form.

## 3 Trefftz methods on a mesh

Based on the effectiveness of the element-level preconditioners $\mathtt{P}^K_{disk}$, $\mathtt{P}^K_{best}$, and $P^K_{reg}$ in terms of both conditioning and error reduction, we propose *block diagonal preconditioners* $\mathtt{P}^{disk}$, $\mathtt{P}^{best}$ and $\mathtt{P}^{reg}$. In Figure 2, we see the performance of these preconditioners on triangular meshes.

The preconditioners are block-wise circulant and these experiments suggest they offer significant reduction in condition numbers.

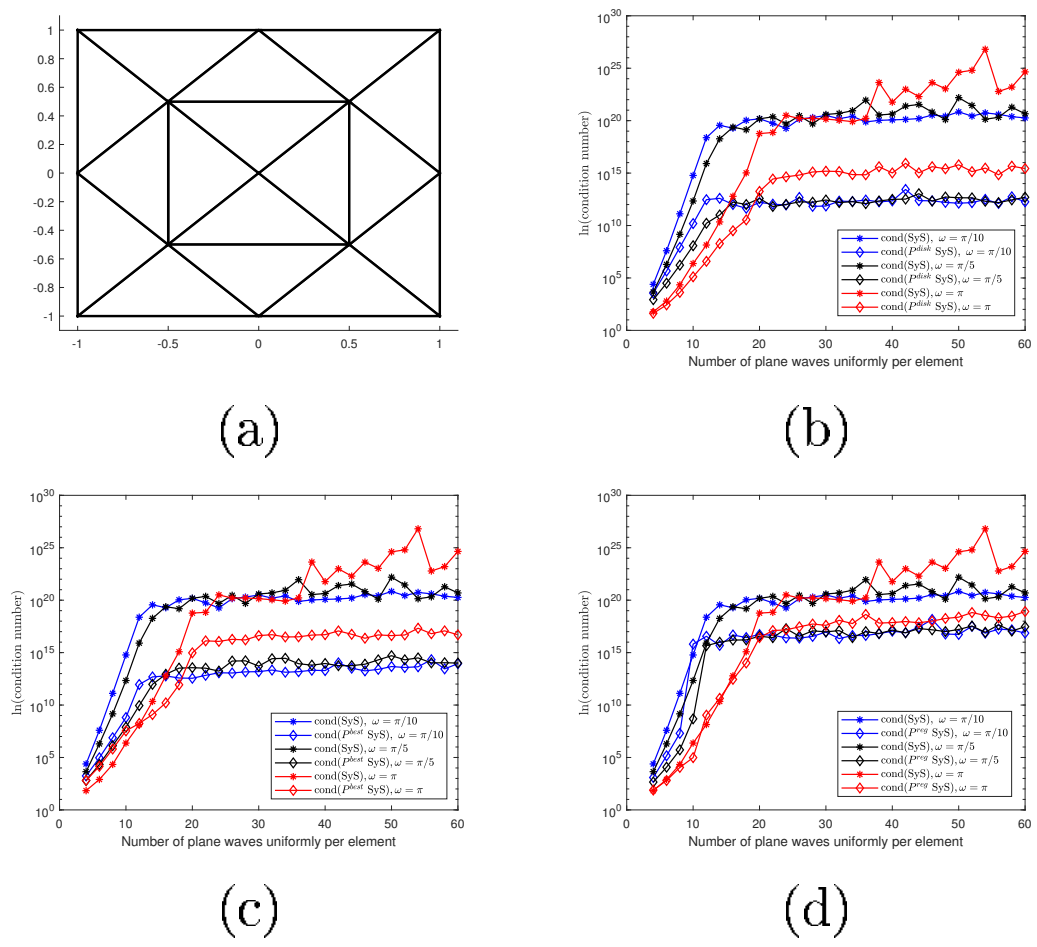

Figure 2: Employing the mesh shown in (a), the condition number of the system matrix $\mathtt{SyS}$ (black curves), and the block-preconditioned matrices $\mathtt{P}^{\mathtt{disk}}\mathtt{SyS}$ (b), $\mathtt{P}^{\mathtt{best}}\mathtt{SyS}$ (c), and $\mathtt{P}^{\mathtt{reg}}\mathtt{SyS}$ (d) with $\delta = 10^{-9}$.

# Space-time BEM based coupling strategies for the analysis of infrastructures solicited by seismic waves

**Giulia Di Credico**[1,*], **Alessandra Aimi**[1], **Heiko Gimperlein**[2], **Chiara Guardasoni**[1], **Francesco Freddi**[3], **Roberto Valentino**[4]

[1]Department of Mathematical, Physical and Computer Sciences, University of Parma, ITALY
[2]Department of Engineering Mathematics, University of Innsbruck, AUSTRIA
[3]Department of Engineering and Architecture, University of Parma, ITALY
[4]Department of Chemistry, Life Sciences and Environmental Sustainability, University of Parma, ITALY

*Email: giulia.dicredico@unipr.it

**Abstract**

We want to develop a numerical coupling method with applications in road infrastructures, specifically concrete bulkheads in contact with the ground and potentially lifted by seismic phenomena. The interaction between materials (concrete and soil) is solicited by the propagation of elastic waves, which, in preexisting engineering analysis literature, are often simplified with a series of elastic springs. To have a realistic inclusion of the elastic solicitations from the ground, the proposed coupling strategy exploits Boundary Element Method for the unbounded soil domain, allowing a natural incorporation of radiation conditions at infinity. Non-penetration and frictional contact constraints are considered at the interface between the concrete and the soil, leading to a reformulation of the time-dependent mathematical problem at hand by a variational inequality solvable with an iterative solver, such as the Uzawa approach.

## 1 The 2D mathematical model problem

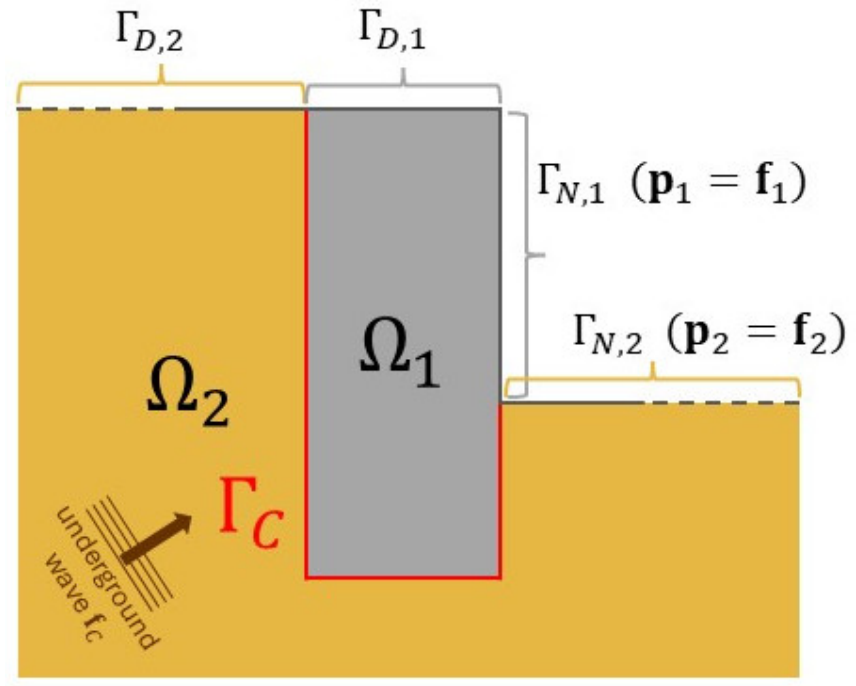

Figure 1: Representation of the geometry of interest

We consider a 2D section of a concrete bulkhead partially embedded in the ground. The geometry is schematically represented, as in Figure 1, by two adjacent domains $\Omega_1, \Omega_2 \in \mathbb{R}^2$ that have boundaries, denoted as $\partial\Omega_1$ and $\partial\Omega_2$, which include a non-trivial contact interface $\Gamma_C = \partial\Omega_1 \cap \partial\Omega_2$. Especially on $\Gamma_C$, an incoming pressure wave $\mathbf{f}_C$ with source in the unbounded domain $\Omega_2$ (the soil) causes elastic deformations and gaps. The dynamic over time is governed by the Navier-Lamé equations for the displacements $\mathbf{u}_1$ and $\mathbf{u}_2$ of the interacting patches:

$$\rho_j\, \partial_t^2 \mathbf{u}_j = \nabla \cdot \sigma_j(\mathbf{u}_j) \quad \text{in } [0,T] \times \Omega_j\,, \tag{1}$$

with $j = 1, 2$. In (1), $\sigma_j$ is the Cauchy stress tensor, incorporating the physical properties of the medium in $\Omega_j$, and $\rho_j$ the corresponding mass density. For each $j = 1, 2$, the boundary of $\Omega_j$ is decomposed as the union $\partial\Omega_j = \Gamma_C \cup \Gamma_{N,j} \cup \Gamma_{D,j}$, where the letters $N$ and $D$ indicate the different kind of boundary constraints (Dirichlet and Neumann), in particular:

$$\begin{aligned} &\mathbf{u}_j = \mathbf{0} \quad \text{in} \quad [0,T] \times \Gamma_{D,j}, \\ &\mathbf{p}_j = \sigma(\mathbf{u}_j) \cdot \mathbf{n}_j = \mathbf{f}_j \quad \text{in} \quad [0,T] \times \Gamma_{N,j}, \end{aligned} \tag{2}$$

being $\mathbf{n}_j$ the outpointing normal direction defined on $\partial\Omega_j$. The first Dirichlet constraint simulates the blockage of soil and concrete (possibly due to the presence of a road above the bulkhead). On $\Gamma_C$ the following non-linear conditions are prescribed to include in the model the non-penetration of the two interacting domains

$$\begin{aligned} &u_{1,\perp} \geq u_{2,\perp},\ p_{1,\perp} \geq 0,\ p_{2,\perp} \leq f_{C,\perp}, \\ &\begin{cases} p_{1,\perp} = f_{C,\perp} - p_{2,\perp} = 0 \text{ if } u_{1,\perp} > u_{2,\perp} \\ p_{1,\perp} + p_{2,\perp} = f_{C,\perp} \text{ if } u_{1,\perp} = u_{2,\perp} \end{cases}. \end{aligned} \tag{3}$$

Frictional forces contrasting the rubbing of the two materials allow the opening of a small gap between the two materials:

$$\begin{aligned} &p_{1,\parallel} + p_{2,\parallel} - f_{C,\parallel} = 0,\ |p_{1,\parallel}| \leq \mathcal{F}, \\ &p_{1,\parallel} \cdot (\dot{u}_{1,\parallel} - \dot{u}_{2,\parallel}) + \mathcal{F}\, |\dot{u}_{1,\parallel} - \dot{u}_{2,\parallel}| = 0. \end{aligned} \tag{4}$$

The symbols $\perp$ and $\parallel$ in (3) and (4) denote the orthogonal and tangential parts of the related vector field, and the constant $\mathcal{F} \geq 0$ (4) represents the friction threshold.

## 2 Weak formulation for a BEM-BEM coupling perspective

At first, we set the above written Neumann data in proper space-time anisotropic Sobolev spaces:

$$\mathbf{f}_i \in H^{\frac{1}{2}}([0,T], \widetilde{H}^{-\frac{1}{2}}(\Gamma_{N,i}))^2 \text{ for } i = 1,2,$$
$$\mathbf{f}_C \in H^{\frac{1}{2}}([0,T], \widetilde{H}^{-\frac{1}{2}}(\Gamma_C))^2.$$

Then, defining the subset $\Gamma_{\Sigma_1} := \Gamma_C \cup \Gamma_{N,1}$ and the tensor product of functional spaces

$$\mathcal{C} := H^{\frac{1}{2}}([0,T], \tilde{H}^{\frac{1}{2}}(\Gamma_{\Sigma_1}))^2 \times H^{\frac{1}{2}}([0,T], \tilde{H}^{\frac{1}{2}}(\Gamma_C))^2 \times H^{\frac{1}{2}}([0,T], \tilde{H}^{\frac{1}{2}}(\Gamma_{N,2}))^2 \times M^+(\mathcal{F}),$$

with $M^+(\mathcal{F})$ containing a suitable Lagrange multiplier $\boldsymbol{\lambda}$ having compact support in $\Gamma_C$, a weak mixed formulation equivalent to the strong problem (1)-(4) reads as:
*find* $(\mathbf{u}_1, \widetilde{\mathbf{u}}, \mathbf{u}_2, \boldsymbol{\lambda}) \in \mathcal{C}$ *such that*

$$\langle \mathcal{S}_1 \mathbf{u}_1, \mathbf{v}_1 \rangle_{\Gamma_{\Sigma_1}} + \langle \mathcal{S}_2 \mathbf{u}_2, \mathbf{v}_2 \rangle_{\Gamma_{N,2}} + \langle \mathcal{S}_2 \mathbf{u}_2, \mathbf{v}_1 - \widetilde{\mathbf{v}} \rangle_{\Gamma_C} = \langle \boldsymbol{\lambda}, \widetilde{\mathbf{v}} \rangle_{\Gamma_C} + \langle \mathbf{f}_1, \mathbf{v}_1 \rangle_{\Gamma_{N,1}} + \langle \mathbf{f}_2, \mathbf{v}_2 \rangle_{\Gamma_{N,2}} + \langle \mathbf{f}_C, \mathbf{v}_1 - \widetilde{\mathbf{v}} \rangle_{\Gamma_C} \tag{5a}$$

$$\langle \boldsymbol{\mu} - \boldsymbol{\lambda}, \partial_{t,\parallel} \widetilde{\mathbf{u}} \rangle_{\Gamma_C} \geq 0, \tag{5b}$$

*for all* $(\mathbf{v}_1, \widetilde{\mathbf{v}}, \mathbf{v}_2, \boldsymbol{\mu}) \in \mathcal{C}$. The function $\widetilde{\mathbf{u}}$ is defined to set $\mathbf{u}_{2|_{\Gamma_C}} = \mathbf{u}_{1|_{\Gamma_C}} - \widetilde{\mathbf{u}}$, i.e. to express the displacement of the second domain in dependence on $\mathbf{u}_{1|_{\Gamma_C}}$ and the gap at the contact interface. The products involved in (5) are the standard integrals $\langle \mathbf{p}, \mathbf{u} \rangle_\Gamma = \int_\Gamma \int_0^T \mathbf{p} \cdot \mathbf{u} \, dt \, d\Gamma_{\mathbf{x}}$, while operators $\mathcal{S}_i$ acting on $\mathbf{u}_i$ and defined as $\mathcal{S}_i \mathbf{u}_i := \mathbf{p}_i$ are the so called Poincarè-Steklov operators, defined at the different boundaries of the interacting domains. For each $i = 1,2$, $\mathcal{S}_i$ is assembled through the composition of the well known boundary integral operators of elastodynamics $\mathcal{V}_i$, $\mathcal{K}_i$ and $\mathcal{W}_i$:

$$\mathcal{S}_i = (\mathcal{K}_i^* + 1/2)\, \mathcal{V}_i^{-1} (\mathcal{K}_i + 1/2) - \mathcal{W}_i. \tag{6}$$

An algebraic form of the weak problem (5) is then constructed on the basis of a space-time discretization of the boundaries $\partial\Omega_1$ and $\partial\Omega_2$ and the time interval $[0,T]$. This induces the approximation of the quartet $(\mathbf{u}_1, \widetilde{\mathbf{u}}, \mathbf{u}_2, \boldsymbol{\lambda}) \in \mathcal{C}$ as a linear combination of space-time piecewise polynomial basis. In particular, time approximation leads to the construction of a linear system with Toeplitz block structure, allowing a block forward time-stepping scheme. The resulting mixed discrete problem can be then solved with a standard solver such as the Uzawa iterative algorithm, as already implemented in [1, 2]. A more sophisticated weak formulation (5) and the treatment of complex geometries, such as half-space domains correlated with mixed boundary conditions and triple nodes as in figure 1, are the improvement with respect to the cited works. The perspective is then the analysis of realistic engineering structures with an efficient BEM based coupled method.

## 3 Preliminary results

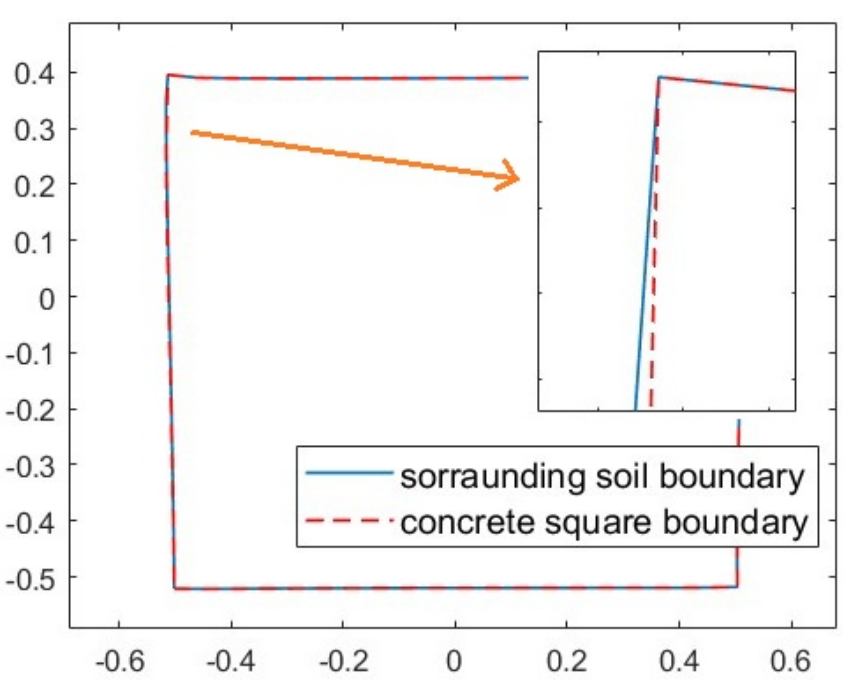

Figure 2: Deformation of a concrete square embedded in a surrounding soil domain

Here, we consider a simplified model, where a square concrete block $\Omega_1$, with sides of length 1 m, is fully embedded in the soil, namely $\Omega_2 = \mathbb{R}^2 \setminus \Omega_1$. The bottom and left sides of the square represent the contact boundary $\Gamma_C$, while the remaining boundary portion is a simple transmission interface, not included in Figure 1. A pure vertical force $\mathbf{f}$ is applied on the top side, compressing the interior domain. Figure 2 shows boundary displacements $\mathbf{u}_1$ and $\mathbf{u}_2$, magnified by a factor of $10^6$. The detail of the opening gap between the two materials is highlighted, confirming the effectiveness of the resolution algorithm applied to the discrete version of the mixed problem (5a)-(5b).

# The largest ever Helmholtz simulations in inhomogeneous media

**Max Cubillos**[1,*], **Edwin Jimenez**[1]

[1]Directed Energy Directorate, Air Force Research Laboratory, Albuquerque, New Mexico, USA

*Email: afrl.rdl.orgbox@us.af.mil

**Abstract**

We report on what are, to our knowledge, the largest ever simulations of time-harmonic wave propagation in three-dimensional heterogeneous media—both in terms of the physical size of the problems in wavelengths and the number of degrees of freedom (some simulations consist of over a trillion unknowns). Our algorithms use an HPC implementation of the Operator Fourier Transform (OFT), which is a functional calculus framework for the representation of pseudodifferential operators. Within the OFT framework, inverting the Helmholtz operator reduces to two sequential solves of a pseudo-time paraxial (or Schrödinger) equation and superimposing these solutions via pseudo-time Fourier integrals. An attractive feature of this approach is that efficient parallel scaling can be achieved with domain decomposition methods for time-dependent partial differential equations; we demonstrate implementations using finite differences and high-order spectral methods as examples.

## 1 Introduction

We are concerned with numerical solution of the Helmholtz equation in an unbounded domain

$$\Delta v(\boldsymbol{x}) + \kappa^2 m(\boldsymbol{x}) v(\boldsymbol{x}) = 0, \tag{1}$$

where $\kappa$ is the wavenumber and $m(\boldsymbol{x})$ is a spatially varying refraction coefficient accounting for the spatial dependence of the wave velocity. The Helmholtz equation is notoriously challenging to solve with direct and iterative methods. However, a scalable direct solver was recently developed using a functional calculus framework we call the Operator Fourier Transform (OFT) [1]. We describe the approach in the following section.

## 2 OFT applied to the Helmholtz equation

Writing the solution $v$ of equation (1) as the sum of a known incident field $v^i$ and scattered field $v^s$ and truncating the domain with first order nonreflecting boundary conditions for simplicity, the equation for the scattered field is then

$$\begin{cases} \big(m + \kappa^{-2}\Delta\big) v^s = g, & \boldsymbol{x} \in D, \\ \dfrac{\partial v^s}{\partial n} - i\kappa v^s = 0, & \boldsymbol{x} \in \partial D, \end{cases} \tag{2}$$

where the source term is given by

$$g(\boldsymbol{x}) = -\big(m(\boldsymbol{x}) - 1\big) v^i(\boldsymbol{x}). \tag{3}$$

Letting $A = m(\boldsymbol{x})I + \kappa^{-2}\Delta$ denote the variable coefficient Helmholtz operator, we may formally write the solution to equation (2) as

$$v^s(\boldsymbol{x}) = \big[f^2(A)g\big](\boldsymbol{x}), \tag{4}$$

where $f(A) = 1/\sqrt{A}$ is a pseudodifferential operator. To evaluate $f(A)$, we use the operator Fourier transform, defined as

$$\begin{aligned} f(A) &= \frac{1}{\sqrt{2\pi}} \int_{-\infty}^{\infty} \widehat{f}(\tau) e^{i\tau A} d\tau \\ &= \sqrt{\frac{-i}{\pi}} \int_0^{\infty} \frac{e^{i\tau}}{\sqrt{\tau}} e^{i\tau(A-I)} d\tau, \end{aligned} \tag{5}$$

where $\widehat{f}(\tau)$ is the Fourier transform of the scalar function $f$ and $u(\boldsymbol{x}, \tau) = e^{i\tau(A-I)} g(\boldsymbol{x})$ is interpreted as the solution at time $t = \tau$ of the initial-boundary-value problem

$$\begin{cases} \dfrac{\partial u}{\partial t} = i(A - I)u, & (\boldsymbol{x}, t) \in D \times (0, \infty), \\ \dfrac{\partial u}{\partial n} - i\kappa u = 0, & (\boldsymbol{x}, t) \in \partial D \times (0, \infty), \\ u = g, & (\boldsymbol{x}, t) \in D \times \{0\}, \end{cases} \tag{6}$$

## 3 Parallel OFT solver

A numerical method for equations (5) and (6) is described in [1]; we briefly review the algorithm together with the differences in the parallel case.

The algorithm in a single domain involves three main elements:

1. Exponentially increasing time steps of the form $t_m = a(b^m - 1)$ for some positive constants $a, b$.

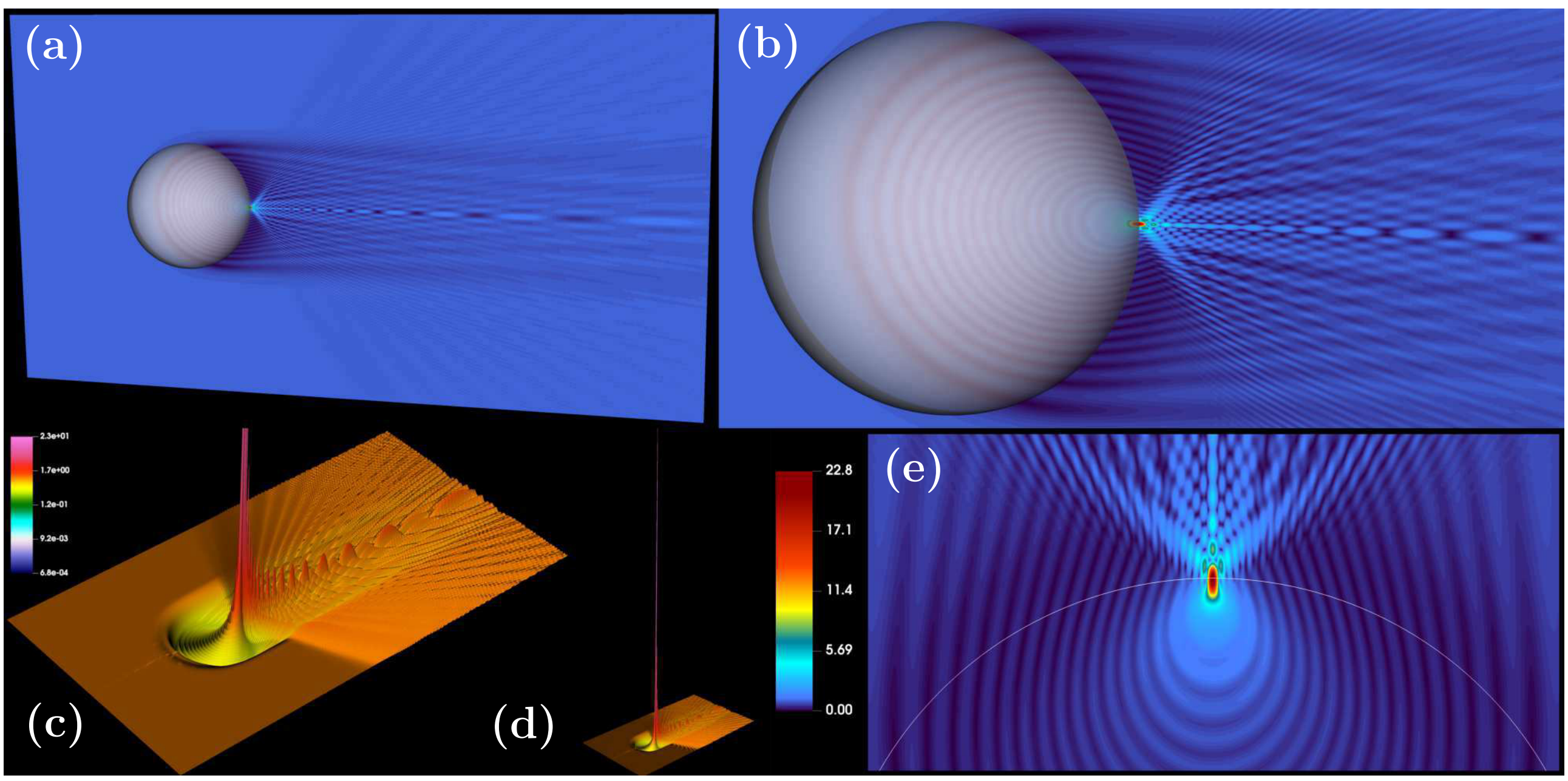

Figure 1: (a), (b) 2D slice of the total field magnitude for plane wave transmission through a Luneburg lens. Elevated profiles of the 2D slice showing (c) a close up of the near field interactions and (d) the large magnitude of the focal spot. (e) Close up of 2D slice in the focal region.

2. Backward Euler alternating direction implicit (BDF1-ADI) scheme, together with Richardson extrapolation to achieve second-order accuracy.

3. Quadrature rule based on exact integration of a piecewise linear temporal approximation of the paraxial solution between time steps.

In the parallel case, the domain is decomposed into overlapping subdomains as described in, e.g., [2]. Within each subdomain, the method is the same as above except that the boundary condition at a subdomain interface is replaced with

$$\partial_n u^{m+1} - i\kappa u^{m+1} = \partial_n \widetilde{u}^m - i\kappa \widetilde{u}^m, \qquad (7)$$

where $u^m$ is the solution at time $t_m$ and the tilde denotes the solution in the neighboring subdomain.

## 4 Numerical results

In one of our examples, we demonstrate the solver by simulating plane wave transmission through a Luneburg lens. The domain is $D = [-3,3] \times [-3,3] \times [-3,8]$ and the wavelength is $\lambda = 0.03$. For the entire domain we use a $8,180 \times 8,180 \times 14,980$ grid divided into $6 \times 6 \times 14$ subdomains, for a total of 1.002 trillion points over 504 subdomains. The refraction coefficient is given by

$$m(\boldsymbol{x}) = \begin{cases} 2 - |\boldsymbol{x}|^2, & |\boldsymbol{x}| \leq 1 \\ 1, & |\boldsymbol{x}| > 1. \end{cases} \qquad (8)$$

For the computation we use 504 nodes with 192 cores per node, for a total of 96,768 cores. In this example we use second-order finite differences for the spatial discretization.

Fig. 1 shows the solution with the expected focusing of the plane wave at the edge of the lens, with an amplitude at the focus of $|v| = 22.8$. The $L^\infty$ residual error for this solution is $5.8 \cdot 10^{-3}$.

## Acknowledgements

Approved for public release; distribution is unlimited. AFRL Public Affairs release approval number AFRL-2026-0372.

# Rendering disordered waveguides transparent to waves

Matthieu Davy[1,*], Isam Ben Soltane[1], Daniel Lenz[2], Stefan Rotter[2]
[1]Université de Rennes, CNRS, IETR - UMR 6164, F-35000 Rennes, France
[2]Institute for Theoretical Physics, Vienna University of Technology (TU Wien), A-1040 Vienna, Austria
*Email: matthieu.davy@univ-rennes.fr

**Abstract**

Wave transmission through disordered media is typically hindered by multiple scattering which strongly reduces transmission. Here, we demonstrate that adding a properly designed scattering medium to a given disordered structure can make it fully transmitting and transparent, preserving the incident modes at the output. This is confirmed numerically in 2D disordered waveguides with the inverse design of antireflection structures being carried out using the coupled dipole approximation.

## 1 Introduction

While wavefront shaping techniques can partially counteract the effect of disorder to obtain perfect transmission through random media, this feature is typically restricted to a few incident wavefronts [1]. In contrast, inverse design of disordered structures greatly broadens the range of achievable functionalities. By optimizing the internal disorder, one can simultaneously control multiple input and output channels. Volumetric metastructures exploit multiple scattering to outperform ordered systems [2].

We investigate transmission through a disordered multichannel waveguide, see Fig.1(a). We have previously demonstrated that such structurally unpatterned medium made of randomly distributed scatterers can be rendered perfectly transmitting for all incoming wavefronts by placing a suitably engineered complementary medium in front of it [3]. However, a fundamental question lies in the degree of control that can simultaneously be exerted on transmitted fields in order to render the medium transparent to incoming waves. Here, we determine and realize the condition of transparency in multimode waveguides using a model based on coupled dipole approximation.

## 2 Model

The scattering matrix $S$ connecting the transmitted and reflected fields, can be partitioned into reflection matrices (RM) $r$ and $r'$ and transmission matrices (TM) $t$ and $t^T$. Our aim is to engineer a scattering medium placed on the left so that the complete system becomes fully transmitting to any incident wave. In this representation, $r$ and $t$ can be expressed in terms of the RM and TM of the two subsystems [3], $r = r_1 + t_1 r_2 (I - r_1' r_2)^{-1} t_1^T$ and $t = t_1 (I - r_1' r_2)^{-1} t_2$.

The condition of vanishing reflection for *all* incident wavefronts, $r = 0$, which extends the notion of critical coupling from a single input channel to the $N$ incoming channels, leads to the relation $r_1' = r_2^*$ [3]. Satisfying this condition enables perfect transmission but without any control on the outgoing wavefront, giving a translucent medium. The condition for transparency involves the joint design of the complete TM $t$. For a desired unitary TM, the condition $r_1' = r_2^*$ together with $t_1 = t_{\text{target}} t_2^\dagger$ finally gives $t = t_{\text{target}}$. For a diagonal TM, incident wavefront are perfectly transmitted and the outgoing wavefront is identical to the incoming wavefront, fullfilling the condition of perfect transparency. While no theoretical guarantee ensures the existence of such a solution, a sufficiently large number of degrees of freedom should enable its implementation.

## 3 Methods and Results

Direct optimization of large disordered waveguides using full-wave solvers is computationally intensive. We employ the coupled dipole approximation (CDA) [4], which enables an efficient computation of $S$ in the basis of flux-normalized waveguide modes. The $N_s$ scatterers, located at $\mathbf{r}_k$, are described by polarizabilities $\alpha_k$ inducing dipole moments $\mathbf{p}_m$ under excitation by an incoming mode $m$. The scattering matrix can be expressed in terms of the vector of Green's functions between modes and points $G_0$

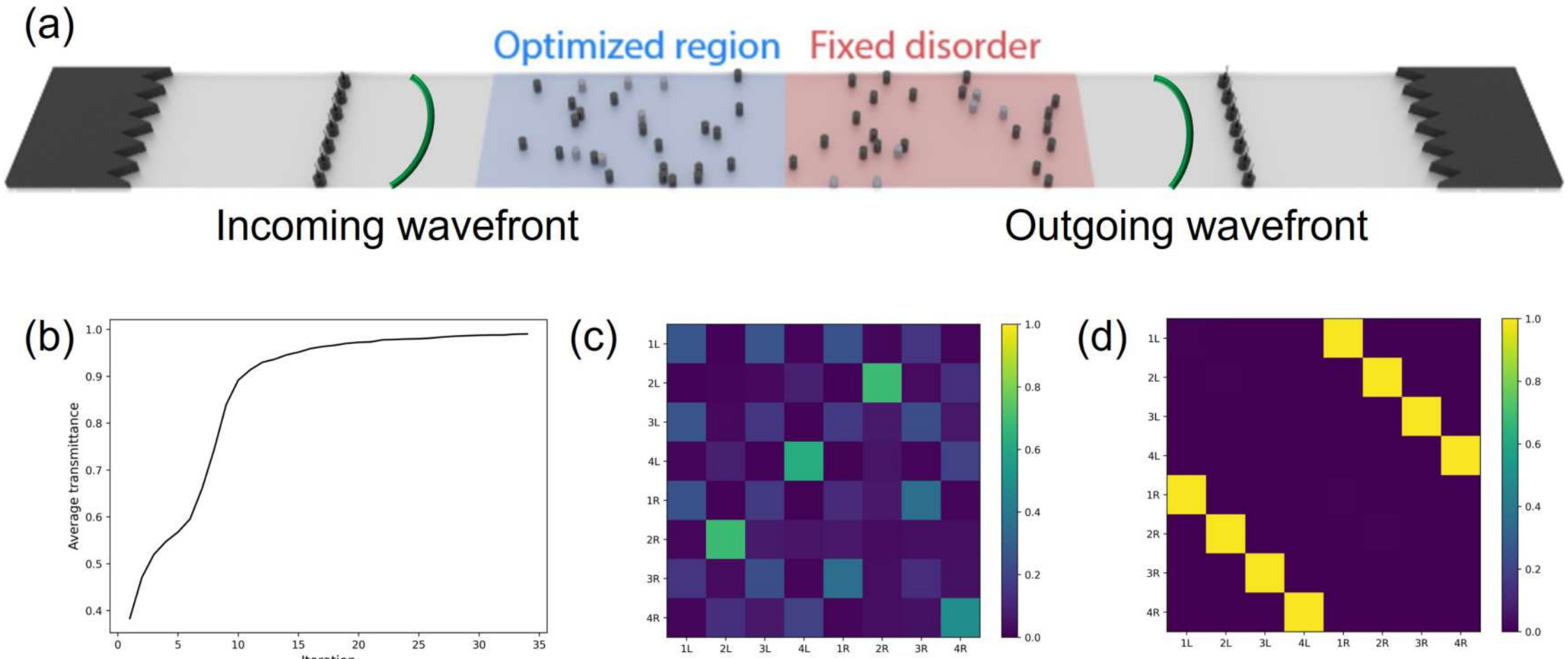

Figure 1: (a) Experimental configuration: a 2D multimode waveguide containing a fixed disordered region, with a complementary optimized medium placed on its left. (b) Average diagonal transmission, defined as $1/N\left(\sum_{n=1}^{N}|t_{nn}|^2\right)$, plotted versus the optimization iterations. A value of unity indicates perfect transmission and a diagonal transmission matrix (TM). (c,d) Intensity of the scattering matrix, $|S|^2$, for the full structure before (c) and after (d) optimization.

and the interaction matrix $W$ with the compact expression [5]

$$S = S_0 + G_0^T W^{-1} G_0. \tag{1}$$

$W$ contains the dipole–dipole Green's functions matrix $G^{dd}$ with $W = \mathrm{diag}(\alpha^{-1}) - G^{dd}$. Due to their large cross-section, metallic cylinders are represented by a few dipoles distributed along its boundary. The values of $\alpha$ are determined through a calibration procedure, by matching the results from the CDA model with full-wave simulations using COMSOL Multiphysics.

We start with a fixed disorder made of 3 metallic and 9 teflon cylinders giving an average transmission of 0.52. We optimize the complementary medium (3 metallic and 50 teflon cylinders) by maximizing the cost function : $f = \mathrm{diag}(|\mathrm{t}|^2)$. At the end of the optimization, the minimal diagonal element of $|t|^2$ is 0.96. Each mode is therefore almost perfectly transmitted with an outgoing wavefront being the same as the incident one. The TM becomes diagonal, see Fig. 1(b-d).

This work, which will be validated experimentally for the conference, opens promising avenues for the realization of antireflection structures with potential applications in wireless communications, filtering, energy harvesting, and imaging. With ongoing progress in computational methods and microfabrication technologies, the concept is expected to scale to systems supporting an ever larger number of modes.

# Spectral properties of a kinetic equation coupled with the linear wave system

**Bruno Després**[1,*]
[1]LJLL, Sorbonne University, Paris, FRANCE
*Email: bruno.despres@sorbonne-universite.fr

**Abstract**

We review the spectral properties of a linear wave equation coupled with a kinetic equation. Such models arise in Vlasov-Euler systems for the modeling of thick sprays encountered in engineering. Under natural conditions, the spectral properties are common to the ones of linearized Landau damping mechanism. The dispersion relations are indeed similar even if the physics is different. This is based on works which started during the PHD thesis of Victor Fournet.

## 1 Introduction

Models for thick sprays [1, 2] are encountered for when a compressible gaseous phase is coupled with particles which may be incompressible in first approximation. A simple example is the non dimensional barotropic compressible Euler equations for the fluid part coupled with a kinetic equation for the particles

$$\left\{ \begin{array}{l} \partial_t(\alpha\rho) + \nabla_x \cdot (\alpha\rho u) = 0, \\ \partial_t(\alpha\rho u) + \nabla_x \cdot (\alpha\rho u \otimes u) + \alpha \nabla_x p(\rho) = 0. \\ \partial_t f + v \cdot \nabla_x f - \nabla_x p(\rho) \cdot \nabla_v f = 0. \end{array} \right. \tag{1}$$

where the two phases are coupled by $\alpha(t,x) = 1 - \kappa \int_{\mathbb{R}^3} f(t,x,v)\, dv$ and by the fact that the particles react to the pressure gradient $\nabla_x p(\rho)$ of the fluid. An additional term, which is important for physical completeness, is the friction between the particles and the gas. This term can be written as $d(v-u)$. This term is discarded here for convenience, that is we take $d = 0$. Many other terms can be written to obtain a more complete model.

## 2 Linearization

It has been proposed in [3] to study such systems using linearization. One starts choosing a meaningful stationary solution. In the simplest case we take

$$\rho_0 > 0,\ u_0 = 0,\ \alpha_0 > 0,\ f_0(v) = e^{-v^2/2}.$$

The Gaussian distribution for particles represents a concentration of particles at low velocity. One rewrites the equation using the specific volume $\tau = 1/\rho$ instead of the density $\rho$. After linearization and convenient simplifications, one gets the integro-differential system in one space dimension and one velocity dimension

$$\left\{ \begin{array}{l} \partial_t \tau = \partial_x u + \partial_x \int_{\mathbb{R}} vg\sqrt{f_0(v)}dv, \\ \partial_t u = \partial_x \tau, \\ \partial_t f + v\partial_x f = \sqrt{f_0(v)}v\partial_x \tau. \end{array} \right. \tag{2}$$

All variables are linearized variables. This system couples a linear wave equation with an integral part. The last equation is a transport equation in the velocity direction.

Well posedness in the quadratic sense of such systems has been shown in [3]. The quadratic well posedness is generic for such systems. In our case it writes

$$\frac{d}{dt}\left( \int (\tau^2 + u^2)dx + \int\int f^2 dxdv \right) = 0.$$

## 3 Dispersion relation

To analyse the dynamics we use the Fourier-Laplace method in the torus $x \in \mathbb{T}$. One gets [4]

$$\widetilde{\tau}(\omega,k) = \int_0^{+\infty} \int_{\mathbb{T}} e^{i\omega t} e^{-ikx} \tau(t,x) dxdt$$

where $k \in \mathbb{Z}$ and $\Im\omega > 0$. It yields the system

$$\left\{ \begin{array}{l} -i\omega\widetilde{\tau} = ik\widetilde{u} + ik\int_{\mathbb{R}} v\sqrt{f_0(v)}\widetilde{f}(v)dv + a(k), \\ -i\omega\widetilde{u} = ik\widetilde{\tau} + b(k), \\ (-i\omega + ikv)\widetilde{f}(v) = ikv\sqrt{f_0(v)}\widetilde{\tau} + c(v,k). \end{array} \right.$$

The terms $\widetilde{\tau}$, $\widetilde{u}$ and $\widetilde{f}(v)$ depend on $\omega$ and $k$. The terms $a$, $b$ and $c$ account for the initial conditions of the dynamical system.

The determinant of the linear system with dimensionless coefficients yields the dispersion $\mathcal{D}(\omega,k) = 0$ where

$$\mathcal{D}(\omega,k) = 1 - \frac{k^2}{\omega^2} - \int \frac{f_0'(v)}{v - \omega/k} dv.$$

This dispersion relation is different but very close to the dispersion relation of the linear Landau damping phenomenon [5] which is $\int \frac{f_0'(v)}{v-\omega/k}dv = k^2$. The physics is different, but the dispersion relations are very similar.

To reconstruct non stationary solutions of the dynamical equations, one needs to discuss the complex roots of the dispersion on the lower part of the complex plane $\Im\omega < 0$. This is performed with analytic continuation of the dispersion relation. For the Gaussian profile $f_0(v) = e^{-v^2/2}$, one gets

$$\mathcal{D}(\omega, k) = 2 - \frac{k^2}{\omega^2} - C\frac{\omega}{k}e^{-\frac{\omega^2}{2k^2}}\left(i - \int_0^{\frac{\omega}{\sqrt{2}k}} e^{v^2}dv\right)$$

where the constant is $C = \sqrt{\pi/2}$. The integral is along any path in the complex plane.

The principal root is calculated with numerical methods. The numerical value is

$$\omega \approx 0.31 - 0.098i. \tag{3}$$

This can checked by means of comparisons with the numerical solution of either the non linear system (1) or the linear system (2). One obtains the results in Figures 1 and 2 where the accordance is remarkable between the linear theory and the non linear numerical simulations.

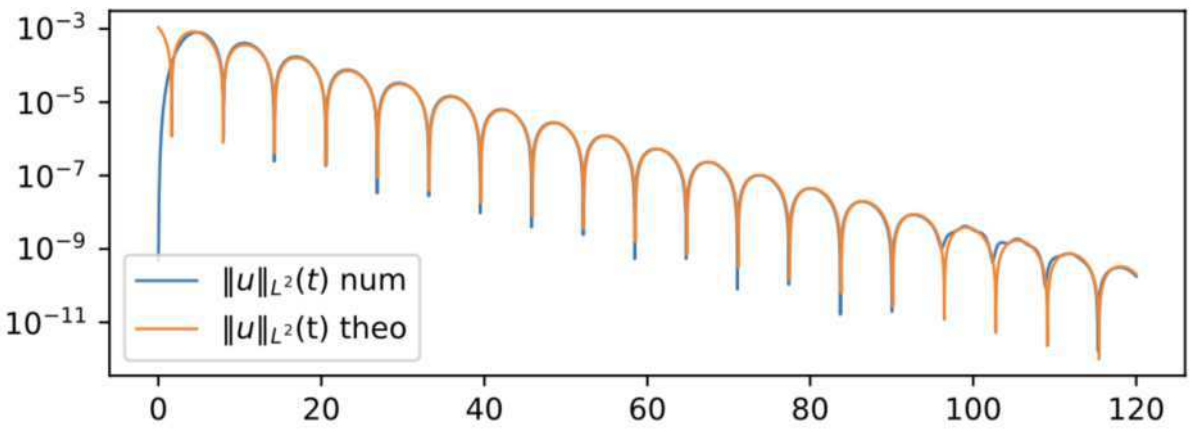

Figure 1: Comparison of the damping rate between the numerical solution $\|u(t)\|$ and the exact value (3).

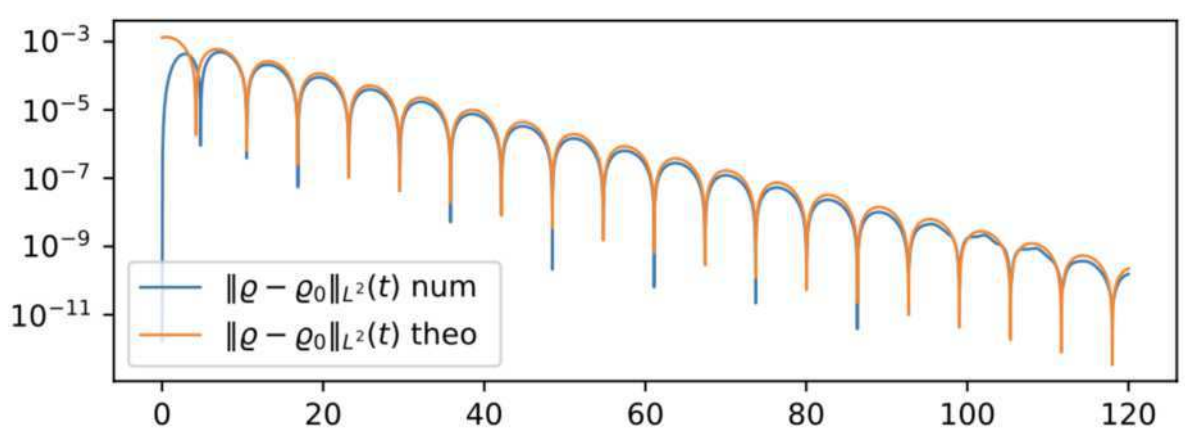

Figure 2: Comparison of the damping rate between the numerical solution $\|\rho(t) - \rho_0\|$ and the exact value (3).

One arrives at an interesting physical principle. *The interaction of a linear acoustic wave and with a cloud of particles (in the thick spray regime, initially close to a Maxwellian distribution) yields a damping phenomenon analogous to the linear Landau damping phenomenon.* More precisely the acoustic energy is damped and transferred to the particles.

For other particles profiles, the situation can be radically different. Linear instabilities are possible and the acoustic energy increases in the linear regime. A recent investigation of the coupling of particles with a more general system of conservation laws yields a linear model which couples particles with a Friedrichs system. A necessary condition for stability is obtained in [6]. The case with non zero friction $d > 0$ is also different as shown by the numerical simulations.

Finally let us mention that the mathematical well posedness of the non linear version of these coupled systems like (1) is a difficult open problem for which only few results exists.

# Transparent boundary conditions in heterogeneous electromagnetic waveguides

**Anne-Sophie Bonnet-Ben Dhia**[1,*], **Lucas Chesnel**[2], **Sonia Fliss**[1], **Eric Lunéville**[1], **Aurélien Parigaux**[1]

[1]POEMS, CNRS, Inria, ENSTA, Institut Polytechnique de Paris, Palaiseau, France
[2]Inria, Unité de Math. Appliquées, ENSTA, Institut Polytechnique de Paris, Palaiseau, France
*Email: anne-sophie.bonnet-bendhia@ensta.fr

**Abstract**

We consider an unbounded domain, with a perfectly conducting boundary, which is the union of several semi-infinite waveguides. The objective is to compute the time-harmonic electromagnetic wavefield generated by some compactly supported source, by a method coupling finite elements in the junction zone and modal expansions in the waveguides.

This can be achieved in homogeneous waveguides [1] by using a modal expression of a so-called Electric-to-Magnetic condition. But the case of transversely heterogeneous waveguides is less obvious, due to the nonselfadjointness of the modal eigenvalue problem. Relying on Kondratiev theory [2], we derive, justify and validate here two types of so called transparent conditions *with overlap* for transversely heterogeneous waveguides.

## 1 Problem setting

To simplify, let us consider here the case of a domain $\Omega$ containing only one semi-infinite waveguide, whose cross-section is a bounded connected domain $S \subset \mathbb{R}^2$. More precisely, we suppose that

$$\Omega \cap \{z < 0\} = \Omega_0 \text{ is bounded}$$
$$\Omega \cap \{z > 0\} = S \times \{z > 0\}$$

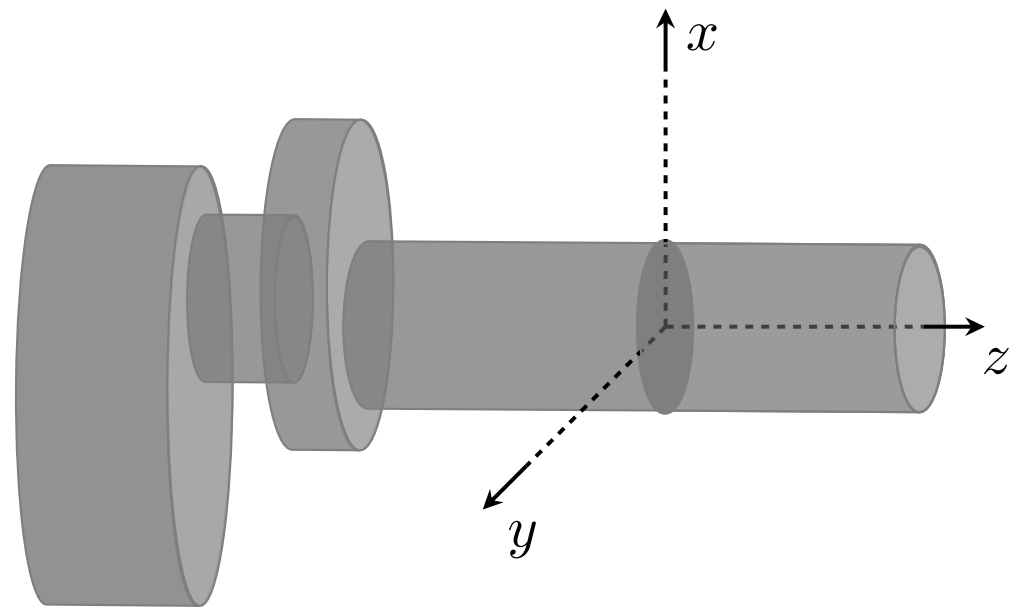

We also suppose that the dielectric permittivity $\varepsilon$ and the magnetic permeability $\mu$ are such that $\varepsilon = \varepsilon_\infty(x, y)$, $\mu = \mu_\infty(x, y)$ for $z > 0$.

Our objective is to compute the *outgoing* solution of the following problem:

$$\left|\begin{aligned} \mathbf{curl}\,\mu^{-1}\mathbf{curl}\,\mathbf{E} - \omega^2\varepsilon\mathbf{E} = i\omega\mathbf{J} \quad &(\Omega) \\ \mathbf{n}\times\mathbf{E} = 0 \qquad &(\partial\Omega) \end{aligned}\right. \tag{1}$$

where the source term $\mathbf{J}$ is supported in $\Omega_0$.

It is proved in [2] that problem (1) is well-posed if the domain $\Omega$ does not support trapped modes at the frequency $\omega$ (trapped modes are non trivial $\mathbf{L}^2$ solutions of (1) with $\mathbf{J} = 0$).

In order to solve this problem numerically, our objective is to derive an equivalent formulation in a bounded domain $\Omega_L := \Omega \cap \{z < L\}$ ($L > 0$), with a transparent boundary condition on the artificial boundary $S_L := S \times \{L\}$.

## 2 Electromagnetic modes

The mathematical definition of the term *outgoing* requires the introduction of the modes of the infinite waveguide $S \times \mathbb{R}$. More precisely, they are solutions of the same equations as (1) with $\mathbf{J} = 0$, replacing $\Omega$ by $S \times \mathbb{R}$, $\varepsilon$ by $\varepsilon_\infty$ and $\mu$ by $\mu_\infty$, that take the form

$$\mathbf{E}(x, y, z) = \mathcal{E}(x, y)e^{i\beta z} \text{ with } \beta \in \mathbb{C}.$$

It is known [2] that there is a finite number $N_p$ of propagating modes ($\beta \in \mathbb{R}$) and an infinite sequence of evanescent ones ($\beta \notin \mathbb{R}$). Assuming that $\omega$ is not a *cut-off frequency*, we denote by $(\beta_n, \mathbf{E}_n^+)_{n\geq 1}$ the family of right-going modes:

- if $n \leq N_p$, the mode is propagating and has a positive group velocity;
- if $n > N_p$, the mode is evanescent and $\mathrm{Im}(\beta_n) > 0$.

A solution $\mathbf{E}$ of (1) is said to be *outgoing* if there exist $a_1, a_2 \cdots a_{N_p} \in \mathbb{C}$ and $\alpha > 0$ such that

$$\mathbf{E} = \sum_{n=1}^{N_p} a_n \mathbf{E}_n^+ + \mathbf{E}_{ev} \text{ in } \Omega \cap \{z > 0\}$$

with $e^{\alpha z}\mathbf{E}_{ev} \in \mathbf{L}^2(\Omega)$.

## 3 Two elementary problems

The derivation of the transparent boundary conditions on $S_L$ will rely on the two following elementary problems.

**The half-guide Dirichlet problem.** Suppose that $\mathbf{n}\times\mathbf{E}$ is known on $S_0$. There exists a unique outgoing solution $\mathbf{E}^\infty$ of

$$\left|\begin{array}{ll}\mathbf{curl}\,\mu_\infty^{-1}\mathbf{curl}\,\mathbf{E}^\infty-\omega^2\varepsilon_\infty\mathbf{E}^\infty=0 & (S\times\mathbb{R}^+)\\ \mathbf{n}\times\mathbf{E}^\infty=0 & (\partial S\times\mathbb{R}^+)\\ \mathbf{n}\times\mathbf{E}^\infty=\mathbf{n}\times\mathbf{E} & (S_0)\end{array}\right.$$

Moreover, as a consequence of the bi-orthogonality between the modes, the coefficients $a_n^\infty$ in the identity $\mathbf{E}^\infty=\sum_{n=1}^{N_p}a_n^\infty\mathbf{E}_n^+ + \mathbf{E}_{ev}^\infty$ have an explicit expression, in function of $\mathbf{n}\times\mathbf{E}|_{S_0}$:

$$a_n^\infty=\int_{S_0}(i\omega\mu_\infty)^{-1}(\mathbf{n}\times\mathbf{E})\cdot\mathbf{curl}\,\mathbf{E}_n^+ \qquad (2)$$

**The full-guide transmission problem.** Suppose now that $\mathbf{n}\times\mathbf{E}$ and $\mathbf{n}\times\mu_\infty^{-1}\mathbf{curl}\,\mathbf{E}$ are known on $S_0=\Omega\cap\{z=0\}$. There exists a unique outgoing solution $\mathbf{E}^\infty$ of

$$\left|\begin{array}{ll}\mathbf{curl}\,\mu_\infty^{-1}\mathbf{curl}\,\mathbf{E}^\infty-\omega^2\varepsilon_\infty\mathbf{E}^\infty=0 & (S\times\mathbb{R})\\ \mathbf{n}\times\mathbf{E}^\infty=0 & (\partial S\times\mathbb{R})\\ [\mathbf{n}\times\mathbf{E}^\infty]=\mathbf{n}\times\mathbf{E} & (S_0)\\ [\mathbf{n}\times\mu_\infty^{-1}\mathbf{curl}\,\mathbf{E}^\infty]=\mathbf{n}\times\mu_\infty^{-1}\mathbf{curl}\,\mathbf{E} & (S_0)\end{array}\right.$$

where the brackets denote the jumps of the quantities through $S_0$. Again, the coefficients $a_n^\infty$ have an explicit expression, as a function of the data on $S_0$:

$$\begin{aligned}a_n^\infty &= \frac{1}{2i\omega}\left(\int_{S_0}\mu_\infty^{-1}\mathbf{n}\times\mathbf{E}\cdot\mathbf{curl}\,\mathbf{E}_n^+\right.\\ &\quad\left.-\int_{S_0}\mu_\infty^{-1}\mathbf{curl}\,\mathbf{E}\cdot\mathbf{n}\times\mathbf{E}_n^+\right)\end{aligned} \qquad (3)$$

## 4 EtM and CtM conditions

The idea is now to consider the problem

$$\left|\begin{array}{ll}\mathbf{curl}\,\mu^{-1}\mathbf{curl}\,\mathbf{E}-\omega^2\varepsilon\mathbf{E}=i\omega\mathbf{J} & (\Omega_L)\\ \mathbf{n}\times\mathbf{E}=0 & (\partial\Omega_L\setminus S_L)\\ \mathbf{n}\times\mu_\infty^{-1}\mathbf{curl}\,\mathbf{E}=\mathbf{n}\times\mu_\infty^{-1}\mathbf{curl}\,\mathbf{E}^\infty & (S_L)\end{array}\right.$$

where $\mathbf{E}^\infty$ is the solution of one of the two elementary problems above.

If $\mathbf{E}^\infty$ is the solution of the half-guide Dirichlet problem, we say that the boundary condition on $S_L$ is an Electric-to-Magnetic (EtM) condition with overlap, because it relates the transverse electric field on $S_0$ to the transverse magnetic field on $S_L$ ($L$ is the overlap). Indeed, up to the multplication by $i\omega$, $\mathbf{n}\times\mu_\infty^{-1}\mathbf{curl}\,\mathbf{E}^\infty$ is equal to the transverse magnetic field $\mathbf{n}\times\mathbf{H}^\infty$.

If $\mathbf{E}^\infty$ is the solution of the full-guide transmission problem, the boundary condition on $S_L$ is called a Current-to-Magnetic (CtM) condition (also with overlap), because it relates both the transverse electric and magnetic fields on $S_0$, taken as electric and magnetic currents, to the transverse magnetic field on $S_L$.

**Theorem 1** *The problem set in $\Omega_L$ is equivalent to the initial problem* (1) *set in $\Omega$*
*(i) except for a discrete sequence of frequencies for the EtM condition,*
*(ii) always for the CtM condition.*

## 5 Numerical implementation

In practice, thanks to the overlap, the evanescent field $\mathbf{E}_{ev}^\infty$ is neglected, so that $\mathbf{E}^\infty$ is replaced by a finite sum of propagating modes. It is also possible to include some evanescent modes in the approximation, as soon as they are not bi-orthogonal to themselves.

Note that the CtM condition requires the knowledge of $\mathbf{n}\times\mu_\infty^{-1}\mathbf{curl}\,\mathbf{E}$ on the internal boundary $S_0$. This quantity is well-defined for a solution $\mathbf{E}$ of Maxwell equations, but not in general for a field $\mathbf{E}$ of the variationnal space

$$\{\mathbf{E}|\mathbf{E},\mathbf{curl}\,\mathbf{E}\in\mathbf{L}^2(\Omega_L),\mathbf{n}\times\mathbf{E}=0\text{ on }\partial\Omega_L\setminus S_L\}$$

The difficulty is overcome by integrating by part the formula (3) for $a_n^\infty$, which amounts to define $\mathbf{n}\times\mu_\infty^{-1}\mathbf{curl}\,\mathbf{E}$ in a weak sense.

Finally, the numerical results are obtained with the library XLiFE++ [3], using Nédéléc finite elements.

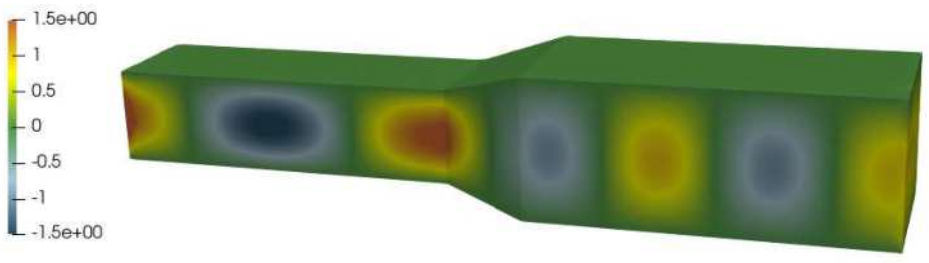

# A THEORY FOR THREE-DIMENSIONAL WAVE BREAKING BASED ON CREST INSTABILITY OF SHORT-CRESTED STOKES WAVES

Frédéric Dias[1,*], Ayoub Mansar[1], Matt Turner[2], Tom Bridges[2]
[1]UCD, Dublin, Ireland and ENS Paris-Saclay, Gif sur Yvette, France
[2]University of Surrey, Guildford, UK
*Email: frederic.dias@ucd.ie

**Abstract**

A theory for three-dimensional wave breaking based on crest instability of short-crested Stokes waves is presented.

## 1 Introduction

In the late 1970s and early 1980s, several papers dealing with instabilities of water waves were published. Well-known papers are the papers by Longuet-Higgins [3,4], McLean et al. [7], McLean [8]. McLean et al. [7] considered a two-dimensional (2D) Stokes wave. They perturbed it by adding three-dimensional (3D) perturbations. For sufficiently steep Stokes waves, it was shown that the dominant instability becomes a 3D perturbation of class II. This instability results from a resonant interaction between five components of the wave field and was first discovered by Longuet-Higgins [3] for purely 2D perturbations. This instability is quite different from the Benjamin-Feir (BF) instability that results from a quartet resonance, that is, a resonant interaction between four components of the wave field. According to the classification of McLean [8], the BF instability belongs to class I interactions that are predominantly 2D. Su et al. [12] and Su [11] performed a series of experiments on instabilities of gravity-wave trains of large steepness in deep water in a long tow tank and a wide basin. They observed 3D structures corresponding to the nonlinear evolution of the dominant class II instability discovered by McLean et al. [7]. The initial 2D wave train of large steepness evolves into a series of 3D spilling breakers, followed by a transition to a more or less 2D wave train. These 3D patterns take the form of crescent-shaped perturbations (horseshoe patterns) riding on the basic waves. It is interesting to note that Melville [10] performed experiments on the instability of 2D waves and saw the 3D instability of 2D waves. The 3D effects are described in an appendix! At the same time, there was an interest in finding steady 3D waves from the analytical point of view (bifurcations) and numerically (see for example [9]). The paper of Longuet-Higgins [3] coins the terminology "superharmonics". The distinguishing feature of that instability is that it is co-propagating. In other words, it has the same wavelength and the same speed as the Stokes wave. This feature is to be contrasted with subharmonic instabilities and the BF instability where the perturbation has a different wavelength and a different speed. The link between wave breaking and instabilities has been made several times. But which instabilities are the most relevant for wave breaking? In his book on wave breaking, Babanin [1] writes that instability is a key word in the breaking process. He also writes (I quote): "Some processes leading to wave breaking, and in particular the most important mechanism of modulational [BF] instability of nonlinear wave trains, may be impaired in three-dimensional natural field circumstances, which fact further depletes the generality of conclusions achieved in two-dimensional wave tanks." This is why we decided to focus on the superharmonic instability, which develops into crest instability. For the case of unstable periodic Stokes waves, the wave-breaking scenario was observed in the numerical simulations of Jillians [2]. Jillians took the eigenfunction as a small perturbation to the Stokes wave and integrated in time. A microbreaker emerged from the superposition of a periodic travelling wave and a superharmonic unstable eigenfunction. The same numerical behaviour was found by Longuet-Higgins and Dommermuth [5]. Mansar et al. [6] checked the robustness of the appearance of this crest instability leading to breaking. They added a perturbation to a large-amplitude unstable Stokes wave, which was then taken as initial data in a direct numerical solution of the Navier-Stokes equations, using the Basilisk numerical software package. Among the effects that were studied were dissipation, approximation error, changes in depth, non-zero air density. The conclusion is that the emergence of breaking waves is quite robust.

## 2 Theory

The next step is to tackle both the 3D instability of 2D waves as well as the instability of 3D waves. We present a theory for breaking of 3D water waves, based on instability of short-crested Stokes waves travelling in deep water. We find that these waves are susceptible to a crest instability at large amplitude, with a dipole structure of the eigenfunctions that varies periodically along the crest. Studying the full 3D problem numerically, with the unstable eigenfunctions as initial data, leads to overturning crests that generate the structure of a micro-breaker, similar to what has been observed in the open ocean. The numerical simulations are tested with the addition of various noise factors, and the theoretical mechanism is found to be robust.

# Decoupling Wave-Heat-Type Problems

Julian Dörner[1,*], Benjamin Dörich[1]
[1]Karlsruhe Institute for Technology
*Email: julian.doerner@kit.edu

**Abstract**

This talk concerns a broad class of time-dependent, coupled wave–heat systems. We propose a decoupled time-stepping scheme where the wave part is advanced explicitly and the heat part implicitly, removing the parabolic CFL restriction and permitting much larger time steps. Under suitable assumptions, we prove stability and derive fully discrete error estimates in a unified framework.

## 1 Introduction

We consider wave problems which are coupled to dissipative heat problems and can be written formally as

$$\partial_t x_1 - S_2 x_2 = -B\theta + f, \tag{1a}$$
$$\partial_t x_2 + S_1 x_1 = 0, \tag{1b}$$
$$\partial_t \theta + L\theta = B^* x_1 + g. \tag{1c}$$

Here $x = (x_1, x_2) \in X$ denotes the wave-variable with adjoint wave-operators $S_1^* = S_2$, and $\theta \in H$ denotes the heat-variable with a positive heat-operator $L$, together with a coupling operator $B$, cast in a dense Gelfand setting

$$Y \overset{d}{\hookrightarrow} X \hookrightarrow Y^*, \quad V \overset{d}{\hookrightarrow} H \hookrightarrow V^*.$$

The most descriptive example fitting in this framework is the thermo-elastic string in $\Omega = [0,1]$ governed by the system

$$\partial_t v - \partial_x^2 u = -\partial_x \theta + f, \tag{2a}$$
$$\partial_t u - v = 0, \tag{2b}$$
$$\partial_t \theta - \partial_x^2 \theta = -\partial_x v + g, \tag{2c}$$

with $x = (v, u)$, and the spaces

$$X = L^2(\Omega) \times H_0^1(\Omega), \qquad Y = H_0^1(\Omega) \times H_0^1(\Omega),$$
$$H = L^2(\Omega), \qquad V = H_0^1(\Omega).$$

Using an explicit time discretization of (2), the heat-part (2c) is subject to a more restrictive CFL condition than the wave-part. Hence, the problem benefits from an implicit treatment of this part.

A more complicated example fitting the assumptions is given by the non-local in time visco-acoustic problem on $\Omega \subset \mathbb{R}^d$ with

$$\rho \partial_t \mathbf{v} - \nabla p_0 = \sum_{i=1}^{G} \nabla p_i + f, \tag{3a}$$
$$\kappa^{-1} \partial_t p_0 - \operatorname{div} \mathbf{v} = 0, \tag{3b}$$
$$\partial_t p_i + \tau_i^{-1} p_i = \kappa \tau_P \operatorname{div} v, \tag{3c}$$

where $x = (\mathbf{v}, p_0)$, $\theta = (p_1, \dots)$, and spaces

$$X = L^2(\Omega)^{d+1}, \qquad Y = H_0(\operatorname{div}, \Omega) \times H^1(\Omega),$$
$$H = L^2(\Omega)^G, \qquad V = H^1(\Omega)^G.$$

Here, the heat-part (3c) stems from a discretization of a non-local in time convolutional operator. Small relaxation times $\tau_i \ll 1$ result in a stiff system of ODEs, which then leads to a severe CFL. Thus, the problem also benefits if we treat this part implicitly.

More models fit into our framework such as Maxwell's equation with non-local in time material and thermo-elastic problems generalizing (2) in higher dimensions.

## 2 Space discretization

In the spirit of [1], we consider a fairly general finite element discretization for the wave-part (possibly non-conforming), and a conforming discretization for the heat part, i.e. there exists discrete spaces

$$Y_h \subset X, \qquad V_h \subset V.$$

Those spaces induce the discrete pivot spaces $X_h$ and $H_h$ with the continuous inner products. Our abstract space discretization of (1) reads

$$\partial_t x_{h1} - S_{h2} x_{h2} + S_{h1}^s x_{h1} = -B_h \theta_h + f_h, \tag{4a}$$
$$\partial_t x_{h2} + S_{h1} x_{h1} + S_{h2}^s x_{h2} = 0, \tag{4b}$$
$$\partial_t \theta_h + L_h \theta_h = B_h^* x_{h1} + g_h, \tag{4c}$$

where we emphasize the possibility of positive stabilization operators $S_{h1}^s$, $S_{h2}^s$, commonly used in discontinuous Galerkin methods. Structurally, we assume that the discrete wave-operators are adjoint and the heat-operator is positive.

## 3 Time stepping

We discretize the wave-part (4a) and (4b) with the leapfrog method, and the heat-part with the Crank-Nicolson method in time. The resulting scheme is given in (5) below. This idea was successfully applied to the Stokes-Darcy problem, see e.g. [2].

As mentioned in the introduction above, we only obtain a CFL condition from the less severe wave-part, i.e. we demand for $\eta \in (0,1)$ that

$$\frac{\tau}{2}\|S_h x_h\|^2 + \frac{\tau}{2}\|B_h\theta_h\|^2 + \tau\|\sqrt{S_h^s}x_h\|^2 \le \eta\|(x_h,\theta_h)\|. \tag{CFL}$$

In our examples the first two norms in (CFL) scale inverse with the minimal mesh size $h_{\text{min}}$, and the third with $h_{\text{min}}^{1/2}$. Hence, we consider a step-size restriction $\tau \approx h_{\text{min}}$.

$$x_{h1}^{n+1/2} = x_{h1}^n + \frac{\tau}{2}S_{h2}x_{h2}^n - \frac{\tau}{2}S_{h1}^s x_{h1}^n - B_h\theta_h^n + \frac{\tau}{2}f_h^n, \tag{5a}$$

$$x_{h2}^{n+1} = x_{h2}^n - \tau S_{h1}x_{h1}^{n+1/2} - \tau S_{h2}^s x_{h1}^n, \tag{5b}$$

$$\theta_h^{n+1} = \theta_h^n - \frac{\tau}{2}L_h(\theta_h^{n+1}+\theta_h^n) + \tau B_h^* x_{h1} + \frac{\tau}{2}(g_h^{n+1}+g_h^n), \tag{5c}$$

$$x_{h2}^{n+1} = x_{h2}^{n+1/2} + \frac{\tau}{2}S_{h1}x_{h1}^{n+1} - \frac{\tau}{2}S_{h2}^s x_{h2}^n - \frac{\tau}{2}B_h\theta_h^{n+1} + \frac{\tau}{2}f_h^{n+1}. \tag{5d}$$

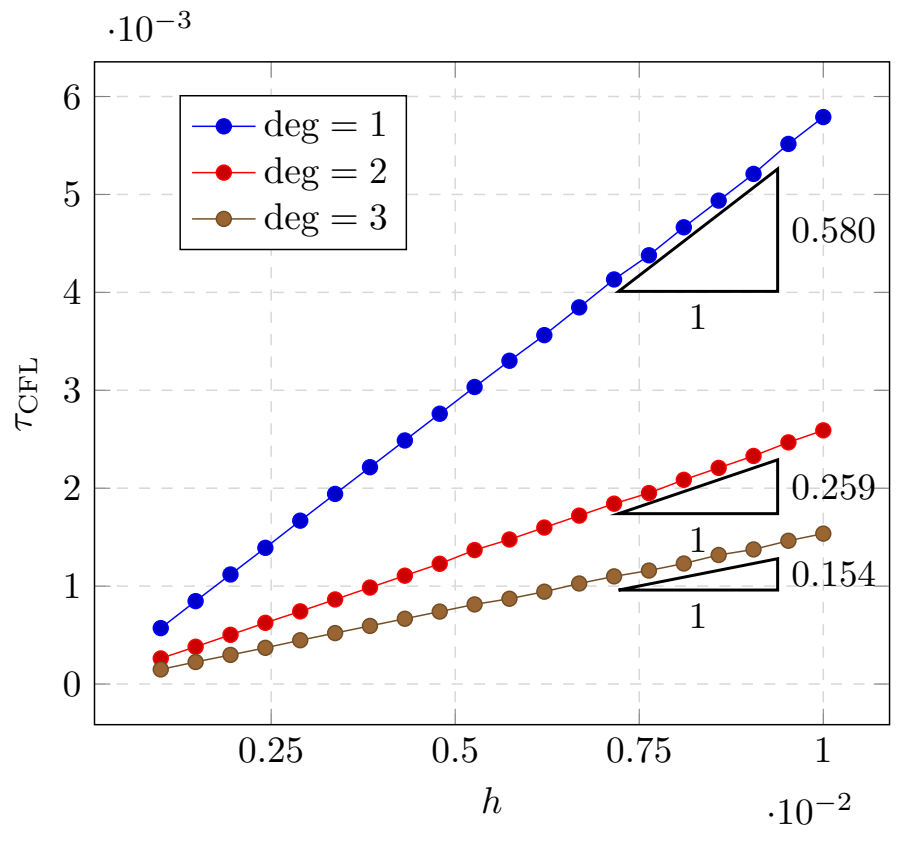

Figure 1: CFL threshold

## 4 Results

Our proposed scheme is stable with constants growing linear in the end-time $T$.

**Lemma 1** *If* (CFL) *holds, then*

$$(1-\eta)\|(x_h^n,\theta_h^n)\|_{X_h\times H_h}^2 \lesssim \|(x_h^0,\theta_h^0)\|_{X_h\times H_h}^2 + \frac{T}{1-\eta}\tau\sum_{\ell\ge 0}^{n}\|(f_h^\ell,g_h^\ell)\|_{X_h\times H_h}^2.$$

Under the assumption that the solution is four times differentiable with respect to time, we show an abstract full-error estimate

$$\|(u,\theta)(t_n) - (u_h^n,\theta_h^n)\|_{X_h\times H_h} \lesssim \tau^2 + \Delta_{\text{space}},$$

where $\Delta_{\text{space}}$ collects defects introduced by the abstract space discretization.

In the talk, we illustrate our theoretical findings with several examples mentioned in the introduction. We present space and time discretization plots that validate (CFL). Figure 1 shows the threshold step size $\tau_{\text{CFL}>0}$ for which the proposed method is stable, plotted against the mesh size, for problem (2) with a conforming space discretization. The figure clearly confirms the linear relation discussed above.

# Comparison of numerical coupling strategies in time-harmonic elasto-acoustic simulations

**Pierre Dubois**[1,2,*], **Hélène Barucq**[2], **Florian Faucher**[2], **Isabelle Ramière**[1], **Raphaël Prat**[1]

[1]CEA/DES/IRESNE/DEC/SESC, Cadarache, F-13108 Saint-Paul-Lez Durance, France

[2]Team Makutu, Inria, University of Pau and Pays de l'Adour, CNRS, France

*Email: pierre.dubois@cea.fr

**Abstract**

This work addresses numerical methods for time-harmonic elasto-acoustic wave problems. We aim at comparing monolithic and partitioned solution strategies, regarding computational cost and scalability. The comparison focuses on the impact of computational domain size and frequency range on solver performance. We also investigate the robustness of partitioned algorithms and possible non-convergence issues.

## 1 Introduction

In many applications, accurately modeling the interaction between acoustic waves propagating in a fluid and elastic waves evolving in a solid is crucial. A representative example is the detection and localization of an object submerged deep in a water body. These phenomena are strongly coupled and their numerical simulation can be computationally expensive, making direct approaches based on a single solver for the fully coupled system difficult to implement efficiently. This so-called monolithic approach is often contrasted with partitioned strategies, in which each physical subsystem is solved independently, while the coupling is enforced through the iterative exchange of appropriate interface quantities until convergence is achieved. Despite their potential advantages, partitioned algorithms have received limited attention in the context of elasto-acoustic wave problems. It is worth noting that recent studies based on domain decomposition methods in the time domain have introduced coupling conditions that ensure convergence and accuracy [1]. In this work, we focus on time-harmonic elasto-acoustic problems and we propose to assess the performance of a monolithic solver based upon a Hybridizable Discontinuous Galerkin (HDG) formulation with a partitioned algorithm using continuous Lagrange finite elements for both subdomains.

## 2 Time-harmonic elasto-acoustic model

Let $\Omega = \Omega_f \cup \Omega_s$ be the computational domain in $\mathbb{R}^d$ $(d = 2, 3)$, each subdomain having a regular boundary $\partial\Omega_f$ and $\partial\Omega_s$ respectively. Both subdomains share a common interface $\Gamma = \partial\Omega_f \cap \partial\Omega_s$. We denote by $\mathbf{n}_f$ and $\mathbf{n}_s$ the outward unit normal vectors to $\Omega_f$ and $\Omega_s$, respectively.

The time-harmonic elasto-acoustic coupling problem reads:

$$\begin{cases} \Delta p + k^2 p = f_f, & \text{in } \Omega_f, \\ \nabla \cdot \boldsymbol{\sigma}(\mathbf{u}) + \omega^2 \rho_s \mathbf{u} = \mathbf{f_s}, & \text{in } \Omega_s, \\ \dfrac{\partial p}{\partial \mathbf{n}_f} = -\omega^2 \rho_f (\mathbf{u} \cdot \mathbf{n}_f), & \text{on } \Gamma, \\ \boldsymbol{\sigma}(\mathbf{u}) \cdot \mathbf{n}_s = -p\, \mathbf{n}_s, & \text{on } \Gamma, \end{cases} \tag{1}$$

with boundary conditions omitted for conciseness. Here, $p$ is the acoustic pressure, $\mathbf{u}$ the elastic displacement, $k = \omega/c_f$ the acoustic wavenumber (with $c_f$ the sound speed in the fluid), $\omega$ the angular frequency, and $\rho_f$, $\rho_s$ the fluid and solid densities, $f_f$ and $\mathbf{f_s}$ are the source terms in the fluid and solid domains. For the Cauchy stress tensor we use the linear elastic isotropic constitutive law:

$$\boldsymbol{\sigma}(\mathbf{u}) = \lambda_s (\nabla \cdot \mathbf{u})\mathbf{I} + 2\mu_s \boldsymbol{\varepsilon}(\mathbf{u}), \tag{2}$$

with Lamé parameters $\lambda_s$, $\mu_s > 0$ and strain tensor $\boldsymbol{\varepsilon}(\mathbf{u}) = \frac{1}{2}(\nabla\mathbf{u} + (\nabla\mathbf{u})^T)$. The interface conditions enforce continuity of normal displacement and balance of normal stresses.

## 3 Coupling strategies

After discretization, the coupled problem reduces to the block linear system:

$$\begin{pmatrix} \mathbf{A}_f & \mathbf{C}_f \\ \mathbf{C}_s & \mathbf{A}_s \end{pmatrix} \begin{pmatrix} \mathfrak{p} \\ \mathfrak{u} \end{pmatrix} = \begin{pmatrix} \mathfrak{f}_\mathfrak{f} \\ \mathfrak{f}_\mathfrak{s} \end{pmatrix}, \tag{3}$$

where $\mathfrak{p}$ and $\mathfrak{u}$ are the discrete pressure and displacement unknowns, $\mathfrak{f}_\mathfrak{f}$ and $\mathfrak{f}_\mathfrak{s}$ are the discrete source term, $\mathbf{A}_f$ and $\mathbf{A}_s$ represent the discrete acoustic and elastic operators, while $\mathbf{C}_f$ and $\mathbf{C}_s$ encode the interface coupling conditions.

There are two main approaches to solve (3): monolithic and partitioned. The monolithic approach solves the entire system simultaneously and yields a robust and accurate solution when convergence is achieved. However, its scalability relies heavily on the performance of the underlying solver, which can be limited for large-scale problems.

The partitioned approach recasts (3) as a fixed-point problem. Starting from an initial guess, the scheme iterates:

$$\begin{cases} \mathbf{A}_f \mathfrak{p}^{n+1} = \mathfrak{f}_\mathfrak{f} - \mathbf{C}_f \mathfrak{u}^n, \\ \mathbf{A}_s \mathfrak{u}^{n+1} = \mathfrak{f}_\mathfrak{s} - \mathbf{C}_s \mathfrak{p}^{n+1}, \end{cases} \tag{4}$$

until a suitable convergence criterion is satisfied. Each iteration of this block Gauss-Seidel scheme solves two smaller systems instead of one large coupled system, reducing memory and computational cost. However, convergence is not guaranteed for all parameter regimes and depends on the problem characteristics, particularly the frequency $\omega$ and the physical properties of the medium [2].

## 4 Numerical experiments

To compare the computational efficiency of the two approaches, our first test considers a three-dimensional domain of size 1km$^3$, with an acoustic medium on the upper half and an elastic medium in the lower half, as illustrated in Fig. 1. For the implementation of the monolithic solver, a HDG method is employed, such that the system (3) is reformulated in terms of numerical traces; see, e.g., [3] All numerical experiments are performed using the software `Hawen`.

A simulation with the full solver is shown in Fig. 1 for the two-medium configuration at frequency $\omega/(2\pi)$ =15Hz where $\lambda_f = c_f/f = 100\,\mathrm{m}$, $\lambda_P = c_P/f = 166\,\mathrm{m}$, and $\lambda_S = c_S/f = 80\,\mathrm{m}$ denote the acoustic and elastic wavelengths. We employ a polynomial approximation of order 4, and the resulting linear system is solved using the direct solver MUMPS. Computational costs are given in Table 4. We will explore the cost difference with the partitioned strategy, for larger test-cases and with multiple right-hand sides.

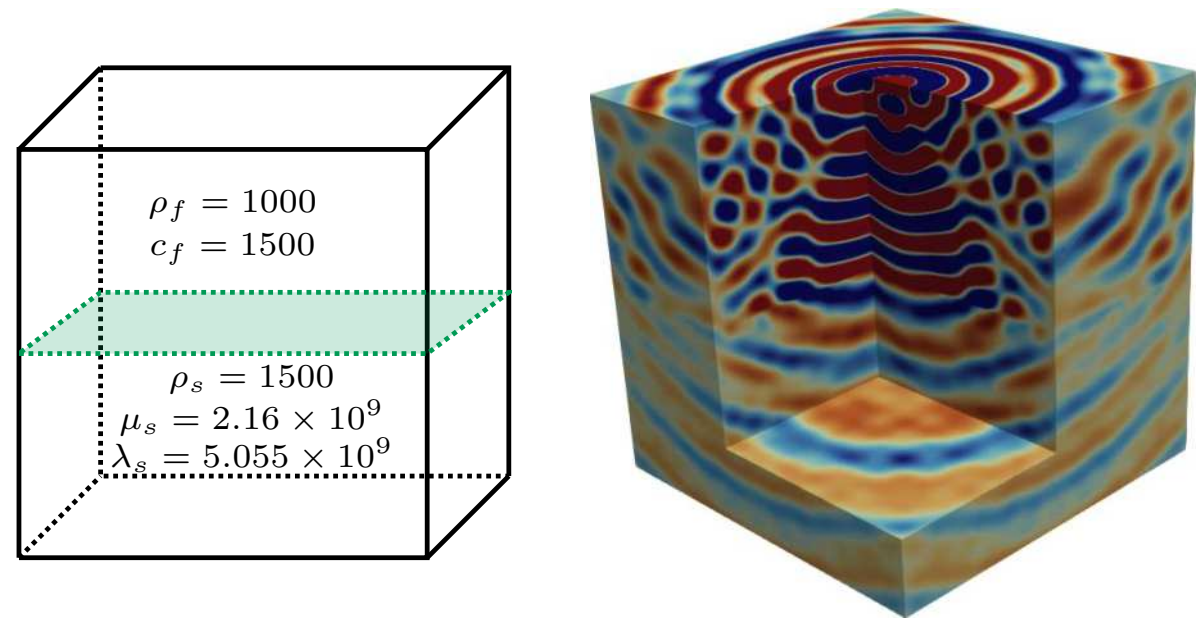

Figure 1: 3D problem of size 1km$^3$ with acoustic medium on the upper half and elastic medium on the lower half. The simulation at frequency 15Hz uses polynomial order 4 and a Gaussian in space source. We show $\mathrm{Re}(u_z)$.

| **system** | **size** | **memory** | **time** |
|---|---|---|---|
| coupled | $6.3 \cdot 10^6$ | 320GB | 9min |
| acoustic (half domain) | $1.1 \cdot 10^6$ | 15GB | 10sec |
| elastic (half domain) | $5.2 \cdot 10^6$ | 260GB | 6min |

Table 1: Computational cost for the matrix factorization with the monolithic solver. With one medium, only half of the domain is considered (see configuration of Fig. 1). The computations are performed on the PlaFRIM cluster using 16 MPI processes with 12 threads.

### Acknowledgements

As part of the "France 2030" initiative, this work has benefited from a national grant managed by the French National Research Agency (Agence Nationale de la Recherche) attributed to the Exa-MA project of the NumPEx PEPR program, under the reference ANR-22-EXNU-0002.

# Nonlinear acoustic inverse problem with cross-correlation data for passive imaging

Jean Dutheil[1], Florian Faucher[1]

[1]Team Makutu, Inria, University of Pau and Pays de l'Adour, CNRS UMR 5142, France.

**Abstract**

We study a quantitative inverse problem in passive imaging, where medium properties are reconstructed using ambient noise signals. We consider the time-harmonic acoustic wave equation, and model the noise as a superposition of wavefields generated by spatially uncorrelated stochastic sources. Under this assumption, the expected value of cross-correlations between two signals is related to the deterministic Green's function and the source covariance. Given these correlations, we aim to reconstruct the density, wavespeed, and covariance of the sources. The inversion is formulated as a nonlinear optimization problem that minimizes a misfit functional, with gradients computed via the adjoint-state method. We carry out numerical experiments to demonstrate the feasibility of the parameter reconstruction.

## 1 Introduction

Passive imaging exploits ambient noise recordings to infer properties of a medium. The recorded signals are modeled as a superposition of waves generated by stochastic sources. Under suitable assumptions on those excitations, the expected cross-correlation of recorded signals can be expressed in terms of the deterministic Green's function and the source covariance [1,4,5]. This relation provides the key ingredient for quantitative parameter identification.

In this work, we consider the inverse wave problem of reconstructing acoustic media from the expected value of cross-correlation data. The problem is formulated as a nonlinear least-squares minimization and solved iteratively, following the approach of Full Waveform Inversion (FWI) used in seismic imaging, [5]. Numerical experiments with synthetic data are performed to demonstrate the accuracy of the proposed methodology for quantitative inversion.

## 2 Time-harmonic modeling

We consider acoustic waves in a bounded domain $\Omega \subset \mathbb{R}^3$, given in terms of the scalar pressure field $p$ and the vector particle velocity $\mathbf{v}$ which are solutions to, at angular frequency $\omega$:

$$\begin{cases} -\dfrac{\mathrm{i}\omega}{\rho(\mathbf{x})c^2(\mathbf{x})}\,p(\mathbf{x}) + \nabla\cdot\mathbf{v}(\mathbf{x}) = f(\mathbf{x}), & \text{(1a)} \\ -\,\mathrm{i}\omega\rho(\mathbf{x})\,\mathbf{v}(\mathbf{x}) + \nabla p(\mathbf{x}) = \mathbf{g}(\mathbf{x}), & \text{(1b)} \end{cases}$$

with absorbing or Dirichlet boundary conditions. The medium is characterized by the density $\rho$ and the wavespeed $c$. We denote Eq. (1) as

$$\mathcal{L}\begin{pmatrix} p \\ \mathbf{v}\end{pmatrix} = \mathcal{F}\,, \quad \text{where} \quad \mathcal{F} := \begin{pmatrix} f \\ \mathbf{g}\end{pmatrix}. \qquad (2)$$

The fundamental Green's kernel $G_{\mathbf{x}_2}$ solves:

$$\mathcal{L}\,G_{\mathbf{x}_2}(\mathbf{x}) = \delta(\mathbf{x}-\mathbf{x}_2)\mathbb{I}_{4\times 4}, \qquad (3)$$

with $\mathbb{I}_{4\times 4}$ the identity matrix and $\delta$ the Dirac distribution. The Green's kernel $G$ has a four-by-four structure. For Eq. (1), $G$ satisfies reciprocity, i.e., $\forall\,\mathbf{x}_1,\,\mathbf{x}_2 \in \Omega$, $G_{\mathbf{x}_2}(\mathbf{x}_1) = G_{\mathbf{x}_1}(\mathbf{x}_2)$.

We model the random source excitations as having zero-mean, and assume that the sources are spatially uncorrelated with covariance $\mathcal{S}$:

$$\begin{aligned}\mathbf{Cov}[\mathcal{F}](\mathbf{x}_1,\mathbf{x}_2) &:= \mathbb{E}[\mathcal{F}(\mathbf{x}_1)\overline{\mathcal{F}^\top(\mathbf{x}_2)}] \\ &= \mathcal{S}(\mathbf{x}_1)\,\delta(\mathbf{x}_1-\mathbf{x}_2),\end{aligned} \qquad (4)$$

This assumption is commonly used, cf. [1,5].

Under this assumption, and using the reciprocity of $G$, the expected value of the cross-correlation of wavefields at positions $\mathbf{x}_1$, $\mathbf{x}_2$, is expressed as, [4,5],

$$\begin{aligned} C(\mathbf{x}_1,\mathbf{x}_2) &:= \mathbf{Cov}[(p,\mathbf{v})^\top](\mathbf{x}_1,\mathbf{x}_2) \\ &= \int_\Omega G_{\mathbf{x}_1}(\mathbf{y})\,\mathcal{S}(\mathbf{y})\,\overline{G_{\mathbf{x}_2}(\mathbf{y})}^\top\,\mathrm{d}\mathbf{y}.\end{aligned}$$

The cross-correlation $C(\mathbf{x}_1,\mathbf{x}_2)$ is a matrix whose entries represent pairwise correlations between the components of $(p,\mathbf{v})$ at $\mathbf{x}_1$ and $\mathbf{x}_2$. For a given pair $(\mathbf{x}_1,\mathbf{x}_2)$ the simulation of $C(\mathbf{x}_1,\mathbf{x}_2)$ is carried out in two stages: first compute a response corresponding to a Dirac source in $\mathbf{x}_2$, Eq. (3). From this solution and the source covariance $\mathcal{S}$, solve:

$$\mathcal{L}\,C(\cdot,\mathbf{x}_2) = \mathcal{S}(\cdot)\overline{G_{\mathbf{x}_2}(\cdot)}^\top\,, \qquad (5)$$

with $C(\mathbf{x}_1,\mathbf{x}_2)$ the expected value of the cross-correlation between $\mathbf{x}_2$ and any position $\mathbf{x}_1$.

We solve the wave equation using the Hybridizable Discontinuous Galerkin (HDG) method, which solves the first-order formulation Eq. (1). The method introduces trace variables defined on the mesh skeleton, and the global system is only written in terms of the trace of the scalar unknown, cf. [3]. This allows us to access all unknowns of the first-order system at a reduced computational cost. The implementation is carried out in the software `Hawen`, [2].

## 3 Nonlinear inversion

The inverse problem consists in reconstructing the model parameters, $\rho$, $c$ and $\mathcal{S}$ from observations of expected cross-correlations. The reconstruction follows an iterative minimization of a misfit functional. Let $\mathbf{d}$ denote the observed cross-correlations between receiver pairs in the discrete set $\mathcal{R} = \{(\mathbf{x}_1, \mathbf{x}_2) : i = 1, \ldots, N_r\}$ ,with $N_r$ the number of receivers. We consider the least-squares misfit functional,

$$J(c, \rho, S) = \frac{1}{2} \sum_{(\mathbf{x}_1,\mathbf{x}_2)\in\mathcal{R}} \|C(\mathbf{x}_1, \mathbf{x}_2) - \mathbf{d}(\mathbf{x}_1, \mathbf{x}_2)\|^2 ,$$

where $C$ is obtained from Eqs. (5), (3).

The minimization of $J$ follows using the adjoint-state method. It requires adapting the standard FWI steps to the computation of cross-correlation, and to derive the adjoint-state method according to the HDG method, [3].

## 4 Numerical results

The forward and inverse algorithms are implemented in the open-source code `Hawen`, [2]. We consider a 2D domain of size 20km × 4.5km. Dirichlet boundary conditions are imposed at the top boundary, and absorbing conditions elsewhere. Synthetic data are generated using the wavespeed model shown in Fig. 1a, together with a structurally similar density model and a constant source covariance $\mathcal{S} = 1$, with noise incorporated to the synthetic data. For $\mathbf{x}_1$ and $\mathbf{x}_2$, $N_r = 400$ positions are taken near the surface.

We perform the reconstruction of wavespeed, starting from the initial model shown in Fig. 1b. The source covariance is assumed to be known, and the density is fixed to the constant value $\rho = 1$. The inversion is performed sequentially over five frequencies from 2 Hz to 6 Hz. The final reconstruction is shown in Fig. 1c.

These preliminary results show promising accuracy of the wavespeed reconstruction from expected cross-correlations, with the stratifications clearly resolved. Our ongoing works investigate the multi-parameter reconstructions for joint inversion of $\mathcal{S}$, $c$ and $\rho$, requiring an appropriate strategy. We further carry out comparisons with (traditional) inversion using direct wavefield, and highlight the impact of working with passive imaging data sets.

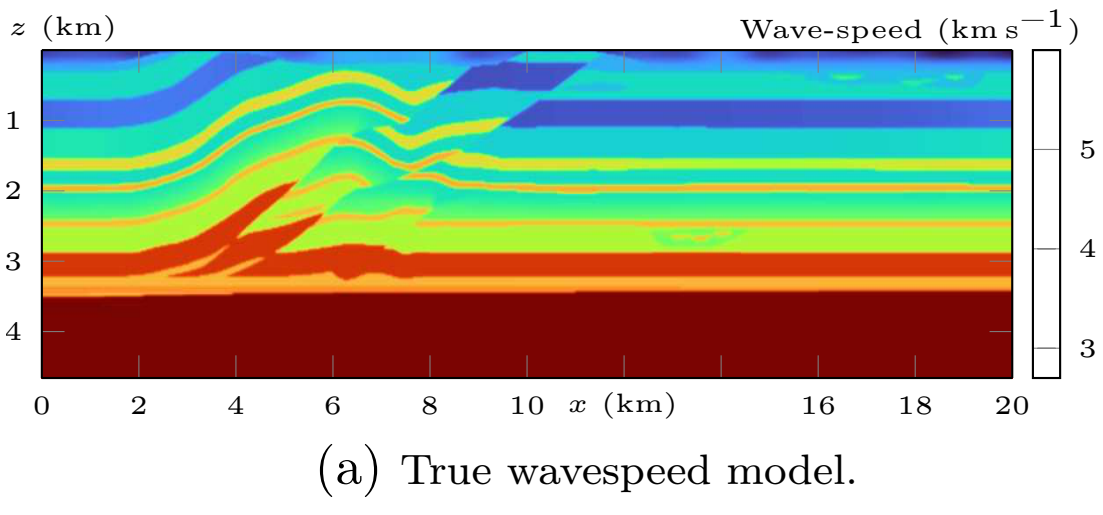

(a) True wavespeed model.

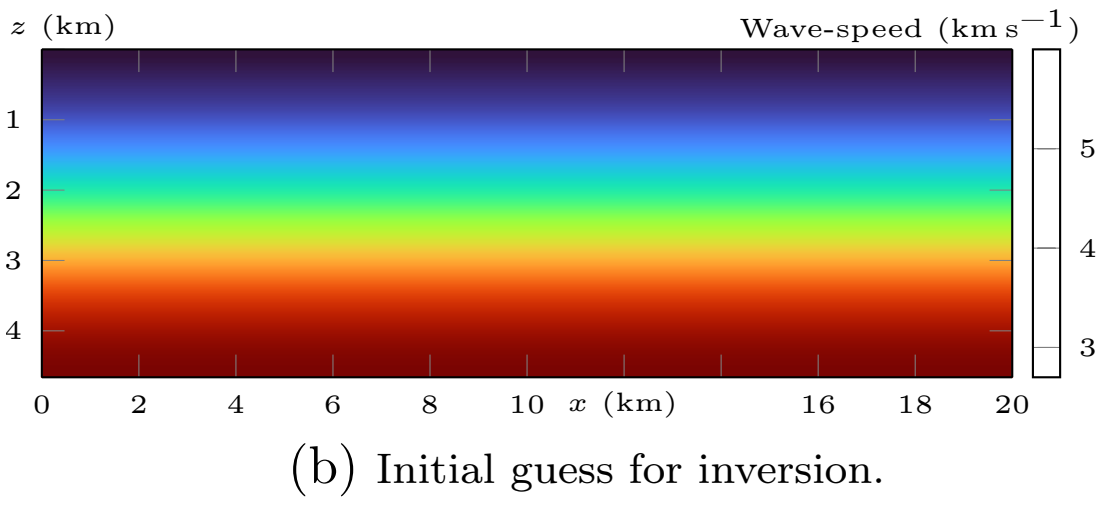

(b) Initial guess for inversion.

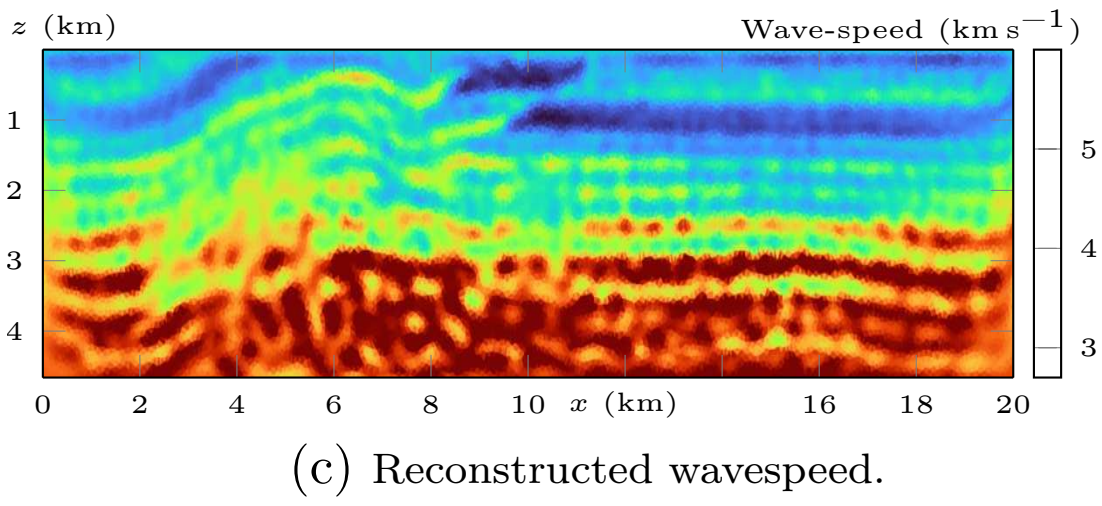

(c) Reconstructed wavespeed.

Figure 1: Experiment of reconstruction using expected values of cross-correlation data.

# Sixth-Order Compact FDM for Helmholtz Interface Equations with Reduced Pollution Effect and Its Extension to Nonlinear Equations

**Qiwei Feng**[1,*], **Bin Han**[2], **Michelle Michelle**[2], **Catalin Trenchea**[3]

[1]Mathematics Department, King Fahd University of Petroleum and Minerals, Dhahran, Saudi Arabia
[2]Department of Mathematical and Statistical Sciences, University of Alberta, Edmonton, Canada
[3]Department of Mathematics, University of Pittsburgh, Pittsburgh, PA, USA

*Email: qiwei.feng@kfupm.edu.sa,qfeng@ualberta.ca

**Abstract**

The equation with the high-frequency solution is numerically challenging to solve. To obtain a reasonable solution, a mesh size that is much smaller than the reciprocal of the wavenumber is typically required (known as the pollution effect). High-order schemes are desirable, because they are better in mitigating the pollution effect. In this talk, we present a high-order compact finite difference method (FDM) for 2D Helmholtz equations with singular sources, which can also handle any possible combinations of boundary conditions (Dirichlet, Neumann, and impedance) on a rectangular domain. To reduce the pollution effect, we propose a new pollution minimization strategy that is based on the average truncation error of plane waves. Precisely, we start with a stencil having a given accuracy order with some free parameters. Then, we plug in the plane wave solution into the scheme and minimize the average truncation error of plane waves over the free parameters of stencils to reduce the pollution effect. Furthermore, we also derive high-order FDMs to solve the Helmholtz interface equation. Finally, we also extend above methods to choose free parameters to propose fourth-order compact FDMs with reduced pollution effects for nonlinear steady and unsteady convection-diffusion equations by applying Crank-Nicolson (CN), BDF3, BDF4 methods with the aid of the iterative method. This is joint work with John Burkardt, Bin Han, Michelle Michelle, and Catalin Trenchea.

## 1 Introduction

The Helmholtz equation is challenging to solve numerically due to several reasons. The first is due to its highly oscillatory solution, which necessitates the use of a very small mesh size $h$ in many discretization methods. Taking a mesh size $h$ proportional to the reciprocal of the wavenumber $\mathsf{k}$ is not enough to guarantee that a reasonable solution is obtained or a convergent behavior is observed. The mesh size $h$ employed in a standard discretization method often has to be much smaller than the reciprocal of the wavenumber $\mathsf{k}$. It is clear that the grid size requirement for high-order schemes is less stringent than low-order ones. This is why we presently focus on the construction of a high-order scheme.

## 2 Model

In [1], we study the 2D Helmholtz equation, which is a time-harmonic wave propagation model, with a singular source term along a smooth interface curve and mixed boundary conditions (Dirichlet, Neumann, and Robin boundary conditions). Let $\Omega := (l_1, l_2) \times (l_3, l_4)$ and $\psi$ be a smooth two-dimensional function. Consider a smooth curve $\Gamma := \{(x, y) \in \Omega : \psi(x, y) = 0\}$, which partitions $\Omega$ into two subregions: $\Omega_+ := \{(x, y) \in \Omega \, : \, \psi(x, y) > 0\}$ and $\Omega_- := \{(x, y) \in \Omega \, : \, \psi(x, y) < 0\}$. The model problem (see Fig. (1) for the illustration) is defined as follows:

$$\begin{cases} \Delta u + \mathsf{k}^2 u = f & \text{in } \Omega \setminus \Gamma, \\ [u] = g, \quad [\nabla u \cdot \vec{n}] = g_\Gamma & \text{on } \Gamma, \\ \mathcal{B}_i u = g_i, & \text{on } \partial\Omega, \end{cases} \qquad (1)$$

where $\mathsf{k}$ is the wavenumber, $f$ is the source term, and for any point $(x_0, y_0) \in \Gamma$,

$$[u](x_0, y_0) := \lim_{(x,y)\in\Omega_+,(x,y)\to(x_0,y_0)} u(x, y) - \lim_{(x,y)\in\Omega_-,(x,y)\to(x_0,y_0)} u(x, y),$$

$[\nabla u \cdot \vec{n}](x_0, y_0)$ can be defined similarly, $\vec{n}$ is the unit normal vector of $\Gamma$ pointing towards $\Omega_+$, and $\mathcal{B}_i \in \{\mathbf{I}_d, \frac{\partial}{\partial \vec{n}}, \frac{\partial}{\partial \vec{n}} - \mathsf{i}\mathsf{k}\mathbf{I}_d\}$ represent Dirichlet, Neumann, Robin boundary conditions.

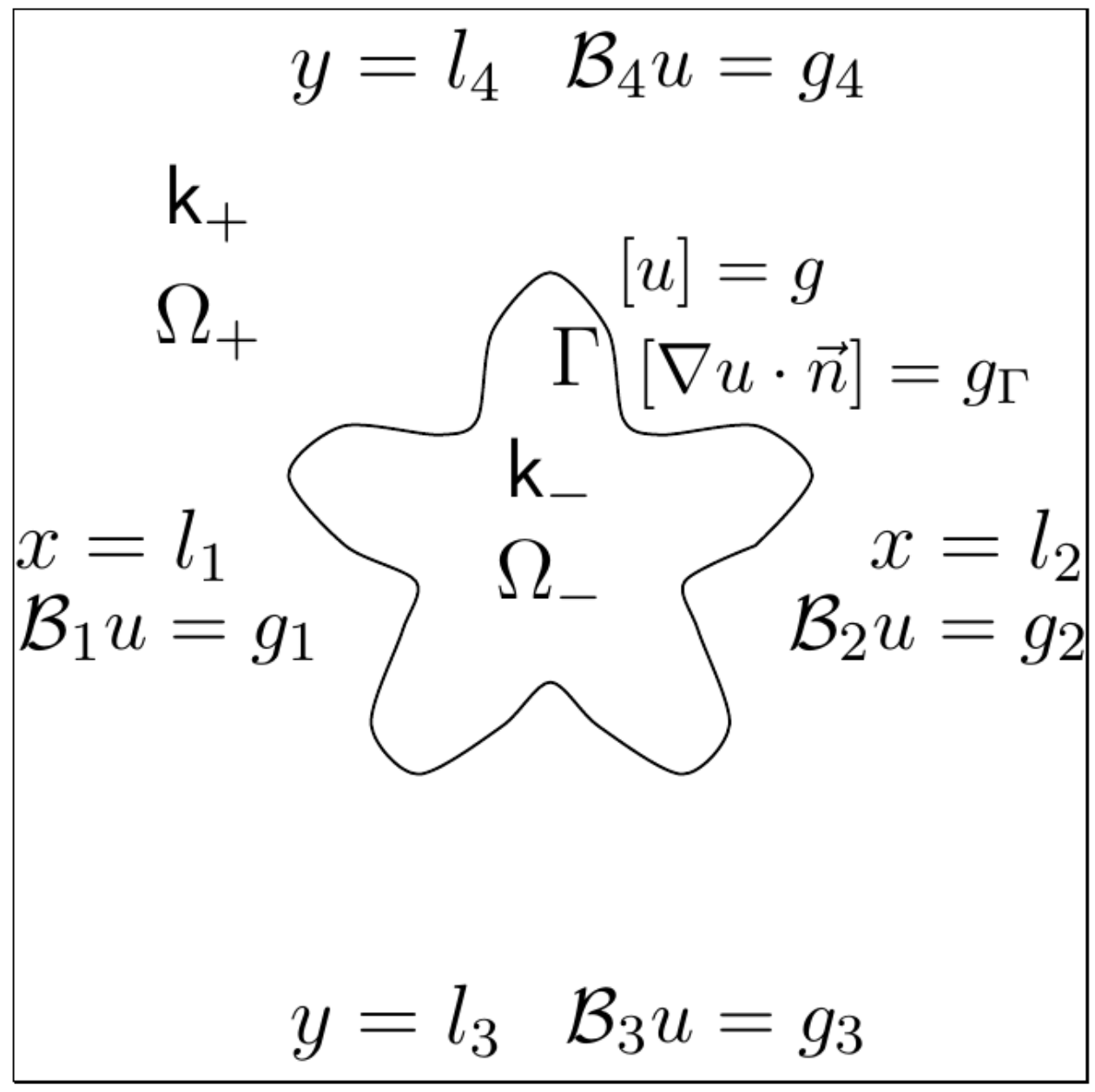

Figure 1: The illustration for the model problem (1), where $\mathcal{B}_1, \ldots, \mathcal{B}_4 \in \{\mathbf{I}_d, \frac{\partial}{\partial \vec{n}}, \frac{\partial}{\partial \vec{n}} - \mathsf{i}\mathsf{k}\mathbf{I}_d\}$, and $\mathsf{k}_\pm$ represents the wavenumber $\mathsf{k}$ in $\Omega_\pm$.

## 3 Methods and Results

We derive new fifth and sixth order compact finite difference schemes to solve Helmholtz equations with interfaces, piecewise constant wavenumbers, and mixed boundary conditions. We present a new pollution minimization strategy, which is based on the average truncation error of plane waves, and is applied to interior, boundary, and corner stencils. The idea is to first construct all possible compact stencils with the maximum consistency order and then to minimize the average truncation error of plane waves over the free parameters of stencils to reduce the pollution effect. For an interior regular point, we use a 9-point stencil with the sixth consistency order. To handle nonzero jump functions, we use a 9-point stencil with the fifth consistency order (for piecewise constant wavenumbers) and the seventh consistency order (for constant wavenumbers). For each corner, we use a 4-point stencil with at least the sixth consistency order. For each side, we use a 6-point stencil with at least the sixth consistency order. Our extensive numerical experiments demonstrate our proposed finite difference scheme with the reduced pollution effect outperforms several state-of-the-art finite difference schemes, particularly in the critical pre-asymptotic region where $\mathsf{k}h$ is near 1 with $\mathsf{k}$ being the wavenumber and $h$ the mesh size.

We also extend the pollution minimization strategy to derive fourth-order compact 9-point FDMs with reduced pollution effects to solve 2D nonlinear convection-diffusion equations in [2] as follows:

$$-\nabla \cdot (\kappa \nabla u) + \nabla \cdot \boldsymbol{F}(u) = f \qquad \text{steady},$$

and

$$u_t - \nabla \cdot (\kappa \nabla u) + \nabla \cdot \boldsymbol{F}(u) = f \qquad \text{unsteady},$$

where $\kappa(x, y)$ (steady) and $\kappa(x, y, t)$ (unsteady) are both smooth variable functions. For the steady nonlinear equation, we use an iteration method to reformulate the nonlinear equation into its linear counterpart, and derive a fourth-order compact 9-point finite difference method (FDM) to solve the reformulated equation on a uniform Cartesian grid. To increase the accuracy, we modify the FDM to reduce the pollution effect. The linear system of the FDM generates an M-matrix, provided the mesh size $h$ is sufficiently small. For the time-dependent nonlinear equation, we discrete the temporal domain using the Crank-Nicolson (CN), BDF3, BDF4 time stepping methods, and apply a similar iterative method to rewrite the nonlinear equation as the same linear convection-diffusion equation. Then we propose the second-order to fourth-order compact 9-point FDMs with reduced pollution effects on a uniform Cartesian grid. We prove that all FDMs satisfy the discrete maximum principle for sufficiently small $h$. Our method is still stable, accurate, and robust for various time at $t = 1, 10, 100, 500$, and for small ($\kappa \leq 10^{-3}$) and high-contrast ($\max \kappa / \min \kappa > 2 \times 10^4$) diffusion coefficients.

# The Helmholtz boundary element method suffers from the pollution effect

**María Ignacia Fierro**[1,*], **Jeffrey Galkowski**[2], **Manas Rachh**[3], **Euan Spence**[1]
[1]Department of Mathematical Sciences, University of Bath, Bath, United Kingdom
[2]Department of Mathematics, University College London, London, United Kingdom
[3]Department of Mathematics, Indian Institute of Technology Bombay, Mumbai, India
*Email: mibfp20@bath.ac.uk

**Abstract**

We consider solving the standard second-kind boundary integral equations for the Helmholtz exterior Dirichlet problem with the $h$-boundary element method (i.e., the Galerkin method with fixed-degree piecewise-polynomial subspaces).

We prove that the $h$-BEM suffers from the pollution effect – i.e., the number of degrees of freedom growing like $k^{d-1}$ is not sufficient for $k$-uniform quasi-optimality – and give numerical experiments demonstrating this.

## 1 Introduction

The standard second-kind boundary integral equations (BIEs) for solving the exterior Dirichlet problem involve the operators

$$A_k := \frac{1}{2}I + K_k - \mathrm{i}kS_k, \quad A'_k := \frac{1}{2}I + K'_k - \mathrm{i}kS_k, \tag{1}$$

where $S_k$ and $D_k$ are the single- and double-layer operators, $D'_k$ is the adjoint of $D_k$ in $L^2(\Gamma)$, and $\Gamma$ is the obstacle boundary. For a scattering problem at frequency $k$, the associated BIE solution oscillates at frequency $\lesssim k$ with little additional structure. The Nyquist–Shannon–Whittaker sampling theorem then implies that the dimension of this space is $\sim k^{d-1}$. Thus, it is not possible to achieve accuracy for data with a space of dimension $\ll k^{d-1}$.

The *pollution effect* is then defined as follows for an integral equation $\mathcal{A}v = f$ posed in $L^2(\Gamma)$.

**Definition 1 (The pollution effect in $L^2(\Gamma)$)** *Let $\mathcal{V}$ be a collection of subspaces of $L^2(\Gamma)$. For $V \in \mathcal{V}$ and $k > 0$, let $\mathcal{A}_V^{-1}(k) : L^2(\Gamma) \to V$ an approximation of $\mathcal{A}^{-1}(k)$.*

*The pair* $(\mathcal{V}, \{\mathcal{A}_V^{-1}\}_{V\in\mathcal{V}})$ does not suffer from the pollution effect in $L^2(\Gamma)$ *if there exists $k_0, \Lambda$, and $C_{\rm qo}$ such that for all $k \geq k_0$, $V \in \mathcal{V}$ with* $\dim V \geq \Lambda k^{d-1}$, *and $f \in L^2(\Gamma)$,*

$$\left\|\mathcal{A}^{-1}f - \mathcal{A}_V^{-1}f\right\|_{L^2(\Gamma)} \leq C_{\rm qo} \min_{w\in V} \left\|\mathcal{A}^{-1}f - w\right\|_{L^2(\Gamma)}.$$

That is $(\mathcal{V}, \{\mathcal{A}_V^{-1}\}_{V\in\mathcal{V}})$ does *not* suffer from the pollution effect if $k$-uniform quasi-optimality is achieved (for all possible data) with a choice of subspace dimension proportional to $k^{d-1}$.

## 2 The pollution effect for $A_k$ and $A'_k$

The paper [1] proves several results showing the the $h$-BEM applied to the integral equations (1) suffers from the pollution effect when the obstacle is trapping. For brevity, we describe here the simplest of these results; the talk will discuss other results in [1], along with new numerical experiments illustrating pollution.

**Definition 2 (Four-diamond domain)** *$\Omega^-$ is a* four-diamond domain *if there exists $0 < \epsilon < \pi/2$ such that $\Omega^- = \cup_{i=1}^4 \Omega_i \Subset \mathbb{R}^2$, where $\Omega_i$ are open, convex, disjoint, and have smooth boundary so that*

$$\partial\big((-\pi,\pi)\times(-\pi,\pi)\big) \setminus \big(\cup_{a,b\in\pm1} B((a\pi,b\pi),\epsilon)\big) \subset \Gamma \subset \mathbb{R}^2 \setminus \big(\cup_{a,b\in\pm1} B((a\pi,b\pi),\epsilon/2)\big),$$
$$\Omega^- \cap (-\pi,\pi)\times(-\pi,\pi) = \emptyset.$$

That is, the exterior region between the components contains a square trapping region with corner exit holes (Figure 1).

In the following result we assume there is $\mathcal{J}$, such that $\mathcal{K} := \mathbb{R} \setminus \mathcal{J}$ as the set of frequencies where $\|(A'_k)^{-1}\|_{L^2(\Gamma)\to L^2(\Gamma)} = \|A_k^{-1}\|_{L^2(\Gamma)\to L^2(\Gamma)}$ is polynomially bounded in $k$, and we denote $p$ the polynomial degree of the BEM space.

**Theorem 3 (Pollution for the four-diamond domain when $p = 0$)** *Let $p = 0$, let $\Omega^-$ be a four-diamond domain, and $k_n := \sqrt{2}n$. Then there is $c > 0$ such that for all $k \in \mathcal{K}$ and $n$ such that $k_n^{-1} \leq \delta \leq 1$, $|k - k_n| < \delta^{-1}k_n^{-1}$, and all $h > 0$ there is $f \in L^2(\Gamma)$ such that the Galerkin* *solution, $v_h$, to $\mathcal{A}v = f$ with $\mathcal{A}$ given by either $A_k$ or $A'_k$, if it exists, satisfies*

$$\frac{\|v_h - v\|_{L^2(\Gamma)}}{\|(I - P^G_{V_h})v\|_{L^2(\Gamma)}} \geq \frac{1}{2} + c \begin{cases} (hk_n)^{-1}, & 1 \leq (hk_n)^2 \delta k_n, \\ (hk_n)\delta k_n, & \delta k_n (hk_n)^2 \leq 1. \end{cases} \tag{2}$$

Note that $k_n = \sqrt{2}n$ correspond to diagonal $(n, n)$ resonant modes of $(-\pi, \pi)^2$.

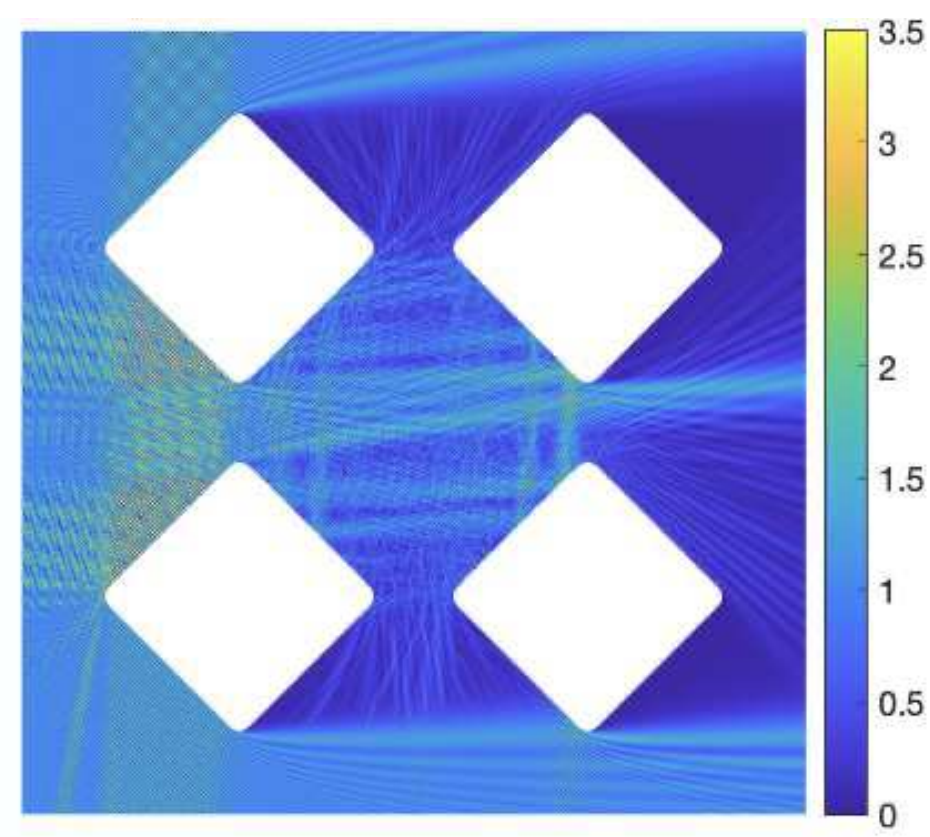

Figure 1: A four-diamond domain.

## 3 Four remarks about Theorem 3

1. Theorem 3 is illustrated by numerical experiments in Figure 2.

2. The bound (2) implies pollution. Indeed, if $hk = \epsilon$ then (2) implies that the quasioptimality constant is *at least* $c\epsilon^{-1}$ and hence that decreasing $\epsilon$ actually makes the quasioptimality constant worse (until $\epsilon \leq (k_n\delta)^{-1/2}$; i.e., until $\epsilon$ is a sufficiently-small $k$-dependent quantity). This "increase then decrease" in the quasioptimality constant for fixed $k$ as the number of points per wavelength increases is shown in Figure 3.

3. Regarding the polynomial boundedness assumption: If $\Omega^-$ consists of two aligned (rounded) squares with parabolic trapping, then [2, Cor. 1.14] shows

$$\|A_k^{-1}\|_{L^2(\Gamma)\to L^2(\Gamma)} + \|(A'_k)^{-1}\|_{L^2(\Gamma)\to L^2(\Gamma)} \leq Ck^2. \tag{3}$$

If this holds here, then $\mathcal{J} = \emptyset$ and (2) is sharp. For analytic $\Gamma$, quasioptimality can be recovered with $p \sim \log k$ [3], but that is not the case here.

4. We emphasise that the phenomenon behind Theorem 3 is *not* sparse in frequency. Indeed, Theorem 3 holds with $\delta = 1$ for a set of $k$ with infinite measure (since $\sum_{n=1}^{\infty} n^{-1} = \infty$). Furthermore, although we prove Theorem 3 with $k_n = \sqrt{2}n$, it is easy to generalize our construction in the proof of Theorem 3 to take $k_n = \sqrt{m_n^2 + \ell_n^2}$ for any $m_n, \ell_n \in \mathbb{Z}$ with $c < |\frac{m_n}{\ell_n}| < C$ and such that, if $m_n/\ell_n = p_n/q_n$ with $p_n$ and $q_n$ relatively prime, then $|q_n| \leq C$.

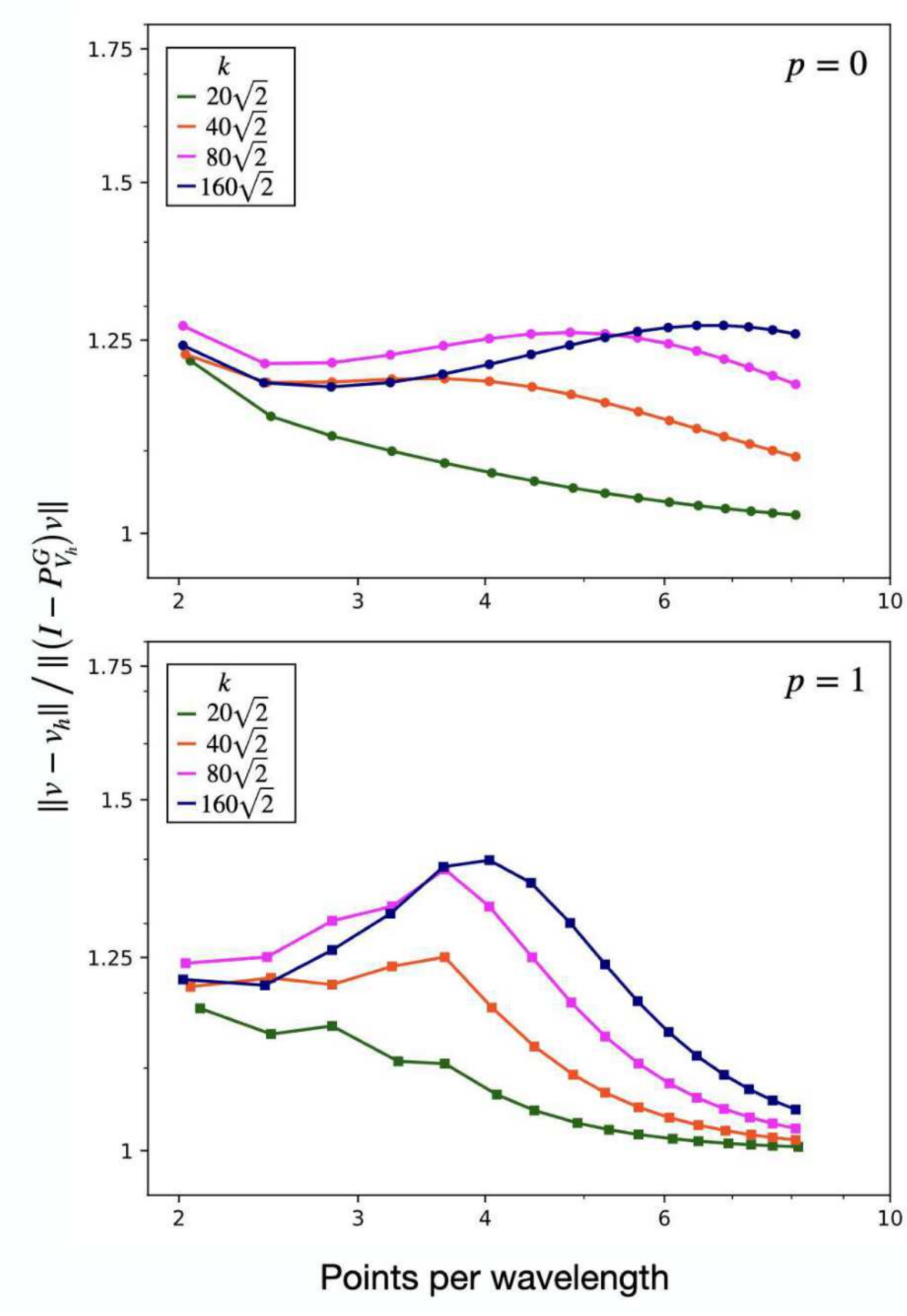

Figure 2: Plots of the ratio of the Galerkin error to the best approximation error for plane-wave incidence on the four diamond domain.

# Scattering problems in junctions of stratified media

**Sarah Al Humaikani**[1], **Anne-Sophie Bonnet-Ben Dhia**[1], **Sonia Fliss**, [,*]
[1]POEMS, CNRS, Inria, ENSTA, Institut Polytechnique de Paris, Palaiseau, France
*Email: sonia.fliss@ensta.fr

**Abstract**

This work deals with the wave scattering by a junction of stratified media, an important case being open waveguides. Our objective is to characterize and compute the outgoing solution using the Complex Scaled Halfspace Matching Method, which has been originally introduced for homogeneous media.

## 1 Introduction

We consider the Helmholtz equation in a domain $\Omega := \mathbb{R}^2 \setminus \mathcal{O}$ ($\mathcal{O}$ being a triangle or a rectangle) that coincides with the union of N (resp. $N = 3$ or $N = 4$) stratified half-planes $\mathrm{H}_n := \{(\mathrm{x}_n, \mathrm{y}_n) \in \mathbb{R}^2 \mid \mathrm{x}_n > 0\}$ for $n = 1, \dots, N$ , where $(\mathrm{x}_n, \mathrm{y}_n)$ is a local coordinate system that is obtained from the global ones $(x, y)$ by a rotation and a translation. Each half-plane $\mathrm{H}_n$ is called stratified in the sense that it is associated with a real-valued wavenumber function $\kappa_n$ that satisfies $\kappa_n = \kappa_n(\mathrm{y}_n)$. In other words, the stratification is orthogonal to the boundary $\Sigma_n$ of $\mathrm{H}_n$. Moreover we assume that $\kappa_n$ is bounded and becomes constant outside a bounded interval. Let $\kappa$ be the wavenumber that coincide with $\kappa_n$ in each half-space. We consider for simplicity an exterior problem and we want to characterize and compute the outgoing solution of

$$- \Delta u - \kappa^2 u = 0 \text{ in } \Omega, \quad u = g \text{ on } \partial\Omega, \quad (1)$$

where $g \in H^{1/2}(\partial\Omega)$.

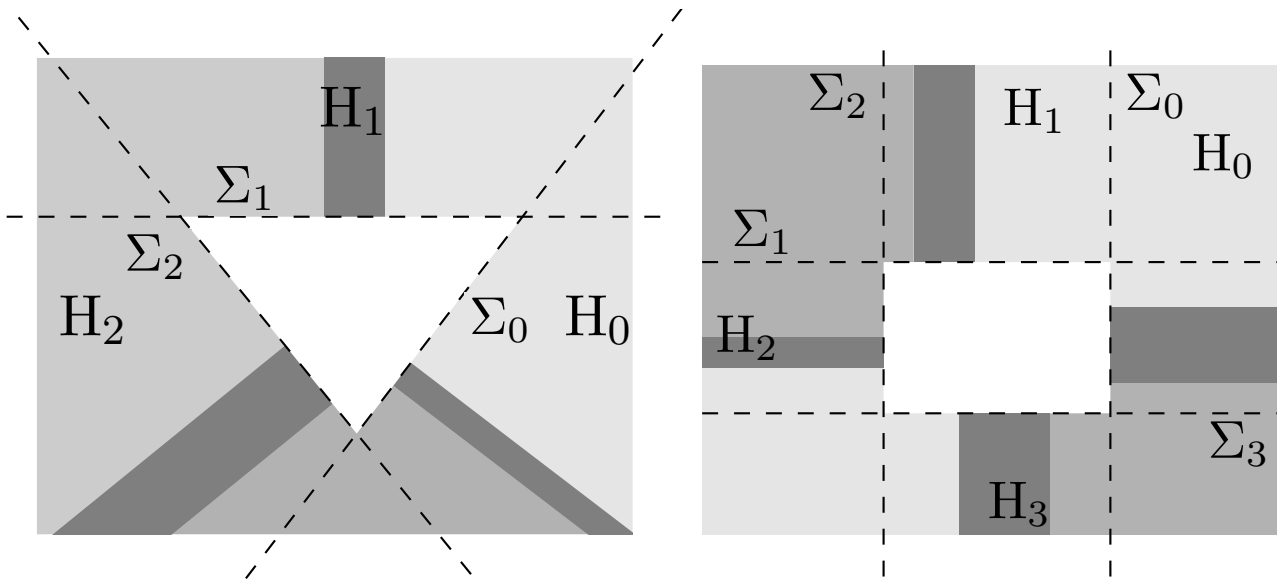

This is a challenging question because, for this class of problems, the definition of the outgoing solution is not clear. This difficulty is linked to the fact that the Green's function of the medium is not known, separation of variables cannot be used to derive either a transparent boundary conditions or an appropriate radiation condition. We propose to use the Halfspace Matching method (HSM). The method is based on half-plane representations of the solution in each $\mathrm{H}_n$ in terms of its trace $\varphi_n$ on the boundary of the half-plane $\Sigma_n$. The problem is then reformulated as a coupled system of integral equations, where the unknowns are precisely the traces $\varphi_n$. For the scattering problem under consideration, the appropriate functional framework for the traces as well as the well-posedness of the HSM system are not clear. Instead, we want to derive a similar system, called Complex Scaled(CS) HSM where the unknowns of the integral equations are exponentially decaying analytic extensions of the traces. Unlike the case of homogeneous media [1], the outgoing solution is not well defined and it is not straightforward to prove that the traces have indeed exponentially decaying analytic extensions. In the present case, the derivation of the CS-HSM system of equations is obtained by using the limiting absorption principle. This leads naturally to the definition of the outgoing solution.

## 2 Half-space (HS) representations

Let consider the problem with dissipation replacing $\kappa^2$ by $\kappa^2 + \mathrm{i}\varepsilon$, $\varepsilon > 0$ and let denote $u^\varepsilon \in H^1(\Omega)$ the associated solution. As explained above, the key tool is the HS representations. For homogeneous media, this representation can be derived by using various expressions of the Dirichlet half-space Green's function: using the Hankel function, the Fourier transform in the $\mathrm{x}_n$-direction or in the $\mathrm{y}_n$-direction. For stratified media, the Green's function is not known in a closed form. Its expressions can be obtained either using the Fourier transform in the $\mathrm{x}_n$-direction (that diagonalizes $-\partial^2_{\mathrm{x}_n}$) or the generalized Fourier transform (that diagonalizes $-\partial^2_{\mathrm{y}_n} - \kappa^2(\mathrm{y}_n)$ ), see for instance [2]. Each representa-tion has its own theoretical and numerical advantages and drawbacks.

## 3 Derivation of the CS-HSM system

**Step 1.** By writing the HS representation

$$u^\varepsilon = P_n^\varepsilon(\varphi_n^\varepsilon) \text{ in } \mathrm{H}_n,$$

one can derive the HSM system of equations, for $n \in \mathbb{Z}/N\mathbb{Z}$,

$$\varphi_n^\varepsilon = D_{n,n\pm1}^\varepsilon(\varphi_{n\pm1}^\varepsilon) \text{ on } \mathrm{H}_{n\pm1} \cap \Sigma_n, \quad (2)$$

that has to be completed with the boundary condition $\varphi_n^\varepsilon = g$ on $\partial\Omega \cap \Sigma_n$. Here the operator

$$D_{n,n\pm1}^\varepsilon : \varphi \mapsto P_{n\pm1}^\varepsilon(\varphi)\big|_{\mathrm{H}_{n\pm1}\cap\Sigma_n},$$

is the sum of an operator of norm $\leq 1/\sqrt{2}$ and a compact operator as soon as $\varepsilon > 0$. This property allows to prove that the HSM system is Fredholm for $\varepsilon > 0$. For this step, either representation may be used interchangeably.

**Step 2.** From (2) and the properties of the HS representations, we prove that for all $n$, $\varphi_n^\varepsilon$, considered as a function of the local coordinate $\mathrm{Y}_n$, has an analytic continuation w.r.t $\mathrm{Y}_n$ and there exists $\theta \in [0, \pi/2[$ such that $\varphi_n^{\varepsilon,\theta} := \varphi_n^\varepsilon \circ \tau_n^\theta$ is in $L^2(\mathbb{R})$ where $\tau_n^\theta$ is the complex path

$$\tau_n^\theta(s) = \pm a_n + (s \mp a_n)e^{\mathrm{i}\theta} \text{ if } \pm s > \pm a_n, \; s \text{ else},$$

where $a_n$ is chosen so that $\Sigma_n \cap \Sigma_{n\pm1} = \{\mathrm{Y}_n = \pm a_n\}$. Note that $\kappa_n$ is constant outside $(-a_n, a_n)$. The $\varphi_n^{\varepsilon,\theta}$'s are called the CS traces. This property is straightforward using the HS representation based on the generalized Fourier transform.

**Step 3.** A key property is that in $\mathrm{H}_n^\theta := \{(\mathrm{X}_n, \mathrm{Y}_n) \in \mathrm{H}_n, \mathrm{X}_n > \max((\mathrm{Y}_n - a)\tan\theta, -(\mathrm{Y}_n + a)\tan\theta)\}$, the solution $u_\varepsilon$ can be reconstructed from the CS trace $\varphi_n^{\varepsilon,\theta}$

$$u^\varepsilon = P_n^{\varepsilon,\theta}(\varphi_n^{\varepsilon,\theta}) \text{ in } \mathrm{H}_n^\theta,$$

thanks to a modified integral representation. This proof requires the use of the HS representation based on the Fourier Transform along $\mathrm{X}_n$. It is based mainly on Cauchy's residue theorem applied to the integral w.r.t. to the Fourier variables. This technique cannot be applied with the other HS representation because of the presence of complex resonances.

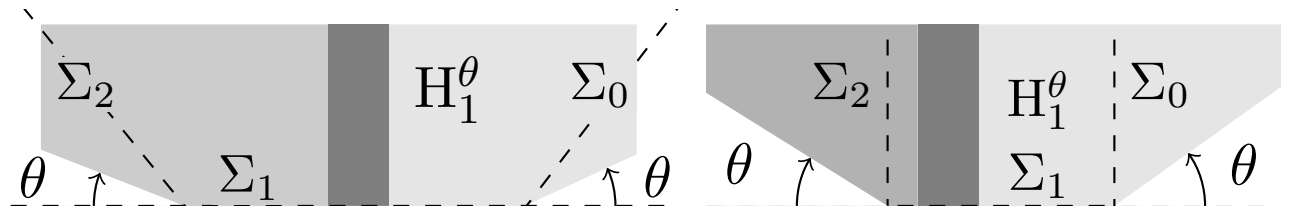

**Step 4.** When $\mathrm{H}_n^\theta \cap \Sigma_{n\pm1} = \mathrm{H}_n \cap \Sigma_{n\pm1}$ (i.e. when $\theta < \Theta$ where $\Theta$ is the maximal angle of the polygon $\mathcal{O}$), we prove that $P_n^{\varepsilon,\theta}(\varphi_n^{\varepsilon,\theta})$ have an analytical continuation along the line $\mathrm{H}_n^\theta \cap \Sigma_{n\pm1}$. This enables to derive the CS-HSM

$$\varphi_n^{\varepsilon,\theta} = D_{n,n\pm1}^{\varepsilon,\theta}(\varphi_{n\pm1}^{\varepsilon,\theta}) \text{ on } \mathrm{H}_{n\pm1} \cap \Sigma_n, \quad (3)$$

that has to be completed with the boundary condition $\varphi_n^{\varepsilon,\theta} = g$ on $\partial\Omega \cap \Sigma_n$. Here the operator is defined as

$$D_{n,n\pm1}^{\varepsilon,\theta} : \varphi \mapsto P_{n\pm1}^{\varepsilon,\theta}(\varphi)\big|_{\mathrm{H}_{n\pm1}\cap\Sigma_n} \circ \tau_n^\theta,$$

where in this expression, $P_{n\pm1}^{\varepsilon,\theta}(\varphi)\big|_{\mathrm{H}_{n\pm1}\cap\Sigma_n}$ is considered as a function of $\mathrm{Y}_n$. Again, the operator $D_{n,n\pm1}^{\varepsilon,\theta}$ is the sum of an operator of norm $\leq 1/\sqrt{2}$ and a compact operator. Now this property does not rely on the presence of dissipation but on $\theta \in (0, \Theta)$.

**Step 5.** We prove that the system of equations (CS-HSM) is Fredholm for $\varepsilon \geq 0$ and $\theta \in (0, \Theta)$. Moreover, whereas uniqueness holds for $\varepsilon > 0$, this is an open question for $\varepsilon = 0$.

**Step 6.** If uniqueness holds, we denote $\{\varphi_n^\theta, n = 1, \ldots, N\}$ the unique solution of the CS-HSM system for $\varepsilon = 0$. We finally prove the limiting absorption principle: namely that for all $n$, $\|\varphi_n^{\varepsilon,\theta} - \varphi_n^\theta\|_{L^2(\mathbb{R})} = \mathcal{O}(\varepsilon)$. Moreover if $\theta < \Theta/2$, we introduce the function $u$ as

$$u = P_n^\theta(\varphi_n^\theta) \text{ in } \mathrm{H}_n^\theta,$$

that is well defined in $\Omega$ when $\theta < \Theta/2$. We show that $u$ satisfies (1) and for all $K \subset \Omega$ bounded,

$$\|u_\varepsilon - u\|_{H^1(K)} = \mathcal{O}(\varepsilon).$$

This method defines as expected the outgoing solution and can be discretized for computations. Numerical results will be presented during the talk.

# Quasi-Trefftz spaces for time-harmonic Maxwell's equations

**Ilaria Fontana**[1,*], **Lise-Marie Imbert-Gérard**[1]
[1]Department of Mathematics, University of Arizona, Tucson, AZ, USA
*Email: ifontana@arizona.edu

**Abstract**

Classical Trefftz methods are high-order numerical schemes for Partial Differential Equations (PDEs) in which the discrete trial and test functions are local solutions of the governing equation. They have been extensively studied for wave propagation problems, for which they achieve high accuracy with reduced computational cost compared to classical finite element or standard Discontinuous Galerkin (DG) methods. However, their application is restricted to linear homogeneous PDEs with piecewise-constant coefficients, since exact local solutions are generally unavailable in more complex settings. Quasi-Trefftz methods overcome this limitation by replacing exact local solutions with approximate ones constructed at the element level. This extends the applicability of Trefftz-type schemes to variable coefficient problems, such as wave propagation in inhomogeneous media. This contribution focuses on extending the quasi-Trefftz approach to the first-order time-harmonic Maxwell's equations. We present recent theoretical results and supporting numerical simulations.

## 1 Introduction

Trefftz and Taylor-based quasi-Trefftz methods have been studied for wave propagation problems. In the context of electromagnetism, classical Trefftz methods have been designed and analyzed for the homogeneous second-order time-harmonic Maxwell's equation [3]. The discrete space is tailored to the equation by taking a set of plane waves as its basis. A first effort to extend this analysis to a polynomial Taylor-based quasi-Trefftz approach can be found in [4]. There, the second-order formulation is considered with possibly varying-in-space electric permittivity. In this work, we introduce, analyze, and implement quasi-Trefftz discrete spaces for the following homogeneous first-order time-harmonic Maxwell's equations [7]:

$$\begin{cases} \nabla \times \mathbf{E} - \mathrm{i}\kappa\varepsilon\mathbf{B} = \mathbf{0} & \text{in } \Omega \\ \nabla \times \mathbf{B} + \mathrm{i}\kappa\mu\mathbf{E} = \mathbf{0} & \text{in } \Omega, \end{cases} \tag{1}$$

where $\Omega \subset \mathbb{R}^3$ is an open set, $\mathbf{E}$ and $\mathbf{B}$ are the unknown electric and magnetic field, i is the imaginary unit, $\kappa > 0$ is a fixed wavenumber, $\varepsilon$ is the electric permittivity and $\mu$ is the magnetic permeability. The latter two coefficients could depend on space, i.e., $\varepsilon = \varepsilon(\mathbf{x})$, $\mu = \mu(\mathbf{x})$.

## 2 Definition of the Quasi-Trefftz Space

Consider a positive integer $p$, an open set $T \subset \Omega$, a point $\mathbf{x}_T \in T$, and $\varepsilon, \mu \in C^p(T;\mathbb{C})$. Let $\mathcal{T}_q^{\mathbf{x}_T}$ and $\boldsymbol{\mathcal{T}}_q^{\mathbf{x}_T}$ denote the Taylor polynomial operator of degree $q$ with center in $\mathbf{x}_T$ for scalar-valued and vector-valued functions, respectively. Then, the *local Taylor-based quasi-Trefftz space* of the order $p$ for the first-order system (1) is the polynomial space

$$\mathbb{QT}_p(T) \subset (\mathbb{P}_p)^3 \times (\mathbb{P}_p)^3$$

whose elements $(\mathbf{Q}^E, \mathbf{Q}^B)$ satisfy the following *quasi-Trefftz properties*:

$$\begin{cases} \boldsymbol{\mathcal{T}}_{p-1}^{\mathbf{x}_T}[\nabla \times \mathbf{Q}^E - \mathrm{i}\kappa\mu\mathbf{Q}^B] = \mathbf{0} \\ \boldsymbol{\mathcal{T}}_{p-1}^{\mathbf{x}_T}[\nabla \times \mathbf{Q}^B + \mathrm{i}\kappa\varepsilon\mathbf{Q}^E] = \mathbf{0} \\ \mathcal{T}_{p-1}^{\mathbf{x}_T}[\nabla \cdot (\varepsilon\mathbf{Q}^E)] = 0 \\ \mathcal{T}_{p-1}^{\mathbf{x}_T}[\nabla \cdot (\mu\mathbf{Q}^B)] = 0 \end{cases} \tag{2}$$

Additionally, given a mesh $\mathcal{T}_h$ covering $\Omega$, the *global Taylor-based polynomial quasi-Trefftz space* of order $p$, $\mathbb{QT}_p(\mathcal{T}_h)$, consists of all couples of polynomials such that $(\mathbf{Q}^E, \mathbf{Q}^B)|_T \in \mathbb{QT}_p(T)$ for all $T \in \mathcal{T}_h$.

## 3 Theoretical and Numerical Results

The following theorem demonstrates that the local quasi-Trefftz space characterized by (2) exhibits the same best approximation properties as the full polynomial space. The proof strongly relies on the properties of the Taylor polynomial operators, as we show in [2].

**Theorem 1 (Approximation properties)** *Let $T \subset \mathbb{R}^3$ be an open set, $\mathbf{x}_T \in T$, $p \in \mathbb{N}$, and $\mathbb{QT}_p(T)$ the local polynomial quasi-Trefftz space characterized by the quasi-Trefftz properties (2), and suppose that $\varepsilon, \mu \in C^p(T;\mathbb{C})$. If $(\mathbf{E},\mathbf{B}) \in C^{p+1}(T;\mathbb{C}^3)\times C^{p+1}(T;\mathbb{C}^3)$ satisfies the Maxwell's system (1), then*

1. $(\boldsymbol{\mathcal{T}}_p^{\mathbf{x}_T}[\mathbf{E}], \boldsymbol{\mathcal{T}}_p^{\mathbf{x}_T}[\mathbf{B}]) \in \mathbb{QT}_p(T)$;
2. *if $T$ is star-shaped about $\mathbf{x}_T$, then*

$$\inf_{(\mathbf{Q}^E,\mathbf{Q}^B)\in\mathbb{QT}_p} \left(\|\mathbf{E}-\mathbf{Q}^E\| + \|\mathbf{B}-\mathbf{Q}^B\|\right) \leq \frac{3^{p+1}}{(p+1)!} h_T^{p+1} \left(\|\mathbf{E}\| + \|\mathbf{B}\|\right)$$

*where $\|\cdot\|$ is the $C^0$-norm on $T$.*

Recent works [1,5,6] have demonstrated that quasi-Trefftz functions and bases can be numerically constructed through an efficient and fully explicit algorithm. This procedure works layer-by-layer over spaces of homogeneous polynomials, and its well-posedness is based on the surjectivity of the PDE's highest-order differential operator. Extending this approach to the time-harmonic Maxwell's system (1) requires new tools, since, in contrast to the scalar case, the highest-order differential operator (i.e., the curl operator) is not surjective. More in detail, the key ingredients of our approach [2] are the Helmholtz decomposition of homogeneous polynomial vector fields and the exactness of the polynomial de Rham sequence. Additionally, we prove that the dimension of the quasi-Trefftz $\mathbb{QT}_p(T)$ is significantly smaller than the full polynomial space.

**Theorem 2 (Quasi-Trefftz dimension)** *Under the hypothesis of Theorem 1, the dimension of the local quasi-Trefftz space is*

$$\dim \mathbb{QT}_p(T) = 2(p+1)(p+3).$$

We conclude by showing a numerical example. Consider the following settings: $\kappa = 1$, $\varepsilon(\mathbf{x}) = 1 - x$, $\mu = \mu(\mathbf{x}) = 1$, and exact solution

$$\begin{cases} \mathbf{E}(\mathbf{x}) = \mathrm{Ai}(x)e^{\mathrm{i}y}\mathbf{e}_3 \\ \mathbf{B}(\mathbf{x}) = \mathrm{Ai}(x)e^{\mathrm{i}y}\mathbf{e}_1 + \mathrm{i}\mathrm{Ai}'(x)e^{\mathrm{i}y}\mathbf{e}_2 \end{cases} \tag{3}$$

Figure 1 illustrates the best approximation properties predicted by Theorem 1. It plots the $C^0$-norm, computed on spheres centered at a point, of the difference between the exact solution (3) and the quasi-Trefftz function matching the coefficients of its Taylor expansion.

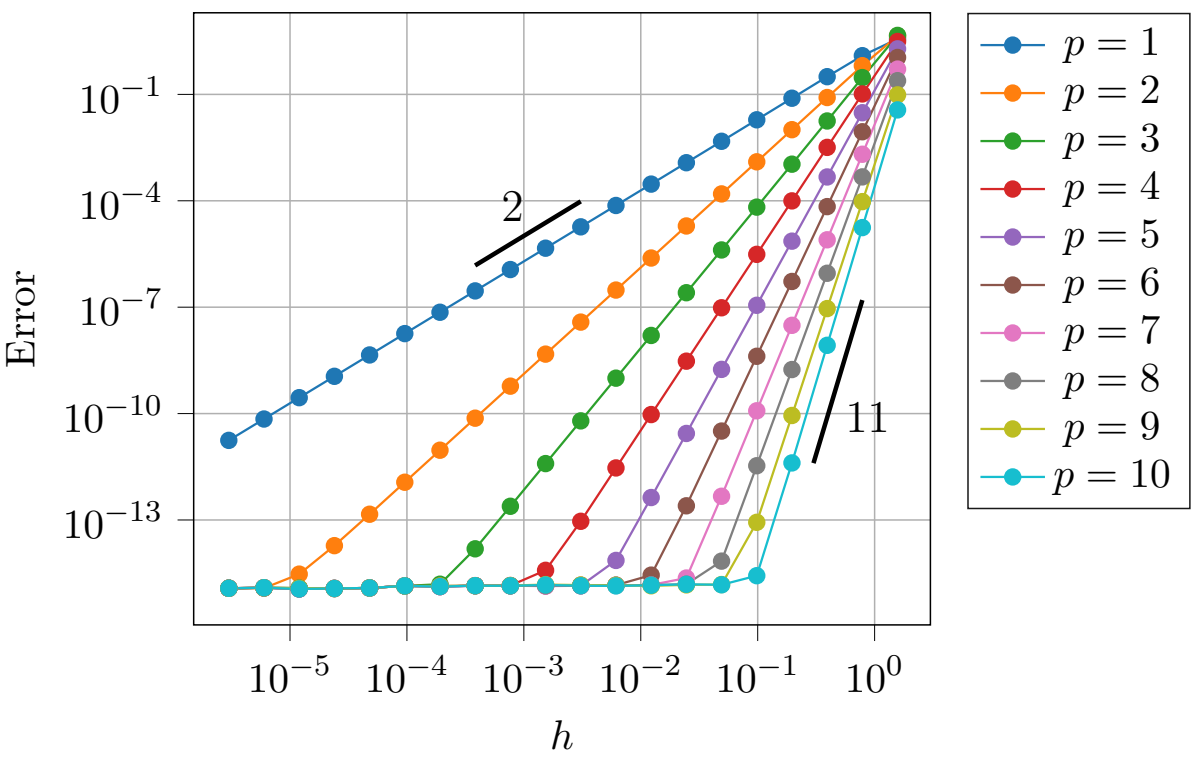

Figure 1: Convergence plots near the point $\mathbf{x}_T = (-4, 1, 0.5)$ for the approximation error. The settings are as follows: $\kappa = 1$, $\varepsilon(\mathbf{x}) = 1 - x$, $\mu = \mu(\mathbf{x}) = 1$, and exact solution defined by (3).

# Reconstructing inertial-mode flows from the correlation function of acoustic modes

**Hélène Barucq**[1], **Florian Faucher**[1], **Damien Fournier**[2,*], **Laurent Gizon**[1], **Zhi-Chao Liang**[2], **Ha Pham**[1]

[1]Project-Team Makutu, Inria, University of Pau and Pays de l'Adour, CNRS UMR 5142, France
[2]Max-Planck-Institut für Sonnensystemforschung, Göttingen, Germany

*Email: fournier@mps.mpg.de

**Abstract**

We consider the inversion of the two-point correlation of acoustic wavefields to infer the 3D structure of flows in the near-surface convection zone. To detect and characterize solar inertial modes, we employ data collected over 15 years to target low-frequency large-scale dynamics. The forward and inverse problems are solved in the frequency-domain with the code `Hawen`.

## 1 Introduction

Correlations of observed acoustic wavefields at the solar surface encode information of 3D flows and thermodynamic perturbations in the solar outer layers. Of current interest in the seismology of stars is the flow associated with inertial modes, which are oscillations restored by the Coriolis force and at frequencies ($\sim$ a few months in period) about three orders of magnitude lower than those of acoustic waves ($\sim$ 5 mins). While evolving on much longer timescales and with small velocity amplitudes [4], inertial modes leave faint but measurable effects on the acoustic waves. Long-time observations allow to measure them, and detect for example Rossby modes using (acoustic) wave travel times [2].

Since the cross-correlation function (CCF) contains more information than travel times [3], we consider inverting observed acoustic CCF, denoted as $\mathcal{C}_{\text{obs}}$, to recover flows associated with the solar inertial modes as shown below.

parameter $\mathbf{u}_m(\eta, z) \xrightarrow{\mathcal{L}}$ Green's kernel $\mathbb{G} \rightsquigarrow v_{\text{los}} \overset{\mathbb{M}}{\rightsquigarrow} \mathcal{C}_{\mathbf{u}}$
Misfit $\updownarrow$
Inversion ← $\mathcal{C}_{\text{obs}}$

The inertial mode is modeled as $\mathbf{u}_m(\eta, z)e^{\mathrm{i}(m\phi - \omega_m t)}$, with $m$ the longitudinal wavenumber, $\omega_m$ the frequency of the mode. Also shown above, the forward problem computes for each model $\mathbf{u}_m(\eta, z)$ the expectation value of CCF, denoted as $\mathcal{C}_{\mathbf{u}}$, from Green's kernel $\mathbb{G}$ and the source covariance matrix $\mathbb{M} = \mathbb{E}[\mathbf{s}(\mathbf{x_1}), \mathbf{s}(\mathbf{x_2})]$.

In dealing with two dynamics, the stated problem is challenging with difficulty ranging from data interpretation and preparation for inversion, forward modeling, to inversion strategy. Furthermore, the last two elements have to be adapted to the current choice of the wave solver constructed in [1].

## 2 Numerical modeling of oscillations

Oscillations result from perturbations of a reference fluid state described by density $\rho_0$, pressure $p_0$, sound speed $c_0$, flow velocity $\mathbf{u}$, and gravity potential $\phi_0$. Below the vector $\mathbf{e}_\eta$, $\mathbf{e}_\phi$, $\mathbf{e}_z$ are basis in the cylindrical coordinates $(\eta, \phi, z)$.

The time-harmonic adiabatic oscillation equation describes the Eulerian perturbations of the fluid states $(\delta_{\mathbf{v}}, \delta_\rho, \delta_p)$ and the Lagrangian displacement $\boldsymbol{\xi}$. For a source $\mathbf{s}$, and in a frame rotating at inertial rate $\Omega_c$, they satisfy, cf. [1],

$$\begin{cases} \rho_0(D_{\mathbf{u}} + \nabla\mathbf{u} + 2\Omega_c\mathbf{e}_z\times)\,\delta_{\mathbf{v}} + \nabla\delta_p + \delta_\rho\mathbf{Y}_\Omega = \mathbf{s}, \\ \delta_{\mathbf{v}} = (D_{\mathbf{u}} - \nabla\mathbf{u})\,\boldsymbol{\xi}, \qquad \delta_\rho = -\nabla\cdot(\rho_0\boldsymbol{\xi}), \\ \delta_p = -\rho_0 c_0^2\nabla\cdot\boldsymbol{\xi} - (\boldsymbol{\xi}\cdot\nabla)p_0, \end{cases}$$

where $\nabla$ is the gradient, flow $\mathbf{u}$ includes differential rotation given by $\mathbf{u}_0 := (\Omega - \Omega_c)r\sin\theta\mathbf{e}_\phi$ in the current frame, $\mathbf{Y}_\Omega = \nabla\phi_0 - \Omega^2\eta\mathbf{e}_\eta$, and $D_{\mathbf{u}} = -\mathrm{i}\omega + \mathbf{u}\cdot\nabla$ the time-harmonic material derivative along $\mathbf{u}$. The value of $\Omega$ is obtained from helioseismic inversions, while $\rho_0$ and $c_0$ follow a standard solar model. Attenuation is added with a term $(-\mathrm{i}\rho_0\gamma_{\text{att}}\boldsymbol{\xi})$, $\gamma_{\text{att}} > 0$ to the system which can be written as, $\mathcal{L}\boldsymbol{\xi} = \mathbf{s}$. The response to a point source $\delta(\cdot - \mathbf{s})\mathbb{I}$ and a boundary condition is written as the tensor $\mathbb{G}$.

When $\mathbf{u}$ only contains differential rotation, i.e. $\mathbf{u} = \mathbf{u}_0$, the operator is written as $\mathcal{L}_0$. In this case, the first-order system is solved in [1], with the Hybridizable Discontinuous Galerkin (HDG) method, for $\boldsymbol{\xi}_m(\eta, z)$ in the cylindrical expansion of $\boldsymbol{\xi}(\eta, z, \phi) = \sum_m \boldsymbol{\xi}_m(\eta, z)\mathrm{e}^{\mathrm{i}m\phi}$. The method reduces the discretized problem to one on the mesh skeleton in terms of the trace of $\delta_p$. Special choices of stabilization is required due to the Tricomi-property of the PDE, cf. [1].

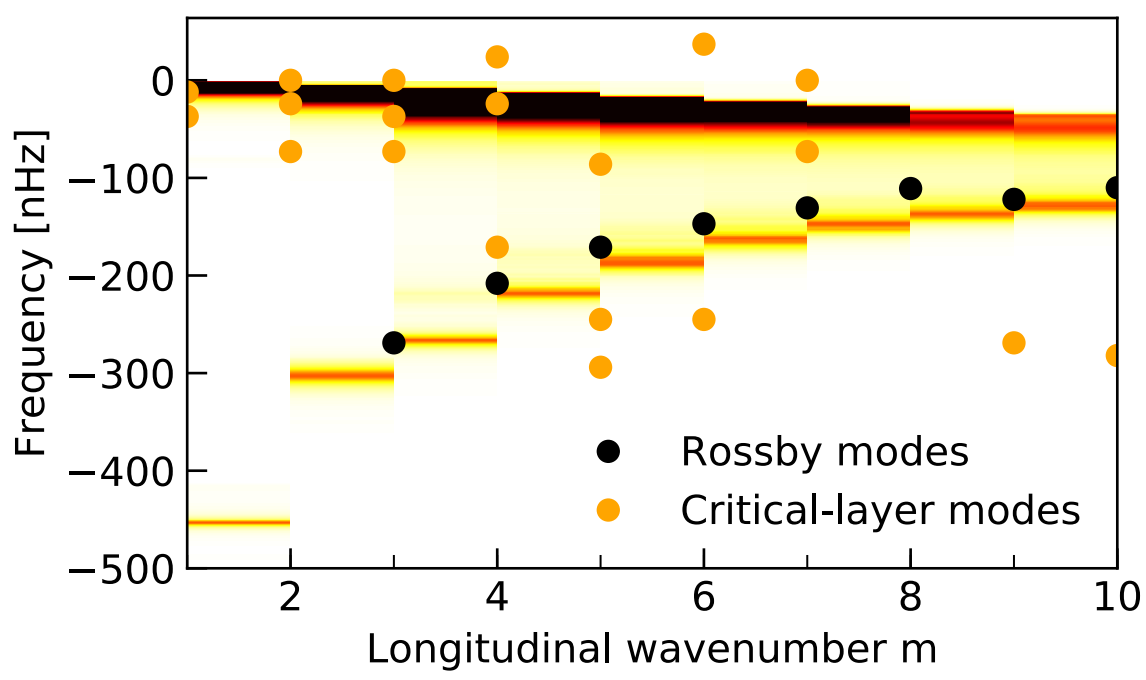

Figure 1: Comparison between the observed frequencies from [4] with the power of the Green kernel $\int_{\mathbb{S}^2} |\mathbb{G}^{yz}_{m,\omega}(||\mathbf{x}|| = 1; \mathbf{s})|^2 \mathrm{d}S$ for a source $\mathbf{s}$ at the equator and 100 km below the surface.

## 3 Spectrum in inertial range.

We consider the Green's kernel $\mathbb{G}_{m,\omega}(\eta, z; \mathbf{s})$ whose elements $\mathbb{G}^{ij}$ corresponds to the $i$-th component of the response to a point source $\delta(.-\mathbf{s})\mathbf{e}_j$. Fig. 1 shows the norm of $\mathbb{G}^{yz}_{m,\omega}(\eta, z; \mathbf{s})$ on $||\mathbf{x}|| = 1$ as a function of $(m, \omega)$ for $\mathbf{s}$ located at the equator at 100 km below the solar surface. It exhibits strong power for certain $(m, \omega)$, particularly those close to the dispersion relation of the observed equatorial Rossby modes from [4]. The Green's kernel for $m = 3$ and $\omega/2\pi = -270$ nHz shown in Fig. 2 also features reasonable agreement with observations in [5]. Our linear adiabatic oscillation model with simplified attenuation is thus able to produce features comparable to the observations in the inertial range.

## 4 Inversion setup

We minimize the cost function

$$\mathcal{J}(\mathbf{u}) = \tfrac{1}{2}||\mathcal{C}_\mathbf{u} - \mathcal{C}_\mathrm{obs}||^2 \,, \qquad (1)$$

for $\mathbf{u}$, fitting between the synthetic and observed CCF. Fig. 3 compares between $\mathcal{C}_{\mathbf{u}_0}$ modeled with differential rotation $\mathbf{u}_0$ and $\mathcal{C}_\mathrm{obs}$ stacked from 10 years of daily acoustic CCFs. As shown there, our forward model can reproduce the main features (skips) of the observed acoustic modes.

We write $\boldsymbol{\xi} = \boldsymbol{\xi}_0 + \boldsymbol{\xi}_\delta$ where $\boldsymbol{\xi}_0$ solves $\mathcal{L}_0\boldsymbol{\xi}_0 = \mathbf{s}$, and only the term $\boldsymbol{\xi}_\delta$ contains the effect of the flow $\mathbf{u}_m(\eta, z)e^{\mathrm{i}(m\phi - \omega_m t)}$. We next decompose $\mathcal{L} = \mathcal{L}_0 + \mathcal{L}_\delta$ to obtain, in first order, $\mathcal{L}_0\boldsymbol{\xi}_\delta = -\mathcal{L}_\delta\boldsymbol{\xi}_0$ where the operator $\mathcal{L}_\delta$ is the perturbation of the forward operator due to $\mathbf{u}_m$. We then follow an iterative minimization approach for the quantitative reconstruction of $\mathbf{u}_m$, using the adjoint-state method to obtain the gradient of $\mathcal{J}$, in the framework of the `Hawen` code [6].

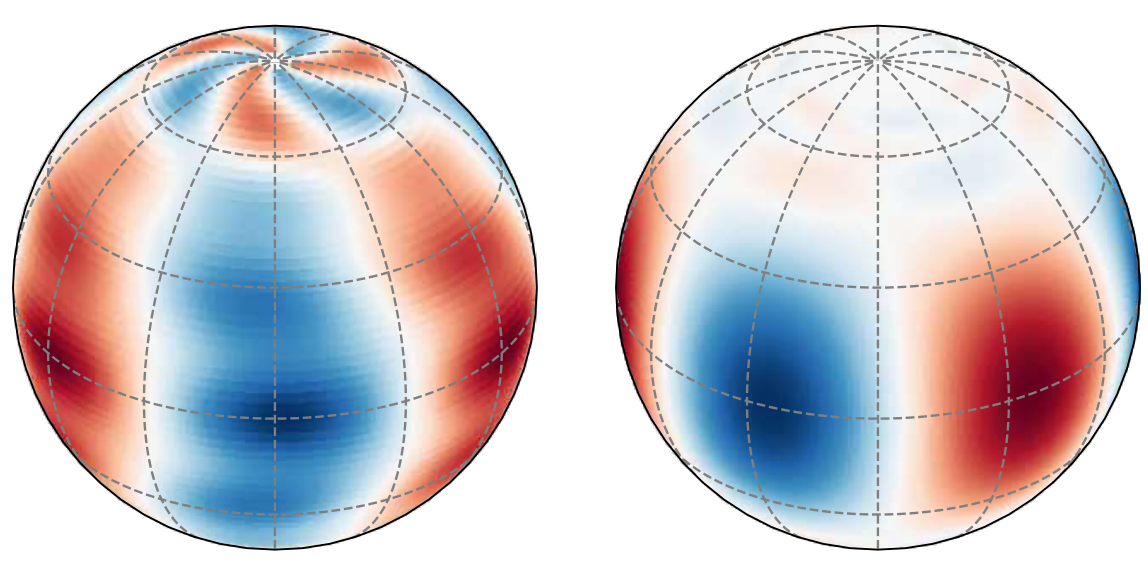

Figure 2: Rossby mode at longitudinal wavenumber $m = 3$. Left: Observations [5]. Right: Green kernel $\mathbb{G}^{yz}_{m,\omega}(\mathbf{x}; \mathbf{s})$ for $m = 3$, $\omega/2\pi = -270$ nHz, $||\mathbf{x}|| = 1$ and $\mathbf{s}$ located at the equator, 100 km below the surface.

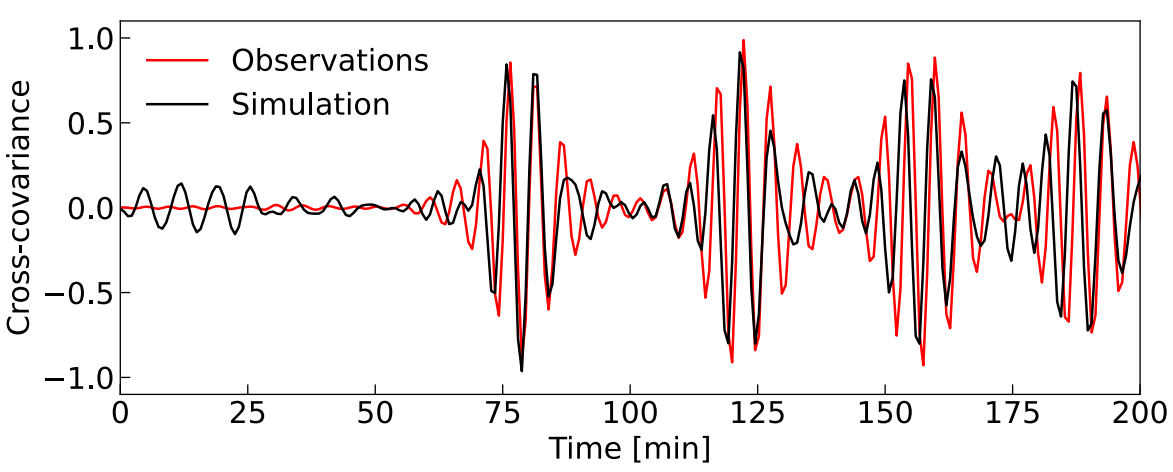

Figure 3: Reference CCF (synthetic with $\mathcal{L}_0$ vs. 10-year stack) for two points separated by 24°.

# A Trefftz Continuous Galerkin method for Helmholtz problems

Nicola Galante[1,*], Bruno Després[2], Emile Parolin[1]
[1]Sorbonne Université, Université Paris Cité, CNRS, INRIA, Laboratoire Jacques-Louis Lions, LJLL, EPC ALPINES, Paris, F-75005, France
[2]Sorbonne Université, Université Paris Cité, CNRS, INRIA, Laboratoire Jacques-Louis Lions, LJLL, Laboratoire de Probabilités, Statistique et Modélisation, EPC MEGAVOLT, Paris, F-75005, France
*Email: nicola.galante@inria.fr

**Abstract**

We present a novel Trefftz continuous Galerkin method for 2D homogeneous Helmholtz problems. Unlike standard discontinuous Trefftz approaches, the method is based on a globally $H^1$-conforming discrete space. We analyze its approximation properties, and show that the resulting discrete representations are stable with bounded coefficients. The best-approximation error is proved to decay algebraically with respect to the number of degrees of freedom for solutions with finite Sobolev regularity, and exponentially for analytic solutions.

## 1 Introduction

We study the approximation properties of a new $H^1$-conforming *Trefftz space* for the Helmholtz equation with wavenumber $\kappa > 0$, i.e. $-\Delta u - \kappa^2 u = 0$, posed in a bounded Lipschitz domain $\Omega \subset \mathbb{R}^2$. The basis functions are defined so that their traces coincide with Laplace–Dirichlet eigenmodes on a mesh edge; they are globally continuous, compactly supported, and admit explicit representations as linear combinations of four plane waves within each element. This ensures they satisfy the Helmholtz equation locally, and allows the system matrix to be assembled in closed form for polygonal domains. The basis functions are either propagative, for small mode numbers relative to $\kappa$, or evanescent for larger ones. Enriching the approximation set with evanescent components corresponding to arbitrary high mode numbers is essential to achieve accurate, bounded-coefficient approximations, and underlies the spectral convergence of the Trefftz space for smooth enough solutions.

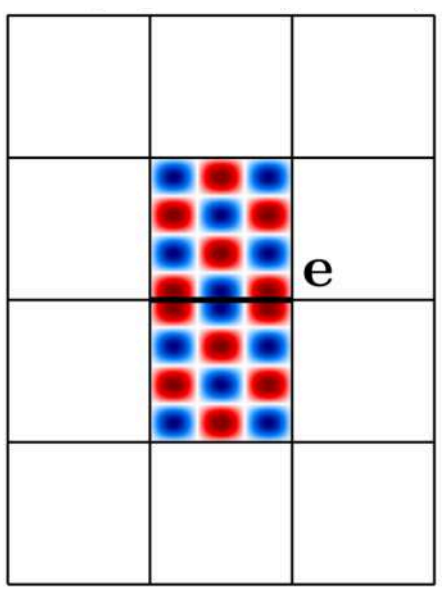

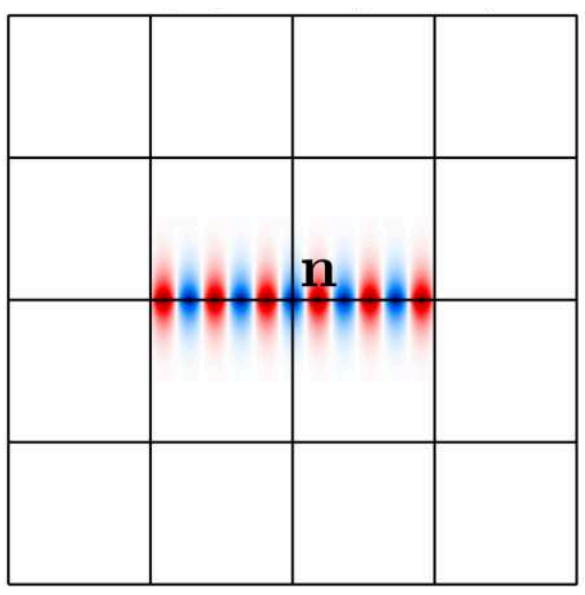

Figure 1: (left) propagative edge basis function, (right) evanescent node basis function.

## 2 A new conforming Trefftz space

The Trefftz space is built on a mesh obtained by intersecting the domain $\Omega$ with a Cartesian grid $\mathcal{T}_h$ of square cells with side length $h > 0$. For any $n \in \mathbb{N}^*$ and $(x, y) \in (0, h)^2$, we define

$$\varphi_n(x,y) := \sin\left(\kappa\nu_n x\right) \frac{\sin(y\kappa\sqrt{1-\nu_n^2})}{\sin(h\kappa\sqrt{1-\nu_n^2})}, \quad (1)$$

where $\nu_n := \frac{n\pi}{\kappa h}$ and $\sqrt{\cdot}$ denotes the principal square root. These functions are simple linear combinations of plane waves with trace supported on the top edge $(0, h) \times \{h\}$. They are well-defined provided that $\kappa^2$ is not a Neumann–Laplace eigenvalue of a cell of $\mathcal{T}_h$.

From these modes, two families of global basis functions are constructed, associated with edges $\mathbf{e}$ and nodes $\mathbf{n}$ of $\mathcal{T}_h$. For any edge $\mathbf{e}$, the *edge functions* are supported on the two cells adjacent to $\mathbf{e}$, and extended by zero elsewhere. On these cells, each function coincides with (1), suitably oriented so that its nonzero trace lies along $\mathbf{e}$; this choice enforces continuity across $\mathbf{e}$, and hence $\Omega$. *Node functions* associated with a node $\mathbf{n}$ are defined analogously on coarser mesh obtained by doubling the cell size in one coordinate direction; their construction similarly requires that $\kappa^2$ is not a Neumann–Laplace eigenvalue for the associated coarse cells. See Figure 1 for an illustration and [1] for details.

The space spanned by $N_e$ edge and $N_n$ node functions is denoted by $V_{\mathbf{N}}$, with $\mathbf{N} = (N_e, N_n)$. By construction, $V_{\mathbf{N}}$ is $H^1$-conforming, and thus suitable for continuous Galerkin discretizations.

## 3 Theoretical results

Best-approximation error estimates are established for the space $V_{\mathbf{N}}$ on domains $\Omega$ tessellated by Cartesian grids $\mathcal{T}_h$. The main result in [1] is:

**Theorem 1** *Let $u \in H^{M+1}(\Omega)$ be a homogeneous Helmholtz solution for some $M \in \mathbb{N}^*$. For any $N_n \in \mathbb{N}^*$, there exists a constant $C > 0$ independent of $N_e$, such that, for any $N_e \gtrsim \kappa h$,*

$$\inf_{v \in V_{\mathbf{N}}} \|u - v\|_{H^1(\Omega)} \leq C\, N_e^{1/2-\mu} \|u\|_{H^{\mu+1}(\Omega)},$$

*where $\mu = \min(2N_n, M)$. Besides, if $u$ is analytic, there exist constants $C, \eta > 0$ independent of $N_e$ and $N_n$, such that, for any $N_e > N_n \gtrsim \kappa h$,*

$$\inf_{v \in V_{\mathbf{N}}} \|u - v\|_{H^1(\Omega)} \leq C\, N_n^{7/2} e^{-\eta N_n} \|u\|_{H^{2N_n+1}(\Omega)}.$$

The dependence of the constants on the wavenumber can be made explicit using $\kappa$-dependent norms; see [1]. The approximation properties of $V_{\mathbf{N}}$ rely on a complementary interplay between edge and node functions. The Trefftz error is controlled by its trace on the mesh skeleton: expanding this trace in the edge Laplace–Dirichlet eigenbasis and applying repeated integration by parts yields increasing algebraic powers of $N_e$, provided the solution is regular enough, but also produces residual terms at edge endpoints. Node functions cancel $N_n$ such nodal contributions, while preserving local Trefftz field regularity.

The analysis in [1] also provides stability estimates, proving that the coefficient norms in the discrete representation are uniformly bounded in $N_e$ and grow at most quadratically in $N_n$. It further addresses the high-frequency regime, showing that, for any fixed mesh $\mathcal{T}_h$, there exists an unbounded sequence of wavenumbers along which the best-approximation error decays exponentially, provided that $N_e > N_n \gtrsim \kappa$.

## 4 Numerical results

We conclude with numerical experiments illustrating the approximation properties of $V_{\mathbf{N}}$. The domain considered is $\Omega = (0,1) \times (-1/2, 1/2)$ with wavenumber $\kappa = 30$. As target functions, circular waves $(r, \theta) \mapsto J_\nu(r)e^{\imath\nu\theta}$ are taken for

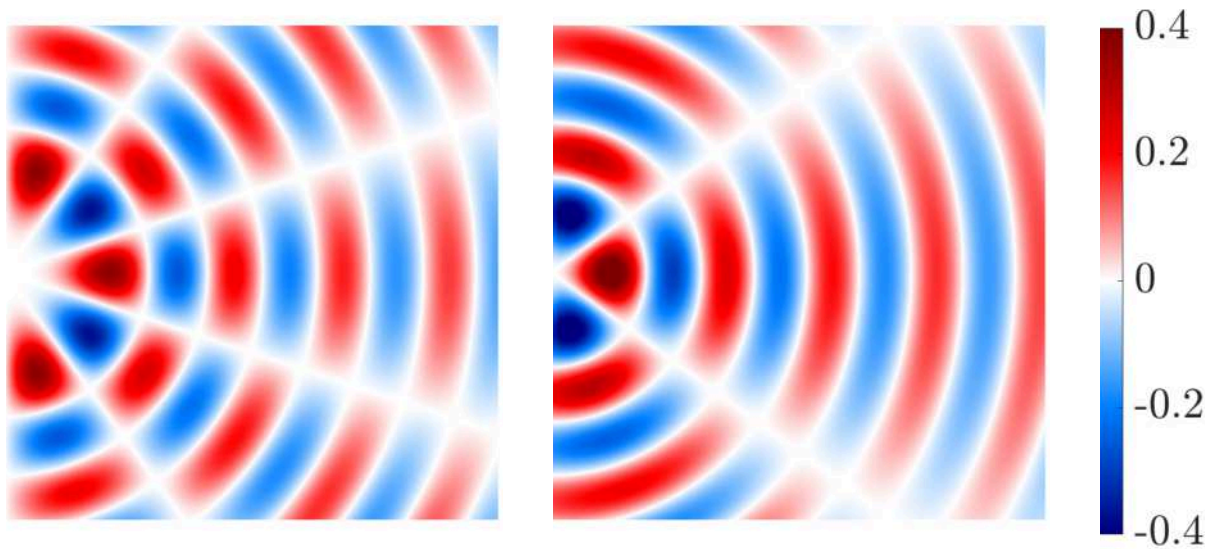

Figure 2: Real part: (left) $\nu = 5$, (right) $\nu = 5/2$.

$\nu = 5$, yielding an analytic solution, and for $\nu = 5/2$, for which the solution has limited Sobolev regularity, as it does not belong to $H^{1+5/2}(\Omega)$; see Figure 2. Due to the redundancy of $V_{\mathbf{N}}$, a least-squares Petrov–Galerkin formulation with regularized SVD stabilization is employed. We consider a mesh $\mathcal{T}_h$ of four square cells with impedance boundary conditions. Figure 3 confirms the rates predicted by Theorem 1 and suggests that the first estimate may be sharpened by a factor of 1/2 in the exponent of $N_e$.

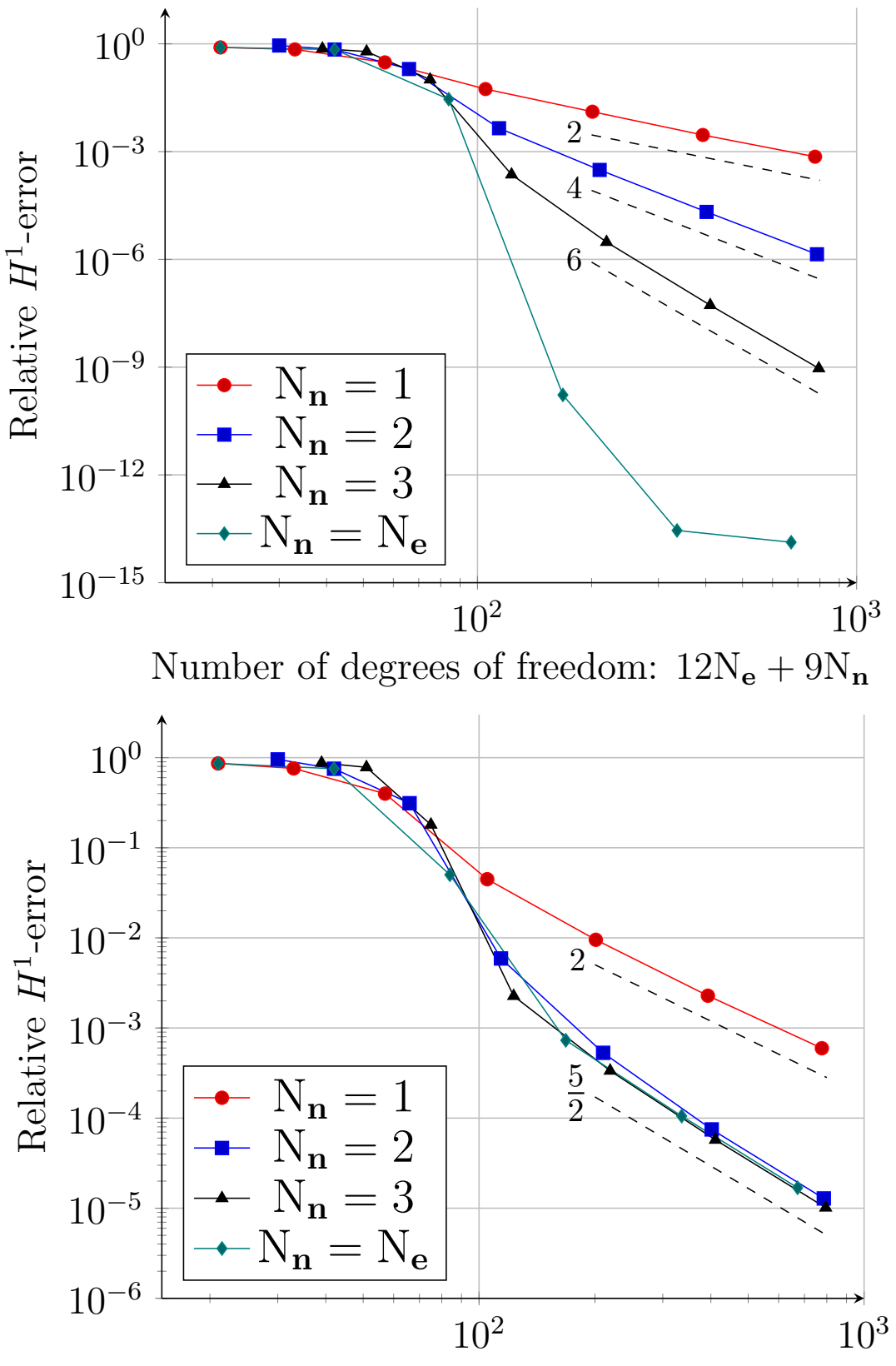

Figure 3: $H^1$-error: (up) $\nu = 5$, (bottom) $\nu = 5/2$.

# A Multi-domain/Multi-method Preconditioner for the Electric Field Integral Equation

**Timothée Galtier**[1,*], **David Levadoux**[1]
[1]DTIS, ONERA, Université de Toulouse, 31000, Toulouse, France
*Email: timothee.galtier@onera.fr

**Abstract**

Wave diffraction problems are commonly formulated as boundary integral equations and solved via the Boundary Element Method (BEM). However, BEM typically leads to large, dense, and ill-conditioned linear systems. While existing preconditioning techniques improve convergence for simplified geometries, industrial aeronautical applications involve complex assemblies (landing gears, air ducts, ...) that necessitate tailored numerical treatments. To address this challenge, we propose a multi-domain/multi-method preconditioner based on an overlapping domain decomposition approach.

## 1 Introduction

In electromagnetism, wave diffraction problems can be cast into the well-known Electric Field Integral Equation (EFIE). While widely used in industry for its accuracy, the EFIE is notoriously ill-conditioned. Various approaches, such as Calderón preconditioners [3], multi-trace integral formulations [2], or the Generalized Combined Sources Integral Equation formalism [1], have been proposed; however, most are tailored for specific geometric configurations. Furthermore, when applied to the entire surface of complex scattering objects, these methods often suffer from reduced efficiency. To address this, we develop a piecewise preconditioner specifically designed for robust performance in large-scale, real-world applications. While our preconditioning method is related to an overlapping Schwarz method through the use of a partition of unity, it differs in that the locally solved problem is not exactly the global one, but is instead intended to synthesize a parametrix.

## 2 A body decomposition preconditioner

The EFIE operator is physically interpreted as an impedance operator $\mathrm{Z} : \mathbf{n} \times \mathbf{H} \mapsto \mathbf{E}_{\tan}$, linking a surface electric current injection to the resulting tangential electric field. To precondition the EFIE, we propose constructing a piecewise approximation $\widetilde{\mathrm{Y}}$ of the admittance operator $\mathrm{Y} = \mathrm{Z}^{-1}$. This method involves partitioning the scattering object into distinct components, or 'bodies', for which a local parametrix for Z can be efficiently obtained by analytical or numerical means. To implement this, the initial surface $\Gamma$ is decomposed into $N$ bodies $\{\Gamma_n\}_{n=1}^N$ and a neighborhood $\Gamma_s$ of $\bigcup_n \partial\Gamma_n$, referred to as the skeleton of $\Gamma$. A quadratic partition of unity, $\{(\Gamma_n, \chi_n)\}_{1\leqslant n\leqslant N} \cup \{(\Gamma_s, \chi_s)\}$, is introduced to ensure a smooth transition between domains. For each body $\Gamma_n$, we define a simplified surrogate geometry $\widetilde{\Gamma}_n$ and a mapping $\psi_n : \Gamma_n \to \widetilde{\Gamma}_n$. This mapping induces a Piola transformation, which maps surface currents from the original geometry to the surrogate one. Finally, the global preconditioner is defined as:

$$\widetilde{\mathrm{Y}} = \chi_s Y_s \chi_s + \sum_n \chi_n \psi_n^{-1} \widetilde{\mathrm{Y}}_n \psi_n \chi_n, \qquad (1)$$

where $Y_s$ and $\widetilde{\mathrm{Y}}_n$ represent local approximations of the admittance for the skeleton and the surrogate bodies, respectively.

We present an example of application of the Body Decomposition Method in Fig. 1 :

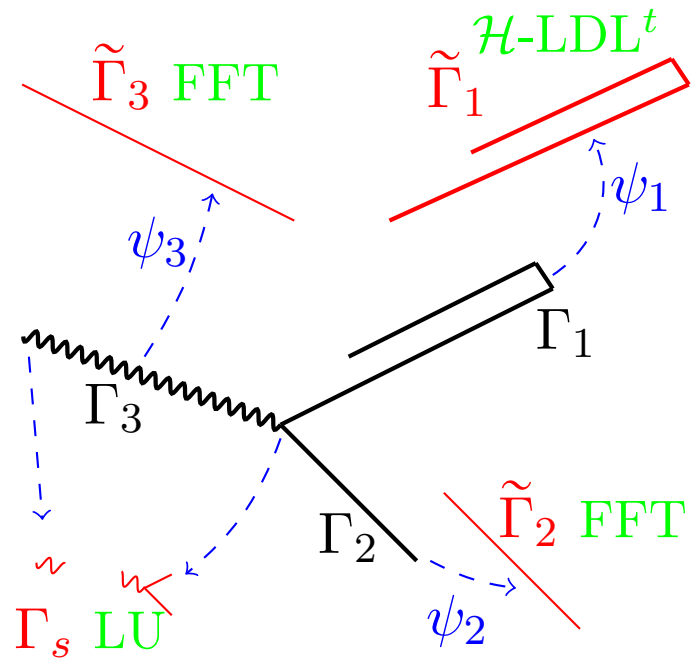

Figure 1: A BDM preconditioner example with local preconditioner solvers (in green).

The local operators $\widetilde{\mathrm{Y}}_n$ and $Y_s$ are chosen according to the geometric nature of each subdomain. Specifically, we consider the following preconditioning techniques:

- FFT refers to FFT-based operators involving the convolution of elementary solutions (or 'elementary currents'), efficiently computed using the Fast Fourier Transform [4].
- $\mathcal{H}$-LDL$^t$ denotes a hierarchical preconditioner obtained by an LDL$^t$ decomposition of Z within the $\mathcal{H}$-matrix framework [5].

It is important to note that the skeleton $\Gamma_s$ and its associated cut-off function $\chi_s$ play a pivotal role in the BDM. Their inclusion ensures that the global preconditioner maintains the correct principal symbol-specifically, that of the admittance operator Y. This consistency is essential for guaranteeing the robustness of the preconditioner across complex geometries.

## 3 A numerical result

We consider a cone-sphere geometry equipped with two wings. The object (in Fig. 2) has an electrical length of approximately $25\lambda$ and is discretized with $273,000$ degrees of freedom.

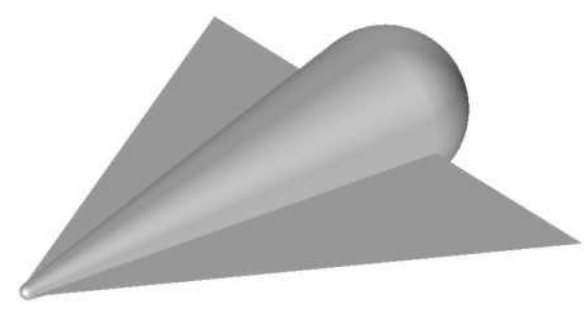

Figure 2: Scattering object

The EFIE is solved using a Generalized Minimal Residual (GMRES) iterative solver, where matrix-vector products are accelerated via the Fast Multipole Method (FMM).

To evaluate the performance of our approach, we compare the BDM preconditioner against $\mathrm{PREC}_{\mathrm{iter}}$, the production preconditioner at ONERA, which is based on an iterative inversion of the FMM near-interaction blocks. For the BDM, the object is decomposed into three bodies: the two wings ($\Gamma_1$ and $\Gamma_2$) and the cone-sphere ($\Gamma_3$). Given their flat geometry, the wings are preconditioned using the FFT-based technique. In contrast, due to its complex shape, the cone-sphere ($\Gamma_3$) is treated with an $\mathcal{H}$-LDL$^t$ decomposition.

As illustrated in Fig. 3, the BDM significantly outperforms the production preconditioner, demonstrating its superior efficiency for complex, multi-component geometries.

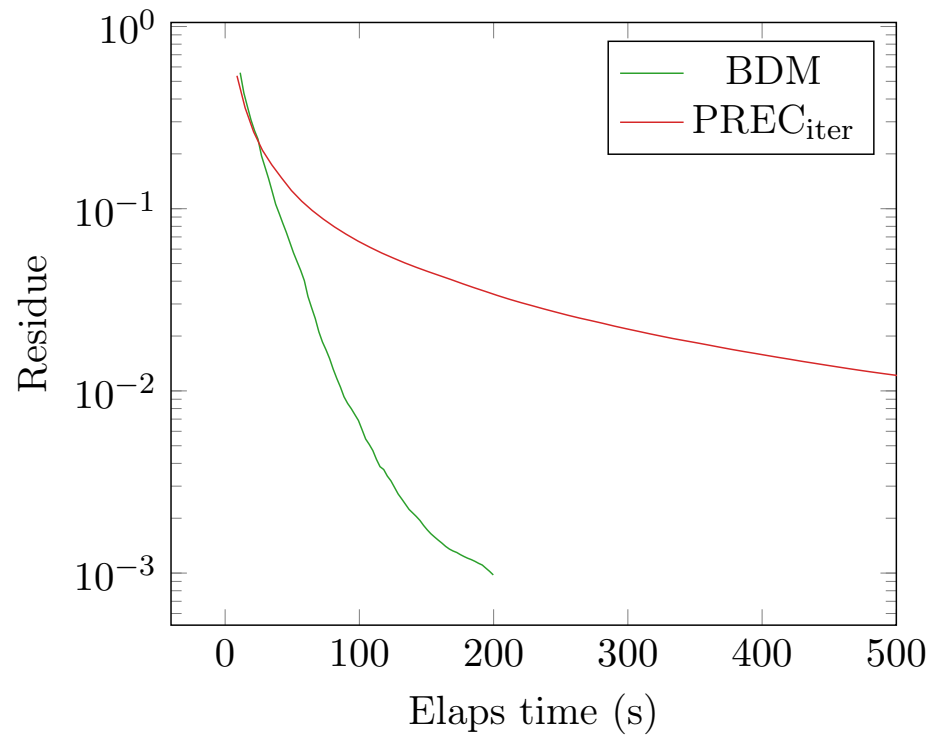

Figure 3: Resolution of the EFIE

The proposed approach is validated through several complex applications, including cavities and intricate plate structures. Furthermore, we investigate optimal BDM configurations (eg. the selection of local solvers) and evaluate the efficiency of the BDM preconditioner against the state of the art.

# Explicit inverse scattering for the one-dimensional Schrödinger equation

**Peter Gibson**[1,*]
[1]Department of Mathematics & Statistics, York University, Toronto, Canada
*Email: pcgibson@yorku.ca

**Abstract**

Working with the one-dimensional Schrödinger equation in impedance form, we derive an exact inverse scattering formula that expresses impedance in terms of the reflection coefficient, and we prove injectivity of the scattering map for impedance functions of lower regularity than previously analyzed. The inverse scattering formula translates directly into an efficient numerical algorithm that accurately transforms digital scattering data into impedance. The results apply to acoustic imaging of layered media, as well as to inverse quantum scattering.

## 1 Introduction

The classical one-dimensional Schrödinger equation comprises one of the most basic and most intensively studied quantum mechanical models. Together with its variants, which are obtained by simple transformations, it provides the building blocks for quantum graphs and models acoustic wave propagation in layered media. Scattering theory for the 1D Schrödinger equation concerns the nonlinear relationship between physical structure, represented by a quantum potential or acoustic impedance, and scattering data, typically a reflection coefficient or scattering matrix. The goal of inverse scattering theory is to elucidate the inverse map going from scattering data to physical structure. Attempts to evaluate the inverse scattering map underlie acoustic imaging of layered media, for instance, including seismic imaging of horizontally stratified sedimentary formations, non-destructive testing of laminated structures in the built environment, and ultrasonic imaging of layered tissue such as skin or the abdominal wall. Famously, the inverse scattering map is also crucial to solving nonlinear equations such as KdV.

Two considerations motivate the work [1] to be presented. Firstly, from a broad theoretical point of view, one would like to formulate the inverse scattering map explicitly, if possible, to reveal the nature of its inherent nonlinearity. We know of no such formulation in the existing literature.

Secondly, from the perspective of acoustic imaging, there is a gap between established mathematical theory concerning inverse scattering and what one would like to know for applications. In particular, it is of interest in acoustics to solve the inverse scattering problem for the Schrödinger equation in impedance form in cases where the impedance function is piecewise absolutely continuous, allowing jump discontinuities. The latter class is of lower regularity than can be handled by known methods.

We present new analytical tools introduced in [2] to address both of these issues. We derive a new formula that expresses impedance explicitly in terms of scattering data. Beyond its theoretical interest, the formula translates readily into a fast and accurate inverse scattering algorithm, applicable to recorded acoustic reflection data. The algorithm is described in detail and turns out to be highly accurate, as illustrated in Figure 1. In addition, we prove injectivity of the forward scattering map for the Schrödinger equation in impedance form in the case where impedance is piecewise absolutely continuous with finitely many jump points, beyond the scope of previous work. The results are primarily motivated by applications to acoustic imaging of layered media, but they apply also to inverse quantum scattering in one spatial dimension. From a purely mathematical perspective, the results illustrate the usefulness for scattering theory of the harmonic exponential operator, a continuum counterpart to orthogonal polynomials that was first described in [2].

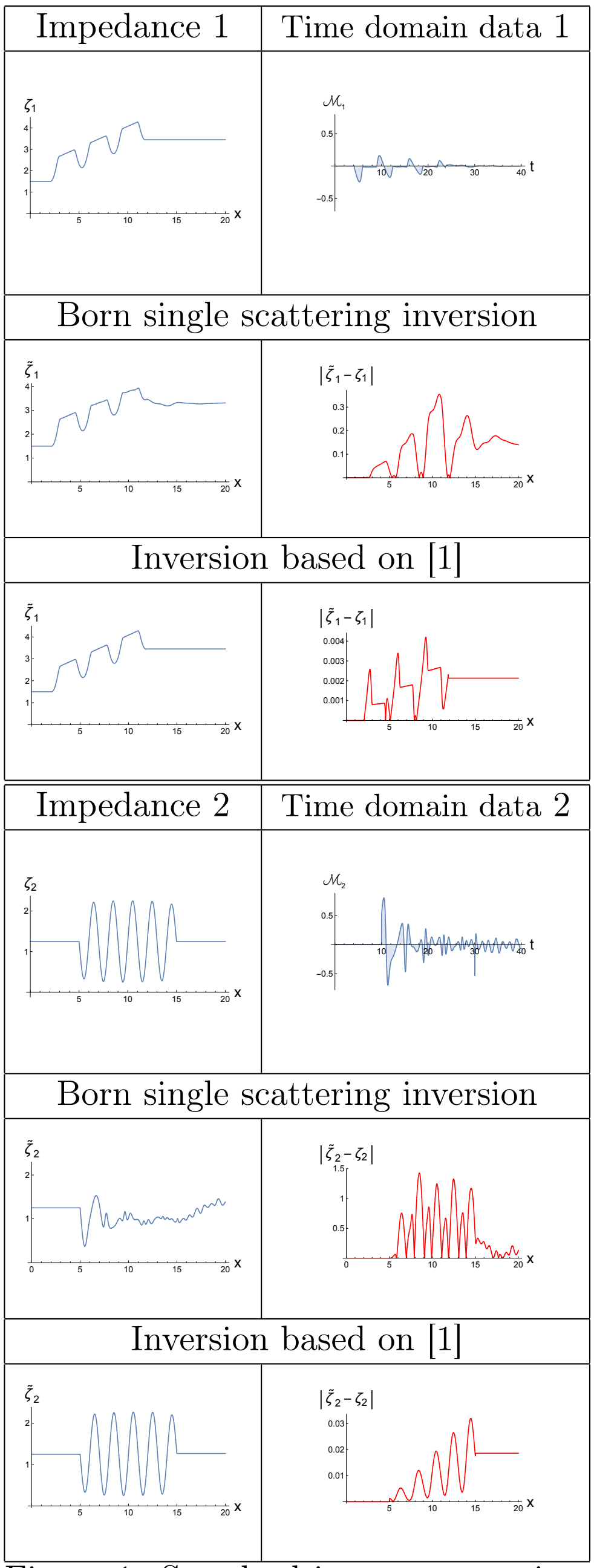

Figure 1: Standard inverse scattering based on the Born approximation compared to the algorithm derived in [1], for two impedance profiles: $\zeta_1$ above, and $\zeta_2$ below. Synthetic time domain data $\mathscr{M}_1(t)$ and $\mathscr{M}_2(t)$ generated at $x_0 = 0$ was sampled on a grid of 10K points for $\zeta_1$ and 20K points for $\zeta_2$. For each example, the computed approximation to impedance $\zeta_j$ is denoted $\tilde{\zeta}_j$, and $|\tilde{\zeta}_j - \zeta_j|$ is the absolute error ($j = 1, 2$).

# A posteriori estimates and adaptive boundary elements for the wave equation

**Heiko Gimperlein**[1,*]
[1]Engineering Mathematics, University of Innsbruck, AUSTRIA
*Email: heiko.gimperlein@uibk.ac.at

**Abstract**

We survey recent work on adaptive mesh refinement procedures for Galerkin and convolution quadrature discretizations of wave equations, formulated as an equivalent boundary integral equation. Reliable a posteriori error estimates of residual type lead to adaptive boundary element methods, based on the four steps: Solve - Estimate - Mark - Refine. We discuss numerical experiments for space-time adaptive Galerkin elements, as well as extensions to convolution quadrature discretizations. The experiments confirm the theoretical results. They study the convergence and computational efficiency of the approaches and the reliability and efficiency of the error estimates.

## 1 Introduction

The efficient numerical solution of boundary integral equations using adaptive mesh refinement procedures has been extensively investigated for homogeneous elliptic problems in unbounded domains. We report on the extension of the a posteriori error estimates and adaptive mesh refinement procedures to the Dirichlet boundary-initial value problem for the wave equation, formulated as a time dependent boundary integral equation.

To be specific, the wave equation

$$\partial_t^2 u - \Delta u = 0 \ , \quad u = 0 \ \text{ for } \ t \leq 0 \ , \qquad (1)$$

is considered in the complement $\mathbb{R}^d \backslash \overline{\Omega}$ of a polyhedral domain or screen, for $d = 2$ or $d = 3$. Dirichlet boundary conditions $u = f$ are imposed on the boundary $\Gamma = \partial\Omega$.

This boundary problem is equivalent to the time dependent boundary integral equation

$$\mathcal{V}\phi(t,x) = f(t,x) \ , \quad (t,x) \in \mathbb{R}^+ \times \Gamma, \qquad (2)$$

involving the weakly singular operator

$$\mathcal{V}\phi(t,x) = \int_{\mathbb{R}^+} \int_\Gamma G(t-\tau,x,y)\phi(\tau,y) \ ds_y \ d\tau.$$

Here, $G$ is a fundamental solution of (1).

The solution $\phi$ of (2) may exhibit singularities due to nonsmooth boundary points (like edges, corners) or nonsmooth data $f$ (like traveling wave fronts), resulting in reduced convergence rates of numerical solutions. Optimal convergence can be recovered for discretizations that are locally refined in space, time or space-time. We here review recent work on a posteriori estimates and the resulting adaptive mesh refinements in space and time, for Galerkin [1, 3] and convolution quadrature discretizations [2] of (2).

## 2 Space-time adaptive boundary elements

An a posteriori estimate for Galerkin discretizations was obtained in [3], in appropriate space-time Sobolev spaces $H^r_\sigma(\mathbb{R}^+, \widetilde{H}^s(\Gamma))$ with weight $\sigma > 0$ (cf. Def. 2.1 in [3]). For $r = s = 0$ they correspond to the weighted $L^2$-space with scalar product $\int_0^\infty e^{-2\sigma t} \int_\Gamma u\overline{v} \, d\Gamma_{\mathbf{x}} \, dt$.

**Theorem A:** *Let $\phi_{h,\Delta t}$ be the solution to a space-time Galerkin discretization of* (2)*, and assume $\mathcal{R} = \partial_t f - \mathcal{V}\partial_t\phi_{h,\Delta t} \in H^0_\sigma(\mathbb{R}^+, H^1(\Gamma))$. Then*

$$\|\phi - \phi_h\|^2_{H^0_\sigma(\mathbb{R}^+, \widetilde{H}^{-\frac{1}{2}}(\Gamma))} \lesssim \sum_{\square} \eta^2_\square \ ,$$

*Here, $\eta^2_\square = h_\square \left( \|\nabla\mathcal{R}\|^2_{L^2(\square)} + \|\partial_t\mathcal{R}\|^2_{L^2(\square)} \right)$, with $\square$ an element (of diameter $h_\square$) of the space-time mesh $\mathbb{R}^+ \times \Gamma = \bigcup \square_{i,\Delta}$, $\square_{i,\Delta} = [t_i, t_{i+1}) \times \Delta$.*

It leads to an adaptive mesh refinement procedure, based on the four steps:

**SOLVE → ESTIMATE → MARK → REFINE** :

**Space–Time Adaptive Algorithm:**
Input: Datum $f$, space-time mesh $\mathcal{T}_0$, refinement parameter $\Theta \in (0,1)$, exit tolerance $\epsilon > 0$, $k = 0$.

1. Solve (2) on $\mathcal{T}_k$.
2. Compute $\eta^2_\square$ in each space-time element $\square \in \mathcal{T}_k$

3. Stop if $\sum_{\square} \eta_{\square}^2 < \epsilon$.
4. Mark all $\square \in \mathcal{T}_k$ with $\eta_{\square}^2 > \Theta\, \eta_{max}^2$.
5. Refine each marked $\square$ to obtain a refined mesh $\mathcal{T}_{k+1}$.
6. Set $k \leftarrow k+1$ and go to 1.

The resulting mesh refinements in space and time were studied in [1].

**Example1:** For a crack $\Gamma = [0,1] \subset \mathbb{R}^2$ and peak-shaped smooth data $f$, Fig. 1 depicts the numerical solution on the adaptively generated mesh, with $\epsilon = 10^{-5}$, $\Theta = 0.5$. Fig. 2 shows the square of the energy error on a sequence of uniform, resp. adaptively generated meshes. Because the solution is smooth, the asymptotic convergence rate of the adaptive and uniform refinements coincide, but the error on an adaptive mesh is a factor $\sim 10$ lower than on a uniform mesh with the same number of DOF. The resulting memory savings are discussed in [1].

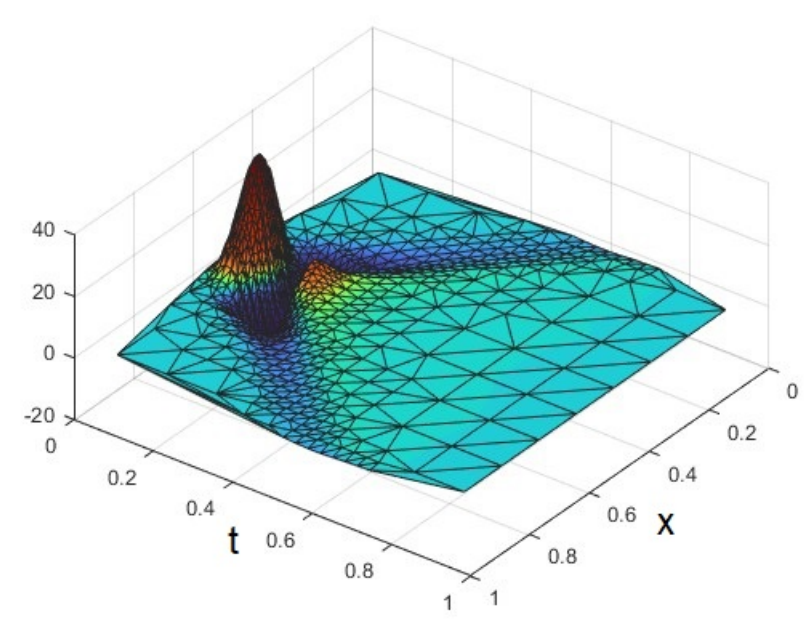

Figure 1: Space-time mesh and solution $\phi_h$ [1]

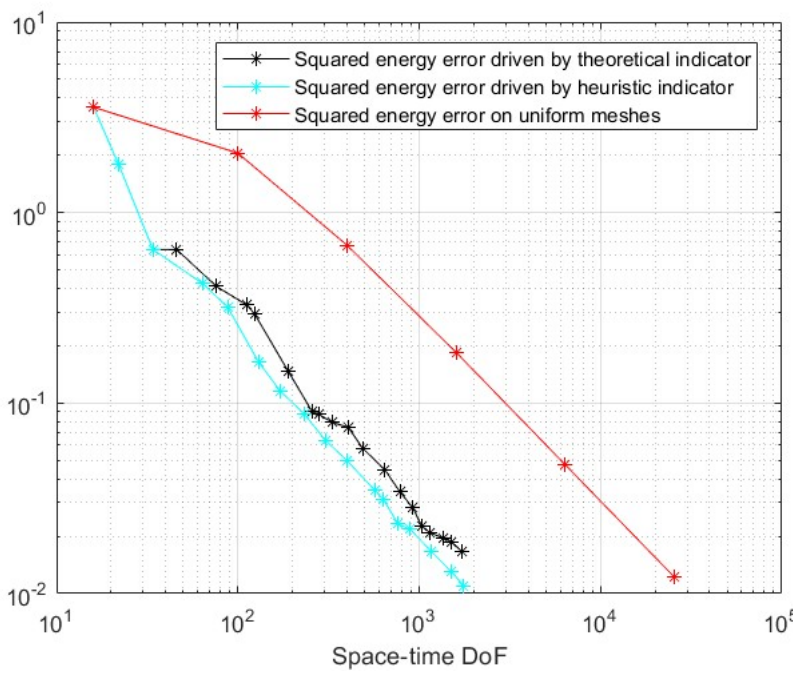

Figure 2: (Squared) energy error on uniform and adaptive space-time meshes [1]

## 3 Space adaptive convolution quadrature

An a posteriori estimate similar to Theorem A for discretizations of (2) by convolution quadrature was obtained in [2]. Such discretizations do not easily allow adaptive mesh refinements in both space and time. Still, [2] showed the relevance of adaptive refinements in space to resolve geometric singularities of the solution.

**Example 2:** For a trapping crack $\Gamma$ and an incident plane wave $f$, Fig. 1 depicts a snapshot of the numerical solution in $\mathbb{R}^2 \setminus \Gamma$. Fig. 4 shows the energy error on a sequence of uniform, $\beta$-graded, resp. adaptively generated meshes. For fixed time step 0.1 the adaptive method converges with the optimal rate 1.5 wrt. the number $M$ of spatial DOF.

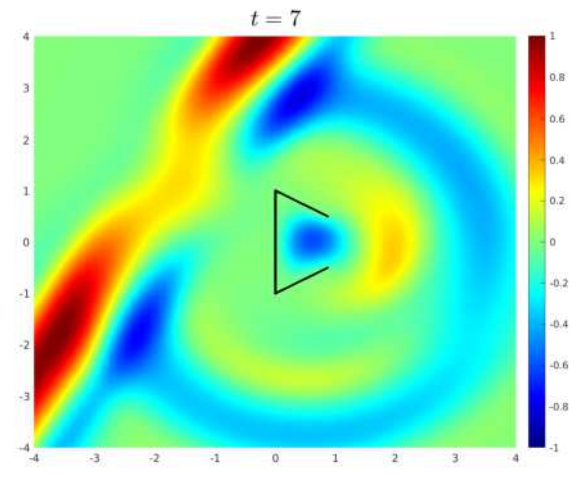

Figure 3: Solution $u_h$ obtained from $\phi_h$ [2]

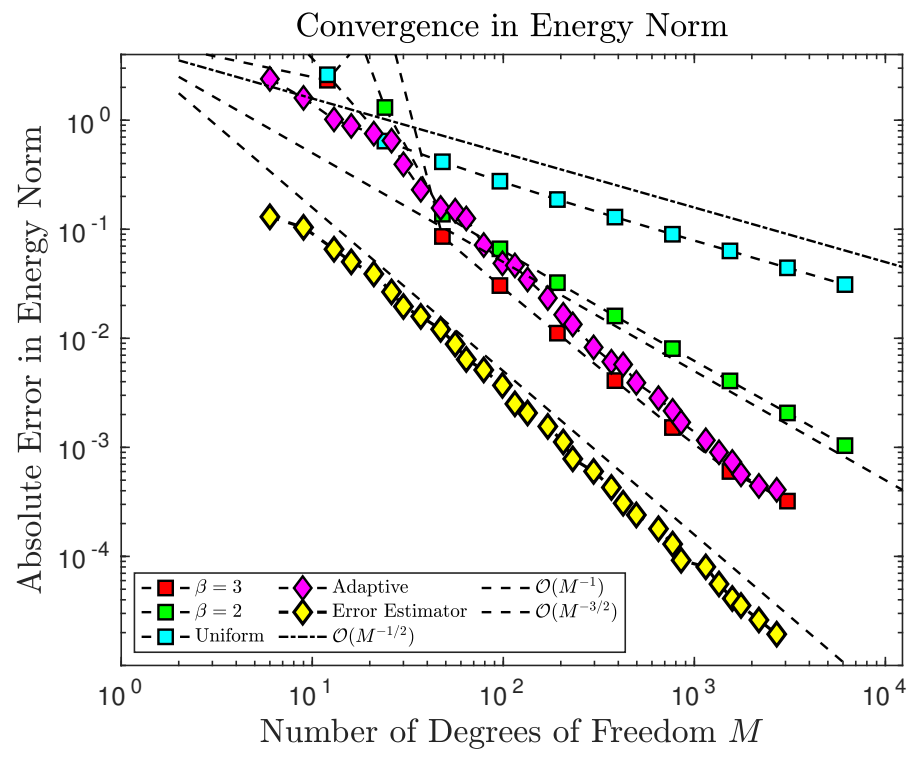

Figure 4: Energy error on uniform, graded and adaptive spatial meshes [2]

# High-order stochastic homogenization of the Helmholtz transmission problem

**Laure Giovangigli**[1,*]
[1]POEMS, CNRS, Inria, ENSTA, Institut Polytechnique de Paris, 91120 Palaiseau, France
*Email: laure.giovangigli@ensta.fr

## Abstract

We consider the time-harmonic scattering of an incident wave by a bounded medium with an heterogeneous random micro-structure. Assuming that the acoustic parameters charaterizing the material are stationary and verify a quantitative mixing assumption, we provide via stochastic homogenization techniques quantitative two-scale expansions of the solution inside and outside the heterogeneous medium that approximate at different orders the true solution.

## 1 Introduction

We are interested in characterizing the field scattered by a bounded heterogeneous object with a random micro-structure embedded in a homogeneous infinite background medium. We derive quantitative two-scale expansions of the field inside and outside the object. In particular we estimate the boundary corrector for the Helmholtz transmission problem in the stochastic setting. Finally we provide a quantitative point-wise expansion with a fluctuation scaling and compare the different expansions. We emphasize that all the error estimates are optimal in terms of scaling with respect to $\varepsilon$.

## 2 Presentation of the problem

Let $D$ be a bounded domain of $\mathbb{R}^d$, $d = 2, 3$. We study the time-harmonic scattering of an incident wave $u^i$ with wavenumber $k$ by a medium occupying $D$. We suppose that the acoustic parameters $a_\varepsilon$ and $n_\varepsilon$ in $D$ can be written as

$$\forall x \in D,\ a_\varepsilon(x) := a\left(\frac{x}{\varepsilon}\right), \text{ and } n_\varepsilon(x) := n\left(\frac{x}{\varepsilon}\right),$$

where $\varepsilon > 0$, $a$ and $n$ are ergodic stationary processes, $a$ is uniformly elliptic and $n$ bounded by positive constants. $\mathbb{R}^d \setminus \overline{D}$ is supposed to be homogeneous with parameters $I$ and $n_0$.

Let $B_R \supset D$ be a ball with radius $R$. The

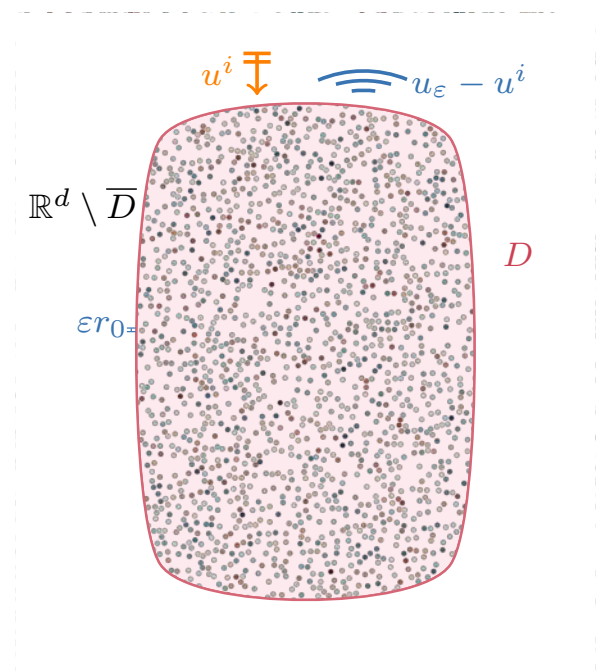

Figure 1: Illustration of the problem

total field $u_\varepsilon$ is the solution in $H^1(B_R)$ of

$$\left|\begin{array}{ll} -\nabla\cdot\left(I\mathbb{1}_{B_R\setminus\overline{D}} + a_\varepsilon\mathbb{1}_D\right)\nabla u_\varepsilon & \\ \quad -k^2\left(n_0\mathbb{1}_{B_R\setminus\overline{D}} + n_\varepsilon\mathbb{1}_D\right)u_\varepsilon = 0 & \text{in } B_R, \\ \partial_\nu(u_\varepsilon - u^i) = \Lambda(u_\varepsilon - u^i) & \text{on } \partial B_R, \end{array}\right. \tag{1}$$

with $\Lambda$ the Dirichlet-to-Neumann operator associated to the exterior Helmholtz problem.

$\mathbb{P}$-a.s. $u_\varepsilon$ converges weakly in $H^1(B_R)$ to $u_0$ the unique solution to

$$\left|\begin{array}{ll} -\nabla\cdot\left(I\mathbb{1}_{B_R\setminus\overline{D}} + a^\star\mathbb{1}_D\right)\nabla u_0 & \\ \quad -k^2\left(n_0\mathbb{1}_{B_R\setminus\overline{D}} + n^\star\mathbb{1}_D\right)u_0 = 0 & \text{in } B_R, \\ \partial_\nu(u_0 - u^i) = \Lambda(u_0 - u^i) & \text{on } \partial B_R, \end{array}\right. \tag{2}$$

where the effective coefficients in $D$ are

$$\forall i,j \in [\![1,d]\!],\quad a^\star_{ij} = \mathbb{E}[e_j \cdot A(e_i + \nabla\phi_i)],$$
$$n^\star = \mathbb{E}[n].$$

Here $(\phi_i)_{i\in[\![1,d]\!]}$ is the standard corrector in stochastic homogenization [2]. We define the first order correction in $D$ $u_1^\varepsilon(\cdot) := \sum_{i=1}^d \phi_i(\frac{\cdot}{\varepsilon})\partial_i u_0(\cdot)$.

Assuming that $a$ and $n$ statisfy a quantitative mixing assumption (Mix) in form of a multiscale weighted variance inequality [1, 2], we use the corrector estimates derived in [2, 3] to establish quantitative two-scale expansions inside and outside $D$. The difficulty is to characterize the error near the boundary.

## 3 Error estimates for the two-scale expansions

We suppose in what follows that $u_{0|_D} \in W^{2,\infty}(D)$ and $(a, n)$ verifies (Mix). The first step is to introduce $v_{1,\varepsilon}$ such that $u_\varepsilon - w_{1,\varepsilon}$ verifies the transmission conditions across $\partial D$, where $w_{1,\varepsilon} := u_0 - \varepsilon u_1^\varepsilon \mathbb{1}_D - v_{1,\varepsilon}$.

**Definition 1 (Boundary corrector)** *Let $v_{1,\varepsilon}$ be the a.s. unique solution in $H^1(B_R \setminus \overline{D}) \times H^1(D)$ of*

$$\begin{cases} -\Delta v_{1,\varepsilon} - k^2 n_0 v_{1,\varepsilon} = 0 \quad in\, B_R \setminus \overline{D}, \\ -\nabla \cdot a_\varepsilon \nabla v_{1,\varepsilon} - k^2 n_\varepsilon v_{1,\varepsilon} = 0 \quad in\, D, \\ [v_{1,\varepsilon}]_{\partial D} = \varepsilon u_1^\varepsilon, \\ [\partial_\nu v_{1,\varepsilon}]_{\partial D} = \varepsilon \sum_{i=1}^d \left(\nabla \cdot (\sigma_i^\varepsilon \partial_i u_0)^+\right) \cdot \nu - k^2 \varepsilon (\beta^\varepsilon u_0)^+ \cdot \nu, \\ \partial_\nu v_{1,\varepsilon} = \Lambda(v_{1,\varepsilon}) \quad on\, \partial B_R. \end{cases}$$

*where $\sigma$ is the standard flux corrector [2] and $\beta$ the corrector associated to $n$ [1].*

$v_{1,\varepsilon}$ verifies the same transmission problem (1) as $u_\varepsilon$ but with highly oscillating boundary data. We prove the following estimate for the two-scale expansion error.

**Proposition 2** *[1] A.s.*

$$\|u_\varepsilon - w_{1,\varepsilon}\|_{H^1(B_R)} \lesssim_{u_0} \varepsilon \mu_d(\frac{1}{\varepsilon}) \chi_\varepsilon^1,$$

*where $\chi_\varepsilon^1$ is a random constant with streched exponential moments, and for $y \in \mathbb{R}^+$,*

$$\mu_d(y) = \begin{cases} |\log(2+y)|^{\frac{1}{2}} & if\, d = 2, \\ 1 & if\, d = 3. \end{cases}$$

Next we show that $v_{1,\varepsilon}$ is of order $\varepsilon^{\frac{1}{2}}$ but only near the boundary.

**Theorem 3** *A.s.*

$$\|v_{1,\varepsilon}\|_{H^1(B_R \setminus \overline{D})} + \|\nabla v_{1,\varepsilon}\|_{H^1(D)} \lesssim_{u_0} \varepsilon^{\frac{1}{2}} \mu_d(\frac{1}{\varepsilon})^{\frac{1}{2}} \chi_\varepsilon^2,$$

$$\|v_{1,\varepsilon}\|_{L^2(B_R)} + \|\nabla v_{1,\varepsilon}\|_{L^2(D_\alpha)} \lesssim_{u_0,\alpha} \varepsilon \mu_d(\frac{1}{\varepsilon}) \chi_\varepsilon^3,$$

*where $D_\alpha := \{x \in B_R \,|\, dist(x, \partial D) > \alpha\}$ for $\alpha > 0$ and $\chi_\varepsilon^{2,3}$ have streched exponential moments.*

For $d = 3$ since $\phi, \sigma$ and $\beta$ can be constructed to be stationary, one can introduce $u_2^\varepsilon$ the second-order correction in $D$ [3] and $v_{2,\varepsilon}$ the associated boundary corrector such that $u_\varepsilon - w_{2,\varepsilon}$ verifies the transmission conditions across $\partial D$, where

$$w_{2,\varepsilon} := u_0 - \varepsilon u_1^\varepsilon \mathbb{1}_D - v_{1,\varepsilon} - \varepsilon^2 u_2^\varepsilon \mathbb{1}_D - v_{2,\varepsilon}.$$

We prove the equivalent of Proposition 2 and Theorem 3 for $u_\varepsilon - w_{2,\varepsilon}$ at order $\varepsilon^{\frac{3}{2}}$ instead of $\varepsilon$. In particular this improves the estimate for $u_\varepsilon - w_{1,\varepsilon}$ in $H^1(D_\alpha)$.

## 4 Point-wise estimate

For $x \in B_R \setminus \overline{D}$, $u_\varepsilon$ admits the following representation formula

$$u_\varepsilon = u_0 + \int_D (a^\star - a_\varepsilon) \nabla u_\varepsilon \cdot \nabla G_0 - k^2 \int_D (n^\star - n_\varepsilon) u_\varepsilon G_0,$$

where $G_0$ is the Green function associated to (2). This together with the homogenization theory leads to considering

$$\tilde{v}_{1,\varepsilon} := \int_D (a^\star - a_\varepsilon) \sum_{i=1}^d (\nabla \phi_i^\varepsilon + e_i) \partial_i u_0 \cdot \nabla G_0 - k^2 \int_D (n^\star - n_\varepsilon) u_0 G_0$$

as a first-order correction in $H^1(B_R \setminus \overline{D})$. $\tilde{v}_{1,\varepsilon}$ verifies the same transmission problem (2) as $u_0$ but with highly oscillating boundary data. We establish the following fluctuation result.

**Theorem 4** *[1] If $G_{0|_D} \in W^{2,\infty}(D)$, then a.s. for $y \in D_\alpha \setminus \overline{D}$, $\alpha > 0$*

$$Var[(u_\varepsilon - u_0 - \tilde{v}_{1,\varepsilon})(y)]^{\frac{1}{2}} \lesssim_{u_0, G_0, \alpha} \varepsilon^{\frac{d+1}{2}} \mu_d(\frac{1}{\varepsilon})^{\frac{1}{2}} \chi_\varepsilon^4,$$

*where $\chi_\varepsilon^4$ has streched exponential moments.*

The difference between $v_{1,\varepsilon}$ and $\tilde{v}_{1,\varepsilon}$ can be proven to be of order $\varepsilon$ in $H^1$-norm.

**Proposition 5** *A.s.*

$$\|\tilde{v}_{1,\varepsilon} - v_{1,\varepsilon}\|_{H^1(D_\alpha)} \lesssim_{u_0} \varepsilon \mu_d(\frac{1}{\varepsilon})^{\frac{1}{2}} \chi_\varepsilon^5,$$

*where $\chi_\varepsilon^5$ has streched exponential moments.*

Numerical simulations illustrate these different expansions.

# Wave propagation and gaps in three-dimensional periodic media

Yuri Godin[1,*], Boris Vainberg[1]
[1]Department of Mathematics and Statistics, University of North Carolina at Charlotte, Charlotte, USA
*Email: ygodin@charlotte.edu

**Abstract**

We investigate the propagation of acoustic waves in an infinite medium containing a periodic array of identical inclusions of arbitrary shape with the Dirichlet or transmission conditions on their interfaces. The inclusion size $a$ is much smaller than the array period while the wavelength is fixed. We derive the dispersion relation and show that there are frequencies for which the solution is a cluster of waves propagating in different directions with different frequencies. The approach is based on reducing the problem to the study of the condition for the existence of a zero eigenvalue (EZE) of a specially constructed Dirichlet-to-Neumann (DtN) operator.

We also introduce and study the concept of local gaps, depending on the choice of the wave vector $\boldsymbol{k}$. The location of local gaps is determined analytically for the Dirichlet and transmission problems.

## 1 Introduction

We are looking for the Bloch solution $u(\boldsymbol{x}) = \varPhi(\boldsymbol{x})\mathrm{e}^{-\mathrm{i}\boldsymbol{k}\cdot\boldsymbol{x}}$ of the equation

$$\Delta u + k^2 u = 0 \tag{1}$$

in the fundamental cell of periodicity $\varPi = [-\pi, \pi]^3$ containing a small inclusion $\varOmega = \varOmega(a)$ which is produced by compressing an $a$-independent region using the factor $a^{-1}$, $0 < a \ll 1$. Here $k = \omega/c$ is the wave number, $k = |\boldsymbol{k}|$, function $\varPhi$ is periodic with the periods of the lattice, $\omega$ is time frequency and $c$ is the wave propagation speed in the host medium. We consider two types of boundary conditions on the surface $\partial\varOmega$ of the inclusion: (i) the Dirichlet boundary condition

$$u|_{\partial\varOmega} = 0 \tag{2}$$

and (ii) the transmission problem with the continuity conditions

$$[\![u(\boldsymbol{x})]\!] = 0, \quad \left[\!\!\left[\frac{1}{\varrho(\boldsymbol{x})}\frac{\partial u(\boldsymbol{x})}{\partial \boldsymbol{n}}\right]\!\!\right] = 0, \tag{3}$$

where $\varrho(\boldsymbol{x})$ is the mass density, $\boldsymbol{n}$ is the outward normal vector, and the brackets $[\![\cdot]\!]$ denote the jump of the enclosed quantity across the boundaries of the inclusions.

Note that in the inclusionless problem, the wave vector $\boldsymbol{k}$ can be multiple if there exist non-trivial vectors $\boldsymbol{m} = \boldsymbol{m}_s \in \mathbb{Z}^3$, $2 \leqslant s \leqslant n$, such that $|\boldsymbol{k}| = |\boldsymbol{k} - \boldsymbol{m}|$. We set $\boldsymbol{m}_1 = \boldsymbol{0}$.

## 2 Outline of the approach

One commonly used technique to study dispersion relations in periodic media analytically is matched asymptotic expansions. Most previous works rely on the construction of waves that approximately satisfy (1) and boundary conditions (2) or (3). However, without an estimate on the inverse operator, it cannot be ensured that the constructed functions are close to the exact solution. We obtained not just a function that almost satisfies the equation, but an asymptotic expansion of the exact solution. In particular, this allows us to see that the clusters of Block waves replace the waves in singular directions when the wave vector is multiple.

We enclose the inclusion in a ball $B_R \subset \varPi$ of radius $R$ and consider two problems in $B_R$ and $\varPi \setminus B_R$:

$$\left(\Delta + k^2\right) u(\boldsymbol{x}) = 0, \ \boldsymbol{x} \in B_R, \ u|_{r=R} = \psi, \tag{4}$$

$$\left(\Delta + k^2\right) v(\boldsymbol{x}) = 0, \ \boldsymbol{x} \in \varPi \setminus B_R, \quad \left]\!\!\right] \mathrm{e}^{\mathrm{i}\boldsymbol{k}\cdot\boldsymbol{x}} v(\boldsymbol{x}) \left[\!\!\left[ = 0, \ v|_{r=R} = \psi, \right.\right. \tag{5}$$

with appropriate boundary conditions on $\partial\varOmega$. Solution of (4),(5) defines two DtN operators: $\boldsymbol{N}^-_{a,\varepsilon}$, $\boldsymbol{N}^+_{\boldsymbol{k},\varepsilon}$

$$\boldsymbol{N}^-_{a,\varepsilon}\psi = \left.\frac{\partial u}{\partial r}\right|_{r=R} \text{ and } \boldsymbol{N}^+_{\boldsymbol{k},\varepsilon}\psi = \left.\frac{\partial v}{\partial r}\right|_{r=R}, \tag{6}$$

where

$$k = (1+\varepsilon)|\boldsymbol{k}|. \tag{7}$$

To obtain the dispersion relation, the following lemma plays a key role.

**Lemma 1** *(i) The operator $\boldsymbol{N}^+_{\boldsymbol{k},\varepsilon}$ is locally analytic in $\boldsymbol{k} \in \mathbb{R}^3$ and $\varepsilon$. The operator $\boldsymbol{N}^-_{a,\varepsilon}$ is locally analytic in $\boldsymbol{k}$ and $\varepsilon$ and infinitely smooth in $a$, $a \ll 1$. (ii) The dispersion relation (7) for which (1–3) have a nontrivial solution is given by the same function $\varepsilon = \varepsilon(a, \boldsymbol{k})$ for which the operator $\boldsymbol{N}^+_{\boldsymbol{k},\varepsilon} - \boldsymbol{N}^-_{a,\varepsilon}$ has a non-trivial kernel.*

The Lemma reduces the study of the dispersion relation to the EZE problem for a smooth operator function. Special matrix representations of operators $\boldsymbol{N}^{+}_{\boldsymbol{k},\varepsilon}$, $\boldsymbol{N}^{-}_{a,\varepsilon}$ were constructed in [1,2] that allowed us to solve the EZE problem and find the asymptotics of the dispersion relation $\varepsilon = \varepsilon(a, \boldsymbol{k})$ as $a \to 0$ which is analytic in $\boldsymbol{k}$ for simple $\boldsymbol{k}$ and is valid for fixed multiple $\boldsymbol{k}$.

## 3 The transmission problem

We illustrate the results for the case of penetrable spherical inclusions.

**Theorem 2** *(i) If $\boldsymbol{k}$ is a simple Bloch vector then the dispersion relation is given by*

$$\overline{\gamma}\omega^2 = \langle \nu \rangle |\boldsymbol{k}|^2 + \mathcal{O}\left(a^4\right), \quad a \to 0. \tag{8}$$

*where*

$$\overline{\gamma} = \gamma_+(1-f) + \gamma_- f \tag{9}$$

*is the average adiabatic bulk compressibility modulus of the medium with $f$ being the volume fraction of the inclusions, and the average specific volume $\langle \nu \rangle$ is determined by Maxwell's formula*

$$\langle \nu \rangle = \nu_+ \left(1 + 3\frac{\nu_- - \nu_+}{\nu_- + 2\nu_+} f\right). \tag{10}$$

*Solution of the problem (1),(3) has the form*

$$u(\boldsymbol{x}) = C\mathrm{e}^{-\mathrm{i}\boldsymbol{k}\cdot\boldsymbol{x}}\left(1 + \tilde{u}(\boldsymbol{x})\right),$$
$$\|\tilde{u}\|_{C^\infty(\Pi\setminus B_R)} = \mathcal{O}(a), \quad a \to 0, \tag{11}$$

*where $C$ is a constant.*
*(ii) Let $\boldsymbol{k}$ be an multiple Bloch vector of order $n \geqslant 2$. Assume that the eigenvalues $\lambda_s$, $1 \leqslant s \leqslant n$, of matrix $M_{i,j} = |\Pi||\boldsymbol{k}_0|^2(\alpha + \beta\hat{\boldsymbol{k}}_i \cdot \hat{\boldsymbol{k}}_j)f$ are distinct and $\boldsymbol{\mu}_s = (\mu_{1,s}, \dots, \mu_{n,s})$ are the corresponding eigenvectors, $\alpha = 1 - \frac{\gamma_-}{\gamma_+}, \beta = 3\frac{\rho_+ - \rho_-}{\rho_+ + 2\rho_-}$. Consider a small neighborhood $I$ of the point $\omega = |\boldsymbol{k}|c_+ = |\boldsymbol{k}|/\sqrt{\rho_+\gamma_+}$. Then, for small $a$, there are exactly $n$ values of the frequencies $\omega = \omega_s \in I$ for which the corresponding problem (1),(3) with $|k_\pm|^2 = \omega_s^2\rho_\pm\gamma_\pm$ has a non-trivial solution. The frequencies $\omega_s$ satisfy the relation*

$$\omega_s^2\rho_+\gamma_+ = |\boldsymbol{k}|^2\left(1 + \lambda_s f\right) + \mathcal{O}(a^4), \; a \to 0. \tag{12}$$

*The corresponding solutions of the problem (1),(3) have the form of clusters*

$$u_s(\boldsymbol{x}) = C_s\left(\sum_{j=1}^{n} \mu_{j,s}\,\mathrm{e}^{-\mathrm{i}(\boldsymbol{k}-\boldsymbol{m}_j)\cdot\boldsymbol{x}} + \tilde{u}_s\right),$$
$$\|\tilde{u}_s\|_{C^\infty(\Pi\setminus B_R)} = \mathcal{O}(a), \; a \to 0, \tag{13}$$

*where $C_s$ are constants and $]\!]\mathrm{e}^{\mathrm{i}\boldsymbol{k}\cdot\boldsymbol{x}}\,\tilde{u}_s[\![= 0$. Here $]\!]f[\![= 0$ denotes the jump of $f$ and its gradient across the opposite sides of the cells of periodicity.*

## 4 The notion of local gaps

Global gaps do not exist in any fixed interval $\varepsilon < \omega < \varepsilon^{-1}$ of the time frequency $\omega$ if the size $a$ of the inclusion is small enough. This result motivates us to introduce local gaps. *The local gap $g(\boldsymbol{k}_0)$* associated with the wave vector $\boldsymbol{k}_0$ consists of frequencies $\omega = \omega_a$ for which Bloch waves with the frequency $\omega_a$ do not exist when $\boldsymbol{k} = (1+\delta)\boldsymbol{k}_0$ for sufficiently small $a$-independent $|\delta|$ and $\omega = \omega_a \to \omega_0$ as $a \to 0$, where $(\omega_0, \boldsymbol{k}_0)$ belongs to the unperturbed dispersion cone $\omega_0 = c|\boldsymbol{k}_0|$. We prove several theorems concerning the existence and location of local gaps.

**Theorem 3** *If $\boldsymbol{k}_0$ is a simple wave vector then there is no local gap for the wave vector $\boldsymbol{k}_0$ and frequencies $\omega$ near $\omega_0$.*

**Theorem 4** *Let the wave vector $\boldsymbol{k}_0$ have a multiplicity of two, i.e. there is an integer-valued vector $\boldsymbol{m}_0$ such that $|\boldsymbol{k}_0| = |\boldsymbol{k}_0 - \boldsymbol{m}_0|$. Then local gap $g(\boldsymbol{k}_0)$ exists in a neighborhood of $\omega_0$ for small $a$ if and only if $\frac{|\boldsymbol{k}_0|}{|\boldsymbol{m}_0|} < \frac{\sqrt{2}}{2}$. In the Dirichlet problem, $g(\boldsymbol{k}_0)$ consists of frequencies $\omega$ such that*

$$|\boldsymbol{k}_0| + a\mu_- + \mathcal{O}(a^2) < \frac{\omega}{c} < |\boldsymbol{k}_0| + a\mu_+ + \mathcal{O}(a^2).$$

*Here*

$$\mu_\pm = \frac{4\pi q}{|\Pi|}\left(1 \pm \frac{4|\boldsymbol{k}_0|}{|\boldsymbol{m}_0|}\sqrt{\frac{1}{2} - \frac{|\boldsymbol{k}_0|^2}{|\boldsymbol{m}_0|^2}}\right),$$

*where $q$ depends on the shape of $\partial\Omega$ and $qa$ is equal to the electrical capacitance of $\partial\Omega$, i.e. its total charge when the potential is unity. In particular, $q = 1$ for a sphere.*

We prove a similar theorem for the transmission problem [3].

# Effective attenuation and Gaussian statistics of the speckle for the paraxial wave equation.

**Christophe Gomez**[1,*], **Olivier Pinaud**[2]
[1]Aix Marseille Univ, CNRS, I2M, Marseille, France
[2]Colorado State University, Department of Mathematics, Fort Collins, CO, USA
*Email: christophe.gomez@univ-amu.fr

**Abstract**

In this talk, we present the asymptotic behavior of solutions to the paraxial wave equation in random media with general mixing assumptions, in the scintillation regime. The coherent part of the wavefield is asymptotically described by the initial condition multiplied by an effective exponential damping term and a phase modulation. This behavior is deterministic and arises from a self-averaging mechanism, whereby high-frequency oscillations suppress the underlying stochasticity. The resulting exponential decay reflects the progressive loss of coherence of the wavefield, with coherent energy being transferred to increasingly small-scale incoherent fluctuations known as speckle. Such incoherent components are typically analyzed via the Wigner transform, yielding deterministic radiative transfer equations. In contrast, our approach focuses directly on the wavefield itself. Although our analysis relies on quantities closely related to the Wigner formalism, our main result provides a precise characterization of the statistics of the speckle. We prove that these fluctuations converge asymptotically to a Gaussian random field, given by an Ornstein–Uhlenbeck-type process, whose covariance structure is explicitly linked to the associated radiative transfer model.

## 1 Introduction

Wave propagation in random media gives rise to speckle patterns whose statistical properties play a central role in optics, acoustics, and imaging. In the scintillation regime (long-distance weak-coupling regime, see (1) and (2)), Gaussianity of the speckle has long been accepted in the physics literature despite the absence of a complete mathematical justification. Early rigorous progress was made within the white-noise paraxial framework (Itô-Schrödinger model), where fourth-order moment calculations were shown to be consistent with Gaussian statistics in the scintillation regime [2]. More recently, a comprehensive analysis of the Itô-Schrödinger model has established convergence of the wave field to a circular complex Gaussian limit in the diffusive scintillation regime [1, 3]. The present work starts from the paraxial wave equation, formulated as a Schrödinger equation with a random potential possessing mixing properties, rather than the Itô-Schrödinger model with explicit white noise, and still establishes Gaussian speckle statistics. We therefore simultaneously consider the scintillation regime and a white-noise limit, leading to a distinct radiative transfer model that characterizes the covariance structure of the speckle, different from the one previously obtained.

## 2 Model

The following paraxial wave equation is considered

$$i\kappa\partial_z\phi + \frac{1}{2}\Delta_{\mathbf{x}}\phi + \kappa^2\sigma V\Big(\frac{z}{\ell_c}, \frac{\mathbf{x}}{\ell_c}\Big)\phi = 0,$$

for $z \geq 0$ and $\mathbf{x} \in \mathbb{R}^2$, and equipped with an initial condition at $z = 0$,

$$\phi(z = 0, \mathbf{x}) = \phi_0\Big(\frac{\mathbf{x}}{r_0}\Big)e^{i\mathbf{k}_0\cdot\mathbf{x}/\lambda}.$$

Here, $\kappa = 2\pi/\lambda$ is the wavenumber ($\lambda$ is the wavelength), $V$ is a random field describing the random fluctuations of the propagation media, $\sigma$ characterizes their order of magnitude, and $\ell_c$ their typical scale of variation (correlation length). $\phi_0$ denotes the envelope of the initial condition, with characteristic width $r_0$, and spatial frequency $\mathbf{k}_0/\lambda$. More precisely, we consider a typical propagation distance of order one, $L \sim 1$, and assume

$$\sigma \ll 1, \qquad \ell_c \sim \lambda \ll 1 \qquad \text{and} \qquad r_0 \gtrsim \lambda.$$

To encode all these assumptions through only one parameter, we consider the dimensionless

parameter

$$\varepsilon := \frac{\lambda}{2\pi L} \ll 1,$$

and

$$\frac{\ell_c}{L} = \varepsilon, \qquad \frac{r_0}{L} = \varepsilon^{1-\alpha}, \qquad \sigma = \sqrt{\varepsilon}, \qquad (1)$$

with $\alpha \in [0,1]$. After nondimensionalizing the paraxial wave equation through the rescaling $z \to zL$, $\mathbf{x} \to \mathbf{x}L$, we obtain

$$i\partial_z \phi^\varepsilon + \frac{1}{2\varepsilon}\Delta_{\mathbf{x}}\phi^\varepsilon + \frac{1}{\sqrt{\varepsilon}}V\Big(\frac{z}{\varepsilon}, \frac{\mathbf{x}}{\varepsilon}\Big)\phi^\varepsilon = 0, \quad (2)$$

where $\phi^\varepsilon(z=0,\mathbf{x}) = \phi_0(\mathbf{x}/\varepsilon^{1-\alpha})e^{i\mathbf{k}_0\cdot\mathbf{x}/\varepsilon}$, with central wavevector $\mathbf{k}_0$. Here $V$ is assumed to have general mixing properties allowing averaging phenomena and exhibit effective wave dynamics. In the following section, we describe two asymptotic behaviors that emerge from the wavefield $\phi^\varepsilon$ in the limit $\varepsilon \to 0$.

## 3 Results

Our first result describes the effective behavior of the coherent part of the wavefield at the scale of the envelope

$$\psi^\varepsilon(z,\mathbf{x}) := \phi^\varepsilon(z, \varepsilon^{1-\alpha}\mathbf{x}).$$

Due to multiple scattering effects, the asymptotic behavior of the coherent part of the wavefield, in the high-frequency limit $\varepsilon \to 0$, is described in terms of an effective attenuation. More precisely, in The Fourier domain, considering the filtered process by the homogeneous Schrödinger semi-group, and properly centered around the central wavevector $\mathbf{k}_0$,

$$\zeta^\varepsilon(z,\mathbf{k}) := \widehat{\psi}^\varepsilon\Big(z, \mathbf{k} + \frac{\mathbf{k}_0}{\varepsilon^\alpha}\Big)e^{it|\mathbf{k}+\mathbf{k}_0/\varepsilon^\alpha|^2/(2\varepsilon^{1-2\alpha})},$$

we obtain

$$\lim_{\varepsilon\to 0}\zeta^\varepsilon(z,\mathbf{k}) = e^{-\Sigma_\alpha(\mathbf{k})z}\widehat{\phi}_0(\mathbf{k}),$$

where the convergence holds in probability in $\mathcal{C}([0,\infty), L^2_w(\mathbb{R}^2))$, and where $w$ stands for the weak topology. Here, we have

$$\Sigma_\alpha(\mathbf{k}) = \begin{cases} \Sigma(\mathbf{k}+\mathbf{k}_0) & \text{if } \alpha = 0, \\ \Sigma(\mathbf{k}_0) & \text{if } \alpha \in (0,1], \end{cases}$$

and

$$\Sigma(\mathbf{k}) = (2\pi)^d \int d\mathbf{p}\int_0^\infty ds\, \widehat{R}(s,\mathbf{p})e^{is(|\mathbf{k}|^2 - |\mathbf{k}-\mathbf{p}|^2)/2}.$$

Here, $\widehat{R}$ is the Fourier transform of the correlation function of the random fluctuations $V$.

Due to multiple scattering and the above effective damping on the coherent part of the wavefield, the missing energy has been converted into small-scale incoherent fluctuations known as speckle. The statistics of this speckle are studied through the locally rescaled wavefield

$$\Phi_{\mathbf{x}}(z,\mathbf{x}') := \phi^\varepsilon(z, \mathbf{x} + \varepsilon\,\mathbf{x}').$$

for a given point $\mathbf{x}$. For $\alpha < 1/2$, the asymptotic behavior of the appropriately compensated wavefield is characterized by an Ornstein-Uhlenbeck-type Gaussian process of the form

$$\xi_{\mathbf{k}}(z,\mathbf{q}) = e^{-\Sigma_\alpha(\mathbf{k})z}\xi_0(\mathbf{k}) + \sqrt{2}\,e^{-i\Sigma^s_\alpha(\mathbf{k})z}\int_0^z\int B_{\mathbf{k}}(ds,d\mathbf{y})e^{-\Sigma^c_\alpha(\mathbf{k})(z-s)}e^{-i\mathbf{q}\cdot\mathbf{y}},$$

with $\xi_0(\mathbf{k}) := \widehat{\phi}_0(\mathbf{k})$ if $\alpha = 0$, and 0 if $\alpha > 0$, where $\Sigma^c_\alpha$ and $\Sigma^s_\alpha$ stand respectively for the real and imaginary part of $\Sigma_\alpha$. Here, $B_{\mathbf{k}}$ is a circular complex Gaussian random field, with a covariance function characterized by the radiative transfer equation

$$\partial_z \mathbf{w} + \mathbf{p}\cdot\nabla_{\mathbf{y}}\mathbf{w} = \int \sigma(\mathbf{p},\mathbf{p}')(\mathbf{w}(\mathbf{p}') - \mathbf{w}(\mathbf{p}))d\mathbf{p}',$$

where

$$\sigma(\mathbf{p},\mathbf{p}') = 2(2\pi)^d\int_0^\infty ds\,\widehat{R}(s,\mathbf{p}'-\mathbf{p}) \times \cos(s(|\mathbf{p}'|^2 - |\mathbf{p}|^2)/2),$$

and with initial conditions

$$\mathbf{w}(z=0) = \begin{cases} |\widehat{\psi}_0(\mathbf{p})|^2\delta(\mathbf{y}-\mathbf{x}) & \text{if } \alpha = 0, \\ (2\pi)^d\|\psi_0\|^2_{L^2_{\mathbf{x}}}\delta(\mathbf{y}-\mathbf{x})\delta(\mathbf{p}) & \text{if } \alpha > 0. \end{cases}$$

For $\alpha \geq 1/2$ it seems that the speckle fluctuations cannot be described in term of a diffusion process.

# Mixed Precision Implicit Runge–Kutta and Explicit Two-Derivative Runge–Kutta Methods

**Sigal Gottlieb**[1,*], **Zachary J. Grant**[1], **César Herrera**[2]
[1]Dept. of Mathematics, UMass Dartmouth, 285 Old Westport Rd, N. Dartmouth, MA 02747
[2]Dept. of Mathematics, Purdue University, 150 North University St. West Lafayette, IN 47907
*Email: sgottlieb@umassd.edu

**Abstract**

In this talk we will present mixed precision time integration for hyperbolic PDEs. We will show how the lower precision perturbations can be controlled and even corrected so that it does not impact the order of accuracy.

## 1 Mixed Precision DIRK Methods

In numerically approximating the solution to a partial differential (PDE) equation

$$U_t = f(U, U_x, U_{xx}), \tag{1}$$

we first approximate the spatial derivatives on the right hand side to obtain the system of ODEs

$$u_t = F(u). \tag{2}$$

We can evolve this using a diagonally implicit Runge–Kutta method, where the implicit component is done in low precision to reduce the computational cost:

$$\begin{aligned} y^{(i)} &= y_n + \Delta t \sum_{j=1}^{i-1} a_{ij} F(y^{(j)}) + \Delta t a_{ii} F_\varepsilon(y^{(j)}) \\ y_{n+1} &= y_n + \Delta t \sum_{i=1}^{s} b_i F(y^{(i)}). \end{aligned} \tag{3}$$

Grant [5] introduced a theoretical framework that describes the perturbation errors introduced by the lower precision computations, which starts with the additive form:

$$\mathbf{y} = u^n \mathbf{e} + \Delta t \tilde{\mathbf{A}} F(\mathbf{y}) + \varepsilon \Delta t \mathbf{A}_\varepsilon \tau(\mathbf{y})$$
$$u^{n+1} = u^n + \Delta t \tilde{\mathbf{b}} F(\mathbf{y}) + \varepsilon \Delta t \mathbf{b}_\varepsilon \tau(\mathbf{y})$$

where $\tilde{A} = \mathbf{A} + \mathbf{A}_\varepsilon$ and $\tilde{\mathbf{b}} = \mathbf{b} + \mathbf{b}_\varepsilon$, and

$$\tau(y) = \frac{1}{\varepsilon}\left(F(y) - F_\varepsilon(y)\right).$$

The second order Taylor expansion

$$u^{n+1} = \underbrace{u^n + \Delta t \tilde{\mathbf{b}} \mathbf{e} F(u^n) + \Delta t^2 \tilde{\mathbf{b}} \tilde{\mathbf{c}} F_y(u^n) F(u^n)}_{scheme}$$
$$+ \underbrace{\varepsilon \left(\Delta t \mathbf{b}_\varepsilon + \Delta t^2 F_y(u^n) \tilde{\mathbf{b}} \mathbf{A}_\varepsilon\right) \tau(\mathbf{y})}_{perturbation}$$
$$+ O(\Delta t^3) + O(\varepsilon \Delta t^3)$$

shows two sources of error, the ones from the scheme and the ones from the perturbation. The order conditions from the scheme are well established. For the scheme to achieve perturbation order $O(\varepsilon \Delta t^{m+1})$ we require :

$$m \geq 1 : |\mathbf{b}_\varepsilon| \mathbf{e} = 0, \quad m \geq 2 : |\tilde{\mathbf{b}}||\mathbf{A}_\varepsilon| \mathbf{e} = 0$$

$$m \geq 3 : |\tilde{\mathbf{b}}||\mathbf{A}_\varepsilon||\tilde{\mathbf{A}}| \mathbf{e} = 0, \quad |\tilde{\mathbf{b}}||\tilde{\mathbf{A}}||\mathbf{A}_\varepsilon| \mathbf{e} = 0.$$

To control the perturbation errors coming from mixed precision, each element of $\mathbf{b}_\varepsilon$ must be zero, and additional conditions depend on which stages are excluded from the final row reconstruction. Grant's framework in [5] shows how to design a method that has higher perturbation order. Explicit [5] or inexpensive implicit [2] corrections can then be used further to improve the perturbation error.

Two-derivative Runge–Kutta schemes may offer enhanced stability and accuracy properties compared to classical one-derivative methods, making them attractive in a wide variety of problems. We extend the numerical analysis mixed precision framework previously developed for Runge–Kutta methods to characterize the propagation of the perturbation errors arising from mixed precision computations in explicit and implicit two-derivative Runge–Kutta methods.

### 1.1 Mixed Precision Runge–Kutta

Two derivative Runge–Kutta (TDRK) methods were applied to the time evolution of PDEs in [1, 3, 6, 7]. When applied to a system of ODEs (2) resulting from a semi-discretization of a PDE

(1), the TDRK where the second derivative term is computed in low precision, can be written in matrix form:

$$\begin{aligned}\mathbf{y} &= u^n + \Delta t \mathbf{A} F(\mathbf{y}) + \Delta t^2 \dot{\mathbf{A}}_\varepsilon \dot{F}_\varepsilon(\mathbf{y}) \\ u^{n+1} &= u^n + \Delta t \mathbf{b} F(\mathbf{y}) + \Delta t^2 \dot{\mathbf{b}}_\varepsilon \dot{F}_\varepsilon(\mathbf{y}),\end{aligned}$$

where $\dot{F}_\varepsilon$ denotes the low precision computation of the second derivative $\dot{F}$.

We can extend the additive mixed precision framework in [5], where we use the notation

$$\dot{F}_\varepsilon = \dot{F} + \varepsilon \dot{\tau}$$

where $\dot{\tau}$ is *not* the derivative of the perturbation $\tau$, but a perturbation of the derivative $\dot{F}$. We write the method as

$$\begin{aligned}\mathbf{y} &= u^n + \Delta t \mathbf{A} F(\mathbf{y}) \\ &\quad + \Delta t^2 \dot{\mathbf{A}} \dot{F}(\mathbf{y}) + \varepsilon \Delta t^2 \dot{\mathbf{A}}_\varepsilon \dot{\tau} \\ u^{n+1} &= u^n + \Delta t \mathbf{b} F(\mathbf{y}) + \Delta t^2 \dot{\mathbf{b}} \dot{F} + \varepsilon \Delta t^2 \dot{\mathbf{b}}_\varepsilon \dot{\tau}.\end{aligned}$$

Taylor expanding this form and matching terms up to third order we obtain:

$O(\Delta t):\ \mathbf{b}\mathbf{e} = 1. \quad O(\Delta t^2):\ \mathbf{b}\mathbf{c} + \dot{\mathbf{b}}\mathbf{e} = \frac{1}{2}$

$O(\Delta t^3):\ \mathbf{b}\mathbf{c}^2 + 2\dot{\mathbf{b}}\mathbf{c} = \frac{1}{3}$

$\quad \mathbf{b}(\mathbf{A}\mathbf{c} + \dot{\mathbf{A}}\mathbf{e}) + \dot{\mathbf{b}}\mathbf{c} = \frac{1}{6}$

$O(\varepsilon \Delta t^2):\ \left|\dot{\mathbf{b}}_\varepsilon\right| \mathbf{e} = 0$

$O(\varepsilon \Delta t^3):\ |\mathbf{b}| \left|\dot{\mathbf{A}}_\varepsilon\right| \mathbf{e} = 0 \quad \& \quad \left|\dot{\mathbf{b}}_\varepsilon\right| |\mathbf{c}| = 0.$

where $\mathbf{c}^2$ is understood as element-wise squaring, and the absolute value $|\cdot|$ is understood element-wise both for matrices and vectors.

**First order perturbation error:** All methods that satisfy the formulation above where $F$ is evaluated in full precision, have a perturbation that enters as $\varepsilon \Delta t^2$ at each time-step, so that at the final time the perturbation error builds up to $E_{pert} = O(\varepsilon \Delta t)$.

**Second order perturbation error:** To improve the errors resulting from the low precision perturbation, we can design a TDRK method that satisfies the additional condition

$$\left|\dot{\mathbf{b}}_\varepsilon\right| \mathbf{e} = 0, \quad \rightarrow \left(\dot{\mathbf{b}}_\varepsilon\right)_j = 0, \quad \forall j,$$

and so the perturbation will result in a final time-error $E_{pert} = O(\varepsilon \Delta t^2)$.

The two-stage third order method in [4] satisfies these conditions

$$\begin{aligned}y^{(1)} &= u^n + \frac{2}{3}\Delta t F(u^n) + \frac{2}{9}\Delta t^2 \dot{F}_\varepsilon(u^n) \\ u^{n+1} &= u^n + \frac{1}{4}\Delta t F(u^n) + \frac{3}{4}\Delta t F(y^{(1)}),\end{aligned}$$

so it has a final time error $E = O(\Delta t^3) + O(\varepsilon \Delta t^2)$.

**Third order perturbation error:** To further improve the perturbation errors we additionally require $\left|\dot{\mathbf{b}}_\varepsilon\right| |\mathbf{A}| \mathbf{e} = 0$ and $|\mathbf{b}| \left|\dot{\mathbf{A}}_\varepsilon\right| \mathbf{e} = 0$. The first of these is automatically satisfied by $\left(\dot{\mathbf{b}}_\varepsilon\right)_j = 0, \ \forall j$, so we need only impose

$$\sum_{i=1}^{s} |b_i| \sum_{j=1}^{s} \left| \left(\dot{\mathbf{A}}_\varepsilon\right)_{ij} \right| = 0.$$

By ensuring that the method does not include in the final reconstruction any stage that has a perturbed value we obtain a final time perturbation error we expect an error of $E_{pert} = O(\varepsilon \Delta t^3)$.

The three stage third order method [4]

$$\begin{aligned}y^{(1)} &= u^n + \frac{2}{3}\Delta t F(u^n) + \frac{2}{9}\Delta t^2 \dot{F}_\varepsilon(u^n) \\ y^{(2)} &= u^n + \frac{1}{3}\Delta t F(u^n) + \frac{1}{3}\Delta t F(y^{(1)}) \\ u^{n+1} &= u^n + \frac{1}{4}\Delta t F(u^n) + \frac{3}{4}\Delta t F(y^{(2)})\end{aligned}$$

satisfies all the conditions above, and so will have a final time error of $E = O(\Delta t^3) + O(\varepsilon \Delta t^3)$.

# Towards kernel-free boundary integral methods for variable coefficients

Benedikt Gräßle[1,*], Stefan A. Sauter[1]

[1]Institute for Mathematics, Universtität Zürich, Zürich, Switzerland

*Email: benedikt.graessle@math.uzh.ch

**Abstract**

Classical boundary integral formulations rely on explicit fundamental solutions and the evaluation of singular kernels—tools that are generally unavailable for variable coefficients. A kernel-free boundary element framework replaces analytic kernels with a discretised boundary operator. This results in a computable algebraic formulation that retains the dimension reduction intrinsic to boundary integral methods and extends quasi-optimal boundary element discretisations beyond constant coefficients.

## 1 Introduction

Boundary integral methods are attractive for elliptic problems because they reduce the effective dimension by one and encode far-field interactions in nonlocal boundary operators. The classical boundary element method (BEM) relies on explicit fundamental solutions and the accurate evaluation of singular kernels. While this paradigm is mature for constant-coefficient operators, for most boundary value problems with varying coefficients, the Schwarz kernel is unknown and its evaluation *not* feasible.

This contribution develops and analyses a kernel-free boundary element formulation. Instead of approximating a known boundary kernel, we construct a *computable approximation* of the boundary integral operator itself from a volume discretisation of the underlying partial differential equation (PDE). The resulting boundary Galerkin system preserves the computational structure of BEM while removing the dependence on analytic kernel formulae. From an algorithmic viewpoint, the resulting algebraic structure is compatible with data-sparse representations in hierarchical matrix formats.

## 2 Model problem

Let $\Omega \subseteq \mathbb{R}^n$ be a Lipschitz domain and consider a strongly elliptic PDE operator

$$\mathcal{L}u := -\operatorname{div}(\mathbb{A}\nabla u) + \mathbf{b}\cdot\nabla u + c\,u \qquad \text{in } \Omega$$

with variable coefficients and suitable boundary conditions on $\partial\Omega$. For an interface (or boundary manifold) $\Sigma$, the single-layer equation seeks a density $b$ such that

$$\mathsf{V}b = f \quad \text{on } \Sigma,$$

where $\mathsf{V}$ denotes the single-layer boundary operator induced by $\mathcal{L}$. Standard BEM is kernel-based: The single-layer operator $\mathsf{V}$ is represented by a Green-kernel integral and then discretised by a Galerkin method. This strategy critically depends on explicit Green function representations, exclusively available in specific situations with the famous example being constant coefficients in the full-space [SS11, Subsec. 3.1] and in the half space $\Omega = \mathbb{R}^n_+$ [CH95, DHN07, LMS25]. Further examples include layered media, wave guides in periodic media, and structured composites [CC12, LLQ18, BY21, LYZZ24].

On general Lipschitz domains (boundary conditions) and for operators with variable coefficients, an explicit Green function is typically unavailable—rendering this strategy infeasible. Instead, we depart from a purely operator theoretic description of the boundary potentials and operators: The single-layer operator

$$\mathsf{V} = \gamma \circ \mathcal{L}^{-1} \circ \gamma' \tag{1}$$

is the composition of the (dual) trace map and the PDE solution operator. This approach extends a recent framework [GHS25, GS25, Grä25] for the wavenumber-explicit analysis of skeleton/boundary equations for general coefficients.

### 2.1 Kernel-free discretisation

Instead of relying on explicit kernels, we discuss a discretisation framework that approximates $\mathcal{L}^{-1}$ in (1) through a one-time preprocessing step from a conforming volume discretisation (e.g., finite elements of mesh-size $d$). This preprocessing step depends only on the coefficients of the differential operator $\mathcal{L}$ and is computed ahead-of-time, independent of the interface location: Its cost is comparable to an (approximate) inversion of a FEM stiffness matrix and does *not* need to be repeated when interfaces change within the same medium.

The resulting discrete single-layer operator $\mathsf{V}_d$ for an interface discretisation $B_h$ (e.g., boundary elements of mesh size $h$) is then employed in a Galerkin formulation: Seek $b_h \in B_h$ with

$$\langle \mathsf{V}_d b_h, g_h \rangle = \langle f, g_h \rangle \quad \forall \psi_h \in B_h.$$

Importantly, the discrete operator is computable without an explicit kernel representation: it is defined purely algebraically as a composition of an (approximate) inverse FEM matrix with sparse transfer operators between the volume and boundary spaces. This factorisation enables a direct implementation and keeps the cost of the boundary discretisation moderate.

### 2.2 Stability and error control

The analysis establishes the quasi-best approximation property of the Galerkin scheme with respect to the volume and boundary discretisations, under a natural compatibility condition on their mesh-sizes. For polynomial finite and boundary elements on quasi-uniform meshes, the theory provides optimal rates in Sobolev norms within the regularity window of the associated transmission problem. Thus the method extends quasi-optimal BEM beyond constant coefficients.

# On the role of regularisation parameters in Adaptive Eigenspace Inversion

**Marie Graff**[1,*], **Nasrin Nikbakht**[2], **Melissa Tacy**[2]

[1]Te Whai Ao — Dodd-Walls Centre for Photonic and Quantum Technologies and Department of Mathematics, The University of Auckland, Auckland, New Zealand

[2]Department of Mathematics, The University of Auckland, Auckland, New Zealand

*Email: marie.graff@auckland.ac.nz

**Abstract**

We analyse two regularisation parameters involved in Adaptive Eigenspace Inversion: the optimal number of eigenfunctions based on sensitivities and the diffusivity parameter based on semi-classical analysis. Through a comprehensive study, we clarify their influence on stability and accuracy. Numerical experiments demonstrate improved robustness compared with fixed-parameter strategies and standard Tikhonov regularisation.

## 1 Introduction

The Adaptive Eigenspace Inversion is a regularisation strategy to solve inverse problems, either for sources or scatterers in the context of wave propagation phenomena. The principle of this method is to represent the profile to recover in a well-chosen eigenfunctions basis, which is updated at each adaptation to enhance the accuracy. Over the last decade [1], the limits of the method have been pushed to test the robustness and accuracy of the reconstructed profiles numerically. However, little theoretical understanding was available to ensure the convergence of the strategy. Recent work [2] gave a first answer to select the number of eigenfunctions, key feature of the method, automatically. In our novel study [3], we validate their approach without reshuffling the function basis, and furthermore, we propose a new strategy to adapt the non-cancelling coefficient in the diffusivity used to build the eigenfunction basis. This new automatic selection of both regularisation parameters (number of eigenfunctions and non-cancelling coefficient) gives an incomparable robustness to the reconstruction algorithm.

## 2 Model and data generation

We consider the Helmholtz equation in an unbounded domain $\Omega \subset \mathbb{R}^2$,

$$\Delta u + k^2 u = f, \qquad \text{in } \Omega, \tag{1}$$

with Sommerfeld radiation conditions

$$iku + \partial_n u = 0, \qquad \text{on } \partial\Omega. \tag{2}$$

Here $f$ is some source term and $k$ represents the wavenumber.

Assume that only partial and corrupted information is known from the wavefield $u$, denoted by $d$, and satisfying

$$d = u_{|\Gamma} + \nu_\delta, \tag{3}$$

where $\Gamma \subset \Omega$ is the partial location of the data collection, that is, where the wavefield is known, and $\nu_\delta$ models white noise with variance $\delta^2$.

The inverse scattering problem is now stated as:

*Given some noisy data $d$, find wavenumber $k$ such that $\|d - G(k)\|^2_{L^2(\Gamma)}$ is minimal.*

Here $G$ represents the nonlinear forward operator involving solving Helmholtz problem (1)-(2) and projecting on $\Gamma$. The above-stated inverse problem is ill-posed and requires regularisation to ensure the stability of the optimisation.

## 3 Adaptive Eigenspace Inversion (AEI)

The principle of AEI, also known as ASI (S for Spectral) is to project the wavenumber to recover into the basis of the eigenfunctions satisfying

$$\begin{aligned} -\nabla \cdot (\mu_\varepsilon[k]\nabla\varphi_\ell) &= \lambda_\ell \varphi_\ell, && \text{in } \Omega, \\ \varphi_\ell &= 0, && \text{on } \partial\Omega, \end{aligned} \tag{4}$$

where $\varphi_\ell$ denotes the $\ell$-th eigenfunctions associated to the positive increasing eigenvalue $\lambda_\ell$, and $\mu_\varepsilon$ is the diffusivity with non-cancelling coefficient $\varepsilon$ defined as

$$\mu_\varepsilon[k] := \frac{1}{\sqrt{|\nabla k|^2 + \varepsilon^2}}. \tag{5}$$

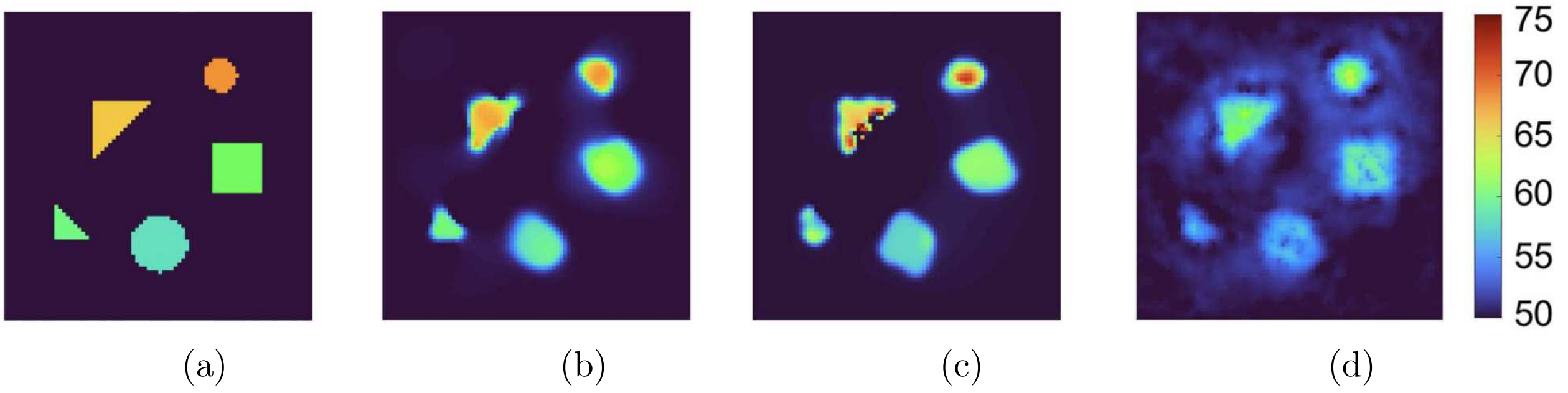

Figure 1: Recovering $k$: (a) reference, (b) using ASI-$\varepsilon$, (c) fixed AEI and (d) Tikhonov regularisation.

Note that the diffusivity relies on the wavenumber $k$, hence yielding to an iterative process: starting from the basis of eigenfunctions of the Laplacian operator, that is, for $k \equiv 1$ and $\varepsilon = 1$, the minimisation is solved a first time to obtain a first estimate. From this estimate, the diffusivity is updated, leading to an improved basis and the minimisation is performed again. Until the misfit reaches the level of noise $\delta$ (also known as Discrepancy Principle of Morozov (DPM) [4]), the adaptation of the basis is repeated. Also, to ensure regularisation, the eigenfunction basis is actually truncated, taking only the first $M$ functions to represent $k$. The choice of $M$ at each adaptation was so far heuristic, but now using sensitivities [2,3], this can be automated for an optimal choice.

## 4 Adapting the non-cancelling coefficient

It has been observed that, in its original form, the AEI algorithm is quite volatile and sensitive to the choice of the regularisation parameters, namely $M$ and $\varepsilon$. In this section, we focus on $\varepsilon$ and an analytic proof [3] led to the ultimate result, at adaptation $j$

$$\varepsilon^{(j)} := \sqrt{\frac{\|\Delta k^{(j-1)}\|_{L^\infty(\Omega)} \|\nabla k^{(j-1)}\|_{L^\infty(\Omega)}}{\tau^{(j)}}}, \quad (6)$$

where $\tau^{(j)}$ is the product of $M^{(j-1)}$, the largest Fourier coefficient in the basis $\{\varphi_\ell\}_{\ell=1}^{M^{(j-1)}}$ and the largest eigenvalue $\lambda_{M^{(j-1)}}$.

Figures 1 and 2 present the results with 2% noise on the data, collected on the boundary. We compare ASI-$\varepsilon$ with standard AEI for arbitrarily chosen, fixed $M = 80$ and $\varepsilon = 1$, and with Tikhonov regularisation, where $\alpha$ is picked so that the resulting misfit satisfies the DPM. In conclusion, adapting $\varepsilon$ together with selecting $M$ from the sensitivities visibly enhances the performance of the algorithm.

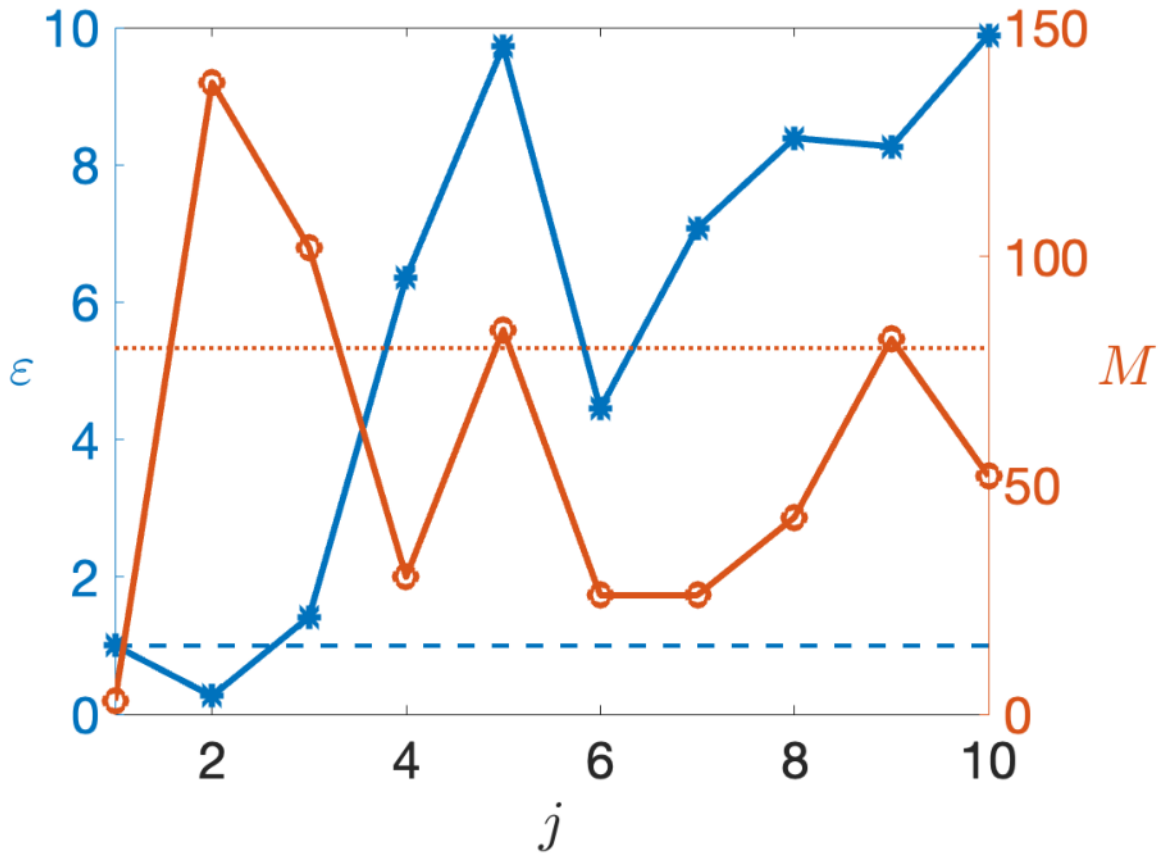

Figure 2: Evolution of $\varepsilon$ and $M$ at each adaptation until convergence (here $j = 10$), starting with $\varepsilon = 1$ and $M = 3$. Dashed line: $\varepsilon = 1$, dotted line: $M = 80$.

# Substructured Optimized Schwarz Domain Decomposition Solver on Multiple GPUs for Helmholtz Problems

R. Greffe[1], A. Chabib[2], A. Modave[2], C. Geuzaine[1]

[1]University of Liège, Department of Electrical Engineering and Computer Science

[2]POEMS, CNRS, Inria, ENSTA, Institut Polytechnique de Paris

**Abstract**

This contribution discusses three strategies to implement a substructured Optimized Schwarz Domain Decomposition Solver for the high-order finite element solution of the Helmholtz Equation on Graphical Processing Units. Advantages and drawbacks of each approach are discussed on a three-dimensional benchmark.

## 1 Introduction

Graphical Processing Units (GPUs) provide most of the computational power in state-of-the-art HPC clusters. The aim of this work is to analyze how to most efficiently implement a substructured Optimized Schwarz Domain Decomposition Method (DDM) on GPUs to solve large scale Helmholtz problems with high-order finite elements. In a substructured DDM, large problems are split into smaller subproblems amenable to direct solvers, and a Krylov method is used to iteratively couple those subproblems through exchange of interface unknowns [1]. The local systems are sparse, which makes them in theory not a great match for GPU architectures. In this contribution, we compare three strategies for implementing the substructured DDM solver on GPUs and compare them on a three-dimensional benchmark.

## 2 Problem statement

We consider a computational domain $\Omega$ with boundary $\Gamma^\infty$, split into $N_{\mathrm{dom}}$ non-overlapping subdomains $\Omega_i$. We define $\Gamma_i^\infty := \partial\Omega_i \cap \Gamma^\infty$ and denote by $\Sigma_{ij}$ the interface between two neighboring subdomains $\Sigma_{ij} := \overline{\partial\Omega_i \bigcap \partial\Omega_j}$.

At iteration $n$, for each $\Omega_i$, the field $u_i^{n+1}$ is the solution the following local problem:

$$\begin{cases} (\Delta + k^2)u_i^{n+1} &= f_i & \text{in } \Omega_i, \\ \partial_{\mathbf{n}_i} u_i^{n+1} + \imath k u_i^{n+1} &= 0 & \text{on } \Gamma_i^\infty, \\ \partial_{\mathbf{n}_i} u_i^{n+1} + \mathcal{S} u_i^{n+1} &= g_{ij}^n & \text{on } \Sigma_{ij}, \end{cases} \quad (1)$$

where $\mathbf{n}_i$ is the unit outward normal on $\partial\Omega_i$, $\mathcal{S}$ is a transmission operator, and $f_i$ is a given source. The unknown $g_{ij}^n$ on the interface $\Sigma_{ij}$ is updated at each iteration using

$$g_{ij}^{n+1} = -g_{ji}^n + 2\mathcal{S} u_j^{n+1}. \quad (2)$$

After discretizing with finite elements and gathering all interface unknowns in a vector $\boldsymbol{g}$, the interface problem can be rewritten as a fixed point iteration $\boldsymbol{g}^{n+1} = \mathcal{F}(\boldsymbol{g}^n) + \boldsymbol{f}$, where $\boldsymbol{f}$ contains the contributions of the sources $f_i$. A Krylov subspace method (e.g. GMRES) is then used to iteratively solve $(I - \mathcal{F})\boldsymbol{g} = \boldsymbol{f}$.

## 3 GPU implementation strategies

The implementations mainly differ in the way the local problems are solved.

**Sparse LU factorization** The first approach consists in using the same strategy on GPUs as on CPUs, i.e. use a sparse direct solver for each local subproblem. While our preliminary results are obtained with the cuDSS library from NVIDIA[2] (which is only compatible with NVIDIA GPUs), other libraries such as SuperLU_dist and PanguLU are being actively ported to GPUs as well. The LU factorizations of the matrices of the local problems are computed in a preprocessing step and reused at each iteration. The main *a priori* disadvantage of this approach is the relatively poor efficiency of sparse factorization on GPUs [2], while its *a priori* advantage is its lower memory requirements compared to dense-algebra approaches (see below), which should make it well suited for a smaller number of large subdomains.

[1]Roland Greffe is a FRIA grantee of the Fonds de la Recherche Scientifique – FNRS

[1]The present research benefited from computational resources made available on Lucia, the Tier-1 supercomputer of the Walloon Region, infrastructure funded by the Walloon Region under the grant agreement n°1910247.

[2]https://developer.nvidia.com/cudss

**Dense LU factorization** In this approach the LU factorizations of the matrices are computed with a batched function for dense matrices from the MAGMA library[3]. At each iteration, the local problems are then solved using batched forward-backward substitutions. The main drawback of this strategy is the memory requirement to store all the dense matrices, which makes it *a priori* better suited for a large number of small subdomains, as this reduces the overall memory footprint and improves the performance of the batched factorizations and forward-backward substitutions.

**Schur Complement** Here, the Schur complement of each local system is explicitly inverted during a preprocessing step in order to only deal with the interface unknowns during iterations. This reduces the size of the local systems and replaces the forward-backward substitution with a batched matrix multiplication which makes each iteration faster. The main *a priori* drawback is the time and memory requirement of the computation of the Schur complements.

## 4 Preliminary results

We consider a single point source in a cubic domain $\Omega$ discretized with P5 finite elements (Figure 1). The transmission operator is $\mathcal{S} = \imath k$. Figure 2 compares the three implementations for cases of increasing sizes on a single GPU.

The Schur complement strategy shows the best performance in terms of time per iteration, closely followed by the sparse LU factorization. As expected, however, the Schur strategy exhibits a much higher preprocessing time compared to the other methods and uses at its peak only just a little bit less memory than the dense LU method.

From preliminary tests, weak scaling on multiple GPUs barely (8.8%) increases the time per iteration on up to 8 GPUs. More details will be shared in the extended paper.

[3] https://icl.utk.edu/magma/

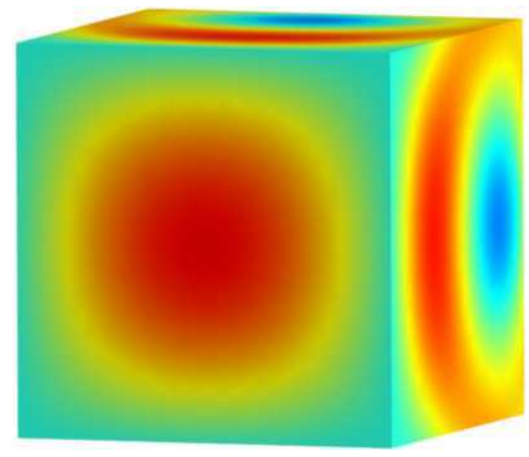

Figure 1: Reference case: point source in a cube.

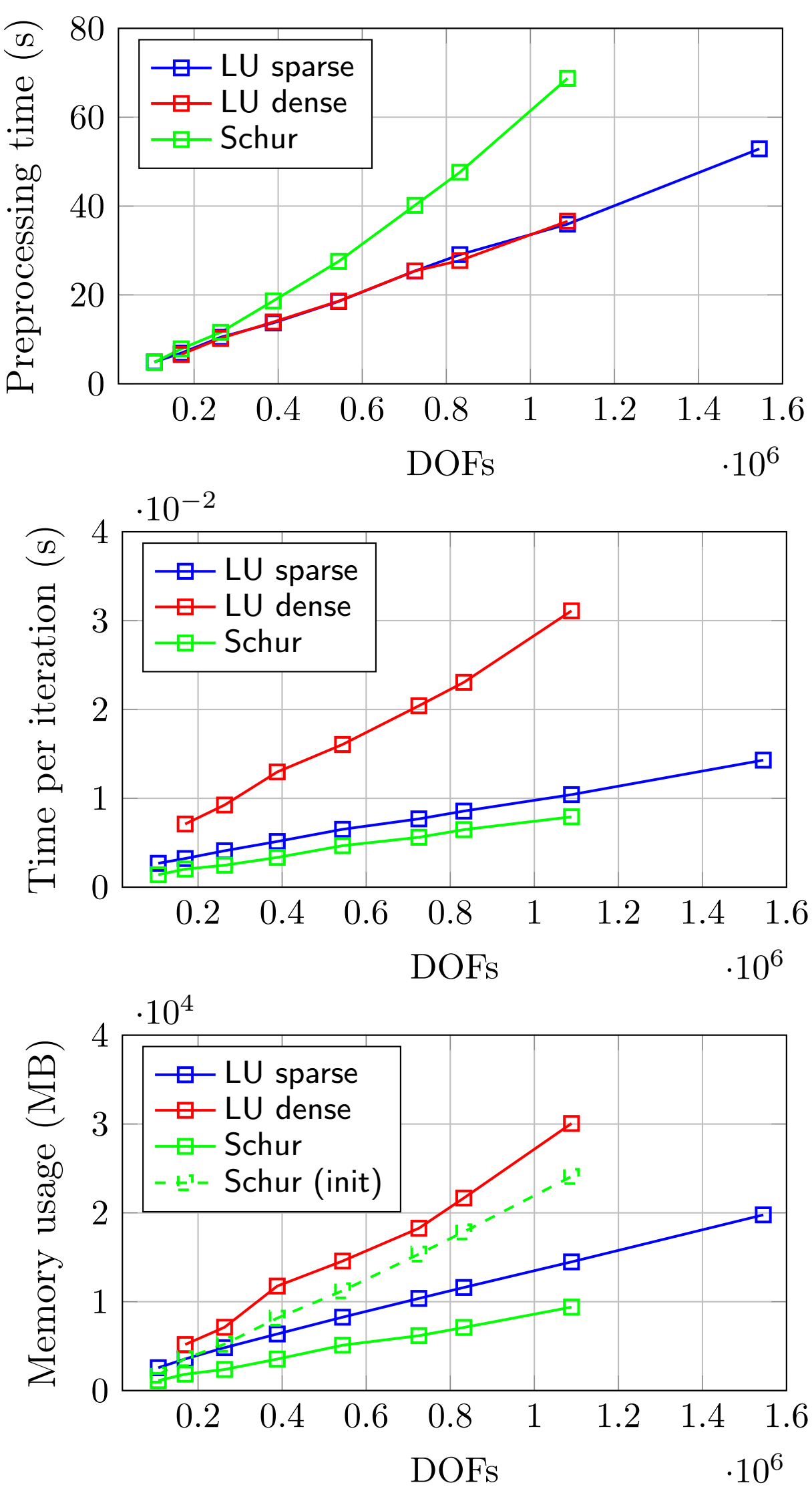

Figure 2: Results for $N_{\text{dom}} = 5000$ on one Nvidia A100 40 GB GPU. Schur (init) refers to the memory usage during the preprocessing step, while Schur refers to the memory usage during the iterations.

# The Inverse Elasto-Acoustic Problem

**Patrick Grice, Yassine Boubendir, Zoi-Heleni Michalopoulou**[1,*]
[1]Department of Mathematical Sciences, New Jersey Institute of Technology, Newark, USA
*Email: pg49@njit.edu

**Abstract**

We investigate inverse obstacle scattering for detecting elastic targets in marine sediments under noisy and uncertain conditions. A boundary integral formulation of the elasto-acoustic problem enables efficient GMRES solution. A domain-derivative Newton type method with multifrequency data is developed for stable reconstruction.

## 1 Introduction

The motivation for this work is the detection and classification of elastic obstacles in marine environments using inverse obstacle scattering techniques. This task becomes particularly challenging when the obstacles are buried within sediment layers. In practice, the reconstructed results are often corrupted by noise and environmental variability. The objective of this work is to develop reconstruction methods that are stable and numerically efficient.

Boundary Integral Equation methods are used to significantly reduce the discretization size of the problem. With numerical results, we demonstrate that by formulating the Elasto-Acoustic problem appropriately, potentially large systems of equations may be solved rapidly by GMRES. We introduce the domain derivative of the Elasto-Acoustic problem and provide examples clarifying its meaning. We describe a method for inverting the domain derivative with regularization and use it to solve the inverse problem by Newton's method using multiple frequencies.

In the scientific literature, Newton's method for inverse obstacle scattering is commonly formulated in free space (e.g., [3]). In this setting, the obstacle is reconstructed from acoustic pressure data on a closed curve that fully encloses the scatterer. Such a configuration is appropriate for certain applications such as medical imaging.

In contrast, for stratified media, acoustic pressure data is measured on curves that do not enclose the obstacle and provide only partial information of the scattered field. Additionally, uncertainty in the geometry and acoustic properties of the water bed further degrades the resolution of the shape reconstruction. Information loss may be mitigated by introducing a perfectly reflecting boundary beneath the scatterer, redirecting downward-propagating waves toward the measurement array. Nevertheless, in many practical scenarios, the obstacle must be reconstructed using incomplete and noisy data.

## 2 Model

The scattering of an incident fluid wave $p_{\text{inc}}$ by an elastic disk is governed by the equations,

$$(\Delta + k^2)p(\mathbf{x}) = 0\,, \qquad \|\,\mathbf{x}\,\| > R\,, \tag{1}$$
$$\nabla \cdot \boldsymbol{\sigma}(\mathbf{x}) + \rho_s \omega^2 \mathbf{u}(\mathbf{x}) = 0\,, \qquad \|\,\mathbf{x}\,\| < R\,, \tag{2}$$

with boundary conditions on $\|\,\mathbf{x}\,\| = R$

$$\boldsymbol{\sigma} \cdot \hat{\mathbf{n}}(\mathbf{x}) + p(\mathbf{x})\hat{\mathbf{n}}(\mathbf{x}) = -p_{\text{inc}}(\mathbf{x})\,, \tag{3}$$
$$\rho_f \omega^2 \mathbf{u} \cdot \hat{\mathbf{n}}(\mathbf{x}) + \frac{\partial p}{\partial n}(\mathbf{x}) = -\frac{\partial p_{\text{inc}}}{\partial n}(\mathbf{x})\,, \tag{4}$$

and the Sommerfeld radiation condition (in $\mathbb{R}^2$),

$$p(\mathbf{x}) - ik\frac{\mathbf{x}}{\|\,\mathbf{x}\,\|} \cdot \nabla p(\mathbf{x}) \to \mathcal{O}\left(\|\,\mathbf{x}\,\|^{-1/2}\right)\,, \tag{5}$$

in the limit $\|\,\mathbf{x}\,\| \to \infty$.

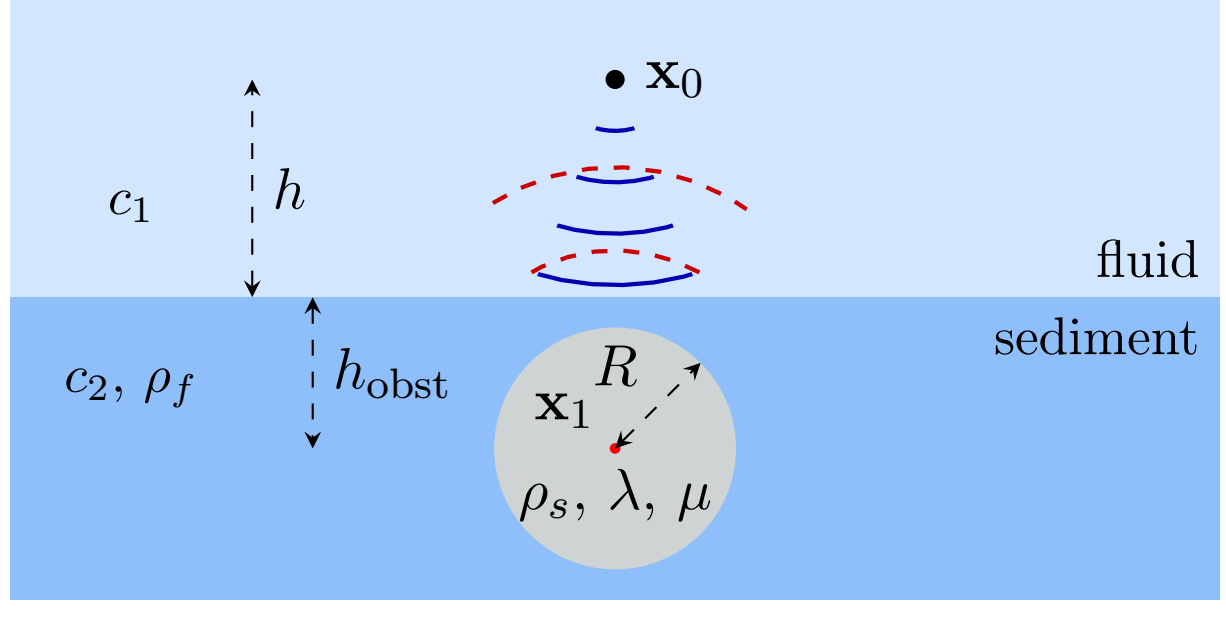

Figure 1: Multiple scattering by a stratified fluid interface and an elastic disk. The sediment is modeled as a dense fluid.

## 3 Methods and Results

The Elasto-Acoustic problem is formulated as a coupled system of Nyström Boundary Integral Equations and discretized using singular Gaussian quadrature rules [2]. By choosing the appropriate the Robin transmission conditions and layer potential representation, the Elasto-Acoustic problem is solved rapidly by GMRES. No expensive matrix factorization or inversion is required to form the block system. The number of GMRES iterations does not increase with the number of discretization points and generlizes well to multiple scattering.

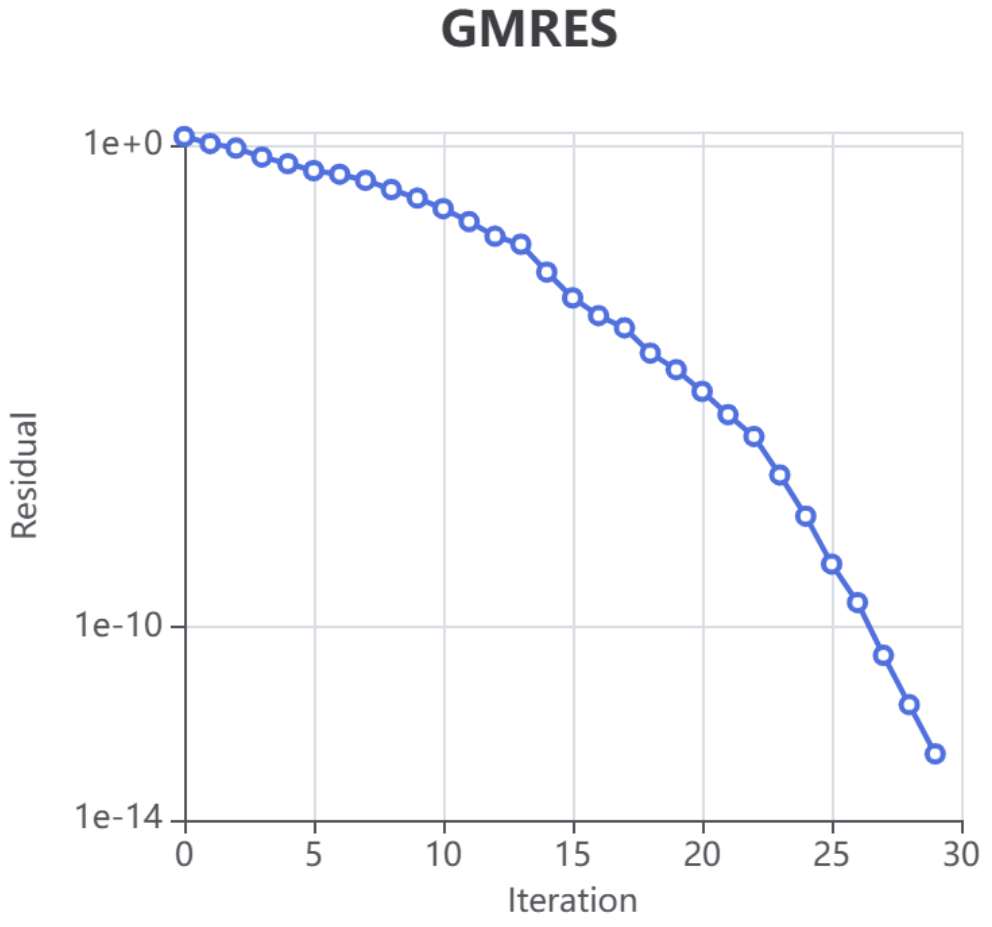

Figure 2: Plane-wave scattering on a disk of radius 2 with 100 discretization points using $\omega = 1$, $\rho_s = 1$, $\rho_f = 1$, $\lambda = 1$, $\mu = 1$, $c = 1$.

To solve the Elasto-Acoustic problem in a stratified environment, constant wave speed is assumed in each layer. The parameterizations of each layer interface is complexified so that the integrands and analytic continuation of the solution decays exponentially [1]. The transmission between layers is formulated as a system of Nyström Boundary Integral Equations and the complexified domain is truncated on a computationally tractable interval.

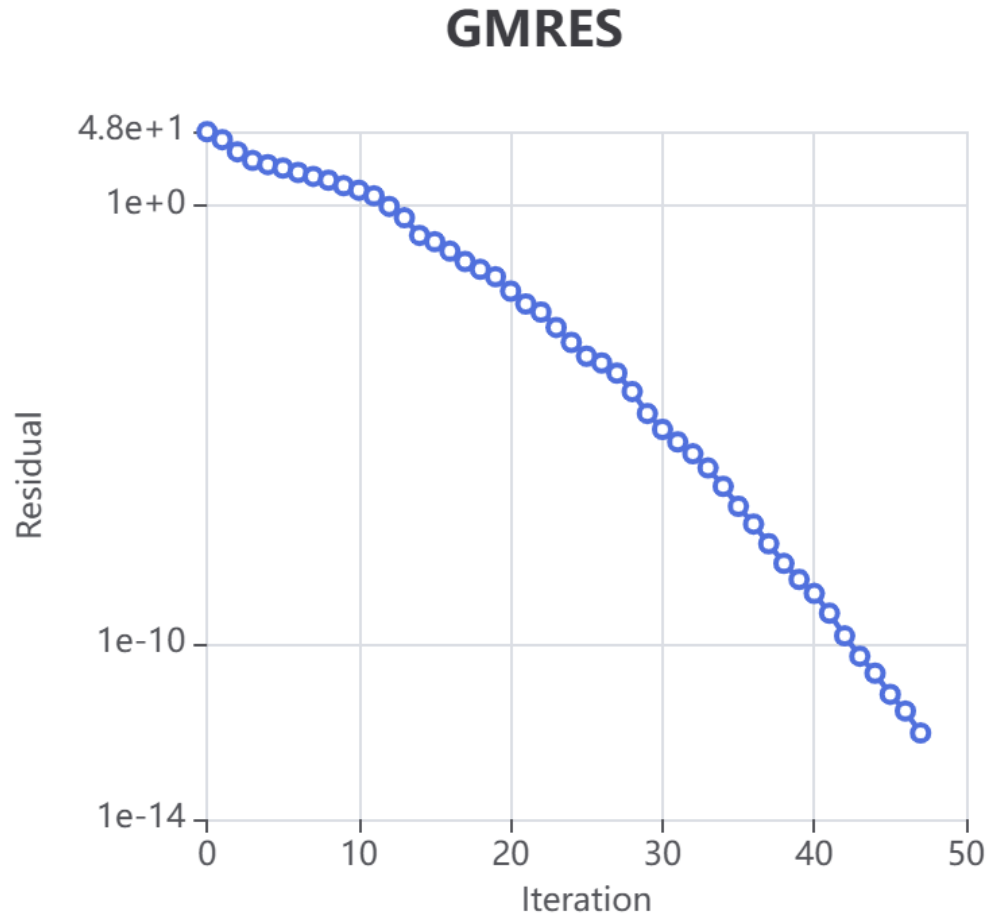

Figure 3: GMRES iterations for multiple scattering by a stratified fluid interface and an elastic disk, see Figure 1. The system solved was 1100 by 1100.

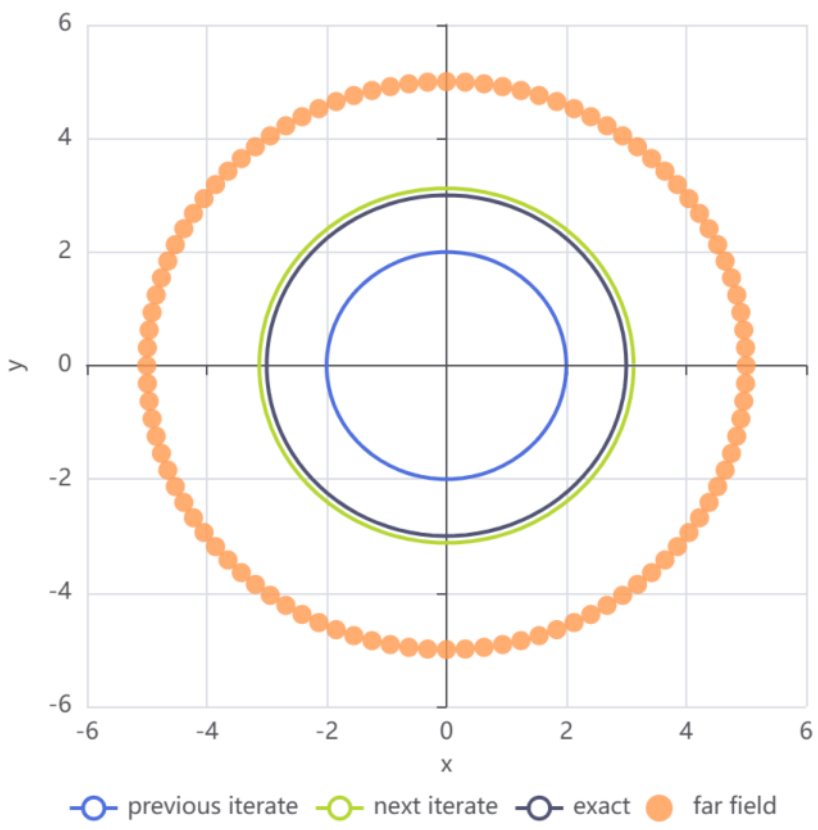

Figure 4: One Newton iteration of the Inverse Elasto Acoustic Problem in free space.

# Adaptive FEM with explicit time integration for the wave equation

Marcus J. Grote[1,*], Omar Lakkis[2], Carina S. Santos[1]
[1]Department of Mathematics and Computer Science, University of Basel, Basel, Switzerland
[2]Department of Mathematics, University of Sussex, Brighton, UK
*Email: marcus.grote@unibas.ch

## Abstract

Starting from a recent a posteriori error estimator for the finite element solution of the wave equation with explicit time-stepping [3], we devise a space-time adaptive strategy which includes both time evolving meshes and local time-stepping [1, 2] to overcome any overly stringent CFL stability restriction on the time-step due to local mesh refinement. Moreover, at each time-step the adaptive algorithm monitors the accuracy thanks to the error indicators and recomputes the current step on a refined mesh until the desired tolerance is met; meanwhile, the mesh is coarsened in regions of smaller errors. Leapfrog based local time-stepping is applied in all regions of local mesh refinement to incorporate adaptivity into fully explicit time integration with mesh change while retaining efficiency. Numerical results illustrate the optimal rate of convergence of the a posteriori error estimators on time evolving meshes [4].

## 1 FEM for the wave equation

We consider the wave equation in $\Omega \subset \mathbb{R}^d$,

$$\begin{aligned} \partial_{tt}u - \nabla\cdot(c^2\,\nabla u) &= f \qquad \text{in } \Omega\times(0,T), \\ u(\cdot,0) = u_0, \quad \partial_t u(\cdot,0) &= v_0 \qquad \text{in } \Omega \end{aligned} \tag{1}$$

with wave speed $c(x) > 0$ and homogeneous Dirichlet boundary conditions.

Spatial discretization by a standard $H^1(\Omega)$-conforming finite element methods (FEM) leads to the system of second-order ODE's:

$$\frac{\mathrm{d}^2}{\mathrm{d}\,t^2}\boldsymbol{M}\boldsymbol{u}(t) + \boldsymbol{K}\boldsymbol{u}(t) = \boldsymbol{f}(t), \tag{2}$$

where $\boldsymbol{M}$ and $\boldsymbol{K}$ denote the mass and stiffness matrix, respectively. Both are symmetric and positive definite while $\boldsymbol{M}$ can be assumed diagonal by applying mass-lumping. For the time discretization we use centered second-order finite differences, which leads to the standard fully explicit leapfrog (LF) method.

## 2 A posteriori error estimates

In [3], we proved the following a posteriori error estimates:

$$\begin{aligned} \max_{0\le n\le N} \|U^n \;&-\; u(\cdot,t_n)\|_{\mathcal{A}} \;\le\; \|\boldsymbol{e}(0)\|_{\mathrm{erg},\mathcal{A}} \\ &+\; C\left\{2\sum_{m=1}^{2N}\zeta^m + \max_{1\le n\le N}\varepsilon_0^n\right\}, \end{aligned}$$

where $\|\cdot\|_{\mathcal{A}}$ denotes the energy norm w.r.t. the spatial elliptic operator $\mathcal{A}$ in (1) whereas $\|\cdot\|_{\mathrm{erg},\mathcal{A}}$ denotes the full (wave) energy norm including the $L^2$-norm of the time derivative. Here $\varepsilon_0^n$ corresponds to a (fully computable) standard residual based error indicator in the energy-norm:

$$\begin{aligned} \varepsilon_0^n := \sum_{K\in\mathcal{M}_{V_n}} &\left\{ h_{\widehat{K}}^2 \left\| A_{V_n}U^n - \nabla\cdot(c(x)\nabla w)\big|_K \right\|_{\mathrm{L}_2(K)}^2 \right. \\ &\left. + \frac{1}{2}h_K \left\| [c(x)\nabla w\big|_K] \right\|_{\mathrm{L}_2(\partial K)}^2 \right\}. \end{aligned}$$

Similarly, the error indicators $\zeta^m$ are local in space and time and fully computable:

$$\begin{aligned} \zeta^m \;&=\; \int_{t_{\frac{m-1}{2}}}^{t_{m/2}} [(\mu_0^n + \vartheta_0^n(t))^2 \\ &+\; (\alpha^n + \mu_1^n + \vartheta_1^n(t) + \delta^n(t))^2]^{1/2}\,\mathrm{d}t, \end{aligned}$$

for $n = \lceil (m+1)/2 \rceil$ and $m = 1,\ldots,2N$. Here $\alpha^n$, $\mu_i^n$, $\vartheta_i^n$ and $\delta^n$ bound at time $t_n$ the errors due to mesh-change, local time-stepping, the time-discretization, and the data (source) approximation, respectively.

## 3 Adaptive FEM and local time-stepping

In the presence of local refinement, the time-step $\Delta t$ is constrained by the smallest elements in the mesh. To overcome that bottleneck, we now apply (stabilized) leapfrog based local time-stepping (LF-LTS) [1,2], which leads to the same LF-type three-term recurrence yet with $\boldsymbol{A} = \boldsymbol{M}^{-1}\boldsymbol{K}$ replaced by a modified matrix $\widetilde{\boldsymbol{A}}$, which involves computing additional but local matrix-vector multiplications in the refined region only.

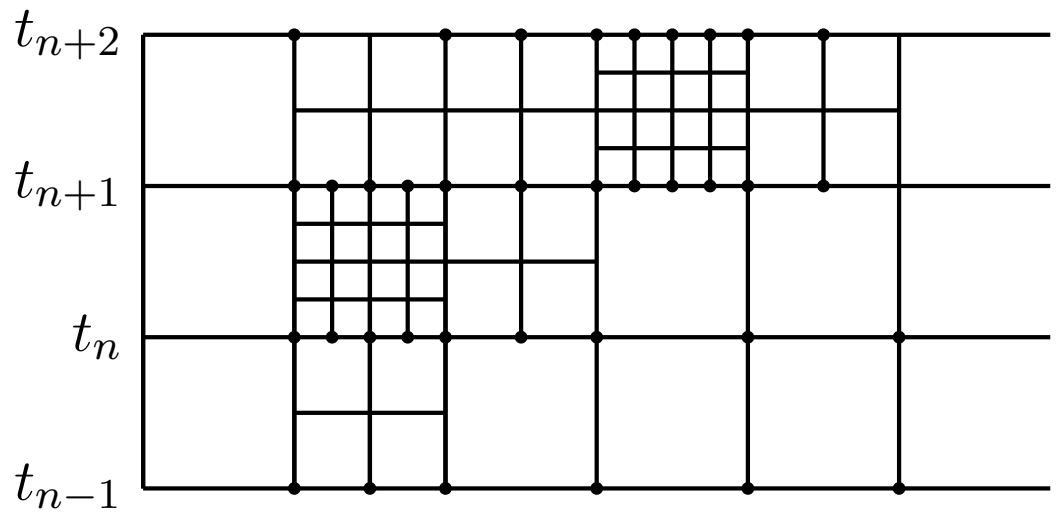

Figure 1: Space-time compatible meshes.

In [2], the stabilized LF-LTS method was shown to converge optimally w.r.t. the $L^2$-norm under a CFL stability condition independent of the coarse-to-fine mesh size ratio.

To include space-time mesh adaptivity, we denote by $\mathbb{V}_n$ the FE space at each time $t_n = n\Delta t$ together with its interpolation (or $L^2$-projection) operator, $\Pi_n$. On a time-evolving mesh, that is, on a sequence of FE spaces $\mathbb{V}_n$, $n = 0, 1, \ldots$ which may change from one time-step to the next, the LF-LTS-FEM method is:

$$\begin{aligned} U^0 &= \Pi_0 u_0 \\ U^1 &= \Pi_1 \left[ U^0 + (\Pi_0 v_0 + (F^0 - \tilde{\boldsymbol{A}}_0 U^0)\Delta t)\Delta t \right] \\ U^{n+1} &= \Pi_{n+1} \left[ 2U^n - \Pi_n U^{n-1} \right. \\ &\qquad \left. + (F^n - \tilde{\boldsymbol{A}}_n U^n)\Delta t^2 \right] \end{aligned}$$

for $n = 1, \ldots, N$. Before advancing with the usual LF step, we thus first project the current and previous FE solutions $U^n$ and $U^{n-1}$ onto the new mesh (or FE space) $V_{n+1}$ at $t_{n+1}$. In the simplest situation with only two local steps per global time-step $\Delta t$ and no stabilization, for instance, the modified stiffness matrix reduces to $\tilde{A}_n = A_n - \frac{\Delta t^2}{16} A_n \Pi_n^f A_n$, where $A_n$ denotes the stiffness matrix associated with $V_n$ and $\Pi_n^f$ restriction to the fine mesh.

During mesh change we restrict ourselves to compatible meshes, that is, we ensure that each element in the new mesh is either present in the old mesh or consists of a union of previous elements, and vice-versa – see Fig. 1. Thus, we ensure that local refinement never increases the numerical error. To minimize the inherent information loss during mesh coarsening, we estimate in advance the potential loss in accuracy by computing *coarsening pre-indicators* [4]

## 4 Numerical results

Here we consider a Gaussian pulse which propagates with constant wave speed $c \equiv 1$ across an L-shaped domain. The initial data is chosen such that at $t = 0$ the solution is a Gaussian centered at $(x, y) = (0.4, 0.6)$:

$$u_0(x, y) = e^{-600\left((x-0.4)^2 + (y-0.6)^2\right)}, \qquad (3)$$

while the initial velocity $v_0 \equiv 0$. At all boundaries we impose homogeneous Dirichlet boundary conditions and choose a uniform initial mesh with $h = 0.08$. We apply the space–time adaptive LF-LTS-FEM algorithm using the iFEM package for mesh refinement and coarsening.

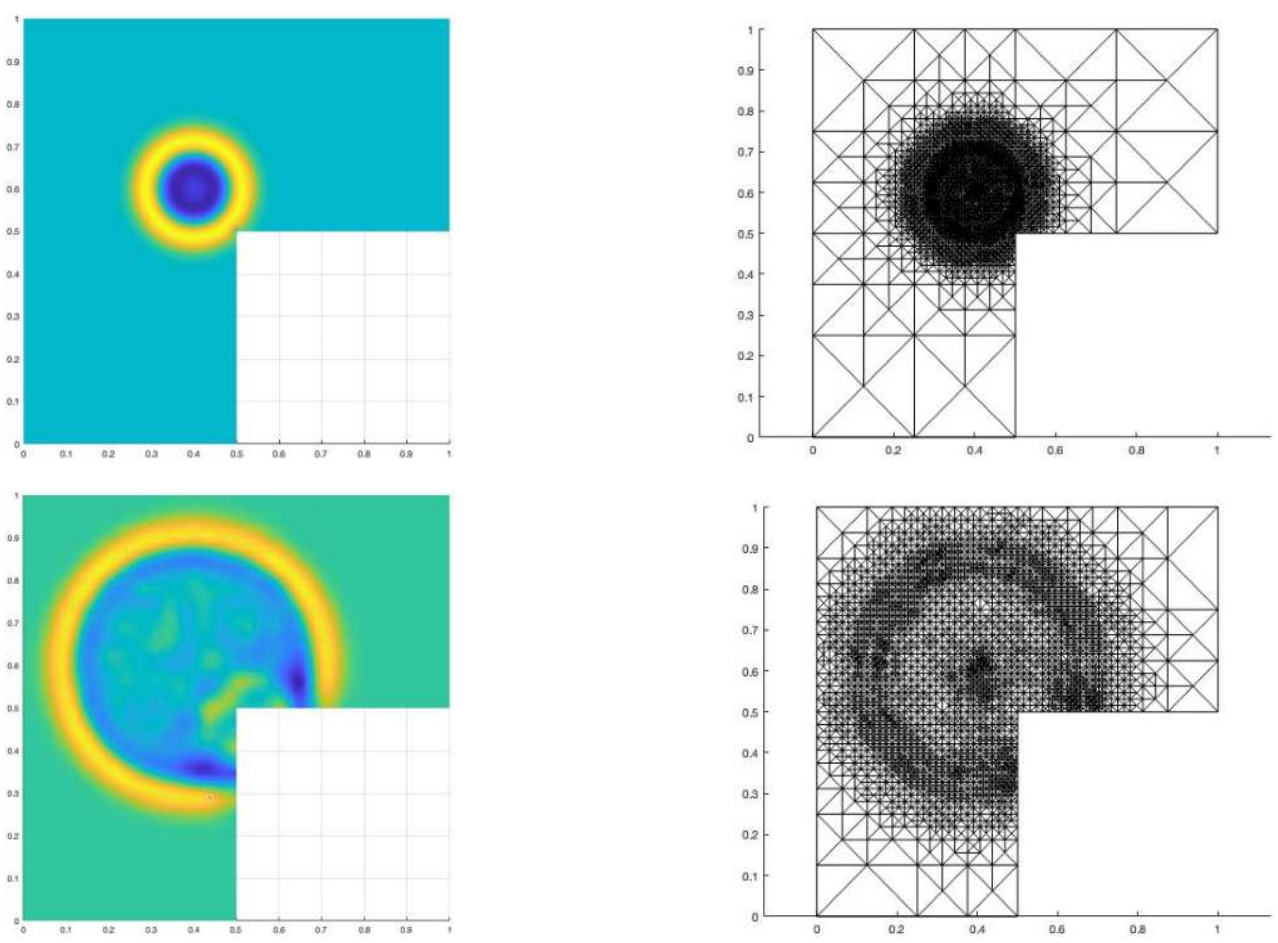

Figure 2: L-shaped domain: solution and adapted mesh at $t = 0.1$ and $t = 0.3$.

As shown in Fig. 2, the refined region of the mesh automatically adapts and tracks the propagating wave front.

# Multiscale analysis of wave propagation in the branched flow regime

**Josselin Garnier**[1], **Natacha Guegan-Fau**[1,*], **Sébastien Imperiale**[1],
[1]CMAP, Ecole polytechnique, Institut Polytechnique de Paris, 91120 Palaiseau, France
*Email: natacha.guegan-fau@polytechnique.edu

## Abstract

The aim of this work is to characterize the properties of a wave field propagating in a complex medium (in underwater acoustics or optics) in a specific regime called the branched flow regime. More precisely we are interested in studying solutions of the wave equation under the paraxial approximation in a random medium with random speckled initial field. We focus on the case where the medium fluctuations are smooth and have weak amplitudes. We also assume that the correlation length of the initial field is larger than the wavelength but smaller than the correlation length of the medium.

In this regime, numerical simulations show that the wave intensity forms branches of large intensity. A main aspect of the work is to justify this phenomenon by means of an asymptotic analysis of the solutions of the random paraxial equation. In particular by using separation of scales techniques and stochastic calculus, we analyze relevant statistical quantities of the random solutions. This approach enables us to derive formulae for the intensity and field correlation functions and to explain the observed phenomena.

## 1 Introduction

The branched flow regime is characterized by the formation of high-intensity branching patterns. These patterns result from smooth random fluctuations of the index of refraction larger than the wavelength (Figure 1). This phenomenon was first observed for electrons in a two-dimensional electron gas before being observed for oceanic waves or light propagation [1]. Our work focuses on the derivation of formulae for the field covariance function and the second-order moment of the intensity to obtain a better understanding of this phenomenon.

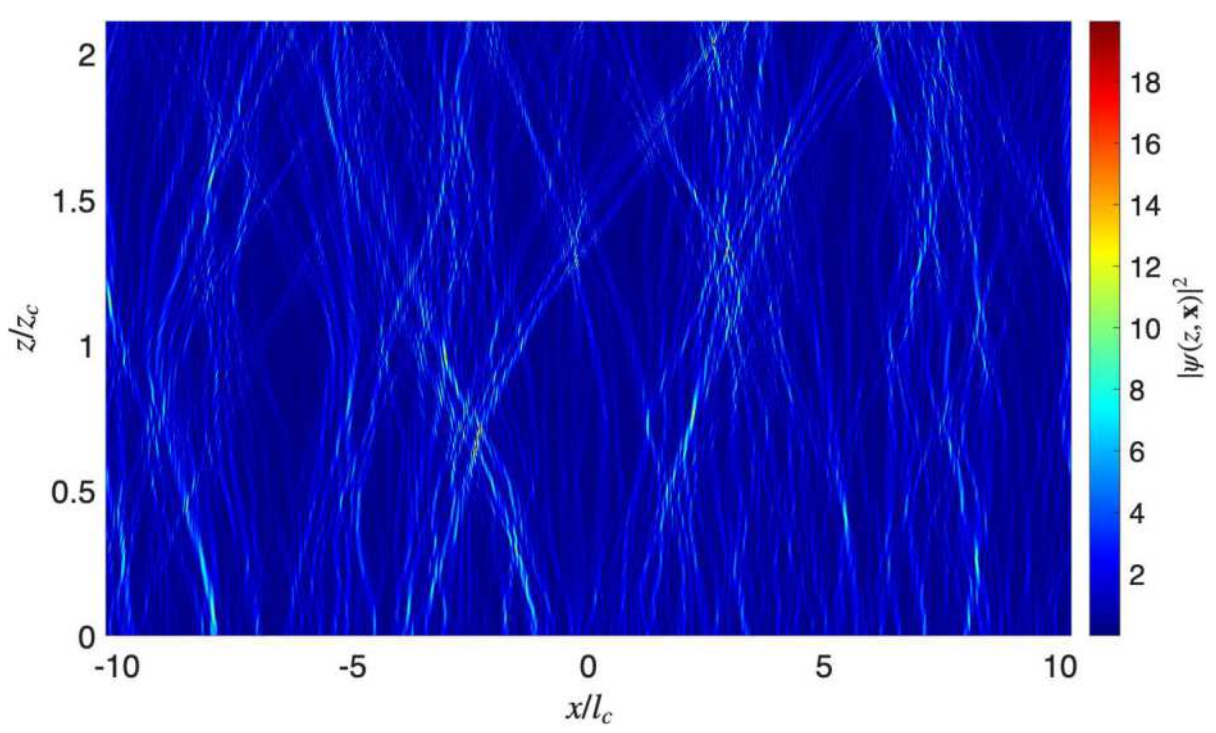

Figure 1: Measured intensity $|\psi(z,\mathbf{x})|^2$ with $\epsilon = \lambda/l_c = 0.005$, where $\lambda$ is the wavelength and $l_c$ the correlation length of the medium. The initial speckled field in the plane $z = 0$ has a correlation radius $\rho_0$ such that $\rho_0/\lambda = 5.85$. The variance of the random potential $\sigma^2$ satifies $\sigma^2\lambda^2 = 10^{-4}$. The characteristic propagation length is $z_c = l_c/(2\sigma^{2/3}\alpha^{2/3}) = 1.17 \times 10^4$. The simulation is done using the Split-Step Fourier Method ([2]).

## 2 Model

We consider a field $\psi$ defined by the relation $u(z,\mathbf{x}) = e^{ikz}\psi(z,\mathbf{x})$, where $u$ is the solution of the Helmholtz equation in a medium with heterogeneous fluctuations of the index of refraction and $k$ is wavenumber in free space. By neglecting the backscattered waves, we obtain the paraxial approximation. In addition, we assume a regime in which: (1) we have $\epsilon = \lambda/l_c \ll 1$ where $\lambda$ is the wavelength and $l_c$ is the correlation radius of the index of refraction. (2) The typical amplitude of the fluctuations of the squared index of refraction is $O(\epsilon^c)$, with $c > 0$. The scaled paraxial wave equation takes the form, for all $(\mathbf{x}, z) \in (\mathbb{R}^n \times \mathbb{R})$, with $n \in \{1, 2\}$,

$$\frac{\partial\psi^\epsilon}{\partial z}(z,\mathbf{x}) = i\epsilon\alpha\Delta_{\mathbf{x}}\psi^\epsilon(z,\mathbf{x}) - i\epsilon^{c-1}\mathcal{V}(z,\mathbf{x})\psi^\epsilon(z,\mathbf{x}).$$

The parameter $\alpha$ depends on the wavenumber in free space and on the homogeneous background index of refraction. The random potential $\mathcal{V}$ is related to the fluctuations of the index of refraction of the medium and is assumed to be a smooth, stationary, real-valued, zero-mean Gaussian process. The initial condition $\psi_{z=0}$ is a speckled field modeled as a complex Gaussian process with Wigner transform $W_0(\mathbf{k}) = \left(4\pi\rho_0^2\right)^{n/2} e^{-\rho_0^2|\mathbf{k}|^2}$, where $\rho_0$ is the correlation radius. It satisfies $\lambda \ll \rho_0 \ll l_c$ such that $\rho_0/l_c$ is of order $\epsilon^d$ for some $d \in (0,1)$. We denote $\mathbb{E}[\cdot]$ the expectation with respect to the distribution of the random process $\mathcal{V}$ and $\langle\cdot\rangle$ the expectation with respect to the distribution of the initial field $\psi_0$.

## 3 Field covariance function

The field covariance function, defined for $(\mathbf{x},\mathbf{y}) \in \mathbb{R}^n$ and $(z,\zeta) \in \mathbb{R}$, is

$$C^\epsilon(z,\mathbf{x},\zeta,\mathbf{y}) = \mathbb{E}\Big[\Big\langle \psi^\epsilon\Big(\frac{z}{\epsilon^b} + \epsilon^e\frac{\zeta}{2}, \mathbf{x} + \epsilon^d\frac{\mathbf{y}}{2}\Big) \times \bar{\psi}^\epsilon\Big(\frac{z}{\epsilon^b} - \epsilon^e\frac{\zeta}{2}, \mathbf{x} - \epsilon^d\frac{\mathbf{y}}{2}\Big)\Big\rangle\Big].$$

where $O(\epsilon^e)$, with $e \geq 0$, characterizes the offset in $z$ relative to the medium correlation length and $O(\epsilon^{-b})$, with $b \geq 0$ denotes the propagation distance. The relevant scaling is $b = 1-d$, $c = \frac{3}{2}(1-d)$, and $e = \frac{1}{2}(3d-1)$, with $d \in (1/3,1)$. The scaling for $e$ is chosen so that the correlation between two points separated by a distance $\epsilon^e\zeta$ in the $z$–direction is not zero (the points are not too far) while still depending on $\zeta$ (the points are not too close). In this regime,

$$\lim_{\epsilon\to 0} C^\epsilon(z,\mathbf{x},\zeta,\mathbf{y}) = C(z,\mathbf{x},0,\mathbf{y})\,\mathbb{E}[e^{-i\mathcal{V}(z,\mathbf{x})\zeta}],$$

where $C(z,\mathbf{x},0,\mathbf{y})$, derived in [3], is the solution of a PDE that depends on the autocorrelation function of the index of the medium. A comparison between theoretical and numerical results shows good agreement (Figure 2).

## 4 Second-order moment of the intensity

The second-order moment of the intensity is defined for $(\mathbf{x},\mathbf{y}) \in \mathbb{R}^n$ and $(z,\zeta) \in \mathbb{R}$ as

$$C_I^\epsilon(z,\mathbf{x},\zeta,\mathbf{y}) = \mathbb{E}\left[\left\langle\left|\psi^\epsilon\Big(\frac{z}{\epsilon^b} + \epsilon^f\frac{\zeta}{2}, \mathbf{x} + \epsilon^d\frac{\mathbf{y}}{2}\Big)\right|^2 \times \left|\psi^\epsilon\Big(\frac{z}{\epsilon^b} - \epsilon^f\frac{\zeta}{2}, \mathbf{x} - \epsilon^d\frac{\mathbf{y}}{2}\Big)\right|^2\right\rangle\right].$$

The relevant scaling is $f = 2d-1$ with $d \in (1/2,1)$. In the limit $\epsilon \to 0$, the correlation is

$$\lim_{\epsilon\to 0} C_I^\epsilon(z,\mathbf{x},\zeta,\mathbf{y}) = \frac{1}{(2\pi)^{2n}}\left[\int_{\mathbb{R}^n} W_0(\mathbf{u})d\mathbf{u}\right]^2 \left[\int_{\mathbb{R}^{2n}} \Pi\left(z,\mathbf{0},\mathbf{u},\mathbf{v}\right)\left(1 + e^{-i(2\alpha\zeta\mathbf{u}-\mathbf{y})\cdot\mathbf{v}}\right) d\mathbf{u}d\mathbf{v}\right]$$

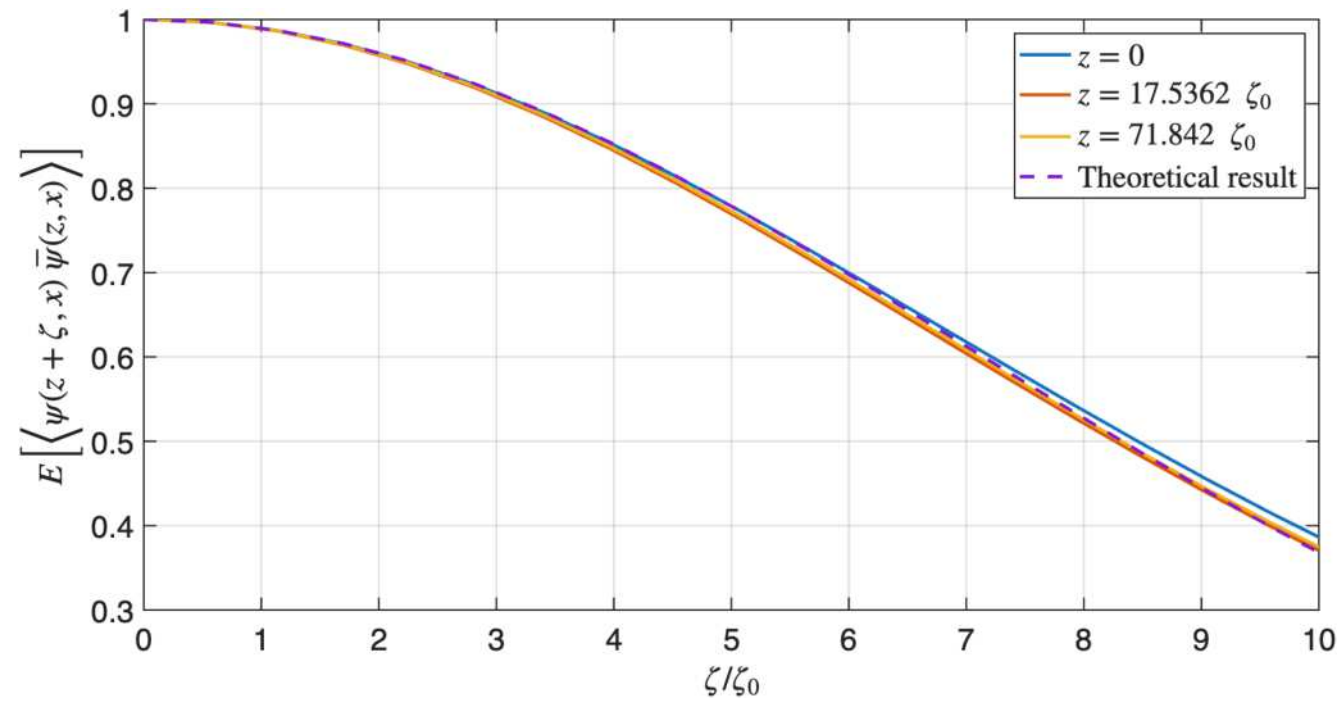

Figure 2: Field Covariance function $\zeta \mapsto \mathbb{E}\big[\big\langle\psi(z+\zeta,x)\,\overline{\psi}(z,x)\big\rangle\big]$ plotted for various values of $z$. The parameters are $\epsilon = 0.005$, $d = 2/3$, $\rho_0/\lambda = \epsilon^{d-1} = 5.85$, $\sigma^2\lambda^2 = 10^{-4}$ and $\zeta_0/l_c = \epsilon^{(3d-1)/2} = 0.07$. The results are obtained by solving the paraxial equation and averaging over 500 realizations of the random medium and 500 realizations of the initial field (250 000 realizations in total). They are compared with the theoretical result derived in Section 3.

with $\Pi\left(z,q,u,v\right)$ being the solution of the equation (for simplicity, $n=1$)

$$\partial_z\Pi = -2\alpha v\partial_q\Pi + \frac{1}{4}\left(\Gamma(0)+\Gamma(q)\right)\partial_u^2\Pi + \left(\Gamma(0)-\Gamma(q)\right)\partial_v^2\Pi$$

where $\Gamma(x) = \int_{\mathbb{R}} \mathbb{E}\left[\partial_x\mathcal{V}(0,0)\partial_x\mathcal{V}(z,x)\right]dz$. The initial condition is of the form

$$\Pi(0,q,u,v) = \frac{\rho_0}{\pi}\exp\left(-2\rho_0^2u^2\right)\exp\left(-\frac{1}{2}\rho_0^2v^2\right).$$

For the second moment of the intensity, the order of $\zeta$ is larger than that used in the field covariance function. In the latter case, the phase modulation term $\mathbb{E}\left[e^{-i\mathcal{V}(z,\mathbf{x})\zeta}\right]$ introduces destructive interference, which produces a loss of coherence. In contrast, this oscillatory term does not appear in the intensity, allowing $\zeta$ to be chosen larger. The theoretical prediction for intensity is validated against numerical simulations.

# On the surface Bloch waves in truncated periodic media

Bojan B Guzina[1,*], Shixu Meng[2], Prasanna Salasiya[1], Long Nguyen[1]
[1]Civil, Environmental, and Geo Engineering, University of Minnesota, Twin Cities, USA
[2]Mathematical Sciences, The University of Texas at Dallas, Richardson, USA
*Email: guzin001@umn.edu

**Abstract**

Much like their counterparts in homogeneous elastic solids, waves in periodic media can be broadly classified into Floquet-Bloch body waves, and evanescent surface waves. Our goal is to elucidate the latter boundary layers, termed surface Bloch (SB) waves, affiliated with rational surface cuts and homogeneous Neumann data. To this end we adopt a 2D scalar wave equation with periodic coefficients as a test bed and develop a unit cell-of-periodicity-based, reduced order model of the SB waves that is capable of describing their dispersion, waveforms, and "skin depth". The centerpiece of our analysis is a quadratic eigenvalue problem (QEP) for the unit cell of periodicity that seeks the complex wavenumber quantifying the evanescence away from the cut plane given (i) the excitation frequency, and (ii) wavenumber in the direction of the cut plane. In this way the sought boundary layer is obtained by superposition of the evanescent QEP eigenstates, whose relative amplitudes are obtained by imposing the homogeneous boundary condition. With the QEP eigenspectrum at hand, evaluation of an SB wave entails only a low-dimensional eigenvalue problem. This feature caters for rapid exploration of the effect of (periodic) surface undulations, and so enables manipulation of SB waves via optimal design of the surface cut.

## 1 Introduction

With reference to a Cartesian coordinate system with an orthonormal vector basis $\boldsymbol{e}_j$ $(j=\overline{1,2})$, we consider the 2D scalar wave equation

$$\nabla\cdot(G(\boldsymbol{x})\nabla u)+\omega^2\rho(\boldsymbol{x})u \;=\; 0, \quad \boldsymbol{x}\in\mathbb{R}^2 \qquad (1)$$

at frequency $\omega$, where $G$ and $\rho$ are $Y$-periodic, and $Y=\{\boldsymbol{x}:\,-\frac{1}{2}<x_i<\frac{1}{2}, j=\overline{1,2}\}$. Recalling the Floquet-Bloch theorem, we write $u(\boldsymbol{x}) = \phi(\boldsymbol{x})e^{i\boldsymbol{k}\cdot\boldsymbol{x}}$ where $\phi(\boldsymbol{x})$ is $Y$-periodic and $\boldsymbol{k} = (k_1,k_2)$ is the wave vector. In this setting we pursue the Bloch-like boundary layers, termed SB waves, propagating in a *semi-infinite periodic medium* bounded by plane $\mathcal{S}$ with unit outward normal $\boldsymbol{\nu}=-\boldsymbol{e}_2$. Assuming homogeneous Neumann data on $\mathcal{S}$, we seek $\mathsf{u}(\boldsymbol{x})$ satisfying

$$\begin{aligned}\nabla\cdot(G(\boldsymbol{x})\nabla\mathsf{u})+\omega^2\rho(\boldsymbol{x})\mathsf{u} &= 0, \quad \boldsymbol{x}\cdot\boldsymbol{e}_2>\mathfrak{d}\\ \boldsymbol{\nu}\cdot\big(G(\boldsymbol{x})\nabla\mathsf{u}\big)|_{\mathcal{S}} &= 0\end{aligned}, \qquad (2)$$

where $\mathfrak{d}\in(-\frac{1}{2},\frac{1}{2})$ signifies the depth of cut within the unit cell.

## 2 Quadratic eigenvalue problem

As examined in [1,2], the QEP relevant to evaluating $\mathsf{u}(x)$ can be formulated as a task of seeking the eigenvalues $\kappa:=k_2$ and eigenfunctions $\phi$ (subject to periodic boundary conditions on $\partial Y$) for given $\omega$ and wavenumber component $k_1$. In weak form, this problem can be written as

$$\begin{aligned}&\kappa^2\!\int_Y G\phi\overline{\psi}\,\mathrm{d}Y_{\boldsymbol{x}}-\kappa\left\{\mathrm{i}\!\int_Y G\left(\frac{\partial\phi}{\partial\xi_2}\overline{\psi}-\phi\frac{\partial\overline{\psi}}{\partial\xi_2}\right)\mathrm{d}Y_{\boldsymbol{x}}\right\}\\ &+\left\{k_1^2\!\int_Y G\phi\overline{\psi}\,\mathrm{d}Y_{\boldsymbol{x}}-\mathrm{i}k_1\!\int_Y G\left(\frac{\partial\phi}{\partial\xi_1}\overline{\psi}-\phi\frac{\partial\overline{\psi}}{\partial\xi_1}\right)\mathrm{d}Y_{\boldsymbol{x}}\right.\\ &\left.+\int_Y\Big(G\nabla\phi\cdot\nabla\overline{\psi}-\omega^2\rho\phi\overline{\psi}\Big)\,\mathrm{d}Y_{\boldsymbol{x}}\right\}=0 \quad \forall\psi\in H^1_p(Y).\end{aligned} \qquad (3)$$

Assuming sufficient conditions for discreteness of the resulting quadratic eigenspectrum [1], we denote the solution of (3) by $\{\kappa_n\in\mathbb{C},\phi_n\in H^1_p(Y)\}_{n=1}^{\infty}$. In [2], it is shown that the eigenvalues $\kappa_n$ relevant to evaluation of $\mathsf{u}(x)$ are contained in the *essential semi-strip*

$$\Re(\kappa_n)\in(-\pi,\pi], \qquad \Im(\kappa_n)>0, \qquad (4)$$

which ensures evanescence of the quadratic eigenstates with depth, $x_2-\mathfrak{d}$, into the half-space. For simplicity of discussion, by $\{\kappa_n,\phi_n\}_{n=1}^{\infty}$ we hereon refer only to the pairs satisfying (4) and we arrange $\kappa_n$ in the order of increasing $\Im(\kappa_n)$.

## 3 Model

QEP (3) allows us to seek a reduced-order model of SB waves via superposition of quadratic eigenstates

$$\tilde{\mathrm{u}}(\boldsymbol{x}) = \sum_{n=1}^{N} a_n \tilde{\phi}_n(\boldsymbol{x}) e^{\mathrm{i}\boldsymbol{k}_n\cdot\boldsymbol{x}}, \quad \boldsymbol{k}_n := (k_1, \kappa_n) \tag{5}$$

that is limited to $N$ deepest-penetrating modes, where $a_n \in \mathbb{C}$ are the modal amplitudes to be computed form the boundary condition. On noting that the *scaled* modal fluxes across $\mathcal{S}$

$$\begin{aligned}\psi_n(x_1) &:= e^{-\mathrm{i}k_1x_1}\, \boldsymbol{\nu}\cdot\big(G(\boldsymbol{x})(\nabla+\mathrm{i}\boldsymbol{k}_n)\phi_n\, e^{\mathrm{i}\boldsymbol{k}_n\cdot\boldsymbol{x}}\big)|_{\mathcal{S}} \\ &= -\Big(G(\boldsymbol{x})\Big(\frac{\partial}{\partial x_2}+\mathrm{i}\kappa_n\Big)\phi_n\, e^{\mathrm{i}\kappa_n x_2}\Big)\Big|_{x_2=\mathfrak{d}}\end{aligned} \tag{6}$$

are $(-\frac{1}{2}, \frac{1}{2})$-periodic and so can be expanded in Fourier series as $\psi_n(x_1) = \sum_{m=-\infty}^{\infty} \hat{\psi}_{nm} e^{\mathrm{i}2\pi m x_1}$, boundary condition $e^{-\mathrm{i}k_1x_1} \times \boldsymbol{\nu}\cdot\big(G(\boldsymbol{x})\nabla\tilde{\mathrm{u}}\big)|_{\mathcal{S}} = 0$ due to (2) can be conveniently restricted to $2M+1$ Fourier modes ($m = \overline{-M, M}$) and so reduced to a low-dimensional algebraic problem

$$\boldsymbol{\Psi}(k_e)\boldsymbol{a} \;=\; \boldsymbol{0}, \quad k_e := k_1 \in (-\pi, \pi] \tag{7}$$

where $\boldsymbol{\Psi}$ is a complex-valued $(2M+1)\times N$ matrix collecting $\hat{\psi}_{mn}$; eigenvalue $k_e$ is a prescribed surface wavenumber providing input to QEP (3), and $\boldsymbol{a} = \boldsymbol{a}(k_e) \in \mathbb{C}^N$ is the vector of SB modal amplitudes $a_n$. In the sequel, we refer to (7) as the *surface-wave eigenproblem* and assume $(2M+1) \geqslant N$ which caters for (but does not guarantee) having an overdetermined system.

To expose the wavenumber $k_e$ and the vector of modal amplitudes $\boldsymbol{a}$ featured by an SB wave, we next consider singular value decomposition $\boldsymbol{\Psi}(k_e) = \boldsymbol{U}\boldsymbol{\Lambda}\overline{\boldsymbol{V}}$ where $\boldsymbol{U}$ and $\boldsymbol{V}$ are the unitary matrices of the left and right singular vectors, and $\boldsymbol{\Lambda}$ is the diagonal matrix of singular values $\sigma_n$ ($n = \overline{1, N}$) arranged in descending order. In this setting, (7) permits a non-trivial solution for $k_e = k_e^\star$ where the condition number of $\boldsymbol{\Psi}$ becomes unbounded, i.e.

$$C_{\boldsymbol{\Psi}}(k_e) := \frac{\sigma_1}{\sigma_N} \to \infty \quad \text{as} \quad k_e \to k_e^\star. \tag{8}$$

In this way, for each $k_e^\star$ one obtains a reduced-order model of the surface Bloch wave as

$$\hat{\mathrm{u}}^\star(\boldsymbol{x}) \;=\; \Big(\sum_{n=1}^{N} a_n^\star \phi_n(\boldsymbol{x})\, e^{\mathrm{i}\kappa_n x_2}\Big) e^{\mathrm{i}k_e^\star x_1}, \tag{9}$$

where $a_n^\star$ are the components of $\boldsymbol{a}(k_e^\star)$, and pairs $\{\kappa_n, \phi_n\}_{n=1}^N$ solve QEP (3) with $\omega$ and $k_1 = k_e^\star$ as input parameters.

## 4 Results

With reference to Fig. 1(a), the surface cut is made at depth $\mathfrak{d} = -\frac{1}{2}$ corresponding to an interface between the two "rows" of unit cells. To compute the dispersion diagram, for each input frequency $\omega$ and a set of trial wavenumbers $k_e^{(q)} \in [0, \pi]$ ($q = 1, 2, \ldots$), we (i) solve QEP (3), see Fig. 1(b) for an example eigenmode; (ii) adopt ROM (5) with $N = 10$; (iii) evaluate the condition numbers $C_{\boldsymbol{\Psi}}(k_e^{(q)})$ with $M = 11$, and (iv) identify the wavenumber $k_e^\star = k_e^\star(\omega)$ of the surface Bloch wave via criterion (8). In Fig. 1(c), points on the first two branches of the *surface-wave* dispersion diagram identified in this way are indicated by bullets. For transparency, we also plot example distributions $C_{\boldsymbol{\Psi}}(k_e)$ versus $k_e \in [0, \pi]$ at frequencies $\omega \in \{0.6, 1.3, 2\}$. Each of these distributions each feature a sharp peak, locating a point on the dispersion diagram via (8).

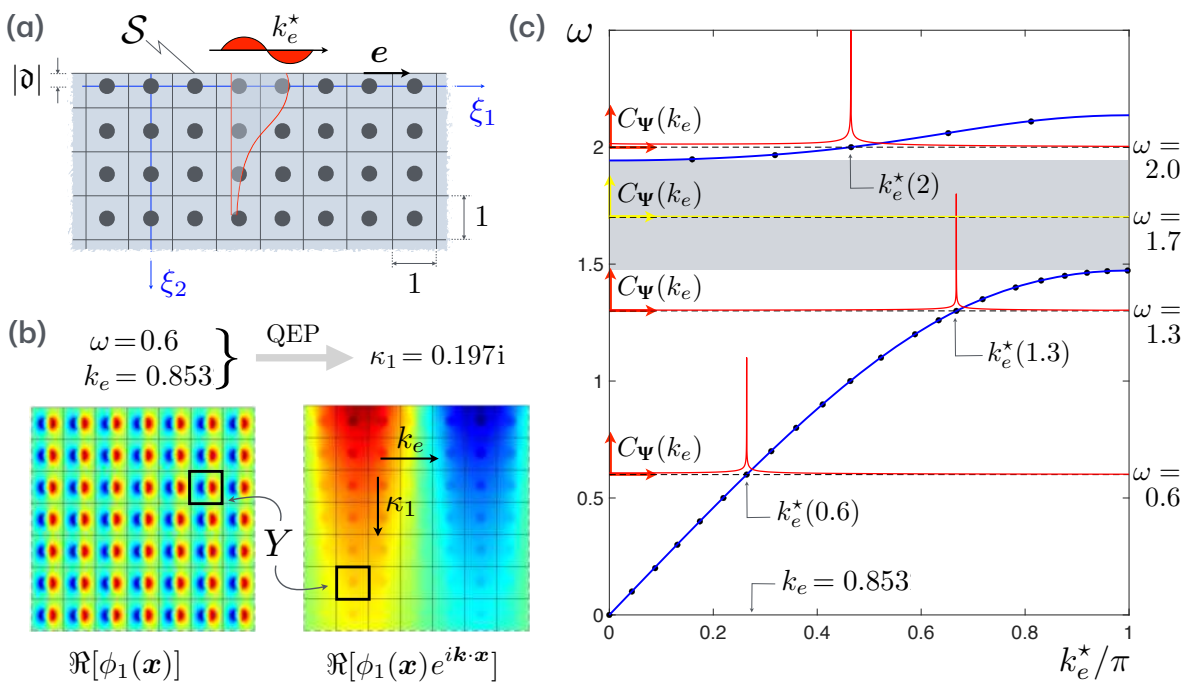

Figure 1: Dispersion of SB waves: (a) evanescent Bloch wave (5) propagating in direction $\boldsymbol{e}_1$ along the traction-free surface $\mathcal{S} = \{\boldsymbol{x} : x_2 = \mathfrak{d}\}$ of a periodic half-space; (b) example QEP solution, fundamental mode: $\phi_1(\boldsymbol{x})$ and $\phi_1(\boldsymbol{x})e^{\boldsymbol{k}\cdot\boldsymbol{x}}$ which illustrate the genesis of SB waves; (c) first two branches of the *surface-wave* dispersion diagram for $\mathfrak{d} = -0.5$. The unit cell features a circular inclusion with $(G, \rho) = (5, 0.1)$ and radius $r = 0.2$ centered in a homogeneous matrix with $(G, \rho) = (1, 1)$.

# Domain Truncation for Quasilinear Wave Equations

**Thomas Hagstrom**[1,*]
[1]Department of Mathematics, Southern Methodist University, Dallas, TX USA
*Email: thagstrom@smu.edu

## Abstract

We consider the problem of developing convergent domain truncation procedures for the numerical simulation of quasilinear wave equations. We focus first on methods which are based on mappings of the unbounded domain to a bounded domain which try to maintain a fixed sampling in the mapped region. These include the use of characteristic coordinates and the introduction of carefully-constructed artificial dissipation. Some preliminary experiments for a $2+1$-dimensional model problem are presented.

## 1 Introduction

Although some open problems remain, provably spectrally-convergent domain truncation procedures are now available for most of the linear hyperbolic equations arising in physical applications. For nonlinear problems, on the other hand, very little has been done to produce methods which can be proved to converge at reasonable cost. Our goal in this work is to critically examine some methods which we believe can work well, at least for small amplitudes, and to take steps towards their rigorous analysis.

## 2 Model and Methods

Our specific interest in this work is second-order scalar hyperbolic problems described by a space-time metric appearing as a perturbation of the flat space Minkowski metric:

$$\frac{\partial^2 u}{\partial t^2} - \Delta u = \sum_{i,j=0}^{d} g_{ij}(Du)\frac{\partial^2 u}{\partial x_i \partial x_j}, \qquad (1)$$

where on the right-hand side we identify $x_0 \equiv t$, $x_1 = r$, a radial variable in an appropriate coordinate system, and $x_j$, $j = 2, \ldots d$ angular variables. We assume an unbounded computational domain which we wish to truncate at $\Gamma$.

### 2.1 Mappings

A straightforward approach to this problem would be to compactify the exterior domain, or more directly discretize the system in the exterior using only finitely many degrees-of freedom. As is well-known, a naive application of this idea fails as waves become underresolved and artificially reflect back through $\Gamma$. However, techniques have been proposed to circumvent this issue. As they work directly with the governing equations, such methods can be applied to nonlinear systems as easily as to linear systems.

**Supergrid**: By far the simplest approach is to rely on carefully-chosen numerical dissipation. A workhorse in computational fluid dynamics, Appelo and Colonius [2] have shown that using high-order methods and high-order dissipation excellent accuracy can be obtained, even for problems where the best linear methods are not applicable due to the presence of reverse modes. Introduce a coordinate stretching in the variable $r$ and high-order finite difference approximations to the spatial derivatives with the addition of a high-order damping term. For linear problems it is straightforward to establish the stability of this method, but for nonlinear problems the creation of large gradients in the layer could render the problem ill-posed.

We have carried out some initial numerical experiments in two space dimensions with the supergrid approach. Here we took a nonlinearity verifying Alinhac's strong null condition [1], which guarantees all-time existence for small data, and an eighth order supergrid layer. Results for a well-refined mesh are shown below - the layers were terminated at $\rho = z_1$ and errors are computed on $(1,3)$ by comparison with a long-domain solution. The results, shown in the figure below, do yield rapid convergence for medium times, but there is late time error growth. The layers are rather thick - instabilities were encountered for smaller layers.

**Characteristic Coordinates**: An alternative approach which can be applied to problems with the specific structure of (1) is based on the change of the time variable to so-called null coordinates, followed by a similar stretching of the $r$ coordinate to a coordinate $\rho$. A motivation for this approach in the linear case follows from the progressive wave expansion:

$$u \sim \sum_j \frac{f_j(t-r,\theta,\phi)}{r^j}.$$

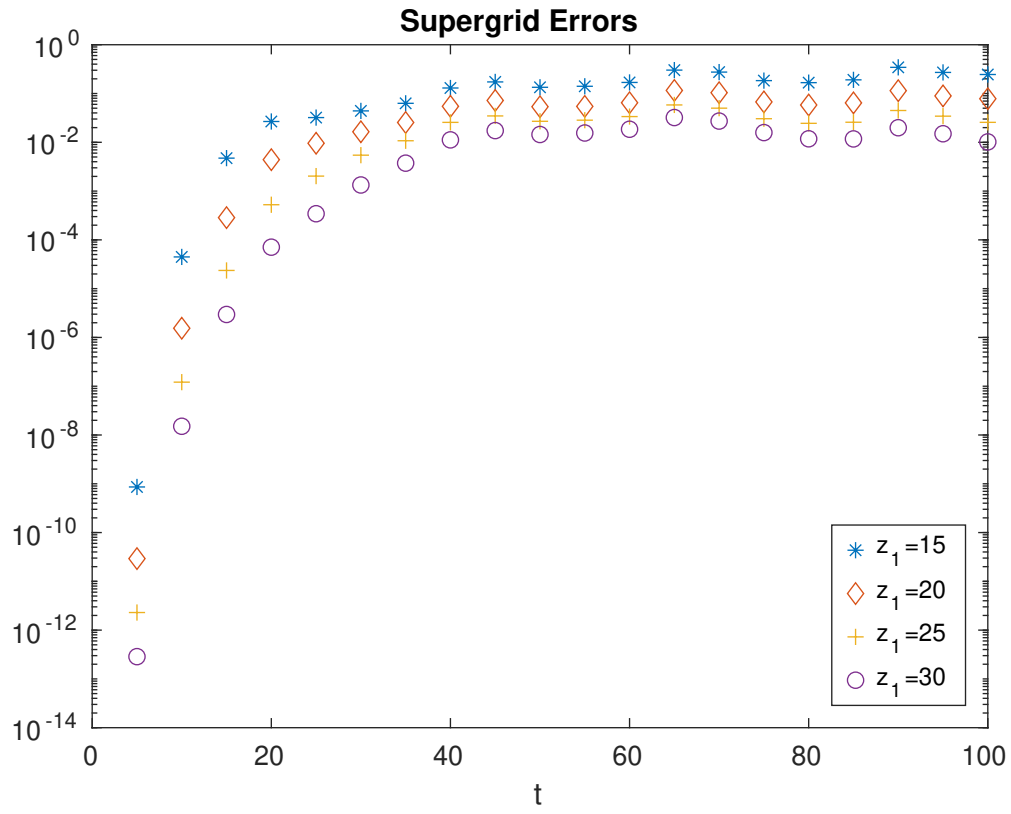

Setting $\tau = t - r$ then the mapped solution $U(\tau, r, \theta, \phi)$ may have very smooth dependence on $r$ out to $\infty$ and so a mapping $r = R(\rho)$, $r \to \infty$ as $\rho \to L$ may yield a smooth solution. An advantage of this approach is that the solution in the far field is directly computed. A disadvantage is that the mapped system is no longer hyperbolic. For the d'Alembertian:

$$\frac{\partial^2 u}{\partial r \partial \tau} + \frac{1}{r}\frac{\partial u}{\partial \tau} = \frac{1}{2}\frac{\partial^2 u}{\partial r^2} + \frac{1}{r}\frac{\partial u}{\partial r} + \frac{1}{2r^2}\Delta_S u.$$

A look at the frozen coefficient problem shows that the phase speeds in the angular directions are unbounded. In [8] a hyperboloidal layer is introduced for which the characteristic form is only approached as $r \to \infty$. Although the system is nonstandard, with an appropriate choice of mappings it seems possible to bound the speeds. There is some discussion in [4] about the numerical issues involved in approximating these systems, but no rigorous theory has appeared.

In [3] a successful application of the use of characteristic coordinates to truncate semilinear problems is presented. For the quasilinear problems considered here, however, the mappings themselves become nonlinear. An experimental demonstration of this idea has appeared in [6]. From the perspective of numerical implementation a nonlinear generalization of the hyperboloidal layer would be of interest - a direct use of the linear mapping would typically lead to instabilities for nonlinear problems.

## 3 Alternative Approaches

We also plan to consider two alternative approaches to the mapping techniques discussed above. The first is the nonlinear extension of high-order radiation conditions based on progressive wave expansions [5]. We believe that for systems obeying the null condition such expansions can be constructed. The second is the use of phase space filters, which have been applied to the nonlinear Schrödinger equation in [7].

## 4 Acknowledgments

This material is based upon work supported by the National Science Foundation under Grant DMS-2309687. Any opinions, findings, and conclusions or recommendations expressed in this material are those of the authors and do not necessarily reflect the views of the NSF.

# Wave-Based Stability Analysis of Couette and Gravity-Driven Bilayer Film Flows

**M. Al Ahmadi Hammou**[1,*], **M. Scholle**[1], **P. H. Gaskell**[2]
[1]Institute for Flow in Additively Manufactured Porous Media (ISAPS), Heilbronn University, Heilbronn, Germany
[2]Department of Engineering, Durham University, Durham, UK
*Email: meriem.al-ahmadi-hammou@hs-heilbronn.de

**Abstract**

The stability of thin liquid films flowing over patterned substrates is considered, focusing on shear- and gravity-driven bilayer flows. Small interfacial disturbances to the base flow are either selectively amplified or damped, depending on their wavelength and the underlying surface geometry. A linear wave analysis is employed, representing perturbations as travelling waves, which allows the spatial structure and temporal growth to be captured by a dispersion relation. Analysis of the resulting dispersion relation demonstrates how substrate geometry can selectively shift instability thresholds, providing a predictive tool for controlling interfacial instabilities in multilayer film flows.

## 1 Introduction

Thin liquid films flowing over patterned substrates are ubiquitous in natural and industrial processes, such as coating, lubrication, and microfluidic devices. In bilayer flows, small disturbances on the layer-separating interface and at a free surface can propagate as waves, whose growth or decay determines the stability of the system. These interfacial waves are influenced by the base flow, layer thicknesses, viscosities, and substrate geometry. A linear wave-based analysis provides a predictive framework by relating wave frequency to wavenumber via a dispersion relation, allowing identification of conditions under which the flow remains stable or becomes unstable.

## 2 Model

Assuming the fluid layers are thin compared to the substrate, $x_2 = b(x_1)$, the flows are essentially two-dimensional (2D), as shown in Fig. 1. Both feature a layer-separating interface with surface tension $\sigma_1$ and layer thicknesses $H_1$ and $H_2$. In the shear-driven case, the upper boundary is a rigid plate moving at speed $U$, while in the gravity-driven case it is replaced by a free surface (surface tension $\sigma_2$). Each is free to deform for flow over the lower fixed rigid substrate. The two layers are incompressible, immiscible Newtonian liquids with viscosity $\eta_1, \eta_2$ and density $\varrho_1, \varrho_2$. The flow structure and interface disturbances are analyzed for both arrangements.

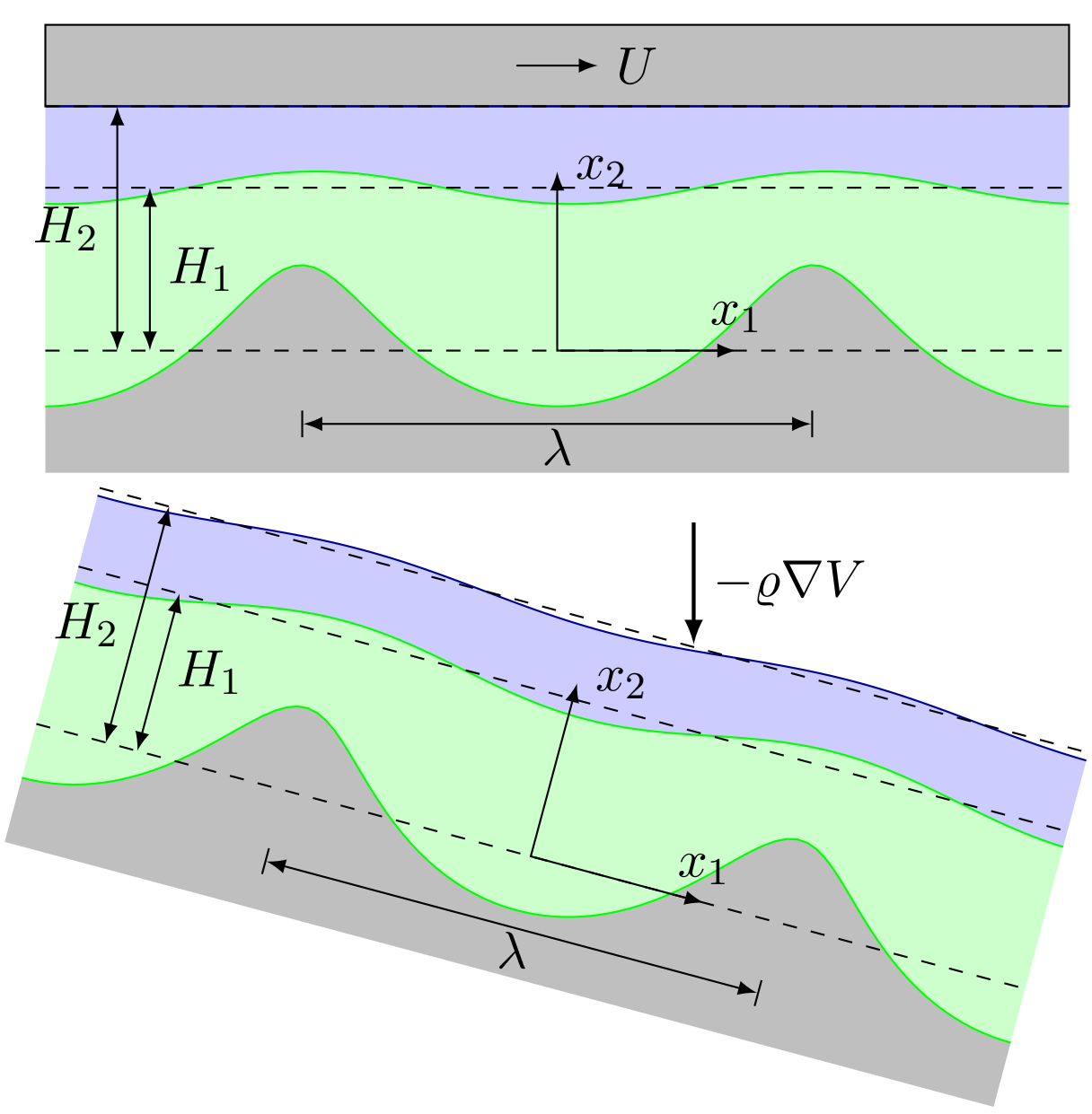

Figure 1: Schematics of bilayer, shear–driven (top) and gravity–driven (bottom) film flow in the presence of a rigid and stationary lower periodically patterned substrate.

An elegant approach exists for exploring such 2D flows [1], based on representing the governing equations in first-integral form:

$$4\frac{\partial^2 X}{\partial \bar{\xi}^2} + \varrho\frac{u^2}{2} = 0\,, \tag{1}$$

$$\frac{P}{\varrho} + \frac{\bar{u}u}{2} + g\mathrm{e}^{\mathrm{i}\alpha}\,\Im\xi - \mathrm{i}\left[\frac{\partial \Psi}{\partial t} - \frac{4\eta}{\varrho}\frac{\partial^2 \Psi}{\partial \bar{\xi}\partial \xi}\right] = \frac{4}{\varrho}\frac{\partial^2 X}{\partial \bar{\xi}\partial \xi},$$

in terms of the complex variable $\xi := x_1 + \mathrm{i}x_2$, its complex conjugate $\bar{\xi}$, the pressure $P$, the complex velocity $u := -2\mathrm{i}\,\frac{\partial \Psi}{\partial \bar{\xi}}$, the streamfunction $\Psi$, and a potential $X$; the latter is a generalisation of Airy's stress function [2] and complex–valued. On neglecting the quadratic terms in (1), steady state solutions can be obtained by twofold integration and determination of the integration functions by the boundary and interface conditions [3].

### 3 Stability Analysis and Results

Having determined steady–state base solutions for the stream function, $\Psi = \Psi_0\left(\xi, \bar{\xi}\right)$, the complex potential, $X = X_0\left(\xi, \bar{\xi}\right)$ and the pressure $P = P_0\left(\xi, \bar{\xi}\right)$, small perturbations $\psi$, $\zeta$ and $p$ are added as follows:

$$\Psi\left(\xi, \bar{\xi}, t\right) = \Psi_0\left(\xi, \bar{\xi}\right) + \psi\left(\xi, \bar{\xi}, t\right) , \quad (2)$$
$$X\left(\xi, \bar{\xi}, t\right) = X_0\left(\xi, \bar{\xi}\right) + \zeta\left(\xi, \bar{\xi}, t\right) , \quad (3)$$
$$P\left(\xi, \bar{\xi}, t\right) = P_0\left(\xi, \bar{\xi}\right) + p\left(\xi, \bar{\xi}, t\right) , \quad (4)$$

the time evolution of which is investigated in terms of a Taylor expansion of the field equations up to and including terms of linear order. Additionally, layer-separating interface and free–surface perturbations have to be considered. After substituting (2)–(4) into (1), the resulting first integral form of the Orr–Sommerfeld equation leads to a second integral form. [4]. For the stability analysis the perturbations are then assumed to possess a wave-like character and are expressed in terms of normal modes. This procedure yields an implicit dispersion relationship relating the wavenumber $k$ to the circular frequency $\omega$. Its roots determine whether the flow is linearly stable or unstable.

Steady base-flow solutions and stability results have been obtained for both flow configurations at various Reynolds numbers:

$$\mathrm{Re} = \max_i \left( \frac{\varrho_i\, u(H_i)\, \lambda}{2\pi\, \eta_i} \right), \quad i = 1, 2.$$

The corresponding stability characteristics for gravity-driven bilayer flow over a moderately sinusoidal substrate, given by:

$$2\pi x_2/\lambda = a\cos(2\pi x_1/\lambda),$$

are shown in Fig. 2. The top figure is the dispersion relation (Re($\omega$) vs. $k$), while the bottom one reveals the corresponding growth-rate curve (Im($\omega$) vs. $k$). Positive growth rates indicate that this flow is unconditionally unstable, even at the smallest Reynolds number considered, Re = 0.02, since stability requires the growth rate to be strictly negative for all $k$.

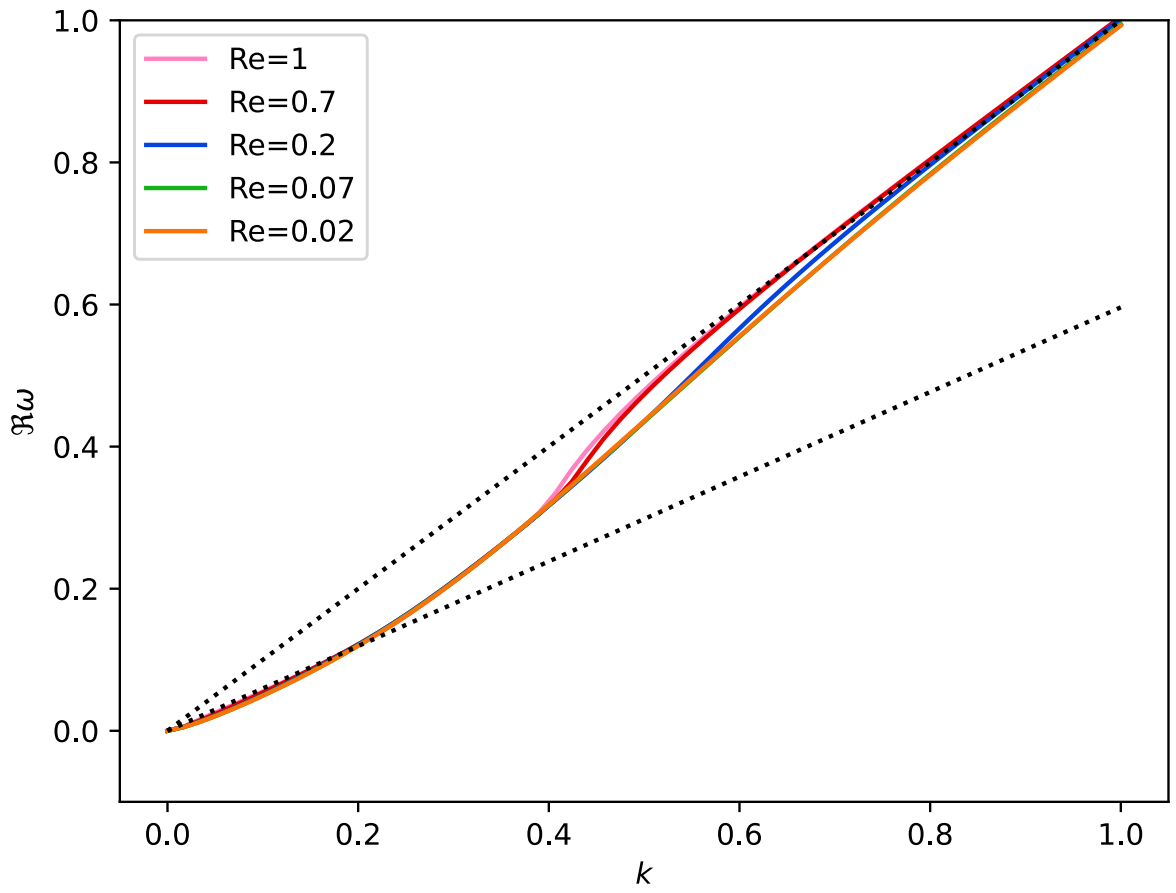

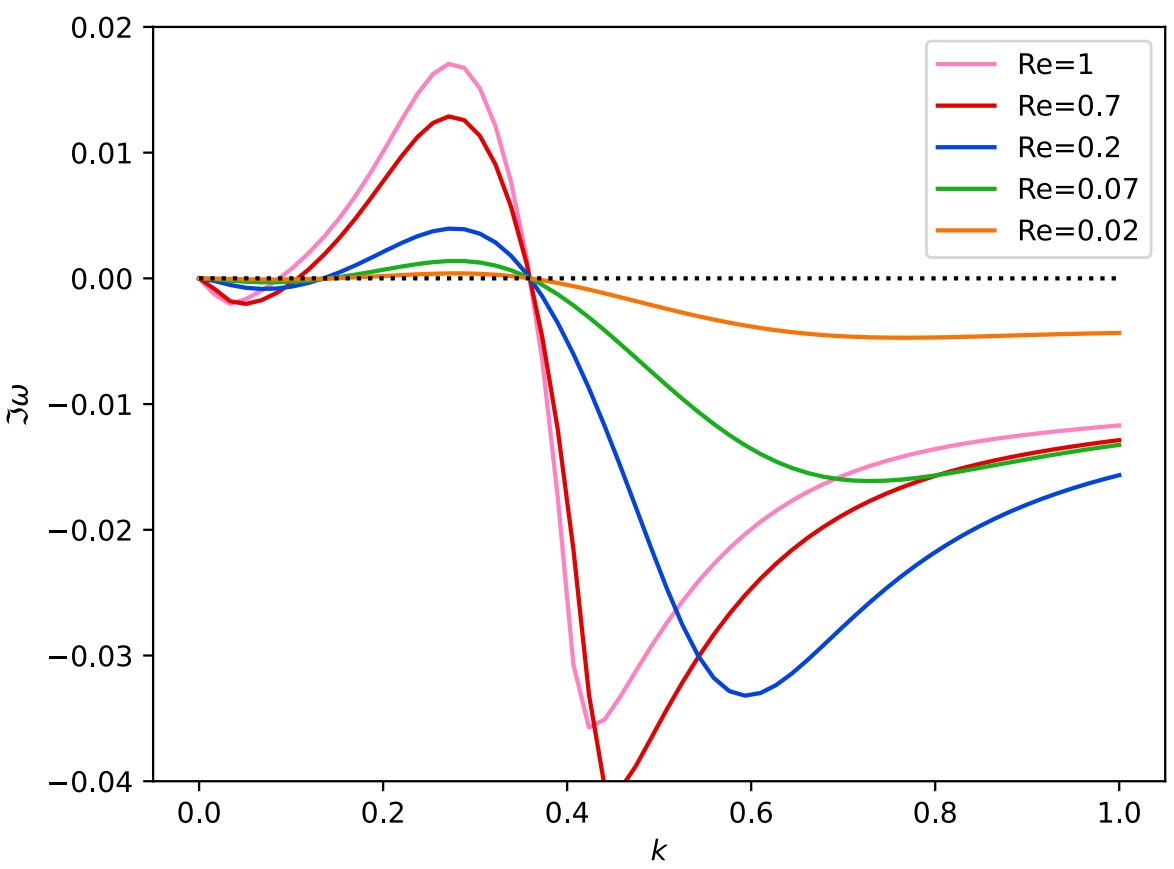

Figure 2: Dispersion relation (top), growth rate (bottom), at the layer-separating interface for gravity driven film flow without surface tension.

# Waves and transport in a circular water wave resonator: experiments and simulations

Uwe Harlander[1,*], Franz-Theo Schön[1], Ion D. Borcia[2], Rodica Borcia[3], Seymen Ilke Kaykanat[3], Michael Bestehorn[4]

[1]Dept. of Aerodynamics and Fluid Mechanics, BTU Cottbus-Senftenberg, Cottbus, Germany
[2]MINT Institute for Physics, BTU Cottbus-Senftenberg, Cottbus, Germany
[3]Dept. of Statistical Physics and Nonlinear Dynamics, BTU Cottbus-Senftenberg, Cottbus, Germany
*Email: uwe.harlander@b-tu.de

**Abstract**

We investigate surface waves in an oscillating *circular* channel that, due to a built-in barrier, is either of finite length or periodic with a submerged local topography mounted at the bottom. Resonance spectra and wave forms depending on the three input parameters, fluid depth, and the tank's oscillation amplitude and frequency, are studied. In contrast to the finite-length resonator, the topographic resonator allows wave scattering, which significantly alters the spectra. Breaking the symmetry of the oscillation or the topography can excite mean fluid transport in the topographic resonator. This transport is tightly connected to the resonant frequency bands observed. A particle image velocimetry system measures the velocity field, and ultrasonic sensors measure the surface displacement. Several numerical models are used to supplement the experimental data.

## 1 Introduction

The experimental study of resonant oscillations in closed containers has a long history, and such oscillations occur in a variety of media. Shocks in closed gas tube resonators have been studied since the 1930s, and nonlinear surface waves in water resonators were investigated in the 1960s. In the present study, surface waves in a circular water wave resonator are studied. The waves are either reflected back and forth by a barrier in the water channel, or they are partly reflected and partly transmitted due to a submerged local topography. Waves are generated by the oscillation of the cylindrical tank (see Fig. 1) with the water channel located at the outer edge.

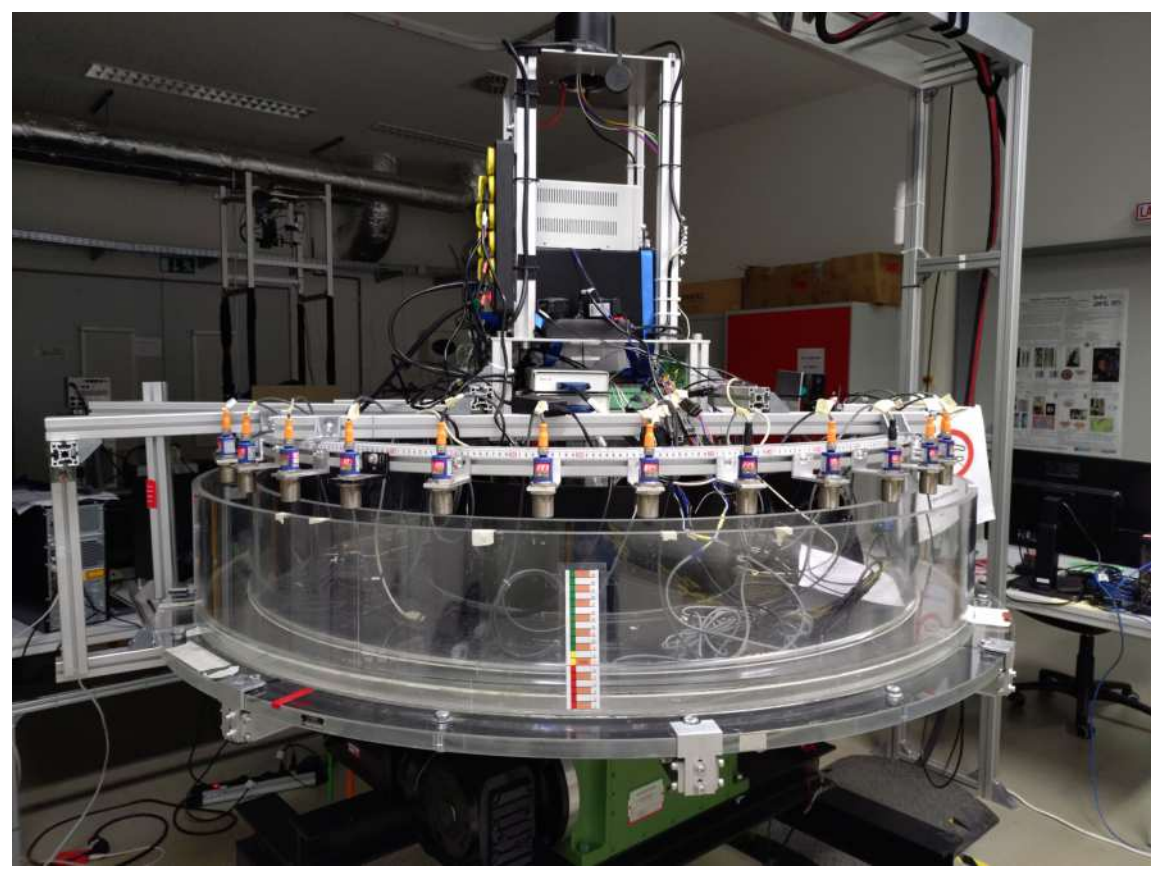

Figure 1: Picture of the experiment. The setup includes two co-rotating computers, 18 ultrasonic sensors, a DC power supply, a green continuous Laser sheet, and a 100 FPS global shutter camera.

## 2 Experiment

The annular channel has a length of 4.76 m and a width of 0.08 m. Water level is flexible and around 0.05 m in most experiments. We measure the height of the water surface by 17 ultrasonic probes covering one half of the channel. For synchronization purposes, an eighteenth sensor is located in the other half. In a small 0.06 m long segment of the channel, we measure the velocity to determine the mean fluid transport.

## 3 Models

We use various numerical models to supplement the experimental data. For studying bores in the closed channel, we solved the Euler equation for an irrotational 2D flow [2]. Figure 2 shows a solution obtained with this model for a two-bore collision. The top figure shows the experimental results, the bottom figure the numerical results. We further use a full Navier-Stokes model with a vertical coordinate transformation to include the rotational part of the flow [3]. To study solitons in different channel geometries, we applied the forced KdV model by [4]. For efficient computations of wave propagation and resonance, we used a long-wave model [5].

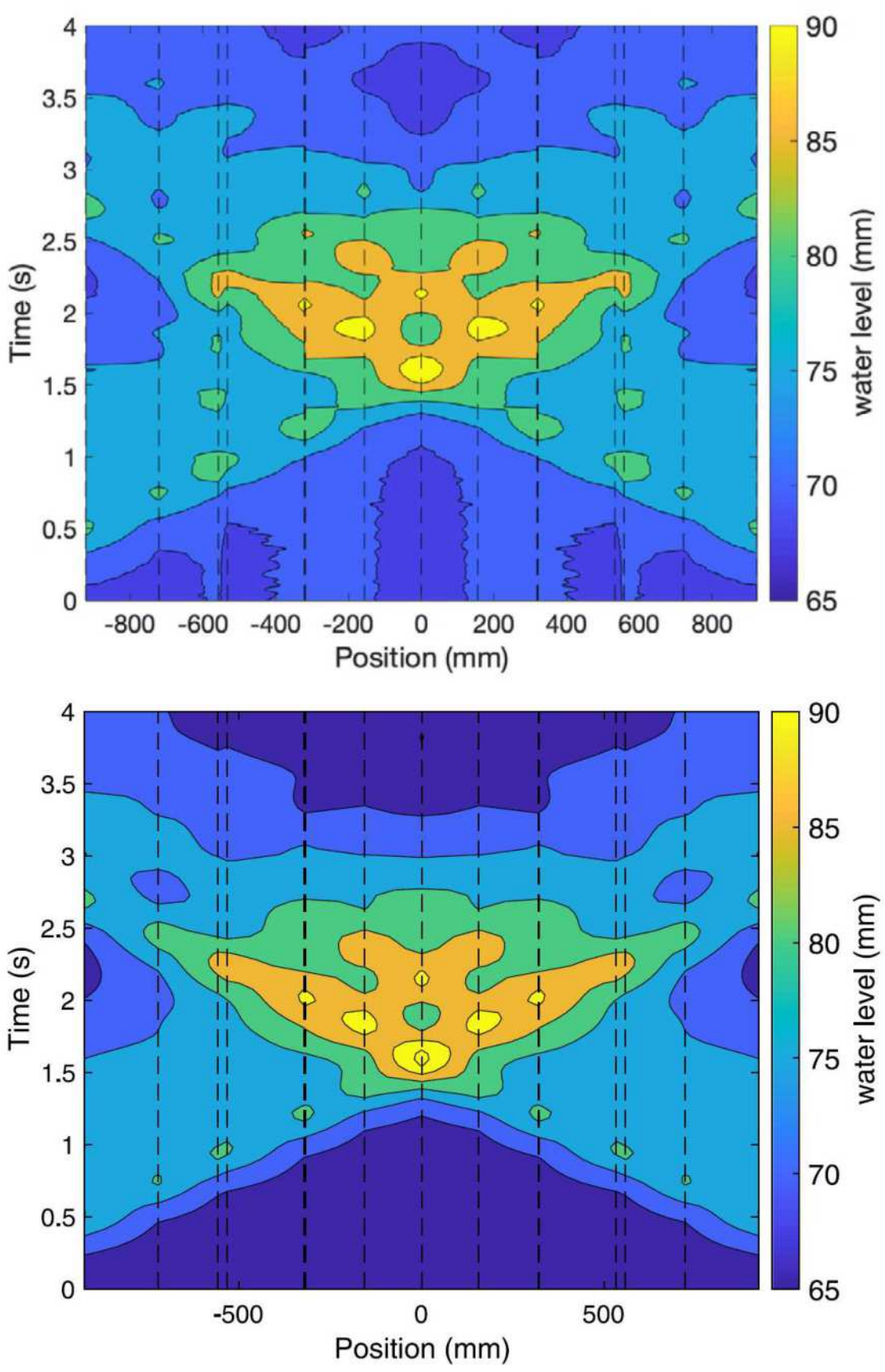

Figure 2: Comparison between experimental data (top) and computer simulations (bottom) for two colliding bores. The position of the sensors is marked with dashed lines. Figure adapted from [1].

## 4 Methods and Results

In contrast to a simple string model, in our finite-length water channel at resonance, we do not find standing waves but strongly localized propagating waves in the form of solitons and bores. Using characteristics of the Poincaré wave equation, we find wave focusing for excitation frequencies close to the linear eigenfrequencies. Normalizing the amplitude resonance peaks with the first linear eigenfrequency $f_0$, we find peaks only near odd integer frequencies. This is due to the anti-symmetric way the barrier excites the waves. In a thought experiment with two barriers moving in opposite directions, we would expect to see the peaks just near the even integer frequencies. Interestingly, replacing the barrier in our periodic channel with a submerged bottom topography, we find peaks at all integers due to wave scattering. The peak heights depend on the ratio $h_b/h_0$ (obstacle height/mean water depth).

We constructed a large number of space-time plots of the fluid surface using the experiment and the models mentioned above. Depending on the excitation frequency, we can observe a wide range of different wave types. In very simple terms, we can say that we usually observe bores, single or a series of solitons, at the resonance peaks, and cnoidal-type waves for the frequencies in between. The transition between the patterns can be quite abrupt.

Breaking the symmetry of the wave forcing, either by asymmetric tank oscillations or an asymmetric profile of the submerged hill, gives rise to mean transport in the periodic channel. We could show experimentally and numerically that the transport strongly depends not only on the forcing amplitude but also on the frequency. There is a large correlation between the resonance peaks and transport maxima. This result is of interest in tidal flow generation and in applications like coastal engineering.

# Truncated kernel windowed Fourier projection: a fast algorithm for the 3D free-space wave equation

**Nour G. Al Hassanieh**[1,*], **Alex H. Barnett**[1], **Leslie Greengard**[1,2]
[1]Center for Computational Mathematics, Flatiron Institute, Simons Foundation, New York, New York 10010, USA
[2]Courant Institute of Mathematical Sciences, New York University, New York, New York 10012, USA
*Email: nalhassanieh@flatironinstitute.org

**Abstract**

We present a spectrally accurate fast algorithm for evaluating the solution to the scalar wave equation in free space driven by a large collection of point sources in a bounded domain. With $M$ sources temporally discretized by $N_t$ time steps of size $\Delta t$, a naive potential evaluation at $M$ targets on the same time grid requires $\mathcal{O}(M^2 N_t)$ work. Our scheme requires $\mathcal{O}((M + N^3 \log N)N_t)$ work, where $N$ scales as $\mathcal{O}(1/\Delta t)$, i.e., the maximum signal frequency. This is achieved by using the recently-proposed windowed Fourier projection (WFP) method to split the potential into a local part, evaluated directly, plus a smooth history part approximated by an $N^3$-point equispaced discretization of the Fourier transform, where each Fourier coefficient obeys a simple recursion relation. The growing oscillations in the spectral representation (which would be present with a naive use of the Fourier transform) are controlled by spatially truncating the hyperbolic Green's function itself. Thus, the method avoids the need for absorbing boundary conditions. We demonstrate the performance of our algorithm with up to a million sources and targets at 6-digit accuracy. We believe it can serve as a key component in addressing time-domain wave equation scattering problems.

## 1 Introduction

We present a fast algorithm to compute free-space hyperbolic potentials in 3D to high accuracy without the need for artificial absorbing boundary conditions. Such potentials are space-time convolutions of the free-space retarded Green's function with a given source distribution. We assume that the source distribution has already been spatially discretized as a large collection of points.

There are many advantages to using integral equation methods. They do not require absorbing boundary conditions, and (in the homogeneous case), they discretize the boundary alone. The rapid evaluation of such potentials accelerates time-domain integral equation-based schemes for the solution of wave scattering problems (for example [1]). The latter are essential in many scientific and engineering applications such as acoustics and, generalizing to vectorial cases, electromagnetic/optical scattering and seismic propagation.

## 2 Problem Statement

We consider the solution of the free-space wave equation

$$\partial_t^2 u(\mathbf{x},t) - \Delta u(\mathbf{x},t) = f(\mathbf{x},t), \tag{1a}$$

$$u(\mathbf{x},0) = \partial_t u(\mathbf{x},0) = 0, \tag{1b}$$

for $\mathbf{x} \in \mathbb{R}^3$, and $t \in (0,T]$. The forcing function $f$ will be assumed to be a sum of $M$ point sources located at $\mathbf{y}_j$, $j = 1,\dots,M$, with smooth time signatures $\sigma_j(t)$, such that $\sigma_j(t) = 0$ for $t \le 0$:

$$f(\mathbf{x},t) = \sum_{j=1}^{M} \delta(\mathbf{x}-\mathbf{y}_j)\sigma_j(t), \tag{2}$$

where $\delta(\mathbf{x})$ denotes the three-dimensional Dirac distribution. We assume that the points $\mathbf{y}_j$ are contained in a finite computational domain $B \equiv [-1,1]^3$ and that the solution is sought within $B$ as well. In particular, we focus on the efficient evaluation of $u(\mathbf{x},t)$ at a large number of spatial targets $\mathbf{x} \in B$, for a set of discrete time slices $t \in (0,T]$.

The problem (1) can be solved exactly using the free-space Green's function for the wave equation (the *wave kernel*) given by

$$G(\mathbf{x},t) = \frac{\delta(t-|\mathbf{x}|)}{4\pi\,|\mathbf{x}|}, \tag{3}$$

where $\delta(t)$ denotes the one-dimensional Dirac distribution. The exact solution to (1) is the space-time convolution

$$u(\mathbf{x}, t) = \int_0^t \int_{\mathbb{R}^3} G(\mathbf{x} - \mathbf{y}, t - \tau) f(\mathbf{y}, \tau) d\mathbf{y} d\tau. \tag{4}$$

Direct evaluation of this retarded potential at $M$ targets at a single time $t$ has a cost of $\mathcal{O}(M^2)$. (In an integral equation setting, such an evaluation would be needed at each time step.) A purely Fourier representation of the solution, on the other hand, is of little numerical interest because (i) the solution is non-smooth spatially due to the non-smooth "switch on" of the response at $\tau = t$, (ii) the solution is dependent on the entire space-time history of $f$, and (iii) the spectral wave kernel becomes more oscillatory as time advances.

## 3 The WFP Method

The first two above concerns are addressed by the windowed Fourier projection (WFP) method, introduced by the authors in [2], which splits the solution (4) into a non-smooth local part, $u_l$, plus a smooth history part, $u_h$, using a temporal blending function, $\phi$, such that $\phi(t) = 0$ for $t \leq 0$, and $\phi(t) = 1$ for $t \geq \delta$. The blending width $\delta > 0$ is chosen to be of the order of the shortest timescale present in the signatures $\sigma_j$. The local part $u_l$ involves interactions only up to distance $\delta$, so that $u_l$ can be computed directly—typically in linear time—using suitable interpolation and quadrature rules. The fast evaluation of $u_h$ relies on the non-uniform fast Fourier transform (NUFFT) and overcomes the history-dependence by exploiting the semigroup property of the governing wave equation. Specifically, a discretized set of Fourier coefficients of $u_h$ will be updated by time-marching to arbitrarily high order accuracy, with a uniform time-step $\Delta t$ that need only be slightly smaller than the time-step needed to resolve the signatures $\sigma_j$ themselves.

## 4 Truncating the Kernel

To eliminate the oscillatory behavior of the spectral wave kernel (concern (iii)), we apply a smooth radial truncation of the kernel $G(\mathbf{x}, t)$ in (3) such that the solution representation remains unchanged in the box $B$. The finite spatial extent of the truncated kernel controls the spectral wave kernel oscillation rate with respect to the wavevector, so that a Nyquist-spaced trapezoid quadrature is exact for a wavenumber grid spacing of only $\mathcal{O}(1)$, no matter how large $t$ is.

## 5 Results

The Truncated Kernel WFP (TK-WFP) algorithm requires $\mathcal{O}\left(M + N^3 \log N\right)$ operations per time-step, where $N = \mathcal{O}(1/\Delta t)$ is the number of quadrature nodes per dimension needed to compute the history part in the Fourier transform domain. For a problem with $10^6$ random sources, whose frequencies vary between 0 and $30\pi$, the TK-WFP algorithm is 75 times faster than a naive direct evaluation. Figure 1 shows the computed solution at time $t = 6$. We aim to present applications of the TK-WFP method in accelerating 2D/3D wave equation scattering schemes [3].

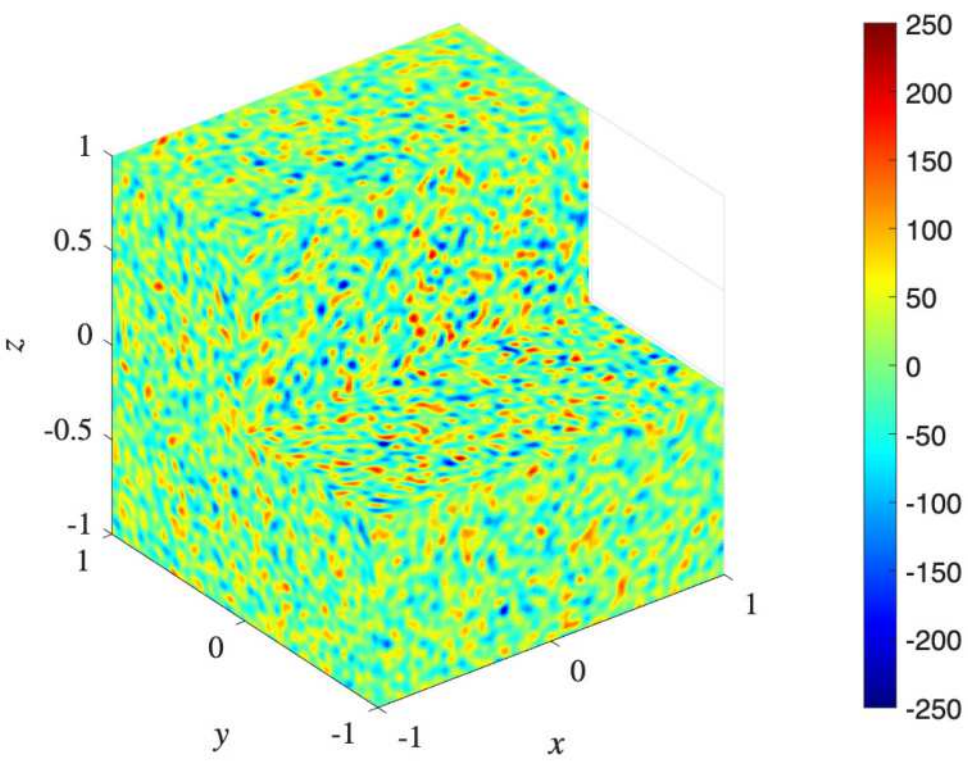

Figure 1: Computed solution at $t = 6$ (right), for a problem with $10^6$ random sources with random frequencies. The solution is shown on a target grid of size $100^3$. The cube has 30 wavelengths per side length at the maximum center frequency of the sources. The relative maximum error at $t = 6$ is $\tilde{\mathcal{E}}_{\Delta t} = 1.4 \times 10^{-6}$.

# Geometry-conforming Galerkin methods for the Lippmann-Schwinger equation on fractal domains

**David P Hewett**[1,*], **Joshua Bannister**[1], **Andrew Gibbs**[1]
[1]Department of Mathematics, University College London
*Email: d.hewett@.ucl.ac.uk

**Abstract**

We propose and analyse a numerical method for acoustic scattering by a class of inhomogeneities (penetrable scatterers) with fractal boundary. Our method is based on a Galerkin discretisation of the Lippmann-Schwinger volume integral equation, using a discontinuous piecewise-polynomial approximation space on a geometry-conforming mesh comprising elements which themselves have fractal boundary.

## 1 Introduction

In applications one often encounters scatterers whose boundaries are complicated and highly non-smooth, with fractal-like multi-scale geometrical features. Examples include trees and other vegetation, building facades, the interfaces between rock strata, and complex ice crystal aggregates in clouds.

In recent years there have been exciting developments in the application of integral equation methods to such problems [1]. But the focus to date has been on impenetrable scatterers. In this work we consider scattering of acoustic waves by a penetrable scatterer (an inhomogeneity) with fractal boundary, as in Figure 1.

## 2 Model

We consider outgoing solutions in $\mathbb{R}^n$, $n = 2, 3$, of the inhomogeneous Helmholtz equation

$$\Delta u^s + k^2(1 + \mathfrak{m})u^s = -f, \tag{1}$$

where $k > 0$ is the wavenumber, $\mathfrak{m}$ is the refractive index perturbation, which we assume to be bounded and supported in a compact set

$$K := \operatorname{supp} \mathfrak{m}, \tag{2}$$

that we refer to as the *inhomogeneity*, and $f$ is a source term (associated with the incident wave), also supported in $K$.

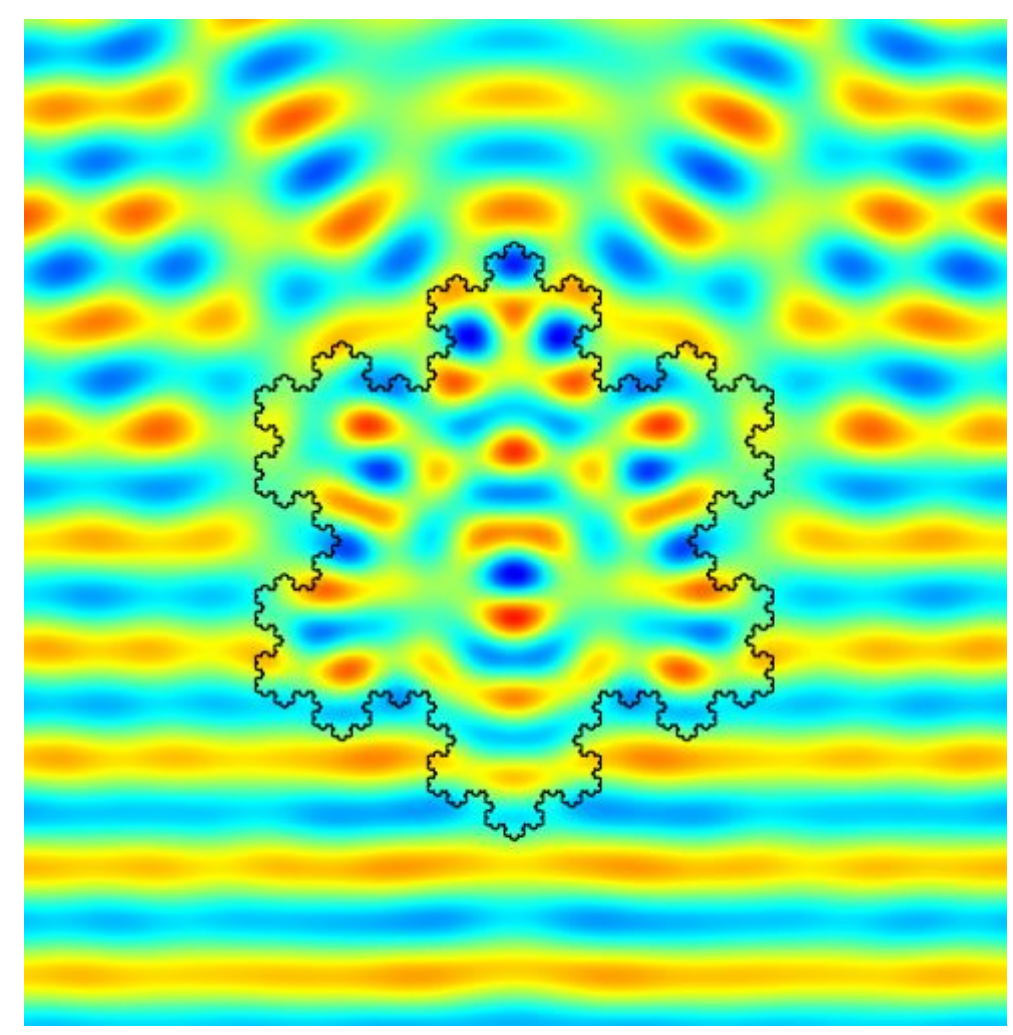

Figure 1: Scattering of a plane wave by the Koch Snowflake. Here $\mathfrak{m} = \chi_K$ so the refractive index is 2 inside $K$ and 1 outside $K$.

The problem can be equivalently formulated as a Lippmann-Schwinger equation

$$u^{\mathrm{s}}(x) = \int_{\mathbb{R}^n} \Big( f(y) + k^2 \mathfrak{m}(y) u^{\mathrm{s}}(y) \Big)\, \Phi(x, y)\, \mathrm{d}y, \tag{3}$$

where $\Phi(x, y)$ is the usual fundamental solution for $\Delta u + k^2 u = 0$.

## 3 Method and Results

To solve the integral equation (3) numerically we restrict it to the compact set $K$, then apply a Galerkin discretization based on a discontinuous piecewise-polynomial approximation on a geometry-conforming mesh. When $\partial K$ is fractal this means the mesh elements necessarily have fractal boundary, as in Figure 2.

We have performed a semi-discrete analysis of the method for completely arbitrary $K$, for both the $h$- and $p$-versions of the method, including proving superconvergence of linear functionals such as scattered field and far-field pattern evaluations, and elucidating how the regularity of $\partial K$ and $\mathfrak{m}$ affects the predicted convergence rates. These results make use of recent developments in the theory of polynomial approximation on non-Lipschitz domains [2].

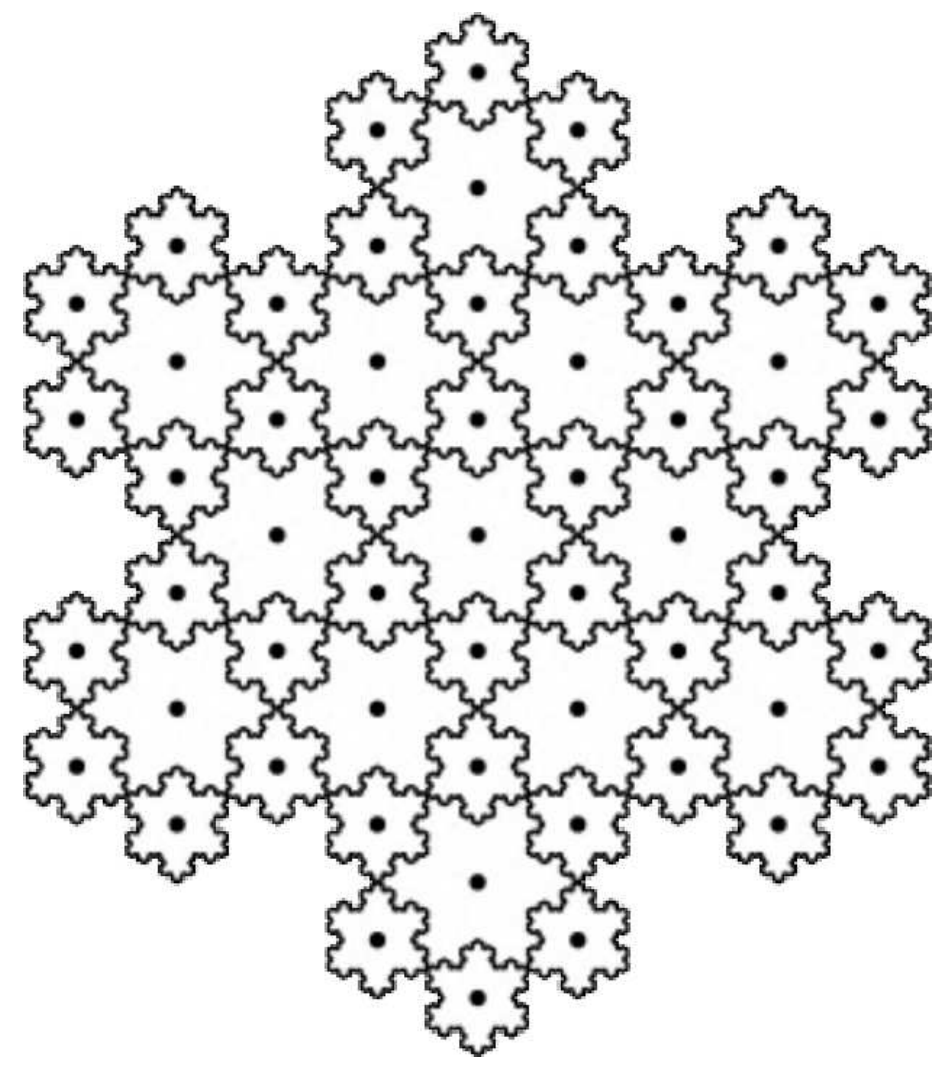

Figure 2: A geometry-conforming mesh.

In the case where $K$ is the attractor of an iterated function system of contracting similarities (such as the Koch Snowflake in Figures 1 and 2) we show how the self-similarity of $K$ can be used to generate geometry-conforming meshes that are closely related to fractal tilings. This self-similarity can also be exploited, using results from [3, 4], to derive accurate singular quadrature rules, permitting practical implementation of the method. We present such rules for the case where $\mathfrak{m}$ is constant on $K$, supported by a fully discrete error analysis.

We present numerical results for piecewise-constant approximation for a number of examples, which validate our theoretical results, and show that our geometry-conforming approach is significantly more accurate than a comparable approach based on classical triangulations of polygonal prefractal approximations to $K$ - see Figures 3 and 4.

Further details can be found in [5].

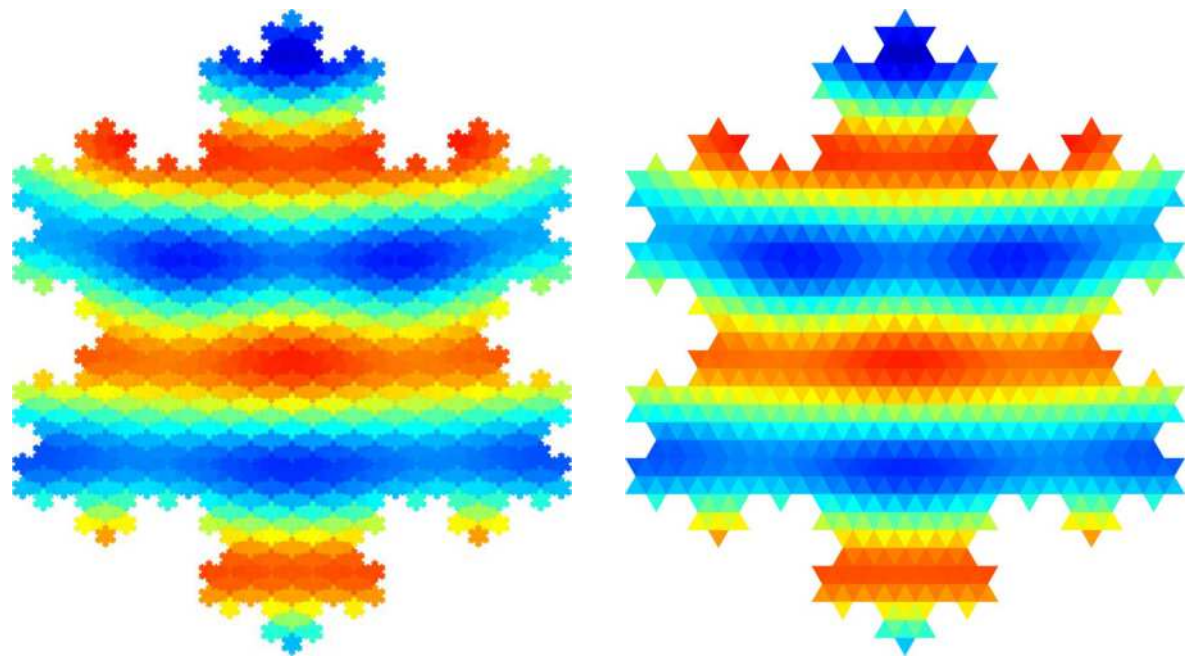

Figure 3: Piecewise-constant approximations on geometry-conforming (left) and prefractal (right) meshes.

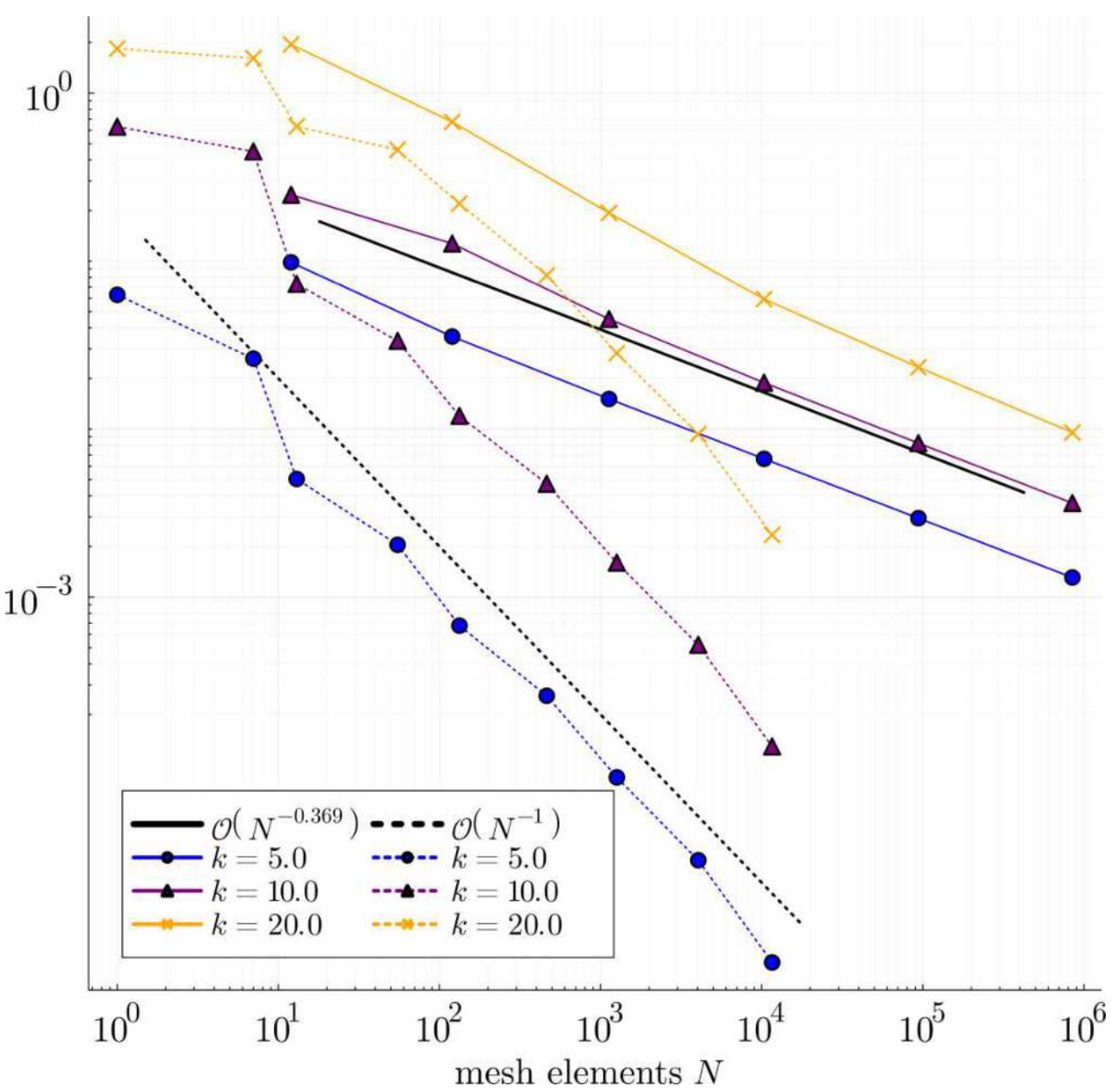

Figure 4: Error in the scattered field for piecewise-constant approximations on the Koch snowflake using our geometry-conforming approach (dashed lines) and a comparable prefractal approach (solid lines), for $\mathfrak{m} = m\chi_K$, with $m = 1.311 + 2.289 \times 10^{-9}$i.

# Uniqueness results in time-harmonic passive imaging with applications to helioseismology

**Thorsten Hohage**[1,*], **Björn Müller**[2]
[1]Institute of Numerical and Applied Mathematics, University of Göttingen, Germany
[2]Max-Planck Institute for Solar System Research, Göttignen, Germany
[*]Email: hohage@math.uni-goettingen.de

**Abstract**

Passive imaging problems involve identifying interior parameters in a wave equation from measurements of the cross-correlations of randomly excited solutions. In this paper, we demonstrate that measurements taken at two distinct hypersurfaces surrounding the interior uniquely determine both the scalar and the vectorial potential of a Schrödinger operator up to a gauge transformation. The gauge transformation can be resolved in helioseismology when measurements are taken at two different frequencies. Our uniqueness result also applies to other passive imaging problems for scalar second-order wave equations.

## 1 Introduction

We consider the inverse problem to identify coefficients $\mathbf{A}, q$ in a linear time-harmonic wave equation

$$\mathcal{L}_{\mathbf{A},q}\psi = s \qquad \text{in } \Omega \tag{1}$$

in a bounded domain $\Omega \subset \mathbb{R}^3$ with a random, unobserved source term $s$ with expectation $\mathbb{E}[s] = 0$ and a differential operator $\mathcal{L}_{\mathbf{A},q} = \mathcal{L}_{\mathbf{A},q,\omega}$ depending on the (known) frequency $\omega > 0$.

The data of the inverse problem are given by the cross-covariance function on a hypersurface $\Gamma$, which is given by

$$C(\mathbf{r}_1, \mathbf{r}_2, \omega) = \mathbb{E}\psi(\mathbf{r}_1, \omega)\overline{\psi(\mathbf{r}_2, \omega)}, \qquad \mathbf{r}_1, \mathbf{r}_2 \in \Gamma$$

since the PDE is linear and $\mathbb{E}[s] = 0$. This cross-covariance function can be estimated by measuring many waves $\psi^{(n)}$, $n = 1, \dots, N$ corresponding to independent realizations $s^{(n)}$ of the source process.

## 2 Model

We consider a differential operator of the form

$$\mathcal{L}_{\mathbf{A},q,\omega} = -\Delta - 2i\mathbf{A}_\omega \cdot \nabla + q_\omega - i(\nabla \cdot \mathbf{A}_\omega) - k_\omega^2 \tag{2}$$

with a real vector potential $\mathbf{A}_\omega$, a complex scalar potential $q$, and a wave number $k_\omega$. On the boundary we impose a general condition

$$\partial_{\mathbf{n}}\psi = B_\omega \psi \quad \text{on } \partial\Omega$$

with an operator $B_\omega$ that may be an exterior Dirichlet-to-Neumann map of an exterior domain or a local operator, e.g., from a transparent boundary condition.

Further standard assumptions are imposed which ensure well-posedness of this forward problem in the sense that there exists a bounded solution operator

$$\mathcal{G}_{\mathbf{A},q} : H_0^{-1}(\Omega) \to H^1(\Omega), \quad \mathcal{G}_{\mathbf{A},q} s := \psi\,.$$

We assume that the solution is measured on two non-intersecting smooth surfaces $\Gamma_a$ and $\Gamma_b$ in $\Omega$ which surround the support of the unknown functions $\mathbf{A}$ and $q$ (see Fig. 1). In helioseismology measurements at different heights are feasible since data for different atomic spectral lines are associated to different heights.

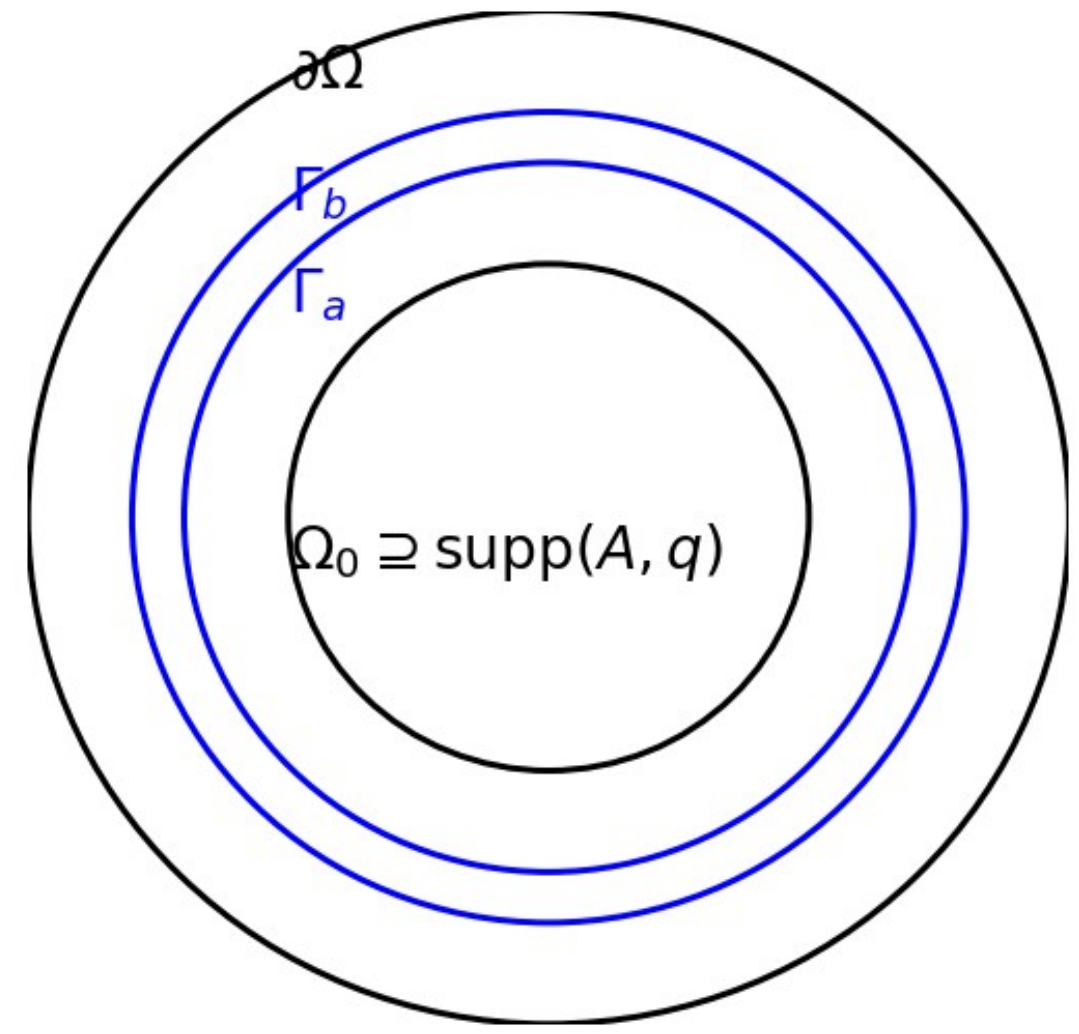

Figure 1: Geometry. $\Gamma_a$ and $\Gamma_b$ are measurement surfaces.

## 3 General uniqueness result

We only state a special case of our main uniqueness theorem here. It is often assumed in helioseismology and passive imaging is that the surface covariance is proportional or equal to the imaginary part of the Green's function:

$$\mathcal{C}^{\Gamma_i\Gamma_j}_{\mathbf{A},q} = \mathrm{Tr}_{\Gamma_i} \Im\left(\mathcal{G}_{\mathbf{A},q}\right) \mathrm{Tr}^*_{\Gamma_j} \tag{3}$$

where $\mathcal{C}^{\Gamma_i\Gamma_j}_{\mathbf{A},q}$ is the integral operator with the cross-covariance function $C(\cdot,\cdot,\omega)$ as kernel, and $\Im\mathcal{G}_{\mathbf{A},q} = \frac{1}{2i}\left[\mathcal{G}_{\mathbf{A},q} - \mathcal{G}^*_{\mathbf{A},q}\right]$. This assumption essentially comes from a particular assumption on the source covariance.

**Theorem 1** *Suppose that* $(\mathbf{A}_1, q_1)$ *and* $(\mathbf{A}_2, q_2)$ *satisfy the assumptions imposed for well-posedness of the forward problem and* $k$ *is not in a some discrete set of exceptional frequencies. Then*

$$\mathrm{Tr}_{\Gamma_i} \Im\mathcal{G}_{\mathbf{A}_1,q_1} \mathrm{Tr}^*_{\Gamma_j} = \mathrm{Tr}_{\Gamma_i} \Im\mathcal{G}_{\mathbf{A}_2,q_2} \mathrm{Tr}^*_{\Gamma_j},$$

*for* $i,j \in \{a,b\}$ *implies that there exists* $\phi \in W^{1,\infty}(\Omega,\mathbb{R})$, $\mathrm{supp}(\phi) \subseteq \Omega_0$ *with*

$$\mathbf{A}_1 = \mathbf{A}_2 + \nabla\phi, \quad q_1 = q_2 + 2\mathbf{A}_2\cdot\nabla\phi + (\nabla\phi)\cdot(\nabla\phi)\,.$$

The gauge transformation in this theorem is an unavoidable source of nonuniqueness.

The core idea of the proof (see [4] for a preliminary version) is to select test functions for the operator-valued data that cause the solutions to vanish in the interior, thereby showing that cross-correlation uniquely determines the Dirichlet-to-Neumann map. This connection enables us to invoke uniqueness theorems for classical inverse problems for the differential equation (1) from [2].

## 4 Application to helioseismology

Line-of-sight velocities of the solar surface caused by acoustic waves are meausred by satellite and ground based instruments using the Doppler shift of certain atomic spectral lines. The wavefields can be approximately described by a scalar wave equation (see [3]), which – after a transformation of variables – takes the form (2) with

$$\mathbf{A}_\omega = \omega\frac{1}{c^2}\mathbf{u}, \qquad q_\omega = k_\omega^2 - \frac{\omega^2 + 2i\gamma(\omega)\omega - \omega_c^2}{c^2},$$
$$k_\omega^2 = \frac{\omega^2 - \omega_c^2}{c^2}, \qquad \omega_c^2 = c^2\rho^{1/2}\Delta(\rho^{-1/2}).$$

Here $c$ is the sound speed, $\rho$ is the density, and $\mathbf{u}$ is the flow field, and $\omega_c$ can be interpreted as acoustic cutoff frequency. Furthermore, $\gamma$ denotes a damping factor, which we assume to be of the form

$$\gamma(\omega,\mathbf{r}) = \gamma_0\left(\frac{\omega}{\omega_0}\right)^{\zeta(\mathbf{r})},$$

where $\gamma_0 \neq 0$ and $\omega_0 \neq 0$ are reference values, and $\zeta$ is real-valued.

**Corollary 2** *Assume that* (3) *holds true and* $\mathrm{Tr}_{\Gamma_i} \Im\mathcal{G}_{\mathbf{A}_1,q_1} \mathrm{Tr}^*_{\Gamma_j} = \mathrm{Tr}_{\Gamma_i} \Im\mathcal{G}_{\mathbf{A}_2,q_2} \mathrm{Tr}^*_{\Gamma_j}$ *for* $i,j \in \{a,b\}$ *for two distinct wavenumbers and* $i,j \in \{a,b\}$. *Then* $(\rho_1, c_1, \gamma_1, \mathbf{u}_1) = (\rho_2, c_2, \gamma_2, \mathbf{u}_2)$.

The result is shown by solving the gauge transformations for the unknowns, see [1] for a similar proof.

We also discuss the following generalizations of this result:

- Measurements only on part of the surfaces $\Gamma_a, \Gamma_b$ if these surfaces are analytic.
- Uniqueness results not assuming (3), but treating the pointwise variance of uncorrelated sources $s$ as additional unknown that can also be identified.

# Domain Decomposition for Elastic Waves in Spatial Networks

Joseph Holten[1,*], Moritz Hauck[1]
[1]Institute for Applied and Numerical Mathematics, Karlsruhe Institute of Technology, DE
*Email: joseph.holten@kit.edu

**Abstract**

We present an algebraic two-level domain decomposition preconditioner for elastic wave propagation in heterogeneous spatial network models discretized by hybrid discontinuous Galerkin (HDG) methods in space and implicit time integration. The spatial network is partitioned algebraically and a coarse space is constructed from the corresponding partition of unity with low-order polynomial enrichment. We prove that, under suitable assumptions, the condition number of the preconditioned system scales as $H/\delta$, independent of the edge lengths and the number of subdomains, where $H$ is the subdomain diameter and $\delta$ the overlap width. Numerical experiments on fiber network models confirm the theoretical predictions.

## 1 Introduction

Many problems in science and engineering modeled by partial differential equations are posed in complex domains, which typically lead to significant computational challenges. When the domain is locally slender, classical simulation techniques, such as the finite element method, require excessive spatial resolution and can result in nearly degenerate elements. As a remedy, one may employ spatial network models, which encode the locally one-dimensional geometry of the underlying domain as a graph where the edges correspond to slender segments of the domain and the nodes represent their connections. The local behavior of the problem is then described by one-dimensional differential equations posed on the edges which are suitably coupled at the nodes. The prominent application of spatial network models considered in this talk are fiber-based materials, which appear in a wide range of engineering applications, including paper, nonwoven fabrics, and sintered metallic fiber networks. Here, the deformation of fiber segments is modeled using one-dimensional beam models that are suitably coupled by rigid joint conditions at the nodes representing the contact regions of fibers [1–3]. Such spatial networks typically exhibit strong heterogeneities in geometry, topology and material parameters and can thus be classified as multiscale problems. Therefore, despite the reduction in dimensionality, they might remain intractable using standard discretization methods that resolve all scales globally.

For elastic wave propagation in these networks, explicit time stepping leads to severely restrictive CFL conditions caused by the presence of very short fibers. Implicit time integration avoids this but requires solving large, ill-conditioned linear systems at each time step. We address this challenge by designing a two-level subspace correction preconditioner whose convergence is robust with respect to the problem size.

Previous work in [4, 5] introduced a domain decomposition method for stationary problems on spatial networks using an artificial Cartesian coarse mesh. We extend this in two directions: (i) we replace the geometric coarse space by an algebraic construction based on graph partitioning, applicable to complex macroscale geometries; (ii) we apply the preconditioner to the wave equation with implicit time integration.

## 2 Model and Discretization

We consider a spatial network $\mathcal{G} = (\mathcal{N}, \mathcal{E})$ with nodes $\mathcal{N} \subset \mathbb{R}^3$ and edges $\mathcal{E}$ identified with line segments. Let $\Sigma := \bigcup \mathcal{E}$ be the union of the line segments. The relevant quantities are displacement $u$, rotation $r$, stress $n$, and moment $m$, considered as functions $\Sigma \to \mathbb{R}^3$. The governing equations in the Timoshenko beam model are posed for an edge $\mathfrak{e} \in \mathcal{E}$ and final time $T > 0$ on $\mathfrak{e} \times (0, T)$ and read as:

$$C_{\mathrm{u}}\ddot{u} + \partial_x n = f, \qquad C_{\mathrm{r}}\ddot{r} + \partial_x m + i \times n = g,$$
$$-C_{\mathrm{n}}(\partial_x u + i \times r) = n, \qquad -C_{\mathrm{m}}\partial_x r = m,$$

where $i$ is the unit basis vector pointing in the direction of the edge $\mathfrak{e}$, and $C_\bullet \in L^\infty(\Sigma; \mathbb{R}^{3\times 3})$ are material coefficient matrices satisfying uniform spectral ellipticity assumptions, and $f, g \in L^2(\Sigma)$ are the distributed forces and moments. The beams are coupled at the nodes by continuity of displacements and rotations and balance of forces and moments: there shall exist $\lambda, \phi : (0,T) \times \mathcal{N} \to \mathbb{R}^3$ satisfying for all $\mathfrak{n} \in \mathcal{N}$ and for all adjacent $\mathfrak{e} \in \mathcal{E}$ the equations

$$u|_{\mathfrak{e}}(t,\mathfrak{n}) = \lambda(t,\mathfrak{n}), \quad r|_{\mathfrak{e}}(t,\mathfrak{n}) = \phi(t,\mathfrak{n}),$$
$$[\![n\nu]\!](t,\mathfrak{n}) = 0, \quad [\![m\nu]\!](t,\mathfrak{n}) = 0,$$

where $\nu = \pm 1$, and $[\![\cdot]\!]$ is the natural jump operator and we assume for simplicity zero concentrated forces and moments.

The system of equations posed on the edges are discretized in space using a hybridizable discontinuous Galerkin (HDG) method [5], yielding local solvers mapping nodal data and distributed forces and moments $(\lambda, \phi, f, g)$ to the interior solution on each edge. Static condensation of the balance conditions produces for all times $t \in (0,T)$ a global system posed on the nodal unknowns: set $V := \{v : \mathcal{N} \to \mathbb{R}^3\}$ and seek $v = (\lambda(t,\cdot), \phi(t,\cdot)) \in V \times V =: V^2$ satisfying

$$\bar{\mathcal{A}}(v,w) = \mathcal{F}(w) \quad \text{for all } w \in V^2. \tag{1}$$

Time integration is performed implicitly (e.g., Crank–Nicolson) so that the large, sparse, ill-conditioned system (1) must be solved at each time step. It can be shown that it is in fact symmetric and positive definite, cf. [5].

## 3 Domain Decomposition Preconditioner

We use the preconditioned conjugate gradient method to efficiently solve (1). As preconditioner, we use a two-level additive Schwarz method based on the following subspace decomposition

$$V^2 = V_0 + V_1 + \ldots + V_p.$$

We algebraically partition the nodes

$$\mathcal{N} = \mathcal{N}_1 \cup \cdots \cup \mathcal{N}_p$$

minimizing cut edges in $\mathcal{E}$, and for $1 \le j \le p$ we define the local spaces as

$$V_j := \{v \in V^2 : d_{\mathcal{G}}(\operatorname{supp} v, \mathcal{N}_j) \le \delta\}$$

where $d_{\mathcal{G}}$ is the edge-length weighted graph distance in $\mathcal{G}$. The coarse space is constructed via the partition of unity $\{\Lambda_j\}_j$ corresponding to the overlapping cover constructed from the partition and define

$$V_0 := (\operatorname{span}\{\Lambda_j\}_j \otimes \mathbb{P}_1(\mathcal{N}))^6.$$

Under natural assumptions on the cover and the underlying graph (bounded multiplicity, homogeneity, connectedness, and an isoperimetric inequality) stability of the decomposition is established using a quasi-interpolation operator into the coarse space $V_0$ and spectral equivalence to a sum of reciprocally edge-length weighted graph Laplacian and edge-length weighted graph mass operators. This yields spectral bounds on the preconditioned operator $\mathcal{P}$ of the form $C_1^{-1} \le \lambda(\mathcal{P}) \le C_2$ with $C_2 \sim 1$ and $C_1 \sim \beta\alpha^{-1}\sigma H\delta^{-1}$, where $H = \max_j \operatorname{diam}_{\mathcal{G}}(\mathcal{N}_j^\delta)$.

Numerical experiments on random fiber networks show that: (i) the number of CG iterations is bounded independently of the network size for fixed $H/\delta$; (ii) algebraic partitioning produces subdomains adapted to the network geometry, outperforming Cartesian decompositions; (iii) the preconditioner is effective for the wave equation with implicit time stepping.

## High-order contour-integral methods for strongly continuous semigroups

**Andrew Horning**[1,*], **Adam R. Gerlach**[2]
[1]Department of Mathematical Sciences, Rensselaer Polytechnic Institute, Troy, USA
[2]Measured Engineering LLC., New Bremen, USA
*Email: hornia3@rpi.edu

**Abstract**

Exponential integrators based on contour integral representations lead to powerful numerical solvers for a variety of ODEs, PDEs, and other time-evolution equations. They are embarrassingly parallelizable and lead to global-in-time approximations that can be efficiently evaluated anywhere within a finite time horizon. However, there are core theoretical challenges that restrict their use cases to analytic semigroups, such as diffusion-dominated processes. In this talk, we show how carefully regularized contour integral representations can be used to construct high-order quadrature schemes for the larger, less regular, class of strongly continuous semigroups, which includes a wide range of wave- and advection-dominated phenomena. Our algorithms are accompanied by explicit high-order error bounds and near-optimal parameter selection. We demonstrate key features of the schemes on wave- and advection-dominated flows associated with some singular PDEs.

### 1 Contour-integral methods

Contour integral methods approximate the action of the matrix exponential with a quadrature rule. Given $\gamma$, a Jordan curve winding once clockwise around the spectrum of a square matrix $A$, quadrature nodes $z_1, \ldots, z_N$, and weights $w_1, \ldots, w_N$, then

$$\begin{aligned} \exp(At)x &= \frac{1}{2\pi i}\int_\gamma e^{zt}(zI-A)^{-1}x\,dz \\ &\approx \frac{1}{2\pi i}\sum_{k=1}^{N} w_k e^{z_k t}(z_k I - A)^{-1}x. \end{aligned} \tag{1}$$

The main cost of the scheme is solving a shifted linear equation at each quadrature node. These can be done in parallel and may be accelerated with iterative methods. After solving these linear systems, the exponential can be evaluated at any time $t > 0$ with only vector operations.

When $A$ is obtained by discretizing a differential operator, its spectrum may fill a large region of the complex plane. As the discretization is refined, this region typically grows and the quadrature approximation in (1) may rapidly deteriorate along with the performance of iterative methods for the linear systems. Famously, this effect can be mitigated through *Talbot quadrature* if the spectrum of $A$ lies in a modest sector along the negative real axis [1]. These techniques have been used to develop highly parallelizable schemes for a broad class of parabolic equations [2]. Recently, Colbrook was awarded a second-place Leslie Fox Prize for extending these techniques to analytic semigroups, which are associated with smooth evolution problems like parabolic and damped wave equations [3].

### 2 Extension to $C^0$-semigroups

Fundamentally, Talbot contour integral schemes exploit *regularity in the operator exponential* for analytic semigroups. Unfortunately, this is not possible for less regular evolution equations associated with challenging causal, singular, or non-normal behavior. Instead, we show how to regularize contour integral schemes to exploit *regularity in the vector* $x$. This is common in time-stepping methods but global-in-time methods require new tools from semigroup theory.

To exploit regularity in $x$ in our scheme and analysis, we develop the algorithm at the infinite-dimensional level and then carefully assess the impact of discretization. Consider a linear operator $\mathcal{A} : D(\mathcal{A}) \to \mathcal{X}$ on a separable Banach space $\mathcal{X}$ and replace $x$ by $u \in \mathcal{X}$ and $A$ by $\mathcal{A}$. If $\mathcal{A}$ generates a strongly continuous semigroup on $\mathcal{X}$, then its spectrum is contained in a left half-plane $\mathrm{Im}(z) < \omega$ and its resolvent decays uniformly in the complimentary right-half plane. We use a regularized analogue of the contour integral in (1) that, for some $\delta > \omega$ and integer

$m \geq 2$ and any $u \in D(\mathcal{A}^2)$, has form

$$\exp(\mathcal{A}t)u = (2\delta\mathcal{I} - \mathcal{A})^m S^{(m)}(t)u, \quad \text{where}$$
$$S^{(m)}(t) = \frac{1}{2\pi i}\int_{\delta - i\infty}^{\delta + i\infty} \frac{e^{zt}}{(2\delta - z)^m}(z\mathcal{I} - \mathcal{A})^{-1}\, dz. \tag{2}$$

To approximate (2), we apply a truncated $N$-point trapezoidal rule with node spacing $h > 0$,

$$S_N^{(m)}(t)u = (2\delta\mathcal{I} - \mathcal{A})^m Q_N^{(m)}(t)u, \quad \text{where}$$
$$Q_N^{(m)}(t) = \frac{h}{2\pi}\sum_{k=-N}^{N} \frac{e^{(\delta + ihk)t}}{(\delta - ihk)^m}((\delta + ihk)\mathcal{I} - \mathcal{A})^{-1}. \tag{3}$$

For a practical computational scheme, the shifted linear equations can be solved with a suitable discretization of smooth functions $u \in D(\mathcal{A}^2) \subset \mathcal{X}$. The main requirement for convergence is that $Ax \to \mathcal{A}u$ and $(z_k I - A)^{-1}x \to (z_k\mathcal{I} - \mathcal{A})^{-1}u$, $k = -N, \ldots, N$, as the discretization is refined.

## 3 Explicit error bounds

At the operator level, the regularized quadrature approximation in (2)-(3) has several key advantages over (1). The integral converges absolutely and uniformly for $u \in D(\mathcal{A}^2)$, while the amplification power of $(2\delta\mathcal{I} - \mathcal{A})^m$ is controlled by the regularity of the vector $u$. Moreover, the integrand decays at a controlled rate along the contour. These features allow us to derive simple explicit error bounds for smooth vectors $u \in D(\mathcal{A}^m)$. In fact, we can derive closed-form expressions for quadrature parameters that, for $u \in D(\mathcal{A}^m)$, achieve [4]

$$\sup_{0 \leq t \leq T} \|\exp(\mathcal{A}t)u - S_N^{(m)}(t)u\| \leq C\left(\frac{\delta}{hN}\right)^{m-1} \|(2\delta\mathcal{I} - \mathcal{A})^m u\|. \tag{4}$$

Here, $C$ is an explicit constant depending on $A$, the contour location $\delta$, and the time horizon $T$.

In practice, one must discretize the infinite-dimensional objects for a numerical algorithm. We show how semigroup theory provides simple, computable bounds for the action of the discretized resolvent that hold for a very broad class of discretization techniques including Galerkin-based schemes and modern spectral methods.

## 4 Applications and outlook

High-order contour integral schemes for strongly continuous semigroups open the door to highly parallelizable methods for traditionally challenging PDEs and related time-evolution processes. They also provide a compact modal-like representation of semigroups that may be useful in operator learning, system identification, and other data-driven modeling techniques. We will illustrate our scheme's performance with several wave- and advection-dominated singular PDEs with variable media (see Figure 1 for a simple example). If time permits, we will also discuss work in progress on applications to Koopman semigroups in data-driven settings.

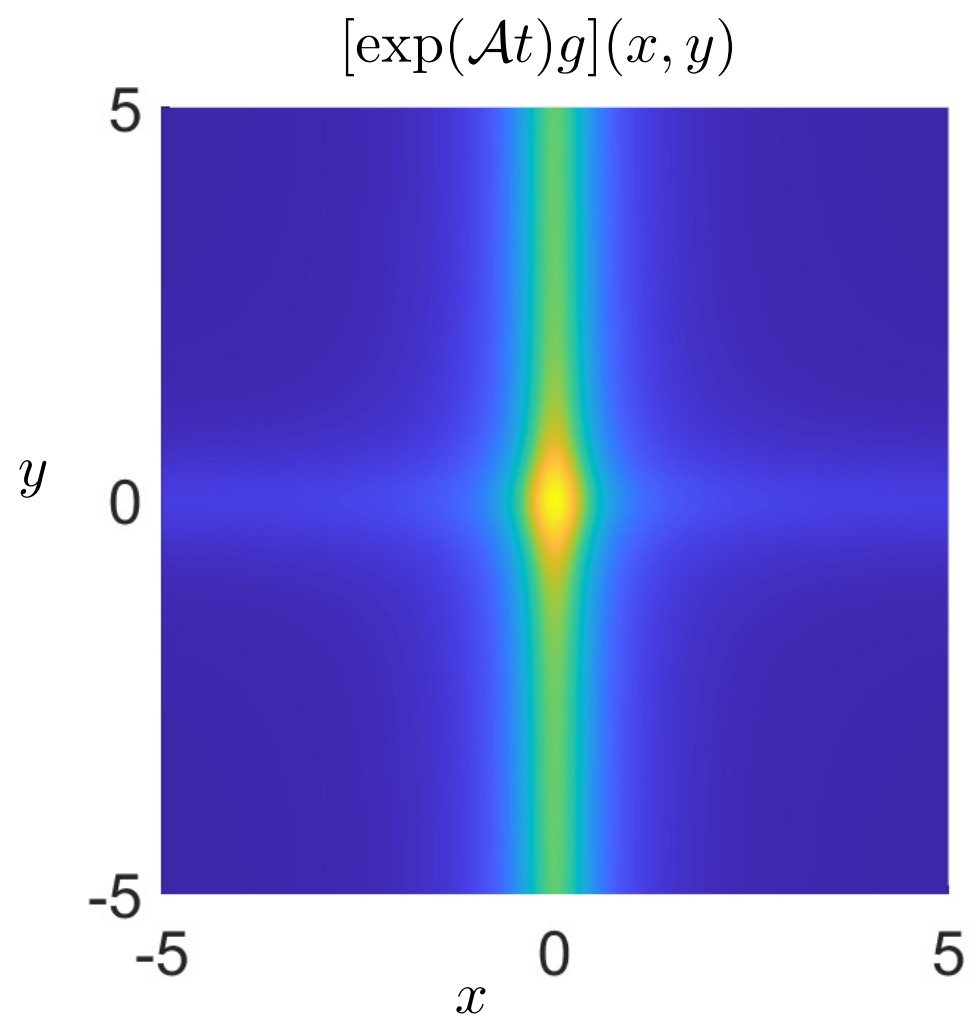

Figure 1: The convection of an anisotropic wave packet, $g(x, y) = \exp(-2x^2 - 0.5y^2)$, in a variable focusing media, governed by a singular differential operator $\mathcal{A}u = F \cdot \nabla u$ with $F(x, y) = (2x - 8x^3, 2y - 8y^3)$, is shown at $t = 1$.

# Bistability of travelling waves and mesa states in a mass-conserved reaction-diffusion system: From bifurcations to implications for actin waves

**Jack M. Hughes**[1,*], **Saar O. Modai**[2], **Leah Edelstein-Keshet**[3], **Arik Yochelis**[2,4]

[1]Department of Mathematics and Statistics, Concordia University, Montreal, Canada
[2]Department of Physics, Ben-Gurion University of the Negev, Be'er Sheva, Israel
[3]Department of Mathematics, University of British Columbia, Vancouver, Canada
[4]Swiss Institute for Dryland Environmental and Energy Research, Blaustein Institutes for Desert Research, Ben-Gurion University of the Negev, Sede Boqer Campus, Midreshet Ben-Gurion, Israel

*Email: jack.hughes@mail.concordia.ca

**Abstract**

The transition between random walk and directional motion is one of the intriguing phenomena observed in eukaryotic cells. Typical theoretical studies distinguish between the two phenomena and focus either on dissipative or gradient flow models, respectively. Using a 3 variable dissipative reaction-diffusion system with mass conservation, we show how pulses, waves, mesas, and fronts generically organize about high codimension bifurcations. Specifically, we demonstrate the novelty of mass conservation, which enters as a persistent large-scale mode possibly leading to a long wavelength bifurcation. Following the biological interest, we address the bistability of travelling wave and mesa solution branches that emerge from a novel codimension-2 long wavelength/ finite wavenumber Hopf bifurcation.

## 1 Introduction

Cell motility is necessary for bodily functions, including wound healing. Filamentous actin (F-actin) is a key driver of motility. Regulators of F-actin formation include GTPases that act as molecular switches, cycling between active and inactive forms. Recent experiments show that stationary patterns that induce directed cell motion and wave patterns that induce cell turning or ruffling can coexist in cells.

To unravel the coexistence of these modes, we focus on the dynamics of F-actin and the GTPases. Typical models of this type focus on either capturing polarization (directed motion) through models with a gradient structure or turning and ruffling with dissipative systems. Models that combine gradient and dissipative structures have been proposed; however, these studies only show patterns in distinct parameters and thus cannot capture coexistence. In this work, we show that solutions indicative of directed cell motion and cell ruffling coexist by analyzing a novel codimension-2 bifurcation (simultaneous long wavelength/ finite wavenumber Hopf). Note that both bifurcations have been studied in isolation.

## 2 Model of GTPase-actin dynamics

The model consists of a three-variable reaction-diffusion system for a GTPase in active ($u$) and inactive ($v$) forms, and F-actin ($F$):

$$\frac{\partial u}{\partial t} = (b+\gamma u^2)v - (1+sF+u^2)u + D_u \Delta u, \tag{1a}$$

$$\frac{\partial v}{\partial t} = -(b+\gamma u^2)v + (1+sF+u^2)u + D_v \Delta v, \tag{1b}$$

$$\frac{\partial F}{\partial t} = \eta(p_0 + p_1 u - F) + D_F \Delta F. \tag{1c}$$

With periodic or no flux boundary conditions, the total mass of GTPase is conserved, i.e.,

$$M := \frac{1}{|\Omega|}\int_\Omega [u(x,t)+v(x,t)]\mathrm{d}x = \text{const.} \tag{2}$$

## 3 Bistability of travelling waves and mesas

The homogeneous steady states (HSSs) $\mathbf{P}^* = (u^*, v^*, F^*)$ of (1), which result in up to three biologically relevant solutions, $\mathbf{P}^*_{0,1,2} > 0$, form an inverse "S", as shown in Fig. 1. Linear stability analysis in 1D leads to solutions $\mathbf{P}(x,t) - \mathbf{P}^* \propto \exp(\sigma t + iqx)$, where $\sigma$ is the growth rate of wavenumber $q$. On an infinite domain, we obtain three dispersion relations for $\sigma(q;s)$. The HSS $\mathbf{P}^*_0$ is linearly stable down to an LW at $s \approx 0.548$, which is after the $\mathbf{P}^*_{0/1}$ saddle node. The HSS $\mathbf{P}^*_2$ exhibits a codimension-2 long wavelength (LW) / finite wavenumber Hopf (WB) instability at $(s, b) = (s_c, b_c) \approx (0.406, 0.067)$, as shown in the inset of Fig. 1. The former instability occurs for $q \to 0$ and leads to the bifurcation of steady states in subcritical direction (towards a stable portion of $\mathbf{P}^*_2$, $s < s_c$). The latter gives rise to both travelling and standing waves (TWs and SWs with $q_c \approx 2.031$, respectively) that bifurcate supercritically (towards the unstable portion of $\mathbf{P}^*_2$, $s > s_c$) and where TWs are linearly stable (here, SWs are ignored and for further details we refer to [1]).

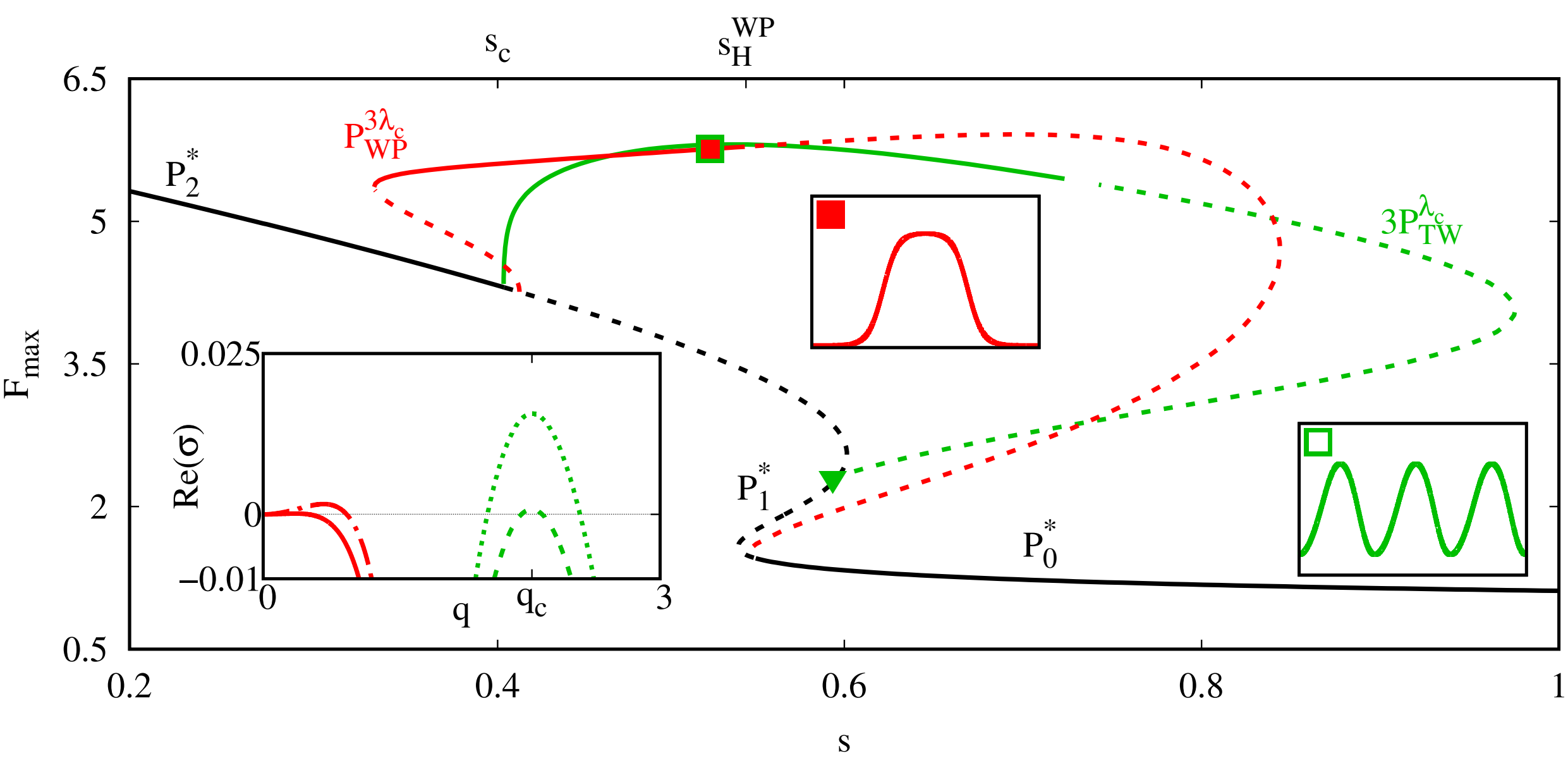

Figure 1: Bifurcation diagram showing the uniform steady states $\mathbf{P}^*_{0,1,2}$ (black), polarized mesas $\mathbf{P}^{3\lambda_c}_{\text{WP}}$ (red), and travelling waves $3\mathbf{P}^{\lambda_c}_{\text{TW}}$ (green), where $\lambda_c \approx 3.093$. The notation $3\mathbf{P}^{\lambda_c}_{\text{TW}}$ denotes the TW branch with wavelength $\lambda_c$ and three spatial copies. Solid lines indicate linear stability and dashed indicate instability. The lower-left inset shows the dispersion relations of the codimension-2 instability onset at $s = s_c \approx 0.406$ and $q_c \approx 2.031$, and after the instability at $s \approx 0.417$; the solid (dashed-dotted) line refers to the real values while the dashed (dotted) are complex conjugates. A profile of a mesa and TW at $s \approx 0.5$ are shown in the remaining insets. The symbol '▼' marks the parity breaking bifurcation onset of $\mathbf{P}^{\lambda_c}_{\text{TW}}$ associated with the $\mathbf{P}^{\lambda_c}_{\text{WP}}$ branch (not shown here). Other parameters: $b = b_c \approx 0.067$, $\gamma = 3.557$, $\eta = 0.6$, $p_0 = 0.8$, $p_1 = 3.8$, $D = 0.1$, $D_v = 1$, $D_F = 0.001$, and $M = 2$.

Fig. 1 shows the resulting bifurcation structure with a finite domain length $L = 3\lambda_c$, where $\lambda_c = 2\pi/q_c \approx 3.093$. The branches were computed using AUTO. We use this domain length to capture the finite size of cells. However, the bifurcation structure is similar for large, finite domain lengths $L$. Stationary solutions emerge for $s \gtrsim s_c$ and stabilize at the initial fold as mesa solutions. The branch loses stability before another fold and terminates at another LW. The travelling waves emerge at $s = s_c$ as stable solutions and lose stability before reaching a fold and then terminating at a parity breaking bifurcation of mesa solutions with wavelength $\lambda_c$. For $s_c < s < s^{\text{WP}}_{\text{H}}$ both the mesas and travelling waves are linearly stable, where $s^{\text{WP}}_{\text{H}}$ denotes the Hopf instability onset of the mesa solutions. Therefore, in this parameter region, switching between these patterns and thus directed motion and cell ruffling can occur. Bifurcation structures with pulses, fronts, and other mesa and wave solutions can be found in [1].

# Reliable and locally space-time efficient *a posteriori* estimates for the transient wave equation

**Nicolas Hugot**[1,2,3,*], **Alexandre Imperiale**[1], **Martin Vohralík**[2,3]
[1]Université Paris-Saclay, CEA, List, F-91120, Palaiseau, France,
[2]Inria, 48 rue Barrault, 75647 Paris, France,
[3]CERMICS, Ecole nationale des ponts et chaussées, IP Paris, 77455 Marne-la-Vallée, France
*Email: nicolas.hugot@cea.fr

**Abstract**

We develop *a posteriori* error estimates for the acoustic wave equation discretised by the spectral finite element method in space and the leapfrog scheme in time. Guaranteed upper bounds are obtained using Prager–Synge-type inequalities, ensuring constant-free reliability of the estimators. Local space-time efficiency is established through the use of appropriately scaled space-time norms and bubble functions. We investigate flux-based and also residual-based estimators. Our approach relies on the space-time $\inf - \sup$ stable formulation introduced by [1]. Then, we investigate the adaptation of these estimators to the elastic wave equation, again discretised by the spectral finite element method in space and the leapfrog scheme in time. We conduct numerical experiments on typical ultrasonic non-destructive testing configurations. This illustrates the practical relevance of the proposed error estimation strategy to validate numerical simulations.

## 1 Introduction

For a final time $T > 0$, $\Omega$ a bounded polytopal open domain of $\mathbb{R}^d$, $d > 0$, $Q := (0, T) \times \Omega$, and $c$ uniformly positive and bounded over $\Omega$, we consider the acoustic wave equation:

$$\partial_{tt} u - \operatorname{div}(c\nabla u) = f \quad \text{on } Q, \tag{1}$$

with zero initial and boundary conditions. We design *a posteriori* estimates of the discretisation error to help decide whether a simulation should be considered as sufficiently precise or not. We also identify the space-time distribution of the error.

## 2 Setting

In [1], the authors introduce an Hilbert space $\mathcal{Y}_0$ and a new $\inf - \sup$-stable weak formulation of (1). The energy norm that equips $\mathcal{Y}_0$ allows us to characterise the discretisation error $||u - u_{h\tau}||_{\mathcal{Y}_0}$ as an appropriate dual norm of the residual.

Let $\mathcal{T}_h(\Omega)$ be a conforming mesh of $\Omega$, $U_h^p \subset L^2(\Omega)$ the space of all element-wise polynomials up to a degree $p > 0$, and $V_h^p := U_h^p \cap H_0^1(\Omega)$. Let $N > 0$ be the number of time steps and $\tau := \frac{T}{N}$. We denote as $\mathcal{T}_\tau(0, T)$ the resulting mesh of the time domain and $\mathcal{T}_{h\tau}(Q) := \mathcal{T}_h(\Omega) \times \mathcal{T}_\tau(0, T)$ the space-time mesh.

The equation is discretised using the spectral finite element method in space [2], and the explicit leapfrog scheme in time. For $n \in \{0..N\}$, let $\ddot{u}_h^n$ be the discrete second time derivative $\frac{1}{\tau^2}(u_h^{n+1} - 2u_h^n + u_h^{n-1})$. The fully discrete equation is: For $n \in \{0..N\}$, find $u_h^{n+1} \in V_h^p$ such that for all $v_h \in V_h^p$:

$$(\ddot{u}_h^n, v_h)_\Omega + (c\nabla u_h^n, \nabla v_h)_\Omega = (f_h^n, v_h)_\Omega. \tag{2}$$

## 3 Discrete solution reconstruction

Since our numerical scheme produces snapshots of the discrete solution at each time step, whereas the variational setting expects full space-time objects, we use a reconstruction procedure to produce a discrete space-time solution $u_{h\tau}$ and source term $f_{h\tau}$. We reconstruct $u_{h\tau}$ in $\mathcal{P}^2(\mathcal{T}_\tau(0, T); V_h^p) \cap \mathcal{C}^0(0, T; V_h^p)$ and $f_{h\tau}$ in $\mathcal{P}^1(\mathcal{T}_\tau(0, T), U_h^p) \cap \mathcal{C}^0(0, T; U_h^p)$ such that the following Galerkin orthogonality holds: For $I_n \in \mathcal{T}_\tau(0, T)$ and $v_h \in V_h^p$,

$$\int_{I_n} (\partial_{tt} u_{h\tau}(t), v_h)_\Omega + (\nabla u_{h\tau}(t), \nabla v_h)_\Omega \mathrm{d}t = \int_{I_n} (f_{h\tau}(t), v_h)_\Omega \mathrm{d}t. \tag{3}$$

## 4 *A posteriori* estimates

Let $\boldsymbol{p}_{h\tau} \in \boldsymbol{L}^2(0, T; \boldsymbol{H}(\operatorname{div}, \Omega))$ and $s_{h\tau} \in H^1(0, T; L^2(\Omega))$ be discrete fluxes produced from the space and time derivatives of $u_{h\tau}$ by the averaging procedure from [3]. For $\star$ either $Q$ or an element of $\mathcal{T}_{h\tau}(Q)$, we introduce:

$$\begin{aligned}\eta_\star^{\mathrm{R}} &:= h\|f_{h\tau} - \partial_t s_{h\tau} - \operatorname{div}\boldsymbol{p}_{h\tau}\|_\star,\\ \eta_\star^{\mathrm{JX}} &:= \|c^{-1/2}\boldsymbol{p}_{h\tau} + c^{1/2}\nabla u_{h\tau}\|_\star,\\ \eta_\star^{\mathrm{JT}} &:= \|s_{h\tau} - \partial_t u_{h\tau}\|_\star,\\ \eta_\star^{\mathrm{J}} &:= \sqrt{(\eta_\star^{\mathrm{JX}})^2 + (\eta_\star^{\mathrm{JT}})^2}.\end{aligned}$$

Our main result is the guaranteed upper bound and local space time efficiency of this estimator.

**Theorem 1** *Assume that $f = f_{h\tau}$ and that $f(0) = 0$. Then there holds:*

- *Guaranteed upper bound:*

$$\|u - u_{h\tau}\|_{\mathcal{Y}_0} \le c_{\mathrm{rel}}^{\mathrm{R}}\eta_Q^{\mathrm{R}} + \eta_Q^{\mathrm{J}}.$$

- *Local space-time efficiency:*

$$\eta_{K\times I_n}^{\mathrm{R}} + \eta_{K\times I_n}^{\mathrm{J}} \le c_{\mathrm{eff}}\left(\|u - u_{h\tau}\|_{\mathcal{Y}_{0,\omega_K\times I_n}}\right).$$

Here, $\|u - u_{h\tau}\|_{\mathcal{Y}_{0,\omega_K\times I_n}}$ is a localised (around $K \times I_n$) version of the full $\mathcal{Y}_0$-norm. $c_{\mathrm{rel}}^{\mathrm{R}}$ and $c_{\mathrm{eff}}$ are constants that can only depend on the space dimension $d$ and the shape-regularity of $\mathcal{T}_h(\Omega)$; $c_{\mathrm{eff}}$ can additionally depend on the polynomial order $p$. We also show similar results for residual-based estimators.

## 5 Application to elastic waves

We heuristically extend the previous estimators to the elastic wave equation. Let $\lambda$ and $\mu$ be the piecewise constant Lamé coefficients from Hook's law and $\boldsymbol{\varepsilon}$ be the symmetric part of $\boldsymbol{\nabla u}$ in:

$$\begin{aligned}\rho\partial_{tt}\boldsymbol{u} - \operatorname{div}(\boldsymbol{\sigma}) &= \rho\boldsymbol{f},\\ \boldsymbol{\sigma} &= \lambda\operatorname{tr}(\boldsymbol{\varepsilon})\boldsymbol{I} + 2\mu\boldsymbol{\varepsilon}.\end{aligned} \tag{4}$$

We directly apply the estimators from above to the elastic case without any structural modification. For the fluxes in particular, we do not reconstruct them in the theoretically prescribed $\boldsymbol{H}(\mathrm{div})$-conform symmetric tensor space. Instead, we process each line independently, using the procedure from the acoustic case for each component. We study experimentally the behaviour of this extended procedure.

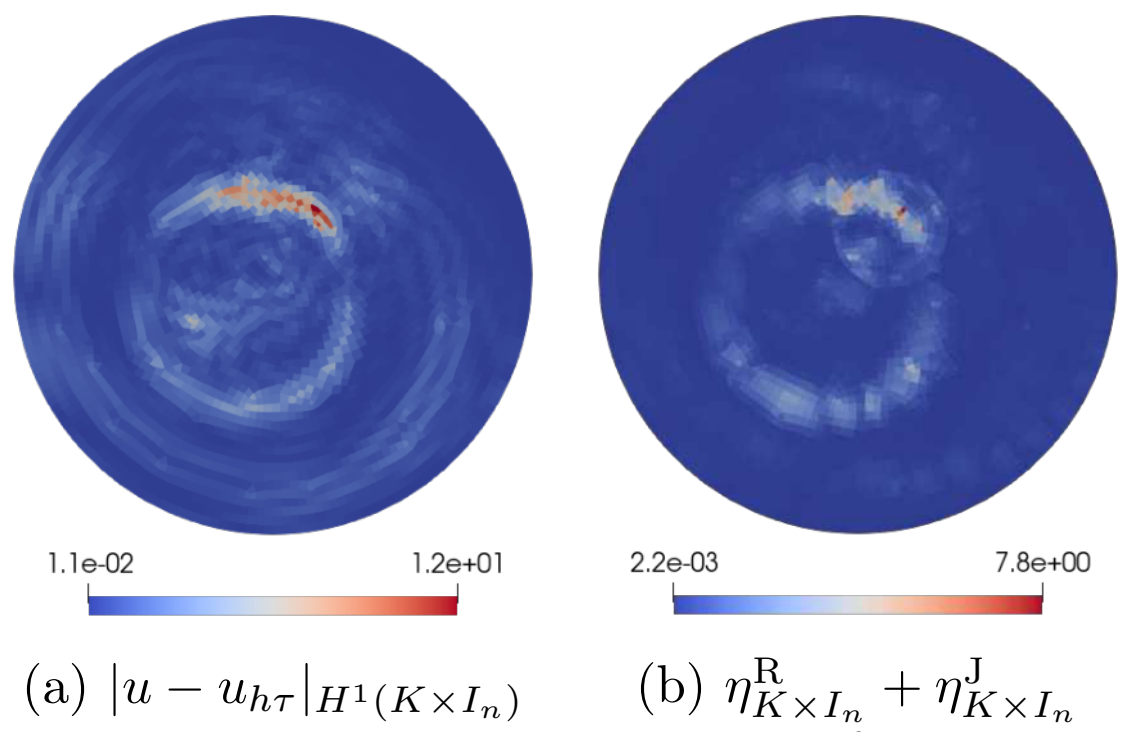

(a) $|u - u_{h\tau}|_{H^1(K\times I_n)}$ (b) $\eta_{K\times I_n}^{\mathrm{R}} + \eta_{K\times I_n}^{\mathrm{J}}$

Figure 1: Error maps at $t = 5.5$ for $T = 10$

## 6 Numerical results

We simulate an elastic wave generated by a source term on a circular domain $\Omega$. Figure 1 shows side-by-side local (space and time) maps of the discretisation error and our estimator. Figure 1 shows that the estimator is able to localise the most important contributions to the discretisation error in the space-time mesh. The global estimators also show predicted error values close to the real one.

## 7 Perspectives

We aim to investigate the adaptation of these estimators to industrial use, where computational cost is a critical constraint. In this context, we expect that the cheapest flux-based estimator previously available is still too expensive compared to the computation of the discrete solution using explicit time stepping. The goal is to progressively simplify the estimators into error indicators that are significantly less expensive to compute but still capture the overall size and distribution of the error in a meaningful way.

# Polynomial quasi-Trefftz spaces for scalar linear equations

**Lise-Marie Imbert-Gérard**[1,*]
[1]Department of Mathematics, University of Arizona, Tucson, USA
*Email: lmig@arizona.edu

## Abstract

Trefftz methods rely on function spaces of exact solutions to the governing equation. In particular, Trefftz Discontinuous Galerkin (DG) methods rely on finite-dimensional discrete spaces of local exact solutions to the governing equation to discretize a DG-type weak formulation. Instead, quasi-Trefftz DG methods rely on finite-dimensional discrete spaces of local exact solutions to the governing equation.

This talk will focus on a systematic study of local Taylor-based polynomial quasi-Trefftz spaces for problems governed by scalar linear equations with smooth variable coefficients, including a definition in terms of a quasi-Trefftz operator, the determination of the spaces dimension, and an explicit algorithm to construct a basis of the spaces.

## 1 Introduction

Trefftz DG methods exploit function spaces that incorporate key structural properties of the governing equation, namely spaces composed of locally exact solutions [2]. Their implementation therefore relies on constructing bases of locally exact solutions, which are generally unavailable for arbitrary scalar linear equations and in particular for variable-coefficient equations. To overcome this limitation, quasi-Trefftz DG methods exploit spaces composed of locally approximate solutions to the governing equation. Such methods rely on the approximation of a weak solution in a finite-dimensional quasi-Trefftz spaces, and their implementation requires to construct quasi-Trefftz bases. This talk, based on [1], will focus on polynomial quasi-Trefftz spaces for scalar linear equations.

In the context of problems governed by scalar linear equations with smooth variable coefficients, we will introduces a notion of quasi-Trefftz solutions based on a Taylor expansion, see Section 2, as well as an associated definition of polynomial quasi-Trefftz space, see Section 3. We will then discuss the dimension of this space, see Section 4, and finally present an explicit algorithm to construct a quasi-Trefftz basis, see Section 5.

## 2 Quasi-Trefftz property

Consider a linear differential operator $\mathcal{L}$ of order $\gamma \in \mathbb{N}$:

$$\mathcal{L} := \sum_{m=0}^{\gamma} \sum_{|\mathbf{j}|=m} c_{\mathbf{j}}(\mathbf{x})\, \partial_{\mathbf{x}}^{\mathbf{j}} \quad \text{for } \in \mathbb{R}^d,$$

where $d \in \mathbb{N}$ refers to the space dimension. The Trefftz property for a function $\varphi$ refers to the fact that $\mathcal{L}\varphi = 0$, in other words the functions belongs to the kernel of the differential operator. By contrast, the quasi-Trefftz property refers to the fact that the Taylor polynomial of $\mathcal{L}\varphi$ is zero up to a desired order. Hence, by essence, this is a local property and, in what follows, every result will be local. Given a point $\mathbf{x}_0 \in \mathbb{R}^d$ and a non-negative integer $q \in \mathbb{N}_0$, we denote by $\mathcal{T}_q$ the operator acting on $\mathcal{C}^q$ functions and computing the Taylor polynomial of degree $q$. Then a function $\varphi \in \mathcal{C}^{\gamma+q}$ satisfies the quasi-Trefftz property if

$$[\mathcal{T}_q \circ \mathcal{L}]\varphi = 0,$$

in other words the function belongs to the kernel of the so-called quasi-Trefftz operator $\mathcal{T}_q \circ \mathcal{L}$.

Thanks to this operator, we can then introduce the notion of polynomial quasi-Trefftz space with the following properties in comparison with standard polynomial spaces:

- it cannot approximate all smooth functions, but can approximate solutions to the PDE $\mathcal{L} = 0$,
- it can reach the same order of approximation as standard polynomial space with significantly less degrees of freedom.

## 3 Quasi-Trefftz space

Polynomial quasi-Trefftz functions are polynomials that satisfy a quasi-Trefftz property. Polynomial quasi-Trefftz spaces are characterized by an adequate choice of polynomial degree with respect to the order of the quasi-Trefftz property. More precisely, given a point $\mathbf{x}_0 \in \mathbb{R}^d$ and for any integer $p \geq \gamma$, we define the polynomial quasi-Trefftz operator

$$\begin{array}{rll} \mathcal{D}_p : & \mathbb{P}_p & \to \mathbb{P}_{p-\gamma} \\ & \Pi & \mapsto \mathcal{T}_{p-\gamma}[\mathcal{L}\Pi] \end{array}$$

so that the proposed quasi-Trefftz space can simply be defined as

$$\mathbb{QT}_p := \ker \mathcal{D}_p.$$

The graded structure of polynomial spaces in terms of homogeneous polynomial spaces, namely

$$\mathbb{P}_p = \bigoplus_{l=0}^{p} \widetilde{\mathbb{P}}_l,$$

will then play a crucial role in the study of $\mathbb{QT}_p$.

## 4 Quasi-Trefftz operator and dimension of the the quasi-Trefftz space

The operator $\mathcal{D}_p$ can be expressed thanks to an ingenuous additive splitting into (1) a "principal part", defined as the highest-order derivative terms with frozen coefficients at $\mathbf{x}_0$, namely

$$\mathcal{L}_* := \sum_{|\mathbf{j}|=\gamma} c_{\mathbf{j}}\left(\mathbf{x}_0\right) \partial_{\mathbf{x}}^{\mathbf{j}},$$

and (2) $\mathcal{R} := \mathcal{D}_p - \mathcal{L}_*$ the corresponding remainder. Leveraging the graded structure of polynomial spaces, we will then highlight specific block structures of these operators restricted to $\bigoplus_{\ell=\gamma}^{p} \widetilde{\mathbb{P}}_\ell$:

- $\mathcal{R}$ is block-strictly lower triangular,
- $\mathcal{L}_*$ is block-diagonal,
- so $\mathcal{L}$ is block-lower triangular.

Besides, unless it is trivial (i.e. $\mathcal{L}$ is of order less than $\gamma$ at $\mathbf{x}_0$), then $\mathcal{L}_*$ is surjective as an operator from $\mathbb{P}_p$ to $\mathbb{P}_{p-\gamma}$. Given a right inverse $\mathcal{S}_*$ of $\mathcal{L}_*$, we will then construct a right inverse $\mathcal{S}$ of $\mathcal{D}_p$ thanks to a block-forward substitution procedure. This can be expressed compactly as an explicit formula for $\mathcal{S}$ as follows:

$$\mathcal{S} := \mathcal{S}_* \circ \sum_{l=0}^{s} [-\mathcal{R} \circ \mathcal{S}_*]^l.$$

Hence the operator $\mathcal{D}_p$ is **surjective**, and we will show that the dimension of the quasi-Trefftz space $\mathbb{QT}_p$ follows from the rank-nullity theorem:

$$\dim \mathbb{QT}_p = \dim \mathbb{P}_p - \dim \mathbb{P}_{p-\gamma}.$$

## 5 Quasi-Trefftz basis

Assume that there exists a subspace $\mathbb{U} \subset \mathbb{P}_p$ such that $\mathbb{P}_p = \mathbb{U} \oplus \mathcal{S}(\mathbb{P}_{p-\gamma})$. Considering the linear equation $\mathcal{D}_p X = Y$ for any right-hand side $Y \in \mathbb{P}_{p-\gamma}$, the set of solutions can be described as $\{U + \mathcal{S}(Y - \mathcal{D}_p(U)), \forall U \in \mathbb{U}\}$. In particular, if $\{X_i, i \in I\}$ is a basis of $\mathbb{U}$, then $\{X_i - \mathcal{S}(\mathcal{D}_p(X_i)), \forall i \in I\}$ is a basis of $\mathbb{QT}_p$.

We will present a choice of space $\mathbb{U}$ and operator $\mathcal{S}_*$ so that the previous assumption is satisfied. Together with an explicit algorithm to implement $\mathcal{S}_*$, this will provide all the building block of the desired explicit algorithm to construct a basis of $\mathbb{QT}_p$.

## 6 Acknowledgement and statement

I have disclosed an outside interest in Airbus Central R&T to the University of Arizona. Conflicts of interest resulting from this interest are being managed by The University of Arizona in accordance with its policies.

I would also like to acknowledge support from the US National Science Foundation (NSF) (Grant No. DMS-2110407).

# Mathematical and numerical modelling of elastic wave propagation with smooth nonlinear contact conditions and friction losses

**Alexandre Imperiale**[1,*], **Bruno Lombard**[2], **Sébastien Imperiale**[3]
[1]Université Paris-Saclay, CEA, List, F-91120, Palaiseau, France
[2]Aix Marseille Univ, CNRS, Centrale Marseille, LMA UMR 7031, Marseille, France
[3]Inria, École polytechnique, CNRS, Palaiseau, France
*Email: alexandre.imperiale@cea.fr

## Abstract

We consider the problem of elastic wave propagation within two domains separated by smooth nonlinear contact conditions. These contact conditions link each component of the continuous normal stress to the jump of the displacement field (across the interface) through a nonlinear constitutive law that includes friction losses. We propose a mathematical analysis of the frictionless case, and prove an existence result for the associated Cauchy problem using the so-called "compactness" method [1]. For this purpose, we use a semi-discrete in time formulation of the problem yielding a bounded sequence (using energy arguments) converging to a solution of the continuous problem. From this approach, we derive a numerical scheme for the complete problem, with friction loss. This numerical scheme is explicit for the linear (volume) terms, and implicit for the nonlinear (interface) terms. We propose a stability analysis for this scheme, which entails a standard CFL condition depending only on the volume operators. We provide validation results against semi-analytical solutions in simplified configurations, and numerical illustrations in 2D.

## 1 Smooth contact with friction

We consider $\Omega_1$ and $\Omega_2$ two open and bounded subsets of $\mathbb{R}^2$ with Lipschitz boundaries $\partial\Omega_1$ and $\partial\Omega_2$ with outer normal fields $\{\boldsymbol{n}_i\}_{i=1}^2$, and $\Omega$ the open domain such that $\overline{\Omega} = \overline{\Omega_1} \cup \overline{\Omega_2}$. We look for $\boldsymbol{u}_1$ and $\boldsymbol{u}_2$ the displacement fields satisfying the linear elastic wave equations

$$\varrho_i \partial^2_{tt} \boldsymbol{u}_i - \operatorname{div} \boldsymbol{\sigma}_i = 0, \quad \boldsymbol{\sigma}_i = \boldsymbol{C}_i : \boldsymbol{\varepsilon}(\boldsymbol{u}_i), \quad \Omega_i,$$

with $\boldsymbol{u}_i(0) = \boldsymbol{u}_{i,0}$ and $\partial_t \boldsymbol{u}_i(0) = \boldsymbol{v}_{i,0}$, and where $\boldsymbol{\sigma}_i$ is the Cauchy stress tensor, $\varrho_i$ is the mass density of the solids, and $\boldsymbol{C}_i$ are uniformly bounded and coercive fourth order tensors. The domains share a common interface $\Gamma = \partial\Omega_1 \cap \partial\Omega_2$. On the exterior boundaries $\partial\Omega_i \backslash \Gamma$, we consider free surface conditions. At the interface, we impose continuity of the normal stress, $[\![\boldsymbol{\sigma}]\!] \cdot \boldsymbol{n} = 0$ where $[\![\cdot]\!]$ is the jump across $\Gamma$. The continuous normal stess is linked to the jump of displacement through the constitutive law [3]

$$\begin{aligned} \boldsymbol{\sigma} \cdot \boldsymbol{n} = K_N \, \mathcal{R}\big( [\![\boldsymbol{u}]\!] \cdot \boldsymbol{n} \big) \boldsymbol{n} + K_T \big( [\![\boldsymbol{u}]\!] \cdot \boldsymbol{t} \big) \boldsymbol{t} \\ + \mathcal{H}\big( [\![\boldsymbol{u}]\!] \cdot \boldsymbol{n} \big) \, \mathcal{G}\big( [\![\partial_t \boldsymbol{u}]\!] \cdot \boldsymbol{t} \big) \boldsymbol{t}, \end{aligned}$$

where $\boldsymbol{t} := \boldsymbol{n}^\perp$, and $K_{N,T} > 0$ are stiffness coefficients. The friction is governed by $\mathcal{G} \in \mathcal{C}^1(\mathbb{R})$, satisfying $\mathcal{G}(0) = 0$ and $\mathcal{G}' > 0$, weighted by a function of the normal displacement $\mathcal{H} > 0$. The nonlinear contact law $\mathcal{R} \in \mathcal{C}^0(\mathbb{R})$ satisfies

$$\forall u_N \in \mathbb{R}, \qquad \int_0^{u_N} \mathcal{R}(\xi) \, \mathrm{d}\xi \geq 0,$$

and the following asymptotic power growth property: $\exists r \geq 1$, $R_0 > 0$, and $\xi_0 > 0$ such that

$$\forall \xi \in \mathbb{R}, \quad |\xi| \geq \xi_0 \quad \Rightarrow \quad |\mathcal{R}(\xi)| \leq R_0 |\xi|^r.$$

## 2 Existence for the frictionless case

Our problem enters the functional framework described hereafter. Let $V$ and $H$ be two real Hilbert spaces such that $V \subset H \subset V'$, with $V$ separable and dense in $H$, with continuous injection. Let $A \in \mathcal{L}(V, V')$ be a linear, positive and symmetric operator satisfying a weak coercivity property. Let $W$ be a Banach space such that $V \subset W \subset H$, and the injection $V \mapsto W$ is continuous and compact. Let $R$ be a nonlinear operator satisfying the following assumptions:

- $R$ is the Fréchet derivative of a positive functional $\mathcal{F} : W \mapsto \mathbb{R}^+$ of class $\mathcal{C}^1$,
- $R$ maps bounded sets in $W$ into bounded sets in $W'$, *i.e.* $\forall c_R > 0$, $\exists C_R > 0$,

$$\forall u \in W, \quad \|u\|_W \leq c_R \;\Rightarrow\; \|R(u)\|_{W'} \leq C_R.$$

For $T > 0$, and $\forall q \in [1, +\infty]$ and any Banach space $E$, we denote $L^q(E) := L^q\big(0, T; E\big)$.

**Theorem 1** *For $u_0 \in V$, $u_1 \in H$, there exists a solution $u \in L^\infty(V)$ satisfying $du/dt \in L^\infty(H)$, and the space-time weak evolution equation*

$$\frac{d^2}{dt^2}u + Au + R(u) = 0, \quad \text{in } L^2(V'),$$

*with $u(0) = u_0$ in $V$, and $du(0)/dt = u_1$ in $H$.*

**Sketch of the proof.** Let $\tau > 0$ be a given time step. We perform Petrov-Galerkin (in time) approximation by looking for an approximate solution $u_\tau \in P^1_\tau(V)$ – where $P^k_\tau(V)$ is the space of (time) polynomial functions of order $k \geq 0$ with values in $V$. Considering $v_\tau \in P^1_\tau(H)$ the velocity variable, and using test functions in $P^0_\tau(V)$ and $P^0_\tau(H)$, we obtain

$$\begin{aligned} 0 = & \frac{v^{n+1} - v^n}{\tau} + A\frac{u^{n+1} + u^n}{2} \\ & + \frac{1}{\tau}\int_{I_n} R\Big(u^n + \frac{t - t^n}{\tau}(u^{n+1} - u^n)\Big)\mathrm{d}t, \end{aligned}$$

satisfied in $V'$, with $u^n := u_\tau(t^n)$ in $V$, $v^n := v_\tau(t^n)$ in $H$, and $t^n := n\tau$. This relation is completed with $(v^{n+1} - v^n)/\tau = (v^{n+1} + v^n)/2$ satisfied in $H$. Using continuity and boundedness of $R$, for each iteration we prove that there exists at least one solution $(u^{n+1}, v^{n+1})$, using Theorem 5.2.3 of [2]. The time discrete sequence satisfies a conservation property $\mathcal{E}_\tau(t^n) = \mathcal{E}_\tau(0)$ where $\mathcal{E}_\tau = \frac{1}{2}(\|v_\tau\|^2_H + \langle Au_\tau, u_\tau\rangle_{V',V}) + \mathcal{F}(u_\tau)$. This property implies boundedness of the sequences $u_\tau$, $v_\tau$, $du_\tau/dt$ and $dv_\tau/dt$ in $L^\infty(V)$, $L^\infty(H)$, $L^\infty(H)$ and $L^\infty(V')$. Hence, it weak* converges in corresponding dual spaces. To show that the limit is indeed a solution of the problem at hand, we use compactness arguments. We first remark that

$$(u_\tau) \subset \Big\{ v \;\Big|\; v \in L^2(V),\ \frac{dv}{dt} \in L^2(H)\Big\}.$$

Since $V \subset W \subset H$ with the compact inclusion $V \subset W$, using Theorem 5.1 in [1], the inclusion of this space in $L^2(W)$ is compact. We conclude by continuity and boundedness of $R$.

**Application to smooth contact.** Let $H := L^2(\Omega)^2$, and $V := \{\boldsymbol{v} \in H \,|\, \boldsymbol{v}|_{\Omega_i} \in H^1(\Omega_i)^2\}$. The linear operator reads

$$\langle A\boldsymbol{u}, \boldsymbol{v}\rangle_{V',V} = \sum_{i=1,2}\int_{\Omega_i}(\boldsymbol{C}_i : \boldsymbol{\varepsilon}(\boldsymbol{u}|_{\Omega_i})) : \boldsymbol{\varepsilon}(\boldsymbol{v}|_{\Omega_i})\,\mathrm{d}\Omega.$$

With $W := \{\boldsymbol{v} \in H \,|\, \boldsymbol{v}|_{\Omega_i} \in H^s(\Omega_i)^2\}$, where $s \in (\frac{1}{2}, 1)$, and $(1 - s)^{-1} \geq r + 1$, the inclusion $V \subset W$ is compact, and the functional

$$\mathcal{F}(\boldsymbol{u}) := \int_\Gamma K_N \int_0^{[\![\boldsymbol{u}]\!]\cdot\boldsymbol{n}} \mathcal{R}(\xi)\,\mathrm{d}\xi + \frac{K_T}{2}\big([\![\boldsymbol{u}]\!]\cdot\boldsymbol{t}\big)^2\,\mathrm{d}\Gamma,$$

is $\mathcal{C}^1$ and positive, with its Fréchet derivative mapping bounded sets into bounded sets, thanks to the assumptions on the contact law $\mathcal{R}$.

## 3 Stability of a fully discrete scheme

We consider $V_h \subset V$ a finite-dimensional approximation space of $V$. For all $(u_h, v_h) \in V_h$, we define $A_h \in \mathcal{L}(V_h)$ and $R_h(u_h) \in V_h$ by $(A_h u_h, v_h)_H := \langle Au_h, v_h\rangle_{V',V}$, $(R_h(u_h), v_h)_H := \langle R(u_h), v_h\rangle_{W',W}$. The fully discrete scheme reads

$$\begin{aligned} 0 = & \frac{u_h^{n+1} - 2u_h^n + u_h^{n-1}}{\tau^2} \\ & + A_h u_h^n + G_h^n\Big(\frac{u_h^{n+1} - u_h^{n-1}}{2\tau}\Big) \\ & + \int_0^1 R_h\Big(\xi\frac{u_h^{n+1} + u_h^n}{2} + (1-\xi)\frac{u_h^n + u_h^{n-1}}{2}\Big)\mathrm{d}\xi, \end{aligned}$$

where $G_h^n$ is a term depending exclusively on the friction. An energy analysis similar to the continuous case, shows stability upon the CFL condition:

$$\exists\alpha \in \mathbb{R}, \qquad \frac{\tau^2}{4}\|A_h\|_{\mathcal{L}(V_h)} \leq \alpha < 1.$$

In practice, $V_h$ is built using spectral finite elements. The mass-lumping technique yields point wise Newton-Raphson iterations. We provide validation results against semi-analytical solutions, and numerical illustrations in 2D.

# Mathematical models for wave propagation phenomena generated by landslides.

**Juliette Dubois**[4], **Natacha Gugan-Fau**[1], **Sebastien Imperiale** [1,*], **Anne Mangeney**[3], **Jacques Sainte-Marie**[2]

[1]Inria, Inria-Saclay - CMAP, Ecole Polytechnique, CNRS, Institut Polytechnique de Paris
[2] Inria, Inria Paris - Sorbonne Université, Paris
[3]Institut de Physique du Globe de Paris [4]TU Berlin, Institute of Mathematics

*Email: sebastien.imperiale@inria.fr

## Abstract

The application context of this work is the early detection of tsunamis generated by seafloor movements. These movements produce several types of waves: seismic waves in the ground, and gravity and acoustic waves in the ocean. Mathematically, the observed phenomena are governed by the prestressed nonlinear Euler equations in the fluid and the prestressed nonlinear elastodynamics equations in the solids. In these nonlinear equations, gravity effects are not neglected (at first). The objective of this work is to derive – through linearization – a hierarchy of linear equations that describe the propagation of these three types of waves in various asymptotic limits or for different types of fluids (e.g., purely acoustic or purely incompressible regimes, neglecting gravity effects in the solid, or barotropic fluids). Particular attention is given to the linearization of the transmission conditions used between the solid and fluid domains. We discuss the well-posedness of the resulting linear wave propagation problems and present a space-time discretization of some of these transient problems. This discretization employs high-order spectral finite elements in space, a leap-frog scheme in time, and perfectly matched layers. Notably, numerical experiments demonstrate that early detection of tsunamis could be achieved by analyzing the acoustic waves at the ocean surface, which are generated by the propagation of seismic waves (the so-called head waves).

## 1 Introduction

Our starting point is the prestressed nonlinear Euler equations in the fluid, $\Omega_F(t)$ and the nonlinear elastodynamics equations in the solid $\Omega_S(t)$, both retaining gravity effects, with free-surface boundary condition at the surface $\Gamma_F(t)$.

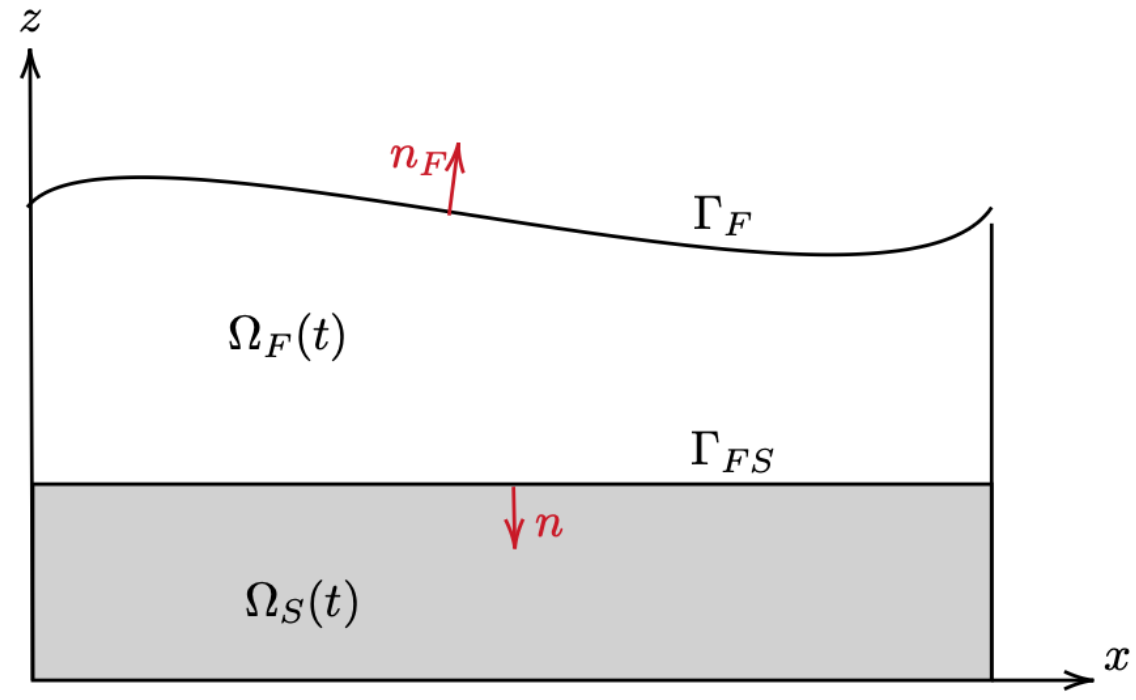

The domains represent the ocean and the Earth beneath, they are supposed unbounded in $\mathbb{R}^d$ and curvature effects of the earth are neglected. These equations are cast in Lagrangian coordinates, the reference domains being $\Omega_S \equiv \Omega_S(0)$ and $\Omega_F \equiv \Omega_F(0)$ (and $\Omega$ is the open union of the closure of these domains). We assume that: i) the configuration is at rest at the initial time; ii) non-deviatoric prestresses are neglected in the earth (see [3]); iii) the source term generates small perturbations; iv) the density $\rho_0$ of matter is smooth and depends only on $z$ across the complete domain; v) the shear in the water is non-zero. The last two assumptions are relaxed in a second step. With these assumptions, we derive the following non-standard variational formulation for the (linearized) displacement field $u(x,t) \in \mathbb{R}^d$, for all smooth test functions $v$,

$$\int_\Omega \rho_0 \ddot{u}\cdot v \,\mathrm{d}x + \int_\Omega \lambda \nabla\cdot u \nabla\cdot v + 2\mu\varepsilon(u) : \varepsilon(v)\,\mathrm{d}x - \int_\Omega \rho_0 g\big(\nabla\cdot u(v\cdot \mathbf{e}_z) + \nabla\cdot v(u\cdot\mathbf{e}_z)\big)\mathrm{d}x$$
$$\int_{\Gamma_F} \rho_0 g(u\cdot\mathbf{e}_z)(v\cdot\mathbf{e}_z)\,\mathrm{d}s - \int_\Omega \rho_0' g(u\cdot\mathbf{e}_z)(v\cdot\mathbf{e}_z)\,\mathrm{d}x = \int_{\Omega_S} f_S\cdot v\,\mathrm{d}x,$$

where $\varepsilon(\cdot)$ is the symmetrized gradient, $\lambda$ and $\mu$ are the Lamé parameters, $f_S$ the source term, $g$ the gravity constant. Vanishing initial data are used.

**Well-posedness** The equation enters a standard variational framework that allows to show existence and uniqueness of a solution $u$ in the space $C^1([0,T];L^2(\Omega)^d)\cap C^0([0,T];H^1(\Omega)^d)$ as soon as $f$ is smooth enough. This is the starting point to several simplifications. Assuming that $\lambda=\rho_0 c_0^2$ and letting $\mu\to 0$ in the fluid we recover a (well-posed) formulation similar to the one obtained in [1] in $\Omega_F$. Note that $u(\cdot,t)_{|\Omega_F}$ no longer belong to $H^1(\Omega_F)$. Finally allowing $\rho$ to lose continuity at the flat interface $\Gamma_{FS}$ and being constant in $\Omega_S$ we obtain the following symmetric variational formulation, for all smooth test functions

$$\int_\Omega \rho_0 \ddot{u}\cdot v\,\mathrm{d}x+\int_{\Omega_S}\lambda\nabla\cdot u\nabla\cdot v+2\mu\varepsilon(u):\varepsilon(v)\,\mathrm{d}x$$
$$+\int_{\Omega_F}\rho_0\big(c_0\nabla\cdot u-\frac{g}{c_0}(u\cdot\mathbf{e}_z)\big)\big(c_0\nabla\cdot v-\frac{g}{c_0}(v\cdot\mathbf{e}_z)\big)\,\mathrm{d}x$$
$$+\int_{\Omega_F}\rho_0 N^2(u\cdot\mathbf{e}_z)(v\cdot\mathbf{e}_z)\,\mathrm{d}x+g\int_{\Gamma_F}\rho_0(u\cdot\mathbf{e}_z)(v\cdot\mathbf{e}_z)\,\mathrm{d}s$$
$$-g\int_{\Gamma_{FS}}[\rho_0]_{\Gamma_{FS}}g(u\cdot\mathbf{e}_z)(v\cdot\mathbf{e}_z)\,\mathrm{d}x=\int_{\Omega_S}f_S\cdot v\,\mathrm{d}x,$$

where $N$ is a function of $z$ and is called the buoyancy frequency and where $[\rho_0]_{\Gamma_{FS}}$ is the difference between the density at the bottom of the ocean and at the surface of the earth (it is negative).

**The barotropic case without gravity effects in the solid.** First, a non-dimensionalization procedure shows that the term over $\Gamma_{FS}$ can be reasonably neglected compared to stiffness terms in $\Omega_S$. Moreover, an important simplification occurs when the fluid is considered barotropic, i.e., $N=0$ and $-\rho_0 g=\rho_0' c_0^2$. We observe that this implies $\rho_0\nabla\cdot u-g\rho_0/c_0^2(u\cdot\mathbf{e}_z)=\nabla\cdot(\rho_0 u)$. We obtain, using standard arguments, for all $t\in[0,T]$, and $u_R:=u_{|\Omega_R}$ with $R\in\{S,F\}$,

$$\begin{cases}\rho_0\ddot{u}_S-\mathrm{Div}\big(\lambda I\nabla\cdot(u_S)+2\mu\varepsilon(u_S)\big)=f_S & \Omega_S,\\ \rho_0\ddot{u}_F-\rho_0\nabla\big(\frac{c_0^2}{\rho_0}\nabla\cdot(\rho_0 u_F)\big)=0 & \Omega_F,\\ c_0^2\nabla\cdot(\rho_0 u_F)+g\rho_0(u\cdot\mathbf{e}_z)=0 & \Gamma_F,\\ c_0^2\nabla\cdot(\rho_0 u_F)\mathbf{e}_z=\lambda\nabla\cdot(u_S)\mathbf{e}_z+\mu\varepsilon(u_S)\mathbf{e}_z & \Gamma_{FS},\\ u_S\cdot\mathbf{e}_z=u_F\cdot\mathbf{e}_z & \Gamma_{FS}.\end{cases}$$

Introducing the scalar field $\dot{\phi}:=\frac{c_0^2}{\rho_0}\nabla\cdot(\rho_0 u_F)$, we derive using the above equations that $\dot{u}_F=\nabla\phi$, moreover

$$\ddot{\phi}-\frac{c_0^2}{\rho_0}\nabla\cdot(\rho_0\nabla\phi)=0\quad\Omega_F.$$

We also obtain the transmission conditions

$$\begin{cases}\rho_0\dot{\phi}\,\mathbf{e}_z=\lambda\nabla\cdot(u_S)\mathbf{e}_z+2\mu\varepsilon(u_S)\mathbf{e}_z & \Gamma_{FS},\\ \dot{u}_S\cdot\mathbf{e}_z=\nabla\phi\cdot\mathbf{e}_z & \Gamma_{FS},\end{cases}$$

and the standard boundary evolution equation

$$\ddot{\phi}+g\nabla\phi\cdot\mathbf{e}_z=0\quad\Gamma_F.$$

Doing so we recover sismo-hydro-acoustics equations. An energy preserving variational formulation can be derived allowing for efficient discretizations routine. Note that the equations without the barotropic assumption will be discussed and treated following the approach of [2].

**Numerical experiments.** The system is discretized using high-order spectral finite elements in space, exploiting the diagonal mass matrix arising from Gauss–Lobatto quadrature, and a second-order leap-frog explicit scheme in time. The Figure below illustrates a snapshot of representative results of an experiment modeling an underwater landslide source (at 1.5 km depth, domain width 48 km). The elastic model produces *head waves* that travel along the fluid-solid interface at the solid P-wave velocity, hence much faster than gravity waves (and around 3x time faster than acoustics waves), and appear as the first arrivals on hydrophone records.

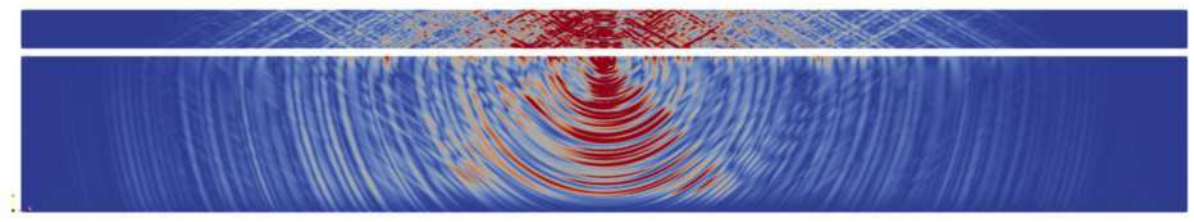

# High-order frequency derivatives of Helmholtz BIE solutions: fast computation and applications

Hiroshi Isakari[1,*]
[1]Faculty of Science and Technology, Keio University, Yokohama, Japan
*Email: isa@keio.jp

**Abstract**

This paper concerns frequency-dependent linear systems arising from boundary integral equations (BIEs) for Helmholtz' equation. We propose a fast method to compute high-order frequency derivatives of the BIE solutions at a reference frequency. The key observation is that all derivative orders can be obtained by solving a sequence of linear systems that share the same coefficient matrix evaluated at the reference frequency; only the right-hand side changes from order to order. To generate these right-hand sides efficiently, we combine forward-mode automatic differentiation with the fast multipole method (FMM), enabling high-order derivative computation with minimal analytic manipulations for kernels and straightforward support for complex frequencies.

## 1 Introduction

We consider the following frequency-dependent linear system

$$\mathsf{A}(\omega)\boldsymbol{u}(\omega) = \boldsymbol{b}(\omega), \tag{1}$$

where $\mathsf{A} : \mathbb{C} \to \mathbb{C}^{n\times n}$ is obtained by discretizing a BIE for a Helmholtz boundary value problem (typically defined in an unbounded domain), and $\omega$ is the (possibly complex-valued) angular frequency. The vectors $\boldsymbol{u}(\omega), \boldsymbol{b}(\omega) \in \mathbb{C}^n$ are the discrete counterparts of the unknown and known boundary densities, respectively, and $n$ denotes the number of degrees of freedom.

This paper develops an efficient method to compute high-order frequency derivatives of the BIE solution at a reference (also possibly complex-valued) frequency $\omega_0$, which shall be denoted as $\boldsymbol{u}^{(m)}$ with the differential order $m$.

Derivative-based rational approximations in BIE date back, at least, to Coyette et al. [1], where a single-frequency BIE solution together with a Padé expansion is used to obtain multi-frequency responses. Since then, fast frequency sweeps and uncertainty quantification using frequency derivatives continue to be actively studied; see, e.g., Marburg [2] for a recent overview. Motivated by recent interest in broadband wave-device design, in this study, we apply derivative-based Padé approximants to topology optimization (TO) of phononic structures (PnS) with wide stopbands, enabling efficient wideband objective and topological sensitivity evaluations.

In parallel, there has been strong recent interest in nonlinear eigenvalue problems (NEPs) induced by BIE for exterior scattering, where resonances are sought in the complex frequency plane. Most successful NEP solvers are based on contour integrals of resolvents and their moments. These approaches repeatedly evaluate solutions at sampling frequencies along a contour in the complex plane, making their dominant computational burden related to wideband frequency sweeps. Recent work has exploited *derivative-free* resolvent surrogates. For example, [3] approximates randomly scalarized resolvents using the AAA rational approximation to locate resonances. Related adaptive rational-approximation and sketching ideas for resonance identification have also been studied in [4].

In this work, we complement these derivative-free approaches by investigating *derivative-based* Padé models for resolvent quantities, enabled by our fast high-order differentiation framework, and demonstrate how the resulting pole information can guide contour selection and preprocessing in contour-integral eigensolvers.

## 2 Method

$m$-th order derivatives of the solution of (1) at $\omega_0$ can naively be evaluated as

$$\mathsf{A}\boldsymbol{u}^{(m)} = \boldsymbol{b}^{(m)} - \sum_{k=0}^{m-1} \binom{m}{k} \mathsf{A}^{(m-k)}\boldsymbol{u}^{(k)}. \tag{2}$$

Instead of expanding Leibniz sums, we here introduce the auxiliary vector-valued function $\boldsymbol{v}(\omega)$

satisfying $\boldsymbol{v}^{(k)}(\omega_0) = \boldsymbol{u}^{(k)}(\omega_0)$ for $k = 0, \ldots, m-1$ and $\boldsymbol{v}^{(m)} = \mathbf{0}$, with which (2) is rewritten as $\mathsf{A}\boldsymbol{u}^{(m)} = \boldsymbol{b}^{(m)} - (\mathsf{A}\boldsymbol{v})^{(m)}$. We obtain its right-hand side by forward-mode automatic differentiation combined with the FMM applied to the matrix-vector multiplication $\mathsf{A}(\omega_0)\boldsymbol{v}(\omega_0)$.

## 3 TO of PnS with a wide stopband

We consider a singly periodic array of rigid inclusions in the homogeneous acoustic medium. For a time-harmonic plane wave at angular frequency $\omega$ incidence, the energy transmission coefficient can be evaluated from far-field coefficients, which can be computed from boundary integrals of the solution of (1). We formulate a topology optimization problem that aims to realize a wide stopband by minimizing the band-averaged transmission under normal incidence. Specifically, the objective is the average transmission over a prescribed frequency interval, expressed as the angular frequency integral.

To evaluate the objective efficiently, we approximate the transmission as a function of $\omega$ by Padé approximants constructed from high-order frequency derivatives of the solution of (1) at a set of adaptively selected expansion points $\{\omega_0\}$ within the target band. The same derivative-based Padé strategy is also used to compute the corresponding topological sensitivity.

Starting from the initial configuration shown in Fig. 1 (left): circular rigid inclusions of radius 0.2 arranged periodically with period $L = 4.0$, we explored an optimized structure. For the Padé approximation, we used cubic numerator and denominator, which requires derivatives up to fifth order. The optimized design (Fig. 1, right) achieves a substantial reduction of the objective value to $9.015 \times 10^{-3}$ times that of the initial configuration.

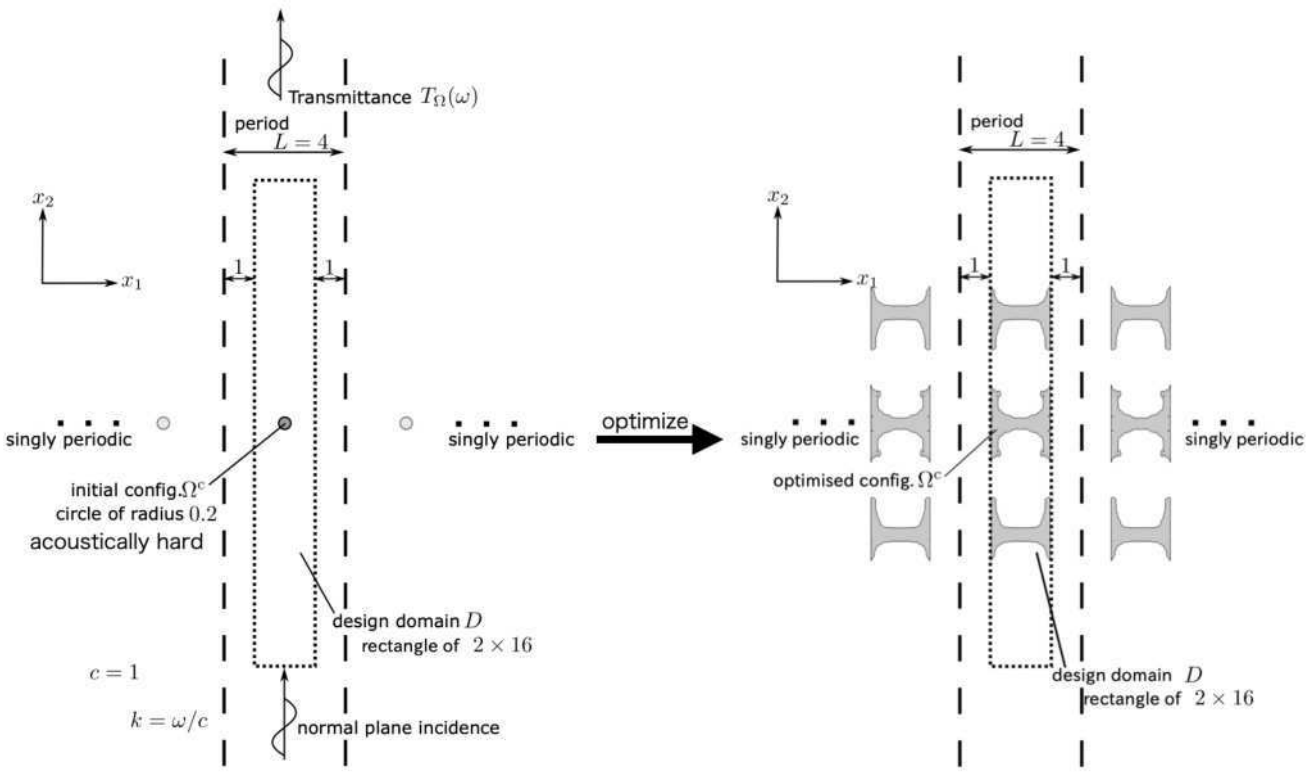

Figure 1: Initial (left) and optimized (right) phononic structures.

## 4 NEP esstimation

We consider the homogeneous BIE (1) arising from a 2D acoustic transmission problem: a unit-disk inclusion (wave speed $\sqrt{2}$) embedded in an infinite background medium (wave speed 1). We search for scattering resonance frequencies, i.e., complex $\omega$ for which (1) admits a nontrivial solution, in the window $0.1 < \Re\omega < 10.1$ and $-2.4 < \Im\omega < 0.1$. We tile the search window into $64 \times 16$ square subdomains and run a Sakurai–Sugiura contour-integral eigensolver on each subdomain (moments by 16-point Gauss–Legendre quadrature per edge), which required $4.37 \times 10^5$ s. Using Padé poles of a (randomized) resolvent surrogate built at each subdomain center to screen subdomains and running the eigensolver only when a pole falls inside, we recovered all resonances in $7.44 \times 10^4$ s.

# Analysis of Wave Propagation in Time-Dependent Media

**Patrick Joly**[1,*], **Marie Touboul**[1]
[1]POEMS, ENSTA, CNRS, INRIA, Institut Polytechnique de Paris, 91120, Palaiseau, France
*Email: patrick.joly@inria.fr

**Abstract**

Recent developments in the theory of metamaterials has raised attention on wave propagation in materials with time-dependent properties since they allow to break some limits of space metamaterials. It seems that the study of fundamental theory of wave equation with space-time dependent coefficents has not received much attention from the applied mathematics community. We aim at contributing to this domain by considering fundamental questions such as existence, uniqueness of the solution, behavior of the corresponding energy, and role of non-smooth coefficients via the study of scattering coefficients. We will provide various results in that direction and illustrate them through time-domain simulations.

## 1 Introduction

We are looking at the wave equation in $\mathbb{R}^d$ :

$$\partial_t\big(\rho(\boldsymbol{x},t)\,\partial_t u\big) - \operatorname{div}\big(\mu(\boldsymbol{x},t)\nabla u\big) = 0, \quad (1)$$

where $\boldsymbol{x} \in \mathbb{R}^d$, $t > 0$, and $\rho$ and $\mu$ are bounded from below and above by strictly positive quantities. We shall be interested by the Cauchy problem, completing (1) by

$$u(\boldsymbol{x},0) = u_0(\boldsymbol{x}), \quad \partial_t u_0(\boldsymbol{x}) = u_1(\boldsymbol{x}). \quad (2)$$

We are looking for finite-energy solutions, i.e. solutions for which the energy defined by

$$E(t) := \int_{\mathbb{R}^m} e(\boldsymbol{x},t)\,d\boldsymbol{x}, \quad (3)$$

with $e := \frac{1}{2}\,\rho\,|\partial_t u|^2 + \frac{1}{2}\,\mu\,|\nabla u|^2$, is finite at any time.

## 2 Energy based analysis

**Differentiable coefficients in time.** The case where $\rho(\boldsymbol{x},t)$ is independent of $\boldsymbol{x}$ and $t$ has been treated in [1]. The energy identity associated to (1) writes

$$\frac{dE}{dt} = -\tfrac{1}{2}\int \partial_t\rho\,|\partial_t u|^2\,d\boldsymbol{x} + \tfrac{1}{2}\int \partial_t\mu\,|\nabla u|^2\,d\boldsymbol{x},$$

from which we deduce that: (i) provided that $\frac{\partial_t\rho}{\rho}$ and $\frac{\partial_t\mu}{\mu}$ are bounded, the problem (1) admits a unique solution whose energy is exponentially bounded in time. (ii) under the assumptions that $\partial_t\rho \geq 0$ and $\partial_t\mu \leq 0$, the energy $E(t)$ is a decreasing function of time.

By the energy approach, we can also establish and describe in detail the finite velocity propagation of the solution. More precisely, for $d = 2$, we assume the initial data satisfy

$$\operatorname{supp}\, u_0 \cup \operatorname{supp}\, u_1 \subset\ B(0,R) \quad (4)$$

and look for $\varphi(\theta,t) : [0,2\pi]\times\mathbb{R}^+ \to \mathbb{R}^+$, positive, periodic in $\theta$, such that if $\Omega(t)$ denotes the moving star shaped domain defined by

$$\Omega(t) = \{\boldsymbol{x} = r\,\boldsymbol{\theta} \in \mathbb{R}^2\ /\ r < \varphi(\theta,t)\}.$$

then

$$\operatorname{supp}\, u(\cdot,t) \subset\ \overline{\Omega(t)}. \quad (5)$$

One shows that it suffices to take $\varphi$ as the solution of the Hamilton-Jacobi type equation (eikonal equation in a moving domain)

$$\varphi\,\partial_t\varphi - c(\varphi,t)\,\big((\partial_\theta\,\varphi)^2 + \varphi^2\big)^{\frac{1}{2}} = 0, \quad (6)$$

with initial condition $\varphi(\theta,0) = R$.

**Medium in subsonic translation.**
For $d = 1$, a particular case, important with respect to some applications, is when the medium is a "moving" static media of the form

$$\rho(x,t) = \rho_0(x-vt), \quad \mu(x,t) = \mu_0(x-vt), \quad (7)$$

where $v$ is a constant velocity. In the subsonic case, when $v < c_0(x) := \sqrt{\mu_0(x)/\rho_0(x)}$, using the change of variable $x \to x - ct$ (moving frame), we prove that there is conservation of an alternative energy. This proves existence and uniqueness of the solution uniformly bounded in time, without any further regularity for the functions contrary to the general case. The analysis will be extended to the case of a non uniform translation (with smooth $x_0(t)$)

$$\rho(x,t) = \rho_0\big(x - x_0(t)\big), \quad \mu(x,t) = \mu_0\big(x - x_0(t)\big).$$

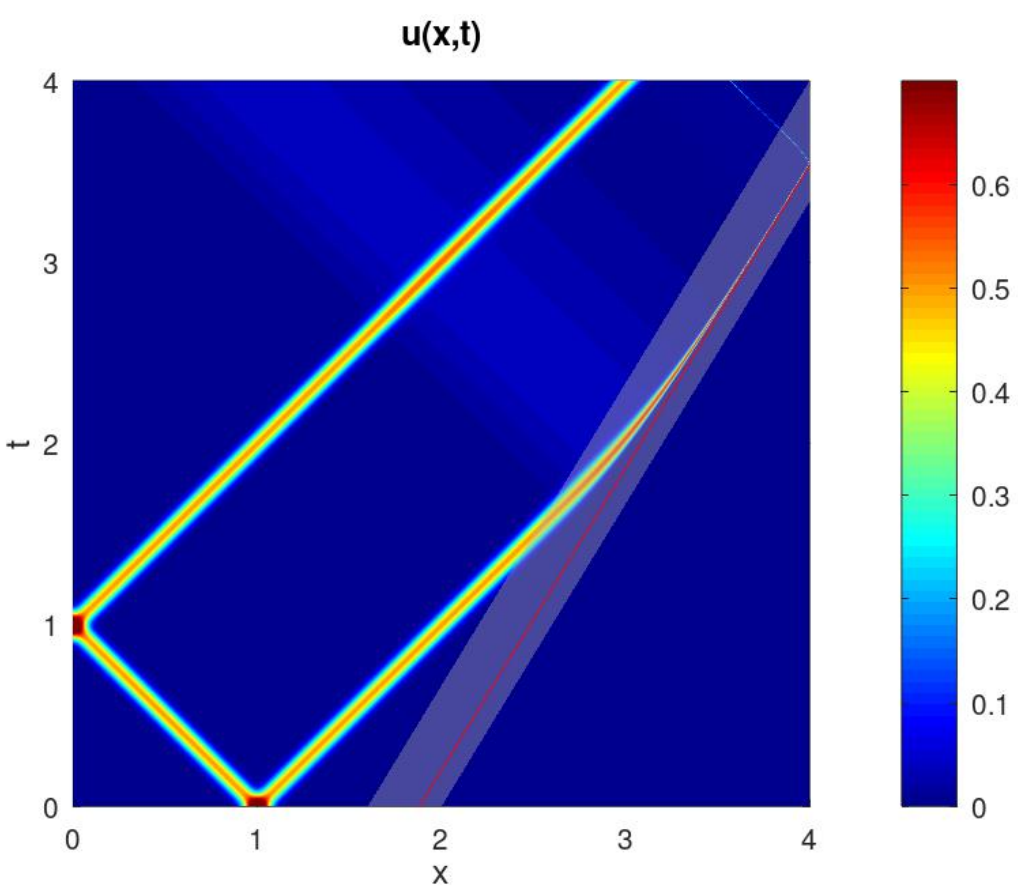

Figure 1: Trans-sonic case of non-existence. The wave does not cross the interface defined by $x = vt$ represented by the colored line.

## 3 Transmission problem across a moving interface in uniform translation

We consider in this section a particular 1D problem where the functions $\rho$ and $\mu$ are piecewise constant and take only two values $(\rho_1, \mu_1)$ for $(x,t) \in \mathcal{Q}_1$ and $(\rho_2, \mu_2)$ for $(x,t) \in \mathcal{Q}_2$, with

$$\begin{aligned} \mathcal{Q}_1 &= \{(x,t) \in \mathbb{R} \times \mathbb{R}_+ \,|\, x - vt < 0\}, \\ \mathcal{Q}_2 &= \{(x,t) \in \mathbb{R} \times \mathbb{R}_+ \,|\, x - vt > 0\}, \\ \Sigma &= \partial\mathcal{Q}_1 \cap \partial\mathcal{Q}_2 \equiv x = vt, \end{aligned} \tag{8}$$

where one can interpret $\Sigma$ as a moving interface. We introduce the velocities $c_j = \sqrt{\mu_j/\rho_j}$ and define $c_- = \min(c_1, c_2)$, $c_+ = \max(c_1, c_2)$. We distinguish three cases:

(i) the subsonic case: $|v| < c_-$,

(ii) the trans-sonic case: $c_- < |v| < c_+$,

(iii) the supersonic case: $|v| > c_+$.

In the subsonic and supersonic case, we compute the analytical solution of this problem, showing in particular existence and uniqueness.

On the other hand, in the trans-sonic case, we get either non-uniqueness or non-existence of the solution (see also [2, 3]). For instance, for $v > 0$ and initial data supported in $\{x < 0\}$, we get non-existence if $c_2 < v < c_1$ and non-uniqueness if $c_1 < v < c_2$. This motivates the introduction of a regularized interface in the trans-sonic case, which is the subject of the last section.

## 4 Regularized interface (trans-sonic case)

We focus on the trans-sonic case with $v > 0$ and initial data supported in $\{x < 0\}$. We regularize the interface by introducing a thickened interface of width $\varepsilon$, within which the speed $c(x,t)$ varies monotonically and continuously from $c_1$ to $c_2$ and moves with a velocity $v$:

$$c_\varepsilon(x) = c_1 + (c_2 - c_1)\,\chi\left(\tfrac{x}{\varepsilon}\right) \tag{9}$$

with $\chi(x) = 1$ if $x > 1$, $\chi(x) = 0$ if $x < -1$. Moreover, if $v = \theta\, c_1 + (1-\theta)\, c_2, \theta \in ]0,1[$, we assume that $\chi(0) = \theta$. Our goal is to understand the limit behaviour of the associated solution $u_\varepsilon$ when $\varepsilon \to 0$, in order to understand the non-existence or non-uniqueness issues raised in section 3. This is still work in progress. We expect that, when $c_2 > c_1$ (non-uniqueness for $\varepsilon = 0$), this procedure will select the natural solution for the limit problem. The case $c_2 < c_1$ is still an open question. However, a first remarkable result is that $u_\varepsilon$ remains supported in the domain $\mathcal{Q}_1$: in other words, the solution never crosses the interface.

# Transmission eigenvalues for biharmonic plate scattering: existence, discreteness, far-field recovery, and numerical computation

**Andreas Kleefeld**[1,*], **Isaac Harris**[2], **Heejin Lee**[3]
[1]Jülich Supercomputing Centre, Forschungszentrum Jülich GmbH, Jülich, Germany
[2]Department of Mathematics, Purdue University, West Lafayette, USA
[3]Zu Chongzhi Center for Mathematics and Computational Sciences, Duke Kunshan University, Kunshan, China
*Email: a.kleefeld@fz-juelich.de

## Abstract

Transmission eigenvalues for biharmonic scattering by a bounded cavity via a coupled interior–exterior Helmholtz formulation posed in all of $\mathbb{R}^2$ is studied. Discreteness and existence of infinitely many positive real eigenvalues are proven and it is shown how they can be identified from far-field data. Boundary integral equations are used to construct a numerical method to efficiently compute the transmission eigenvalues.

## 1 Introduction

Scattering by an impenetrable clamped cavity in a thin elastic plate, modeled by the Kirchhoff–Love (biharmonic) equation is studied. By splitting the biharmonic operator, a new transmission eigenvalue problem coupling an interior Helmholtz field in $D$ with an exterior modified Helmholtz field in $\mathbb{R}^2 \setminus D$, linked by transmission boundary conditions on $\partial D$ is obtained. Fundamental spectral results, including discreteness and existence of infinitely many positive real transmission eigenvalues via variational methods is shown (inspired by [2]). Further, these eigenvalues can be recovered from far-field data. Finally, a numerical method based on boundary integral equations is developed to efficiently compute transmission eigenvalues and used to support the theoretical findings such as interlacing with Dirichlet spectra. Further, they can be obtained from far-field data. Open problems that merit further investigation are given. Precisely, an interlacing property with Dirichlet and Neumann eigenvalues for all transmission eigenvalues is observed. Future work can go in the direction of surface concentration of the transmission eigenfunctions [3].

## 2 The model

The direct scattering model is biharmonic scattering in the exterior of a bounded clamped cavity $D \subset \mathbb{R}^2$. Given an incident field $u^i$, find the scattered field $u^s$ such that $\Delta^2 u^s - k^4 u^s = 0$ in $\mathbb{R}^2 \setminus D$, with transmission boundary conditions $u^s|_{\partial D} = -u^i$ and $\partial_\nu u^s|_{\partial D} = -\partial_\nu u^i$, plus a Sommerfeld-type radiation condition at infinity. Using operator splitting, write $u^s = u_H + u_M$ where $u_H$ solves $\Delta u_H + k^2 u_H = 0$ and $u_M$ solves $\Delta u_M - k^2 u_M = 0$ in $\mathbb{R}^2 \setminus D$, coupled on $\partial D$ by $u_H + u_M = -u^i$ and $\partial_\nu(u_H + u_M) = -\partial_\nu u^i$, with $u_H$ radiating (and $u_M$ exponentially decaying). The associated clamped transmission eigenvalue problem is: find $k > 0$ and a nontrivial pair $(v, w)$ such that

$$\Delta v + k^2 v = 0 \ \text{ in } D\,, \qquad (1)$$

$$\Delta w - k^2 w = 0 \ \text{ in } \mathbb{R}^2 \setminus D\,, \qquad (2)$$

and the transmission conditions

$$v + w = 0 \text{ on } \partial D\,, \qquad (3)$$

$$\partial_\nu(v + w) = 0 \text{ on } \partial D\,. \qquad (4)$$

This problem is derived when assuming that the corresponding far field pattern is zero.

## 3 Theoretical results

The set of transmission eigenvalues $k > 0$ satisfying (1)–(4) is at most a discrete set of the real line [5, Theorem 3.3]. In addition, there exist infinitely many positive transmission eigenvalues [4, Theorem 3.7]. Further, it is proven that the first transmission eigenvalue is smaller or equal to the first Dirichlet eigenvalue of the negative Laplacian [4, Theorem 3.8]. Moreover, it is shown that transmission eigenvalues can be recovered from far-field data [5, Theorem 4.1].

## 4 Numerical results

Boundary integral equations are used to solve the problem, since the domain of interest is unbounded. The resulting boundary integral equations is discretized with the boundary element collocation method (BEM) and the non-linear eigenvalue problem is solved with the Beyn algprithm [1]. In Table 1, the first four transmission eigenvalues (TE) $k > 0$ for a unit disk computed with BEM are listed and compared with the zeros of the determinant for $\ell \in \mathbb{N}_0$

$$\det \begin{pmatrix} H_\ell^{(1)}(\mathrm{i}k) & J_\ell(k) \\ \mathrm{i}H_\ell^{(1)'}(\mathrm{i}k) & J_\ell'(k) \end{pmatrix} = 0$$

to show that five digits accuracy are obtained when using 120 collocation nodes. Here, $J_\ell$ and $H_\ell^{(1)}$ denote the first kind Bessel and Hankel of order $\ell \in \mathbb{N}_0$, respectively. Additionally, the real-part of the eigenfunctions $w$ within $\mathbb{R}^2\backslash\overline{D}$ and $v$ within $D$ are numerically computed for a unit disk as shown in Figure 1 for $k_2 \approx 3.05164$ and $k_4 \approx 4.36453$. Four TE for an ellipse with semi-axis $(1, b)$ for various choices of $b$ are given in Table 2. Interestingly, we observe a monotonicity behavior with respect to all TEs when decreasing the area. The interlacing of the first four TEs with the Dirichlet (DE) and Neumann (NE) eigenvalues for a unit disk are given in Table 3, although only the interlacing between the first TE and first DE is proven. The interlacing property also holds for other scatterers. More numerical results for different scatterers as well as reconstructions from far-field data can be found in [4, Section 4] and [5, Section 5].

| TE | determinant | BEM |
|---|---|---|
| $k_1$ | 1.614634999563989 ($\ell = 0$) | 1.61464 |
| $k_2$ | 3.051633502815541 ($\ell = 1$) | 3.05164 |
| $k_3$ | 3.051633502815541 ($\ell = 1$) | 3.05164 |
| $k_4$ | 4.364516909785722 ($\ell = 2$) | 4.36453 |

Table 1: First four TE for a unit disk.

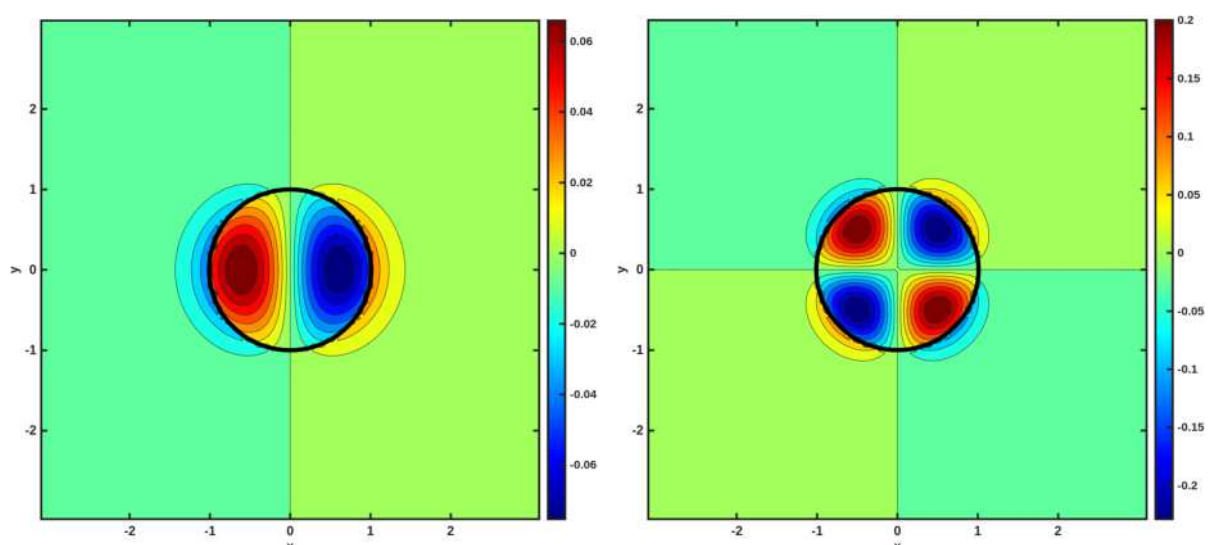

Figure 1: Real part of the eigenfunctions $w$ within $\mathbb{R}^2\backslash\overline{D}$ and $v$ within $D$ for the unit disk corresponding to $k_2 \approx 3.0564$ and $k_4 \approx 4.36453$.

| TE | $(1,1)$ | $(1,0.9)$ | $(1,0.8)$ | $(1,0.7)$ |
|---|---|---|---|---|
| $k_1$ | 1.61464 | 1.70401 | 1.81492 | 1.95646 |
| $k_2$ | 3.05164 | 3.13673 | 3.24382 | 3.38233 |
| $k_3$ | 3.05164 | 3.30629 | 3.62554 | 4.03737 |
| $k_4$ | 4.36453 | 4.54060 | 4.67571 | 4.82125 |

Table 2: First four TE for various ellipses.

| NE | 0.00000 | 1.84119 | 1.84119 | 3.05424 |
|---|---|---|---|---|
| TE | 1.61464 | 3.05164 | 3.05164 | 4.36453 |
| DE | 2.40483 | 3.83171 | 3.83171 | 5.13563 |

Table 3: First four DE, TE, NE for a unit disk.

## 5 Open questions

Can one proof that the interlacing property of the TE with DE and NE holds for any scatterer and all eigenvalues? Numerical results indicate a positive answer. Do complex-valued TE exist? Numerical results indicate a negative answer. What about giving a proof of monotonicity with respect to meas($D$), see [4, p. 16] or showing spectral patterns of transmission eigenfunctions (as in [3])?

# Domain derivative based shape reconstruction for time-dependent Maxwell's equations

Marvin Knöller[1,*]

[1]Department of Mathematics and Statistics, University of Helsinki, Helsinki, Finland

*Email: marvin.knoller@helsinki.fi

**Abstract**

In this talk we study the inverse scattering problem, which is to reconstruct a perfectly conducting scattering obstacle from a known incident field and measurements of the electric field at some measurement positions located away from the scattering obstacle. Our approach for this reconstruction is the use of a temporal domain derivative, which is applicable to perfectly conducting obstacles in free space. We first provide a short overview on Maxwell's equations and give a detailed introduction to the shape reconstruction problem. Subsequently we derive the temporal domain derivative for perfect conductors by proceeding through the Laplace domain. We apply the convolution quadrature method to the temporal domain derivative, which yields a provably convergent semi-discretization in time. We present several numerical examples for three-dimensional shape reconstructions and highlight algorithmic aspects and implementational difficulties.

## 1 Scattering from perfect conductors

Let the pair $(\boldsymbol{E}^i, \boldsymbol{H}^i)$ denote a pair of incident fields satisfying

$$\begin{aligned} \varepsilon_0 \partial_t \boldsymbol{E}^i - \mathbf{curl}\, \boldsymbol{H}^i &= 0 \quad \text{in } \mathbb{R}^3 \times (0,\infty)\,, \\ \mu_0 \partial_t \boldsymbol{H}^i + \mathbf{curl}\, \boldsymbol{E}^i &= 0 \quad \text{in } \mathbb{R}^3 \times (0,\infty)\,, \end{aligned}$$

where $\varepsilon_0$ and $\mu_0$ denote the electric permittivity and magnetic permeability of free space. As time proceeds $(\boldsymbol{E}^i, \boldsymbol{H}^i)$ propagate towards a perfectly conducting bounded and sufficiently regular scattering object $D$. On $\Gamma = \partial D$ they induce a current, which generates the scattered fields $(\boldsymbol{E}^s, \boldsymbol{H}^s)$ that satisfy

$$\varepsilon_0 \partial_t \boldsymbol{E}^s - \mathbf{curl}\, \boldsymbol{H}^s = 0 \quad \text{in } \Omega \times (0,\infty)\,, \tag{1a}$$

$$\mu_0 \partial_t \boldsymbol{H}^s + \mathbf{curl}\, \boldsymbol{E}^s = 0 \quad \text{in } \Omega \times (0,\infty)\,, \tag{1b}$$

$$\boldsymbol{\nu} \times \boldsymbol{E}^s = -\boldsymbol{\nu} \times \boldsymbol{E}^i \quad \text{on } \Gamma \times (0,\infty)\,. \tag{1c}$$

In our inverse scattering problem we assume to have knowledge on $(\boldsymbol{E}^i, \boldsymbol{H}^i)$ and that we can measure the time-dependent scattered electric field $\boldsymbol{E}^s$ at some receivers $\mathbf{z}_j$, $j = 1, \ldots, M$, which are located away from $D$.

## 2 A domain derivative in the time-domain

Let $\boldsymbol{h} \in C^1(\Gamma, \mathbb{R}^3)$ be a vector field on $\Gamma$ that is supposed to deform $D$. Let $R > 0$ be so large that $B_R(0)$ includes all scatterers under consideration. The vector field $\boldsymbol{h}$ can be extended to $B_R(0)$. We use this extended $\boldsymbol{h} \in C_0^1(B_R(0), \mathbb{R}^3)$ to define a diffeomorphism

$$\varphi : B_R(0) \to B_R(0)\,, \qquad \varphi(\boldsymbol{x}) = \boldsymbol{x} + \boldsymbol{h}(\boldsymbol{x})$$

and let $D_{\boldsymbol{h}} = \varphi(D)$ and $\Gamma_{\boldsymbol{h}} = \varphi(\Gamma)$. For the points $\mathbf{z}_j$, $j = 1, \ldots, M$ and for a fixed incident wave denote by $X$ the operator that maps the boundary $\Gamma$ to $(\boldsymbol{E}^s(\mathbf{z}_j))_{j=1,\ldots,M}$, i.e., $X(\Gamma) = (\boldsymbol{E}^s(\mathbf{z}_j))_{j=1,\ldots,M} \in \mathbb{R}^{3M} \times (0,\infty)$. Moreover, we define the operator

$$F_\Gamma(\boldsymbol{h}) = X(\Gamma_{\boldsymbol{h}}) \tag{2}$$

that maps from a sufficiently small neighborhood of the zero function in $C^1(\Gamma, \mathbb{R}^3)$ to $\mathbb{R}^{3M}$. We aim to determine the Fréchet derivative of $F_\Gamma$ at 0, i.e., the linear operator $F_\Gamma'(0)$ satisfying

$$\frac{1}{\|\boldsymbol{h}\|_{C^1(\Gamma)}} \|F_\Gamma(\boldsymbol{h}) - F_\Gamma(0) - F_\Gamma'(0)\boldsymbol{h}\|_Y \to 0$$

as $\|\boldsymbol{h}\|_{C^1(\Gamma)} \to 0$ for some appropriate norm $\|\cdot\|_Y$. This Fréchet derivative is characterized through the following theorem (see also the time-harmonic case in [2]).

**Theorem 1** *Under some regularity assumptions on $D$ and $\boldsymbol{E}^i$, the Fréchet derivative $F_\Gamma'(0)$ exists and is given by the solution $\boldsymbol{E}'$, which is a causal and tempered distribution, that satisfies*

$$\partial_t^2 \boldsymbol{E}' + \mathbf{curl}^2\, \boldsymbol{E}' = 0 \quad \textit{in } \Omega \times [0,\infty)\,, \tag{3a}$$

$$\boldsymbol{\nu} \times \boldsymbol{E}' = \mathrm{rhs}(\boldsymbol{E}, \boldsymbol{\nu}, \boldsymbol{h}_{\boldsymbol{\nu}}) \quad \textit{on } \Gamma \times [0,\infty) \tag{3b}$$

*evaluated at $\mathbf{z}_j$, i.e., $F_\Gamma'(0)\boldsymbol{h} = (\boldsymbol{E}'(\mathbf{z}_j))_{j=1,\ldots,M}$. The right hand side in* (3b) *depends on the total* *field $\boldsymbol{E} = \boldsymbol{E}^s + \boldsymbol{E}^i$ the exterior normal $\boldsymbol{\nu}$ to $D$ and on the normal component of the perturbation $\boldsymbol{h}_{\boldsymbol{\nu}}$.*

Using the Maxwell single layer potential $S(\partial_t)$ and the corresponding single layer operator $V(\partial_t)$ (see, e.g., [1]), the scattered field $\boldsymbol{E}^s$ from (1) and the temporal domain derivative $\boldsymbol{E}'$ from (3) can be written as

$$\begin{aligned}\boldsymbol{E}^s &= -S(\partial_t)V^{-1}(\partial_t)(\boldsymbol{\nu} \times \boldsymbol{E}^i)\\ \boldsymbol{E}' &= -S(\partial_t)V^{-1}(\partial_t)L(\partial_t)(\boldsymbol{\nu} \times \boldsymbol{E}^i)\,,\end{aligned}$$

where $L(\partial_t)$ arises from the boundary condition in (3b).

## 3 Reconstructing perfect conductors

We restrict ourselves to the reconstruction of star-shaped domains that are centered at the origin. Let measurements of the scattered waves at the receivers $\mathbf{z}_j$, $j = 1, \ldots, M$ for $N_t \in \mathbb{N}$ time steps be denoted by $\mathbf{g} \in \ell^2(0, N_t, \mathbb{R}^{3M})$. In our shape reconstruction algorithm we develop the radius $r = r(\alpha_n^m, \beta_n^m)$ of star-shaped objects into real-valued spherical harmonics, i.e.,

$$r = \sum_{n=0}^{N}\sum_{m=0}^{n} \alpha_n^m \operatorname{Re} Y_n^m + \sum_{n=1}^{N}\sum_{m=1}^{n} \beta_n^m \operatorname{Im} Y_n^m$$

Any star-shaped object's boundary, whose radius $r$ is parametrized by $\alpha_n^m, \beta_n^m$ is denoted by $\Gamma_{\alpha_n^m,\beta_n^m}$. We reformulate the inverse problem into a minimization problem, in which we aim to find $(\alpha_N^*, \beta_N^*) = \arg\min_{(\alpha_N,\beta_N)} f(\alpha_N, \beta_N)$, where

$$f(\alpha_N, \beta_N) = \frac{\left\| F_{\Gamma_{\alpha_n^m,\beta_n^m}}(0) - \mathbf{g} \right\|^2_{\ell^2_\tau(0,N_t;\mathbb{R}^{3M})}}{\|\mathbf{g}\|^2_{\ell^2_\tau(0,N_t;\mathbb{R}^{3M})}}$$

We introduce a regularization term $\psi$ , that aims at penalizing high oscillations. Finally, we minimize the regularized function

$$\gamma(\alpha_n^m, \beta_n^m) = f(\alpha_n^m, \beta_n^m) + \alpha\psi(\alpha_n^m, \beta_n^m)$$

using the Gauß–Newton algorithm. The approximation of both $F_{\Gamma_{\alpha_n^m,\beta_n^m}}(0)$ and $F'_{\Gamma_{\alpha_n^m,\beta_n^m}}(0)$ is done by using Rao–Wilton–Glisson and scaled Nédelec functions in a Galerkin formulation in space and a Runge–Kutta convolution quadrature method in time. Using the all-at-once formulation for convolution quadrature from [1], the numerical approximation in time reduces to several frequency domain problems, which can be solved in parallel on a high-performance cluster.

## 4 A numerical example

We study the situation, in which the exact scattering object is a cube that is aligned with the coordinate axes and centered at the origin with side length 3m. As an incident wave we pick a Gaussian beam with incident direction $[1, 0, 0]$, polarization $[0, 1, 1]$, a pulse width of 20ns and center frequency 30MHz. We pick 5 receivers located to the front, left, right, top and bottom and 6m away from the cube and measure the scattered electric field up to the final time $T = 600$ns. For the inverse problem we start with a ball centered at the origin with the radius 0.5m. In our implementation we use a spatial discretization of star-shaped scatterers with 3072 degrees of freedom and a temporal discretization of 400 time steps. The regularization parameter is chosen to be $\alpha = 0.01$. The reconstruction is found in Figure 1.

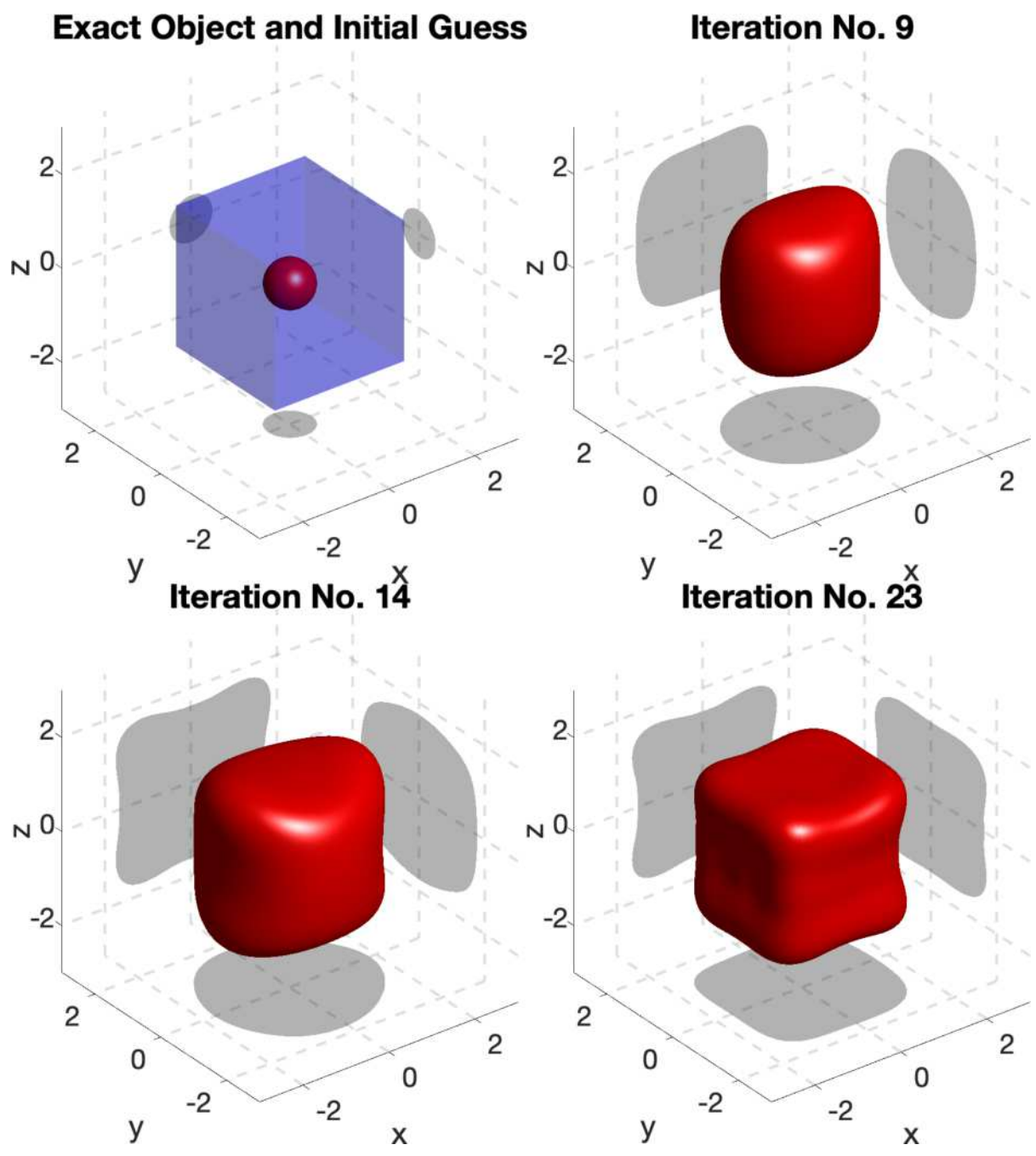

Figure 1: Reconstruction of the cube using 5 measurement points located around it.

# Optimized Schwarz methods in time for discrete transport control

**Duc Quang Bui**[1], **Bérangère Delourme**[2], **Laurence Halpern**[2], **Felix Kwok**[3,*]
[1]Laboratoire Manceau de Mathématiques, Le Mans Université, Le Mans, France
[2]LAGA, Université Sorbonne Paris Nord, Villetaneuse, France
[3]Department of Mathematics and Statistics, Université Laval, Québec, Canada
*Email: felix.kwok@mat.ulaval.ca

**Abstract**

Finding an optimal control numerically for a system governed by hyperbolic partial differential equations is a subtle problem, even for something as simple as the linear transport equation. In this paper, we consider the parallel-in-time solution to the discrete transport control problem: the time horizon is subdivided and the sub-interval problems are solved independently. Consistency is enforced by exchanging interface traces that are linear combinations of the state and adjoint unknowns and iterating to convergence. Using Fourier analysis, we analyze three different iterations: the fixed-point iteration, the relaxed iteration, and preconditioned GMRES. For each case, we propose optimized parameters for the transmission conditions that lead to fast convergence of the method. We illustrate our results by numerical examples.

## 1 Time-domain decomposition

In this work, we study the distributed control of a transport equation over the time horizon $[0,T]$. Let $y := y(x,t;v,y_{\text{ini}})$ be the solution of

$$\begin{cases} \partial_t y + \partial_x y = v & \text{in } \mathbb{R}\times(0,T), \\ y(\cdot,0) = y_{\text{ini}}, \end{cases}$$

where the initial conditions $y_{\text{ini}}$ is a continuous 1-periodic function in space, and the source term $v$ is in the *space of controls* $\mathscr{V}$, consisting of 1-periodic functions in space such that

$$\|v\|^2_{\mathscr{V}} = \int_0^T \int_0^1 |v(x,t)|^2\,dx\,dt < \infty.$$

Let the target state $y_{\text{tar}}$ be another continuous 1-periodic function in space, and consider the set of admissible controls defined by

$$\mathscr{K} = \{v \in \mathscr{V} : y = y(v, y_{\text{ini}}) \implies y(\cdot,T) = y_{\text{tar}}\}.$$

Over all admissible controls $v \in \mathscr{K}$, we seek the one with minimal $L^2$-norm; the unique minimizer $u \in \mathscr{V}\cap\mathcal{C}^1([0,T]\times\mathbb{R})$ is characterized by the optimality system [1, Proposition 2.1]

$$\begin{cases} \partial_t y + \partial_x y = u \text{ in } \mathbb{R}\times(0,T), \\ \partial_t u + \partial_x u = 0 \text{ in } \mathbb{R}\times(0,T), \\ y(\cdot,0) = y_{\text{ini}},\ y(\cdot,T) = y_{\text{tar}}. \end{cases} \tag{1}$$

The problem is then discretized using first-order upwind in space and semi-implicit Euler in time (explicit time-stepping and implicit source term):

$$\begin{cases} \dfrac{Y_j^m - Y_j^{m-1}}{\Delta t} + \dfrac{Y_j^{m-1} - Y_{j-1}^{m-1}}{\Delta x} = U_j^m, \\ \dfrac{U_j^m - U_j^{m-1}}{\Delta t} + \dfrac{U_{j+1}^m - U_j^m}{\Delta x} = 0, \end{cases} \tag{2}$$

with constant CFL number $\Delta t/\Delta x \in (0,1]$. To parallelize in time, we partition $[0,T]$ into $L$ disjoint subintervals $I_0,\dots,I_{L-1}$ of equal length $\Delta T = M\Delta t$. We then solve (2) for the local unknowns $(Y_{i,j}^m)$ and $(U_{i,j}^m)$ on each subinterval $I_i$ in parallel, with initial and final conditions

$$q_i Y_{i,j}^M + U_{i,j}^M = \xi_{i,j}^+, \qquad -p_i Y_{i,j}^0 + U_{i,j}^0 = \xi_{i,j}^-.$$

If one of the subinterval end points is $t=0$ or $t=T$, we impose $Y_{0,j}^0 = \mathbf{Y}_{\text{ini},j}$ or $Y_{L-1,j}^M = \mathbf{Y}_{\text{tar},j}$ instead. Since the values of $\xi_{i,j}^+$, $\xi_{i,j}^-$ that ensure continuity of $y$ and $u$ across subdomain interfaces are not known in advance, we must update these traces after each solve using neighbouring subdomain values:

$$\begin{aligned} \xi_{i,j}^+ &\leftarrow (1-\theta)\big(q_i Y_{i+1,j}^0 + U_{i+1,j}^0\big) + \theta\xi_{i,j}^+, \\ \xi_{i,j}^- &\leftarrow (1-\theta)\big(-p_i Y_{i-1,j}^M + U_{i-1,j}^M\big) + \theta\xi_{i,j}^-. \end{aligned}$$

When the $\xi_{i,j}^\pm$ converge, the solution becomes continuous and coincides with the global solution of (2), see [2,3]. The goal is to choose $p_i$, $q_i$ and $\theta \in [0,1)$ to optimize convergence. The case of $p_i = q_i$ has been analyzed in [1]. It was also shown there that the continuous variant of the algorithm, where the subinterval problems are based on (1) rather than (2), exhibits exact convergence after $L$ iterations if one chooses $p_i = \frac{1}{i\Delta T}$, $q_i = \frac{1}{(L-1-i)\Delta T}$ and $\theta = 0$. However, this no longer holds after discretization.

## 2 Optimizing interface conditions for two sub-intervals

For simplicity, we consider the two subinterval case and analyze the error system by letting $y_{\text{ini}} = y_{\text{tar}} = 0$. Here, we only have two sets of unknown traces[1] $(\xi_{0,j})$ and $(\xi_{1,j})$ and two interface parameters $p$ and $q$ to determine. Periodicity in space allows us to take the discrete Fourier transform of (2) in space. The transformed variables satisfy the recurrence

$$\begin{bmatrix} \hat{\xi}^k_{1,\ell} \\ \hat{\xi}^k_{0,\ell} \end{bmatrix} = \underbrace{\begin{bmatrix} \theta & (1-\theta)m_{12} \\ (1-\theta)m_{21} & \theta \end{bmatrix}}_{M_\theta(p,q,\ell)} \begin{bmatrix} \hat{\xi}^{k-1}_{1,\ell} \\ \hat{\xi}^{k-1}_{0,\ell} \end{bmatrix},$$

where $k$ is the iteration index, $\ell$ is the Fourier mode and $m_{12}$, $m_{21}$ depend on $p$, $q$ and $\ell$. In the unrelaxed case ($\theta = 0$), doing a double step leads to a scalar recurrence

$$\hat{\xi}^k_{i,\ell} = \rho(p,q,\ell)\,\hat{\xi}^{k-2}_{i,\ell}, \qquad i = 0, 1.$$

We seek $p, q$ that solve the min-max problem

$$\rho^* = \min_{p,q>0} \max_{\ell\in[0,N-1]} |\rho(p,q,\ell)|, \tag{3}$$

where $N = \frac{1}{\Delta x}$. Our main result is as follows:

**Theorem 1** *There exists $\Delta t_0$ such that, for any $\Delta t < \Delta t_0$, the maximum of $\ell \mapsto \rho(p,q,\ell)$ on $[0, N-1]$ is attained at $\ell^* = \ell^*(p,q) \in (0, N-1)$. The triple alternation problem*

$$-\rho(p,q,0) = -\rho(p,q,N-1) = \rho(p,q,\ell^*(p,q))$$

*has a unique solution $(p^*, q^*)$; this solution also solves the min-max problem (3). Moreover, as $\Delta t \to 0$, $(p^*, q^*)$ admit the asymptotic formulas*

$$p^* = \frac{c_p}{\Delta T} + o(1), \quad q^* = \frac{c_q}{\Delta T} + o(1) \tag{4}$$

*with $c_p \simeq 1.1993$ and $c_q \simeq 0.0906$. For these values, we have $\rho^* = c_\rho + o(1)$, where $c_\rho \simeq 0.0755$.*

The estimates of $\rho(p,q,\ell)$ required for the proof are also needed to analyze GMRES applied to the fixed-point equation $\boldsymbol{\xi} - \mathbf{M}_\theta(p,q)\boldsymbol{\xi} = 0$. Using the same optimization strategy as in [1], we find that for $(c_p, c_q) \simeq (2.1367, 0.3623)$, GMRES converges with the factor $\rho \simeq 0.0392$, which is faster than the fixed-point iteration.

[1] We omit $\pm$ from the superscript, since only one of them is valid for $\xi_{0,j}$ and $\xi_{1,j}$ in the two-subdomain case.

## 3 Numerical illustrations

We illustrate the performance of the unrelaxed algorithm for $y_{\text{ini}} = \sin(1-\pi x) + \sin(2\pi x)$ and $y_{\text{tar}} = \cos(14\pi x)$. We split the time horizon $[0,1]$ into two equal subdomains and discretize using $\Delta t = \frac{1}{240}$, $\Delta x = \frac{1}{120}$ (CFL $= \frac{1}{2}$). To show that our method converges even for unoptimized parameters, we choose arbitrarily $p = \frac{3}{\Delta T}$ and $q = \frac{1}{2\Delta T}$: Figure 1 shows the iterate $y^k$ along the interface $t = \Delta T$ after $k = 1$ and 3 iterations. Figure 2 shows the contraction factors for various $p$ and $q$, obtained numerically by running the algorithm for many pairs of $(p,q)$ and calculating the mean contraction factor in each case. The red dots, which show the asymptotic values $(p^*, q^*)$ found in Theorem 1, are close to the experimentally found minimum, as expected.

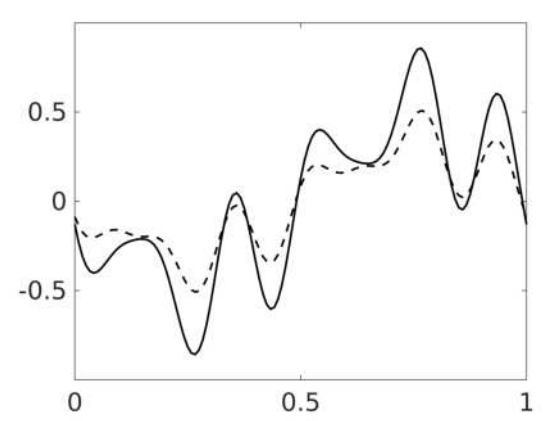

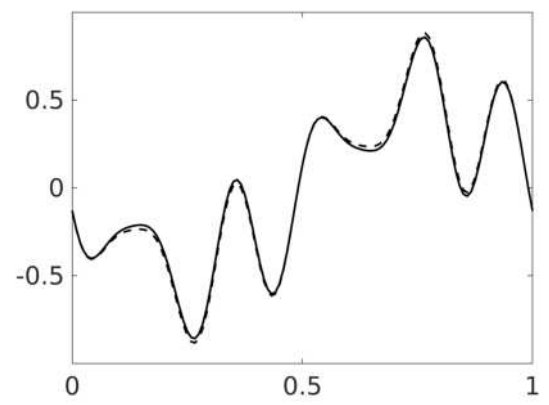

Figure 1: Iterate $y^k(\cdot, \Delta T)$ (dotted) versus the exact solution (solid). Left: $k = 1$, right: $k = 3$.

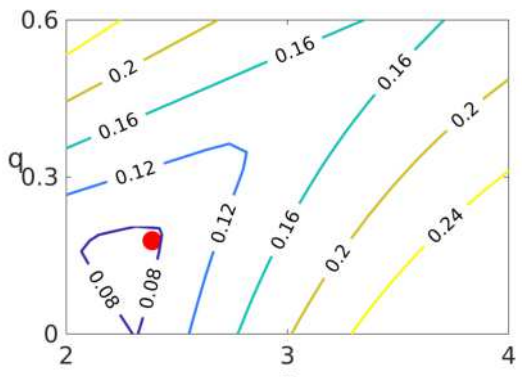

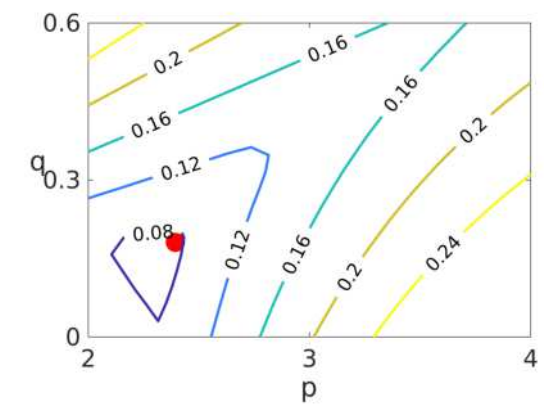

Figure 2: Discrete convergence factor as a function of $p$ and $q$. CFL $= \frac{1}{2}$ (left), CFL $= \frac{2}{3}$ (right).

**Propagation of the heave in a harbor: geometrical optics approach.**

__Olivier Lafitte__[1,*], **Emmanuel Audusse**[2], **Catherine Sulem**[3]
[1]LAGA, Université Sorbonne Paris Nord, Villetaneuse, France
[2]LAGA, Université Sorbonne Paris Nord, Villetaneuse, France
[3]Department of Mathematics, University of Toronto, Toronto, Canada
[*]Email: lafitte@math.univ-paris13.fr

**Abstract**

We study the waves in a harbor, with a given frequency $\omega$. The linearized system approximating the solution leads to a Dirichlet to Neumann operator for which $\frac{\omega^2}{g}$ is an element of the spectrum. We show that the WKB method that is classically used for wave asymptotics can be applied for finding the pseudo eigenvector (which does not belong to $L^2(\mathbb{R}^2)$) associated with a given element of the spectrum by writing a cascade of equations similar to the equations obtained for an hyperbolic system, without a high frequency parameter. This cascade of equations can be solved recursively and yields approximate solutions to the wave problem, for which we give the leading term.

## 1 Introduction

**Laplace equation yields waves.** Indeed, consider $\Omega = \{(x, y, z), x < 0, -h_0(x, y) < z < 0\}$, with $0 < \alpha \leq h_0 \leq M$ on $\mathbb{R}^2$ modeling a harbor. The linearized system on $\Psi(x, y, z, t)$ potential, $\eta(x, y, t)$ elevation of the surface is

$$\begin{cases} \partial_t \Psi(., 0) + g\eta = 0, & \text{linearized Bernouilli} \\ \partial_t \eta - \partial_z \Psi(., 0) = 0, & \text{displacement of the surface} \\ \Delta\Psi = 0 \text{ in } \Omega, & \text{incompressible irrotational} \\ \partial_n \Psi(., -h_0(.)) = 0 & \text{no penetration.} \end{cases} \tag{1}$$

If $h_0$ is constant, Airy [1] derived a solution of this system of frequency $\omega$ and of wave number $\mathbf{k}$ where $g|\mathbf{k}| \tanh |\mathbf{k}|h_0 = \omega^2$:
$\Psi(x, y, z, t) = \cos(\omega t) \cos(\mathbf{k}.(x, y)) \frac{\cosh(|\mathbf{k}|(z+h_0))}{\cosh(|\mathbf{k}|h_0)}$
$\eta(x, y, t) = \frac{\omega}{g} \sin(\omega t) \cos(\mathbf{k}.(x, y))$
How can one generalize this analysis to a variable bottom? The Dirichlet problem

$$\begin{cases} \Phi(., 0) = \phi^0(.), \\ \Delta\Phi = 0 \text{ in } \Omega, \\ \partial_n \Phi(., h_0(.)) = 0 \end{cases} \tag{2}$$

has a unique solution, which defines the associated Dirichlet to Neumann operator (DTN) $L(h_0) : \phi^0 \to \partial_z \Phi(., 0) = L(h_0)(\phi^0)$. This operator acts from $\dot{H}^{\frac{1}{2}}(\mathbb{R}^2)$ to $\dot{H}^{-\frac{1}{2}}(\mathbb{R}^2)$ [4] and is represented by a pseudodifferential operator [3]. Finding an oscillating in time solution of the water wave problem (1) is thus, for a given frequency $\omega$, finding a (complex) pseudomode $\phi_\omega$ associated with the eigenvalue $\lambda = \frac{\omega^2}{g}$ of the operator $L(h_0)$, deduce $\eta_\omega = -\frac{i\omega}{g}\phi_\omega$ and deduce the solution $\Phi_\omega$ of the Dirichlet problem (2).

We propose in this work a calculation of the solution of the water waves in $\Omega$ by using the methods of geometrical optics, and obtain the principal term of the solution (and define why it can be called a principal term). This renders the comparison of the Dirichlet to Neumann solution with the numerical solution of the linearized system with Robin and Neumann conditions (1) and with the WKB approx. solution of the mild slope equation ( [2]) on $\phi(x, y)$, $\Re(e^{-i\omega t}\phi(x, y)w(x, y, z))$ approx. solution of (1).

## 2 WKB analysis of the water wave system without asymptotic parameter

Using the formulation of Craig-Sulem [3], it is proven that there exists three pseudodifferential operators $Op(A), Op(C), Op(G_0)$ such that $Op(C)$ is invertible of inverse $Op(B)$ and $L(h_0) = Op(G_0) - Op(|\xi|)Op(B)Op(A)$. Consider $H_0 > 0$ a given height, write the solution $\Phi$ through its Fourier transform in $(x, y)$:

$$\begin{aligned} \Phi(x, y, z) &= Op(\tfrac{\cosh |\xi|(z+H_0)}{\cosh |\xi| H_0})(\phi) \\ &\quad + Op(\sinh |\xi| z)(\psi)(., z), \end{aligned}$$

$\phi$ and $\psi$ being the Dirichlet and Neumann traces. Operators $Op(A)$ and $Op(C)$ are given by

$$Op(A)(f(\mathbf{x}) = \int e^{i\xi.\mathbf{x}} \frac{\sinh |\xi|(H_0 - h_0(\mathbf{x}))}{\cosh(|\xi| H_0)} \hat{f}(\xi) d\xi,$$

$$Op(C)(\phi)(\mathbf{x}) = \int e^{i\xi.\mathbf{x}} \cosh |\xi|(h_0(\mathbf{x}))\hat{f}(\xi) d\xi,$$

whereas $Op(G_0)$ is the Fourier multiplier of symbol $|\xi| \tanh(H_0|\xi|)$. The symbols of $Op(A), Op(C)$ are denoted by $A(h_0, \xi), C(h_0, \xi)$.

Let us define now the action of an operator on $\mathbf{x} \to a(\mathbf{x})e^{i\theta(\mathbf{x})}$. We consider the clas-sical asymptotic expansion of a high frequency oscillating function $a(\mathbf{x},\epsilon)e^{i\theta(\mathbf{x})/\epsilon}$, namely with $a(\mathbf{x},\epsilon) \simeq \sum a_j(\mathbf{x})\epsilon^j$, as $Op(A(\mathbf{x},\xi))[a(\mathbf{x},\epsilon)e^{i\theta(\mathbf{x})/\epsilon}] = A(\mathbf{x},\frac{\nabla\theta}{\epsilon})a(\mathbf{x},\epsilon) + \epsilon L^1_{A,\theta}(a(\mathbf{x},\epsilon)) + \epsilon^2 L^2_{A,\theta}(a(\mathbf{x},\epsilon))\cdots$, where each operator $L^p_{A,\theta}$ involve derivatives of $a$ and of $\nabla\theta$. The action of $Op(A(\mathbf{x},\xi))$ on $[a_0(\mathbf{x}) + a_1(\mathbf{x}) + \cdots]e^{i\theta(\mathbf{x})}$ is **structured as**:

order 0: $A(h_0(\mathbf{x}),\nabla\theta)a_0(\mathbf{x})$,

order 1: $A(h_0(\mathbf{x}),\nabla\theta)a_1(\mathbf{x}) + \frac{1}{i}[\partial_\xi A.\nabla_{\mathbf{x}}a_0 + \frac{1}{2}\nabla^2_{\xi^2}A\nabla^2_{x^2}\theta a_0]$,

higher order terms: $\cdots$.

Finding a pseudo eigenmode of the DTN amounts to finding $a = (a_0, a_1, \cdots), b = (b_0, b_1, \cdots)$ and $\theta$ such that

$$\begin{pmatrix} Op(A) & Op(C) \\ Op(G_0) & Op(|\xi|) \end{pmatrix} \begin{pmatrix} ae^{i\theta} \\ be^{i\theta} \end{pmatrix} = \begin{pmatrix} 0 \\ \lambda ae^{i\theta} \end{pmatrix}.$$

This equality yields the cascade of equations:

$$\begin{cases} A(h_0,\nabla_x\theta)a + C(h_0,\nabla_x\theta)b & +\frac{1}{i}[\nabla_\xi A\nabla_x a + \frac{1}{2}\nabla^2_{\xi^2}A\nabla^2_{x^2}\theta a \\ & +\nabla_\xi C\nabla_x b + \frac{1}{2}\nabla^2_{\xi^2}A\nabla^2_{x^2}\theta b] = l.o.t. \\ G_0(|\nabla_x\theta|)a + |\nabla_x\theta|b & +\frac{1}{i}[\nabla_\xi G_0\nabla_x a + \frac{1}{2}\nabla^2_{\xi^2}A\nabla^2_{x^2}\theta a \\ & +\frac{\nabla_x\theta}{|\nabla_x\theta|}\nabla_x b] - \lambda a = l.o.t \end{cases}$$

We **impose** the following **hierarchical system** (similar to the high frequency expansion)

$$\begin{cases} Aa_0 + Cb_0 = 0 \\ G_0a_0 + |\nabla\theta|b_0 - \lambda a_0 = 0 \end{cases}$$

$$\begin{cases} Aa_1 + Cb_1 + & \frac{1}{i}[\nabla_\xi A\nabla_x a_0 + \frac{1}{2}\nabla^2_{\xi^2}A\nabla^2_{x^2}\theta a_0 \\ & +\nabla_\xi C\nabla_x c_0 + \frac{1}{2}\nabla^2_{\xi^2}C\nabla^2_{x^2}\theta b_0] = 0. \\ (G_0-\lambda)a_1 + |\nabla\theta|b_1 & +\frac{1}{i}[\nabla_\xi G_0\nabla_x a_0 + \frac{1}{2}\nabla^2_{\xi^2}G_0\nabla^2_{x^2}\theta a_0 \\ & +\frac{\nabla_x\theta.\nabla_x b_0}{|\nabla_x\theta|} + \frac{1}{2}\nabla^2_{\xi^2}|\xi|\nabla^2_{x^2}\theta b_0] = 0, \\ \cdots \end{cases}$$

The first relation is the existence of a nontrivial solution of $M_0(h_0,\nabla\theta)\begin{pmatrix} a_0 \\ b_0 \end{pmatrix} := \begin{pmatrix} Aa_0 + Cb_0 \\ (G_0-\lambda)a_0 + |\nabla\theta|b_0 \end{pmatrix} = 0$. This yields, after calculations:

$$|\nabla\theta(\mathbf{x})|\tanh h_0(\mathbf{x})|\nabla\theta(\mathbf{x})| = \lambda, \tag{3}$$

called the **eikonal equation**.

As $\begin{pmatrix} Aa_1 + Cb_1 \\ (G_0-\lambda)a_1 + |\nabla\theta|b_1 \end{pmatrix} \in \mathrm{Im}(M_0(h_0,\nabla\theta))$ for all $a_1, b_1$, one obtains a PDE on $a_0, b_0$. Using the relation $\begin{pmatrix} a_0 \\ b_0 \end{pmatrix} \in \mathrm{Ker}(M_0(h_0,\nabla\theta))$, one deduces a PDE on $a_0$, called the **homogeneous transport equation**:

$$E(\mathbf{x})\nabla_x\theta\nabla_x a_0 + K(\mathbf{x})a_0 = 0, \tag{4}$$

which translates on $a_0$ known on all characteristic curves of $\nabla\theta$, these curves can be interpreted as rays associated with the phase $\theta$.

The approximate solution of the water wave problem then reads

$$\begin{aligned} &\Psi(\mathbf{x},z,t) = \cos(\sqrt{g\lambda}t)Op(\sinh|\xi|z)(b_0(.)\cos(\theta(.)))(\mathbf{x},z) \\ &+(\cos(\sqrt{g\lambda}t)Op(\tfrac{\cosh|\xi|(z+H_0)}{\cosh|\xi|H_0})(a_0(.)\cos(\theta(.)))(\mathbf{x},z), \\ &\eta(\mathbf{x},t) = \sqrt{\tfrac{\lambda}{g}}\sin(\sqrt{g\lambda}t)a_0(\mathbf{x})\cos\theta(\mathbf{x}). \end{aligned} \tag{5}$$

## 3 Application to the 2d case

### 3.1 Exact principal solution

In the 2d case, one can integrate the homogeneous transport equation (4), the solution of the eikonal equation (3) being explicit, and one obtains rigorously the leading order terms $a_0$ and $b_0$ through the equalities (with $y = \frac{\lambda}{\theta'(x)}$)

$$\begin{aligned} &\theta' h_0\tanh\theta' h_0 = \tfrac{\omega^2 h_0}{g}, \theta(0) = 0, \\ &b_0(x) = -\tfrac{\cosh h_0(x)\theta'}{\cosh h_0(x)\theta'\cosh H_0\theta'}a_0(x), \\ &a_0(x) = [(1-y^2)\tanh^{-1}y + y]^{-\frac{1}{2}}, \end{aligned} \tag{6}$$

### 3.2 Action of the DTN on the GO formal approximation

Assume that $h_0(x) - H_0$ is compactly supported, of class $C^2$, and $\theta, a_0, b_0$ given in (6). We prove the following approximation result

**Proposition 1** *There exist $C_0 > 0$ such that, in the max norm,*
$|Op\begin{pmatrix} A & C \\ G_0-\lambda & k \end{pmatrix}[\begin{pmatrix} a_0 \\ b_0 \end{pmatrix}e^{i\theta}]| \le C_0(|h_0''| + |h_0'|^2)$.

The proof relies on the stationary phase theorem in its exact formulation, using $\theta''$ and $a_0'$ of the form $T(h_0)h_0'$ hence $\theta''', a_0''$ contains only terms proportional to $h_0''$ and $(h_0')^2$.

# Quasi-Trefftz Discontinuous Galerkin Method for the Convected Helmholtz Equation with Smooth Coefficients

Andréa Lagardère[1,2,*], Lise-Marie Imbert-Gerard[3], Guillaume Sylvand[1], Sébastien Tordeux[2]

[1]Airbus Central Research and Technology, Airbus, Toulouse, France
[2]Project-team Makutu, Inria, Pau, France
[3]Department of Mathematics, University of Arizona, Tucson, USA

*Email: andrea.lagardere@airbus.com

**Abstract**

This work presents a Quasi-Trefftz Discontinuous Galerkin (QTDG) method for the convected Helmholtz equation with spatially varying coefficients. The problem is formulated as a first-order system, and local basis functions are constructed using Taylor expansions to satisfy the governing equation up to a prescribed order. Numerical experiments in heterogeneous flows demonstrate that the method achieves the expected convergence rates with significantly fewer degrees of freedom than standard high-order polynomial methods.

## 1 Introduction

The numerical simulation of acoustic wave propagation with background flows in subsonic regime is essential for aircraft noise reduction. In these applications, the medium is inherently heterogeneous and in this case the convected Helmholtz equation has space-dependent coefficients. Standard high-order methods, such as the Finite Element Method or the Discontinuous Galerkin (DG) Method, often require a large number of degrees of freedom (DoFs) to accurately resolve wave phenomena, leading to significant computational costs.

By contrast, Trefftz methods [4] use basis functions satisfying locally the governing equation. In heterogeneous media, analytical solutions are usually unavailable. In this case, Quasi-Trefftz methods build basis functions satisfying the governing equation locally up to a prescribed order [1, 2].

In this work, we generalize the QTDG framework for first-order systems, originally developed for the Helmholtz problem [3], to the convected Helmholtz problem. This choice of a first-order system formulation is motivated by its natural compatibility with DG numerical fluxes, providing a general and efficient framework for complex aeroacoustic models.

## 2 Mathematical Model

We consider the variable density $\rho$, the speed of sound $c$, and the Mach vector $\boldsymbol{M} = \boldsymbol{v}/c$. Given the pulsation $\omega$ and the wavenumber $\kappa = \omega/c$, the velocity potential $\varphi$ satisfies the convected Helmholtz equation:

$$\begin{aligned}&\rho(\kappa^2\varphi + \mathrm{i}\kappa\boldsymbol{M}\cdot\nabla\varphi) + \nabla\cdot(\rho\nabla\varphi)\\ &-\nabla\cdot(\rho((\boldsymbol{M}\cdot\nabla\varphi)\boldsymbol{M} - \mathrm{i}\kappa\varphi\boldsymbol{M})) = 0.\end{aligned} \tag{1}$$

We introduce a non-physical auxiliary vector variable $\boldsymbol{\sigma}$ designed to rewrite (1) as a first-order system denoted by $\mathcal{S}(\varphi, \boldsymbol{\sigma}) = \boldsymbol{0}$:

$$\begin{cases} -\dfrac{\mathrm{i}\omega\rho}{c^2}\varphi - \nabla\cdot\left(\dfrac{\rho}{c}\boldsymbol{M}\right)\varphi \\ \qquad +\nabla\cdot\boldsymbol{\sigma} + \nabla\cdot\left(\dfrac{2\rho}{c}\boldsymbol{M}\varphi\right) &= 0,\\ \dfrac{-\mathrm{i}\omega}{\rho}(\mathcal{I} - \boldsymbol{M}\boldsymbol{M}^T)^{-1}\boldsymbol{\sigma} + \nabla\varphi &= 0.\end{cases}$$

Note that this system is physically relevant in the subsonic regime $|\boldsymbol{M}| < 1$. The problem is closed by boundary conditions on two separate parts of the boundary, denoted $\Gamma_{rig}$ and $\Gamma_{abc}$, given by:

$$\begin{cases}\boldsymbol{\sigma}\cdot\boldsymbol{n} = 0 & \text{on } \Gamma_{rig},\\ \boldsymbol{\sigma}\cdot\boldsymbol{n} + \dfrac{\rho(\boldsymbol{v}\cdot\boldsymbol{n} - c)}{c^2}\varphi = g & \text{on } \Gamma_{abc},\end{cases}$$

where $\boldsymbol{n}$ denotes the outward unit normal vector and $g$ is a source.

## 3 Numerical Method

The boundary value problem is approximated using a QTDG method. Specifically, the QT

basis functions are polynomials $(\hat{\varphi}, \hat{\boldsymbol{\sigma}}) \in \mathbb{P}_p \times (\mathbb{P}_{p-1})^d$ that satisfy locally

$$\mathcal{S}(\hat{\varphi}, \hat{\boldsymbol{\sigma}}) = \left(\mathcal{O}(h^{p-1}), (\mathcal{O}(h^{p-2}))^d\right)^T. \quad (2)$$

The corresponding discrete QT spaces are better suited to represent the solution than standard polynomial spaces as they incorporate the local physics of the heterogeneous convected Helmholtz equation within each basis function: the QTDG method requires significantly fewer DoFs than a standard DG method, as shown in Table 1.

Table 1: Number of local DoFs per element in 2D between standard DG and QTDG for a velocity potential $\varphi$ of degree $p$.

| Degree | Standard DG | QTDG |
|---|---|---|
| $p$ | $(p+1)(p+2)/2$ | $2p+1$ |
| 2 | 6 | 5 |
| 3 | 10 | 7 |
| 4 | 15 | 9 |
| 5 | 21 | 11 |
| 6 | 28 | 13 |

## 4 Numerical Results

The QTDG method has been implemented in 2D for the convected Helmholtz equation. To assess the method's performance in a realistic aeroacoustic configuration, we consider the scattering of a plane wave by a rigid disk of radius $R = 1$ immersed in a potential background flow. The flow velocity $\boldsymbol{v}$ varies spatially, as illustrated in Fig. 1, creating a heterogeneous medium for the acoustic waves.

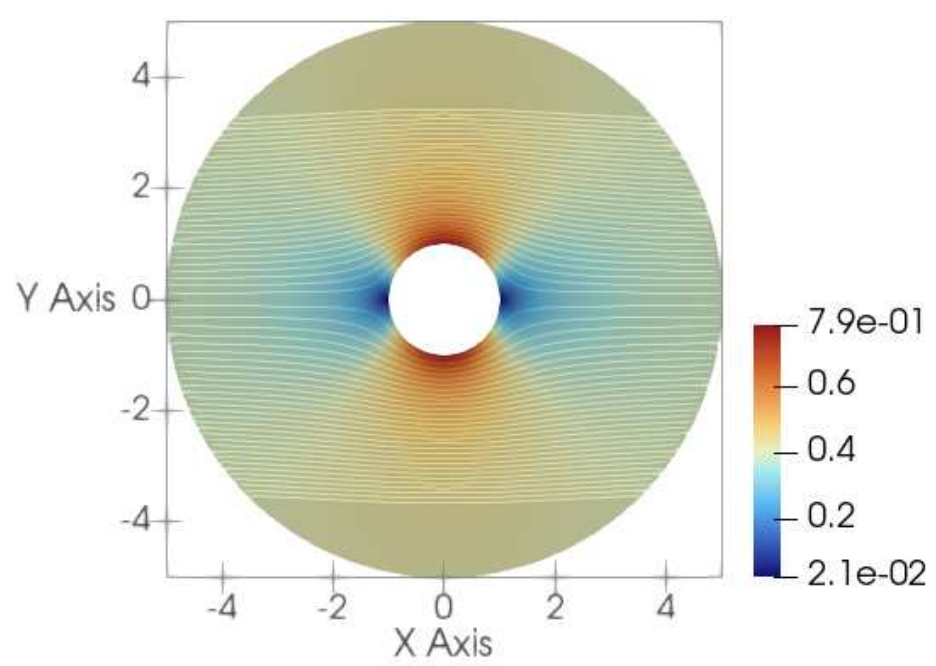

Figure 1: Distribution of the background Mach number $|\mathbf{M}|$ for a potential flow around the rigid disk ($R = 1$). Streamlines indicate the flow direction.

To validate the construction of the basis functions, we first verified that all local basis functions $(\hat{\varphi}, \hat{\boldsymbol{\sigma}})$ satisfy the theoretical property (2) at the desired order $q$.

Fig. 2 represents the QTDG numerical solution of the velocity potential $\varphi$. Numerical convergence is then evaluated by comparing the results to a reference solution computed on a very fine mesh. The results confirm that the QTDG method achieves the expected high-order convergence rates.

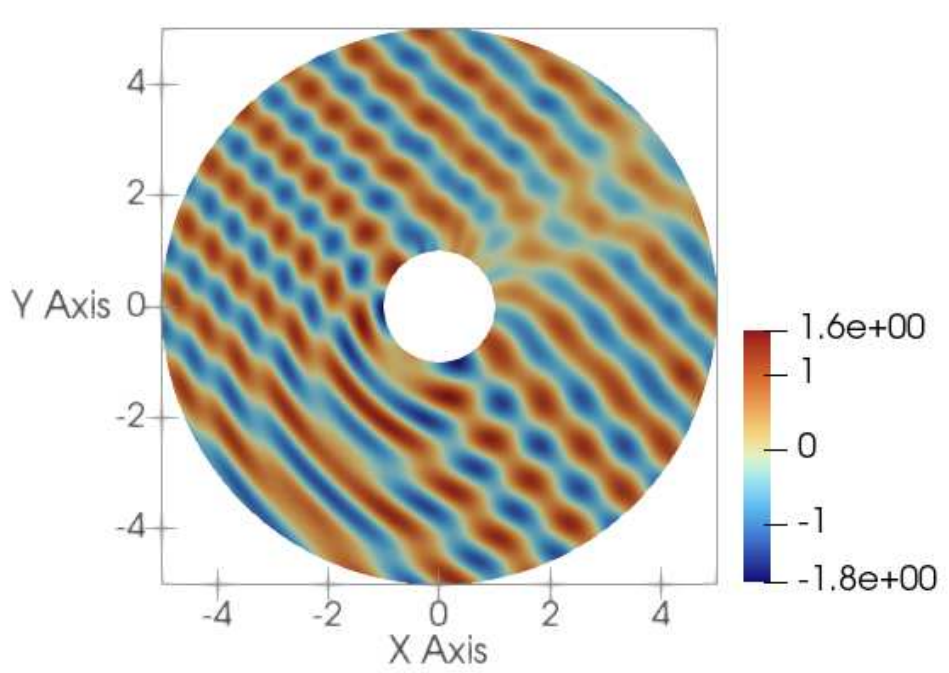

Figure 2: Numerical solution of the velocity potential $\varphi$ obtained with the QTDG method.

# The Hermite-Taylor Discrete Correction Function Method for Surface-Driven Electromagnetic Problems

Yann-Meing Law[1,*]
[1]Department of Mathematics and Industrial Engineering, Polytechnique Montréal, Montréal, CAN
*Email: yann-meing.law@polymtl.ca

**Abstract**

Hermite-Taylor methods are efficient high-order methods for hyperbolic problems. These methods evolve the variables and their spatial derivatives through order $m$ to achieve a $2m+1$ order of accuracy. Unfortunately, the enforcement of boundary conditions is challenging due to the amount of data required at the boundary. In this work, a novel correction function method is proposed to enforce boundary and interface conditions for Hermite-Taylor methods. The key idea is to directly constrain the coefficients of the correction function polynomials to enforce the governing equations, the surface conditions and the correction functions to match the Hermite method's solution. Numerical examples in 2-D are performed and the expected rate of convergence is obtained.

## 1 Introduction

Hermite-Taylor methods are based on two processes: a Hermite interpolation method in space and a Taylor method in time. They achieve a $2m+1$ order of accuracy by evolving the variables and their spatial derivatives through order $m$. There are $(m+1)^d$ data for each variable at each mesh node in $d$ space dimensions. These high-order methods are especially well-suited for linear hyperbolic problems since the stability condition on the time step size depends only on the largest wave speed in the system and is independent of the order [1]. The enforcement of boundary conditions is challenging in this setting because Hermite methods require $(m+1)^d$ data on the boundary, which are usually not available. Many approaches were developed to overcome this challenge, such as hybrid DG-Hermite [2], compatibility boundary conditions [3] and correction function methods. Here we propose a novel discrete correction function method [4]. The key idea is to directly constrain the coefficients of the correction function polynomials to enforce the governing equations. The proposed method avoids any numerical integration, reduces the condition number of the matrices used in the minimization procedure, and can be straightforward extended to 3-D problems.

For simplicity, we are seeking numerical solutions of Maxwell's equations in 1-D

$$\partial_t H - \partial_x E = 0, \quad \partial_t E - \partial_x H = 0, \qquad (1)$$

in a domain $\Omega = [x_{\Gamma_\ell}, x_{\Gamma_r}]$, enclosed within the computational domain $\Omega_c = [x_\ell, x_r]$, and a time interval $I = [t_0, t_f]$. Here $H$ is the magnetic field and $E$ is the electric field. The initial condition is given by $H(x,t_0) = H_0(x)$ and $E(x,t_0) = E_0(x)$. We focus on enforcing the perfect electric conductor (PEC) boundary condition $E = 0$ at $x_{\Gamma_1}$ and $x_{\Gamma_2}$. Extension to other types of boundary and interface conditions can be done.

## 2 Hermite-Taylor Method

Two staggered grids in space and time are required for these methods. The primal grid nodes are defined at $(x_i, t_n)$ with $x_i = x_\ell + ih$, $i = 0,\dots,N$, and $t_n = t_0 + n\Delta t$, $n = 0,\dots,N_t$. Here $h = (x_r - x_\ell)/N$ and $\Delta t = (t_f - t_0)/N_t$ for which $N$ and $N_t$ are respectively the number of cells on the primal grid and the number of time steps. The second grid, called the dual grid, has its nodes defined at the cell centers of the primal grid, that is, $(x_{i+1/2}, t_{n+1/2}) = (x_i + h/2, t_n + \Delta t/2)$ for $i = 0,\dots,N-1$ and $n = 0,\dots,N_t - 1$.

Assuming that the data required for Hermite-Taylor methods are available on the primal grid at the initial time. For each cell and for each variable, we compute the Hermite interpolant of degree $2m+1$, centered at $x_{i+1/2}$. Expanding the polynomial coefficients in time and enforcing Maxwell's equations allow us to compute a space-time polynomial and evolve the data on the dual grid at $(x_{i+1/2}, t_{1/2})$. The procedure is repeated from the dual grid to the primal grid to complete a time step.

## 3 Discrete Correction Function Method

Focusing on the left boundary at $x_{\Gamma_\ell}$, assume a node $(x^*, t^*)$ in its vicinity where the Hermite-Taylor method cannot be applied. In this situation, we provide the data at $(x^*, t^*)$ using the discrete correction function method. Defining a small space-time domain $\Omega_c \times I_c$, called a local patch, that includes the boundary, $(x^*, t^*)$ and some nodes where the Hermite method's solution is available. We are seeking polynomials, called correction functions, approximating the electromagnetic fields of the form

$$P_d^f(x,t) = \sum_{k=0}^{d}\sum_{s=0}^{d} c_{k,s}^f \left(\tfrac{x-x_c}{h}\right)^k \left(\tfrac{t-t^*}{\Delta t}\right)^s.$$

Here $f$ is either $H$ or $E$, and $x_c$ is the spatial center of the local patch.

Considering Faraday's law, substituting the electromagnetic fields by their correction functions and matching the coefficients leads to

$$\begin{aligned}
&\tfrac{L_0(s+1)}{\Delta t} c_{k,s+1}^H = \tfrac{L_0(k+1)}{h} c_{k+1,s}^E, \quad k,s = 0,\dots,d-1,\\
&\tfrac{L_0(s+1)}{\Delta t} c_{d,s+1}^H = 0, \quad s = 0,\dots,d-1,\\
&\tfrac{L_0(k+1)}{h} c_{k+1,d}^E = 0, \quad k = 0,\dots,d-1.
\end{aligned}$$

Here $L_0$ is a scaling factor to ensure that all the terms share the same dimensional unit. The procedure is repeated for Maxwell-Ampère's law and lead to a linear system $\mathcal{G}\mathbf{c} = 0$.

The boundary condition is enforced by constraining

$$\left(\tfrac{L_0}{c_0}\right)^j \partial_t^j P_d^E(x_{\Gamma_1}, t_i) = 0, \tag{2}$$

$i = 1,\dots,N_s$ and $j = 0,\dots,d$. Here $c_0$ is the reference wave speed, $t_i \in I_c$ and the time derivatives are converted into spatial derivatives using Maxwell's equations. This leads to a linear system $\mathcal{S}\mathbf{c} = 0$.

We match the correction functions with the Hermite method's solution pointwise within the local patch, leading to $\mathcal{B}\mathbf{c} = \mathbf{b}_{\mathcal{B}}$. We therefore obtain an overdetermined system $M\mathbf{c} = \mathbf{b}$, where $M = [\mathcal{G}; \mathcal{B}; \mathcal{S}]$ and $\mathbf{b} = [0; \mathbf{b}_{\mathcal{B}}, 0]$, which is solved in the least square sense using a pivoted QR factorization. Finally, we compute the $2(m+1)$ data required by the Hermite method at $(x^*, t^*)$ using the correction functions of $H$ and $E$.

## 4 Numerical Examples

Consider Maxwell's equations in 2-D using the TM$_z$ mode. Convergence studies show a $2m+1$ rate of convergence for $m = 1-3$, leading up to a seventh-order method, using a standing-mode solution and a smooth tri-star boundary. As a second example, we consider a smooth tri-star PEC illuminated by an incident wave and perform self-convergence studies. The domain is $\Omega = [-1,1]^2$ and $I = [0,1]$. We set $\Delta t = h/2$ and the reference solution is computed using $m = 3$ and $h = 1/1280$. Fig. 1 and Fig. 2 show the reference solution for $E$ and the self-convergence plot for $m = 1-3$. The expected $2m+1$ rate of convergence is obtained.

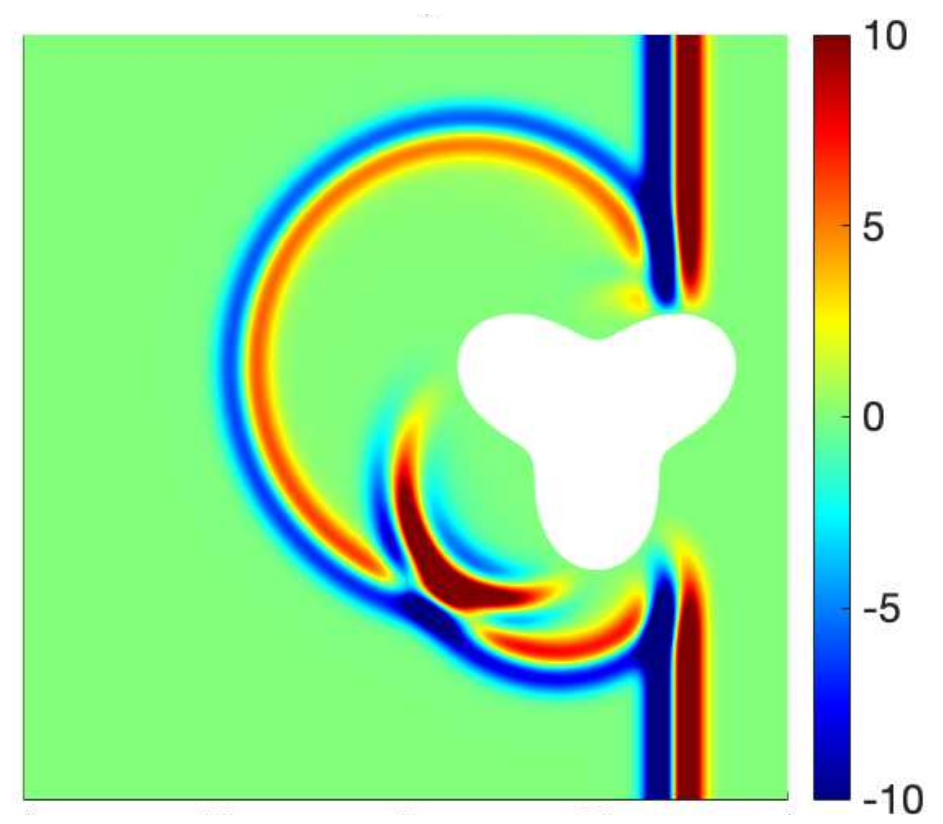

Figure 1: The reference solution for $E$ at $t_f$.

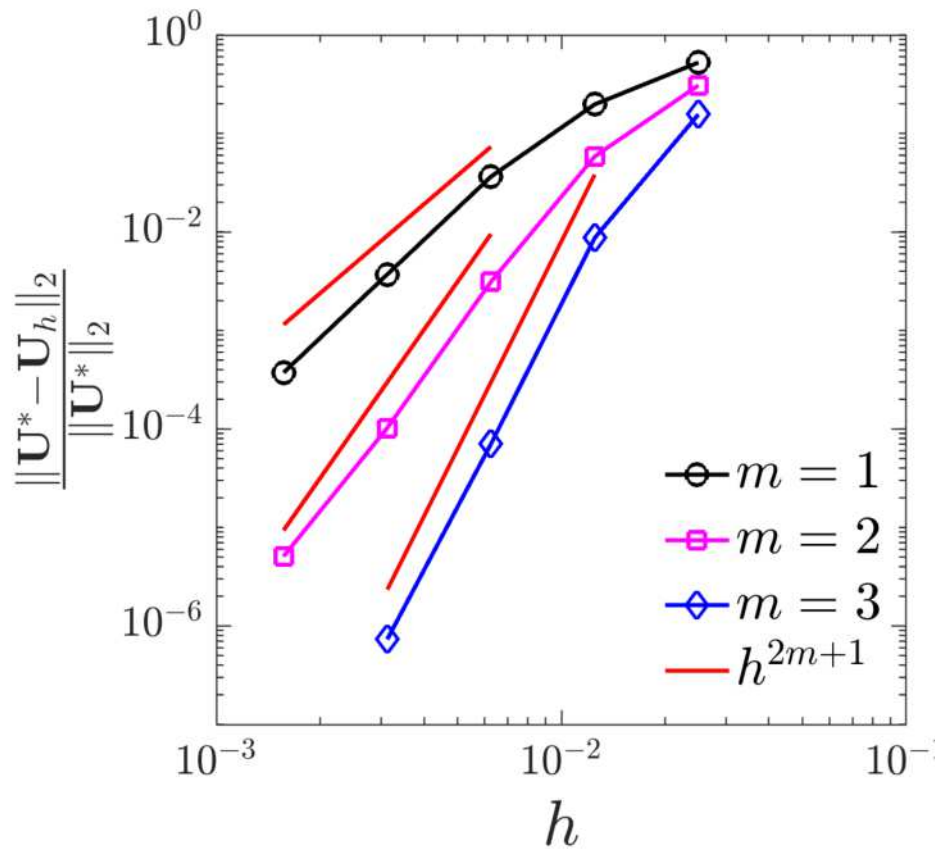

Figure 2: Self-convergence plot ($\mathbf{U} = [H, E]^T$).

# A Helmholtz Solver for Open Domains with Complex Geometry using WaveHoltz and Overset Grids

**Jeffrey Banks**[1], **William Henshaw,**[1] **Ngan Le** [1,*], **Donald Schwendeman**[1]
[1]Department of Mathematical Sciences, Rensselaer Polytechnic Institute, Troy, NY 12180, USA
*Email: leb2@rpi.edu

## Abstract

We develop an efficient Helmholtz solver for domains with open boundaries where radiation boundary conditions are applied. The schemes are based on the WaveHoltz algorithm, which computes solutions of the Helmholtz equation by time-filtering solutions of the wave equation. Complex geometry is treated with overset grids. The solution of the wave equation is solved efficiently with implicit time-stepping. Fully discrete analysis based on an eigenvector expansion provides insight into the convergence rate of the WaveHoltz iteration.

## 1 Introduction

Designing an effective Helmholtz solver is a challenging task for a variety of reasons [3]. One of the challenges is the highly indefinite character of the discretized linear system, causing traditional iterative algorithms, such as preconditioned GMRES or multigrid methods, to converge slowly or fail to converge. Another challenge is the resolution requirements to manage pollution errors at high frequencies. A wide range of solvers have been developed to address these problems, and good overviews of these methods are included in [4]. The method presented here is based on the WaveHoltz [1], which avoids inverting an indefinite matrix. The basic WaveHoltz fixed-point iteration is accelerated using a matrix-free GMRES method. The time-domain solver is adjusted to remove dispersion errors in time, enabling the use of large time-steps without degrading the accuracy. Pollution errors are addressed through the use of high-order accurate spatial approximations. Numerical results are presented demonstrating the algorithm for both simple and complex geometries.

## 2 Model

We are interested in finding real-valued time-periodic solutions to the time-domain problem.

$$\partial_t^2 w = \mathcal{L}(\mathbf{x})w - \gamma(\mathbf{x})\partial_t w - f(\mathbf{x})\cos(\omega t). \tag{1}$$

Where $\mathcal{L}(\mathbf{x})$ is a variable-coefficient elliptic operator, e.g. $\mathcal{L}(\mathbf{x}) = \nabla(c^2(\mathbf{x})\nabla)$, $\gamma(\mathbf{x})$ is a damping coefficient. For numerical simulation, we truncate the unbounded domain and restrict to $\Omega$, where $\Omega$ is a bounded domain in $\mathbb{R}^d$, $d = 1, 2, 3$, and introduce radiation boundary conditions (RBCs). Second-order Engquist-Majda RBCs are used throughout. Figure 1 shows an example of the domain $\Omega$.
The (real-valued) time-periodic solution to (1) is $w(\mathbf{x}, t) = \mathrm{Re}\left(u(\mathbf{x})e^{-i\omega t}\right) = u_r(\mathbf{x})\cos(\omega t) + u_i(x)\sin(\omega t)$, where $u = u_r + iu_i$ solves the Helmholtz problem

$$\mathcal{L}(\mathbf{x})u + \omega^2 u + i\omega\gamma(\mathbf{x})u = f(\mathbf{x}). \tag{2}$$

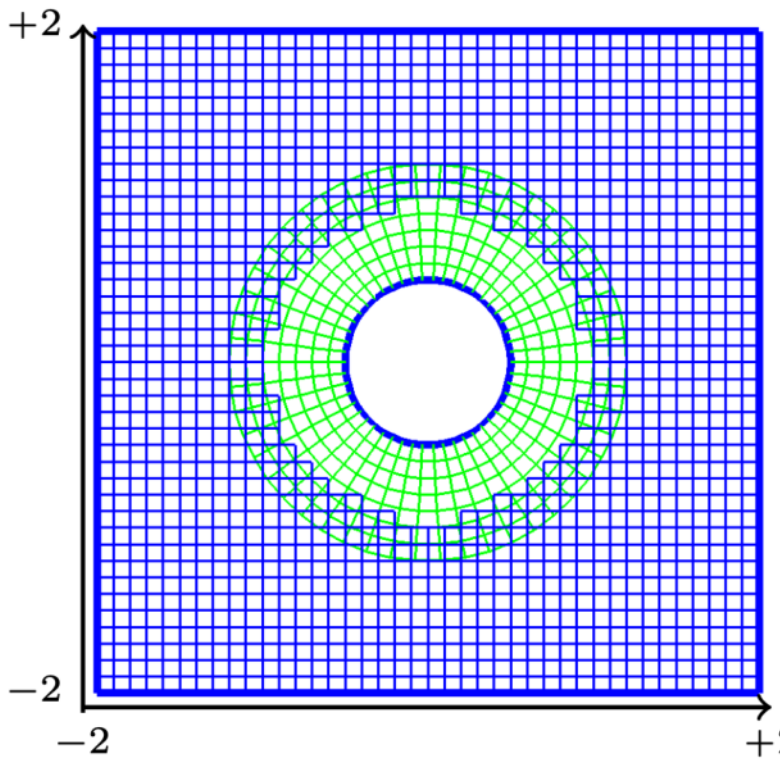

Figure 1: Overset grid with Dirichlet boundary condition imposed on the cylinder (blue circle) and RBCs on all boundaries of the square.

## 3 Methods and Results

The WaveHoltz is used to solve the Helmholtz problem in (2). The iteration first solves the wave equation with two initial conditions over time $0 < t \leq T$, where $T = 2\pi N_p/\omega$, $N_p$ is the number of periods.

$$\partial_t^2 w = \mathcal{L}(\mathbf{x})w - \gamma(\mathbf{x})\partial_t w - f(\mathbf{x})\cos(\omega t), \tag{3a}$$
$$w(\mathbf{x}, 0) = v^{(k)}(\mathbf{x}), \quad \partial_t w(\mathbf{x}, 0) = \dot{v}^{(k)}(\mathbf{x}), \tag{3b}$$
$$\mathcal{B}[\partial_t, \nabla]w = 0 \tag{3c}$$

Then filter the solution in time to get the next iterates

$$v^{(k+1)} = \frac{2}{T}\int_0^T \left(\cos(\omega t) - \frac{1}{4}\right) w(\mathbf{x}, t)\, dt, \tag{4a}$$
$$\dot{v}^{(k+1)} = \frac{2}{T}\int_0^T \left(\cos(\omega t) - \frac{1}{4}\right) \partial_t w(\mathbf{x}, t)\, dt. \tag{4b}$$

The wave equation (3) is solved efficiently using an implicit scheme, with only ten time-steps per period needed ($N_t = 10$)

$$D_{+t}D_{-t}W_{\mathbf{j}}^n = \frac{L_{ph}}{2}(W_{\mathbf{j}}^{n+1} + W_{\mathbf{j}}^{n-1}) - \tilde{\gamma}_{\mathbf{j}} D_{0t} W_{\mathbf{j}}^n - f_{\mathbf{j}}\cos(\tilde{\omega} t^n)\cos(\tilde{\omega}\Delta t), \tag{5a}$$
$$W_{\mathbf{j}}^0 = V_{\mathbf{j}}^{(k)}, \quad D_{0t}W_{\mathbf{j}}^0 = \dot{V}_{\mathbf{j}}^{(k)}, \tag{5b}$$

where $D_{+t}$, $D_{-t}$, and $D_{0t}$ are the forward, backward, and central divided difference operators respectively; $L_{ph}$ is the $p^{\text{th}}$−order accurate spatial approximation. The time filters are approximated by the composite trapezoidal rule

$$V_{\mathbf{j}}^{(k+1)} = \frac{2}{T}\sum_{n=0}^{N_t} \sigma_n(\cos(\tilde{\omega} t^n) - \frac{\alpha_d}{2})W_{\mathbf{j}}^n, \tag{6a}$$
$$\dot{V}_{\mathbf{j}}^{(k+1)} = \frac{2}{T}\sum_{n=0}^{N_t} \sigma_n(\cos(\tilde{\omega} t^n) - \frac{\alpha_d}{2})D_{0t}W_{\mathbf{j}}^n. \tag{6b}$$

Here $\tilde{\omega}$, $\tilde{\gamma}_{\mathbf{j}}$, and $\alpha_d$ are corrections due to time discretization errors (see [2]) so that

$$V_{\mathbf{j}}^{(k)} \to U_{r,\mathbf{j}}, \quad \dot{V}_{\mathbf{j}}^{(k)} \to \frac{\sin(\tilde{\omega}\Delta t)}{\Delta t} U_{i,\mathbf{j}}, \tag{7}$$

where $U_{\mathbf{j}} = U_{r,\mathbf{j}} + U_{i,\mathbf{j}}$ solves the discretized problem of (2)

$$L_{ph}U_{\mathbf{j}} + \omega^2 U_{\mathbf{j}} + i\omega\gamma_{\mathbf{j}} U_{\mathbf{j}} = f_{\mathbf{j}}. \tag{8}$$

The convergence of the fully-discrete scheme is demonstrated for scattering problem from an incident plane wave

$$w^I(\mathbf{x}, t) = e^{i(\mathbf{k}\cdot\mathbf{x} - \omega t)}, \tag{9}$$

where $\mathbf{k} = [k_x, k_y]$, and $\omega = c|\mathbf{k}|$. We solve for the scattered field $w^S = w - w^I$, with the homogeneous Dirichlet boundary condition replaced by $w^S(\mathbf{x}, t) = -w^I(\mathbf{x}, t)$, while the RBCs remain unchanged. Figure 2 shows the total field $|w|$ and convergence rates (theoretical, fixed-point, and GMRES) for the cylinder and split-annulus resonator. A fourth-order accurate approximation in space is used. Results show good convergence for GMRES, demonstrating the potential of WaveHoltz for complex Helmholtz problems.

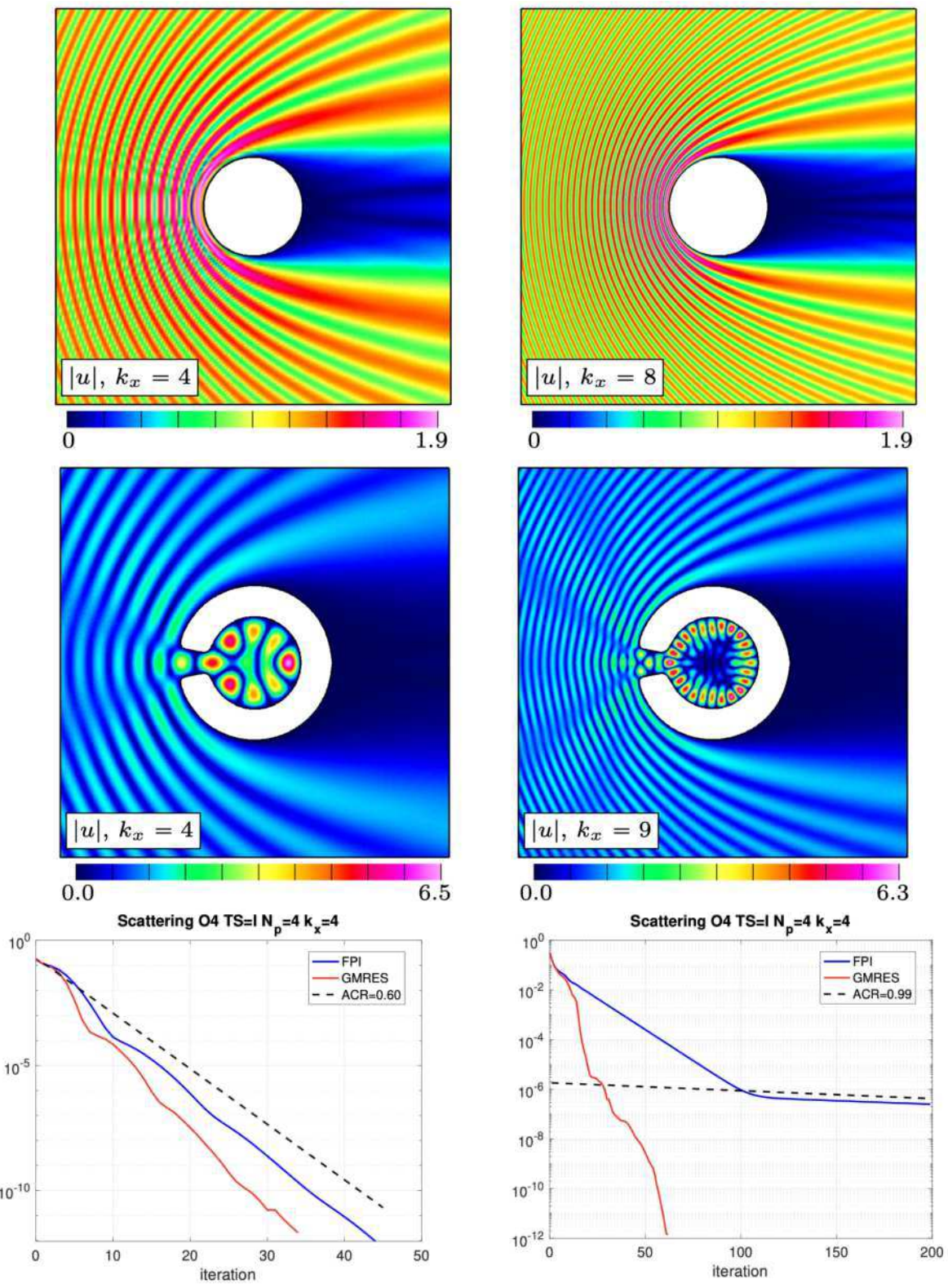

Figure 2: WaveHoltz solutions from scattering problems for a cylinder (top) and split annulus (middle). The bottom row shows the convergence rates ($k_x = 4$) for the cylinder (left) and split annulus (right).

# Model reduction for a wave propagation problem in a slotted structure

Augustin Leclerc[1,*], Pierre Benjamin[1,2], Luc Giraud[1], Sofiane Haddad[1,2], Paul Mycek[1], Guillaume Sylvand[1,2], Isabelle Terrasse[2]
[1]CONCACE Team, joint work with Inria, Airbus and Cerfacs, France.
[2]Airbus Group Innovations, France.
*Email: augustin.leclerc@inria.fr

**Abstract:** In this work, we present a method to learn a reduced model for a wave propagation problem in a slotted structure, using the Boundary Element Method (BEM).

## 1 Introduction

Numerical simulation of wave propagation in multi-scale structures is often compromised by field singularities occurring near fine geometric features like slots [1]. These singularities require an extremely refined local mesh, leading to huge linear systems that may become computationally intractable.

To alleviate this computational cost, model order reduction techniques offer a framework to project complex systems onto low-dimensional subspaces. Substructuring approaches especially allow the condensation of particular components. For example, it is shown in [2] that truncated Singular Value Decomposition (SVD) can effectively reduce coupling operators, enabling the efficient handling of large-scale models without sacrificing accuracy. This concept shares connections with Proper Orthogonal Decomposition techniques, which similarly project models onto low-dimensional subspaces.

In this work, we propose a reduced model based on a data-driven approximation of the Schur complement of the singular sub-domain (with datas from lightweight numerical simulations). By replacing the fine mesh discretization of the slot with a reduced operator, we preserve the physical influence of the singularity while reducing the number of degrees of freedom.

## 2 Model reduction

The training of our reduced model is based on solving a series of linear systems, each parameterized by a mesh geometry $g \in \mathcal{G}_{Train}$ such that

$$\boldsymbol{M_g X_g} = \boldsymbol{Y_g}, \tag{1}$$

where $\boldsymbol{Y_g}$ and $\boldsymbol{X_g}$ are multicolumn vectors (in our case, we make the incidence angle vary to get several solutions), and $\boldsymbol{M_g}$ is the BEM matrix, which is complex symmetric. Moreover, all the meshes of our dataset contain a subset which is invariant with respect to $g$. Hence, we can identify different blocks in the matrices:

$$\begin{gathered} \boldsymbol{M_g} = \begin{pmatrix} \boldsymbol{A_g} & \boldsymbol{B_g} \\ \boldsymbol{B_g}^T & \boldsymbol{D} \end{pmatrix}, \\ \boldsymbol{X_g} = \begin{pmatrix} \boldsymbol{X_{A,g}} \\ \boldsymbol{X_{D,g}} \end{pmatrix}, \quad \boldsymbol{Y_g} = \begin{pmatrix} \boldsymbol{Y_{A,g}} \\ \boldsymbol{Y_D} \end{pmatrix}, \end{gathered} \tag{2}$$

where the block $\boldsymbol{D} \in \mathbb{C}^{n_D \times n_D}$ contains the interactions of the degrees of freedom of the invariant subset.

Based on this linear system we intend to define a new one of lower dimension where the number of unknowns in the block $\boldsymbol{D}$ is reduced while the solution $\boldsymbol{X_{A,g}}$ remains almost unchanged. To best represent the influence of the unknowns $\boldsymbol{X_{D,g}}$ on $\boldsymbol{X_{A,g}}$, the governing idea is to construct a new linear system such that the Schur complement $\boldsymbol{A_g} - \boldsymbol{B_g D}^{-1} \boldsymbol{B_g}^T$ is a good approximation of the original one. More specifically, we approximate $\boldsymbol{S_g} = \boldsymbol{B_g D}^{-1} \boldsymbol{B_g}^T$ by its truncated SVD of a given rank $n_{VS}$

$$\widetilde{\boldsymbol{S}}_{\boldsymbol{g}} = \boldsymbol{U}(:, 1:n_{VS})\, \boldsymbol{\Sigma_{n_{VS}}} \boldsymbol{V}^*(1:n_{VS}, :), \tag{3}$$

such that $\boldsymbol{S_g} = \boldsymbol{U\Sigma V}^*$. It leads to the computation of the following system

$$\begin{pmatrix} \boldsymbol{A_g} & \widetilde{\boldsymbol{B}}_{\boldsymbol{g}} \\ \widetilde{\boldsymbol{B}}_{\boldsymbol{g}}^T & \widetilde{\boldsymbol{D}}_{\boldsymbol{g}} \end{pmatrix} \begin{pmatrix} \boldsymbol{X^{\Sigma}_{A,g}} \\ \widetilde{\boldsymbol{X}}^{\boldsymbol{\Sigma}}_{\boldsymbol{D,g}} \end{pmatrix} = \begin{pmatrix} \boldsymbol{Y_{A,g}} \\ \widetilde{\boldsymbol{Y}}_{\boldsymbol{D,g}} \end{pmatrix}, \tag{4}$$

where $\widetilde{\boldsymbol{B}}_{\boldsymbol{g}} = \boldsymbol{U}(:, 1:n_{VS})$, $\widetilde{\boldsymbol{D}}_{\boldsymbol{g}} = \left( \widetilde{\boldsymbol{B}}_{\boldsymbol{g}}^* \boldsymbol{S_g} \overline{\widetilde{\boldsymbol{B}}_{\boldsymbol{g}}} \right)^{-1}$ and $\widetilde{\boldsymbol{Y}}_{\boldsymbol{D,g}} = \widetilde{\boldsymbol{D}}_{\boldsymbol{g}} \widetilde{\boldsymbol{B}}_{\boldsymbol{g}}^* \boldsymbol{B_g D}^{-1} \boldsymbol{Y_D}$. It can be shown that the Schur complement of the approximate matrix in (4) with respect to $\widetilde{\boldsymbol{D}}_{\boldsymbol{g}}$ is $\boldsymbol{A_g} - \widetilde{\boldsymbol{S}}_{\boldsymbol{g}}$.

In order to generalize the model, we look for $n_{VS}$ vectors that generate the subspace in which $\widetilde{\boldsymbol{X}}^{\boldsymbol{\Sigma}}_{\boldsymbol{D}}$ is the projection of $\boldsymbol{X_D}$ (respectively the column concatenations of $\{\widetilde{\boldsymbol{X}}^{\boldsymbol{\Sigma}}_{\boldsymbol{D,g}}\}_{g \in \mathcal{G}_{Train}}$ and $\{\boldsymbol{X_{D,g}}\}_{g \in \mathcal{G}_{Train}}$) with respect to the Frobenius norm $\| \cdot \|_F$. This is equivalent to looking for

$$\widehat{\boldsymbol{Q}} = \underset{\boldsymbol{Q} \,\in\, \mathbb{C}^{n_D \times n_{VS}}}{\operatorname{argmin}} \|\boldsymbol{Q} \widetilde{\boldsymbol{X}}^{\boldsymbol{\Sigma}}_{\boldsymbol{D}} - \boldsymbol{X_D}\|_F. \tag{5}$$

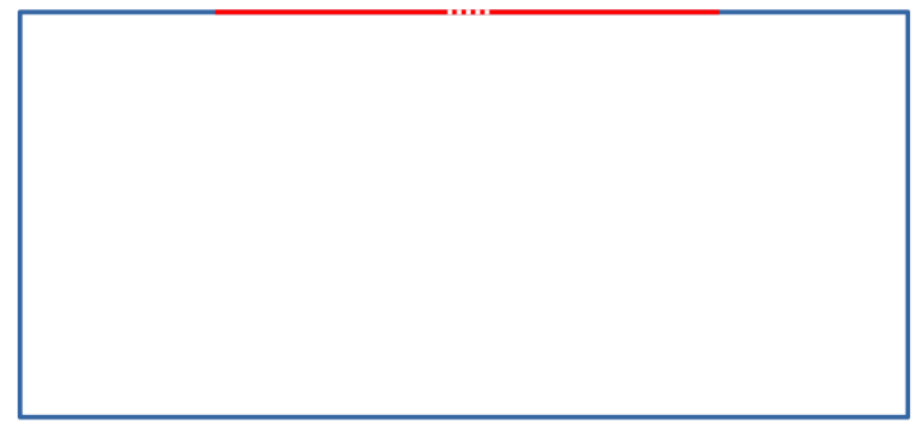

Figure 1: Cross-section of the structure. The red set (whose dotted part is a slot) is the one we aim to reduce.

Hence, considering that $\widehat{Q}\widetilde{X}^{\Sigma}_{D} \approx X_D$, we construct another reduced model for the mesh $g'$ picked from a train or test dataset:

$$\begin{pmatrix} A_{g'} & B_{g'}\widehat{Q} \\ \left(B_{g'}\widehat{Q}\right)^T & \widehat{Q}^T D\widehat{Q} \end{pmatrix} \begin{pmatrix} X^{Q}_{A,g'} \\ \widetilde{X}^{Q}_{D,g'} \end{pmatrix} = \begin{pmatrix} Y_{A,g'} \\ \widehat{Q}^T Y_D \end{pmatrix}. \quad (6)$$

## 3 Application on a 2D wave propagation problem

We consider the problem of electromagnetic wave propagation in a medium with an invariant thin metal slotted structure with respect to $z$ (see its cross-section on Figure 1), excited by an electromagnetic wave of the form

$$(\boldsymbol{E}(x,y,z), \boldsymbol{H}(x,y,z)) = e^{i\gamma z}\left(\begin{pmatrix} \boldsymbol{E_T} \\ E_z \end{pmatrix}, \begin{pmatrix} \boldsymbol{H_T} \\ H_z \end{pmatrix}\right)(x,y). \quad (7)$$

We use the integral representation of Maxwell's equations, so that the unknowns are the electric and magnetic current density $\boldsymbol{j} = \boldsymbol{H} \times \boldsymbol{n}$ and $\boldsymbol{m} = \frac{1}{Z_0}\boldsymbol{n} \times \boldsymbol{E}$, with $\boldsymbol{n}$ is the outward unit normal to the structure and $Z_0$ is the impedance of free space. We get this problem by coupling the interior and exterior problems via the degrees of freedom located on the slot. Hence, if we only consider the Transverse Electric polarization, the unknowns are reduced to the tangential field $\boldsymbol{j_T}$ which is confined to the cross-sectional plane of the slot and the structure walls, and $m_z$ (the component of $m$ along the direction of invariance) located exclusively on the slot. The problem is then discretized via the BEM with first-order elements.

However, the slot width can be very small in relation to the wavelength. This results in a singularity in the electric field, which requires the area near the slot to be meshed very finely. The number of unknowns in this area is therefore very large, and we would like to reduce this part of the system by learning an equivalent model. Using the reduction method presented above, $g$ being the mesh of the structure, we will consider that the unknowns in this area form the block $X_{D,g}$, their intern interactions are in the matrix $D$, their interactions with the other degrees of freedom are in $B_g$, and their internal interactions of the other unknowns are in $A_g$. If we consider various structures, the mesh that corresponds to $D$ has to be invariant. In this context, the columns of $\widehat{Q}$ turn out to be the coefficients in the basis of first-order elements of a new basis of functions representing the slot area, and containing the singularity.

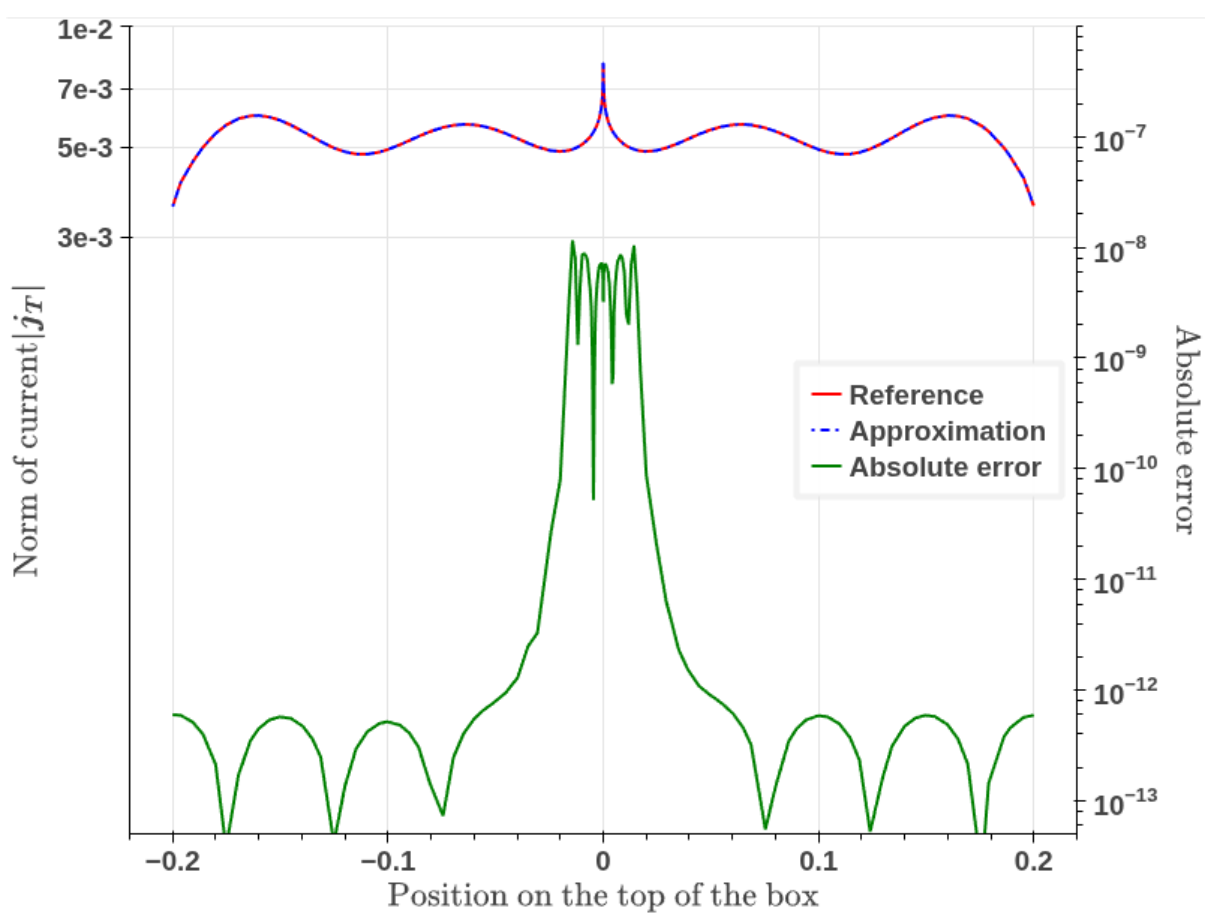

Figure 2: Norm of field $\boldsymbol{j_T}$ onto the top of the box, with a normal incident wave, for a test case.

We train our model on a set of rectangular boxes with centered slots but asymmetrical meshes (see in Figure 2 the result for a test case). This allows us to test the model on more general cases, such as those with off-center slots, those with a different overall structure geometry, or even those with multiple slots, enabling us to test the interactions of one reduced model on another. Numerical experiments will be presented to illustrate our purpose and to show more general cases.

# A parallel space-time $p$-adaptive discontinuous Galerkin method for nonlinear acoustics

Pascal Lehner[1,*], Christian Wieners[2], Danielle Corallo[2]
[1]Department of Mathematics, Klagenfurt, Austria
[2]KIT, Karlsruhe, Germany
*Email: pascal.lehner@aau.at

## Abstract

In this paper, we introduce and analyze a space-time $p$-adaptive discontinuous Galerkin method for nonlinear acoustics. We first present the underlying mathematical model, which is based on a recently derived formulation involving, in particular, only first order in time derivatives. We then propose a spacetime discontinuous Galerkin discretization of this model, combining a symmetric Friedrichs systems discretization for symmetric hyperbolic systems with an interior penalty discretization for damping terms. The resulting nonlinear system is solved using Newton's method. Next, we present a well-posedness analysis of the discrete problem. The analysis begins with a linearized system, for which stability is shown. Using a fixed point argument, these results are extended to the fully discrete nonlinear system, yielding a priori error estimates in a natural discontinuous Galerkin norm. Finally, we present numerical experiments demonstrating the parallel solvability of the space-time formulation and the effectiveness of $p$-adaptivity. The results confirm the theoretical convergence rates and show that adaptive refinement can reduce the number of degrees of freedom required to accurately approximate selected goal functionals. Moreover, the experiments demonstrate that the model reproduces characteristic phenomena of nonlinear acoustics, such as harmonic generation, thereby validating the proposed model.

## 1 Introduction

Efficient numerical methods for nonlinear acoustics are crucial in engineering applications, as well as in medical fields including imaging, tissue ablation, lithotripsy, and drug delivery. Simulating nonlinear acoustic wave propagation is challenging due to the interaction of high frequency waves and nonlinear effects, which can produce steep gradients or shock like structures. Classical models of nonlinear acoustics are Westervelt's and Kuznetsov's equations, which involve second-order time derivatives acting on nonlinear terms. Various finite element and discontinuous Galerkin (DG) methods with error analyses have been developed for these models. In this work, the authors propose a novel space-time DG method for a nonlinear acoustic model related to Kuznetsov's equation that uses only first-order time derivatives and spatial nonlinearities.

## 2 Model

Let $\Omega \subset \mathbb{R}^d$ be an open, bounded space domain with dimension $d \in \{2, 3\}$, and let $\Gamma := \partial\Omega$ denote the boundary of $\Omega$, which is assumed to be Lipschitz. Furthermore, let $I := (0, T) \subset \mathbb{R}$ denote a time interval with final time $T > 0$, and define the corresponding space-time domain $Q := I \times \Omega$. The model is governed by parameters satisfying

$$\alpha,\, \beta,\, \gamma,\, \delta,\, \varepsilon,\, \zeta,\, \eta,\, \theta > 0,\; \gamma \neq \delta,$$

and takes the form

$$\begin{aligned} \alpha\partial_t p + \nabla\cdot v - \beta\Delta p + \gamma p\nabla\cdot v + \delta\nabla p\cdot v &= p_S \\ \varepsilon\partial_t v + \nabla p - \zeta\Delta v + \frac{1}{2}\nabla\big(\eta v^2 - \theta p^2\big) &= v_S \\ p = p_0, \quad v = v_0 \quad &\text{on} \quad \{0\}\times\Omega \\ p = p_D, \quad v = v_D \quad &\text{on} \quad I\times\Gamma_{\mathrm{D}} \\ \beta\nabla p\cdot n = p_N, \quad v = v_N \quad &\text{on} \quad I\times\Gamma_{\mathrm{N}} \end{aligned}$$

with $\Gamma_{\mathrm{D}} \cup \Gamma_{\mathrm{N}} = \Gamma$, where $\Gamma_{\mathrm{D}} \cap \Gamma_{\mathrm{N}} = \emptyset$, and $\Gamma_{\mathrm{D}}$ has positive $(d-1)$-dimensional Lebesgue measure. In this model of nonlinear acoustics $p : Q \to \mathbb{R}$ can be understood as the fluctuation of pressure and $v : Q \to \mathbb{R}^d$ as the fluctuation of the velocity field, respectively, cf. [2, Section 2] for details. Well posedness of the model has been shown in [2].

## 3 Space-time method

The proposed DG method for a suitable discontinuous space-time polynomial space $Z_h$ is:

Solve for $u_h \in Z_h$ such that $b_h(u_h, z_h) = \ell_h(z_h)$

for all $z_h \in Z_h$ with

$$b_h(u_h, z_h) := m_h(u_h, z_h) + a_h(u_h, z_h) \\ +d_h(u_h, z_h) + n_h(u_h, z_h)$$

and a suitable right hand side $\ell_h$ reflecting the boundary conditions. Here $m_h$ stands for a space-time discretiation of the time derivatives, $a_h$ for the of the linear wave operator $(\nabla, \nabla\cdot)$ as done in space-time methods for symetric Friedrich systems [3]. $d_h$ stands for the discretization of the Laplacian terms with an interior penalty method, and $n_h$ for the nonlinear terms.

The analysis is carried out with respect to the $h$-dependent norm

$$\|u_h\|_{Z_h} := \left(\frac{1}{2}\sum_{j=0}^{k} \|M^{\frac{1}{2}}[u_h]_j\|^2_{L^2(\Omega_h)} \right. \\ \left. + \int_0^T c_A \|A_n[u_h]\|^2_{L^2(\partial\Omega_h)} + c_\Theta \|u_h\|^2_h \, \mathrm{d}t \right)^{\frac{1}{2}},$$

where $\Omega_h$ is the discrete space domain, $c_A, c_\theta > 0$, $M, A_n$ are suitable matrices, $[\cdot]$ denotes a jump in space, $[\cdot]_j$ denotes a jump in time at timestep $j$, $k$ is the number of time steps, and $\|\cdot\|_h$ is the usual norm of the broken Sobolev space $H^1(\Omega_h)$. After showing well-posedness of a linearized system by means of standard methods, using a generalized Newton-Kantorvich argument, one can show the following well-posedness and convergence result.

**Theorem 1** *Let $u = (p, v) \in U$ be a regular enough solution to the continuous system and assume that $\|u\|_U, \|\ell_h\|_{Z_h'}$ are sufficiently small. Let $\chi \in (0,1)$ and $\xi > 0$ be determined by space-time mesh and approximation space properties. For $\bar{h} > 0$ define*

$$D_{\bar{h},*} := \{u_h \in Z_h : \|u_h - u\|_{Z_h} < C_* \bar{h}^{\xi}\}$$

*with $C_* > 0$ independent of $\bar{h}$. Then, there exists a $\tilde{h} < \bar{h}$ such that for all $h < \tilde{h}$ the discrete system has a unique solution $u_h \in Z_h$ and $(1+\chi)$-order of convergence in Newton's method with initial guess $u_{h,*} \in D_{h,*}$ holds. In particular, we have the following error estimate*

$$\|u_h - u\|_{Z_h} \le C_* h^{\xi}.$$

## 4 Numerical Experiments

For the numerical experiments, we aim to model sound waves propagating through soft tissue. To this end, we choose $\Omega = (-\frac{1}{2}, \frac{1}{2})^2$, $I = (0, 1)$,

$$\alpha = 1,\ \beta = 10^{-4},\ \gamma = 6,\ \delta = 1, \\ \varepsilon = 1,\ \zeta = 10^{-4},\ \eta = 1,\ \theta = 1$$

and consider a harmonically time varying source. The results are compared with those obtained from the linear wave equation, revealing the expected local sound speed shift characteristic of nonlinear acoustics, see below.

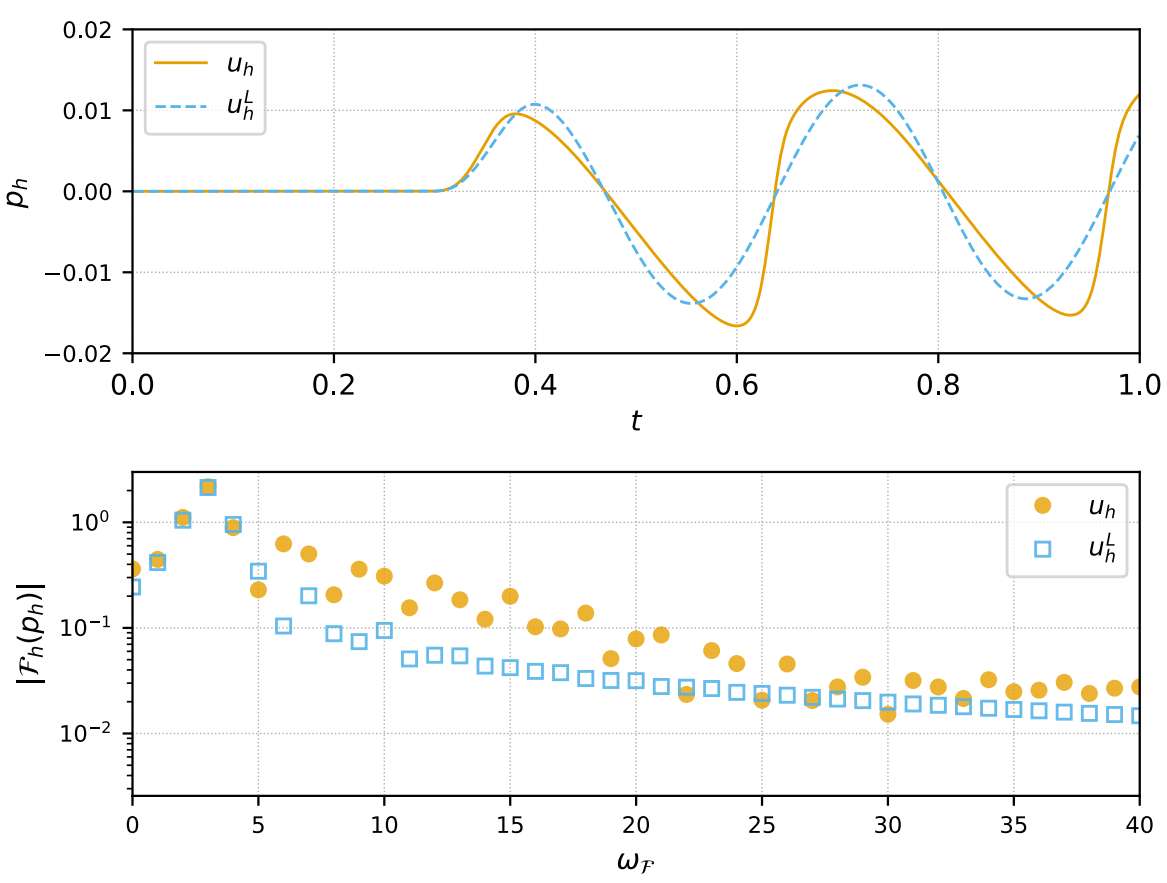

Figure 1: At the top, the time evolution of the pressure at $(x, y) = (0.246094, 0.246094)$ is shown. The blue dotted line corresponds to the pressure component obtained with $\gamma = \delta = \eta = \theta = 0$, while the orange solid line represents the pressure from the full nonlinear model. At the bottom, the corresponding magnitudes of the discrete Fourier transforms of $p_h$ and $p_h^L$ are displayed.

# Asymptotic Solutions to the Radial Wave Equation in Kerr Spacetime

**John Levis**[1,*], **Thomas Hangstrom**[1]
[1]Department of Computational and Applied Mathematics, SMU, Dallas, USA
*Email: jlevis@smu.edu

**Abstract**

We explore new analytical approaches for solving the radial component of the scalar wave equation in Kerr spacetime, with a focus on modeling outgoing gravitational radiation from rotating black holes. Beginning with the Teukolsky Master Equation in Eddington–Finkelstein coordinates, we carry out a full separation of variables and derive the resulting radial ordinary differential equation, emphasizing its complex, non–self-adjoint structure and the presence of an irregular singularity at infinity. Several complementary solution strategies are explored. First, we obtain far-field asymptotics by reducing the radial equation to a Whittaker form, yielding explicit outgoing-wave solutions with logarithmic phase corrections. The central contribution of the work is a detailed construction of power-series expansions using two different techniques resulting in general recurrence relations for the coefficients. We then explore various characteristics of these expansions including their asymptoticity, convergence and the impact of the eigenvalue to gain a deeper understanding of their validity and utility.

## 1 Introduction

As an outcome of the separation of variables, we derive the resulting radial wave equation (RWE)

$$R_{rr} + \frac{2r-2M+2aim+4Mr\sigma}{\Delta} R_r - \frac{\sigma^2\Gamma+2M\sigma+\lambda}{\Delta} R = 0$$

where

$$\Delta = r^2 - 2Mr + a^2 \text{ and } \Gamma = r^2 + 2Mr + a^2$$

The physical parameters are the black hole mass M, spin a, azimuthal number m, separation constant $\lambda$ and $\sigma$ as the dual Laplace variable. As such, the RWE is non-self-adjoint, contains complex coefficients, and exhibits an irregular singularity at infinity—features that complicate both direct analysis and numerical treatment.

Our objective is to construct asymptotic and semi-analytic approximations that illuminate the structure of the radial solutions and quantify their accuracy. The analysis begins with far-field asymptotics, leading to a power-series expansion with terms satisfying a short recurrence relation. These tools provide insight into the oscillatory behavior, decay rates, and divergence properties of the solutions, ultimately yielding a practical approximation strategy for large-radius behavior.

## 2 Far-Field Asymptotics

To understand the behavior of outgoing waves at large radius, we simplify the radial equation by retaining only the dominant terms as r $\to \infty$. After appropriate substitutions, the resulting approximation yields the outgoing-wave expression

$$R(r) \approx \frac{C}{r} e^{-\sigma(r+2Mlnr)}$$

The logarithmic phase shift is a direct manifestation of the long-range nature of the gravitational potential and corresponds to the tortoise coordinate $r_* = r + 2Mlnr$.

This form is consistent with the Sommerfeld radiation condition and provides a physically meaningful leading-order approximation. Substituting this expression back into the full radial equation shows that the residual error is bounded by 4CM$\omega^2$ confirming its asymptotic validity.

## 3 Power-Series Ansatz and Recurrence Structure

While the far-field approximation captures the leading behavior, many applications require higher-order corrections. To obtain these, we introduce the enhanced ansatz

$$R(r) = A(r)B(r)$$

with $A(r) = e^{(\alpha r)} r^{\beta}$ and $B(r) = \sum_{n=0}^{\infty} b_n r^{-n}$

Matching the highest powers of r results in $\alpha^2 = \sigma^2$ and $\beta = -(1 + 4M\sigma)$. Thus, the out-going solution takes the form

$$R(r) = e^{-\sigma r} r^{-(1+4M\sigma)} \sum_{n=0}^{\infty} b_n r^{-n}$$

We also deduce the $b_0$ is arbitrary and can be derived from some near-field boundary condition.

From here, we take two approaches to deriving the appropriate recurrence relation. With the direct method we posit that for $n \geq 1$ and $b_{-1} = 0 = b_{-2}$

$$b_n = \frac{-1}{F_n}\left(G_n b_{n-1} + H_n b_{n-2} + J_n b_{n-3}\right)$$

where $F_n$, $G_n$, $H_n$, $J_n$ come directly from the RWE.

In the second approach, we expand the damping and potential coefficients of the radial equation as follows

$$F(r) = \textstyle\sum_{n=0}^{\infty} f_n r^{-n} \quad G(r) = \sum_{n=0}^{\infty} g_n r^{-n}$$

and substitute the ansatz into the full RWE. Leveraging first order terms for F(r) and G(r) allows us to deduce a general expression for the recurrence as

$$b_n = \frac{-1}{2n\sigma}\Big[(n+1+4M\sigma)(n+4M\sigma)b_{n-1} - \sum_{j=1}^{n}\Big(\sigma f_{j+1} - g_{j+1} + (n+1+4M\sigma-j)f_j\Big)b_{n-j}\Big]$$

## 4 Asymptoticity, Convergence and Eigenvalues

It is shown that both expansions are indeed asymptotic in r . However, both recurrences demonstrate a factorial nature (at least) leading to a conclusion that they are not convergent for all r. That said, we can determine that there exists N = N(r) where

$$R(r) \approx e^{-\sigma r} r^{-(1+4M\sigma)} \sum_{n=0}^{N-1} b_n r^{-n}$$

is convergent and furthermore that the error function $E^N(r) = \sum_{n=N}^{\infty} b_n r^{-n}$ extracted from

$$R(r) = e^{-\sigma r} r^{-(1+4M\sigma)} \left(\sum_{n=0}^{N-1} b_n r^{-n} + \sum_{n=N}^{\infty} b_n r^{-n}\right)$$

is bounded.

The key observations of the eigenvalues impact on the solutions include that as the ratio $\lambda r^{-2}$ moves from negative to positive, the solutions evolve from locally oscillatory to locally exponentially dominant.

## 5 Conclusions and Future Work

This work provides a unified asymptotic and semi-analytic framework for the scalar radial equation in Kerr spacetime. The far-field reduction yields explicit outgoing-wave behavior with logarithmic phase corrections, while the enhanced exponential–power-series ansatz produces a general recurrence relation that clarifies the divergence structure of the asymptotic expansion. Although the resulting series is non-convergent, optimal truncation yields accurate approximations over large radial domains.

This paper is intended to serve as the foundation for further investigation into the analytical modeling of gravitational waves in a variety of increasingly complicated regimes including solutions that involve the expansion in the frequency variable, analysis of non-linear terms in the metric perturbation, and generalization to the Teukolsky equation for gravitational fields. These developments will support more accurate modeling of gravitational radiation, both by the development of approximate formulas for special cases and for the analysis of radiation conditions employed in numerical studies.

**Acknowledgment:** The work of the second author was supported in part by NSF Grant DMS-2309687. Any opinions, findings, and conclusions or recommendations expressed in this material are those of the authors and do not necessarily reflect the views of the NSF.

# Radiation Boundary Conditions with Rational Function Approximations

Steve Li[1,*], Thomas Hagstrom[1]
[1]Department of Mathematics, Southern Methodist University, Dallas TX, USA
*Email: steveli@smu.edu

**Abstract**

We are concerned with computing the discrete spectrum and associated bound states of semi-infinite Hamiltonians. The applications include modeling electronic edge states localized at boundaries of crystalline materials. The problem is non-trivial since arbitrary "hard truncations" produce spurious bound states. Currently, the improved method involves computing the Green's function of the semi- infinite Hamiltonian by imposing appropriate boundary conditions, where the spectral data is subsequently recovered via Riesz projection [1]. However, we will attempt to apply rational function approximation at the radiation boundary in order to transform it into an eigenvalue problem.

## 1 Introduction

The accurate modeling of wave propagation in unbounded domains remains a fundamental challenge in computational physics and applied mathematics. To overcome this difficulty, radiation boundary conditions (RBCs) are introduced to emulate open-domain behavior and allow outgoing waves to leave the computational domain without reflection. When studying Hamiltonian systems on a lattice model, we are concerned with finding bound states for decaying modes by solving a non-standard eigenvalue problem under RBC. In this project, we explore two options for implementing the rational function approximation that transforms the problem into a standard eigenvalue problem: AAA rational function approximation [2] and contour integral eigen-projection [3].

## 2 Model

The semi-infinite operator $H$ is given by

$$[H\psi]_k = A^*(k-1)\psi_{k-1} + V(k)\psi_k + A(k)\psi_{k+1}, \quad (1)$$

where

$$V(k) := \begin{pmatrix} 0 & t_1^k \\ t_1^k & 0 \end{pmatrix}, \; A(k) := \begin{pmatrix} 0 & 0 \\ t_2^k & 0 \end{pmatrix}, \quad k \in \mathbb{N}^+$$

The potentials at each stage $t_{1,2}^k$ are randomly perturbed values of final stage potentials $t_{1,2}^\infty \in \mathbb{R}$. In our example, we have $t_1^\infty = 1, t_2^\infty = 1.6$, and $t_{1,2}^k = t_{1,2}^\infty + \mathcal{N}(0,1)$. We also define

$$V^\infty := \begin{pmatrix} 0 & t_1^\infty \\ t_1^\infty & 0 \end{pmatrix}, \; A^\infty := \begin{pmatrix} 0 & 0 \\ t_2^\infty & 0 \end{pmatrix}.$$

In block-matrix form we may express the eigenvalue problem as $H\psi = \lambda\psi$, with the truncation condition that

$$t_1^k = t_1^\infty, \; t_2^k = t_2^\infty, \; \psi_{k+1} = \kappa\psi_k, \quad \forall k \geq m-1.$$

The truncation is thus given by

$$\left[\begin{array}{cccc|c} V(1) & A(1) & \dots & \dots & \dots \\ A^*(1) & V(2) & A(2) & \dots & \dots \\ \dots & \dots & \dots & \dots & \dots \\ \dots & \dots & (A^\infty)^* & V^\infty & A^\infty \\ \hline \dots & \dots & \dots & \dots & \dots \end{array}\right] \left[\begin{array}{c} \psi_1 \\ \psi_2 \\ \dots \\ \psi_m \\ \hline \dots \end{array}\right] = \lambda \left[\begin{array}{c} \psi_1 \\ \psi_2 \\ \dots \\ \psi_m \\ \hline \dots \end{array}\right], \quad (2)$$

where the last line before the truncation is

$$(A^\infty)^*\psi_{m-1} + V^\infty\psi_m + A^\infty\psi_{m+1} = \lambda\psi_m. \quad (3)$$

To enforce such condition for truncation, we need

$$\det((A^\infty)^* + \kappa(V^\infty - \lambda I) + \kappa^2 A^\infty) = 0,$$

which is equivalent to

$$\kappa^2 t_1^\infty t_2^\infty + \kappa((t_1^\infty)^2 + (t_2^\infty)^2 - \lambda^2) + t_1^\infty t_2^\infty = 0. \quad (4)$$

## 3 Methods and Results

We are currently facing a non-standard eigenvalue problem. Although we have an explicit form for $\kappa(\lambda)$, simply plugging in $\psi_{m+1} = \kappa\psi_m$ back into equation (2) will result in a nonlinear $\lambda$-dependency on the LHS. We propose two separate methods to implement a rational function approximation for the equation (4), namely the AAA algorithm [2] and contour integral eigen-projection [3].

For the AAA algorithm, we have the analytical solution to the quadratic equation (4) for which we may obtain the following rational function approximation:

$$\kappa(\lambda) \approx r_\infty + \sum_{j=1}^{l} \frac{r_j}{\lambda - p_j}, \quad r_\infty, r_j, p_j \in \mathbb{C} \quad (5)$$

where the number of approximation terms $l$ is given by the error tolerance. We are now able to convert the previously nonlinear $\lambda$-dependency in our eigenvalue problem to a linear one, using some algebraic manipulations. The resulting linear system will be re-integrated into the truncation formulation (2) and solved by a function call for EIG in MATLAB.
The AAA algorithm works for the one dimensional case where explicit solutions are available for a single $\kappa$. In multi-dimensional problems, we often encounter multiple values of $\kappa$ which we would like our method to capture. All values of $\kappa$ lie within the unit circle, and the idea is to discretize the contour integral formulation for the eigen-projection

$$\begin{pmatrix} \psi_m \\ \psi_{m+1} \end{pmatrix} = \frac{1}{2\pi i} \oint \begin{pmatrix} \kappa I & -I \\ (A^\infty)^* & (V-\lambda I)+\kappa A^\infty \end{pmatrix}^{-1} \begin{pmatrix} \psi_m \\ \psi_{m+1} \end{pmatrix},$$

which can be discretized using an arbitrary number of terms $N$:

$$\begin{pmatrix} \psi_m \\ \psi_{m+1} \end{pmatrix} \approx \frac{1}{N} \sum_{j=0}^{N-1} \begin{pmatrix} re^{i2\pi/Nj} I & -I \\ (A^\infty)^* & (V-\lambda I)+re^{i2\pi/Nj} A^\infty \end{pmatrix}^{-1} \begin{pmatrix} \psi_m \\ \psi_{m+1} \end{pmatrix} * re^{i2\pi/Nj}. \quad (6)$$

We may define the term inside the summation:

$$\begin{pmatrix} \phi_{j,m} \\ \phi_{j,m+1} \end{pmatrix} = \frac{1}{N} \begin{pmatrix} re^{i2\pi/Nj} I & -I \\ (A^\infty)^* & (V-\lambda I)+re^{i2\pi/Nj} A^\infty \end{pmatrix}^{-1} \begin{pmatrix} \psi_m \\ \psi_{m+1} \end{pmatrix} * re^{i2\pi/Nj}. \quad (7)$$

Noting that $\psi_{m+1} = \sum_{j=0}^{N-1} \phi_{j,m+1}$, we eventually arrive at the standardized form with linear dependency on $\lambda$:

$$(A^\infty)^* \psi_{m-1} + V^\infty \psi_m + A^\infty (\sum_{j=0}^{N-1} \phi_{j,m+1}) = \lambda \psi_m \quad (8)$$

which eliminates the $\psi_{m+1}$ term from the original problem given by equation (3). The eigenvalues can be subsequently computed in MATLAB with a standard eigensolve.
In both scenarios, we are able to obtain the same eigenvalues within our desired range of bound state eigenvalues, as shown in Table (1). The behavior of the computational results match those obtained from the Green's function method i.e. $\lambda = 0$ is always a bound state eigenvalue and all other eigenvalues are symmetrical across 0 [1]. However, we note that one of the zero eigenvalues experiences a relatively large erroneous perturbation. We will examine the cause of this error in future work.

| AAA | Contour Integral |
|---|---|
| -0.3723 | -0.3723 |
| -0.3152 | -0.3152 |
| -0.2827 | -0.2827 |
| -1.0673e-13 | -0.0426 |
| 1.0695e-13 | -1.1894e-15 |
| 0.2827 | 0.2827 |
| 0.3152 | 0.3152 |
| 0.3723 | 0.3723 |

Table 1: Computed eigenvalues from both methods within the discrete spectrum.

In conclusion, both the AAA and contour integration methods are able to transform the truncation problem into a standard eigensolve. The AAA algorithm works for the one-dimensional lattice problem, whereas the contour integration method works in both one-dimensional and multi-dimensional cases. Numerical results from both algorithms match, given the same parameters. Upcoming steps include convergence studies on the contour integration method and criteria for bound state eigenvalue selection.

## 4 Acknowledgment

Research is supported by NSF DMS-2309687. Any opinions, findings, and conclusions or recommendations expressed in this material are those of the author(s) and do not necessarily reflect the views of the National Science Foundation.

# A discontinuous Galerkin approach for the radiative transfer approximation of high-frequency vibrational energy flow on curved surfaces

Robert J. Lockett[1,*], David J. Chappell[1]
[1]School of Science and Technology, Nottingham Trent University, Nottingham, United Kingdom
*Email: rob.lockett2017@my.ntu.ac.uk

## Abstract

The effective modelling of high-frequency wave propagation is important across a range of settings, yet conventional finite and boundary element methods become prohibitively expensive in the high-frequency regime. We present a transport based framework formulated via the stationary Radiative Transfer Equation (RTE) on curved surfaces. The approach builds on high-frequency energy transport modelling ideas, most notably Dynamical Energy Analysis, retaining improved fidelity over Statistical Energy Analysis while relaxing key restrictions. Our method is compatible with finite element method-style unstructured meshes with arbitrary polygonal elements and achieves competitive simulation times for large-scale discretisations. To compute numerical solutions, we discretise the RTE in space using the discontinuous Galerkin method on general polygonal meshes and in direction using a discrete ordinates scheme with directions projected into the plane of the local finite elements. Numerical results are presented to validate the formulation and to demonstrate its effectiveness as an efficient high-frequency wave transport simulation method on curved surfaces.

## 1 Introduction

High-frequency wave phenomena arise routinely in noise and vibration applications, yet their numerical simulation in complex built-up structures remains difficult. A key reason is that shorter wavelengths require finer resolution of both intricate geometry and rapidly oscillatory fields, so the computational cost increases steeply with frequency [1].

In this work we present a transport-based alternative, modelling wave energy propagation via the stationary RTE for a directional energy density. Building on our existing two-dimensional discontinuous Galerkin (DG) formulation [2], we extend the method to curved three-dimensional surfaces by discretising the geometry with two-dimensional surface elements. We describe the resulting DG–discrete ordinates discretisation and provide validation and illustrative results demonstrating the accuracy of the approach.

## 2 Model

Let $\Omega \subset \mathbb{R}^3$ be a curved surface with boundary $\Gamma = \partial\Omega$. We denote the wave energy density at $\boldsymbol{x} \in \Omega$ by $u(\boldsymbol{x}, \theta)$, propagating through $\Omega$ in the direction $\hat{\boldsymbol{s}}(\theta)$. Although $\Omega \subset \mathbb{R}^3$, the transport is intrinsically two-dimensional and so, on each flat surface element we work in local coordinates induced by the element plane. We model $u$ using the stationary scattering free RTE given by:

$$\hat{\boldsymbol{s}}(\theta)\cdot\nabla_\Omega u(\boldsymbol{x},\theta) + \mu(\boldsymbol{x})\, u(\boldsymbol{x},\theta) = f(\boldsymbol{x},\theta), \quad (1)$$

where $\nabla_\Omega$ is the surface (tangential) gradient, $\mu(\boldsymbol{x})$ is the attenuation coefficient, and $f(\boldsymbol{x},\theta)$ represents a prescribed source.

On $\Gamma$ we impose either transparent conditions or specular reflection. For $(\boldsymbol{x},\theta) \in \Gamma_+$, i.e. $\hat{\boldsymbol{s}}(\theta)\cdot\hat{\boldsymbol{n}}(\boldsymbol{x}) > 0$, where $\hat{\boldsymbol{n}}(\boldsymbol{x})$ is the outward normal to $\Gamma$ lying in the tangent plane of $\Omega$, we enforce

$$\big(u(\boldsymbol{x},\theta) - u_0(\boldsymbol{x},\theta)\big)\big|_{\Gamma_+} = u(\boldsymbol{x},\, \hat{\boldsymbol{s}} - 2\,(\hat{\boldsymbol{s}}\cdot\hat{\boldsymbol{n}})\,\hat{\boldsymbol{n}})\,.$$

Here $u_0$ prescribes a boundary source (and we take $u_0 \equiv 0$ when absent). Finally, the fluence (or total energy density) is obtained by integrating over all directions,

$$\Phi(\boldsymbol{x}) = \int_{S^1} u(\boldsymbol{x},\theta)\,\mathrm{d}\hat{\boldsymbol{s}} = \int_{-\pi}^{\pi} u(\boldsymbol{x},\theta)\,\mathrm{d}\theta, \quad (2)$$

where $S^1$ denotes the unit circle of directions in a flat polygonal element.

## 3 Methods

We discretise (1) using a DG method in space coupled with a discrete ordinates approximation in direction, on a polygonal surface discretisation of $\Omega$. Let $\Omega = \bigcup_{j=1}^{N_E} E_j$ be a partition into non-overlapping (planar) surface elements, with element edges $\Gamma_{i,j}$, $i = 1, \dots, \kappa$. On each $E_j$ we approximate $u$ using a Legendre polynomial basis of total degree $k$, where $\tilde{P}_m$ denotes the $m$th Legendre polynomial on $[-1, 1]$. To define local coordinates $(\hat{x}, \hat{y})$, we introduce an orthonormal tangent frame $(\boldsymbol{e}_{x,j}, \boldsymbol{e}_{y,j})$ spanning the plane of $E_j$, and represent the propagation direction as $\hat{\boldsymbol{s}}_j(\theta) = \cos\theta\, \boldsymbol{e}_{x,j} + \sin\theta\, \boldsymbol{e}_{y,j}$, so that $\hat{\boldsymbol{s}}_j(\theta)$ is tangential and constant on $E_j$.
Element and edge integrals are evaluated via reference-element mappings, to allow exact integration on general geometries as described in [3]. For the directional discretisation, we sample $u$ at $N_\theta$ uniformly spaced directions $\{\theta_n\}_{n=1}^{N_\theta}$. For each $\theta_n$ we expand

$$u(\hat{\boldsymbol{x}}, \theta_n) \approx \sum_{\alpha_1=0}^{k} \sum_{\alpha_2=0}^{k-\alpha_1} \hat{u}^n_{j,\alpha_1,\alpha_2}\, \tilde{P}_{\alpha_1}(\hat{x})\, \tilde{P}_{\alpha_2}(\hat{y}),$$

with $\alpha_1 + \alpha_2 \leq k$.
Applying Green's identity on each element yields the standard DG formulation (see [2]). On inter-element edges we use an upwind numerical flux. Writing $\hat{\boldsymbol{n}}_{i,j}$ for the outward normal to $\Gamma_{i,j}$ in the plane of $E_j$, the flux contribution takes the form $\gamma_{i,j}(\theta) = \big(\hat{\boldsymbol{s}}_j(\theta) \cdot \hat{\boldsymbol{n}}_{i,j}\big)\, u^{\mathrm{p}}_{i,j}(\theta)$, where the upwind trace $u^{\mathrm{p}}_{i,j}$ is defined by

$$u^{\mathrm{p}}_{i,j}(\hat{\boldsymbol{x}}, \theta) = \begin{cases} u_j(\hat{\boldsymbol{x}}_i, \theta), & \text{if } \hat{\boldsymbol{s}}_j(\theta) \cdot \hat{\boldsymbol{n}}_{i,j} > 0, \\ u_{j+}(\hat{\boldsymbol{x}}_{i+}, \theta), & \text{otherwise.} \end{cases}$$

Here $E_{j+}$ denotes the neighbouring element across $\Gamma_{i,j} = \Gamma_{i+,j+}$, and $\hat{\boldsymbol{x}}_{i+}$ is the corresponding boundary point on $E_{j+}$. Across $\Gamma_{i,j}$, we map directions between the two tangent planes via a dihedral rotation and perform interpolation between angles in the discrete ordinates set.
This discretisation yields a linear system for the coefficients $\hat{u}^n_{j,\alpha_1,\alpha_2}$. After solving, we compute $\Phi$ by applying a quadrature rule for (2).

## 4 Results

As an initial validation test, we consider a half-cylinder of radius $\pi^{-1}$ parameterised by the arc-length coordinate $s \in [0, 1]$, with perfectly reflective boundary conditions. At the boundary $s = 0$, we prescribe a boundary emission concentrated in direction perpendicular to the boundary edge which is modelled by a Dirac delta. We fix the angular discretisation with $N_\theta = 2$ having ordinates $\theta \in \{0, -\pi\}$ (corresponding to propagation in the $\pm s$ directions). In this effectively one-dimensional setting in $s$, the fluence admits the closed-form solution

$$\Phi_{\mathrm{sol}}(s) = \frac{\cosh\big(\mu(s-1)\big)}{\sinh(\mu)}, \qquad 0 \leq s \leq 1.$$

Table 1 reports the $\ell_2$ errors in the computed $\Phi$ versus $\Phi_{\mathrm{sol}}$ for increasing mesh resolution and polynomial order. Figure 1 illustrates the specific case where $k = 4$ and $N_E = 8192$.

| $N_E$ | $k=0$ | $k=2$ | $k=4$ |
|---|---|---|---|
| 128 | 1.0418e-3 | 1.3651e-3 | 9.2065e-4 |
| 512 | 6.4460e-4 | 4.8462e-4 | 4.0821e-4 |
| 2048 | 3.5413e-4 | 2.0949e-4 | 1.9698e-4 |
| 8192 | 1.8505e-4 | 9.9859e-5 | 9.7698e-5 |

Table 1: $\ell_2$ errors for varying $k$ and $N_E$.

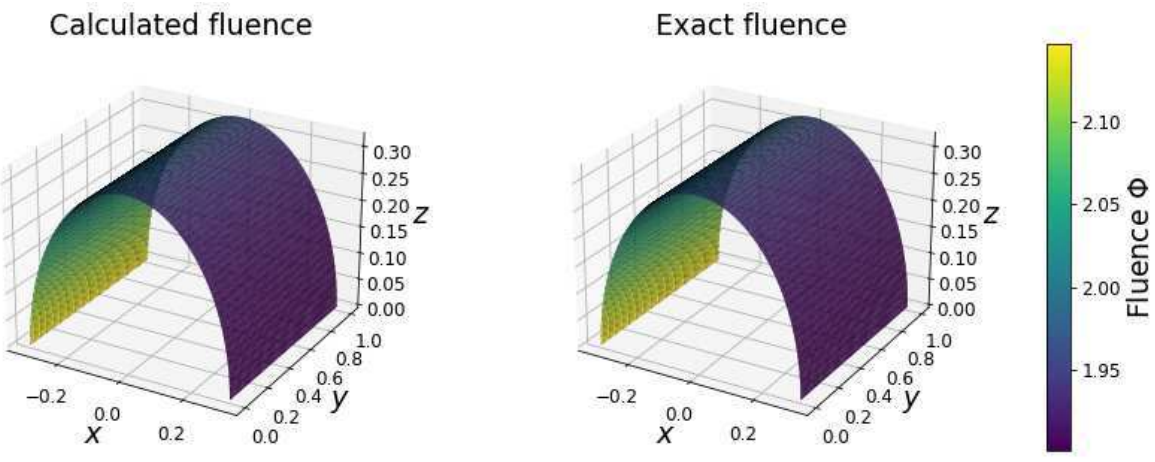

Figure 1: Fluence on half-cylinder versus solution for $k = 4$ and $N_E = 8192$.

Table 1 shows that the error decreases under mesh refinement, but the convergence is first order for any choice of $k$, which is slower than in a strictly planar test. This is primarily due to geometric error from the piecewise planar surface DG discretisation compared to the curved exact solution. Consequently, the most efficient strategy for this example is to take $k = 0$ and reduce the error via mesh refinement.
In the accompanying talk, I will present additional results to showcase our method, including cases where increasing $k$ provides a significant advantage.

**Leafnet algorithm for $N$-dimensional nonlinear hyperbolic systems**

**Emmanuel Lorin**[1,*], **Arian Novruzi**[2]
[1]School of Mathematics and Statistics, Carleton University, Ottawa, Canada
[2]Department of Mathematics and Statistics, University of Ottawa, Ottawa, Canada
*Email: elorin@math.carleton.ca

**Abstract**

In this presentation, we introduce *leafnet*, a neural network-based algorithm for the diffusion-free approximation of weak solutions to $N$-dimensional hyperbolic conservation laws. A mathematical analysis of this algorithm, as well as some numerical experiments will be proposed.

## 1 Introduction

The use of neural networks for solving (parameterized) PDE has recently become very popular, as beyond the use of automatic differentiation, they allow for the inclusion of external (experimental or numerical) data. Among the existing techniques the most popular are the Physics Informed Neural Network (PINN) methods [1].

This presentation is devoted to neural network-based solutions of $N$-dimensional nonlinear hyperbolic conservation laws (HCL). PINN-like methods allow for the approximation of smooth solutions to HCL, avoiding to a certain extent stability issues usually due to time-stepping. Indeed, PINN algorithms allow for formulating and solving IVP problems as a BVP, where the initial condition is treated as a boundary condition. Stability issues may however be "transferred" to the optimization algorithm. As for the shock waves, usually their approximation is known to be very inaccurate. This is due to the fact that shock waves are non-smooth solutions, defined from the weak formulation of the HCL, while PINN methods approximate PDEs in their strong form. As a consequence, a direct use of PINN algorithms for approximating shock waves usually requires the addition of artificial diffusion.

We develop the Leafnet algorithm (see [2]) allowing for accurately approximating shock waves using neural networks, by i) a direct use of the PINN method where solutions are smooth, ii) approximating Rankine-Hugoniot jump conditions, iii) reconstructing overall piecewise smooth solutions. The proposed method can be seen as a multidimensional extension of the algorithm which was developed in [3] for one-dimensional systems, and which largely relies on the structure of solutions as well as mathematical properties of HCL.

We will focus on the approximation of scalar shock waves for $N$-dimensional hyperbolic conservation laws, including their generation from smooth Cauchy data. We also will design self-similar neural networks allowing for an accurate approximation of rarefaction waves. A detailed mathematical analysis of the error for shock wave solutions will be presented, as well as some numerical experiments.

## 2 Methodology for scalar shock waves in $N$ dimensions

Let us summarize the strategy in the case of scalar shock wave, solution to nonlinear scalar hyperbolic conservation law in $Q := \Omega \times [0,T]$ where $\Omega = \prod_{i=1}^{N}(a_i,b_i)$, or $\Omega = \mathbb{R}^N$ with $N \geq 1$:

$$\partial_t u + f(u)\cdot\nabla u = 0, \qquad \text{in } Q, \tag{1a}$$

$$u(\cdot,0) = u_0(\cdot), \qquad \text{on } \Omega, \tag{1b}$$

where $f(u) = (f_1(u),\ldots,f_N(u)) =: F'(u)$, $F \in C^2(\mathbb{R};\mathbb{R}^N)$, $\nabla = (\partial_1,\ldots,\partial_N)$, and $u : (x,t) := (x_1,\ldots,x_N,t) \in Q \to u(x,t) \in \mathbb{R}$. Given $u_0^{\pm} \in C^1(\Omega)$, and a $C^1$-surface $\Gamma_0 \subset \Omega$, parameterized by $g_0 \in C^1(\Omega')$, $\Omega' := \prod_{i=1}^{N-1}(a_i,b_i)$,

$$\Gamma_0 = \left\{(x',x_N) \in \Omega,\ x_N = g_0(x'), x' \in \Omega'\right\},$$

we define $\Omega^{\pm} \subset \Omega$ by

$$\Omega^+ = \{(x',x_N) \in \Omega,\ x_N > g_0(x')\},$$
$$\Omega^- = \{(x',x_N) \in \Omega,\ x_N < g_0(x')\},$$

and the $C^1(\Omega)$-piecewise function $u_0$ by

$$u_0(x) = u_0(x',x_N) = \begin{cases} u_0^+(x), & x_N > g_0(x'), \\ u_0^-(x), & x_N < g_0(x'). \end{cases} \tag{2}$$

We assume that the shock wave is parameterized by $x_N = g(x',t)$

$$\Gamma = \left\{(x',g(x',t),t) \in Q \ :\ (x',t) \in Q'\right\},$$

where $Q' := \Omega' \times (0,T)$ and $g$ has to be identified. Let $\nu = (\nu_1, \dots, \nu_N, \nu_t)$ be the unit normal vector on $\Gamma$ oriented on the direction of $x_N$

$$\nu = \frac{\nabla_{(x,t)}(x_N - g(x',t))}{|\nabla_{(x,t)}(x_N - g(x',t))|}\,.$$

The Rankine-Hugoniot jump condition on $\Gamma$ in terms of $u^\pm$ reads

$$\nu_t [\![u]\!] + [\![F(u)]\!] \cdot (\nu_1, \dots, \nu_N) \;=\; 0, \quad (3)$$

where $u^\pm$ satisfy

$$\partial_t u^\pm + f(u^\pm) \cdot \nabla u^\pm = 0, \quad \text{in } Q, \quad (4a)$$
$$u^\pm(\cdot, 0) = u_0^\pm, \quad \text{on } \Omega. \quad (4b)$$

Subject to taking $T$ smaller, we assume that there is no shock formation for both IBVPs (4). Ranking-Hugoniot (3) can then be rewritten

$$\begin{array}{rcl} \partial_t g + \sum_{i=1}^{N-1} a_i(g, x', t)\partial_i g & = & b(g, x', t) \\ g(\cdot, 0) & = & g_0(\cdot)\,, \end{array} \quad (5)$$

where

$$\begin{array}{rcl} a_i(g, x', t) & = & \int_0^1 f_i(u^-(x', g, t) \\ & & +s(u^+(x', g, t) - u^-(x', g, t)))ds, \\ b(g, x', t) & = & \int_0^1 f_N(u^-(x', g, t) \\ & & +s(u^+(x', g, t) - u^-(x', g, t)))ds. \end{array}$$

Hence, "solving" the shock $\Gamma$ of (1) is equivalent to solve the $N-1$ dimensional IVP (5). The *leafnet method* consists in using physics-informed neural networks with randomly chosen collocation points, for:

1. Solving 2 IBVP (4): $u^-$ and $u^+$ in $Q$.
2. Solving Rankine-Hugoniot (5): $g$ in $Q'$.
3. Reconstructing solution from IBVP (4) for $(x', x_N, t) \in Q$

$$u(x', x_N, t) = \left\{ \begin{array}{ll} u^-(x', x_N, t), & \text{if } x_N < g(x', t) \\ u^+(x', x_N, t), & \text{if } x_N > g(x', t) \end{array} \right.$$

The method can also handle shock formation. In this goal, we consider a virtual discontinuity curve $\Gamma_0$ with equations $x_N = g_0(x')$, which separates $\Omega$ in $\Omega^\pm$, and generates two initial functions $u_0^\pm$ as in (2), with naturally $u_0^- = u_0^+$ on $\Gamma_0$. Then, we solve $u^\pm$ and $g$ as above. The time formation $t_*$ is smallest time where the jump $[\![u]\!]$ on $\Gamma$ becomes discontinuous. For $t \in (t_*, T)$, $g$ must satisfy the Rankine-Hugoniot condition and therefore $g$ will parameterize the shock surface $\Gamma$.

## 3 Example

Let us illustrate the methodology in the case of the formation of a scalar shock wave from a smooth initial condition. We set $Q = \Omega \times [0,T]$ with $\Omega = (-1,1) \times (-2,2)$ and $T = 0.8$. The flux function $F = (F_1, F_2)$ is $F_1(u) = u^2$, $F_2(u) = u^2/4$ and $u_0$ is given by

$$u_0(x_1, x_2) = \frac{1}{2}(1 - \tanh(5(x_1 + (x_2+2)^2/20)))\,.$$

The solution is smooth up to time $t = t^*$, corresponding to the generation of a shock wave, here equal to 0.2. The neural networks have 2 hidden layers with 15 neurons each. In Fig. 1 (Left) we report $u_0$ and the solution at $t = t^*$ (Middle), and the leafnet solution (Right) at $t = T$. We report in Fig. 2 (Left) the surface of discontinuity, and the loss functions for the leafnet algorithm in Fig. 2 (Right).

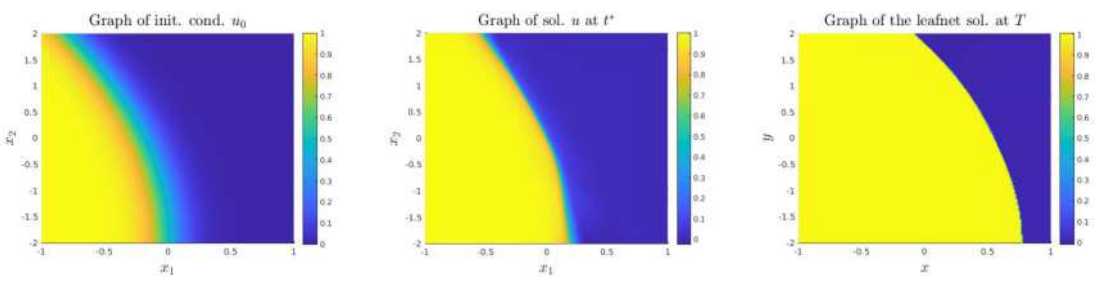

Figure 1: (Left) Graph of $u_0$. (Middle) Graph of the solution at $t = t^*$. (Right) Graph of the leafnet solution at $t = T$.

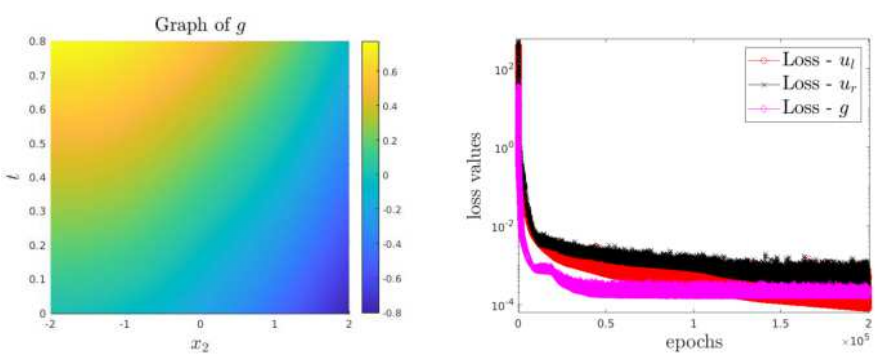

Figure 2: (Left) Graph of the surface of discontinuity. (Right) Loss functions.

# Optical wave reconstruction in highly turbulent random media

**Kevin Luna**[1,*], **Max Cubillos** [1]
[1]Air Force Research Lab, Albuquerque, United States
*Email: afrl.rdl.orgbox@us.af.mil

**Abstract**

The problem of reconstructing a wavefront from noisy phase gradient data is central to adaptive optics. While standard algorithms work for smooth, single-valued phase functions, they fail when waves propagate through strongly turbulent media, where the wave's intensity develops zeros and the phase becomes singular. Our framework overcomes this limitation by reconstructing the entire complex wave field, which is unaffected by phase singularities. It uses a variational approach to minimize the error between the field's gradient and the measured data, resulting in significant accuracy improvements over traditional methods in strong turbulence.

## 1 Introduction

A Shack-Hartmann sensor, using an array of lenses, produces noisy pointwise intensity $I$ and phase gradient $\boldsymbol{g} = (g_1, g_2)$ measurements over a region of the plane $\Omega \subset \mathbb{R}^2$. We associate $I$ with $\rho^2$ and $\boldsymbol{g}$ with $\nabla\varphi$. In the optics context, the wave field is almost exclusively considered in polar form, $u(\boldsymbol{x}) = \rho(\boldsymbol{x}) \exp\big(i\varphi(\boldsymbol{x})\big)$, where $\rho$ and $\varphi$ are real scalar functions of $\boldsymbol{x} = (x, y) \in \mathbb{R}^2$. In practice, subsequent processing steps and hardware in imaging systems exclusively require $\phi(x, y)$. Consequently, optical scientists and engineers heavily focus on the following linear inverse problem for the phase of the wave.

*Provided noisy phase gradient measurements* $\mathbf{g}$*, find an estimate* $\widetilde{\phi}$ *for the phase of the incoming wave* $u$ *such that*

$$\widetilde{\phi} = \underset{\phi}{\operatorname{argmin}} \frac{1}{2} \|\nabla\varphi - \boldsymbol{g}\|^2. \tag{1}$$

The "phase reconstruction" problem is central to many applications, but traditional solution methods fail when the phase function develops branch point singularities. These singularities form when strong turbulence in the propagation media causes substantial fluctuations in the field's intensity. Standard approaches are ill-posed in these conditions as they amount to reconstructing a singular field directly.

We propose a reformulation of the problem that circumvents this issue. By focusing on the full optical field, which remains well-behaved, we can reconstruct it accurately. The singular phase, if desired, can then be recovered from the imaginary part of the logarithm of the reconstructed complex field. This approach avoids the inherent difficulties of directly reconstructing a singular function.

## 2 Linear Least Squares for Full Field reconstruction

We start by writing the optical field in polar form, $u(\boldsymbol{x}) = \rho(\boldsymbol{x}) \exp\big(i\varphi(\boldsymbol{x})\big)$, where $\rho$ and $\varphi$ are real scalar functions of $\boldsymbol{x} = (x, y) \in \mathbb{R}^2$. The gradient of $u$ can then be written as

$$\nabla u = u\Big(\frac{\nabla(\rho^2)}{2\rho^2} + i\nabla\varphi\Big), \tag{2}$$

where we use the fact that $\nabla(\rho^2) = 2\rho\nabla\rho$. Therefore, we would like to find $u$ that satisfies (2) with $I$ and $\boldsymbol{g}$ substituted into the right-hand side. In other words, we would like to find $u$ such that

$$\nabla u - \Big(\frac{\nabla I}{2I} + i\boldsymbol{g}\Big)u \approx 0. \tag{3}$$

Therefore, to obtain the solution to this problem, we formulate a minimization problem in a weighted least-squares sense by finding the minimizer of the functional

$$\widetilde{J}[u] = \frac{1}{2}\|I^p \nabla u - \boldsymbol{h}u\|^2. \tag{4}$$

where, $p \geq 0$ is a weighting parameter and the vector-valued function $\boldsymbol{h}(\boldsymbol{x})$ is given in terms of the data,

$$\boldsymbol{h} = I^{p-1}\boldsymbol{h}_0, \quad \boldsymbol{h}_0 = \frac{1}{2}\nabla I + iI\boldsymbol{g}. \tag{5}$$

Minimization of the functional $\widetilde{J}[u]$ is an ill-posed problem due to non-uniqueness; $u = 0$ and any scaled multiple $cu$ are also always solutions. To remedy this, we introduce an additional constraint. Since applications typically only require phase differences, we can fix the solution at a single point $\boldsymbol{x}_0$. To account for measurement noise, this condition is enforced weakly by adding a penalty term, yielding the modified functional:

$$J[u] = \frac{1}{2}\|I^p \nabla u - \boldsymbol{h} u\|^2 + \frac{c_0}{2}\left|u(\boldsymbol{x}_0) - \sqrt{I}(\boldsymbol{x}_0)\right|^2, \tag{6}$$

where $c_0 > 0$ is a penalty constant. The reconstructed field $u_{\text{rec}}$ is found by solving the following linear least squares problem:

$$u_{\text{rec}} = \underset{u}{\operatorname{argmin}}\, J\,[u]\,. \tag{7}$$

Depending on the application and hardware either direct or iterative solvers are used in practice for the phase reconstruction problem. The direct algorithm is the version we examine in this work. However, an iterative solver can be developed using the Euler-Lagrange equation that extremizes 6. This PDE can be found in [1]

## 3 Solution Algorithm

Without loss of generality, we assume the Shack-Hartmann sensor is given by an $M \times M$ lenslet array whose centers are located on equispaced points in $[-1, 1]^2$. The computational grid is given by an $N \times N$ Chebyshev-Gauss-Lobatto grid on $[-1, 1]^2$. The data we interpolate onto the computational grid is the intensity $I$ and the product of the intensity and the phase gradients $I\boldsymbol{g}$. Note that we apply spectral filtering in the presence of noise.

Let $u_{jk}$, $I_{jk}$, $\boldsymbol{h}_{jk}$ denote, respectively, the values of the optical field, intensity, and the function $\boldsymbol{h}(x)$ in (5) at the Chebyshev nodes. Let $D_{lm}$ denote the Chebyshev collocation derivative matrix, and let $(j_0, k_0)$ denote the location of the grid point corresponding to that of the penalty term $\boldsymbol{x}_0$. We may write the least squares problem with the penalty term added only to the first component as

$$I^p_{jk} \sum_{l=1}^{N} D_{jl} u_{lk} - h_{1,jk} u_{jk} + c_0 \delta_{jj_0} \delta_{kk_0} u_{j_0 k_0} = c_0 \delta_{jj_0} \delta_{kk_0} \sqrt{I_{j_0 k_0}},$$

$$I^p_{jk} \sum_{l=1}^{N} D_{kl} u_{jnaul} - h_{2,jk} u_{jk} = 0,$$

which is clearly overdetermined. Let $v = \text{vec}(u)$ denote the vectorization of $u_{jk}$, $A$ the matrix corresponding to the vectorized left-hand-side operator, and $f$ the vectorization of the right-hand side. Then we may write the least squares problem succinctly as $Av = f$. These equations may be solved using standard methods, based off the SVD and QR decomposition.

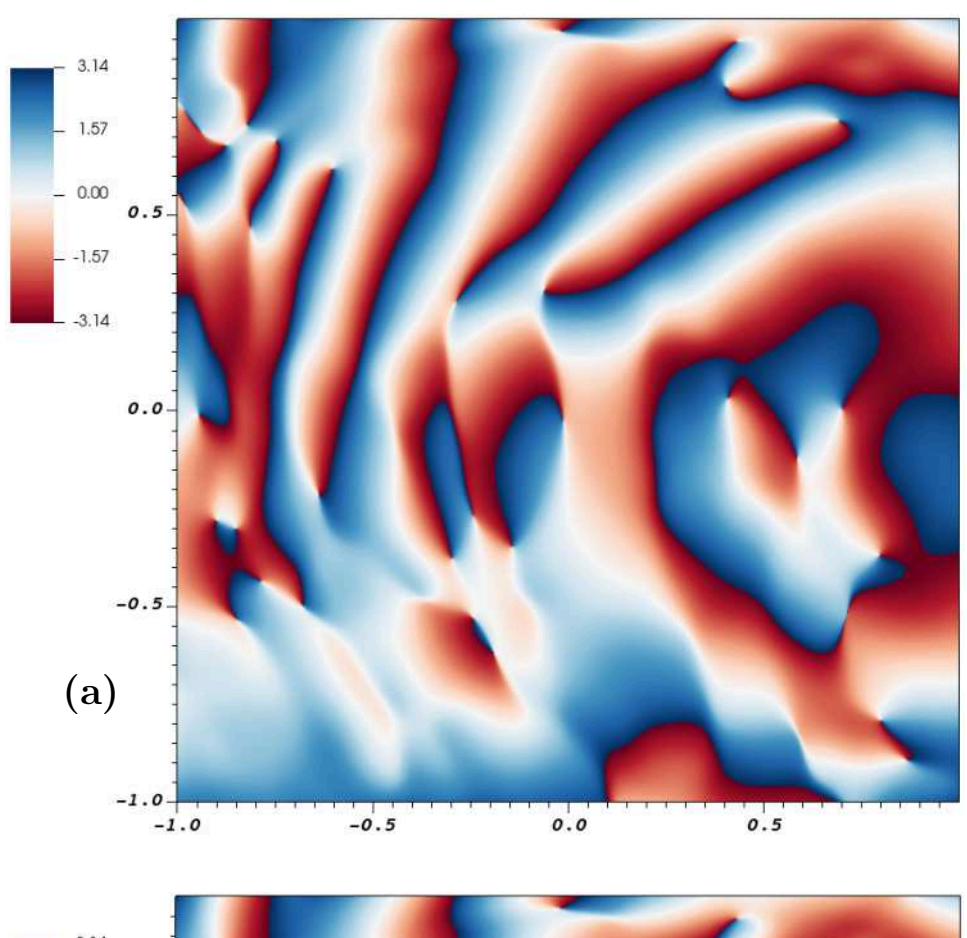
(a)

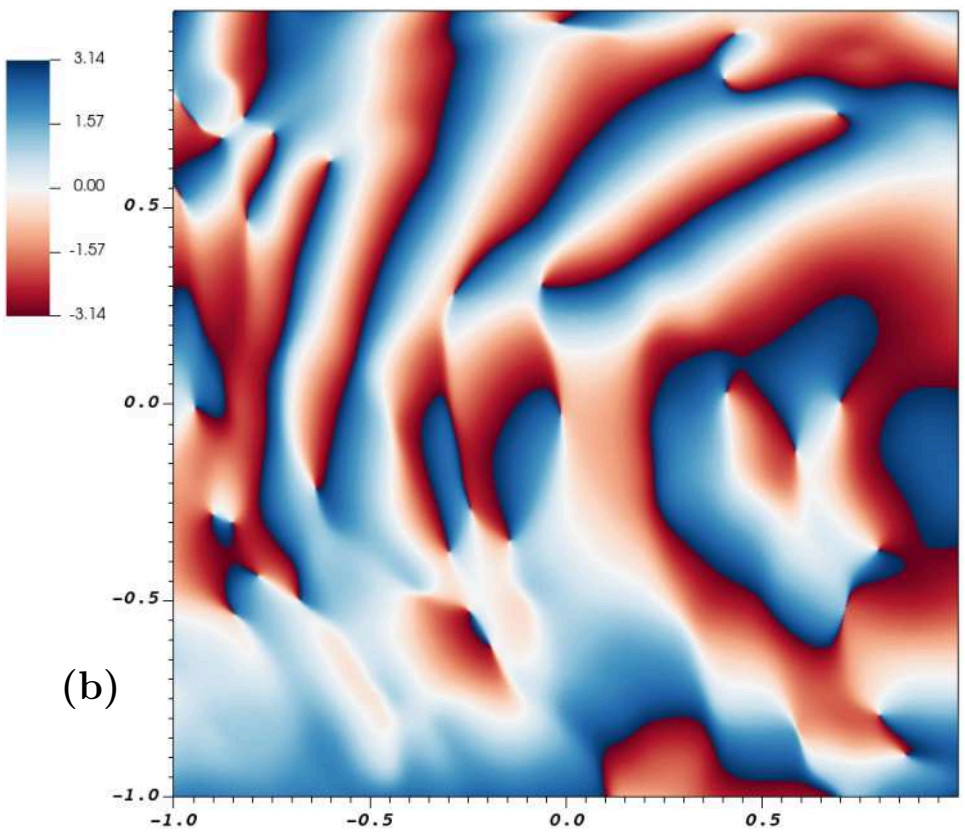
(b)

Figure 1: The (a) exact and (b) reconstructed phase profiles for simulation of a Gaussian beam propagated through strong turbulence.

## 4 Results

As an example, we show a simulation of a Gaussian beam propagated through strong turbulence in figure 1. The parameters for the simulation are: wavelength $\lambda = 1\mu$m, propagation distance $z = 1.2$km, width of the beam $w = 1$cm, and the turbulent refractive index is modeled using a Von-Karman spectrum with inner scale $\ell_0 = 1$cm, outer scale $L_0 = 0.5$m, and structure parameter $C_n^2 = 3.5 \times 10^{-13}\text{m}^{-2/3}$. We use $100 \times 100$ subapertures and a $64 \times 64$ Chebyshev grid for the reconstruction. The relative $L^\infty$ error for the optical field is $E = 2.6 \times 10^{-3}$.

## 5 Acknowledgements

Approved for public release; distribution is unlimited. Public Affairs release approval AFRL-2026-0373.

# Spectral analysis of a closed perturbed waveguide filled with a hyperbolic medium

**Maryna Kachanovska**[1], **Dylan Machado**[1,*]
[1]POEMS, CNRS, Inria, ENSTA, Institut Polytechnique de Paris, 91120 Palaiseau, France.
*Email: dylan.machado@ensta.fr

**Abstract**

We investigate the spectrum of the first-order operator that governs the propagation of electromagnetic waves in the case of a closed locally perturbed waveguide filled with a hyperbolic medium. We prove a limiting absorption principle (LAP) for a range of frequencies in an interval defined by a certain geometric condition on the perturbation. In addition, we study the well-posedness of the limit problem arising from the LAP.

## 1 Introduction and motivation

Hyperbolic electromagnetic metamaterials are anistropic artificial media which exhibit a lot of interesting physics phenomena such as the negative refraction or the rainbow-trapping [1].
The propagation of electromagnetic waves is governed by Maxwell's equations which depend on the constitutive laws of the considered environment. In the case of hyperbolic mediums, the simplest mathematical model stems from the physics of plasma, in particular considering strongly magnetised cold plasmas. In the time-harmonic regime, 2D Maxwell's equations read

$$\left| \begin{aligned} & i\omega\underline{\underline{\varepsilon}}(\omega)\hat{\mathbf{E}}(\omega,\cdot) - \mathbf{curl}\,\hat{\mathrm{H}}_z(\omega,\cdot) = 0 \\ & i\omega\hat{\mathrm{H}}_z(\omega,\cdot) + \mathrm{curl}\,\hat{\mathbf{E}}(\omega,\cdot) = 0, \qquad \omega\in\mathbb{R}. \end{aligned} \right. \tag{1}$$

The dielectric tensor is given by

$$\underline{\underline{\varepsilon}}(\omega) = \begin{pmatrix} 1 & 0 \\ 0 & \varepsilon(\omega) \end{pmatrix}, \; \varepsilon(\omega) = 1 - \frac{\omega_p^2}{\omega^2} \tag{2}$$

where $\omega_p > 0$ is the plasma frequency. Eliminating the electric field $\hat{\mathbf{E}}$ yields

$$\varepsilon(\omega)^{-1}\partial_x^2\hat{\mathrm{H}}_z + \partial_y^2\hat{\mathrm{H}}_z + \omega^2\hat{\mathrm{H}}_z = 0. \tag{3}$$

In the frequency range $\mathcal{I}_\mathrm{p} := (-\omega_p,\omega_p)\backslash\{0\}$, $\varepsilon(\omega)$ becomes negative so that the equation (3) becomes hyperbolic [2].
We consider the equations (1) posed in a straight waveguide $\Omega\subset\mathbb{R}^2$ with a locally varying boundary, see Figure 1a. In addition, we require that $\hat{\mathrm{H}}_z(\omega,\cdot)$ satisfies the homogeneous Dirichlet condition on the boundary $\Gamma = \partial\Omega$. The temporal form of (1) can be written in terms of a Schrödinger-like equation

$$\left| \begin{aligned} & \mathbf{u}'(t) + i\mathbb{A}\mathbf{u}(t) = 0 \quad t\in\mathbb{R} \\ & \mathbf{u}(0) = \mathbf{u}_0 \end{aligned} \right. \tag{4}$$

where $(\mathbb{A}, D(\mathbb{A}))$ is a self-adjoint operator on the energy space $\mathcal{H} := \mathbb{L}^2(\Omega)^4$. The goal of this work in progress is to describe the spectrum of the first-order operator $\mathbb{A}$ leading to the nonlinear eigenvalue problem (5). Let us finally point out that a related problem has been studied in [3]. Here, we present an alternative approach to the LAP.

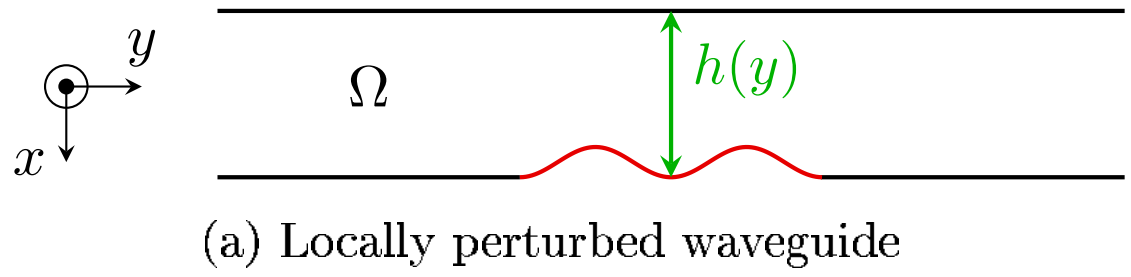

(a) Locally perturbed waveguide

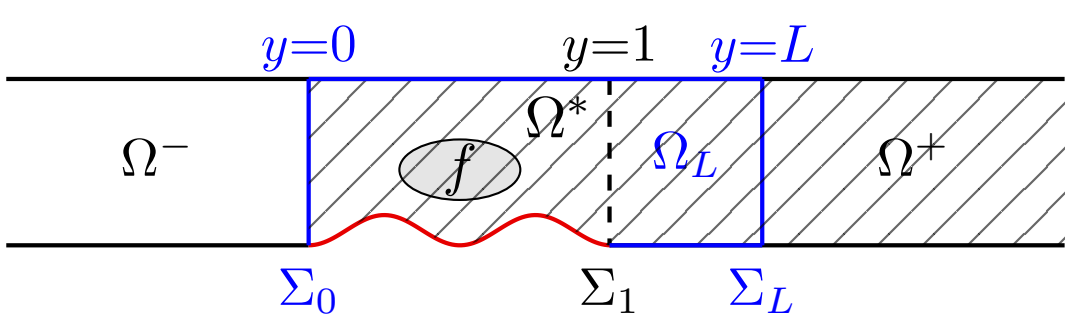

(b) Truncation of the waveguide

## 2 Limiting absorption principle

Let us consider the hyperbolic boundary value problem

$$\left| \begin{aligned} & \varepsilon(\omega)^{-1}\partial_x^2 u + \partial_y^2 u + \omega^2 u = f \quad \text{in } \Omega \\ & u\in\mathbb{H}_0^1(\Omega) \end{aligned} \right. \tag{5}$$

for $\omega\in\mathcal{I}_\mathrm{p}$ and $f\in\mathbb{L}^2_{\mathrm{comp}}(\Omega)$. We introduce the map $h:\mathbb{R}\to\mathbb{R}$ that parametrizes the perturbation, see Figure 1a. We asssume that $h$ is a positive diffeomorphism of class $\mathcal{C}^2$ which is constant outside the interval $(0,1)$. The perturbed domain is denoted by $\Omega^*$ and we assume that the support of $f$ is localized inside $\Omega^*$, see Figure 1b. Adding absorption $\eta>0$ to the frequency, $\omega_\eta := \omega + i\eta$, the problem (5) becomes elliptic. Thus, there is a unique solution $u_\eta\in\mathbb{H}_0^1(\Omega)$ of (5). The main result is the

**Theorem 1.** *Let $f \in \mathbb{L}^2_{\text{comp}}(\Omega)$ and $\omega \in \mathcal{I}_{\text{p}}$ a frequency satisfying the subsonic condition*

$$\sup_{y\in\mathbb{R}} |h'(y)| < -\varepsilon(\omega)^{-1}. \tag{6}$$

*For $L > 1$, there are $C_L > 0$ and $\eta_L > 0$ s.t*

$$\forall 0 < \eta < \eta_L, \quad \|u_\eta\|_{\mathbb{L}^2(\Omega_L)} \leq C_L \|f\|_{\mathbb{L}^2(\Omega)}. \tag{7}$$

*Proof.* Let us give a concise sketch of the proof. Without loss of generality, we assume that $\varepsilon(\omega) = -1$ and $\omega = 1$. The geometric setting is displayed in the figure 1b. With the absorption $\eta > 0$, there is a unique $u_\eta \in \mathbb{H}^1_0(\Omega)$ verifying in $\Omega$

$$-(1+i\eta)\partial_x^2 u_\eta + \partial_y^2 u_\eta + (1+i\eta)u_\eta = f. \tag{8}$$

Introducing now the Dirichlet-to-Neumann (DtN) operator $\Lambda^-_\eta$ associated to the exterior problem in $\Omega^- = \Omega\backslash\Omega^+$, one can prove that the original problem is equivalent to the truncated problem

$$\left| \begin{array}{ll} -(1+i\eta)\partial_x^2 u_\eta + \partial_y^2 u_\eta + (1+i\eta)u_\eta = f & \Omega^+ \\ \gamma_0 u_\eta = 0 & \Gamma \\ \partial_y u_\eta = -\Lambda^-_\eta \gamma_0 u_\eta & \Sigma_0. \end{array} \right. \tag{9}$$

In particular, indentifying the $y$ variable with the time, we re-interpret the latter as the wave equation posed in a time-dependent domain with a special condition linking the initial conditions. In moving domains whose boundary satisfies the subsonic condition, the well-posedness of the Cauchy problem for the wave equation justifies this interpretation. Hence, one can use the energy estimates up to a time $y = L$. To conclude the proof, it suffices to control uniformly in $\eta$ the data (initial conditions and right-hand side) using the ellipticity nature of the problem with the absorption and also the uniform positivity of the DtN maps. □

Finally, we underline the fact that our proof of the LAP can be adapted to $\mathbb{R}^d$.

## 3 The spectrum of $\mathbb{A}$

In the same spirit as [4], the resolvent of $\mathbb{A}$, see (4), is given by

$$\mathcal{R}(z) = \mathrm{T}_z + \mathrm{V}_z\mathcal{S}_z\mathrm{F}_z, \ z \in \mathbb{C}\backslash\mathbb{R}$$

where $\mathrm{T}_z$, $\mathrm{V}_z$ and $\mathrm{F}_z$ are known operators of differential type, and $\mathcal{S}_z : \mathbb{H}^{-1}(\Omega) \to \mathbb{H}^1_0(\Omega)$ is the solution operator of (5). Let be $\Lambda_{\text{sub}}$ the set of frequencies $\omega \in \mathcal{I}_{\text{p}}$ satisfying the subsonic condition (6). On can prove that $\Lambda_{\text{sub}} = (-\omega_p, -\omega^*) \cup (\omega^*, \omega_p)$ where $\omega^* \in \mathcal{I}_{\text{p}}$ is a critical frequency. A consequence of the theorem 1 is the following

**Corollary 2.** *The spectrum of $\mathbb{A}$ on $\Lambda_{\text{sub}}$ is purely absolutely continuous, i.e*

$$\Lambda_{\text{sub}} \cap \sigma(\mathbb{A}) \subset \sigma_{\text{ac}}(\mathbb{A}).$$

## 4 Well-posedness of the limit problem

Let us study the limit problem obtained when the absorption $\eta \searrow 0$ in (9). We introduce the exterior DtN operators $\Lambda^\pm$ without the absorption on the boundaries $\Sigma_0$ and $\Sigma_1$, see Figure 1b. Also, let $\varphi : Q \to \Omega$ be the waveguide diffeomorphism where $Q = (0,1)\times\mathbb{R}$ is the straight waveguide.

**Theorem 3.** *Let $f \in \mathbb{L}^2(\Omega^*)$, $\omega \in \mathcal{I}_{\text{p}}$ verifying (6) and $I = [0,1]$. There is a unique $u \in \mathbb{L}^2(\Omega^*)$ such that $u\circ\varphi \in \mathcal{C}^0(I, \mathbb{H}^{1/2}_{00}(I)) \cap \mathcal{C}^1(I, \mathbb{H}^{-1/2}(I))$ and satisfies the following two-point boundary-value problem for the wave equation*

$$\left| \begin{array}{ll} (\Box + 1)\, u = f & \Omega^* \\ \gamma_0 u = 0 & \Gamma \\ \partial_y u = \pm\Lambda^\pm \gamma_0 u & \Sigma_0 \cup \Sigma_1 \end{array} \right. \tag{10}$$

The proof of this result uses the well-posedness of the Cauchy problem and the shooting method like ideas.

# A Posteriori Estimates for Space Discretizations of the Wave Equation with Dynamic Meshes

**Théophile Chaumont-Frelet**[1], **Sumit Mahajan**[1,*], **Martin Vohralík**[2]

[1]Inria, Univ. Lille, CNRS, UMR 8524 – Laboratoire Paul Painlevé, 40 Av. Halley, 59650 Villeneuve-d'Ascq, France

[2]Inria, 48 rue Barrault, 75647 Paris, France & CERMICS, Ecole nationale des ponts et chaussées, IP Paris, 77455 Marne-la-Vallée, France

*Email: sumit.mahajan@inria.fr

**Abstract**

In this work, we derive a posteriori error estimates for the semidiscrete finite element approximation of the scalar wave equation on dynamically evolving spatial meshes. The analysis requires the right-hand side to possess finite temporal regularity (up to a finite number of time derivatives), while only $L^2(\Omega)$ spatial regularity is assumed in space, allowing for discontinuities such as piecewise polynomial data with jumps across mesh elements. The error is measured in a damped energy norm, and the proposed estimator is proven to be reliable and globally efficient in this norm, with constants independent of the mesh size and polynomial degree. Numerical experiments illustrate the theory and demonstrate the robustness of the estimator for long-time simulations with dynamic mesh adaptation.

## 1 Introduction

An accurate numerical simulation of the wave equation is challenging due to the presence of propagating wave fronts and localized transient phenomena. Classical finite element methods on fixed spatial meshes are often inefficient for such problems, as the relevant features evolve in space and time and require globally fine discretizations. Dynamic mesh strategies address this limitation by allowing the spatial discretization to change at selected time instances, concentrating degrees of freedom where they are most needed and thereby reducing computational cost for a prescribed accuracy. However, mesh transitions introduce additional difficulties, including projection-induced jump terms and potential loss of energy consistency, which must be carefully controlled to ensure stability and to design reliable and efficient a posteriori error estimators.

## 2 Model

In this work, we consider the scalar wave equation posed on a bounded Lipschitz polyhedral domain $\Omega \subset \mathbb{R}^d$ $(d = 2, 3)$, and over the semi-infinite time interval $\mathcal{I} = [0, +\infty)$. This consists of finding $u : \mathcal{I} \times \Omega \to \mathbb{R}$ such that

$$\begin{cases} \mu\ddot{u} - \boldsymbol{\nabla} \cdot (\boldsymbol{A}\boldsymbol{\nabla} u) = \mu f & \text{in } \mathcal{I} \times \Omega, \\ u = 0 & \text{on } \mathcal{I} \times \partial\Omega, \\ u|_{t=0} = \dot{u}|_{t=0} = 0 & \text{in } \Omega, \end{cases} \quad (1)$$

where $\mu : \Omega \to \mathbb{R}$ and $\boldsymbol{A} : \Omega \to \mathbb{R}^{d\times d}$ are spatially varying coefficients and $f : \mathcal{I} \times \Omega \to \mathbb{R}$ is a given source term, homogeneous Dirichlet boundary conditions and zero initial data are assumed. Furthermore, the source term $f$ is assumed to be smooth in time and compactly supported away from $t = 0$.

## 3 Methodology

Our analysis extends the framework introduced in [1, 2] for static meshes based on frequency-splitting argument to dynamically evolving spatial meshes. In the dynamic setting, changes in the spatial discretization at temporal interfaces induce discontinuities in the global discrete solution across time slabs. These interfacial jumps, defined in (4) and quantified by the functionals $\mathcal{J}^{\ell}_{\rho}(\cdot)$ in (3), are intrinsic to the numerical method and give rise to additional error contributions that must be explicitly accounted for in the a posteriori analysis.

A central ingredient of our approach is the measurement of the error in a *damped energy norm*, defined by

$$|||u|||^2_{\rho} = \int_0^{+\infty} \Big(\|\dot{u}(t)\|^2_{\mu,\Omega} + \|\boldsymbol{\nabla} u(t)\|^2_{\boldsymbol{A},\Omega}\Big) e^{-2\rho t}\,\mathrm{d}t,$$

where $\rho > 0$ is a damping parameter that prevents the accumulation of errors over long time interval by selecting $\rho = 1/T$.

To accommodate mesh-induced discontinuities, the error is estimated through the residual of a suitably reconstructed discrete solution $\tilde{u}_h \in C^1(\mathcal{I}; X)$, which restores sufficient temporal regularity across mesh transitions. This reconstruction removes the temporal jump discontinuities induced by mesh changes while preserving Galerkin orthogonality on each time slab.

This leads naturally to a fully computable time-dependent a posteriori estimator $\eta(t)$, defined via the instantaneous residual

$$v \;\mapsto\; (\mu f(t), v)_\Omega - (\mu \ddot{\tilde{u}}_h(t), v)_\Omega - (\boldsymbol{A} \boldsymbol{\nabla} \tilde{u}_h(t), \boldsymbol{\nabla} v)_\Omega,$$

for all admissible test functions $v$. Integrating the estimator in time yields

$$\Lambda_\rho^2 := \int_0^\infty \eta^2(t)\, e^{-2\rho t}\, \mathrm{d}t,$$

which forms the core quantity in the error bounds.

## 4 Results

We partition the time domain $\mathcal{I} = [0, \infty)$ into a finite sequence of subintervals as

$$\mathcal{I} = \bigcup_{\ell=0}^{N-1} \mathcal{I}^\ell \cup \mathcal{I}^N, \quad \text{where } \mathcal{I}^\ell := (T_\ell, T_{\ell+1}],$$

for $\ell = 0, \dots, N-1$, and $\mathcal{I}^N := (T_N, \infty)$, with $0 = T_0 < T_1 < \cdots < T_N < \infty$. On each time interval $\mathcal{I}^\ell$, the spatial domain $\Omega$ is discretized by a conforming simplicial mesh $\mathcal{T}_h^\ell$.

For the error $\tilde{\xi}_h := u - \tilde{u}_h$, the reliability estimate separates the contributions of the time-dependent residual, data oscillation, and mesh-change induced jump terms. Our key results read as follows: We have the upper bound

$$\left|\!\left|\!\left| \tilde{\xi}_h \right|\!\right|\!\right|_\rho^2 \le \; \Lambda_\rho^2 + C \sum_{\ell=1}^{N} \mathcal{J}_\rho^\ell(\tilde{u}_h) + h.o.t. \qquad (2)$$

The final term in the bound arises from the dynamic nature of the mesh given as

$$\mathcal{J}_\rho^\ell(v) := \frac{1}{2} \sum_{k=2}^{M} \left( \|\mathfrak{J}_\ell^{(k)}(v)\|_{\mu,\Omega}^2 + \|\mathfrak{J}_\ell^{(k-1)}(\boldsymbol{\nabla} v)\|_{\boldsymbol{A},\Omega}^2 \right) e^{-2\rho T_\ell}, \qquad (3)$$

for $\ell = 1, \dots, N$, where $M \ge 2$ is a fixed integer denoting the highest order of temporal derivatives included in the jump contributions and the temporal jump (higher derivative) of order $k = 2, \dots, M$, is

$$\mathfrak{J}_\ell^{(k)}(v) := v_\ell^{(k)}(T_\ell) - v_{\ell-1}^{(k)}(T_\ell). \qquad (4)$$

These terms are unavoidable in the presence of evolving finite element spaces, and the present estimate makes their influence explicit and fully controlled.

We also show the lower bound as

$$\Lambda_\rho^2 \le C \left( \left|\!\left|\!\left| \tilde{\xi}_h \right|\!\right|\!\right|_\rho^2 + \sum_{\ell=1}^{N} \mathcal{J}_\rho^\ell(\tilde{u}_h) + h.o.t. \right). \qquad (5)$$

All constants in the lower bound are independent of the polynomial degree $p$ of the finite element approximation. Thus, the estimator is efficient and $p$–robust, up to the unavoidable contribution of the jump terms, which characterize precisely the effect of the dynamic mesh on the a posteriori error control.

We note that, in the case of static meshes, the jump contributions vanish, i.e. $\mathcal{J}_\rho^\ell(\cdot) = 0$, and, under the assumed temporal smoothness of the data, the oscillation terms are of higher order. Consequently, the estimates reduce to

$$C^{-1} \Lambda_\rho^2 - \text{h.o.t.} \;\le\; \left|\!\left|\!\left| \tilde{\xi}_h \right|\!\right|\!\right|_\rho^2 \;\le\; C \left( \Lambda_\rho^2 + \text{h.o.t.} \right),$$

which corresponds to the residual-based a posteriori error bounds obtained in [1].

These results indicate that dynamic mesh transitions can be rigorously incorporated into a posteriori error analysis without sacrificing robustness or asymptotic accuracy. The proposed estimator therefore provides a reliable tool for adaptive wave simulations on evolving meshes and opens the way toward fully adaptive space-time strategies for hyperbolic problems.

# Adapting Full-Waveform Inversion for Subsurface $CO_2$ Sequestration

Alison Malcolm[1,*], Vahid Entezaar-Saadat, Xingda Jiang, Abolfazl Khan Mohammadi, Colin Farquharson[1]

[1]Department of Earth Sciences, Memorial University of Newfoundland, St John's, Canada

*Email: amalcolm@mun.ca

**Abstract**

To help mitigate climate change, one strategy is to inject $CO_2$ into underground systems. Ensuring that the $CO_2$ has remained in place is an important aspect of doing this, and is typically done using seismic waves. This monitoring sets up a fundamentally different inverse problem in which we would like to verify that something is there, but are not particularly concerned with its specific properties or boundaries. This presentation will discuss ways in which we modify the full-waveform inversion problem to answer this different question efficiently. It will include discussion of ways to parameterize the $CO_2$ plume to allow for fast uncertainty quantification as well as ways to plan surveys to place instruments to answer our question of interest with as little data acquisition cost as possible.

## 1 Introduction

Full-waveform inversion (FWI) has been well-studied [1] to optimize the recovery of the Earth's velocity field pixel-by-pixel throughout a subsurface model. This is ideal for many applications, including reservoir characterization as well as understanding tectonics and other Earth processes. For $CO_2$ sequestration however, the goal is different in that we only want to ensure that $CO_2$ remains in place but do not need a pixel-by-pixel image of exactly where it is.

We will discuss two ways to adapt FWI to these requirements. The first is to exploit interrogation theory [2, 3], which allows one to determine the answer to a question about the data, rather than to form a pixel-by-pixel image. This theory works within a framework similar to that of uncertainty quantification, relying on sampling many models to determine which data best constrain the result of interest.

A key challenge in uncertainty quantification for FWI is in the number of degrees of freedom. One way to address this is to incorporate methods like Stein-Variational Gradient Descent (SVGD) or Normalizing Flows into the inverse problem [4]. Another more direct method is to parameterize the key parts of the model using only a few degrees of freedom. This latter method is well-suited for $CO_2$ sequestration where we have reliable a priori information collected before the injection of $CO_2$.

Our goal is to setup the tools to combine interrogation theory and a low-parameterization method for FWI to determine which data we need to constrain subsurface $CO_2$.

## 2 Model

Our model is the acoustic wave equation

$$\frac{\partial^2 u}{\partial t^2} - c^2 \nabla^2 u = f \tag{1}$$

where $u$ is the pressure, $c$ is the velocity, and $f$ is the source. Within our model, shown in Figure 1, we assume that we have a known background velocity, $c_0$, and that there is an unknown perturbation $\Delta c$ that represents a plume of $CO_2$. This models 4D seismic imaging in which an initial dataset is recorded before injection of $CO_2$ into the subsurface. This initial survey is used to determine the background structure of the subsurface. A second (or more) dataset is then acquired after the injection of $CO_2$ – it is this second dataset that we use to determine the location and outline of the plume.

## 3 Methods and Results

### 3.1 Measurement Location Optimization

Optimizing the location of receivers for seismic imaging is done using interrogation theory [2], which allows us to answer a question of interest (e.g. what is the volume of subsurface $CO_2$?) and to determine what data are necessary to answer this question. The basic idea is to maximize $U$ defined by ( [2, eq. 32]

$$U(d) = E\left[a(y_d, d)^2 | d\right] \tag{2}$$

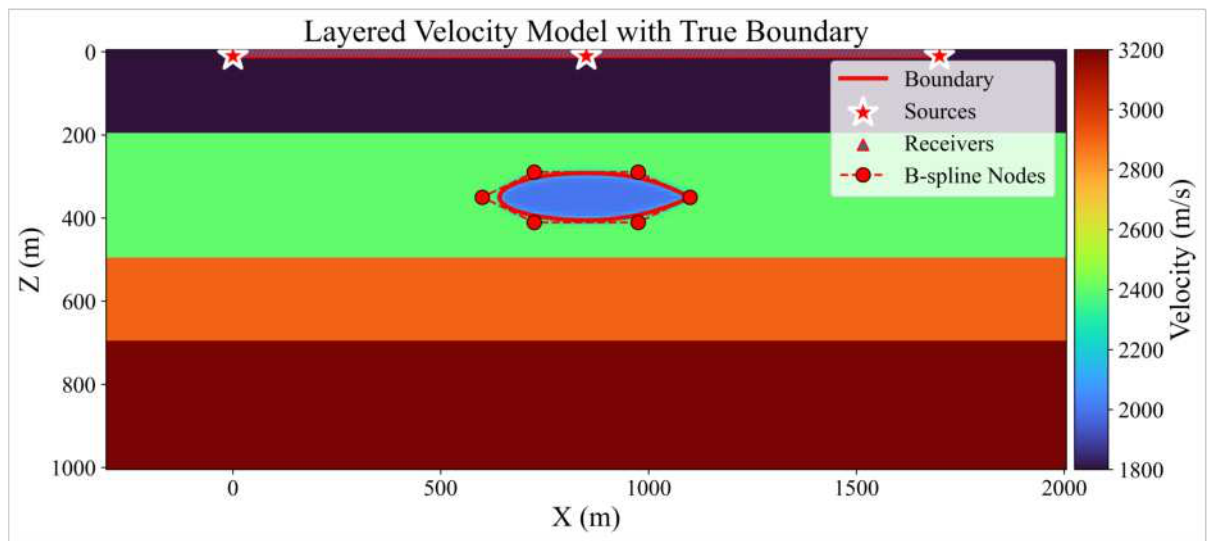

Figure 1: Velocity model. We assume the background (here layered) model is known and are trying to recover the location of the $CO_2$ plume.

where $d$ is the experimental design (set of observation locations), $a$ is the answer to our question and $y_d$ is observed data within our experimental design, $d$. The computation of $U$ requires the integration of the likelihood, $p(m|y_d, d)$. We have tested this with gravity modelling, where the solution can be determined analytically. Work to compute or approximate this likelihood for FWI is ongoing, as described below.

### 3.2 Shape Inversion

There are many ways to estimate the likelihood in seismic inversion. Here we choose to use normalizing flows [4–6] because they run quickly and are not limited to small numbers of degrees of freedom. To step towards answering our questions, rather than simply generating a pixel-by-pixel image, we parameterize the model with six nodes, as illustrated in Figure 1 recalling that we are interesting in recovering the blue inclusion within the second layer. With this parameterization, we estimate the posterior for the location of each node, from which we can calculate the uncertainty on the boundary, as shown in Figure 2. We see that with a full dataset, we can recover the shape adequately (the wavelength, $\lambda \approx 85$ m) and that the quality of the recovery varies with position as expected.

Tying together the measurement location optimization and the recovery of the shape will allow us to estimate which data are most important to constrain the shape of the $CO_2$ plume without requiring a pixel-by-pixel reconstruction.

**Acknowledgements:** Funding: Total Energies, NSERC Alliance and Discovery Grant programs. Compute: Digital Research Alliance of Canada. We thank Paul Barbone of Boston University and Yen Sun of Total for useful discussions.

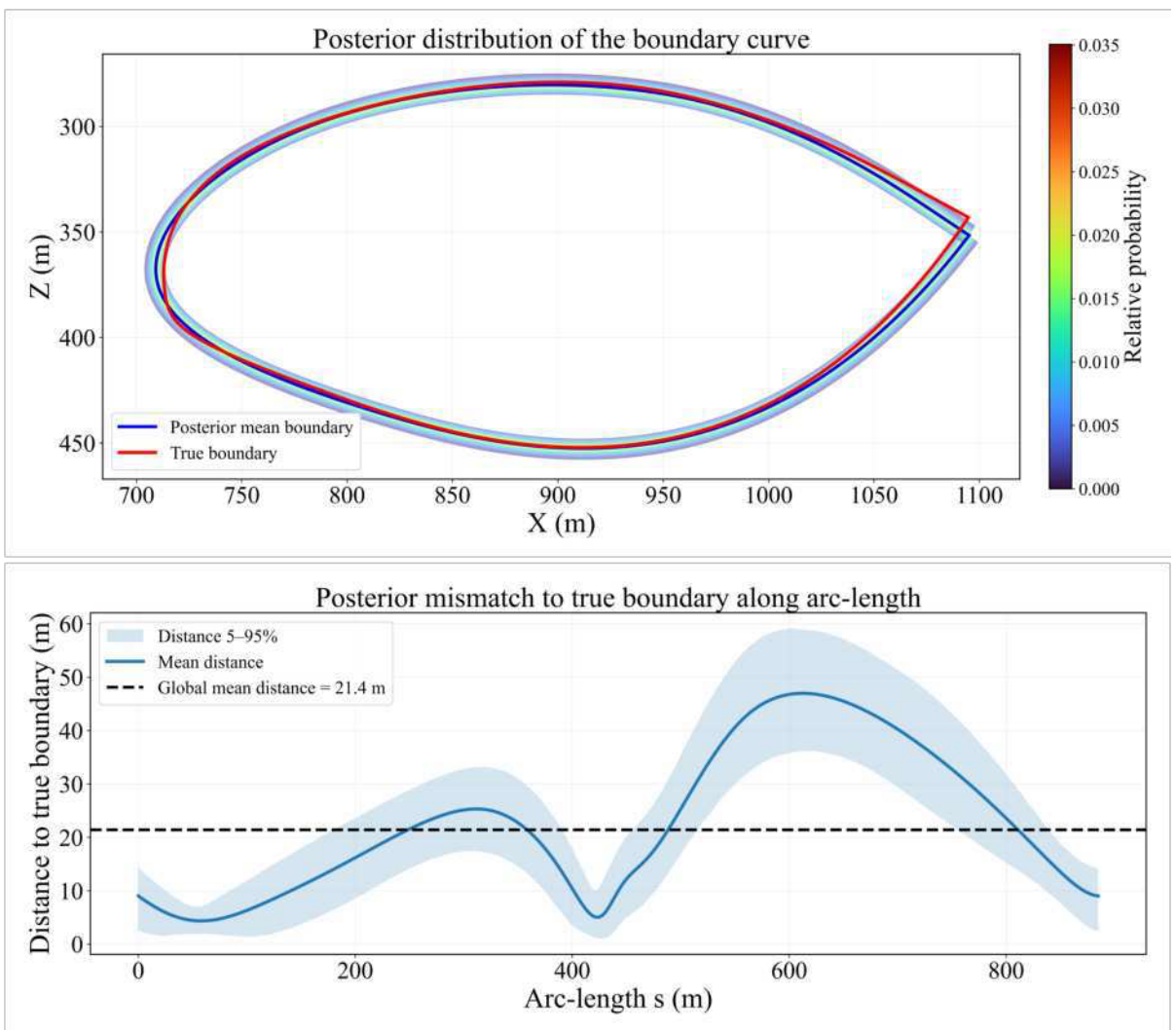

Figure 2: Top: Recovered shape of the plume using a parameterization of six nodes and normalizing flows to recover the posterior distribution. Bottom: Error in the recovered shape.

# A Hermite-Taylor Compact Discrete Correction Function Method for Maxwell's Equations

Seth Malonson[1,*], Yann-Meing Law[1], Jean-Christophe Nave[2]
[1]Department of Mathematics and Industrial Engineering, Polytechnique Montreal, Montreal, CA
[2]Department of Mathematics and Statistics, McGill University, Montreal, CA
*Email: seth.malonson@polymtl.ca

**Abstract**

High-order accurate numerical solutions to hyperbolic problems can be obtained from Hermite-Taylor methods. Function values and derivatives through order $m$ of each variable are evolved in time to attain a $2m + 1$ order convergence. This large data requirement creates a challenge in handling boundary conditions such that the accuracy of the Hermite-Taylor method is preserved. The discrete correction function method (DCFM) overcomes this challenge by seeking a polynomial approximation for each field in the vicinity of the boundary through a least squares approach. Here, we propose a modified DCFM which leverages all available data evolved by the Hermite-Taylor solution, that is, field values and all derivatives through order $m$. This framework results in a compact stencil for the DCFM, enabling larger values of $m$ and an order-adaptive algorithm. Numerical examples in 1-D suggest that the modified DCFM is stable for $m = 1-6$ and the expected rate of convergence is observed for these values of $m$, thus providing up to a thirteenth-order accurate scheme.

## 1 Introduction

Consider a domain $\Omega = [x_l, x_r]$, a time interval $I = [t_0, t_f]$, and a 1-D simplification of Maxwell's equations using the $x$-directed, $z$-polarized transverse electromagnetic (TEM) mode. We are seeking a numerical solution to

$$\begin{aligned} \mu\partial_t H &= \partial_x E \quad \text{in } \Omega \times I, \\ \varepsilon\partial_t E &= \partial_x H \quad \text{in } \Omega \times I, \\ H(x, 0) &= H_0(x) \quad \text{for } x \in \Omega, \\ E(x, 0) &= E_0(x) \quad \text{for } x \in \Omega, \\ E(x_l, t) &= E(x_r, t) = 0 \quad \text{for } t \in I, \end{aligned}$$

where $H$ is the magnetic field, $E$ is the electric field, $\mu$ is the magnetic permeability, $\varepsilon$ is the electrical permittivity, and $H_0(x)$ and $E_0(x)$ are known functions.

## 2 Hermite-Taylor

Hermite methods require two staggered meshes, a primal mesh and a dual mesh. The primal mesh consists of $N + 1$ evenly spaced nodes located at $(x_i, t_n) = (x_l + i\Delta x, t_0 + n\Delta t)$ for $i = 0, \ldots, N_x$ and $n = 0, \ldots, N_t$, with $\Delta x = (x_r - x_l)/N_x$ and $\Delta t = (t_f - t_0)/N_t$. Here, $N_x$ is the number of cells into which the primal mesh divides our domain, and $N_t$ is the number of time steps required to reach $t_f$. For the dual mesh, the nodes are centered between each adjacent pair of primal nodes, located at $(x_{i+1/2}, t_{n+1/2}) = (x_l + (i + 1/2)\Delta x, t_0 + (n + 1/2)\Delta t)$ for $i = 0, \ldots, N_x - 1$ and $n = 0, \ldots, N_t - 1$.

For each field and for each node of the primal mesh, we assume to know the function and derivative values up to order $m$ at time $t_n$. On each cell, we construct the unique Hermite polynomial in $\mathbb{P}_{2m+1}$ which matches the known function and derivative values at each endpoint. We center this Hermite interpolant at the center of the cell. We then expand each coefficient of the Hermite polynomial in time and use a recursion relation derived from Maxwell's equations to determine the additional coefficients from the expansion in time. Finally, we evaluate the resulting space-time polynomials and their derivatives through order $m$ at the dual mesh nodes at time $t_{n+1/2}$. We repeat this process to evolve the data from the dual mesh at $t_{n+1/2}$ to the primal mesh at $t_{n+1}$, thus completing a time step. The Hermite-Taylor method is known to be $2m + 1$ order accurate for linear hyperbolic problems with periodic boundary conditions [1].

## 3 Compact DCFM

Imposing boundary conditions in the Hermite-Taylor method is inherently difficult due to the necessity to evolve field values and their first $m$ derivatives at nodes immediately adjacent to the boundary. The DCFM overcomes this challenge by seeking a polynomial approximation of each field in the vicinity of the boundary. To do so, for all of the nodes where the Hermite method cannot be applied, referred to as correction function (CF) nodes, we create a local patch of length $\beta\Delta x$ in space, where $\beta$ is a positive constant. The time domain of the local patch is defined as $[t_{n-1/2}, t_n]$. On this patch, for each field, we seek a space-time polynomial of degree $d$ in each dimension, centered in space at the CF node and in time at $t_n$. This results in $2(d+1)^2$ coefficients to determine in 1-D. The coefficients of the polynomial approximations of the electromagnetic fields are computed by minimizing an overdetermined linear system, formed by equations that enforce Maxwell's equations, the Hermite-Taylor solution at nodes within the local patch, and the boundary conditions [2]. Finally, knowing the polynomial coefficients, we evaluate the correction functions and their derivatives through order $m$ at the CF nodes.

In addition to constraining the correction functions to match the field values of the Hermite solution, the modified DCFM also constrains their derivatives through order $m$, leveraging all available data from the Hermite-Taylor method. By including these additional constraints, we are able to reduce the length of each local patch to $\Delta x$ by choosing $\beta = 1$. Furthermore, this enables an order-adaptive algorithm for the proposed compact DCFM by dynamically changing the polynomial degree of the correction functions according to the maximum value of $m$ carried at each node within the local patch. Given a tolerance, we truncate the correction function polynomials at the degree where all higher derivative values are smaller than the tolerance in absolute value, as it is done for $p$-adaptive Hermite methods [3]. Current work is ongoing to extend the compact DCFM to 2-D.

## 4 Numerical Examples

We consider a domain $\Omega = [0,1]$ and a time interval $I = [0,100]$. We set $\mu = 1$, $\varepsilon = 1$, and $\Delta t = 0.9\Delta x$. The initial conditions are $E_0(x) = H_0(x) = \exp((x-x_0)^2/s^2)$, where $x_0 = 0.5, s = \sqrt{0.0005}$. Fig. 1 shows convergence plots for the electromagnetic fields at time $T = 100$. We observe the expected $2m+1$ rate of convergence. Fig. 2 illustrates the value of $m$ at each primal node and the numerical solution $H$ at time $T = 2$.

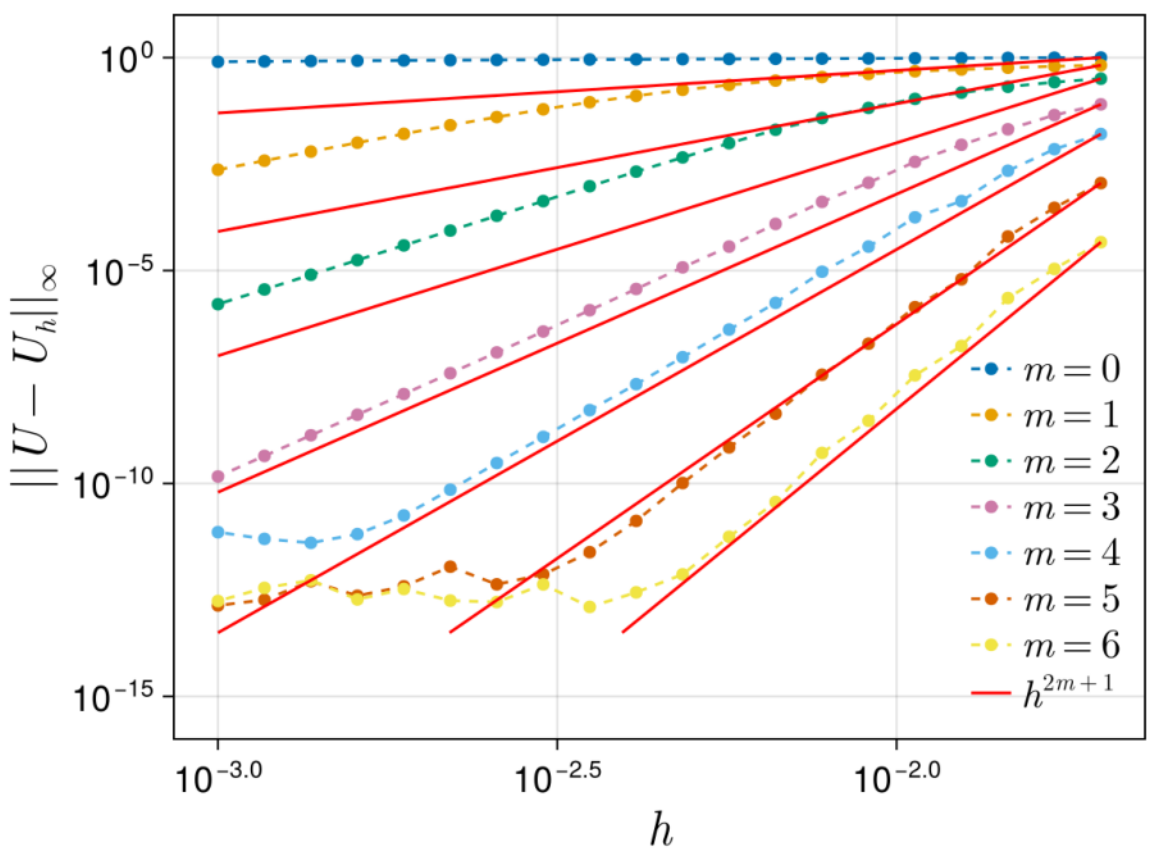

Figure 1: Convergence plot in the maximum norm at $T = 100$. Here, $U = [E, H]$.

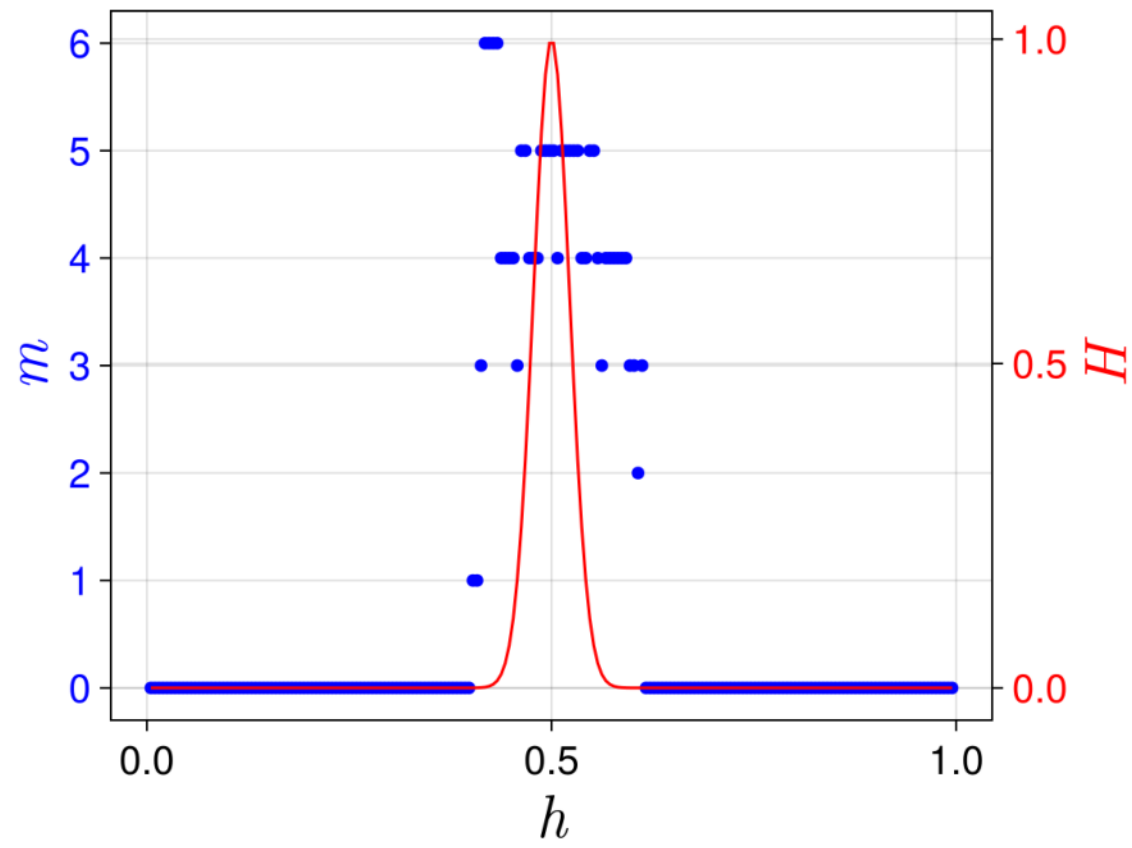

Figure 2: Values of $m$ and the magnetic field $H$ computed at $T = 2$.

# A stochastic boundary element method for electromagnetic scattering in random media

**Titouan Marquaille**[1,*], **Thomas Bonnafont**[1,2], **Arnaud Coatanhay**[1]
[1]Lab-STICC, UMR CNRS 6285, ENSTA, Institut Polytechnique de Paris, 29806 Brest, France
[2]DGA, 60 Bd. du général Martial Valin, Paris 75509, France
*Email: titouan.marquaille@ensta.fr

## Abstract

The boundary element method (BEM) is standard for wave scattering, yet classical deterministic formulations neglect critical medium uncertainties. While Monte Carlo methods can quantify these effects, their slow statistical convergence leads to prohibitive costs. We propose a Stochastic BEM based on non-intrusive polynomial chaos expansion (PCE), which transforms the uncertainty quantification problem into independent deterministic BEM solves, enabling rapid computation of statistical moments with spectral convergence.

## 1 Introduction

Electromagnetic scattering problems arise in antenna design or radar cross-section prediction where accurate field predictions are critical for system performance. The BEM is the standard numerical tool for exterior problems, as it discretizes only the scatterer's surface and naturally satisfies the Sommerfeld radiation condition at infinity.

However, deterministic models rely on idealized assumptions that neglect real-world variability. While geometric uncertainties (shape variations, surface roughness) have been extensively studied [3], material uncertainties, arising from manufacturing tolerances or environmental variations—pose distinct challenges, as they directly modify the governing wave equation. Variations in permittivity and permeability directly alter the effective wave number, shifting resonance frequencies and scattering signatures. Although Monte Carlo (MC) methods can quantify these effects, its slow convergence rate ($\mathcal{O}(N^{-1/2})$) renders it unsuitable for high-fidelity solvers.

To overcome this bottleneck, we employ Polynomial Chaos Expansion (PCE). This method approximates the random output as a series expansion in orthogonal stochastic polynomials, achieving spectral convergence for smooth parameter dependencies. Its non-intrusive formulation relies on regression and requires no modification of the deterministic BEM solver.

## 2 Problem Formulation

We consider the scattering of a time-harmonic ($e^{j\omega t}$) transverse magnetic (TM) plane wave by a perfectly conducting (PEC) obstacle $\Omega \subset \mathbb{R}^2$ embedded in a stochastic background medium. Let $(\Theta, \mathcal{F}, \mathbb{P})$ be a probability space. The background relative permittivity $\varepsilon_r(\theta)$ and permeability $\mu_r(\theta)$ are modeled as independent random variables. Consequently, the total scalar electric field $u(\mathbf{r}, \theta)$ satisfies :

$$\Delta u(\mathbf{r}, \theta) + k^2(\theta) u(\mathbf{r}, \theta) = 0 \quad \text{in } \mathbb{R}^2 \setminus \overline{\Omega}, \quad (1)$$

where $k(\theta) = k_0\sqrt{\varepsilon_r(\theta)\mu_r(\theta)}$ is the random wave number. The field corresponds $u = u^{\text{inc}} + u^{\text{sc}}$, where $u^{\text{inc}}(\mathbf{r}, \theta) = E_0 e^{jk(\theta)\mathbf{d}\cdot\mathbf{r}}$, with $\mathbf{d}$ the unit propagation vector. The scattered field $u^{\text{sc}}$ satisfies the Sommerfeld radiation condition.

On the boundary $\Gamma = \partial\Omega$, the PEC condition $u = 0$ leads to the electric field integral equation (EFIE). For a realization $\theta$, the induced surface current $J$ satisfies:

$$\frac{j\omega\mu_r(\theta)}{4} \int_\Gamma H_0^{(1)}(k(\theta)|\mathbf{r} - \mathbf{s}|) J(\mathbf{s}, \theta)\, d\Gamma(\mathbf{s}) = -\, u^{\text{inc}}(\mathbf{r}, \theta), \quad \forall \mathbf{r} \in \Gamma, \quad (2)$$

where $H_0^{(1)}(z)$ denotes the zeroth-order Hankel function of the first kind.

Applying the BEM with pulse basis functions and point matching yields the linear system $\mathbf{Z}(\theta)\mathbf{J}(\theta) = \mathbf{V}(\theta)$, where the impedance matrix $\mathbf{Z}$ depends on the random parameter $\theta$. The quantity of interest is the complex surface current density $\mathbf{J}(\theta) \in \mathbb{C}^N$, essential for deriving the RCS through near-to-far-field transformation.

## 3 Stochastic Framework (PCE)

We place ourselves within the Hilbert space $L^2(\Theta, \mathcal{F}, \mathbb{P})$ and define the vector of uncertain inputs as $\mathbf{X}(\theta) = [\varepsilon_r(\theta), \mu_r(\theta)]^\top$ modeled as log-normal variables.

After an isoprobabilistic transformation of the log-normal inputs into a standard Gaussian space, a Hermite polynomial chaos expansion is used to represent the stochastic solution $J \in L^2(\Theta)$. We approximate $J$ by a truncated expansion of order $p$:

$$J(\mathbf{r}, \theta) \approx \sum_{\boldsymbol{\alpha} \in \mathcal{A}_p} J_{\boldsymbol{\alpha}}(\mathbf{r}) \Psi_{\boldsymbol{\alpha}}(\boldsymbol{\xi}(\theta)), \tag{3}$$

where $\mathcal{A}_p = \{\boldsymbol{\alpha} \in \mathbb{N}^2 : |\boldsymbol{\alpha}| \leq p\}$. The family $\{\Psi_{\boldsymbol{\alpha}}\}$ forms an orthogonal basis of multivariate polynomials [2]. To accommodate the complex-valued nature of the surface current and material parameters within the standard PCE framework, the solution $J$ is decomposed into real and imaginary components, which are approximated as independent real-valued fields over the same polynomial basis.

Directly computing the expansion coefficients for each degree of freedom of the BEM mesh is computationally expensive. To remedy this, we employ proper orthogonal decomposition (POD) [1] to identify a reduced spatial basis $\{\phi_k\}_{k=1}^K$ that optimally captures the energy of the system. The stochastic solution is then approximated as a linear combination of these few deterministic modes:

$$J(\mathbf{r}, \theta) \approx \bar{J}(\mathbf{r}) + \sum_{k=1}^{K} \beta_k(\theta) \phi_k(\mathbf{r}), \tag{4}$$

where $\bar{J}$ is the mean field and $\beta_k(\theta)$ are the stochastic reduced coordinates. Consequently, the PCE regression is performed solely on the $K$ scalar coefficients $\beta_k$ rather than the full discretized field, drastically reducing the computational burden while filtering spatial noise.

## 4 Numerical Results

The framework is validated on a PEC circular cylinder ($R = 0.5\,\text{m}$) illuminated at $f = 300\,\text{MHz}$. The background parameters have mean values $\mu_\varepsilon = 2.5$ and $\mu_\mu = 1$.

To assess the robustness of the method, we test three distinct coefficients of variation $CV = \frac{\text{std}}{\text{mean}} \in \{5\%, 20\%, 80\%\}$. The Monte Carlo reference is computed using $10^6$ realizations. A key challenge is the nonlinear dependence of stochastic parameters in $H_0^{(1)}(k(\theta)r)$, which may affect convergence at high coefficient of variation. The error analysis for the complex surface current $\mathbf{J}$ is summarized in Table 1.

| Moment | $L^2$ Error (%) | Max Abs. (%) |
|---|---|---|
| **CV = 5%** (Speed-up ×1833) | | |
| Mean | $1.11 \times 10^{-6}$ | $2.64 \times 10^{-6}$ |
| Variance | $2.87 \times 10^{-6}$ | $1.67 \times 10^{-5}$ |
| **CV = 20%** (Speed-up ×457) | | |
| Mean | $9.47 \times 10^{-5}$ | $2.32 \times 10^{-4}$ |
| Variance | $4.58 \times 10^{-4}$ | $3.15 \times 10^{-3}$ |
| **CV = 80%** (Speed-up ×77) | | |
| Mean | $3.49 \times 10^{-3}$ | $9.42 \times 10^{-3}$ |
| Variance | $8.66 \times 10^{-2}$ | $1.03 \times 10^{-1}$ |

Table 1: Relative $L^2$ and maximum absolute errors (in %) of the PCE statistical moments compared to the Monte-Carlo reference.

The results show that the proposed PCE–POD framework remains accurate despite the strong kernel non-linearity. Although the regression error naturally scales with the input dispersion, precision remains excellent: even at $CV = 80\%$, the variance error is negligible ($< 0.01\%$). This accuracy is obtained with speed-up factors ranging from $\approx 1833$ to over 77 compared to the Monte Carlo reference. Consequently, the method offers a highly favorable cost-accuracy trade-off, validating its robustness for stochastic scattering problems involving large parameter variations.

# A Discontinuous Petrov–Galerkin Framework for Time-Harmonic Maxwell Equations in Microwave-Heated Flows

**Oreste Marquis**[1,*], **Bruno Blais**[1], **Matthias Maier**[2]
[1]Department of Chemical Engineering, Polytechnique Montréal, Montreal, Canada
[2] Department of Mathematics Texas A&M University, College Station, USA
*Email: oreste.marquis@etud.polymtl.ca

**Abstract**

In this work, we present an open-source Discontinuous Petrov–Galerkin (DPG) solver for the time-harmonic Maxwell equations, implemented using the finite element library `deal.II` and coupled with heat transfer and incompressible flow solvers to simulate microwave-heated flows.

## 1 Introduction

Microwave-assisted reactors are widely used in chemical, pharmaceutical, and materials processing to overcome conventional heat and mass-transfer limitations, thereby enhancing yield, selectivity, and energy efficiency. However, predictive simulation tools for their design remain limited due to the mathematical challenges posed by the time-harmonic Maxwell equations. While direct solvers can handle these indefinite systems, they scale as $\mathcal{O}(n^2)$ in memory, making them impractical for realistic three-dimensional reactor geometries, which require $\mathcal{O}(10^5)$ cells. To address this, we develop a DPG framework [1] for the time-harmonic Maxwell equations and couple it with heat-transfer and incompressible Navier–Stokes solvers within the open-source software `Lethe`. The DPG methodology always delivers a symmetric positive-definite system matrix, enabling the use of scalable iterative solvers despite the indefiniteness of the underlying operator.

## 2 Method

The electromagnetic problem is treated as quasi-static: the characteristic propagation time scale is orders of magnitude smaller than the heating time scale,

$$\tau_{\mathrm{em}} \sim L\sqrt{\varepsilon\mu} \ll \tau_{\mathrm{heating}} \sim \frac{\rho C_p \Delta T}{Q_{\mathrm{em}}}, \qquad (1)$$

where $Q_{\mathrm{em}}$ denotes the volumetric electromagnetic power absorbed by the medium. Consequently, the electromagnetic fields are computed in the frequency domain with frozen material properties and updated only when temperature-dependent properties change beyond a prescribed threshold. The cycle-averaged electromagnetic power dissipation,

$$Q_{\mathrm{em}} = \frac{1}{2}\sigma|\mathbf{E}|^2 + \frac{1}{2}\omega\varepsilon_0\varepsilon''|\mathbf{E}|^2 + \frac{1}{2}\omega\mu_0\mu''|\mathbf{H}|^2, \qquad (2)$$

where $|\cdot|^2$ denotes the squared modulus of the complex field amplitudes, enters the transient heat equation as a source term,

$$\rho C_p\left(\frac{\partial T}{\partial t} + \mathbf{u}\cdot\nabla T\right) = \nabla\cdot(\kappa\nabla T) + Q_{\mathrm{em}}, \qquad (3)$$

which is solved in time together with the incompressible Navier–Stokes equations. The resulting temperature field updates the material properties for the next electromagnetic solve, thus enabling fully coupled multiphysics simulations.

To solve the time-harmonic Maxwell equations, we use the DPG method. It is based on the idea that the test space $V$ can be chosen in a way that guarantees uniform and unconditional discrete stability. In practice, the optimal test space is approximated by a broken enriched space whose polynomial degree exceeds that of the trial space $U_h$. Because the broken test functions lack inter-element continuity, additional interface unknowns in the trace space $\hat{U}_h$ are introduced to enforce conformity weakly. The abstract DPG formulation on a mesh $\Omega_h$ reads:

$$\begin{cases} \text{Find } u_h \in U_h, \hat{u}_h \in \hat{U}_h \text{ and } \Psi \in V_h \text{ so that} \\ (v,\Psi)_V + b_h(v,u_h) + \langle v,\hat{u}_h\rangle = l(v), \forall v \in V_h, \\ b_h(w_h,\Psi) = 0, \quad \forall w_h \in U_h, \\ \langle \hat{w}_h,\Psi\rangle = 0, \quad \forall \hat{w}_h \in \hat{U}_h, \end{cases}$$

where $(\cdot,\cdot)_V$ denotes the sesquilinear inner product inducing the test-space norm and $\Psi$ is the Riesz representative of the residual, which serves as a computable a posteriori error estimator.

Applying the method in the electromagnetic frequency-domain setting under the $e^{-i\omega t}$ convention, and adopting an ultraweak variational formulation in which both Ampère–Maxwell's and Faraday's laws are weakened using broken test functions $\mathbf{F}, \mathbf{I} \in \prod_{K \in \Omega_h} H(\text{curl}, K)$, leads to the following formulation:

$$(\nabla \times \mathbf{F}, \mathbf{H})_{\Omega_h} + (\mathbf{F}, i\omega\varepsilon_r \mathbf{E})_{\Omega_h} + \langle \mathbf{F}, \mathbf{n} \times \hat{\mathbf{H}} \rangle_{\partial\Omega_h} = (\mathbf{F}, \mathbf{J})_{\Omega_h}, \tag{4}$$

$$(\nabla \times \mathbf{I}, \mathbf{E})_{\Omega_h} - (\mathbf{I}, i\omega\mu_r \mathbf{H})_{\Omega_h} + \langle \mathbf{I}, \mathbf{n} \times \hat{\mathbf{E}} \rangle_{\partial\Omega_h} = 0, \tag{5}$$

where all inner products are conjugate-linear in the first argument. In addition to the above equations, the test-space norm is chosen to control both field variables and their residuals, including boundary contributions from impedance conditions on $\Gamma_R$:

$$\|(\mathbf{F}, \mathbf{I})\|_V^2 = \|\mathbf{F}\|^2 + \|\mathbf{I}\|^2 + \|\nabla \times \mathbf{F} - i\omega\mu_r \mathbf{I}\|^2 + \|\nabla \times \mathbf{I} + i\omega\varepsilon_r \mathbf{F}\|^2 + \|\mathbf{n} \times \mathbf{I}_\parallel + \frac{\mathbf{k} \cdot \mathbf{n}}{\omega\mu_r} \mathbf{F}_\parallel\|_{\Gamma_R}^2, \tag{6}$$

where $\mathbf{k}$ is the wave vector associated with the imposed impedance condition.

## 3 Preliminary results and outlook

The implementation is verified on two benchmark problems. The first is a rectangular waveguide (0.25×0.25×1 m) with perfectly electrically conducting walls, filled with air and excited by the fundamental $TE_{10}$ mode at 2.45 GHz—a configuration representative of lab-scale microwave-assisted reactors. The observed convergence of the error in the $L^2$ norm (Fig. 1) matches theoretical expectations confirming the correct implementation of the DPG formulation.

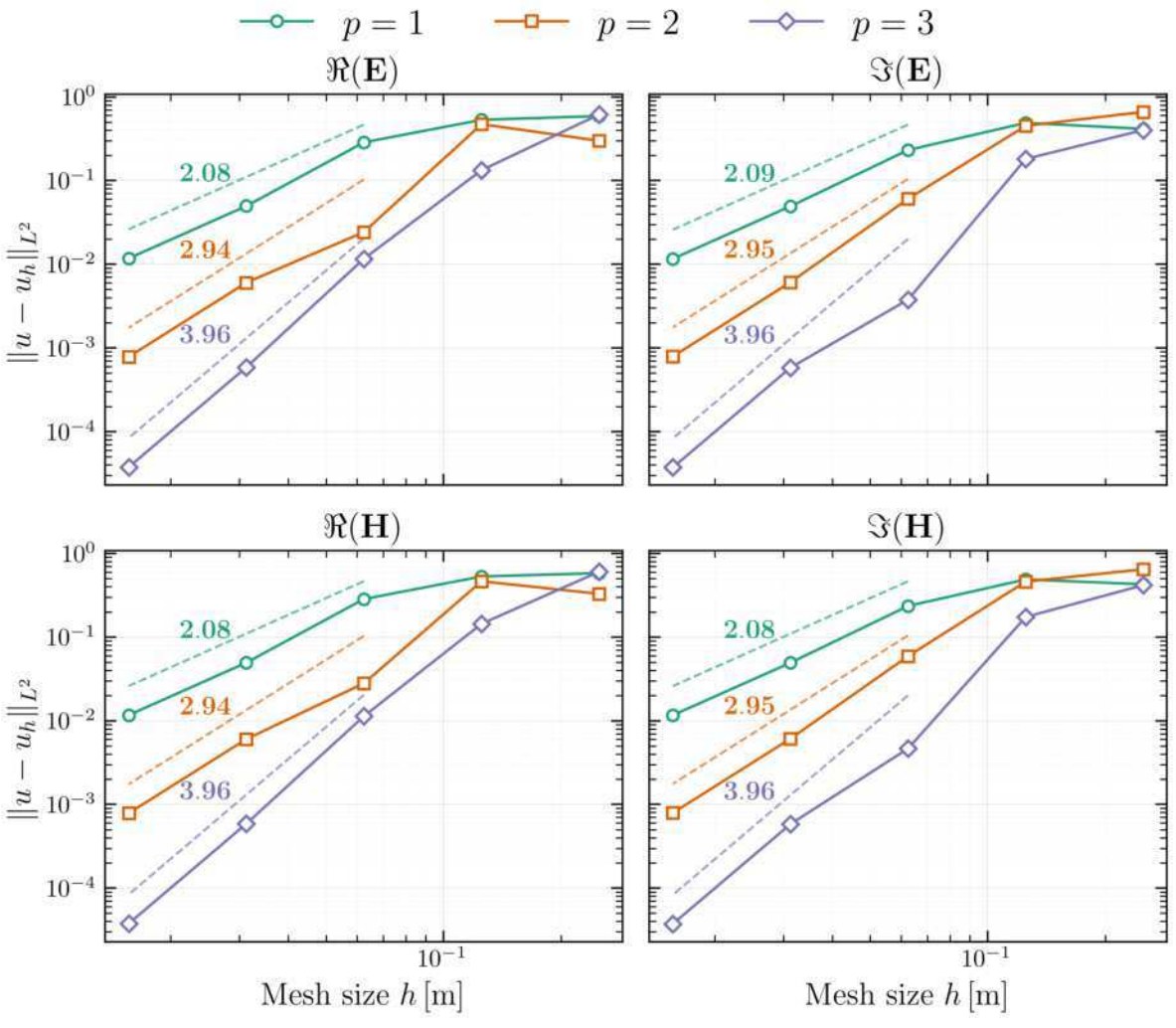

Figure 1: Convergence plot for the waveguide test case showing the expected rate of $\mathcal{O}(h^{p+1})$.

The second benchmark is the Fichera "oven" problem, a canonical test case featuring field singularities at re-entrant corners [2]. This problem is well-suited to verify the implementation of the residual-based error estimator and the $h$-adaptivity strategy. The computed solution (Fig. 2) accurately captures the singular field behaviour by refining near edges and corners and is consistent with results reported in the literature.

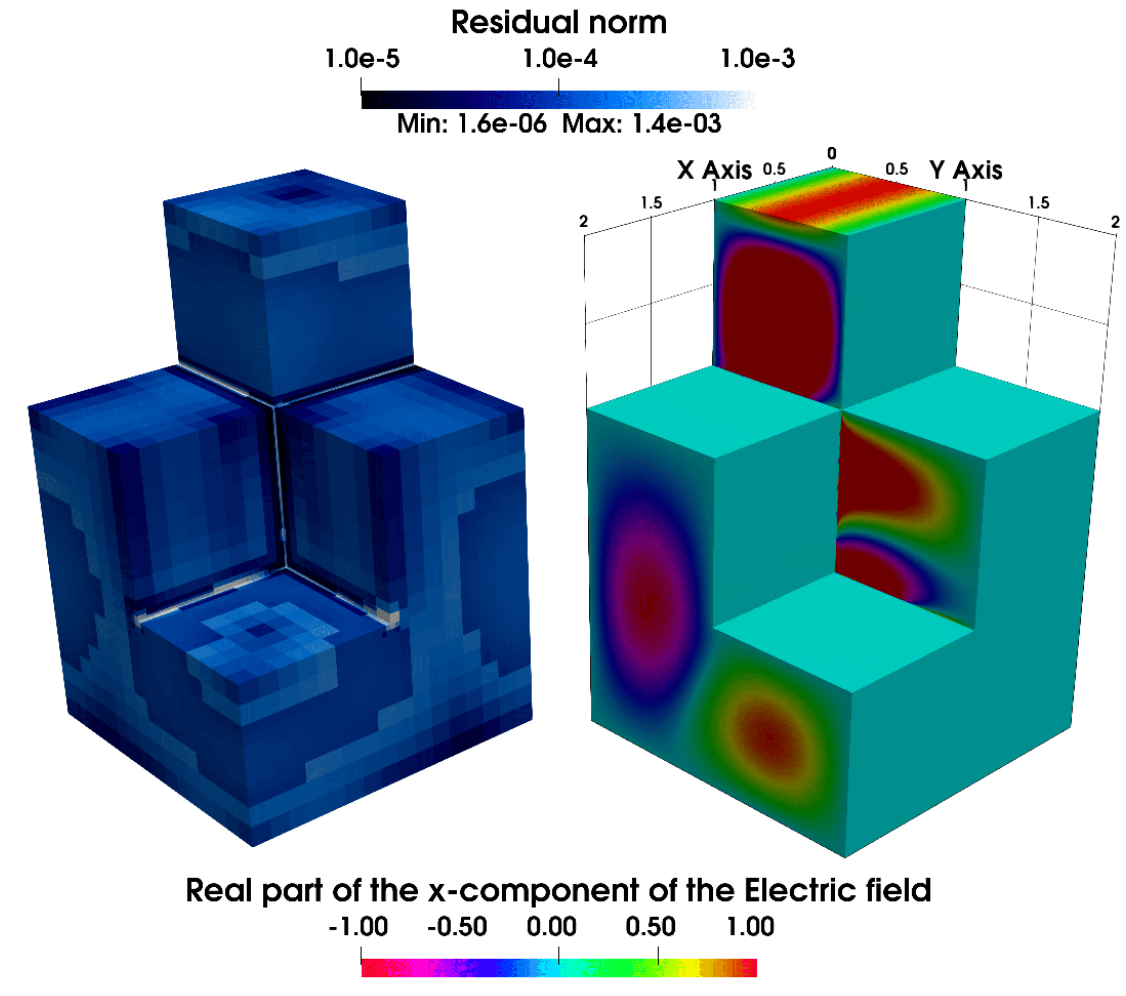

Figure 2: Fichera "oven" solution and residual norm $\|\Psi\|_V$ after nine adaptive refinements.

Ongoing work targets fully coupled multiphysics simulations. As a representative test case, we consider a lossy dielectric sample placed inside a waveguide and subjected to microwave irradiation, while a forced air flow extracts the generated heat. This configuration will enable the investigation of the nonlinear interplay between electromagnetic energy deposition, temperature-dependent material properties, and convective heat-transfer dynamics: central elements to the design of microwave-assisted reactors.

# Full Waveform Inversion Using Second Order Optimizers and Domain Decomposition with Spectral Coarse Space

**Boris Martin**[1,*], **Christophe Geuzaine**[1]

[1]Montefiore Institute, University of Liège, Belgium

*Email: boris.martin@uliege.be

## Abstract

Full Waveform Inversion uses gradient-based or second-order optimizers, the latter improving convergence at the cost of additional solutions of the underlying PDE. We explore combining second-order optimizers with recent spectral coarse grid Domain Decomposition preconditioners, whose high setup costs can be amortized over multiple solves.

## 1 Introduction

Large-scale 3D Full Waveform Inversion (FWI), i.e., the reconstruction of the heterogeneous medium's properties through full wave simulation, requires solving the underlying PDE (e.g., Helmholtz equation) repeatedly and with multiple right-hand sides (RHS). In the frequency domain, the discretized PDE takes the form of a sparse linear system, which can be solved using a sparse direct solver or iterative solvers, such as GMRES preconditioned by a Domain Decomposition Method (DDM). Iterative methods however require iterating anew for each RHS.

Different preconditioners offer different tradeoffs in terms of convergence, iteration cost and setup cost. In particular, two-level DDMs using spectral coarse spaces such as $H_k$-GenEO [1] or MS-GFEM [2] can accelerate convergence, at the cost of an elaborate setup involving Generalized Eigenvalue Problems on each subdomain. Each subdomain provides a few modes to deflate in a *coarse problem*. Depending on the optimizer chosen for FWI, the number of times a given operator must be inverted can vary; impacting whether a good preconditioner with high setup cost might beat a mediocre preconditioner with fast setup.

Our goal in this contribution is to identify synergies between optimizer choice and linear solver selection.

## 2 Model problem

We focus on the heterogeneous Helmholtz equation with impedance boundary conditions on the Overthrust 2D model (Figure 1). The goal of FWI is to recover the wave velocity $c(x)$ from data at 40 source-receiver locations.

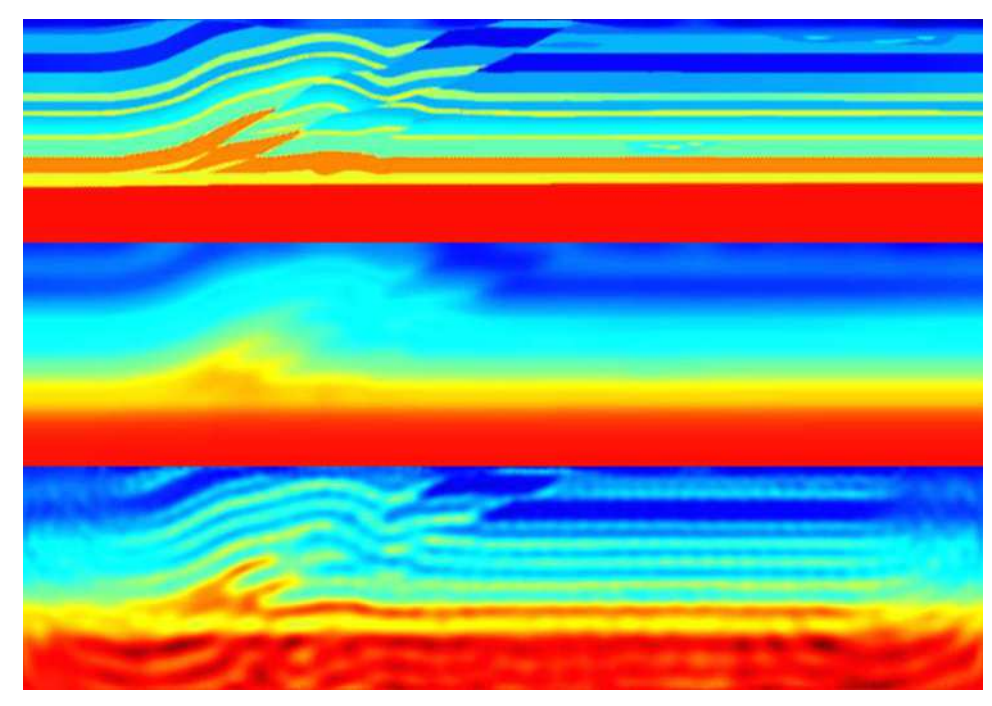

Figure 1: Exact (top), initial (middle) and reconstructed (bottom) wave velocity field $c(x)$ with Truncated Gauss-Newton.

## 3 Optimization algorithm

Let $J(m)$ be the cost function with $m \in \mathbb{R}^n$ representing degrees of freedom of $c(x)$. A second order expansion gives

$$J(m + \delta m) \approx J(m) + \nabla J \cdot \delta m + \frac{1}{2}\delta m^T H \, \delta m,$$

where $H$ is the Hessian. A Newton approach proceeds by choosing $\delta m = -H^{-1}\nabla J(m)$ and iterating until convergence. Computing, storing and inverting $H$ is however computationally prohibitive for large $n$, so approximate Newton methods are preferred. One can either build an approximation of $H^{-1}$ and iteratively refine it using previous gradients, leading to the family of Quasi-Newton methods (including L-BFGS), or, if the action of $H$ can be computed at a reasonable cost, one can approximately solve $H\delta m = -\nabla J$ using e.g. Conjugate Gradients, leading to a Truncated Newton (TN) method. The TGN variant uses the Gauss-Newton approximation of the Hessian [3], which is cheaper to evaluate and less accurate but can lead to a more stable optimization process. Table 1 illustrates the difference in computational structure between L-BFGS, TN and TGN on the Overthrust 2D model, and Figure 2 presents the resulting computational breakdown when a direct sparse solver is used (MUMPS).

| | L-BFGS | TN | TGN |
|---|---|---|---|
| PDE solves | 80 | 160 | 119 |
| Operator updates | 39 | 13 | 7 |

Table 1: Number of iterations to reduce $J$ by a factor $10^3$ (using frequencies of 1, 2, 3, 4, and 5 Hz in the FWI).

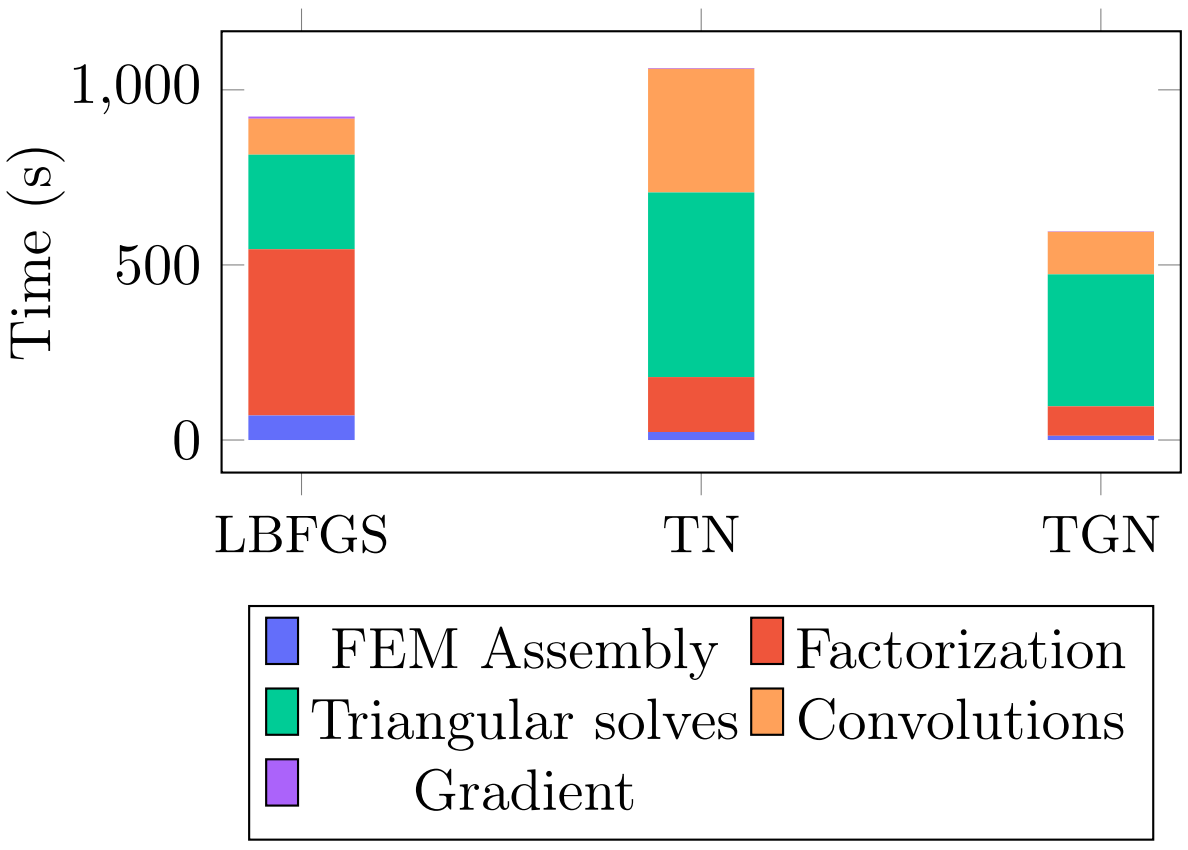

Figure 2: Computational time breakdown with a direct solver on an AMD Ryzen 7 CPU: the higher iteration count of TGN is offset by the reuse of the factorizations. Convolutions refer to adjoint-forward correlation to compute gradients and gradient perturbations, as well as computing the RHS of the perturbed problems in TN/TGN cases.

## 4 Iterative solver

While sparse direct solvers for the Helmholtz equation are a viable choice in 2D, iterative methods will become necessary for high frequency 3D cases with hundreds of millions of DoFs [4]. We compare the number of GMRES iterations on the 2D Overthrust case between a one-level ORAS DDM and a two-level ORAS DDM with MS-GFEM coarse space in Table 2.

When a few large subdomains are used, hundreds of RHS are needed to amortize the GEVP cost. However for intermediate size (e.g. 32 subdomains), building a good coarse space is already effective for only 10 RHS. For 30 modes per subdomain, we even observe superlinear convergence: from 9s on 8 domains to 0.5s on 64 domains (setup not included).

| | $n_{\text{ev}}$ | | | |
|---|---|---|---|---|
| $N_{\text{dom}}$ | 0 | 10 | 20 | 30 |
| 8 | 79 | 66 | 46 | 25 |
| 16 | 120 | 113 | 57 | 20 |
| 32 | 318 | 206 | 50 | 14 |
| 64 | 647 | 264 | 23 | 10 |
| 128 | $> 1000$ | 137 | 11 | 7 |
| 256 | $> 1000$ | 44 | 7 | 6 |

Table 2: Number of GMRES iterations as a function of the number of subdomains $N_{\text{dom}}$ for the Overthrust 2D model at 6 Hz, for single-level ORAS ($n_{\text{ev}} = 0$) and for two-level ORAS with MS-GFEM coarse space with $n_{\text{ev}} = 10$, 20 and 30. Simulation for 10 RHS.

In the full contribution we will explore the respective computational costs with multiple RHS and each optimizer/iterative solver combination, for both 2D and 3D models. Switching from L-BFGS to Truncated Newton methods is expected to increase the optimal coarse space size.

# Exceptional points and defective resonances in a sound-soft scattering problem

Kei Matsushima[1,*]

[1]Informatics and Data Science Program, Graduate School of Advanced Science and Engineering, Hiroshima University, Japan

*Email: matsushima@acs.hiroshima-u.ac.jp

**Abstract**

This study investigates non-Hermitian degeneracy and exceptional points associated with resonances in a two-dimensional scattering problem with sound-soft obstacles. Our aim is to identify non-Hermitian degenerate (defective) resonances in a sound-soft scattering setting using numerical methods based on layer potentials and boundary integral operators. For this purpose, we formulate resonances in the sound-soft system as eigenvalues of a holomorphic Fredholm operator-valued function. This allows us to define the defectiveness of resonances and associated exceptional points based on their geometric and algebraic multiplicities. Subsequently, it is shown that the mismatch of the two multiplicities leads to fractional-power sensitivity of defective resonances with respect to operator perturbations, known as exceptional points. We present some numerical evidence for the existence of a defective resonance.

## 1 Sound-soft scattering problem

Let $\Omega \subset \mathbb{R}^2$ be an open and bounded subset with smooth boundary. We assume that the exterior $\mathbb{R}^2 \setminus \overline{\Omega}$ is connected. In this study, we consider the following Dirichlet problem:

$$(P_k) \quad \begin{cases} -\Delta u - k^2 u = 0 & \text{in } \mathbb{R}^2 \setminus \overline{\Omega}, \\ u = g & \text{on } \partial\Omega, \end{cases}$$

where $g$ denotes prescribed Dirichlet data on the boundary $\partial\Omega$, and $k \in \mathbb{C}$ is a given constant.

### 1.1 Scattering resonances

Using layer potential theory, for given $g \in C(\partial\Omega)$, we see that

$$u_k(x) := \int_{\partial\Omega} \frac{\partial G_k}{\partial\nu(y)}(x-y)\varphi_k(y)\mathrm{d}s(y) \quad \text{in } \mathbb{R}^2 \setminus \overline{\Omega}$$

belongs to $C^2(\mathbb{R}^2\setminus\overline{\Omega}) \cap C(\mathbb{R}^2\setminus\Omega)$ and solves $(P_k)$, where $\varphi_k := 2(I+K(k))^{-1}g$, provided that the linear operator $I + K(k)$ is invertible on $C(\partial\Omega)$. Here, the fundamental solution $G_k$ is given by

$$G_k(x) := \frac{\mathrm{i}}{4} H_0^{(1)}(k|x|) \quad \text{in } \mathbb{R}^2,$$

and $\nu$ denotes the unit normal vector pointing outward from $\Omega$. The boundary integral operator $K(k) : C(\partial\Omega) \to C(\partial\Omega)$ is defined by

$$(K(k)\varphi)(x) = 2\int_{\partial\Omega} \frac{\partial G_k}{\partial\nu(y)}(x-y)\varphi(y)\mathrm{d}s(y).$$

See, e.g., [1] for details.

Let $\Lambda = \{z \in \mathbb{C} : \mathrm{Im}[z] < 0\}$ denote the lower half-plane. We say that $k \in \Lambda$ is a *scattering resonance* (or scattering pole) if $I + K(k)$ is not invertible.

### 1.2 Multiplicity

Using the compactness of $K(k)$ and its analyticity with respect to $k$, we obtain several important properties of scattering resonances. More precisely, the mapping $k \mapsto I + K(k)$ is a holomorphic Fredholm operator-valued function on $\Lambda$, in the sense that

1. $k \mapsto I + K(k)$ is holomorphic on $\Lambda$,
2. $I + K(k)$ is a Fredholm operator for all $k \in \Lambda$, and
3. there exists at least one $k \in \Lambda$ such that $I + K(k)$ is invertible.

From the general theory of holomorphic Fredholm operator-valued functions [2], the set of scattering resonances is discrete in $\Lambda$, and *geometric multiplicity* $m_{\mathrm{g}} := \dim N(I + K(k))$ is finite for each scattering resonance $k$, where $N$ denotes the null space.

The algebraic multiplicity satisfies a similar finiteness property. An ordered collection $\varphi_0, \varphi_1, \dots$ in $C(\partial\Omega)$ is called a *Jordan chain* associated with a scattering resonance $k$ if $\varphi_0 \neq 0$ and

$$\sum_{l=0}^{j} \frac{1}{l!}\frac{\mathrm{d}^l}{\mathrm{d}k^l}(I + K(k))\varphi_{j-l} = 0 \quad \text{for all } j = 0, 1, \dots.$$

Any Jordan chain has finite length. An element $\varphi \in C(\partial\Omega)$ is called a generalized eigenvector if it belongs to a Jordan chain. We define the *algebraic multiplicity* $m_{\mathrm{a}}$ as the dimension of the linear hull spanned by all the generalized eigenvectors. One can show that $1 \le m_{\mathrm{g}} \le m_{\mathrm{a}} < \infty$.

### 1.3 Defective scattering resonances

Our goal is to identify a *defective scattering resonance*, i.e., a resonance whose geometric and algebraic multiplicities differ ($m_{\mathrm{g}} < m_{\mathrm{a}}$). To this end, we need a characterization of defective scattering resonances that can be checked numerically. We use the following result from [3]:

**Theorem 1** *Let $k_0 \in \Lambda$ be a scattering resonance with $m_{\mathrm{g}} = 1$ and let $B : C(\partial\Omega) \to C(\partial\Omega)$ be a nonzero bounded linear operator. Then there exist open neighborhoods $V$ of $0$ in $\mathbb{C}$ and $\Lambda_0$ of $k_0$ and holomorphic functions $c_j : V \to \mathbb{C}$ ($j = 0, \dots, m_{\mathrm{a}} - 1$) such that*

1. $c_0(0) = \cdots = c_{m_{\mathrm{a}}-1}(0) = 0$,
2. *for every $\varepsilon \in V$, the polynomial*

$$k \mapsto W(\varepsilon, k) := (k - k_0)^{m_{\mathrm{a}}} + c_{m_{\mathrm{a}}-1}(\varepsilon)(k - k_0)^{m_{\mathrm{a}}-1} + \cdots + c_0(\varepsilon)$$

*has $m_{\mathrm{a}}$ zeros (counting multiplicity) in $\Lambda_0$, and*

3. *for every $\varepsilon \in V$, $k \in \Lambda_0$ satisfies $N(I + K(k) + \varepsilon B) \neq \{0\}$ if and only if $k \in \Lambda_0$ is a zero of $W(\varepsilon, \cdot)$.*

## 2 Numerical result

We present a numerical result suggesting the existence of a defective resonance. As shown in Figure 1, we take the obstacle $\Omega$ to be the union of

- a disk of radius 0.137 centered at $(0.7495, 0)$, and
- an ellipse with semi-major axis 1.2 and semi-minor axis 0.8 centered at $(-0.7495, 0)$.

A Nyström method is used to find scattering resonances numerically. See [3] for the details of the numerical scheme. We found two scattering resonances near $k = 5.0 - 2.47\mathrm{i}$. In view of Theorem 1, we introduce the perturbation $B : C(\partial\Omega) \to C(\partial\Omega)$ defined by

$$(B\varphi)(x) = \int_{\partial\Omega} \sin|x - y|\varphi(y)\mathrm{d}s(y),$$

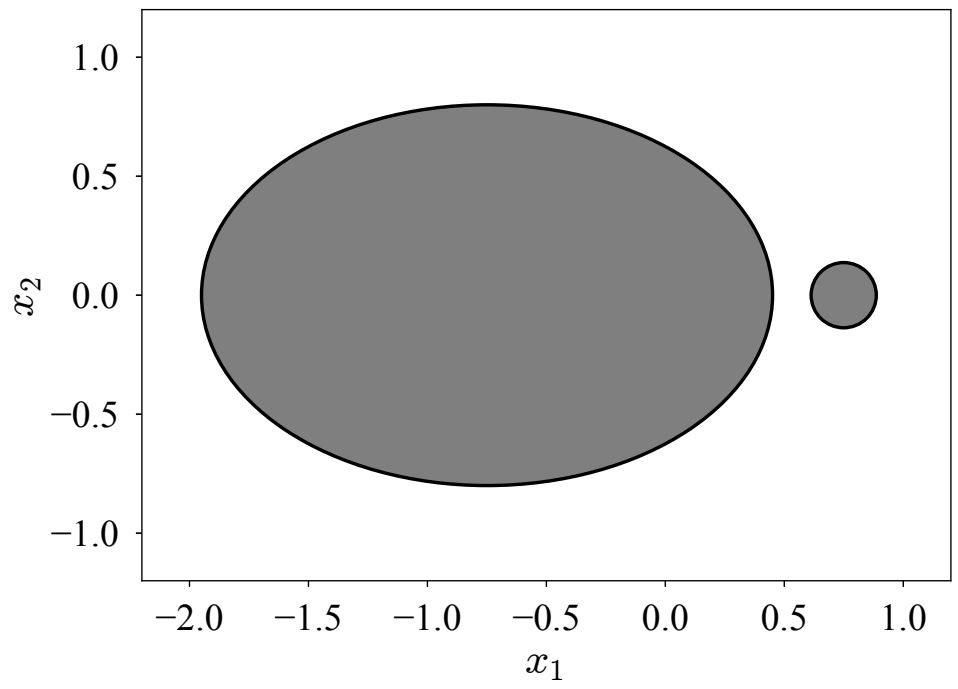

Figure 1: Obstacle $\Omega$ used in the numerical example.

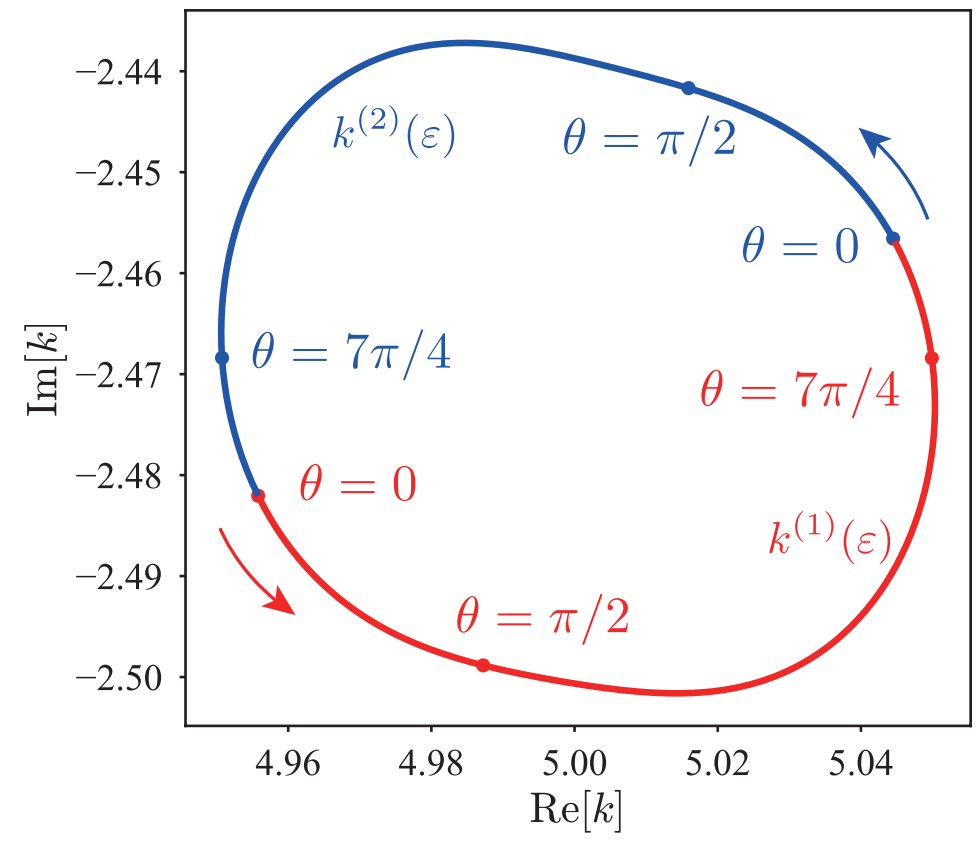

Figure 2: Trajectorie of two zeros of $k \mapsto I + K(k) + \varepsilon B$.

and numerically seek the zeros of $k \mapsto I + K(k) + \varepsilon B$. Figure 2 shows the trajectory of the two zeros $k_1(\varepsilon)$ and $k_2(\varepsilon)$ as the perturbation parameter $\varepsilon$ moves along the small circle $10^{-3}\mathrm{e}^{\mathrm{i}\theta}$ from $\theta = 0$ to $2\pi$. The result strongly suggests that $k_1(\varepsilon)$ and $k_2(\varepsilon)$ obey the fractional-power asymptotics $O(\varepsilon^{1/2})$, consistent with Theorem 1 with $m_{\mathrm{a}} = 2$.

# Interfacial Radiation Boundary Conditions for Topologically Protected Edge Modes

**Mason McCallum**[1,*], **Thomas Hagstrom**[1]
[1]Department of Mathematics, Southern Methodist University, Dallas, USA
*Email: mmccallum@smu.edu

## Abstract

We study time-domain numerical simulation of continuum models that support topologically protected edge states. Our high-order time-domain method requires nonreflecting boundary conditions, but existing analyses of complete radiation boundary conditions (CRBC) assume smooth coefficients and therefore do not apply to discontinuous media. We apply a high-order discontinuous Galerkin method coupled with CRBC and observe stable optimal convergence without any special interface treatment. Numerical evidence indicates that interface-localized modes lie in the kernel of the limiting CRBC operator even in the presence of coefficient discontinuities. These results indicate CRBC may remain valid for a broader class of hyperbolic interface problems than is currently proven.

## 1 Introduction

Topologically Invarient metamaterials (TI) have garnered high interest in engineering and physics communities. These materials support robust topologically protected intefacial waves with nonreciprocal transport. We consider continuum models for TI. Much of the existing numerical schemes have been in a band-limited frequency domain regime. We seek a full broadband time domain solution. High-order time-domain solutions of such models require careful treatment of a radiation boundary conditions. We show below that the standard CRBC formulation [3] remains sufficient in the case of the following model problem.

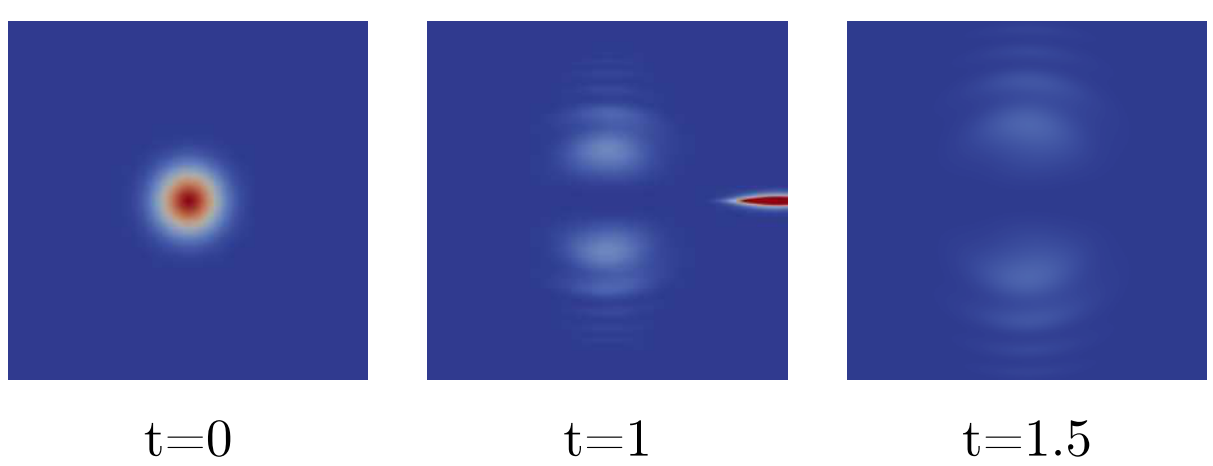

Figure 1: Evolution of gaussian initial condition along the interface

## 2 Model

As a starting point, we will consider the general model proposed by Bal [2]. In future work, we will explore other models and relax geometric assumptions. In the following, the vector field $\psi$ represents a spinor field within the material. $\sigma_{1-3}$ are the standard Pauli matrices, and m(y) is the mass potential of the medium.

$$\frac{\partial \psi}{\partial t} + \sigma_1 \frac{\partial \psi}{\partial x} + \sigma_2 \frac{\partial \psi}{\partial y} + im(y)\sigma_3 \psi = 0 \quad (1)$$

We consider the following discontinuous material parameter. $m \in \mathbb{R}^+$.

$$m(y) = \begin{cases} m, & y > 0 \\ -m, & otherwise \end{cases}$$

$m_\pm$ will be used as shorthand for the two half spaces $\pm m$. In each half space above and below the interface, we take the ansatz

$$\psi = e^{st+\lambda x+\eta_\pm y}\phi_\pm .$$

We arrive at the corresponding dispersion relation.

$$\lambda = \pm\left(s^2 + m_\pm^2 - \eta_\pm^2\right)^{1/2}$$

From the continuity of $\lambda$ accross the interface, we conjecture bound states of this interface problem lie within the of kernel the limiting CRBC boundary operator, without modification around the interface. In this work, we show numerical convergence supporting this hypothesis.

## 3 Methods and Results

There have been other attempts to solve these systems numerically, such as the integral equation formulation in [1]. Our time-domain numerical scheme is constructed to be extensible to more general interface geometries such as branching interfaces and smooth material transitions. Our mesh is aligned to the discontinuous interface. The discretization is a high order nodal Discontinuous Galerkin scheme on mass-lumped Gauss-Lobatto-Legendre elements with an upwind flux combined with a standard low storage Runge-Kutta 4/5 scheme.

### 3.1 CRBC implementation

In the following, we will consider the formulation of the CRBC layer along the east boundary without loss of generality. Let $\phi = (e \quad w)^T$ be the characteristic variables along the east face of the computational domain. $e$ will denote the outgoing characteristic variables from the bulk and $w$ the incoming, with $Q$ denoting the number of auxiliary variables normal to the boundary surface. The auxiliary variables in the radiation boundary obey the same evolution equations as the bulk problem posed on a $(d-1)$-dimensional boundary surface. This boundary system has a tensor-product structure which means the function space in the tangential direction is conforming to the bulk problem. What we lack are the appropriate normal derivatives for the function space in the normal direction. These come from a recurrence relation that specifies how the auxiliary equations are coupled. These recurrence relations are constructed such that outgoing characteristics lie within the operators null space and maintain the well-posedness of the bulk problem.

$$\mathcal{L}^+_{j-1}\phi_{j-1} = \mathcal{L}^-_j\phi_j, \quad j \in 1, 2, ..., Q$$
$$\mathcal{L}^+_{j-1} = a_{2j-1}\frac{\partial}{\partial t} + \frac{\partial}{\partial y} + s_{2j-1}$$
$$\mathcal{L}^-_j = a_{2j}\frac{\partial}{\partial t} - \frac{\partial}{\partial y} + s_{2j}$$

Without termination, this is an exact characterization of the boundary condition. In practice termination achieves exponential convergence in the number of auxiliary recursion variables. To terminate the recursion we insist that the incoming characteristic variable at level $Q$ is zero. The trace of the bulk problem is supplied as data for the boundary operator.

$$\frac{\partial e_0}{\partial x} = e^{\text{volume}}, \quad \text{(bulk crbc coupling)}$$
$$\frac{\partial w_Q}{\partial x} = 0, \quad \text{(termination conditions)}$$

The operators $\mathcal{L}^{\pm}_j$ are parameterized by scalars $a_j$ and $s_j$. These parameters are found by optimizing an expression of the numerical reflection error off a radiation boundary. What we show here is numerical evidence that across a discontinuous interface these parameters $a_j$ and $s_j$ are sufficient to acheive the expected convergence. An analytical proof is left for future work.

The convergence result below was computed for an initial Gaussian pulse centered at the interface. Complete radiation boundary conditions are implemented along all faces of the computational boundary with the appropriate corner closures as in [3]. These results were computed with ten auxiliary CRBC variables ($Q = 10$). The result below is the self-convergence of the DG scheme. A future presentation of this work will show convergence to the free space solution via long domain comparisons.

**Convergence Result**

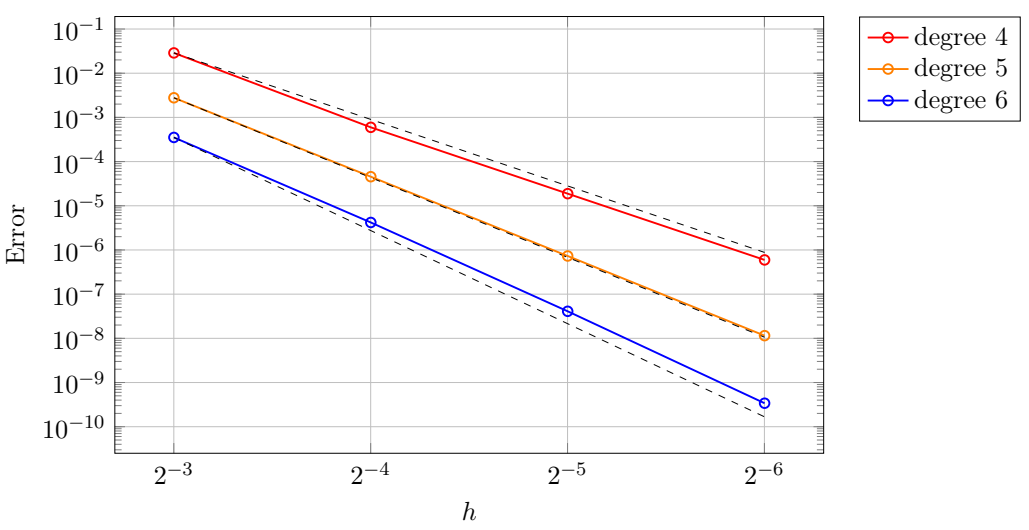

#### Acknowledgments

This material is based upon work supported by the National Science Foundation under Grants DMS-1840260 and DMS-2309687. Any opinions, findings, and conclusions or recommendations expressed in this material are those of the authors and do not necessarily reflect the views of the NSF.

# Flow rate measurement in a heterogeneous fluid with acoustic waves

**Alexander McSweeney-Davis**[1,*], **Lorenzo Audibert**[2], **Houssem Haddar**[1]

[1] Inria, Unité de Mathématiques Appliquées, ENSTA, Institut Polytechnique de Paris, 91120 Palaiseau, France

[2]Departement PRISME, EDF R&D, 6 quai Watier BP 49 Chatou, 78401 Cedex, France

*Email: alex.mcsweeney@ensta.fr

**Abstract**

This work studies the propogation of acoustic waves in a waveguide with a non-uniform flow and impenetrable obstacles in the context of modelling of ultrasonic flowmeters. The basis of our model is the Convected Helmholtz equation with obstacles and a non-uniform flow. The Fredholm nature of this problem is established, and the potential is then decomposed into an incident field (for which we have an analytical expression for a wide variety of ultrasonic flowmeters), and a diffracted field. Numerical results are presented. We also study a homogenisation of this system, where we assume particles are periodically distributed and take the limit as the size of the periodic cells goes to 0.

## 1 Introduction

Ultrasonic flowmeters are a non-intrusive non-invasive tool for measuring the flow rate of fluids through a pipe. The velocity of the fluid is found by measuring the difference in time taken for ultrasonic waves to travel upstream and downstream the flow. However, water extracted from natural bodies of water contains particles which disturb the propagation of these ultrasonic waves. Our aim is to construct a model for the ultrasonic waves which accounts for the particles and permits the rigorous study of the inverse problem. In this paper we consider a two-dimensional waveguide containing $N$ particles $\{\mathcal{O}_i\}_{i=1}^N$. We denote by $\mathcal{O} = \cup_i \mathcal{O}_i$ the collection of all the particles and $\Omega = ([-L, L] \times [0, h]) \setminus \mathcal{O}$ the fluid part of the domain. We denote by $\partial\Theta$ the exterior wall of the waveguide and $\Sigma_{\pm L}$ the outflow/inflow interfaces, see Figure 1 for an overview.

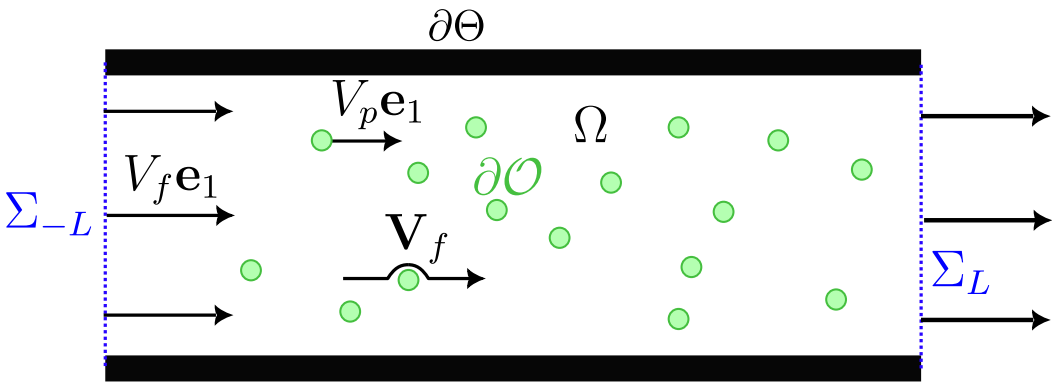

Figure 1: The two-dimensional waveguide.

## 2 Modelling

### 2.1 General model results

We model the fluid as inviscid isentropic homogeneous fluid moving with an incompressible potential flow and carrying impenetrable particles moving at a uniform speed $V_p$ in the $x$-direction. Using $V_f$ to denote the speed of the fluid in the absence of particle, we first perform a change in reference frame, to the frame moving with the particles. The relevant flow is this frame is the flow of the fluid relative to the particles, which we denote by $\mathbf{M}$. We arrive at a model with two sequential problems, a Laplace problem for the flow potential (Equation (1)) and a Convected Helmholtz problem for the propagation of the acoustic waves in the fluid (Equation (2)). First we solve Equation (1) from which we obtain the fluid Mach velocity $\mathbf{M} = \nabla u$. This is substituted into Equation (2) in which we denote by $D_k$ the convective derivative $D_k = -ik + \mathbf{M} \cdot \nabla$, where $k$ is the wave number.

$$\begin{cases} \Delta u &= 0, \text{ in } \Omega \\ \partial_{\mathbf{n}} u &= 0, \text{ on } \partial\Theta \cup \partial\mathcal{O} \\ \partial_{\mathbf{n}} u &= \frac{(V_f - V_p)}{c}\mathbf{e}_1, \text{ on } \Sigma_{\pm L} \end{cases} \tag{1}$$

$$\begin{cases} \nabla \cdot (\nabla\varphi - D_k\varphi\mathbf{M}) + ikD_k\varphi &= 0 \text{ in } \Omega \\ \partial_{\mathbf{n}}\varphi &= 0 \text{ on } \partial\mathcal{O} \\ \partial_{\mathbf{n}}\varphi &= f \text{ on } \partial\Theta \\ (\nabla\varphi - D_k\varphi\mathbf{M}) \cdot \mathbf{n} &= -T_{\pm}\varphi \text{ on } \Sigma_{\pm L} \end{cases} \tag{2}$$

The Dirichlet-to-Neumann operators $T_{\pm}$ are constructed from the homogenous problem of a fluid moving uniformly at a velocity $(V_f - V_p)\mathbf{e}_1$. We first prove that the acoustic wave problem (2) is of Fredholm type, for a subsonic incompressible flow. The proof for this is based on the Analytical Fredholm theorem and is very similar to [1], the main difficulty being to prove that the DtN maps $T_\pm$ are analytic in $k$.

**Theorem 1** *Suppose that* $\mathbf{M}$ *is an incompressible vector field tangential to the particles, pipe wall and satisfying* $\mathbf{M}\cdot\mathbf{n} = \pm(V_f - V_p)$ *on* $\Sigma_{\pm L}$ *and is uniformly subsonic* $||\mathbf{M}||_{L^\infty(\Omega)} < 1$. *Then problem* (2) *is well-posed in* $H^1(\Omega)$ *except for a countable collection of wavenumbers.*

This reduces the problem to proving that the flow model (1) yields a subsonic flow, which we prove for sufficiently small $V_f - V_p$ for circular particles. In (2) the Neumann boundary data $f$ represents the signal from the ultrasonic flowmeter. To better quantify the impact of the particles on the system, we study Equation (3), the problem without particles, where the flow is uniform with the Mach speed $\frac{V_f - V_p}{c}$ denoted by $M$.

$$\begin{cases} (1-M^2)\partial_{x_1x_1}\varphi^{\text{inc}} + \partial_{x_2x_2}\varphi^{\text{inc}} \\ +2ikM\partial_{x_1}\varphi^{\text{inc}} + k^2\varphi^{\text{inc}} = 0 \\ \partial_{\mathbf{n}}\varphi^{\text{inc}} = f \text{ on } \partial\Theta \end{cases} \quad (3)$$

Using a Green's function, we are able to compute analytically the solution to (3) for a wide range of realistic force signals. Applying a Fourier Transform allows us to obtain the corresponding time domain signal, shown in Figure 2, allowing us to measure the travel time of the wave packet and simulate an ultrasonic flowmeter reading.

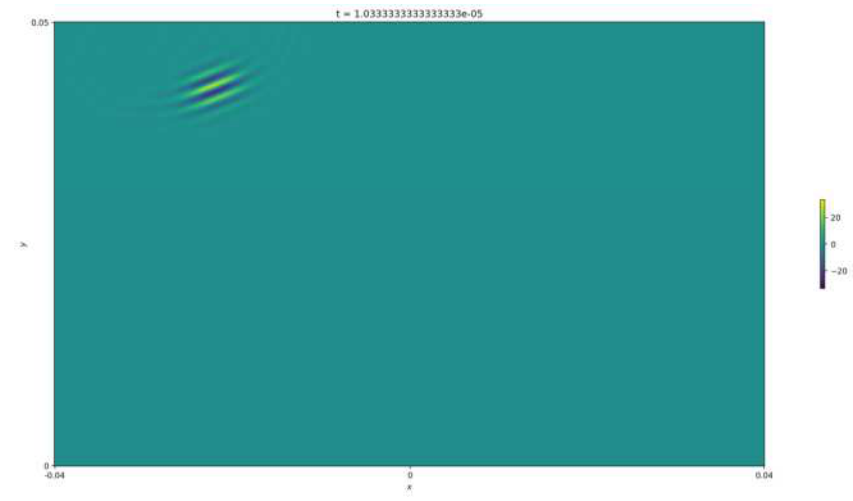

Figure 2: $\varphi^{\text{inc}}$ in the time domain.

We further demonstrate how (2) and (3) can be implemented numerically to obtain the acoustic wave for a given waveguide, particle configuration, and input signal from the ultrasonic flowmeter, shown in Figure 3.

### 2.2 Homogenisation of general model

We consider the problem with a periodic array of particles, with the size of the periodic cells being denoted by $\varepsilon$. We consider various homogenisation regimes, where the size of the particles $r(\varepsilon)$ is given by $r_0\varepsilon^m$, and the velocity of the particles $V_p(\varepsilon) = V_f - (V_f - V_p^0)\varepsilon^n$ converges to the velocity of the fluid. In this case formal calculations yield the following homogenised equation:

$$\nabla\cdot(\mathbf{A}^*\nabla\varphi_0) + 2ik\mathbf{M}^*\cdot\nabla\varphi_0 + k^2\rho^*\varphi_0 = 0 \quad (4)$$

where $\mathbf{A}^*, \mathbf{M}^*, \rho^*$ are effective tensors independent of $\varepsilon$ computed from the periodic cell problems. We approach this problem using the Two-scale convergence method [2]. An interesting scaling is $m = n = 1$ and we aim to present theoretical and numerical convergence properties of this system. In this case the limit problem is the classical homogenised problem and higher-order terms are needed to capture the impact of the particles on the flow. We will also discuss other possible scaling.

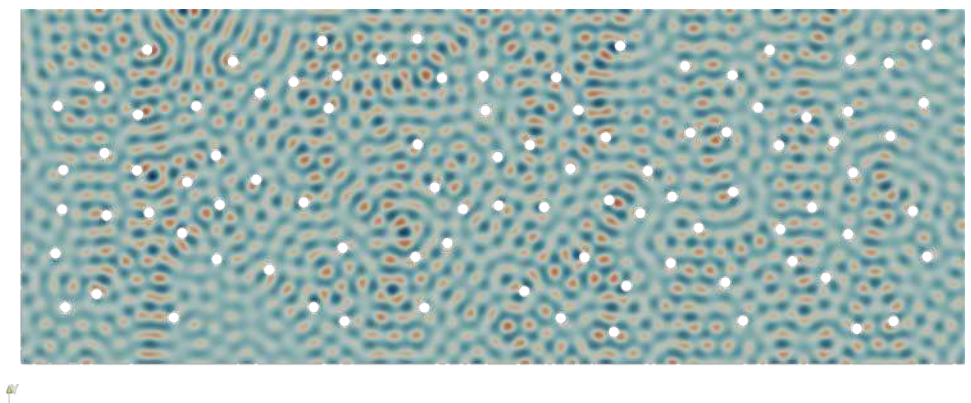

Figure 3: $\varphi - \varphi^{\text{inc}}$ in the frequency domain.

The challenge for proving convergence comes from the non-linear term $\nabla\cdot([\mathbf{M}^\varepsilon\cdot\nabla\varphi^\varepsilon]\,\mathbf{M}^\varepsilon)$ appearing in the divergence part of Equation (2), which requires careful uniform estimates for the flow potential $u^\varepsilon$. The natural functional framework of Equation (1) yields only $L^2$ bounds for $\mathbf{M}^\varepsilon$, while for this non-linear term we will need $L^\infty$ bounds. While such homogenisation problems have been well-studied for the Helmholtz equation, to our knowledge this has not been presented for the Convected Helmholtz equation with non-uniform flow.

# The $hp$-FEM does not suffer from the pollution effect when applied to the Helmholtz equation with piecewise Gevrey coefficients and Gevrey obstacles

**Jeffrey Galkowski**[1], **Mostafa Meliani**[2], [*], **Euan Spence**[1]
[1]Department of Mathematics, University College London, London, UK
[2]Department of Mathematical Sciences, University of Bath, Bath, UK
*Email: mm4138@bath.ac.uk

**Abstract**

We present new results that weaken the regularity conditions under which the $hp$-FEM is proved not to suffer from the pollution effect.

## 1 Context

The $h$-version of the finite-element method (FEM) famously suffers from the pollution effect; namely, $\gg k^d$ degrees of freedom are required to maintain accuracy as $k \to \infty$ [1].

It was first shown in [7, 8] that, if

(i) the solution operator is polynomially bounded in $k$ (which holds for "most" $k$ by [5]), and

(ii) the domain and the coefficients satisfy appropriate regularity assumptions,

then the $hp$-FEM does not suffer from the pollution effect, i.e., the Galerkin solution is quasi-optimal, uniformly in $k$, if

$$\frac{hk}{p} \leq c \text{ and } p \geq C \log k. \qquad (1)$$

These results were then generalised to wider classes of Helmholtz problems in [2–4, 6].

Regarding the regularity assumptions (ii): the current state of the art is that $k$-uniform quasi-optimality of the $hp$-FEM under (1) is proved for the variable-coefficient Helmholtz Dirichlet problem when *either*

(a) the obstacle is analytic and the coefficients are smooth and analytic in a neighbourhood of the obstacle [3, 4, 6], *or*

(b) the obstacle is analytic and the coefficients are piecewise analytic, with the interfaces across which the coefficients jump also analytic [2].

## 2 Summary of our contribution

We consider the $hp$-FEM applied to the PML approximation of the variable coefficient Helmholtz equation outside a Dirichlet obstacle (see Sections 4 and 5 below).

We prove $k$-uniform quasi-optimality of the $hp$-FEM under the assumptions that (i) the solution operator is polynomially bounded in $k$ and (ii) *the obstacle is Gevrey and the coefficients are Gevrey in a neighbourhood of the obstacle, piecewise smooth elsewhere, and Gevrey in one-sided neighbourhoods of all the interfaces across which the coefficients jump.*

Since the Gevery space of functions is an intermediate space between $C^\infty$ and analytic functions, this result generalises both sets of results (a) and (b) above.

## 3 Gevrey space of functions

**Definition.** *Given $\sigma \geq 1$, a function $g : \Omega \to \mathbb{C}$ is $\sigma$-Gevrey if given $K \Subset \Omega$ there exists $M$ such that*

$$\sup_{x \in K} |\partial^\gamma g(x)| \leq M^{|\gamma|+1}(\gamma!)^\sigma.$$

Note that a 1-Gevrey function is analytic.

## 4 The Helmholtz scattering problem we consider

To avoid introducing lots of notation, we only state here our results for globally smooth coefficients that are Gevrey in the neighourhood of the obstacle. The talk will present, in addition, results for the piecewise smooth coefficients described above.

Let $\sigma \geq 1$. Let $\Omega_- \subset B_{R_0} := \{x : |x| < R_0\} \subset \mathbb{R}^d$, $d = 2, 3$, be a bounded open set with $\sigma$-Gevrey boundary $\Gamma$ and set $\Omega_+ := \mathbb{R}^d \setminus \overline{\Omega_-}$. Let $A_{\rm scat} \in C^\infty(\Omega_+; \mathbb{R}^{d\times d})$ be symmetric positive definite, let $c_{\rm scat} \in C^\infty(\Omega_+; \mathbb{R})$ be strictly positive and bounded and let $A_{\rm scat}$ $c_{\rm scat}$ be such that there exist $R_* > R_{\rm scat} > R_0 > 0$ such that

$$\overline{\Omega_-} \cup \operatorname{supp}(I - A_{\rm scat}) \cup \operatorname{supp}(1 - c_{\rm scat}) \Subset B_{R_{\rm scat}},$$

and $A_{\rm scat}$, $c_{\rm scat}$ are $\sigma$-Gevrey in $B_{R_*}$. Let $g \in L^2(\Omega_+)$ with $\operatorname{supp} g \Subset \mathbb{R}^d$ and $k > 0$, $u \in H^1_{\rm loc}(\Omega_+)$ satisfies the exterior Dirichlet problem if

$$\begin{aligned} -c_{\rm scat}^2 \nabla \cdot (A_{\rm scat} \nabla u) - k^2 u = g \quad & \text{in } \Omega_+ \\ u = 0 \quad & \text{on } \Gamma \end{aligned} \tag{2}$$

and $u$ is outgoing in the sense that it satisfies the Sommerfeld radiation condition

$$\partial_r u(x) - iku(x) = o\left(r^{-(d-1)/2}\right) \tag{3}$$

as $r := |x| \to \infty$, uniformly in $\hat{x} := x/r$.

## 5 The PML problem

In order to approximate the solution to the problem stated above, we use a PML truncation. Let $R_{\rm tr} > R_1 > R_* > R_{\rm scat}$ and let $\omega_{\rm tr} \subset \mathbb{R}^d$ be a bounded Lipschitz domain with $B_{\rm tr} \subset \Omega_{\rm tr}$. Let $\Omega := \Omega_{\rm tr} \cap \Omega_+$, $\Gamma_{\rm tr} := \partial\Omega_{\rm tr}$, and $0 \geq \theta < \pi/2$. Let $P := -c_{\rm scat}^2 \nabla \cdot (A_{\rm scat} \nabla)$. The PML problem is then

$$P_\theta v - k^2 v = g \quad \text{in } \Omega_+, \quad v = 0 \quad \text{on } \Gamma \cup \Gamma_{\rm tr} \tag{4}$$

where $P_\theta := P$ for $r \leq R_1$ and $P_\theta := -\Delta_\theta$ for $r \geq R_1$, where $-\Delta_\theta$ is a second-order differential operator defined in spherical coordinates as in, e.g., [4, §1.3.1].

## 6 Statement of the main result

**Theorem** (Quasioptimality of *hp*-FEM). *Suppose that $\Omega_-$, $A_{\rm scat}$, $c_{\rm scat}$, and $\Omega_{\rm tr}$ are as in §4. Suppose further that $\Omega_-$, $A_{\rm scat}$, $c_{\rm scat}$ and $K \subset [k_0, \infty)$ are such that the solution operator of the exterior Dirichlet problem is polynomially bounded. Let $\{V_N\}_{N=0}^\infty$ be the piecewise-polynomial approximation spaces described in [7, §5.1.1] (with the requirement that the element maps be analytic relaxed to that they be Gevrey). Given $\varepsilon > 0$, there exist $k_1$, $C$, $C_{\rm qo} > 0$ such that the following is true. Given $G \in (H^1(\Omega))^*$, for all $k \in K \cap [k_1, \infty)$, $\varepsilon \leq \theta \leq \pi/2 - \varepsilon$, and $R_{\rm tr} \geq R_1(1+\varepsilon)$, the solution $v$ to the PML problem (4) exists and is unique. Furthermore, if*

$$\frac{hk}{p} \leq C \qquad \text{and} \qquad p \geq 1 + \varepsilon \log k,$$

*then the Galerkin solution of the PML problem (4) satisfies*

$$\|v - v_n\|_{H^1_k(\Omega)} \leq C_{\rm qo} \min_{w_N \in V_N} \|v - w_n\|_{H^1_k(\Omega)},$$

*where*

$$\|u\|^2_{H^1_k(\Omega)} := k^{-2} \|\nabla u\|^2_{L^2(\Omega)} + \|u\|^2_{L^2(\Omega)}.$$

# Characterizing the imaging indicator of the linear sampling method for the inverse medium scattering problem

**Lorenzo Audibert**[1], **Shixu Meng**[2,*]
[1]PRISME, EDF R&D, 6 Quai Watier, 78400 Chatou, France and UMA, Inria, ENSTA Paris, Institut Polytechnique de Paris, 91120 Palaiseau, France
[2]Department of Mathematical Sciences, The University of Texas at Dallas, Richardson, USA
*Email: smeng@utdallas.edu

**Abstract**

We propose a new characterization of an imaging indicator in the class of linear sampling method for the inverse medium scattering problem. We show that the imaging indicator represents a weighted average of a certain nonlinear transformation of the unknown.

## 1 Introduction

Let $k > 0$ be the wave number. A plane wave takes the following form

$$e^{ikx\cdot\hat{\theta}}, \quad \hat{\theta} \in \mathbb{S}^{d-1} := \{x \in \mathbb{R}^d : |x| = 1\},$$

where $\hat{\theta}$ is the direction of propagation, $d = 2$ or 3. Let $\Omega \subset \mathbb{R}^d$ be an open and bounded set with Lipschitz boundary $\partial\Omega$ such that $\mathbb{R}^d\backslash\overline{\Omega}$ is connected. Let the real-valued function $q(x) \in L^\infty(\mathbb{R}^d)$ be the contrast of the medium, $q > 0$ on $\Omega$ and $q = 0$ on $\mathbb{R}^d\backslash\overline{\Omega}$ such that the support of the contrast is $\Omega$. We assume that $q_{\sup} \geq q(y) \geq q_{\inf}$, a.e., $y \in \Omega$, for some positive constants $q_{\sup}$ and $q_{\inf}$; the case when $q$ is complex-valued can be treated similarly. The forward scattering problem due to a plane wave $e^{ikx\cdot\hat{\theta}}$ is to find total wave filed $e^{ikx\cdot\hat{\theta}} + u^s(x;\hat{\theta};k)$ belonging to $H^1_{loc}(\mathbb{R}^d)$ such that

$$\Delta u^s + k^2(1+q)u^s = -k^2 q e^{ikx\cdot\hat{\theta}} \text{ in } \mathbb{R}^d, \quad (1)$$

$$\lim_{r:=|x|\to\infty} r^{\frac{d-1}{2}}\left(\frac{\partial u^s}{\partial r} - iku^s\right) = 0. \quad (2)$$

A solution is called radiating if it satisfies the Sommerfeld radiation condition (2) uniformly for all directions. There exists a unique radiating solution to (1)–(2). Note that as $r = |x| \to \infty$,

$$u^s(x) = \begin{cases} \frac{e^{i\frac{\pi}{4}}}{\sqrt{8k\pi}}\frac{e^{ikr}}{\sqrt{r}}\left\{u^\infty(\hat{x}) + \mathcal{O}\left(\frac{1}{r}\right)\right\} & d=2 \\ \frac{e^{ik|x|}}{4\pi|x|}\left\{u^\infty(\hat{x}) + \mathcal{O}\left(\frac{1}{r}\right)\right\} & d=3 \end{cases}$$

uniformly with respect to all directions $\hat{x} := x/|x| \in \mathbb{S}^{d-1}$, we introduce $u^\infty(\hat{x};\hat{\theta};k)$ which is known as the far-field pattern due to incident plane wave $e^{ikx\cdot\hat{\theta}}$ at wave number $k$. The multi-static data at a fixed frequency are given by

$$\{u^\infty(\hat{x};\hat{\theta};k) : \hat{x} \in \mathbb{S}^{d-1}, \hat{\theta} \in \mathbb{S}^{d-1}\}. \quad (3)$$

The goal of linear sampling type methods is to determine $\Omega$ and other possible information. The work [1] showed that an imaging indicator in the class of linear sampling method represents the average of the unknown for a linear problem. This work extends [1] to the nonlinear inverse medium problem and demonstrates that the imaging indicator represents a weighted average of the a nonlinear transformation of the unknown.

## 2 Far-field operator and its factorization

The far-field operator $\mathcal{F} : L^2(\mathbb{S}^{d-1}) \to L^2(\mathbb{S}^{d-1})$ is defined by

$$(\mathcal{F}g)(\hat{x}) := \int_{\mathbb{S}^{d-1}} u^\infty(\hat{x};\hat{\theta};k)g(\hat{\theta})\,\mathrm{d}s(\hat{\theta}), \quad (4)$$

where the kernel is given by the data (3). Introduce $\mathcal{H} : L^2(\mathbb{S}^{d-1}) \to L^2(\Omega)$ by

$$(\mathcal{H}g)(y) := \int_{\mathbb{S}^{d-1}} e^{iky\cdot\hat{\theta}}g(\hat{\theta})\,\mathrm{d}s(\hat{\theta}), \ y \in \Omega,$$

$\forall g \in L^2(\mathbb{S}^{d-1})$, and it follows directly that its adjoint $\mathcal{H}^* : L^2(\Omega) \to L^2(\mathbb{S}^{d-1})$ is given by

$$(\mathcal{H}^*h)(\hat{x}) := \int_\Omega e^{-ik\hat{x}\cdot y}h(y)\,\mathrm{d}y, \ \hat{x} \in \mathbb{S}^{d-1},$$

$\forall h \in L^2(\Omega)$. Another operator $\mathcal{T}$ is needed for the factorization, namely $\mathcal{T} : L^2(\Omega) \to L^2(\Omega)$ which is given by

$$\mathcal{T}f := k^2 qf + k^2 qv|_\Omega, \quad \forall f \in L^2(\Omega),$$

where $v \in H^1_{loc}(\mathbb{R}^d)$ is the unique radiating solution to $\Delta v + k^2(1+q)v = -k^2 qf$. Then it holds that [3]

$$\mathcal{F} = \mathcal{H}^*\mathcal{T}\mathcal{H}.$$

For later purposes, let $\phi_z \in L^2(\mathbb{S}^{d-1})$ be given by

$$\phi_z(\hat{x}) := e^{-ik\hat{x}\cdot z}, \quad \hat{x} \in \mathbb{S}^{d-1} \tag{5}$$

It holds that $z \in \Omega$ if, and only if, $\phi_z \in \text{Range}(\mathcal{H}^*)$, cf. [3]. Moreover let $B(z,\epsilon) := \{x \in \mathbb{R}^d : |x-z| < \epsilon\}$ for some fixed small $\epsilon > 0$, if $\overline{B(z,\epsilon)} \subset \Omega$, then it holds that $\phi_z = \mathcal{H}^* E_z$, where $E_z \in L^2(\Omega)$ is given by

$$E_z := \begin{cases} -(\Delta w_z + k^2 w_z) & \text{in } B(z,\epsilon) \\ 0 & \text{otherwise} \end{cases}, \tag{6}$$

here $w_z(x) := \chi_{|x-z|}\Phi(x,z)$ in $\mathbb{R}^d$, and $\chi \in C^\infty(\mathbb{R})$ is a cut off function with $\chi(t) = 1$ for $|t| \geq \epsilon$ and $\chi(t) = 0$ for $|t| \leq \epsilon/2$, cf. [3].

## 3 New imaging indicator and physical interpretation

Let $R(\mathcal{H})$ denote the range of $\mathcal{H}$ and $Y_\Omega := \text{closure}(R(\mathcal{H}))$. Note that $R(\mathcal{H})$ is dense in $\{v \in H^1(\Omega) : \Delta v + k^2 v = 0 \text{ in } \Omega \text{ variationally}\}$. Let $\mathcal{P}_\Omega : L^2(\Omega) \to L^2(\Omega)$ be the projection operator given by $\mathcal{P}_\Omega x = x|_{Y_\Omega}$ for any $x \in L^2(\Omega)$. We also recall [2] that $k \in \mathbb{C}$ is called an *interior transmission eigenvalue* if there exists nontrivial $(w,v) \in L^2(\Omega) \times L^2(\Omega)$ such that $w - v \in H^2(\Omega)$ and

$$\begin{aligned} \Delta w + k^2(1+q)w = 0 &\quad \text{in} \quad \Omega, \\ \Delta v + k^2 v = 0 &\quad \text{in} \quad \Omega, \\ w = v &\quad \text{on} \quad \Omega, \\ \frac{\partial w}{\partial \nu} = \frac{\partial v}{\partial \nu} &\quad \text{on} \quad \partial\Omega. \end{aligned}$$

**Lemma 1** *Assume that $k$ is not an interior transmission eigenvalue, then $\mathcal{P}_\Omega \mathcal{T} \mathcal{P}_\Omega : Y_\Omega \to Y_\Omega$ is coercive and there exists a positive constant $T_{\inf} > 0$ such that, for any $f \in Y_\Omega$,*

$$|\langle \mathcal{P}_\Omega \mathcal{T} \mathcal{P}_\Omega f, f \rangle_{L^2(\Omega)}| \geq T_{\inf} \|f\|^2_{L^2(\Omega)}. \tag{7}$$

We introduce a family of regularization schemes $\{\mathcal{R}_\alpha\}_{\alpha>0}$ by

$$\mathcal{R}_\alpha h := \sum_{n=0}^{\infty} f_\alpha(|\mu_n|^2)\overline{\mu_n}\langle h, \zeta_n\rangle \zeta_n, \tag{8}$$

where $f_\alpha$ is a regularizing filter that is a bounded, real-valued, and piecewise continuous function $f_\alpha : (0,\infty) \to \mathbb{R}$ such that $\lim_{\alpha\to 0} f_\alpha(\mu) = \frac{1}{\mu}$ for all $\mu > 0$, $|\mu f_\alpha(\mu)| \leq d_0$ for all $\alpha \geq 0$ and $\mu > 0$, here $d_0 > 0$ is a constant. Now we are ready to obtain the first main result.

**Theorem 2** *Assume that $k$ is not an interior transmission eigenvalue. Suppose that $\{\mathcal{R}_\alpha\}_{\alpha>0}$ is a family of regularization schemes given by (8) and set $g_{z,\alpha} = \mathcal{R}_\alpha \phi_z$. Then it holds for any $\overline{B(z,\epsilon)} \subset \Omega$ that*

$$\lim_{\alpha\to 0}\langle g_{z,\alpha}, \phi_z\rangle = \langle (\mathcal{P}_\Omega \mathcal{T} \mathcal{P}_\Omega)^{-1}\mathcal{P}_\Omega E_z, \mathcal{P}_\Omega E_z\rangle. \tag{9}$$

*where the left (resp. right) inner product is with respect to $L^2(\mathbb{S}^{d-1})$ (resp. $L^2(\Omega)$).*

The following lemma constructs $\mathcal{P}_\Omega E_z$ explicitly.

**Lemma 3** *Let $E_z^{\mathcal{P}} := \mathcal{P}_\Omega E_z$. Then $E_z^{\mathcal{P}} = -(\Delta + k^2)w_z^{\mathcal{P}}$ in $\Omega$, where $w_z^{\mathcal{P}} \in H^2(\Omega)$ is the unique solution to*

$$\begin{aligned} (\Delta + k^2)(\Delta + k^2)w_z^{\mathcal{P}} = 0 &\quad \text{in} \quad \Omega, \quad (10)\\ w_z^{\mathcal{P}} = \Phi(\cdot, z) &\quad \text{on} \quad \partial\Omega, \quad (11)\\ \frac{\partial w_z^{\mathcal{P}}}{\partial \nu} = \frac{\partial \Phi(\cdot,z)}{\partial \nu} &\quad \text{on} \quad \partial\Omega. \quad (12) \end{aligned}$$

Now we are ready to obtain the main theorem. Note that $E_z^{\mathcal{P}}$ is independent of the unknown and the information about the unknown is represented by $\mathcal{P}_\Omega \mathcal{T} \mathcal{P}_\Omega$.

**Theorem 4** *Assume that $k$ is not an interior transmission eigenvalue. Suppose that $\{\mathcal{R}_\alpha\}_{\alpha>0}$ is a family of regularization schemes given by (8) and set $g_{z,\alpha} = \mathcal{R}_\alpha \phi_z$.*

- *For any $z \in \Omega$, it holds that*

  $$\lim_{\alpha\to 0}\langle g_{z,\alpha}, \phi_z\rangle = \langle (\mathcal{P}_\Omega \mathcal{T} \mathcal{P}_\Omega)^{-1} E_z^{\mathcal{P}}, E_z^{\mathcal{P}}\rangle \tag{13}$$

  *where $E_z^{\mathcal{P}} = -(\Delta + k^2)w_z^{\mathcal{P}}$ with $w_z^{\mathcal{P}} \in H^2(\Omega)$ satisfying (10)–(12).*

- *For any $z \in \mathbb{R}^d \backslash \overline{\Omega}$, it holds that*

  $$\lim_{\alpha\to 0}\langle g_{z,\alpha}, \phi_z\rangle = \infty. \tag{14}$$

# Error estimation for a homogenized model of a space-time modulated metamaterial

Jean-Francois Mercier[1,*], Marie Touboul[1]
[1]POEMS, CNRS, Inria, ENSTA Paris, Institut Polytechnique de Paris, 91120 Palaiseau, France
*Email: jean-francois.mercier@ensta.fr

**Abstract**

We study metamaterials governed by a 1D wave equation with periodic coefficients modulated as a traveling wave, i.e. they depend on $X - cT$, where $c$ is the modulation velocity. Using time-domain periodic homogenization, we replace the infinite periodic medium by an equivalent homogeneous effective model. We estimate the error between the exact solution and the homogenized model for a medium with small parameter $\varepsilon$ over an observation time $T$. Theoretical error bounds of order $\varepsilon T^3$ are derived using energy estimates.

## 1 Setting of the problem

The medium of periodicity $h$ satisfies the non-dimensionalized wave equation for $U$

$$\frac{\partial}{\partial T}\left(\alpha\left(\frac{X-cT}{\varepsilon}\right)\frac{\partial U}{\partial T}\right) = F(X,T)$$

$$+\frac{\partial}{\partial X}\left(\beta\left(\frac{X-cT}{\varepsilon}\right)\frac{\partial U}{\partial X}\right), \quad (1)$$

with zero initial conditions. $X$ and $T$ denote the non-dimensionalized space coordinate and time, respectively. The parameters, $\alpha$ and $\gamma$ are 1-periodic functions assumed to be piecewise continuous, $F(X,T)$ is a source term and $\varepsilon = \lambda^\star/h$ with $\lambda^\star$ the characteristic source wavelength.

To simplify homogenization and analysis, we reformulate the problem in the moving frame defined by the coordinates $(x,t) = (X - cT, T)$. The fields in this frame are denoted $u(x,t) = U(X,T)$, $f(x,t) = F(x,t)$ and (1) yields

$$L^\varepsilon u^\varepsilon = -f, \; u^\varepsilon(x,0) = 0 = (\partial u^\varepsilon/\partial t)(x,0), \; (2)$$

with $L^\varepsilon u^\varepsilon = \partial_x\left(\gamma\left(x/\varepsilon\right)\partial_x u^\varepsilon\right) - \alpha\left(x/\varepsilon\right)\partial_t^2 u^\varepsilon + c\alpha\left(x/\varepsilon\right)\partial_x\partial_t u^\varepsilon + c\partial_x\left(\alpha\left(x/\varepsilon\right)\partial_t u^\varepsilon\right)$. We suppose that $\alpha$ and $\gamma = \beta - \alpha c^2 \in L^\infty(Y)$ where $Y = (0,1)$ is the unit cell, and are strictly positive.

## 2 Two-scale analysis

We now make the assumption that $\varepsilon \ll 1$ (low-frequency regime) and we note $y = x/\varepsilon$ the fast variable. In [1] the two-scale asymptotic expansion method combined to the notion of slow and fast variables are used with the field $u^\varepsilon$ expanded using the following ansatz

$$u^\varepsilon(x,t) = \sum_{j\geq 0} \varepsilon^j u_j(x, x/\varepsilon, t), \quad (3)$$

with $y \to u_j(x,y,t)$ 1-periodic functions. Now we sum up the results in [1].

The $0^{th}$ order field reads $u_0(x,y,t) = U_0(x,t)$ satisfying the effective medium equation

$$\mathcal{L}_0 U_0 = f, \;\; U_0(x,0) = 0 = (\partial U_0/\partial t)(x,0), \; (4)$$

where $\mathcal{L}_0 U_0 = \alpha_0 \partial_t^2 U_0 - 2W_0 \partial_x \partial_t U_0 - \beta_0 \partial_x^2 U_0$. The $1^{th}$ order field reads

$$u_1 = U_1 + P_0(\partial U_0/\partial t) + P_1(\partial U_0/\partial x), \quad (5)$$

with $U_1(x,t)$ satisfying

$$\mathcal{L}_0 U_1 = \sum_{j=0}^{3} B_j(\partial^3 U_0/\partial x^j \partial t^{3-j}), \quad (6)$$

with $B_j$ some constants explicitly determined.

The $2^{nd}$ order field reads

$$u_2 = U_2 + \sum_{j=0}^{1} P_j(\partial U_0/\partial x^j \partial t^{1-j})$$

$$+\sum_{j=0}^{2} Q_j(\partial^2 U_0/\partial x^j \partial t^{2-j})$$

The functions $P_j(y)$ and $Q_j(y)$ are solutions of the so-called cell problems and we suppose that for all $j$, $P_j \in L^\infty(Y)$, $Q_j \in L^\infty(Y)$ and $dQ_j/dy \in L^\infty(Y)$. This is naturally satisfied if $\alpha$ and $\gamma$ are regular functions.

Now we show our main result, that $u_0(x, t \leq T)$ is a good approximation of $u^\varepsilon(x, t \leq T)$ for $\varepsilon$ and $T$ small enough, summarized in the

**Theorem 1** *A constant $C > 0$ exists such that $\sup_{t\in[0,T]} \|(u^\varepsilon - u_0)(x,t)\| \leq C\varepsilon T^3$.*

Here $\|.\|$ designs the $L^2$-norm in space. To prove this theorem, we need first to derive estimates for the homogenized solutions $u_0$, $u_1$, $u_2$.

## 3 Estimates of the homogenized solutions

The crucial step to estimate $u_0 = U_0$ is to establish an energy control. Multiplying (4) by $\partial U_0/\partial t$ leads, after manipulations, to the

**Lemma 2** *For all $t > 0$, the energy $E(t) = \int_{\mathbb{R}} \left[\alpha_0 \left[(\partial U_0/\partial t)(x,t)\right]^2 + \beta_0 \left[(\partial U_0/\partial x)(x,t)\right]^2\right] dx$, satisfies $E^{1/2}(t) \leq (1/\sqrt{\alpha_0}) \int_0^t \|f(\cdot,s)\| ds$.*

Multiplying (6) by $\partial U_1/\partial t$, a similar energy estimation is found for $U_1$. Then introducing for any $U(x,t)$ the useful compact notations $\|\partial^n U(\cdot,t)\| = \sum_{j=0}^n \left\|(\partial^n U/\partial x^j \partial t^{n-j})(\cdot,t)\right\|$ and $M_n(t) = \int_0^t \|\partial^n f(\cdot,s)\| \, ds$, we deduce the

**Lemma 3** *A constant $C(\alpha,\gamma) > 0$ exists such that for all $t > 0$ and $1 \leq n \leq 3$, $\|\partial^n U_0(\cdot,t)\| \leq C(\alpha,\gamma) M_{n-1}(t)$ and $\|u_1(x, x/\varepsilon, t)\| \leq L_\infty C(\alpha,\gamma) \left(M_0(t) + \int_0^t M_1(s)\, ds + \int_0^t (t-s) M_2(s)\, ds\right)$, where $L_\infty = \max\left(\max_{0\leq j\leq 1} \|P_j\|_{L^\infty(Y)}, \max_{0\leq j\leq 2} \|Q_j\|_{L^\infty(Y)}, \max_{0\leq j\leq 2} \left\|Q'_j\right\|_{L^\infty(Y)}\right)$.*

Similar estimations are obtained for $u_2$. Thanks to these three estimates, we have now all the ingredients to prove the main theorem.

## 4 Convergence of the homogenized solution

The error is defined by $\tilde{e}^\varepsilon = u^\varepsilon - u_0$ and to estimate it we introduce the intermediate quantity $e^\varepsilon(x,t) = u^\varepsilon(x,t) - [u_0(x,t) + \varepsilon u_1(x, x/\varepsilon, t)]$. We first establish the equation satisfied by $e^\varepsilon$.

**Lemma 4** *For all $t > 0$, $e^\varepsilon$ satisfies*

$$-(L^\varepsilon e^\varepsilon)(x,t) = \varepsilon s(x,t), \qquad (7)$$

*with $s(x,t) = -\partial_x \tilde{v}_2(x,t) + \partial_x v_2(x, x/\varepsilon, t) + (L_0 u_1)(x,t)$ where $\tilde{v}_2(x,t) = v_2(x, x/\varepsilon, t)$ with $v_2(x,y,t) = \gamma(y)\partial_y u_2(x,y,t)$ and where $(L_0 u_1)(x,t) = [\gamma(y)\partial_x^2 u_1(x,y,t) - \alpha(y)\partial_t^2 u_1(x,y,t) + 2c\alpha(y)\, \partial_x \partial_t u_1(x,y,t)]\ (x, y = x/\varepsilon, t)$.*

**Remark 5** *It is crucial to have a right hand side of order $\varepsilon$ which implied to introduce the function $\tilde{v}_2$. Note that the introduction of $e^\varepsilon$ is necessary since the quantity of interest $\tilde{e}^\varepsilon$ satisfies $L^\varepsilon \tilde{e}^\varepsilon = g$ with a right hand side $g$ of order $1/\varepsilon$, thus not tending to zero when $\varepsilon \to 0$.*

Now we are able to prove the main result, that $\|\tilde{e}^\varepsilon\| = O(\varepsilon T^3)$. To do so we need several steps. First we will prove that $\|\partial e^\varepsilon/\partial t\| = O(\varepsilon T^2)$ and this is once again done by using energetic arguments, multiplying (7) by $\partial e^\varepsilon/\partial t$

**Lemma 6** *The energy $E(t) = \|(\partial e^\varepsilon/\partial x)(\cdot,t)\|^2 + \|(\partial e^\varepsilon/\partial t)(\cdot,t)\|^2$ satisfies for all $t \in [0,T]$, $E(t)^{1/2} \leq 2\varepsilon L_\infty C(\alpha,\gamma)\ (M_1(T) + T M_2(T) + \int_0^t M_2(s)\, ds + \int_0^t (t-s) M_3(s)\, ds)$.*

Integrating over time, we deduce the estimate for $\|e^\varepsilon\|$. To establish the final result, we separate the source time $T_f$ and the observation time $T > T_f$. We consider that the source $f(x,t)$ is switched off at a time called $T_f$. Then using $u^\varepsilon - u_0 = \tilde{e}^\varepsilon = e^\varepsilon + \varepsilon u_1$, we can prove the

**Theorem 7** *If Supp$(f(t)) \subset [0, T_f]$ for some final source time $T_f > 0$, then for all time $T$ such that $T > T_f$ we have the estimate $\sup_{t\in[0,T]} \|(u^\varepsilon - u_0)(x,t)\| \leq L_\infty C(\alpha,\gamma) M(T_f)\, \varepsilon T^3$ with $M(T_f) = \max_{j\in\{1,2,3,4\}} M_j(T_f)$.*

Note that the $T^3$ dependence is not usual. In the non space-time modulated case, corresponding to $c = 0$, thus parameters $\alpha$ and $\beta$ depending only on $X$, we get $U_1 = 0$ from which we recover the usual result of an error $O(\varepsilon T^2)$ [2]. Note that a result of an error $O(\varepsilon T)$ can be found in the litterature but it is valid only for bounded domains, with rigid boundaries (like in [2]) or with periodic boundary conditions to represent infinite domains [3,4].
Time-domain simulations to compare the exact solution and the homogenized solution have been developed and are currently performed to confirm the theoretical error estimates.

# High-order Wavelet Galerkin Methods for Scattering Problems

**Michelle Michelle**[1,*], **Bin Han**[1]

[1]Department of Mathematical and Statistical Sciences, University of Alberta, Edmonton, AB, Canada

*Email: mmichell@ualberta.ca

## Abstract

High-order schemes are known to be effective in mitigating the pollution effect that arises when solving the Helmholtz equation. We present a high-order wavelet method for solving certain scattering problems modeled by the Helmholtz equation. We explore how these wavelets, which constitute a sparse multiscale representation system, can be constructed to handle rectangular and some non-rectangular domains. Additionally, we discuss several advantages of this wavelet method.

## 1 Introduction

The Helmholtz equation is challenging to solve numerically due to its sign-indefinite standard weak formulation and its highly oscillating solution. We address this issue by designing a high-order discretization scheme in the form of a wavelet Galerkin method. Wavelets on the real line form sparse multiscale representation systems and constitute Riesz bases in Sobolev spaces. See [3] for a comprehensive treatment on the subject.

Wavelets come in many forms, and we have great flexibility in their constructions. We are especially interested in high-order spline wavelets due to their approximation power, explicit expressions, and ease of implementation. One compelling reason for choosing wavelets as a basis in the Galerkin method is that they produce a sparse coefficient matrix comparable to that of a standard finite element method (FEM), while being much better conditioned [5–7]. In particular, the condition number of the wavelet coefficient matrix is uniformly bounded irrespective of its size. This property follows from the fact that wavelets are Riesz bases in Sobolev spaces. Numerically, this translates to significantly fewer iterations for an iterative scheme, such as GMRES, to achieve a prescribed tolerance (e.g., see [6, Table 5]).

## 2 Models

We consider two model problems. The first one is modeled by the following equation [1]:

$$\begin{aligned} &\Delta u + \kappa^2 u = f \quad \text{in} \quad \Omega := (0,1)^2, \\ &u = 0 \quad \text{on} \quad \partial\Omega \backslash \Gamma, \\ &\frac{\partial u}{\partial \nu} = \mathcal{T}(u) + g \quad \text{on} \quad \Gamma, \end{aligned} \tag{1}$$

where $\kappa$ is a variable wavenumber, $\Gamma := (0,1) \times \{1\}$, $f$ is a variable source term, $g$ is a boundary condition, $\nu$ is the unit outward normal, and

$$\mathcal{T}(u) := \frac{\mathrm{i}\kappa_0}{2} ⨍_0^1 \frac{1}{|x-x'|} H_1^{(1)}(\kappa_0|x-x'|) u(x',1) dx',$$

where $\kappa_0 > 0$, $⨍$ stands for the Hadamard finite part integral, and $H_1^{(1)}$ is the Hankel function of the first kind of degree 1.

The second one is modeled by the following equation [8]:

$$\begin{aligned} &-\nabla \cdot (A \nabla u) - \kappa^2 u = f \text{ in } \Omega := B_R \backslash \overline{S}, \\ &u = g \text{ on } \partial S, \quad u = 0 \text{ on } \partial B_R, \end{aligned} \tag{2}$$

where $B_R := \{x \in \mathbb{R}^2 : |x| < R\}$, $S$ is a star-shaped domain (also called a scatterer) with smooth boundary such that $\overline{S} \subsetneq B_R$ and $A$ is a bounded spatially varying scalar or symmetric matrix-valued function containing information on the media and the perfectly matched layer (PML).

## 3 Methods and Results

Since our method is based on the Galerkin framework, we naturally start by constructing a Riesz wavelet basis on the domain of interest. Consider the model problem in (1) and define

$$\mathcal{H}(\Omega) := \{u \in H^1(\Omega) : u = 0 \text{ on } \partial\Omega \backslash \Gamma\}$$

Its weak formulation is to find $u \in \mathcal{H}(\Omega)$ such that for all $v \in \mathcal{H}(\Omega)$

$$\begin{aligned} \langle \nabla u, \nabla v \rangle_\Omega - \langle \kappa^2 u, v \rangle_\Omega - \langle \mathcal{T}(u), v \rangle_\Gamma \\ = \langle g, v \rangle_\Gamma - \langle f, v \rangle_\Omega. \end{aligned} \tag{3}$$

Now, starting with a compactly supported biorthogonal multiwavelet in $L^2(\mathbb{R})$, we can adapt it to a bounded interval with our method in [4], and form a Riesz basis of the Sobolev space $\mathcal{H}(\Omega)$, denoted by $\mathcal{B}^{\mathcal{H}(\Omega)}$, via tensor product. The latter is theoretically guaranteed by [6, Theorem 1.2], whose proof is direct in that we avoid the use of the classical Jackson and Bernstein inequalities from approximation theory (e.g., [2]), which often involve technical assumptions.

In the actual computation, we take a finite subset of $\mathcal{B}^{\mathcal{H}(\Omega)}$ truncated at scale level $J$, denoted by $\mathcal{B}_J^{\mathcal{H}(\Omega)}$; i.e., we only keep basis functions with a dilation factor less than or equal to $J-1$. The approximate solution takes the form of $u_J = \sum_{\eta \in \mathcal{B}_J^{\mathcal{H}(\Omega)}} c_\eta \eta$, where $c_\eta$ can be found from (3) for all $v \in \mathcal{B}_J^{\mathcal{H}(\Omega)}$ with $u$ replaced by $u_J$. Below, we reproduce a portion of [6, Table 3] using the spline wavelet of order 4 in [6, Example 3.6]. For the following numerical results, we let $\kappa = 16\pi$ and the PDE data be chosen such that the true solution $u = \exp(xy)\sin(\kappa x)\sin((\kappa+\pi/2)y)$. We can see that even though the matrix size quadruples, the condition number plateaus, because the smallest singular values are bounded away from zero. The $L^2$ convergence rates agree with the approximation order of the underlying wavelet basis.

| $N_J$ | $\sigma_{\max}$ | $\sigma_{\min}$ | $\frac{\sigma_{\max}}{\sigma_{\min}}$ | $\frac{\|u-u_J\|_2}{\|u\|_2}$ | Order |
|---|---|---|---|---|---|
| 2256 | 5.43 | 8.68E-3 | 6.26E+2 | 5.15E-2 | |
| 9120 | 5.68 | 6.26E-3 | 9.08E+2 | 2.97E-3 | 4.08 |
| 36672 | 5.87 | 6.06E-3 | 9.68E+2 | 2.09E-4 | 3.81 |

Table 1: $\sigma_{\max}$, $\sigma_{\min}$ are largest and smallest singular values respectively, order is the convergence rate in the $L^2(\Omega)$-norm, and $N_J$ is the cardinality of $\mathcal{B}_J^{\mathcal{H}(\Omega)}$.

Figure 1 shows a wavelet-based approximate solution of (1) with a variable wavenumber. Refer to [6] for other examples involving variable wavenumbers or unknown solutions.

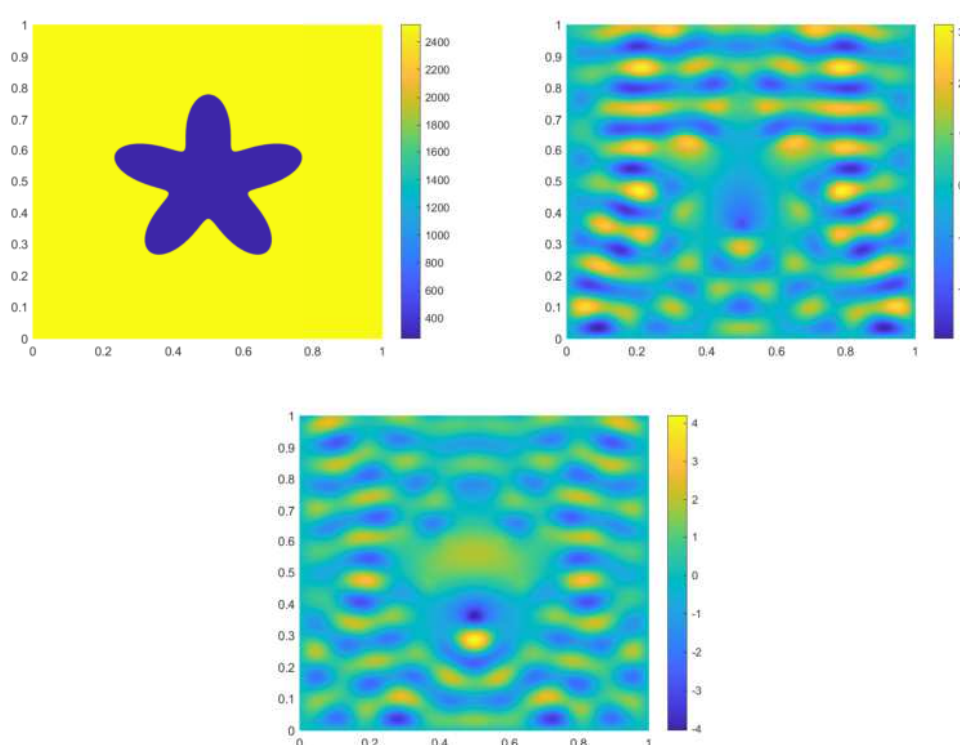

Figure 1: Model problem (1) with a variable $\kappa$.

## 4 In-Progress Work

We also report our progress on numerically solving the model problem in (2) by means of a wavelet Galerkin method. A direct construction of a Riesz wavelet basis on a general domain is known to be challenging both in theory and in practice. However, for special domains such as the one in (2), we can circumvent the need for such a construction by applying a domain transformation technique combined with recasting the problem in polar coordinates. This enables us to use the approach developed in [4] together with a tensor product. Moreover, the overall framework discussed in the previous section carries over. In a sense, our current work takes a step toward using wavelet Galerkin methods for the Helmholtz equation on nonrectangular domains, while keeping the construction of the Riesz wavelet basis on the bounded domain as simple as possible.

# High-order simulation of elasto-acoustic wave propagation using hybrid nonconforming methods (HHO) on general meshes

**Romain Mottier**[1,*], **Alexandre Ern**[2], **Laurent Guillot**[3]
[1]Department of Mathematics and Computer Science, University of Basel, Switzerland
[2]CERMICS, ENPC, Institut Polytechnique de Paris, and SERENA Project-Team, INRIA
[3]CEA, DAM, DASE
*Email: romain.mottier@unibas.ch

## Abstract

Owing to the complexity of geological structures, the simulation of elasto-acoustic wave propagation in the Earth remains a challenging task. To address this difficulty efficiently, we consider hybrid nonconforming methods (HHO methods), which offer several attractive features, including local conservativity, geometric flexibility through the support of general polyhedral meshes, and high-order accuracy [1].

## 1 Introduction and model problem

Let us consider a domain $\Omega$ decomposed into two disjoint polyhedral subdomains, the fluid domain $\Omega^{\mathrm{S}}$ and the solid domain $\Omega^{\mathrm{F}}$, sharing the interface $\Gamma := \partial\Omega^{\mathrm{F}} \cap \partial\Omega^{\mathrm{S}}$.

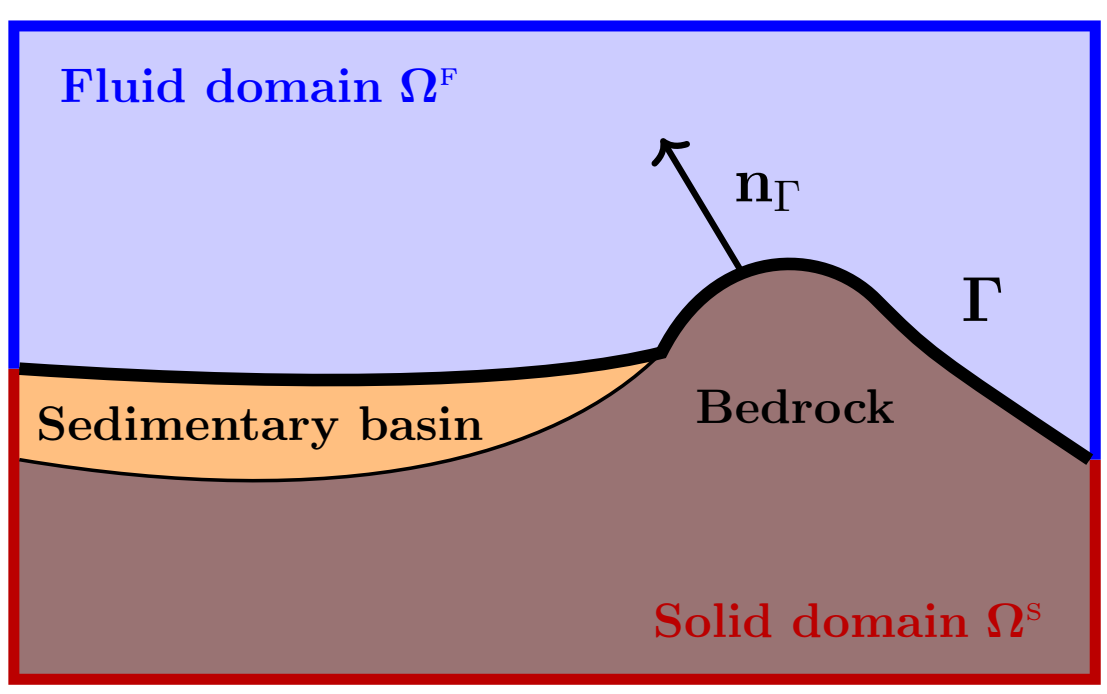

Figure 1: Example of computational domain: solid subdomain $\Omega^{\mathrm{S}}$ (sedimentary basin and bedrock) and fluid subdomain $\Omega^{\mathrm{F}}$.

We consider a first-order formulation in time and space, following the structure of a time-dependent Friedrichs system.

- Acoustic waves are governed by the pressure $p$ and velocity $\boldsymbol{v}^{\mathrm{F}}$:

$$\rho^{\mathrm{F}}\partial_t \boldsymbol{v}^{\mathrm{F}} - \nabla p = \mathbf{0}, \tag{1.1a}$$
$$\frac{1}{\kappa}\partial_t p - \nabla\cdot \boldsymbol{v}^{\mathrm{F}} = f^{\mathrm{F}}. \tag{1.1b}$$

- Elastic waves are described by the stress $\mathbb{s}$ and velocity $\boldsymbol{v}$:

$$\mathbb{C}^{-1}\partial_t \mathbb{s} - \nabla_{sym}\boldsymbol{v}^{\mathrm{S}} = \mathbb{0}, \tag{1.2a}$$
$$\rho^{\mathrm{S}}\partial_t \boldsymbol{v}^{\mathrm{S}} - \nabla\cdot\mathbb{s} = \boldsymbol{f}^{\mathrm{S}}. \tag{1.2b}$$

Here, $\rho^{\mathrm{F}}$ and $\rho^{\mathrm{S}}$ denote the fluid and solid densities, $\kappa$ the fluid bulk modulus, and $\mathbb{C}$ the Hooke tensor. The material parameters are assumed piecewise constant in $\Omega^{\mathrm{S}}$.

This formulation enables high-order discretizations in space (HHO) and in time (Runge–Kutta), leading to fully discrete schemes that are efficient in both explicit and implicit settings.

## 2 General meshes

In the context of general meshes, a cell is characterized by its number of faces rather than its geometrical shape. HHO methods naturally accommodate with this type of, as each cell is treated as a polygon with an arbitrary number of faces.

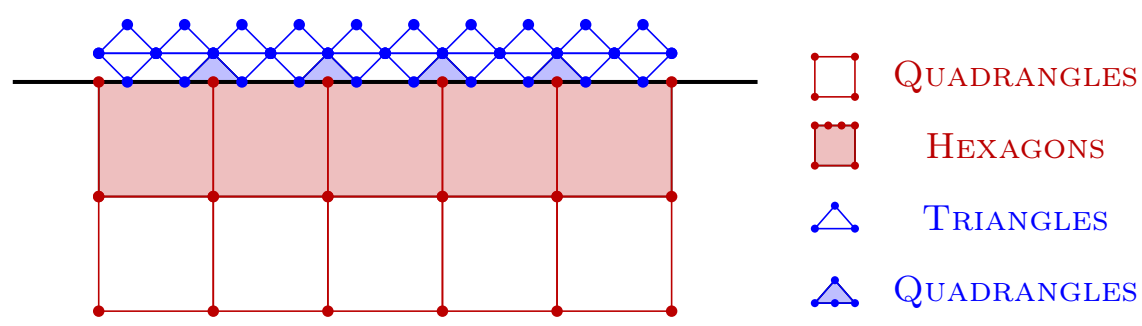

Figure 2: Example of general meshes

For instance, in Figure 2 triangles with a hanging node on one edge are treated as quadrangles (four faces), and quadrangles with two hanging nodes as hexagons (six faces).

## 3 HHO space semi-discretization

The coupled problem studied here involves primal variables ($p$, $\boldsymbol{v}$) discretized with HHO and dual variables ($\boldsymbol{m}$, $\mathbb{s}$) discretized with a standard dG aproach. HHO unknowns combine a cell polynomial of degree $k' \in \{k, k+1\}$ and a face polynomial of degree $k$, while dG unknowns are cell-based polynomials of degree $k$.

Then, the HHO discretization is formulated locally using the following two key operators:

- a local gradient reconstruction operator whose design mimics an integration by part (approximates the continuous gradient),

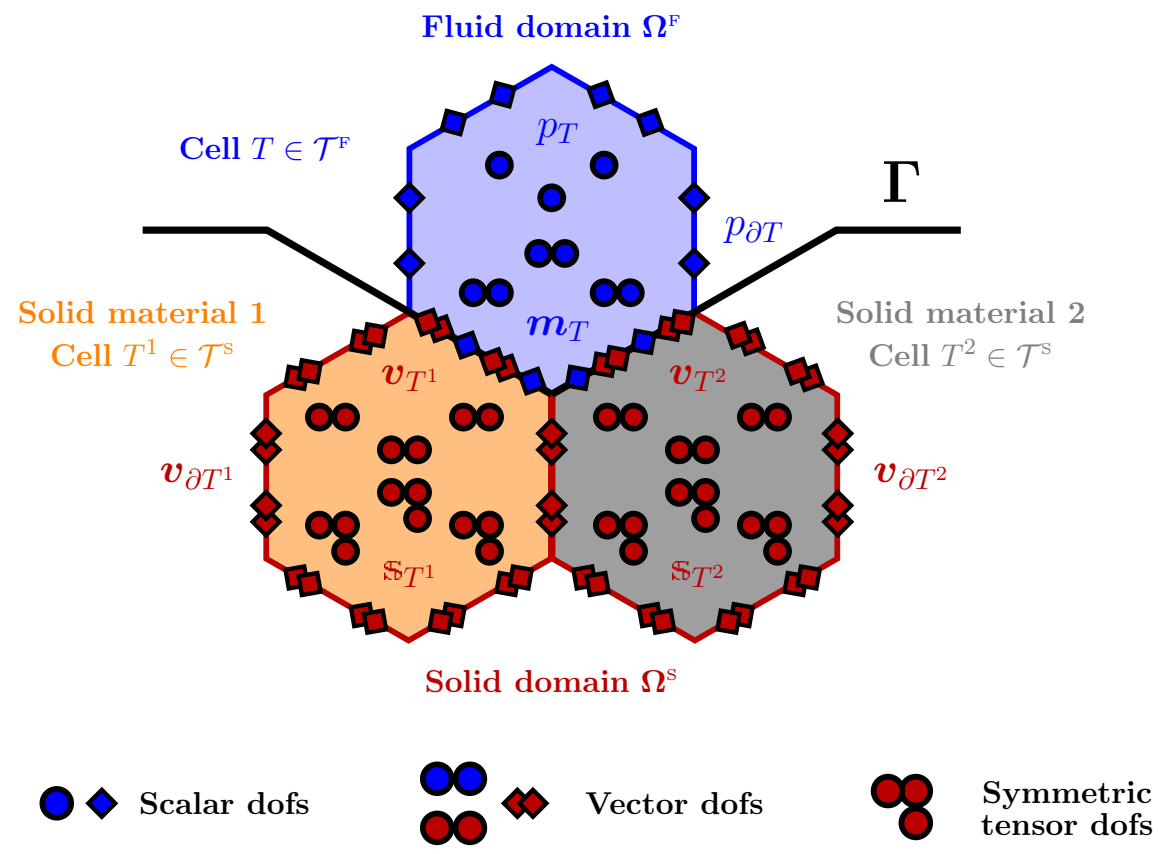

Figure 3: Elasto-acoustic dofs along a solid-fluid interface, with hexagonal mesh cells and lowest equal-order discretization ($k' = k = 1$).

- a stabilization operator which weakly imposes the continuity between the trace of the cell unknowns and the face unknowns.

## 4 Fully discrete problem & static condensation

In the explicit setting (ERK), the face unknowns are eliminated through a dG rewriting, whereas in the implicit time discretization (SDIRK) static condensation is performed to eliminate the cell unknowns, leading to increased computational efficiency compared to classical dG methods.

Both approaches are presented in detail and enriched with a local time-stepping version.

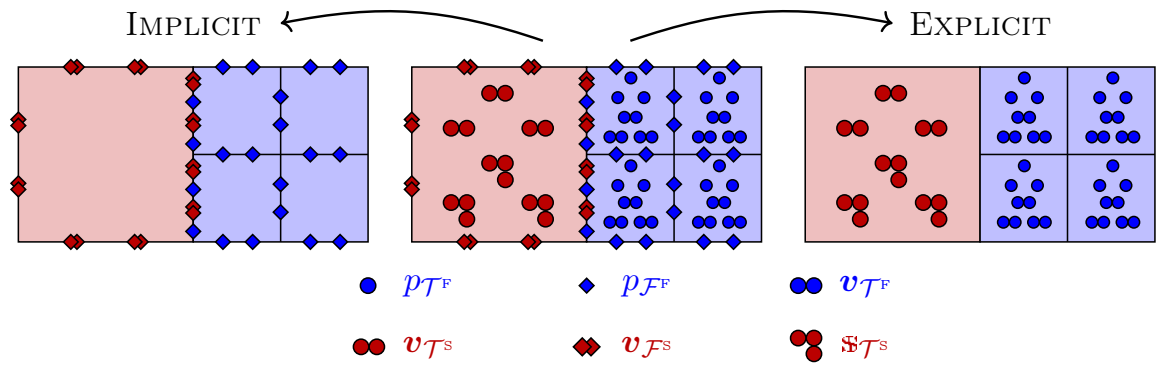

Figure 4: Fully coupled unknowns after static condensation for SDIRK schemes (left) and for ERK schemes (right) for the lowest equal-order discretization ($k' = k = 1$).

## 5 Results

Numerical experiments on test cases with analytical and semi-analytical solutions show that the method delivers optimal convergence rates [4].

Efficiency test cases indicate that the implicit discretization remains competitive with respect to the explicit one.

A numerical stability study of the local time-stepping approach shows that the method remains stable independently of the size of the local time steps used in the fine region, provided that the global time step satisfies $\Delta t \propto h_{max}$.

Finally, a more realistic case involving the propagation of an elastic pulse in a sedimentary basin coupled to the atmosphere is presented [3].

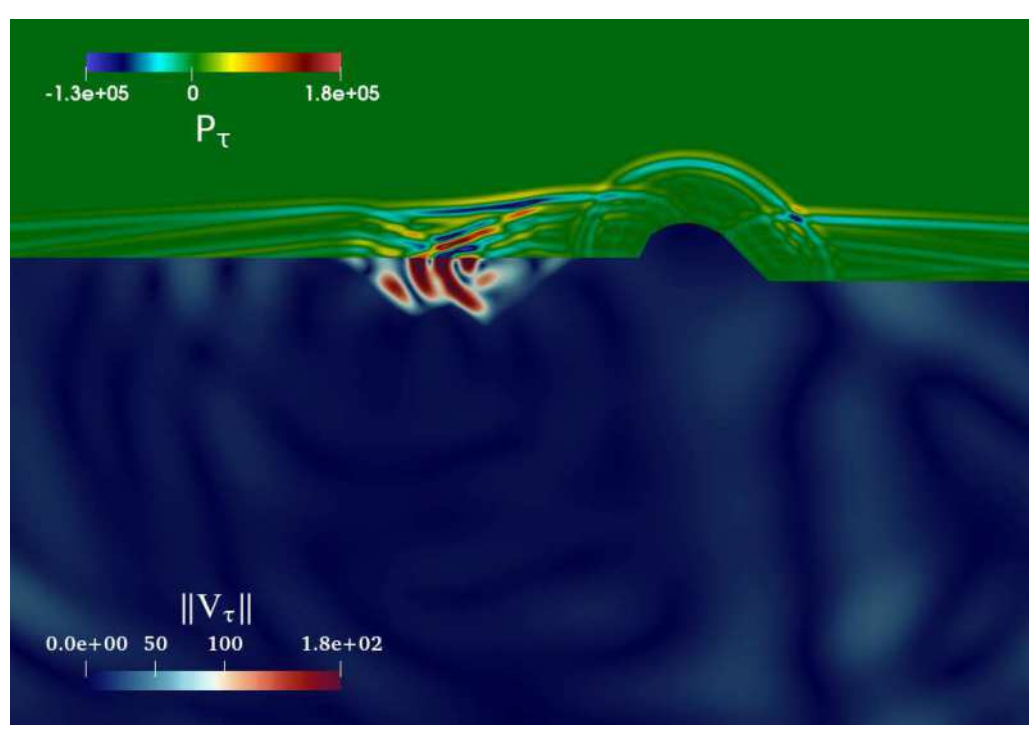

Figure 5: Spatial distribution of the acoustic pressure (atmosphere) and the elastic velocity norm (solid subdomain).

A part of the seismic energy leaks into the atmosphere where the mountain is located; another part is trapped in the sedimentary basin, which continuously emits waves to the atmosphere.

# Towards a Scale-Bridging Time Discretization Scheme for the Acoustic Wave Equation

**Dustin Mühlhäuser**[1,*], **Roland Maier**[1]
[1]Institute for Applied and Numerical Mathematics, Karlsruhe Institute of Technology, Karlsruhe, Germany
*Email: dustin.muehlhaeuser@kit.edu

## Abstract

We propose a multiscale method for the simulation of wave propagation through in space and time strongly heterogeneous media. Based on a variational ansatz in time, coarse-scale basis functions in time are enriched by fine-scale corrections, yielding an adapted coarse-scale solution that is exact at all coarse time steps w.r.t. a fine-scale reference solution. The method exploits that the corrections in time are by construction local, enabling parallelization and making otherwise possibly infeasible computations tractable. Further, we use a classical multiscale method for the spatial discretization which yields a reduction of dimension in space. Additionally, we outline a preliminary integration of deep learning to efficiently approximate spatial variations across time steps.

## Setup

We consider the acoustic wave equation with heterogeneous coefficients in a space-time cylinder $Q = \Omega \times [0, T]$ with zero initial conditions. That is, we seek a solution $u : Q \to \mathbb{R}$ such that for a given $g \in L^2(0, T; L^2(\Omega))$

$$\left.\begin{aligned} \partial_t(B(x,t)\partial_t u(x,t)) - \nabla\cdot(A(x,t)\nabla u(x,t)) & \\ = g(x,t) \quad &, (x,t) \in Q, \\ u(x,0) = 0, \quad &, x \in \Omega, \\ \partial_t u(x,0) = 0, \quad &, x \in \Omega. \end{aligned}\right\} \quad (1)$$

We assume that the coefficients $A, B : Q \to \mathbb{R}_{>0}$ vary strongly in space and time, which is why a suitable multiscale method is in demand.

We combine a *Localized Orthogonal Decomposition* (LOD)-based spatial multiscale approach, see [3], with a suitable variational multiscale method in time. For the variational formulation of (1) we introduce the following Bochner spaces. We define $H^1_{0,}(0, T; L^2(\Omega)) := \{u \in L^2(Q) : \partial_t u \in L^2(Q)$ and $u(x, 0) = 0$ for $x \in \Omega\}$ and $\mathcal{U} := L^2(0, T; H^1_0(\Omega)) \cap H^1_{0,}(0, T; L^2(\Omega))$, and, analogously, $\mathcal{V} := L^2(0, T; H^1_0(\Omega)) \cap H^1_{,0}(0, T; L^2(\Omega))$ for functions with $v(x, T) = 0$ for $x \in \Omega$. These are indeed Hilbert spaces, see [4] for further details. Now, the weak formulation of (1) reads as follows: seek $u \in \mathcal{U}$ such that

$$c(u, v) = \int_{[0,T]} (g, v)_{L^2(\Omega)} \,\mathrm{d}t \quad (2)$$

for all $v \in \mathcal{V}$, where the bilinear form is given by

$$c(u, v) := \int_{[0,T]} \int_{\Omega} -B\partial_t u \partial_t v + (A\nabla u)\cdot\nabla v \,\mathrm{d}x\,\mathrm{d}t.$$

One can show as in [2, Chapter IV, Thm. 3.1 and Thm. 3.2], that there exists a unique solution $u \in \mathcal{U}$ to (2) for sufficiently smooth coefficients $A, B$, which are uniformly bounded from above and below in space and time.

## Model and Ideas

**Spatial discretization.** We use the above mentioned LOD framework based on a Galerkin *finite element* (FE) spatial discretization of (2). Let $V_H = V_H(\Omega) \subset H^1_0(\Omega)$ denote the coarse-scale first-order FE space of dimension $N = \dim V_H$, and let $W_H \subsetneq H^1_0(\Omega)$ denote the, w.r.t. the spatial part of the bilinear form, orthogonal complement of $V_H$ that contains the fine-scale functions, called *corrections*. The idea of the LOD method is to compute the corrections via a, in this case time-dependent, (localized) correction operator $\mathcal{C}(t) : H^1_0(\Omega) \to W_H$. The complement $\tilde{V}_H(t) := (1 - \mathcal{C}(t))V_H$ of $W_H$ is spanned by the adapted LOD-basis functions $\{\Psi_1(t; x), \ldots, \Psi_N(t; x)\}$, which operate on the coarse scale but contain by construction the fine-scale information. Applied to (2), this leads to the semi-discrete problem of seeking a time-dependent vector $\mathbf{u}(t) \in U_N := \{\mathbf{u}(t) \in H^1((0, T; \mathbb{R}^N) : \mathbf{u}(0) = \mathbf{0}_N\}$ such that

$$c_N(\mathbf{u}, \mathbf{v}) = \int_{[0,T]} [\mathbf{g}(t)]^\mathsf{T} \mathbf{v}(t) \,\mathrm{d}t \quad (3)$$

for all $\mathbf{v}(t) \in V_N \subset H^1(0,T;\mathbb{R}^N)$, where the semi-discrete bilinear form is given by

$$c_N(\mathbf{u},\mathbf{v}) := \int_{[0,T]} -\partial_t \mathbf{u}(t)^\mathsf{T}\mathbb{M}(t)\partial_t\mathbf{v}(t) + \mathbf{u}(t)^\mathsf{T}\mathbb{S}(t)\mathbf{v}(t)\,\mathrm{d}t,$$

whereas the entries of the spatial mass and stiffness matrix are given for all $i,j = 1,\ldots,N$, by

$$[\mathbb{M}(t)]_{i,j} = \int_\Omega B(x,t)\Psi_j(t;x)\Psi_i(t;x)\,\mathrm{d}x,$$

$$[\mathbb{S}(t)]_{i,j} = \int_\Omega (A(x,t)\nabla\Psi_j(t;x))\cdot\nabla\Psi_i(t;x)\,\mathrm{d}x.$$

Note that both matrices possibly vary strongly in time due to the assumed high variations of the coefficients $A$ and $B$. On one hand, there is the direct dependency on the coefficients, and on the other hand this results from time-dependent spatial corrector $\mathcal{C}(t)$ of $V_H$ which yields for every point in time possibly different adapted spatial basis-functions $\Psi_1(t;x),\ldots,\Psi_N(t;x)$ of $\tilde{V}_H(t)$. Therefore, a fine temporal discretization is required to resolve the variations and, theoretically, at each fine time step an entirely new set of spatial corrections must be computed. One possibility to alleviate the computational bottleneck is to replace the expensive spatial corrector problems with a pre-trained neural network [1].

**Time discretization.** Let $\mathcal{I}_\tau$ and $\mathcal{I}_\mathcal{T}$ denote uniform fine- and coarse-scale partitions of the interval $[0,T]$ with step size $\tau \ll \mathcal{T}$, where the fine step-size is chosen such that it resolves the fine variations in time. Further, let $I_\mathcal{T}^{\mathrm{nod}}$ denote the componentwise, coarse-scale nodal interpolation operator and let $U_\mathcal{T} := \mathcal{P}_1(\mathcal{I}_\mathcal{T})^N \cap U_N$, and analogously $V_\mathcal{T}$, denote the coarse-scale first-oder FE ansatz- respectively test space in time to discretize (3).

**Adapted time-stepping method.** The idea of our method is to construct an adapted coarse-scale solution to (3) with the property that it is nodally exact in time, i.e., it coincides in all coarse time-steps with a fine-scale reference solution. We enrich the test-space $V_N$ with temporal fine-scale information in the spirit of the, originally spatial, LOD framework as follows. We define the temporal LOD-corrector $\mathcal{D} : H^1(0,T;\mathbb{R}^N) \to W_N$ as the $c_N$-orthogonal projection onto $W_N := \ker(I_\mathcal{T}^{\mathrm{nod}}|_{V_N}) \subsetneq V_N$. If we now choose as test space

$$\tilde{V}_\mathcal{T} := (1-\mathcal{D})V_\mathcal{T},$$

the solution $u_\mathcal{T} \in U_\mathcal{T}$ of the discretized version of (3) is by construction nodally exact in time. As above, the corrector $\mathcal{D}$ is computed numerically by solving temporal corrector problems. These are by construction local in time. However, we observe that the corrections are generally fully coupled across the spatial discretization points. Therefore, further modifications that employ redundancies in the construction are desirable. This, as well as the above mentioned replacement of the spatial corrector problems with a pre-trained neural network, is subject to current studies.

**Numerical example.** To illustrate the desired nodal exactness of our method, we consider a simple academic example of (3) where the dimension of the spatial discretization is reduced to 1 and the coefficients are chosen as $A = B = 1$. For the right-hand side $g(t) = 100\cos(6\pi t)$ and a time-step size of $\tau = 1/256$ and $\mathcal{T} = 1/16$, one can see the solution on the whole time interval $[0,1]$ (left) and a zoomed-in part around two coarse steps to highlight the effect. The observed nodal exactness transfers to higher-dimensional spatial discretization spaces.

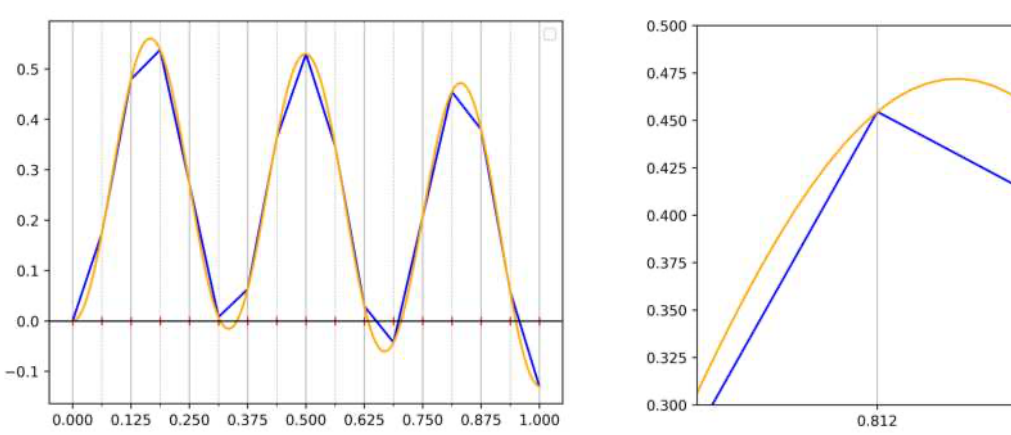

Figure 1: The solution of the adapted method (blue) compared to the reference solution (orange) are identical at the coarse time-steps.

## Stability of a passive inverse parameter problem with applications to helioseismology

**Björn Müller**[1,*], **Thorsten Hohage**[2]
[1]Max-Planck Institute for Solar System Research, Göttingen, Germany
[2]Institut für Numerische und Angewandte Mathematik, University of Göttingen, Göttingen, Germany
*Email: muellerb@mps.mpg.de

### Abstract

We discuss the stability of a passive imaging problem in a spherically symmetric background medium occurring in global helioseismology. The data of this inverse problem consists of correlations of randomly excited acoustic waves recorded on a hypersurface enclosing the solar interior. The uniqueness of interior parameters and the source strength, up to an unavoidable Gauge transformation, will be discussed in a separate talk by Thorsten Hohage.

### 1 Introduction

Acoustic waves in the solar interior are described by a time-harmonic wave equation and are randomly excited. Helioseismology is the scientific discipline that uses observations of these oscillations at different observation heights to identify coefficients of the wave equation from independent measurements on the visible part of the solar surface. Despite its great success in inferring the solar interior, the uniqueness and the stability of this passive imaging problem remain open questions.

### 2 Forward and inverse problem

Let $\Omega_0 \subset \Omega \subset \mathbb{R}^3$ be smooth and bounded domains (e.g., solar interior and the end of the solar photosphere), and let $\Gamma \subset \Omega \setminus \overline{\Omega_0}$ be the hypersurface on which measurements are available. The wavefield $\psi$ is excited by random sources $s$ with $\mathbb{E}[s] = 0$ and finite second moments, according to

$$\mathcal{L}_q\psi := \left(-\Delta + q - k^2\right)\psi = s \qquad \text{in } \Omega, \tag{1a}$$

$$\partial_{\boldsymbol{n}}\psi = \mathcal{B}\psi \qquad \text{on } \partial\Omega, \tag{1b}$$

where $q \in L^\infty(\Omega, \mathbb{C})$ with $\mathrm{supp}(q) \subset \Omega_0$, $\boldsymbol{n}$ is the outer normal on $\partial\Omega$, $k \in \mathbb{C}$, and $\mathcal{B} : H^{1/2}(\partial\Omega) \to H^{-1/2}(\partial\Omega)$ is the boundary condition satisfying $\Im(\mathcal{B}) > 0$ and $\Re(\mathcal{B}) \geq 0$. The boundary condition ensures unique solvability of (1) for all $q \in \mathbb{X}$, where

$$\mathbb{X} := \{\psi \in L^\infty(\Omega, \mathbb{C}), \mathrm{supp}(\psi) \subset \Omega_0, \Im(\psi) < 0\}.$$

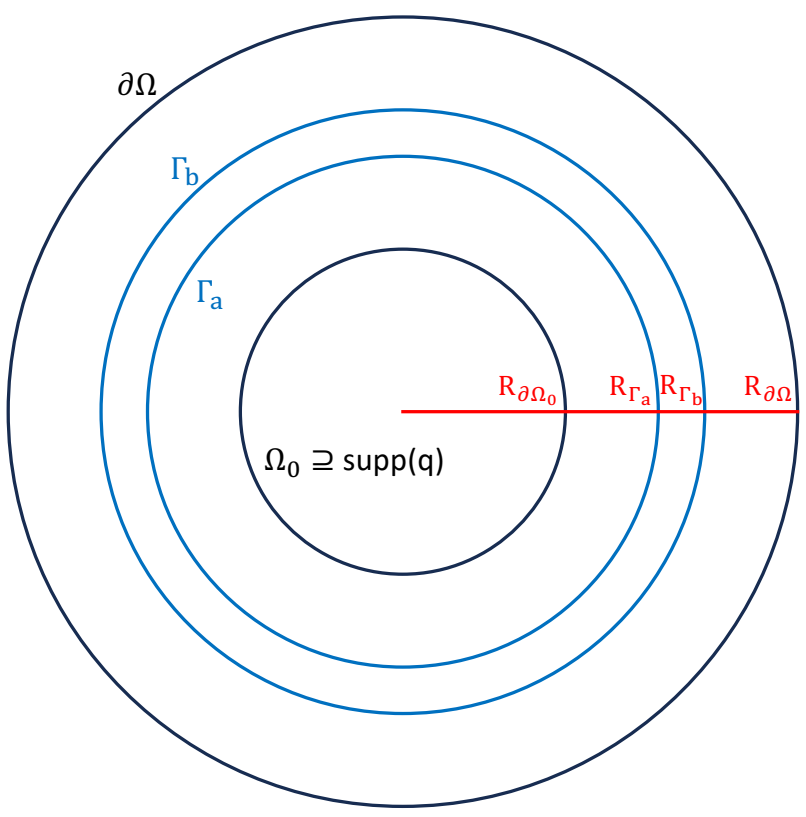

Figure 1: Sketch of the geometry.

Solutions to (1) may be written in the form

$$\psi = \mathcal{G}_q s = \int_\Omega G_q(\cdot, \boldsymbol{y}) s(\boldsymbol{y})\, \mathrm{d}\boldsymbol{y},$$

where $G_q$ denotes the Green's function to (1) and $\mathcal{G}_q$ the corresponding volume potential operator. Since $\mathbb{E}[s] = 0$, we study the cross-correlation of solutions $\psi$ to the wave equation (1) between two points $\boldsymbol{r}_1, \boldsymbol{r}_2 \in \Gamma$ for independent samples $s_1, \ldots, s_N$

$$\mathrm{Corr}^{\Gamma\Gamma}(\boldsymbol{r}_1, \boldsymbol{r}_2) := \frac{1}{N}\sum_{k=1}^{N} \mathrm{Tr}_\Gamma\, \psi_k(\boldsymbol{r}_1) \cdot \overline{\mathrm{Tr}_\Gamma\, \psi_k(\boldsymbol{r}_2)},$$

where $\psi_k = \mathcal{G}_q s_k$ and $\mathrm{Tr}_\Gamma\, \psi := \psi|_\Gamma$. The cross-covariance $\mathbb{E}\left[\mathrm{Corr}^{\Gamma\Gamma}\right]$ is the integral kernel of the covariance operator

$$\mathcal{C}^{\Gamma\Gamma}_{q,\boldsymbol{Cov}[s]} := \mathrm{Tr}_\Gamma\, \mathcal{G}_q\, \boldsymbol{Cov}[s] \mathcal{G}_q^*\, \mathrm{Tr}_\Gamma^*,$$

where the source covariance $\boldsymbol{Cov}[s] : H^1(\Omega) \to H_0^{-1}(\Omega)$ is Hilbert-Schmidt. The inverse parameter problem consists of determining $q \in \mathbb{X}$ from observations of $\mathcal{C}^{\Gamma\Gamma}_{q,\boldsymbol{Cov}[s]}$.

### 3 Stability in global helioseismology

Suppose that measurements are available on two separated hypersurfaces in the solar photosphere,

meaning that $\Gamma = \Gamma_a \cup \Gamma_b$ with disjoint $\Gamma_a = \partial\Omega_a$ and $\Gamma_b = \partial\Omega_b$. Furthermore, we assume spherical symmetry, meaning that $\Omega$, $\Omega_b$, $\Omega_a$ and $\Omega_0$ are balls centered at 0 with radii $R_{\partial\Omega} > R_{\Gamma_b} > R_{\Gamma_a} > R_{\partial\Omega_0}$ (see Fig. 1). In helioseismology, measurements at two observation heights can be obtained by analyzing the motion of different spectral lines formed at different heights. We assume that there are separated sources in the solar interior and the atmosphere, and the source covariance takes the form

$$\boldsymbol{Cov}[s] = M_S + \boldsymbol{Cov}[s^{\Omega\setminus\Omega_b}], \quad M_S\psi = S\cdot\psi,$$

where $S \in L^\infty(\Omega, \mathbb{R}_{\geq 0})$ with $\mathrm{supp}(S) \subset \Omega_0$ and $\boldsymbol{Cov}[s^{\Omega\setminus\Omega_b}]$ is positive definite on $\{\psi \in H^1(\Omega \setminus \overline{\Omega_b}) : (\Delta + k^2)\psi = 0\}$.

In global helioseismology, one is interested in the global (spherically symmetric) structure of the solar interior such that $q(\boldsymbol{r}) = q(|\boldsymbol{r}|)$, $S(\boldsymbol{r}) = S(|\boldsymbol{r}|)$ and the exterior covariance operator $\boldsymbol{Cov}[s^{\Omega\setminus\Omega_b}]$ is spherically symmetric. The stability result will be based on the decay of the eigenvalues $(\lambda_\ell)_{\ell\in\mathbb{N}}$ of $\mathcal{N}^*\,\boldsymbol{Cov}[s^{\Omega\setminus\Omega_b}]\mathcal{N}$, where

$$\mathcal{N} : H^{1/2}(\Gamma_b) \to H^1(\Omega \setminus \overline{\Omega_b}), \quad g \mapsto u$$

and $u$ is the unique solution to

$$\begin{cases} (\Delta + k^2)u = 0, & \text{in } \Omega \setminus \overline{\Omega_b}, \\ \partial_{\boldsymbol{n}} u = \mathcal{B}u, & \text{on } \partial\Omega, \\ u = g & \text{on } \Gamma_b\,. \end{cases}$$

For some proper $C_N > 0$, we define the function:

$$F(x) \mapsto x \cdot \inf_{i \leq N(x)} \lambda_i, \quad N(x) := C_N x^{-\frac{1}{2}}\,.$$

It can be shown that $F : [0, \infty] \to [0, \infty]$ is invertible and strictly monotonically increasing.

**Theorem 1** *Suppose that* $0 < s < \frac{3}{2}$. *There exists constants* $C_\omega, C_F > 0$ *and* $0 < \nu^* < 1$ *dependent on* $s$, *and it yields*

$$\|q_1 - q_2\|_{H^{-s}(\Omega_0)} \leq C_\omega \left|\ln\left[F^{-1}\left(C_F\nu\right)\right]\right|^{\frac{-2s}{3}}, \tag{2}$$

*where* $\nu = \|\mathcal{C}^{\Gamma\Gamma}_{q_1,\boldsymbol{Cov}[s_1]} - \mathcal{C}^{\Gamma\Gamma}_{q_2,\boldsymbol{Cov}[s_2]}\| \leq \nu^*$.

The crucial step is to prove the estimate

$$\|\mathcal{S}^{\Gamma_b\Gamma_b}_{q_1} - \mathcal{S}^{\Gamma_b\Gamma_b}_{q_2}\| \leq F^{-1}\left(C_F\|\mathcal{C}^{\Gamma\Gamma}_{q_1,\boldsymbol{Cov}[s_1]} - \mathcal{C}^{\Gamma\Gamma}_{q_2,\boldsymbol{Cov}[s_2]}\|\right),$$

where $\mathcal{S}^{\Gamma_b\Gamma_b}_q$ is the single layer integral operator. Importantly, this inequality connects the passive problem with the active problem. The proof is completed by the well-established theory on the Calderón problem (e.g., [2]).

## 4 Convergence rates

Theorem 1 allows us to quantify the stability of the solution to the inverse problem against the data noise level.

**Theorem 2** *Let* $\tilde{\mathrm{Corr}}^{\Gamma\Gamma}$ *be the integral operator with integral kernel* $\mathrm{Corr}^{\Gamma\Gamma}$ *and* $\psi_s$ *the concave envelope in a neighborhood at* 0 *for the right-hand side of* (2). *Let* $0 < s < \frac{3}{2}$, $\delta > 0$ *and suppose that* $\|\mathcal{C}^{\Gamma\Gamma}_{q^\dagger,\boldsymbol{Cov}[s]} - \tilde{\mathrm{Corr}}^{\Gamma\Gamma}\| \leq \delta$ *for the ground truth* $q^\dagger$ *satisfying* $q^\dagger \in H^s(\Omega)$. *For the parameter choice* $\alpha \sim \frac{\delta}{\psi_s(s\delta)}$, *it yields*

$$\|q^\dagger - q^\delta_\alpha\|_{L^2(\Omega_0)} = \mathcal{O}\left(\psi_s(\delta)\right), \tag{3}$$

*where* $q^\delta_\alpha$ *is the unique minimizer to the Tikhonov functional:*

$$\|\mathcal{C}^{\Gamma\Gamma}_{q,\boldsymbol{Cov}[s]} - \tilde{\mathrm{Corr}}^{\Gamma\Gamma}\|^2 + \alpha\|q\|^2_{L^2(\Omega_0)}\,.$$

## 5 Application to helioseismology

Motivated by observations, one often assumes that $\mathcal{C}^{\Gamma\Gamma}_{q,\boldsymbol{Cov}[s]} = \frac{1}{2\mathrm{i}}\left[\mathcal{S}^{\Gamma\Gamma}_q - \overline{\mathcal{S}^{\Gamma\Gamma}_q}\right]$, realized by sources on the computational domain balancing the incoming and outgoing energy flux. In this case, it follows from Theorem 2, that there exists constants $C_1 > 0$ and $C_2 \in \mathbb{R}$, such that

$$\psi_s(\delta) = C_1\left|\ln\left[C_2 - \ln(\delta)\right]\right|^{-\frac{2s}{3}}\,.$$

A different popular source model extends turbulent motions to the solar photosphere. The source covariance can be modeled by $\boldsymbol{Cov}[s^{\Omega\setminus\Omega_b}] = M_{S^{\mathrm{ext}}}$ with $S^{\mathrm{ext}} \in C^\infty\left((R_{\Gamma_b}, R_{\partial\Omega})\right)$. If we furthermore assume that $S^{\mathrm{ext}}(R_{\Gamma_b}) \neq 0$, it follows that there exists $C > 0$, such that

$$\psi_s(\delta) = C\left|\ln(\delta)\right|^{-\frac{2s}{3}}\,.$$

The noise level scales with the observation time $T_{\mathrm{obs}}$ like $\delta \propto T_{\mathrm{obs}}^{-1/2}$. Therefore, the results allow us to quantify the achievable temporal resolution for inversions of different quantities within different solar background models.

# Imaging in two-dimensional waveguides with complex geometry

**Gregory Mwamba**[1,*], **Chrysoula Tsogka**[2], **Symeon Papadimitropoulos**[3]
[1]Department of Applied Mathematics, University of California, Merced, Merced, USA
[2]Department of Applied Mathematics, University of California, Merced, Merced, USA
[3]Department of Applied Mathematics, University of California, Merced, Merced, USA
*Email: gmwamba@ucmerced.edu

**Abstract**

We consider the problem of imaging extended scatterers in waveguides with complex geometries. The data consist of the array response matrix, obtained by transmitting signals from each array element and recording the responses at all array elements. For imaging, we propose a method based on back-propagating the data after projecting it onto the propagating modes of an infinite waveguide of constant depth, equal to the depth at the array location. We show that this approach allows us to recover information about both the scatterer and the waveguide geometry. The recovered geometric information can then be incorporated to improve the imaging of the scatterer.

## 1 Formulation of the problem

We consider a two-dimensional (2D) Cartesian coordinate system $\vec{\boldsymbol{x}} = (z, x) \in \mathbb{R}^2$, where $z$ denotes the range coordinate and $x$ the cross-range coordinate, with the downward direction taken to be positive. The computational domain $\Omega$, shown in Figure 1, is a two-dimensional waveguide that is unbounded in the range direction and has variable depth. Specifically, the waveguide consists of three regions: $\Omega_{L^-}$, which has constant depth $D$ for $z \leq L$; a smoothly varying section, denoted by $\Omega_L$, where the depth transitions from $D$ to $D_m > D$ for $L < z < L_m$; and $\Omega_{L_m^+}$, which has constant depth $D_m$ for $z \geq L_m$.

Wave propagation is governed by the Helmholtz equation,

$$-\Delta u(\vec{\boldsymbol{x}}, \omega) - k^2 \eta(\vec{\boldsymbol{x}}) u(\vec{\boldsymbol{x}}, \omega) = f(\vec{\boldsymbol{x}} - \vec{\boldsymbol{x}}_s, \omega), \quad (1)$$

subject to homogeneous Dirichlet boundary conditions at the waveguide's boundaries and appropriate radiation conditions at infinity. Here, $\eta(\vec{\boldsymbol{x}})$ is the index of refraction that is equal to one in $\Omega \backslash \mathcal{O}$ and takes a different value inside the scatterer $\mathcal{O}$.

We place an active array, denoted by $\mathcal{A}$, parallel to the waveguide's cross-section at range $z = z_a > L_m$. The array consists of $N$ transducers and spans the entire cross-section of the waveguide at that range, corresponding to the local depth $D_m$ of $\Omega_{L_m^+}$. We assume that the waveguide geometry is otherwise unknown, apart from this depth at the array location.

Our objective is to locate and image an extended scatterer contained within the waveguide. Here, the term *extended* refers to a scatterer whose size is comparable to a reference wavelength $\lambda_0$.

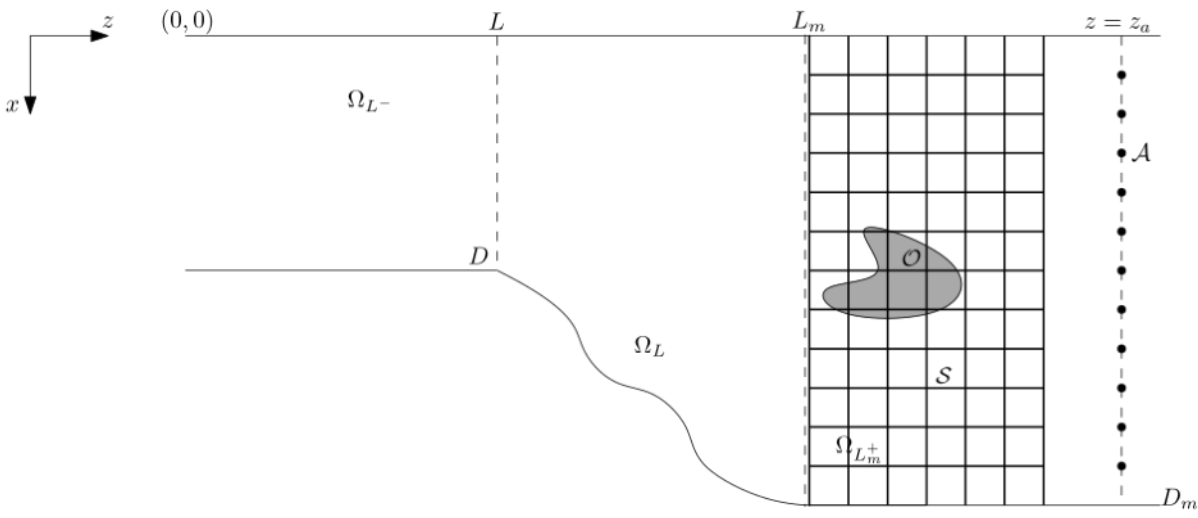

Figure 1: Imaging setup in a 2D waveguide with complex geometry

The data used for imaging consist of recordings of the total field, i.e. the field recorded at the array location when a scatterer is present. These are collected in the *array response matrix*, denoted by $\Pi(\omega) \in \mathbb{C}^{N \times N}$. The $(r, s)$-th entry of $\Pi(\omega)$ corresponds to the signal recorded at the $r$-th receiver when the $s$-th source emits a signal at a frequency $\omega$, that is the solution of (1) evaluated at $\vec{\boldsymbol{x}} = \vec{\boldsymbol{x}}_r$.

## 2 Imaging

To generate an image, we define a search domain $\mathcal{S} \subset \Omega$ and evaluate an appropriate imaging functional. In [1], an imaging approach based on projecting the response matrix onto the propagating modes was proposed for a two-dimensional infinite waveguide of constant depth. We adopt the same imaging functional in the present setting, using the propagating modes for an infinite waveguide with depth $D_m$, the depth where the array is located. Let $\{\mu_n, \Phi_n\}_{n=1}^{\infty}$ be the eigenvalues and corresponding orthonormal eigenfunctions of the Dirichlet transverse Laplacian $-\Delta$ in $[0, D_m]$. It is well-known that $\mu_n \in \mathbb{R}$ and are positive, while $\{\Phi_n\}_{n=1}^{\infty}$ forms an orthonormal basis of $L^2[0, D_m]$. Let also $M$ be the number of propagating modes in $\Omega_{L_m^+}$, i.e. $M$ is an index that satisfies $\mu_M < k^2 < \mu_{M+1}$, where $k$ is the constant wavenumber in $\Omega_{L_m^+}$. Finally, let

$$\beta_n = \begin{cases} \sqrt{k^2 - \mu_n}, & 1 \le n \le M, \\ \mathrm{i}\sqrt{\mu_n - k^2}, & n \ge M+1, \end{cases}$$

be the horizontal wavenumbers.

The imaging functional is given by

$$\mathcal{I}(\vec{\boldsymbol{y}}^{\mathrm{s}}) = \sum_{i=1}^{M}\sum_{j=1}^{M} \overline{\mathbb{Q}_{ij}}\, \mathcal{G}_i(z_a, \vec{\boldsymbol{y}}^{\mathrm{s}})\, \mathcal{G}_j(z_a, \vec{\boldsymbol{y}}^{\mathrm{s}}),$$

where $\vec{\boldsymbol{y}}^{\mathrm{s}} \in \mathcal{S}$ and $\mathcal{G}_i(z_a, \cdot)$ are the Fourier coefficients of the Green's function with respect to the orthonormal basis $\{\Phi_n\}_{n=1}^{\infty}$ of $L^2[0, D_m]$, i.e.,

$$\mathcal{G}_i(z_a, \cdot) = \beta_i \int_0^{D_m} G\big((z_a, x'), \cdot\big)\Phi_i(x')\, dx'.$$

The matrix $\mathbb{Q}$ is an $M \times M$ matrix, whose entries are the projections of the response matrix on the propagating modes.

$$\mathbb{Q}_{ij} = \beta_i\beta_j \int_{\mathcal{A}}\int_{\mathcal{A}} \Pi_{rs}(\omega)\Phi_i(x_s')\Phi_j(x_r')\, dx_s'\, dx_r'.$$

## 3 Numerical Results

We consider a 2D homogeneous waveguide with complex geometry that has constant depth $D = 10\lambda_0$, for $z \in [0, 10]\lambda_0$, a variable depth described by the equation $x = 3z - 400$ for $z \in (10, 20)\lambda_0$ and $D_m = 20\lambda_0$ for $z > 20\lambda_0$; here $\lambda_0 = 20$ m is a reference wavelength. We consider a disc-shaped scatterer $\mathcal{O}$ with radius $\lambda_0$, centered at $\vec{\boldsymbol{x}}^{\star} = (27, 16)\lambda_0$ and $\eta(\vec{\boldsymbol{x}}) = 0.36$ in $\mathcal{O}$. We place an array with $N = 81$ transducers spanning the entire depth $D_m$ at $z_a = 38\lambda_0$. Our data are for frequencies ranging from 40 Hz to 106.5 Hz with an increment of 0.83125 Hz.

For the images presented in Figure 2, the only prior information available is the waveguide's depth at the array location. Using this depth, we compute analytically the Green's function of an infinite waveguide with the corresponding constant depth. Using this background model, the imaging functional produces reasonably accurate reconstructions of both the extended scatterer shown in the left subplot and the shape of $\Omega_L$ shown in the right subplot. While the resulting image reconstructs only part of the scatterer, its location and the waveguide geometry are recovered with great accuracy.

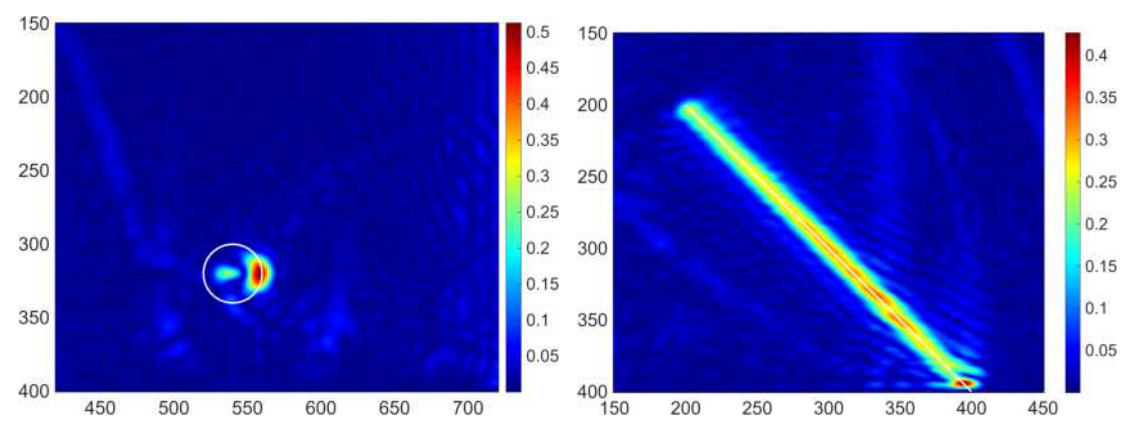

Figure 2: Imaging with the total field and the Green's function for a waveguide with constant depth $D_m$. Left: $\mathcal{I}(\vec{\boldsymbol{y}}^{\mathrm{s}})$ in the region that contains the disc-shaped scatterer plotted in white. Right: $\mathcal{I}(\vec{\boldsymbol{y}}^{\mathrm{s}})$ focused on the topography.

The estimated topography is then used to update the background Green's function and compute the incident field, which is subtracted from the total field to obtain the scattered field. This yields a noticeably improved reconstruction of the scatterer in Figure 3, with a shape that more closely matches the true object. These results illustrate that even limited initial knowledge of the waveguide geometry can be leveraged to improve both geometric estimation and scatterer imaging.

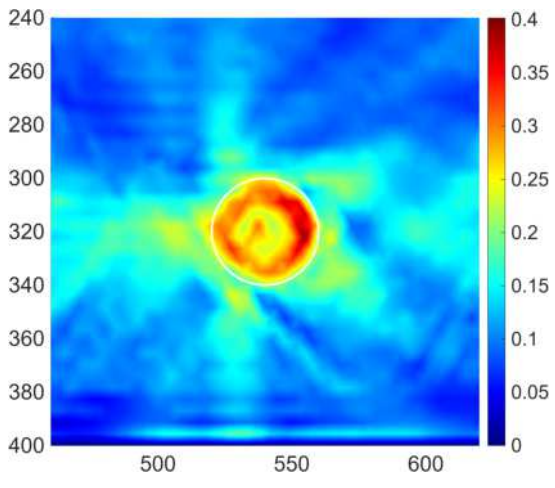

Figure 3: $\mathcal{I}(\vec{\boldsymbol{y}}^{\mathrm{s}})$ using the scattered field as data and the Green's function in the true waveguide geometry. We observe that we get a better image of the extended scatterer.

# Complex scaling for the homogeneous wave equation with scaling direction non-normal to the transparent boundary

**Lothar Nannen**[1,*], **Constantin Pierer**[1], **Markus Wess**[1]
[1]Institute of Analysis and Scientific Computing, TU Wien, Austria
*Email: lothar.nannen@asc.tuwien.ac.at

**Abstract**

In this talk, we investigate instabilities in a generalized radial perfectly matched layer (PML) for the homogeneous wave equation when the direction of the complex scaling is chosen to be non-normal to the transparent boundary. To this end, we study the spectrum both analytically and with the aid of numerical tests. It turns out that the dependency of the complex scaling on the frequency can be modified in a way that instabilities are eliminated.

## 1 Introduction

For the numerical calculation of solutions to the scalar wave equation in unbounded domains, the domain is often truncated and surrounded by a perfectly matched layer (PML). This is intended to absorb waves radiating from the interior and thus ensure that no artificial reflections occur at the outer boundary. The problem on the now restricted domain can be approximated e.g. using finite elements.

For the basic construction of a PML for time-dependent problems, please refer to [1]. Numerical simulation is performed in the time domain, while a PML is formulated and analysed in the frequency domain. There, it depends on three things: the direction of the complex scaling, the choice of scaling profile as a function of frequency, and the discretisation of the complex-scaled problem.

## 2 Direction of complex scaling

We use generalised radial coordinates. For simplicity, let $\Omega_{\text{int}} \subset \mathbb{R}^d$ with $d = 2, 3$, be a bounded, open, convex domain with boundary $\Gamma$, and let $m \in \Omega_{\text{int}}$ be a centre point. Then, $x \in \Omega_{\text{ext}} := \mathbb{R}^d \backslash \overline{\Omega_{\text{int}}}$ has a representation of the form $x(\xi, \hat{x}) = \hat{x} + \xi(\hat{x} - m)$ with a surface variable $\hat{x} \in \Gamma$ and a radial variable $\xi > 0$. We use a linear complex scaling in radial direction, i.e. for $\sigma = \sigma(\mathrm{i}\omega) \in \mathbb{C}$ we define the complex scaled coordinate by

$$x_\sigma(\xi, \hat{x}) = \hat{x} + \sigma(\mathrm{i}\omega)\xi(\hat{x} - m).$$

Note, that this scaling is non-normal to $\Gamma$ if $\Gamma$ is not a sphere with centre point $m$.

## 3 Scaling profile

We chose

$$\sigma(\mathrm{i}\omega) := \beta + \frac{\alpha}{\gamma - \mathrm{i}\omega} \tag{1}$$

with $\alpha > 0$, $\beta, \gamma \geq 0$. The following serves as a brief motivation: if a wave contains the exponential factor $\exp(\mathrm{i}\omega\xi)$, then for $\omega \in \mathbb{R}\backslash\{0\}$ there holds for the complex scaled wave

$$|\exp(\sigma(\mathrm{i}\omega)\mathrm{i}\omega\xi)| = \exp\left(-\frac{\alpha}{1 + \gamma^2/\omega^2}\xi\right). \tag{2}$$

Hence, the complex scaled wave decays exponentially for $\xi \to \infty$.

## 4 Discretisation

In frequency domain, the complex scaled differential equation contains fractions and higher order polynomials in $\mathrm{i}\omega$. Introducing auxiliary functions leads us to a system of differential equations which is linear in $\mathrm{i}\omega$. This results in time domain into a linear, real-valued system of differential equations, which is first order in time and second order in space. In time we use the Crank-Nicolson scheme. For the spatial discretisation in $\Omega_{\text{ext}}$ one could either use infinite elements as in [2, Sec.6] or a classical PML based on finite elements for a truncated problem.

## 5 Spectrum

We study in this section the spectrum of the non-discretised problem in frequency domain. If there exists $\omega \in \mathbb{C}$ with $\Im(\omega) > 0$ such that $\mathrm{i}\omega$ is part of the spectrum, then we have to expect instabilities in time-domain due to the exponentially in time $t$ increasing factor $\exp(-\mathrm{i}\omega t)$. There are two sources for such problematic frequencies.

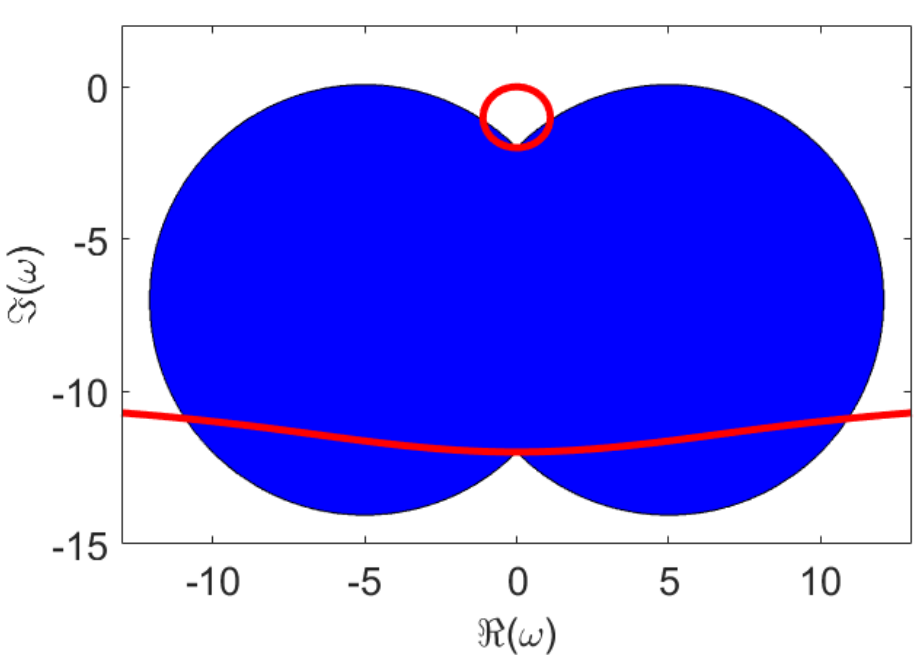

Figure 1: solutions $\omega$ to (3) (red) and (4) (blue) for $\mu = 3/4\pi$, $\alpha = 10$, $\beta = 1$, $\gamma = 2$ .

### 5.1 Spectrum related to the damping

In (2) we have noted the exponential decay in space of the complex scaled solution. For $\omega \in \mathbb{C}$ with

$$\Re\left(\sigma(\mathrm{i}\omega)\mathrm{i}\omega\right) = 0 \tag{3}$$

the complex scaled function belongs no longer to $L^2((0,\infty))$ resulting in an essential spectrum. Most optimal for this part of the spectrum would be the choice $\beta = \gamma = 0$, since in this case no solutions $\omega$ to (3) exist. $\beta > 0$ and $\gamma = 0$ leads to solutions $\omega = -\mathrm{i}\alpha/\beta$. Hence, $\Im(\omega) < 0$ and the system remains stable. If $\gamma > 0$ it can be shown that $\omega = 0$ is the only solution to (3) with non-negative imaginary part. This choice therefore results in a borderline case, but one that is just about stable, see Fig.1 for an example.

### 5.2 Spectrum related to the scaling direction

This part of the spectrum was studied in [2, Sec. 5.3.1] in a slightly different setting. In our setting, it should be possible to construct a Weyl sequence to a $\omega$ if

$$\arg(\sigma(\mathrm{i}\omega)) \in [\pi - \mu, \pi + \mu]\,, \tag{4}$$

whereas $\mu \in [\pi/2, \pi)$ represents the largest angle between the scaling direction (in our case $\hat{x} - m$) and the interface $\Gamma$. A scaling direction normal to $\Gamma$ gives the smallest angle $\pi/2$.

If $\beta = \gamma = 0$, solutions $\omega$ to (4) fulfil $\arg(\omega) \in [-\pi/2-\mu, -\pi/2+\mu]$. In other words, if the scaling direction is normal to $\Gamma$, then $\mu = \pi/2$ and the solutions to (4) are all $\omega$ with $\Im(\omega) \leq 0$. This is again a borderline case. As soon as $\Gamma$ is no longer a perfect circle with centre point $m$, $\mu > \pi/2$ and there exist solutions $\omega$ to (4) with positive imaginary part.

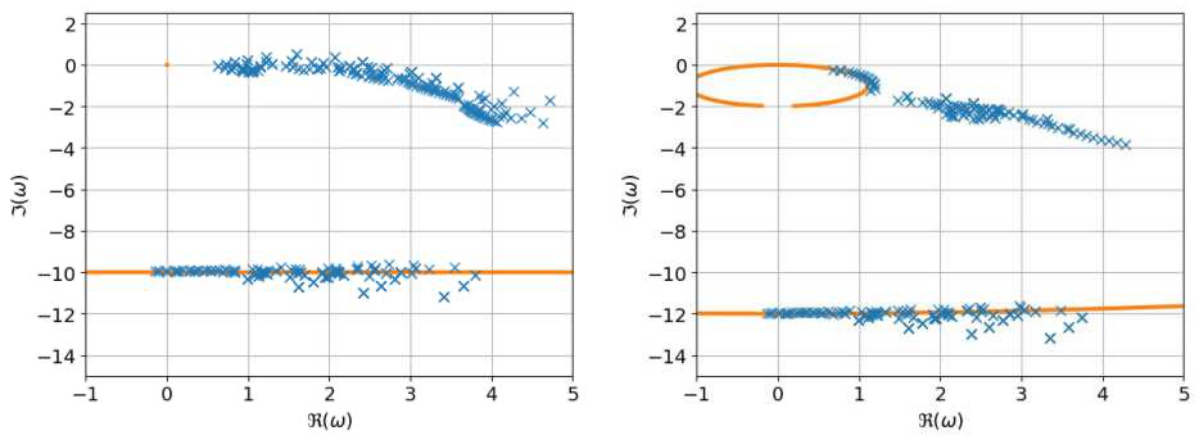

Figure 2: Some computed resonances for $\Omega_{\text{int}} = [0,1]^2$ and $m = (0.95, 0.95)$. Left $\gamma = 0$ and right $\gamma = 2$. The orange lines are the solutions to (3).

If $\beta > 0$, the solutions to (4) are a union of two disks in the complex plane, see Fig.1. This set might contain $\omega$ with positive imaginary part in particular for $\mu > \pi/2$. For $\gamma > 0$ this set is shifted by $-\mathrm{i}\gamma$ such that for sufficiently large $\gamma$ all solutions to (4) have strictly negative imaginary part.

## 6 Numerical test and conclusion

In the last section we studied the continuous setting only. Nevertheless, for sufficiently fine discretisations we should expect instabilities in time domain if the spectrum of the continuous problem contains problematic frequencies.

Fig.2 shows computed resonances $\omega$ in an extremal case, where the centre point is very close to the transparent boundary leading to $\mu$ close to $\pi$. We have used the scaling profile (1) with $\alpha = 10$, $\beta = 1$ and $\gamma = 0$ (left) and $\gamma = 2$ (right). Clearly, a discretisation of the essential spectrum related to (3) can be seen. The other resonances are related to (4). For $\gamma = 0$ there exist resonances with positive imaginary part. In this case a time simulation would be instable. For $\gamma = 2$ the resonances related to (4) move down such that all computed resonances have non-positive imaginary part.

Hence, our studies indicate that the scaling profile (1) with positive $\alpha$ and $\beta$ and sufficiently large $\gamma$ might be a good choice in particular if the scaling direction is chosen to be non-normal to the transparent boundary.

# Identifying defective units in infinite periodic arrays of point sources

**Dinh-Liem Nguyen**[1,*], **Nhung H. Nguyen**[1], **Thi-Phong Nguyen**[2]
[1]Department of Mathematics, Kansas State University, Manhattan, USA
[2]Department of Mathematical Sciences, New Jersey Institute of Technology, Newark, USA
*Email: dlnguyen@ksu.edu

**Abstract**

This work focuses on identifying defective units in unbounded periodic arrays of point sources using boundary data. Unlike classical inverse source problems in free space, the key challenge here lies in the disruption of periodicity caused by defective sources in the infinite array. To address this, we employ the Floquet–Bloch transform to reformulate the original inverse source problem as a quasi-periodic inverse source problem. We first establish uniqueness theorems for both the original and the quasi-periodic formulations. Then, we develop a numerical method for identifying the number of defective sources, their locations and intensities.

## 1 Introduction

We consider a one-dimensional array of point sources that is unbounded, periodic along the $x_1$-direction with period $L > 0$, and bounded in the $x_2$-direction between $-h$ and $h$, for some $h > 0$. Each unit cell of length $L$ contains one point source. We denote by

$$\Omega := \mathbb{R} \times (-h, h), \quad \Omega_0 := \left(-\frac{L}{2}, \frac{L}{2}\right) \times (-h, h),$$

and let $\boldsymbol{x}_0$ be the location of the source in $\Omega_0$. The locations of all other sources in the array are then given by $\boldsymbol{x}_j := \boldsymbol{x}_0 + (jL, 0)^\top$ for $j \in \mathbb{Z}$. These sources are represented by the delta distribution $\delta(x; \boldsymbol{x}_0 + (jL, 0)^\top)$ for $x \in \Omega$ and have corresponding intensities $\gamma_j$.

In this work, we assume that some sources are not functioning properly, in the sense that their intensities differ from 1. In particular, there are $N \in \mathbb{N}^*$ defective sources located at $\boldsymbol{x}_{m_\ell}$, each having intensity $\gamma_{m_\ell} \neq 1$, for $\ell = 1, \dots, N$. We denote by $\mathbf{I}$ the set of indices of all defective sources, i.e., $\mathbf{I} := \{\mathbf{m} : \gamma_\mathbf{m} \neq \mathbf{1}\}$ and assume that $\mathbf{I}$ is a finite set.

The wave generated by all point sources in a homogeneous background with a wave number $k > 0$ satisfies the equation

$$\Delta v + k^2 v = -\left(\delta_{\text{per}}(\cdot; \boldsymbol{x}_0) + \sum_{m \in \mathbf{I}} \sigma_m \delta(\cdot; \boldsymbol{x}_m)\right), \tag{1}$$

where $\sigma_m := \gamma_m - 1$ and

$$\delta_{\text{per}}(\cdot; \boldsymbol{x}_0) := \sum_{m \in \mathbb{Z}} \delta(\cdot; \boldsymbol{x}_0 + (mL, 0)^\top), \qquad \text{in } \Omega.$$

We denote by $G_{\text{per}}(\cdot, \boldsymbol{x}_0)$ the outgoing periodic Green's function with period $L$ associated with the point source at $\boldsymbol{x}_0$. Then we impose the outgoing behavior of $v$ by requiring that $v - G_{\text{per}}(\cdot, \boldsymbol{x}_0)$ satisfies the Sommerfeld radiation condition. Set $\Gamma_{\pm t} := \mathbb{R} \times \{\pm t\}$ for some $t > h$. We consider the following inverse problem below.

**Inverse problem:** Determine the indices $m \in \mathbf{I}$ corresponding to the defective sources located at $\boldsymbol{x}_m = \boldsymbol{x}_0 + (mL, 0)^\top$ and their intensities $\gamma_m$ using measurements of the radiated waves $v$ on $\Gamma_t \cup \Gamma_{-t}$ with $t > h$.

The uniqueness of the solution to this inverse problem can be proved using the unique continuation and linear independence of Delta distributions, as in the case of free space [2]. Over the past two decades, inverse source problems in free space or bounded domains have been studied extensively from both theoretical and computational perspectives, see for instance [1, 2] and references therein. Unlike the classical inverse source problems, the problem of identifying defective sources within an infinite periodic array appears, to the best of our knowledge, not to have been studied before.

## 2 Reconstruction method

We first reformulate the inverse problem using the Floquet–Bloch transform. Recall that the Floquet-Bloch transform of a function $\varphi \in \S(\mathbb{R}^2)$ in $x_1$-direction, with period $L > 0$, is defined by

$$\mathcal{F}\varphi(x; \xi) := \sum_{m \in \mathbb{Z}} \varphi(x_1 + mL, x_2) e^{-i\xi mL}, \tag{2}$$

for $\xi \in (-\pi/L, \pi/L)\,, x \in \mathbb{R}^2$. We are now ready to reformulate the inverse problem of interest as a quasi-periodic problem. To this end, we first define $u := v - G_{\text{per}}(\cdot, \boldsymbol{x}_0)$ from (1). Then $u$ satisfies

$$\Delta u + k^2 u = -\sum_{m \in \mathbf{I}} \sigma_m \delta(\cdot; \boldsymbol{x}_m)$$

and the Sommerfeld radiation condition. Next, we denote by $u_\xi := \mathcal{F}u(\cdot;\xi)$, $\xi \in \left(-\frac{\pi}{L}, \frac{\pi}{L}\right)$, then $u_\xi$ is $\xi-$quasi-periodic, satisfies

$$\Delta u_\xi + k^2 u_\xi = -\sum_{m \in \mathbf{I}} \sigma_m e^{-i\xi m L} \delta(\cdot; \boldsymbol{x}_0), \quad \text{in } \Omega_0, \tag{3}$$

and the Rayleigh expansion radiation condition. We can state the quasiperiodic inverse problem as follows

**Quasi-periodic inverse problem:** Given $u_\xi$ on

$$\Lambda_{\pm t} := \left(-\frac{L}{2}, \frac{L}{2}\right) \times \{\pm t\},$$

for a finite number of $\xi \in \left(-\frac{\pi}{L}, \frac{\pi}{L}\right)$, determine the number of defective sources $N$, the indices $m \in \mathbf{I}$ corresponding to their locations $\boldsymbol{x}_0 + (mL, 0)^\top$, along with their intensities $\gamma_m$.

Next, we propose an imaging function to determine $\boldsymbol{x}_0$. Let $\xi \in \left(-\frac{\pi}{L}, \frac{\pi}{L}\right)$. For $z \in \Omega_0$, we define the following imaging function

$$I_\xi(z) = \int_{\Lambda_t \cup \Lambda_{-t}} \Phi_\xi(x, z) u_\xi(x) ds(x), \tag{4}$$

where

$$\Phi_\xi(x, z) = \frac{-i}{2L} \sum_{\substack{j \in \mathbb{Z} \\ \beta_j > 0}} e^{-i\xi_{(j)}(x_1 - z_1) - i\beta_j |x_2 - z_2|}.$$

The behavior of $I_\xi(z)$ is analyzed in the following theorem.

**Theorem 1** *The imaging function satisfies*

$$I_\xi(z) = \sum_{n \in \mathbf{I}} \sigma_n e^{-i\xi n L} \mathcal{K}_\xi(z, \boldsymbol{x}_0), \tag{5}$$

*where*

$$\mathcal{K}_\xi(z, \boldsymbol{x}_0) = \sum_{\substack{j \in \mathbb{Z} \\ \beta_j > 0}} \frac{e^{i\xi_{(j)}\left(z_1 - \boldsymbol{x}_0^{(1)}\right)}}{2L\beta_j} \cos\left(\beta_j \left(z_2 - \boldsymbol{x}_0^{(2)}\right)\right).$$

The function $\mathcal{K}_\xi(z, \boldsymbol{x}_0)$ can be shown to have a significant peak as $z$ approaches $\boldsymbol{x}_0$ and thus allows us to locate $\boldsymbol{x}_0$ by evaluating $I_\xi(z)$ in some search domain.

Now computing $\Theta_\xi = I_\xi(\boldsymbol{x}_0)/\mathcal{K}_\xi(\boldsymbol{x}_0, \boldsymbol{x}_0)$ for multiple $\xi$, we can determine the number of defective sources from the following theorem.

**Theorem 2** *For $M \in \mathbb{N}^*$, define an $M \times M$ matrix*

$$D_M := \begin{bmatrix} \Theta_{\xi_{M-1}} & \dots & \Theta_{\xi_0} \\ \Theta_{\xi_M} & \dots & \Theta_{\xi_1} \\ \vdots & \ddots & \vdots \\ \Theta_{\xi_{2M-2}} & \dots & \Theta_{\xi_{M-1}} \end{bmatrix}.$$

*Then $D_M$ are invertible for all $M \leq N$ and $D_M$ are singular otherwise if and only if there are $N$ defective sources.*

We now proceed to present the method for determining the indices $m$ and intensities $\gamma_m$ of the defective sources.

**Theorem 3** *Let $N^* \in \mathbb{N}$ such that $\mathbf{I} \subset [-\mathbf{N}^*, \mathbf{N}^*]$. Then the linear system*

$$\begin{bmatrix} e^{-i\xi_0(-N^*)L} & \dots & e^{-i\xi_0 N^* L} \\ \vdots & \ddots & \vdots \\ e^{-i\xi_{2N^*}(-N^*)L} & \dots & e^{-i\xi_{2N^*} N^* L} \end{bmatrix} \begin{bmatrix} \varkappa_1 \\ \vdots \\ \varkappa_{2N^*+1} \end{bmatrix} = \begin{bmatrix} \Theta_{\xi_0} \\ \vdots \\ \Theta_{\xi_{2N^*}} \end{bmatrix}$$

*is uniquely solvable. Furthermore, the indices $m$ for the location of defective sources and their intensities $\gamma_m$ are determined by*

$$m = -N^* - 1 + j, \qquad \gamma_m = \varkappa_j + 1,$$

*for $j \in \{n \in \mathbb{N} : 1 \leq n \leq 2N^* + 1, \varkappa_n \neq 0\}$.*

Here we present a numerical example demonstrating the performance of the method. The integer indices $m$ of the defective sources are all accurately reconstructed.

| True $\boldsymbol{x}_0$ | Comp. $\boldsymbol{x}_0$ | Comp. $m$ | True $\gamma$ | Comp. $\gamma$ |
|---|---|---|---|---|
| $(0.61, -1.59)$ | $(0.60, -1.59)$ | 3 | 0.5 | 0.52 |
| $(-0.70, -1.15)$ | $(0.70, -1.14)$ | $-2$ | 0.9 | 0.89 |
| | | 5 | 0.1 | 0.09 |
| $(-0.52, -0.48)$ | $(0.51, -0.48)$ | $-1$ | 0.3 | 0.29 |
| | | 3 | 0.5 | 0.49 |
| | | 10 | 0.7 | 0.69 |

Table 1: Reconstruction results for the indices and intensities of defective sources (noise level $\delta = 10\%$)

# Computational spectral geometry for the Maxwell cavity problem.

Nilima Nigam[1,*], Kewen Bu[1], Jayda Grant[1], Rodrigo Stehling[1]
[1]Department of Mathematics, Simon Fraser U.Burnaby,BC
*Email: nigam@math.sfu.ca

**Abstract**

## 1 Introduction

In this talk we report on recent results in spectral optimization for the Maxwell eigenvalue problem. Let $\Omega \in \mathbb{R}^3$ be a bounded domain with Lipschitz boundary, the permittivity $\epsilon$ be a symmetric positive-definite matrix and $\lambda, \mathbf{E}$ be an eigenpair of the cavity problem with perfectly conducting boundary conditions $n \times \mathbf{E} = 0$ with $\mathbf{E} \in X_N := \{H(\nabla\times)|E\times n = 0\,\text{on the boundary}\} \cap H(div)$ and

$$\nabla \times (\nabla \times \mathbf{E}) = \lambda\epsilon\mathbf{E}, \quad \nabla \cdot (\epsilon\mathbf{E}) = 0, \qquad in\,\Omega. \tag{1}$$

Clearly the eigenvalues $\lambda_k = \lambda_k(\Omega, \epsilon)$ depend on the domain as well as the permittivity $\epsilon$. We suppress the dependence on these quantities if the context is clear. A *scale invariant* spectral quantity for varying $\Omega$ is $F_k(\Omega) := |\partial\Omega|\lambda_k(\Omega, \epsilon)$. Theoretical results on optimization wrt. $\epsilon \in (W^{1,\infty}(\Omega))^{3\times 3}$ were recently studied [4]. It is known that $\lambda_k(\Omega, \epsilon)$ are not generically Frechét differentiable as a function of the parameters if the eigenvalue possesses multiplicity, [5]. In [4], the authors consider instead a constrained optimization problem for symmetric functions of the eigenvalues, which bypass this issue, and show there are *no* critical permittivities. However, this leaves open the question: in a narrower class of permittivities, are there local extrema of $F_k(\epsilon)$?

In another direction, one may fix the permittivity, and vary $\Omega$. Spectral shape optimization is of extensive interest in engineering applications, where cavities may be designed to optimize either the eigenfrequency, or the field distribution. A standard approach is consider perturbations $\Omega_t$ of a reference domain $\Omega$, and use standard pullback arguments to calculate both shape derivatives and eigenpairs (for a related problem) on $\Omega$ [3]. Once again, these approaches can fail in the presence of multiple eigenvalues.

Our starting point is [2] which collects results on the spectrum of 1 in product domains $\Omega := \omega \times (a, b)$ where $\omega \in \mathbb{R}^2$ is a bounded, Lipschitz domain. The Maxwell eigenvalues can be represented as sums of Dirichlet (dir)or Neumann (neu) eigenvalues $\Delta$ on $\omega$ and $(a, b)$. The lowest nonzero eigenvalue on a product-domain $\Omega$ is

$$\lambda_1(\Omega) = \min\{\lambda_1^{neu}(\omega) + (\frac{\pi}{|b-a|})^2, \lambda_1^{dir}(\omega)\}.$$

Simple arguments show that $|\Omega|\lambda_1(\Omega)$ is minimized when $\Omega$ is a flat (two-sided) disk, [1].

## 2 Optimization over domains

If $\omega = (0, \pi)^2$, then $F_k(\Omega)$ for the domain $\Omega := \omega\times(0, \pi l_3)$ takes the form $J_k := F_k(\omega\times(0, \ell_3)) = 2\pi^2(1+2\ell_3)(k_1^2+k_2^2+(k_3/\ell_3)^2)$ for integers $k_1, k_2, k_3 \geq 0$ (with atmost 1 of $k_i$ zero); the eigenvalues have multiplicity 2 if $k_i > 0$. Ordering the eigenvalues (and hence the spectral invariants) is non-trivial and depends on $\ell_3$ and $k_i$. For instance,

$$J_1(\Omega) := \begin{cases} 2\pi^2(1 + 2\ell_3), & 0 \leq \ell_3 \leq 1 \\ 2\pi^2(1 + 2\ell_3)(1 + \ell_3^{-2}), & \ell_3 > 1. \end{cases}$$

In 1 we plot $J_1..J_5$. The computations are performed using quadratic edge elements in Netgen; we enforce the divergence constraint by projecting out the gradients, following the thesis of S. Zaglmayr (2006). We see the flat two-sided square ($\ell = 3$) is a minimizer of $J_1$ in this class of domains. We see the multiplicity of eigenvalues makes $J_1 = J_2$ after $\ell_3 = 1$. The first local maximizer of $J_k$ is the first local minimizer of $J_{k+1}$. More detailed expressions and proofs will be presented in forthcoming work. We also plot the same spectral invariants for $\Omega$ which are ellipsoids $x^2 + y^2 + (\ell z)^2 = 1$. These are not product-form domains. Once again, though, the flat two-sided disk $\ell_3 = 0$ minimizes $F_1(\Omega)$. The sphere $\ell_3 = 1$ is a local maximizer. We note that if topology changes are permitted, then $F_k(\Omega)$ will be minimized over a domain of genus $k$.

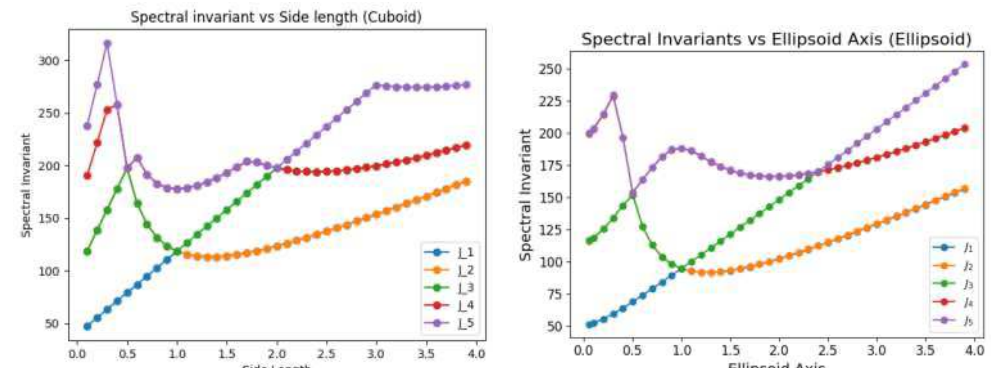

Figure 1: Left: Spectral invariants $J_k(\ell_3)$ v/s cuboid height $\ell_3$.

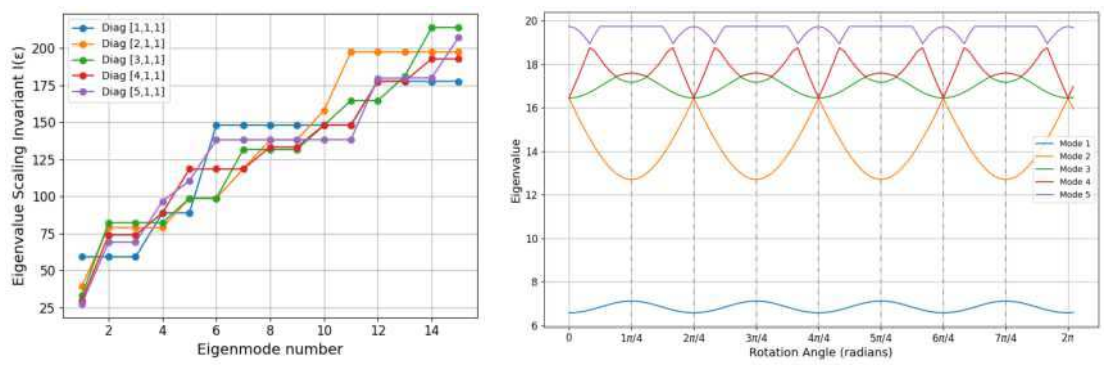

Figure 2: Left: Computed spectral invariants on a cube, $I_k(\epsilon)$ v/s k, for $\epsilon = diag(A, 1, 1)$. Right: $I_k(\alpha)$ v/s $\alpha$, $\epsilon = R^T diag(3, 1, 1)R$, $R$ a rotation around z axis.

## 3 Optimization over permittivity

In the case of permittivity optimization over the class of real spd matrices, a scale-invariant spectral quantity is $I_k(\epsilon) := tr(\epsilon)\lambda_k(\epsilon)$. We numerically investigated the behaviour of $I_k(\epsilon)$ for permittivities in highly restricted classes.

In figure 2, we plot $I_k(\epsilon)$ for permittivities of the form $\epsilon = diag(A, 1, 1)$. The domain $\Omega$ is a cube. Perhaps unsurprisingly, the breaking of metric symmetry leads to changes in eigenvalue multiplicity. We also consider permittivities of the form $\epsilon = R^T DR$ where $R$ is a rotation and $D$ is a diagonal matrix.We note that if $\lambda_1$ is a simple eigenvalue, then the Frechet derivative $\lambda'(\alpha) = (\int_\Omega \epsilon'(\alpha)\mathbf{E}(\alpha \cdot \mathbf{E}(\alpha)) / \int_\Omega \epsilon(\alpha)\mathbf{E}(\alpha) \cdot \mathbf{E}(\alpha)$. We see rich (periodic) structure; one also sees the optima for $k = 2, 3...$ may occur at locations of high multiplicity, posing challenges for optimization. The domain $\Omega$ plays a key role; in 3, the same computation on a tetrahedral region (which is also highly symmetric) shows a completely different multiplicity and optimization structure.

## 4 Conclusion

Computations allow us to investigate the dependance of the Maxwell spectrum on both domain shape and permittivity. Even for simple situations, rich spectral structure is revealed. We see that while global critical permittivities in the wide class $(W^{1,\infty})^{3\times 3}$ do not exist, in more constrained admissable classes the optimization problem has solutions. In forthcoming work, we present some conjectures concerning the symmetries of shape optimizers.

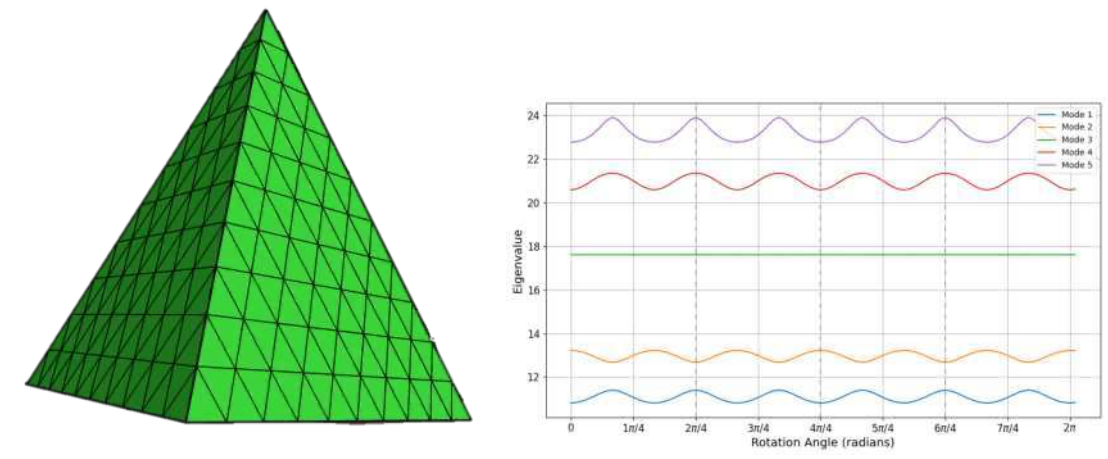

Figure 3: Left: Computed spectral invariants on a tetrahedron Right: $I_k(\alpha)$ v/s $\alpha$, $\epsilon = R^T diag(3, 1, 1)R$, $R$ a rotation around z axis.

# A Factorization Method Approach to the Biharmonic Transmission Problem in Absorbing Media

General Ozochiawaeze[1,*], Rafael Ceja Ayala[2], Isaac Harris[1]

[1]Department of Mathematics, Purdue University, West Lafayette, Indiana, USA

[2]School of Mathematical and Statistical Sciences, Arizona State University, Tempe, Arizona, USA

*Email: oozochia@purdue.edu

**Abstract**

In this talk, I will present an inverse shape reconstruction problem arising in biharmonic wave scattering, motivated by flexural wave propagation in thin elastic plates. I will focus on a transmission problem, where the goal is to recover a penetrable, absorbing inclusion from far-field data at a fixed frequency. Using Kirsch's factorization method, we will show one can directly recover the shape, size, and location of the unknown medium from spectral information of the far-field operator. I will also discuss the Born approximation for weak scatterers in this setting and compare its accuracy against exact far-field data for simple weak scatterers, illustrating both its usefulness and limitations.

## 1 Introduction

We study the inverse problem of imaging unknown inclusions in thin elastic (Kirchhoff–Love) plates using biharmonic waves, where a *penetrable, absorbing* inclusion allows wave entry while dissipating energy. Applications include platonic gratings [3], structural health monitoring of thin-walled aircraft [6], and ice floes embedded in elastic sheets [1]. Furthermore, biharmonic scattering has a broader relevance in ultra-broadband elastic cloaks and earthquake-resistant passive systems [4]. For thin plates relative to the wavelength, propagation is governed by a two-dimensional biharmonic equation, whose separation into radiating and evanescent components complicates geometric reconstruction from far-field data. We extend the factorization method [5] to this fourth-order setting under a natural absorption assumption, allowing a spectral characterization of the far-field operator and a binary test for interior sampling points.

## 2 Model

We consider time-harmonic scattering in a thin elastic (Kirchhoff–Love) plate containing a penetrable, absorbing inclusion $D \subset \mathbb{R}^2$ with $C^2$ boundary, where the exterior $\mathbb{R}^2 \setminus \overline{D}$ is simply connected. For an incident plane wave

$$u^i(x,d) = e^{i\kappa x\cdot d}, \quad \kappa > 0,\ d \in \mathbb{S}^1,$$

the total field $u = u^i + u^s$ satisfies

$$(\Delta^2 - \kappa^4 n(x))\, u^s = \kappa^4 (n(x) - 1)\, u^i \quad \text{in } \mathbb{R}^2,$$

with transmission conditions on $\partial D$:

$$[\![u]\!] = [\![\partial_\nu u]\!] = [\![\Delta u]\!] = [\![\partial_\nu(\Delta u)]\!] = 0,$$

and Sommerfeld radiation conditions for $u^s$ and $\Delta u^s$ at infinity. The refractive index $n \in L^\infty(D)$ satisfies

$$\operatorname{supp}(n-1) = \overline{D}, \quad \Im(n) \geq \alpha > 0 \quad [2],$$

ensuring well-posedness and uniqueness of $u^s$. The inverse problem is to reconstruct the support $D$ from far-field measurements

$$\{u^\infty(\hat{x}, d) : \hat{x} \in \mathbb{S}^1,\ d \in \mathbb{S}^1\},$$

determining for each sampling point $z \in \mathbb{R}^2$ whether $z$ lies inside or outside the scatterer without assuming prior knowledge of $n(x)$.

## 3 Methods and Results

Let $F : L^2(\mathbb{S}^1) \to L^2(\mathbb{S}^1)$ denote the far-field operator mapping a weight function $g \in L^2(\mathbb{S}^1)$ to the superposition of the biharmonic far-field data. Under the absorption assumption, we managed to show $F$ admits a factorization of the form

$$\Im(F) = H^* \Im(T) H,$$

where $H$ is the Herglotz operator and $\Im(T)$ is a coercive, self-adjoint operator on $L^2(D)$. We prove a rigorous range characterization: a point $z \in \mathbb{R}^2$ lies inside the scatterer if and only if a test function $\phi_z := e^{-ik\hat{x}\cdot z}$ belongs to the range of $\Im(F)^{1/2}$, providing a computable binary criterion for qualitative reconstruction. This yields a computable indicator based solely on the spectral data of $\Im(F)$,

$$W(z) := \left[ \sum_{n=1}^{\infty} \frac{|(\phi_z, \psi_n)_{L^2(\mathbb{S}^1)}|^2}{\lambda_n} \right]^{-1},$$

where $(\lambda_n, \psi_n)$ are the eigenpairs of $\Im(F)$. Large indicator values correspond to points inside $D$, while small values indicate points outside. For weak inhomogeneities, the Born approximation provides an efficient and accurate alternative, as confirmed by numerical experiments. In this framework, the scattered field is linearized with respect to the contrast $(n-1)$, yielding the approximation

$$u_{\mathrm{B}}^{\infty}(\theta, \phi) = \frac{\kappa^2}{2} \int_{B_\epsilon} (n-1)\, e^{-i\kappa y\cdot(\hat{x}-d)}\, dy.$$

For small, weak inhomogeneities, this approximation provides a simple and computationally efficient estimate of the far-field pattern. In fact, up to a constant, it coincides with the Born approximation used in classical acoustic scattering problems, motivating a numerical justification for its application to biharmonic waves in the thin-plate regime.

We consider small circular scatterers $B_\epsilon$, for which separation of variables allows explicit computation of the biharmonic far-field pattern via a Fourier–Bessel decomposition and a linear system for the boundary-matched coefficients. For weak scatterers, the Born approximation provides a simple alternative, and comparison of the resulting far-field matrices $\mathbf{F}_{\text{exact}}$ and $\mathbf{F}_{\text{Born}}$ shows decreasing errors with smaller $\epsilon$, as illustrated in Table 1 and Figure 1.

| $\epsilon$ | $\|\mathbf{F}_{\text{exact}} - \mathbf{F}_{\text{Born}}\|_\infty$ |
|---|---|
| 1.00 | 0.6480 |
| 0.90 | 0.5008 |
| 0.80 | 0.3784 |
| 0.70 | 0.2774 |
| 0.60 | 0.1935 |
| 0.50 | 0.2134 |

Table 1: Maximum entry-wise difference between exact and Born far-field matrices for small disks $B_\epsilon$, i.e., $\|\mathbf{F}_{\text{exact}} - \mathbf{F}_{\text{Born}}\|_\infty = \max_{i,j} \left|(\mathbf{F}_{\text{exact}})_{ij} - (\mathbf{F}_{\text{Born}})_{ij}\right|$.

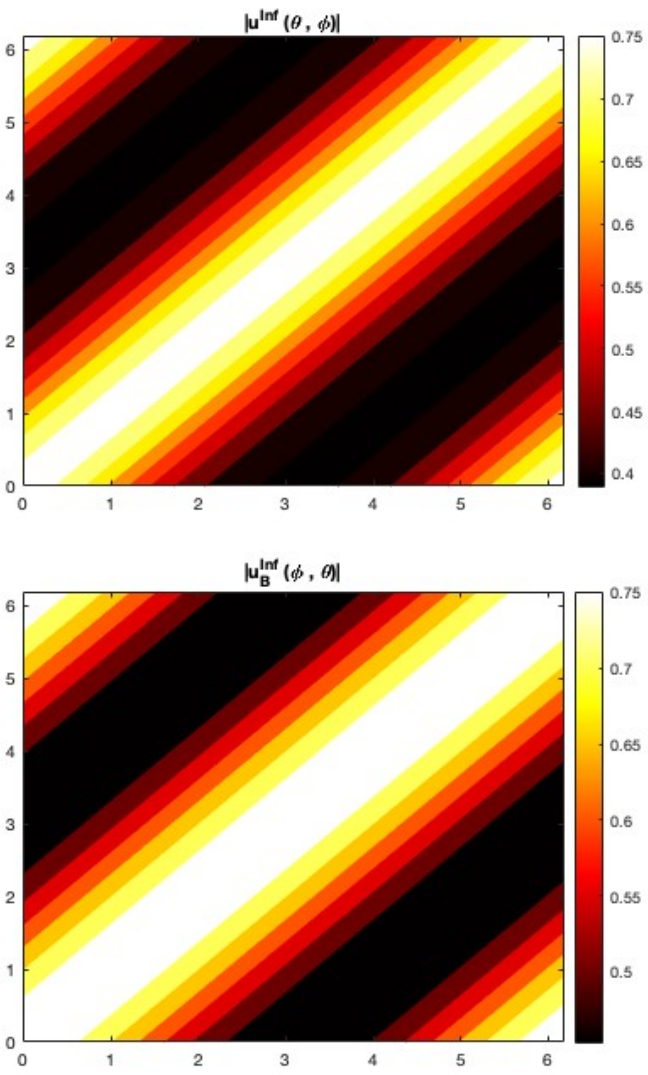

Figure 1: Exact (left) and Born-approximated (right) far-field matrices for a disk of radius $\epsilon = 0.5$, $n = 1.0 + 0.5\mathrm{i}$, $\kappa = 2$.

# High-order Nyström method for scattering by penetrable inhomogeneous media with discontinuous refractive index

**Ambuj Pandey**[1,*], **Krishna Yamanappa Poojara**[2]

[1]Department of Mathematics, Indian Institute of Science Education and Research Bhopal, Madhya Pradesh 462066, India

[2]Institute of Computational Mathematics and Scientific/Engineering Computing, Academy of Mathematics and Systems Science, Chinese Academy of Sciences, Beijing 100190, China

*Email: ambuj@iiserb.ac.in

**Abstract**

We present a high-order Nyström method for two-dimensional wave scattering by penetrable inhomogeneous media with discontinuous refractive indices, based on the Lippmann-Schwinger volume integral formulation. Most popular, fast, high-order methods for the Lippmann-Schwinger equation perform well for smooth media but lose accuracy and exhibit low-order convergence when the material contrast is discontinuous. The proposed approach evaluates the volume potential using a patchwise quadrature that separates regular and singular interactions: smooth contributions are treated with standard high-order quadrature, while singular interactions are handled by a Chebyshev expansion of density, a singularity-regularizing change of variables, and precomputed integral moments. The resulting discretization retains high-order convergence without the need for interface-fitted meshing. Numerical experiments demonstrate near machine-precision accuracy (errors of order $10^{-13}$) even for scattering by media with sharp material interfaces.

## 1 Introduction

We consider time-harmonic wave scattering by a bounded penetrable inhomogeneity $\Omega \subset \mathbb{R}^2$ with refractive index $n(\mathbf{x})$, where $n(\mathbf{x}) = 1$ outside $\Omega$ and may be discontinuous across $\partial\Omega$. For an incident field $u^i$ satisfying $\Delta u^i + \kappa^2 u^i = 0$ in $\mathbb{R}^2$, the total field $u$ satisfies $\Delta u + \kappa^2 n^2(\mathbf{x}) u = 0$ together with the Sommerfeld radiation condition [2] for the scattered field $u^s = u - u^i$. Discontinuities in $n(\mathbf{x})$ reduce solution regularity and make high-order discretization challenging. To avoid truncation of the unbounded domain, we use an equivalent Lippmann-Schwinger formulation [2]

$$u(\mathbf{x}) + \kappa^2 \mathcal{A}[u](\mathbf{x}) = u^i(\mathbf{x}), \qquad (1)$$

where $\mathcal{A}[u](\mathbf{x}) = \int_\Omega G_\kappa(\mathbf{x} - \mathbf{y})\, m(\mathbf{y})\, u(\mathbf{y})\, d\mathbf{y}$, $m(\mathbf{x}) = 1 - n^2(\mathbf{x})$ and $G_\kappa(\mathbf{x}) = \frac{i}{4} H_0^{(1)}(\kappa|\mathbf{x}|)$ is the free-space Green's function. This formulation automatically enforces the radiation condition and restricts computation to $\Omega$, but the weak logarithmic singularity of $G_\kappa$ and discontinuities in $m(\mathbf{x})$ pose challenges for constructing high-order solver. While fast, high-order solvers for (1) are well developed for smooth media, their accuracy deteriorates in the presence of discontinuous material contrast because the density $m(\mathbf{x})u(\mathbf{x})$ is only piecewise smooth. The proposed method retains high-order convergence and yields significantly improved accuracy for such problems without requiring interface-fitted meshes. Numerical experiments demonstrate solutions accurate to near machine precision even for discontinuous media.

## 2 Numerical Method

The Nyström discretization is obtained by evaluating the convolution operator $\mathcal{A}[u](\mathbf{x})$ at the quadrature (Nyström) nodes; thus discretizing (1) reduces to computing the volume potential $\mathcal{A}[u]$ at these points. The scattering region $\Omega$ is decomposed into non-overlapping patches $\Omega_\ell$ such that $\Omega = \bigcup_{\ell=1}^P \Omega_\ell$, where each patch is obtained from a smooth parametrization $\xi_\ell : [-1,1]^2 \to \Omega_\ell$. Under this mapping the operator $\mathcal{A}[u]$ can be written as

$$\mathcal{A}[u](\mathbf{x}) = \sum_{\ell=1}^{P} \int_{-1}^{1}\int_{-1}^{1} G_\kappa(\mathbf{x} - \xi_\ell(s,\tau))\varphi_\ell(s,\tau)\, ds d\tau$$

with $\varphi_\ell(s,\tau) = m(\xi_\ell(s,\tau))\, u(\xi_\ell(s,\tau))\, J_\ell(s,\tau)$, where $J_\ell$ denotes the Jacobian of the parametrization.

The treatment of each patch integral depends on the location of the target point $\mathbf{x}$. If $\mathbf{x}$ is well separated from $\Omega_\ell$, the integrand is smooth and the integral is evaluated using a high-order tensor-product quadrature (in our implementation, we use the Fejér first quadrature rule). If $\mathbf{x} \in \Omega_\ell$, the integral kernel $G_\kappa$ is singular. In this case the smooth density $\varphi_\ell$ is approximated by a truncated Chebyshev expansion, which produces a collection of singular integral moments. To compute these accurately, we apply a rectangular-polar change of variables [1] centered at the target point. This transformation regularizes the integrand, rendering the integrand smooth, after which the integrals are evaluated to high order using the Fejér first quadrature. Note that, these singular integral moments depend only on the geometry and the kernel and therefore precomputed and stored for the repeated use. For near-singular interactions, when $\mathbf{x}$ lies close to $\Omega_\ell$, the point is projected onto the integration patch $\Omega_\ell$ and evaluated to high-order by employing the same strategy used for singular integrals. Having obtained an accurate evaluation of $\mathcal{A}[u](\mathbf{x})$, the Nyström system follows directly. Let $\mathcal{N} = \{\mathbf{x}_i\}_{i=1}^N \subset \Omega$ denote the quadrature nodes. Enforcing the equation at each $\mathbf{x}_i \in \mathcal{N}$ gives $u(\mathbf{x}_i) + \kappa^2\mathcal{A}[u](\mathbf{x}_i) = u^i(\mathbf{x}_i), \quad i = 1, \ldots, N$, which determines the unknown values $u(\mathbf{x}_i)$ and is solved iteratively using GMRES.

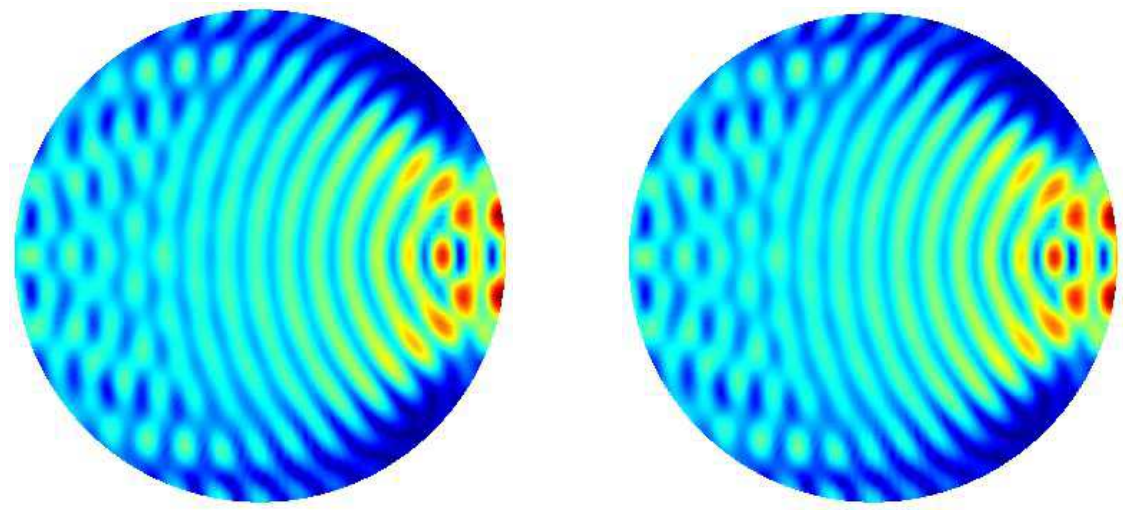

Figure 1: Plane-wave scattering with $\kappa = 18$ for a unit disc, where $n(\mathbf{x}) = \sqrt{2}$ in $\Omega$ and $n(\mathbf{x}) = 1$ outside. Left: $|u|$ computed using the exact Mie series solution. Right: $|u|$ computed using the proposed method.

## 3 Computational Results

To demonstrate the performance of the proposed method, we present scattering simulations for unit-disc and kite-shaped obstacles under two types of incidence: a plane wave (PW) $u^i(\mathbf{x}) = e^{i\kappa\mathbf{x}\cdot\mathbf{d}}$ (where $\mathbf{d}$ is the incidence direction) and a partial plane wave (PPW) $u^i(\mathbf{x}) = J_0(\kappa|\mathbf{x}|)$. The numerical performance is studied in two settings: (i) increasing the number of patches $P$ while keeping the number of quadrature nodes per patch fixed, and (ii) fixing the patch resolution while varying the number of nodes.

Tables 1 and 2, report convergence results for the unit disc, where the numerical solution is compared with the known analytical solution. The results clearly confirm the high-order convergence of the proposed scheme. In the tables, $P, n_1$, and $n_2$ denote the number of patches and the number of quadrature nodes in the parametric variables, respectively; $\epsilon_\infty$ is the maximum relative error and noc the numerical order of convergence. The refractive index is taken as $n(\mathbf{x}) = \sqrt{2}$ inside $\Omega$ and $n(\mathbf{x}) = 1$ outside, with wavenumber $\kappa = 10$. Similar behavior is observed for the kite-shaped obstacle, where the solution is compared against a reference solution computed on a refined grid. For illustration, Figure 1 depicts the total field for plane-wave scattering by a unit disk with wavenumber $\kappa = 18$, where the refractive index is $n(\mathbf{x}) = \sqrt{2}$ in $\Omega$ and $n(\mathbf{x}) = 1$ outside.

| $P \times n_1 \times n_2$ | PW | | PPW | |
|---|---|---|---|---|
| | $\epsilon_\infty$ | noc | $\epsilon_\infty$ | noc |
| $80 \times 2 \times 2$ | $1.41{\times}10^{00}$ | - | $9.94{\times}10^{-01}$ | - |
| $80 \times 4 \times 4$ | $2.10{\times}10^{-02}$ | 6.07 | $8.68{\times}10^{-03}$ | 6.84 |
| $80 \times 8 \times 8$ | $3.56{\times}10^{-05}$ | 9.20 | $8.11{\times}10^{-06}$ | 10.06 |
| $80 \times 16 \times 16$ | $2.12{\times}10^{-09}$ | 14.04 | $7.99{\times}10^{-11}$ | 16.63 |
| $80 \times 32 \times 32$ | $2.74{\times}10^{-13}$ | 12.92 | $1.78{\times}10^{-13}$ | 8.81 |

Table 1: Convergence for the fixed patches while varying the grids in each patch.

| $P \times n_1 \times n_2$ | PW | | PPW | |
|---|---|---|---|---|
| | $\epsilon_\infty$ | noc | $\epsilon_\infty$ | noc |
| $5 \times 10 \times 10$ | $9.82{\times}10^{-02}$ | - | $2.94{\times}10^{-03}$ | - |
| $20 \times 10 \times 10$ | $5.18{\times}10^{-04}$ | 7.57 | $2.89{\times}10^{-05}$ | 6.67 |
| $80 \times 10 \times 10$ | $1.01{\times}10^{-06}$ | 9.00 | $6.28{\times}10^{-08}$ | 8.85 |
| $320 \times 10 \times 10$ | $1.07{\times}10^{-09}$ | 9.88 | $2.25{\times}10^{-10}$ | 8.12 |

Table 2: Convergence for the varying patches.

# Maxwell à la Helmholtz: Frequency-Robust Direct Boundary Integral Equations for 3D PEC Scattering via Helmholtz Integral Operators

Carlos Pérez-Arancibia[1,*], Catalin Turc[2]
[1]Department of Applied Mathematics, University of Twente, Enschede, the Netherlands
[2]Department of Mathematical Sciences, New Jersey Institute of Technology, Newark, USA
*Email: c.a.perezarancibia@utwente.nl

**Abstract**

We introduce direct boundary integral equation (BIE) formulations for time-harmonic electromagnetic scattering by smooth perfectly electrically conducting (PEC) obstacles, based exclusively on Helmholtz boundary integral operators. The approach exploits an equivalence between the Maxwell PEC scattering problem and vector Helmholtz equations with boundary conditions expressed in terms of the Dirichlet and Neumann traces of the electromagnetic fields. The resulting combined-field formulations are proved uniquely solvable at all frequencies and free from spurious resonances. A Calderón-type regularization yields Fredholm second-kind equations amenable to efficient iterative solution. Numerical experiments using a Nyström discretization based on the Density Interpolation Method, implemented in `Inti.jl`, demonstrate the accuracy and frequency robustness of the proposed approach.

## 1 Introduction

The numerical solution of time-harmonic electromagnetic scattering problems via BIEs is significantly more challenging than its acoustic counterpart governed by the Helmholtz equation. Classical electromagnetic BIE formulations—such as the EFIE, MFIE, and their combined-field variant CFIE—involve surface-tangential vector unknowns with subtle regularity properties, complicated kernel structures, and are prone to low-frequency breakdown.

In this work, we introduce direct BIE formulations for three-dimensional PEC scattering that bypass these difficulties by relying exclusively on classical Helmholtz integral operators. The key ingredient, rigorously established in the companion paper [1], is an equivalence between the Maxwell PEC scattering problem and vector Helmholtz boundary value problems whose boundary conditions involve the Dirichlet and Neumann traces of the fields and the curvature of the obstacle. This allows any Helmholtz BIE solver to be adapted for PEC scattering, provided a sufficiently accurate surface geometry representation is available. Calderón-type regularization yields Fredholm second-kind equations, and numerical results using the Density Interpolation Method [2] in `Inti.jl` [3] validate the approach.

## 2 Problem formulation

Let $\Omega \subset \mathbb{R}^3$ be a bounded domain occupied by a smooth PEC obstacle, with boundary $\Gamma = \partial\Omega$ of class $C^{2,\alpha}$, $\alpha \in (0,1]$, and outward unit normal $\nu$. We denote by $\mathrm{P}_\nu$ and $\mathrm{P}_t$ the normal and tangential projectors onto $\Gamma$, and by $\gamma$ and $\partial_\nu$ the exterior Dirichlet and Neumann traces on $\Gamma$.

An incident time-harmonic electromagnetic field $(E^i, H^i)$, satisfying Maxwell's equations in an open set containing $\Omega$, impinges on the obstacle and generates scattered fields $(E^s, H^s)$ governed by:

$$\operatorname{curl} E^s - i\omega\mu H^s = 0 \quad \text{in } \mathbb{R}^3 \setminus \Omega, \tag{1a}$$

$$\operatorname{curl} H^s + i\omega\epsilon E^s = 0 \quad \text{in } \mathbb{R}^3 \setminus \Omega, \tag{1b}$$

$$\mathrm{P}_t\, \gamma(E^s + E^i) = 0 \quad \text{on } \Gamma, \tag{1c}$$

$$\lim_{|x|\to\infty} |x| \left( \sqrt{\mu}\, E^s \times \tfrac{x}{|x|} - \sqrt{\epsilon}\, H^s \right) = 0, \tag{1d}$$

where $\omega > 0$ is the angular frequency, $\epsilon, \mu > 0$ are the dielectric and magnetic constants, and $k = \omega\sqrt{\epsilon\mu}$ is the wavenumber. Condition (1c) is the PEC boundary condition and (1d) is the Silver–Müller radiation condition. Problem (1) admits a unique solution $E^s, H^s \in C^2(\mathbb{R}^3 \setminus \Omega, \mathbb{C}^3) \cap C^{1,\alpha}(\mathbb{R}^3 \setminus \Omega, \mathbb{C}^3)$ for all $\omega > 0$ [1].

## 3 Maxwell–Helmholtz equivalence

The following result, established in [1], is the cornerstone of our approach. Here, $\mathcal{R} = \mathrm{D}_\Gamma\,\nu$ and $\mathcal{H} = \frac{1}{2}\mathrm{tr}\,\mathcal{R}$ denote the shape operator and mean curvature of $\Gamma$, respectively.

**Theorem 1** *The fields $E^s, H^s$ solve (1) if and only if either:*

**(E)** *$H^s = (i\omega\mu)^{-1}\,\mathrm{curl}\,E^s$ and $E^s$ satisfies $\Delta E^s + k^2 E^s = 0$ in $\mathbb{R}^3 \setminus \Omega$, the Sommerfeld radiation condition, and*

$$\mathrm{P}_t\,\gamma(E^s + E^i) = 0, \quad (2)$$

$$\mathrm{P}_\nu[\partial_\nu(E^s + E^i) + 2\mathcal{H}\gamma(E^s + E^i)] = 0; \quad (3)$$

**(M)** *$E^s = -(i\omega\epsilon)^{-1}\,\mathrm{curl}\,H^s$ and $H^s$ satisfies $\Delta H^s + k^2 H^s = 0$ in $\mathbb{R}^3 \setminus \Omega$, the Sommerfeld radiatoin condition, and*

$$\mathrm{P}_\nu\,\gamma(H^s + H^i) = 0, \quad (4)$$

$$\mathrm{P}_t[\partial_\nu(H^s + H^i) + \mathcal{R}\gamma(H^s + H^i)] = 0. \quad (5)$$

## 4 Boundary integral formulations

The equivalences **(E)** and **(M)** in Theorem 1 each give rise to a direct BIE formulation for the Maxwell PEC scattering problem using only Helmholtz integral operators—one for the electric field and one for the magnetic field. Due to space constraints, we focus here on **(E)**; the magnetic formulation follows analogously. Applying Green's representation formula componentwise to both $E^s$ and $E^i$, and combining the results, yields the representation

$$E^s(x) = \mathcal{D}[\gamma E](x) - \mathcal{S}[\partial_\nu E](x), \quad (6)$$

for $x \in \mathbb{R}^3 \setminus \Omega$, where $E = E^s + E^i$ is the total electric field, and $\mathcal{S}$, $\mathcal{D}$ are the componentwise Helmholtz single- and double-layer potentials.

Substituting the boundary conditions (2)–(3) into the representation formula (6), taking traces on $\Gamma$, and combining the resulting Dirichlet and Neumann trace equations with coupling parameter $\eta \in \mathbb{R}\setminus\{0\}$ yields the *direct electric combined-field-only integral equation* (D-ECFOIE):

$$\mathcal{L}_e\,\phi = f_e, \qquad f_e := -\partial_\nu E^i - i\eta\gamma E^i, \quad (7)$$

where $\phi = \mathrm{P}_\nu\,\gamma E + \mathrm{P}_t\,\partial_\nu E \in C^{1,0}_{\nu,t}(\Gamma)$, with $C^{p,q}_{\nu,t}(\Gamma) := \{\phi : \Gamma \to \mathbb{C}^3 \mid \mathrm{P}_\nu\,\phi \in C^{p,\alpha}(\Gamma, \mathbb{C}^3),\ \mathrm{P}_t\,\phi \in C^{q,\alpha}(\Gamma, \mathbb{C}^3)\}$, and

$$\begin{aligned}\mathcal{L}_e\phi := &-\tfrac{1}{2}\big\{\mathrm{P}_t\,\phi + (i\eta I - 2\mathcal{H})\,\mathrm{P}_\nu\,\phi\big\}\\ &+ \big\{T + i\eta K + 2(K' + i\eta S)\mathcal{H}\big\}\,\mathrm{P}_\nu\,\phi\\ &- (K' + i\eta S)\,\mathrm{P}_t\,\phi. \end{aligned} \quad (8)$$

Here $S$, $K$, $K'$, and $T$ are the componentwise Helmholtz single-layer, double-layer, adjoint double-layer, and hypersingular boundary integral operators, respectively. The scattered field is then recovered via (6).

We then have the following theorem:

**Theorem 2** *For all $k > 0$ and $\eta \in \mathbb{R} \setminus \{0\}$, the D-ECFOIE (7) admits a unique solution $\phi \in C^{1,0}_{\nu,t}(\Gamma)$.*

The proof proceeds via a Calderón-type regularization: substituting $\phi = \mathcal{R}_\nu\tilde{\phi}$ into (7) and exploiting the Calderón identity $T_0 S_0 = K_0'^2 - \frac{1}{4}I$, where the subscript 0 denotes the $k = 0$ (Laplace) counterparts of the BIE operators, recasts the equation as $\tilde{\phi} + \mathcal{A}\tilde{\phi} = \tilde{f}$, with $\mathcal{A}$ compact on $C^{0,\alpha}(\Gamma, \mathbb{C}^3)$, so existence follows from uniqueness via the Fredholm alternative. The resulting second-kind equation, discretized via a Nyström method, yields GMRES-favorable linear systems.

Analogous results hold for the magnetic field, leading to the D-MCFOIE with unknown $\psi = \mathrm{P}_t\,\gamma H + \mathrm{P}_\nu\,\partial_\nu H \in C^{0,1}_{\nu,t}(\Gamma)$, whose solution gives direct access to the surface currents $J = \nu \times P_t\psi$.

Regarding low-frequency behavior, the D-ECFOIE suffers from breakdown as $k \to 0$, which is mitigated by augmenting (7) with a finite-rank correction enforcing the charge conditions $\int_\Gamma \nu \cdot \gamma E^s\,\mathrm{d}s = 0$; the D-MCFOIE, by contrast, remains well-conditioned for all $k \geq 0$ when $\Gamma$ is simply connected.

Numerical results for both formulations, implemented via the Nyström-DIM solver in `Inti.jl`, closely match those of the indirect counterparts in [1] and are omitted due to space constraints.

# A HDG method with transmission variables for time-harmonic wave propagation problems with constant coefficients

Simone Pescuma[1,*], Gwénaël Gabard[2], Théophile Chaumont-Frelet[3], Axel Modave[4]
[1]CERMICS, ENPC, Institut Polytechnique de Paris, Inria, CNRS, Marne-la-Vallée, France
[2]LAUM, CNRS, IA-GS, Le Mans Université, Le Mans, France
[3]Inria, Laboratoire Paul Painlevé, Université de Lille, Lille, France
[4]POEMS, CNRS, Inria, ENSTA, Institut Polytechnique de Paris, Palaiseau, France
*Email: simone.pescuma@enpc.fr

## Abstract

We consider the iterative finite element solution of time-harmonic problems with a hybridized discontinuous Galerkin (HDG) method. While standard HDG methods use numerical traces in the hybridization process, our approach, called CHDG, use transmission variables. For acoustic problems, this modification accelerates the convergence of standard iterative procedures. Here, CHDG is extended to general problems with constant coefficients. We prove that the hybridized system can be solved with a fixed point. The method is applied the propagation of acoustic waves in a constant subsonic mean flow.

## 1 General model

Let an open bounded domain $\Omega \subset \mathbb{R}^3$ and a non-overlapping partition of its boundary $\partial\Omega = \cup_{i=1\dots N_{\text{bnd}}} \overline{\Gamma_i}$ (with $N_{\text{bnd}} \geq 1$). The $l$-dimensional complex unknown field $\boldsymbol{U}(\boldsymbol{x})$ verifies

$$\begin{cases} -\imath\omega \boldsymbol{U} + \displaystyle\sum_{j=1}^{3} \boldsymbol{A}_j \frac{\partial \boldsymbol{U}}{\partial x_j} = \boldsymbol{S}_{\text{vol}}, & \text{in } \Omega, \\ \boldsymbol{G}_i^- - \boldsymbol{\alpha}_i \boldsymbol{G}_i^+ = \boldsymbol{S}_{\text{sur},i}, & \text{on each } \Gamma_i, \end{cases} \tag{1}$$

with real, constant, symmetric square matrices $\{\boldsymbol{A}_j\}_j$ and complex data $\boldsymbol{S}_{\text{vol}}$ and $\{\boldsymbol{S}_{\text{sur},i}\}_i$. For each region $\Gamma_i$, the vectors $\boldsymbol{G}_i^+(\boldsymbol{x})$ and $\boldsymbol{G}_i^-(\boldsymbol{x})$ correspond to transmission variables *(see below)*, and $\boldsymbol{\alpha}_i$ is a constant rectangular coefficient matrix that verifies

$$\rho(\boldsymbol{\alpha}_i^\dagger \boldsymbol{\alpha}_i) \leq 1. \tag{2}$$

This assumption ensures that the boundary condition is passive.

## 2 DG scheme and transmission variables

Let $\boldsymbol{\mathcal{V}}_h^l := \prod_{K\in\mathcal{T}_h} \boldsymbol{\mathcal{P}}_{\text{p}}^l(K)$ be the vector broken polynomial space defined on a mesh $\mathcal{T}_h = \bigcup K$. Defining $(\cdot,\cdot)_{\mathcal{T}_h} := \sum_K (\cdot,\cdot)_K$ and $\langle\cdot,\cdot\rangle_{\partial\mathcal{T}_h} := \sum_K \sum_{F\subset\mathcal{F}_K} (\cdot,\cdot)_F$, where $\mathcal{F}_K$ is the collection of faces of $K$, the DG formulation of (1) reads

**Problem 1** *Find* $\boldsymbol{U}_h \in \boldsymbol{\mathcal{V}}_h^l$ *s.t.,* $\forall \boldsymbol{V}_h \in \boldsymbol{\mathcal{V}}_h^l$,

$$-\imath\omega \big(\boldsymbol{U}_h, \boldsymbol{V}_h\big)_{\mathcal{T}_h} - \sum_{j=1}^{3} \Big(\boldsymbol{A}_j \boldsymbol{U}_h, \frac{\partial \boldsymbol{V}_h}{\partial x_j}\Big)_{\mathcal{T}_h} + \big\langle \boldsymbol{F}\widehat{\boldsymbol{U}}_h, \boldsymbol{V}_h \big\rangle_{\partial\mathcal{T}_h} = 0,$$

with the flux matrix $\boldsymbol{F} := \sum_j n_j \boldsymbol{A}_j$, the outward unit normal vecror $\boldsymbol{n}$ and the numerical flux $\boldsymbol{F}\widehat{\boldsymbol{U}}_h$, which are defined face by face.

The *upwind flux* is defined thanks to the flux vector splitting method. Let $\boldsymbol{\Lambda}_{K,F}^\pm$ and $\boldsymbol{W}_{K,F}^\pm$ be the eigenvalues and eigenvectors matrices associated to the strictly positive/negative eigenvalues of $\boldsymbol{F}$. For an interior face $F \not\subset \partial\Omega$ of an element $K$ shared with $K'$, the upwind flux is

$$\boldsymbol{F}_{K,F}\widehat{\boldsymbol{U}}_F := \boldsymbol{W}_{K,F}^+ \sqrt{|\boldsymbol{\Lambda}_{K,F}^+|}\, \boldsymbol{G}_{K,F}^+ - \boldsymbol{W}_{K,F}^- \sqrt{|\boldsymbol{\Lambda}_{K,F}^-|}\, \boldsymbol{G}_{K,F}^-,$$

with the *outgoing* (+) and *incoming* (−) *transmission variables* defined as

$$\boldsymbol{G}_{K,F}^\pm := \sqrt{|\boldsymbol{\Lambda}_{K,F}^\pm|}\, [\boldsymbol{W}_{K,F}^\pm]^\dagger \boldsymbol{U}_{K,F},$$

which are associated to traveling information. Their sizes correspond to the number of strictly positive/negative eigenvalues, denoted by $n_{K,F}^\pm$. The definition is adapted for boundary faces.

## 3 Hybridized formulation

The DG problem is hybridized using the incoming variable at each face. These variables are recast as the global variable $\boldsymbol{G}_h^-$ belonging to

$$\boldsymbol{\mathcal{G}}_h^- := \prod_{K\in\mathcal{T}_h}\prod_{F\in\mathcal{F}_K}\boldsymbol{\mathcal{P}}_{\mathrm{p}}^{n_{K,F}^-}(F).$$

Eliminating the physical unknown $\boldsymbol{U}(\boldsymbol{x})$ by solving element-wise systems then leads to

**Problem 2** *Find* $\boldsymbol{G}_h^- \in \boldsymbol{\mathcal{G}}_h^-$ *s.t.,* $\forall \boldsymbol{\xi}_h \in \boldsymbol{\mathcal{G}}_h^-$,

$$\langle \boldsymbol{G}_h^-, \boldsymbol{\xi}_h\rangle_{\partial\mathcal{T}_h} - \langle \Pi(\mathrm{S}(\boldsymbol{G}_h^-)), \boldsymbol{\xi}_h\rangle_{\partial\mathcal{T}_h} = \langle \boldsymbol{b}, \boldsymbol{\xi}_h\rangle_{\partial\mathcal{T}_h},$$

with a *global exchange operator* $\Pi$ that swaps the transmission variables at the interior faces and weakly enforces the boundary conditions at the boundary faces, a *global scattering operator* S that computes outgoing variables for given incoming variables by solving all the local systems, and a *global right-hand side* $\boldsymbol{b}$.

The local problem is well-posed and the following theorem holds

**Theorem 3** *The operator* $\Pi\mathrm{S}$ *is a strict contraction:* $\|\Pi\mathrm{S}(\boldsymbol{G}_h^-)\| < \|\boldsymbol{G}_h^-\|$, $\forall \boldsymbol{G}_h^- \in \boldsymbol{\mathcal{G}}_h^- \setminus \{\boldsymbol{0}\}$.

Consequently, the fixed-point iteration applied to the CHDG system converge. This system is well-posed and equivalent to the DG system.

## 4 Numerical results in aeroacoustics

We denote by $p$ and $\boldsymbol{u}$ the fluctuations of pressure and velocity, and by $\rho_0$, $c_0$ and $\boldsymbol{u}_0$ the constant density, sound speed and subsonic mean flow. Wave propagation is modeled using

$$\boldsymbol{U} := \begin{bmatrix} p \\ \rho_0 c_0 \boldsymbol{u}\end{bmatrix}, \quad \boldsymbol{A}_j := \begin{bmatrix} \boldsymbol{e}_j\cdot\boldsymbol{u}_0 & c_0\boldsymbol{e}_j^\dagger \\ c_0\boldsymbol{e}_j & (\boldsymbol{e}_j\cdot\boldsymbol{u}_0)\boldsymbol{I}_3\end{bmatrix}$$

in (1), with proper data and boundary conditions. The passivity condition (2) holds if the impedance condition $p - \rho_0 c_0 \boldsymbol{n}\cdot\boldsymbol{u} = s_{\mathrm{imp}}$ is prescribed on the part of $\partial\Omega$ where the background flow is grazing or outgoing ($\boldsymbol{u}_0\cdot\boldsymbol{n} \geq 0$). When the background flow is incoming ($\boldsymbol{u}_0\cdot\boldsymbol{n} < 0$), a boundary condition on the tangential component of $\boldsymbol{u}$ must be prescribed. The CHDG formulation can then be applied. By contrast with formulations proposed for the convected Helmholtz equation (e.g. [3]), this formulation can deal with vorticity waves.

To validate the method, we compute the field generated by a Dirac point source at position $\boldsymbol{x}_s$ in a uniform horizontal mean flow $\boldsymbol{u}_0 = u_0\boldsymbol{e}_x$ within the unit disc in $\mathbb{R}^2$ (see figure). The impedance condition is prescribed on $\partial\Omega$, and

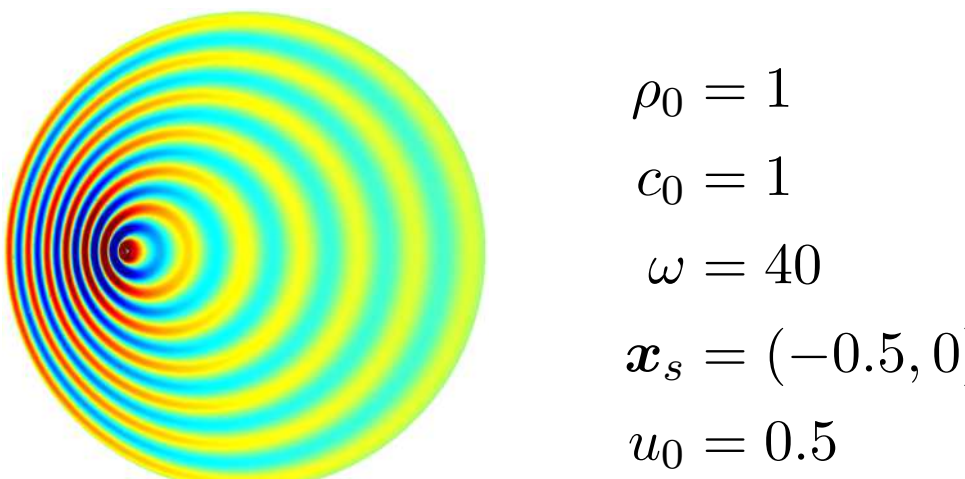

$\rho_0 = 1$
$c_0 = 1$
$\omega = 40$
$\boldsymbol{x}_s = (-0.5, 0)$
$u_0 = 0.5$

a condition on the tangential components of $\boldsymbol{u}$ is prescribed at the inflow boundary, with data computed with a reference solution.

We compare the relative error (computed with the analytical solution) with GMRES and CGNR (i.e. the conjugate gradient method applied to the normal system) for both DG and CHDG, and the fixed-point iteration for CHDG. We pick $h = 1/25$ and $\mathrm{p} = 3$.

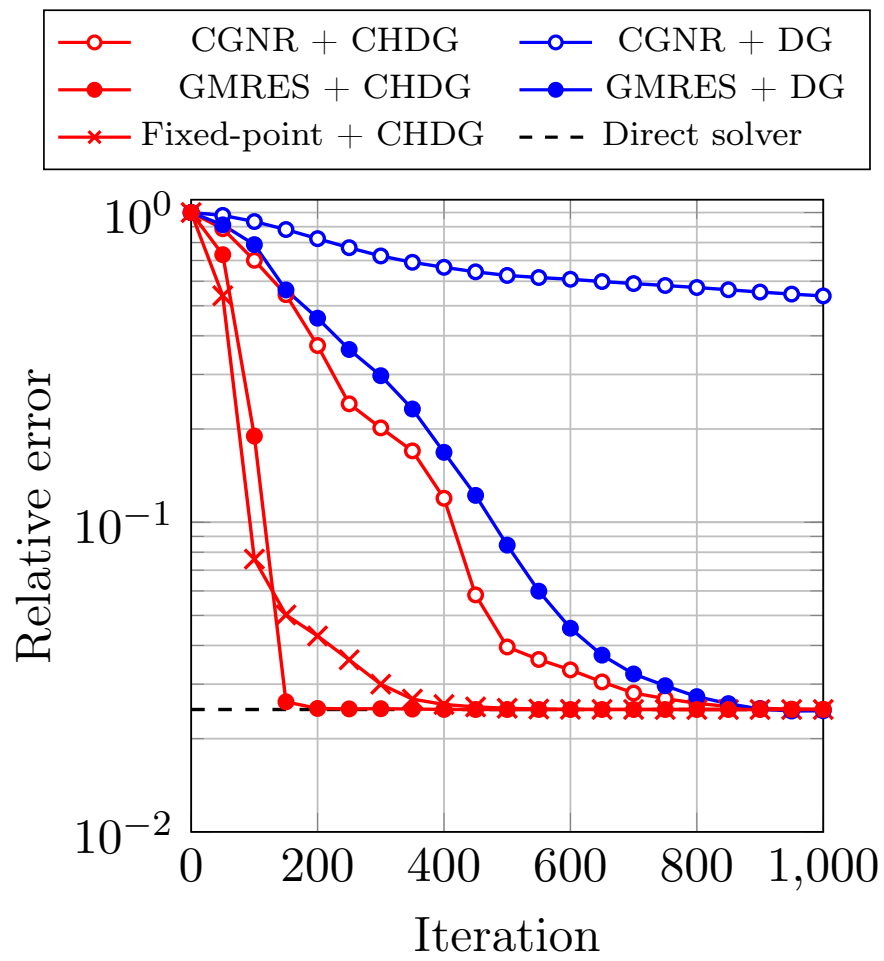

As expected, the fixed point applied to the CHDG system converges. Additionally, CHDG yields faster iterative solutions than DG. These results are consistent with those obtained for acoustic waves, e.g. in [1, 2]. Additional results confirm the method's robustness with respect to boundary conditions and parameters. More details will be shared in the extended paper.

*This work was partially supported by the ANR JCJC project WavesDG (ANR-21-CE46-0010).*

# Towards a fast conservative methods for the numerical resolution of Maxwell's equations in multipatch domains

**Ángel Pita-da-Veiga**[1,*], **Jerónimo Rodríguez**[1,2], **Rafael Vázquez**[1,2]

[1]Department of Applied Mathematics, Universidade de Santiago de Compostela, 15782 Santiago de Compostela, A Coruña (Spain)

[2]Galician Centre for Mathematical Research and Technology (CITMAga), 15782, Santiago de Compostela, A Coruña (Spain)

*Email: angel.pitadaveiga@usc.es

## Abstract

This talk aims to show how the numerical method presented in [1] can be generalized to complex geometries. The starting point is a fast, efficient, and conservative method for the numerical solution of Maxwell's equations that employs isogeometric analysis; that is, it constructs discrete spaces using splines to discretize the primal and dual de Rham sequences. Building on this formulation for single-patch geometries, we show how it can be extended to more complex domains using discontinuous Galerkin methodology to couple between patches. To achieve this, it is necessary to introduce a modification in the dual sequence, with the drawback that spurious eigenvalues are introduced in the associated eigenvalue problem. The overall method was implemented in GeoPDEs library. The numerical results show optimal order of convergence (resp. subobtimal orders in half order) for upwind fluxes (resp. for centered fluxes).

## 1 Singlepatch case

We consider the time-domain Maxwell equations in terms of differential forms:

$$\left|\begin{aligned} \mathrm{d}^2 D &= Q, \\ \mathrm{d}^2 B &= 0, \\ \partial_t D - \mathrm{d}^1 H &= -J, \\ \partial_t B + \mathrm{d}^1 E &= 0, \end{aligned}\right. \tag{1}$$

where the unknowns are $E \in \mathrm{H}_0\Lambda^1(\Omega)$, $H \in \mathrm{H}\Lambda^1(\Omega)$, $B \in \mathrm{H}_0\Lambda^2(\Omega)$ and $D \in \mathrm{H}\Lambda^2(\Omega)$. The model is complemented with the constitutive laws $B = \star^1_\mu H$ and $D = \star^1_\varepsilon E$, or their equivalent reformulation $H = \star^2_{1/\mu} B$ and $E = \star^2_{1/\varepsilon} D$. The first two equations of (1) can be discarded if the compatibility condition $\mathrm{d}^2 J + \partial_t Q = 0$ is assumed. We consider initially a singlepatch domain $\Omega$, that is, a domain that is the image of the unit cube by a spline transformation. Then, on each patch, to discretize the spaces of the de Rham sequences, we consider discrete spaces $X^k_h$ and $\tilde{X}^k_h$ defined by tensor-product constructions, one starting with degree $p$ and the other with degree $p-1$, which form the following discrete sequences:

$$X^0_{h,0} \xrightarrow{\mathrm{d}^0} X^1_{h,0} \xrightarrow{\mathrm{d}^1} X^2_{h,0} \xrightarrow{\mathrm{d}^2} X^3_{h,0}$$

$$\tilde{X}^3_h \xleftarrow{\mathrm{d}^2} \tilde{X}^2_h \xleftarrow{\mathrm{d}^1} \tilde{X}^1_h \xleftarrow{\mathrm{d}^0} \tilde{X}^0_h$$

To proceed, we must provide a definition a the discrete Hodge operator. One possible choice is the following: if $\omega \in \mathrm{H}\Lambda^k(\Omega)$, then there exists a unique $\star^k_{h,\gamma}\omega_h \in \tilde{X}^{3-k}_h$ such that

$$\langle \omega, \mu_h \rangle_{\mathrm{L}^2_\gamma \lambda^k(\Omega)} = \int_\Omega \star^k_{h,\gamma}\omega \wedge \mu_h, \quad \forall \mu_h \in X^k_{h,0}.$$

It is analogous for $\tilde{\star}^k_{h,\gamma}\omega \in X^{3-k}_{h,0}$. Once these discretizations of both the differential form spaces and the discrete Hodge operator are combined, the following matrix formulations are obtained:

$$\left|\begin{array}{rl} (\mathrm{I}) & \partial_t \mathbf{D}_h = \tilde{\mathbb{D}}_1 \mathbf{H}_h - \mathbf{J}_h, \\ (\mathrm{II}) & \mathbb{M}_{1,\varepsilon}\mathbf{E}_h = \tilde{\mathbb{K}}_2 \mathbf{D}_h \text{ or } \mathbb{K}_1 \mathbf{E}_h = \tilde{\mathbb{M}}_{2,1/\varepsilon}\mathbf{D}_h, \\ (\mathrm{III}) & \partial_t \mathbf{B}_h = -\mathbb{D}_1 \mathbf{E}_h, \\ (\mathrm{IV}) & \tilde{\mathbb{M}}_{1,\mu}\mathbf{H}_h = \mathbb{K}_2 \mathbf{B}_h \text{ or } \tilde{\mathbb{K}}_1 \mathbf{H}_h = \mathbb{M}_{2,1/\mu}\mathbf{B}_h, \end{array}\right. \tag{2}$$

where the $\mathbb{M}$ matrices are mass matrices, the $\mathbb{K}$ matrices are pairing matrices, and the $\mathbb{D}$ matrices represent the exterior derivative on finite-dimensional subspaces. There are two (non-equivalent) possible options in equations (II) and (IV), depending on the discrete constitutive laws considered. If $B = \star^1_\mu H$ and $D = \star^1_\varepsilon E$ are adopted, then it is necessary to invert the mass matrices in order to compute $E$ and $H$ as functions of $D$ and $B$, respectively. Conversely, if equalities $H = \star^2_{1/\mu} B$ and $E = \star^2_{1/\varepsilon} D$ are used, pairing matrices must be inverted instead. The second choice is significantly more efficient in terms of computational cost, since the coupling matrices are obtained as Kronecker products of matrices defined on the parametric domain. By applying the leapfrog method to either choice of (2), two numerical schemes, that are of high-order in space, are obtained and they preserve exactly both electromagnetic energy and charges. For more details, we refer to [1].

## 2 Multipatch case

We consider now a multipatch geometry, that is, of the form $\Omega = \cup_{\ell=1,\dots,N} K_\ell$, where each patch $K_\ell$ is the image of the unit cube under a spline mapping. We are going to take the following global discrete spaces:

$$X_h^k(\mathcal{T}) := \left\{ \omega_h \colon \ \omega_h|_K \in X_h^k(K), \forall K \in \mathcal{T} \right\},$$

$$\tilde{X}_h^k(\mathcal{T}) := \left\{ \omega_h \colon \ \omega_h|_K \in \tilde{X}_h^k(K), \forall K \in \mathcal{T} \right\},$$

where $X_h^k(K)$ and $\tilde{X}_h^k(K)$ are singlepatch spaces. These singlepatch spaces are different from the ones used in [1], since now we are not strongly imposing essential boundary conditions in the exterior of each patch. They will be imposed weakly in the exterior of the global boundary through numerical fluxes. A modification of the dual sequence is also needed to recover square pairing matrices in the second scheme. A possible choice for the DG fluxes is:

$$\begin{aligned} \mathcal{F}_K^{\mathrm{upw},\alpha}(\omega_1, \omega_2, \mu_1, \mu_2) &= \frac{1}{2} \left( \mathrm{tr}_K \, \omega_1 + \mathrm{tr}_{K'} \, \omega_2 \right) \\ &\quad + \frac{\alpha}{2} \star^{1,p-1}_{h,\sqrt{\mu/\varepsilon}} \left( \star^1 \mathrm{tr}_K \, \mu_1 - \star^1 \mathrm{tr}_{K'} \, \mu_2 \right), \\ \tilde{\mathcal{F}}_K^{\mathrm{upw},\alpha}(\omega_1, \omega_2, \mu_1, \mu_2) &= \frac{1}{2} \left( \mathrm{tr}_K \, \mu_1 + \mathrm{tr}_{K'} \, \mu_2 \right) \\ &\quad - \frac{\alpha}{2} \star^{1,p-1}_{h,\sqrt{\varepsilon/\mu}} \left( \star^1 \mathrm{tr}_K \, \omega_1 - \star^1 \mathrm{tr}_{K'} \, \omega_2 \right). \end{aligned}$$

For $\alpha = 0$, they are centered fluxes and for $\alpha = 1$ they will be called upwind fluxes, even though values of $\alpha \in (0, 1)$ are also allowed. The results obtained with these fluxes are shown in Figure 1, for three different methods: inverting mass matrices and inverting mass and pairing matrices with the modified dual sequence. We note that convergence is suboptimal for the centered flux, and optimal for the upwind flux with high degree. The modification of the dual sequence introduces spurious eigenvalues in the approximated solution. We are currently working on filtering or reducing these spurious results.

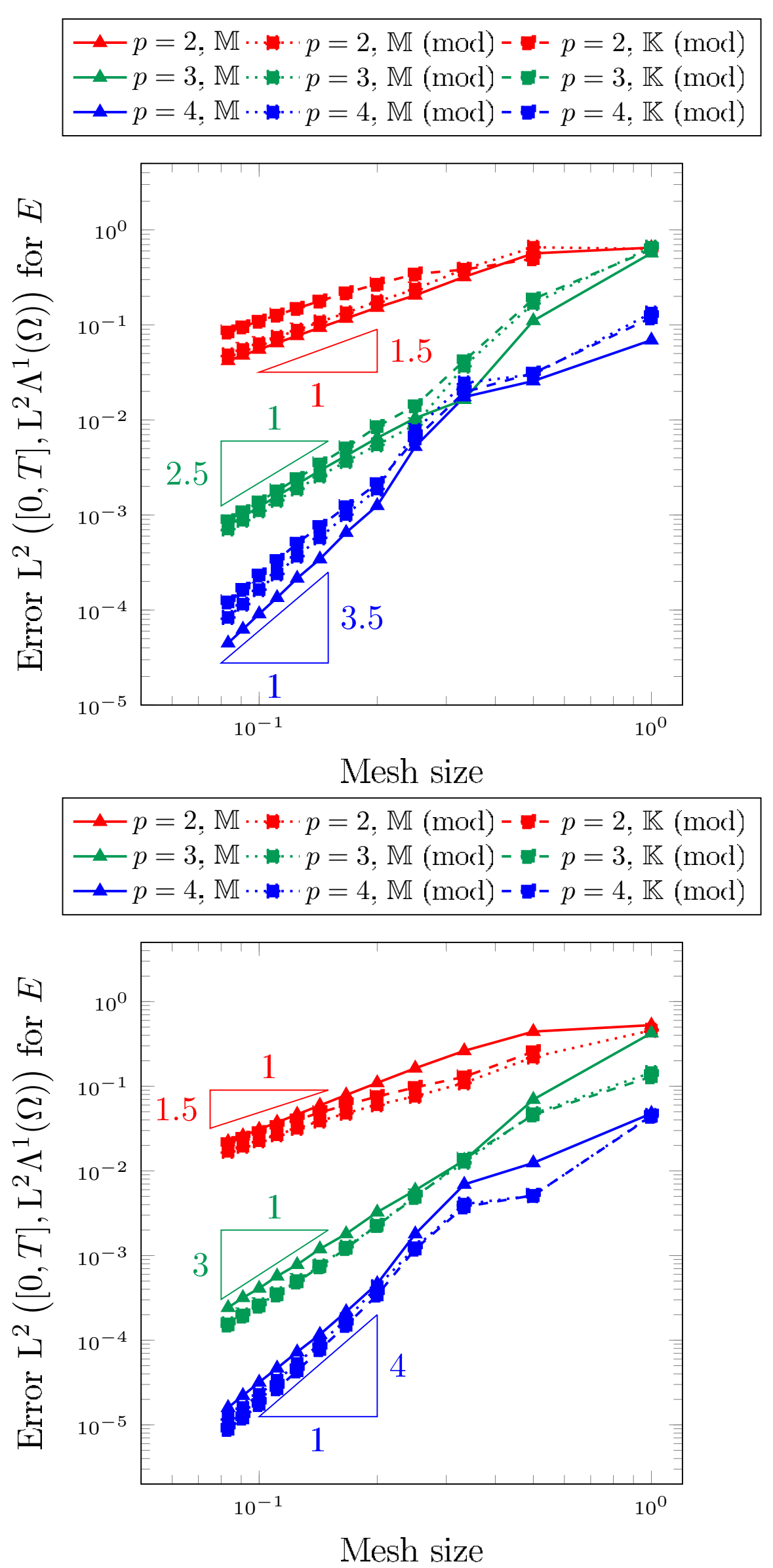

Figure 1: Curves of convergence of the error of the electric 1-form $E$ for the centered (top figure) and upwind flux (low figure).

# Some approaches for the field extrapolation in the context of inverse magnetisation problem and beyond

**Dmitry Ponomarev**[1,*]
[1]Centre Inria d'Université Côte d'Azur
*Email: dmitry.ponomarev@inria.fr

## Abstract

We review several methods of the data extrapolation in the inverse magnetisation problem arising when modern magnetometric techniques are used in paleomagnetism aiming to estimate remanent rock magnetisation. While some of these approaches have recently been developed, some are totally new and some are still under development. The overall goal is to replace the original inverse problem with the one that is less ill-posed, namely, that having a unique solution and less severe instability. We illustrate the results numerically and discuss implications and possible generalisations.

## 1 Introduction

The size of the region of available data plays a crucial role in solving inverse problems. Motivated by a particular paleomagnetic context of the inverse magnetisation problem (i.e., an inverse source problem for the Laplace equation with a source in the divergence form), where a magnetic distribution is to be reconstructed from partial measurements of its magnetic field, we consider the data extrapolation problem. That is, we study the mapping of the measured data on a limited region into the data on a larger region which is unaccessible for reliable measurements. While such a problem is still ill-posed, it can be shown to be less unstable and, at least for the ideal set of data, having a unique solution, unlike the full inverse magnetisation problem which is known to admit silent sources. The magnetic field extrapolation problem can thus be viewed as a preliminary data processing step whose reliable solution should improve solution of the inverse magnetisation problem, regardless of the choice of a particular approach for the latter.

## 2 Model

A SQUID (Superconducting QUantum Interference Device) magnetometer used for measuring rock magnetisation such as that in the Paleomagnetism lab of Earth, Atmospheric & Planetary Sciences department at Massachusetts Institute of Technology (USA), measures one component of the magnetic field which is vertical to the measurement plane. For planar magnetisation distribution (corresponding to a typical thinly sliced magnetised rock sample), we have the following explicit relation between the magnetisation vector $\vec{M} \equiv (M_1, M_2, M_3)^T$ and the magnetic field component $B_3$ measured at height $h > 0$ on a patch $Q \subset \mathbb{R}^2$ corresponding to the localisation region of the sample but on the plane parallel to it

$$B_3(\mathbf{x}, h) = \sum_{j=1}^{3} (K_j \star_Q M_j)(\mathbf{x}), \quad \mathbf{x} \in Q.$$

Here, the auxiliary kernel functions $K_{1,2}(\mathbf{x}) := -\frac{1}{2}\partial_{x_{1,2}} p_h(\mathbf{x})$, $K_3(\mathbf{x}) := -\frac{1}{2}\partial_h p_h(\mathbf{x})$ are defined in terms of the 3-dimensional Poisson kernel for a half-space $p_h(\mathbf{x}) := \frac{h}{2\pi\left(|\mathbf{x}|^2+h^2\right)^{3/2}}$, and we employ the notations $\mathbf{x} \equiv (x_1, x_2)^T$, $(K \star_Q f)(\mathbf{x}) \equiv \int_Q K(\mathbf{x}-\mathbf{t}) f(\mathbf{t})\, d\mathbf{t}$.

## 3 Field extrapolation

We are interested in obtaining the extrapolated field $B_3^{\text{ext}} = B_3|_{\mathbb{R}^2}$ based on the knowledge of the measured field $B_3^{\text{meas}} = B_3|_Q$. The possibility of unique extrapolation was shown in [2] and the so-called double-spectral approach was considered. In the present work, we propose and compare several extrapolation methods which are simpler to implement and which all can be put in the framework where the extrapolated field is computed, for sufficiently large $N \in \mathbb{N}$, as

$$B_3^{\text{ext}}(\mathbf{x}) \simeq \sum_{n=1}^{N} \langle B_3^{\text{meas}}, \varphi_n \rangle_Q \, \widetilde{\varphi_n}(\mathbf{x}), \quad \mathbf{x} \in \mathbb{R}^2,$$

where $\langle \cdot, \cdot \rangle_Q$ is the inner product in $L^2(Q)$, $\varphi_n$ is an eigenfunction of some self-adjoint positive integral operator $\mathcal{K} : L^2(Q) \to L^2(Q)$, $f \mapsto \mathcal{K}[f] = \int_Q K(\cdot, \mathbf{t}) f(\mathbf{t}) \, \mathrm{d}\mathbf{t}$, with the kernel function $K(\mathbf{x}, \mathbf{t})$ defined for $\mathbf{t} \in Q$, $\mathbf{x} \in \mathbb{R}^2$, and

$$\widetilde{\varphi_n}(\mathbf{x}) := \frac{1}{\lambda_n} \int_Q K(\mathbf{x}, \mathbf{t}) \varphi_n(\mathbf{t}) \, \mathrm{d}\mathbf{t}, \quad \mathbf{x} \in \mathbb{R}^2,$$

can be viewed as the continuation of $\varphi_n$ from $Q$ to $\mathbb{R}^2$, with $\lambda_n > 0$ the corresponding eigenvalue.

*Method 1:*
$K(\mathbf{x}, \mathbf{t}) = K_3(\mathbf{x} - \mathbf{t}) := -\frac{1}{2} \partial_h p_h(\mathbf{x} - \mathbf{t})$.

*Method 2:* $K(\mathbf{x}, \mathbf{t}) = K_{12}(\mathbf{x}, \mathbf{t}) := \frac{1}{2} p_h(\mathbf{x} - \mathbf{t})$.

*Method 3:* $K(\mathbf{x}, \mathbf{t}) = K_{12}^{[2]}(\mathbf{x}, \mathbf{t}) + K_3^{[2]}(\mathbf{x}, \mathbf{t})$,
$K_{12}^{[2]}(\mathbf{x}, \mathbf{t}) := \frac{1}{4} \langle p_h(\mathbf{x} - \cdot), p_h(\cdot - \mathbf{t}) \rangle_Q$,
$K_3^{[2]}(\mathbf{x}, \mathbf{t}) := \frac{1}{4} \langle \partial_h p_h(\mathbf{x} - \cdot), \partial_h p_h(\cdot - \mathbf{t}) \rangle_Q$.

*Method 4:* $K(\mathbf{x}, \mathbf{t}) = \sum_{j=1}^{3} K_j^{[2]}(\mathbf{x}, \mathbf{t})$,
$K_{1,2}^{[2]}(\mathbf{x}, \mathbf{t}) := \frac{1}{4} \langle \partial_{x_{1,2}} p_h(\mathbf{x} - \cdot), \partial_{x_{1,2}} p_h(\cdot - \mathbf{t}) \rangle_Q$.

*Method 5:* $K(\mathbf{x}, \mathbf{t}) = \frac{1}{16} \partial_h^2 p_h(\mathbf{x} - \mathbf{t})$.

Methods 1–3 assume additional regularity on the magnetisation: $\partial_{x_1} M_1 + \partial_{x_2} M_2 \in L^2(Q)$ in addition to $M_3 \in L^2(Q)$, whereas Methods 4 and 5 require only $L^2$ regularity of the magnetisation components. Moreover, Method 5, being an alternative formulation of the approach proposed in [1], is actually applicable to volumetric magnetisations through an equivalent planar magnetisation distribution with an infinite support. This method is, however, sub-optimal for planar samples (see Fig. 1) as it fails to incorporate the known support information.

Furthermore, we note that Methods 1 and 2 are simple but limited to extrapolation only in a sense that their by-products cannot be used to further solve the full inverse magnetisation problem. In contrast to that, the other proposed methods are based on reconstruction of some magnetisation-dependent quantities from which, as we show, it is possible to advance towards the full magnetisation inversion, hence these methods can be seen as a partial inversion.

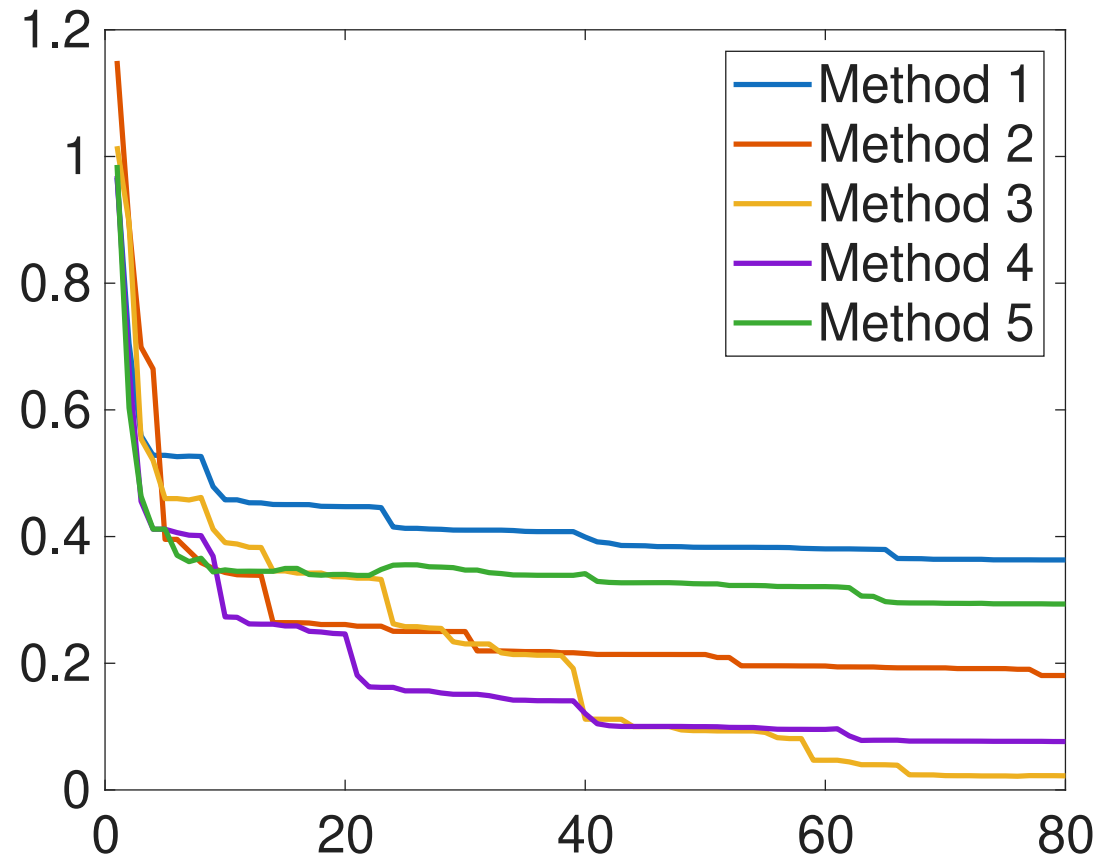

Figure 1: Relative error in the extrapolation from $Q = (-1, 1)^2$ to $(-10, 10)^2$ at height $h = 1$ for the field of a planar synthetic magnetisation supported on $Q$ versus the number of eigenfunctions $N$.

## 4 Discussion and extensions

One of the initial motivation for the field extrapolation was the asymptotic formulas for the net magnetisation (moment) reconstruction from a large region with known magnetic field, see [3]. However, we show that Methods 3–5 allow direct alternative expressions for the magnetisation moment in terms of the auxiliary quantities which are by-products of the field extrapolation.

Finally, we mention that our results are expected to be extendable to other inverse source problems in potential theory where the data of only one field component is practically available.

# Transparent boundary conditions for the stellar oscillation equations without Cowling approximation

Janosch Preuss[1,*], Hélène Barucq[1], Florian Faucher[1], Ha Pham[1], Sébastien Tordeux[1]

[1]Project-Team Makutu, Inria, Université de Pau et des Pays de l'Adour, France

*Email: janosch.preuss@inria.fr

**Abstract**

We present a transparent boundary condition (tbc) for the gravity potential equation appearing in the equations modeling stellar oscillations. This tbc is constructed with the learned finite element method which exploits the separability of the Dirichlet-to-Neumann ($\mathcal{DtN}$) operator in specific settings and formulates its approximation as a scalar optimization problem.

## 1 Introduction

The system of equations modeling stellar oscillations (Unno et al, 1979) consists of the equation of motion coupled with the Poisson equation. The latter which models perturbations to the gravitational potential is often ignored in the Cowling's approximation (Cowling, 1941). This however is only valid for high frequencies and modes, and its removal is particularly important for distant stars due to low-resolution observations.

The Poisson equation is posed on the whole of $\mathbb{R}^3$, however finite element simulations require a bounded domain. In helioseismology truncation is usually carried out with a tbc placed slightly above observation height. In recent years, under Cowling's approximation, [1, 2] have developed tbcs for the equation of motion to model wave propagation in and interaction with the atmosphere. In this work, with the eventual goal of solving the full system, we develop a tbc for the Poisson equation.

There is a vast literature on solving the Poisson equation on an unbouded domain. Some tbcs are based on a direct discretization of nonlocal operators on the trunction boundary like boundary elements or the DtN-FE method which leads to dense submatrices. Other approaches like classical infinite elements methods (CIE) avoid this by approximating the exact tbc $\mathcal{DtN}$ by sparse matrices whose entries are determined analytically. This however requires homogeneous exterior.

In [3] learned infinite elements (LIE) have been introduced for time-harmonic wave equations on stratified exterior domains with separable geometry. Unlike in CIE, the matrix entries are obtained with a small minimization problem fitting the action of $\mathcal{DtN}$ on a set of eigenmodes. This feature allows LIE to be applicable on a domain with stratified inhomogeneity. Additionally, it allows to reduce the dimension and optimize the sparsity of these matrices rendering the method very efficient. This motivates us to apply LIE to the gravity potential equation.

## 2 Learned infinite elements for the Laplace equation

Consider the Poisson equation in $\mathbb{R}^3$ with source $f$ compactly supported in $B_a := \{|x| < a\}$

$$-\Delta\phi = f, \quad \text{in } \mathbb{R}^3, \qquad \lim_{|x|\to\infty} \phi(|x|) = 0,$$

Denote by $\Gamma := \{|x| = a\}$. The exact exterior $\mathcal{DtN}$ truncates the problem exactly to $B_a$

$$-\Delta\phi = f, \quad \text{in } B_a, \qquad -\partial_{\mathbf{n}}\phi = \mathcal{DtN}\phi \text{ on } \Gamma.$$

Denote by $(w_\ell, \lambda_\ell)$ the eigenpairs of the Laplace-Beltrami operator $\Delta_\Gamma$ on $\Gamma$

$$-\Delta_\Gamma w_\ell = \lambda_\ell w_\ell, \text{ with } \lambda_\ell = \ell(\ell+1)/a^2.$$

In this basis the $\mathcal{DtN}$ is diagonal with eigenvalue $dtn(\lambda_\ell)$:

$$\mathcal{DtN} w_\ell = dtn(\lambda_\ell) w_\ell, \quad dtn(\lambda) = \frac{1}{2a} + \sqrt{\frac{1}{4a^2} + \lambda}. \tag{1}$$

Due to its non-localness, naive discretizations of $\mathcal{DtN}$ in finite element methods result in dense matrices. Instead of discretizing the $\mathcal{DtN}$ directly on $\Gamma$, many tbc methods solve a sparse augmented system by adding an auxiliary block involving exterior degrees of freedom (dofs) so that the Schur complement of this block results in a good approximation of $\mathcal{DtN}$.

In [3] it is shown that if the augmented system is of the form $A \otimes M_\Gamma + B \otimes K_\Gamma$ with $A, B \in \mathbb{C}^{(N+1)\times(N+1)}$ for $N \in \mathbb{N}_0$ and stiffness matrix $K_\Gamma$ discretizing $-\Delta_\Gamma$ and mass matrix $M_\Gamma$ the identity on $\Gamma$, then the corresponding approximation $\mathrm{DtN}_N$ of $\mathcal{DtN}$ (obtained after taking the Schur complement) is diagonal in the discrete eigenbasis $\underline{w}_\ell$ of the matrix pair $(K_\Gamma, M_\Gamma)$, $K_\Gamma \underline{w}_\ell = \underline{\lambda}_\ell M_\Gamma \underline{w}_\ell$. That is, we have

$$\mathrm{DtN}_N\, \underline{w}_\ell = \mathrm{dtn}_N(\underline{\lambda}_\ell)\underline{w}_\ell, \tag{2}$$

where $\mathrm{dtn}_N$ is a rational function of $\underline{\lambda}_\ell$ depending on $A$ and $B$ only. Additionally, assuming a certain sparse and symmetric ansatz for $A$ and $B$ we have that $\mathrm{dtn}_N$ simplifies to

$$\mathrm{dtn}_N(\lambda) = A_{00} + \lambda B_{00} - \sum_{j=1}^{N} \frac{(A_{0j} + \lambda B_{0j})^2}{A_{jj} + \lambda}.$$

In view of (1) and (2), the question of finding a good approximation $\mathrm{DtN}_N$ of $\mathcal{DtN}$ is reduced to finding entries $A_{0j}$, $A_{jj}$ and $B_{0j}$ so that $\mathrm{dtn}_N$ best approximates $\mathcal{dtn}$. This is phrased as a weighted nonlinear optimization problem fitting the action of $\mathcal{DtN}$ and $\mathrm{DtN}_N$ on $L$ modes. The number of $L$ depends on the smoothness of the exact solution. Figure 1 shows that only very few terms $N$ are needed to obtain an excellent fit up to $\lambda_\ell = 210$, i.e. $\ell = 14$. The relative $L^2$- error when solving the potential equation with $\mathrm{DtN}_N$ imposed at $a = 1$ is shown in Table 1. As reference we impose on $\Gamma$ the value of the exact solution and solve a Dirichlet boundary value problem (denoted 'exact bc'). The reference accuracy is attained by $\mathrm{DtN}_N$ with $N = 2$. Since the convergence is exponential in $N$, the corresponding increase in number of dofs 'ndof' and number of non-zero entries 'nze' in the linear systems remains small.

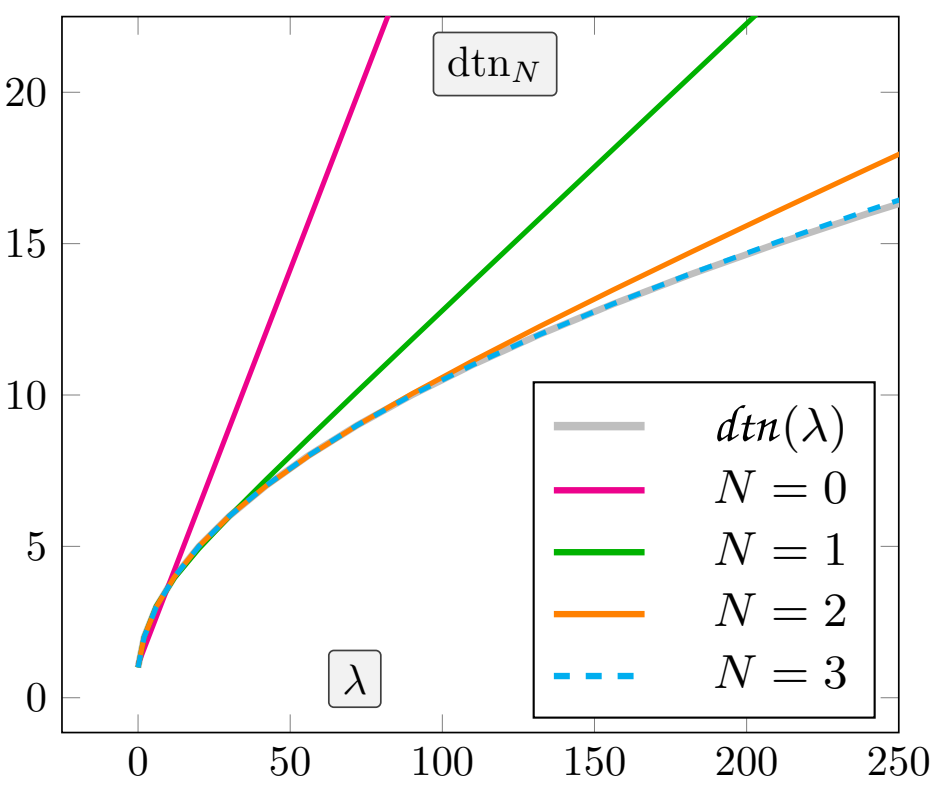

Figure 1: Fitting of $\lambda \mapsto \mathcal{dtn}(\lambda)$ by $\lambda \mapsto \mathrm{dtn}_N(\lambda)$.

| Method | ndof | nze | $L^2$-error |
|---|---|---|---|
| exact bc. | $2.4 \cdot 10^5$ | $9.35 \cdot 10^6$ | $2.52 \cdot 10^{-6}$ |
| $N = 0$ | $2.6 \cdot 10^5$ | $1.07 \cdot 10^7$ | $4.0 \cdot 10^{-2}$ |
| $N = 1$ | $2.8 \cdot 10^5$ | $1.13 \cdot 10^7$ | $1.3 \cdot 10^{-4}$ |
| $N = 2$ | $3.0 \cdot 10^5$ | $1.21 \cdot 10^7$ | $2.54 \cdot 10^{-6}$ |

Table 1: Result for solving the potential equation with our tbc.

## 3 Perspectives

The promising preliminary results encourage to applying LIE to the Poisson equation in the full system of stellar oscillation equations analogous to [4] for the Earth. This should not pose problems since these equations are not coupled at the level of the tbc. In this way we could extend the work in [5] which considers differential rotation and gravity, however in the Cowling approximation. The method also shows promises in computing eigenvalues for fast-rotating stars and on ellipse geometries [6].

# Geometry-dependent well-posedness and finite elements for shear Alfvén waves

M. Renoldner[1,*], A. Buffa[1], T. Miehling[1], M. Picasso[1]
[1]Department of Mathematics, EPFL, Lausanne, Switzerland
*Email: markus.renoldner@epfl.ch

## Abstract

Alfvén waves form the fastest modes of the drift-reduced Braginskii system and determine time-stepping constraints in plasma turbulence simulations [4]. We prove well-posedness of the Alfvén system in anisotropic Sobolev spaces under a geometric condition expressed by a perpendicular Poincaré inequality. A conservative finite-element scheme with a priori error bounds is presented. We show that the system's condition number depends on the same assumption as the perpendicular Poincaré constant.

## 1 The Alfven wave equations

Let $\Omega \subset \mathbb{R}^n$ with $n > 1$ be a bounded domain, $T$ the final time, and let $\mathbf{B}$ be the magnetic field, smooth and nonvanishing. We further set $\mathbf{b} := \frac{\mathbf{B}}{\|\mathbf{B}\|}$. The Alfvén wave equations [1–3] describe the $\mathbf{b}$-parallel electron velocity and the electric potential $u, \phi : \Omega\times[0, T] \to \mathbb{R}$, respectively, such that

$$\mu\,\partial_t u - \mathbf{b}\cdot\nabla\phi = f, \quad (1)$$

$$\partial_t \operatorname{div}(\nabla_\perp \phi) + \operatorname{div}(\mathbf{b}u) = g, \quad (2)$$

where $\mu > 0$ is a given, fixed, small parameter describing the ratio between the masses of different particle species. To simplify the presentation, we will assume $\mu = 1$ in the following. $f, g$ are given functions on $\Omega \times (0, T)$. We set homogeneous Dirichlet boundary condition for $\phi$. The operator in the first term in equation (2), is the *perpendicular gradient*

$$\nabla_\perp\phi := (I - \mathbf{b}\mathbf{b}^\top)\nabla\phi. \quad (3)$$

One mechanism guiding the analysis in sections 2 and 3 is the energy conservation law,

$$\frac{1}{2}\frac{\mathrm{d}}{\mathrm{d}t}\left(\|u\|^2_{L^2(\Omega)} + \|\nabla_\perp\phi\|^2_{L^2(\Omega)}\right) = \int_\Omega fu + \int_\Omega g\phi.$$

For $f = g = 0$ the energy is conserved. Controlling the energy for non-vanishing data requires a perpendicular Poincaré inequality for $\nabla_\perp$. This constitutes the main analytical obstacle.

## 2 Theoretical results

We introduce an orthonormal frame $\mathbf{b}_1, \dots, \mathbf{b}_n \in \mathrm{Lip}(\Omega, \mathbb{R}^n)$ with $\mathbf{b} := \mathbf{b}_1$, where the perpendicular gradient is realized as the collection of weak derivatives along the perpendicular frame vectors $\mathbf{b}_2, \dots, \mathbf{b}_n$. The corresponding anisotropic Sobolev space is then

$$H^1_\perp(\Omega) := \{\psi \in L^2(\Omega) \mid \nabla_\perp\psi \in L^2(\Omega)\},$$

and $H^1_{\perp,0}(\Omega)$ the closure of $C^\infty_c(\Omega)$. A key ingredient of the analysis is a perpendicular Poincaré inequality requiring one *globally directed* perpendicular frame vector which extends Azerad's result [5] from divergence-free to general fields.

**Theorem 1 (Poincaré-type inequality)** *Suppose that one frame vector $\mathbf{b}_j$ with $j \neq 1$ has bounded differential and is* globally directed *in the sense that $\mathbf{b}_j(x)\cdot\mathbf{k} > 0$ for some fixed $\mathbf{k} \in \mathbb{R}^n$.Then there exists a constant $C_P > 0$, only depending on $\Omega$ and $\mathbf{b}_j$, such that for all $\phi \in H^1_{\perp,0}(\Omega)$*

$$\|\phi\|_{L^2(\Omega)} \le C_P\|\nabla_\perp\phi\|_{L^2(\Omega)}. \quad (4)$$

We will now introduce a weak form of (1)-(2). We seek $u$ and $\phi$ such that for all $v \in L^2(\Omega)$, $\psi \in H^1_0(\Omega)$ we have

$$\int_\Omega \partial_t uv - \int_\Omega v\mathbf{b}\cdot\nabla\phi = \int_\Omega fv, \quad (5)$$

$$-\int_\Omega \partial_t\nabla_\perp\phi\cdot\nabla_\perp\psi - \int_\Omega u\mathbf{b}\cdot\nabla\psi = \int_\Omega g\psi, \quad (6)$$

and $u(0) = u_0$, $\phi(0) = \phi_0$ in $H^1(\Omega)$.

**Theorem 2 (Well-posedness)** *Under the perpendicular Poincaré inequality and standard regularity assumptions on the data, the Alfvén wave system* (5)-(6) *admits a unique solution*
$u \in H^1(0, T; L^2(\Omega))$
*and* $\phi \in L^2(0, T; H^1_0(\Omega)) \cap H^1(0, T; H^1_{\perp,0}(\Omega))$.

The standard domain of Alfvén waves is the *tokamak*, modelled by

$$\mathbb{T} := \left\{r^2_{\min} < x^2 + y^2 < r^2_{\max},\ 0 < z < z_{\max}\right\},$$

with $0 < r_{\min} < r_{\max}$ and $0 < z_{\max}$. In this case, the magnetic field is usually oriented towards the *toroidal* direction, meaning that $\mathbf{b}(x)\cdot\mathbf{e}_\varphi(x) > 0$, where $\{\mathbf{e}_\rho, \mathbf{e}_\varphi, \mathbf{e}_z\}$ is the cylindrical basis. One can show that Theorem 2 applies to this setting.

## 3 Numerical analysis

We introduce now a conservative Crank Nicolson in time and FEM in space discretisation of the weak form (5)-(6). Let $X_h^r(\Omega)$ be the usual Lagrange FEM space of order $r$. Let $\Delta t$ denote the time step size. We solve

$$\int_\Omega \frac{u_h^{n+1} - u_h^n}{\Delta t} v_h - \int_\Omega v_h \mathbf{b}\cdot\nabla \frac{\phi_h^{n+1} + \phi_h^n}{2} = \int_\Omega f^{n+\frac{1}{2}} v_h,$$

$$-\int_\Omega \frac{\nabla_\perp \phi_h^{n+1} - \nabla_\perp \phi_h^n}{\Delta t}\cdot\nabla_\perp \psi_h - \int_\Omega \frac{u_h^{n+1} + u_h^n}{2}\mathbf{b}\cdot\nabla\psi_h = \int_\Omega g^{n+\frac{1}{2}}\psi_h,$$

for all $v_h \in X_h^r(\Omega), \psi_h \in X_{h,0}^r(\Omega)$. The scheme preserves a discrete analogue of the continuous energy and is numerically stable.

**Theorem 3 (A priori error estimate)** *Suppose that the exact solutions $u, \phi$ from Theorem 2 are in $C^3(0,T;H^{r+1}(\Omega))$. Let $(u_h^M, \phi_h^M)$ be the solution of the FEM scheme. Then, there exists a constant $C > 0$, independent of $h$ and $\Delta t$ such that*

$$\|\epsilon_u^M\|_{L^2(\Omega)} + \|\nabla_\perp \epsilon_\phi^M\|_{L^2(\Omega)} \le C(h^r + \Delta t^2).$$

The discrete perpendicular Laplacian is invertible precisely when the perpendicular Poincaré inequality holds, i.e. when at least one perpendicular frame vector is globally directed. On $\Omega = \mathbb{R}/\mathbb{Z}\times(0,1)$ with Dirichlet boundaries in $y$, consider

$$\mathbf{b}_\epsilon : x \mapsto \left(\sqrt{1-\epsilon^2}, \epsilon\right), \quad \epsilon \in [0,1].$$

As $\epsilon \to 1$ the matrix condition number diverges, as shown in Fig. 1. This degeneration explains conditioning issues in Alfvén wave solvers and indicates necessity of the assumption on $\mathbf{b}_j$ stated in Theorem 1.

## 4 Conclusion

We establish well-posedness of Alfven waves and an energy-stable finite-element scheme whose conditioning reflects the magnetic geometry.

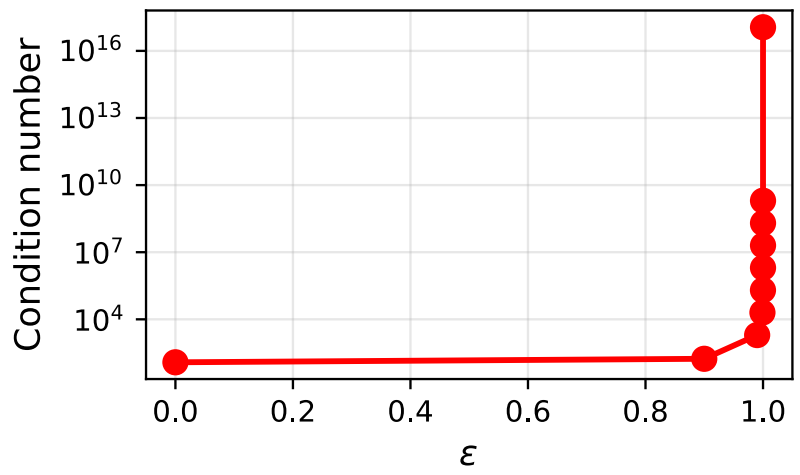

Figure 1: Condition number vs. $\epsilon$.

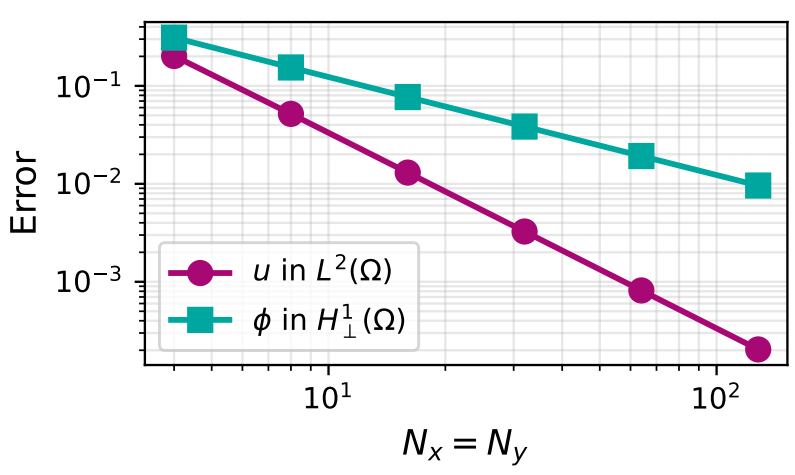

Figure 2: Convergence test with $r = 1$.

# Inversion-free binary segmentation for computed tomography and inverse scattering

**Jacob D. Rezac**[1,*], **Zachary J. Grey**[1], **Newell H. Moser**[1,2], **Jason P. Killgore**[1]
[1]National Institute of Standards and Technology, Boulder, USA
[2]Department of Physics, University of Colorado, Boulder, USA
*Email: jacob.rezac@nist.gov

## Abstract

Quantifying geometric information contained in objects from indirect measurements is a fundamental problem in non-destructive testing, e.g., inverse scattering problems or X-Ray computed tomography. A standard approach to identifying key geometric features – inversion followed by segmentation – can result in noisy, low-resolution estimates due to the ill-posed nature of corresponding inverse problems. In this work, we instead solve directly for binary segmentations of unknown object geometries. This approach identifies the support of internal features in an object by testing an indicator function against shapes from a given design dictionary. The technique is demonstrated on an example in X-Ray computed tomography[1].

## 1 Introduction

Segmenting the structures in an object's interior from indirect measurements is a fundamental challenge of non-destructive imaging. Our motivation is to segment without needing to explicitly calculate any inverse mappings between object interiors and measurements. In this work, we assume that object interiors are composed of simple and single-material shapes residing in a constant background material (e.g., air). An example of such an object is a semiconductor – see the X-Ray Computed Tomography (CT) scan shown in Figure 1. Semiconductors are composed primarily of well-defined shapes with sharp material transitions. The specific arrangement of these shapes depends on design and manufacturing details. We assume knowledge of the shapes, but not their location, orientation, or size within the object-of-interest.

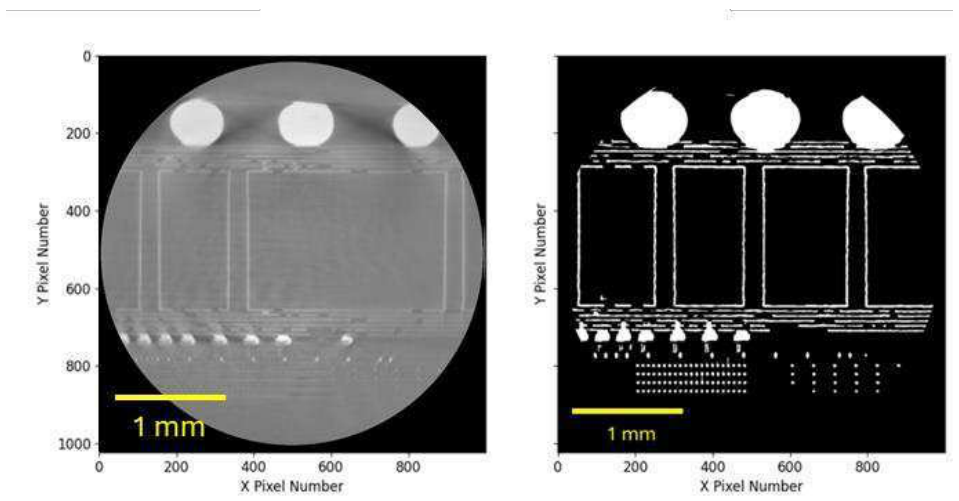

Figure 1: Reconstructed slice (left) and binary segmentation (right) from an X-Ray CT scan.

The technique we introduce below relies on a shape dictionary, $\mathcal{S}$, that describes probable combinations of these geometric features (e.g., translations, rotations, and scalings of simple shapes). Given $\mathcal{S}$, the goal of the work is to identify which elements in the shape dictionary are present in the unknown object *without solving any measurement-dependent problems.* Instead, the technique approximates a binary segmentation of the object with a quick-to-calculate indicator function of measured data and pre-computed functions of elements of $\mathcal{S}$.

If the unknown object's interior is described by a compactly supported function $m \in L^2(\mathbb{R}^d)$, $d = 2$ or 3, we consider noisy measurements from the physical model $y(\mathbf{x}) = (Km)(\mathbf{x})$ where $K : L^2(\mathbb{R}^d) \to L^2(\Gamma_m)$ is a linear operator that maps from $m$ to measurement locations $\Gamma_m \subset \mathbb{R}^d$. This work's applications-of-interest can be represented by

$$(Km)(\mathbf{x}) := \int_{\mathbb{R}^d} k(\mathbf{x}, \mathbf{z}) m(\mathbf{z}) \, \mathrm{d}V(\mathbf{x}) \qquad (1)$$

with a modality-specific kernel $k(\cdot, \cdot)$. This model includes, e.g., acoustic scattering under the Born approximation where $K$ is the Lippmann-Schwinger integral operator [1] and parallel beam X-Ray CT with a Radon transform operator denoted by $R$ [2]. For the latter, we use the fact that $(R^* Rm)(\mathbf{z})$ has a representation as (1) with kernel $k(\mathbf{x}, \mathbf{z}) = \|\mathbf{x} - \mathbf{z}\|^{-1}$, $\mathbf{x}, \mathbf{z} \in \mathbb{R}^d$ [2].

[1]Official contribution of the National Institute of Standards and Technology; not subject to copyright in the United States. This work was performed with funding from the CHIPS Metrology Program, part of CHIPS for America, National Institute of Standards and Technology, U.S. Department of Commerce.

## 2 Theoretical Results

Our method is qualitative in nature [1]; we seek an indicator function $I : \mathcal{S} \to \mathbb{R}$ for each element of the shape dictionary $\mathcal{S}$ that is large unless a test shape falls within an unknown object's true interior. The indicator function is defined via a generalized Rayleigh quotient and test shape $D_\tau \in \mathcal{S}$,

$$Q(\phi; D_\tau) = \frac{\int_{D_\tau} |(K^*\phi)(\mathbf{z})|^2 \, \mathrm{d}V(\mathbf{z})}{\int_{\mathbb{R}^d} |(K^*\phi)(\mathbf{z})|^2 \, \mathrm{d}V(\mathbf{z})}.$$

If $\phi_\tau \in \mathrm{argmax}_{\varphi \in L^2(\Gamma_m)} Q(\varphi; D_\tau)$, we define

$$I(D_\tau) = (\phi_\tau, y) - \mathrm{Vol}(D_\tau), \quad D_\tau \in \mathcal{S}$$

where $(u, v) = \int_{\Gamma_m} u(\mathbf{x})\overline{v(\mathbf{x})} \, \mathrm{d}V(\mathbf{x})$ and $\mathrm{Vol}(\cdot)$ is the volume of test shape. The following theorem demonstrates that *in the noise-free case when $m$ is an binary function*, $\mathrm{supp}(m)$ is approximately the union of $D_\tau \in \mathcal{S}$ where $I(D_\tau) = 0$.

**Theorem 1** *Assume $\mathbb{1}_{D_\tau} \in \mathcal{R}(K^*)$, $m$ is a binary function with compact support, and let $\varphi_\tau \in \arg\max_{\varphi \in L^2(\Gamma_m)} Q(\varphi; D_\tau)$. Then, $I(D_\tau) = 0$ if and only if $D_\tau \subseteq \mathrm{supp}(m)$.*

***Proof.*** Letting $\psi_\varphi = K^*\varphi$,

$$\varphi_\tau = \underset{\varphi \in L^2(\Gamma_m)}{\mathrm{argmax}} \frac{(\mathbb{1}_{D_\tau}\psi_\varphi, \psi_\varphi)}{\|\psi_\varphi\|^2}.$$

Thus, $\psi_\varphi(\mathbf{z}) = (K^*\varphi_\tau)(\mathbf{z}) = \mathbb{1}_{D_\tau}(\mathbf{z})$ and

$$(\varphi_\tau, Km) = (K^*\varphi_\tau, m) = \mathrm{Vol}(D_\tau \cap \mathrm{supp}(m)).$$

Hence, $I(D_\tau) = 0$ when $\mathrm{Vol}(D_\tau \cap \mathrm{supp}(m)) = \mathrm{Vol}(D_\tau)$, which is equivalent to $D_\tau \subseteq \mathrm{supp}(m)$. □

To account for errors from noise and for if $\mathbb{1}_{D_\tau} \notin \mathcal{R}(K^*)$, we binarize by plotting the smallest 10% of values of $|I(D_\tau)|$, $D_\tau \in \mathcal{S}$. Though the theorem above requires $m$ to be a binary function, we see in the example below that this may not be needed in practice.

## 3 Simulations and Commentary

We demonstrate the new procedure with simulations of a noisy 2-dimensional Radon transform using the X-Ray CT ASTRA toolbox [3]. We build $\mathcal{S}$ from all translations of a $10 \times 10$ pixel square across the $256 \times 256$ pixel image in Figure 2. The simulation includes Poisson noise with approximately 4 dB signal-to-noise ratio.

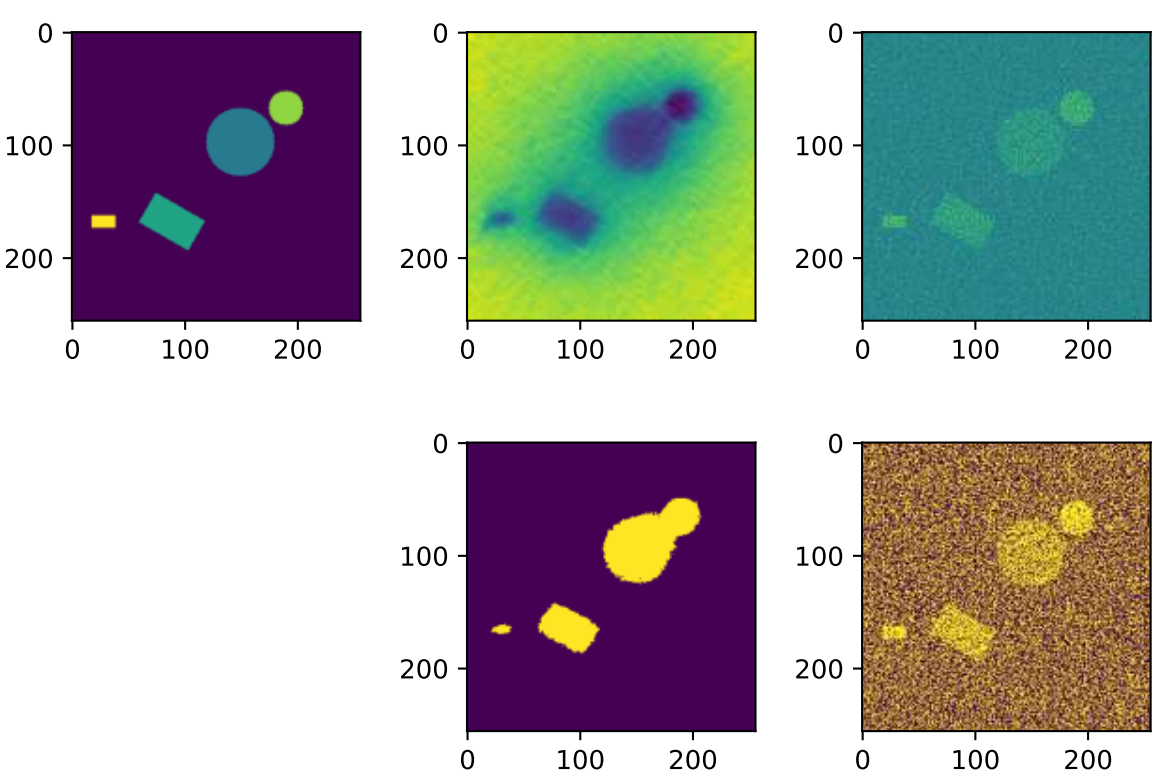

Figure 2: Simulation geometry (top-left); $I(D)$ (top-middle) and binarization (bottom-middle); filtered backprojection (top-right) and binarization (bottom-right). Colorbars differ between images except for the binarized bottom row.

There were 90 measurement angles in $[0°, 180°]$ and 256 detector locations.

Figure 2 (middle) shows the values of $I(D)$ (top) and its binarization (bottom). Figure 2 (right) also shows a filtered backprojection (FBP) reconstruction and its Otsu binarization. The FBP segmentation is more impacted by noise than our method.

Unlike standard qualitative methods that test against infinitesimally-small shapes, restricting to a pre-build shape dictionary further regularizes the problem. However, a poorly built shape dictionary reduces the method's effectiveness. To address this, segmentations can be used as initial guesses to shape-based optimization problems that further subselect the elements of $S$ that are found in $\mathrm{supp}(m)$, e.g., [4].

# Interactions between water waves and floating structures, a numerical method

**Geoffrey Beck**[1], **Alan Riquier** [2,*]

[1]Univ Rennes, IRMAR UMR 6625, Centre Inria de l'Université de Rennes (MINGuS), France

[2]Fields Institute for Research in Mathematical Sciences, Toronto, Canada

*Email: ariquier@fieldsinstitute.ca

**Abstract**

We propose a FEM-based numerical approach for wave–floating-structure interactions. The evolution of the 2D fluid is governed by the incompressible free-surface Navier-Stokes equations. The motion of the rigid body is prescribed by Newton's equations, accounting for both the gravitational and buoyancy forces.

## 1 Floating body and buoyancy

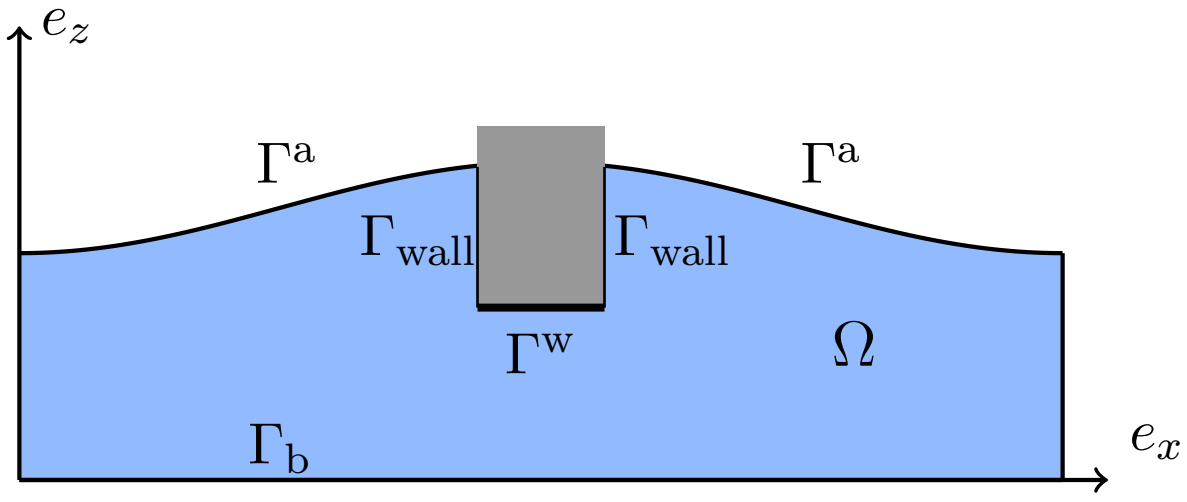

Figure 1: Schematic representation of the configuration.

Consider a freely-floating rigid body in a viscous two-dimensional fluid of density $\rho$ and viscosity $\nu$. The fluid domain is denoted $\Omega$. Its boundary is decomposed as $\partial\Omega = \Gamma^a \cup \Gamma^w \cup \Gamma_b \cup \Gamma_{wall}$ (Fig. 1). Due to the displacement of both the free surface and the floating body, the quantities $\Omega$, $\Gamma^a$, $\Gamma^w$ and $\Gamma_{wall}$ are, in fact, time-dependent. For the sake of simplicity and brevity, we shall only consider the case of a rectangular body moving in the mere vertical direction.

The distance $\delta$ between the center of mass and the body equilibrium position is given by the Newton equation

$$m\ddot{\delta} + \rho g|\Gamma^w|\delta = \int_{\Gamma^w} (\Pi - 2\nu\partial_z u_z)$$

where the fluid velocity $u$ and the hydrodynamical pressure $\Pi = P - P_{atm} + \rho g z$ are computed through the incompressible free-surface Navier-Stokes equations on $\Gamma^a$, homogeneous Dirichlet on $\Gamma^w \cup \Gamma_b$ and $u_{|\Gamma^w} = \dot{\delta} e_z$ at the wetted line $\Gamma^w$. The equation on the vertical displacement $\delta$ is in fact a combination of Newton's equations, accounting for both the gravitational and the stress exerted on the solid by the fluid (including the viscosity), and the Archimedes equilibrium equation.

## 2 Added mass effect

When a partially submerged structure accelerates in a fluid, it does not move alone: it sets a portion of the surrounding fluid into motion too. The added mass corresponds precisely to the inertial component of this supplementary fluid motion. This inertial contribution is contained in the pressure forces acting on the body's wetted interface. Let us first consider an ideal, incompressible, and irrotational fluid. In this setting, the velocity field is the gradient of a scalar harmonic potential $\phi$. This potential depends linearly on the velocity of the solid. More precisely $\phi(t,x) = \dot{\delta}(t)\phi(t,x) + \phi_{fluid}(t,x)$ where $\phi$ is the harmonic potential which satisfies at the boundaries:

$$(\partial_n \phi)_{|\Gamma^w} = 1, \quad (\partial_n \phi)_{|\Gamma_b \cup \Gamma_{wall}} = 0, \quad \phi_{|\Gamma^a} = 0.$$

Note that since the fluid domain moves, this harmonic potential $\phi$ consequently depends on time.

Within the framework of irrotational fluids, the pressure can be computed from the potential through Bernoulli's equation, relating the evolution of $\phi$ to the velocity's squared Euclidean norm and the pressure. When the latter is integrated over the surface of the body to evaluate the hydrodynamic force, the term involving the time derivative of the potential produces a contribution proportional to the acceleration of the solid. We can interpret this inertial contribution as an added mass. Within the Navier-Stokes framework, the velocity field can be decomposed into the gradient of a potential responsible for the added mass effects, plus a non-inertial remainder that can no longer be expressed as a gradient due to viscous effects. Then the hydrodynamical pressure is decomposed into an added-mass term and a non-inertial hydrodynamical pressure $\Pi = -\rho\ddot{\delta}\varphi + \widetilde{\Pi}$. Thus, the Newton equation becomes

$$(m + m_a)\ddot{\delta} + \rho g|\Gamma^{\mathrm{w}}|\delta = \widetilde{f}_{\mathrm{hyd}} := \int_{\Gamma^{\mathrm{w}}} (\widetilde{\Pi} - \nu\partial_z u_z)$$

where the added mass is given by

$$m_a := \int_{\Omega} \rho|\nabla\varphi|^2.$$

The evolution of the fluid's quantities $(u, \widetilde{\Pi})$ is achieved through the Navier-Stokes equations:

$$\begin{cases} \rho D_t u + \nabla\widetilde{\Pi} - \rho\nu\Delta u = a_\delta \nabla\varphi & \text{in } \Omega(t), \\ \nabla \cdot u = 0 & \text{in } \Omega(t), \\ \widetilde{\Pi} n + 2\nu D[u]n = a_\delta \varphi n + \rho g z n & \text{on } \Gamma^{\mathrm{a}}(t), \\ u = 0 & \text{on } \Gamma_{\mathrm{b}} \cup \Gamma_{\mathrm{wall}}(t), \\ u = \dot{\delta} e_z & \text{on } \Gamma^{\mathrm{w}}(t), \end{cases}$$

where

$$a_\delta := \rho \frac{\widetilde{f}_{\mathrm{hyd}} - \rho g|\Gamma^{\mathrm{w}}|\delta}{m + m_a}, \qquad D_t := \partial_t + (u\nabla).$$

We denoted as $D[u]$ is the symmetric part of the gradient tensor and $n$ the unit vector perpendicular to the free surface.

## 3 Numerical method

Numerical approximation of the solution to this problem is achieved extending the scheme introduced in [3] to handle the movements of the solid body. At each time iteration, the harmonic potential $\phi$, the velocity $u$ and the pressure $p$ are approximated using the Finite Element Method, implemented using the FreeFEM library [1]. Taylor-Hood $(\mathbb{P}^2, \mathbb{P}^1)$ elements are used for the Navier-Stokes solver. Advection of the free-surface $\Gamma^{\mathrm{a}}(t)$ and of the wetted boundary $\Gamma^{\mathrm{w}}(t) \cup \Gamma_{\mathrm{wall}}(t)$ is achieved transporting the whole triangulation with the fluid (the ALE method of [2]). This ensures an accurate depiction of the movement of both the free surface and the floating solid body. The fluid's pressure and velocity, needed for Newton's equation, are first computed during a predictor step. The solid's velocity $\dot{\delta}$ is then computed from these first estimates. From this, the fluid's pressure and velocity are recomputed (corrector step) and then used to advect the mesh.

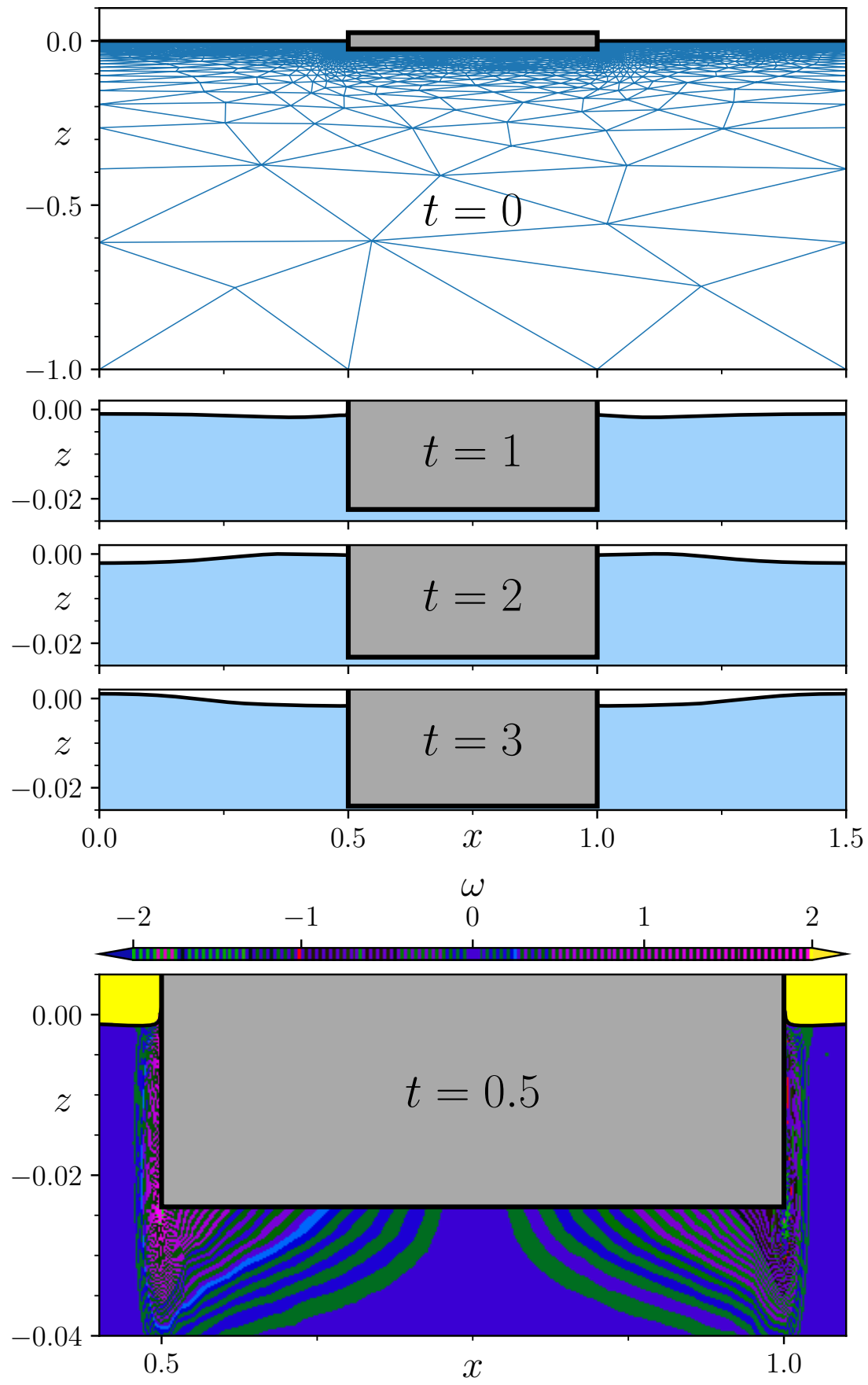

Figure 2: Simulation of an immersed solid moving upward due to buoyancy: initial mesh (top), free surface evolution (middle) and vorticity $\omega$ (bottom)

### Data accessibility and reproducibility

The FreeFEM scripts along with the post-processing python notebook can be found here.

### Acknowledgement

GB is supported by the BOURGEONS project, grant ANR-23-CE40-0014-01.

# Scalable WaveHoltz in MFEM

**Amit Rotem**[1,*], **Daniel Appelö**[1], **Will Pazner**[2]
[1]Department of Mathematics, Virginia Tech., Blacksburg, USA
[2]Department of Mathematics and Statistics, Portland State University, Portland, USA
*Email: arotem@vt.edu

**Abstract**

We present a scalable WaveHoltz solver for the Helmholtz equation in MFEM, leveraging optimized finite element infrastructure to achieve performance portability on modern GPU architectures. Our formulation combines off-the-shelf MFEM components with a novel modified midpoint time-stepping scheme accelerated by algebraic multigrid resulting in a robust and performant method providing a practical pathway to scalable Helmholtz solvers. We compare against alternative off-the-shelf solvers and analyze key design choices, highlighting tradeoffs impacting robustness, efficiency, and implementation simplicity.

## 1 Introduction

The Helmholtz equation is challenging for iterative solvers due to its indefiniteness: general-purpose techniques such as multigrid and domain decomposition, effective for definite problems, must be non-trivially specialized for the Helmholtz equation [4]. We consider the WaveHoltz iteration, which filters the solution of the wave equation over a short time interval [2]. We develop a scalable implementation within MFEM and study its behavior for high-order continuous (CG) and discontinuous (DG) Galerkin discretizations. High-order methods mitigate pollution error and reduce degrees of freedom at high frequencies [5]. We identify design choices that balance robustness and efficiency and find that CG with implicit time stepping is most efficient.

## 2 Model

Consider the Helmholtz equation:

$$-\nabla \cdot (c^2 \nabla u) - \omega^2 \rho^2 u = f, \quad x \in \Omega, \tag{1a}$$

$$\alpha \mathbf{n} \cdot (c^2 \nabla u) - i\beta\omega\rho u = g, \quad x \in \partial\Omega. \tag{1b}$$

Discretizing (1) with CG or DG yields:

$$(A - i\omega H - \omega^2 M)\mathbf{u} = \mathbf{f}. \tag{2}$$

Here $A$ discretizes the negative Laplacian, $M$ is the mass matrix, and $H$ is the boundary mass matrix from the impedance boundary conditions.

## 3 Methods

### 3.1 WaveHoltz Iteration

The WaveHoltz fixed-point iteration is:

$$\mathbf{u}^{(k+1)} = \Pi \mathbf{u}^{(k)} = \frac{2}{N} \sum_{n=0}^{N-1} \kappa_n (\mathbf{p}^n - \frac{1}{i\omega}\mathbf{q}^n). \tag{3}$$

Here $\mathbf{p}^n, \mathbf{q}^n$ are time discrete solutions of the real valued wave equation:

$$\dot{\mathbf{p}} = \mathbf{q}, \tag{4a}$$

$$M\dot{\mathbf{q}} + H\mathbf{q} + A\mathbf{p} = \Re\{\mathbf{f}e^{-i\omega t}\}, \tag{4b}$$

$$\mathbf{p}(0) = \Re\{\mathbf{u}^{(k)}\}, \quad \mathbf{q}(0) = \omega\Im\{\mathbf{u}^{(k)}\}. \tag{4c}$$

The time-filter $\kappa_n = \cos\frac{2\pi n}{N} - \frac{1}{4} + \frac{1}{4}\tan^2\frac{\pi}{N}$ over $N$ time steps is designed to remove time-stepping errors and ensure convergence to the solution of (2).

### 3.2 Time Stepping

In each WaveHoltz iteration we evolve (4) over $P$ periods of $e^{-i\omega t}$ (larger $P$ improves convergence). We present a novel modified implicit midpoint scheme which exactly evolves $e^{-i\omega t}$. Define $\theta = \frac{\tan(\pi/N)}{\omega}$ and $\sigma = \frac{\sin(\pi/N)}{\omega/2}$, then

$$\mathbf{p}^{n+1} - \theta\mathbf{q}^{n+1} = \mathbf{p}^n + \theta\mathbf{q}^n, \tag{5a}$$

$$\begin{aligned} \theta A\mathbf{p}^{n+1} &+ (M + \theta H)\mathbf{q}^{n+1} \\ &= -\theta A\mathbf{p}^n + (M - \theta H)\mathbf{q}^n \\ &+ \sigma\Re\left\{\mathbf{f}\, e^{-2\pi i(n+\frac{1}{2})/N}\right\}. \end{aligned} \tag{5b}$$

Since the time stepping scheme is A-stable and introduces no time-stepping errors to the converged solution, the number of time steps can be kept small and independent of mesh resolution.

### 3.3 Scalable High Order Finite Element Methods

To take advantage of modern parallel computing environments, and especially many GPU environments, finite element methods are implemented as partially assembled (PA) operators in MFEM [1] without the use of sparse matrices. Gradients and integrals are computed on the fly, reducing memory footprint, and exposing finer grain parallelism, particularly on GPUs.

The block system (5) reduces, by Schur complement, to solving $(\theta^2 A + \theta H + M)\mathbf{x} = \mathbf{y}$ which is positive definite and can be solved by algebraic multigrid (AMG). For high order methods, AMG is more expensive than the same system for a low order method, so we use the low order reduction method of [3] coupled with AMG to precondition the linear system. In practice, AMG is robust to problem size, making the WaveHoltz iteration scalable.

## 4 Results

For this abstract we limit the presentation to one compelling result. We solve (1) in the unit cube with constant coefficients, $\omega = 10$, fourth order elements, and increasing mesh resolution. We prescribe the Dirichlet condition $u = e^{i(x+2y+3z)}$ on the sides and impedance conditions on the top and bottom. The WaveHoltz solver uses only five time steps per period for three periods and accelerated with GMRES(20). For reference, we compare against MINRES. For both solvers we stop at relative residual $10^{-6}$. All computations performed on a NVIDIA A100 GPU.

In Figure 1 we observe that the WaveHoltz iteration count is effectively constant whereas the MINRES iteration count grows with the problem size. Although each WaveHoltz iteration is more expensive, both methods have asymptotically linear per-iteration cost, making WaveHoltz overall linear while MINRES is superlinear. For the largest CG problem, WaveHoltz was more than twenty times faster than MINRES; the largest DG problem was solved in just eleven minutes with WaveHoltz, but did not converge in twenty hours with MINRES.

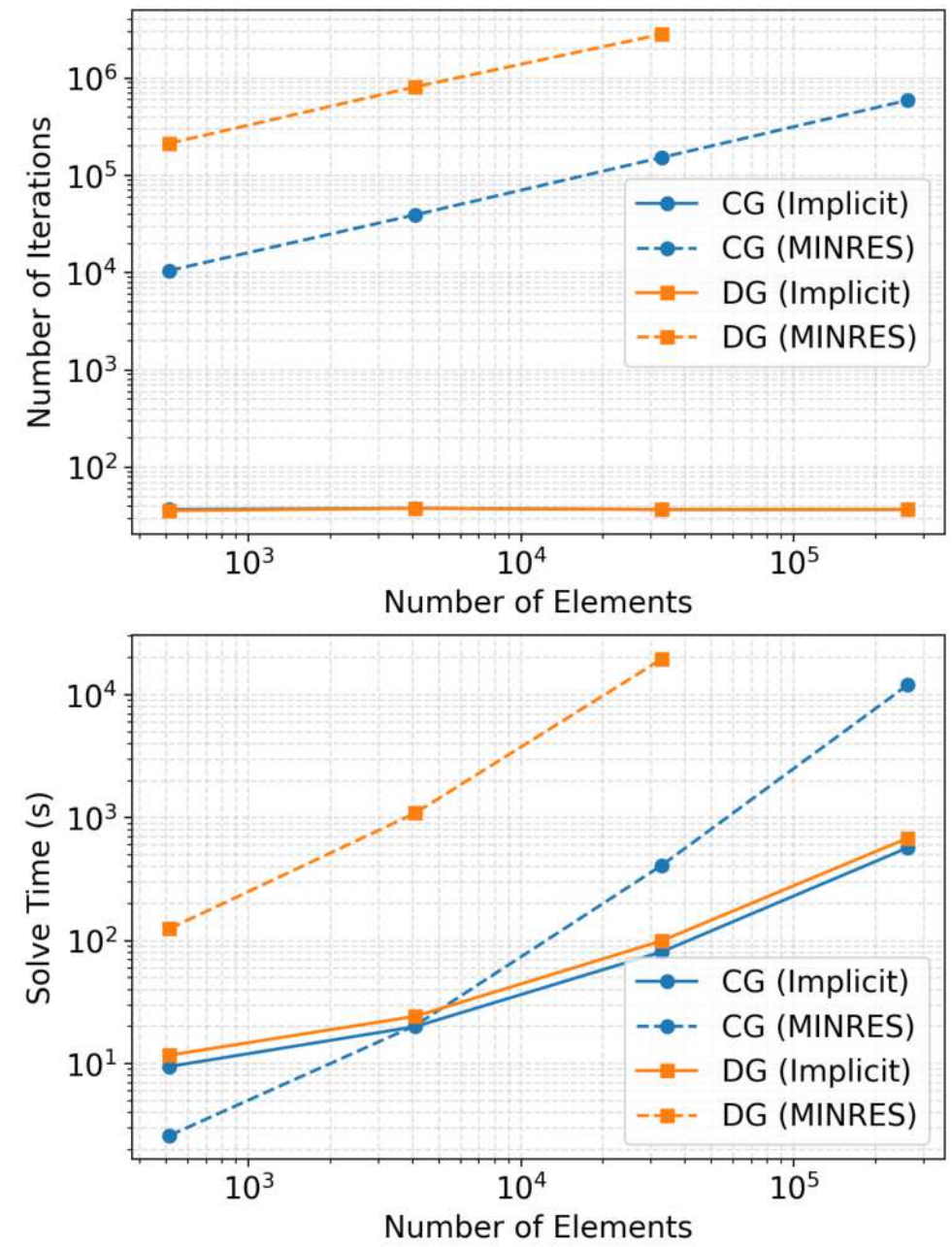

Figure 1: In the top figure, the number of iterations vs. number of elements, and in the bottom figure, the solve time vs. the number of elements

# The cut-off resolvent can grow arbitrarily fast in obstacle scattering

Siavash Sadeghi[1,*], Simon N. Chandler-Wilde[1]
[1]Department of Mathematics and Statistics, University of Reading, UK
*Email: s.sadeghi@pgr.reading.ac.uk

## Abstract

We consider time-harmonic acoustic scattering by a compact sound-soft obstacle $\Gamma \subset \mathbb{R}^n$ ($n \geq 2$) that has connected complement $\Omega := \mathbb{R}^n \setminus \Gamma$. This scattering problem is modelled by the inhomogeneous Helmholtz equation $\Delta u + k^2 u = -f$ in $\Omega$, the boundary condition that $u = 0$ on $\partial\Omega$, and the standard Sommerfeld radiation condition. It is well-known that, if the boundary $\partial\Omega$ is smooth, then the norm of the cut-off resolvent of the Laplacian, that maps the compactly supported inhomogeneous term $f$ to the solution $u$ restricted to some ball, grows at worst exponentially with $k$. In this paper we show that, if no smoothness of $\Gamma$ is imposed, then the growth can be arbitrarily fast.

## 1 Introduction

It follows from the well-posedness of the above scattering problem that, given any $R > R_\Gamma := \max_{x\in\Gamma} |x|$ and any $k > 0$, there exists some minimal $C_{k,R} > 0$ such that

$$\|u\|_{L^2(\Omega_R)} \leq C_{k,R}\|f\|_{L^2(\Omega_R)}, \tag{1}$$

whenever $\mathrm{supp}(f) \subset \Omega_R := \{x \in \Omega : |x| < R\}$. This paper is concerned with answering the following question: For fixed $R > R_\Gamma$, how fast can $C_{k,R}$ grow as $k \to \infty$?

For each $k > 0$, we have a well-defined resolvent

$$R(k) := (-\Delta - k^2)^{-1} : L^2_{\mathrm{comp}}(\Omega) \to H_0^{1,\mathrm{loc}}(\Delta, \Omega), \tag{2}$$

that maps $f \in L^2_{\mathrm{comp}}(\Omega)$ continuously to the solution $u \in H_0^{1,\mathrm{loc}}(\Delta, \Omega)$ of the above scattering problem. For $R > R_\Gamma$ let $\chi_R$ denote the characteristic function of $\Omega_R$. Then, for $k > 0$ and $R > R_\Gamma$, the minimal $C_{k,R}$ in (1) is precisely

$$C_{k,R} := \|\chi_R R(k)\chi_R\|_{L^2(\Omega)\to L^2(\Omega)}. \tag{3}$$

From [1, Thm. 2], when $\Gamma$ is the closure of a $C^\infty$ domain with connected complement $\Omega$, the cut-off resolvent can grow at most exponentially fast, i.e., for every $R > R_\Gamma$, there exists $c, C > 0$ such that

$$C_{k,R} \leq Ce^{ck}, \quad k \geq k_0. \tag{4}$$

Furthermore, this bound is sharp in the sense that, if $\Gamma$ has an ellipse shaped cavity, then there exists an unbounded sequence $0 < k_1 < k_2 < ...$ such that the cut-off resolvent achieves exponential growth along this sequence.

## 2 Main result

We show that, in the absence of smoothness of the obstacle $\Gamma$, we can achieve arbitrary growth of the cut-off resolvent through a moderately increasing sequence of wavenumbers

**Theorem 1** [2, Thm. 1.4] *Suppose the sequences $(k_j)_{j\in\mathbb{N}}$ and $(a_j)_{j\in\mathbb{N}}$ are positive, increasing and unbounded, and that, for some $c > 0$,*

$$k_j \geq c(j \log(j + \mathrm{e}))^{1/n} \log^2(\log(j + \mathrm{e}^{\mathrm{e}})), \quad j \in \mathbb{N}. \tag{5}$$

*Then there exists a compact set $\Gamma \subset \mathbb{R}^n$, with $\Omega := \mathbb{R}^n \setminus \Gamma$ connected, such that, if $R > R_\Gamma$, then, for each $j \in \mathbb{N}$,*

$$C_{k_j,R} \geq a_j. \tag{6}$$

## 3 Methods and Significance

The construction of the compact set $\Gamma$ is explicit, see figure 1 below and [2, Thms. 1.5,1.9]. The bounded, open set $\mathcal{O}_\infty$ is the union of open cubes $\mathcal{O}_j$ of sidelength $\ell_j = \pi\sqrt{n}k_j^{-1}$. The set $\Gamma$ is the closure of the union of "slit boxes" $\partial\mathcal{O}_j \setminus S_j$, where $S_j \subset \partial\mathcal{O}_j$ is an $(n-1)$-dimensional cube of sidelength proportional to $(k_j a_j)^{-\frac{2}{3(n-1)}}$. We show that $\mathcal{O}_\infty$ is bounded and $\Omega = \mathbb{R}^n \setminus \Gamma$ is connected. By construction, we have that $k_j^2$ is the first eigenvalue of the Dirichlet Laplacian on $\mathcal{O}_j$ and hence a Dirichlet eigenvalue of $\mathcal{O}_\infty$. The associated eigenfunction $u_j$, extended to $\mathbb{R}^n \setminus \mathcal{O}_j$ by zero, is a sufficiently good quasimode to imply that (6) holds.

Extra care is taken in how we arrange the slit boxes in order to ensure that the sequence

$(k_j^2)_{j\in\mathbb{N}}$ need not grow much faster than a sequence of Dirichlet eigenvalues on a bounded domain (the Weyl law has $k_j$ proportional to $j^{\frac{1}{n}}$) while at the same time ensuring that $\mathcal{O}_\infty$ is bounded. The same cut-off resolvent growth can be achieved with a much simpler configuration of boxes which are placed in decreasing order of size along some line, but the sequence $k_j$ will then need to grow much faster.

It was shown in [4], for $\Gamma$ which is the closure of a bounded Lipschitz domain with $\Omega = \mathbb{R}^n \setminus \Gamma$ connected, that if one excludes a set of wavenumbers of arbitrarily small Lebesgue measure from $[k_0, \infty)$ one can ensure at most polynomial growth of the cut-off resolvent outside this excluded set. Of course the Lipschitz condition does not hold for the set $\Gamma$ in Theorem 1. Theorem 2 below is a refinement of the main result (Theorem 1.1) of [4]; our refinement is to prove the bound (7) for arbitrary compact $\Gamma$ and to sharpen the exponent (the proof in [4] establishes (7) only for Lipschitz $\Gamma$, and with our exponent $2n+\delta$ replaced by the larger exponent $5n/2+\delta$).

**Theorem 2** [3] *Let $\Gamma \subset \mathbb{R}^n$ be a compact set such that $\Omega = \mathbb{R}^n\setminus\Gamma$ is connected. Then, for any $R > R_\Gamma$, $\epsilon > 0, k_0 > 0$ and $\delta > 0$, there exists $J \subset [k_0, \infty)$, with Lebesgue measure $|J| \le \epsilon$, and a constant $C = C(\epsilon, \delta, k_0, R, n)$ such that*

$$C_{k,R} \le Ck^{2n+\delta}, \quad k \in [k_0, \infty) \setminus J. \tag{7}$$

We note the large contrast between Theorems 1 and 2; for any modestly increasing sequence $(k_j)_{j\in\mathbb{N}}$ there exists a compact set $\Gamma$ (with connected complement), such that $C_{k_j,R}$ grows arbitrarily fast, yet $C_{k,R}$ grows at worst at a polynomial rate outside some set of wavenumbers of arbitrarily small measure.

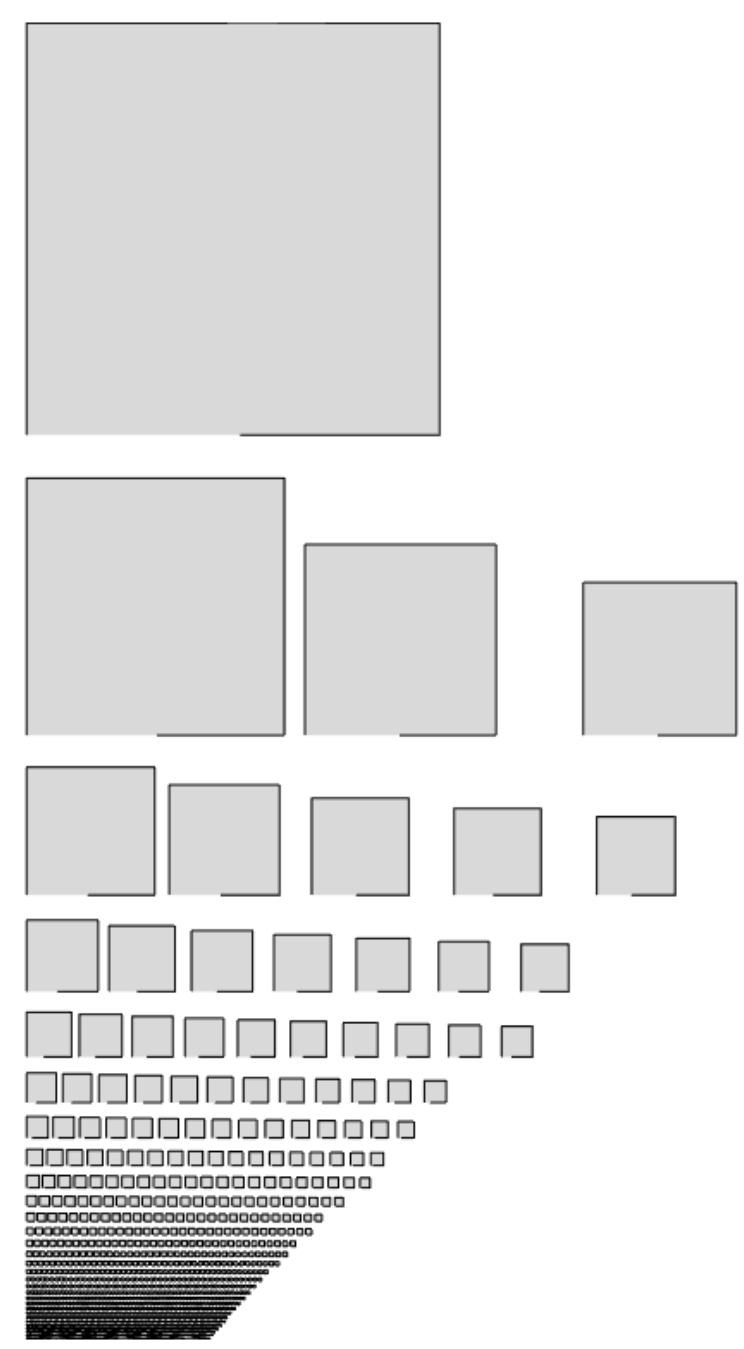

Figure 1: Plot when $n = 2$. $\Gamma$ is the closure of the union of $\Gamma_j$, where $\Gamma_j \subset \partial\mathcal{O}_j$ consists of the solid black lines on the boundary of $\mathcal{O}_j$ which is a box of sidelength $\pi\sqrt{n}k_j^{-1}$ (shaded in gray).

# Addressing Seismic Acquisition Repeatability in $CO_2$ Storage Monitoring Using Machine Learning

Manon Sarrouilhe[1,*], Mauricio Araya-Polo[2], Hélène Barucq[1], Henri Calandra[3], Stefano Frambati[3]

[1]Makutu, Inria-TotalEnergies, Université de Pau et des Pays de l'Adour, France
[2]TotalEnergies, Houston, USA
[3]TotalEnergies, Pau, France

[*]Email: manon.sarrouilhe@inria.fr

## Abstract

$CO_2$ geological storage requires reliable seismic monitoring to ensure safety and effectiveness. Repeatable seismic acquisitions are challenging due to source variability and reservoir changes. We address this by generating seismograms for corrected source positions under unknown reservoir velocities. We propose an attention-based architecture that we compare with symmetric convolutional autoencoders formerly developed by Bharadwaj, Li, and Demanet for tomography.

## Introduction

$CO_2$ geological storage is critical for reaching long-term carbon neutrality and seismic monitoring will play a crucial role to ensure that storage is safe, effective and permanent. Traditional technologies for the monitoring of $CO_2$ injection rely on geophysical techniques such as seismic surveying, gravimetry and electromagnetic methods. However, achieving repeatable seismic acquisitions remains challenging due to variations in source locations between surveys and possible changes in reservoir velocities caused by natural phenomena, including seismic events. In this context, numerical methods can serve as effective predictive tools. Nevertheless, their applicability is often limited by the difficulty of balancing high accuracy with affordable computational costs. Recent advances in machine learning have shown strong potential to address these challenges by complementing numerical-analysis-based methods. In this work, we present the challenge of acquisition repeatability by creating seismograms corresponding to corrected source positions and an unknown reservoir velocity. The proposed approach builds on the method introduced by Bharadwaj, Li and Demanet [1] for tomography problems. We extend their work which employs symmetric convolutional autoencoders (SymAE) by proposing a new architecture based on attention mechanisms for seismic surveying. We illustrate our study using a 2D model inspired by SEAMCO2, enabling a direct comparison between both architectures.

## 1 General setting

Seismic acquisition consists in measuring signals at sensors recording reflected waves generated by sources. Signals are solution to the acoustic wave equation:

$$\begin{cases} \dfrac{1}{c^2}\dfrac{\partial^2 p}{\partial t^2}(x,t) - \Delta p(x,t) = f(x,t) \text{ in } \Omega_T \\ p(x,t) = 0 \text{ on } \Gamma_{0,T} \\ \dfrac{1}{c}\dfrac{\partial p}{\partial t}(x,t) + \dfrac{\partial p}{\partial n}(x,t) = 0 \text{ on } \Gamma_{1,T} \\ p(x,0) = 0, \dfrac{\partial p}{\partial t}(x,0) = 0 \text{ in } \Omega \end{cases} \tag{1}$$

The propagation domain is denoted by $\Omega_T = \Omega \times (0,T)$. The space boundary $\partial\Omega = \Gamma_0 \cup \Gamma_1$ is composed of a free surface $\Gamma_0$ on the top and on $\Gamma_1$ we impose an absorbing boundary condition. It contains a reservoir characterized by the propagation speed $c$. The pressure is computed for different positions of sources denoted by $f$ and velocities in the reservoir using a Spectral Element Method (SEM) in space and a Leap-Frog scheme in time. The receivers are static from one model to another while sources can move randomly. Figure 1 represents a simplified SEAMCO2 model.

Monitoring $CO_2$ storage relies on seismic acquisition. This process is repeated over many years, and during this time, the reservoir may change shape or properties due to natural interference. In addition, the position of the sources used during the first acquisition may have changed. Due to these uncertainties, a naive comparison between two acquisitions cannot be made. Our goal is to apply neural network architectures to learn seismic acquisitions in the form of images in order to understand how the interior of the reservoir has evolved over time. It is worth noting that the velocity in the reservoir is unknown.

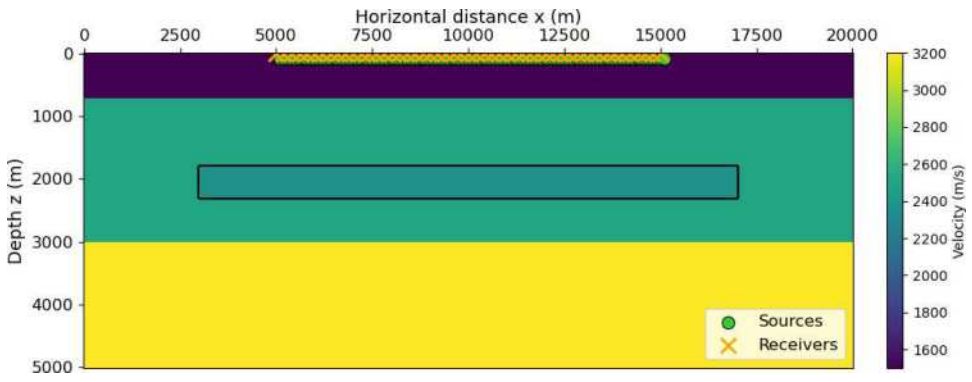

Figure 1: Simplified SEAMCO2

## 2 Data preparation and Method

Initial data is a set of $N$ velocity models denoted by $M_i$ where $i = 1, \cdots, N$. They represent $i$ successive seismic campaigns on an evolving reservoir model, $M_0$ representing the initial state (i.e., no reservoir) and $M_i$ representing the current state of the reservoir, with an unknown fluid distribution. For each model $M_i$, we associate an acquisition geometry composed of a set of $N_i$ sources $\mathcal{S}_i = (s^i_1, \cdots, s^i_{N_i})$ and $P$ receivers. We denote $\sigma_{i,j}$ the seismic campaign associated with model $M_i$ and sources $\mathcal{S}_j$. We generate each $\sigma_{i,i}$ with a SEM code. SymAE is trained in autoencoder mode by randomly sampling a velocity model $M_i$ and drawing $n_\tau$ seismograms from $\mathcal{S}_i$. The sampled seismograms share the same latent coherent content (the velocity model) but differ through nuisance factors related to source location. The $n_\tau$ seismograms are aggregated and treated as a single training instance. As in [1], the autoencoder architecture is decomposed into two branches to promote disentanglement between coherent and nuisance components. One branch uses a symmetrization/sum on the $n_\tau$ instances, reducing the dependency on the shots and the other masks randomly some neurons to eliminate the coherent information. Unlike [1], we replace convolutional operators with stacked attention modules [2]. First, each trace is embedded and time-projected in a smaller space permitting to work on a representation of the trace. Then, we add a positional encoding for helping the neural network to learn the position/index of receivers with one another. We further employ a combination of L2 and L1 loss to address large amplitude variability caused by source–receiver geometry and reflection strength. This enables training on unnormalized data after removing the direct arrival while preserving subtle reflection details. Once the network has been trained, it can be used in inference mode to create an artificial set of seismograms $\sigma_{1,0}$ which predicts the seismograms that would have been produced by simulating the sources $\mathcal{S}_0$ with the (unknown) velocity model $M_1$. Note that the network does not need to know the location of sources, receivers or the specific models $M_0$, $M_j$, and thus can be applied when source locations are unknown (natural sources), and/or when the model itself is unknown.

## 3 Results

Table 1 compares [1] approach and ours.

| Threshold | $< 5\%$ | $< 2\%$ | $< 1\%$ |
|---|---|---|---|
| Using [1] | 100% | 99.21% | 95.17% |
| Using transformer | 100% | 100% | 99.98% |

Table 1: Accuracy for two shots without direct arrival for different threshold and approaches.

Trace correlation learning has been eased by applying the method introduced in [2] so that the neural network learns long-term dependencies between receivers thanks to the embedding and positional encoding. Both methods were trained for the same wall-clock time (1h30 on an NVIDIA H200 GPU). One training iteration takes 9 s with our approach, compared to 40 s for [1], allowing for substantially more optimization steps within the same computational budget. The number of parameters is comparable (23M for ours vs 32M for theirs). While their method uses slightly more parameters and achieves slightly lower accuracy, both approaches operate at a similar computational scale. The next step is to incorporate an inversion procedure to recover the velocity field.

# Efficient evaluation of the two-layered acoustic Green's function (extended abstract)

**Stefan.A. Sauter**[1,*], **Laura Depedrini**[1], **Benedikt Grässle**[1], **J.Markus. Melenk**[2]
[1]Institut für Mathematik, Universität Zürich, Winterthurerstrasse 190, CH-8057 Zürich, Switzerland
[2]Institut für Analysis und Scientific Computing, Technische Universität Wien, Wiedner Hauptstrasse 8–10, A-1040 Wien, Austria
*Email: (stas@math.uzh.ch)

**Abstract**

In this paper we consider the numerical solution of high-frequency wave propagation from a scatterer within a two-layered ambient medium by boundary element methods. The Green's function can be represented as a Sommerfeld-type integral but is not available explicitly.

We introduce parameter-dependent integral transformations and integration contours such that the integrand becomes non-oscillatory and suitable for fast numerical quadrature. This paper is the abridged version of [1].

## 1 Introduction

The two-layered medium is described by the disjoint decomposition of the full space $\mathbb{R}^d$ into $H_+ \cup H_- \cup H_0$, $d \in \mathbb{N} = \{1, 2, \ldots\}$, for the upper/lower half-space

$$H_\pm := \left\{ \mathbf{x} = (x_j)_{j=1}^d \in \mathbb{R}^d \mid \pm x_d > 0 \right\}$$

and their common interface $H_0 := \overline{H_+} \cap \overline{H_-}$ with normal vector $\mathbf{n} = (0, \ldots, 0, -1)^T$. For $\mathbf{z} = (z_j)_{j=1}^d \in \mathbb{R}^d$, let $\mathbf{z}' = (z_j)_{j=1}^{d-1}$ so that $\mathbf{z} = (\mathbf{z}', z_d)$.

The (reciprocal) material density coefficient is assumed to be a piecewise constant function $a : \mathbb{R}^d \backslash H_0 \to \mathbb{R}_{>0}$ with $a(\mathbf{x}) = a_\pm$ if $\pm x_d > 0$ for given constants $a_+, a_- > 0$. We assume for the wavenumber $s = -\mathrm{i}\, k$, $k \in \mathbb{R} \backslash \{0\}$, and define the effective wavenumber for each half-space $H_\pm$ by $s_\pm := s/\sqrt{a_\pm}$.

Let $\Omega \subset\subset H_+$ be a bounded Lipschitz domain with boundary $\Gamma$. Its exterior complement is denoted by $\Omega^c := \mathbb{R}^d \backslash \overline{\Omega}$. We consider the exterior Dirichlet problem: For given $\psi_D \in H^{1/2}(\Gamma)$, find $u \in H^1_{\mathrm{loc}}(\Omega^c)$ such that

$$\begin{aligned} &-\operatorname{div}(a \nabla u) + s^2 u = 0, \text{ in } (\Omega^c \cap H_+) \cup H_-, \\ &u\big|_\Gamma = \psi_D, \text{ (Dirichlet boundary condition)}, \\ &[u]_{H_0} = [a \partial_{\mathbf{n}} u]_{H_0} = 0, \text{ (interface condition)} \end{aligned}$$

with appropriate radiation conditions imposed at infinity. Green's representation formula for this problem yields (see, e.g., [3, §3.4.2.2])

$$u(\mathbf{x}) = -\int_\Gamma G(\mathbf{x}, \mathbf{y}) \frac{\partial u(\mathbf{y})}{\partial \mathbf{n}_+} \, ds_{\mathbf{y}} + \int_\Gamma \frac{\partial G(\mathbf{x}, \mathbf{y})}{\partial \mathbf{n}_+(\mathbf{y})} u(\mathbf{y}) \, ds_{\mathbf{y}}, \qquad \forall \mathbf{x} \in \Omega^c.$$

Here the normal vector $\mathbf{n}$ at $\Gamma$ points into $\Omega^c$ and $\mathbf{n}_+ = a_+ \mathbf{n}$. Applying the trace operator for the boundary $\Gamma$ to this equation we obtain a relation between the Dirichlet trace of $u$ and the co-normal derivative of $u$ on $\Gamma$:

$$\int_\Gamma G(\mathbf{x}, \mathbf{y}) \frac{\partial u(\mathbf{y})}{\partial \mathbf{n}_+} \, ds_{\mathbf{y}} = -\frac{1}{2} \psi_D(\mathbf{x}) + \int_\Gamma \frac{\partial G(\mathbf{x}, \mathbf{y})}{\partial \mathbf{n}_+(\mathbf{y})} \psi_D(\mathbf{y}) \, ds_{\mathbf{y}}, \quad \forall \mathbf{x} \in \Gamma \tag{1}$$

Since $\Omega \subset H_+$, both variables $\mathbf{x}, \mathbf{y}$ belong to $H_+$ and the Green's function in (1) is given by

$$G(\mathbf{x}, \mathbf{y}) = g^{\mathrm{full}}_{++}(\|\mathbf{x} - \mathbf{y}\|) + g^{\mathrm{imp}}_{++}(\mathbf{x}, \mathbf{y}),$$

$$\text{where} \quad g^{\mathrm{imp}}_{++}(\mathbf{x}, \mathbf{y}) := \frac{1}{2\pi} \frac{1}{(2\pi\omega)^\nu} \times \int_0^\infty \frac{a_+ \hat{\beta}_+ - a_- \hat{\beta}_-}{a_+ \hat{\beta}_+ + a_- \hat{\beta}_-} \frac{e^{-\hat{\beta}_+ z_d}}{2 a_+ \hat{\beta}_+} J_\nu(\omega\rho) \, \rho^{\nu+1} \, d\rho \tag{2}$$

and $g^{\mathrm{full}}_{++}$ is the full space Green's function for constant wavenumber $s_+$ with $\omega := \|\mathbf{x}' - \mathbf{y}'\|$, $z_d := x_d + y_d$, and $\hat{\beta}_\pm(\rho) := \sqrt{\rho^2 + s_\pm^2}$. Since $g^{\mathrm{full}}_{++}$ is given explicitly, it is key for an efficient numerical discretization of (1), e.g., by Galerkin boundary elements, that a fast and stable method is available to approximate the integral in the second summand of (2) and this is the topic of this paper.

## 2 Parameter-dependent integral transforms

The integrand in (2) is oscillatory and, in addition, contains two singular points $\rho_\pm = k/\sqrt{s_\pm}$ where we split the integral so that $g^{\mathrm{imp}}_{++}(\mathbf{x}, \mathbf{y}) := \int_0^{\rho_+} \ldots + \int_{\rho_+}^{\rho_-} \ldots \int_{\rho_-}^\infty \ldots$. (To fix ideas we assume here $\rho_+ \leq \rho_-$). The essential step is to find integration contours in the analyticity regions of the parameter-dependent integrand such that the oscillatory part of the integrand becomes constant along the contour and the real part decays exponentially. To understand the oscillatory behaviour we write the Bessel function in the form

$$J_\nu(z) = \chi_+(z)\, e^{-\mathrm{i}z} + \chi_-(z)\, e^{\mathrm{i}z} \qquad (3)$$

for slowly varying functions $\chi_\pm$ which are decaying algebraically as $|z| \to \infty$ and are holomorphic in $\mathbb{C}\backslash\,]-\infty, 0]$.

Hence, the oscillatory and decay behaviour of the integrands along a contour is determined by the behaviour of the terms (in scaled coordinates)

$$\begin{aligned} &\mathrm{e}^{-|s_+|p_1^\pm(x)} \text{ for } p_1^\pm(x) := -\mathrm{i}xz_d \pm \mathrm{i}\omega\sqrt{1-x^2},\ x \in (0,1), \\ &\mathrm{e}^{-|s_+|p_2^\pm(x)} \text{ for } p_2^\pm(x) := xz_d \pm \mathrm{i}\omega\sqrt{1+x^2},\ x > 0. \end{aligned} \qquad (4)$$

Here we consider exemplarily the integral related to $p_1^+$ for the interval $(0,1)$, i.e., $a_+ > a_-$. By a simple scaling transform we obtain $\rho_+ = 1$. Then,

$$\int_0^1 \cdots = \int_{\Gamma_0} \cdots - \int_{\Gamma_1} \cdots$$

where

$$\Gamma_{u_0} = \left\{ \begin{pmatrix} u \\ |w_1|\sqrt{\dfrac{u^2-1+w_1^2}{u^2+w_1^2}} \end{pmatrix} : u > u_0 \right\} \text{ and}$$

$$w_1 = -\frac{z_d}{\omega}(u-u_0) + \sqrt{1-u_0^2}. \qquad (5)$$

**Theorem 1** *The transformed integrand along the contours $\Gamma_0$ and $\Gamma_1$ has the property that $\operatorname{Im} p_1^- = \text{const}$ and $\operatorname{Re} p_1$ is positive and linearly growing for increasing argument.*

## 3 Numerical integration

The integrands along the transformed contours are non-oscillatory but contain an algebraic endpoint singularity. Since our contour integral representations render the integrals exponentially decaying towards infinity we applied scaled version of Gauss-Laguerre quadrature to these integrals. The problem parameters are chosen as $a_+ = 2$, $a_- = 1$, $z_d = \omega = 1$, the integral is over the interval $(0,1)$ and related to $p_1^-$. We consider the contour which starts from the real axes at $u = 1$.

The well-known ( [2]) root-exponential convergence behaviour for smooth functions $f$ in $\int_0^\infty \mathrm{e}^{-x}\, x^\alpha f(x)\, dx$ is clearly visible for the considered problem.

In summary we have presented an efficient method for approximate point evaluations of the acoustic Green's function in a two-layered medium.

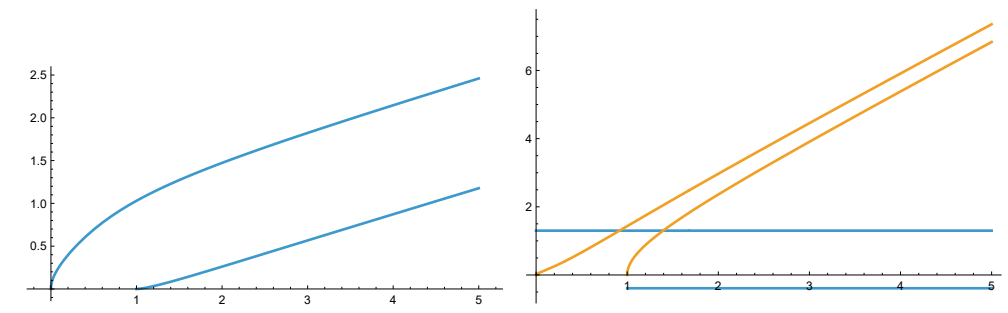

Table 1: Left: contour consisting of two infinite pieces, starting from the real axis at 0 and coming back from $+\infty$ on the real axis to the interval endpoint 1. Right: Imaginary parts of $p_1^-$ along the contour (blue) and real parts (orange). The real part is always non-negative and increasing towards infinity while the imaginary part is constant.

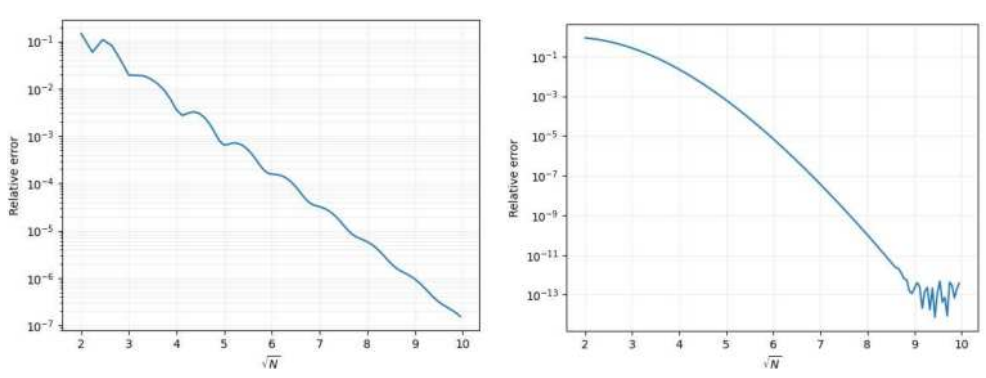

Table 2: Root-exponential convergence of Gauss-Laguerre quadrature along the contour (left) for small imaginary wavenumber $s = \mathrm{i}$ and (right) for larger wavenumber $s = 100\,\mathrm{i}$

# Fast 3D Time Domain Boundary Element Approach for Elastodynamics based on Multivariate Adaptive Cross Approximation

**Vibudha Lakshmi Keshava**[1], **Martin Schanz**[1,*]

[1]Institute of Applied Mechanics, Graz University of Technology, Graz, Austria

*Email: m.schanz@tugraz.at

### Abstract

A time domain boundary element formulation for elastodynamics is proposed. The respective integral equation is discretised in space with low order shape functions and in time with the generalized Convolution Quadrature (gCQ) method. The fast and data sparse representation is based on the multivariate adaptive cross approximation, enabling fast methods for both spatial and temporal discretisation. Combined with gCQ, the algorithm adaptively selects a minimal set of complex frequencies for temporal discretisation. For spatial degrees of freedom, $\mathcal{H}$-matrices with classical ACA or a fast multipole method based on a kernel interpolation are applied. The combination of spatial and temporal fast methods yields a data-sparse, efficient solver for elastodynamics. Parametric studies demonstrate the technique's performance.

## 1 Introduction

Time-domain boundary integral equations are well known for solving many engineering problems, e.g., in acoustics and elastodynamics. The corresponding boundary element (BE) formulations are also well-known. Here, the method based on convolution quadrature is considered. Formulations for a non-uniform time discretisation also exist, which utilise the generalised Convolution Quadrature (gCQ) [3]. However, all time-domain BE formulations require fast methods to handle real-world three-dimensional problems. Here the so-called multivariate adaptive cross approximation [1] will be used based either on a fast multipole method (3D-FMM) or on adaptive cross approximation (3D-ACA) to establish a time domain elastodynamic BE formulation. The respective formulation for acoustics can be found in [2, 4].

## 2 Model and Algorithm

Let $\Omega \subset \mathbb{R}^3$ be a bounded Lipschitz domain and $\Gamma = \partial\Omega$ its boundary with the outward normal $\mathbf{n}$ and the time interval $t \in [0, T]$. Denoting with $u$ the displacement vector and with $q$ the traction vector the first boundary integral equation for elastodynamics

$$\mathcal{C}(\mathbf{x})\mathbf{u}(\mathbf{x}, t) = (\mathcal{V} * \mathbf{q})(\mathbf{x}, t) - (\mathcal{K} * \mathbf{u})(\mathbf{x}, t) \quad (1)$$

for all $(\mathbf{x}, t) \in \Gamma \times (0, T)$ can be derived. The integral free term is denoted with $\mathcal{C}$, the single layer potential with $\mathcal{V}$, and the double layer potential with $\mathcal{K}$. Both last mentioned operators contain the fundamental solution of elastodynamics or its co-normal derivative the tractions. Here, this integral equation is sufficient as a collocation approach is used with a simple splitting of the boundary with respect to given boundary conditions. This step is performed after a spatial discretisation with linear shape functions for the displacements and constant shape functions for the tractions.

The time discretisation with gCQ and $[0, T] = [0, t_1, t_2, \dots, t_N], \Delta t_n = t_n - t_{n-1}$ results in a discrete integral equation on the boundary $\Gamma$ at discrete times $t_n$. The respective expression for the single layer potential has the form

$$\begin{aligned}(\mathsf{V} * \mathsf{q}^h)(t) \approx & \hat{\mathsf{V}}\left((\Delta t_n \mathsf{A})^{-1}\right)(\mathsf{q}^h)_n + \\ & \sum_{\ell=1}^{N_Q} \hat{\mathsf{V}}(s_\ell)\, \mathsf{W}_\ell^{\Delta t_n}\Big((\mathsf{q}^h)_{n-1}\Big)\end{aligned} \quad (2)$$

with $\mathsf{W}_\ell^{\Delta t_n}\big((\mathsf{q}^h)_{n-1}\big) = \omega_\ell \left(\mathsf{b}^\mathsf{T}\mathsf{A}^{-1} \cdot \mathsf{x}_{n-1}(s_\ell)\right)$ $(\mathsf{I} - \Delta t_n s_\ell \mathsf{A})^{-1}\,\mathbb{1}$ where $\mathsf{A} \in \mathbb{R}^{m\times m}$, $\mathsf{b}, \mathsf{c} \in \mathbb{R}^m$ describes the underlying A- and L-stable Runge-Kutta method of order $m$. As usual the superscript $h$ indicates a discretised quantity and $\hat{()}$ denotes that the elements in this matrix are Laplace transformed quantities and $s_\ell$ are the complex Laplace variables at the $N_Q$ discrete complex frequencies (details can be found in [2, 3]).

The multivariate adaptive cross approximation is applied on the sum on the right hand side of (2) and can be represented as the tensor product

$$\sum_{\ell=1}^{N_Q} \hat{\mathsf{V}}\left(s_\ell\right) \mathsf{W}_\ell^{\Delta t_n} \approx \sum_{k=1}^{r} \mathbf{H}_k \left(\hat{\mathsf{V}}\right) \otimes \sum_{\ell=1}^{N_Q} \mathbf{f}_k[\ell] \mathsf{W}_\ell^{\Delta t_n} , \tag{3}$$

where the rank $r$ is signifanct smaller than $N_Q$. Essentially, the matrix $\mathbf{H}_k$ contains the spatial discretisation and the vector $\mathbf{f}_k$ the frequency information. The former can be approximated by either a classical FMM (3D-FMM) or an adaptive cross approximation (3D-ACA) approach. Hence, the complexity of the original operation $\mathcal{O}\left(N_Q M^2\right)$ for $M$ spatial unknowns is reduced to $\mathcal{O}\left(r(M \log M + N_Q)\right)$. The essential difference of the proposed algorithm between the published acoustic formulation [4] and elastodynamics is the selection for the adaptively determined complex frequencies. In elastodynamics $\mathbf{f}_k$ is a vector of $3 \times 3$-blocks and not a scalar as in acoustics. The algorithm requires a pivot block to be selected at each iteration, which is determined by sorting the minimal singular values of each block and then selecting the largest.

## 3 Numerical Results

The proposed method is tested by applying a smooth function, which solves the governing equation in the full space, on a unit cube. The triangulation and the time steps are refined, where the overall observation time $T$ is kept constant. Based on the known solution, an $L_2$-error in space and time can be defined to determine the convergence rate. Using a two stage Radau IIA method for the gCQ and linear/constant shape function in space, the Dirichlet data converge with order two and the Neumann data with order one. Applying as well 3D-ACA or 3D-FMM to obtain a data sparse method the compression rates given in Fig. 1 can be obtained. The compression is defined as the ratio between the storage necessary for a dense computation, i.e., without any fast methods, and the storage used in the 3D-ACA or 3D-FMM. In Fig. 1a, a pure Dirichlet problem is shown and in Fig. 1b, a mixed problem where one half of the cube has Dirichlet boundary conditions and the other half Neumann boundary conditions. It is obvious that a dramatic reduction in storage is possible, which is mainly driven by a very low rank $r$ and not that much by the ACA or FMM for the matrices $\mathbf{H}_k$. For the latter, larger amount of spatial degrees of freedom are necessary.

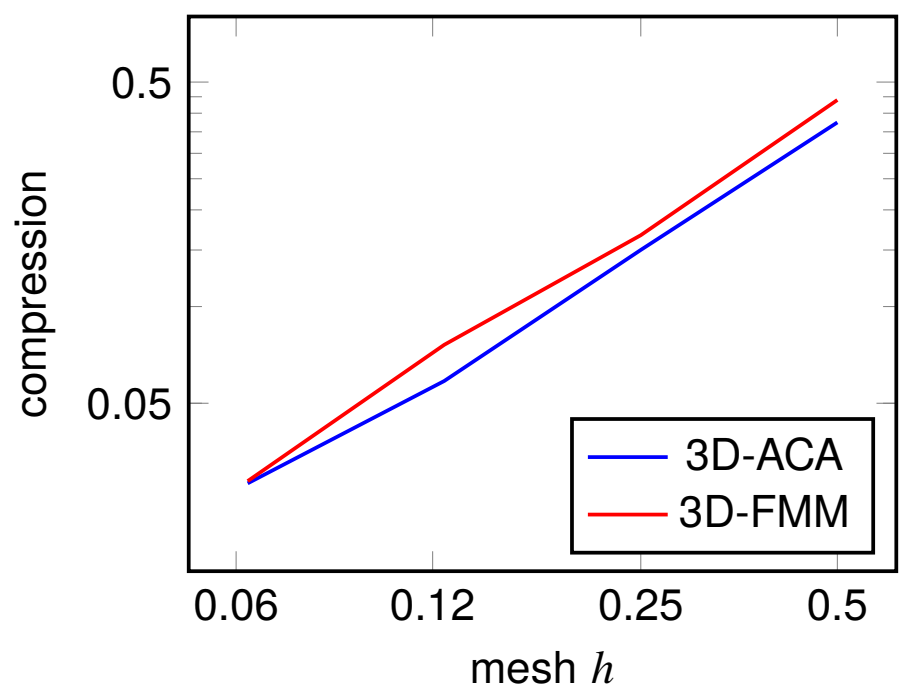

(a) Dirichlet problem

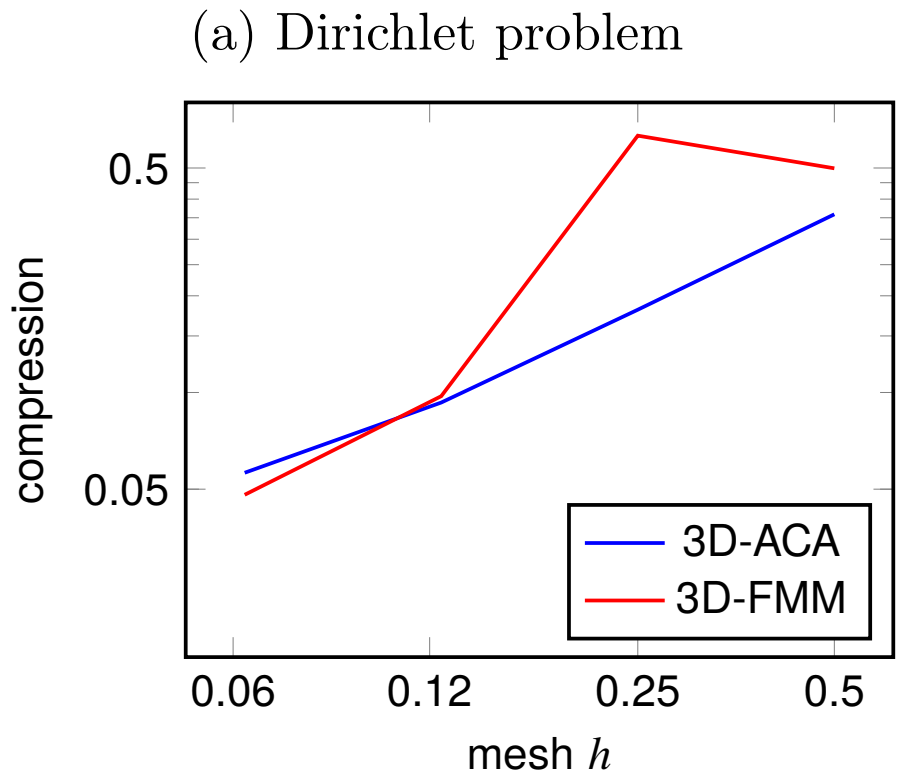

(b) Mixed problem

Figure 1: Compression versus refinement in space $h$ and time $\Delta t$

# Direct reconstruction for acoustic inverse Born scattering

**Lisa Schätzle**[1,*], **Nuutti Hyvönen**[1]
[1]Department of Mathematics and Systems Analysis, Aalto University, Finland
*Email: lisa.schatzle@aalto.fi

**Abstract**

In this talk, we derive a reconstruction formula for the inverse medium scattering problem for the Helmholtz equation in two dimensions under the Born regime. Through an appropriate choice of expansion functions for both the unknown contrast and the observed far field data, the resulting system matrix decouples in angular direction and attains, for each angular frequency separately, a triangular structure that allows for efficient inversion. We show numerical examples and demonstrate that the resulting reconstruction formula together with an adequate regularization method remains effective even when applied to full nonlinear far field data beyond the Born regime.

## 1 Introduction

We consider the inverse medium scattering problem for the Helmholtz equation in two dimensions, i.e., the task to recover a compactly supported contrast function from full knowledge of the associated far field data or, equivalently, of the associated far field operator. Although this problem is uniquely solvable, it is (i) severely ill-posed as small perturbations in the observed far field data may lead to large reconstruction errors and (ii) nonlinear. In the regime of weak scattering, the Born approximation yields a linearized relation between the contrast and the far field data and thus overcomes difficulty (ii). In this situation, the inverse scattering problem reduces to inverting non-uniformly sampled restricted Fourier data of the contrast. The specific sampling structure of these Fourier data enables us, building on recent work [1, 2] for linearized EIT, to derive an explicit reconstruction formula that expresses the expansion coefficients of the contrast in terms of those of the far field data.

## 2 Inverse medium scattering

Let $\kappa > 0$ denote the wave number and $\Omega \subset \mathbb{R}^2$ a penetrable scatterer modelled by the compact support of a contrast function $q \in L^\infty(\Omega)$ satisfying $q > -\infty$ a.e. in $\mathbb{R}^2$, where $\partial\Omega$ is assumed to be of Lipschitz class. This scatterer is illuminated by a plane wave incident field

$$u^i(x, d) := \mathrm{e}^{\mathrm{i}\kappa x\cdot d}\,, \qquad x \in \mathbb{R}^2\,,$$

at illumination direction $d \in S^1$, which produces a scattered field $u^s$ such that the associated total field $u := u^i + u^s \in H^1_{\mathrm{loc}}(\mathbb{R}^2)$ satisfies the scattering problem for the Helmholtz equation

$$\Delta u + \kappa^2(1 + q)u = 0 \qquad \text{in } \mathbb{R}^2\,, \tag{1a}$$

$$\lim_{|x|\to\infty} \sqrt{|x|}\Big(\frac{\partial u^s}{\partial |x|}(x, d) - \mathrm{i}\kappa u^s(x, d)\Big) = 0 \tag{1b}$$

uniformly w.r.t. $\widehat{x} := x/|x| \in S^1$. The unique weak solution $u \in H^1_{\mathrm{loc}}(\mathbb{R}^2)$ of (1) can be characterized by the Lippmann-Schwinger equation

$$\begin{aligned} u(\cdot, d) &= u^i(\cdot, d) + \frac{\mathrm{i}\kappa^2}{4}\int_\Omega q(y)u(y, d)H_0^{(1)}(\kappa|\cdot - y|)\,\mathrm{d}y \\ &=: u^i(\cdot, d) + L_q u(\cdot, d) \qquad \text{in } \Omega \end{aligned} \tag{2}$$

with $L_q : L^2(\Omega) \to L^2(\Omega)$ and $H_0^{(1)}$ denoting the Hankel function of first kind and order zero. For $|x| \to \infty$ the scattered field $u^s$ fulfills the asymptotic expansion

$$u^s(x, d) = \frac{\mathrm{e}^{\mathrm{i}\pi/4}}{\sqrt{8\pi}}\frac{\mathrm{e}^{\mathrm{i}\kappa|x|}}{\sqrt{\kappa|x|}}u^\infty(\widehat{x}, d) + O\big(|x|^{-\frac{3}{2}}\big)\,,$$

uniformly w.r.t. $\widehat{x} \in S^1$, where the far field pattern $u^\infty \in L^2(S^1 \times S^1)$ is given by

$$u^\infty(\widehat{x}, d) \;=\; \kappa^2 \int_\Omega q(y)u(y, d)\mathrm{e}^{-\mathrm{i}\kappa\widehat{x}\cdot y}\,\mathrm{d}y\,.$$

Knowledge of $u^\infty(\widehat{x}, d)$ for all $\widehat{x}, d \in S^1$ uniquely determines the contrast $q$, but this inversion is severely ill-posed and nonlinear.

As soon as $\|L_q\|_{L^2(\Omega)\leftarrow L^2(\Omega)} \ll 1$ with $L_q$ given as in (2), the far field pattern $u^\infty$ can well-approximated by its associated Born far field pattern

$$u_B^\infty(\widehat{x}, d) := \kappa^2 \int_{B_R(\boldsymbol{c})} q(y) \mathrm{e}^{\mathrm{i}\kappa(d-\widehat{x})\cdot y} \, \mathrm{d}y \,,$$

which linearizes $u^\infty$ with respect to $q$.

## 3 Matrix representation in Born regime

Assuming $\Omega \subset B_R(c)$, we aim to construct orthonormal systems (ONS) $(f_{m,n})_{m,n}$ of $L^2(S^1 \times S^1)$ and $(\Psi_{j,k})_{j,k}$ of $L^2(B_R(c))$ such that the system matrix describing the mapping of the expansion coefficients of $q$ w.r.t. $(\Psi_{j,k})_{j,k}$ onto the expansion coefficients of $u_B^\infty$ w.r.t. $(f_{m,n})_{m,n}$

(i) decouples in $j$ direction,

(ii) exhibits a triangular structure for each $j$.

Denoting by

$$f_m(\widehat{x}) := \frac{1}{\sqrt{2\pi}} \mathrm{e}^{\mathrm{i}m \arg \widehat{x}} \,, \quad m \in \mathbb{Z} \,, \; \widehat{x} \in S^1 \,,$$

the standard Fourier basis of $L^2(S^1)$ and building on the findings from [3], we choose the modulated Fourier basis

$$f_{m,n}(\widehat{x}, d) := \mathrm{e}^{\mathrm{i}\kappa c\cdot(d-\widehat{x})} f_m(\widehat{x}) f_{-n}(d) \,,$$

$m, n \in \mathbb{Z}$, $\widehat{x}, d \in S^1$, which forms an ONS of $L^2(S^1 \times S^1)$.

From [3] we know that the expansion coefficients of $u_B^\infty$ w.r.t. $(f_{m,n})_{m,n}$ are given by

$$a_{m,n} := (2\pi)^{3/2} (\kappa R)^2 \mathrm{i}^{n-m} \int_{B_1(0)} q(Ry + c) \times f_{n-m}(\widehat{y}) J_m(\kappa R|y|) J_n(\kappa R|y|) \, \mathrm{d}y \,, \quad (3)$$

$m, n \in \mathbb{Z}$, and that these are, caused by the decay behavior of the Bessel function $J_m$, essentially supported on the index set $[-N, N]^2$ with $N \in \mathbb{N}$ slightly larger than $\kappa R$.

A suitable choice of an ONS of $L^2(B_1(0))$ for expanding $q(R(\cdot) + c)$ in (3) is suggested by the following lemma.

**Lemma 1** *There exists an ONS* $(\Psi_{j,k})_{j\in\mathbb{Z},k\in\mathbb{N}_0}$ *of* $L^2(B_1(0))$ *such that*

*(i)* $\Psi_{j,k}(y) = f_j(\widehat{y}) R_k^{|j|}(|y|)$ *for* $y \in B_1(0)$ *and*

*(ii) for each* $j \in \mathbb{Z}$ *and* $K \in \mathbb{N}$ *there holds*

$$\operatorname{span}\{R_k^{|j|} : k \in \{0, \dots, K - \lceil \tfrac{|j|}{2} \rceil\}\} = \operatorname{span}\{P_k^{|j|} : k \in \{\lceil \tfrac{|j|}{2} \rceil, \dots, K\}\} \,,$$

*where* $P_k^{|j|}(t) := J_k(\kappa R t) J_{k-|j|}(\kappa R t)$ *for* $t \in (0,1)$.

Unfortunately, the functions $(R_k^{|j|})_{j,k}$ are not known in explicit form but, however, they can for each $|j|$ be constructed numerically by a Gram-Schmidt orthonormalization of $(P_k^{|j|})_k$ w.r.t. the weighted inner product $\langle \cdot, \cdot \rangle_{L_t^2(0,1)}$.

As our main result, we obtain the following formula for the expansion coefficients of the rescaled and shifted contrast $q(R(\cdot) + c)$.

**Theorem 2** *Denote by* $(\Psi_k^{|j|})_{j,k}$ *the ONS from Lemma 1. Then, for each* $j \in \mathbb{Z}$ *the expansion coefficients* $(c_{j,k})_k$ *of* $q(R(\cdot)+c)$ *w.r.t. this ONS can be computed successively by the formula*

$$c_{j,k+\lceil j/2 \rceil} = \frac{a_{k+\lceil j/2 \rceil, k-\lfloor j/2 \rfloor}}{(2\pi)^{3/2} (\kappa R)^2 (-\mathrm{i})^j \langle R_k^{|j|}, P_{k+\lceil j/2 \rceil}^{|j|} \rangle_{L_t^2(0,R)}} - \sum_{i=0}^{k-1} \frac{\langle R_i^{|j|}, P_{k+\lceil j/2 \rceil}^{|j|} \rangle_{L_t^2(0,R)}}{\langle R_k^{|j|}, P_{k+\lceil j/2 \rceil}^{|j|} \rangle_{L_t^2(0,R)}} c_{j,i+\lceil j/2 \rceil} \,. \quad (4)$$

The inner products in (4) are known from the Gram-Schmidt orthonormalization. Additional stabilization of the reconstruction formula (4) is achieved by restricting the computation to the essential index support $[-N, N]^2$ of the given data $(a_{m,n})_{m,n}$ and by employing a truncated singular value decomposition as a regularization method.

# Boundary determination in anisotropic elasticity

**Joonas Ilmavirta[1], Mikko Salo[1], Hjørdis Schlüter[2],*, Daniel Windisch[3]**

[1]Department of Mathematics and Statistics, University of Jyväskylä, Jyväskylä, Finland

[2]Department of Mathematics and Statistics, University of Helsinki, Helsinki, Finland

[3]Department of Computer Science, KU Leuven, 3001 Leuven, Belgium

*Email: hjordis.schluter@helsinki.fi

**Abstract**

We address the inverse problem of recovering an anisotropic stiffness tensor $c$ and the density of mass $\rho$ from the associated Dirichlet-to-Neumann map in two- and three-dimensional domains. We show that for generic stiffness tensors the normalized stiffness tensor $a = \rho^{-1}c$ is uniquely determined up to zeroth order by the associated Dirichlet-to-Neumann map at the boundary. This result builds on a construction of special solutions that concentrate near a boundary point and are highly oscillatory. The generic assumption is needed in order to ensure that the slowness surface corresponding to the stiffness tensor is square free and to ensure that eigenvectors corresponding to different roots in the slowness surface are linearly independent.

## 1 Introduction

The elastic medium is modeled by $\Omega \subset \mathbb{R}^n$ with $n \geq 2$. The stiffness of the material is modeled as a rank 4 tensor $c \in \mathbb{R}^{n\times n\times n\times n}$ and the density of mass is modeled as a scalar-valued function $\rho \in \mathbb{R}$. The stiffness tensor links the stress and strain tensors of the material and following their symmetries and their relationship, $c$ has the following minor and major symmetries

$$c_{ijk\ell} = c_{jik\ell} = c_{k\ell ij}, \tag{1}$$

and is positive in the sense that

$$c^{ijk\ell} A_{ij} A_{k\ell} \geq \delta A_{ij} A_{k\ell} \tag{2}$$

for some $\delta > 0$ and for all symmetric matrices $A$. The waves are modelled by the following linear elastic wave equation in $\Omega \times (0, T)$:

$$\begin{aligned}&(Pu(x,t))_i := \\ &\rho(x)\partial_t^2 u_i(x,t) - \partial_{x_j}(c_{ijk\ell}(x)\partial_{x_k} u_\ell(x,t)) \\ &= 0,\end{aligned} \tag{3}$$

for $i = 1, ..., n$ and where we denote by $P$ the elastic wave operator that acts on a displacement $u \in \mathbb{R}^n$. Imposing a displacement at the bottom of the space time cylinder $\Omega\times[0,T]$ gives rise to the following boundary value problem for $u$:

$$\begin{cases} Pu = 0 & \text{in } \Omega \times (0,T) \\ u = f & \text{on } \partial\Omega \times (0,T) \\ u(x,0) = \partial_t u(x,0) = 0 & \text{in } \Omega, \end{cases}$$

with $f(x,0) = 0$ and $\partial_t f(x,0) = 0$ for $x \in \partial\Omega$. The boundary measurements that we consider is the Dirichlet-to-Neumann (DN) map $\Lambda$ that maps the Dirichlet boundary function $f$ to the corresponding Neumann boundary function and which is defined as follows:

$$(\Lambda_{(\rho,c)} f)_i = \nu_j\, c_{ijk\ell}\, \partial_{x_k} u_\ell \Big|_{(0,T)\times\partial\Omega} \tag{4}$$

for $i = 1, ..., d$ and where $\nu$ denotes the unit outward normal to $\partial\Omega$.

If we expand the elastic wave operator $P$ from above, we obtain:

$$P_{i\ell} = \delta_{i\ell}\partial_t^2 - a_{ijk\ell}\partial_j\partial_k - \rho^{-1}(\partial_j c_{ijk\ell})\partial_k.$$

This yields the principal symbol:

$$(p_2)_{i\ell}(x,t,q,\omega) = -\delta_{i\ell}\omega^2 + a_{ijk\ell} q_j q_k$$

We introduce the slowness vector $\xi = \omega^{-1} q$ and define the Christoffel matrix:

$$\Gamma_{i\ell} := a_{ijk\ell}\xi_j\xi_k$$

With these definitions the principal symbol is on the form:

$$p_2(x,t,\xi,\omega) = \omega^2 \left[\Gamma(x,\xi) - I\right].$$

The slowness surface $\sigma$ is defined by the kernel of the principal symbol. It is defined as a polynomial in $\xi$:

$$\sigma(x,\xi) := \det\left(\Gamma(x,\xi) - I\right) = 0.$$

In the following we consider polynomials $\sigma$ that are square free in the sense that $\sigma$ cannot be factored into polynomials that have power larger than 1.

The inverse problem that we address is to recover $c$ and $\rho$ from $\Lambda_{(c,\rho)}$ at the boundary $\partial\Omega$. This problem has been studied for isotropic stiffness tensors in [1] and been generalized to transversely isotropic and orthorhombic stiffness tensors [2]. We obtain boundary determination results for the normalized stiffness tensor $a$ with $a_{ijk\ell} = \rho^{-1}c_{ijk\ell}$. This is a first step towards proving [3, Conjecture 2.4] that conjectures that $\Lambda_{(c,\rho)}$ determines $c$ and $\rho$ uniquely in $\Omega$ (in dimensions $n = 2$ and $n = 3$).

## 2 Results

Let us now formulate the main result that we have obtained.

**Theorem 1** *Let $\Omega \subset \mathbb{R}^n$ with $n = 2$ or $n = 3$ be a bounded and connected domain with a smooth boundary $\partial\Omega$. There is an open and dense set of stiffness tensors that satisfy the symmetry and positivity conditions in (1)-(2) for $i = 1, 2$. If the DN maps defined in (4) satisfy $\Lambda_{\rho_1,c_1} = \Lambda_{\rho_2,c_2}$ on $\partial\Omega \times (0, T)$ for all $T > 0$ then*

$$\rho_1^{-1}c_1 = \rho_2^{-1}c_2$$

*on $\partial\Omega$.*

This result builds on the construction of special solutions to (3) near a boundary point, following the procedure in [4, Thm 5.1] for real principal type operators with scalar-valued principal symbols. We extend this construction to real principal type operators with matrix-valued principal symbol. For the argument to apply, the slowness polynomial $\sigma(x, \xi)$ associated with the stiffness tensor must be square free. Moreover, to obtain the conclusion of Theorem 1, the eigenvectors of $\Gamma(x, \xi_1, \ldots, z)$ corresponding to distinct roots $z_1$ and $z_2$ of $\sigma(x, \xi_1, \ldots, z)$ must be linearly independent. Both squarefreeness and linear independence hold for most, but not all, stiffness tensors; consequently, Theorem 1 is valid for an open and dense subset of stiffness tensors. Establishing the existence of such a set requires tools from algebraic geometry. To show that it is non-empty, one must construct, in each dimension $n$, an example for which the slowness polynomial is square free and the corresponding eigenvectors are linearly independent. In two dimensions, this can be achieved by the isotropic stiffness tensor

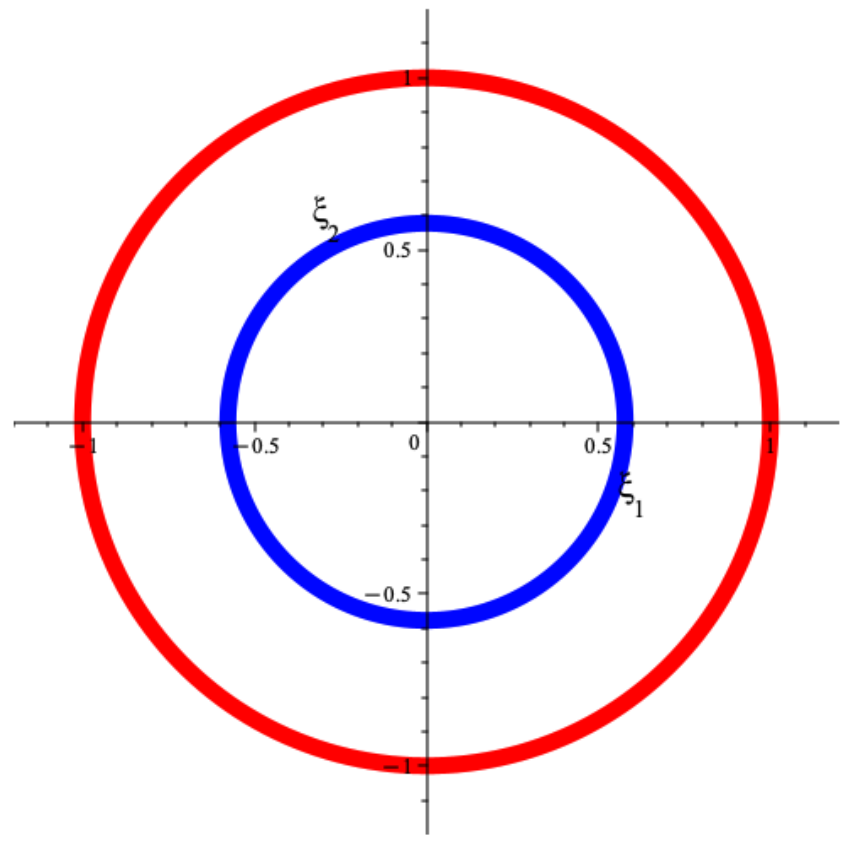

Figure 1: Slowness surface $\sigma(x, \xi) = 0$ for the isotropic stiffness tensor in equation (5).

$$c = \begin{bmatrix} 3 & 1 & 0 \\ 1 & 3 & 0 \\ 0 & 0 & 1 \end{bmatrix}, \tag{5}$$

whose slowness surface is illustrated in Figure 1.

### Acknowledgements

The present work is supported by the Research Council of Finland (Flagship of Advanced Mathematics for Sensing Imaging and Modelling grant 359208).

# Broken-FEEC on multipatch domains with local refinements

Frederik Schnack[1,*], Martin Campos Pinto[1]
[1]Max-Planck-Institut for Plasmaphysics, Garching b. München, Germany
*Email: frederik.schnack@ipp.mpg.de

**Abstract**

We discuss recent advancements in broken-FEEC (Finite Element Exterior Calculus) methods in two dimensions and their extension to locally refined multipatch spline spaces with non-matching interfaces.

## 1 Introduction

Broken-FEEC methods like [1] allow for discontinuous discretizations at patch interfaces, leading to block-diagonal mass matrices, while preserving the de Rham structure. In the case of local refinements, where the two discretizations on shared patch-interfaces might differ from each other, this method leads to numerical artifacts, for example spurious reflections. There are many areas of research on how to deal with such effects, for example by employing artificial dissipation. In order to resolve this issue, we developed discrete conforming projection operators, that preserve polynomial moments. These operators are metric-independent, localized, and explicit, ensuring continuity across non-matching interfaces while preserving high-order polynomial moments. This results in broken-FEEC de Rham sequences with accurate strong and weak derivatives, leading to energy-preserving Maxwell solvers that are explicit and virtually free of spurious modes. Numerical simulations confirm the efficacy of these methods, for example, in eliminating spurious waves without any dissipation mechanism.

## 2 (Non-)conforming de Rham diagrams

We consider a multipatch domain of the form

$$\Omega = \operatorname{int}\Big(\bigcup_{k\in\mathcal{K}} \overline{\Omega_k}\Big) \quad \text{with} \quad \Omega_k = F_k(\hat\Omega), \quad (1)$$

with disjoint, geometrically conforming subdomains $\Omega_k$ associated with smooth mappings $F_k$ defined on a reference domain $\hat\Omega = ]0,1[^2$.

At the heart of our discretization lies the following de Rham diagram:

$$\begin{array}{ccccc}
V^0 & \xrightarrow{\nabla} & V^1 & \xrightarrow{\operatorname{curl}} & V^2 \\
\downarrow \Pi^0 & & \downarrow \Pi^1 & & \downarrow \Pi^2 \\
V_h^0 & \xrightarrow{\nabla} & V_h^1 & \xrightarrow{\operatorname{curl}} & V_h^2 \\
\uparrow P^0 & & \uparrow P^1 & & \uparrow P^2 \\
V_{\mathrm{pw}}^0 & \xrightarrow{\nabla P^0} & V_{\mathrm{pw}}^1 & \xrightarrow{\operatorname{curl} P^1} & V_{\mathrm{pw}}^2
\end{array} \quad (2)$$

In the first line of the diagram (2), we have the continuous $\nabla$-curl de Rham sequence, consisting of the continuous spaces $V^0 = H^1(\Omega)$, $V^1 = H(\operatorname{curl};\Omega)$ and $V^2 = L^2(\Omega)$.

In the bottom line, there are the (discontinuous) patch-wise discrete spaces $V^\ell_{\mathrm{pw}}$ for $\ell = 0,1,2$, which are the juxtaposition of the finite element spaces $V_k^\ell$ defined on each patch, with possibly different resolution, as the pull-back $V_k^\ell = \mathcal{F}_k^\ell(\hat V_k^\ell)$ of the logical spaces $\hat V_k^\ell$, which we choose to be the usual tensor-product spline spaces of [2]. They also form a discrete de Rham sequence with the patch-wise differential operators.

FEEC is based on a sequence of conforming finite element spaces

$$V_h^\ell = \{v \in V^\ell : v|_{\Omega_k} \in V_k^\ell\} \subset V^\ell, \quad (3)$$

that form discrete de Rham sequences, as seen in the second line of the diagram (2). The existence of bounded projections $\Pi^\ell : V^\ell \to V_h^\ell$ that commute with the differential operators, i.e.

$$\nabla \Pi^0 \phi = \Pi^1 \nabla \phi \qquad \forall \phi \in V^0, \quad (4)$$
$$\operatorname{curl} \Pi^1 \boldsymbol{u} = \Pi^2 \operatorname{curl} \boldsymbol{u} \qquad \forall \boldsymbol{u} \in V^1, \quad (5)$$

is the key tool for obtaining stable discretizations, for which we refer to [3].

In order to relate the patch-wise spaces $V^\ell_{\mathrm{pw}}$ to the conforming spaces $V_h^\ell$, we introduce the conforming projection operators $P^\ell$ and discuss their properties in the next section.

## 3 Conforming projection operators

The conforming projection operators

$$P^\ell : V^\ell_{\mathrm{pw}} \to V^\ell_h \subset V^\ell_{\mathrm{pw}}, \tag{6}$$

are the key ingredients for our discrete schemes. They establish vertex- and edge-conformity sequentially, for which they involve averaging operations along non-matching interfaces in order to project onto a shared interface space. In particular, we define extension and restriction operators mapping from the coarse to the fine one-dimensional space and vice versa.

In addition, we impose the preservation of polynomial moments of degree $r > 0$, which is defined by

$$\int P^\ell \phi q_j = \int \phi q_j, \quad \forall \phi \in V^\ell_{\mathrm{pw}}, \tag{7}$$

where the polynomials $q_j$, for $j = 0, \ldots, r$, form a basis of the polynomial space of degree $r$. This introduces a set of linear constraints on the projectors $P^\ell$, but ensures high-order approximation of their adjoint operator $(P^\ell)^*$. In turn, this enables high-order approximation of weak differential operators in our schemes and controls the formation of spurious reflections.

## 4 Time-domain Maxwell problem

We consider the time-domain Maxwell problem on a square domain: For $\boldsymbol{J} \in L^2(\Omega)$, solve

$$\begin{aligned} \partial_t \boldsymbol{E} - \mathbf{curl}\ B &= -\boldsymbol{J}, \\ \partial_t B + \mathrm{curl}\ \boldsymbol{E} &= 0, \end{aligned} \tag{8}$$

where $\boldsymbol{E} \in H_0(\mathrm{curl}, \Omega)$ and $B \in L^2(\Omega)$.

The discrete scheme in the broken-FEEC framework, using a simple leap-frog time-stepping with $\Delta t > 0$, then reads

$$\begin{aligned} \mathbf{b}^{n+\frac{1}{2}} &= \mathbf{b}^n - \frac{\Delta t}{2}\mathbf{C}\mathbf{P}^1\mathbf{e}^n, \\ \mathbf{M}^1\mathbf{e}^{n+1} &= \mathbf{M}^1\mathbf{e}^n + \Delta t(\mathbf{C}\mathbf{P}^1)^T\mathbf{M}^2\mathbf{b}^{n+\frac{1}{2}} \\ \mathbf{b}^{n+1} &= \mathbf{b}^{n+\frac{1}{2}} - \frac{\Delta t}{2}\mathbf{C}\mathbf{P}^1\mathbf{e}^{n+1}, \end{aligned} \tag{9}$$

with coefficients $\mathbf{b}$, $\mathbf{e}$ and $\mathbf{j}$ and where $\mathbf{M}^\ell$ are block-diagonal mass-matrices, $\mathbf{C}$ is the discrete patch-wise curl operator and $\mathbf{P}^1$ the conforming projection operator.

We discretize the square domain $\Omega = [0, \pi]^2$ by $3 \times 3$-patches, where the enclosed patch in the middle is of a coarse resolution, while the other surrounding patches are refined. In order to compare a discretization with and without the preservation of polynomial moments, we consider the initial condition

$$\begin{aligned} \boldsymbol{E}(t=0,\boldsymbol{x}) &= \frac{1}{\sigma^2}\begin{pmatrix} \boldsymbol{x}_2 - \bar{\boldsymbol{x}}_2 \\ -(\boldsymbol{x}_1 - \bar{\boldsymbol{x}}_1) \end{pmatrix} e^{\frac{-\|\boldsymbol{x}-\bar{\boldsymbol{x}}\|_2^2}{2\sigma^2}}, \\ \boldsymbol{B}(t=0,\boldsymbol{x}) &= \mathrm{curl}\,\boldsymbol{E}, \quad \boldsymbol{J}(t,\boldsymbol{x}) = 0, \end{aligned} \tag{10}$$

where $\bar{\boldsymbol{x}} = (0.5\pi, 0.5\pi)$ and $\sigma = 0.1$, and show the resulting amplitude of the electric field in Figure 1.

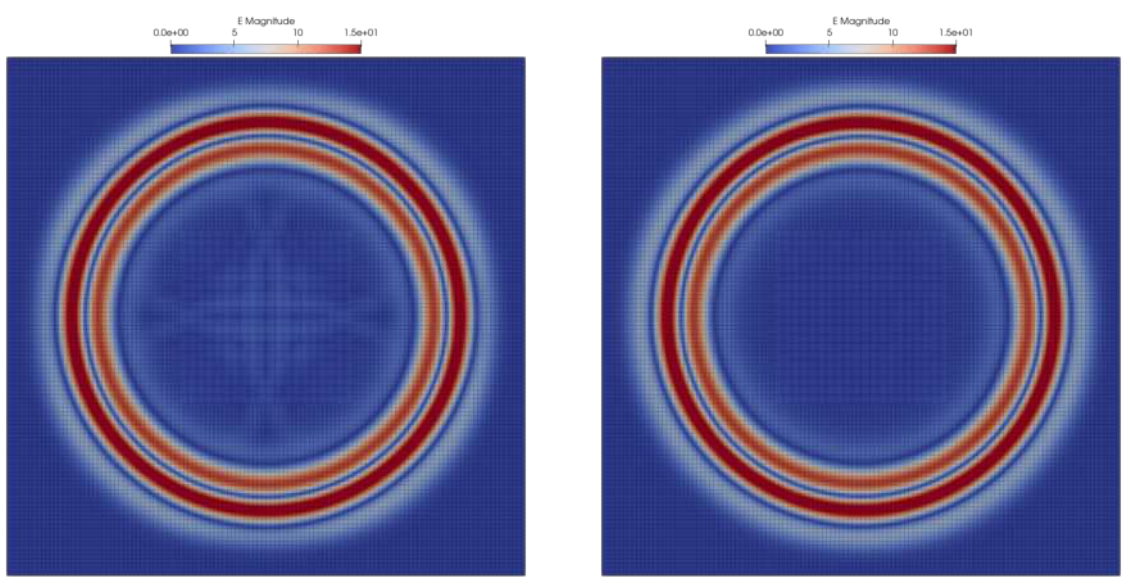

Figure 1: The amplitude of the electric field after leaving the interior (coarse) patch. The discretization on the right involves a conforming projection that preserves polynomial moments, while the discretization on the left does not.

After the initial Gaussian pulse leaves the interior patch, we observe spurious reflections if the conforming projection operators do not preserve polynomial moments. These reflections are significantly reduced and virtually vanishing if the conforming projection operators preserve polynomial moments, and we are able to recover the optimal convergence order.

# Nonlinear Acoustic Model Equations resulting from Non-Relativistic and Relativistic Lagrangian Pairing Inspired by Quantum Theory

**Markus Scholle**[1,*], **Philip H. Gaskell**[2], **Meriem Al Ahmadi Hammou**[1]

[1]Institute for Flow in Additively Manufactured Porous Media (ISAPS), Heilbronn Univ., Germany.
[2]Department of Engineering, Durham University, Durham, UK

*Email: markus.scholle@hs-heilbronn.de

**Abstract**

A method of obtaining acoustic model equations (linear and nonlinear) directly from the Lagrangian of the respective continuum is established, and the question of how the corresponding Lagrangian densities can be obtained is examined. The pairing of non–relativistic and relativistic Lagrangians, originally discovered in quantum theory, plays a crucial role, ultimately enabling the derivation of acoustic model equations in pairs of both non-relativistic and associated relativistic forms. The procedure is exemplified using the Lagrangian of a potential flow, with the newly developed systematic approach used to first obtain its relativistic counterpart, followed by the Kuznetsov equation and its relativistic analogue as associated acoustic equations.

## 1 From Lagrangians to acoustic models

For a continuum whose dynamics is determined by a Lagrangian $\ell\left(\psi_n, \partial_t\psi_n, \nabla\psi_n\right)$ depending on $N$ independent fields, $\psi_n$, it is standard practice to first obtain the Euler–Lagrange equations by variation, leading to the equations of motion which, in turn, serve as the starting point for deriving acoustic model equations by considering small or moderate deviations, $\tilde{\psi}_n$, from a given state of equilibrium, $\psi_{0n}$, i.e.:

$$\psi_n = \psi_{0n} + \tilde{\psi}_n\,. \tag{1}$$

Applying a Taylor expansion, followed by elimination of unknowns, leads to a wave equation for a single unknown, as per the left hand path of Fig. 1. An alternative approach proposed by [1], based on applying the Taylor expansion directly to the Lagrangian and defining a perturbation Lagrangian, $\Omega\left(\tilde{\psi}_n, \partial_t\tilde{\psi}_n, \nabla\tilde{\psi}_n\right)$, consisting of all quadratic terms of the expansion, follows the right hand path of Fig. 1. This can be extended to a weakly nonlinear analysis by additionally considering the cubic terms.

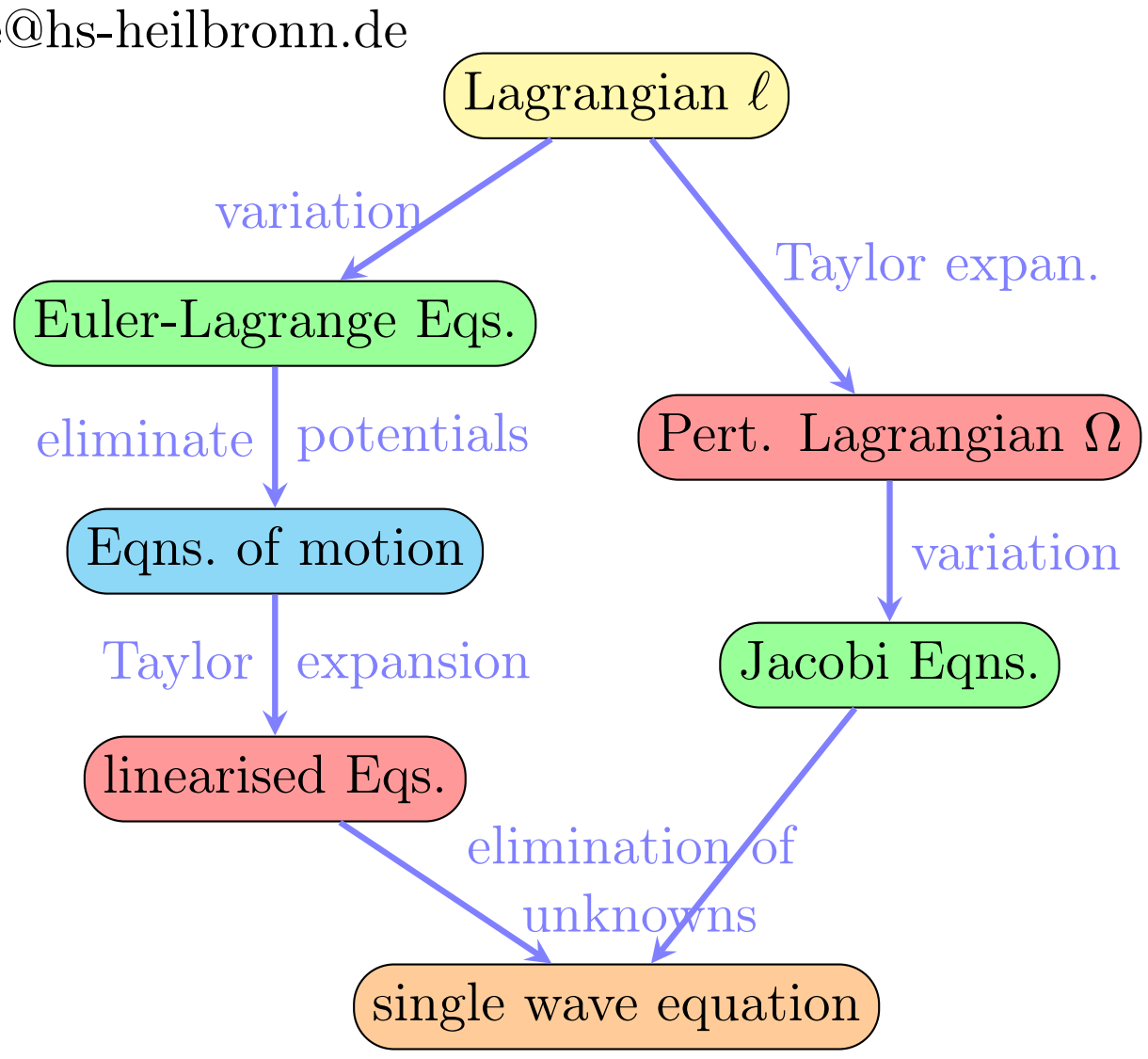

Figure 1: Obtaining a wave equation from a Lagrangian: proposed procedure (right hand path), standard way (left hand path).

## 2 Construction of Lagrangians

Despite the elegance of the left-hand path, it requires a Lagrangian for the respective continuum as a starting point. In order to address this issue, [2] provides the following general construction scheme for Lagrangians:

$$\ell = L\left(\omega, \vartheta_i, \overset{\odot}{\vartheta}_i, \nabla\vartheta_i\right)\,, \tag{2}$$

with:

$$\omega = \partial_t\phi + \frac{(\nabla\phi)^2}{2}\,, \quad \overset{\odot}{\vartheta}_i = \partial_t\vartheta_i + \nabla\phi\cdot\nabla\vartheta_i\,, \tag{3}$$

which results from inverse treatment of Noether currents related to universal symmetries such as spatial/temporal translations and Galilei boosts, with invariant fields $\vartheta_i = \psi_{i+1}$, $i = 1, \cdots, N-1$ and a potential field $\phi = \psi_1$ that transforms w.r.t. Galilei boosts according to $\phi \to \phi - \vec{u}\cdot\vec{x} + \vec{u}^2 t/2$. Furthermore, the scheme (2) ensures the equality $\vec{p} = \varrho\vec{u}$ between momentum density and mass flux density. As an attempt toward a relativistic generalisation, [3] reveal, having replaced Galilei boosts by Lorentz boosts, the same scheme (2) results but with modifications:

$$\omega = \frac{c^2}{2} - \frac{1}{2}\partial^\mu\phi\partial_\mu\phi\,, \quad \overset{\odot}{\vartheta}_i = -\partial^\mu\phi\partial_\mu\vartheta_i\,, \tag{4}$$

suggesting the concept of *Lagrangian pairing* between a non–relativistic and an associated relativistic Lagrangian, both resulting from one 'blueprint form' (2) by applying either the substitution rules (3) for Galilean invariance or (4) for Poincaré invariance. Serving as a prominent example, the Schrödinger Lagrangian and the Klein–Gordon Lagrangian prove to be such a pair. In their follow–up article [4], normalisation of velocity $u_\mu u^\mu = 1$ and the eigenvalue relation $T^{\mu\nu}u_\nu = mn\left(c^2+e\right)u^\mu$ are also taken into account so as to restrict the 'blueprint Lagrangian' to the form:

$$\frac{L}{mn} = c\sqrt{c^2 - 2\omega + F\left(\theta^j, \overset{\odot}{\theta^j}, \partial\theta^j\right)} - c^2 - e \quad (5)$$

with particle density $n$, specific inner energy $e$ and a function $F$ that takes into account the individual continuum. Lagrangian pairing, illustrated in Fig. 2, allows determination of relativistic pendants to given classical Lagrangians.

$$\frac{L}{mn} = c\sqrt{c^2 - 2\omega + F\left(\theta^j, \overset{\odot}{\theta^j}, \partial\theta^j\right)} - c^2 - e$$

non-relativ. $(c \to \infty)$
$\omega \to \partial_t\phi + \frac{1}{2}(\nabla\phi)^2$
$\overset{\odot}{\vartheta} \to \partial_t\vartheta + \nabla\phi\cdot\nabla\vartheta$

relativistic
$\omega \to \frac{c^2}{2} - \frac{1}{2}\partial^\mu\phi\partial_\mu\phi$
$\overset{\odot}{\vartheta} \to -\partial^\mu\phi\partial_\mu\vartheta$

Figure 2: Lagrangian pairing schematic.

## 3 Resulting model equations

As an introductory example, $F=0$, related to a relativistic barotropic potential flow with fields $n$ and $\phi$, yields the perturbation Lagrangian:

$$\Omega = -m\tilde{n}\partial_t\phi - \frac{K_0\tilde{n}^2}{2n_0^2} - \left(\frac{B}{A}-1\right)\frac{K_0\tilde{n}^3}{6n_0^3} - \frac{m}{2}\left[n_0 + \tilde{n} + \frac{n_0\partial_t\phi}{c^2(1+\hat{h}_0)}\right]\frac{(\nabla\phi)^2}{1+\hat{h}_0}, \quad (6)$$

with respective relativistic specific enthalpy, compression modulus and nonlinearity parameter:

$$\hat{h}_0 := \frac{1}{c^2}\frac{\mathrm{d}}{\mathrm{d}n}(ne)\bigg|_{n=n_0}, \quad (7)$$

$$K_0 := mn_0^2\frac{\mathrm{d}^2(ne)}{\mathrm{d}n^2}\bigg|_{n=n_0}, \quad (8)$$

$$\frac{B}{A} - 1 := \frac{mn_0^3}{K_0}\frac{\mathrm{d}^3(ne)}{\mathrm{d}n^3}\bigg|_{n=n_0}. \quad (9)$$

A linear analysis yields D'Alembert's equation with relativistically corrected small–signal speed:

$$a_0 = \sqrt{\frac{K_0}{mn_0\left(1+\hat{h}_0\right)}}, \quad (10)$$

and weakly nonlinear analysis the relativistic analogue of frictionless Kuznetsov's equation:

$$-\partial_t^2\phi + a_0^2\nabla^2\phi = \partial_t\left[\frac{1-a_0^2/c^2}{1+\hat{h}_0}(\nabla\phi)^2 + \frac{B/A - a_0^2/c^2}{2a_0^2\left(1+\hat{h}_0\right)}(\partial_t\phi)^2\right], \quad (11)$$

where additional contributions due to relativistic effects are highlighted in blue.

Although less relevant to acoustic phenomena under terrestrial conditions, relativistic acoustic model equations are relevant to astrophysical problems, such as those involving neutron star cores, where the speed of sound can even exceed the conformal limit of $c/\sqrt{3}$ [5].

Via further examples of increasing complexity, the inclusion of thermodynamic degrees of freedom and irreversible effects — such as viscosity and heat conduction — into the Lagrangian framework is readily demonstrable.

# Homogenization of wave equations and dispersive profiles for ring solutions

**Ben Schweizer**[1,*], **Grégoire Allaire**[2], **Agnes Lamacz**[1]
[1]Department of Mathematics, TU Dortmund, Dortmund, Germany
[2]CMAP, Ecole Polytechnique, Palaiseau, France
*Email: ben.schweizer@tu-dortmund.de

## Abstract

When waves travel over large distances and for long time, one can typically observe ring structures. When waves travel through a medium with microstructure, one can additionally observe dispersive effects. The dispersive effects are small when the microstructure is small, but over long time spans, they create contributions of leading order. We present an overview over known results, including stochastic media, and present new contributions that cover systems such as elasticity. We derive a dispersive effective system with the Bloch method of homogenization. The approximate representation formulas for solutions in Fourier space allow additionally to describe solutions as superpositions of ring waves. The rings expand with constant speed, their profiles change on a slow time scale according to the dispersive terms.

## 1 Introduction

We study wave systems of the form

$$\partial_t^2 u^\varepsilon(x,t) = \nabla \cdot (A(x/\varepsilon)\nabla u^\varepsilon(x,t)) \,. \quad (1)$$

The independent variables are $(x,t) \in \mathbb{R}^d \times [0,\infty)$ for space dimension $d \in \mathbb{N}$, the solution is a map $u^\varepsilon : \mathbb{R}^d \times [0,\infty) \to \mathbb{R}^p$ for some $p \in \mathbb{N}$. We are interested in the description of solutions in the limit $\varepsilon \to 0$. When solutions are studied on fixed time intervals $t \in [0,T_0]$, this is the classical task of homogenization theory. The result is essentially that solutions $u^\varepsilon$ converge to the solution $u$ of the limit problem which is defined with the classically homogenized elliptic operator. One has to be careful with initial conditions, we mention [3] for such an analysis.

New effects appear when large time intervals $t \in [0,T_0\varepsilon^{-2}]$ are studied, the effective system then is a dispersive wave equation. While such large time effects are well understood for scalar equations, systems are more complicated.

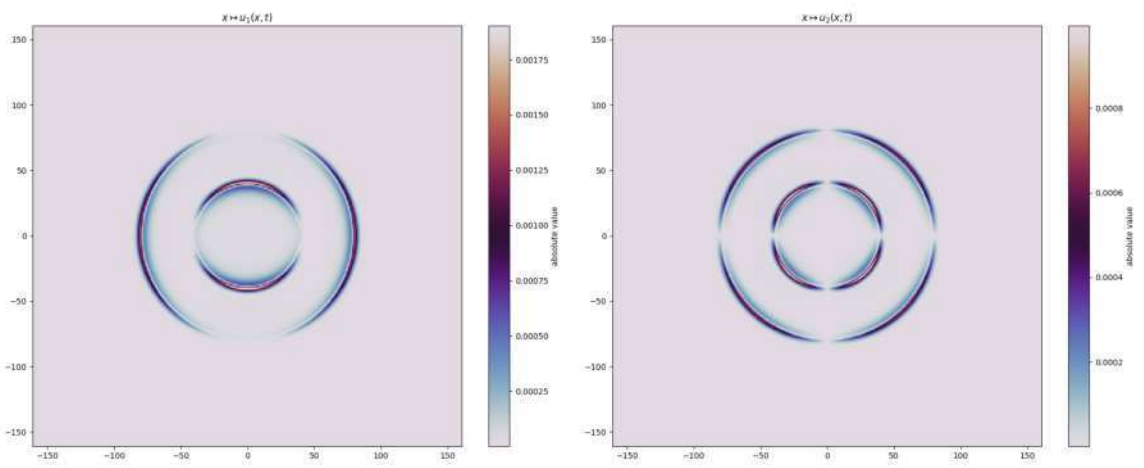

Figure 1: The solution $u$ to the Lamé system, plotted are the absolute values of the two real components $u_1$ (left) and $u_2$ (right) of the solution $x \mapsto u(x,t_0)$ in physical space, both $x_1$ and $x_2$ run from $-150$ to 150 and the values of $u$ range from 0 to $10^{-3}$. The pressure wave is concentrated along the outer ring, the shear wave along the inner ring.

To simplify formulas, we restrict ourselves to initial data

$$u^\varepsilon(x,0) = u_0(x)\,, \quad \partial_t u^\varepsilon(x,0) = 0\,, \quad (2)$$

the elliptic operator is

$$L_\varepsilon u^\varepsilon(x,t) := -\nabla \cdot (A(x/\varepsilon)\nabla u^\varepsilon(x,t)) \,, \quad (3)$$

where $A : \mathbb{R}^d \ni y \mapsto A(y) \in \mathbb{R}^{d\times d\times p\times p}$ is $2\pi$-periodic in every direction and $A(y) : \mathbb{R}^{p\times d} \to \mathbb{R}^{p\times d}$ is self-adjoint for every $y$. We study two coercivity assumptions. One is that $\zeta : A(y)\zeta \geq \gamma\|\zeta\|^2$ holds with $\gamma > 0$ independent of $y \in \mathbb{R}^d$ and $\zeta \in \mathbb{R}^{d\times p}$. The other one is adapted to the system of elasticity; essentially, we consider $d = p$ and demand coercivity on symmetric matrices combined with the property that $A(y)\zeta$ vanishes for skew-symmetric matrices $\zeta$.

## 2 Homogenization

With $\hat{k} = k/|k|$, we define $U^\varepsilon$ as the inverse Fourier transform of

$$\hat{U}^\varepsilon(k,t) = \sum_{m=1}^{p} B_m(k)\xi_m(\hat{k}) \exp\left(i\omega_m^*(\varepsilon k)t/\varepsilon\right) \,. \quad (4)$$

In this expression, the terms $\xi_m$ and $\omega_m^*$ are obtained from the coefficient function $A$ of the elliptic operator $L_\varepsilon$ of (3). The $p$ functions $k \mapsto$

$\omega_m^*(k)$ are constructed from the Taylor expansion of the dispersion relation,

$$\omega_m^*(\varepsilon k) = c_m(\hat{k})|\varepsilon k| + b_m(\hat{k})|\varepsilon k|^3 \,. \tag{5}$$

The $p$ functions $(\xi_m)_{m\leq p}$ are approximations of Bloch functions. The coefficients $(B_m(k))_{m\leq p}$ are constructed from the initial data $u_0$ using projections that depend on the vectors $\xi_m$.

**Theorem 1 (Homogenization)** *We consider initial data* $u_0 \in L^1(\mathbb{R}^d, \mathbb{C}^p) \cap L^2(\mathbb{R}^d, \mathbb{C}^p)$ *with compact support in Fourier space and* $T_0 > 0$. *Let* $L_\varepsilon$ *be as in* (3) *and let* $u^\varepsilon$ *be the solution of* $\partial_t^2 u^\varepsilon = -L_\varepsilon u^\varepsilon$ *with initial condition* (2). *Let* $U^\varepsilon$ *be defined by the Fourier representation formula* (4). *Then there holds, as* $\varepsilon \to 0$,

$$\sup_{t\in[0,T_0/\varepsilon^2]} \|u^\varepsilon(\cdot,t) - U^\varepsilon(\cdot,t)\|_{L^2(\mathbb{R}^d)} \to 0 \,. \tag{6}$$

## 3 Ring solutions

A further aim is to show that $U^\varepsilon$ of (4) is close to a superposition of $p$ rings. Here, we can assume that the functions $B_m$, $\xi_m$ and $b_m$ are smooth functions, that $B_m$ has a compact support and that 0 is not contained in the support of $B_m$. Closeness is measured uniformly in time on intervals $[t_0/\varepsilon^2, T_0/\varepsilon^2]$ with $0 < t_0 < T_0$ and in the sense of $L^2(\mathbb{R}^d)$ in space as in (6).

Functions $\hat{G}_m$ are defined from the $p$ summands in the representation of $\hat{U}^\varepsilon(k,t)$ in (4). Profile functions $V_m$

$$V_m : \mathbb{R} \times S^{d-1} \times (0,\infty) \to \mathbb{C}^p \tag{7}$$

are defined using slices of $\hat{G}_m$ and a one-dimensional inverse Fourier transform thereof, examples for $V_m$ are shown in Fig. 2. The single ring function $w_m$ is defined by setting $w_m^\varepsilon(x,t)$ equal to

$$\frac{1}{(c_m t)^{(d-1)/2}} V_m\left(|x| - c_m t, \frac{x}{|x|}, \varepsilon^2 t\right) \,, \tag{8}$$

the reconstruction with $p$ rings is

$$w^\varepsilon(x,t) := \sum_{m=1}^{p} w_m^\varepsilon(x,t) \,. \tag{9}$$

Some regularity and growth assumptions must be made for the subsequent convergence result.

**Theorem 2 (Ring approximation)** *Let* $p \in \mathbb{N}$ *be the number of rings and let the dimension be* $d = 2$ *or* $d = 3$ *and let* $T_0 > 0$. *We assume that the speeds* $c_m(q) = c_m$ *are independent of the direction* $q$ *and different from each other,* $c_m \neq c_\ell$ *for* $m \neq \ell$. *Let* $U^\varepsilon : \mathbb{R}^d \times (0,\infty) \to \mathbb{C}^p$ *have the Fourier-transform* (4) *with* $\omega_m^*$ *of* (5). *Let* $w^\varepsilon : \mathbb{R}^d \times (0,\infty) \to \mathbb{C}^p$ *be given by* (8)–(9). *Then, for an arbitrary sequence* $\varepsilon \to 0$ *and an arbitrary time sequence* $t_\varepsilon \in (0, T_0/\varepsilon^2)$ *with* $\liminf_{\varepsilon\to 0} \varepsilon^2 t_\varepsilon > 0$, *there holds*

$$\|w^\varepsilon(\cdot,t_\varepsilon) - U^\varepsilon(\cdot,t_\varepsilon)\|_{L^2(\mathbb{R}^d)} \to 0 \,. \tag{10}$$

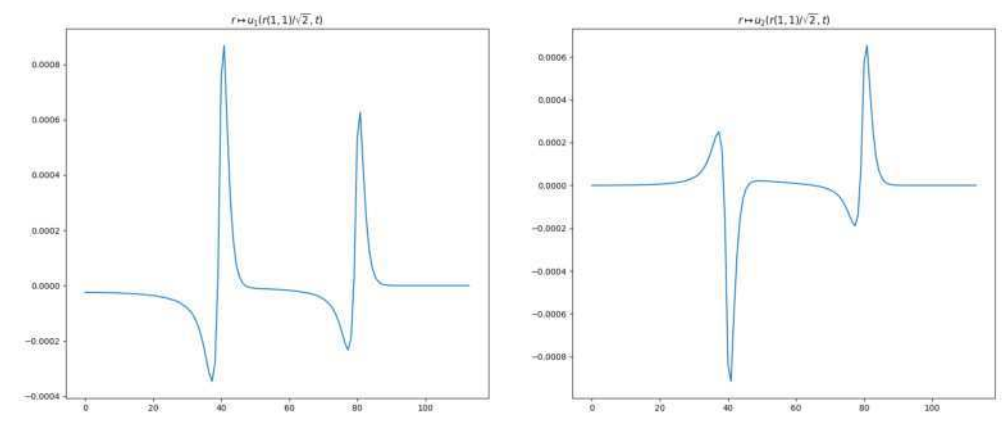

Figure 2: Numerically determined profiles in a 45°-direction ($x_1 = x_2$), first and second component of the solution $u$. The horizontal axis shows the radial variable (from 0 to 120), the vertical axis shows $u_1$ and $u_2$ in the range from about $-4\cdot 10^{-4}$ to $8\cdot 10^{-4}$. The profile of the shear-wave is centered at $c_s t_0 = 40$, the pressure-wave at $c_p t_0 = 80$.

# Application of the multivariate Loewner Framework to the modelling of acoustic propagation in rigid-frame porous materials

Lucie Schwoebel[2,3,*], Charles Poussot-Vassal[1], Denis Matignon[2], Rémi Roncen[3], Estelle Piot[3], Pierre Vuillemin[1]

[1]DTIS, ONERA, Université de Toulouse, 31000, Toulouse, France
[2]Fédération ENAC ISAE-SUPAERO ONERA, Université de Toulouse, 31000, Toulouse, France
[3]DMPE, ONERA, Université de Toulouse, 31000, Toulouse, France

*Email: lucie.schwoebel@onera.fr

**Abstract**

The propagation of acoustic waves in rigid-frame porous materials is commonly described using the Equivalent Fluid (EF) model, which depends on multiple parameters such as porosity and permeability. These porous properties are usually not accurately known, thus affecting the predicted acoustic response. In this work, reduced-order modelling is employed to approximate with controlled accuracy the parametrized EF model with a rational multivariate representation, using an extension of the Loewner Framework (LF), namely the multivariate LF (mLF) [1]. The latter is a data-driven approach allowing for constructing rational $n$-variable models by means of interpolation. It is adapted for constructing of models that depend on parameters, describing either uncertainties or parametric dependences. The main result of [1] concerns the variable decoupling achieved in the LF, providing a solution to the Kolmogorov Superposition Theorem restricted to rational functions. As a byproduct, it tames the curse of dimensionality in storage, computational complexity and accuracy, allowing dealing with a high number of parameters.

## 1 Multivariate Loewner Framework

The mLF rationale consists in re-arranging the data, computing the barycentric coefficients via the Loewner matrix, and constructing the rational interpolant [1].

**Tensor grid data.** The data-set (tensor grid) is the primary ingredient. It is obtained by evaluating the $n$-variable function $\mathbf{H}(x_1, \cdots, x_n)$, either through simulations or experiments along the discretized grid $\mathbf{x}_1, \cdots, \mathbf{x}_n$ [1]. In the context of the LF, grid points $\mathbf{x}_l$ are split into columns $\lambda_l(j_l)$ and rows $\mu_l(i_l)$ interpolation, or support points[2]. Evaluating the multivariate function $\mathbf{H}$ along with the combinations of the support points $\lambda_l(j_l), \mu_l(i_l) \in \mathbb{C}$ thus forms an $n$-dimensional tensor, denoted $\mathcal{T}_n^{\otimes}$ [3]. Following the Loewner philosophy detailed in [1], we define $P_c^{(n)}$, the column data, and $P_r^{(n)}$, the row data, being two subsets of $\mathcal{T}_n^{\otimes}$, leading to $\mathbf{w}_{j_1,j_2,\cdots,j_n}$ and $\mathbf{v}_{i_1,i_2,\cdots,i_n}$:

$$\begin{aligned} P_c^{(n)} &: \{(\lambda_1(j_1), \lambda_2(j_2), \cdots, \lambda_n(j_n); \mathbf{w}_{j_1,j_2,\cdots,j_n})\} \\ P_r^{(n)} &: \{(\mu_1(i_1), \mu_2(i_2), \cdots, \mu_n(i_n); \mathbf{v}_{i_1,i_2,\cdots,i_n})\}. \end{aligned}$$

**Multivariate Loewner matrix.** The tensor data now serves for constructing the $n$-D Loewner matrix $\mathbb{L}_n$, which may be viewed as an operator mapping the interpolation points and $\mathcal{T}_n^{\otimes}$ onto a $Q \times K$ matrix, with $Q = q_1 q_2 \dots q_n$ (rows) and $K = k_1 k_2 \dots k_n$ (columns), where each entry of the $\mathbb{L}_n$ matrix is given by

$$\ell_{j_1,j_2,\cdots,j_n}^{i_1,i_2,\cdots,i_n} = \frac{\mathbf{v}_{i_1,i_2,\cdots,i_n} - \mathbf{w}_{j_1,j_2,\cdots,j_n}}{(\mu_1(i_1) - \lambda_1(j_1)) \cdots (\mu_n(i_n) - \lambda_n(j_n))}.$$

**Lagrangian (barycentric) rational model.** By considering an appropriate number of column interpolation points $k_l$ $(l = 1, \cdots, n)$, one can compute $\mathbb{L}_n \mathbf{c}_n = 0$, the right null space of $\mathbb{L}_n$, which contains the so-called barycentric coefficients, $\mathbf{c}_n^\top = [c_{1,\cdots,1}, \cdots, c_{1,\cdots,k_n} | \cdots, c_{k_1,\cdots,k_n}] \in \mathbb{C}^K$. Then, the multivariate Lagrangian (or barycentric) form $\mathbf{G}_{\text{lag}}(x_1, \cdots, x_n) =$

$$\frac{\sum_{j_1=1}^{k_1} \cdots \sum_{j_n=1}^{k_n} \frac{c_{j_1,\cdots,j_n} \mathbf{w}_{j_1,\cdots,j_n}}{(x_1 - \lambda_1(j_1)) \cdots (x_n - \lambda_n(j_n))}}{\sum_{j_1=1}^{k_1} \cdots \sum_{j_n=1}^{k_n} \frac{c_{j_1,\cdots,j_n}}{(x_1 - \lambda_1(j_1)) \cdots (x_n - \lambda_n(j_n))}},$$

interpolates the $n$-D data tensor and eventually reveals the original underlying function $\mathbf{H}$ (if rational). Reducing either $k_l$, or directly $K$ (the null space entries), reduces the complexity of $\mathbf{G}_{\text{lag}}$ and leads to tensor approximation.

[1] $\mathbf{x}_l = [x_l(1), \cdots, x_l(N_l)]$ and $x_l(i)$ is the $l$-th variable evaluated at the $i$-th element (with $l = 1, \cdots, n$).

[2] We also denote $j_l = 1, \cdots, k_l$ and $i_l = 1, \cdots, q_l$ with $l = 1, \cdots, n$ and $N_l = k_l + q_l$.

[3] We define $\mathbf{x}_1 = [\lambda_1(j_1), \mu_1(i_1)]$, $\mathbf{x}_2 = [\lambda_2(j_2), \mu_2(i_2)]$, ..., $\mathbf{x}_n = [\lambda_n(j_n), \mu_n(i_n)]$.

**Decoupling of variables.** The principal bottleneck of the above Loewner-driven model construction is the need for constructing a $Q \times K$ matrix and computing a $K$ dimensional null space. In [1, Theorem 5.8], the $n$-D Loewner null space $\mathbf{c}_n$ is expressed as a linear combination of a collection of 1-D Loewner matrix null spaces. A direct consequence is [1, Theorem 5.9] which states how to decouple the evaluation of the barycentric coefficients, reducing complexity and memory cost [ [1, Theorems 10 & 13], allowing the framework to be applied to large multi-dimensional tensors.

## 2 Porous media modelling

In this study, the mLF is applied to the EF model of rigid-frame porous media, that describes the homogenised behaviour of the acoustic pressure and velocity fields, governed by the coupled PDEs given in the Laplace domain by:

$$\begin{cases} \rho_0\, \alpha(s,\mathbf{x}) \cdot s\widehat{\mathbf{u}} &= -\nabla \widehat{p}, \\ \chi_0\, \beta(s,\mathbf{x}) s\widehat{p} &= -\nabla \cdot \widehat{\mathbf{u}}, \end{cases}$$

where $\rho_0$, $\chi_0 \in \mathbb{R}^+$ are the ambient fluid density and compressibility, $p$ the spatial average of the acoustic pressure in the porous domain and $\mathbf{u}$ the spatial average of the particle velocity in the pore, normalised by the total volume. The symbol $\widehat{\cdot}$ denotes the Laplace transform and $s = \mathrm{j}\omega$ the associated complex variable.

The functions $\alpha(s,\mathbf{x})$ and $\beta(s,\mathbf{x})$ are irrational multivariate functions of $s$ and of the parameter vector $\mathbf{x} = [\phi, \overline{r}, \sigma_r]$, respectively pore mean size, pore size standard deviation and porosity, all independent positive real numbers. They are defined as [2]:

$$\begin{cases} \alpha(s,\mathbf{x}) = \frac{e^{4S_2}}{\phi}\left[1 + \frac{a_1}{\mathrm{j}\omega\overline{r}^2} e^{6S_2} + \frac{a_1 e^{6S_2}\sqrt{1+\frac{\mathrm{j}\omega\overline{r}^2}{2a_1}e^{-7S_2}} - 1}{\mathrm{j}\omega\overline{r}^2}\right], \\ \beta(s,\mathbf{x}) = \phi\gamma - \phi(\gamma-1)\left[1 + \frac{a_2 e^{-2S_2}}{\mathrm{j}\omega\overline{r}^2} + \frac{a_2 e^{-2S_2}\sqrt{1+\frac{\mathrm{j}\omega\overline{r}^2 e^{3S_2}}{a_3}} - 1}{\mathrm{j}\omega\overline{r}^2}\right]^{-1}. \end{cases}$$

Noting $S_2 = (\sigma_r \log_{10} 2)^2$ and $\gamma$, $a_{1,2,3}$ constant positive values associated to the ambient fluid.

The acoustic behaviour can then be described by the material impedance $Z$ or absorption coefficient $A$, with $h$ the material thickness and $Z_0$, $c_{1,2}$ constant positive values associated to the ambient fluid:

$$Z = -\mathrm{j}\sqrt{\frac{c_1\alpha}{\beta}}\cot\left(\omega h\sqrt{c_2\alpha\beta}\right) \text{ and } A = 1 - \left|\frac{Z - Z_0}{Z + Z_0}\right|^2.$$

## 3 Application

The mLF is then applied on the porous media EF model constitutive functions $\alpha(s,\mathbf{x})$ and $\beta(s,\mathbf{x})$. This method identifies the optimal approximation order for each parameter, following a user-defined threshold $\epsilon$ depending on the required approximation fidelity (Fig. 1).

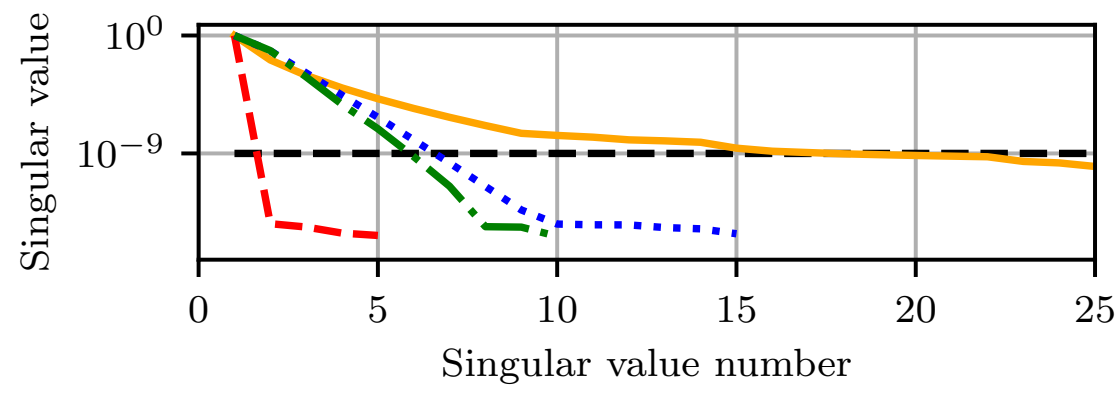

Figure 1: *Singular values for $\omega$ (—), $\phi$ (- -), $\overline{r}$ (···) and $\sigma_r$ (-·-·). Threshold $\epsilon$ (- -).*

From the approximated functions, we can then get the absorption coefficient for different parameter configurations to validate the approximation (Fig. 2).

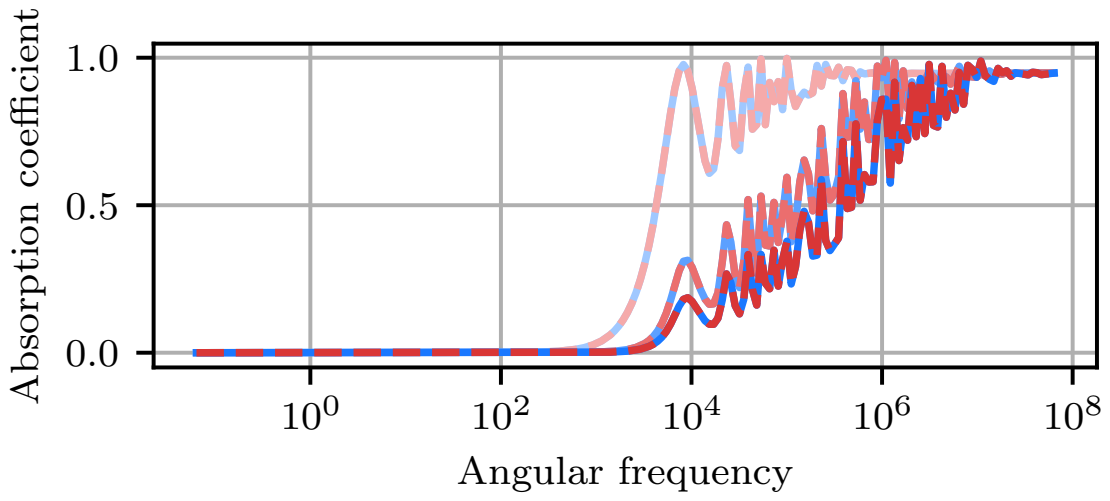

Figure 2: *Absorption coefficient from frequency model (—) and from its approximation (- -) for three different pore mean sizes (darker shade for bigger pores). Mean absolute error of order* $10^{-5}$.

By partial fraction decomposition with fixed parameters $\mathbf{x} = [x_1 = \phi, x_2 = \overline{r}, x_3 = \sigma_r]$, the multivariate Lagrangian form $\mathbf{G}_{\text{lag}}(s, x_1, \cdots, x_n)$ can then be expressed in the Laplace domain, $G(s) = \sum_{k=1}^{k_1} \frac{r_k(\mathbf{x})}{s+\xi_k(\mathbf{x})}$, or in the time-domain, $g(t) = \sum_{k=1}^{k_1} r_k(\mathbf{x})\, e^{-\xi_k(\mathbf{x})t} 1_{t>0}$. This allows for an $\mathbf{x}$-dependent time-domain expression of the acoustic behaviour of porous media.

# Passivity enforcing multipole based method for time domain modelling of acoustic wave propagation in rigid frame porous media

**Lucie Schwoebel[1,2,*], Denis Matignon[1], Remi Roncen[2], Estelle Piot[2]**
[1]Fédération ENAC ISAE-SUPAERO ONERA, Université de Toulouse, 31000, Toulouse, France
[2]DMPE, ONERA, Université de Toulouse, 31000, Toulouse, France
*Email: lucie.schwoebel@onera.fr

**Abstract**

Simulation of the acoustic wave propagation in rigid-frame porous media in the time domain involves multipole approximations of the frequency-dependent functions which model the homogenised behaviour of the acoustic pressure and velocity fields within the material. Their irrational components can be approximated through a diffusive representation framework. The behaviour of the homogenised porous model is then expressed as a sum of first-order low-pass filters, and previous work has proved that the passivity condition leads to stability in those augmented systems [1]. However conventional pole-fitting methods may compromise the passivity, particularly for more complex media, e.g., the double porosity media, formed by the Boolean intersection of two solid pore networks, showcasing a diffusion effect at their interface. To address the passivity issue, a new passivity-enforcing identification strategy is introduced, providing a robust framework for porous material time-domain acoustic simulations [2].

## 1 Introduction

The homogenised behaviour of the acoustic pressure and velocity fields within rigid-frame porous media is described in the Laplace domain by the equivalent fluid model (EFM):

$$\begin{cases} \rho_0\, \alpha(s) \cdot \widehat{\partial_t \mathbf{u}} &= -\nabla \widehat{p}, \\ \chi_0\, \beta(s) \widehat{\partial_t p} &= -\nabla \cdot \widehat{\mathbf{u}}, \end{cases}$$

where $\rho_0\,,\, \chi_0 \in \mathbb{R}^+$ are the ambient fluid density and compressibility, $p$ the spatial average of the acoustic pressure and $\mathbf{u}$ the spatial average of the particle velocity in the pore, normalised by the total volume. The symbol $\widehat{\cdot}$ denotes the Laplace transform and $s = \mathrm{j}\omega$ is the complex variable of the Laplace domain.

The porous medium is represented by two frequency-dependent irrational functions, $\alpha(s)$ - dynamic viscous tortuosity - and $\beta(s)$- dynamic thermal compressibility [1]:

$$\begin{cases} \alpha(s) &= \frac{\alpha_\infty}{\phi} \left[ 1 + \frac{M}{s} + N \frac{\sqrt{1+\frac{s}{L}}-1}{s} \right], \\ \beta(s) &= \phi\gamma - \phi(\gamma-1) \left[ 1 + \frac{M'}{s} + N' \frac{\sqrt{1+\frac{s}{L'}}-1}{s} \right]^{-1} \end{cases}$$

with $\alpha_\infty$, $\phi$, $M$, $N$, $L$, $M'$, $N'$, $L' > 0$ the material properties and $\gamma$ the air heat capacity ratio.

These transfer functions are first expressed as diffusive representations in the frequency domain:

$$G(s) = \frac{c_0}{s} + \int_0^{+\infty} \frac{\mu(\xi)}{s+\xi}\, d\xi + c_\infty\,, \qquad (1)$$

with $\mu$ the diffusive weights ($\int_0^{+\infty} \frac{\mu(\xi)}{1+\xi}\, d\xi < \infty$), $c_0 = \lim_{s\to 0} sG(s)$ and $c_\infty = \lim_{s\to+\infty} G(s)$. In the time domain, the corresponding impulse response is:

$$g(t) = c_0 1_{t>0} + \int_0^{\infty} \mu(\xi)\, e^{-\xi t} 1_{t>0}\, d\xi + c_\infty \delta_0(t)\,.$$

Then, they are approximated in the frequency domain by $G_{\text{approx}}(s) = \frac{c_0}{s} + \sum_{q=1}^{Q} \frac{r_q}{s+\xi_q} + c_\infty$, a rational function.

## 2 Passivity of the representation

Previous work showed that the EFM system is stable provided that $\xi_q > 0\,,\, \forall q$ and $\forall\, \Re(s) \geq 0$, $\Re\left(G_{\text{approx}}(s)\right) \geq 0$, see [1][4].

When modelling more complex media, conventional pole-fitting methods may compromise passivity. To address this issue, a new passivity-enforcing identification strategy is introduced [2]. The method proceeds in two stages, see Fig. 1: poles are first distributed on a logarithmic grid and weights computed analytically under a reality constraint; then both poles and weights are refined jointly in log space to enforce positivity.

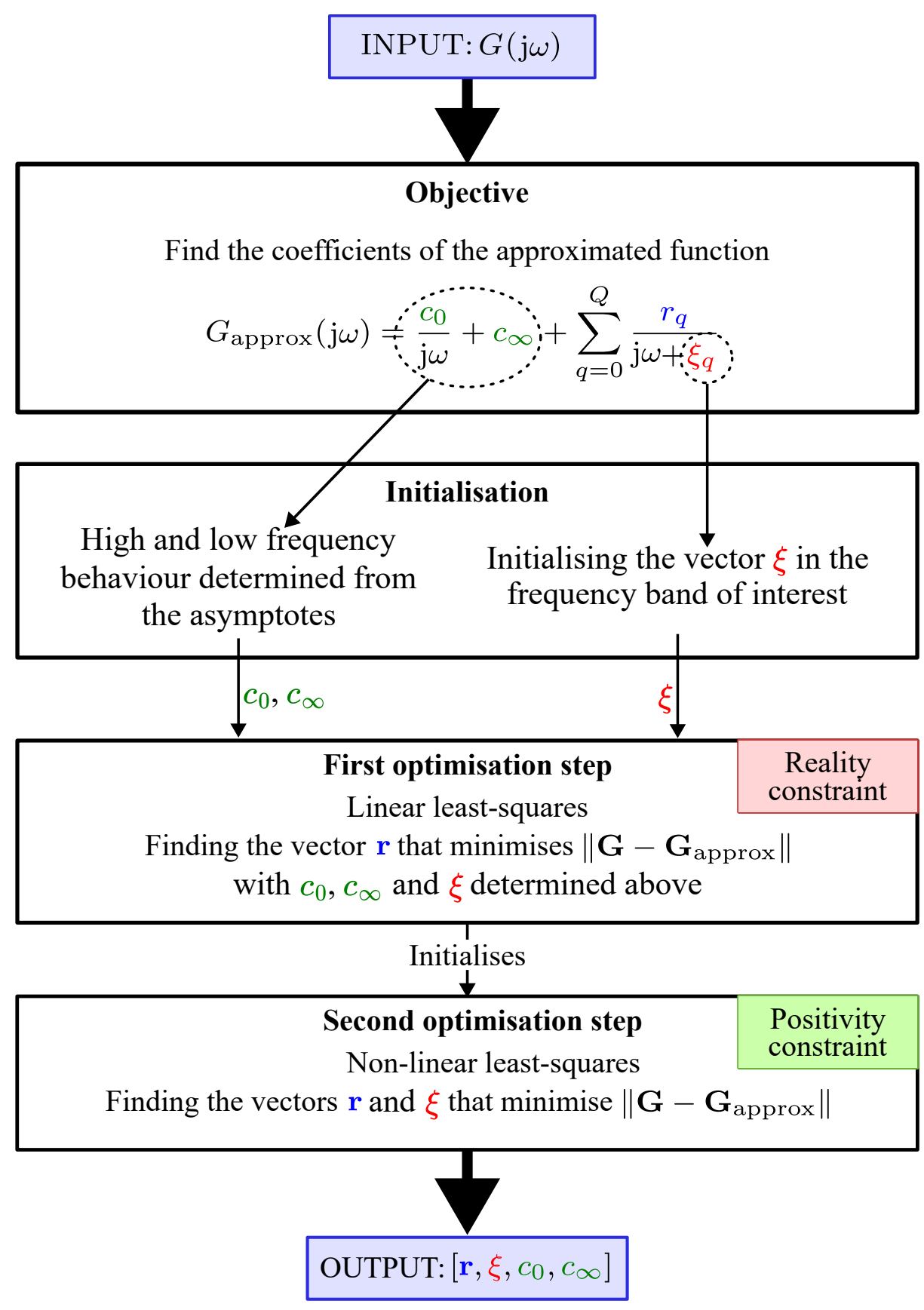

Figure 1: *Identification process (from [2]).*

## 3 Extension to double porosity media

The analysis is now extended to double-porosity media, defined as the Boolean intersection of two solid pore networks. The EFM of this medium is obtained by combining the EFMs of the two constituent single-porosity materials, taking into account of pressure-diffusion effects at the interface between the networks [5]. The resulting $\alpha(s)$ and $\beta(s)$ are:

$$\begin{cases} \alpha_{\text{dp}}(s) & = & \left[\alpha_1(s)^{-1} + (1-\phi_1)\alpha_2(s)^{-1}\right]^{-1}, \\ \beta_{\text{dp}}(s) & = & \beta_1(s) + (1-\phi_1)F_d(s)\beta_2(s), \end{cases}$$

with $\phi_{\text{dp}} = \phi_1 + (1-\phi_1)\phi_2$ the total porosity, and $F_d$ the pressure diffusion function is:

$$F_d(s) = 1 - C(s)\left(C(s) + \sqrt{1 + C(s)\frac{M_d}{2}}\right)^{-1},$$

with $C(s) = sC_0\beta_2(s)$ and $C_0$ , $M_d \in \mathbb{R}^+$.

In the sequel, we make the hypothesis that $F_d$ and both $\beta_{\text{dp}}$ and $\alpha_{\text{dp}}$ are diffusive transfer functions (Eq. 1). Note that some combinations of diffusive functions are proved to remain diffusive (see [2] appendix B for examples).

## 4 Results

The identification algorithm was applied on $\alpha_{\text{dp}}(s)$ and $\beta_{\text{dp}}(s)$ of double porosity media with $Q = 5$ pairs of poles and residues (Fig. 2).

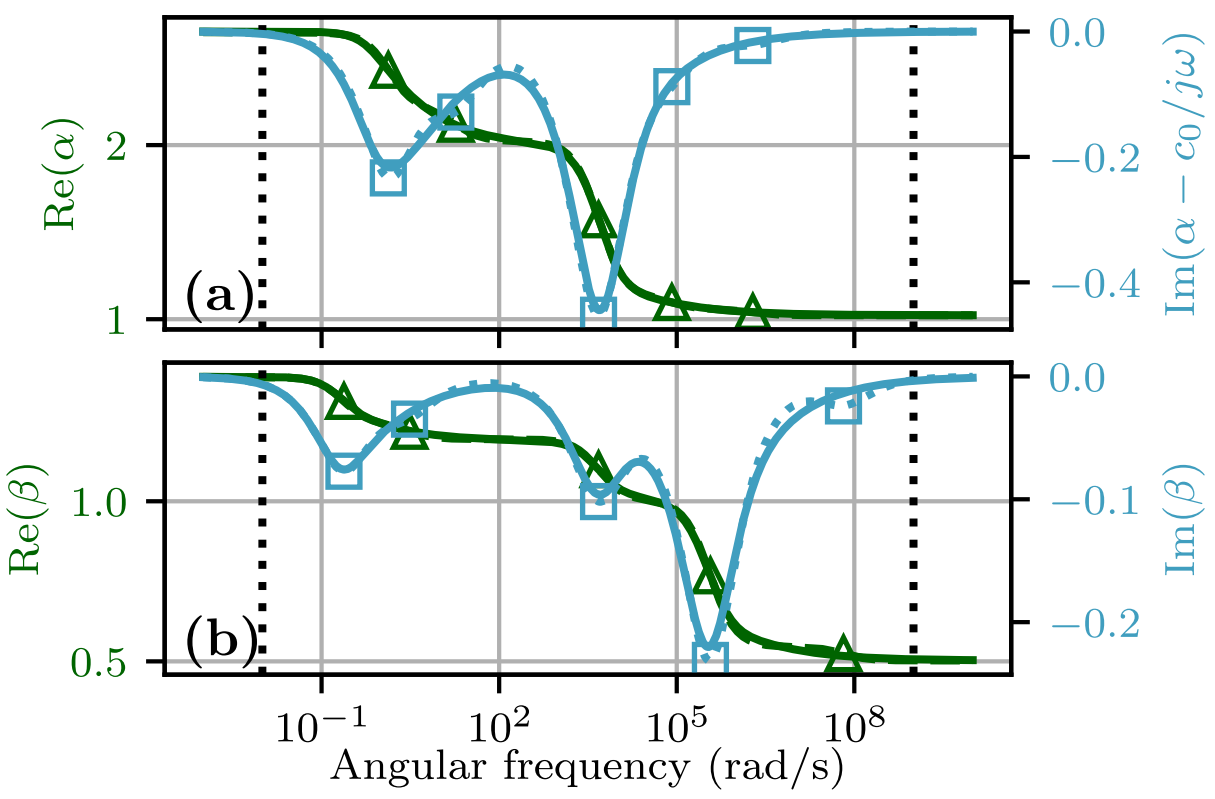

Figure 2: *Theoretical (solid line) and approximated functions (dashed line). The markers abscissae correspond to the $\xi_q$ values (example taken from [2]).*

Then, numerical time-domain simulations were performed with the approximated model to validate the representation (Fig. 3).

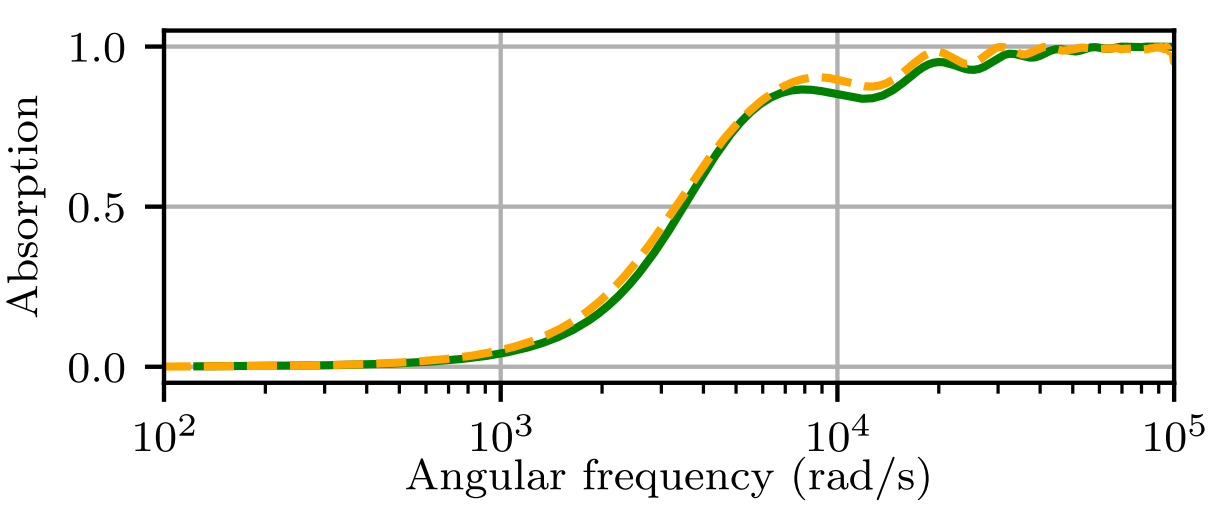

Figure 3: *Absorption coefficient from frequency model (green solid line) and obtained from simulations (orange dashed line) (example taken from [2]).*

# An efficient numerical solver for non-periodic scattering problems in locally perturbed periodic structures

**Nasim Shafieeabyaneh**[1,*], **Tilo Arens**[1], **Ruming Zhang**[2]
[1]Department of Mathematics, Karlsruhe Institute of Technology, Karlsruhe, Germany
[2]Department of Mathematics, Technical University of Berlin, Berlin, Germany
*Email: nasim.shafieeabyaneh@kit.edu

**Abstract**

We study scattering of a time-harmonic acoustic wave by a locally perturbed periodic surface, where the incident field need not be periodic. This scattering problem is challenging due to the unboundedness of the domain and the lack of periodicity. Using a coordinate transformation, we remove the perturbation and obtain an equation with non-constant coefficients in a periodic domain. This allows application of the Floquet–Bloch (FB) transform, reducing the problem to a coupled family (indexed by the Floquet parameter) of periodic problems in a bounded cell. We approximate the solution with the perfectly matched layer (PML) method and prove its exponential convergence on compact sets. Moreover, we show that the FB transform of the solution is analytic with respect to the Floquet parameter. This enables us to accurately discretize the inverse FB transform using only a few values of the Floquet parameter. Finally, by exploiting the system structure via the Schur complement, we propose a fast parallel solver, validated by numerical examples.

## 1 Introduction

We consider $\Gamma := \{(x_1, \zeta(x_1)) : x_1 \in \mathbb{R}\}$, where $\zeta$ is a $2\pi$-periodic Lipschitz continuous function. Moreover, we set $\Gamma_\delta := \{(x_1, \zeta_\delta(x_1)) : x_1 \in \mathbb{R}\}$, where $\zeta_\delta - \zeta$ is compactly supported in $(-\pi, \pi)$. The domains lying above these surfaces are denoted by $\Omega$ and $\Omega_\delta$, respectively (see Fig. 1).

For a given non-periodic incident field $u^i$, we aim to find the scattered field $u^s$ that satisfies

$$\begin{aligned} \Delta u^s + k^2 u^s = 0 \quad &\text{in } \Omega_\delta, \\ u^s = -u^i \text{ on } \Gamma_\delta \end{aligned} \tag{1}$$

with the wave number $k > 0$.

For $\mathrm{H} > \|\zeta\|_\infty$, we define $\Gamma_\mathrm{H} := \mathbb{R} \times \{\mathrm{H}\}$ and denote by $\Omega^\mathrm{H}$ the domain between $\Gamma$ and $\Gamma_\mathrm{H}$.

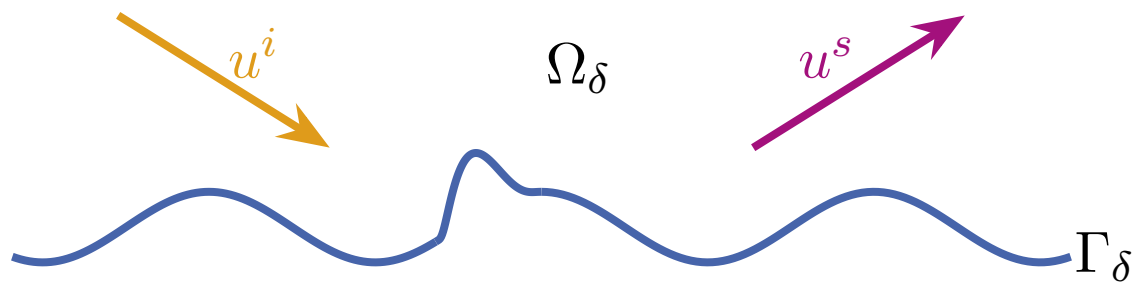

Figure 1: Unbounded domain $\Omega_\delta$

Using the upward propagating radiation condition to $u^s$ leads to $\partial_{x_2} u^s = T^+ u^s$ on $\Gamma_\mathrm{H}$, where the Dirichlet-to-Neumann (DtN) map is defined by $T^+ := \mathcal{F}^{-1}\sqrt{k^2 - \xi^2}\mathcal{F}$ with the Fourier operator $\mathcal{F}$.

To transform the problem into a periodic domain, we define a diffeomorphism $\psi$ from $\Omega$ to $\Omega_\delta$ (see [2, Sec. 4]). Using $u_\mathrm{tra} := u^s \circ \psi$ in (1), it follows that

$$\nabla \cdot (A \nabla u_\mathrm{tra}) + k^2 c\, u_\mathrm{tra} = 0 \quad \text{in } \Omega^\mathrm{H},$$

with $u_\mathrm{tra} = -u^i \circ \psi$ on $\Gamma$, $\partial_{x_2} u_\mathrm{tra} = T^+ u_\mathrm{tra}$ on $\Gamma_\mathrm{H}$, $c := |\det \nabla\psi|$ and $A := c\,(\nabla\psi)^{-1}(\nabla\psi)^{-\top}$.

## 2 The PML Method

The PML method adds a layer above $\Gamma_\mathrm{H}$ to damp outgoing waves (see [1, Sec. 2]). Solving a Dirichlet problem analytically in the PML provides us the PML approximation of the DtN map, $T^+_\rho := \mathcal{F}^{-1}\sqrt{k^2 - \xi^2}\coth(-\mathrm{i}\rho\sqrt{k^2 - \xi^2})\mathcal{F}$, where $\rho \in \mathbb{C}$ denotes the PML damping parameter. The PML approximation of $u_\mathrm{tra}$, denoted by $u_\rho$, satisfies

$$\nabla \cdot (A \nabla u_\rho) + k^2 c u_\rho = 0 \quad \text{in } \Omega^\mathrm{H}, \tag{2}$$

with $u_\rho = -u^i \circ \psi$ on $\Gamma$, $\partial_{x_2} u_\rho = T^+_\rho u_\rho$ on $\Gamma_\mathrm{H}$.

**Theorem 1** *Let $u^i = G(\cdot, \mathbf{y})$, where $G$ is the Dirichlet Green's function for the upper half-space containing $\Omega^\mathrm{H}$ and $\mathbf{y} \in \Omega^\mathrm{H}$. Then, for every compact subset $K \subseteq \Omega^\mathrm{H}$ and sufficiently large $\rho$, there exist $C$, $c > 0$ such that*

$$\|u_\rho - u_\mathrm{tra}\|_{H^1(K)} \le C\,\mathrm{e}^{-c|\rho|}.$$

As Problem (2) is posed in the periodic domain $\Omega^{\mathrm{H}}$, we can apply the FB transform to decompose it into a family of periodic problems in a bounded cell.

## 3 The FB Transform

For $f \in C_0^\infty(\Omega^{\mathrm{H}})$, the FB transform of $f$, denoted by $\mathcal{J}f$, is defined by

$$(\mathcal{J}f)(\alpha, \mathbf{x}) := \sum_{j\in\mathbb{Z}} f(x_1 + 2\pi j, x_2) \mathrm{e}^{-\mathrm{i}\alpha(x_1+2\pi j)}$$

with $(\alpha, \mathbf{x}) \in \mathbb{R} \times \Omega^{\mathrm{H}}$ (see [2, Sec. A]). For fixed $\alpha \in [-1/2, 1/2]$, $\mathcal{J}f$ is $2\pi$-periodic with respect to $x_1$. Accordingly, it suffices to restrict the analysis to the fundamental cell of periodicity $\Omega^{\mathrm{H}}_{2\pi} := \{\mathbf{x} \in \Omega^{\mathrm{H}} : -\pi < x_1 < \pi\}$.

Adding and subtracting $\Delta u_\rho$ and $k^2 u_\rho$ to (2) and then applying the FB transform, we find that $w := (\mathcal{J}u_\rho)(\alpha, \mathbf{x})$ satisfies

$$\begin{aligned} &\Delta w + 2\mathrm{i}\alpha\partial_{x_1} w + (k^2 - \alpha^2)w \\ &\quad + \nabla \cdot ((A - I)\nabla u_\rho) + k^2(c - 1)u_\rho = 0 \text{ in } \Omega^{\mathrm{H}}_{2\pi} \end{aligned} \tag{3}$$

together with $w = -(\mathcal{J}u^i)(\alpha, \psi(\mathbf{x}))$ on $\Gamma^{2\pi}$ and $\partial_{x_2} w = T^+_{\alpha,\rho} w$ on $\Gamma^{2\pi}_{\mathrm{H}}$, where $T^+_{\alpha,\rho}$ is the periodized version of $T^+_\rho$. The PML approximation of the solution is computed by the inverse FB transform, i.e., $u_\rho = \mathcal{J}^{-1} w = \int_{-1/2}^{1/2} w\, \mathrm{e}^{\mathrm{i}\alpha x_1}\, \mathrm{d}\alpha$.

**Theorem 2** *Let $\mathcal{J}u^i$ be analytic with respect to $\alpha$. Then, for sufficiently large $\rho$, the function $w$ that solves (3) is also analytic with respect to $\alpha$.*

To discretize (3), we consider quadrature nodes and weights $\{(\alpha_j, \eta_j)\}_{j=1}^N$ to compute the inverse FB transform. Afterward, for each $\alpha_j$, $w$ is approximated by the piecewise linear functions $\{\varphi_m\}_{m=1}^M$. Defining the vectors of unknowns $W_j = [w_{j,1}, \ldots, w_{j,M}]^\top$ and $U = [u_1, \ldots, u_M]^\top$, we obtain the following linear system

$$\begin{bmatrix} A_1 & 0 & \cdots & 0 & B_1 \\ 0 & \ddots & \ddots & \vdots & \vdots \\ \vdots & \ddots & \ddots & 0 & \vdots \\ 0 & \cdots & 0 & A_N & B_N \\ C_1 & \cdots & \cdots & C_N & I \end{bmatrix} \begin{bmatrix} W_1 \\ \vdots \\ \vdots \\ W_N \\ U \end{bmatrix} = \begin{bmatrix} F_1 \\ \vdots \\ \vdots \\ F_N \\ 0 \end{bmatrix}, \tag{4}$$

where the complex $M \times M$ matrices $A_j$ and $B_j$ correspond to the sesquilinear forms that arise from the variational formulation of (3), the matrix $C_j := -\eta_j \operatorname{diag}(\mathrm{e}^{\mathrm{i}\alpha_j x_{1,1}}, \ldots, \mathrm{e}^{\mathrm{i}\alpha_j x_{M,1}})$ and the vector $F_j(m) := \big(\mathcal{J}u^i\big)(\alpha_j, \psi(\mathbf{x}_m))$ for all $j \in \{1, \ldots, N\}$ and $\mathbf{x}_m = (x_{m,1}, x_{m,2})$ are the nodes of the spatial discretization for all $m \in \{1, \ldots, M\}$.

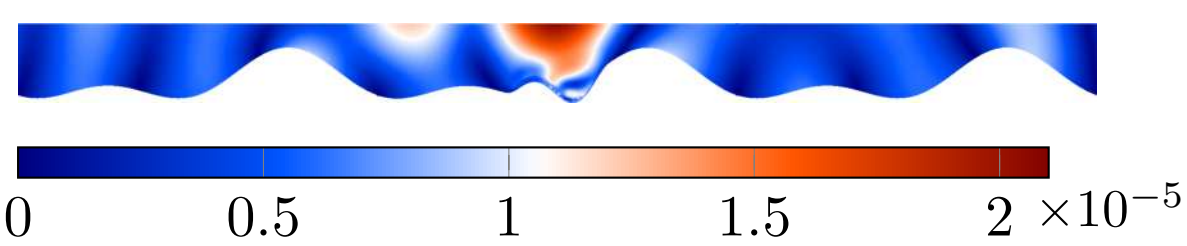

Figure 2: Numerical error for $k = 1.5, |\rho| = 20$

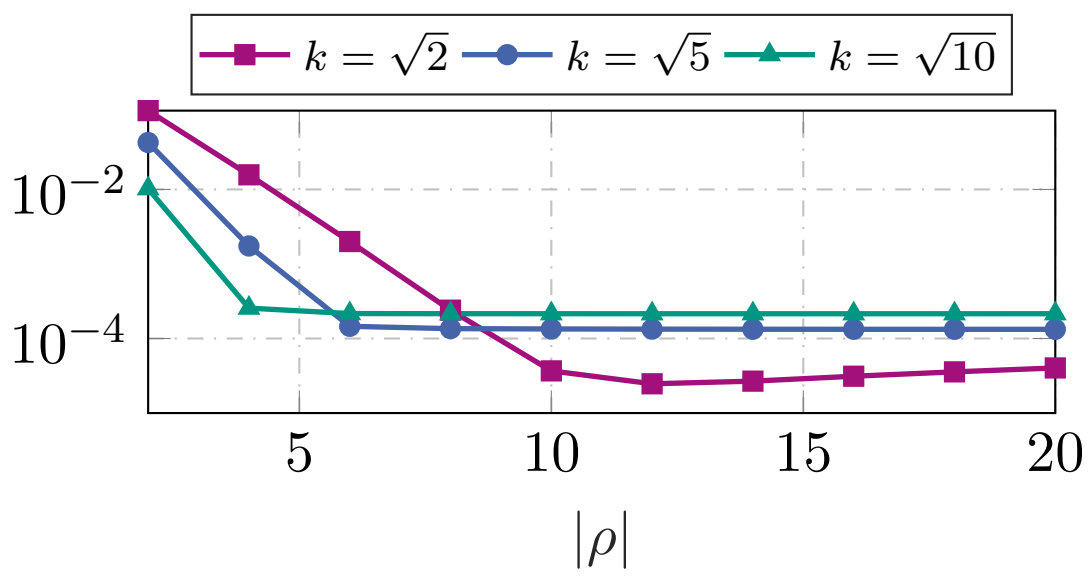

Figure 3: Relative error with respect to $|\rho|$

To propose an efficient method, we first recursively apply the Schur complement to reformulate the above linear system into

$$\left( I + \sum_{\ell=1}^{N} C_\ell A_\ell^{-1} B_\ell \right) U = -\sum_{\ell=1}^{N} C_\ell A_\ell^{-1} F_\ell\,.$$

The summands on the right-hand side and matrix–vector multiplications by the summands on the left-hand side are all independent of each other. Hence, they can be carried out in parallel. As a result, this approach significantly reduces the computational time.

## 4 Numerical Results

To show the efficiency of the method, we consider $u^i = G(\cdot, \mathbf{y})$ with $\mathbf{y}$ below $\Gamma_\delta$.

In Fig. 2, we plot the numerical error in the main bounded cell and two neighboring cells below the PML. In Fig. 3, we illustrate the dependence of the relative error on $|\rho|$. This result shows that the error decays exponentially until the error of the spatial discretization dominates.

# An Unconditionally Stable Energetic Galerkin Scheme for Acoustic Penetrable Objects Using Time-Domain Surface Integral Equations

**Junze Shao**[1,*], **Rui Chen**[2], **Hakan Bagci**[1]
[1]Electrical and Computer Engineering (ECE) Program
Computer, Electrical, and Mathematical Sciences and Engineering (CEMSE) Division
King Abdullah University of Science and Technology (KAUST), Thuwal 23955, Saudi Arabia
[2]School of Microelectronics, Nanjing University of Science and Technology, Nanjing 210094, China
*Email: junze.shao@kaust.edu.sa

## Abstract

An energetic Galerkin weak formulation is developed for acoustic scattering from homogeneous penetrable objects using time-domain surface integral equations. The formulation is obtained by incorporating potential integral equations and their normal derivatives into the global energy identity. The resulting equations are discretized using higher-order Galerkin methods in space and time. Numerical results demonstrate high accuracy and long-time stability.

## 1 Introduction

The energetic Galerkin weak form, proven to possess coercivity and continuity properties, achieves stability when used with conforming Galerkin approximation schemes [1]. This energetic formulation has been successfully applied to interior wave propagation problems with mixed boundary conditions and exterior Neumann problems in two- and three-dimensional settings, demonstrating its robustness, accuracy, and long-time stability.

In this work, the energetic Galerkin formulation is extended to model acoustic wave scattering from penetrable objects. By incorporating time-domain potential surface integral equations (TDPIEs) into the energy identity, an unconditionally stable energetic Galerkin weak form is obtained. The resulting bilinear form is discretized using higher-order Galerkin methods in both space and time [2]. Numerical results, including comparisons with the original TDPIE method [3], demonstrate that the proposed formulation is free from spurious resonances and remains stable over long simulation times.

## 2 Formulation

Consider a penetrable object $\Omega_2$ with mass density $\rho_2$ and wave speed $c_2$, embedded in a background medium $\Omega_1$ with mass density $\rho_1$ and wave speed $c_1$. Under excitation by an incident field $\varphi^{\text{inc}}(\mathbf{r},t)$, scattered fields $\varphi_1^{\text{scat}}(\mathbf{r},t)$ and $\varphi_2^{\text{scat}}(\mathbf{r},t)$ are generated in $\Omega_1$ and $\Omega_2$, respectively. The total fields are given by $\varphi_1(\mathbf{r},t) = \varphi_1^{\text{scat}}(\mathbf{r},t) + \varphi^{\text{inc}}(\mathbf{r},t)$ in $\Omega_1$ and $\varphi_2(\mathbf{r},t) = \varphi_2^{\text{scat}}(\mathbf{r},t)$ in $\Omega_2$. According to the Kirchhoff-Helmholtz representation theorem [4], the fields in $\Omega_1$ and $\Omega_2$ satisfy the following set of TDPIEs on surface $\Gamma = \partial\Omega_2$:

$$\begin{aligned}
&\partial_t\varphi_1(\mathbf{r},t)+\partial_t S_1[\partial_n\varphi_1](\mathbf{r},t) \\
&\qquad -\partial_t D_1[\varphi_1](\mathbf{r},t) = \partial_t\varphi^{\text{inc}}(\mathbf{r},t)\\
&\partial_n\varphi_1(\mathbf{r},t)+D_1'[\partial_n\varphi_1](\mathbf{r},t)\\
&\qquad -N_1[\varphi_1](\mathbf{r},t) = \partial_n\varphi^{\text{inc}}(\mathbf{r},t)\\
&\partial_t\varphi_2(\mathbf{r},t)-\partial_t S_2[\partial_n\varphi_2](\mathbf{r},t)+\partial_t D_2[\varphi_2](\mathbf{r},t) = 0\\
&\partial_n\varphi_2(\mathbf{r},t)-D_2'[\partial_n\varphi_2](\mathbf{r},t)+N_2[\varphi_2](\mathbf{r},t) = 0
\end{aligned} \tag{1}$$

Here, $S_k[\cdot]$, $D_k[\cdot]$, $D_k'[\cdot]$ and $N_k[\cdot]$ denote the time-domain single-layer, double-layer, adjoint double-layer, and hypersingular integral operators associated with medium $\Omega_k$, respectively. The velocity potential satisfies the following transmission conditions on $\Gamma$:

$$\begin{aligned}
\partial_n\varphi_1(\mathbf{r},t) &= \partial_n\varphi_2(\mathbf{r},t)\\
\rho_1\partial_t\varphi_1(\mathbf{r},t) &= \rho_2\partial_t\varphi_2(\mathbf{r},t).
\end{aligned} \tag{2}$$

The total energy of the coupled system at time $T$ is given by the following energy identity:

$$\begin{aligned}
\mathcal{E}_{\text{tot}}(T) = &-\rho_1\int_\Gamma\int_0^T \partial_n\varphi_1(\mathbf{r},t)\partial_t\varphi_1(\mathbf{r},t)\,ds\,dt\\
&+\rho_2\int_\Gamma\int_0^T \partial_n\varphi_2(\mathbf{r},t)\partial_t\varphi_2(\mathbf{r},t)\,ds\,dt.
\end{aligned}$$

Substituting the TDPIEs (1) into the energy identity and enforcing the transmission conditions (2) yields the following energetic Galerkin

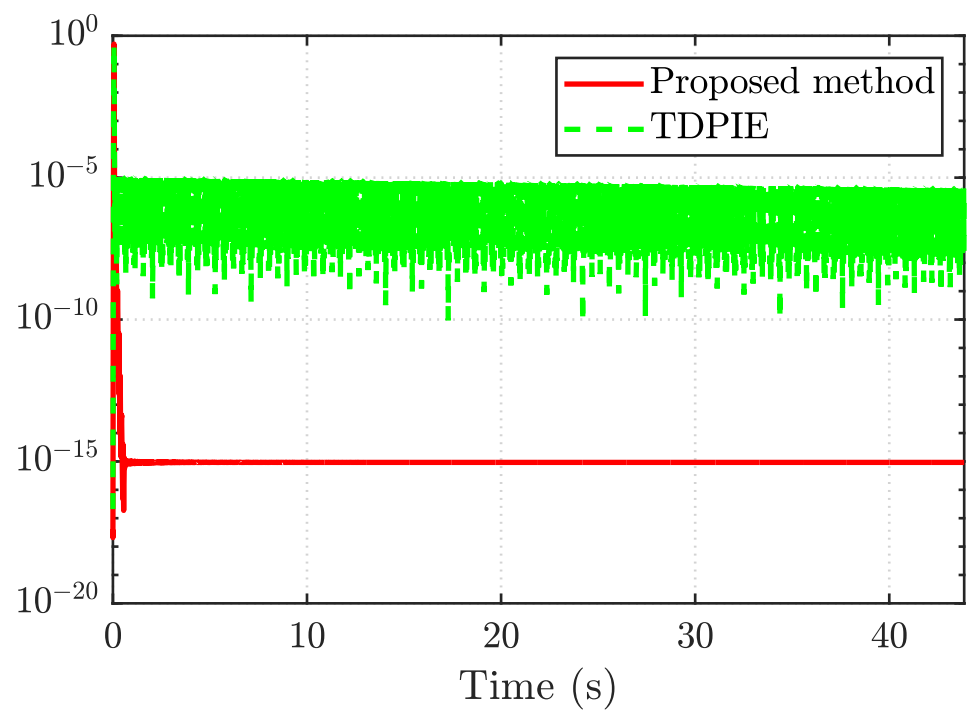

Figure 1: Representative expansion coefficient of $\varphi_1(\mathbf{r}, t)$ computed using the proposed scheme and the original TDPIE formulation [3].

weak formulation:

$$\begin{aligned}&\langle \partial_n \varphi_1, \partial_t S_1[\partial_n \varphi_1] + (\rho_2/\rho_1)\partial_t S_2[\partial_n \varphi_1] \\ &\quad - \partial_t \tilde{D}_1[\varphi_1] - \partial_t \tilde{D}_2[\varphi_1] \rangle = \langle \partial_n \varphi_1, \partial_t \varphi^{\mathrm{inc}} \rangle \\ &\langle \partial_t \varphi_1, \tilde{D}'_1[\partial_n \varphi_1] + \tilde{D}'_2[\partial_n \varphi_1] \\ &\quad - N_1[\varphi_1] - (\rho_1/\rho_2) N_2[\varphi_1] \rangle = \langle \partial_t \varphi_1, \partial_n \varphi^{\mathrm{inc}} \rangle\end{aligned}$$

Here, '$\sim$' on the surface integral operators $D_k[\cdot]$ and $D'_k[\cdot]$ indicates that the integrals are evaluated in the Cauchy principal value sense, and the space–time inner product over $\Gamma \times [0, T]$ is defined as

$$\langle a, b \rangle = \int_\Gamma \int_0^T a(\mathbf{r}, t) b(\mathbf{r}, t) \, dt \, ds.$$

The resulting bilinear form is discretized using higher-order Galerkin methods in both space and time [2], yielding a causal discrete system that is solved using a marching-on-in-time (MOT) procedure.

## 3 Numerical Results

A sphere of radius $r = 1\,\mathrm{m}$ with $\rho_2 = 0.1\,\mathrm{kg\,m^{-3}}$ and $c_2 = 1500\,\mathrm{m\,s^{-1}}$, embedded in air ($\rho_1 = 1.0\,\mathrm{kg\,m^{-3}}$ and $c_1 = 343\,\mathrm{m\,s^{-1}}$) is considered. The sphere is excited by a Gaussian-modulated plane wave with central frequency 120 Hz and bandwidth 70 Hz. The surface $\Gamma$ is discretized using curvilinear patches with average edge length 0.15 m. First-order spatial and second-order temporal Galerkin basis functions [2] are used, resulting in 1 130 unknowns. The simulation is performed for 100 000 time steps with step size $4.368 \times 10^{-4}\,\mathrm{s}$. Fig. 1 compares a representative coefficient of $\varphi_1(\mathbf{r}, t)$ obtained using the proposed scheme and the original TDPIE formulation [3], demonstrating absence of spurious resonances and long-time stability. Fig. 2 shows the eigenvalues of the companion matrix of the MOT system, all of which lie within the unit circle, confirming late-time stability. The accuracy of the proposed solver is validated by comparing scattering cross section results with the analytical Mie series solution. A relative $L_2$ norm error of $2.5178 \times 10^{-5}$ is achieved at the central frequency.

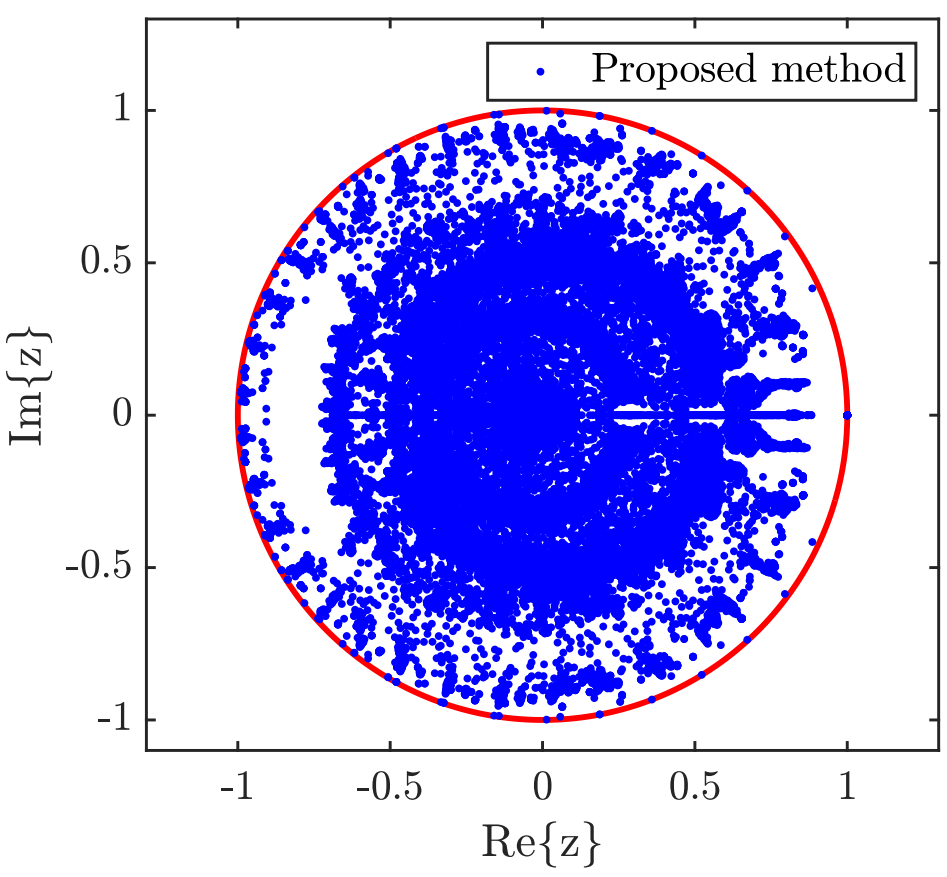

Figure 2: Eigenvalues of the companion matrix of the MOT system.

# Spectrally characterizing subsurface targets with GPR-SAR

**Arnold D. Kim**[1], **Symeon Papadimitropoulos**[1], **Joseph Simpson**[1,*], **Chrysoula Tsogka**[1]
[1]Department of Applied Mathematics, University of California, Merced, Merced, USA
*Email: jsimpson9@ucmerced.edu

## Abstract

When target sizes are comparable to the resolution of a synthetic aperture radar (SAR) imaging system, images only produce representative points which can be used to detect and locate different targets but not distinguish any further between them. Using the representative points reconstructed by Kirchhoff migration (KM), we can recover sets of frequency spectra from the same measurements. These spectra are related to the radar cross-section (RCS) of the target and provide valuable information allowing for differentiating targets that would otherwise be indistinguishable [1]. We apply this approach to image subsurface targets with ground-penetrating radar SAR (GPR-SAR) with the aim of reducing false positives in landmine detection.

## 1 Problem Description

A diagram of the setup for GPR-SAR with an unmanned aerial vehicle (UAV) is shown in Fig. 1. As a UAV travels along a prescribed flight path, at locations $\boldsymbol{x}_n$ for $n = 1, \dots, N$, which make up the synthetic aperture, it emits a multifrequency signal at frequencies $f_m$ for $m = 1, \dots, M$ that propagates downward. Some of this signal is reflected by the air-soil interface (ground bounce) and the rest transmits into the soil and scatters off of subsurface targets. The receiver on the UAV measures these signals. The measurements are stored in the matrix $D \in \mathbb{C}^{M\times N}$ whose entry $d_{mn}$ is the signal measured at $\boldsymbol{x}_n$ for $f_m$.

Due to the computational cost of simulating 3D rough surface scattering, we restrict our attention here to the 2D problem, which contains the essential features needed to test and validate our imaging method such as clutter due to the rough air-soil interface and additive measurement noise. For a Gaussian pulse, the recorded scattered signals are $d_{mn} = u_1(f_m;\boldsymbol{x}_n) - u_2(f_m;\boldsymbol{x}_n)$,

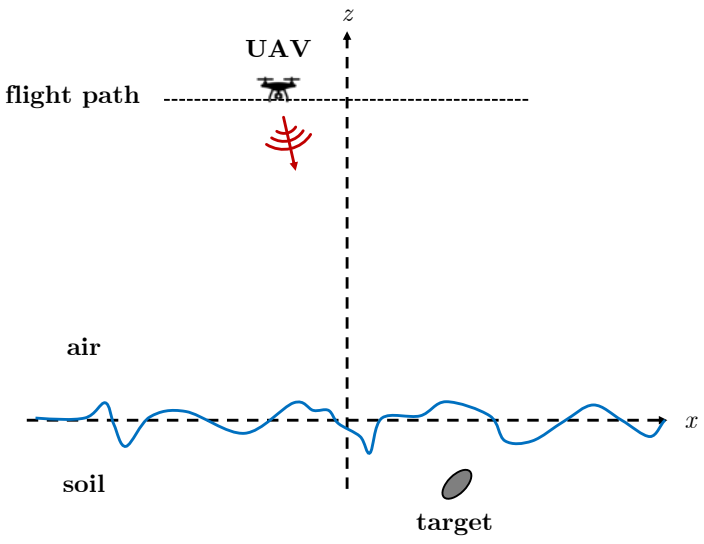

Figure 1: A sketch of the UAV-based GPR-SAR imaging system.

where $u_1$ and $u_2$ satisfy

$$\nabla \cdot (\epsilon_i^{-1} \nabla u_i) + k_m^2 u_i = -\frac{\alpha}{\pi} e^{-\alpha|\boldsymbol{x}-\boldsymbol{x}_n|^2} \tag{1}$$

in the weak sense for $i = 1, 2$ where $k_m = 2\pi f_m/c$, $c = 2.997 \times 10^8$ m/s, $\epsilon_1 = 1$ above the air-soil interface, $\epsilon_1 = \epsilon_{\text{soil}}$ below the air-soil interface outside of targets, $\epsilon_1 = \epsilon_{\text{tar},q}$ inside target $q$, and $\epsilon_2 = 1$ everywhere. We model a rough air-soil interface using a mean-zero Gaussian correlated random height function parameterized by its correlation length and root-mean-square (rms) height [2]. We solve (1) with appropriate outgoing conditions in a truncated domain using perfectly matched layers and Finite Element Methods (FEM).

## 2 Imaging and Spectral Recovery

For $Q$ scatterers indexed by $q = 1, \dots, Q$, we can decompose the data matrix $D$ into

$$D = R + \sum_{q=1}^{Q} F^{(q)} \odot S^{(q)} + E \tag{2}$$

where $R$ contains the ground-bounce signals, $F^{(q)}$ is a frequency- and direction-dependent reflectivity matrix, $S^{(q)}$ is the signal from a subsurface point target located at the center of the target, and $E$ contains additive measurement noise and other small modeling errors not accounted for in the other matrices.

Here $\|R\|_F \gg \|F^{(q)} \odot S^{(q)}\|_F$. To remove the dominant contribution from $R$, we apply PCA and project off the first $j^* - 1$ rank-1 matrices

from the singular value decomposition of $D$. We choose $j^*$ empirically to remove as much ground-bounce signal as possible without also removing target signals. We apply Kirchhoff migration (KM) to the projected data matrix $\tilde{D}$ using the phase of the leading asymptotic behavior for a single-scattering point-scatterer located below a flat air-soil interface. These KM images give us representative points for the targets.

We can also use the decomposition in (2) together with the projection matrix from PCA to setup linear systems to approximately recover the entries of $F^{(q)}$ using the peaks from the KM images in a model for the entries of $S^{(q)}$. These systems will be underdetermined due to the projection, which we address with a piecewise constant approximation for the spatial component of $F^{(q)}$, effectively breaking the synthetic aperture into sub-apertures and treating the reflectivities as constant in space over each sub-aperture. Furthermore, the entries of $F^{(q)}$ can be interpreted as a scaled version of the scattering amplitude measured in the direction of backscattering in the soil for target $q$, $f^{\text{sca},(q)}_{\text{B,soil}}$, which can be computed using boundary integral equation (BIE) methods [3]. Without knowledge of the locations of target centers, the phase of $f^{\text{sca},(q)}_{\text{B,soil}}$ is unknown, so we compute the RCS, which is proportional to $|f^{\text{sca},(q)}_{\text{B,soil}}|^2$. Since these RCS spectra depend on target shape, size, and material properties, they open opportunities to distinguish false positives from true targets.

## 3 Results

For simulation results we use 41 frequencies uniformly sampling the band from 3.1 to 5.1 GHz with 80 spatial measurements uniformly sampling a linear synthetic aperture of length 102 cm at an average elevation of 75 cm. We use $\epsilon_{\text{soil}} = 9$, and $\alpha = 10$ for the pulse. We include additive measurement noise so that SNR≈20 dB. The rough air-soil interface has correlation length 8 cm and rms height 0.5 cm. For PCA, we use $j^* = 3$. For the RCS we use 10 sub-apertures.

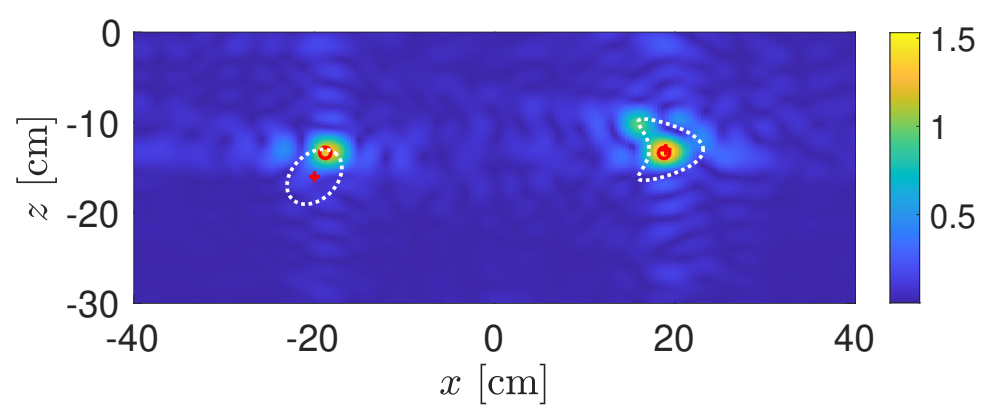

Figure 2: KM image after applying PCA.

Fig. 2 shows a plot of $|I^{\text{KM}}|$ after PCA for a kite and rotated ellipse, both with $\epsilon_{\text{tar}} = 2.3$. Target boundaries are dashed white lines, centers are red '+' symbols, and recovered locations are red 'o' symbols. The image only exhibits peaks at representative points, preventing further identification of targets.

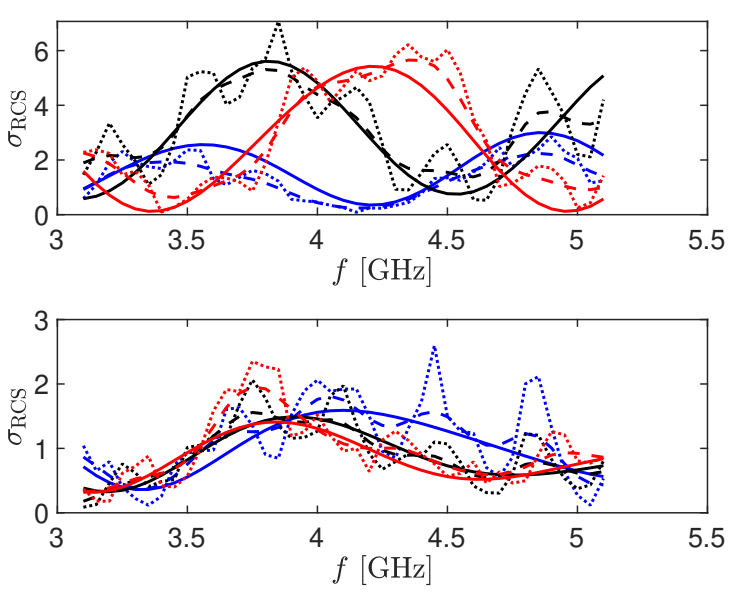

Figure 3: Sample RCS spectra for the kite (top) and ellipse (bottom).

Fig. 3 shows plots of sample recovered RCS curves for the kite (top) and ellipse (bottom). The blue curves are the second sub-aperture, red the fifth sub-aperture, and black the ninth sub-aperture. Dotted lines are recovered RCS curves, dashed lines are smoothed recovered RCS curves, and solid lines are predicted RCS curves computed with BIE. We can see clear differences between the RCS spectra of the two targets. Furthermore, we see relatively close agreement between the recovered curves and the predicted curves, potentially providing opportunities for future classification schemes.

## Acknowledgements

This work was supported by the AFOSR Grant FA9550-24-1-0191 and NSF Grant DMS-1840265.

# Recovery-Based Error Indicator for Finite Difference Methods

**Ferhat Sindy**[1,*], **Annalisa Buffa**[1], **Marco Picasso**[1]

[1] École Polytechnique Fédérale de Lausanne (EPFL), Institute of Mathematics, CH-1015 Lausanne, Switzerland

*Email: ferhat.sindy@epfl.ch

**Abstract**

A recovery-based error indicator for high-order Finite Difference Methods, based on post-processing of the Finite Difference values, is presented. The performance and accuracy of the proposed error indicator are demonstrated through several numerical experiments.

## 1 Introduction

The theory and practice of a posteriori error estimation are well established for the Finite Element Method (FEM), see for instance [1]. In contrast, error estimators for high-order Finite Difference Method (FDM) remain scarce. Several works have investigated a posteriori error estimation for FDMs, often by reformulating the FDM as an equivalent FEM. This reformulation enables the use of standard residual-based estimators. We propose a novel recovery-based indicator that employs a finite element reconstruction of the numerical solution and its gradient, and uses the Polynomial Preserving Recovery (PPR) method [2] to estimate the gradient error of the reconstructed solution. The primary focus of this work is on error evaluation and solution assessment.

## 2 Recovery-Based Error Indicator

We consider the PPR recovery-based error indicator [2]. We let $\Omega \subset \mathbb{R}^d$ be an open, $d$-dimensional polygonal domain. Given $h > 0$, let $\mathcal{T}_h$ be a partition of $\Omega$ into elements $K$ with diameter less or equal $h$. We define the continuous piecewise polynomial Finite Element space of degree r:

$$\mathcal{V}_h^r = \left\{ v_h \in C^0(\overline{\Omega}) : v_h|_K \in \mathbb{P}_r(K), \quad \forall K \in \mathcal{T}_h \right\}.$$

A basis of $\mathcal{V}_h^r$ is the standard Lagrange basis functions $\{\phi_i\}_{i=1}^N$, where $N$ is the number of degrees of freedom. Let $u_h \in \mathcal{V}_h^r$ be an approximation of $u$ obtained, for instance, by solving a linear elliptic problem with the FEM. The PPR gradient recovery operator $\mathrm{G}_h$ is a linear mapping from $\mathcal{V}_h^r$ to $(\mathcal{V}_h^r)^d$, mimicking the $\nabla u$, and defined on piecewise polynomials of degree $r$. Thus:

$$\mathrm{G}_h u_h(x) = \sum_{i=1}^N (\mathrm{G}_h u_h)(z_i)\phi_i(x), \quad x \in \overline{\Omega},$$

where $z_i$ is the nodal coordinate and $\phi_i(x)$ the Lagrange basis function at nodal coordinate $z_i$. Note that $\mathrm{G}_h u_h$ is fully defined by its values at the nodes, which are determined through constructing local polynomials in a discrete least squares sense on patches [2]. We define the local error indicators $\{\eta_K : K \in \mathcal{T}_h\}$ and the global indicator $\eta_h$ as:

$$\eta_K := \|\mathrm{G}_\mathrm{h} u_h - \nabla u_h\|_{L^2(K)} \text{ and } \eta_h := \left( \sum_{K \in \mathcal{T}_h} \eta_K^2 \right)^{\frac{1}{2}}.$$

It turns out that $\eta_h$ is a good approximation of $\|\nabla(u - u_h)\|_{L^2(\Omega)}$ if $\mathrm{G}_\mathrm{h} u_h$ is a better ("super converged") approximation of $\nabla u$ in $\Omega$.

## 3 Example: The Wave Equation

We define $\Omega = (0,1)^2 \times (0,T)$, with $d = 3$. We consider the 2D wave equation:

$$\begin{aligned} &\frac{\partial^2 u}{\partial t^2}(\mathbf{x},t) - \Delta u(\mathbf{x},t) = 0, && (\mathbf{x},t) \in \Omega, \\ &u(\mathbf{x},t) = 0, && (\mathbf{x},t) \in \partial(0,1)^2 \times [0,T], \\ &u(\mathbf{x},0) = f(\mathbf{x}), && \mathbf{x} \in (0,1)^2, \\ &\frac{\partial u}{\partial t}(\mathbf{x},0) = g(\mathbf{x}), && \mathbf{x} \in (0,1)^2. \end{aligned} \tag{1}$$

where $f, g$ are sufficiently smooth and $\mathbf{x} = (x_1, x_2)$. Given positive integers $N_x$ and $N_t$, we define the grid spacing $\frac{1}{N_x+1}$ and time stepsize $\frac{T}{N_t+1}$. We have a uniform Cartesian grid in space and time consisting of FD grid points $(x_{1,i}, x_{2,j}, t^n)$, where $0 \leq i,j \leq N_x + 1$, and $0 \leq n \leq N_t + 1$. We solve problem (1) with a $r+1$ order

FDM to obtain a FD solution $u_{ij}^n$ approximating $u(x_{1,i}, x_{2,j}, t^n)$. We assume that $N_x + 1$ and $N_t + 1$ are multiples of $r$ and construct an "interpolant mesh" in space and time, which consists of cuboids, $Q_{ij}^n$, with sides $\frac{r}{N_x+1}$, $\frac{r}{N_x+1}$ and $\frac{rT}{N_t+1}$ and forms a partition of $\Omega$. We define the mesh parameter of the interpolant mesh to be $h := \max\left\{\frac{r}{N_x+1}, \frac{rT}{N_t+1}\right\}$. The mesh nodes associated with the interpolant mesh coincide with the FD grid points, allowing us to define for $(\mathbf{x}, t) \in \Omega$:

$$u_h(\mathbf{x}, t) = \sum_{i=0}^{N_x+1} \sum_{j=0}^{N_x+1} \sum_{n=0}^{N_t+1} u_{ij}^n \phi_i(x_1)\phi_j(x_2)\phi^n(t).$$

We denote $e_h = u - u_h$ and we indicate the following error:

$$\overline{|||e_h|||}^2 := \overline{\left\|\frac{\partial e_h}{\partial t}\right\|}^2_{L^2((0,1)^2)} + \overline{\|\nabla e_h\|}^2_{L^2((0,1)^2)},$$

where the overline denotes the root mean square in time over the interval $[T - r\Delta t, T]$, by

$$\overline{\eta_h}^2 := \sum_{i=0}^{\frac{N_x+1}{r}-1} \sum_{j=0}^{\frac{N_x+1}{r}-1} \overline{\eta_{Q_{ij}}^t}^2 + \sum_{i=0}^{\frac{N_x+1}{r}-1} \sum_{j=0}^{\frac{N_x+1}{r}-1} \overline{\eta_{Q_{ij}}^x}^2,$$

where $\eta_{Q_{ij}}^x := \|\mathrm{G}_h^x u_h - \nabla u_h\|_{L^2(Q_{ij})}$, $\eta_{Q_{ij}}^t := \left\|\mathrm{G}_h^t u_h - \frac{\partial u_h}{\partial t}\right\|_{L^2(Q_{ij})}$, $Q_{ij}^n$ are squares with sides $\frac{r}{N_x+1}$ and form a partition of $(0,1)^2$, $\mathrm{G}_\mathrm{h}^x u_h$ is a better approximation of $\nabla u$ in $\Omega$ and $\mathrm{G}_\mathrm{h}^t u_h$ is a better approximation of $\frac{\partial u}{\partial t}$ in $\Omega$.

## 4 Numerical Experiments

We consider the 2D acoustic wave equation in a discontinuous, heterogeneous medium and the 3D ($\Omega = (0,1)^3 \times (0,T)$) acoustic wave equation in a homogeneous medium, with a Gaussian pulse as an initial condition. Local error features are accurately captured by the indicators, as illustrated in Figures 1 and 2.

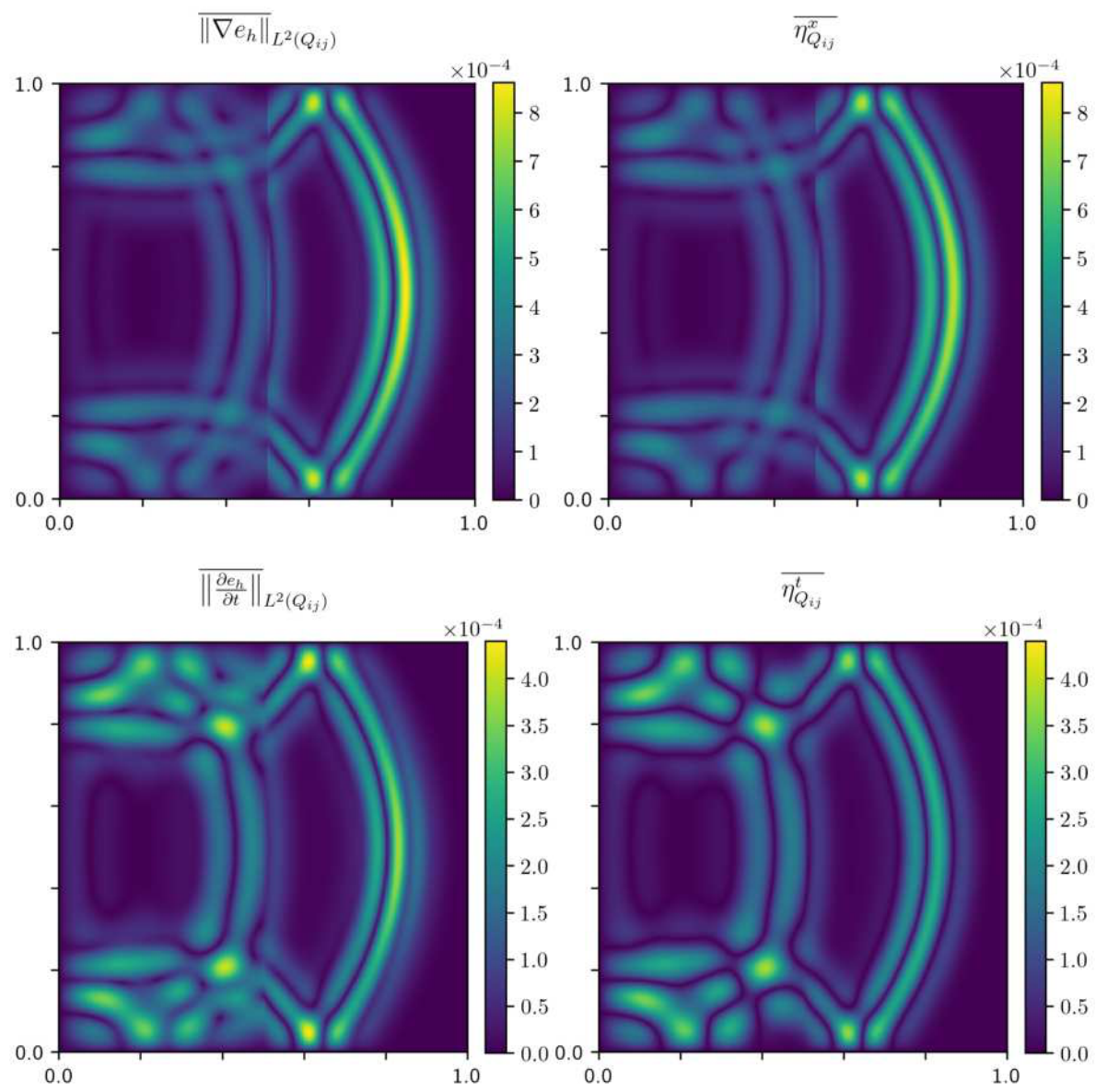

Figure 1: Left, the local error and right, the local indicator.

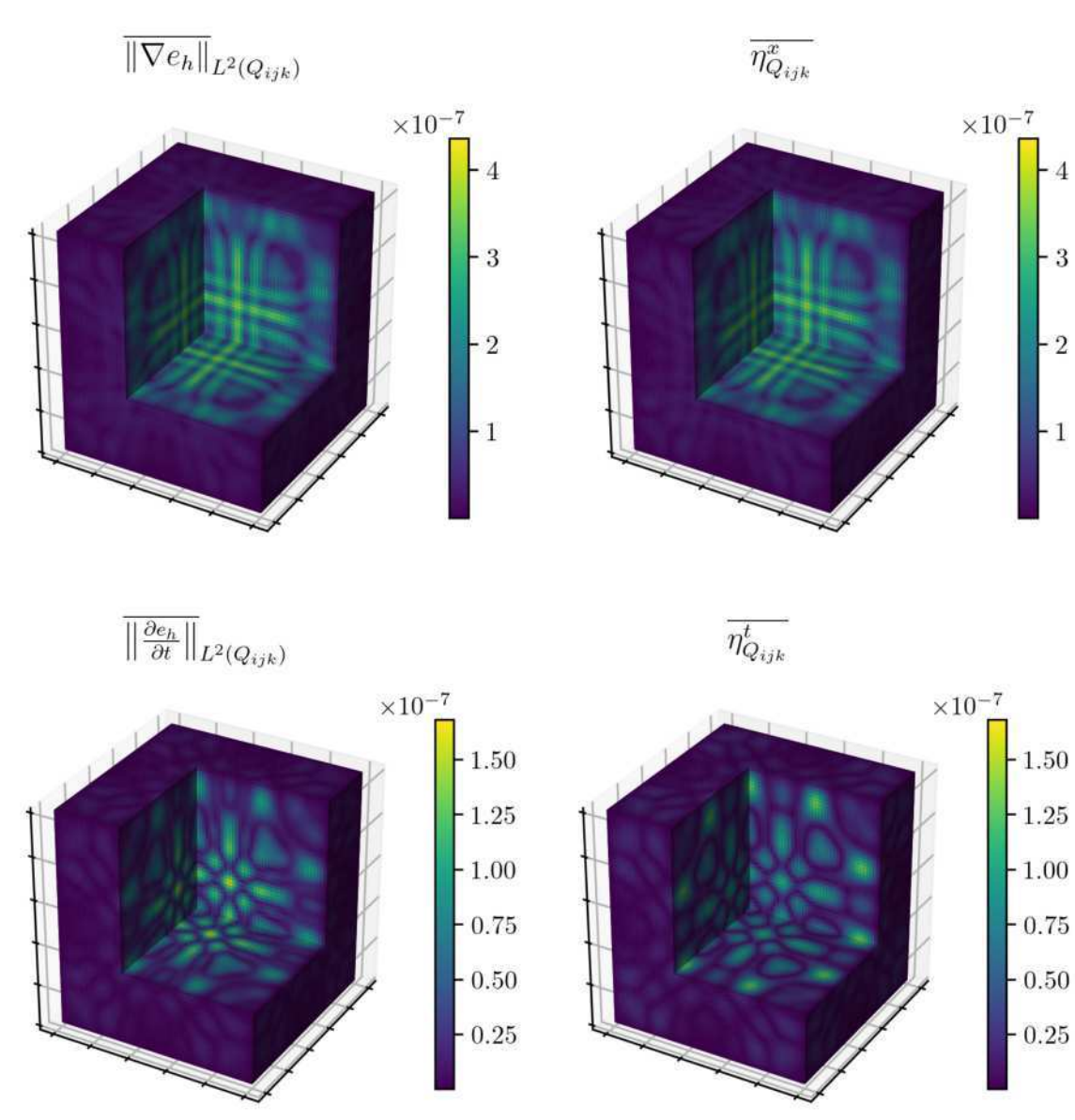

Figure 2: Left, the local error and right, the local indicator.

# Seismic Reflection of Three-Dimensional Plane Waves in Functionally Graded Triclinic Media with Interfacial Imperfections

Akanksha Srivastava[1,*], Federico J Sabina[1]
[1]Instituto de Investigaciones en Matemáticas Aplicadas y en Sistemas, Universidad Nacional Autónoma de México, Ciudad de México, México
*Email: mathsakanksha@gmail.com

**Abstract**

Understanding seismic wave behavior in complex media is crucial for both geophysical exploration and the development of advanced engineered materials. This study investigates the reflection characteristics of three-dimensional plane waves propagating through a functionally graded anisotropic medium bounded by an imperfect contact. By integrating concepts from reflection seismology and anisotropic elastodynamics, we derive analytical expressions for amplitude ratios and energy partitioning as functions of the incident angle using the generalized Snell's law. Additionally, we compute slowness surfaces to visualize wavefront geometry in anisotropic graded media. The influence of material anisotropy, functional grading, and interfacial imperfections on amplitude–versus–angle (AVA) behavior is examined through detailed graphical demonstrations. This work contributes a rigorous analytical framework for interpreting seismic responses in graded anisotropic medium and offers insights relevant to both seismic exploration and in advancing our comprehension of the Earth's interior structure.

## 1 Introduction

The existence of reflection and transmission phenomenona are proven and considered useful for exploration of the earth's internal structure as well as exploration of valuable materials buried inside the earth, for example, metals, minerals, petroleum and hydrocarbons. For this purpose, AVO analysis is useful, and this technique is based on the relationship between amplitude ratios and angle of incidence. This analysis also aids in estimating P and S-wave velocities as well as it infers the petro physical properties of the formation, which is useful in determining the source location. Ruger [1] studied the reflection coefficients and azimuthal AVO analysis of anisotropic modelling. Consideration of the effects of anisotropy on reflection and transmission coefficients, especially the functional grading and imperfect boundary, is of utmost importance.

Functionally graded materials (FGM) are composite materials which are designed to present a particular spatial variation of their properties [Somiya [2]]. This is usually achieved by forming a compound of two components whose volume fraction is changed across a certain direction. These materials have multifunctional properties in which the properties may vary gradually with position or direction [Mandal and Acikbas [3]]. In general, the inclusion of functionally graded and anisotropy are two attributes that introduce complexity into problem as well as help to simulate a realistic elastodynamic model. A very hard layer present inside the earth can be modelle as rigid boundary [Kumar et al. [4]]. The earth's composition is sometimes inhomogeneous, including a particularly hard stratum. The medium functionally graded anisotropic medium and rigid boundary play significant role in the analysis of reflection phenomena. The development of analytical study for the different types of earth's structure in three-dimensional model has been the subject of a continuous effort in the last few decades, e.g., Chatterjee et al. [5], Chattopadhyay et al. [6].

## 2 Model

This study considers a three-dimensional plane wave propagating in functionally graded triclinic crystalline semi-infinite medium occupying the region, $M$ $(y < 0)$. The depth of the semi-infinite medium is represented by the $y$-axis, while the imperfect boundary condition is acquired along the $x$-axis. The incident wave encounters the boundary with an oblique angle, defined by $\theta_0$, as depicted in the schematic rep-resentation of the problem illustrated in Fig. 1.

The stress-strain relationships for triclinic crys-

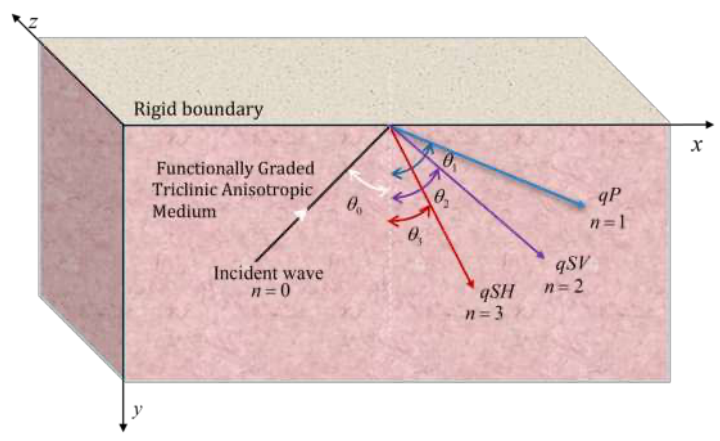

Figure 1: Diagrammatic representation of the present problem.

talline medium are governed by the generalized Hooke's law, represented as follows:

$$\underbrace{\begin{bmatrix}\sigma_1\\ \sigma_2\\ \vdots\\ \sigma_6\end{bmatrix}}_{\sigma} = \begin{bmatrix}c_{11} & c_{12} & \cdots & c_{16}\\ c_{12} & c_{22} & \cdots & c_{26}\\ \vdots & \vdots & \ddots & \vdots\\ c_{16} & c_{26} & \cdots & c_{66}\end{bmatrix} \underbrace{\begin{bmatrix}\epsilon_1\\ \epsilon_2\\ \vdots\\ \epsilon_6\end{bmatrix}}_{\epsilon}, \quad (1)$$

where $\sigma = [\sigma_{11} \;\; \sigma_{22} \;\; \sigma_{33} \;\; \sigma_{23} \;\; \sigma_{13} \;\; \sigma_{12}]^T$, $\epsilon = [\epsilon_{11} \;\; \epsilon_{22} \;\; \epsilon_{33} \;\; 2\epsilon_{23} \;\; 2\epsilon_{13} \;\; 2\epsilon_{12}]^T$, $c_{ij}$ are elastic stiffness coefficients, and $\epsilon_{ij}$ represent strain tensor. Due to the absence of planes of symmetry, triclinic materials are categorized as anisotropic elastic materials with 21 independent elastic constants.

The functional gradient property holds in the FGTM substrate, therefore the elastic coefficients and density are assumed to have exponential variations along with $y$ direction, which are expressed as follows:

$$c_{ij} = M_{ij}e^{\gamma y}, \;\; \rho = \rho^{'} e^{\gamma y}. \quad (2)$$

where $\gamma$ is the exponential coefficient indicating the profile of the material gradient along the x-axis, and the quantities $(M_{ij}, \rho^{'})$ are attached to indicate the $x$-direction factors in the material coefficients.

## 3 Methods and Results

The constitutive equations for functionally graded triclinic medium are employed to derive the governing equations of motion and the associated Christoffel matrices for aformentioned semi-infinite media. The dispersion relations are obtained by setting the determinant of the Christoffel matrix to zero, leading to distinct cubic characteristic equation. This cubic equation is solved analytically using Cardano's method to obtain explicit expressions for the phase velocities of quasi-longitudinal (qP), quasi-shear vertical (qSV), and quasi-shear horizontal (qSH) waves. The derived phase velocities depend explicitly on the propagation vector, reflecting the anisotropic nature of the materials.

The generalized Snell's law is then applied to determine the amplitude ratios of reflected waves for the considered model. The results show that, unlike isotropic media, the functionally graded triclinic model permits mode conversion.

# Exploring Explicit Formulas and Multi-solitons for Integrable Systems

**Ruoci Sun**[1,*]

[1]School of Mathematical Sciences, Laboratory of Mathematics and Complex Systems, MOE, Beijing Normal University, 100875 Beijing, China.

*Email: ruoci.sun.16@normalesup.org

**Abstract**

An integrable system is a mathematical model that can be solved exactly. Beyond their pivotal role in the study of nonlinear partial differential equations, these systems offer a plethora of applications and powerful tools across diverse real-world fields such as fluid mechanics, statistical physics, quantum mechanics, etc. In this presentation, I will introduce several methods aimed at establishing explicit formulas for general solutions of some integrable systems, using the Lax pair structure, spectral theory, Fourier analysis, and tools from Hardy spaces. These explicit formulas not only enable the extension of the flow map, but also offer deeper insights into the long time behavior and perturbations of these equations. Furthermore, I will highlight the multi-soliton dynamics arising from these explicit representations, illustrating the rich wave phenomena captured by such systems.

## 1 Introduction

An evolution equation is *integrable* if its general solutions can be expressed explicitly in terms of its initial datum and the time variable. Let $L$ be some $\mathbb{C}$-linear operator. The Linear equation

$$\partial_t w(t) + Lw(t) = f(t), \quad w(t_0) = \varphi \tag{1}$$

is usually integrable. The explicit formula of its general solutions can be established by using the following Duhamel's formula

$$w(t) = e^{-(t-t_0)L}\varphi + \int_{t_0}^{t} e^{-(s-t_0)L} f(s)\mathrm{d}s. \tag{2}$$

The goal of my presentation is to generalize the idea of Duhamel formula (2) to some nonlinear evolution equation via the Lax pair structure and tools of Hardy spaces.

## 2 Model and Motivation

Let $\mathbb{M}$ be either $\mathbb{R}$ or $\mathbb{T} := \mathbb{R}/2\pi\mathbb{Z}$. The intertwined derivative nonlinear Schrödinger system on $\mathbb{M}$ reads as

$$\begin{cases} \partial_t u = i\partial_x^2 u + u\partial_x \Pi(\overline{v}u) + v\partial_x\Pi(|u|^2), \\ \partial_t v = i\partial_x^2 v + v\partial_x \Pi(\overline{u}v) + u\partial_x\Pi(|v|^2), \end{cases} \tag{3}$$

where $(t,x) \in \mathbb{R}\times\mathbb{M}$ and $\Pi : L^2(\mathbb{M};\mathbb{C}) \to L^2(\mathbb{M};\mathbb{C})$ denotes the Szegő projector that cancels all negative Fourier modes and preserves all nonnegative Fourier modes, i.e.

$$\begin{cases} \widehat{\Pi f}(\xi_1) = \hat{f}(\xi_1), & \forall \xi_1 \geq 0; \\ \widehat{\Pi f}(\xi_2) = 0, & \forall \xi_2 < 0. \end{cases} \tag{4}$$

for any $f \in L^2(\mathbb{M};\mathbb{C})$. System (3) is locally well-posed on every high regularity Hardy–Sobolev space

$$H^s_+(\mathbb{M};\mathbb{C}) := \Pi\left(H^s(\mathbb{M};\mathbb{C})\right), \quad \forall s > \tfrac{3}{2}. \tag{5}$$

System (3) is $L^2$-critical (i.e., its scaling leaves the $L^2$-norm invariant) and is invariant under phase rotations, spatial translations, scaling, and Galilean boosts. It can be regarded as the unification of the following three classical integrable models on the manifold $\mathbb{M}$:

1. the focusing Calogero–Moser–Sutherland derivative nonlinear Schrödinger equation
$$\partial_t u = i\partial_x^2 u + 2u\partial_x\Pi(|u|^2); \tag{6}$$
2. the defocusing Calogero–Moser–Sutherland derivative nonlinear Schrödinger equation
$$\partial_t u = i\partial_x^2 u - 2u\partial_x\Pi(|u|^2); \tag{7}$$
3. the linear Schrödinger equation
$$\partial_t u = i\partial_x^2 u. \tag{8}$$

According to [1], the focusing equation (6) can be interpreted as the thermodynamic continuum limit of the Calogero–Moser–Sutherland model. According to [3] and [4], the defocusing equation (7) is the deep-water limit of the intermediate nonlinear Schrödinger equation introduced by D. Pelinovsky. Physicists aim to find a single equation that can encapsulate multiple principles, and unify different theories or phenomena under a single overarching framework. The intertwined system (3) unifies equations (6), (7), (8), preserving their complete integrability. It offers a potential unified framework for understanding these distinct physical phenomena. According to [2], multi-soliton solutions of (6) can be interpreted as a complexified version of the Calogero–Moser system

$$\partial_t^2 z_j = \sum_{k\neq j} \frac{8}{(z_j - z_k)^3}, \quad 1 \leq j \leq N \qquad (9)$$

Exploring the multi-soliton dynamics of system (3) with two variables could potentially yield more controllable Calogero–Moser type systems and facilitate the study of perturbation theories for the original systems.

## 3 Methods and Results

First, I establish the Lax pair structure of (3)

$$\tfrac{\mathrm{d}}{\mathrm{d}t}\mathbf{L}_{\Psi(t)} = [\mathbf{B}_{\Psi(t)}, \mathbf{L}_{\Psi(t)}]. \qquad (10)$$

The operators $\mathbf{L}_\Psi$ and $\mathbf{B}_\Psi$ consist of differential operators and Toeplitz operators. In [5], I also express every general solution of (3) in terms of its initial datum and the time variable, leading to the explicit formulas both on $\mathbb{T}$ and on $\mathbb{R}$.

Secondly, I have classified all multi-soliton potentials for system (3) on the real line by using the spectral analysis of the Lax–Beurling shift-semigroups and the Lax operator $\mathbf{L}_\Psi$ given by (1.18) of [5].

By using an auxiliary identity, I obtain the *simplicity* and *finiteness* of eigenvalues of the Lax operator $\mathbf{L}_\Psi = \mathbf{L}_\Psi^* : H_+^1(\mathbb{R};\mathbb{C}) \to L_+^2(\mathbb{R};\mathbb{C})$, whenever $\psi = (u,v) \in H_+^1(\mathbb{R};\mathbb{C}^{1\times 2})$. My *spectral characterization theorems* indicate that if $\psi = (u,v)$ is a multi-soliton potential of system (3), then $u$ and $v$ both belong to the purely point spectral subspace $\mathcal{H}_{\mathrm{pp}}(\mathbf{L}_\Psi)$. The simplicity of eigenvalues allows for a better description for each multi-soliton potential.

Finally, I discovered a new phenomenon for system (3) by classifying all its traveling solitary waves. This phenomenon is entirely different from Gérard–Lenzmann's result [2] on the focusing equation (6) on $\mathbb{R}$. According to Theorem 1.2 of [2], every traveling solitary wave $\tilde{u}(t) \in H_+^1(\mathbb{R})$ for (6) has the same $L^2$-mass, i.e. $\|\tilde{u}(t)\|_{L_+^2}^2 = 2\pi$. It cannot be changed via scaling in the focusing case (6). However, this property does not hold for system (3). Even though (3) on $\mathbb{R}$ is $L^2$-critical, the $L^2$-norm of its traveling solitary waves can be arbitrarily large. Precisely, for any $\mathbf{m} \geq 4\pi$, there exists a solitary wave $\Psi = (u,v)$ such that

$$\|\Psi(t)\|_{L^2}^2 := \|u(t)\|_{L^2}^2 + \|v(t)\|_{L^2}^2 = \mathbf{m}. \qquad (11)$$

Moreover, if $\|u(t)\|_{L^2} \to 0^+$, then $\|v(t)\|_{L^2} \to +\infty$ necessarily. Any solitary wave $\Psi$ satisfies $\|\Psi(t)\|_{L_+^2(\mathbb{R};\mathbb{C}^{1\times 2})} \geq 2\sqrt{\pi}$. If $\mathrm{Re}\left(\frac{v(0,x)}{u(0,x)}\right) \to 0^+$, then $\|u(t)\|_{L^2} \to +\infty$. Compared with the focusing equation (6), system (3) exhibits much more traveling wave solutions whose $L^2$-masses are controllable. The classification of traveling solitary waves for (3) relies on the analysis for a non local Euler–Lagrange equation and the use of some Pohozaev-type identity.

To obtain these results, we combine the Lax pair structure (10) with tools from Hardy space including Toeplitz operators, Poisson integrals, and the Lax–Beurling shift semigroups. Furthermore, we employ both direct and inverse spectral theory for the Lax operator $\mathbf{L}_\Psi$ to study the long-time behavior of pure multi-soliton solutions of system (3).

# Quadrature for Space-Time Galerkin BEM of the Wave Equation

**Johannes Tausch**[1,*]

[1]Department of Mathematics, Southern Methodist University, Dallas, TX, USA

*Email: tausch@smu.edu

**Abstract**

A new approach to computing integrals that arise in the Galerkin discretization of retarded layer potentials is introduced. The kernel presents two difficulties (a) a singularity in the origin and (b) support on the light cone, a lower dimensional manifold in space-time. The quadrature method is based on decomposing the integration domain (tensor product of two elements) into the convex hulls of simpler polytopes where the singularity occurs in a single variable.

## 1 Space Time Galerkin Discretization

A popular solution approach to the Dirichlet problem of the wave equation $u_{tt} = \Delta u$ is to represent $u$ by a retarded single layer potential

$$\mathcal{V}q(\mathbf{x},t) = \int_\Gamma \frac{1}{4\pi r} q(\mathbf{y}, t-r)\, ds(\mathbf{y}),$$

where $r = |\mathbf{x}-\mathbf{y}|$ and $q$ is an unknown density on the space time cylinder $\Gamma \times [0,T]$. The jump relations of the single layer potential lead to the boundary integral equation

$$\mathcal{V}q(\mathbf{x},t) = f(\mathbf{x},t), \qquad (\mathbf{x},t) \in \Gamma,$$

where the right hand side is the Dirichlet datum. Several discretization methods exist, see, e.g., [1]. Here we consider a time-domain BEM with tensor products of standard finite elements in space and time. For space, the surface $\Gamma$ is subdivided into a conforming triangulation $\Gamma = \bigcup_{j=1}^{N_f} P_j$, with corresponding nodal basis functions $\phi_j$. The time interval is subdivided uniformly with resulting basis functions $\varphi_k$.

The Galerkin projection is obtained by multiplying with the time derivative of the test functions. Expanding $q = \sum_{\ell j} \varphi_\ell \phi_j q_{\ell,j}$ leads to the block lower triangular linear system

$$\sum_{\ell=1}^{k} A_{k-\ell}\, q_\ell = f_k, \quad k = 0, 1, \dots N_t,$$

where the matrices $A_d$ are defined by

$$A_{d,ij} = \int_\Gamma \int_\Gamma \frac{g_d(r)}{4\pi r} \phi_i(\mathbf{x}) \phi_j(\mathbf{y})\, ds(\mathbf{y}) ds(\mathbf{x}),$$

and $r = |\mathbf{x}-\mathbf{y}|$. The kernel $g_d$ accounts for the time integration of the $\varphi$'s

$$g_d(r) = \int_0^T \varphi_d'(t) \varphi_0(t-r)\, dt.$$

Since the $\varphi_d$'s are piecewise polynomials with finite support it follows that

$$g_d(r) = \begin{cases} 0 & \text{for } r < 0, \\ g(r/h - d) & \text{for } r \geq 0. \end{cases}$$

Here $g$ is another piecewise polynomial with support in $[-M, M]$, where $M$ is an integer. This function can be written as a telescoping sum of polynomials $\mathrm{k}_j$ such that

$$g(\rho) = \mathrm{k}_{-M}(\rho) + \cdots + \mathrm{k}_j(\rho), \text{ for } \rho \leq j$$

and $j \in \{-M, \dots, M\}$. This makes clear that the matrix coefficients $A_{d,ij}$ are sums of integrals of the form

$$I = \int_{P_i} \int_{P_j} \frac{\mathrm{k}(r)}{4\pi r} \phi(\mathbf{x},\mathbf{y})\, ds(\mathbf{y}) ds(\mathbf{x}). \qquad (1)$$

where $\mathrm{k}(\cdot)$ is one of the functions $\mathrm{k}_j$ the interval $[0, j]$ which is extended by zero. Further, $\phi(\cdot,\cdot)$ is a smooth function that represents the contributions of the shape functions. The main difference of this integral to its elliptic counterpart is the fact the kernel k in (1) has finite support and has discontinuous derivatives at $r = R$.

## 2 Quadrature for the non-smooth Kernel

It is well known how to handle integrals in (1) when the kernel is singular in the origin and smooth otherwise. Our approach follows [2], with the novelty that the integrand has low regularity when $r = R$. We briefly recall the main concepts from polytope theory that are necessary to understand this construction.

**Polytopes, Pyramids and Convex Hulls.** The triangles $P_j$ are a special case of convex polytopes. In general, a convex polytope in $\mathbb{R}^d$ is the convex hull of its vertices $v_i \in \mathbb{R}^d$,

$$\mathrm{conv}(v_0, \ldots, v_n) = \left\{ \sum_{i=0}^{n} \lambda_i v_i \colon \lambda_i \geq 0, \sum_{i=0}^{n} \lambda_i = 1 \right\}.$$

It is important to distinguish between the space dimension $d$ and the dimension of the polytope. For instance, the triangles $P_j$ are two dimensional polytopes in $\mathbb{R}^3$. Likewise, the Cartesian product $P_i \times P_j$ in (1) is a four dimensional polytope in $\mathbb{R}^6$. Note that the faces of $P_i \times P_j$ are the Cartesian products of $P_i$'s and $P_j$'s faces.

An important class of polytopes are pyramids. These are polytopes where all but one vertex are in an affine plane. The vertex that is not in the plane is the apex of the pyramid. The convex hull of the other vertices is the base. The notation is $\mathrm{pyr}(v_0, B_0)$ where $v_0$ is the apex and $B_0$ is the base. If the base itself is a pyramid, $B_0 = \mathrm{pyr}(v_1, B_1)$, then we write $\mathrm{pyr}(v_0, v_1, B_1)$. For iterated pyramids one can see that see that $\mathrm{pyr}(v_0, \ldots, v_s, B) = \mathrm{conv}(A, B)$, where $A$ is the $s$-dimensional simplex formed by the apices.

**Pyramidal decomposition.** $P_i \times P_j$ is not a pyramid, but can be decomposed into two pyramids as follows: Select one vertex $v_x$ of $P_x$ and one vertex $v_y$ of $P_y$. Let $E_x$ be the edge of $P_x$ that does not contain $v_x$ and $E_y$ be the edge of $P_y$ that does not contain $v_y$. Further, let $\mathbf{v} = (v_x, v_y)$ be the apex. One can see that

$$P_x \times P_y = \mathrm{pyr}(\mathbf{v}, E_x \times P_y) \cup \mathrm{pyr}(\mathbf{v}, P_x \times E_y).$$

If desired, the decomposition can be repeated for the bases, generating iterated pyramids. For instance, selecting vertices $w_x$ of $E_x$ and $w_y$ of $P_y$ leads to

$$E_x \times P_y = \mathrm{pyr}(\mathbf{w}, u_x \times P_y) \cup \mathrm{pyr}(\mathbf{w}, E_x \times F_y).$$

where $\mathbf{w} = (w_x, w_y)$, $u_x$ is the other vertex of $E_x$ and $F_y$ is the edge of $P_y$ that does not contain $w_y$. This process leads to the decomposition

$$P_x \times P_y = \bigcup_k \mathrm{conv}(A_k, F_{k,x} \times F_{k,y}), \qquad (2)$$

where $F_{k,x}$ and $F_{k,y}$ are faces of $P_x$ and $P_y$ and $A_k$ is the simplex of the apices that appeared in the decomposition that lead to $F_{k,x} \times F_{k,y}$.

**Integrals over Convex Hulls** With the parameterization $\mathbf{z} = (1-\lambda)\mathbf{z}_A + \lambda \mathbf{z}_B$, where $\mathbf{z}_A = (\mathbf{x}_A, \mathbf{y}_A) \in A$, $\mathbf{z}_B = (\mathbf{x}_B, \mathbf{y}_B) \in B$, an integral in (2) is

$$\int_{\mathrm{conv}(A,B)} f(\mathbf{z}) d\mathbf{z} = \delta_{AB} \int_A \int_B I(\mathbf{z}_B, \mathbf{z}_A) d\mathbf{z}_B d\mathbf{z}_A,$$

where the integrand contains the $\lambda$-integration

$$I(\mathbf{z}_B, \mathbf{z}_A) = \int_0^1 f\big((1-\lambda)\mathbf{z}_A + \lambda \mathbf{z}_B\big)(1-\lambda)^s \lambda^q d\lambda.$$

Here, $\delta_{AB}$ is a constant, $s = \dim(A)$ and $q = \dim(B)$. In the integrand of (1)

$$r = \rho(\lambda) := |(1-\lambda)(\mathbf{x}_A - \mathbf{y}_A) + \lambda(\mathbf{x}_B - \mathbf{y}_B)|$$

is the root of a parabola that opens upward. The $\lambda$-integral can be reduced to the part of $[0, 1]$ for which $\rho(\lambda) \leq R$.

We will show how to perform the pyramidal decomposition such that the resulting iterated pyramids fall only into two distinct groups.

1. $A$ and $B$ are such that $|\mathbf{x}_A - \mathbf{y}_A| \leq R$ and $|\mathbf{x}_B - \mathbf{y}_B| > R$. In this case the equation $\rho(\lambda) = R$ has a unique solution $\lambda^* \in [0, 1]$, that is a smooth function of $\mathbf{z}_A$ and $\mathbf{z}_B$. Hence the $\lambda$-integral reduces to $[0, \lambda^*]$.

2. $A$ and $B$ are such that $|\mathbf{x}_A - \mathbf{y}_A| \leq R$ and $|\mathbf{x}_B - \mathbf{y}_B| \leq R$. In this case the convexity of $\rho^2(\lambda)$ implies that $\rho(\lambda) = R$ has no solution in $[0, 1)$ and the $\lambda$-integral is smooth.

In both cases $I(\mathbf{z}_A, \mathbf{z}_B)$ is a smooth function in $A \times B$. This enables the use of high order Gauss- or Gauss-like quadrature to compute all pyramids in the decomposition (2). The significance of the discussed decomposition is that it avoids all other cases that would lead to non-smooth integrands.

# A tsunami inverse problem in the time domain: numerical analysis of an optimal control approach

**Laurent Bourgeois**[1], **Philippe Moireau**[2,*], **Raphaël Terrine**[1]
[1]Poems, CNRS, Inria, Ensta Paris, Institut Polytechnique de Paris
[2]CMAP, Inria, CNRS, École polytechnique, Institut Polytechnique de Paris
*Email: philippe.moireau@polytechnique.edu

**Abstract**

This work addresses the inverse problem of reconstructing the seabed displacement, the source of a tsunami, from surface measurements of ocean waves. Assuming the measurements are modeled as boundary observations of a linearized model derived from the Euler equations, the study formulates a data assimilation procedure that combines a data fidelity term and a model constraint. Absorbing boundary conditions and parabolic regularization are incorporated to facilitate numerical resolution and analysis. The associated optimality adjoint system is then derived, and a time-discrete iterative scheme based on finite differences is proposed. Its consistency and convergence are rigorously analyzed by recasting the discretization as a mixed space-time finite element formulation. This enables an efficient inversion method based on an iterative gradient descent with convergence guarantees.

## 1 Direct problem setting

Following [1], tsunami propagation is originally modeled by an acoustic and gravity wave system of the form

$$\begin{cases} c^{-2}\partial_t^2\varphi - \Delta\varphi = 0 & \text{in } \Omega\times(0,T), \\ \partial_t\varphi + g\, s = 0 & \text{on } \Gamma_0\times(0,T), \\ \partial_{x_2}\varphi - \partial_t s = 0 & \text{on } \Gamma_0\times(0,T), \\ \partial_\nu\varphi = \chi & \text{on } \Gamma_{-H}\times(0,T), \\ \varphi = \partial_t\varphi = 0 & \text{on } \Omega\times\{0\}, \\ s = 0 & \text{on } \Gamma_0\times\{0\}, \end{cases} \tag{1}$$

where $\Gamma_0$ is the the free surface at rest and $\Gamma_{-H}$ is the the straight bottom of the ocean $\Omega$. We assume we have a noisy measurement $t \mapsto y(t)$ of the sea height perturbation $s$ for an unknown source $\chi$. Our objective is to reconstruct $\chi$ from $y \in L^2(0,T;\Gamma_0)$, an ill-posed inverse problem for which we can justify a uniqueness result using a Holmgren-type theorem.

## 2 Inverse problem setting

To reconstruct the source in practice, we must bound the domain in the $x_2$ direction at $\pm R$:

$$\partial_x\varphi + \eta\,\partial_{x_1 t}^2\varphi \pm \frac{1}{c}\,\partial_t\varphi = 0 \text{ on } \Sigma_{\pm R}\times(0,T),$$

and we also add parabolic regularization $\eta > 0$, to the now bounded domain equation:

$$\frac{1}{c^2}\partial_t^2\varphi - \Delta\varphi - \eta\,\Delta\partial_t\varphi = 0 \text{ in } \Omega^R\times(0,T).$$

Then, we introduce the spaces

$$H = H^1(\Omega^R)\times L^2(\Omega^R)\times L^2(\Gamma_0^R),$$

and

$$V = \{z = (\varphi,\psi,s) \in H^1(\Omega^R)\times H^1(\Omega^R)\times L^2(\Gamma_0^R) \text{ such that } \psi|_{\Gamma_0^R} + g\,s = 0\},$$

equipped with their usual norms. $V$ is dense in $H$ with continuous injection. Therefore, we can recast the model with parabolic regularization and absorbing boundary conditions in the abstract form: For all $v \in V$,

$$\begin{cases} \frac{\mathrm{d}}{\mathrm{d}t}(z,v)_H + a_\eta(z,v) = b(\chi,v), \\ z(0) = 0. \end{cases} \tag{2}$$

where for all $z = (\varphi,\psi,s)$, $\tilde z = (\tilde\varphi,\tilde\psi,\tilde s)$ in $V$,

$$a_\eta(z,\tilde z) = -(\psi,\tilde\varphi)_{H^1(\Omega^R)} + (\nabla\varphi,\nabla\tilde\psi)_{L^2(\Omega^R)} + \eta\,(\nabla\psi,\nabla\tilde\psi)_{L^2(\Omega^R)} + \frac{1}{c}\sum_{\pm}(\psi,\tilde\psi)_{L^2(\Sigma_{\pm R})}.$$

and $b(\chi,v) = (\chi,v)_{L^2(\Gamma_{-H}^R)}$

Our reconstruction strategy consists in solving the optimal control problem

$$\bar\chi = \operatorname*{argmin}_{\chi\in L^2(0,T;\Gamma_{-H}^R)} \left\{ J(\chi) = \frac{\varepsilon}{2}\|\chi\|^2_{L^2((0,T)\times\Gamma_{-H}^R)} + \frac{1}{2}\|y - s_\chi\|^2_{L^2((0,T)\times\Gamma_0^R)} \right\}, \tag{3}$$

for a well-chosen $\varepsilon > 0$, where $s_\chi$ is the last component of $z_\chi \in W_0(0,T) = \{z \in L^2(0,T;V) \cap H^1(0,T;V') \mid z(0) = 0\}$ solution of (2).

Such an optimal control problem can, under inf-sup conditions, be shown to be equivalent to solving a mixed formulation involving the adjoint equation as the Lagrange multiplier. Find $\bar{Z} = (\bar{z}, \bar{\chi}) \in W_0(0,T) \times L^2((0,T) \times \Gamma^R_{-H})$ and $\bar{p} \in L^2(0,T;V)$ such that for all $\tilde{Z} = (\tilde{z}, \tilde{\chi}) \in W_0(0,T) \times L^2((0,T) \times \Gamma^R_{-H})$ and $\tilde{p} \in L^2(0,T;V)$,

$$\begin{cases} j(\bar{Z}, \tilde{Z}) + d(\tilde{Z}, \bar{p}) = \ell(\tilde{Z}), \\ d(\bar{Z}, \tilde{p}) = 0, \end{cases} \tag{4}$$

with for all $Z = (z, \chi) = (\varphi, \psi, s, \chi)$,

$$j(Z, \tilde{Z}) = \int_0^T [\varepsilon(\chi, \tilde{\chi})_{L^2(\Gamma^R_{-H})} + (s, \tilde{s})_{L^2(\Gamma^R_0)}\, dt,$$

$$d(Z, p) = \int_0^T [\langle \dot{z}, p \rangle_{V',V} + a_\eta(z, p) - b(\chi, p)\, dt,$$

$$\text{and } \ell(Z) = \int_0^T (y, s)_{L^2(\Gamma^R_0)}\, dt.$$

## 3 Optimize then discretize or discretize then optimize

Returning to the optimal control approach, classical numerical methods for solving such variational data assimilation problem typically use a discretize-then-optimize strategy. First, the direct problem is discretized in space and time. The cost functional is then discretized in time, and the adjoint is defined from the discretized problem and functional. This approach computes an exact gradient, which is advantageous for gradient descent, typically based on BFGS algorithm. However, analyzing convergence, especially convergence rates, can be challenging.

In contrast, space-time domain data assimilation formulations directly discretize the mixed formulation. This optimize-then-discretize path often requires inf-sup stabilization terms and necessitates solving the full space-time formulation.

In this work, we show that, in the context of our regularized model, we can propose a time discretization that corresponds exactly to the natural 4D-Var, but can also be analyzed through the lens of a discrete mixed formulation derived from (4) and associated with a discrete inf-sup condition.

The strategy to reconcile the two points of view is to introduce in (4) a time and space discretization of the form $P_1 \times P_k$ for the state variable, and $P_0 \times P_k$ for the input variable, while the adjoint variable also follows a $P_0 \times P_k$ time-space discretization. Moreover, the time integrals in the bilinear forms are discretized using a one-point quadrature rule, allowing, as in discontinuous Galerkin methods, the recovery of finite-difference Euler schemes [2]. This allows us to reinterpret the discretization of the (4) problem as the mixed formulation associated with a discretized criterion under time-discretized dynamics. We can then demonstrate, after proving a discrete inf-sup condition and using the classical Strang lemma, that under a quasi-uniformity assumption of the meshes, we obtain the following convergence theorem.

**Theorem 1** *Assume a sufficiently regular $\bar{\chi}$, the unique solution to the optimal control problem (3), where $\bar{z}$ is the solution to (2), and $\bar{p}$ is the corresponding adjoint variable in (4). Let $\bar{\chi}^\tau_h$, $\bar{z}^\tau_h$, and $\bar{p}^\tau_h$ be the corresponding solutions to the discretized mixed formulation. Under the quasi-uniformity assumptions for the meshes,*

$$\|\bar{z}^\tau_h - \bar{z}\|_{W_0(0,T)} + \|\bar{\chi}^\tau_h - \bar{\chi}\|_{L^2((0,T)\times\Gamma^R_{-H})} \lesssim \frac{\tau + h^k}{h^{\frac{3}{2}}},$$

$$\text{and } \|\bar{p}^\tau_h - \bar{p}\|_{L^2(0,T;V)} \lesssim \frac{\tau + h^k}{h^2}.$$

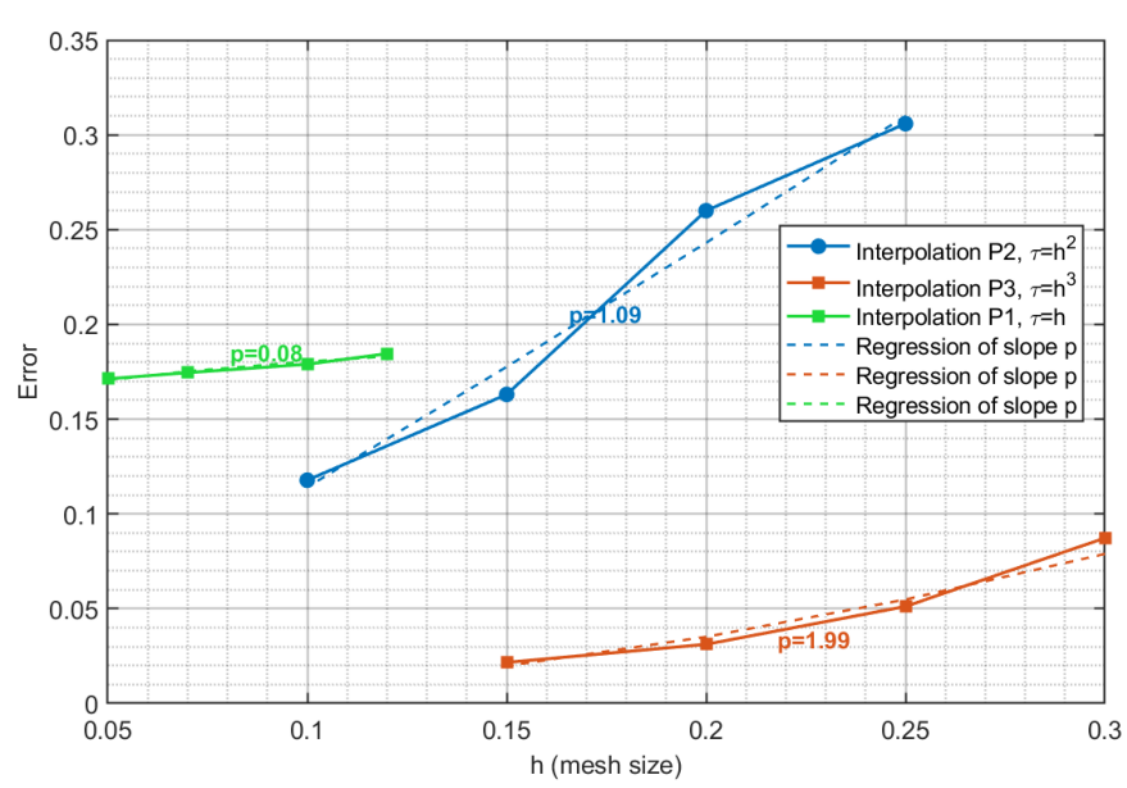

Figure 1: $\chi^\tau_h$ convergence for $k = 1, 2, 3$

# Variational Quasi-Trefftz Method for Time-Harmonic First-Order Wave Propagation Systems

**Matthias Rivet**[1], **Sébastien Pernet**[2], **Sébastien Tordeux**[2], [,*]

[1]ONERA-Toulouse, Project-team Makutu, Inria, Pau, France

[2]Project-team Makutu, Inria, Pau, France

[*]Email: sebastien.tordeux@inria.fr

**Abstract**

This work presents a quasi-Trefftz framework aimed at improving the conditioning of classical Trefftz methods. The proposed approach relies on the introduction of an auxiliary solver used to construct basis functions that are significantly less correlated.

## 1 Introduction

In the frequency domain, the numerical solution of wave propagation problems using classical discretization methods requires the LU factorization of the matrix appearing on the left-hand side. Unfortunately, this leads to extremely high memory costs due to fill-in in the matrix sparsity pattern, even when reordering strategies are used.

However, domain decomposition strategies can be employed to enable the use of iterative solvers together with matrix-free techniques, which drastically reduce the computational cost.

Variational Trefftz methods can be interpreted as an algebraic reformulation of domain decomposition methods, making them particularly well suited for high-performance computing (HPC). However, when the Galerkin space is constructed from sums of plane waves, it is well known that these methods suffer from severe round-off errors, which significantly deteriorate the quality of the numerical results.

In this talk, we present a unified framework for Trefftz methods together with the definition of a new quasi-Trefftz approximation space that significantly improves the robustness and accuracy of the numerical method.

## 2 Time-harmonic first-order systems

Symmetric Friedrichs hyperbolic systems provide a unified framework for the modeling of wave propagation phenomena. They encompass a wide range of physical models, including acoustic waves in fluids, Maxwell's equations in vacuum, and the equations of linear elasticity.

These systems can be written in the form

$$M\,\partial_t u + \sum_{i=1}^{d} F_i\,\partial_{x_i} u = 0 \qquad \text{in } \Omega \subset \mathbb{R}^d, \quad (1)$$

where $u : \Omega \times \mathbb{R} \to \mathbb{R}^m$ is the unknown field, $M \in \mathbb{R}^{m\times m}$ is a symmetric positive definite matrix, and $(F_i)_{1\le i\le d}$ are symmetric matrices of size $m \times m$.

In the time-harmonic regime we seek solutions of the form

$$u(x,t) = \Re\big(u(x)e^{-i\omega t}\big)\,,$$

which leads to the time-harmonic first-order system

$$-i\omega M u + \sum_{i=1}^{d} F_i\,\partial_{x_i} u = 0 \qquad \text{in } \Omega. \quad (2)$$

Boundary conditions are expressed in terms of the normal flux matrix

$$F(\mathbf{n}) = \sum_{i=1}^{d} n_i F_i,$$

where $\mathbf{n}$ denotes the unit outward normal to $\partial\Omega$. The matrix $F(\mathbf{n})$ admits the spectral decomposition

$$F(\mathbf{n}) = F^+(\mathbf{n}) + F^-(\mathbf{n}),$$

where $F^+(\mathbf{n})$ and $F^-(\mathbf{n})$ correspond to the positive and negative parts of the spectrum.

Boundary conditions are imposed on the incoming characteristic variables associated with the negative eigenvalues of $F(\mathbf{n})$. For instance, an absorbing boundary condition can be written as

$$F^-(\mathbf{n})u = F^-(\mathbf{n})u_{\text{inc}},$$

where $u_{\text{inc}}$ denotes an incident field.

This unified structure turns out to be particularly well suited for Trefftz discretization. In particular, we observed that it yields an algebraic framework tailored for HPC, enabling the systematic and largely automated design of Trefftz methods.

## 3 Trefftz variational formulation

Let $\mathcal{T}_h$ be a mesh of $\Omega$ and denote by $\mathcal{F}_h$ the skeleton of the mesh. For each element $K \in \mathcal{T}_h$, we introduce the local Trefftz space

$$T(K) = \Big\{ v \in (L^2(K))^m \ : \ - i\omega M v + \sum_{i=1}^{d} F_i \, \partial_{x_i} v = 0 \quad \text{in } K, \quad F(\mathbf{n}) v \in L^2(\partial K) \Big\}.$$

The global Trefftz space is then defined as

$$T_h = \{ v \in L^2(\Omega)^m \ : \ v_{|K} \in T(K), \ \forall K \in \mathcal{T}_h \}.$$

Since the functions in $T_h$ satisfy the differential equation inside each element, the volume contributions vanish and the variational formulation involves only interface contributions over the mesh skeleton.

The Trefftz variational formulation reads: find $u_h \in T_h$ such that for all $v_h \in T_h$

$$\sum_{K \in \mathcal{T}_h} \int_{\partial K} \big( F^+(\mathbf{n}_K) \, u_K + F^-(\mathbf{n}_K) \, u_V \big) \cdot \overline{v_K} \, ds = 0.$$

Here $u_K$ denotes the restriction of $u_h$ to the element $K$, while $u_V$ denotes the trace coming from the neighbouring element across the considered interior face. On boundary faces, $u_V$ is replaced by the incident field $u_{\text{inc}}$.

## 4 Quasi-Trefftz discrete space

Classically, discrete Trefftz spaces are constructed from exact local solutions of the governing equation (plane waves, Bessel functions, evanescent modes). However, this approach may lead to large numerical errors due to ill-conditioning, which in three dimensions may severely deteriorate the quality of the approximation or even cause the numerical method to fail.

Our main idea is to define the functions of the Trefftz space not as exact solutions of the first-order system, but as highly accurate approximations obtained by solving a local first-order system numerically in each element. The numerical solver must be sufficiently accurate in order to avoid introducing additional consistency errors.

More precisely, for each element $K$, we consider the local problem

$$-i\omega M v + \sum_{i=1}^{d} F_i \, \partial_{x_i} v = 0 \qquad \text{in } K, \tag{3}$$

together with the boundary condition

$$F^-(\mathbf{n}) \, v = \varphi \qquad \text{on } \partial K, \tag{4}$$

where $\varphi$ belongs to a suitable piecewise polynomial space defined on $\partial K$.

## 5 Content of the talk

In this talk, we present the construction of the quasi-Trefftz formulation based on a Flux Reconstruction method [2]. This approach is particularly well suited for the type of local problems considered here and enables the construction of highly accurate quasi-Trefftz basis functions.

Moreover, the impact on the conditioning of the method and on the quality of the approximation will be investigated. In particular, we establish a clear relationship between the degree of the polynomial trace space and the degree of the Flux Reconstruction method.

Finally, an optimization of the decomposition of the flux operator $F(\mathbf{n})$ into $F^+(\mathbf{n})$ and $F^-(\mathbf{n})$ using artificial intelligence techniques will be performed in order to further improve the quality of the Trefftz method.

More details and references can be found in the PhD thesis of Matthias Rivet [3].

# Optimization-based transionospheric SAR autofocus with screen projection precursor

**Mikhail Gilman**[1] **and Semyon Tsynkov**[1,*]
[1]Department of Mathematics, North Carolina State University, Raleigh, NC 27695, USA
*Email: tsynkov@math.ncsu.edu

**Abstract**

We propose a novel optimization-based autofocus algorithm for spaceborne synthetic aperture radar (SAR) that corrects the distortions of radar images due the ionospheric turbulence. Those are the distortions that cannot be corrected by traditional autofocus techniques.

## 1 Introduction

Traditional autofocus algorithms for synthetic aperture radar (SAR) are designed for correcting the distortions of SAR images caused by uncertainties in the radar antenna trajectory. They are not well suited for correcting the distortions due to the ionospheric turbulence. The reason is that in conventional imaging, the phase perturbations of SAR signals that cause the distortions depend only on the antenna position while in the case of propagation through the turbulent plasma, for each antenna position they additionally depend on the target coordinates. The latter dependence cannot be adequately addressed by traditional autofocus [1].

Previously, we have built a variational transionospheric SAR autofocus algorithm that circumvents the deficiencies of conventional techniques [2]. It is applied concurrently with SAR inversion per se, and demonstrates superior performance for a range of imaging scenarios. It also outperforms its nearest alternative, which is a two stage non-variational procedure that involves traditional autofocus preceded by screen projection for approximating the ionospheric effects [3]. Yet the performance of our optimization-based autofocus is not uniform. For high levels of turbulent fluctuations and clutter, its focusing capacity appears not as robust because the optimization formulation is not convex [4].

In this work [5], we enhance the optimization performance in the non-convex case by carefully selecting the initial guess for gradient descent. We choose it as the outcome of the two stage screen projection algorithm [3]. This leads to a drastic improvement in the quality of focusing while adding little to the computational cost.

## 2 Model, analysis, and results

The range-compressed signal due to ground reflectivity $\mu$ received by the radar antenna at the orbit admits the following representation [5]:

$$u(x) = \int_{x-L_{\mathrm{SA}}/2}^{x+L_{\mathrm{SA}}/2} e^{\mathrm{i}\frac{\omega_0(x-z)^2}{Rc}} \cdot e^{-\mathrm{i}\Psi(s(x,z))}\mu(z)\,\mathrm{d}z,$$
$$s(x,z) = \xi x + (1-\xi)z, \quad 0<\xi<1, \qquad (1)$$

where $z$ and $x$ denote the azimuthal coordinates at the ground level and orbit level, respectively, $L_{\mathrm{SA}}$ is the synthetic aperture, $R$ is the distance from the orbit to the target area on the ground, $c$ is the speed of light, $\omega_0$ is the radar carrier frequency, and $\Psi$ is the density of the phase screen at relative elevation $\xi$ between the ground and orbit that accounts for the scintillation phase error (SPE) due to the ionospheric turbulence.

SAR image approximates the unknown $\mu$:

$$\mathcal{I}(y) = \frac{\omega_0}{\pi Rc}\int_{y-L_{\mathrm{SA}}/2}^{y+L_{\mathrm{SA}}/2} e^{-\mathrm{i}\frac{\omega_0(x-y)^2}{Rc}} \cdot e^{\mathrm{i}\Psi^{\mathrm{rec}}(s(x,y))}u(x)\,\mathrm{d}x, \quad (2)$$

where $\Psi^{\mathrm{rec}}$ is supposed to reconstruct the unknown SPE $\Psi$. We represent the SPE as

$$\Psi(s) = \sum\nolimits_{n=1}^{N}\big(p_n\cos(k_n s) + q_n\sin(k_n s)\big),$$

where $k_n>0$, $p_n$, $q_n\in\mathbb{R}$, and the various spatial scales $k_n$ characterize the turbulence. Then, we seek the reconstruction SPE in the form

$$\Psi^{\mathrm{rec}}(s) = \sum_{n=1}^{N^{\mathrm{rec}}}\big(p_n^{\mathrm{rec}}\cos(k_n^{\mathrm{rec}}s) + q_n^{\mathrm{rec}}\sin(k_n^{\mathrm{rec}}s)\big).$$

To determine $\Psi^{\mathrm{rec}}$, we make a common assumption that the ground reflectivity $\mu(z)$ contains bright point scatterers, define a sparsity promoting cost function [6] with the penalty term that facilitates spectral decay for higher frequencies:

$$\text{Cost}[\mathcal{I}, \Psi^{\text{rec}}] = -\frac{1}{L_{\text{SA}}} \int |\mathcal{I}[\Psi^{\text{rec}}](y)|^4 \text{d}y + \zeta \frac{\omega_0}{Rc} \sum_{n=1}^{N^{\text{rec}}} k_n^2 \big( (p_n^{\text{rec}})^2 + (q_n^{\text{rec}})^2 \big), \quad (3)$$

and minimize it with respect to $p_n^{\text{rec}}$ and $q_n^{\text{rec}}$. Gradient-based minimization can be used since the cost function (3) is differentiable. The analysis of [6] also suggests that the first term in (3), which is $\propto \| \cdot \|_4^4$, offers a balanced performance between promoting sparsity in the case of bright point targets and applicability to lower contrast scenes as well, see [2]. It is important to emphasize that to obtain $\Psi^{\text{rec}}$ by minimizing (3), one does not need to know the ground truth.

The resulting image (2) shows superior quality of focusing [2], but it may deteriorate for high levels of turbulent perturbations $\Psi$ and clutter (the white noise component of $\mu$). The reason is that the cost function (3) is non-convex, and the performance of gradient descent appears very sensitive to the choice of the initial guess.

To supply a good initial guess for the minimization of (3), we employ the two stage screen projection autofocus of [3]. Compared to the optimization-based autofocus of [2], it underperforms as a standalone technique [7], but it is a simple and numerically inexpensive algorithm that does not involve any optimization. The single stage inversion (2) is replaced with

$$p(s) = \frac{1}{\eta L_{\text{SA}}} \int_{s-\eta L_{\text{SA}}/2}^{s+\eta L_{\text{SA}}/2} e^{-\text{i}\frac{\omega_0 (x-s)^2}{\eta Rc}} u(x)\, \text{d}x, \quad (4\text{a})$$

$$\mathcal{I}(y) = \frac{C\omega_0}{\xi \pi Rc} \int_{y-\xi L_{\text{SA}}/2}^{y+\xi L_{\text{SA}}/2} e^{-\text{i}\frac{\omega_0 (y-s)^2}{\xi Rc}} e^{\text{i}\Psi^{\text{rec}}(s)} p(s)\, \text{d}s, \quad (4\text{b})$$

where $\eta = 1 - \xi$ and $C = \left( \frac{\pi Rc}{\omega_0 \xi \eta L_{\text{SA}}^2} \right)^{-1/2} e^{\text{i}\frac{\pi}{4}}$. The first stage (4a) produces a partially focused image $p(s)$ at screen elevation and does not require any ionospheric correction at all. The second stage (4b) merges the partially focused image and ionospheric correction into one function $e^{\text{i}\Psi^{\text{rec}}(s)} p(s)$ of a single argument $s$ that does not depend $x$, $z$, or $y$ [cf. equations (1)–(2)]. This allows one to apply a conventional SAR autofocus technique, such as phase curvature autofocus (PCA) for stripmap imaging, to determine the reconstruction SPE $\Psi^{\text{rec}}$ in (4b). The latter is then used as initial guess for the gradient descent minimization of the cost function (3).

In Figure 1, we are showing an example of SAR inversion (2) with reconstruction SPE obtained by minimizing (3) starting from initial guess (4). The specific parameters of simulation include the magnitude of SPE $\|\Psi\|_2$, magnitude of clutter, and ratio of synthetic aperture to outer scale of turbulence $L_{\text{SA}}/l_{\max}$. The performance of autofocus is on par with ideal reconstruction, where the true $\Psi$ (unknown in practice) is substituted for $\Psi^{\text{rec}}$. In [5], we considered multiple realizations of SPE and clutter and conducted a comprehensive statistical analysis of autofocus performance that corroborated its superior focusing capacity for a broad range of imaging scenarios and ionospheric conditions.

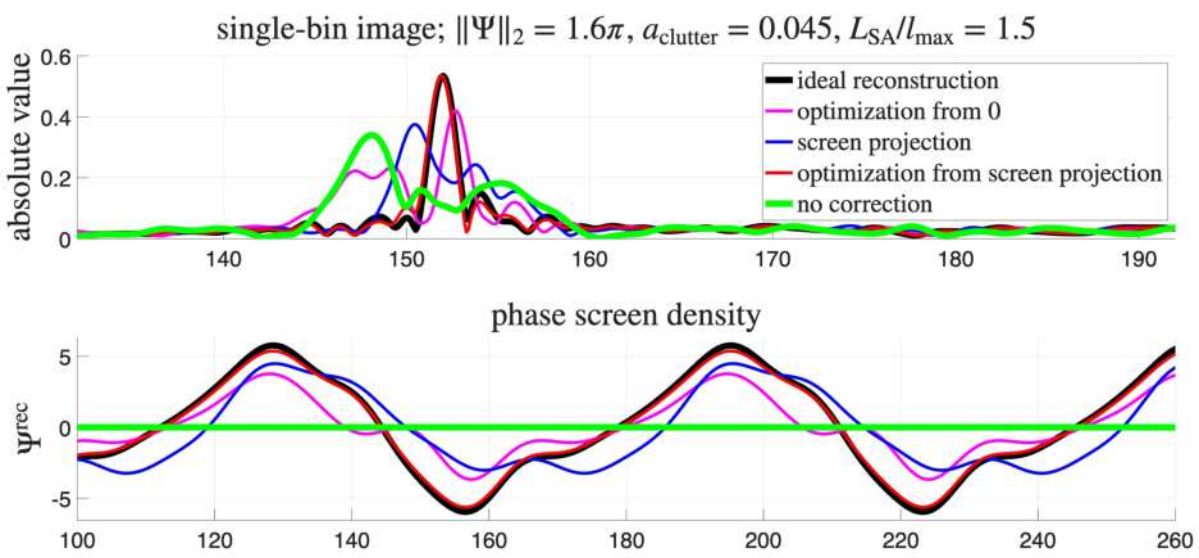

Figure 1: Reconstruction example: $\mathcal{I}$ and $\Psi$.

**Acknowledgement.** Supported by AFOSR, grants FA9550-23-1-0101 and FA9550-24-1-0172.

# Topological solutions and accelerated waves in electromagnetic theory

**Radmila Sazdanovic**[1], **Paul Spears**[1], **and Semyon Tsynkov**[1,*]
[1]Department of Mathematics, North Carolina State University, Raleigh, NC 27695, USA
*Email: tsynkov@math.ncsu.edu

## Abstract

We discuss three classes of special electromagnetic solutions that may hold a considerable potential for advancing various applications ranging from communications theory to directed energy to remote sensing. Those are (1) solutions to unsteady Maxwell's equations where the field lines form specific topological knots or links; (2) time-harmonic topological vortex knots; and (3) bending beams (also known as Airy beams or accelerated waves), including the non-paraxial beams. Vortex knots and accelerated beams have been observed experimentally in optics. Their laboratory realization using terahertz radio waves is of substantial future interest.

## 1 Time-dependent topological solutions

"Marrying" topology with physics is not trivial. These are two vastly different disciplines, and it is not immediately obvious how one can embed a topological structure into the electromagnetic field. First and foremost, what physical quantity(ies) and/or their attributes shall the topological properties be assigned to?

The original idea proposed by Rañada [1] was to have the electric, magnetic, and Poynting field lines mutually orthogonal and linked, yielding a Hopf fibration. Consider

$$\eta(x,y,z,t) = \frac{(Az + t(A-1)) + i(tx - Ay)}{(Ax + ty) + i(A(A-1) - tz)},$$

where $A = (r^2 - t^2 + 1)/2$, and

$$\zeta(x,y,z,t) = \frac{(Ax + ty) + i(Az + t(A-1))}{(tx - Ay) + i(A(A-1) - tz)}.$$

Then, the vector fields $\boldsymbol{E}$ and $\boldsymbol{B}$ given by the following expressions solve the homogeneous time-dependent Maxwell's equations in vacuum normalized to the speed of light $c = 1$:

$$\boldsymbol{E} = \frac{1}{4\pi i}\frac{\nabla\zeta \times \nabla\bar{\zeta}}{(1+\zeta\bar{\zeta})^2}, \quad \boldsymbol{B} = \frac{1}{4\pi i}\frac{\nabla\eta \times \nabla\bar{\eta}}{(1+\eta\bar{\eta})^2}.$$

The vectors $\boldsymbol{E}$ and $\boldsymbol{B}$ are orthogonal; the tangent lines to $\boldsymbol{E}$, $\boldsymbol{B}$ and $\boldsymbol{E}\times\boldsymbol{B}$ form Hopf links. The electromagnetic field is null: $\boldsymbol{E}\cdot\boldsymbol{B} = 0$ and $|\boldsymbol{E}|^2 - |\boldsymbol{B}|^2 = 0$, and decays sufficiently fast toward infinity, so that the total energy is finite.

Work [1] has been generalized by various authors using the Bateman's method [2] of generating the null fields, to include the entire family of torus knots. The design of the corresponding solutions is quite elaborate though, and most importantly, their experimental realization remains elusive, even for the simplest topology.

## 2 Topological vortex knots

An alternative idea proposed by Dennis, et. al. [3] was to rather have the vortex lines of the field (lines of zero intensity and singular phase) knotted or linked. Consider $\mathbb{C} \ni \psi = \psi(x,y,z)$, such that $\psi = 0$ defines a knot in $\mathbb{R}^3$. $\psi$ could be a Milnor polynomial or a Bode-Hirasawa polynomial [4]. Then, $\psi$ is interpreted as a physical quantity of interest, e.g., linearly polarized electric field, and propagated paraxially:

$$-2ik\frac{\partial u}{\partial z} = \frac{\partial^2 u}{\partial x^2} + \frac{\partial^2 u}{\partial y^2},$$

with propagation started off by a cross-section of the polynomial:

$$u(x,y,z)\big|_{z=0} = \psi(x,y,0).$$

The solution $u = u(x,y,z)$ is also a polynomial, not the same as the original $\psi(x,y,z)$. Its vortex lines may have the same topology as that of the original knot, in which case the knot is called stable. Otherwise, it is called unstable. No general theory is available that would render a priori classification; the torus knots $3_1$, $5_1$, $7_1$, as well as some other knots, are stable. Stability is important for future applications in telecommunications, where the transmitted knots will be used in the capacity of alphabet characters.

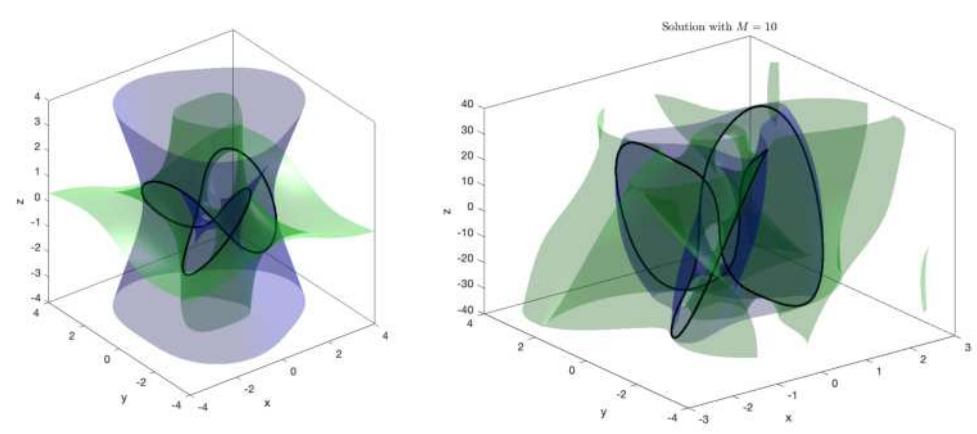

Figure 1: Original and propagated trefoil knot.

Vortex knots have been successfully realized in optics experiments. One important component that has been missing though is their automatic recognition/identification in the received field. The development of a methodology for automated knot recognition involves the analysis of sampling strategies, curve reconstruction, knot identification using the appropriate invariants, and the analysis of sensitivity to noise. The latter will allow one to obtain the probabilities of receiving the same knots as those transmitted, which is needed for computing the key characteristics of the resulting noisy communication channel, such as entropy and capacity.

## 3 Accelerated beams

Airy packets were first introduced in quantum mechanics by Berry and Balazs [5]. Those are solutions to the linear Schrödinger equation: $-i\frac{\partial\psi}{\partial t} = \frac{\hbar}{2m}\frac{\partial^2\psi}{\partial x^2}$ having the following form:

$$\psi(x,t) = \mathrm{Ai}\left[\frac{B}{\hbar^{2/3}}\left(x - \frac{B^3t^2}{4m^2}\right)\right]e^{i\frac{B^3t}{2m\hbar}\left[x-\frac{B^3t^2}{6m^2}\right]}.$$

The quantity $|\psi|^2$ "slides" undistorted along the parabola $x = \frac{B^3t^2}{4m^2}$, which is interpreted as acceleration. Yet it is shown in [5] that particle orbits remain straight lines, and what accelerates in an Airy packet is not any individual particle, but the caustic, i.e., envelope of the family of orbits.

Airy packets were later reinvented and generalized in optics, becoming Airy beams. Generalizations include truncation that makes the total energy finite, as well as transition to non-paraxial regime, where radiation solution to the Helmholtz equation $\frac{\partial^2\psi}{\partial x^2} + \frac{\partial^2\psi}{\partial z^2} + k^2\psi = 0$, $x \in \mathbb{R}$, $z \geqslant 0$, subject to $\psi(x,z)|_{z=0} = \psi_0(x)$, is given by the Sommerfeld integral:

$$\psi(x,z) = 2\int_{-\infty}^{\infty}\psi_0(\xi)\frac{1}{4i}\frac{\partial H_0^{(1)}\left(k\sqrt{(x-\xi)^2+z^2}\right)}{\partial z}d\xi.$$

In doing so, the light rays in the solution $\psi(x,z)$ are still straight lines (can be derived using stationary phase argument), while the envelope of the light rays (i.e., the caustic) may be bent.

The shape of the envelope is completely controlled by the phase $\phi(\xi)$ of the data $\psi_0(\xi) = A(\xi)e^{i\phi(\xi)}$, see [6]. The amplitude $A(\xi)$ can be varied independently to achieve the additional desired properties of an accelerated beam.

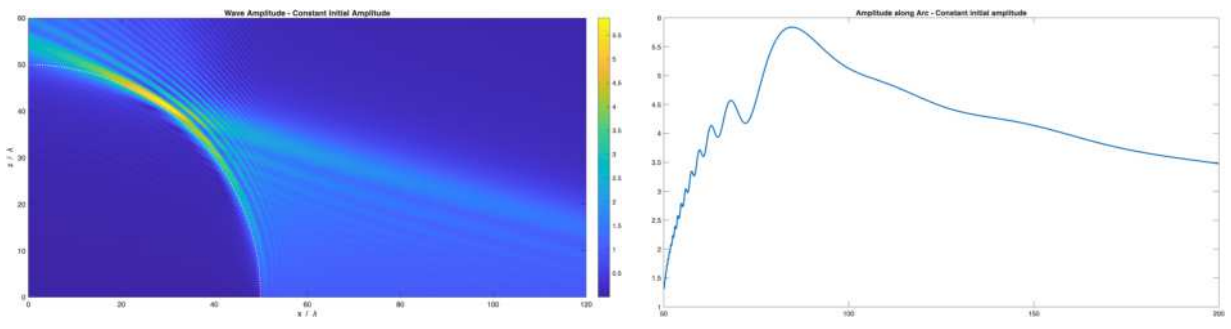

Figure 2: Circular beam for $A(\xi) = \text{const}$.

In Figure 2, we are showing a 2D plot and the intensity profile along the caustic shaped as a circle for an accelerated beam that corresponds to $A(\xi) = \text{const}$. In Figure 3, we are showing the results of optimization, where the objective was to achieve a constant beam intensity along a given portion of the arc. To perform optimization, we represented $A(\xi)$ as a truncated Chebyshev expansion, used its coefficients as control variables, and minimized the norm of the difference between the beam amplitude and a given constant value.

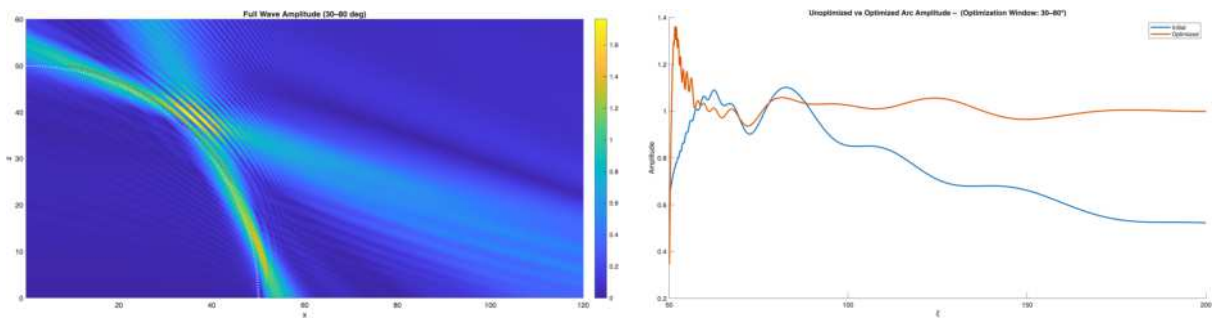

Figure 3: Beam optimized for constant intensity.

**Acknowledgement.** Work supported by US ARO, grant W911NF2410111.

## Spectral properties of waves in high-contrast media

**Yuri A. Godin**[1], **Leonid Koralov**[2] **Boris Vainberg**[1,*]
[1]Department of Mathematics and Statistics, Univ. of North Carolina at Charlotte, Charlotte, USA
[2]Department of Mathematics, University of Maryland, College Park, USA
*Email: brvainbe@charlotte.edu

**Abstract**

We study the spectral properties of elliptic operators with highly inhomogeneous coefficients and related issues concerning wave propagation in high-contrast media. A unified approach to solving problems in bounded domains with Dirichlet or Neumann boundary conditions is proposed, as well as in infinite periodic media. For a small parameter $\varepsilon > 0$ characterizing the contrast of the medium components, the analyticity of the eigenvalues and eigenfunctions is established in a neighborhood of $\varepsilon = 0$ despite the presence of the factor $1/\varepsilon$ in the operator under consideration. The effective operators corresponding to $\varepsilon = 0$ are described.

### 1 Introduction

Problems involving media with high-contrast inclusions arise in many areas of physics and engineering: modeling flow in a porous medium, wave propagation in perforated elastic media, homogenization of composites, description of the bandgap structure of photonic and phononic crystals, the calculation of thermal fields in materials with highly conductive components, impedance tomography, and so on. In such problems, a small parameter $\varepsilon$ naturally arises, equal to the ratio of the material parameters of the constituent components. For example, $\varepsilon$ could be the square of the ratio of the wave propagation speeds in the host material and the inclusion (or the ratio of the diffusion coefficients).

The presence of a small parameter $\varepsilon$ in the problem suggests using asymptotic methods; see [1–7, 10] for recent results on the asymptotic analysis of solutions in high-contrast media.

### 2 Model

For simplicity, we restrict ourselves here by considering a Dirichlet problem in a bounded domain $\Omega \subset \mathbb{R}^d$ containing a connected subdomain (inclusion) $\Omega_-$ (or a set of subdomains). Let $\Omega_+ = \Omega \setminus \overline{\Omega}_-$. We will write functions $u$ in $\Omega$ in the form $u = (u^+, u^-)$, where $u^\pm$ are the restrictions of $u = u(\boldsymbol{x})$ on $\Omega_\pm$. Consider the negative Laplacian in $\Omega$ with a piecewise-constant coefficient equal to $1, \varepsilon^{-1}$, $0 < \varepsilon \ll 1$, in $\Omega_\pm$, respectively.

Let $A_\varepsilon$ be the following operator in $L_2(\Omega)$:

$$A_\varepsilon u = (-\Delta u^+, -\frac{1}{\varepsilon}\Delta u^-), \quad \varepsilon > 0, \quad u \in \mathcal{H}_\varepsilon,$$

where the domain $\mathcal{H}_\varepsilon = \mathcal{H}_\varepsilon(\Omega)$ of $A_\varepsilon$ consists of functions $u$ such that $u^\pm$ belongs to Sobolev spaces $H^2(\Omega_\pm)$, respectively, and the following boundary conditions hold:

$$u^+ = 0, \quad \boldsymbol{x} \in \partial\Omega;$$
$$u^+ = u^-, \quad \frac{\partial u^-}{\partial \boldsymbol{n}} = \varepsilon \frac{\partial u^+}{\partial \boldsymbol{n}}, \quad \boldsymbol{x} \in \partial\Omega_-,$$

where $\boldsymbol{n}$ is outward with respect to $\Omega_-$, unit normal vector. Note that the boundary conditions on $\partial\Omega_-$ appear naturally when $A_\varepsilon$ is defined by the quadratic form involving the Laplacian with a piecewise-constant coefficient.

### 3 Results

The operator $A_\varepsilon$ has a discrete spectrum, its eigenvalues are squares of the time frequencies of the corresponding harmonic motions, i.e., the asymptotic analysis of the problem as $\varepsilon \downarrow 0$ is important. All available results (usually on the main term of asymptotics with an estimate of the reminder) are not simple since $A_\varepsilon$ does not have a limit as $\varepsilon \downarrow 0$: the domain $\mathcal{H}_\varepsilon$ of the operator $A_\varepsilon$ depends on $\varepsilon$ and the equation corresponding to $A_\varepsilon$ contains a factor of $1/\varepsilon$.

Despite this difficulty, we suggest a simple approach that allows us to prove the existence of the analytic extension of the eigenvalues and eigenfunctions of $A_\varepsilon$ from $\varepsilon > 0$ to a neighborhood of $\varepsilon = 0$ and provide an explicit description of the limit eigenvalues and eigenfunctions at $\varepsilon = 0$. We propose a unified approach that allows us to solve different high-contrast problems with the Dirichlet or Neumann boundary conditions, as well as periodic problems in $\mathbb{R}^d$ with the Bloch conditions, where we study dispersion relations.

We employ two simple tricks to solve high-contrast problems with any boundary conditions and obtain the result using fairly simple arguments. First, we study the eigenvalues and eigenfunctions not for $A_\varepsilon$ but for the inverse operator $B_\varepsilon = A_\varepsilon^{-1}$. Second, we study $B_\varepsilon$ by reducing the equation $A_\varepsilon u = f$ to an equation for the trace of $u$ on the boundary between the matrix and the inclusion. The latter equation is analytic in $\varepsilon$, and therefore only general simple statements from functional analysis are needed.

The characteristic equation for the limit eigenvalues (as $\varepsilon \downarrow 0$) is described using a non-local eigenvalue problem with Dirichlet data on the boundary of the inclusion equal to a constant determined by the equation containing the integral of the normal derivative of the solution. This type of problem has appeared in different publications; see, for example, [1, 6, 8, 9]. Here, we obtain this problem as one of the results of our general approach. The characteristic equation for the limit eigenvalues (as $\varepsilon \downarrow 0$) has the following form. We show that a function $u$ is a limit eigenfunction of the operator $A_\varepsilon$, as $\varepsilon \downarrow 0$, with the limit eigenvalue $\lambda$, if

(i) $u^- \equiv c_0$, $\boldsymbol{x} \in \Omega_-$, with some constant $c_0$,
(ii) $u^+$ is a solution of the problem

$$-\Delta u^+ = \lambda u^+, \quad u^+ \in H^2(\Omega_+), \tag{1}$$

$$u^+|_{\partial\Omega} = 0, \quad u^+|_{\partial\Omega_-} = c_0, \tag{2}$$

(iii) the following relation holds

$$\int_{\partial\Omega_-} \frac{\partial u^+}{\partial \boldsymbol{n}} \, \mathrm{d}S + c_0 \lambda |\Omega_-| = 0. \tag{3}$$

If $c_0 = 0$, then $\lambda$ is an eigenvalue of the Dirichlet problem (1), (2) with an eigenfunction $u$ such that $\int_{\partial\Omega_-} \frac{\partial u}{\partial \boldsymbol{n}} \, \mathrm{d}S = 0$. Thus, not an arbitrary eigenvalue of the Dirichlet problem (1), (2) (with $c_0 = 0$) is a limit eigenvalue of $A_\varepsilon$. If $c_0 \neq 0$, then $c_0$ can be canceled in (3), i.e., one can put $c_0 = 1$ in (i)-(iii), and equation (3) becomes a characteristic equation for additional limit eigenvalues $\lambda$.

All the above mentioned results can be found in [10].

# A potential theory on weighted graphs

**Fernando Guevara Vasquez**[1,*], **Trent DeGiovanni**[2]
[1]Department of Mathematics, University of Utah, Salt Lake City, UT, USA
[2]Currently: IDeaS, a SAS Company. Formerly: Mathematics Department, Dartmouth College, Hanover, NH, USA
*Email: fguevara@math.utah.edu

**Abstract**

We present a potential theory on weighted graphs analogous to the continuum potential theory, relying only on knowledge of the graph Laplacian and a partition of the nodes into interior, exterior and boundary. Our approach brings classic results to weighted graphs including a Calderón projector, Dirichlet to Neumann map identities and a spectral result for the Neumann Poincaré operator.

## 1 Introduction

Potential theory is an important tool to study and solve time harmonic wave problems in a homogeneous medium [5]. Although there are several potential theories on graphs [1, 2] and lattices [2, 4], few allow to use identities and results similar to those in the continuum. The approach we summarize here was introduced in [3] for weighted graph Laplacians and mimicks the continuum approach, operator-to-operator.

## 2 The Laplacian and Green function

We consider a finite undirected graph $\mathcal{G} = (\mathcal{V}, \mathcal{E})$ with vertices $\mathcal{V}$ and edges $\mathcal{E} \subset \mathcal{V} \times \mathcal{V}$. We work with a partition of the nodes into an external boundary $B$ (analogous to the radiation boundary condition at infinity in the continuum), an exterior $\Omega^+$, an interior $\Omega^-$ and a boundary $\partial\Omega$, as in Figure 1. The connectivity of the graph is such that the whole graph is connected, and the only connections between the node sets above are between $B$ and $\Omega^+$; between $\Omega^+$ and $\partial\Omega$; and between $\partial\Omega$ and $\Omega^-$. Consider a conductivity $\sigma : \mathcal{E} \to (0, \infty)$. The weighted graph Laplacian is $L(\sigma) = \nabla^T \text{diag}\,(\sigma)\nabla$, where $\text{diag}\,(\sigma)$ is the diagonal matrix representing multiplication by $\sigma$ and $\nabla$ is the discrete gradient operator that takes a function defined on the nodes $\mathcal{V}$ and returns the function defined on the edges $\mathcal{E}$ given by the difference in nodal values spanning an edge. The edge orientation is fixed ahead of time and does not affect the definition of $L(\sigma)$. The Green function on this graph is a matrix such that the $j-$th column is zero if $j \in B$ and the voltage obtained by injecting a unit current at the $j-$th node with homogeneous Dirichlet boundary conditions on the nodes $B$ if $j \notin B$.

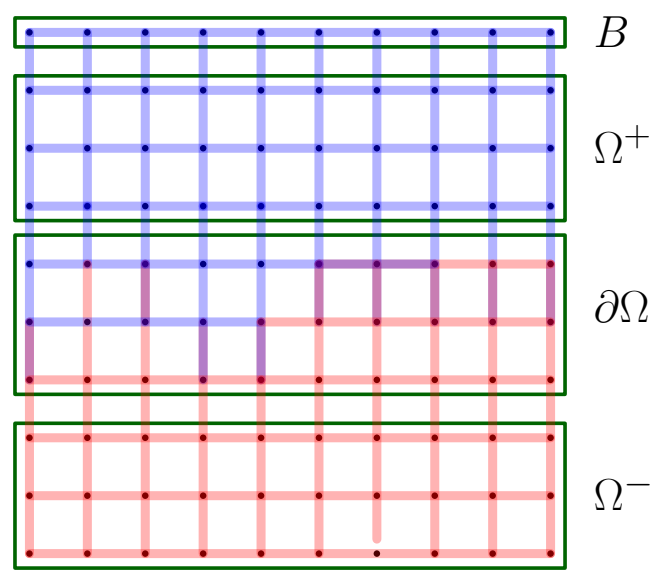

Figure 1: A graph showing the partition of the nodes and of the conductivity (edge colors).

## 3 Trace operators

In the continuum, Dirichlet $\gamma_0^\pm$ and Neumann $\gamma_1^\pm$ trace operators are defined by taking limits as one approaches the boundary from either the exterior (+) or the interior (−). In a weighted graph, it is natural to define the Dirichlet trace operators $\gamma_0^\pm = \gamma_0$ as the restriction to $\partial\Omega$ of a function on the nodes. To define the Neumann trace operators, we need to distinguish between currents (fluxes) flowing from $\Omega^-$ to $\partial\Omega$ and from $\partial\Omega$ to $\Omega^+$. We may choose to assign, in whole or in part, the fluxes within $\partial\Omega$ to either. To codify this choice, we introduce a partition of the conductivity, i.e. two edge functions $\sigma^\pm \geq 0$ such that $\sigma^+ + \sigma^- = \sigma$, $\sigma^+ = \sigma$ on all edges with a node in $B \cup \Omega^+$ and $\sigma^- = \sigma$ on all edges with a node in $\Omega^-$. In this way, $\gamma_1^- u = -(L(\sigma^+)u)|_{\partial\Omega}$ are the currents flowing from $\partial\Omega$ to $\Omega^+$ resulting from a voltage $u$ and similarly $\gamma_1^- u = (L(\sigma^-)u)|_{\partial\Omega}$ are the currents flowing from $\Omega^-$ to $\partial\Omega$. Unless otherwise noted, the results presented here are independent of the conductivity partition, provided it is fixed a priori.

## 4 Reproduction formulas

In the continuum, the single (resp. double layer) potential operator maps a density defined on the boundary to the superposition of monopoles (resp. dipoles) weighted by this density. On graphs, we define the single layer (resp. double layer) potential operator by $\mathcal{S} = G\gamma_0^T$ (resp. $\mathcal{D}^\pm = G[\gamma_1^\pm]^T$). We point out an important difference with the continuum: we use two double layer operators instead of a single one, each associated with the Neumann trace operators $\gamma_1^\pm$. These operators allow to reproduce a function from its Dirichlet and Neumann traces, provided the function is harmonic (i.e. the net currents are zero) on a particular set. See [3] for a proof.

**Theorem 1** *Let $u$ be a function defined on the nodes. Depending on the set $X$ where $(L(\sigma)u)|_X = 0$ we have the following.*

| $X$ | *reproduction formula* |
|---|---|
| $\Omega^- \cup \partial\Omega$ | $\mathcal{S}(\gamma_1^+ u) - \mathcal{D}^+(\gamma_0 u) = u\mathbb{1}_{\partial\Omega\cup\Omega^-}$ |
| $\Omega^-$ | $\mathcal{S}(\gamma_1^- u) - \mathcal{D}^-(\gamma_0 u) = u\mathbb{1}_{\Omega^-}$ |
| $\Omega^+$ | $-\mathcal{S}(\gamma_1^+ u) + \mathcal{D}^+(\gamma_0 u) = u\mathbb{1}_{B\cup\Omega^+}$ |
| $\Omega^+ \cup \partial\Omega$ | $-\mathcal{S}(\gamma_1^- u) + \mathcal{D}^-(\gamma_0 u) = u\mathbb{1}_{B\cup\Omega^+\cup\partial\Omega}$ |

Here we used $\mathbb{1}_X$ for the indicator function of a set and multiplication is componentwise.

## 5 Jump formulas

We define the boundary layer operators by taking Dirichlet or Neumann traces of $\mathcal{S}$ and $\mathcal{D}^\pm$. In particular, the single boundary layer operator is $S = \gamma_0\mathcal{S}$, the double boundary layer operator is $D = \frac{1}{2}\gamma_0(\mathcal{D}^+ + \mathcal{D}^-)$ (with adjoint $D' = D^T$) and the hypersingular operator $H = -\gamma_1^+\mathcal{D}^- = -\gamma_1^-\mathcal{D}^+$. The last equality is proven in [3]. These operators obey jump formulas, similar to the ones in the continuum, see [3].

**Theorem 2** *The following jump formulas hold*

$$\gamma_1^\pm\mathcal{S} = \mp\frac{I}{2} + D', \gamma_0\mathcal{D}^\pm = \mp\frac{I}{2} + D.$$

It is also possible to show that the matrices

$$P^\pm = \mp\frac{I}{2} + C, \text{ where } C = \begin{bmatrix} -D & S \\ H & D' \end{bmatrix},$$

are projectors, i.e. $(P^\pm)^2 = P^\pm$, analogous to the Calderón projectors in the continuum [3].

## 6 Other identities and results

Since $S$ is positive definite [3], we may consider two Dirichlet to Neumann maps: the map $\Lambda_\sigma^+$ (resp. $\Lambda_\sigma^-$) that associates the Dirichlet trace to the exterior (resp. interior) Neumann trace. These obey similar identities as in the continuum,

$$\begin{aligned}\Lambda_\sigma^\pm &= (\mp\frac{I}{2} + D')S^{-1} = S^{-1}(\mp\frac{I}{2} + D) \\ &= \mp H \mp (\mp\frac{I}{2} + D')S^{-1}(\mp\frac{I}{2} + D).\end{aligned}$$

and also agree with the common Schur complement definition [3]. Another result in [3] is that the spectrum of the Neumann-Poincaré operator $D'$ must be in $[-1/2, 1/2)$, under suitable assumptions on the graph connectivity.

## 7 Future work

The present approach is limited to finite graphs and the Laplacian. Computationally it only requires the Green function at the boundary and neighboring nodes. We are working on extending this approach to lattices and the discrete Helmholtz equation. We envision this technique could help formulate and solve a variety of scattering problems in lattices.

# Energy-preserving simulation of electromagnetic waves at oblique incidence using Split-Field FETD method

**Valentin Wasquel**[1,2,*], **Paula Aguilera**[1], **Michel Fournié**[2], **Denis Matignon**[2]
[1]CEA, Bordeaux, France
[2]Fédération ENAC ISAE-SUPAERO ONERA, Université de Toulouse, Toulouse, France
*Email: valentin.wasquel@cea.fr

## Abstract

We define a new formulation of the Maxwell's equations for oblique waves under the port-Hamiltonian system (pHs) framework, using a field transformation method. This new system was implemented using a mixed finite element method and validated against analytical solutions.

## 1 Introduction

Compared to spectral methods, time-domain methods allow to simulate broadband incident waves and non-linear systems. When simulating normal incident wave in a periodic medium, the numerical domain can be restricted to one unit cell by applying Periodic Boundary Conditions (PBCs) in the truncated direction. However, when the object is illuminated at an oblique incidence, PBCs cannot be implemented as some borders require knowledge of the future value of the field on the opposite border.

In 1998, Roden et al. [1] presented the Split-Field (SF) Finite-Difference Time-Domain method which involves a space-varying time-shift to remove the delay between borders. In this method, auxiliary fields are solved, but the classical electromagnetic energy is not directly accessible. Because the energy is the quantity that is conserved, better results can be expected from a method that guarantees the conservative property of the simulated system.

The pHs framework is a way of describing physical system, centered around the conservation and transfer of energy between its components. Since the early 2000, it has been used in numerous domains [2], and when coupled with a structure-preserving method [3], a mixed finite element method conserving the passivity of the original system can be obtained.

## 2 SF-PHS formulation

We consider a passive medium, characterized by its relative permittivity and permeability $(\varepsilon_r, \mu_r)$, periodic in the $y$ and $z$ directions. For simplicity, we consider an oblique incident plane wave, with wave vector $\vec{k}$ lying in the $(x, y)$-plane, and making an angle $\theta$ with the $x$-axis : $\vec{k} = \cos(\theta)k\vec{e}_x + \sin(\theta)k\vec{e}_y$.

The SF method [1] consists of applying a time offset to the considered field depending on the position. Instead of solving for $\vec{E}(\vec{x}, t)$ and $\vec{H}(\vec{x}, t)$, we define four auxiliary fields: the time shifted electric and magnetic fields : $\vec{P}(\vec{x}, t) := \vec{E}(\vec{x}, t+\alpha y)$ and $\vec{Q}(\vec{x}, t) := \vec{H}(\vec{x}, t+\alpha y)$ with $\alpha := \frac{n\sin(\theta)}{c_0}$, for $n$ the refractive index, $c_0$ the speed of light in vacuum; and the fields :

$$\begin{bmatrix}\vec{\mathcal{D}}\\ \vec{\mathcal{B}}\end{bmatrix} = \begin{bmatrix}\varepsilon_0\varepsilon_r & \alpha\vec{e}_y\times.\\ -\alpha\vec{e}_y\times. & \mu_0\mu_r\end{bmatrix}\begin{bmatrix}\vec{P}\\ \vec{Q}\end{bmatrix} = A\begin{bmatrix}\vec{P}\\ \vec{Q}\end{bmatrix}. \quad (1)$$

Although these fields are not the time-shifted electric displacement and magnetic flux, they satisfy equations which are formally identical to Maxwell's equations, namely:

$$\frac{\partial}{\partial t}\begin{bmatrix}\vec{\mathcal{D}}\\ \vec{\mathcal{B}}\end{bmatrix} = \begin{bmatrix}0 & \vec{\nabla}_\times.\\ -\vec{\nabla}_\times. & 0\end{bmatrix}\begin{bmatrix}\vec{P}\\ \vec{Q}\end{bmatrix}. \quad (2)$$

This system follows the formulation of a conservative pHs with the flow vector $f = \left[\vec{\mathcal{D}}^\top \ \ \vec{\mathcal{B}}^\top\right]^\top$ and the effort vector $e = \left[\vec{P}^\top \ \ \vec{Q}^\top\right]^\top = A^{-1}f$. It preserves a new energy defined as the Hamiltonian:

$$\mathcal{H}(t) := \frac{1}{2}\int_\Omega \left(\vec{P}\cdot\vec{\mathcal{D}} + \vec{Q}\cdot\vec{\mathcal{B}}\right)\mathrm{d}v = \frac{1}{2}\int_\Omega \begin{bmatrix}\vec{P}\\ \vec{Q}\end{bmatrix}^\top A\begin{bmatrix}\vec{P}\\ \vec{Q}\end{bmatrix}\mathrm{d}v$$

$$= \int_\Omega \left[\frac{\varepsilon_0\varepsilon_r}{2}\left\|\vec{P}\right\|^2 + \frac{\mu_0\mu_r}{2}\left\|\vec{Q}\right\|^2 - \alpha\left(\vec{P}\times\vec{Q}\right)\cdot\vec{e}_y\right]\mathrm{d}v.$$

The cross-product operator being skew-symmetric, the matrix $A$ is symmetric and positive definite; hence $\mathcal{H}(t)$ defines a positive energy. In the expanded expression of that energy, the first two terms are the electric and magnetic energies of the time-shifted field we expect to see in the canonical problem. The third term is proportional to a component of the time-shifted Poynting vector. This reasoning can be extended to dielectric media and Perfectly Matched Layers (PML).

## 3 SF-FETD method

To validate this methodology, a mixed finite element method is implemented. Based on the weak formulation of the system (2), we integrated by part one of the equations, see [3]. Integrating the equation on $\vec{P}$ gives:

*Find* $(\vec{\mathcal{D}}, \vec{P}, \vec{\mathcal{B}}, \vec{Q}) \in V$ *such that, for all test functions* $(\vec{\phi}_\mathcal{D}, \vec{\phi}_P, \vec{\phi}_\mathcal{B}, \vec{\phi}_Q) \in V$, *we have :*

$$\begin{cases} \displaystyle\int_\Omega \frac{\partial}{\partial t}\vec{\mathcal{D}} \cdot \vec{\phi}_\mathcal{D} = \int_\Omega \vec{\nabla}_\times \vec{Q} \cdot \vec{\phi}_\mathcal{D}, \\ \displaystyle\int_\Omega \frac{\partial}{\partial t}\vec{\mathcal{B}} \cdot \vec{\phi}_\mathcal{B} = -\int_\Omega \vec{\nabla}_\times \vec{\phi}_\mathcal{B} \cdot \vec{P} + \int_{\partial\Omega} \vec{u}_\partial \cdot \vec{\phi}_{\mathcal{P}_{//}}, \\ \displaystyle\int_\Omega \begin{bmatrix} \vec{\phi}_P & \vec{\phi}_Q \end{bmatrix} \begin{bmatrix} \vec{P} \\ \vec{Q} \end{bmatrix} = \int_\Omega \begin{bmatrix} \vec{\phi}_P & \vec{\phi}_Q \end{bmatrix} A^{-1} \begin{bmatrix} \vec{\mathcal{D}} \\ \vec{\mathcal{B}} \end{bmatrix}, \end{cases}$$

with $V := [H(div)]^2 \times [H(rot)]^2$, $\partial\Omega$ is the boundary of $\Omega$, $\vec{F}_{//}$ corresponds to the tangential component of the trace of $\vec{F}$ onto $\partial\Omega$ and $\vec{u}_\partial$ is a control we can use to impose the value of $\vec{P}$ on $\partial\Omega$. To preserve the pHs structure and to be able to monitor energy flux though the boundaries, we introduce a corresponding observable variable $\vec{y}_\partial$ that respects:
$\forall \vec{\phi}_y \in \left[L^2(\partial\Omega)\right]^2, \int_{\partial\Omega} \vec{y}_\partial \cdot \vec{\phi}_y = \int_{\partial\Omega} \vec{Q}_{//} \cdot \vec{\phi}_y$.

This system to consider is spatially discretized using Nédélec elements (N1curl) of degrees 0 for $\vec{Q}$ and $\vec{\mathcal{B}}$, Raviart–Thomas elements (RT0) for $\vec{P}$ and $\vec{\mathcal{D}}$, and the boundary variables are discretized using P1-Lagrange elements. Finally, the Crank–Nicolson method is used in time.

## 4 Numerical results

We consider a problem with 2 media : a small passive medium with $\varepsilon_r = 2.1$ and $\mu_r = 1$ and a layer of PLA-carbon, a dielectric material with loss. These two layers were placed on top of a PEC and an incident TE polarized planar-wave was simulated at various angles. We monitor the evolution of the energy in the simulation and compute the ratio of the reflected energy over the incident energy. This result was compared with solutions obtained with an analytical method (Leontovich impedance) and the commercial tool HFSS (an FEM-FD tool developed by ANSYS). The comparison is presented for $\theta = 20°$ in Fig. 1. As we can see, the SF-FETD method gives the same result as the other methods. However, this method also allows to precisely monitor the evolution of the energy in the system at each time step. The evolution of the energy is reported in Fig. 2. Thanks to the pHs formulation, we can compute the energy balance at each time-step and observe that once the source has been turned off, the energy in the system can only be dissipated.

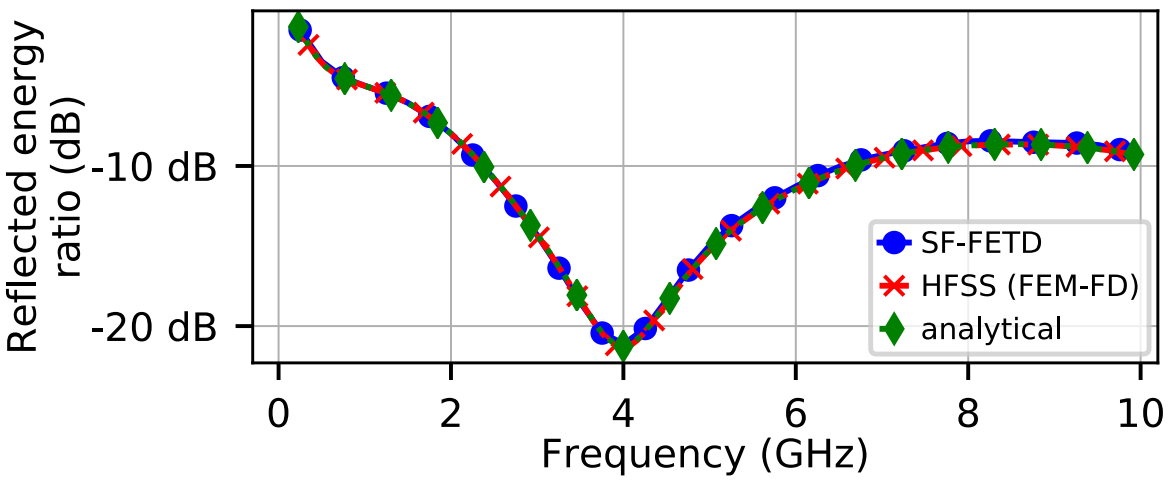

Figure 1: *Reflected energy ratio.*

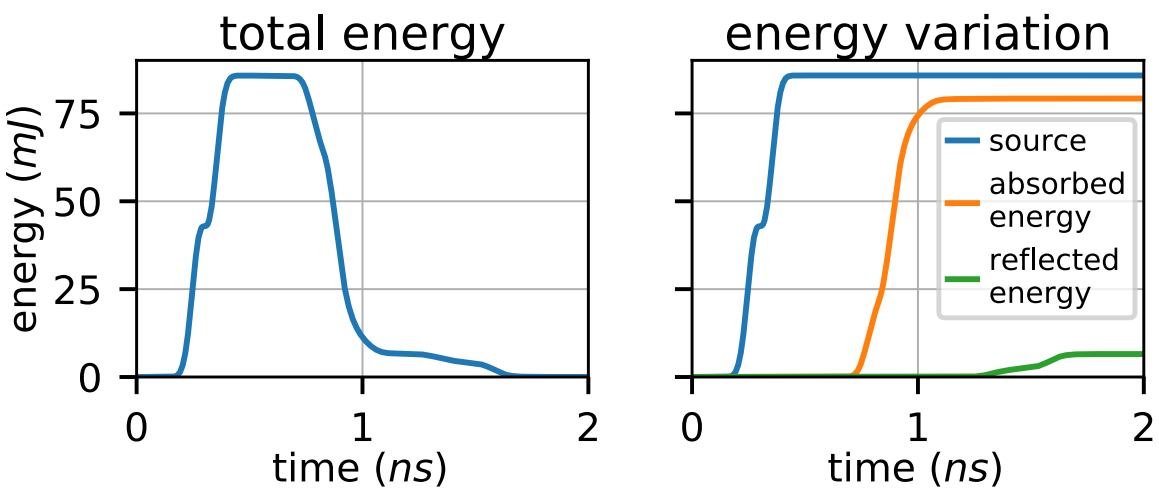

Figure 2: *Quantity of energy in the system at each time-step (L) and the contribution of the source or the dissipative terms (R).*

## 5 Conclusion

Using the Split-Field transformation method, we were able to write the Maxwell's equations for an incident electromagnetic wave at oblique angle as a port-Hamiltonian system. This formulation leads to the definition of a new energy useful to preserve the conservative property at the continuous and discrete levels.

# Finite element discretizations of hyperboloidal compactification layers for unbounded time-domain wave problems

Markus Wess[1,*], Anıl Zenginoğlu[2]
[1]Institute of Analysis and Scientific Computing, TU Wien, Vienna, Austria
[2]Institute for Physical Science and Technology, University of Maryland, Maryland, USA
*Email: markus.wess@tuwien.ac.at

**Abstract**

Time-domain simulations of wave propagation on unbounded domains require artificial boundary treatment to prevent spurious reflections. We explore finite element discretizations of hyperboloidal compactification layers for treating open boundaries in time-domain wave problems. By reformulating the governing equations on hyperboloidal spacelike slices that extend to future null infinity in an exterior layer, outgoing radiation is naturally transported off the computational domain. We review the underlying geometric and analytical ideas, discuss practical discretization strategies, and present numerical results demonstrating stable long-time evolution. Compared to perfectly matched layers (PMLs) this approach leads to systems with less auxiliary unknowns and straightforward energy decay without artificial absorption. Our findings suggest that hyperboloidal compactification offers a robust and conceptually clean alternative to PML methods for a broad class of wave propagation problems on unbounded domains.

## 1 Introduction

Hyperboloidal compactification layers were introduced in [1] for radially symmetric problems and discretized using one-dimensional spectral methods. In this work we focus on extending the method to general geometries using discretizations based on unstructured meshes.

In general we study time-domain acoustic or electromagnetic waves on $\tilde{\mathcal{D}} = \mathbb{R}^d \setminus \overline{\mathcal{O}}$ where $d \in \{1,2,3\}$ and the obstacle $\mathcal{O}$ is a bounded, sufficiently regular open set. We decompose the domain $\tilde{\mathcal{D}}$ disjointly into $\tilde{\mathcal{D}} = \tilde{\mathcal{D}}_{int} \cup \Gamma \cup \tilde{\mathcal{D}}_{ext}$ where $\tilde{\mathcal{D}}_{ext} = \{\tilde{x} \in \mathbb{R}^d : \|\tilde{x}\| > R\}$ for some $R > 0$ s.t. $\mathcal{O} \subset \mathbb{R}^d \setminus \tilde{\mathcal{D}}_{ext}$, $\tilde{\mathcal{D}}_{int} = \tilde{\mathcal{D}} \setminus \overline{\tilde{\mathcal{D}}_{ext}}$ and $\Gamma = \overline{\tilde{\mathcal{D}}_{int}} \cap \overline{\tilde{\mathcal{D}}_{ext}}$.

The method of hyperboloidal compactification layers now consists of the following three steps:

- **compactification**: the unbounded exterior $\tilde{\mathcal{D}}_{ext}$ is mapped via the mapping $\tilde{x} \mapsto x$ to the bounded layer $\mathcal{D}_{ext} = \{x \in \mathbb{R}^d : R < \|x\| < S\}$ for a suitable $S > R$.
- **time transformation**: a space dependent time transformation $t(\tau, \tilde{x}) = \tau + h(\|\tilde{x}\|)$ for a suitable positive function $h$ is introduced, in a way that the transformed equal-time slices (i.e., the level sets of $\tau$) approach the outgoing characteristics.
- **rescaling of the solution**: the solution is rescaled by a scalar function $\omega(x)$ which takes out the algebraic decay of the outgoing solution.

The three free parameters of the method, i.e., the compactification $\tilde{x}(x)$, the time transformation $h$ and the rescaling $\omega$ have to be chosen such that the resulting weak formulation and its solutions are sufficiently regular and feasible for numerical approximation.

## 2 Weak formulation of the transformed wave equation in compactified space

For simplicity we describe the method for the homogeneous, sourceless wave equation. We apply the three steps outlined above to

$$\partial_t^2 U(t, \tilde{x}) - \Delta U(t, \tilde{x}) = 0, \quad \tilde{x} \in \tilde{\mathcal{D}}$$

with appropriate initial data and homogeneous Neumann boundary data on $\partial\tilde{\mathcal{D}}$.

Now the rescaled and time transformed solution

$$u(\tau, x) = \frac{1}{\omega(x)} U(t(\tau, \tilde{x}(x)), \tilde{x}(x))$$

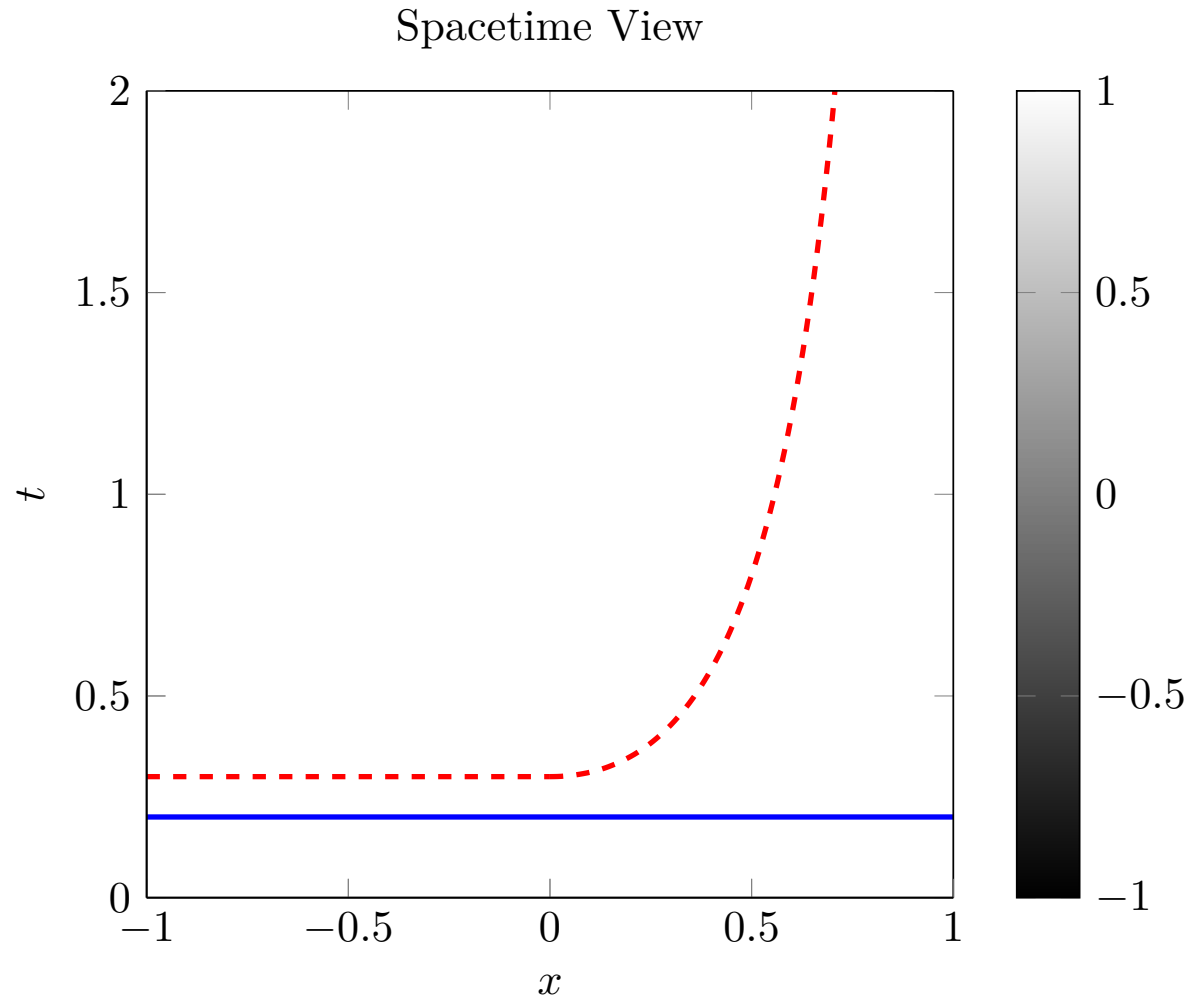

Figure 1: One-dimensional space-time wave in compactified coordinates with $\mathcal{D}_{int} = (-1,0), \mathcal{D}_{ext} = (0,1)$. The solution along a standard time slice (blue line) is highly oscillatory at mapped infinity ($x = 1$). The solution along a time slice with respect to the time transformed variable (dashed red line) is regular.

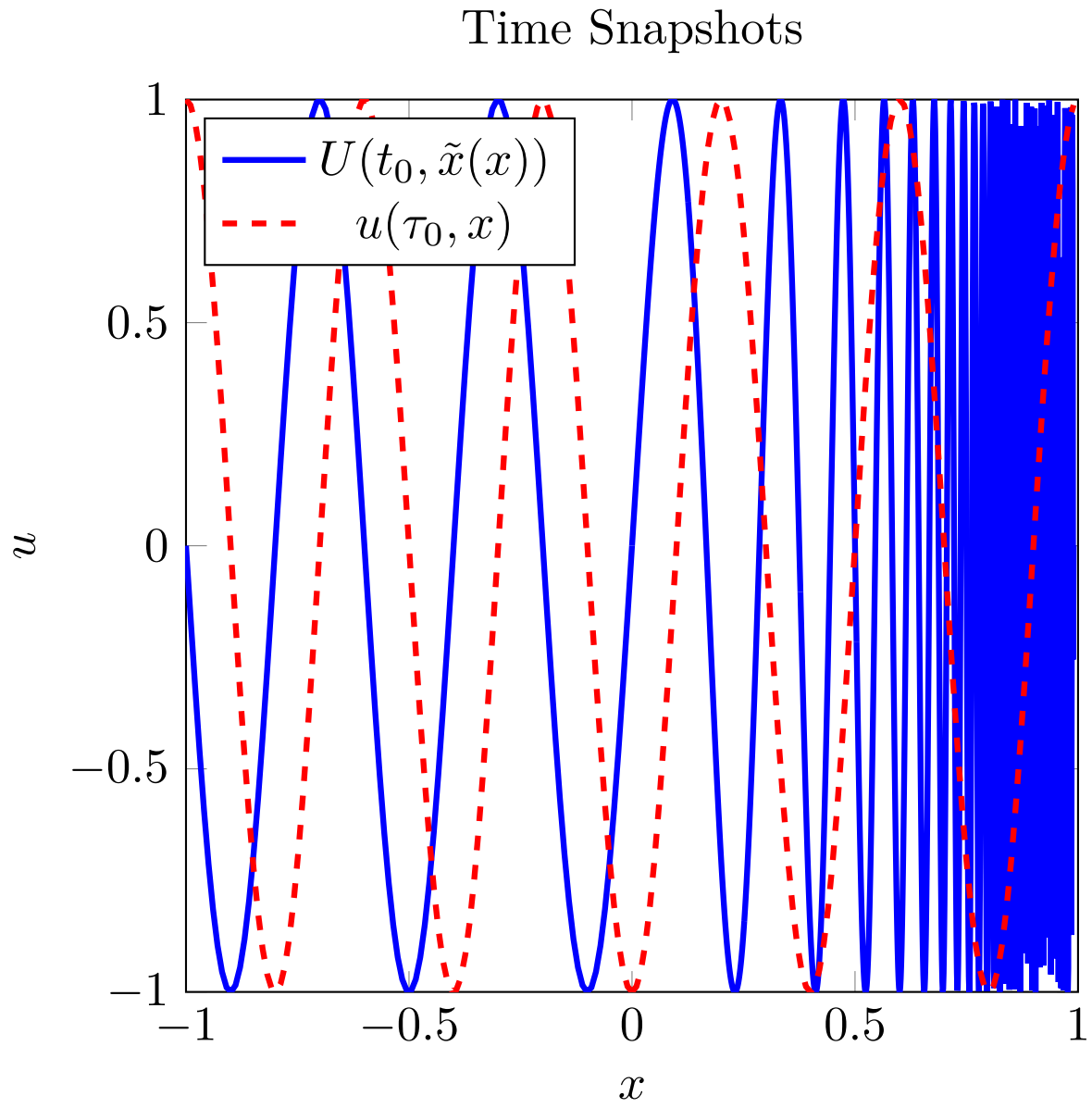

Figure 2: Solutions along an original time-slice (blue line) and transformed time-slice (dashed red line) from the example from Figure 1.

on the compactified domain $\mathcal{D}$ satisfies the weak formulation

$$\begin{aligned}&\int_{\mathcal{D}} a\partial_\tau^2 uv + \nabla(\omega u)^\top B\nabla(\omega v)\\ &= \int_{\mathcal{D}} C^\top \left(\nabla\partial_\tau uv - \partial_\tau u\nabla v\right) - \int_{\partial\mathcal{D}\setminus\partial\mathcal{O}} d\partial_\tau uv\end{aligned} \tag{1}$$

for all $v \in H^1(\mathcal{D})$. The coefficients $a, d \in \mathbb{R}, B \in \mathbb{R}^{d\times d}, C \in \mathbb{R}^d$ depend on the spacial variable, as well as the choices of the compactification, time-transformation and rescaling and under suitable conditions can be shown to be positive and guarantee regular and well-behaved solutions (cf. Figures 1 and 2 for a one-dimensional example).

## 3 Energy decay

Testing (1) with $\partial_\tau u$ immediately yields

$$\begin{aligned}\partial_\tau \mathcal{E}(u)^2 := \partial_\tau \frac{1}{2}\int_{\mathcal{D}} a|\partial_\tau u|^2 + \nabla(\omega u)^\top B\nabla(\omega u)\\ = -\int_{\partial\mathcal{D}\setminus\partial\mathcal{O}} d|\partial_\tau u|^2\end{aligned} \tag{2}$$

which, given the correct signs of the coefficients $a, B, d$ gives that energy $\mathcal{E}$ decays and flows out at mapped infinity $\partial\mathcal{D} \setminus \partial\mathcal{O}$.

## 4 Space and time discretization

The resulting weak formulation (1) is posed on the bounded compactified spacial domain $\mathcal{D}$ and thus can be discretized in space using standard conforming finite elements. We choose a suitable symplectic explicit or implicit time-stepping for the second order in time equation to guarantee discrete stability.

Alternatively one may introduce an auxiliary field to obtain a first order system in time, which can be discretized using various choices of mixed finite elements and again appropriate time-stepping.

We show the applicability of both approaches using the finite element library Netgen/NGSolve [2] and conduct numerical experiments to study stability and performance compared to perfectly matched layers.

# Derivation of polarization dependent interface conditions for Maxwell's equations

David Wiedemann[1,*], Ben Schweizer[1], Bérangère Delourme[2]
[1]Department of Mathematics, TU Dortmund University, Dortmund, Germany
[2]LAGA, Université Sorbonne Paris Nord, Villetaneuse, France
*Email: david.wiedemann@tu-dortmund.de

**Abstract**

We investigate the time-harmonic Maxwell equations in complex geometries, with particular emphasis on configurations that model polarization filters or Faraday cages. The domains under consideration contain perfectly conducting inclusions that are periodically distributed along a surface. We analyze the asymptotic regime in which both the size of the inclusions and the distance between them tend to zero, and we derive the resulting effective interface conditions. Depending on the geometric properties of the inclusions, the limiting system exhibits markedly different behaviors: perfect transmission, perfect reflection or polarization effects.

## 1 Introduction

The transmission of electromagentic (EM) waves through thin periodic metal structures that are smaller than the wavelength depends crucially on the topological properties of the metal obstacle: the metal grid in the window of a microwave oven reflects EM waves ($\lambda \approx 12$cm). Vice versa tiny slits ($\approx 0.05$mm) in a metal coating yield an EM wave transmission for radio waves ($\lambda \approx 8$–43cm). A metal microstructure which is only connected in one tangential direction leads to polarization dependent reflection.

## 2 Microscopic model

We consider the time-harmonic Maxwell equations with unknowns $E^\eta, H^\eta$ in a domain $\Omega_\eta = \Omega \setminus \Sigma_\eta$ for a small parameter $\eta > 0$, which scales the microstructure:

$$\operatorname{curl} E^\eta = i\omega\mu H^\eta + f_h \quad \text{in } \Omega_\eta\,, \tag{1a}$$

$$\operatorname{curl} H^\eta = -i\omega\varepsilon E^\eta + f_e \quad \text{in } \Omega_\eta\,, \tag{1b}$$

$$E^\eta \times \nu = 0 \quad \text{on } \partial\Omega_\eta\,. \tag{1c}$$

The boundary condition (1c) models a perfect conductor outside of $\Omega$ and in $\Sigma_\eta$.

Macroscopic geometry:

$$\Omega = (0,1)^2 \times (-1,1)\,,$$
$$\Gamma = (0,1)^2 \times \{0\}\,,$$

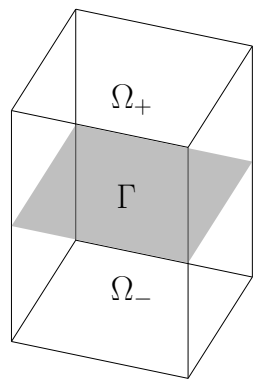

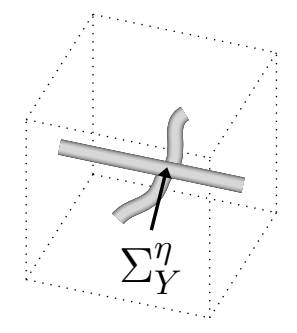

Cell geometry:

$$\Sigma_Y^\eta \subset Y = [0,1]^3\,.$$

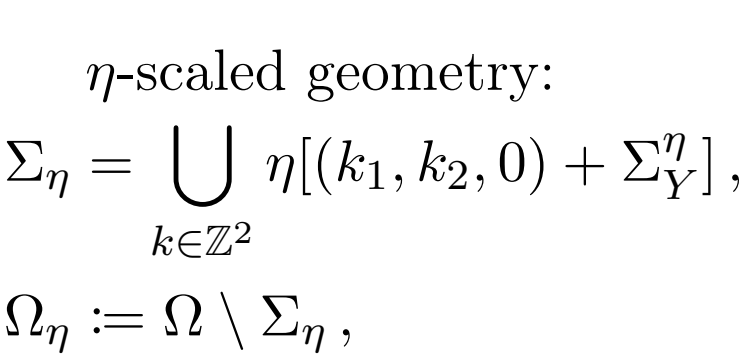

$\eta$-scaled geometry:

$$\Sigma_\eta = \bigcup_{k\in\mathbb{Z}^2} \eta[(k_1,k_2,0) + \Sigma_Y^\eta]\,,$$
$$\Omega_\eta := \Omega \setminus \Sigma_\eta\,,$$

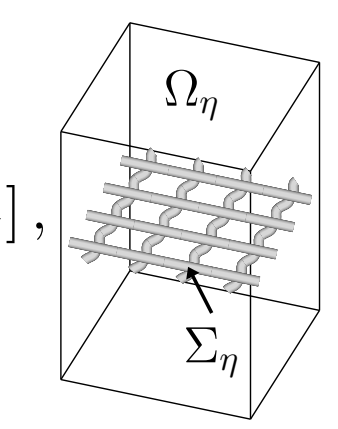

## 3 Asymptotic connectedness

The interaction of EM waves with thin obstacles is determined by the connectivity of the obstacles: EM waves with electric fields orientated parallel to the directions in which the obstacles are **connected** are **reflected**. EM waves with electric fields parallel to the directions in which the obstacles are **disconnected** are **transmitted**.

To define an analytic concept for the connectivity of a sequence of obstacles $\Sigma_Y^\eta \subset Y$, we use the following domains and function space:

$$Z := (0,1)^2 \times \mathbb{R}\,,$$
$$Z^\eta := Z \setminus \Sigma_Y^\eta\,,$$
$$\overset{\bullet}{H}(\operatorname{curl}, Z) := \Big\{ u|_Z \;\Big|\; u \in H_{\mathrm{loc}}(\operatorname{curl}, \mathbb{R}^3) \text{ is } e_1\text{- and } e_2\text{-periodic} \Big\}$$

**Definition 1** *Let $i \in \{1,2\}$ and $j = 3 - i$. A sequence of obstacles $\Sigma_Y^\eta$ is called asymptotically connected in the direction $e_i$ if there exists a sequence $\Psi_\eta^{(i)} \in \overset{\bullet}{H}(\operatorname{curl}, Z)$ such that:*

$$\eta^{1/2} \| \Psi_\eta^{(i)} - \operatorname{sgn}(x_3) e_j \|_{L^2(Z^\eta)} \to 0$$
$$\eta^{-1/2} \| \operatorname{curl} \Psi_\eta^{(i)} \|_{L^2(Z^\eta)} \to 0$$

The notion of *asymptotically disconnected* in a tangential direction can be defined in a conceptually analogous way.

## 4 Effective interface conditions

Independent on the geometry of the obstacles $\Sigma_Y^\eta$, one obtains the limit system:

$$\operatorname{curl} E^{\text{hom}} = i\omega\mu H^{\text{hom}} + f_h \quad \text{in } \Omega\,, \tag{2a}$$

$$\operatorname{curl} H^{\text{hom}} = -i\omega\varepsilon E^{\text{hom}} + f_e \quad \text{in } \Omega\setminus\Gamma\,, \tag{2b}$$

$$E^{\text{hom}}\times\nu = 0 \quad \text{on } \partial\Omega\,. \tag{2c}$$

System (2) is incomplete; while (2a) implies the continuity of the tangential components of $E^{\text{hom}}$ across $\Gamma$, the tangential components of $H^{\text{hom}}$ can jump across $\Gamma$. The remaining interface conditions, which complete (2), are derived using using the cell functions of the definition of asymptotic (dis-)connectedness. We denote the jump across $\Gamma$ by $[\![\cdot]\!]_\Gamma$.

**Theorem 2 (Homogenization result [1])** *Let $(E^\eta, H^\eta)\in L^2(\Omega_\eta,\mathbb{C}^3)\times L^2(\Omega_\eta,\mathbb{C}^3)$ be solutions of (1) such that their trivial extension to $\Omega$ converge weakly to some $(E^{\text{hom}}, H^{\text{hom}})\in L^2(\Omega,\mathbb{C}^3)\times L^2(\Omega,\mathbb{C}^3)$. Then, $(E^{\text{hom}}, H^{\text{hom}})$ solves the limit system (2).*

*Case 1: **Reflecting**. If $\Sigma_Y^\eta$ is asymptotically connected in directions $e_1$ and $e_2$*

$$\Rightarrow \quad E_1^{\text{hom}}|_\Gamma = E_2^{\text{hom}}|_\Gamma = 0\,.$$

*Case 2: **Inactive**. If $\Sigma_Y^\eta$ is asymptotically disconnected in directions $e_1$ and $e_2$*

$$\Rightarrow [\![H_1^{\text{hom}}]\!]_\Gamma = [\![H_2^{\text{hom}}]\!]_\Gamma = 0\,.$$

*Case 3: **Polarizing**. Let $i\in\{1,2\}$ and $j = 3-i$. If $\Sigma_Y^\eta$ is asymptotically connected in the direction $e_i$ and asymptotically disconnected in the direction $e_j$*

$$\Rightarrow E_i^{\text{hom}}|_\Gamma = [\![H_i^{\text{hom}}]\!]_\Gamma = 0\,.$$

## 5 Wire geometries

We study the asymptotic (dis-) connectedness of thin wires with small gaps, that are characterized by a sequence of radii $r_\eta$ and intervals $I_\eta$ representing the gaps:

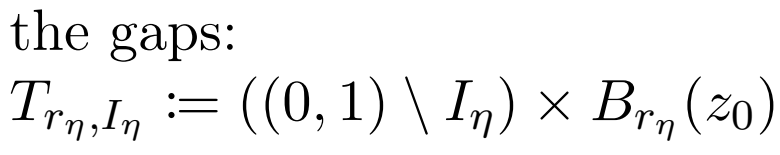

$T_{r_\eta,I_\eta} := ((0,1)\setminus I_\eta)\times B_{r_\eta}(z_0)$

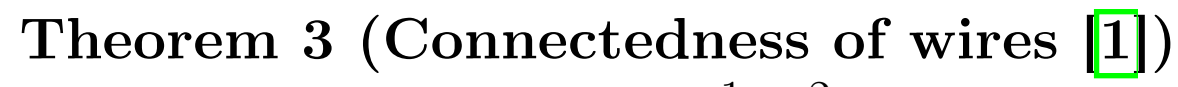

**Theorem 3 (Connectedness of wires [1])**

- *If $\eta|\ln(r_\eta)|\to 0$ and $\eta^{-1}r_\eta^{-2}|I_\eta|\to 0$, then $T_{r_\eta,I_\eta}$ is asymptotically connected in $e_1$.*
- *Wires aligned parallel to the $e_1$ direction are asymptotically disconnected in $e_2$.*
- *If $\eta|\ln(r_\eta)|\to\infty$, then $T_{r_\eta} = T_{r_\eta,\emptyset}$ is asymptotically disconnected in $e_1$.*

## 6 The polarization filter interface

The limit systems of Case 1 and Case 2 in Theorem 2 are classical cavity problems in $\Omega_+$, $\Omega_-$ and $\Omega$ and, thus, well-understood. The polarization interface condition of Case 3 cannot be addressed by classical methods. Nevertheless, a reformulation in a system of 3D- and 1D-Helmholtz equations yields the following Fredholm alternative:

**Theorem 4 (Polarization interface [2])** *Let $\Omega = (0,l_1)\times(0,l_2)\times(0,l_3)$ with an interface $\Gamma$ in the $x_1$–$x_2$-plane. For the system (2) with $E_1^{\text{hom}}|_\Gamma = [\![H_1^{\text{hom}}]\!]_\Gamma = 0$ and lateral periodic boundary conditions, $\omega^2\notin\{4\pi^2k^2(\varepsilon\mu l_1^2)^{-1}\mid k\in\mathbb{N}\}$, exactly one of the following statements hold:*

- *for $f_h = f_e = 0$ there exists a non-trivial solution.*
- *for every $f_h, f_e\in L^2(\Omega,\mathbb{C}^3)$ there exists a solution.*

So far this approach is restricted to scalar-valued and constant coefficients $\varepsilon$ and $\mu$ as well as cuboid domains $\Omega$ with a flat interface $\Gamma$.

Related difficulties appear at polarization dependent boundary conditions. For an impedance boundary condition with a vanishing impedance for one polarization, a Fredholm alternative is shown in [3].

# Transcranial Ultrasound Simulations using Volume-Surface Integral Equations

**Elwin van 't Wout**[1,*], **Alberto Almuna-Morales**[1], **Danilo Aballay**[1], **Ignacio Labarca-Figueroa**[1], **Pierre Gélat**[2], **Reza Haqshenas**[3]

[1]Institute for Mathematical and Computational Engineering, School of Engineering and Faculty of Mathematics, Pontificia Universidad Católica de Chile, Santiago, Chile

[2]Department of Surgical Biotechnology, Division of Surgery and Interventional Science, University College London, London, United Kingdom

[3]Department of Mechanical Engineering, University College London, London, United Kingdom

*Email: e.wout@uc.cl

**Abstract**

Transcranial ultrasound is a non-invasive therapy that uses low-intensity ultrasound to target the brain. The acoustic transmission through the skull may shift the focus and attenuate the field. We designed a volume-surface integral equation that models acoustic propagation through a heterogeneous medium in an unbounded domain. The simulations on CT-derived skull geometries quantify the focal aberrations due to material heterogeneity at high frequencies.

## 1 Introduction

Therapeutic ultrasound is a non-invasive treatment that has been approved for various health conditions, including neuromodulation. Acoustic fields are transmitted through the human body, delivering sufficient energy to the target region to produce the desired mechanical or thermal effects in support of the treatment's clinical objectives. These acoustic fields must be sharply localized to avoid damage to healthy tissue or organs surrounding the lesion. A bone's presence can significantly impact treatment by attenuating the acoustic field or aberrating the focal region. Hence, each patient's treatment must be planned carefully to comply with safety guidelines. Computational simulations could play a pivotal role in assessing focused ultrasound, but only when they deliver accurate predictions of the acoustic energy field fast enough on available computer hardware [1].

## 2 Methodology

The Helmholtz equation is the most widely used mathematical model to simulate therapeutic ultrasound, especially in scenarios with relatively low acoustic intensities where linear propagation models suffice. Most numerical methods to solve this partial differential equation are variants of ray tracing, volumetric grid-based solvers (finite differences, finite elements, pseudospectral, etc.), or boundary integral equations [2].

This study considers an integral equation approach. In general, the boundary element method (BEM) efficiently solves high-frequency scattering at high-contrast interfaces, handling unbounded domains naturally while providing accurate solutions on meshes with relatively few elements per wavelength. However, this approach has been restricted to piecewise constant domains so far, severely limiting its scope in therapeutic ultrasound, where material parameters are heterogeneous and derived from medical images. This study introduces a volume-surface integral equation (VSIE) that allows for local heterogeneous domains, and maintains BEM's high accuracy for acoustic scattering in unbounded domains [3].

## 3 Formulation

The Helmholtz equation reads

$$-\rho(\mathbf{x})\,\nabla\cdot\left(\frac{1}{\rho(\mathbf{x})}\nabla p(\mathbf{x})\right)-\left(\frac{\omega}{c(\mathbf{x})}\right)^2 p(\mathbf{x})=s(\mathbf{x})$$

for the unknown pressure $p(\mathbf{x})$ and known source $s(\mathbf{x})$, density $\rho(\mathbf{x})>0$, speed of sound $c(\mathbf{x})>0$, and angular frequency $\omega$. Let us consider a bounded region $\Omega$ embedded in an unbounded domain $\Omega_{\text{ext}}=\mathbb{R}^3\setminus\Omega$, representing bone and soft tissue, respectively. The material interface $\Gamma$ is Lipschitz continuous with the unit normal vector $\hat{\mathbf{n}}$ pointing outwards. In this study, we assume $\rho(\mathbf{x})=\rho_0$ and $c(\mathbf{x})=c_0$ to be constant in $\Omega_{\text{ext}}$ and $s(\mathbf{x})=0$ in $\Omega$. Hence, we only have external acoustic sources, and the bone region $\Omega$ may be heterogeneous. We impose continuity of the pressure field and normal particle velocity at the bone-tissue interface, and close the system with the Sommerfeld radiation condition.

The Green's function

$$G_0(\mathbf{x},\mathbf{y}) = \frac{e^{\imath(\omega/c_0)|\mathbf{x}-\mathbf{y}|}}{4\pi|\mathbf{x}-\mathbf{y}|} \quad \text{for} \quad \mathbf{x} \neq \mathbf{y}$$

is valid for the homogeneous exterior Helmholtz equation. The main idea behind the VSIE is to apply the representation theorem to the Helmholtz transmission problem and to retain the Newton potential in the heterogeneous interior. This yields a coupled system of volume integral equations in the heterogeneous bone and a surface integral at the high-contrast bone-tissue interface [4]. Precisely,

$$\begin{bmatrix} M_\Gamma + K_{\Gamma,\Gamma} & -V_{\Gamma,\Omega} + T_{\Gamma,\Omega} \\ K_{\Omega,\Gamma} & M_\Omega - V_{\Omega,\Omega} + T_{\Omega,\Omega} \end{bmatrix} \begin{bmatrix} p_\Gamma \\ p_\Omega \end{bmatrix} = \begin{bmatrix} f_\Gamma \\ f_\Omega \end{bmatrix},$$

must be solved for

$$M_\Sigma[\phi](\mathbf{x}) = \tfrac{\rho_0}{\rho(\mathbf{x})}\phi(\mathbf{x}),$$

$$V_{\Sigma,\Omega}[\phi](\mathbf{x}) = \int_\Omega G_0(\mathbf{x},\mathbf{y})\tfrac{\rho_0}{\rho(\mathbf{y})}\big((\tfrac{\omega}{c(\mathbf{y})})^2 - (\tfrac{\omega}{c_0})^2\big)\phi(\mathbf{y})\mathrm{d}\mathbf{y},$$

$$K_{\Sigma,\Gamma}[\phi](\mathbf{x}) = \int_\Gamma \big(\partial_{\hat{\mathbf{n}}(\mathbf{y})}G_0(\mathbf{x},\mathbf{y})\big)\big(\tfrac{\rho_0}{\rho(\mathbf{y})} - 1\big)\phi(\mathbf{y})\mathrm{d}\mathbf{y},$$

$$T_{\Sigma,\Omega}[\phi](\mathbf{x}) = \nabla_\mathbf{x} \cdot \int_\Omega G_0(\mathbf{x},\mathbf{y})\,\nabla_\mathbf{y}\big(\tfrac{\rho_0}{\rho(\mathbf{y})}\big)\phi(\mathbf{y})\mathrm{d}\mathbf{y},$$

$$f_\Sigma(\mathbf{x}) = \int_{\Omega_{\text{ext}}} G_0(\mathbf{x},\mathbf{y})s(\mathbf{y})\mathrm{d}\mathbf{y}$$

at $\mathbf{x} \in \Sigma$ and $\Sigma \in \{\Gamma, \Omega\}$. These operators denote the scaled mass, single-layer, double-layer, adjoint double-layer integral operators and incident wave field, respectively.

The VSIE is discretized using a collocation method with piecewise-constant basis functions and solved iteratively with GMRES. We also designed a hierarchical matrix compression algorithm to limit the memory footprint.

## 4 Results

We calculated the density and speed of sound from CT images of a human skull [4]. The hexahedral mesh has 199 693 voxels of size $0.489 \times 0.489 \times 0.5$ mm. The density and speed of sound are on average 1788 kg/m$^3$ and 2648 m/s in the heterogeneous bone, while $\rho_0 = 1000$ kg/m$^3$ and $c_0 = 1500$ m/s are typical values for soft tissue. The bowl transducer emits a focused field at 500 kHz, yielding a 3 mm wavelength in the exterior and a longer wavelength in the heterogeneous interior. Hence, we have at least 6 elements per wavelength.

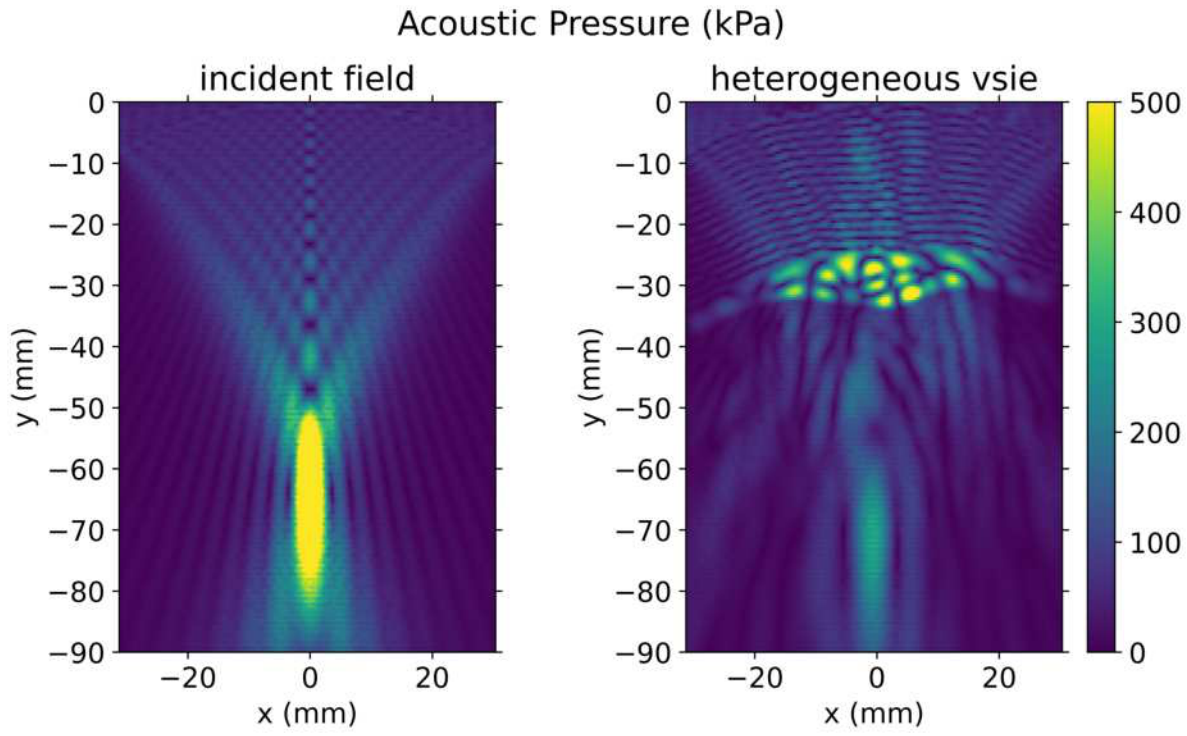

Figure 1: The focused acoustic field, propagating from top to bottom: incident field (left) and total field with presence of skull bone (right).

The VSIE-calculated fields in Figure 1 show the skull's acoustic influence: attenuation in the brain region, a shift in focal peak location, and high pressure in the skull bone. Additional grid-convergence studies and comparisons with BEM confirm the accuracy of the VSIE. Finally, hierarchical matrix compression reduces the memory footprint to a fraction of the original dense matrix, making ultrasound simulations on standard computers feasible.

# Destabilization of traveling waves in a nonlinear lattice with nonreciprocal coupling

**Andrus Giraldo**[1], **Stefan Ruschel**[2], **Behrooz Yousefzadeh**[3], [,*]
[1]School of Computational Sciences, Korea Institute for Advanced Study, Seoul, South Korea
[2]School of Mathematical Sciences, University of Nottingham, Nottingham, United Kingdom
[3]Department of Mechanical, Industrial and Aerospace Engineering, Concordia University, Montreal, Canada
[*]Email: behrooz.yousefzadeh@concordia.ca

## Abstract

We explore the propagation and stability of traveling waves in a ring of nonlinear oscillators with local nonreciprocal coupling. We develop a numerical framework based on delay-advance differential equations that allows for formulating a two-point boundary-value problem to compute the spectrum of a traveling wave independent of the lattice size. We present certain wave features that emerge from the interplay between nonreciprocal coupling and nonlinearity, along with the relevant instability mechanisms.

## 1 Introduction

Materials with odd elasticity possess anti-symmetric components in their elastic tensor, which account for unconventional wave propagation properties [1]. A balanced interplay between odd elasticity and mechanical damping structures allows for sustained free wave propagation in a very limited parameter range [2]. Here, we explore the combined effects of nonlinearity and odd elasticity on the existence and stability of traveling waves through a discrete model where odd elasticity manifests as local nonreciprocal interaction.

## 2 Model

We consider a discrete model of a mechanical system with nonreciprocal coupling springs: an arrangement of masses on a ring. In terms of non-dimensional parameters, the governing equations are

$$\ddot{x}_n + 2\zeta\dot{x}_n + x_n + \kappa(2x_n - x_{n+1} - x_{n-1}) + \alpha\kappa(x_{n+1} - x_{n-1}) + \beta x_n^3 = 0 \quad (1)$$

where $x(t)$ denotes the displacement of each degree of freedom from its static equilibrium position, overdot denotes differentiation with respect to time, $\zeta$ is the damping ratio, $\kappa$ represents the regular coupling force between adjacent cite (springs connecting masses), $\alpha$ is the nonreciprocity parameter, and $\beta$ is the coefficient of the nonlinear spring.

Parameter $\alpha$ is a measure of the degree of nonreciprocity in coupling (odd elasticity), with $\alpha = 0$ corresponding to a reciprocal connection (regular mechanical spring) and $0 < \alpha < 1$ corresponding to a nonreciprocal connection. When $\alpha = 1$, each mass only feels the spring on its right-hand side, and the spring constant on the left becomes negative when $\alpha > 1$. We only consider $0 \leq \alpha \leq 1$ in this work.

## 3 Methods and Results

We determine the stability of a traveling wave by computing its master stability curve [3]. This curve contains the complete spectrum of the Floquet exponents $\lambda$ of a periodic traveling wave in a lattice, as well as the spectra of all the larger networks embedding the original lattice. The master stability curves are computed using numerical continuation of a suitable two-point boundary value problem of a delay differential equation of the mixed type. We show the results for a lattice of 5 oscillators with $\zeta = 0.05$, $\kappa = 0.1$ and $\beta = 0.1$.

Fig. 1 shows the master stability curves associated with traveling waves at $\alpha = 0.8$, with modes $k = 1/5$ (wave number $q = 2\pi/5$). This traveling wave with $k = 1/5$ can be embedded in networks of size $5m$, with $m = 1, 2, \ldots$. In particular, the black circle denotes the Floquet exponent $\lambda$ of the traveling wave for a network of size 5 (its minimum admissible embedding), while the solid dots represent the corresponding Floquet exponents when the traveling wave is embedded in a network of size 25, obtained by concatenating five copies of the wave. These Floquet exponents – and those of any higher-dimensional embedding of the same wave – lie on the same master stability curve. For the example in Fig. 1, the traveling wave is unstable in its smallest embedding (5 oscillators) because two of its Floquet exponents have positive real parts. As a result, all the higher embeddings of this wave remain unstable.

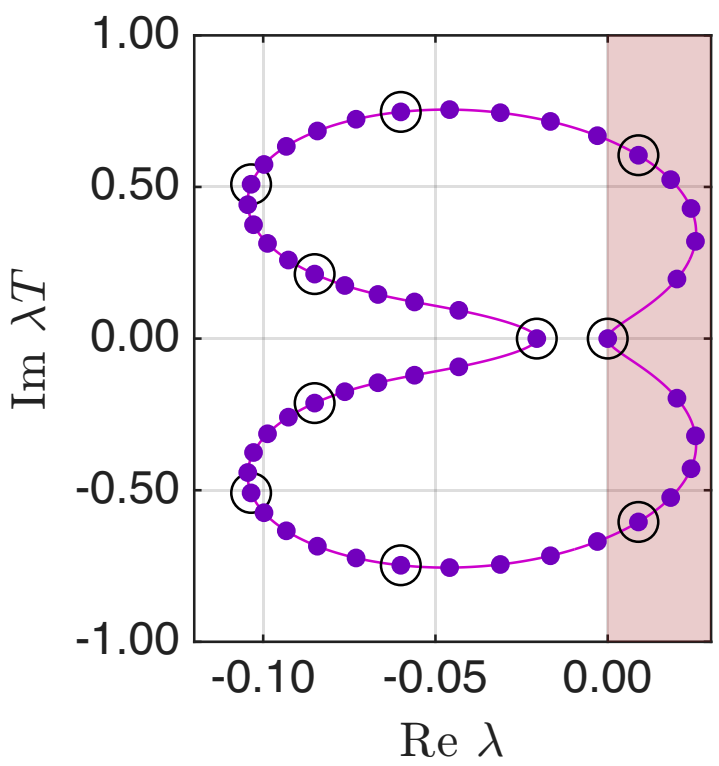

Figure 1: Master stability curves of traveling waves for $\alpha = 0.8$.

Because the master stability curves always pass through $\lambda = 0$, its zero-curvature loci will delimit a region in which the Floquet multipliers have strictly non-positive real parts. By computing these loci, we can determine regions in the parameter space that correspond to a change in stability. An example of this transition is shown in Fig. 2, where the master stability curve in the middle panel lies entirely in the negative half-plane. The associated traveling waves are stable for all lattice sizes at this wavenumber. This perspective toward the stability of traveling waves is analogous to destabilization of traveling waves via an *Eckhaus instability* in partial differential equations [4].

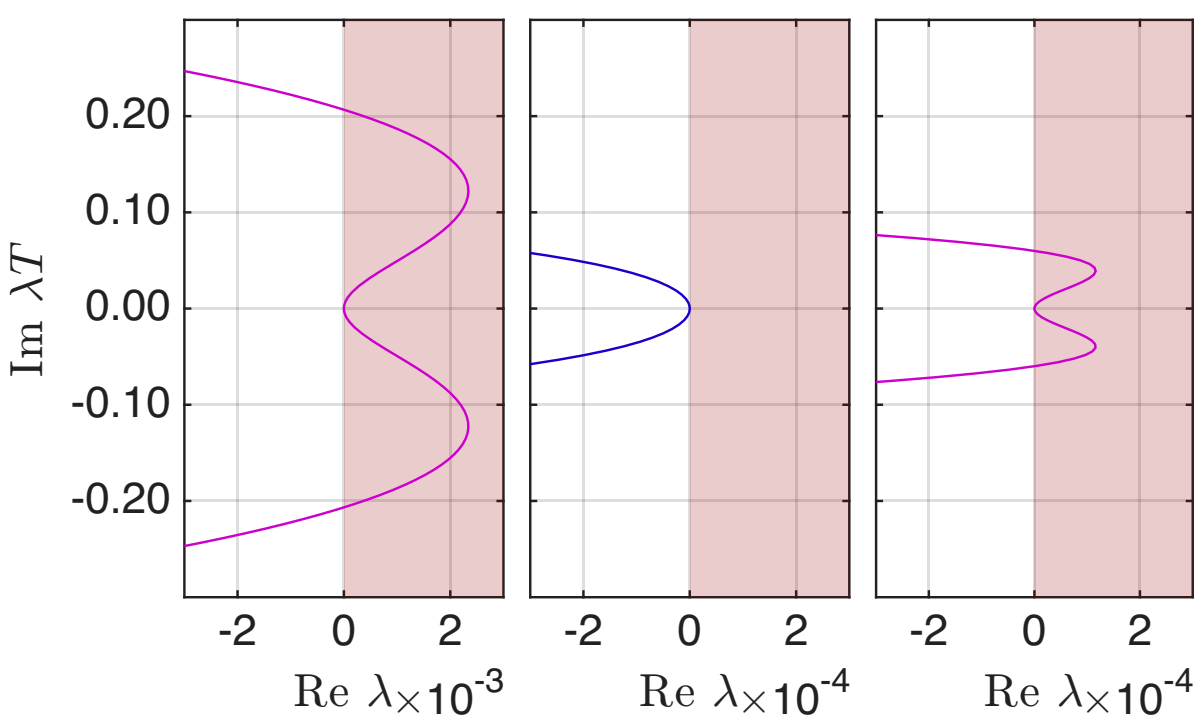

Figure 2: Master stability curves for $\alpha = 0.8$ and three different values of the wavenumber.

In summary, we have extended the master stability framework [3] to investigate the stability of periodic traveling waves in discrete models of materials with odd elasticity. The existence of finite-amplitude periodic traveling waves in these materials, despite the absence of external forcing, is due to a balance between nonreciprocal springs (an active element) and energy loss (damping).

# Viscoelastic gravity wave force on a circular cylinder in open channel

**Boyuan Yu**[1,*], **Vincent H. Chu**[1]

[1]Department of Civil Engineering, McGill University, Montreal, Canada

*Email: boyuan.yu@mail.mcgill.ca

## Abstract

Free-surface viscoelastic flow past obstacles is numerically examined using three-dimensional direct numerical simulation of the full Navier-Stokes equations. The current work includes the previously neglected surface-wave effect and investigates its influence on flow patterns and impact force. The flow around the obstacle evolves from a smooth standing wave to a breaking bow shock wave as the Froude number increases. For a fixed Froude number, the drag coefficient increases with the Weissenberg number. The drag enhancement due to viscoelasticity declines with Froude number. However, as the Weissenberg number becomes sufficiently large, the drag coefficient approaches an asymptote that depends only on the Froude number.

## 1 Introduction and Motivation

The impact force of viscoelastic flows on an obstacle is a problem of fundamental interest, which dates back to the 1940s [1, 2]. However, despite the ubiquity of free-surface viscoelastic (FSVE) flows, the role of gravity waves in their impact on obstacles has received much less attention.

This study focuses on FSVE flows past a circular cylinder in an open channel on a mild slope, highlighting the effect of surface waves on flow pattern and impact force. Fig. 1 defines the problem and key parameters to be determined by numerical simulation.

## 2 Formulation and Numerical Methods

Three-dimensional direct numerical simulation with volume-of-fluid free-surface capturing technique combined with Oldroyd-B rheology [3] is conducted for the studied problem. The dimensionless governing equations using the scales of steady-uniform solution of FSVE flows down inclined slopes ($\overline{U}$ and $H$ denote the depth-averaged velocity and fluid depth for steady-uniform solution of FSVE, respectively) [4] are:

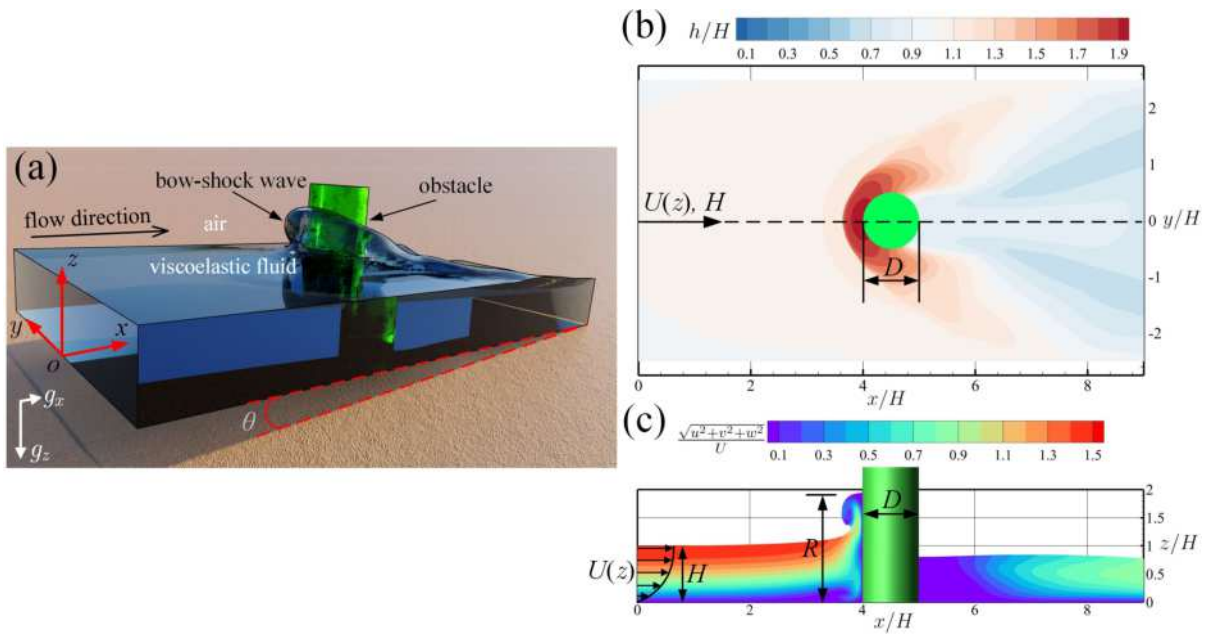

Figure 1: Schematics of free-surface viscoelastic flow past a circular cylinder on an inclined plane. (a) The three-dimensional view of the studied problem rendered by Blender. (b) The plane view of the free-surface around the circular cylinder. The centerline is denoted as the dashed line. (c) The vertical view of velocity distribution at the centerline.

$$\nabla \cdot \tilde{\boldsymbol{u}} = 0 \tag{1}$$

$$\frac{\partial \tilde{\boldsymbol{u}}}{\partial \tilde{t}} + \tilde{\boldsymbol{u}} \cdot \nabla \tilde{\boldsymbol{u}} = -\nabla \tilde{p} + \nabla \cdot (\frac{1}{\mathrm{Re}} \tilde{\boldsymbol{\tau}}) + \frac{1}{\mathrm{Fr}^2} \boldsymbol{f} \tag{2}$$

$$\tilde{\boldsymbol{\tau}} + \mathrm{Wi}\, \overset{\triangledown}{\tilde{\boldsymbol{\tau}}} = 2(\tilde{\boldsymbol{D}} + \beta \mathrm{Wi}\, \overset{\triangledown}{\tilde{\boldsymbol{D}}}) \tag{3}$$

$\tilde{\boldsymbol{u}}$, $\tilde{\boldsymbol{\tau}}$ and $\tilde{\boldsymbol{D}}$ are velocity vector, extra-stress tensor and deformation tensor following the coordinate system shown in Fig. 1, respectively. Therefore, the six dimensionless parameters/factors defining the problem are: (1) The Froude number $\mathrm{Fr} \equiv \overline{U}/\sqrt{g'H}$. (2) The Weissenberg number $\mathrm{Wi} \equiv \lambda_1 \overline{U}/H$. (3) The solvent viscosity ratio $\beta \equiv \lambda_2/\lambda_1$. (4) The slope inclination $S_o \equiv \tan\theta$. (5) The shape of obstacles. (6) The ratio of obstacle width to steady-uniform-flow depth $r \equiv D/H$. This study concentrates on the effect of Fr and Wi. We use a circular cylinder with fixed normalized diameter $D/H = 1$ on a mild slope $S_o = 0.05$ and fixed value of $\beta = 0.5$ is used. The Reynolds number is related to Fr and $S_o$ by $\mathrm{Re} = 3\mathrm{Fr}^2/S_o$ [4].

The numerical simulation is conducted using `rheoInterFoam` solver of RheoTool toolbox [5] which includes improved numerical schemes for high Weissenberg number problems.

## 3 Combined Effect of Surface-wave and Viscoelasticity on Drag Force

Fig. 2 shows the surface-wave pattern around the circular cylinder for Fr = 0.5, 1.5 and 2.5. Surface wave around the circular cylinder is a smooth standing wave for Fr = 0.5. In comparison, for larger Froude numbers 1.5 and 2.5, the runup heights increase and surface waves become bow shock waves with steep wave fronts. Different from Newtonian fluids, the first-normal-stress difference is notably large, which is due to the effect of viscoelasticity.

The normalized drag force by inertia force $\frac{1}{2}\rho\overline{U}^2 HD$ and viscosity force $(\mu_s + \mu_p)\overline{U}^2 H$ are denoted as $C_D$ and $C_\mu$, respectively. Fig. 3 shows $C_D$, $C_\mu$ and $R/H$ as a function of Wi for Fr in the range of $[0.5, 2.5]$. Drag enhancement due to viscoelasticity is observed for all the Froude numbers in the range and runup height slightly decreases with Wi. However, $C_D$ and $C_\mu$ reaches asymptotes as Wi is sufficiently large and the asymptotic values only depends on Fr.

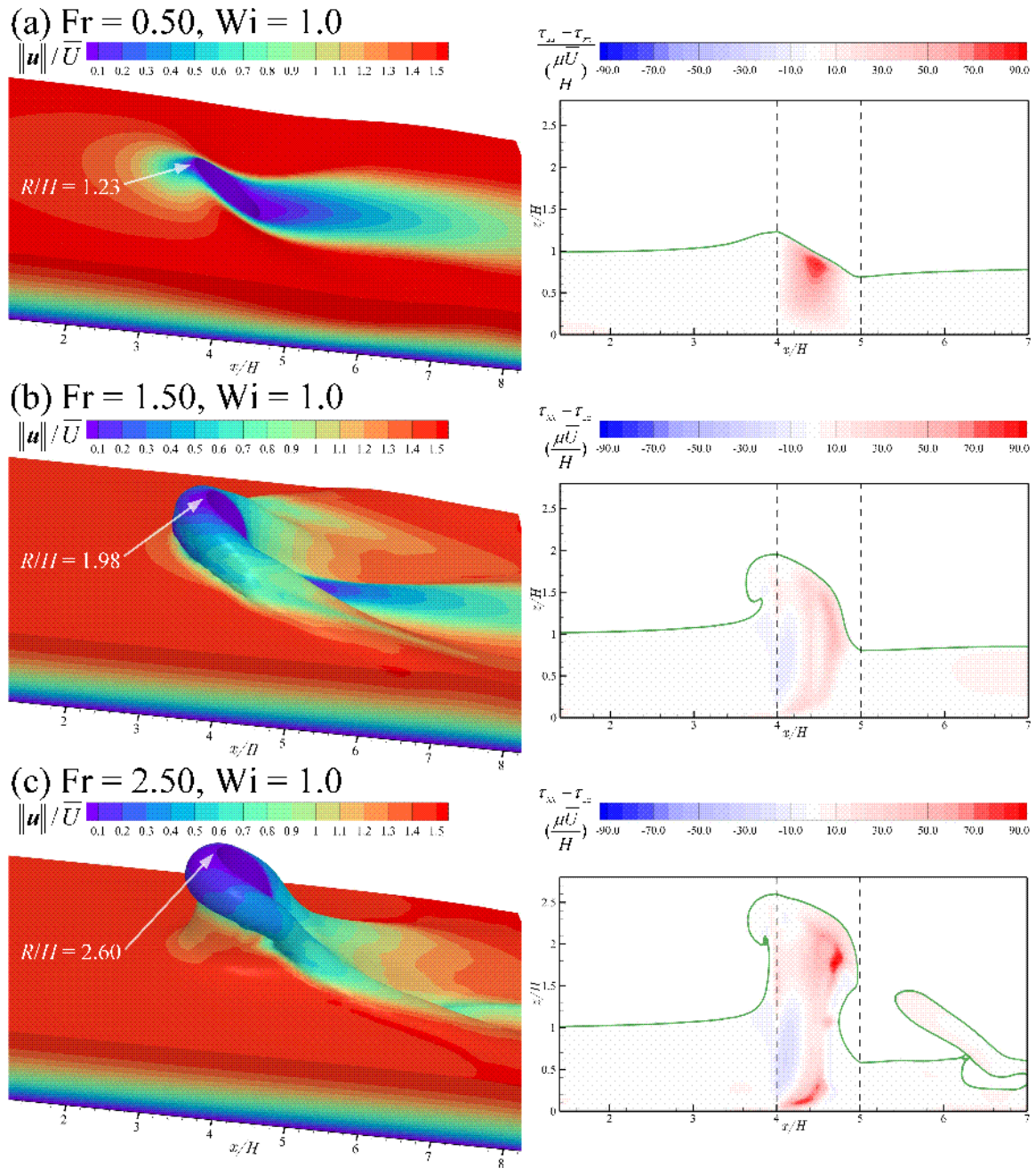

Figure 2: The influence of free-surface on flow pattern. Left column: velocity contours and 3D view of free-surfaces; right column: contours of first-normal-stress difference at the centerline slice. The dashed lines mark the extent of the circular cylinder. The green curves denote the fluid-air interface. (a) Fr = 0.50, Wi = 1.0 (Re = 15.00). (b) Fr = 1.50, Wi = 1.0 (Re = 135.0). (c) Fr = 2.50, Wi = 1.0 (Re = 375.0).

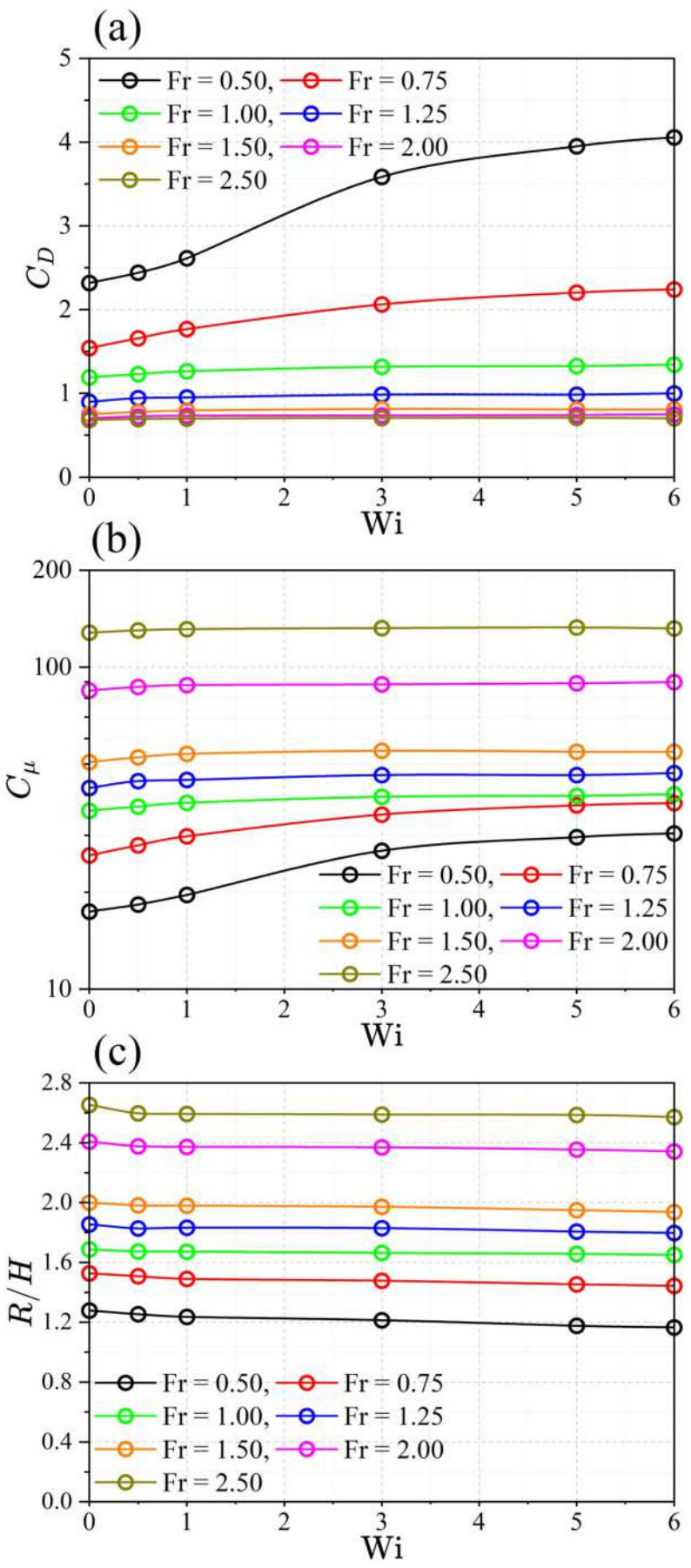

Figure 3: The viscoelastic effect on drag coefficient $C_D$, $C_\mu$ and runup height $R/H$ for fixed slope $S_o = 0.05$ and varying Froude number. (a) $C_D$ as a function of Weissenberg number Wi. (b) $C_\mu$ as a function of Wi. (c) $R/H$ as a function of Wi.

# Finite Element Formulation with Integral Dielectric Kernel for ICRH Plasma Heating on Non-Field-Aligned Meshes

**Jules Zaleski**[1,*], **Bernard Reman**[2], **Philippe Lamalle**[2], **Christophe Geuzaine**[1]
[1]Department of Electrical Engineering and Computer Science, University of Liège, Belgium
[2]Laboratory for Plasma Physics, Royal Military Academy, Brussels, Belgium
*Email: j.zaleski@uliege.be

### Abstract

This work presents a finite element formulation to model Ion Cyclotron Resonance heating (ICRH) of fusion plasmas in the frequency domain. In the limit of vanishing Larmor radius, wave-particle interactions are modelled by a 1D integral dielectric kernel, integrated over static magnetic field lines. The approach is validated on general (non-field-aligned) 2D unstructured meshes.

## 1 Introduction

To simulate the impact of warm plasma on ICRH, the main approach for decades has been to use spectral solvers with a finite number of excitation modes. Recent reformulations in configuration space based on an integral representation of the RF interactions for Maxwellian plasma population [1] open the door to discretization using the finite element method (FEM). This would greatly facilitate the modelling of more complex geometries using locally adapted meshes.

Initial numerical investigation of finite element discretization of the integral dielectric kernel have focused on the particular case of field-aligned meshes, where the elements are locally aligned with the static magnetic field lines. A fully analytical integration of the kernel on such meshes is proposed in [2,3]. In this contribution, we extend this methodology to general unstructured meshes, where we integrate the dielectric kernel semi-analytically.

## 2 Integral dielectric kernel formulation

In the frequency domain, the FEM formulation of the Maxwell-Vlasov equations consists in finding the electric field $\boldsymbol{E}$ such that

$$\frac{i}{2}\int_{\mathcal{V}}\left[\frac{1}{\omega\mu_0}\left(\nabla\times\mathbf{F}\right)^*\cdot\left(\nabla\times\mathbf{E}\right)-\omega\epsilon_0\mathbf{F}^*\cdot\mathbf{E}\right]d\mathcal{V} + \sum_{\beta}\mathcal{W}_\beta = -\frac{1}{2}\int_{\mathcal{V}}\mathbf{F}^*\cdot\mathbf{j}_s dS, \quad (1)$$

holds for suitable test functions $\boldsymbol{F}$, at a given angular frequency $\omega$ and for a given source current density $\boldsymbol{j}_s$. The plasma contribution is taken into account for each particle species $\beta$ and wave polarization $E_L$ (such that $E_{-1} := E_-$, $E_0 := E_\parallel$, $E_{+1} := E_+$) through the nonlocal integral dielectric kernel $\mathcal{W}_\beta$ (see [1]). For a 2D problem:

$$\mathcal{W}_\beta = -\frac{i\epsilon_0}{2}\int dr\int ds\sum_{L=-1}^{+1}2^{\alpha/2}\frac{\omega_{p_\beta}^2}{v_{\mathrm{T}_\beta}} \times\int ds' F_L^*(s)\,\Upsilon_\alpha\left(\frac{\omega-L\omega_{c_\beta}}{v_{T_\beta}}\left|s-s'\right|\right)E_L(s'), \quad (2)$$

with $\alpha = 2\delta_{L,0}$, $s$ and $s'$ curvilinear coordinates along the magnetic field line, $\omega_p$ the plasma frequency, $v_T$ the thermal velocity, $\omega_c$ the resonant frequency, and $\Upsilon_\alpha$ the so-called kernel dispersion function. The latter exhibits a logarithmic singularity at the origin for $\alpha = 0$ and for all $\alpha$, the function tends asymptotically to 0 on both sides of the domain.

To integrate the local terms in (1), a standard Gauss quadrature is used. The integration of the non-local terms in (2) requires computing a piece-wise linear path in the unstructured mesh along the static field line: see Figure 1.

For each quadrature point $g$ in element $T_i$, both the self-interaction and the interaction with element $T_j$ are computed analytically over the corresponding line segments in the path, exploiting the recursive properties of the kernel and the values of the finite element basis functions along the segment. This allows an efficient and accurate handling of the logarithmic singularity. The plasma parameters $\frac{\omega-L\omega_{c_\beta}}{v_{T_\beta}}$ govern

[1]Jules Zaleski is a FRIA grantee of the Fonds de la Recherche Scientifique – FNRS.

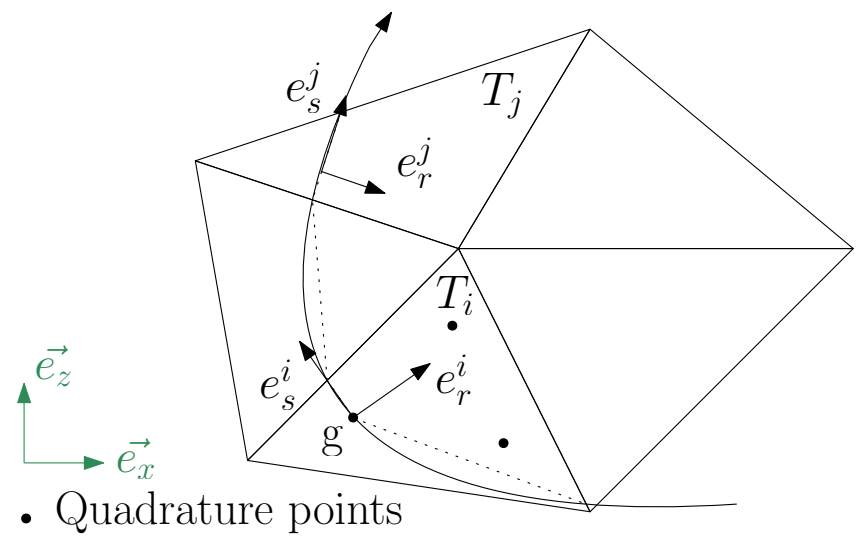

Figure 1: Non-local coupling between one quadrature point in element $T_i$ and distant element $T_j$ along the static field line.

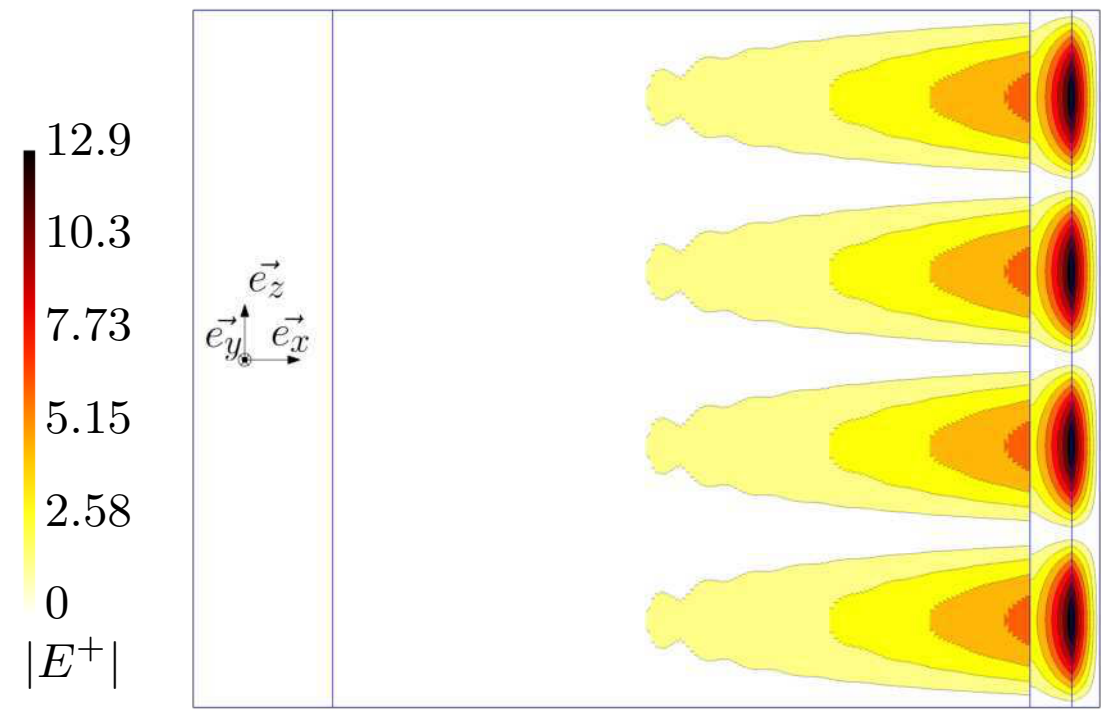

Figure 2: Norm of the solution $E_+ := (E_x + iE_y)/\sqrt{2}$ at $f = 51.8$ MHz with an in-plane static magnetic flux density with constant amplitude $B^y_{\text{stat}} = 3.45$ T, with three particle species: electrons (0 keV, $n_{\text{elec}} = 3 \cdot 10^{19}$ m$^3$), hydrogen ($T_H = 2$ keV, 5 % $\cdot\, n_{\text{elec}}$) and deuterium (0 keV, 95 % $\cdot\, n_{\text{elec}}$).

the rate of decay of $\Upsilon_\alpha$, which allows to reduce the integration support.

## 3 Preliminary Results

Figure 2 demonstrates the new approach in a 2D setting, where the electric field is discretized using high-order Whitney edge elements for the in-plane components, and high-order Lagrange elements for the out-of-plane component[2]. An out-of-plane current density $j^z_s = \sin(k_\parallel y)$ is imposed on a line strap antenna on the right, which leads to a single antenna mode $(k_\parallel)$. In this configuration, the results can be validated against those obtained using a complex permittivity tensor obtained through a spectral approach [4]: the excellent agreement is confirmed by the slope of the absorption profile shown in Figure 3.

Figure 4 illustrates the new approach by modelling a more realistic antenna setting. The current density is applied on four radiating ICRH antenna straps, which excites a broad spectrum of spatial modes.

[2]https://gitlab.onelab.info/gmsh/fem

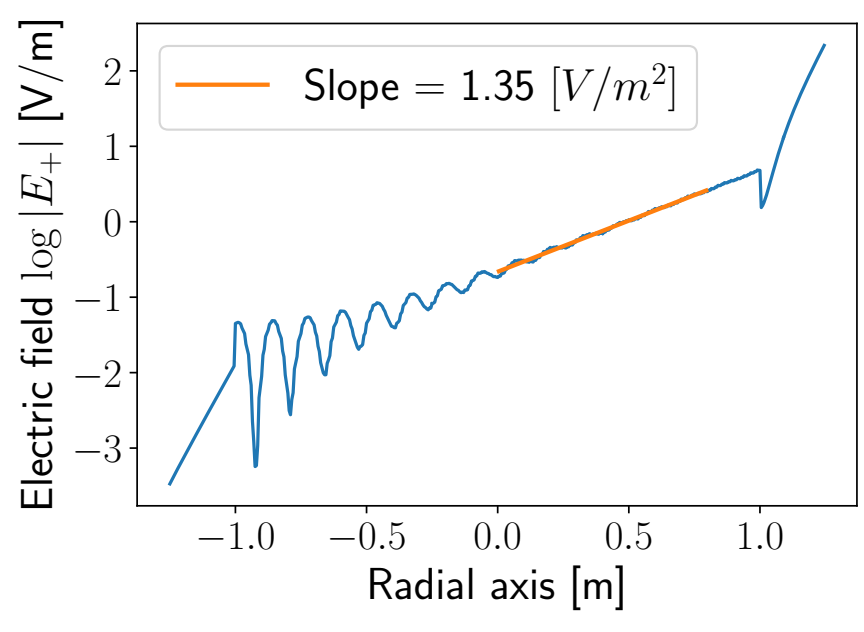

Figure 3: Horizontal absorption profile of the solution of Figure 2. The slope of the absorption profile is validated against the analytical dispersion relation.

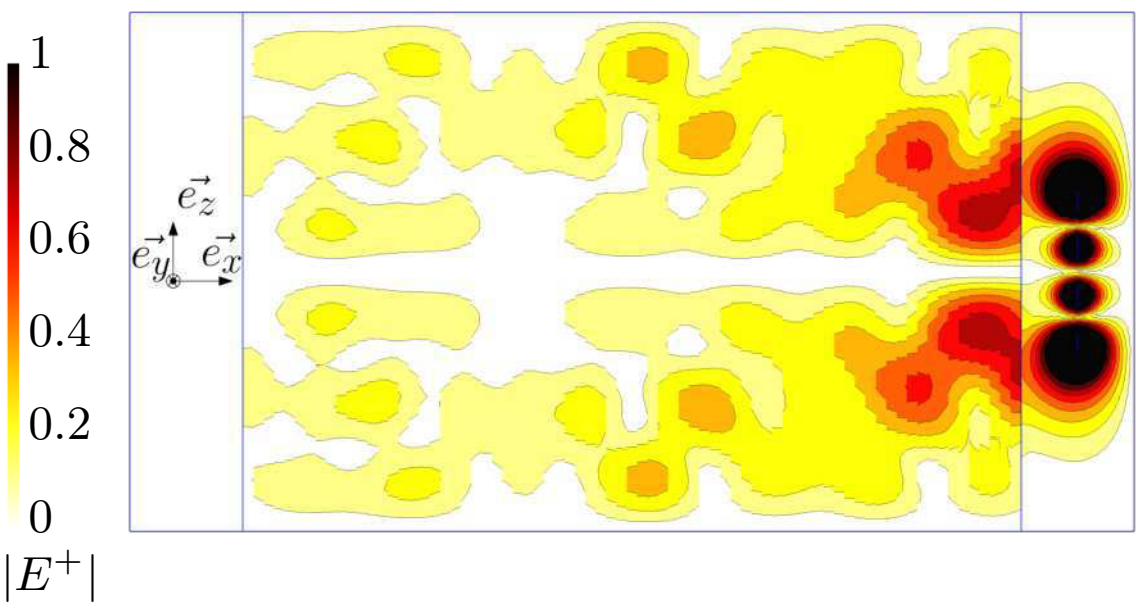

Figure 4: Norm of the solution $E_+$ at the same frequency and with same particle species as for Figure 2. The field is saturated at 1V/m near the antenna straps for visualization purposes.

# Preserving multiple conserved quantities for the Korteweg–de Vries equation

**Andy T.S. Wan**[1], **Tareq Uz Zaman**[2],

[1]Department of Applied Mathematics, University of California, Merced, USA
[2]Department of Mathematics, Simon Fraser University, Burnaby, Canada

## Abstract

Many evolutionary partial differential equations possesses multiple conservation laws and preserving such quantities is essential for accurate long-time predictions. The Korteweg–de Vries (KdV) equation is one such nonlinear dispersive model possessing an infinite hierarchy of conservation laws. We investigate a class of structure-preserving semi-discretizations in space for the KdV equation which inherit a finite number of discrete conservation laws. To preserve multiple conserved quantities across time, we employ the Minimal $\ell^2$-Norm Discrete Multiplier Method (MN-DMM) [1] for the proposed semi-discretizations. Numerical results on solitons illustrate accurate preservation of conserved quantities, as well as stable and accurate long-time dynamics over nonconservative schemes.

## 1 Introduction

The KdV equation is a cornerstone in the study of nonlinear wave phenomena, arising in shallow water dynamics, ion-acoustic plasmas, and optical fibers. A defining feature of the KdV equation is its infinite hierarchy of conserved quantities [2], which encode the integrable structure of the equation and are critical for the physical fidelity of solutions. Thus, preserving invariants at the discrete level is important for stable and accurate long-time predictions [3].

Conservative numerical schemes for the KdV equation have been developed in many settings, including Hamiltonian-preserving discontinuous Galerkin methods [4] and linearly-fitted energy-mass-preserving schemes [5]. However, preserving multiple invariants simultaneously for general classes of structure-preserving spatial discretizations remains challenging. In this work, we develop a structure-preserving FE discretiztion by combining *compatible* finite elements spatial discretization with the MN-DMM [1] to preserve multiple invariants in time, resulting in a fully conservative discrete scheme.

## 2 Model

We consider the periodic KdV equation on the domain $x \in [0, L]$:

$$u_t + u\, u_x + u_{xxx} = 0, \qquad (1)$$

with periodic boundary conditions. It is well-known the first four conserved quantities (with the first three invariants being the mass, momentum, and energy) are given by

$$M = \int_0^L u\, dx, P = \int_0^L u^2\, dx, E = \int_0^L \left(u^3 - 3u_x^2\right) dx,$$

$$H = \int_0^L \left(5u^4 - 10u\, u_x^2 + u_{xx}^2\right) dx.$$

These belong to an infinite hierarchy of conserved quantities generated by the recursion operator $\mathcal{R} = \partial_{xx} + \frac{2}{3}u + \frac{1}{3}u_x\partial_x^{-1}$ associated with the integrable structure of the KdV equation [2]. Preserving these quantities discretely is thus vital for capturing the integrable structure and long-time dynamics of soliton solutions.

## 3 Methodologies and Results

**Compatible spatial discretization.** In principle, our approach applies to a broader class of structure-preserving spatial semi-discretizations, including finite difference, spectral, and finite element methods. As a concrete example, we employ $C^k$-continuous Hermite finite elements on a uniform periodic mesh with $N$ elements. For $k = 1$, each element carries two nodes with degrees of freedom $(u_i,\, u_x|_i)$, yielding $2N$ global unknowns; for $k = 2$, the additional degree of freedom $u_{xx}|_i$ gives $3N$ unknowns and enables preservation of higher-order invariants. Specifically, the Hermite basis enforces the prequisite smoothness for the third-order dispersive term, which enables exact discrete integration-by-parts identities, mirroring the skew-symmetry of the continuous operator. This spatial semi-discretization of (1) results in the abstract ODE of the form, $\mathbf{M}\,\dot{\mathbf{U}} = \mathbf{F}(\mathbf{U})$ where $\mathbf{M}$ is the FEM mass matrix, $\mathbf{U}$ is a vector for nodal degrees of freedom, and $\mathbf{F}(\mathbf{U})$ is the discrete spatial operator corresponding to the nonlinear and dispersive terms. Thanks to the compatible discretization, the discrete conserved invariants inherit the structure of the continuous invariants; for instance, the discrete mass and momentum are $M_h = \mathbf{1}^\top \mathbf{M}\mathbf{U}$ and $P_h = \mathbf{U}^\top \mathbf{M}\mathbf{U}$.

**Conservative temporal discretization.** The above semi-discrete system is advanced in time using the MN-DMM [1]. Specifically, let $\boldsymbol{f}^\tau$ be any consistent one-step scheme to $\mathbf{M}^{-1}\mathbf{F}(\mathbf{U})$, viewed as a predictor. For a time-independent conserved vector $\boldsymbol{\Psi}$, the MN-DMM corrects $\boldsymbol{f}^\tau$ by the conservative scheme (generally implicit),

$$\boldsymbol{f}^\tau_{\mathrm{MN}} = \boldsymbol{f}^\tau - (\Lambda^\tau)^+ \Lambda^\tau \boldsymbol{f}^\tau,$$

with the fully discrete system to be solved,

$$\mathbf{U}^{n+1} = \mathbf{U}^n + \tau\, \boldsymbol{f}^\tau_{\mathrm{MN}}(t^n, \mathbf{U}^{n+1}, \mathbf{U}^n).$$

Here $(\Lambda^\tau)^+$ is the right Moore–Penrose pseudoinverse of the discrete multiplier matrix $\Lambda^\tau \in \mathbb{R}^{m\times n}$, whose $j$-th row is $\frac{\Delta}{\Delta \boldsymbol{U}}\, \Psi_j\big(\mathbf{U}^n, \mathbf{U}^{n+1}\big)$, the divided difference of the $j$-th component of $\boldsymbol{\Psi}$.

The main idea of MN-DMM is to construct a unique conservative scheme $\boldsymbol{f}^\tau_{\mathrm{MN}}$ which is closest to $\boldsymbol{f}^\tau$ in $\ell^2$ norm; for time-independent invariants this is equivalent to projecting $\boldsymbol{f}^\tau$ onto the kernel of $\Lambda^\tau$, which automatically satisfies the discrete multiplier condition $\Lambda^\tau \boldsymbol{f}^\tau_{\mathrm{MN}} = \mathbf{0}$ [1]. The minimum $\ell^2$-norm correction is computed via the Moore–Penrose pseudoinverse and applied through Picard fixed-point iterations until the invariants are restored to machine precision.

**Preliminary numerical results.** We initialize a soliton $u(x, 0) = A\,\mathrm{sech}^2\big(\sqrt{A/12}\,(x - x_0)\big)$ with $A = 1$, $x_0 = L/3$, $L = 40$, using $C^1$ Hermite elements with $N = 100$ and $\Delta t = 0.001$ over $t \in [0, 30]$. Table 1 reports the absolute errors in the four conserved quantities at $t = 30$, comparing the Euler method versus the MN-DMM scheme using Picard iterations. MN-DMM achieves near machine-precision conservation with all four invariants, improving the fourth-invariant error by nine orders of magnitude. Figure 1 further illustrates the evolution of conservation errors over time, confirming that MN-DMM effectively suppresses the numerical dissipation and dispersion that otherwise accumulate in long-time simulations.

Table 1: Absolute errors in conserved quantities at $t = 30$ for the four KdV invariants using $C^1$ Hermite FEM ($N = 100$, $\Delta t = 0.001$).

| Invariant | Euler | MN-DMM (Picard) |
|---|---|---|
| $\lvert M(T) - M(0)\rvert$ | $3.55 \times 10^{-15}$ | $1.78 \times 10^{-15}$ |
| $\lvert P(T) - P(0)\rvert$ | $6.84 \times 10^{-12}$ | $5.33 \times 10^{-15}$ |
| $\lvert E(T) - E(0)\rvert$ | $2.40 \times 10^{-6}$ | $3.02 \times 10^{-14}$ |
| $\lvert H(T) - H(0)\rvert$ | $3.27 \times 10^{-3}$ | $1.02 \times 10^{-12}$ |

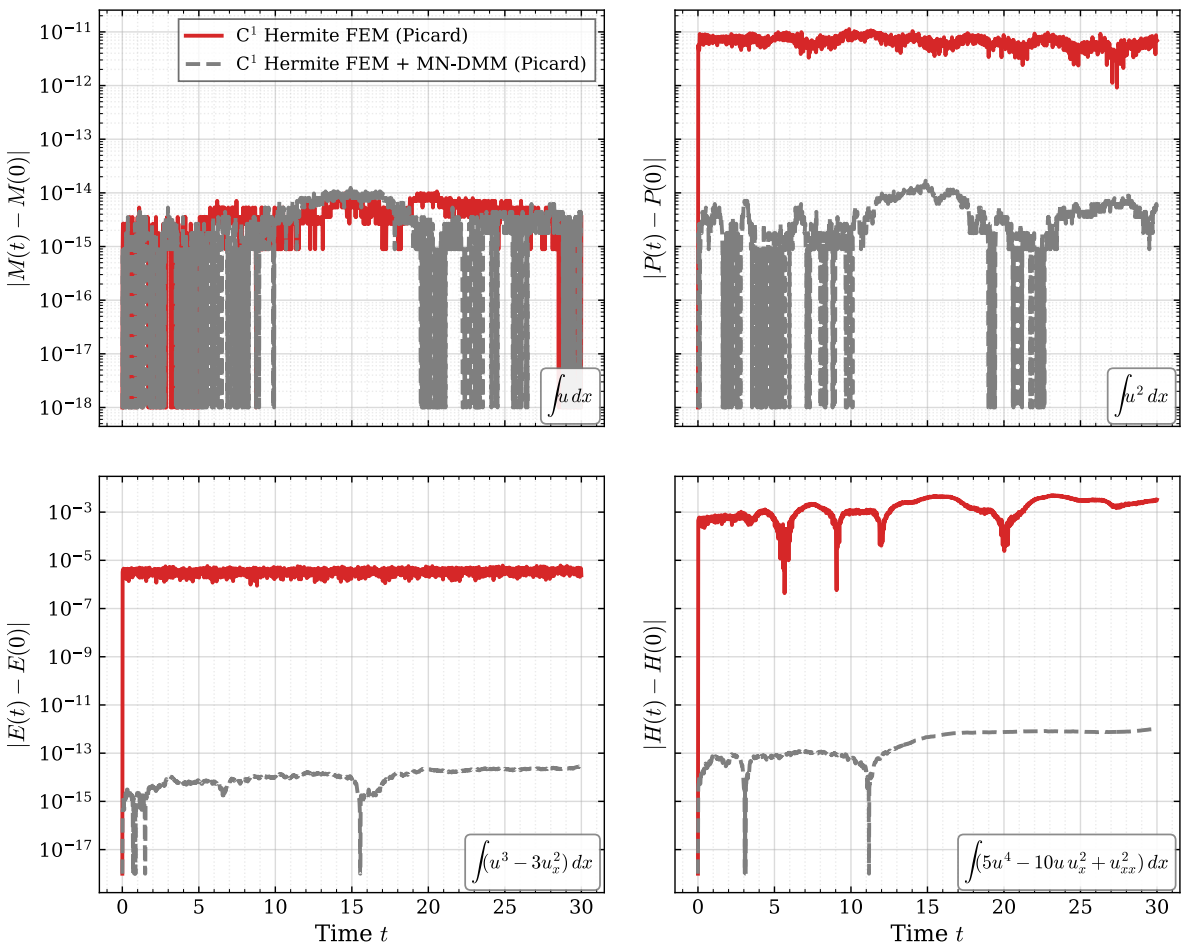

Figure 1: Errors in conserved quantities for the four KdV invariants using $C^1$ Hermite FEM ($N = 100$, $\Delta t = 0.001$).

# An inverse problem for estimating effective parameters in asymptotic elastodynamics with non-periodic roughness

**Rodrigo Zelada**[1,2,*], **Alexandre Imperiale**[1], **Philippe Moireau**[2]
[1]Université Paris-Saclay, CEA, List, F-91120, Palaiseau, France
[2]CMAP, Inria, CNRS, École polytechnique, Institut Polytechnique de Paris
*Email: rodrigo.zeladamancini@cea.fr

**Abstract**

We consider the problem of elastic wave propagation in a two-dimensional domain where part of the boundary is defined by a non-periodic rough perturbation. We assume that the maximum amplitude of this perturbation is small compared to the minimum wavelength. Standard discretization procedures in this context can be computationally intensive. Therefore, we seek an approximate model whose numerical complexity is robust with respect to this roughness. In this context, asymptotic expansions are well-known tools for providing effective models. Although a priori estimates exist for certain cases, most asymptotic models suffer from a modeling error proportional to the amplitude of the perturbation. In our work, we aim to reduce this modeling error by introducing a modification of the effective coefficients in the boundary condition. This modification is automatically adjusted by solving an inverse problem formulated as a minimization problem. In this minimization procedure, the data fidelity term is informed by synthetic data generated from a single run of the complete model. We illustrate our approach with two-dimensional numerical simulations.

## 1 Presentation of the problem

Let $\eta$ be a function taking values between 0 and $\bar{\eta} > 0$, where $\bar{\eta}$ is a small parameter. Let $\Omega^\eta$ be an open bounded connected subset of $\mathbb{R}^d$. We suppose that $\Omega^\eta$ has a Lipschitz boundary, composed by four parts : $\partial\Omega^\eta = \Gamma_N \cup \Gamma_L \cup \Gamma_R \cup \Gamma_\eta$, where $\Gamma_N$ is the free surface and $\Gamma_\eta$ is a non-periodic rough boundary and parametrized by $\eta$ representing a variable thickness. We consider $H^1_{LR}(\Omega^\eta) := \{v \in H^1(\Omega^\eta) \mid v \text{ periodic on } \Gamma_L, \Gamma_R\}$. The displacement $\boldsymbol{u}^\eta$ in an homogeneous isotropic material satisfies the elastic wave equation

$$\frac{d^2}{dt^2}\int_{\Omega^\eta} \rho \boldsymbol{u}^\eta \cdot \boldsymbol{v} \, \mathrm{d}x + \int_{\Omega^\eta} \sigma(\boldsymbol{u}^\eta) : \varepsilon(\boldsymbol{v}) \, \mathrm{d}x = \ell(\boldsymbol{v}), \quad (1)$$

where $\rho > 0$ is the mass density, and $\ell(\boldsymbol{v}) := \int_{\Gamma_N} \boldsymbol{g} \cdot \boldsymbol{v} \, \mathrm{d}s$, $\forall \boldsymbol{v} \in H^1_{LR}(\Omega^\eta)^d$, with $\boldsymbol{g}$ a known regular (in time and space) load on $\Gamma_N$ and initial displacement and velocity equal to zero. The problem (1) is very expensive to solve (since we need to mesh very thin around $\Gamma_\eta$) and even more if we plan to use time explicit schemes. Instead, we split the domain $\Omega^\eta$ into two open bounded disjoints domains $\Omega$ and $\Omega^\eta_-$, where only $\Omega^\eta_-$ contains the rough boundary $\Gamma_\eta$ and $\Omega$ is independent of $\eta$. We denote $\Gamma = \partial\Omega \cap \partial\Omega^\eta_-$. We solve the elastic wave equation in $\Omega$ by including an effective boundary condition on $\Gamma$, that represents the propagation within the rough domain $\Omega^\eta_-$.

## 2 Effective model

We suppose that $\eta$ is small with respect to the minimum wavelength. There are several ways to perform the asymptoptic development and rigorously justify the convergence, e.g. in the case of a constant thickness [1], smooth variable or even non smooth but periodic. In order to justify the approximate model, we consider $\eta \in \mathcal{C}^1(\Gamma)$.

For the sake of simplicity, let us consider $d = 2$, $\Omega^\eta$ to be divided into two open bounded connected subdomains $\Omega = (0, a) \times (0, b)$ and $\Omega^\eta_- = \{(x_1, x_2) \in (0, a) \times (-\bar{\eta}, 0) / -\eta(x_1) < x_2 < 0\}$, which are separated by an interface $\Gamma = \partial\Omega \cap \partial\Omega^\eta_- = \{0\} \times (0, a)$, with $a, b > 0$. Then, we obtain the following *effective boundary conditions* (EBC)

$$\sigma(\boldsymbol{u})\boldsymbol{n} = -\rho\eta(x_1)\partial^2_{tt}\boldsymbol{u} + c\partial_1 \begin{pmatrix} \eta(x_1)\partial_1 u_1 \\ 0 \end{pmatrix} \text{ on } \Gamma,$$

where $c$ is a positive constant depending on the Lamé coefficients. We define the space $U := \{\boldsymbol{v} \in H^1_{LR}(\Omega)^2 \mid v_1|_\Gamma \in H^1(\Gamma)\}$. The variational formulation of the approximate problem consists in finding $t \mapsto \boldsymbol{u}(t) \in U$, s.t. $\forall \boldsymbol{v} \in U$,

$$\frac{d^2}{dt^2}\left(\int_\Omega \rho \boldsymbol{u}\cdot\boldsymbol{v}\mathrm{d}x + \int_\Gamma \rho\eta\boldsymbol{u}\cdot\boldsymbol{v}\mathrm{d}s\right) + \int_\Omega \sigma(\boldsymbol{u}) : \varepsilon(\boldsymbol{v})\mathrm{d}x + \int_\Gamma c\eta\partial_1 u_1 \partial_1 v_1 \mathrm{d}s = \ell(\boldsymbol{v}). \tag{2}$$

Defining $V := \{(\boldsymbol{u}', \boldsymbol{\gamma}') \in H \mid \boldsymbol{u}' = (u_1', u_2') \in H^1_{LR}(\Omega)^2, \boldsymbol{\gamma}' = (\gamma_1', \gamma_2') \in H^1(\Gamma)\times L^2(\Gamma),\ \boldsymbol{u}'|_\Gamma = \boldsymbol{\gamma}'\}$, dense in $H := L^2(\Omega)^2 \times L^2(\Gamma)^2$, problem (2) admits a unique solution $(\boldsymbol{u}, \boldsymbol{\gamma}) \in \mathcal{C}^0(0,T;V) \cap \mathcal{C}^1(0,T;H)$. Furthermore, it can be shown that this approximation is of order 3 w.r.t $\eta$.

## 3 Inverse problem

We consider the non-smooth and non-periodic case where $\eta \in L^\infty(\Gamma)$. In this case, the asymptotic development is more delicate since singularities appear. In particular, the convergence rate of the smooth case is lost. In [2], a modification of the standard EBC is proposed, in order to improve the rate of convergence. We add a perturbation to the EBC, which we reconstruct by solving an inverse problem, expecting to decrease the error. Let $\alpha, \beta \in H^1(\Gamma)$. We consider $\phi$ a $C^1$ positive increasing function in $\mathbb{R}$. We denote $\theta := (\alpha, \beta) \in H^1(\Gamma)^2$. We consider the problem: find $t \mapsto \boldsymbol{u}_\theta(t) \in U$, $\forall \boldsymbol{v} \in U$,

$$\frac{d^2}{dt^2}\Big(\int_\Omega \rho \boldsymbol{u}_\theta \cdot \boldsymbol{v}\ \mathrm{d}x + \int_\Gamma \phi(\alpha)\boldsymbol{u}_\theta\cdot\boldsymbol{v}\ \mathrm{d}s\Big) + \int_\Omega \sigma(\boldsymbol{u}_\theta) : \varepsilon(\boldsymbol{v})\mathrm{d}x + \int_\Gamma \phi(\beta)\partial_1 u_{\theta,1}\partial_1 v_1 \mathrm{d}s = \ell(\boldsymbol{v}).$$

The function $\phi$ ensures by design the coercivity of the associated bilinear form. Indeed, we prefer to reparameterize the optimization problem rather than add a constraint of positivness to the parameter space.

We aim at reconstructing the parameter $\theta$. To this end, we introduce the observation region $\omega \subset \Omega$, whose distance to $\Gamma$ is at least one wavelength. We compute the solution $\mathbf{u}^\eta$ of the reference problem (1), we define the observations as $\boldsymbol{y} := \boldsymbol{u}^\eta|_\omega$ and consider the problem:

$$\inf_{\theta \in H^1(\Gamma)^2} \Big\{ \|\boldsymbol{y} - \boldsymbol{u}_\theta|_\omega\|^2_{L^2(0,T;\omega)^2} + \gamma\|\theta\|^2_{H^1(\Gamma)^2} \Big\},$$

where $\gamma > 0$ is the regularization parameter. To solve this problem, we rely on a BFGS based gradient descent method based on solving iteratively a succesion of direct and adjoint equations.

## 4 Numerical results and perspectives

We use an explicit leapfrog scheme with spectral finite elements. We start from the parameters given by the asymptotic model (2) and we iterate until convergence. We are able to reconstruct effective parameters whose associated solution better approximates the reference solution. Qualitatively, we observe improvements in the reflected wave near the top right corner, as depicted in Figure 1b, where we plot the displacement magnitudes over a horizontal line at fixed $x_2$ (across the reflected wave): the reference in red, the initial in black, and the optimized in blue.

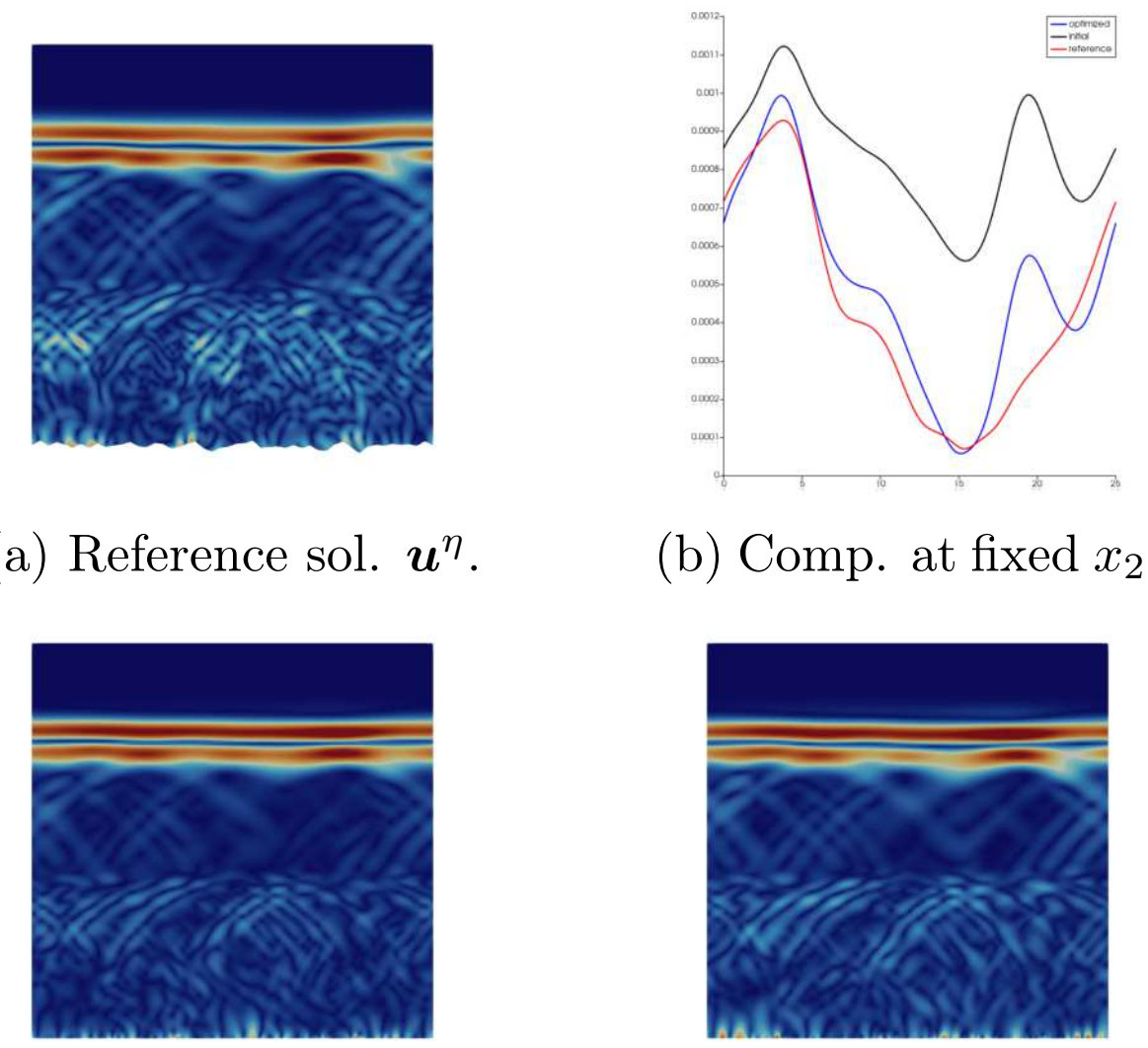

(a) Reference sol. $\boldsymbol{u}^\eta$. (b) Comp. at fixed $x_2$.

(c) Initial guess $\boldsymbol{u}_{\theta_0}$. (d) Optimal sol. $\boldsymbol{u}_{\bar{\theta}}$.

Figure 1: Snapshots after reflection (fixed $t$).

Next, we plan to consider a set of incident waves with continuously varying angles. After sampling the incident angles, solving the reference equation for each of them, and reconstructing the parameters, we hope to build an approximate model precise enough for any incident wave, in particular those that were not used to generate the synthetic data.

# Exponential and harmonic regimes of anti-plane shear wave in a piezoelectric thin layer over half-space with surface elasticity

**Anna Zemlyanova[1,*], Gennadi Mikhasev[2], Arindam Nath[2], Jia Fei[2]**

[1]Department of Mathematics, Kansas State University, Manhattan, KS, USA

[2]International Center for Applied Mechanics, Department of Astronautical Science and Mechanics, Harbin Institute of Technology, Harbin, China

*Email: azem@ksu.edu

## Abstract

This paper delineates an analytical framework for examining the transference of anti-plane shear waves in a piezoelectric thin film that is imperfectly bonded to a half-space medium, based within the context of surface-interface elasticity theory. Piezoelectric materials with layered design are extensively utilized in smart devices owing to their intrinsic abilities in sensing, actuation, and control. Two distinct wave propagation regimes can be identified: one exhibiting exponentially decaying amplitudes (TE) and the other displaying harmonically varying amplitudes (TH). Using the Gurtin–Murdoch surface elasticity model, we find the dispersion relation for waves propagating through the piezoelectric layered structure. The study also looks at interfacial imperfections, for which the displacement and electric potential are not continuous because the two media are imperfectly bonded. The derived dispersion relations are compared with established results in the literature and the numerical simulations are performed.

## 1 Introduction

Piezoelectric materials have emerged as innovative smart materials with extensive applicability across several technologies. A key characteristic of piezoelectricity is the production of mechanical strain in response to external voltage. Piezoelectric materials are extensively utilized in smart buildings owing to their intrinsic abilities in sensing, actuation, and control. Surface-interface effects are regarded as the principal cause of the extraordinary behaviors shown by nanostructures, attributable to high surface-to-volume ratio. Gurtin and Murdoch introduced foundational ideas concerning surface-free energy and surface stress in elastic media. In our study, we were the first to examine the simultaneous jump in both mechanical displacements and electric field potential at the interface between different piezoelectric materials and study their influence on the propagation of anti-plane shear waves in two different anti-plane wave regimes, namely, transverse exponential (TE) and transverse harmonic (TH) waves, in a piezoelectric material, taking into account the influence of surface-interface elasticity theory.

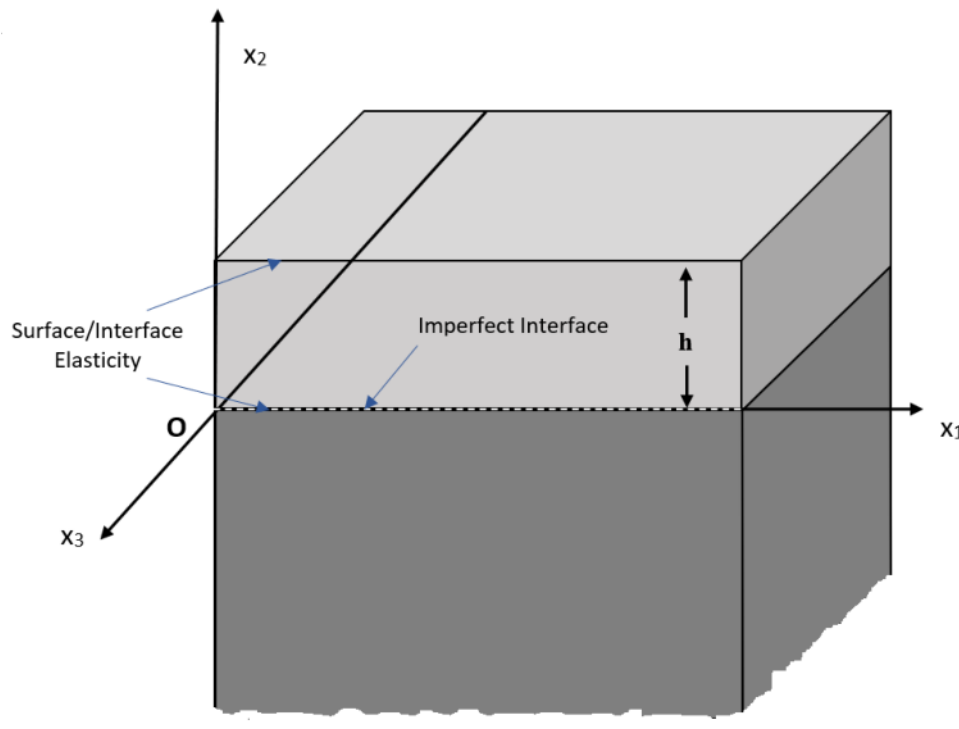

Figure 1: Thin infinite piezoelectric plate-like domain lying on half-space media

## 2 Formulation of the problem

Consider a three-dimensional piezoelectric thin plate of thickness $h$ occupying a three-dimensional domain $D = (-\infty < x_1 < \infty, 0 \leqslant x_2 \leqslant h, -\infty < x_3 < \infty)$ as shown in Fig. 1. For the propagation of anti-plane shear wave in the piezoelectric layer and semi-space, the components of displacement $\mathbf{u}^{(j)} = (u_1^{(j)}, u_2^{(j)}, u_3^{(j)})$ and electric potential $\phi^{(j)}$ are taken in the form

$$u_3^{(j)} = u_3^{(j)}(x_1, x_2, t),\; u_1^{(j)} = u_2^{(j)} = 0,\\ \phi^{(j)} = \phi(x_1, x_2, t), \tag{1}$$

where $j = 1$ corresponds to the layer, and $j = 2$ corresponds to the semi-space.

The shear stress and electric displacement in the piezoelectric medium are

$$\begin{aligned} \tau_{23}^{(j)} &= c_{44}^{(j)} u_{3,2}^{(j)} + e_{15}^{(j)} \phi_{,2}^{(j)}, \\ D_2^{(j)} &= e_{15}^{(j)} u_{3,2}^{(j)} - \varepsilon_{11}^{(j)} \phi_{,2}^{(j)}. \end{aligned} \tag{2}$$

The governing equations can be written in terms of $u_3^{(j)}$ and $\psi^{(j)} = \phi^{(j)} - \frac{e_{15}^{(j)}}{\varepsilon_{11}^{(j)}} u_3^{(j)}$ as

$$\bar{c}_{44}^{(j)} \nabla^2 u_3^{(j)} = \rho^{(j)} \ddot{u}_3^{(j)}, \tag{3}$$

$$\nabla^2 \psi^{(j)} = 0, \tag{4}$$

where $\bar{c}_{44}^{(j)} = c_{44}^{(j)} + \frac{\left(e_{15}^{(j)}\right)^2}{\varepsilon_{11}^{(j)}}$ is the stiffened elastic constant due to the piezoelectric effect. Boundary conditions:

1. The mechanical boundary condition pertinent to the stress-free boundary of the layer with surface energy, i.e., at $x_2 = h$:

$$\tau_{31,1}^{(s_1)} - \tau_{23}^{(1)} - \rho^{(s_1)} \ddot{u}_3^{(1)} = 0 \tag{5}$$

2. The electric boundary conditions with surface energy, i.e., at $x_2 = h$:

$$D_{1,1}^{(s_1)} - D_2^{(1)} = 0 \tag{6}$$

3. The mechanical boundary condition at the interface of piezoelectric/ dielectric medium, considering interface energy, i.e. at $x_2 = 0$:

$$\tau_{31,1}^{(s_2)} + \tau_{23}^{(1)} - \rho^{(s_2)} \ddot{u}_3^{(1)} = \tau_{23}^{(2)} \tag{7}$$

4. The electric boundary condition at the interface of piezoelectric/ dielectric medium considering interface energy, i.e., at $x_2 = 0$:

$$D_{1,1}^{(s_1)} + D_2^{(1)} = D_2^{(2)} \tag{8}$$

5. At the interface $x_2 = 0$, we assume a sliding represented by the bond stiffness $K_1$ in the $x_3$-direction:

$$\tau_{23}^{(2)} = K_1 \left( u_3^{(1)} - u_3^{(2)} \right). \tag{9}$$

6. A dielectrically weakly conducting interface represented $K_2$, therefore i.e., at $x_2 = 0$:

$$D_2^{(2)} = K_2 \left( \phi^{(1)} - \phi^{(2)} \right). \tag{10}$$

## 3 Conclusions

The propagation of anti-plane shear wave in a thin piezoelectric layer weakly bonded to a half-space structure is examined taking into account surface-interface elasticity. The dispersion equations were derived analytically. The phase velocity curves for TE-TE regime and TH-TE regime are analyzed to study the effects of surface-interface parameters and interfacial imperfections. The following are the main conclusions:

1. The dimensionless dispersion equations for the transference of shear wave in piezoelectric thin layer considering size-dependent surface–interface effects have been derived analytically. The findings have been reduced to specific cases in order to compare with previous findings.

2. The existence of surface/interface energy introduces further typical length-scale factors into the model. Consideration for the latter facilitated the identification of two unique regimes of anti-plane surface wave: TE-TE regime and TH-TE regime.

3. The TE-TE regime contains two modes due to the consideration of both surface and interface elasticity parameters at the two boundaries of the layer.

4. The analysis of the dispersion curves showed that the interfacial imperfect bond stiffness significantly affected the phase velocities of the propagating wave.

5. The weak conduction at the interface only affects the wave dispersion when the permittivities of the layers are quite comparable to each other.

6. The effect of imperfect interface parameters is more pronounced for the second mode of dispersion in the TE-TE regime when taking into account surface-interface elasticity.

7. The phase velocity for TE-TE regime in the piezoelectric plate attached to piezoelectric substrate is more influenced by the jump in the electric potential at the interface surface than the phase velocity for TH-TE regime in the same structure.

# Infinity as a finite boundary for Helmholtz equations

Anıl Zenginoğlu[1,*], Markus Wess[2]
[1]Institute for Physical Science and Technology, University of Maryland, College Park, USA
[2]Institute of Analysis and Scientific Computing, TU Wien, Vienna, Austria
*Email: anil@umd.edu

## Abstract

Accurate simulation of time-harmonic waves in unbounded domains relies on a proper treatment of "infinity." Truncation, absorbing layers, and approximate radiation conditions can pollute near-field solutions and make far-field extraction expensive, especially for high frequencies, trapping geometries, and variable exterior media. We introduce a compactification method that maps the exterior region to a bounded computational domain by combining a real coordinate compactification with a rescaling that factors out the known oscillatory decay. In the transformed problem, the far-field limit becomes a finite outer boundary where the outgoing Sommerfeld condition is enforced automatically by regularity, and the radiation field is obtained as boundary data. We demonstrate efficient spectral and finite-element solutions and benchmark against perfectly matched layers (PML), showing improved accuracy at fixed cost while delivering scattering amplitudes without Green's functions or post-processing. The approach is applicable to a wide range of wave propagation problems where far-field accuracy is critical.

## 1 Introduction

We study time-harmonic scattering by an obstacle embedded in a medium. We consider an obstacle $\mathcal{O} \subset \mathbb{R}^d$ in dimension $d = 2, 3$ with boundary $\Gamma = \partial\mathcal{O}$, the exterior domain $\mathcal{D} = \mathbb{R}^d \setminus \overline{\mathcal{O}}$, and a real-valued refractive index $n : \mathbb{R}^d \to \mathbb{R}$. A scattering solution at frequency $k > 0$ is a function $U : \mathcal{D} \to \mathbb{C}$ satisfying

$$(\Delta + k^2 n(x)^2)\, U = 0 \quad \text{in } \mathcal{D}, \tag{1}$$

with a boundary condition on $\Gamma$, and a radiation condition at infinity. We consider media with

$$n(x)^2 = 1 + b(x), \quad b(x) \to 0 \quad \text{as } r = |x| \to \infty, \tag{2}$$

uniformly in the angular variable $\omega = x/|x| \in \mathbb{S}^{d-1}$. At large radius $r = |x|$, the general solution admits an incoming/outgoing decomposition. For short-range media, the decomposition can be written as

$$U(r,\omega) \sim r^{-\frac{d-1}{2}} \left( e^{-ikr}\, u^-_\infty(\omega) + e^{+ikr}\, u^+_\infty(\omega) \right). \tag{3}$$

The problem we are interested in is the direct, numerical computation of the far-field $u^+_\infty(\omega)$ for a prescribed incident field $u^-_\infty(\omega)$.

## 2 Model

A generic solution to the Helmholtz scattering problem oscillates infinitely many times on an unbounded domain, making a naive compactification of space unsuitable for numerics. To make the problem numerically tractable, we split the total field into an incident and scattered part as $U = U^- + U^+$, and take out the asymptotic oscillatory decay by introducing rescaled unknowns,

$$U^\pm(r,\omega) = \frac{e^{\pm ik\, h(r)}}{r^{\frac{d-1}{2}}} u^\pm(r,\omega), \tag{4}$$

where $h$ is called the height function and satisfies $h'(r) \to 1$ as $r \to \infty$. The rescaled fields $u^\pm$ do not oscillate toward infinity and can be computed on the unbounded domain by compactification. We map infinity to a finite boundary at $\rho = S$ with $S \in \mathbb{R}_+$ by choosing a compactification $r = r(\rho) : (0, S) \to (0, \infty)$ such that

$$r(0) = 0, \quad \lim_{\rho \to S} r(\rho) = \infty, \quad \lim_{\rho \to S} G(\rho) r^2(\rho) = \eta \tag{5}$$

for some fixed $\eta > 0$, where $G(\rho) := 1/r'(\rho)$.
**Key point:** In the compactified variables, the far-field is obtained by direct evaluation at the finite outer boundary

$$u^+_\infty(\omega) = \lim_{\rho \to S} u^+(\rho,\omega).$$

We denote $H(\rho) := h'(r(\rho))$ and assume the following regularity condition,

$$n^2(\rho,\omega) - H^2(\rho) \sim G(\rho) \quad \text{as } \rho \to S. \tag{6}$$

The Helmholtz equation (1) can be expressed in terms of the rescaled field $u^{\pm}(\rho,\omega)$ as

$$0 = \partial_\rho\big(G\,\partial_\rho u^{\pm}\big) \pm 2ik\,H\,\partial_\rho u^{\pm} + \frac{1}{G\,r^2}\Delta_{\mathbb{S}^{d-1}}u^{\pm} + \left[k^2\frac{n^2-H^2}{G} \pm ikH' + \frac{(d-1)(3-d)}{4\,G\,r^2}\right]u^{\pm}. \tag{7}$$

Using this formulation, we solve scattering problems for obstacles in unbounded media, including constant, short-range, and long-range cases.

## 3 Methods and Results

We demonstrate the numerical solution of (7) using spectral and finite-element methods in two dimensions, $d = 2$. The spectral method is implemented in the Dedalus framework [4] which represents the solution in a Fourier–Chebyshev basis,

$$u^{\pm}(\rho,\theta) \approx \sum_{n=0}^{N}\sum_{m=-M}^{M}\hat{u}^{\pm}_{n,m}T_n(\rho)e^{im\theta},$$

where $T_n(\rho)$ are Chebyshev polynomials and $e^{im\theta}$ are Fourier modes. We show spectrally accurate solutions for scattering at the unit disk. The transformed scattered solution oscillates primarily in angular directions resulting in a highly accurate solver on the full unbounded domain (see Fig. 1). The total solution can be computed by inverting the explicitly known transformations as shown in Fig. 2.

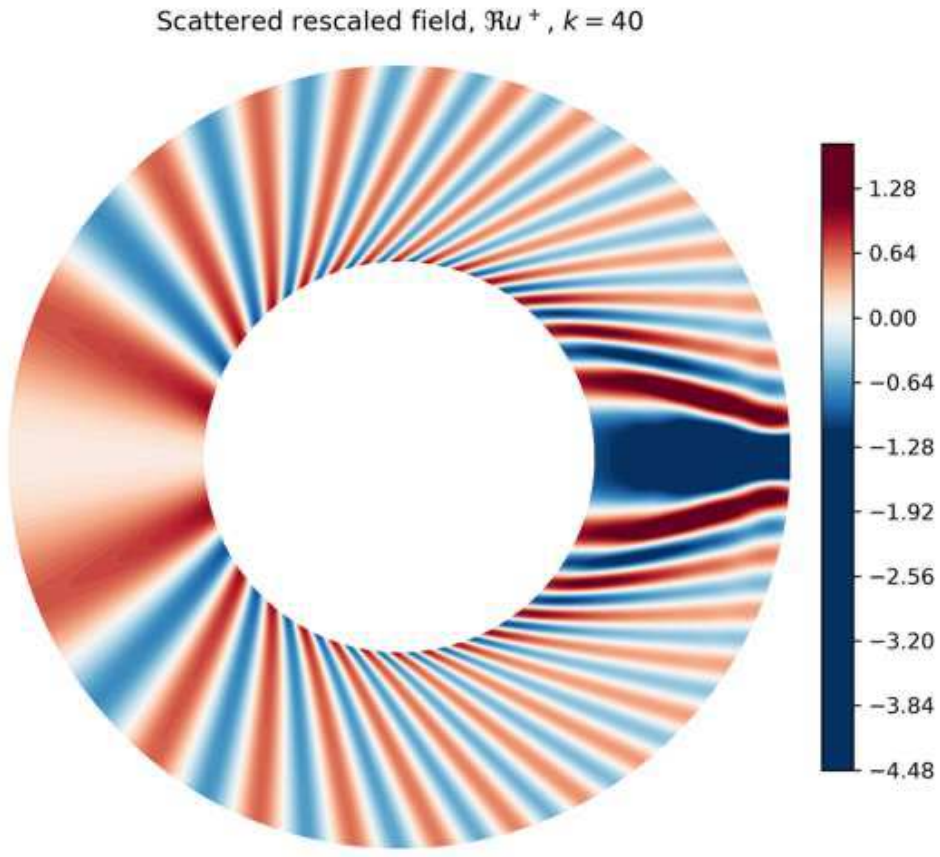

Figure 1: Real part of the rescaled field $u^+(\rho,\theta)$ scattered at a unit disk on a compactified mesh. The outer boundary represents infinity.

We perform efficient and highly accurate computations of scattering amplitudes in short-range, and long-range variable media without relying on Green's functions or explicit far-field representations, and with spectral convergence in representative experiments.

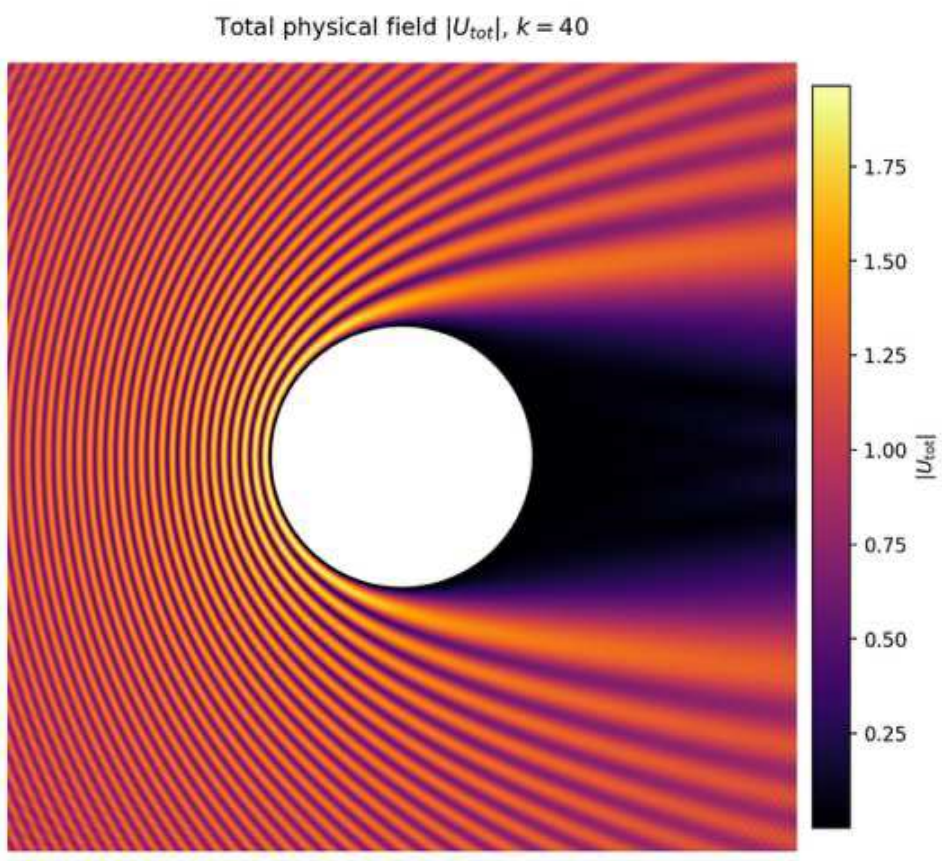

Figure 2: Magnitude of the physical total field $U = U^- + U^+$ represented on a finite box recovered from the compactified solution by inverting the rescaling and compactification.

We also develop and analyze finite element discretizations of the compactified problem (7). The finite-element implementation is based on the NGSolve package [5]. We benchmark the compactification approach against PML for various obstacle shapes including trapping geometries. The compactified formulation delivers improved accuracy at comparable computational cost, and provides direct access to radiative fields at infinity.